\documentclass[oneside,11pt,dvipsnames]{book}
\usepackage{jmlr2e}

\usepackage[noframe]{showframe}
\usepackage{framed}
\usepackage{lipsum}

\makeatletter
\newcommand{\paragrapharrow}{%
	\@startsection{paragraph}{4}{\z@}%
	{3.25ex \@plus1ex \@minus.2ex}%
	{-1em}%
	{\normalfont\normalsize\bfseries$\blacktriangleright$\ }}
\makeatother

\newcommand{\roundcornertheorem}{0pt}
\newcommand{\linewidththeorem}{0.1pt}
\newcommand{\frametitlerulewidththeorem}{0.1pt}
\newcommand{\innerbottommargintheorem}{2pt}
\newcommand{\innerleftmargintheorem}{2pt}
\newcommand{\innerrightmargintheorem}{2pt}
\newcommand{\innertopmargintheorem}{2pt}
\newcommand{\outerlinewidththeorem}{1pt}

\renewcommand{\thesection}{\@arabic\c@section}
\newcounter{theo}[chapter] 
\renewcommand{\thetheo}{\thechapter.\arabic{theo}}

\newenvironment{theoremHigh}[1][]{%
\refstepcounter{theo}%
\ifstrempty{#1}%
{\newcommand{\theoName}{}}
{\newcommand{\theoName}{:~(#1)}}
\mdfsetup{
	backgroundcolor=structurecolorHighTheorem,
	linecolor=structurecolorHighTheorem,
	frametitlerulewidth=\frametitlerulewidththeorem,
	roundcorner=\roundcornertheorem,
	linewidth=\linewidththeorem,
	innerbottommargin=\innerbottommargintheorem,
	innerleftmargin=\innerleftmargintheorem,
	innerrightmargin=\innerrightmargintheorem,
	innertopmargin=\innertopmargintheorem,
	outerlinewidth=\outerlinewidththeorem,
	topline=false,
	innertopmargin=-5pt,
	innerbottommargin=1pt,
	linewidth=0,
	startinnercode=\paragraph{{\strut Theorem~\thetheo\theoName}}
}
\begin{mdframed}[]\relax%
}{\end{mdframed}}

\newenvironment{lemma}[1][]{%
	\refstepcounter{theo}%
\ifstrempty{#1}%
{\newcommand{\theoName}{}}
{\newcommand{\theoName}{:~(#1)}}
\mdfsetup{
	backgroundcolor=\mdframecolorTheorem,
	linecolor=\mdframecolorTheorem,
	frametitlerulewidth=\frametitlerulewidththeorem,
	roundcorner=\roundcornertheorem,
	linewidth=\linewidththeorem,
	innerbottommargin=\innerbottommargintheorem,
	innerleftmargin=\innerleftmargintheorem,
	innerrightmargin=\innerrightmargintheorem,
	innertopmargin=\innertopmargintheorem,
	outerlinewidth=\outerlinewidththeorem,
	topline=false,
	innertopmargin=-5pt,
	innerbottommargin=1pt,
	linewidth=0,
	startinnercode=\paragraph{{\strut Lemma~\thetheo\theoName}}
}
\begin{mdframed}[]\relax%
}{\end{mdframed}}

\usepackage{array,multirow,graphicx}
\usepackage{bm}
\usepackage{cancel}

\usepackage{extarrows}

\usepackage[english]{babel}

\usepackage{imakeidx}
\makeindex[columns=2, title=Alphabetical Index]

\usepackage{stackengine}

\usepackage{enumitem}
\setenumerate[1]{itemsep=0pt,partopsep=0pt,parsep=\parskip,topsep=5pt}
\setitemize[1]{itemsep=0pt,partopsep=0pt,parsep=\parskip,topsep=5pt}
\setdescription{itemsep=0pt,partopsep=0pt,parsep=\parskip,topsep=5pt}

\usepackage{dsfont}

\usepackage[framemethod=tikz]{mdframed}

\usepackage{arydshln}

\numberwithin{equation}{chapter}
\usepackage[absolute,overlay]{textpos}

\usepackage[final]{pdfpages}

\usepackage[noframe]{showframe}
\usepackage{framed}
\usepackage{lipsum}

{\endMakeFramed}
\definecolor{shadecolor}{gray}{0.75}

\usepackage{xcolor}

\usepackage{tcolorbox}

\usepackage{graphicx}    
\usepackage[hang]{subfigure}

\usepackage{algorithm}
\usepackage{algpseudocode}
\usepackage{amsmath}
\usepackage{graphics}
\usepackage{epsfig}

\usepackage{blkarray}

\usepackage{listings}

\usepackage[Bjornstrup]{fncychap}

\newcommand{\colortitlechap}{\color[RGB]{60,113,183}} 
\newcommand{\colornumberchap}{\color[RGB]{60,113,183}} 
\newcommand{\colorbackchap}{\colorbox[RGB]{200,200,200}} 

\makeatletter
\renewcommand{\DOCH}{%
	\settowidth{\py}{\CNoV\thechapter}
	\addtolength{\py}{-10pt}%
	\fboxsep=0pt%
	\colorbackchap{\rule{0pt}{40pt}\parbox[b]{\textwidth}{\hfill}}%
	\kern-\py\raise20pt%
	\hbox{\colornumberchap\CNoV\thechapter}\\%
}

\renewcommand{\DOTI}[1]{%
	\nointerlineskip\raggedright%
	\fboxsep=\myhi%
	\vskip-1ex%
	\colorbackchap{\parbox[t]{\mylen}{\CTV\FmTi{\colortitlechap#1}}}\par\nobreak%
	\vskip 40\p@%
}

\renewcommand{\DOTIS}[1]{%
	\fboxsep=0pt%
	\colorbackchap{\rule{0pt}{40pt}\parbox[b]{\textwidth}{\hfill}}\\%
	\nointerlineskip\raggedright%
	\fboxsep=\myhi%
	\colorbackchap{\parbox[t]{\mylen}{\CTV\FmTi{\colortitlechap#1}}}\par\nobreak%
	\vskip 40\p@%
}
\makeatother

\usepackage[colorlinks]{hyperref}

\usepackage{color}   

\usepackage{hyperref}
\hypersetup{
colorlinks=true, 
linktoc=all,     
linkcolor=mydarkblue,  
anchorcolor=blue,
citecolor=mydarkgreen,
}

\usepackage{booktabs}  

\usepackage{tikz}
\usetikzlibrary{mindmap,trees,backgrounds}
\usetikzlibrary{arrows.meta}

\usetikzlibrary{decorations.text}

\usetikzlibrary{calc,backgrounds}

\newcommand*\circled[1]{\tikz[baseline=(char.base)]{
		\node[shape=circle,draw,inner sep=2pt] (char) {#1};}}

\usepackage{sidecap}
\sidecaptionvpos{figure}{t}

\usepackage{verbatimbox}

\usepackage[mathscr]{euscript}

\usepackage{stmaryrd}

\usepackage{minitoc}   
\let\cleardoublepage\clearpage
\usepackage[titletoc]{appendix}

\usepackage[labelfont=bf]{caption}
\newcommand\myhrulefill[1]{\leavevmode\leaders\hrule height#1\hfill\kern0pt}
\DeclareCaptionFormat{myformat}{
	{\color[RGB]{0,0,0}\myhrulefill{0.12em}}
	\\#1#2#3
}
\renewcommand\thesection{\thechapter.\arabic{section}}

\makeatletter
\renewcommand{\l@section}{\@dottedtocline{1}{1.5em}{2.6em}}
\renewcommand{\l@subsection}{\@dottedtocline{2}{4.0em}{3.6em}}
\renewcommand{\l@subsubsection}{\@dottedtocline{3}{7.4em}{4.5em}}
\makeatother

\usepackage{enumitem}
\newcommand*{\eitemi}{\tikz \draw [baseline, ball color=structurecolor,draw=none] circle (2pt);}
\newcommand*{\eitemii}{\tikz \draw [baseline, fill=structurecolor,draw=none,circular drop shadow] circle (2pt);}
\newcommand*{\eitemiii}{\tikz \draw [baseline, fill=structurecolor,draw=none] circle (2pt);}
\setlist[enumerate,1]{label=\color{black}\arabic*.,itemsep=0pt,partopsep=0pt,parsep=\parskip,topsep=5pt}
\setlist[enumerate,2]{label=\color{black}(\alph*).,itemsep=0pt,partopsep=0pt,parsep=\parskip,topsep=5pt}
\setlist[enumerate,3]{label=\color{black}\Roman*.,itemsep=0pt,partopsep=0pt,parsep=\parskip,topsep=5pt}
\setlist[enumerate,4]{label=\color{black}\Alph*.,itemsep=0pt,partopsep=0pt,parsep=\parskip,topsep=5pt}
\setlist[itemize,1]{label={\eitemi},itemsep=0pt,partopsep=0pt,parsep=\parskip,topsep=5pt}
\setlist[itemize,2]{label={\eitemii},itemsep=0pt,partopsep=0pt,parsep=\parskip,topsep=5pt}
\setlist[itemize,3]{label={\eitemiii},itemsep=0pt,partopsep=0pt,parsep=\parskip,topsep=5pt}

\usepackage{fancyhdr, blindtext}

\newcommand{\changefonts}{%
	\fontsize{9}{11}\selectfont
}

\fancypagestyle{plain}{%

	\fancyhf{}  
}

\usepackage{afterpage}

\newcommand{\natu}{\mathbb{N}}
\definecolor{blue-violet}{rgb}{0.54, 0.17, 0.89}
\definecolor{brightlavender}{rgb}{0.44, 0.16, 0.39}

\usepackage{amsthm}
\newtheoremstyle{normalfontstyle} 
{3pt}                           
{3pt}                           
{\normalfont}                   
{}                              
{\bfseries}                     
{}                             
{ }                             
{}                              

\usepackage{thmtools}
\declaretheoremstyle[
spaceabove=3pt,
spacebelow=3pt,
headfont=\bfseries,
notefont=\bfseries, 
notebraces={(}{)}, 
bodyfont=\normalfont,
postheadspace=1em, 
]{normalfontboldhead}

\theoremstyle{normalfontstyle}
\newcommand{\BlackBox}{\rule{1.5ex}{1.5ex}}  
\renewenvironment{proof}{\par\noindent{\bf Proof\ }}{\hfill\BlackBox\\[2mm]}

\declaretheorem[style=normalfontboldhead, name=Definition, numberlike=theo]{definitionT}
\newmdenv[skipabove=7pt,
skipbelow=7pt,
rightline=false,
leftline=true,
topline=false,
bottomline=false,
linecolor=mydarkblue,
innerleftmargin=5pt,
innerrightmargin=5pt,
innertopmargin=0pt,
leftmargin=2cm,
rightmargin=0cm,
linewidth=4pt,
innerbottommargin=0pt]{dBox}
\newenvironment{definition}{\begin{dBox}\begin{definitionT}}{\end{definitionT}\end{dBox}}

\declaretheorem[style=normalfontboldhead, name=Exercise, numberlike=theo]{exerciseC}
\newmdenv[skipabove=7pt,
skipbelow=7pt,
rightline=false,
leftline=true,
topline=false,
bottomline=false,
linecolor=mydarkgreen,
innerleftmargin=5pt,
innerrightmargin=5pt,
innertopmargin=0pt,
leftmargin=2cm,
rightmargin=0cm,
linewidth=4pt,
innerbottommargin=0pt]{eBox}
\newenvironment{exercise}{\begin{eBox}\begin{exerciseC}}{\end{exerciseC}\end{eBox}}

\declaretheorem[style=normalfontboldhead, name=Remark, numberlike=theo]{remarekC}
\newmdenv[skipabove=7pt,
skipbelow=7pt,
rightline=false,
leftline=true,
topline=false,
bottomline=false,
linecolor=mydarkpurple,
innerleftmargin=5pt,
innerrightmargin=5pt,
innertopmargin=0pt,
leftmargin=2cm,
rightmargin=0cm,
linewidth=4pt,
innerbottommargin=0pt]{rBox}
\newenvironment{remark}{\begin{rBox}\begin{remarekC}}{\end{remarekC}\end{rBox}}

\declaretheorem[style=normalfontboldhead, name=Assumption, numberlike=theo]{assumptionC}
\newmdenv[skipabove=7pt,
skipbelow=7pt,
rightline=false,
leftline=true,
topline=false,
bottomline=false,
linecolor=mydarkpurple,
innerleftmargin=5pt,
innerrightmargin=5pt,
innertopmargin=0pt,
leftmargin=2cm,
rightmargin=0cm,
linewidth=4pt,
innerbottommargin=0pt]{asBox}

\declaretheorem[style=normalfontboldhead, name=Condition, numberlike=theo]{conditionC}
\newmdenv[skipabove=7pt,
skipbelow=7pt,
rightline=false,
leftline=true,
topline=false,
bottomline=false,
linecolor=mydarkpurple,
innerleftmargin=5pt,
innerrightmargin=5pt,
innertopmargin=0pt,
leftmargin=2cm,
rightmargin=0cm,
linewidth=4pt,
innerbottommargin=0pt]{cdBox}

\declaretheorem[style=normalfontboldhead, name=Example, numberlike=theo]{exampleC}
\newmdenv[skipabove=7pt,
skipbelow=7pt,
rightline=false,
leftline=false,
topline=false,
bottomline=false,
linecolor=mydarkgreen,
innerleftmargin=1pt,
innerrightmargin=5pt,
innertopmargin=0pt,
leftmargin=2cm,
rightmargin=0cm,
linewidth=4pt,
innerbottommargin=0pt]{xBox}
\newenvironment{example}{\begin{xBox}\begin{exampleC}}{\exampbar\end{exampleC}\end{xBox}}

\usepackage{adforn}  
\newcommand{\xchaptertitle}{Chapter~\thechapter~}
\newcommand{\problemname}{Problems}
\newenvironment{problemset}[1][\xchaptertitle~\problemname]{
	\vspace*{10pt}
	\begin{center}
		\phantomsection\addcontentsline{toc}{section}{\texorpdfstring{\xchaptertitle~\problemname}{\problemname}}
		\markright{#1}
		\textcolor{structurecolor}{\Large\bfseries\adftripleflourishleft~#1~\adftripleflourishright}
	\end{center}
	\begin{enumerate}[ref=\thechapter.\theenumi]}{
\end{enumerate}}

\newcommand{\topone}{{(1)}}

\newcommand{\toptzero}{{(t)}}
\newcommand{\toptzeroTOP}{{(t)\top}}

\newcommand{\toptone}{{(t+1)}}
\newcommand{\toptoneTOP}{{(t+1)\top}}

\newenvironment{bmatrixfoot}
{\footnotesize\begin{bmatrix}}
	{\end{bmatrix}\normalsize}

\newcommand{\holders}{\text{H{\"o}lder's} }

\let\oldforall\forall
\renewcommand{\forall}{\oldforall\, }

\usepackage{color}   

\newcommand{\mdframecolor}{gray!10}

\newcommand{\mdframehidelineNote}{true}

\definecolor{mylightbluetitle}{RGB}{60,113,183}
\definecolor{mylightbluetext}{RGB}{27,54,189}
\definecolor{structurecolorblue}{RGB}{60,113,183}
\definecolor{structurecolorgreen}{RGB}{63,145,182}
\colorlet{structurecolor}{structurecolorblue}
\definecolor{structurecolorelegant}{RGB}{60,113,183}
\definecolor{structurecolorlt}{RGB}{31,119,185}

\definecolor{structurecolorHighTheoremBlue}{RGB}{220,227,248}
\definecolor{structurecolorHighTheoremGreen}{RGB}{188,222,231}
\colorlet{structurecolorHighTheorem}{structurecolorHighTheoremBlue}
\definecolor{colorBlue1}  {RGB}{220,227,248}
\newcommand{\mdframecolorTheorem}{gray!35}  

\definecolor{winestain}{rgb}{0.5,0,0}
\definecolor{mydarkblue}{rgb}{0,0.08,0.45}
\definecolor{mydarkred}{rgb}{0.70,0.00,0.00}
\definecolor{mydarkgreen}{rgb}{0.00,0.30,0.00}
\definecolor{mydarkyellow}{RGB}{197,151,13}
\definecolor{mydarkpurple}{RGB}{90,35,140}

\definecolor{color0}  {RGB}{174,225,254} 
\definecolor{color1}  {RGB}{220,227,248} 
\definecolor{color2}  {RGB}{28,130,185} 
\definecolor{color3}  {RGB}{255,253,250} 
\definecolor{colormiddleright}  {RGB}{245,253,250} 
\definecolor{colorbottomleft}  {RGB}{255,243,250} 
\definecolor{coloruppermiddle}  {RGB}{255,253,230} 
\definecolor{colormiddleleft}  {RGB}{255,244,237}
\definecolor{colorcr}  {RGB}{249,253,232} 
\definecolor{colorreduction}  {RGB}{255,235,254} 
\definecolor{colorqr}  {RGB}{254,221,199} 
\definecolor{colorbiconjugate}  {RGB}{251,149,161} 
\definecolor{colorsvd}  {RGB}{215,247,235} 
\definecolor{colorupperright}  {RGB}{239,246,251} 
\definecolor{colorspectral}  {RGB}{206,226,243} 
\definecolor{colorbottomright}  {RGB}{220,224,236} 
\definecolor{coloreigenvalue}  {RGB}{197,203,224} 
\definecolor{colorcp} {RGB}{217, 234, 186} 
\definecolor{colorcpborder} {RGB}{233, 243, 216} 
\definecolor{colorupperleft}  {RGB}{235,243,240} 
\definecolor{colorsemidefinite}  {RGB}{217,232,226} 
\definecolor{colormiddle} {RGB}{235, 240,255}
\definecolor{colorlu}  {RGB}{220,227,255} 
\definecolor{colorals}  {RGB}{240,230,255} 
\definecolor{coloralsbkg}  {RGB}{248,243,255} 
\definecolor{canaryyellow}{rgb}{1.0, 0.75, 0.0}
\definecolor{bluepigment}{rgb}{0.0, 0.0, 1.0}
\definecolor{canarypurple}{RGB}{208, 13, 241}
\definecolor{ocre}{RGB}{51,102,0} 
\definecolor{colorGreenOcre}{RGB}{51,102,0} 
\definecolor{colorBlue2}{RGB}{200,207,248}
\definecolor{shadecolor}{gray}{0.75}
\definecolor{darkblue}{rgb}{0.0, 0.0, 0.55}
\definecolor{citrine}{rgb}{0.89, 0.82, 0.04}

\newcommand{\hadaprod}{\circ}

\newcommand{\leadto}{\qquad\underrightarrow{ \text{leads to} }\qquad}

\newcommand{\gapthree}{\,\,\,}

\newcommand{\gap}{\,\,\,\,\,\,\,\,}  
\mathchardef\mhyphen="2D

\newcommand{\real}{\mathbb{R}}
\newcommand{\naturalset}{\mathbb{N}}
\newcommand{\complex}{\mathbb{C}}
\newcommand{\integer}{\mathbb{Z}}
\newcommand{\exampbar}{\hfill $\square$\par}
\newcommand{\argmax}{\text{arg max}}
\newcommand{\argmin}{\text{arg min}}

\newcommand\mathopmax[1]{\mathop{\max}_{#1}}
\newcommand\mathopmin[1]{\mathop{\min}_{#1}}
\newcommand\mathoplim[1]{\mathop{\lim}_{#1}}

\newcommand\abs[1]{\left\lvert#1\right\rvert}
\newcommand\absbig[1]{\big\lvert#1\big\rvert}

\newcommand\norm[1]{\left\lVert#1\right\rVert}

\newcommand\normzero[1]{\left\lVert#1\right\rVert_0}

\newcommand\normone[1]{\left\lVert#1\right\rVert_1}

\newcommand\normmone[1]{\left\lVert#1\right\rVert_{m_1}}
\newcommand\normtwo[1]{\left\lVert#1\right\rVert_2}
\newcommand\normtwobig[1]{\big\lVert#1\big\rVert_2}
\newcommand\normp[1]{\left\lVert#1\right\rVert_p}

\newcommand\norma[1]{\left\lVert#1\right\rVert_a}
\newcommand\normb[1]{\left\lVert#1\right\rVert_b}

\newcommand\normf[1]{\left\lVert#1\right\rVert_F}

\newcommand\norminf[1]{\left\lVert#1\right\rVert_{\infty}}

\newcommand\normminf[1]{\left\lVert#1\right\rVert_{m_\infty}}

\newcommand{\convd}{\breve d}
\newcommand{\concd}{\widehat d}
\newcommand{\cnstd}{\overline d}

\newcommand{\concG}{\widehat G}

\newcommand{\Cov}{\mathbb{C}\mathrm{ov}}

\newcommand{\Exp}{\mathbb{E}}
\newcommand{\Var}{\mathbb{V}\mathrm{ar}}
\newcommand{\indicator}{\mathds{1}}
\newcommand{\prob}{\Pr}

\newcommand{\bpi}{\bm{\pi}}
\newcommand{\exponential}{\mathcal{E}}
\newcommand{\cauchydist}{\mathcal{C}}
\newcommand{\gammadist}{\mathcal{G}}
\newcommand{\inversegammadist}{\mathcal{G}^{-1}}
\newcommand{\inversenormaldist}{\mathcal{N}^{-1}}
\newcommand{\rectifieddist}{\mathcal{RN}}
\newcommand{\studentt}{\tau}
\newcommand{\standardstudentt}{\mathrm{t}}
\newcommand{\normal}{\mathcal{N}}
\newcommand{\laplacedist}{\mathcal{L}}
\newcommand{\skewlaplacedist}{\mathcal{SL}}

\newcommand{\truncatednormal}{\mathcal{TN}}
\newcommand{\generaltruncatednormal}{\mathcal{GTN}}
\newcommand{\gtnsng}{\mathcal{GTNSNG}}
\newcommand{\halfnormal}{\mathcal{HN}}
\newcommand{\betadist}{\mathrm{Beta}}
\newcommand{\wishartdist}{\mathrm{Wi}}
\newcommand{\inversewishart}{\mathrm{IW}}
\newcommand{\inversechidist}{\mathrm{\chi^{-2}}}
\newcommand{\normalinversegamma}{\mathcal{NIG}}
\newcommand{\normalgamma}{\mathcal{NG}}
\newcommand{\tnsng}{\mathcal{TNSNG}}
\newcommand{\rnsng}{\mathcal{RNSNG}}
\newcommand{\niw}{\mathcal{NIW}}
\newcommand{\nix}{\mathcal{NIX}}
\newcommand{\nig}{\mathcal{NIG}}
\newcommand{\chisquared}{\chi^2}
\newcommand{\nonchisquared}{\mathcal{NX}^2}
\newcommand{\poissondist}{\mathcal{P}}

\newcommand{\bxn}{\bx_n}
\newcommand{\bmo}{\bm{m}_0}
\newcommand{\bso}{\bm{S}_0}
\newcommand{\bxstar}{\bx^{\star}}
\newcommand{\xstar}{x^{\star}}
\newcommand{\Nstar}{{N^{\star}}}

\newcommand{\smu}{\mu}
\newcommand{\ssigma}{\sigma}

\newcommand{\dirichlet}{\mathrm{Dirichlet}}
\newcommand{\multinomial}{\mathrm{Multi}}
\newcommand{\binomialdist}{\mathrm{Binom}}
\newcommand{\multinoulli}{\mathrm{Multin}}  
\newcommand{\bernoulli}{\mathrm{Bern}}  
\newcommand{\bernoullidist}{\mathrm{Bern}}

\newcommand{\tr}{\mathrm{tr}}

\newcommand\comple[1]{#1^c}
\newcommand{\cspace}{\mathcal{C}}
\newcommand{\nspace}{\mathcal{N}}
\newcommand{\bone}{\mathbf{1}}
\newcommand{\diag}{\mathrm{diag}}
\newcommand{\rank}{\mathrm{rank}}
\newcommand{\adjugate}{\mathrm{adj}}
\newcommand{\trace}{\mathrm{tr}}

\newcommand{\ba}{\bm{a}}
\newcommand{\bA}{\bm{A}}
\newcommand{\bb}{\bm{b}}
\newcommand{\bB}{\bm{B}}
\newcommand{\bc}{\bm{c}}
\newcommand{\bC}{\bm{C}}
\newcommand{\bd}{\bm{d}}
\newcommand{\bD}{\bm{D}}
\newcommand{\be}{\bm{e}}
\newcommand{\bE}{\bm{E}}
\newcommand{\bff}{\bm{f}}
\newcommand{\bF}{\bm{F}}
\newcommand{\bg}{\bm{g}}
\newcommand{\bG}{\bm{G}}

\newcommand{\bH}{\bm{H}}

\newcommand{\bI}{\bm{I}}

\newcommand{\bK}{\bm{K}}
\newcommand{\bl}{\bm{l}}
\newcommand{\bL}{\bm{L}}
\newcommand{\bmm}{\bm{m}}
\newcommand{\bM}{\bm{M}}
\newcommand{\bn}{\bm{n}}
\newcommand{\bN}{\bm{N}}
\newcommand{\bo}{\bm{o}}

\newcommand{\bp}{\bm{p}}
\newcommand{\bP}{\bm{P}}
\newcommand{\bq}{\bm{q}}
\newcommand{\bQ}{\bm{Q}}
\newcommand{\br}{\bm{r}}
\newcommand{\bR}{\bm{R}}
\newcommand{\bs}{\bm{s}}
\newcommand{\bS}{\bm{S}}

\newcommand{\bT}{\bm{T}}
\newcommand{\bu}{\bm{u}}
\newcommand{\bU}{\bm{U}}
\newcommand{\bv}{\bm{v}}
\newcommand{\bV}{\bm{V}}
\newcommand{\bw}{\bm{w}}
\newcommand{\bW}{\bm{W}}
\newcommand{\bx}{\bm{x}}
\newcommand{\bX}{\bm{X}}
\newcommand{\by}{\bm{y}}
\newcommand{\bY}{\bm{Y}}
\newcommand{\bz}{\bm{z}}
\newcommand{\bZ}{\bm{Z}}

\def\vmu{{\bm{\mu}}}
\def\vtheta{{\bm{\theta}}}
\def\va{{\bm{a}}}

\def\ve{{\bm{e}}}

\def\vx{{\bm{x}}}

\def\mA{{\bm{A}}}
\def\mB{{\bm{B}}}

\def\mH{{\bm{H}}}

\def\mJ{{\bm{J}}}

\def\mX{{\bm{X}}}

\def\mSigma{{\bm{\Sigma}}}

\DeclareMathAlphabet{\mathsfit}{\encodingdefault}{\sfdefault}{m}{sl}
\SetMathAlphabet{\mathsfit}{bold}{\encodingdefault}{\sfdefault}{bx}{n}

\def\sA{{\mathbb{A}}}
\def\sB{{\mathbb{B}}}

\def\sG{{\mathbb{G}}}

\def\sI{{\mathbb{I}}}
\def\sJ{{\mathbb{J}}}
\def\sK{{\mathbb{K}}}
\def\sL{{\mathbb{L}}}

\def\sS{{\mathbb{S}}}
\def\sT{{\mathbb{T}}}

\def\sX{{\mathbb{X}}}

\def\sZ{{\mathbb{Z}}}

\def\1{\bm{1}}

\def\rS{{\textnormal{S}}}

\def\ra{{\textnormal{a}}}
\def\rb{{\textnormal{b}}}
\def\rc{{\textnormal{c}}}

\def\rt{{\textnormal{t}}}

\def\rw{{\textnormal{w}}}
\def\rx{{\textnormal{x}}}
\def\ry{{\textnormal{y}}}
\def\rz{{\textnormal{z}}}

\def\rva{{\mathbf{a}}}
\def\rvb{{\mathbf{b}}}

\def\rvv{{\mathbf{v}}}

\def\rvx{{\mathbf{x}}}
\def\rvy{{\mathbf{y}}}
\def\rvz{{\mathbf{z}}}

\def\rmA{{\mathbf{A}}}
\def\rmB{{\mathbf{B}}}

\def\rmY{{\mathbf{Y}}}

\def\eva{{a}}

\def\gD{{\mathcal{D}}}
\def\gE{{\mathcal{E}}}

\newcommand{\R}{\mathbb{R}}

\newcommand{\entropy}{\mathrm{H}}
\newcommand{\KL}{D_{\mathrm{KL}}}

\newcommand{\bzero}{\boldsymbol{0}}

\newcommand{\balpha}{{\boldsymbol\alpha}}
\newcommand{\bbeta}{{\boldsymbol\beta}}

\newcommand{\bepsilon}{{\boldsymbol\epsilon}}

\newcommand{\bgamma}{{\boldsymbol\gamma}}

\newcommand{\blambda}{{\boldsymbol\lambda}}
\newcommand{\bmu}{{\boldsymbol\mu}}
\newcommand{\bnu}{{\boldsymbol\nu}}
\newcommand{\bomega}{{\boldsymbol\omega}}

\newcommand{\bsigma}{{\boldsymbol\sigma}}
\newcommand{\btau}{{\boldsymbol\tau}}
\newcommand{\btheta}{{\boldsymbol\theta}}

\newcommand{\bLambda}{{\boldsymbol\Lambda}}
\newcommand{\bOmega}{{\boldsymbol\Omega}}

\newcommand{\bSigma}{{\boldsymbol\Sigma}}
\newcommand{\bTheta}{{\boldsymbol\Theta}}

\newcommand{\widehattheta}{{\widehat\theta}}

\newcommand{\widetildealpha}{{\widetilde\alpha}}
\newcommand{\widetildebeta}{{\widetilde\beta}}

\newcommand{\widehatbbeta}{{\widehat\bbeta}}

\newcommand{\widehatbtheta}{{\widehat\btheta}}

\newcommand{\widehatbSigma}{{\widehat\bSigma}}

\newcommand{\widebarblambda}{{\overline\blambda}}

\newcommand{\widebarx}{\overline{{x}}}

\newcommand{\widebarbW}{\overline{\bm{W}}}

\newcommand{\widebarbZ}{\overline{\bm{Z}}}
\newcommand{\widebarba}{\overline{\bm{a}}}
\newcommand{\widebarbb}{\overline{\bm{b}}}

\newcommand{\widebarbw}{\overline{\bm{w}}}
\newcommand{\widebarbx}{\overline{\bm{x}}}

\newcommand{\widebarbz}{\overline{\bm{z}}}

\newcommand{\widehata}{\widehat{a}}

\newcommand{\widehatp}{\widehat{p}}

\newcommand{\widehatx}{\widehat{x}}

\newcommand{\widehatbx}{\widehat{\bm{x}}}

\newcommand{\widehatbz}{\widehat{\bm{z}}}

\newcommand{\mathcalC}{\mathcal{C}}

\newcommand{\mathcalF}{\mathcal{F}}

\newcommand{\mathcalH}{\mathcal{H}}

\newcommand{\mathcalL}{\mathcal{L}}

\newcommand{\mathcalO}{\mathcal{O}}
\newcommand{\mathcalP}{\mathcal{P}}
\newcommand{\mathcalQ}{\mathcal{Q}}

\newcommand{\mathcalV}{\mathcal{V}}

\newcommand{\mathcalX}{\mathcal{X}}
\newcommand{\mathcalY}{\mathcal{Y}}
\newcommand{\mathcalZ}{\mathcal{Z}}

\newcommand{\widetildebA}{\widetilde{\bm{A}}}
\newcommand{\widetildebB}{\widetilde{\bm{B}}}
\newcommand{\widetildebC}{\widetilde{\bm{C}}}
\newcommand{\widetildebD}{\widetilde{\bm{D}}}

\newcommand{\widetildebS}{\widetilde{\bm{S}}}

\newcommand{\widetildebW}{\widetilde{\bm{W}}}
\newcommand{\widetildebX}{\widetilde{\bm{X}}}

\newcommand{\widetildebZ}{\widetilde{\bm{Z}}}
\newcommand{\widetildeba}{\widetilde{\bm{a}}}

\newcommand{\widetildebl}{\widetilde{\bm{l}}}

\newcommand{\widetildebw}{\widetilde{\bm{w}}}
\newcommand{\widetildebx}{\widetilde{\bm{x}}}
\newcommand{\widetildeby}{\widetilde{\bm{y}}}
\newcommand{\widetildebz}{\widetilde{\bm{z}}}

\newcommand{\widetildea}{\widetilde{a}}

\newcommand{\widetildew}{\widetilde{w}}

\newcommand{\widetildey}{\widetilde{y}}
\newcommand{\widetildez}{\widetilde{z}}

\definecolor{titlepagecolor}{cmyk}{75,68,67,90}
\definecolor{titlepagecolor2}{rgb}{1.0, 0.08, 0.58}
\definecolor{emerald}{rgb}{0.31, 0.78, 0.47}
\definecolor{deeppink}{HTML}{D14064}
\definecolor{lowpink}{HTML}{ffe6ec}
\newcommand{\partcolor}{gray!65} 
\definecolor{lowblue}{HTML}{E1EBFE}

\makeatletter \renewcommand*\cleardoublepage{
	\clearpage
	\if@twoside   
	\ifodd\c@page 
	\hbox{}\newpage
	\if@twocolumn\hbox{}   
	\newpage
	\fi
	\fi
	\fi
} \makeatother
\let\originalpart=\part
\def\part#1{\cleardoublepage\clearpage \pagecolor{\partcolor} \originalpart{#1}\nopagecolor }

\usepackage{yfonts,color,lettrine}
\definecolor{caligraphcolor}{HTML}{74AECB}

\usepackage{setspace}

\AtBeginDocument{\setlength{\DefaultFindent}{0.5em}}
\def\algoalign#1{\parbox[t]{\dimexpr\linewidth-\algorithmicindent}{#1}}

\jmlrheading{1}{2023--}{1-48}{4/00}{10/00}{Jun Lu}
\firstpageno{1}
\begin{document}

\frontmatter
\pagecolor{lowblue}
\newpage
\thispagestyle{empty}
\includegraphics{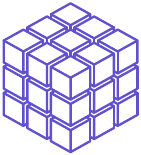}
%
%
%
\clearpage

\clearpage
\nopagecolor 

\thispagestyle{empty}  
\title{Bayesian Matrix Decomposition and Applications}

\author{
\begin{center}
\name Jun Lu \\ 
\email jun.lu.locky@gmail.com \\
\copyright \,\, 2023\~ \,\, Jun Lu\\
\end{center}
}

\maketitle

\chapter*{\centering \begin{normalsize}Preface\end{normalsize}}

In 1954, Alston S. Householder published \textit{Principles of Numerical Analysis}, an early and influential work on matrix decomposition, specifically favoring (block) LU decomposition---a method that factorizes a matrix into the product of lower and upper triangular matrices.
Today, matrix decomposition is a foundational tool in machine learning, statistics and many other fields. Its prominence has been greatly amplified by the development of the backpropagation algorithm for training neural networks. Matrix decomposition plays a crucial role in reducing data dimensionality and transforming data into representations that are better suited for machine learning algorithms.

Bayesian matrix decomposition is a relatively recent subfield within the broader landscape of matrix decomposition and machine learning. It combines the principles of Bayesian statistics with classical matrix factorization techniques to perform probabilistic matrix decomposition. The use of Bayesian methods in this context was first introduced in the early 2000s, primarily to address limitations of traditional matrix factorization approaches---such as poor interpretability, limited predictive performance, and difficulties in handling uncertainty.

The primary goal of this book is to provide a self-contained introduction to the core concepts and mathematical tools of Bayesian matrix decomposition, thereby laying a smooth foundation for discussing specific decomposition techniques and their applications in later sections.
That said, due to space constraints, we acknowledge that we cannot cover all the rich and valuable developments in this area---for example, a detailed treatment of variational inference for optimization is beyond our scope. Readers seeking deeper coverage are encouraged to consult specialized literature on Bayesian analysis and probabilistic modeling.

This book primarily serves as a concise overview of the purpose, significance, and practical utility of key Bayesian matrix decomposition methods, including real-valued decomposition, nonnegative matrix factorization (NMF), and Bayesian interpolative decomposition. It explores the origins and computational characteristics of these methods and highlights their real-world applications. The only prerequisites are a basic course in linear algebra and introductory statistics, ensuring accessibility while maintaining mathematical rigor through carefully presented proofs throughout the text.

\section*{\centering \begin{normalsize}Keywords\end{normalsize}}
Bayesian inference, Gibbs sampling, Conjugate model,
Alternating least squares (ALS), Multiplicative update, Real-valued matrix decomposition, Low-rank approximation,  Nonnegative matrix decomposition, Autoencoder, Principal component analysis, Interpolative decomposition, Ordinal matrix decomposition, Poisson matrix decomposition.

\vspace{5em}
\noindent
\section*{\centering \begin{normalsize}Acknowledgment\end{normalsize}}
We gratefully acknowledge Ulrich Paquet and Josh Chang for sharing valuable insights on Bayesian ordinal and Poisson matrix decompositions, and Gilbert Strang for his helpful discussion regarding the proof of the CUR decomposition.
The author is deeply grateful to Joerg Osterrieder, Christine P. Chai, and Xuanyu Ye for their collaboration on the Bayesian approach to nonnegative matrix factorization and interpolative decomposition. Their contributions have greatly enriched many of the discussions presented in this book.
Finally, the author thanks Nicolas P. Rougier for open-sourcing the poster design \citep{nicolas2015nn}.

\newpage
\begingroup
\hypersetup{
linkcolor=structurecolor,
linktoc=page,  
}
\dominitoc
\pdfbookmark{\contentsname}{toc} 
\tableofcontents 
\endgroup

\chapter*{Notation}\label{notation}
\index{Notation}


This section provides a concise reference describing notation used throughout this
book.
If you are unfamiliar with any of the corresponding mathematical concepts,
the book describes most of these ideas in Chapter~\ref{chapter_introduction}.

\vspace{0.4in}
\begin{minipage}{\textwidth}
\centerline{\bf Numbers and Arrays}
\bgroup
\def\arraystretch{1.5}
\begin{tabular}{cp{4.25in}}
$\displaystyle a$   & A scalar (integer or real)\\
$\displaystyle \ba$ & A vector\\
$\displaystyle \bA$ & A matrix\\
$\displaystyle \bI_n$ & Identity matrix with $n$ rows and $n$ columns\\
$\displaystyle \bI$   & Identity matrix with dimensionality implied by context\\
$\displaystyle \ve_i$ & Standard basis vector $[0,\dots,0,1,0,\dots,0]$ with a 1 at position $i$\\
$\displaystyle \text{diag}(\va)$ & A square, diagonal matrix with diagonal entries given by $\va$\\
$\displaystyle \ra$   & A scalar random variable\\
$\displaystyle \rva$  & A vector-valued random variable\\
$\displaystyle \rmA$  & A matrix-valued random variable\\
\end{tabular}
\egroup
\index{Scalar}
\index{Vector}
\index{Matrix}
\index{Tensor}
\end{minipage}

\index{Sets}
\vspace{0.2in}
\begin{minipage}{\textwidth}
\centerline{\bf Sets}
\bgroup
\def\arraystretch{1.5}
\begin{tabular}{cp{4.25in}}
$\displaystyle \sA$ & A set\\
$\displaystyle \varnothing$ & The null set \\
$\displaystyle \real$ & The set of real numbers \\
$\displaystyle \naturalset$ & The set of natural numbers \\
$\displaystyle \complex$ & The set of complex numbers \\
$\displaystyle \{0, 1\}$ & The set containing 0 and 1 \\
$\displaystyle \{0, 1, \dots, n \}$ & The set of all integers between $0$ and $n$\\
$\displaystyle [a, b]$ & The real interval including $a$ and $b$\\
$\displaystyle (a, b]$ & The real interval excluding $a$ but including $b$\\
$\displaystyle \sA \backslash \sB$ & Set subtraction, i.e., the set containing the elements of $\sA$ that are not in $\sB$\\
\end{tabular}
\egroup
\index{Scalar}
\index{Vector}
\index{Matrix}
\index{Tensor}
\index{Graph}
\index{Set}
\end{minipage}

\index{Matrix indexing}
\vspace{0.2in}
\begin{minipage}{\textwidth}
\centerline{\bf Indexing}
\bgroup
\def\arraystretch{1.5}
\begin{tabular}{cp{4.25in}}
$\displaystyle \eva_i$ & Element $i$ of vector $\va$, with indexing starting at 1 \\
$\displaystyle \ba_{-i}$ & All elements of vector $\va$ except for element $i$ \\
$\displaystyle  a_{ij}$ & Element $i, j$ of matrix $\mA$ \\
$\displaystyle \mA_{i, :}=\mA[i, :]$ & Row $i$ of matrix $\mA$ \\
$\displaystyle \mA_{:, i}=\mA[:, i],\, \ba_i$ & Column $i$ of matrix $\mA$ \\
\end{tabular}
\egroup
\end{minipage}

\vspace{0.2in}
\begin{minipage}{\textwidth}
\centerline{\bf Linear Algebra Operations}
\bgroup
\def\arraystretch{1.5}
\begin{tabular}{cp{4.25in}}
$\displaystyle \bA^\top$ & Transpose of matrix $\mA$ \\
$\displaystyle \bA^+$ & Moore-Penrose pseudoinverse of $\mA$\\
$\displaystyle \bA \hadaprod \bB $ & Element-wise (Hadamard) product of $\mA$ and $\mB$ \\
$\displaystyle \det(\bA)$ & Determinant of $\bA$ \\
$\displaystyle \mathrm{rref}(\bA)$ & Reduced row echelon form of $\bA$ \\
$\displaystyle \cspace(\bA)$ & Column space of $\bA$ \\
$\displaystyle \nspace(\bA)$ & Null space of $\bA$ \\
$\displaystyle \mathcalV$ & A general subspace \\
$\displaystyle \rank(\bA)$ & Rank of $\bA$ \\
$\displaystyle \trace(\bA)$ & Trace of $\bA$ \\
\end{tabular}
\egroup
\index{Transpose}
\index{Element-wise product|see {Hadamard product}}
\index{Hadamard product}
\index{Determinant}
\end{minipage}

\vspace{0.4in}
\begin{minipage}{\textwidth}
\centerline{\bf Calculus}
\bgroup
\def\arraystretch{1.5}
\begin{tabular}{cp{4.25in}}
$\displaystyle\frac{d y} {d x}$ & Derivative of $y$ with respect to $x$\\ [2ex]
$\displaystyle \frac{\partial y} {\partial x} $ & Partial derivative of $y$ with respect to $x$ \\
$\displaystyle \nabla_{\bx} y $ & Gradient of $y$ with respect to $\bx$ \\
$\displaystyle \nabla_{\bX} y $ & Matrix derivatives of $y$ with respect to $\bX$ \\
$\displaystyle \frac{\partial f}{\partial \vx} $ & Jacobian matrix $\mJ \in \R^{m\times n}$ of $f: \R^n \rightarrow \R^m$\\
$\displaystyle \nabla_\vx^2 f(\vx)\text{ or }\mH( f)(\vx)$ & The Hessian matrix of $f$ at input point $\vx$\\
$\displaystyle \int f(\vx) d\vx $ & Definite integral over the entire domain of $\vx$ \\
$\displaystyle \int_\sS f(\vx) d\vx$ & Definite integral with respect to $\vx$ over the set $\sS$ \\
\end{tabular}
\egroup
\index{Derivative}
\index{Integral}
\index{Jacobian matrix}
\index{Hessian matrix}
\end{minipage}

\vspace{0.4in}
\begin{minipage}{\textwidth}
\centerline{\bf Probability and Information Theory}
\bgroup
\def\arraystretch{1.5}
\begin{tabular}{cp{4.25in}}
$\displaystyle \ra \bot \rb$ & The random variables $\ra$ and $\rb$ are independent\\
$\displaystyle \ra \bot \rb \mid \rc $ & They are conditionally independent given $\rc$\\
$\displaystyle P(\ra)$ & A probability distribution over a discrete variable\\
$\displaystyle p(\ra)$ & A probability distribution over a continuous variable, or over
a variable whose type has not been specified\\
$\displaystyle \ra \sim P$ & Random variable $\ra$ has distribution $P$\\
$\displaystyle  \Exp_{\rx\sim P} [ f(x) ]\text{ or } \Exp [f(x)]$ & Expectation of $f(x)$ with respect to $P(\rx)$ \\
$\displaystyle \Var[f(x)] $ &  Variance of $f(x)$ under $P(\rx)$ \\
$\displaystyle \Cov[f(x),g(x)] $ & Covariance of $f(x)$ and $g(x)$ under $P(\rx)$\\
$\displaystyle \entropy(\rx) $ & Shannon entropy of the random variable $\rx$\\
$\displaystyle \KL [P \parallel Q ] $ & Kullback--Leibler divergence of P and Q \\
$\displaystyle \mathcal{N} ( \vx \mid \vmu , \mSigma)$ & Gaussian distribution %
over $\vx$ with mean $\vmu$ and covariance $\mSigma$ \\
\end{tabular}
\egroup
\index{Independence}
\index{Conditional independence}
\index{Variance}
\index{Covariance}
\index{Kullback--Leibler divergence}
\index{Shannon entropy}
\end{minipage}

\vspace{0.4in}
\begin{minipage}{\textwidth}
\centerline{\bf Functions}
\bgroup
\def\arraystretch{1.5}
\begin{tabular}{cp{4.25in}}
$\displaystyle f: \sA \rightarrow \sB$ & The function $f$ with domain $\sA$ and range $\sB$\\
$\displaystyle f \circ g $ & Composition of the functions $f$ and $g$ \\
  $\displaystyle f(\vx ; \vtheta) $ & A function of $\vx$ parametrized by $\vtheta$.
  (Sometimes we write $f(\vx)$ and omit the argument $\vtheta$ to lighten notation) \\
$\displaystyle \log(x), \ln(x)$ & Natural logarithm of $x$ \\
$\displaystyle \sigma(x),\, \text{Sigmoid}(x)$ & Logistic sigmoid, i.e., $\displaystyle \frac{1} {1 + \exp\{-x\}}$ \\
$\displaystyle \zeta(x)$ & Softplus, $\log(1 + \exp\{x\})$ \\
$\displaystyle \norm{\bx}_p $ & $\ell_p$-norm of $\vx$ \\
$\displaystyle \norm{\bx}=\norm{\bx}_2 $ & $\ell_2$-norm of $\vx$ \\
$\displaystyle \norm{\bx}=\norm{\bx}_1 $ & $\ell_1$-norm of $\vx$ \\
$\displaystyle \norm{\bx}=\norm{\bx}_\infty $ & $\ell_\infty$-norm of $\vx$ \\
$\displaystyle x^+$ & Positive part of $x$, i.e., $\max(0,x)$\\
$\displaystyle u(x)$ & Step function with value 1 when $x\geq0$ and value 0 otherwise\\
$\displaystyle \indicator\{\mathrm{condition}\}$ & is 1 if the condition is true, 0 otherwise\\
\end{tabular}
\egroup
\index{Sigmoid}
\index{Softplus}
\index{Norm}
\end{minipage}

Sometimes we use a function $f$ whose argument is a scalar but apply
it to a vector, matrix: $f(\vx)$, $f(\mX)$.
This denotes the application of $f$ to the
array element-wise. For example, if $\bC = \sigma(\bX)$, then $c_{ij} = \sigma(x_{ij})$
for all valid values of $i$ and  $j$.

\vspace{0.4in}
\begin{minipage}{\textwidth}
\centerline{\bf Abbreviations}
\bgroup
\def\arraystretch{1.5}
\begin{tabular}{cp{4.25in}}
EM & Expectation maximization \\
LVM & Latent variable model\\
ELBO & Evidence lower-bound \\
VFE & Variational free-energy\\
GMM & Gaussian mixture model\\
VI & Variational inference\\
PD & Positive definite  \\
PSD & Positive semidefinite \\
MCMC & Markov chain Monte Carlo \\
i.i.d. & Independently and identically distributed \\
p.d.f., PDF & Probability density function \\
p.m.f., PMF & Probability mass function \\
OLS & Ordinary least squares\\
NG & Normal-Gamma distribution \\
NIG & Normal-inverse-Gamma  distribution \\
NIX & Normal-inverse-Chi-squared distribution\\
TN & Truncated-normal distribution \\
GTN & General-truncated-normal distribution\\
RN & Rectified-normal distribution \\
IW & Inverse-Wishart distribution \\
NIW & Normal-inverse-Wishart distribution \\
ARD & automatic relevance determination\\
ALS & Alternating least squares \\
GD, SGD & Gradient descent, stochastic gradient descent\\
MU & Multiplicative update \\
MSE & Mean squared error\\
NMF & Nonnegative matrix factorization\\
ID & Interpolative decomposition\\
IID & Intervened interpolative decomposition \\
BID & Bayesian interpolative decomposition \\

\end{tabular}
\egroup
\end{minipage}


\clearpage


\mainmatter

\part{Backgrounds}

\newpage
\chapter{Introduction and Background}\label{chapter_introduction}
\begingroup
\hypersetup{
	linkcolor=structurecolor,
	linktoc=page,  
}
\minitoc \newpage
\endgroup
\section*{Introduction and Background}
\addcontentsline{toc}{section}{Introduction and Background}
Matrix decomposition has become a core technology across numerous fields, including statistics \citep{banerjee2014linear, gentle1998numerical}, optimization \citep{gill2021numerical}, clustering and classification \citep{li2009non, wang2013non, lu2021survey}, computer vision \citep{goel2020survey}, and recommender system \citep{symeonidis2016matrix}. 
Its importance stems largely from its integration into machine learning applications \citep{goodfellow2016deep, bishop2006pattern}.
Machine learning algorithms are designed to uncover hidden patterns and relationships in data, yet they often face challenges when dealing with high-dimensional and complex datasets. Matrix decomposition techniques address this by reducing data dimensionality and representing it in a form that is more amenable to processing by machine learning models.

Algorithms such as QR decomposition, singular value decomposition (SVD), alternating least squares (ALS), and nonnegative matrix factorization (NMF) decompose a matrix into a smaller set of constituent matrices that capture the underlying structure of the data. These factor matrices can then serve as features for machine learning algorithms, enabling them to learn patterns and relationships more effectively. Additionally, this process helps reduce noise and redundancy in the data, making it easier to identify meaningful structures.

Beyond feature extraction, matrix decomposition is widely applied to various matrix-based problems, such as collaborative filtering and link prediction \citep{marlin2003modeling, lim2007variational, mnih2007probabilistic, raiko2007principal, chen2009collaborative}:
\begin{itemize}
\item In collaborative filtering, matrix decomposition reveals latent patterns in user-item interaction matrices. For instance, given a matrix of user ratings for items, decomposition methods can infer latent user preferences and item attributes. These latent factors can then predict ratings for unseen items, enabling personalized recommendations (e.g., articles, movies, or music).
	
\item In link prediction problems, matrix decomposition algorithms can be used to uncover hidden patterns in networks. For example, given a network of users and their connections, matrix decomposition methods can be used to infer latent user preferences, which can then be used to predict new links in the network.
\end{itemize}

\index{Collaborative filtering}
\index{Link prediction}

A (bilinear) matrix decomposition expresses a complex matrix as the product of two (or more) simpler factor matrices. The fundamental idea behind this decompositional approach is not to solve specific problems directly, but rather to simplify complex matrix operations by performing them on the decomposed components instead of the original matrix. Although computing a decomposition can be computationally expensive, once obtained---say, a factorization of $\bA$---it can be reused efficiently across multiple tasks. For example, it enables solving an entire set of linear systems $\{\bb_1=\bA\bx_1, \bb_2=\bA\bx_2, \ldots, \bb_K=\bA\bx_K\}$
without recomputing the factorization each time.

There are two main approaches that have been applied to inference in the (low-rank) matrix factorization task. The first is to define a loss function and optimize over
the factored components using alternating updates \citep{comon2009tensor, lee1999learning}. 
The second is to build a probabilistic model representing the matrix factorization and then to perform
statistical inference to compute any desired components \citep{salakhutdinov2008bayesian, ari2012probabilistic}.

Another core application of  matrix decomposition  in machine learning is that it provides a way of incorporating prior knowledge into the machine learning algorithms. For example, in \textit{Bayesian matrix decomposition (BMD, or Bayesian matrix factorization, BMF)}, assumptions about sparsity, nonnegativity, or other structural properties can be encoded through prior distributions. This allows the model to make more informed inferences and often yields improved results.

Traditional methods like SVD and NMF have proven effective at capturing intrinsic data structures. However, they often struggle with missing or incomplete data, modeling uncertainty, and integrating prior information. Bayesian matrix factorization addresses these limitations by embedding Bayesian principles into the factorization framework, offering a more flexible and robust approach.

Bayesian matrix decomposition is a probabilistic model that factorizes a matrix into latent components. It was initially explored in the contexts of factor analysis \citep{canny2004gap, dunson2005bayesian} and matrix completion \citep{zhou2010nonparametric}. As a generative graphical model, BMD infers low-dimensional latent factors---such as user preferences and item attributes---from observed data. It has been successfully applied to tasks including matrix completion, image inpainting, denoising, and super-resolution. Inference in BMD is typically performed using Markov chain Monte Carlo (MCMC) methods within a Bayesian framework.

Given a data matrix $\bA$, bilinear matrix decomposition seeks factors $\bW$ and $\bZ$ such that  $\bA=\bW\bZ$ or $\bA\approx\bW\bZ$.
In the Bayesian setting, this becomes a problem of inferring the posterior distributions over the latent variables $\bW$ and $\bZ$ given the observed data $\bA$.
Priors are placed on  $\bW$ and $\bZ$, in which case 
we can either try to infer a \textit{point estimate} of a \textit{maximum likelihood estimator} by 
maximizing the likelihood $ \mathop{\max}_{\bW,\bZ} p(\bA\mid\bW,\bZ)$; or of a \textit{maximum a posteriori estimators} by $ \mathop{\max}_{\bW,\bZ} p(\bW,\bZ \mid \bA)$; or compute  the full posterior distribution $p(\bW,\bZ \mid\bA)$.

The choice of likelihood reflects the nature of the data: for example, a Poisson likelihood is suitable for count data, while a Gaussian likelihood is appropriate for real-valued or nonnegative observations. Similarly, priors are selected based on structural constraints---such as nonnegativity, sparsity, count support, or ordinality---and are placed on the entries of $\bW$ and $\bZ$. From a non-probabilistic perspective, the likelihood corresponds to a cost function (e.g., mean squared error), and the priors act as regularization terms that encourage desired properties like sparsity.
Within the scope of this book, our goal is to expand the repertoire of Bayesian matrix factorization algorithms.

\index{Numerical rank}
In numerical matrix decomposition methods \citep{lu2021numerical}, a matrix decomposition task on matrix $\bA$ can be cast as,
\begin{itemize}
\item $\bA=\bQ\bU$: where $\bQ$ is an orthogonal matrix spanning the same column space as $\bA$, and $\bU$ is a simpler, often sparse matrix used to reconstruct $\bA$.
\item $\bA=\bQ\bT\bQ^\top$: where $\bQ$ is orthogonal such that $\bA$ and $\bT$ are \textit{similar matrices} that share the same properties (e.g., eigenvalues, sparsity). Moreover, working on $\bT$ is an easier task compared to that on $\bA$.
\item $\bA=\bU\bT\bV$: where $\bU$ and $\bV$ are orthogonal matrices such that the columns of $\bU$ and the rows of $\bV$ constitute an orthonormal basis for the column space and row space of $\bA$, respectively.
\item $\underset{M\times N}{\bA}=\underset{M\times R}{\bB}\gapthree \underset{R\times N}{\bC}$: where $\bB$ and $\bC$ are full rank matrices that can reduce the memory storage of $\bA$. In practice, a low-rank approximation $\underset{M\times N}{\bA}\approx \underset{M\times K}{\bD}\gapthree \underset{K\times N}{\bF}$ can be employed, where $K<R$ is called the \textit{numerical rank} of the matrix such that  the matrix can be stored much more inexpensively and can be multiplied rapidly with
vectors or other matrices.
An approximation of the form $\bA=\bD\bF$ is useful for storing the matrix $\bA$ more frugally (we can store $\bD$ and $\bF$ using $K(M+N)$ floating-point numbers, as opposed to $MN$ numbers for storing $\bA$), for efficiently computing a matrix-vector product $\bb = \bA\bx$ (via $\bc = \bF\bx$ and $\bb = \bD\bc$), for data interpretation, and much more.
\end{itemize}
In contrast, Bayesian matrix decomposition typically focuses on simpler factorization forms, such as low-rank real-valued factorization, nonnegative factorization, models tailored for count or ordinal data, and Bayesian interpolative decomposition (ID).

The primary aim of this book is to provide a self-contained introduction to the foundational concepts and mathematical tools in Bayesian inference and matrix analysis, thereby paving the way for a seamless presentation of matrix decomposition (or factorization) techniques and their applications in later sections. That said, we acknowledge that it is impossible to cover all relevant and interesting developments in Bayesian matrix decomposition within this volume. Due to space constraints, topics such as high-order Bayesian tensor decomposition, nonparametric matrix factorization, and detailed treatments of variational inference for Bayesian models are not included. Readers seeking deeper coverage are encouraged to consult specialized literature in Bayesian analysis; notable references include \citet{rai2015leveraging, qian2016bayesian, lu2021numerical, takayama2022bayesian}.

\paragraph{Notation and preliminaries.} 
In the remainder of this section, we introduce and review some fundamental concepts from linear algebra that are relevant to matrix factorization. Additional important concepts will be defined and discussed as needed for clarity throughout the text.
Readers who already have a solid background in matrix analysis may choose to skip this section.
For simplicity, we restrict our discussion to real-valued matrices unless otherwise stated. 

In all cases, scalars will be denoted in a non-bold font possibly with subscripts (e.g., $a$, $\alpha$, $\alpha_i$). We will use \textbf{boldface} lowercase letters possibly with subscripts to denote vectors (e.g., $\bmu$, $\bx$, $\bx_n$, $\bz$) and
\textbf{boldface} uppercase letters possibly with subscripts to denote matrices (e.g., $\bA$, $\bL_j$). The $i$-th element of a vector $\bz$ will be denoted by $z_i$ in the non-bold font.
In the meantime, the \textit{normal fonts} of scalars denote  \textbf{random variables} (e.g., $\textnormal{a}$ and $\textnormal{b}_1$ are random variables, while italics $a$ and $b_1$ are scalars); 
the normal fonts of \textbf{boldface} lowercase letters possibly with subscripts denote \textbf{random vectors} (e.g., $\rva$ and $\rvb_1$ are random vectors, while italics $\ba$ and $\bb_1$ are vectors); 
and the normal fonts of \textbf{boldface} uppercase letters possibly with subscripts denote \textbf{random matrices} (e.g., $\rmA$ and $\rmB_1$ are random matrices, while italics $\bA$ and $\bB_1$ are matrices).

\index{Matlab-style notation}
Subarrays are formed when a subset of the indices is fixed.
{The $i$-th row and $j$-th column value of matrix $\bA$ (i.e., entry ($i,j$) of $\bA$) will be denoted by $a_{ij}$}. Furthermore, it will be helpful to utilize the \textbf{Matlab-style notation}, the $i$-th row to the $j$-th row and the $k$-th column to the $m$-th column submatrix of the matrix $\bA$ will be denoted by $\bA_{i:j,k:m} \equiv \bA[i:j,k:m]$. 
A colon is used to indicate all elements of a dimension, e.g., $\bA_{:,k:m}\equiv \bA[:,k:m]$ denotes the $k$-th column to the $m$-th column of the matrix $\bA$, and $\bA_{:,k}\equiv \bA[:,k]$ denotes the $k$-th column of $\bA$. Alternatively, the $k$-th column of $\bA$ may be denoted more compactly by $\ba_k$. 

When the index is not continuous, given ordered subindex sets $\sI$ and $\sJ$, $\bA[\sI, \sJ]$ denotes the submatrix of $\bA$ obtained by extracting the rows and columns of $\bA$ indexed by $\sI$ and $\sJ$, respectively; and $\bA[:, \sJ]$ denotes the submatrix of $\bA$ obtained by extracting the columns of $\bA$ indexed by $\sJ$, where the $[:, \sJ]$ syntax in this expression selects all rows from $\bA$ and only the columns specified by the indices in $\sJ$.

\begin{definition}[Matlab Notation]\label{definition:matlabnotation}
Suppose $\bA\in \real^{M\times N}$, and $\sI=[i_1, i_2, \ldots, i_K]$ and $\sJ=[j_1, j_2, \ldots, j_L]$ are two index vectors. 
Then $\bA[\sI,\sJ]$ denotes the $K\times L$ submatrix
$$
\bA[\sI,\sJ]=
\begin{bmatrix}
a_{i_1,j_1} & a_{i_1,j_2} &\ldots & a_{i_1,j_L}\\
a_{i_2,j_1} & a_{i_2,j_2} &\ldots & a_{i_2,j_L}\\
\vdots & \vdots&\ddots & \vdots\\
a_{i_K,j_1} & a_{i_K,j_2} &\ldots & a_{i_K,j_L}\\
\end{bmatrix}.
$$
Whilst, $\bA[\sI,:]$ denotes a $K\times N$ submatrix, and $\bA[:,\sJ]$ denotes a $M\times L$ submatrix analogously. We should also notice the range of the index:
$$
\left\{
\begin{aligned}
0&\leq \min(\sI) \leq \max(\sI)\leq M;\\
0&\leq \min(\sJ) \leq \max(\sJ)\leq N.
\end{aligned}
\right.
$$
\end{definition}

And in all cases, vectors are formulated in a column rather than in a row. A row vector will be denoted by a transpose of a column vector such as $\ba^\top$. A specific column vector with values is split by the semicolon symbol $``;"$, e.g., $\bx=[1;2;3]$ is a column vector in $\real^3$. Similarly, a specific row vector with values is split by the comma symbol $``,"$, e.g., $\by=[1,2,3]$ is a row vector with 3 values. Furthermore, a column vector can be denoted by the transpose of a row vector e.g., $\by=[1,2,3]^\top$ is a column vector.

The transpose of a matrix $\bA$ will be denoted by $\bA^\top$, and its inverse will be denoted by $\bA^{-1}$. 
We will denote the $P \times P$ identity matrix by $\bI_P$ (or simply $\bI$ when the dimension is clear from context). 
A vector or matrix of all zeros will be denoted by a \textbf{boldface} zero $\bzero$, whose size should be clear from context, or we denote $\bzero_P$ to be the vector of all zeros with $P$ entries.
Similarly, a vector or matrix of all ones will be denoted by a \textbf{boldface} one $\bone$, whose size is clear from  context, or we denote $\bone_P$ to be the vector of all ones with $P$ entries.
Subscripts on $\bI$, $\bzero$, and $\bone$ are often omitted when the dimensions are unambiguous.

\subsection*{Linear Algebra}

\begin{definition}[Eigenvalue, Eigenvector\index{Eigenvalue}\index{Eigenvector}\index{Eigenpair}]
Let $\bA\in\complex^{N\times N}$. A scalar $\lambda \in \complex$ is called a \textit{(right) eigenvalue} (also known as a  \textit{proper value}, or \textit{characteristic value}) of $\bA$ if there exists  a nonzero vector $\bu \in \complex^N$ such that
\begin{equation*}
\bA \bu = \lambda \bu.
\end{equation*}
In this case, $\bu$ is called a \textit{(right) eigenvector} of $\bA$ associated with $\lambda$.
\end{definition}
For simplicity, we restrict our attention to real-valued matrices unless otherwise stated. Unless explicitly noted, all eigenvalues discussed are assumed to be real as well.  

Intuitively, an eigenvector $\bu$ of a matrix $\bA$ represents a direction in $\real^N$ that remains  invariant under the linear transformation defined by $\bA$: applying $\bA$ to $\bu$ does not change its direction---it only scales it by the corresponding eigenvalue $\lambda$.
Note that while real matrices can have complex eigenvalues in general, symmetric real matrices always have real eigenvalues (see Theorem~\ref{theorem:spectral_theorem}).

The pair $(\lambda, \bu)$  is commonly  referred to as an \textit{eigenpair}. 
Importantly, eigenvectors are not unique:  if  $\bu$ is an eigenvector, then so is any nonzero scalar multiple $\eta\bu$ (with $\eta\in\real\setminus \{\bzero\}$).
To resolve this ambiguity, eigenvectors are typically normalized---for example, by requiring $\normtwo{\bu}=1$ (unit $\ell_2$-norm; see Definition~\ref{definition:vec_l2_norm}). 
Moreover, since both $\bu$ and $-\bu$ correspond to the same eigenvalue and represent the same direction up to sign, it is common to fix the sign by convention (e.g., by requiring the first nonzero entry of $\bu$ to be positive).

In fact, real-valued matrices may have complex eigenvalues. However, all  eigenvalues of symmetric matrices are real (a consequence of the spectral theorem; see see Theorem~\ref{theorem:spectral_theorem}).

\begin{definition}[Spectrum and Spectral Radius\index{Spectrum}\index{Spectral radius}]\label{definition:spectrum}
The set of all eigenvalues of $\bA$ is called the \textit{spectrum} of $\bA$ and is denoted by $\Lambda(\bA)$. The largest magnitude of the eigenvalues is known as the \textit{spectral radius}, denoted $\rho(\bA)$:
$$
\rho(\bA) = \mathop{\max}_{\lambda\in \Lambda(\bA)}  \abs{\lambda}.
$$
\end{definition}

In linear algebra, every vector space admits a basis, and any vector in the space can be expressed as a linear combination of basis vectors. We now define key concepts related to subspaces.

\begin{definition}[Subspace and Span\index{Subspace}\index{Span}]
A nonempty subset $\mathcalV\subseteq \real^N$ is called a \textit{subspace} if  for all $\ba,\bb\in \mathcalV$ and all scalars $x,y\in \real$, the linear combination $x\ba+y\bb$ also belongs to $\mathcalV$.
Moreover, a set of vectors  $\{\ba_1, \ba_2, \ldots,$ $\ba_M\}$ is said to \textit{span} a subspace $\mathcalV$ if every vector $\bv\in\mathcalV $ can be written as a linear combination of these vectors.
\end{definition}

In this context, we will often use the idea of the linear independence of a set of vectors. Two equivalent definitions are given as follows.
\begin{definition}[Linearly Independent\index{Linearly independent}]
A set of vectors $\{\ba_1, \ba_2, \ldots, \ba_M\}$ is  \textit{linearly independent} if  the only solution to $x_1\ba_1+x_2\ba_2+\ldots+x_M\ba_M=\bzero $ is $x_1=x_2=\ldots=x_M=0$. 
Equivalently, the set is linearly independent if $\ba_1\neq \bzero$, and for each $k>1$, the vector $\ba_k$ does not lie in the span of $\{\ba_1, \ba_2, \ldots, \ba_{k-1}\}$.
\end{definition}

\begin{definition}[Basis and Dimension\index{Basis}\index{Dimension}]
A set of vectors $\{\ba_1, \ba_2, \ldots, \ba_M\}$ is  a \textit{basis} for a subspace  $\mathcalV$ if the vectors are linearly independent and span $\mathcalV$. 
Every basis of a given subspace contains the same number of vectors; this number is called the \textit{dimension} of $\mathcalV$. 
By convention, the trivial subspace $\{\bzero\}$ has dimension zero.
Furthermore, any subspace of positive dimension admits an orthogonal basis---that is, a basis in which all vectors are mutually orthogonal (see Definition~\ref{definition:orthogo_mat}).
\end{definition}

\begin{definition}[Column Space (Range)\index{Column space}]
For an $M\times N$ real matrix $\bA$, the \textit{column space} (or \textit{range}) of $\bA$ is the subspace of $\real^M$ spanned by its columns:
\begin{equation*}
\cspace (\bA) = \{ \by\in \real^M: \exists\, \bx \in \real^N, \, \by = \bA \bx \}.
\end{equation*}
The \textit{row space} of  $\bA$ is the column space of  $\bA^\top$:
\begin{equation*}
\cspace (\bA^\top) = \{ \bx\in \real^N: \exists\, \by \in \real^M, \, \bx = \bA^\top \by \}.
\end{equation*}
\end{definition}

\begin{definition}[Null Space (Nullspace, Kernel)\index{Null space}]
For an $M\times N$ real matrix $\bA$,  the \textit{null space} (or \textit{nullspace}, \textit{kernel}) of  $\bA$ is defined as
\begin{equation*}
\nspace (\bA) = \{\by \in \real^N:  \, \bA \by = \bzero \}.
\end{equation*}
Similarly, the null space of  $\bA^\top$ is
\begin{equation*}
\nspace (\bA^\top) = \{\bx \in \real^M:  \, \bA^\top \bx = \bzero \}.
\end{equation*}
\end{definition}

Both the column space of $\bA$ and the null space of $\bA^\top$ are subspaces of $\real^M$. In fact, every vector in $\nspace(\bA^\top)$ is orthogonal  to vectors in $\cspace(\bA)$, and vice versa.
Similarly, every vector in $\nspace(\bA)$ is also orthogonal to vectors in $\cspace(\bA^\top)$, and vice versa.

\begin{definition}[Rank\index{Rank}]
The \textit{rank} of a matrix $\bA\in \real^{M\times N}$ is the dimension of its column space. 
That is, the rank of $\bA$ is equal to the maximum number of linearly independent columns of $\bA$, and is also the maximum number of linearly independent rows of $\bA$. The matrix $\bA$ and its transpose $\bA^\top$ have the same rank. 
A matrix is said to have \textit{full rank} if its rank equals $\min\{M,N\}$.
As a special case, for any nonzero vectors  $\bu \in \real^M$ and  $\bv \in \real^N$, the outer product  $\bu\bv^\top\in \real^{M\times N}$ is an $M\times N$ matrix of rank 1. 
In summary, the rank of $\bA$ equals:
\begin{itemize}
\item number of linearly independent columns;
\item number of linearly independent rows.
\end{itemize}
Remarkably, these two quantities are always equal (see Theorem~\ref{lemma:equal-dimension-rank}).
\end{definition}

\begin{definition}[Orthogonal Complement in General]
The \textit{orthogonal complement} $\mathcalV^\perp$ of a subspace $\mathcalV\subseteq\real^N$ is the set of all vectors orthogonal to every vector in $\mathcalV$. That is,
$$
\mathcalV^\perp = \{\bv :  \bv^\top\bu=0, \,\,\, \forall \bu\in \mathcalV  \}.
$$
The two subspaces are disjoint and together span the entire space $\real^N$. 
Their dimensions satisfy $\dim(\mathcalV)+\dim(\mathcalV^\perp)=N$, and
 $(\mathcalV^\perp)^\perp=\mathcalV$.
\end{definition}

\begin{definition}[Orthogonal Complement of Column Space]
For an $M\times N$ real matrix $\bA$, the orthogonal complement of $\cspace(\bA)$, denoted $\cspace^{\bot}(\bA)$, is the subspace defined as:
\begin{equation*}
\begin{aligned}
\cspace^{\bot}(\bA) &= \{\by\in \real^M: \, \by^\top \bA \bx=\bzero, \, \forall \bx \in \real^N \} \\
&=\{\by\in \real^M: \, \by^\top \bv = \bzero, \, \forall \bv \in \cspace(\bA) \}.
\end{aligned}
\end{equation*}
\end{definition}

These ideas lead to the four fundamental subspaces associated with any matrix $\bA\in \real^{M\times N}$ of rank $R$, as described in Theorem~\ref{theorem:fundamental-linear-algebra}.
To establish the fundamental theorem of linear algebra, we first prove a key result: the equality of row rank and column rank.

\begin{theoremHigh}[Row Rank Equals Column Rank\index{Rank}\index{Matrix rank}]\label{lemma:equal-dimension-rank}
For any matrix $\bA\in \real^{M\times N}$, the dimension of its column space equals the dimension of its row space. That is, the row rank and column rank of $\bA$ are equal.
\end{theoremHigh}

\begin{proof}[{of Theorem~\ref{lemma:equal-dimension-rank}}]
We first observe  that the null space of $\bA$ is orthogonal complementary to the row space of $\bA$: $\nspace(\bA) \bot \cspace(\bA^\top)$ (where the row space of $\bA$ is exactly the column space of $\bA^\top$), that is, vectors in the null space of $\bA$ are orthogonal to vectors in the row space of $\bA$. To see this, suppose $\bA$ has rows $\{\ba_1^\top, \ba_2^\top, \ldots, \ba_M^\top\}$ and $\bA=[\ba_1^\top; \ba_2^\top; \ldots; \ba_M^\top]$ is the row partition. For any vector $\bx\in \nspace(\bA)$, we have $\bA\bx = \bzero$, that is, $[\ba_1^\top\bx; \ba_2^\top\bx; \ldots; \ba_M^\top\bx]=\bzero$. 
Since the row space of $\bA$ is spanned by $\{\ba_1^\top, \ba_2^\top, \ldots, \ba_M^\top\}$, then $\bx$ is perpendicular to any vectors from $\cspace(\bA^\top)$, which means $\nspace(\bA) \bot \cspace(\bA^\top)$.

Now suppose  the dimension of row space of $\bA$ is $R$. \textbf{{Let $\{\br_1, \br_2, \ldots, \br_R\}$ be a set of vectors in $\real^N$ and form a basis for the row space}}. Then the $R$ vectors $\{\bA\br_1, \bA\br_2, \ldots, \bA\br_R\}$ are in the column space of $\bA$, which are linearly independent. To see this, suppose we have a linear combination of the $R$ vectors: $x_1\bA\br_1 + x_2\bA\br_2+ \ldots+ x_R\bA\br_R=\bzero$, that is, $\bA(x_1\br_1 + x_2\br_2+ \ldots+ x_R\br_R)=\bzero$, and the vector $\bv=x_1\br_1 + x_2\br_2+ \ldots+ x_R\br_R$ lies in null space of $\bA$. 
Moreover, since $\{\br_1, \br_2, \ldots, \br_R\}$ is a basis for the row space of $\bA$, $\bv$ is thus also in the row space of $\bA$. 
We have shown that vectors from the null space of $\bA$ is perpendicular to vectors from the row space of $\bA$; thus, it holds that $\bv^\top\bv=0$ and $x_1=x_2=\ldots=x_R=0$. Then \textbf{$\bA\br_1, \bA\br_2, \ldots, \bA\br_R$ are in the column space of $\bA$, and they are linearly independent}, which means the dimension of the column space of $\bA$ is larger than $R$. This result shows that \textbf{row rank of $\bA\leq $ column rank of $\bA$}. 

Applying the same argument to $\bA^\top$ yields \textbf{column rank of $\bA\leq $ row rank of $\bA$}. 
Hence, the two ranks are equal.
\end{proof}

An important consequence of this proof is that if $\{\br_1, \br_2, \ldots, \br_R\}$ is a basis for the row space of $\bA$, then \textbf{$\{\bA\br_1, \bA\br_2, \ldots, \bA\br_R\}$} forms a basis for the column space. We state this formally:

\begin{lemma}[Column Basis from Row Basis]\label{lemma:column-basis-from-row-basis}
Let $\bA\in \real^{M\times N}$. 
If $\{\br_1, \br_2, \ldots, \br_R\}$ is a basis for the row space of $\bA$, then $\{\bA\br_1, \bA\br_2, \ldots, \bA\br_R\}$ is a basis for the column space of $\bA$.
\end{lemma}

It is straightforward to verify that any vector in the row space of $\bA$ is orthogonal to any vector in  $\nspace(\bA)$. 
Indeed, if  $\bx_n \in \nspace(\bA)$, then $\bA\bx_n = \bzero$, meaning  $\bx_n$ is orthogonal to every row of $\bA$, and hence to the entire row space.

Similarly, any vector in  $\cspace(\bA)$ is orthogonal to any vector in $\nspace(\bA^\top)$. Moreover, $\cspace(\bA)$ and $\nspace(\bA^\top)$ together span $\real^M$. 
This is the essence of the \textit{fundamental theorem of linear algebra}, which consists of two parts: orthogonality and dimensionality.

While orthogonality follows directly from the definitions, the dimensional relationship requires proof. Specifically, if the row space has dimension $R$, then the null space has dimension $N-R$. This is formalized below.

\index{Fundamental spaces}
\index{Fundamental theorem of linear algebra}
\begin{figure}[h!]
\centering
\includegraphics[width=0.98\textwidth]{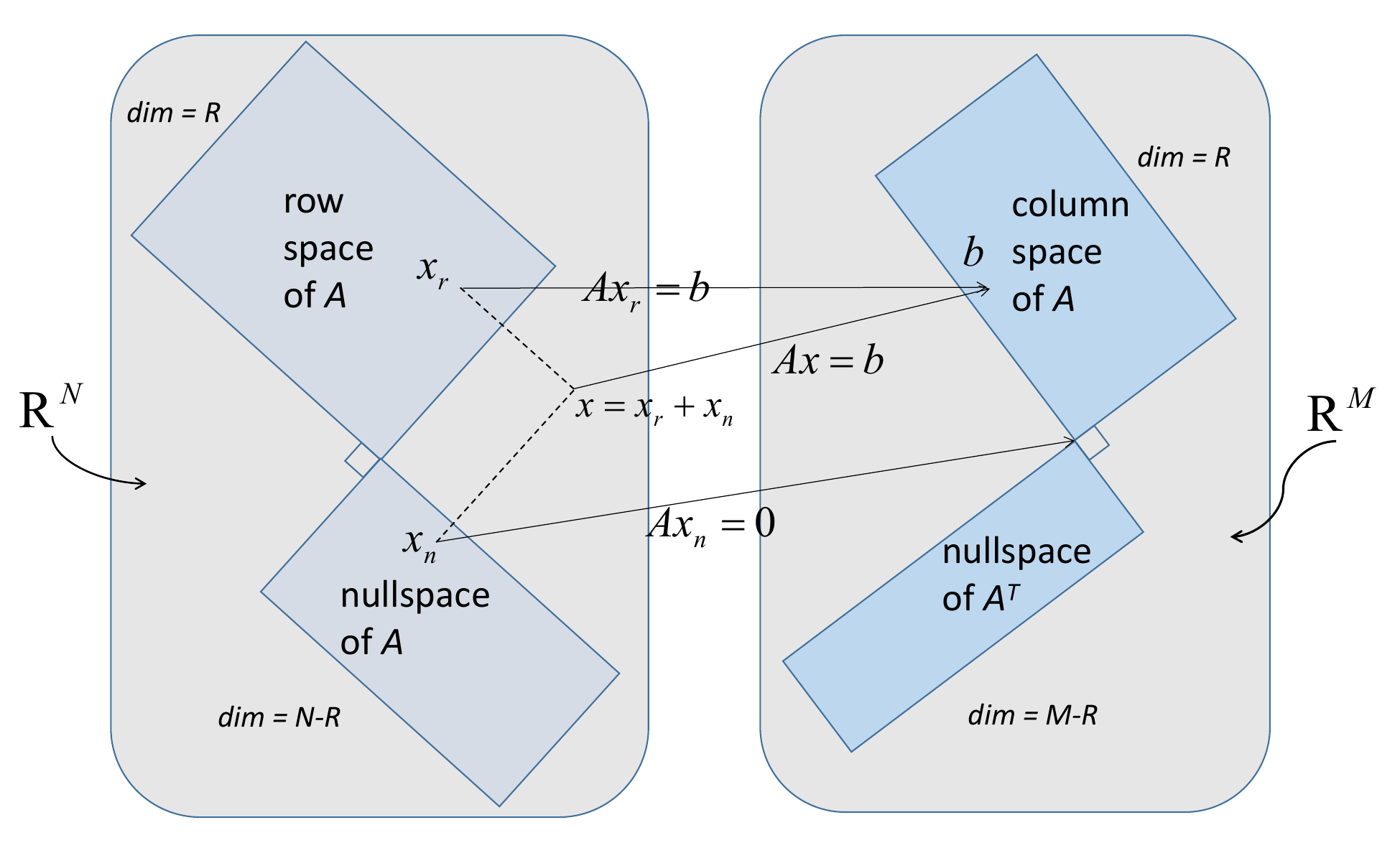}
\caption{Two pairs of orthogonal subspaces in $\real^N$ and $\real^M$. $\dim(\cspace(\bA^\top)) + \dim(\nspace(\bA))=N$ and $\dim(\nspace(\bA^\top)) + \dim(\cspace(\bA))=M$. The null space component maps to zero as $\bA\bx_n = \bzero \in \real^M$. The row space component maps into the column space as $\bA\bx_R = \bA(\bx_R+\bx_n)=\bb \in \cspace(\bA)$.}
\label{fig:lafundamental}
\end{figure}
\begin{theoremHigh}[The Fundamental Theorem of Linear Algebra]\label{theorem:fundamental-linear-algebra}

\textit{Orthogonal Complement} and \textit{Rank-Nullity Theorem}: for any matrix $\bA\in \real^{M\times N}$ of rank $R$, the following hold:
\begin{itemize}
\item $\nspace(\bA)$ is orthogonal complement to the row space $\cspace(\bA^\top)$ in $\real^N$: $\dim(\nspace(\bA))+\dim(\cspace(\bA^\top))=N$;

\item $\nspace(\bA^\top)$ is orthogonal complement to the column space $\cspace(\bA)$ in $\real^M$: $\dim(\nspace(\bA^\top))+\dim(\cspace(\bA))=M$;

\item $\dim(\cspace(\bA^\top)) = \dim(\cspace(\bA)) = R$, that is, $\dim(\nspace(\bA)) = N-R$ and $\dim(\nspace(\bA^\top))=M-R$.
\end{itemize}
\end{theoremHigh}

\begin{proof}[of Theorem~\ref{theorem:fundamental-linear-algebra}]
Following from the proof of Theorem~\ref{lemma:equal-dimension-rank}. Let $\br_1, \br_2, \ldots, \br_R$ be a set of vectors in $\real^N$ that form a basis for the row space, then \textbf{$\bA\br_1, \bA\br_2, \ldots, \bA\br_R$ is a basis for the column space of $\bA$}. Let $\bn_1, \bn_2, \ldots, \bn_K \in \real^N$ form a basis for the null space of $\bA$. Following again from the proof of Theorem~\ref{lemma:equal-dimension-rank}, $\nspace(\bA) \bot \cspace(\bA^\top)$, thus, $\br_1, \br_2, \ldots, \br_R$ are perpendicular to $\bn_1, \bn_2, \ldots, \bn_K$. Then, $\{\br_1, \br_2, \ldots, \br_R, \bn_1, \bn_2, \ldots, \bn_K\}$ is linearly independent in $\real^N$.

For any vector $\bx\in \real^N $, $\bA\bx$ is in the column space of $\bA$. Then it can be represented as a linear combination of $\bA\br_1, \bA\br_2, \ldots, \bA\br_R$: $\bA\bx = \sum_{i=1}^{R}a_i\bA\br_i$ which states that $\bA(\bx-\sum_{i=1}^{R}a_i\br_i) = \bzero$, and $\bx-\sum_{i=1}^{R}a_i\br_i$ is thus in $\nspace(\bA)$. Since $\{\bn_1, \bn_2, \ldots, \bn_K\}$ is a basis for the null space of $\bA$, $\bx-\sum_{i=1}^{R}a_i\br_i$ can be represented by a linear combination of $\bn_1, \bn_2, \ldots, \bn_K$: $\bx-\sum_{i=1}^{R}a_i\br_i = \sum_{j=1}^{k}b_j \bn_j$, i.e., $\bx=\sum_{i=1}^{R}a_i\br_i + \sum_{j=1}^{k}b_j \bn_j$. That is, any vector $\bx\in \real^N$ can be represented by $\{\br_1, \br_2, \ldots, \br_R, \bn_1, \bn_2, \ldots, \bn_K\}$, and the set forms a basis for $\real^N$. Thus, the dimension sum to $N$: $R+K=N$, i.e., $\dim(\nspace(\bA))+\dim(\cspace(\bA^\top))=N$. Similarly, we can prove $\dim(\nspace(\bA^\top))+\dim(\cspace(\bA))=M$.
\end{proof}

Figure~\ref{fig:lafundamental} demonstrates two pairs of such orthogonal subspaces and shows how $\bA$ takes $\bx$ into the column space. The dimensions of the row space of $\bA$ and the null space of $\bA$ add to $N$. And the dimensions of the column space of $\bA$ and the null space of $\bA^\top$ add to $M$. The null space component goes to zero as $\bA\bx_{\bn} = \bzero \in \real^M$, which is the intersection of the column space of $\bA$ and the null space of $\bA^\top$. 
Conversely, the row space component goes to the column space as $\bA\bx_{\br} = \bA(\bx_{\br} + \bx_{\bn})=\bb\in \real^M$.

\begin{definition}[Orthogonal and Semi-Orthogonal Matrix\index{Orthogonal matrix}\index{Semi-orthogonal matrix}]\label{definition:orthogo_mat}
A real square matrix $\bQ$ is called an \textit{orthogonal matrix} if its inverse equals its transpose, i.e., $\bQ^{-1}=\bQ^\top$ and $\bQ\bQ^\top = \bQ^\top\bQ = \bI$. In other words, suppose $\bQ=[\bq_1, \bq_2, \ldots, \bq_N]$, where $\bq_i \in \real^N$ for all $i \in \{1, 2, \ldots, N\}$, then $\bq_i^\top \bq_j = \delta(i,j)$ with $\delta(i,j)$ denoting the Kronecker delta function. If $\bQ$ contains only $\gamma$ of these columns with $\gamma<N$ (called a \textit{semi-orghogonal matrix}), then $\bQ^\top\bQ = \bI_\gamma$ stills holds with $\bI_\gamma$ denoting the $\gamma\times \gamma$ identity matrix. But $\bQ\bQ^\top=\bI$ will not be true. For any vector $\bx$, the orthogonal matrix will preserve the length: $\norm{\bQ\bx} = \norm{\bx}$.
\end{definition}

\begin{definition}[Permutation Matrix]\label{definition:permutation-matrix}
A \textit{permutation matrix} $\bP$ is a square binary matrix with exactly one entry equal to 1 in each row and each column, and zeros elsewhere.
\paragraph{Row Point.} That is, the permutation matrix $\bP$ has the rows of the identity $\bI$ in any order and the order decides the sequence of the row permutation: the product $\bP\bA$ rearranges the rows of $\bA$ according to the order of the rows in $\bP$.
\paragraph{Column Point.} Or, equivalently, the permutation matrix $\bP$ has the columns of the identity $\bI$ in any order and the order decides the sequence of the column permutation: the product $\bA\bP$ rearranges the columns of $\bA$ according to the order of the columns in $\bP$.
\end{definition}

A permutation matrix $\bP\in\real^{N\times N}$ can be compactly represented by an index vector $\sJ \in \integer_+^N$ where $\sJ$ is a permutation of $\{1,2,\ldots,N\}$, such that $\bP = \bI[:, \sJ]$, where $\bI$ is the $N\times N$ identity matrix. 
Note that the entries of $\sJ$ are distinct integers from 1 to $N$, so their sum is always $1+2+\ldots+N= {(N^2+N)}/{2}$.

\begin{example}[Permutation]
Suppose, 
$$\bA=\begin{bmatrix}
1 & 2&3\\
4&5&6\\
7&8&9
\end{bmatrix}
\qquad \text{and} \qquad
\bP=\begin{bmatrix}
	&1&\\
	&&1\\
	1&&
\end{bmatrix}.
$$
The row and columns permutations are given, respectively, by 
$$
\bP\bA = \begin{bmatrix}
4&5&6\\
7&8&9\\
1 & 2&3\\
\end{bmatrix}
\qquad\text{and}\qquad
\bA\bP = \begin{bmatrix}
3 & 1 & 2 \\
6 & 4 & 5\\
9 & 7 & 8
\end{bmatrix},
$$
where the order of the rows of $\bA$ appearing in $\bP\bA$ matches the order of the rows of $\bI$ in $\bP$, and the order of the columns of $\bA$ appearing in $\bA\bP$ matches the order of the columns of $\bI$ in $\bP$.
\end{example}

\begin{definition}[Positive Definite and Positive Semidefinite\index{Positive definite}\index{Positive semidefinite}]\label{definition:psd-pd-defini}
A matrix $\bA\in \real^{N\times N}$ is called \textit{positive definite (PD)} if $\bx^\top\bA\bx>0$ for all nonzero $\bx\in \real^N$.
And a matrix $\bA\in \real^{N\times N}$ is called \textit{positive semidefinite (PSD)} if $\bx^\top\bA\bx \geq 0$ for all $\bx\in \real^N$. 
\footnote{
In this book, a positive definite or a semidefinite matrix is always assumed to be symmetric, i.e., the notion of a positive definite matrix or semidefinite matrix is only interesting for symmetric matrices.
}
\footnote{A symmetric matrix $\bA\in\real^{N\times N}$ is called \textit{negative definite} (ND) if $\bx^\top\bA\bx<0$ for all nonzero $\bx\in\real^N$; 
a symmetric matrix $\bA\in\real^{N\times N}$ is called \textit{negative semidefinite} (NSD) if $\bx^\top\bA\bx\leq 0$ for all $\bx\in\real^N$;
and a symmetric matrix $\bA\in\real^{N\times N}$ is called \textit{indefinite} (ID) if there exist $\bx$ and $\by\in\real^N$ such that $\bx^\top\bA\bx<0$ and $\by^\top\bA\by>0$.
}
\end{definition}

We can establish that a matrix $\bA$ is positive definite if and only if it possesses exclusively \textit{positive eigenvalues}. Similarly, a matrix $\bA$ is positive semidefinite if and only if it exhibits solely \textit{nonnegative eigenvalues}. See Problem~\ref{problem:psd_eigen}.

In introductory linear algebra, the following equivalence is fundamental:
\begin{remark}[List of Equivalence of Nonsingularity for a Matrix]\label{remark:equiva_nonsingular}
For a square matrix $\bA\in \real^{N\times N}$, the following claims are equivalent:
\begin{itemize}
\item $\bA$ is nonsingular;
\item $\bA$ is invertible, i.e., $\bA^{-1}$ exists;
\item $\bA\bx=\bb$ has a unique solution $\bx = \bA^{-1}\bb$;
\item $\bA\bx = \bzero$ has a unique, trivial solution: $\bx=\bzero$;
\item Columns of $\bA$ are linearly independent;
\item Rows of $\bA$ are linearly independent;
\item $\det(\bA) \neq 0$; 
\item $\dim(\nspace(\bA))=0$;
\item $\nspace(\bA) = \{\bzero\}$, i.e., the null space is trivial;
\item $\cspace(\bA)=\cspace(\bA^\top) = \real^N$, i.e., the column space or row space span the entire $\real^N$;
\item $\bA$ has full rank $R=N$;
\item $\bA^\top\bA$ is symmetric positive definite;
\item $\bA$ has $N$ nonzero (positive) singular values;
\item All eigenvalues of $\bA$ are nonzero;
\end{itemize}
\end{remark}

Keeping these equivalences in mind is essential to avoid confusion in theoretical and computational contexts.
Similarly, the following remark characterizes singular matrices---particularly in the context of eigenvalues.
\begin{remark}[List of Equivalence of Singularity for a Matrix]
For a square matrix $\bA\in \real^{N\times N}$ with eigenpair $(\lambda, \bu)$, the following claims are equivalent:
\begin{itemize}\label{remark:equiva_singular}
\item $(\bA-\lambda\bI)$ is singular;
\item $(\bA-\lambda\bI)$ is not invertible;
\item $(\bA-\lambda\bI)\bx = \bzero$ has nonzero $\bx\neq \bzero$ solutions, and $\bx=\bu$ is one of such solutions;
\item $(\bA-\lambda\bI)$ has linearly dependent columns;
\item $\det(\bA-\lambda\bI) = 0$; 
\item $\dim(\nspace(\bA-\lambda\bI))>0$;
\item Null space of $(\bA-\lambda\bI)$ is nontrivial;
\item Columns of $(\bA-\lambda \bI)$ are linearly dependent;
\item Rows of $(\bA-\lambda \bI)$ are linearly dependent;
\item $(\bA-\lambda \bI)$ has rank $R<N$;
\item Dimension of column space = dimension of row space = $R<N$;
\item $(\bA-\lambda \bI)^\top(\bA-\lambda \bI)$ is symmetric semidefinite;
\item $(\bA-\lambda \bI)$ has fewer than $N$ nonzero singular values;
\item Zero is an eigenvalue of $(\bA-\lambda \bI)$.
\end{itemize}
\end{remark}

\index{Vector norm}
\index{Matrix norm}
\begin{definition}[Vector $\ell_1, \ell_2, \ell_\infty$, $\ell_p$-Norms]\label{definition:vec_l2_norm}
For a vector $\bx\in\real^N$, the \textit{$\ell_2$-norm} is defined as $\normtwo{\bx} = \sqrt{x_1^2+x_2^2+\ldots+x_N^2}$.
Similarly, the \textit{$\ell_1$-norm} can be obtained by 
$
\normone{\bx} = \sum_{n=1}^{N} \abs{x_n} .
$
And the \textit{$\ell_\infty$-norm} can be obtained by 
$
\norminf{\bx} = \mathop{\max}_{n=1,2,\ldots,N} \abs{x_n} .
$
More generally, the $\ell_p$-norm is defined as $\normp{\bx}=\sqrt[p]{ \sum_{n=1}^{N}\abs{x_n}^p  }$ for $p\geq 1$.
\end{definition}

For matrices, two common norms are the Frobenius norm and the spectral norm.
\begin{definition}[Matrix Frobenius Norm\index{Frobenius norm}]\label{definition:frobernius-in-svd}
The \textit{Frobenius norm} of a matrix $\bA\in \real^{M\times N}$ is defined as 
$$
\norm{\bA}_F = \sqrt{\sum_{m=1,n=1}^{M,N} (a_{mn})^2}=\sqrt{\trace(\bA\bA^\top)}=\sqrt{\trace(\bA^\top\bA)} = \sqrt{\sigma_1^2+\sigma_2^2+\ldots+\sigma_R^2}, 
$$
where $\sigma_1, \sigma_2, \ldots, \sigma_R$ are the nonzero singular values of $\bA$; see Theorem~\ref{theorem:reduced_svd_rectangular}.
\end{definition}

\begin{definition}[Matrix Spectral Norm]\label{definition:spectral_norm}
The \textit{spectral norm} of a matrix $\bA\in \real^{M\times N}$ is defined as 
$$
\norm{\bA}_2 = \mathop{\max}_{\bx\neq\bzero} \frac{\norm{\bA\bx}_2}{\norm{\bx}_2}  =\mathop{\max}_{\bu\in \real^N: \norm{\bu}_2=1}  \norm{\bA\bx}_2 ,
$$
which equals the largest singular value of  $\bA$, i.e., $\norm{\bA}_2 = \sigma_1(\bA)$.
\end{definition}

The Frobenius norm is the natural matrix analogue of the vector $\ell_2$-norm.
For simplicity, when the context is clear, we often omit subscripts for these two norms:: $\norm{\bA}=\normf{\bA}$ and $\normtwo{\bx}=\norm{\bx}$.
However, the subscript must be retained for the spectral norm: $\normtwo{\bA}$ should never be written as $\norm{\bA}$ to avoid ambiguity.

Finally, given a norm, we define open and closed balls as follows:
\begin{definition}[Open Ball, Closed Ball]
The \textit{open ball} centered at $\bc\in\real^N$ with radius $r$  is defined as 
$$
B(\bc, r) =\{\bx\in\real^N\mid \normtwo{\bx-\bc} <r\}.
$$
Similarly, the \textit{closed ball} with center $\bc\in\real^N$ and radius $r$  is defined as 
$$
B[\bc,r] =\{\bx\in\real^N\mid \normtwo{\bx-\bc} \leq r\}.
$$
\end{definition}

\subsection*{Singular Value Decomposition (SVD)}\label{section:SVD}

We introduce the \textit{singular value decomposition (SVD)} of a matrix in this section.
Before presenting the general form of the SVD, we first recall the spectral decomposition of a symmetric matrix.
The \textit{spectral theorem}---also known as the spectral decomposition for symmetric matrices---states that any real symmetric matrix has real eigenvalues and can be diagonalized using a real orthonormal basis of eigenvectors.
\footnote{For Hermitian matrices (the complex analogue), a similar result holds: eigenvalues are real, and diagonalization is possible via a complex orthonormal basis.}

\begin{theoremHigh}[Spectral Decomposition\index{Spectral decomposition}\index{Spectral theorem}]\label{theorem:spectral_theorem}
A real matrix $\bA \in \real^{N\times N}$ is symmetric if and only if there exist an orthogonal matrix $\bQ$ and a diagonal matrix $\bLambda$ such that
\begin{equation*}
\bA = \bQ \bLambda \bQ^\top,
\end{equation*}
where the columns of $\bQ = [\bq_1, \bq_2, \ldots, \bq_n]$ are  mutually orthonormal eigenvectors of $\bA$, and the diagonal entries of $\bLambda=\diag(\lambda_1, \lambda_2, \ldots, \lambda_n)$ are the corresponding eigenvalues of $\bA$, which are real. 
In particular, the following properties hold:
\begin{enumerate}
\item All eigenvalues of a symmetric matrix are \textbf{real}.
\item The eigenvectors can be chosen to form an \textbf{orthonormal} set.
\item The rank of $\bA$ equals the number of its nonzero eigenvalues.
\item If all eigenvalues are distinct, then the corresponding eigenvectors are linearly independent.
\end{enumerate}
\end{theoremHigh}
\begin{proof}
See \citet{lu2022matrix}.
\end{proof}

Using spectral decomposition, we can factor a symmetric matrix into a diagonal form. However, this approach does not extend to non-symmetric or non-square matrices, for which such a diagonalization is generally impossible.
The SVD overcomes this limitation. Instead of relying on a single orthogonal matrix of eigenvectors, the SVD expresses any real matrix as a product of two orthogonal matrices and a diagonal (or diagonal-like) matrix of singular values. We now state the SVD theorem.
\begin{theoremHigh}[Full Singular Value Decomposition\index{Singular value decomposition}]\label{theorem:full_svd_rectangular}\label{theorem:reduced_svd_rectangular}
Every real $M\times N$ matrix $\bA$ of rank $R$  admits a decomposition of the form
$$
\bA = \bU \bSigma \bV^\top,
$$ 
where $\bSigma\in \real^{M\times N}$ has the block structure $\bSigma=\footnotesize\begin{bmatrix}
\bSigma_R & \bzero \\
\bzero & \bzero
\end{bmatrix}$ with $\bSigma_R=\diag(\sigma_1, \sigma_2 \ldots, \sigma_R)\in \real^{R\times R}$, and the singular values satisfy $\sigma_1 \geq \sigma_2 \geq \ldots \geq \sigma_R$.
\begin{itemize}
\item The scalars  $\sigma_i$  are the nonzero \textit{singular values} of  $\bA$. Each $\sigma_i$ equals the positive square root of a nonzero eigenvalue of both $\bA^\top \bA$ and $\bA \bA^\top$.

\item $\bU\in \textcolor{black}{\real^{M\times M}}$ is an orthogonal matrix whose first $R$ columns are eigenvectors of $\bA \bA^\top$ corresponding to its $R$ nonzero eigenvalues; the remaining $M - R$ columns form an orthonormal basis for $\nspace(\bA^\top)$.

\item $\bV\in \textcolor{black}{\real^{N\times N}}$ is an orthogonal matrix whose first $R$ columns are eigenvectors of $\bA^\top \bA$ corresponding to its $R$ nonzero eigenvalues; the remaining $N - R$ columns form an orthonormal basis for $\nspace(\bA)$.

\item The columns of $\bU$ and $\bV$ are called the \textit{left and right singular vectors} of $\bA$, respectively. 
\end{itemize}
Moreover, the decomposition can be expressed as a sum of rank-one matrices: $ \bA = \bU \bSigma \bV^\top = \sum_{i=1}^R \sigma_i \bu_i \bv_i^\top$.
\end{theoremHigh}

\begin{figure}[h!]
\centering  
\vspace{-0.35cm}  
\subfigtopskip=2pt  
\subfigbottomskip=2pt  
\subfigcapskip=-5pt  
\subfigure[Reduced SVD decomposition.]{\label{fig:svdhalf}
\includegraphics[width=0.47\linewidth]{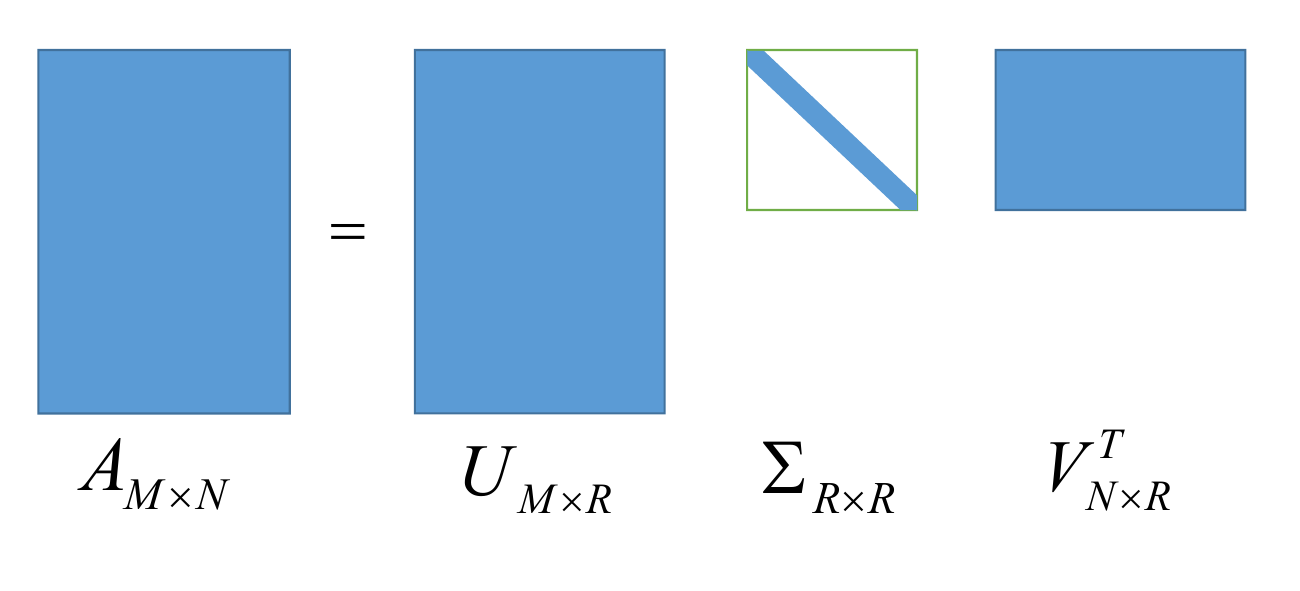}}
\quad 
\subfigure[Full SVD decomposition.]{\label{fig:svdall}
\includegraphics[width=0.47\linewidth]{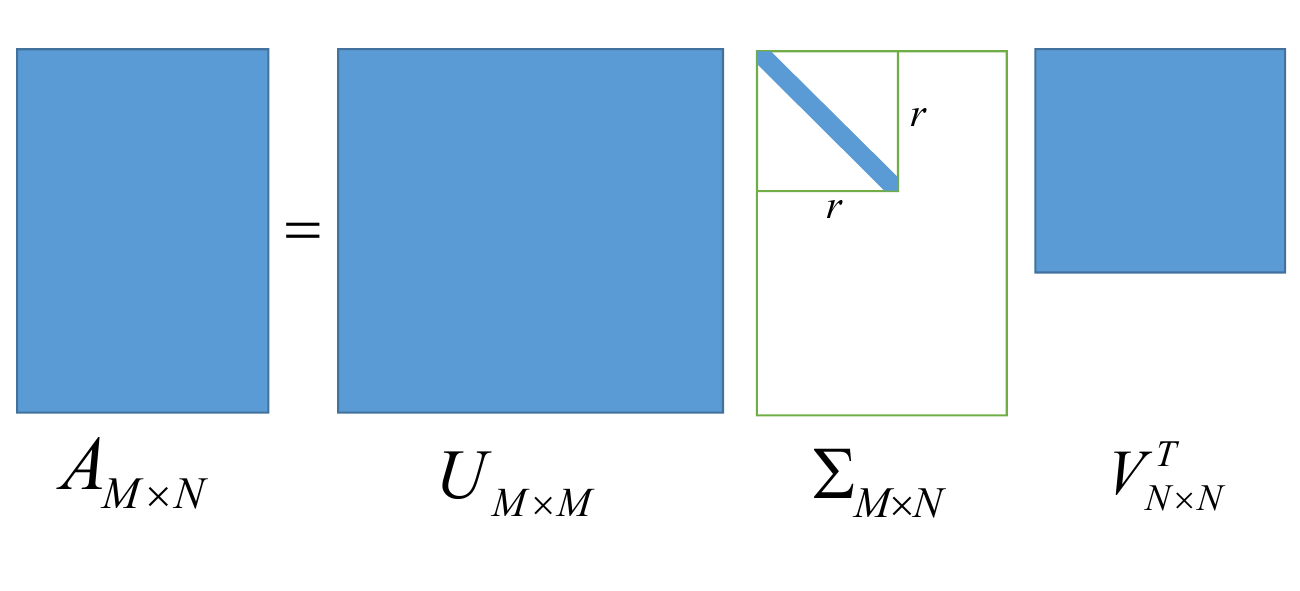}}
\caption{Comparison between the reduced and full SVD. White entries are zero, and blue entries are not necessarily zero}
\label{fig:svd-comparison}
\end{figure}

If $\bA$ has rank $R$, then its singular values satisfy $\sigma_1,  \sigma_2, \ldots, \sigma_R > 0$, and all remaining singular values (if any) are zero
In many applications, it is more convenient to work with the \textit{reduced singular value decomposition (reduced SVD)}. 
For $\bA$ of rank $R$ with (full) SVD $\bA = \bU\bSigma \bV^\top$, we take the submatrices $\bU_R \in \real^{M\times R}$, $\bV_R \in \real^{N\times R}$ such that $\bU = [\bU_R, \bu_{R+1}, \ldots, \bu_M]$, $\bV = [\bV_R, \bv_{R+1}, \ldots,\bv_N]$, and $\bSigma_R = \diag([\sigma_1,\ldots,\sigma_R]) \in \real^{R\times R}$.
Therefore, we obtain the reduced SVD of $\bA$:
$$
\bA=\bU_R\bSigma_R\bV_R^\top.
$$
The comparison between the reduced and full SVD is shown in Figure~\ref{fig:svd-comparison}, where white entries denote zeros and blue entries are possibly nonzero.

Given $\bA\in\real^{M\times N}$ with reduced SVD $\bA = \bU_R\bSigma_R\bV_R^\top$, we observe that
\begin{align*}
\bA^\top\bA &= \bV_R\bSigma_R^\top\bU_R^\top\bU_R\bSigma_R\bV_R^\top = \bV_R\bSigma_R^2\bV_R^\top;\\
\bA\bA^\top &= \bU_R\bSigma_R\bV_R^\top\bV_R\bSigma_R^\top\bU_R^\top = \bU_R\bSigma_R^2\bU_R^\top.
\end{align*}
Thus, we recover the (reduced) spectral decompositions of symmetric positive semidefinite matrices $\bA^\top\bA$ and $\bA\bA^\top$, respectively. 
In particular, the singular values $\sigma_i = \sigma_i(\bA)$ satisfy
\begin{equation}\label{equation:sigbd_nearortho_eqp}
\sigma_i(\bA) = \sqrt{\lambda_i(\bA^\top\bA)} = \sqrt{\lambda_i(\bA\bA^\top)}, \quad i = 1,\ldots,\min\{M,N\}, 
\end{equation}
where $\lambda_1(\bA^\top\bA) \geq \lambda_2(\bA^\top\bA) \geq \ldots$ are the eigenvalues of $\bA^\top\bA$ in nonincreasing order. 
Moreover, the left and right singular vectors listed in $\bU,\bV$ can be obtained by the spectral decomposition of the positive semidefinite matrices $\bA^\top\bA$ and $\bA\bA^\top$.
Consequently, one can construct the SVD of $\bA$ directly from the spectral decompositions of $\bA^\top\bA$ and $\bA\bA^\top$---which also provides an alternative proof of the existence of the SVD.


\begin{exercise}[Proof of SVD]
Using the discussion above and the spectral decomposition, prove the existence of the singular value decomposition of a matrix $\bA\in\real^{M\times N}$.
\end{exercise}

\subsection*{Four Orthonormal Bases in SVD}

\index{Fundamental theorem of linear algebra}
For any matrix $\bA\in\real^{M\times N}$, the following properties hold:
\begin{itemize}
\item The null space $\nspace(\bA)$ is the orthogonal complement of the row space $\cspace(\bA^\top)$ in $\real^N$: $\dim(\nspace(\bA))+\dim(\cspace(\bA^\top))=N$.

\item The left null space $\nspace(\bA^\top)$ is the orthogonal complement of the column space $\cspace(\bA)$ in $\real^M$: $\dim(\nspace(\bA^\top))+\dim(\cspace(\bA))=M$.
\end{itemize}
This result is known as  the fundamental theorem of linear algebra  (see Theorem~\ref{theorem:fundamental-linear-algebra}).
In particular, the construction of the SVD provides orthonormal bases for all four fundamental subspaces described in this theorem. To establish this, we first require the following lemma.
\begin{lemma}[Subspace of $\bA^\top \bA$ and $\bA\bA^\top$]\label{lemma:rank-of-ttt}
Let $\bA\in \real^{M\times N}$ be given. 
Then,
\begin{itemize}
\item The column space of $\bA^\top \bA$ is identical  to the column space of $\bA^\top$ (i.e., row space of $\bA$): $\cspace(\bA^\top\bA)=\cspace(\bA^\top)$; this also shows $\nspace(\bA^\top\bA)=\nspace(\bA)$ by fundamental theorem of linear algebra.
\item The column space of $\bA\bA^\top$ is identical  to the column space of $\bA$: $\cspace(\bA\bA^\top)=\cspace(\bA)$. 
Hence, this also shows $\nspace(\bA\bA^\top)=\nspace(\bA^\top)$.
\end{itemize}
\end{lemma}
\begin{proof}[of Lemma~\ref{lemma:rank-of-ttt}]
Let $\bx\in \nspace(\bA)$, we have 
$
\bA\bx  = \bzero \implies \bA^\top\bA \bx =\bzero, 
$
i.e., $\bx\in \nspace(\bA) \implies \bx \in \nspace(\bA^\top \bA)$. Therefore, $\nspace(\bA) \subseteq \nspace(\bA^\top\bA)$. 
Furthermore, let $\bx \in \nspace(\bA^\top\bA)$, we have 
$$
\bA^\top \bA\bx = \bzero\implies \bx^\top \bA^\top \bA\bx = 0\implies \normtwo{\bA\bx}^2 = 0 \implies \bA\bx=\bzero, 
$$
i.e., $\bx\in \nspace(\bA^\top \bA) \implies \bx\in \nspace(\bA)$. Therefore, $\nspace(\bA^\top\bA) \subseteq\nspace(\bA) $. 
Combining both inclusions yields
$
\nspace(\bA) = \nspace(\bA^\top\bA).
$
By the fundamental theorem of linear algebra, it follows that
$$
\cspace(\bA^\top)=\cspace(\bA^\top\bA).
$$
Applying the same argument to $\bA^\top$  establishes the second claim.
\end{proof}

\index{Orthonormal basis}
\index{Fundamental theorem}
\begin{theoremHigh}[Four orthonormal bases in SVD]\label{theorem:svd-four-orthonormal-Basis}
Let $\bA = \bU \bSigma \bV^\top$ be the full SVD of $\bA\in\real^{M\times N}$, where $\bU=[\bu_1, \bu_2, \ldots,\bu_M]$ and $\bV=[\bv_1, \bv_2, \ldots, \bv_N]$ are the column partitions of $\bU$ and $\bV$, respectively. 
Then:
\begin{itemize}
\item $\{\bv_1, \bv_2, \ldots, \bv_R\} $ is an orthonormal basis for the row space  $\cspace(\bA^\top)$;

\item $\{\bv_{R+1},\bv_{R+2}, \ldots, \bv_N\}$ is an orthonormal basis for the null space $\nspace(\bA)$;

\item $\{\bu_1,\bu_2, \ldots,\bu_R\}$ is an orthonormal basis for the column space $\cspace(\bA)$;

\item $\{\bu_{R+1}, \bu_{R+2},\ldots,\bu_M\}$ is an orthonormal basis for the left null space $\nspace(\bA^\top)$. 
\end{itemize}
\end{theoremHigh}
\begin{proof}[of Theorem~\ref{theorem:svd-four-orthonormal-Basis}]
By the spectral decomposition, for the symmetric matrix $\bA^\top\bA$, its column space $\cspace(\bA^\top\bA)$ is spanned by the eigenvectors. 
Therefore, the set $\{\bv_1,\bv_2, \ldots, \bv_R\}$ forms an orthonormal basis for $\cspace(\bA^\top\bA)$.
Thus, $\{\bv_1, \bv_2,\ldots, \bv_R\}$ also serves as an orthonormal basis for $\cspace(\bA^\top)$ by Lemma~\ref{lemma:rank-of-ttt}. 

Since $\bV$ is orthogonal, the space spanned by the remaining vectors $\{\bv_{R+1}, \bv_{R+2},\ldots, \bv_N\}$ is the orthogonal complement to the space spanned by $\{\bv_1,\bv_2, \ldots, \bv_R\}$, which is precisely $\nspace(\bA)$.
Thus, $\{\bv_{R+1},\bv_{R+2}, \ldots, \bv_N\}$ constitutes an orthonormal basis for $\nspace(\bA)$. 

A similar argument applied to $\bA \bA^\top$ shows that $\{\bu_1,\bu_2, \ldots, \bu_R\}$ spans $\cspace(\bA)$ and the set  $\{\bu_{R+1},\bu_{R+2}, \ldots, \bu_M\}$ spans $\nspace(\bA^\top)$.
Alternatively, we can see that $\{\bu_1,\bu_2, \ldots,\bu_R\}$ forms a basis for the column space of $\bA$ by Lemma~\ref{lemma:column-basis-from-row-basis}, since $\bu_i = \frac{\bA\bv_i}{\sigma_i},\, \forall\, i \in\{1, 2, \ldots, R\}$ (which is a key property of SVD). 
\end{proof}

\begin{figure}[h!]
\centering
\includegraphics[width=0.8\textwidth]{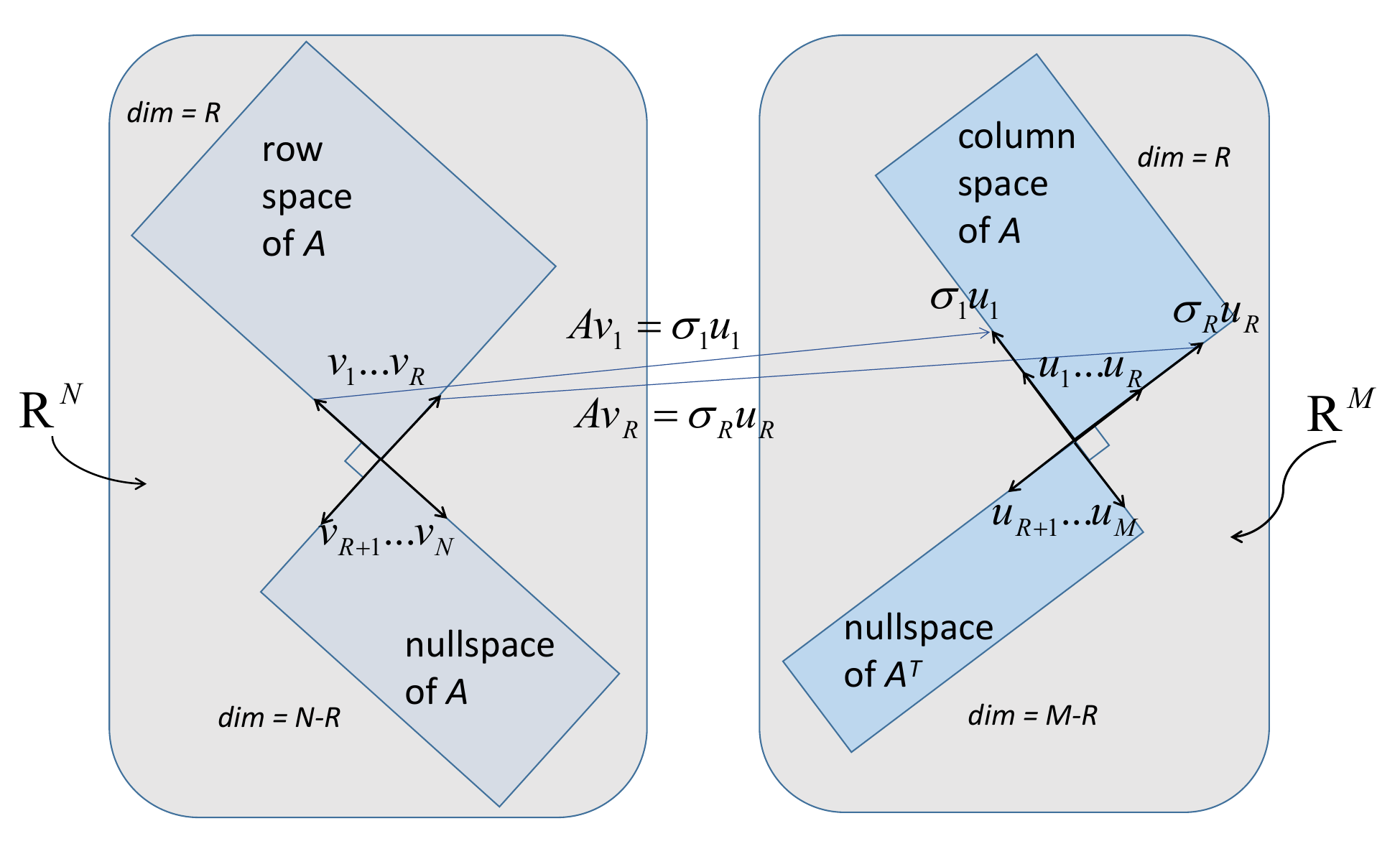}
\caption{Orthonormal bases that diagonalize $\bA$ via the SVD. The set $\{\bv_1, \bv_2, \ldots, \bv_R\} $ forms an orthonormal basis for the row space  $\cspace(\bA^\top)$, and $\{\bu_1,\bu_2, \ldots,\bu_R\}$ forms an orthonormal basis for the column space $\cspace(\bA)$.
The action of $\bA$ links these bases: for each $i \in \{1, 2, \ldots, R\}$, it transforms the row-space basis vector $\bv_i$ into the column-space basis vector $\bu_i$ scaled by the singular value $\sigma_i$, i.e., $\bA \bv_i = \sigma_i \bu_i$.
}
\label{fig:lafundamental3-SVD}
\end{figure}

The relationship among the four subspaces is illustrated  in Figure~\ref{fig:lafundamental3-SVD}.
Specifically, for each $i \in \{1, 2, \ldots, R\}$, the matrix $\bA$ maps the row-space basis vector $\bv_i$ to the column-space basis vector $\bu_i$ according to the relation $\sigma_i\bu_i=\bA\bv_i$.

\subsection*{Differentiability and Differential Calculus}

\index{Continuously differentiability}
\index{Second-order partial derivative}
\begin{definition}[Directional Derivative, Partial Derivative]\label{definition:partial_deri}
Let $f$ be a real-valued function defined on a set$\sS\subseteq \real^N$, and let  $\bd\in\real^N$ be a nonzero vector. Then the \textit{directional derivative} of $f$ at $\bx$ w.r.t. the direction $\bd$ is given by, if the limit exists, 
$$
\mathop{\lim}_{t\rightarrow 0^+}
\frac{f(\bx+t\bd) - f(\bx)}{t}.
$$
And it is denoted by $f^\prime(\bx; \bd)$ or $D_{\bd}f(\bx)$. 
The directional derivative is sometimes called the \textit{G\^ateaux derivative}.

For any $n\in\{1,2,\ldots,N\}$, the directional derivative at $\bx$ w.r.t. the direction of the $n$-th standard basis $\be_n$ is called the $n$-th \textit{partial derivative} and is denoted by $\frac{\partial f}{\partial x_n} (\bx)$, $D_{\be_n}f(\bx)$, or $\partial_n f(\bx)$.
\end{definition}

If all  partial derivatives of a function $f$ exist at a point $\bx\in\real^N$, then the \textit{gradient} of $f$ at $\bx$, denoted $\nabla f(\bx)$, is defined as the column vector containing all the partial derivatives:
$$
\nabla f(\bx)
\triangleq 
\begin{bmatrix}
\frac{\partial f}{\partial x_1} (\bx)\\
\frac{\partial f}{\partial x_2} (\bx)\\
\vdots \\
\frac{\partial f}{\partial x_N} (\bx)	
\end{bmatrix}
\in \real^N.
$$
A function $f$ defined over an open set $\sS\subseteq \real^N$ is called \textit{continuously differentiable} over $\sS$ if all the partial derivatives exist and are continuous on $\sS$.
Under this assumption of continuous differentiability, the directional derivative and the gradient are related by
\begin{equation}
f^\prime(\bx; \bd) = \nabla f(\bx)^\top \bd, \gap \text{for all }\bx\in\sS \text{ and }\bd\in\real^N.
\end{equation} 
Moreover, continuous differentiability implies that $f$ is well-approximated by its linearization:
\begin{equation}
\mathop{\lim}_{\bd\rightarrow \bzero}
\frac{f(\bx+\bd) - f(\bx) - \nabla f(\bx)^\top \bd}{\norm{\bd}} = 0\gap 
\text{for all }\bx\in\sS,
\end{equation}
or equivalently,
\begin{equation}
f(\by) = f(\bx)+\nabla f(\bx)^\top (\by-\bx) + o(\norm{\by-\bx}),
\end{equation}
where $o(\cdot): \real_+\rightarrow \real$ is a one-dimensional function satisfying $\frac{o(t)}{t}\rightarrow 0$ as $t\rightarrow 0^+$.

The partial derivative $\frac{\partial f}{\partial x_i} (\bx)$ is also a real-valued function of $\bx\in\sS$ that can be partially differentiated. The $j$-th partial derivative of $\frac{\partial f}{\partial x_i} (\bx)$ is defined as 
$$
\frac{\partial^2 f}{\partial x_j\partial x_i} (\bx)=
\frac{\partial \left(\frac{\partial f}{\partial x_i} (\bx)\right)}{\partial x_j} (\bx).
$$
This is called the ($j,i$)-th \textit{second-order partial derivative} of function $f$.
A function $f$ defined over an open set $\sS\subseteq$ is called \textit{twice continuously differentiable} over $\sS$ if all the second-order partial derivatives exist and are continuous over $\sS$. In the setting of twice continuously differentiability, the second-order partial derivative are symmetric:
$$
\frac{\partial^2 f}{\partial x_j\partial x_i} (\bx)=
\frac{\partial^2 f}{\partial x_i\partial x_j} (\bx).
$$
The \textit{Hessian} of the function $f$ at a point $\bx\in\sS$, denoted $\nabla^2f(\bx)$, is defined as the symmetric $N\times N$ matrix 
$$
\nabla^2f(\bx)
\triangleq 
\begin{bmatrix}
\frac{\partial^2 f}{\partial x_1^2} (\bx) & 
\frac{\partial^2 f}{\partial x_1\partial x_2} (\bx) & \ldots &
\frac{\partial^2 f}{\partial x_1\partial x_N} (\bx)\\
\frac{\partial^2 f}{\partial x_2\partial x_1} (\bx) & 
\frac{\partial^2 f}{\partial x_2\partial x_2} (\bx) & \ldots &
\frac{\partial^2 f}{\partial x_2\partial x_N} (\bx)\\
\vdots & 
\vdots & \ddots &
\vdots\\
\frac{\partial^2 f}{\partial x_N\partial x_1} (\bx) & 
\frac{\partial^2 f}{\partial x_N\partial x_2} (\bx) & \ldots &
\frac{\partial^2 f}{\partial x_N^2} (\bx)
\end{bmatrix}.
$$

We now present a classical result from calculus.
\begin{theoremHigh}[Taylor's Expansion with Lagrange Remainder] \index{Taylor's formula} \index{Taylor's expansion}
Let $f(x): \real\rightarrow \real$ be $k$-times continuously differentiable on the closed interval $\sI$ with endpoints $x$ and $y$, for some $k\geq 0$. If $f^{(k+1)}$ exists on the interval $\sI$, then there exists a $x^\star \in (x,y)$ such that 
$$
\begin{aligned}
&\gap f(x) \\
&= f(y) + f^\prime(y)(x-y) + \frac{f^{\prime\prime}(y)}{2!}(x-y)^2+ \ldots + \frac{f^{(k)}(y)}{k!}(x-y)^k
+ \frac{f^{(k+1)}(x^\star)}{(k+1)!}(x-y)^{k+1}\\
&=\sum_{i=0}^{k} \frac{f^{(i)}(y)}{i!} (x-y)^i + \frac{f^{(k+1)}(x^\star)}{(k+1)!}(x-y)^{k+1}.
\end{aligned}
$$ 
The Taylor's expansion can be extended to a function of vector $f(\bx):\real^N\rightarrow \real$ or a function of matrix $f(\bX): \real^{M\times N}\rightarrow \real$.
\end{theoremHigh}
Taylor's expansion---also called the \textit{Taylor series}---approximates a function near a point using a polynomial. 
For example, from elementary calculus, we know that near  $\theta=0$,   
$$
\cos (\theta) \approx 1-\frac{\theta^2}{2}.
$$
This is a second-degree polynomial approximation. 
To see how such an approximation arises, suppose we seek a quadratic polynomial $ f(\theta) = c_1+c_2 \theta+ c_3 \theta^2$ that matches $\cos(\theta)$ and its first two derivatives at $\theta=0$. 
Imposing the conditions
$$\left\{
\begin{aligned}
\cos(0) &= f(0); \\
\cos^\prime(0) &= f^\prime(0);\\
\cos^{\prime\prime}(0) &= f^{\prime\prime}(0);\\
\end{aligned}
\right.
\quad\implies\quad  
\left\{
\begin{aligned}
1 &= c_1; \\
-\sin(0) &=0= c_2;\\
-\cos(0) &=-1= 2c_3.\\
\end{aligned}
\right.
$$
This makes $f(\theta) = c_1+c_2 \theta+ c_3 \theta^2 = 1-\frac{\theta^2}{2}$, which agrees with the known second-order Taylor approximation of $\cos(\theta)$ at $0$. We omit the full proof of the general Taylor theorem here.
For multivariate functions, analogous approximations hold.

\begin{theoremHigh}[Linear Approximation Theorem]\label{theorem:linear_approx}
Let $f(\bx):\sS\rightarrow \real$ be a twice continuously differentiable function on an open set $\sS\subseteq\real^N$, and given two points $\bx, \by\in\sS$. Then there exists a point $\bx^\star\in[\bx,\by]$ such that 
$$
f(\by) = f(\bx)+ \nabla f(\bx)^\top (\by-\bx) + \frac{1}{2} (\by-\bx)^\top \nabla^2 f(\bx^\star) (\by-\bx).
$$ 
\end{theoremHigh}

\begin{theoremHigh}[Quadratic Approximation Theorem]\label{theorem:quadratic_app_theo}
Let $f(\bx):\sS\rightarrow \real$ be a twice continuously differentiable function on an open set $\sS\subseteq\real^N$, and given two points $\bx, \by\in\sS$. Then it follows that 
$$
f(\by) = f(\bx)+ \nabla f(\bx)^\top (\by-\bx) + \frac{1}{2} (\by-\bx)^\top \nabla^2 f(\bx) (\by-\bx)
+
o(\norm{\by-\bx}^2).
$$ 
\end{theoremHigh}

\begin{problemset}
\item \label{problem:psd_eigen} \textbf{Eigenvalue characterization theorem.} 
Prove that a symmetric matrix $\bA$ is positive definite if and only if all its eigenvalues are strictly positive. Similarly, $\bA$ is positive semidefinite if and only if all its eigenvalues are nonnegative.

\item \label{prob:tr_de_pd} \textbf{Trace, det of PD/PSD/ND matrices.} Let $\bA$ be positive definite (resp., positive semidefinite). Show that $\trace(\bA)$ and $\det(\bA)$  are all positive (resp., nonnegative). Moreover, show that $\trace(\bA)=0$ if and only if $\bA=\bzero$. 
Let $\bB\in\real^{N\times N}$ be negative definite. Show that $\trace(\bB)$ is negative; $\det(\bB)$ is negative for odd $N$ and is positive for even $N$.
\textit{Hint: Use Problem~\ref{problem:psd_eigen}.}

\item \textbf{Submultiplicativity.}
Prove that both the Frobenius norm and the spectral norm are submultiplicative; that is, for any matrices
$\bA,\bB\in\real^{N\times N}$,
$$
\normf{\bA\bB}\leq \normf{\bA}\normf{\bB}
\qquad\text{and}\qquad 
\normtwo{\bA\bB}\leq \normtwo{\bA}\normtwo{\bB}.
$$
\item \label{prob:cs_matvec}  \textbf{Cauchy--Schwarz inequality.}
For any vectors $\bu,\bv\in\real^N$, show that 
$$
\abs{\bu^\top\bv}
\leq \normtwo{\bu}\normtwo{\bv}
$$

\item \label{theorem:holder-inequality}\textbf{\holders  inequality\index{\holders  inequality}.}
Let $p,q>1$ satisfy $\frac{1}{p}+\frac{1}{q} = 1$. 
Show that  for any vector $\bx,\by\in \real^N$, we have
$$
\sum_{n=1}^{N}x_n y_n
\leq 
\abs{\sum_{n=1}^{N}x_n y_n}
\leq \sum_{n=1}^{N}\abs{x_n} \abs{y_n} \leq \left(\sum_{n=1}^{N}  \abs{x_n}^p\right)^{1/p}  \left(\sum_{n=1}^{N} \abs{y_n}^q\right)^{1/q}=\normp{\bx}\norm{\by}_q,
$$
where $\normp{\bx} = \left(\sum_{n=1}^{N}  \abs{x_n}^p\right)^{1/p} $ denotes the \textit{$\ell_p$-norm}  of $\bx$. 
Show that the equality holds if the two sequences $\{\abs{x_n}^p\}$ and $\{\abs{y_n}^q\}$ are linearly dependent.
When $p=q=2$, this reduces to the vector Cauchy--Schwarz inequality.

\item \label{prob:cauch_sc_l1l2}\textbf{Standard bounds on vector norms}\index{Standard bounds on vector norms}
Let $\bx\in\real^N$. Use the Cauchy--Schwarz inequality (see Problem~\ref{prob:cs_matvec}) to prove the following relationships:
$$
\begin{aligned}
\norminf{\bx} &\leq \normone{\bx} \leq N\norminf{\bx}; \\
\norminf{\bx} &\leq \normtwo{\bx} \leq \sqrt{N}\norminf{\bx}; \\
\normtwo{\bx} &\leq \normone{\bx} \leq \sqrt{N}\normtwo{\bx}. \\
\end{aligned}
$$

\item Prove Remark~\ref{remark:equiva_nonsingular} and Remark~\ref{remark:equiva_singular}.
\end{problemset}

\chapter{Bayesian Inference}\label{chapter:mcmcs}
\begingroup
\hypersetup{
	linkcolor=structurecolor,
	linktoc=page,  
}
\minitoc \newpage
\endgroup
\index{Markov chain Monte Carlo}

\lettrine{\color{caligraphcolor}I}
In statistics, the complete set of data that needs to be investigated or studied for a certain phenomenon or entity is referred to as the \textit{statistical population}, or simply the \textit{population}. For example,  if we wish to study the fuel efficiency (measured in miles per gallon, or MPG) of a specific model of car, then the fuel efficiency data of all  cars of that model ever produced constitute the population. The distribution of these fuel efficiency data points is called the \textit{population distribution}. Each car's fuel efficiency, or each single data point, is referred to as an \textit{individual}.

However, in practice, it is usually impossible to collect data from every single car of that model; thus, the full population is typically unknown. 
In such cases, we resort to \textit{sampling}: randomly selecting a subset of individuals (cars) from the population and measuring their fuel efficiency
The resulting collection of measurements is called a \textit{sample}. For instance, if 500 cars are randomly selected from the entire production run, and their fuel efficiencies are recorded, we obtain a sample consisting of 500 data points. 
The number of data points in the sample is called the \textit{sample size}, which in this case is 500. It is important to distinguish between a \textit{sample} and an \textit{individual}: a sample results from one sampling process and contains multiple individual data points, and its size is the count of those individuals.

We often assume that the population distribution follows a known probability distribution, though some of its parameters remain unknown. For example, fuel efficiency might follow a normal distribution, but its mean and variance are not known in advance. In such situations, our goal is to infer (or estimate) these unknown parameters---such as the mean and variance---using the observed sample.

\textit{Statistical inference} refers to the set of methods used to deduce characteristics (typically parameters) of a population based on sample data. For example, we might use the sample average to estimate the population mean. More broadly, statistical inference involves making probabilistic statements about unknown quantitative features of a population based on a limited set of observations. There are many established techniques for parameter estimation from samples, including the method of moments, maximum likelihood estimation (MLE), and Bayesian estimation.

This book focuses specifically on Markov chain Monte Carlo (MCMC) methods for probabilistic and statistical inference, which aim to draw conclusions from probabilistic models.
This chapter provides a concise overview of the mathematical foundations of probabilistic inference, emphasizing concepts that will serve as the basis for the material in subsequent chapters. Our hope is that it offers a solid stepping stone toward understanding the Bayesian inference algorithms presented later and used throughout the rest of the book.

\section{The Bayesian Approach}\label{section:bayes_approach}
Over the past decade, the Bayesian approach has been widely applied across diverse areas of data analysis, including economic forecasting, medical imaging, and population studies \citep{besag1986statistical, hill1994bayesian, marseille1996bayesian}. In modern statistics, Bayesian methods have become increasingly important and prevalent.
The foundational idea is attributed to \textit{Thomas Bayes}, who conceived it but died before publishing his work. Fortunately, his friend \textit{Richard Price} edited and published Bayes' findings in 1764. The same principle was later independently rediscovered by \textit{Pierre-Simon Laplace} at the end of the 18-th century.
In this section, we introduce the core ideas of the Bayesian approach and illustrate them using two simple models---the Beta-Bernoulli model and the Bayesian linear model---as an appetizer to highlight the advantages and role of prior information in Bayesian modeling.

Bayesian modeling and statistics are fundamentally grounded in Bayes' theorem, which is formally stated as follows:
\begin{theoremHigh}[Bayes' Theorem \citep{bayes1958essay, laplace1820theorie}\index{Bayes' theorem}]\label{theorem:bayes_theo}
Let $\sS$ be a sample space and let $\sB_1, \sB_2, \ldots, \sB_K$ be a partition of $\sS$ such that (1). $\cup_k \sB_k=\sS$ and (2). $\sB_i \cap \sB_j=\varnothing$ for all $i\neq j$. 
Let further $\sA$ be any event. Then it follows that 
$$
\prob(\sB_k \mid \sA) = \frac{\prob(\sA \mid \sB_k)\prob(\sB_k)}{\prob(\sA)} = \frac{\prob(\sA\mid \sB_k)\prob(\sB_k)}{\sum_{i=1}^{K}\prob(\sA\mid \sB_k)\prob(\sB_k)}.
$$
\end{theoremHigh}
In Bayesian statistics, Bayes' theorem provides a principled rule for updating probabilities when new information---such as observed data---becomes available. This allows us to refine our prior beliefs about parameters of interest.

More concretely, let $\mathcalX=\mathcalX (\bx_{1:N})= \{\bx_1, \bx_2, \ldots, \bx_N\}$ denote a set of $N$ observed data points, assumed to be independent and identically distributed (i.i.d.) according to a distribution parameterized by $\btheta\in\bTheta$. 
Note that the parameters $\btheta$ may include  hidden or latent variables, such as   latent variables in a mixture model indicating the cluster to which a data point belongs. 
One common approach  to learning the model parameters involves finding the best-fit parameters $\widehat{\btheta}_{\text{MLE}}$ that maximize the likelihood (hence the name \textit{maximum likelihood estimation}, MLE or ML estimation~\footnote{\textbf{Estimation method vs. estimator vs. estimate}: \textbf{Estimation method} is a general algorithm to produce the estimator. \textbf{An estimate} is the specific value that \textbf{an estimator} takes when observing the specific value, i.e., an estimator is a random variable and the realization of this random variable is called an estimate.}):
$$
\widehat{\btheta}_{\text{MLE}} = \mathop{\argmax}_{\btheta\in\bTheta} p(\mathcalX \mid \btheta).
$$

In contrast, the Bayesian approach treats parameters as random variables, reflecting our uncertainty about their true values. This aligns with the broader Bayesian philosophy: all uncertain quantities are modeled as random variables, and probability theory is used to reason about them. Rather than seeking a single best-fitting parameter (as in MLE), Bayesian inference accounts for all plausible values of $\btheta$ through integration.

Naturally, we want $p(\mathcalX \mid \btheta)$ to be flexible enough to adapt to the data, giving us the opportunity to develop a sufficiently accurate model. At the same time, we aim to integrate any prior knowledge about the data distribution into the model.
The idea of the Bayesian approach involves assuming a \textit{prior} probability distribution for $\btheta$ with hyper-parameters $\balpha$ (i.e., $p(\btheta\mid \balpha)$, also known as the probability of the model). 
This  distribution represents the plausibility of each possible value of $\btheta$ before observing the data, and it captures our prior uncertainty regarding $\btheta$. 
The joint distribution of $\btheta$ and $\mathcalX$ is given by
$$
p(\btheta, \mathcalX) = p(\btheta \mid \balpha) p(\mathcalX \mid \btheta).
$$
And we can integrate out $\btheta$ to obtain the \textit{marginal distribution (marginal likelihood)} of $\mathcalX$, 
$$
p(\mathcalX \mid \balpha)= \int p(\btheta \mid \balpha)p(\mathcalX \mid \btheta) \,d\btheta.~\footnote{We assume $\mathcalX$ and $\balpha$ are conditionally independent given $\btheta$. Otherwise, the marginal likelihood can be represented by $p(\mathcalX \mid \balpha)= \int p(\btheta \mid \balpha)p(\mathcalX \mid \btheta, \balpha) \,d\btheta$.}
$$
In the machine learning community, this particular measure is occasionally termed the ``\textit{evidence}" for the model under hyper-parameters $\balpha$, since it represents the element of the posterior distribution across models that is influenced by the data. 

Since probabilistic models include elements that are unknown and the available data seldom provides a comprehensive view of these unknowns, we usually have to incorporate a certain degree of uncertainty regarding various aspects of the model. This uncertainty is defined through (conditional) probability distributions, which characterize both the extent and the type of uncertainty involved.
In the model described above, then, to make inferences about $\btheta$, one simply considers the conditional distribution of $\btheta$ given the observed data. This is referred to as the \textit{posterior} distribution, since it represents the plausibility of each possible value of $\btheta$ after seeing the data.
The posterior distribution is the solution space for given problems and  allows us to quantify our \textit{uncertainty} about parameter values after observing the data, since it measures the probability of the present model in light of the data.
Mathematically, this relationship is expressed via Bayes' theorem,
\begin{equation}\label{equation:posterior_abstract_for_mcmc}
\begin{aligned}
p(\btheta \mid  \mathcalX, \balpha) 
&= \frac{p(\mathcalX \mid  \btheta ) p(\btheta \mid  \balpha)}{p(\mathcalX \mid \balpha)} \\
&= \frac{p(\mathcalX \mid  \btheta ) p(\btheta \mid  \balpha)}{\int_{\btheta}  p(\mathcalX, \btheta \mid \balpha) }  
= \frac{p(\mathcalX \mid \btheta ) p(\btheta \mid \balpha)}{\int_{\btheta}  p(\mathcalX \mid  \btheta ) p(\btheta \mid  \balpha) }  \propto p(\mathcalX \mid \btheta ) p(\btheta \mid \balpha),
\end{aligned}
\end{equation}
where ``$\propto$" means ``proportional to" (see Problem~\ref{problem:proptionto}), $\mathcalX$ is the observed data set, and the marginal distribution $p(\mathcalX \mid \balpha )$ can be disregarded in this case since it acts as a scaling parameter (and we shall see the MCMC algorithm only needs relative probabilities; as a matter of fact, the marginal likelihood is usually impossible to compute). 
In other words, we say the posterior is proportional to the product of the likelihood and the prior. 
This means that the relative probability at a point in the solution space is determined completely by the likelihood,
which is easily determined by comparing the model to the data, and the prior, which is the probability of the model independent of the data. 
The prior encodes any a priori knowledge of the solution irrespective  of the observed data. 
For example, a prior for a system reducing over-clustering might assign a higher probability to a larger cluster than to a small cluster \citep{lu2021survey}.

The elegance of Bayes' theorem becomes apparent as it distinguishes inference from modeling. The model, encompassing the prior distribution and the likelihood, fully dictates the posterior distribution, leaving the computation of the inference as the only remaining step.
More generally, the Bayesian approach---in a nutshell---is to assume a prior distribution for any unknowns ($\btheta$ in our case), and then just follow the rules of probability to answer any questions of interest. For example, when we find the parameter based on the maximum posterior probability of $\btheta$, we turn to the \textit{maximum a posteriori (MAP)} estimation.

\begin{definition}[Maximum a Posterior Estimator]
Maximum a posteriori (MAP) estimate is the parameter value that maximizes the posterior distribution.
The MAP estimate balances information from the prior distribution with information from the likelihood. 
The influence of the prior is stronger when the likelihood provides less information, and vice versa.
\end{definition}

Other than the MAP estimator, this posterior distribution alllows us to compute the density at a new coming data point $\bx^\prime$, called the \textit{posterior predictive distribution}, by averaging over both the uncertainty in the model and in the parameters:
$$
p(\bx^\prime \mid \mathcalX) = \int p(\bx^\prime \mid \btheta) p(\btheta \mid \mathcalX, \balpha) \,d\btheta.~\footnote{If the problem  follows from a generative process $\by\sim p(\by\mid \bx,\btheta)$. Then the predictive distribution is $p(\by^\prime\mid\bx^\prime,\mathcalX, \mathcalY)=\int p(\by^\prime \mid \bx^\prime, \btheta)p(\btheta \mid \mathcalX,\mathcalY, \balpha)\,d\btheta$.}
$$
The posterior predictive distribution can be employed to design test statistics of interest and then compare the posterior predictive distributions to the test statistics of observed values so as to determine the best model among several candidates. This process is known as \textit{model checking or selection} \citep{haugh2021tutorial}.

\index{Model checking}
\index{Model selection}

\paragrapharrow{Frequentists V.S. Bayesian.}
The \textit{frequentist approach} to statistics, developed by Neyman, evaluates statistical procedures based on a probability distribution over all possible datasets.
To be more specific, frequentists consider the parameter vector $\btheta$ to be fixed (albeit unknown), while introducing uncertainty over possible datasets $\mathcalX$. Frequentist methods are often considered more objective as they avoid incorporating subjective prior information. 
In contrast, Bayesian methods allow for the incorporation of prior beliefs. The Bayesian approach treats the data set $\mathcalX$ as given, while introducing uncertainty over $\btheta$. 
However, statisticians nowadays tend to move comfortably between these approaches and popular statistical procedures often combine both of them,  incorporating Bayesian methods for certain aspects of the analysis while using frequentist methods for others. 
For instance, empirical Bayesian methods  have a Bayesian spirit but are not strictly Bayesian; their analysis is frequently frequentist \citep{haugh2021tutorial}.

\subsection{Laplace Approximation}
We mentioned earlier that the posterior distribution can be used to answer any question of interest, including MAP estimation:
$$
\widehat{\btheta}_{\text{MAP}} 
= \mathop{\argmax}_{\btheta\in\bTheta} p(\btheta \mid  \mathcalX, \balpha) 
= \mathop{\argmax}_{\btheta\in\bTheta} p(\mathcalX \mid \btheta ) p(\btheta \mid \balpha).
$$
The motivation for the \textit{Laplace approximation} stems from a fundamental limitation of relying solely on the MAP estimate in Bayesian inference: the MAP provides only a point estimate and discards all information about uncertainty and shape of the posterior distribution:
\begin{itemize}
\item \textbf{Loss of uncertainty quantification:}
The MAP estimate, $\widehat{\btheta}_{\text{MAP}} 
= \mathop{\argmax}_{\btheta\in\bTheta} p(\btheta \mid  \mathcalX, \balpha) $, identifies the single most probable parameter value under the posterior. However, it tells us nothing about:
(i) how confident we are in that estimate, (ii) the spread or variability of plausible parameter values, and (iii) whether there are multiple distinct regions of high posterior density (e.g., multimodality).
In many applications---such as model comparison, prediction with calibrated confidence intervals, or decision-making under uncertainty---this full posterior information is essential.
\item \textbf{Inadequate for marginal likelihood estimation:}
To perform Bayesian model selection or hyper-parameter learning, we need the marginal likelihood $p(\mathcalX\mid \balpha)=\int p(\mathcalX\mid \btheta)p(\btheta\mid \balpha )\,d\btheta$.
The MAP alone cannot provide this integral---it only gives the integrand's maximum, not its volume. Ignoring the width/curvature of the posterior leads to overconfident model comparisons.
\item \textbf{No basis for propagating uncertainty:}
In downstream tasks (e.g., predictive distributions  
$p(\bx'\mid \mathcalX)=\int p(\bx'\mid \btheta)p(\btheta\mid \mathcalX)\,d\btheta$, using only the MAP amounts to plug-in approximation, which underestimates predictive variance and fails to account for parameter uncertainty.
\end{itemize}
The {Laplace approximation} involves approximating the posterior with a Gaussian distribution centered at the mode of the posterior (i.e., the MAP estimate $\widehat{\btheta}_{\text{MAP}} $). 
This provides a practical way to approximate the posterior when its exact form is analytically intractable or computationally expensive to evaluate \citep{kass1995bayes, mackay1998choice, friston2007variational}. 
Define the logarithm of the posterior distribution as 
$$
\mathcalL(\btheta) =\ln p(\mathcalX \mid \btheta ) p(\btheta \mid \balpha)
=\ln p(\mathcalX \mid \btheta ) + \ln p(\btheta \mid \balpha) =\ln p(\btheta \mid  \mathcalX, \balpha) + \mathcalC,
$$
where \texttt{ln} denotes  the natural logarithm (to base $e$),  and $\mathcalC$ is a constant with respect to $\btheta$.
According to the quadratic approximation theorem (Theorem~\ref{theorem:quadratic_app_theo}) and assuming that the parameter space $\bTheta$ is an open set (the gradient of the MAP has vanished gradient),
we can approximate $L(\btheta)$ using a second-order Taylor expansion around $\widehat{\btheta}_{\text{MAP}}$:
$$
\begin{aligned}
\mathcalL(\btheta) &\approx 
\mathcalL(\widehat{\btheta} ) +\nabla \mathcalL(\widehat{\btheta} )^\top (\btheta - \widehat{\btheta} )
+\frac{1}{2} (\btheta - \widehat{\btheta} )^\top \nabla^2 \mathcalL(\widehat{\btheta} ) (\btheta - \widehat{\btheta} )\\
&=\mathcalL(\widehat{\btheta} ) 
+\frac{1}{2} (\btheta - \widehat{\btheta} )^\top \nabla^2 \mathcalL(\widehat{\btheta} ) (\btheta - \widehat{\btheta} ),
\end{aligned}
$$
where we let $\widehat{\btheta}=\widehat{\btheta}_{\text{MAP}} $ for brevity, and we use the fact that since the first-order term is zero at the mode.
Using this quadratic approximation, the log marginal likelihood (also known as the log \textit{evidence}) can be approximated as:
$$
\begin{aligned}
\begin{aligned}
\ln p(\mathcalX \mid \balpha) 
&= \ln \int p(\mathcalX \mid \btheta) p(\btheta \mid \balpha) \,d\btheta 
=\ln  \int \exp\{\mathcalL(\btheta) \}\,d\btheta\\
&\approx \ln p(\mathcalX \mid \widehat{\btheta} ) + \ln p(\widehat{\btheta} \mid \balpha) + \frac{D}{2} \ln (2\pi) - \frac{1}{2}\ln \abs{\nabla^2 \mathcalL(\widehat{\btheta} )},
\end{aligned}
\end{aligned}
$$
where  the last approximation follows  from the definition of the multivariate Gaussian distribution (see, for example, Section~\ref{section:multi_gaussian_dist}), and $D$ is the dimension of the parameter space: $\btheta\in\real^D$.
Exponentiating both sides gives the Laplace approximation to the marginal likelihood:
$$
p(\mathcalX \mid \balpha)_{\text{Lap}} = 
\underbrace{p(\mathcalX \mid \widehat{\btheta} )}_{\text{data likelihood under MAP}} \,
\underbrace{p(\widehat{\btheta} \mid \balpha)}_{\text{penalty from prior}} 
\underbrace{\abs{2\pi (\nabla^2 \mathcalL(\widehat{\btheta} ))^{-1}}}_{\text{local curvature}}.
$$
Thus, the Laplace approximation decomposes into three interpretable components: (i) the likelihood of the data evaluated at the MAP estimate, (ii) a penalty (or regularization) term from the prior, and (iii) a volume correction that accounts for the local curvature of the log-posterior around the mode.

Although the Laplace approximation is computationally efficient and often useful, it has several notable limitations:
\begin{itemize}
\item \textbf{Gaussian assumption}: 
The method assumes the posterior is approximately Gaussian. This assumption may fail badly for multimodal, skewed, or heavy-tailed posteriors, leading to poor uncertainty quantification. Moreover, the Gaussian approximation assigns nonzero probability to invalid parameter values---for example, negative precisions or mixing proportions outside $[0,1]$. While reparameterization (e.g., using log or logit transforms) can help \citep{mackay1998choice}, the approximation is generally not invariant under reparameterization in finite samples---a significant drawback.

\item \textbf{Mode dependence}: The quality of the approximation depends entirely on a single posterior mode. If the posterior is multimodal or the mode is flat or poorly defined, the Laplace approximation can be highly inaccurate.

\item \textbf{Curvature assumption}: 
The method assumes the posterior curvature is well captured by a constant Hessian near the mode. In complex models---especially those with strong nonlinearities or varying curvature---this assumption often fails, resulting in poor global approximation.

\item \textbf{Computation of Hessian}: Computing the Hessian matrix, which is required to determine the variance of the Gaussian approximation, can be computationally expensive and unstable, particularly for models with many parameters or non-smooth likelihood functions. The computation of the volume term, which depends on the determinant of the Hessian matrix ($\absbig{\nabla^2 \mathcalL(\widehat{\btheta} )}$), poses another challenge. Calculating the derivatives within the Hessian requires  $\mathcalO(ND^2)$ operations, followed by $\mathcalO(D^3)$ operations to find the determinant, making it computationally intensive for high-dimensional problems. To simplify this process, approximations often ignore off-diagonal elements or assume a block-diagonal structure for the Hessian, effectively disregarding interdependencies among parameters. 

\item \textbf{Sensitivity to priors}: The approximation can be overly sensitive to the choice of prior, especially when the prior is strongly informative. In such cases, the Gaussian fit may reflect prior assumptions more than the actual data-informed posterior shape.

\item \textbf{Dimensionality issues}: As the number of parameters increases, the Laplace approximation becomes less reliable due to the curse of dimensionality, where the volume of the parameter space grows exponentially and the Gaussian approximation becomes increasingly poor.
\end{itemize}
Despite these shortcomings, the Laplace approximation remains a valuable tool---particularly as a fast initial approximation or when combined with more robust methods such as Markov chain Monte Carlo (MCMC) or variational inference, which are discussed in later sections and can provide more accurate characterizations of complex posterior distributions.

\index{Bayesian information criterion}
\index{Schwarz criterion}
\subsection{Bayesian Information Criterion}
We can express the Laplace approximation of the marginal likelihood together with its computational complexity in terms of the data size $N$:
$$
\ln p(\mathcalX \mid \balpha)_{\text{Lap}}  = 
\underbrace{\ln p(\mathcalX \mid \widehat{\btheta} )}_{\mathcalO(N)} + \underbrace{\ln p(\widehat{\btheta} \mid \balpha)}_{\mathcalO(1)} + 
\underbrace{\frac{D}{2} \ln (2\pi) }_{\mathcalO(1)}
- 
\underbrace{\frac{1}{2}\ln \abs{\nabla^2 \mathcalL(\widehat{\btheta} )}}_{\mathcalO(D\ln N)}.
$$
The \textit{Bayesian information criterion (BIC)}, also known as the  \textit{Schwarz criterion},  retains only those terms that grow with the sample size $N$.
Since the entries of the Hessian scale linearly with $N$ \citep{schwarz1978estimating}, we have:
$$
\begin{aligned}
\ln p(\mathcalX \mid \balpha)_{\text{Lap}} 
&\approx 
\ln p(\mathcalX \mid \widehat{\btheta} ) -\frac{1}{2}\abs{\nabla^2 \mathcalL(\widehat{\btheta} )}
\overset{N\rightarrow \infty}{\approx}
\ln p(\mathcalX \mid \widehat{\btheta} ) - \lim_{N\rightarrow \infty }\frac{1}{2}\abs{\nabla^2 \mathcalL(\widehat{\btheta} )}\\
&= \ln p(\mathcalX \mid \widehat{\btheta} ) -\frac{1}{2}\abs{N \bH_0} = \ln p(\mathcalX \mid \widehat{\btheta} ) - \frac{D}{2}\ln N - \underbrace{\frac{1}{2}\ln \abs{\bH_0}}_{\mathcalO(1)}. 
\end{aligned}
$$
Therefore, the BIC score becomes
$$
\ln p(\mathcalX \mid \balpha)_{\text{BIC}} 
=
\ln p(\mathcalX \mid \widehat{\btheta} ) - \frac{D}{2}\ln N.
$$
The BIC has several appealing properties:
\begin{itemize}
\item It includes a penalty term proportional to the number of parameters $D$, which discourages overfitting by favoring simpler models.
\item It is straightforward to compute and interpret, requiring only the maximized log-likelihood and the model dimension. It does not require any additional assumptions beyond those inherent in the models being compared.
\item It does not require specifying prior distributions beyond what is needed to define the model itself, making it accessible even to practitioners unfamiliar with Bayesian methods.
\item Under standard regularity conditions, BIC is a consistent model selection criterion: as $N\rightarrow \infty$, it selects the true model with probability approaching one, provided the true model is among the candidates.
\end{itemize}
However, from a fully Bayesian perspective, the lack of explicit prior dependence may be viewed as a limitation, as it discards potentially useful prior information.

On the other hand, BIC is invariant to reparameterization---a desirable property. Because it depends only on the maximized likelihood and the parameter count (not on how parameters are expressed and  the local geometry of the parameter space), BIC yields consistent results regardless of the chosen parameterization. This aligns with a core principle of Bayesian inference: the posterior should be invariant under smooth reparameterizations \citep{hoff2009first}. Such invariance enhances the reliability and fairness of model comparisons, avoiding biases introduced by arbitrary choices in model formulation \citep{beal2003variational}.

\index{Occam's razor}
\subsection{Occam's Razor and Occam Factor}
In the BIC framework, model complexity is equated with the number of parameters, and overly complex models are penalized to avoid overfitting.
However, this view can be misleading: a model with many parameters might still be highly constrained and capable of explaining only a narrow range of data, while a model with a single parameter might be flexible enough to fit a wide variety of datasets.
A more principled approach uses the marginal likelihood (or evidence):
$$
p(\mathcalX \mid \balpha)= \int p(\btheta \mid \balpha)p(\mathcalX \mid \btheta) \,d\btheta,
$$
which integrates out the parameters $\btheta$. 
This automatically penalizes models with excessive degrees of freedom, as such models spread their predictive probability mass thinly over many possible datasets.
This built-in trade-off is known as \textit{Bayesian Occam's razor}: given equal fit to the observed data, simpler models are preferred because they concentrate their predictive mass more sharply \citep{mackay1995probable, beal2003variational}.

\begin{figure}[h!]
\centering
\includegraphics[width=0.65\textwidth]{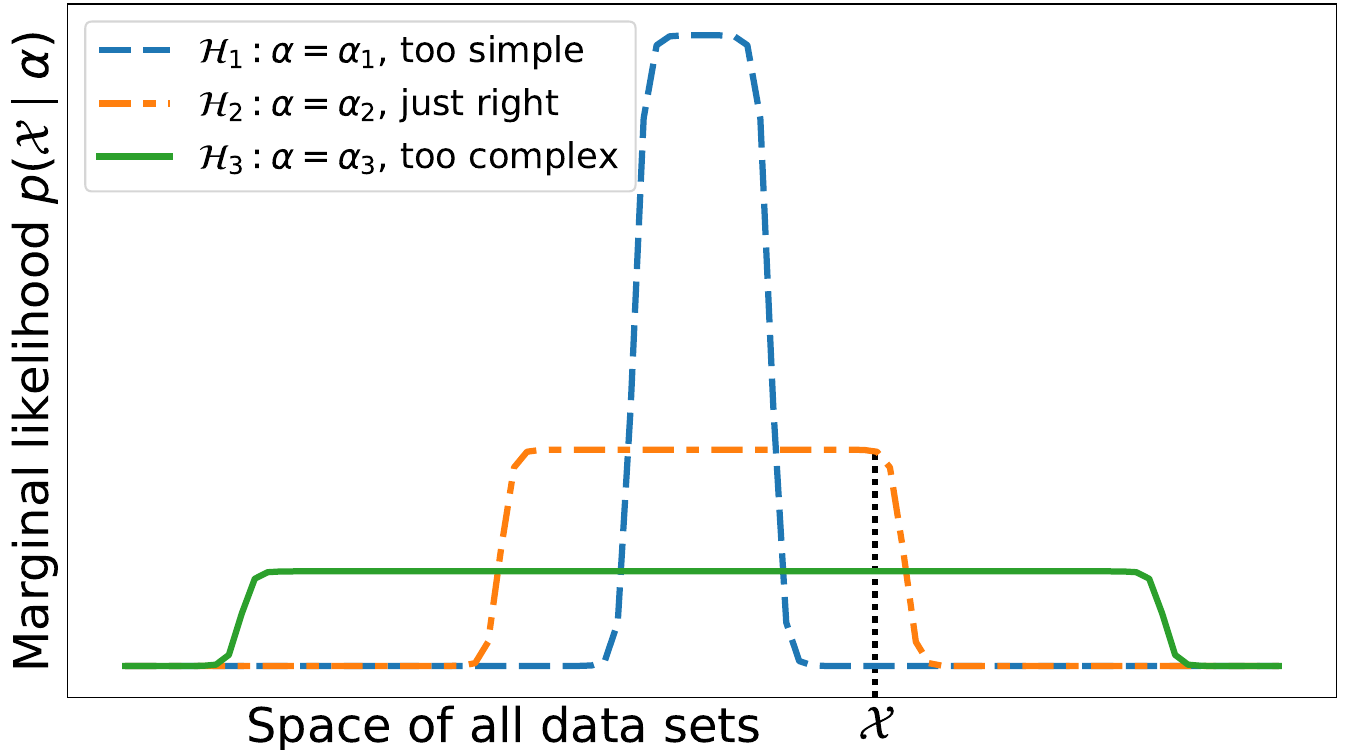}
\caption{Bayesian inference embodies Occam's razor. This figure provides the fundamental intuition for why more complex models tend to be less probable. The horizontal axis represents the space of all possible datasets, $\mathcalX$. 
According to Bayes' theorem, models are favored in proportion to how well they predicted the observed data. These predictions are represented by a marginal probability distribution over $\mathcalX$. 
A simple model  makes only a limited range of predictions; while a more powerful model is capable of predicting a greater variety of datasets.}
\label{fig:occam_razor}
\end{figure}

This principle Occam's razor is illustrated in   Figure~\ref{fig:occam_razor}.
Since the marginal likelihood $p(\mathcalX \mid \balpha)$ defines a probability distribution over all possible datasets $\mathcalX$, its total integral equals one. 
If the the model is overly complex such that it can model a vast variety of datasets, the probability value for each data set can be reduced (the ``too complex" case in the figure with hypothesis \{$\mathcalH_3: \balpha=\balpha_3$\}). 
When the model is too simple, it might not cover the observed data set, rendering a small marginal probability (the ``too simple" case in the figure with hypothesis \{$\mathcalH_1: \balpha=\balpha_1$\}).

In Figure~\ref{fig:occam_razor}, the model hypothesis $\mathcalH_1$ is not compatible with the observed data set $\mathcalX$.
However, in the case where the data are compatible with both theories $\mathcalH_2$ and $\mathcalH_3$, the simpler model $\mathcalH_2$ will turn out to be more probable than the more complex model $\mathcalH_3$, without us having to express any subjective bias against complex models. Our subjective prior should simply assign equal probabilities to the possibilities of simplicity and complexity.
Therefore, given a data set $\mathcalX$, it is possible to discard both models that are too complex and those that are too simple, based on their marginal likelihood.

\begin{figure}[h!]
\centering
\includegraphics[width=0.65\textwidth]{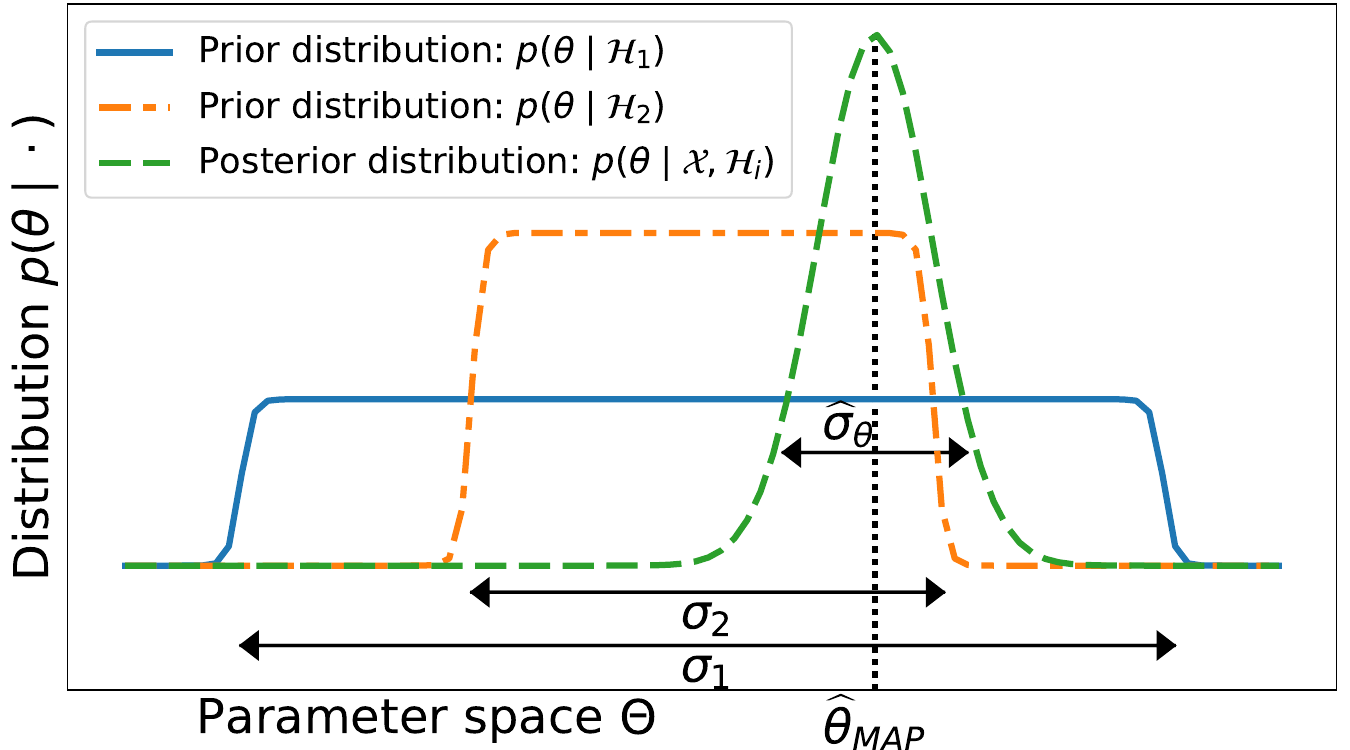}
\caption{Occam factor. The prior distribution $p(\btheta \mid \mathcalH_1)$ for the
	parameter has width $\sigma_{1}$, and the prior distribution $p(\btheta \mid \mathcalH_2)$ for the
	parameter has width $\sigma_{2}$ ($\sigma_2<\sigma_1$). The posterior distribution  has a single peak at $\widehat{\btheta}_{\text{MAP}}$ with width $\widehat{\sigma}_{\btheta}$. }
\label{fig:occam_factor}
\end{figure}

As mentioned previously, the marginal likelihood or evidence is often  intractable or impossible to compute. Bayesian Occam's razor provides a way to approximate the marginal likelihood \citep{mackay1995probable}. As a recap, the marginal likelihood under a hypothesis $\{\mathcalH_1: \balpha=\balpha_1\}$ is 
$$
p(\mathcalX \mid \mathcalH_1)= \int p(\btheta \mid \mathcalH_1)p(\mathcalX \mid \btheta) \,d\btheta.
$$
For many problems, it is not uncommon that the posterior distribution $p(\btheta \mid \mathcalX, \mathcalH_1)= \frac{p(\mathcalX \mid \btheta)p(\btheta \mid \mathcalH_1)}{\text{marginal likelihood}}$ has a strong peak at the most probable parameter $\widehat{\btheta}_{\text{MAP}}$, i.e., the MAP estimate (see Figure~\ref{fig:occam_factor}).
Therefore, the marginal likelihood can be approximated by the height of the peak of the integrand $p(\btheta \mid \mathcalH_1)p(\mathcalX \mid \btheta)$ times its width, denoted by $\widehat{\sigma}_{\btheta}$ (see Figure~\ref{fig:occam_factor}):
$$
\underbrace{p(\mathcalX \mid \mathcalH_1)}_{\text{marginal likelihood}} \approx 
\underbrace{p(\mathcalX \mid \widehat{\btheta}_{\text{MAP}})}_{\text{MAP fit likelihood}} 
\underbrace{p(\widehat{\btheta}_{\text{MAP}} \mid \mathcalH_1) \cdot\widehat{\sigma}_{\btheta}}_{\text{Occam factor}},
$$
where $p(\widehat{\btheta}_{\text{MAP}} \mid \mathcalH_1) \cdot \widehat{\sigma}_{\btheta}$ is defined as the \textit{Occam factor}~\footnote{When the posterior is approximated by a Gaussian, then the width is obtained by the determinant of the covariance matrix: $\widehat{\sigma}_{\btheta}=\det^{-1/2} \left(-\frac{1}{2\pi} \nabla^2 \ln p(\widehat{\btheta}_{\text{MAP}} \mid \mathcalX, \mathcalH_1)\right)$. See \citet{mackay1995probable} for more details.}.
The Occam factor is a value smaller than one  if $\widehat{\sigma}_{\btheta}< \sigma_{1}$, where the latter is the width of the prior distribution $p(\btheta \mid \mathcalH_1)$ (see Figure~\ref{fig:occam_factor}), and acts as a regularization that penalizes the parameter $\btheta$.

The width of the  posterior distribution signifies the uncertainty in parameter $\btheta$; while the width of the prior distribution represents the range of  values that were possible a priori.
Suppose the prior $p(\btheta \mid \mathcalH_1)$ is uniform. Then $p(\btheta \mid \mathcalH_1)=\frac{1}{\sigma_{1}}$, and the Occam factor is 
$$
\mathcalO_1=\frac{\widehat{\sigma}_{\btheta}}{\sigma_{1}},
$$
which quantifies how much the hypothesis space collapses upon observing the data.
The model $\mathcalH_1$ can be viewed as consisting of a certain number of exclusive submodels, of which only one remains viable upon receiving the data. (The Occam factor is the fraction that remains viable after seeing the data.) The logarithm of the Occam factor measures the amount of \textbf{information we gain} about the model's parameters when the data become available \citep{mackay1995probable}.

Now consider a second hypothesis $\{\mathcalH_2:\balpha=\balpha_2\}$ with a smaller width $\sigma_2<\sigma_1$.
And assume the posterior distribution under $\mathcalH_2$ and $\mathcalH_1$ are the same: $p(\btheta \mid \mathcalX, \mathcalH_i)$ with the same width $\widehat{\sigma}_{\btheta}$ (this is a strong assumption for ease of evaluation; see Figure~\ref{fig:occam_factor}). 
The corresponding Occam factors satisfy the following relationship:
$$
\mathcalO_1=\frac{\widehat{\sigma}_{\btheta}}{\sigma_{1}} < \mathcalO_2=\frac{\widehat{\sigma}_{\btheta}}{\sigma_{2}}.
$$
Although $\mathcalH_2$ has the same number of parameters as $\mathcalH_1$, it is more informative a priori---it commits more strongly to a specific region of parameter space. Thus, the Occam factor reveals that model complexity depends not only on the number of parameters but also on the prior distribution over those parameters. A model with strong prior constraints may be effectively simpler than one with vague priors, even if both have identical parameter counts.

\index{Graphical model representation}
\subsection{Graphical Model Representation}
We will explore  \textit{latent variable model}  in greater depth in Section~\ref{section:lvm}. For now, we provide a brief overview. 
A latent variable model extends the standard statistical framework by introducing two random vectors: an observed vector $\rvx\in \sX$ and an unobserved (latent) vector $\rvz\in \sZ$. 
These are jointly generated from a parametric family of distributions $\mathcalF = \{f_{\btheta} = f(\cdot , \cdot \mid \btheta) :\btheta\in\bTheta\} $. 
In other words, the pair $(\rvx, \rvz)$  is drawn from $f_{\btheta^*}$ for some true (but unknown) parameter $\btheta^*\in \bTheta$. 
However, only realizations of $\rvx$ are observed; the corresponding $\rvz$ values remain hidden.
Concretely, although the data-generating process produces paired samples $(\bx_1, \bz_1), (\bx_2, \bz_2), \ldots , (\bx_N, \bz_N)$, we only observe  $\bx_1, \bx_2, \ldots, \bx_N$. 
The unobserved components $\rvz$ are therefore called  \textit{latent variables} or \textit{hidden variables}. 

\begin{figure}[h]
\centering  
\vspace{-0.35cm} 
\subfigtopskip=2pt 
\subfigbottomskip=2pt 
\subfigcapskip=-5pt 
\subfigure[Latent variable model.]{\label{fig:lvm}
\includegraphics[width=0.431\linewidth]{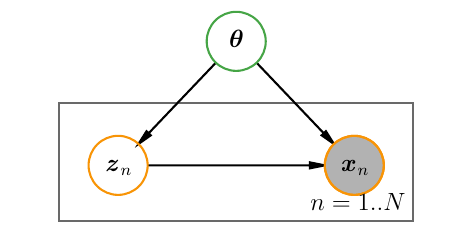}}
\subfigure[Latent variable model with hyperprior.]{\label{fig:lvm_hyper}
\includegraphics[width=0.431\linewidth]{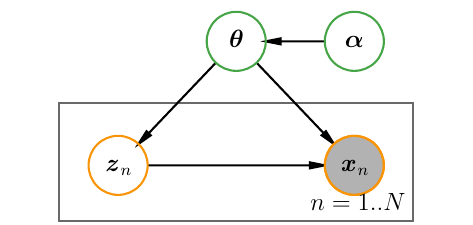}}
\caption{Graphical model representation of latent variable models. Green circles denote prior or hyperprior variables, orange circles represent observed (shaded) and latent (unshaded) variables, and plates indicate repeated structures across data points.}
\label{fig:lvm_and_hyp}
\end{figure}

To represent such models visually and reason about their structure, we use graphical models---a formalism that encodes probabilistic relationships and conditional dependencies among random variables using a graph. In this representation: 
(i) Nodes correspond to random variables or deterministic parameters;
(ii) Observed random variables are typically shown as \colorbox{\mdframecolor}{shaded} circles, while unobserved variables and parameters appear as unshaded circles;
(iii) Directed edges indicate direct probabilistic dependence (e.g., a parent node influences its child);
(iv) A plate (a rectangle enclosing a set of nodes) denotes replication: variables inside the plate are repeated independently across multiple instances (e.g., over $N$ data points).
Figure~\ref{fig:lvm}
illustrates the graphical model for the latent variable setup described above. 
For example, the edge from $\btheta$ to $\bz_n$ indicates that each latent variable $\bz_n$ depends on the global parameter $\btheta$. 
Similarly, $\bx_n$ depends on both $\bz_n$ and $\btheta$ (depending on the specific model).
The plate around $\bx_n$ and $\bz_n$ signifies that these variables are replicated for $n=1,2,\ldots,N$.

Furthermore, in Section~\ref{section:vb_inference}, we will treat the model parameter $\btheta$ itself as a random variable governed by a hyperprior $p(\btheta\mid\balpha)$, where $\balpha$ is a fixed hyper-parameter. The corresponding graphical model is shown in Figure~\ref{fig:lvm_hyper}.

\index{Approximate Bayesian inference}
\section{Approximate Bayesian Inference}
In this book, we focus on approximate probabilistic inference methods.
In some cases, it is computationally feasible to compute the posterior distribution exactly---for example, when using exponential family distributions with conjugate priors, which often admit closed-form solutions.
However, while exact inference is precise and useful for certain model classes, it becomes intractable in complex models. This is because exact methods typically rely on integrals, summations, or intermediate representations whose computational cost grows rapidly with the size of the state space, quickly becoming impractical. For instance, even in a Gaussian mixture model---where conjugate priors are available---the hierarchical structure often renders exact posterior computation infeasible.
Consequently, approximate inference methods are not only useful but often necessary in such settings.

Broadly speaking, there are two main families of approximate inference techniques:
(i) \textit{Variational methods} (also known as \textit{variational inference} or \textit{ensemble learning}), and 
(ii) \textit{Monte Carlo methods} (or  Monte Carlo approximations).
We now provide a brief comparison of these approaches \citep{bonawitz2008composable}.

In variational inference, we approximate the true posterior with a simpler, tractable distribution---typically chosen from a parametric family $q(\btheta\mid \blambda)$, where $\blambda$ denotes \textit{variational parameters}. 
The goal is to find the setting of $\blambda$ that makes $q$ as close as possible to the true posterior $p(\btheta\mid \mathcalX,\balpha)$, usually by minimizing a divergence measure (e.g., the Kullback--Leibler divergence). This optimization is typically performed deterministically, using gradient-based or other numerical methods. Once fitted, inference queries (e.g., expectations, marginals) are computed under the simplified distribution $q$.
The main advantage of variational methods is their computational efficiency and determinism. However, they provide only a lower bound on the marginal likelihood (the evidence), and the quality of the approximation depends critically on how well the chosen family $q$ can capture the structure of the true posterior. Despite this limitation, variational inference has become a cornerstone of Bayesian deep learning \citep{jordan1999introduction, graves2011practical, hoffman2013stochastic, ranganath2014black, mandt2014smoothed}. For a detailed example, see \citet{ma2014bayesian}.

In contrast, Monte Carlo methods sample directly from the target posterior distribution (or an approximation thereof). A set of samples $\{\btheta^{(1)},\btheta^{(2)},\ldots,\btheta^{(N)} \}$ is drawn, and inference is performed by treating this empirical collection as a proxy for the full distribution. Crucially, Monte Carlo estimators are asymptotically exact: as the number of samples increases, the approximation converges (almost surely) to the true target distribution. Thus, higher accuracy can always be achieved by running the algorithm longer---a property not shared by variational methods.
However, direct sampling from the posterior $p(\btheta\mid \mathcalX,\balpha)$ is rarely feasible in practice. Although the unnormalized posterior density $p(\mathcalX\mid \btheta)p(\btheta\mid \balpha)$ is often easy to evaluate, the normalizing constant $p(\mathcalX\mid \balpha)$ (the marginal likelihood) involves an intractable integral. This is where Markov chain Monte Carlo (MCMC) becomes essential; see Section~\ref{section:mcmc}.
Moreover, it is possible to design hybrid inference algorithms that combine variational inference for certain parts of the model with Monte Carlo methods for the remaining components.

\section{Monte Carlo (MC) Methods}\label{sec:monte_carlo_methods}
In Monte Carlo methods, we begin by drawing a sequence of $N$ samples $\{\btheta^{(1)}, \btheta^{(2)}, \ldots, \btheta^{(N)}\}$ from the posterior distribution $p(\btheta \mid  \mathcalX, \balpha)$, as defined in Equation~\eqref{equation:posterior_abstract_for_mcmc}.
We then approximate the target distribution by the empirical measure
\begin{equation}
p(\btheta \mid  -) \approx \overset{\sim}{p}(\btheta \mid  -) = \frac{1}{N} \sum_{n=1}^N \delta_{\btheta^{(n)}}(\btheta),
\end{equation}
where $\delta_{\btheta^{(n)}}(\btheta)$ denotes the Dirac delta function~\footnote{The Dirac delta function $\delta_{\bx_0}
	(\bx)$ satisfies $\int f(\bx)\delta_{\bx_0}\,d\bx = f(\bx_0)$;  it is zero everywhere except at $\bx = \bx_0$, where it is ``infinite" in such a way that its integral is 1.}. 
As the number of samples increases, the approximation
(almost surely) converges to the true target distribution, i.e., $\overset{\sim}{p}(\btheta) \overset{\overset{a.s.}{N\rightarrow \infty} }{\longrightarrow} p(\btheta)$.

Sampling-based methods like this are widely used in modern statistics due to their simplicity and broad applicability. Their core purpose is to approximate expectations of the form
\begin{equation}
\Exp [{h(\bTheta)}] = \int_{\btheta} h(\btheta) p(\btheta) d \btheta,
\end{equation}
when $\bTheta$ is a continuous random variable with probability density function (p.d.f.) $p$, or
\begin{equation}
\Exp[{h(\bTheta)}] = \sum_{\btheta} h(\btheta) p(\btheta),
\end{equation}
when $\bTheta$ is discrete with probability mass function (p.m.f.) $p$. 
In both cases, the Monte Carlo principle replaces the expectation with an empirical average:
\begin{equation}
\Exp[ {h(\bTheta)}] \approx \frac{1}{N}\sum_{n=1}^N h(\btheta^{(n)}).
\end{equation}

If it were generally feasible to sample directly from $p(\btheta \mid  \mathcalX, \balpha)$, Monte Carlo inference would be straightforward. Unfortunately, this is rarely possible in practice. 
We can consider the posterior form $p(\btheta \mid  \mathcalX, \balpha) = \frac{p(\mathcalX \mid  \btheta ) p(\btheta \mid  \balpha)}{p(\mathcalX \mid  \balpha)}$, where the unnormalized posterior $p(\mathcalX \mid  \btheta ) p(\btheta \mid  \balpha)$ can often be evaluated pointwise, but the normalizing constant $p(\mathcalX \mid  \balpha)$ is typically intractable due to high-dimensional integrals or sums. 
In such cases, Markov chain Monte Carlo (MCMC) provides a powerful alternative.

\subsection{Markov Chain Monte Carlo (MCMC)}\label{section:mcmc}

\textit{Markov chain Monte Carlo (MCMC)} algorithms, also called \textit{MCMC samplers}, are numerical methods that generate samples from a target distribution by constructing a Markov chain whose stationary distribution is the desired posterior $p(\btheta\mid \mathcalX, \balpha)$. 
Because MCMC explores the full solution space stochastically, it simultaneously yields both a ``best" estimate (e.g., a posterior mode or mean) and a quantification of uncertainty. Moreover, when supported by the data, the method can reveal multiple plausible solutions.
Intuitively, MCMC performs a stochastic hill-climbing search over the entire parameter space, allocating more computational effort to regions of high posterior probability \citep{andrieu2003introduction, bonawitz2008composable, hoff2009first, geyer2011introduction}.

The algorithm executes a stochastic walk through the state space $\bTheta$, designed so that, in the long run, the probability of visiting any state $\btheta$ matches its posterior probability. Samples from the true posterior are then approximated by recording the states visited during this walk, possibly after discarding initial samples (\textit{burn-in}) or applying \textit{thinning} to reduce autocorrelation.
This walk follows the Markov property: the next state depends only on the current one, not on the full history. Formally, if $\btheta^\toptzero$ denotes the state at iteration $t$, then
This ``history-free" property (i.e., $p(\btheta^\toptone\mid \btheta^{(1)}, \btheta^{(2)},\ldots, \btheta^\toptzero) = p(\btheta^\toptone\mid  \btheta^\toptzero)$) offers two key advantages:
\begin{itemize}
\item Memory requirements remain constant regardless of chain length.
\item The history-free property also indicates that the MCMC stochastic walk can be completely characterized by $p(\btheta^\toptone\mid  \btheta^\toptzero)$, known as the \textit{transition kernel}.
\end{itemize}

We then focus on the discussion of the transition kernel. The transition kernel $\bK$ can also be expressed as a linear transform. 
If $p_t = p_t (\btheta)$ is a row vector that encodes the probability of the walk being in state $\btheta$ at time $t$, then $p_{t+1} = p_t \bK$. If the stochastic walk starts from state $\btheta^{(0)}$, then the distribution from this initial state is the delta distribution $p_0 = \delta_{\btheta^{(0)}} (\btheta)$, and the state distribution for the chain after step $t$ is $p_t = p_0\bK^t$. We can easily find that the key to Markov chain Monte Carlo lies in choosing a kernel $\bK$ such that $\underset{t\rightarrow \infty}{\mathrm{lim}} p_t = p(\btheta \mid  \mathcalX, \balpha)$, independent of the choice of $\btheta^{(0)}$. Kernels exhibiting this property are said to converge to an \textit{equilibrium distribution} $p_{eq} = p(\btheta \mid  \mathcalX)$. Convergence is guaranteed if both of the following criteria are met (see, for example, \citet{bonawitz2008composable} for more details):
\begin{itemize}
\item \textit{Stationarity.} $p_{eq}$ is an invariant (or stationary) distribution for $\bK$. A distribution $p_{inv}$ is considered an invariant distribution for $\bK$ if $p_{inv} = p_{inv} \bK$;
\item \textit{Ergodicity.} $\bK$ is \textit{ergodic}. A kernel is called ergodic if it is \textit{irreducible} (meaning that any state can be reached from any other state) and \textit{aperiodic} (indicating that the stochastic walk never gets stuck in cycles).
\end{itemize}

Numerous MCMC algorithms exist---too many to cover here in full. Common examples include \textit{Gibbs sampling}, \textit{Metropolis--Hastings (MH)}, \textit{slice sampling}, \textit{Hamiltonian Monte Carlo}, and \textit{adaptive rejection sampling (ARS)}.
Though the name is potentially  misleading, Metropolis-within-Gibbs (MWG) was initially developed by \citet{metropolis1953equation}, and MH subsequently emerged as a generalization of MWG \citep{hastings1970monte}. All MCMC algorithms are recognized as special instances of the MH algorithm. Regardless of the specific method, the aim of Bayesian inference is to draw samples from the (unnormalized) joint posterior and use them to approximate marginal posterior distributions for downstream inference tasks.

The most generalizable MCMC algorithm is the MH generalization \citep{metropolis1953equation, hastings1970monte} of the MWG algorithm. The MH algorithm extended MWG to accommodate  asymmetric proposal distributions. 
In this method, it converts an arbitrary \textit{proposal kernel} $q(\btheta_{\star} \mid  \btheta^\toptzero  )$ into a transition kernel with the desired invariant distribution $p_{eq}(\btheta)$. 
In order to generate a sample from a MH transition kernel, the process involves drawing a proposal $\btheta_{\star} \sim q( \btheta_{\star} \mid  \btheta^\toptzero )$ and subsequently evaluating the \textit{MH acceptance probability} by 
\begin{equation}
\prob[A(\btheta_{\star}\mid  \btheta^\toptzero )]  
\triangleq  \min \left(1, \frac{p(\btheta_{\star} \mid  \balpha)q(\btheta^\toptzero \mid  \btheta_{\star}  )}{p(\btheta^\toptzero \mid  \balpha)q(\btheta_{\star} \mid  \btheta^\toptzero  )} \right),
\end{equation}
with probability $\prob[A(\btheta_{\star} \mid  \btheta^\toptzero  )] $ the proposal is accepted and we set $\btheta^\toptone =  \btheta_{\star}$; otherwise the
proposal is rejected and we set $\btheta^\toptone =  \btheta^\toptzero$. That is, 
\begin{equation}
\btheta^\toptone =\left\{
\begin{array}{ll}
\btheta_\star,  \text{ with probability } \prob[A(\btheta_{\star} \mid  \btheta^\toptzero  )]; \\
\btheta^\toptzero,  \text{ with probability } 1 - \prob[A( \btheta_{\star} \mid  \btheta^\toptzero  )] .
\end{array}
\right.
\end{equation}
Intuitively, the ratio $\frac{p(\btheta_{\star} \mid  \balpha)}{p(\btheta^\toptzero \mid  \balpha)}$ encourages moves toward higher-probability regions of the state space, while  the $\frac{q( \btheta^\toptzero \mid  \btheta_{\star} )}{q(\btheta_{\star} \mid  \btheta^\toptzero  )}$ term tends to accept moves that are easy to undo (corrects for asymmetry in the proposal mechanism). 
Since in MH, we only evaluate $p(\btheta)$ as a part of the ratio $\frac{p(\btheta_{\star} \mid  \balpha)}{p(\btheta^\toptzero \mid  \balpha)}$ (because the acceptance ratio depends only on the unnormalized posterior, the intractable marginal likelihood $p(\mathcalX\mid \balpha)$ cancels out), we do not need to compute the intractable normalization constant $p(\mathcalX \mid  \balpha)$ as mentioned in Section~\ref{sec:monte_carlo_methods}. 

Although the proposal kernel $q$ drives candidate generation, the actual transition kernel of MH is more complex. Informally, it combines the probability of accepting a move with the probability of staying in place:
$$
p (\btheta^\toptone \mid  \mathrm{accept}) \prob[\mathrm{accept}] + p (\btheta^\toptone\mid \mathrm{reject}) \prob[\mathrm{reject}].
$$ 
Formally, as shown by \citet{tierney1998note}, the transition kernel is:
\begin{equation}
\begin{aligned}
K(&\btheta^\toptzero \rightarrow \btheta^\toptone) 
= p(\btheta^\toptone \mid  \btheta^\toptzero) \\
&= q(\btheta^\toptone \mid  \btheta^\toptzero  ) A(\btheta^\toptone \mid  \btheta^\toptzero  ) + \delta_{\btheta^\toptzero}(\btheta^\toptone) \int_{\btheta_\star} q(\btheta_{\star} \mid  \btheta^\toptzero  ) (1 - A(\btheta_{\star} \mid  \btheta^\toptzero  )) .
\end{aligned}
\end{equation}

%
%

\subsection{MC vs. MCMC}
Both MC and MCMC aim to produce a sequence $\{\btheta^{(1)},\btheta^{(2)}, \ldots, \btheta^{(N)}\}$ such that, for any integrable function $h$,
\begin{equation}
\frac{1}{N} \sum_{n=1}^N h(\btheta^{(n)}) \approx \int_{\btheta} h(\btheta) p(\btheta) \,d\btheta, 
\end{equation}
(in the case of continuous random variables). 
In other words, the empirical average of $\{h(\btheta^{(1)}),h(\btheta^{(2)}), \ldots, h(\btheta^{(N)})\}$ should approximate the expectation of $h(\btheta)$ under the target distribution $p(\btheta)$. 
For this to hold reliably across a wide class of functions $h$, the empirical distribution of the samples  $\{\btheta^{(1)},\btheta^{(2)}, \ldots, \btheta^{(N)}\}$ must closely resemble $p(\btheta)$. 
MC and MCMC are two ways of generating such a sequence. 
MC simulation, in which we generate independent samples from the target distribution, is in some sense the ``true situation." Independent MC samples automatically create a sequence that is representative of $p(\btheta)$, which means the probability that $\btheta^{(n)} \in \sA, \,\forall\,n \in \{1,2,\ldots, N\}$ for any measurable set  $\sA$ is
\begin{equation}
\prob(\btheta^{(n)} \in \sA) =\int_{\sA} p(\btheta) \,d\btheta,
\end{equation}
In MCMC, samples are correlated due to the Markovian dependence. We only guarantee asymptotic correctness:
\begin{equation}
\lim_{n \rightarrow \infty} \prob(\btheta^{(n)} \in \sA) = \int_\sA p(\btheta) \,d\btheta.
\end{equation}
Thus, while MCMC is indispensable when direct sampling is infeasible, its finite-sample performance depends on chain mixing, autocorrelation, and convergence diagnostics (considerations absent in independent MC).

\subsection{Gibbs Sampler}\label{section:gibbs-sampler}

Gibbs sampling was first introduced by \citet{turchin1971computation} and later popularized by the brothers Geman and Geman \citep{geman1984stochastic} in the context of image restoration. They named the algorithm after the physicist \textit{J. Willard Gibbs}---roughly eight decades after his death---as an homage to the analogy between the sampling procedure and concepts in statistical physics.

Gibbs sampling is particularly useful when the \textbf{joint posterior distribution} $p(\btheta\mid \mathcalX,\balpha)$ is not known explicitly or is difficult to sample from directly. Instead, the method relies on the availability of \textbf{full conditional distributions} for each parameter, which are assumed to be tractable and easy to sample from.
The algorithm proceeds by iteratively sampling each component of the parameter vector $\btheta = \{\theta_1, \theta_2, \ldots, \theta_D\}$ from its conditional distribution given the current values of all other components. Formally, at iteration $t$, we update each $\theta_i$ as:
\begin{equation}\label{equation:gibbs_thetai_t}
\theta_i^\toptzero \sim p(\theta_i \mid  \btheta_{-i}^{(t-1)}, \mathcalX, \balpha),
\end{equation}
where $\btheta_{-i}^{(t-1)}$ denotes all parameters except $\theta_i$ (using their most recently updated values (from iteration $t-1$ or earlier in the same sweep).
Because each step conditions on the latest available values, Gibbs sampling is a componentwise MCMC algorithm. 
Under mild regularity conditions, the sequence of samples $\{\btheta^\toptzero\}_{t=1}^\infty$
converges in distribution to the target posterior $p(\btheta\mid \mathcalX,\balpha)$.

In practice, the initial portion of the chain---known as the \textit{burn-in period}---is discarded because the sampler has not yet reached its stationary distribution. Additionally, due to autocorrelation between successive draws, it is common to apply \textit{thinning}, i.e., retaining only every $k$-th sample, to reduce dependence and storage requirements.

A key insight that simplifies the derivation of full conditionals is that
\begin{equation}
p(\theta_i  \mid   \btheta_{- i}, \mathcalX)
= \frac{
p(\theta_1, \theta_2, \ldots,\theta_D, \mathcalX)
}{
p(\btheta_{- i}, \mathcalX)
} \propto p(\theta_1, \theta_2, \ldots,\theta_D, \mathcalX),
\end{equation}
since the denominator does not depend on $\theta_i$.
Thus, the conditional distribution is proportional to the joint distribution. This allows us to ignore terms constant with respect to $\theta_i$ when deriving sampling steps---a significant practical advantage.

As a simple illustration, consider a bivariate posterior $p(\theta_1,\theta_2\mid \mathcalX)$.
A Gibbs sampler alternates between: 
$$
\theta_1^\toptzero\sim p(\theta_1\mid \theta_2^{(t-1)},\mathcalX)
\qquad\text{and}\qquad 
\theta_2^\toptzero\sim p(\theta_2\mid \theta_1^\toptzero,\mathcalX),
$$
producing a sequence of states:
\begin{equation}
(\theta_1^{(0)}, \theta_2^{(0)}),\,\,\, 
(\theta_1^{(1)}, \theta_2^{(1)}), 
\,\,\,(\theta_1^{(2)}, \theta_2^{(2)}),\,\,\, \cdots, \nonumber
\end{equation}
which, under suitable conditions, converges to the joint distribution $p(\theta_1, \theta_2\mid \mathcalX)$. For further reading, see \citet{turchin1971computation, geman1984stochastic, hoff2009first, gelman2013bayesian}.

\subsection{Adaptive Rejection Sampling (ARS)}
\textit{Adaptive rejection sampling (ARS)} provides an efficient method for sampling from log-concave probability densities \citep{gilks1992adaptive, wild1993algorithm}. We offer a concise overview here; more details can be found in the original papers.



\begin{figure}[h!]
\centering
\includegraphics[width=0.5\textwidth]{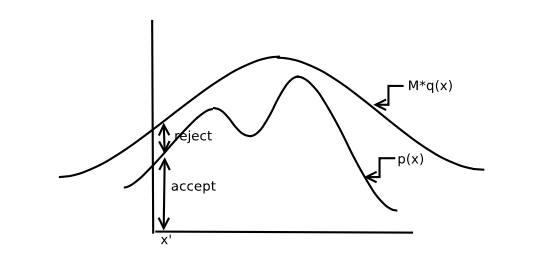}
\caption{Rejection sampling. Figure from Michael I. Jordan's lecture notes.}
\label{fig:rejection_sampling}
\end{figure}

\index{Rejection sampling}
\index{Adaptive rejection sampling}
\subsubsection{Rejection Sampling}


In rejection sampling, the objective is to sample from a target probability density function $p(x)$, given that we can sample from a probability density function $q(x)$ easily (known as the \textit{proposal density}). 
Although the target density $p(x)$ is unknown, the approach relies on establishing an envelop by considering $M \times q(x)$ such that it covers $p(x)$ for some $M>1$, as illustrated in Figure~\ref{fig:rejection_sampling}.
This is expressed as:
\begin{equation}
\frac{p(x)}{q(x)} < M, \text{ for all $x$.}
\end{equation}
Subsequently, when sampling $x_i$ from $q(x)$, and if $y_i=u \times M\times q(x_i)$ lies below the region under $p(x)$ for some $u \sim \mathrm{Uniform}(0,1)$, then we accept $x_i$; otherwise, it is rejected.

Accepted samples follow the distribution $p(x)$.
In essence, the method involves sampling $x_i$ from a distribution and making an acceptance or rejection decision based on the comparison with the envelope.
However, the efficiency of this method depends critically on how tightly $Mq(x)$ envelopes $p(x)$: a large $M$ leads to high rejection rates.

\subsubsection{Adaptive Rejection Sampling}\label{section:ars-sampling}
\begin{figure}[h!]
\centering
\includegraphics[width=0.5\textwidth]{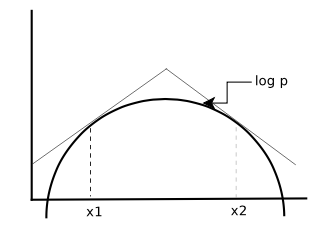}
\caption{Adaptive rejection sampling. Figure from Michael I. Jordan's lecture notes.}
\label{fig:adaptive_rejection_sampling}
\end{figure}

ARS improves upon standard rejection sampling by adaptively tightening the envelope around the target density---specifically when $p(x)$ is log-concave, meaning $\log p(x)$ is a concave function.
The basic idea involves dynamically constructing an upper envelope (the upper bound on $p(x)$), serving as an adaptive replacement for  $M\times q(x)$ in rejection sampling.

As shown in Figure~\ref{fig:adaptive_rejection_sampling}, the logarithm of the density, $\log p(x)$, is considered. We then sample $x_i$ from the upper envelope, and the sample is  either accepted or rejected akin to rejection sampling. 
In case of rejection, a tangent is drawn passing through $x = x_i$ and $y = \log(p)$; and the tangent is used to reduce the upper envelope so as to decrease the number of rejected samples. The intersections of these tangent planes enable the formation of an envelope adaptively. 
This adaptive refinement reduces the rejection rate over time. Because the envelope is built from tangents to a concave function, it always lies above $\log p(x)$, ensuring validity. 
To sample from the upper envelope, we need to transform from log space by exponentiating and using properties of the exponential distribution.

\index{Beta-Bernoulli model}
\index{Beta distribution}
\index{Bernoulli distribution}
\section{Bayesian Appetizers}\label{section:bayes_appetizers}
This section explores Bayesian inference using semi-conjugate priors with the Gibbs sampler and fully conjugate priors that do not require approximate inference. It provides an in-depth look at how Bayesian methods work.
Readers who already have a basic understanding of Bayesian inference may choose to skip this section.
\subsection{Beta-Bernoulli Model}\label{sec:beta-bernoulli}
We formally introduce a \textit{Beta-Bernoulli} model to illustrate the core ideas of Bayesian inference.
The \textit{Bernoulli distribution} models binary outcomes---i.e., random variables that take one of two possible values (typically 0 or 1). Its probability mass function, parameterized by $\theta$, is given by:
\begin{equation}\label{equation:bernoulli_distribution}
\bernoulli(x\mid \theta) = p(x\mid \theta) = \theta^x (1-\theta)^{1-x} \indicator(x\in \{0,1\}),
\end{equation}
(where $\indicator(\cdot)$ is the indicator function, equal to 1 when its argument is true and 0 otherwise)
or equivalently,
$$
\bernoulli(x\mid \theta)=p(x\mid \theta)=\left\{
\begin{aligned}
	&1-\theta ,& \mathrm{\,\,if\,\,} x = 0;  \\
	&\theta , &\mathrm{\,\,if\,\,} x =1,
\end{aligned}
\right.
$$
where $\theta$ is the probability of observing a 1 (success), and  $1-\theta$ is the probability of observing  a 0 (failure).
The mean of the  distribution is $\theta$. 

Suppose we observe a dataset $\mathcalX=\{x_1, x_2, \ldots , x_N\}$, where each $x_n$ is drawn independently from $\bernoulli( \theta)$. 
The likelihood of the data under this model is:
$$
\begin{aligned}
\text{likelihood} =  p(\mathcalX \mid \theta) &= \theta^{\sum x_n} (1-\theta)^{N-\sum x_n}.
\end{aligned}
$$
This expression, viewed as a function of $\theta$, is called the \textit{likelihood function} on $\mathcalX$.

\begin{SCfigure}
\centering
\includegraphics[width=0.55\textwidth]{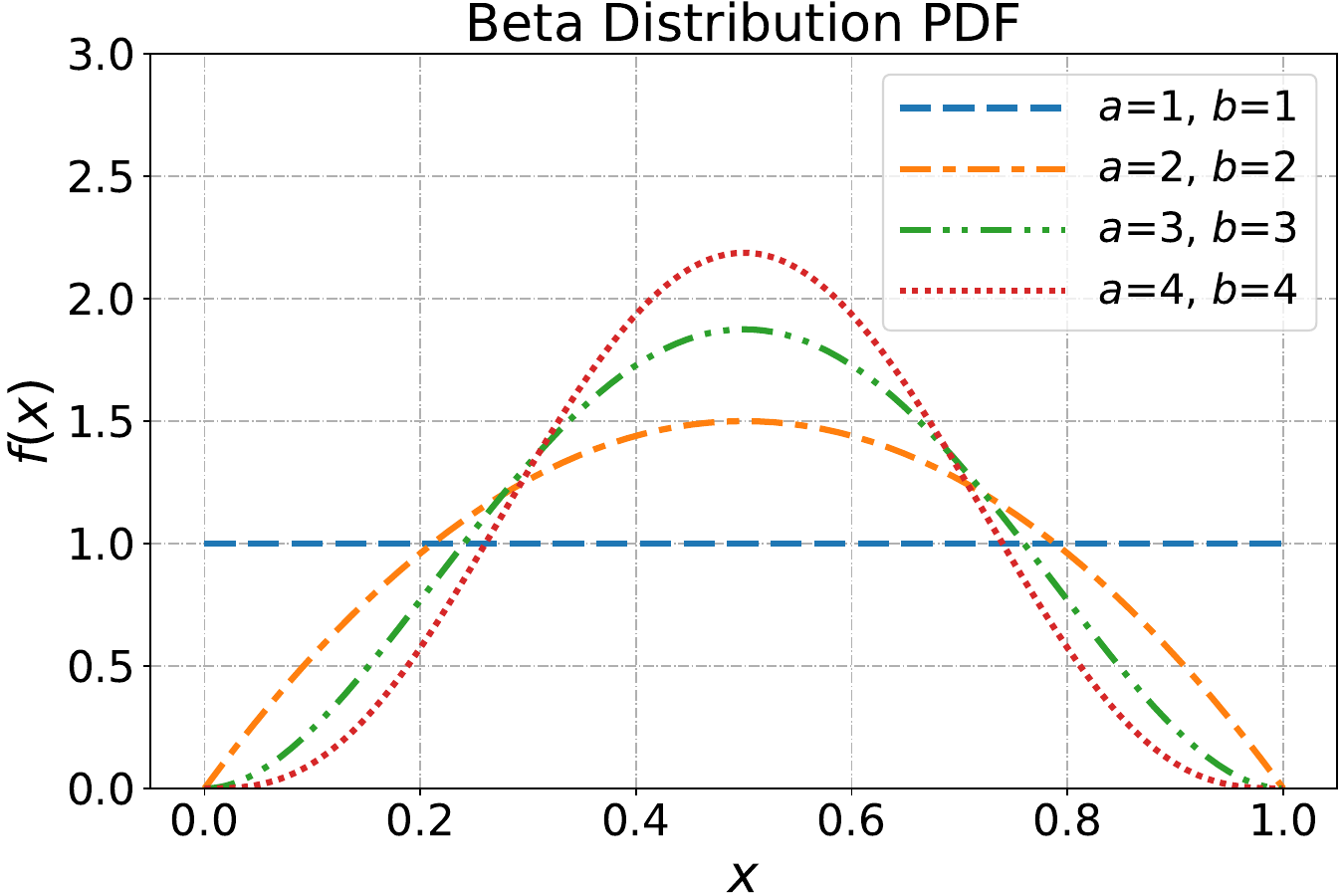}
\caption{Beta probability density functions for different values of the parameters $a$ and $b$. 
When $a=b=1$, the Beta distribution reduces to a \textit{uniform distribution} in the support of $[0,1]$.
The mean, variance, and mode  of the Beta distribution are $\Exp[\rx]=\frac{a}{a+b}$, 
$\Var[\rx]=\frac{ab}{(a+b+1)(a+b)^2}$, and $\mathrm{mode}[\rx]=\frac{a-1}{(a-1)+(b-1)}$ if $a>1, b>1$, respectively. }
\label{fig:dists_beta}
\end{SCfigure}

In the Bayesian framework, we place a prior distribution over the unknown parameter $\theta$. For the Bernoulli model, the natural choice is the \textit{Beta distribution}, whose probability density function is:
\begin{equation}
\mathrm{prior} = \betadist(\theta\mid a, b)=	p(\theta\mid a, b) =\frac{1}{B(a,b)} \theta^{a-1}(1-\theta)^{b-1} \indicator(0\leq \theta\leq 1), \nonumber
\end{equation}
where $B(a,b)$ denotes the \textit{Euler's Beta function} (a normalizing constant ensuring the density integrates to 1).
Figure~\ref{fig:dists_beta} shows Beta densities for different values of $a$ and $b$.
Notably, when $a=b=1$, the Beta distribution becomes the \textit{uniform distribution} on $[0,1]$.

Assigning a Beta prior to $\theta$, the posterior distribution is proportional to the product of the likelihood and the prior:
\begin{equation}
\begin{aligned}
\mathrm{posterior} = p(\theta\mid \mathcalX) 
&\propto p(\mathcalX \mid \theta) p(\theta\mid a,b) \\
&=\theta^{\sum x_n} (1-\theta)^{N-\sum x_n} \times \frac{1}{B(a,b)} \theta^{a-1}(1-\theta)^{b-1}\cdot \indicator(0\leq\theta\leq1) \\
&\propto \theta^{a+\sum x_n-1}(1-\theta)^{b+N-\sum x_n-1}\cdot \indicator(0\leq\theta\leq 1) \\
&\propto \betadist\left(\theta \,\bigg|\, a+\sum_{n=1}^N x_n, b+N-\sum_{n=1}^N x_n\right). \nonumber
\end{aligned}
\end{equation}
The posterior has the same functional form as the prior---both are Beta distributions, differing only in their parameters. When this occurs, the prior is called a \textit{conjugate prior}. Conjugacy offers significant computational advantages: it yields closed-form expressions for the posterior, simplifies differentiation, and makes sampling straightforward---eliminating the need for numerical approximation methods.

\begin{remark}[Prior Information in Beta-Bernoulli Model]
Comparing the prior and posterior reveals an intuitive interpretation: the hyper-parameter $a$ can be viewed as the prior count of successes (1s), and $b$ as the prior count of failures (0s). Their sum, $a+b$, reflects the effective prior sample size. An uninformative prior corresponds to $a=b=1$, which yields a uniform distribution over $[0,1]$.
\end{remark}

\index{Bayesian estimator}
\begin{remark}[Bayesian Estimator]
Like maximum likelihood estimation (MLE) or the \textit{method of moments (MoM)}, Bayesian inference provides a form of point estimation---but instead of returning a single best estimate, it yields a full posterior distribution $p(\theta \mid \mathcalX)$ over the parameter.

For prediction on a new observation $x'$, Bayesian inference integrates over the uncertainty in $\theta$:
\begin{equation}
p(x' \mid  \mathcalX) = \int p(x' \mid  \theta) p(\theta \mid  \mathcalX) d\theta.\nonumber
\end{equation}
Thus, predictions depend on the observed data $\mathcalX$ only through their influence on $\theta$:  $\mathcalX \rightarrow \theta \rightarrow x'$.
\end{remark}

\begin{SCfigure}
\centering
\includegraphics[width=0.55\textwidth]{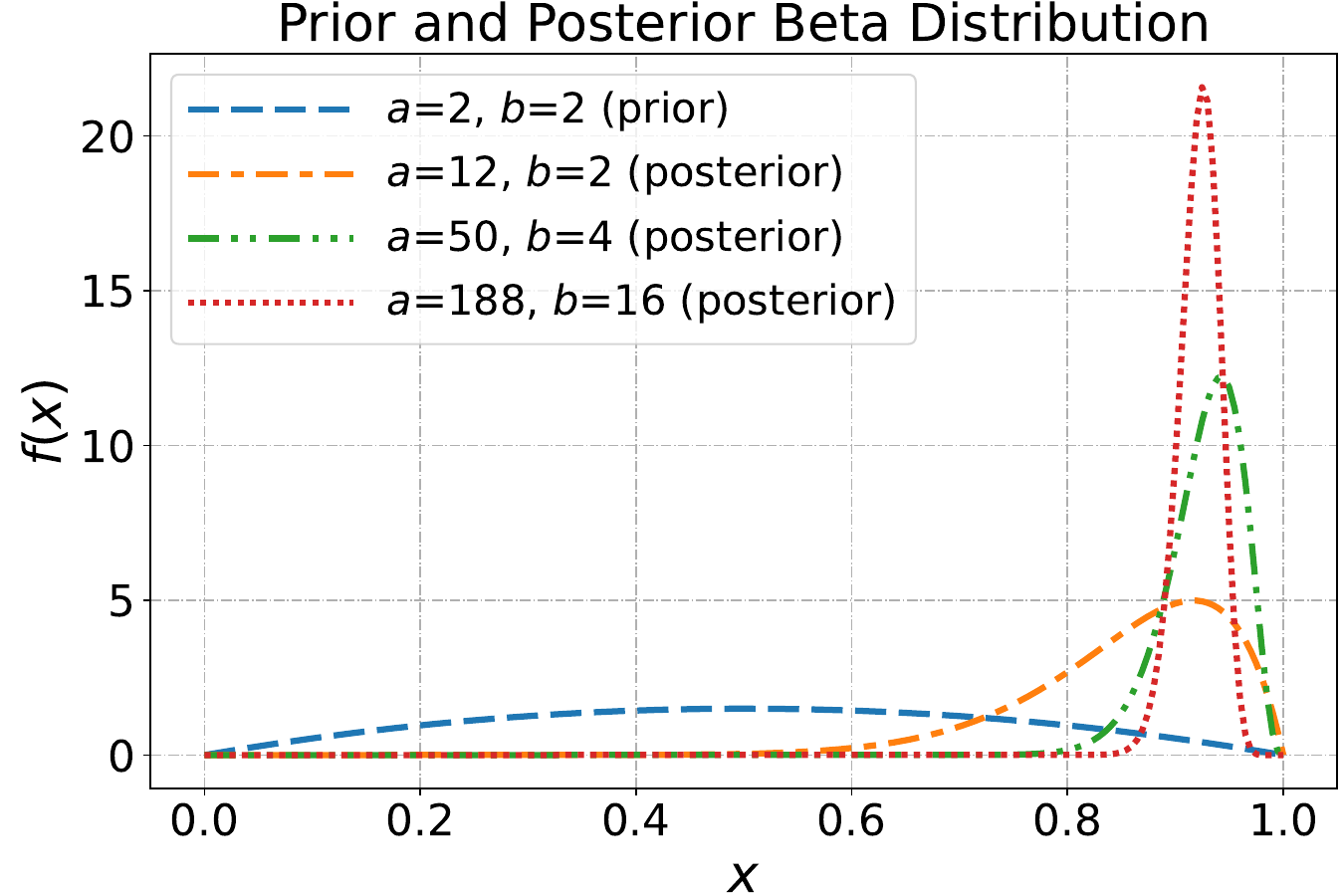}
\caption{The prior distribution is $\betadist(x\mid 2,2)$. The posterior distributions for the three cases in Example~\ref{example:amountofdata} are $\betadist(x\mid 12,2)$, $\betadist(x\mid 50,4)$, and $\betadist(x\mid 188,16)$, respectively.}
\label{fig:dists_beta_posterior}
\end{SCfigure}

\begin{example}[Amount of Data Matters\index{Add-one rule}]\label{example:amountofdata}
Bayesian methods are especially valuable with small or sparse datasets, where frequentist approaches may produce unreliable estimates.
Consider three scenarios involving observed successes in Bernoulli trials:
\begin{enumerate}
\item 10 successes out of 10 trials;
\item 48 successes out of 50 trials;
\item 186 successes out of 200 trials.
\end{enumerate}
The empirical success rates are 100\%, 96\%, and 93\%, respectively. However, the first case is based on very little data, so its estimate may be overly optimistic due to sampling noise.

Using a $\betadist(1,1)$ prior (uniform), the posterior mean success probabilities become: $\frac{11}{12}=91.6\%$, $\frac{49}{52}=94.2\%$, and $\frac{187}{202}=92.6\%$, respectively. 
Now, case 1 no longer appears more certain than case 2---a more reasonable conclusion given the limited data.

This adjustment is known as \textit{Laplace's rule of succession} \citep{ollivier2015laplace}. The ``add-one" rule (using a $\betadist(1,1)$ prior) adds one pseudo-count to both successes and failures, preventing zero-probability estimates and reflecting a uniform prior belief.
If we instead use a $\betadist(2,2)$ prior, Figure~\ref{fig:dists_beta_posterior} compares the prior and posterior distributions for the three cases.
\end{example}

\begin{mdframed}[hidealllines=\mdframehidelineNote,backgroundcolor=\mdframecolor,frametitle={Why Bayes?}]
This example highlights a key strength of Bayesian modeling: it naturally incorporates prior knowledge about parameters, which is especially useful for regularization when data are scarce. This property has contributed significantly to the widespread adoption of Bayesian methods.

In this framework, the prior $p(\theta)$ and likelihood $p(x\mid \theta)$ jointly encode a rational agent's beliefs. Bayes' rule then provides the optimal mechanism for updating those beliefs in light of observed data \citep{fahrmeir2007regression, hoff2009first}.

Of course, the prior $p(\theta)$ might not perfectly reflect true prior knowledge. But as the famous saying goes: ``\textit{All models are wrong, but some are useful}" \citep{box1987empirical}. If the prior reasonably approximates our beliefs, then the resulting posterior $p(\theta \mid  x)$ serves as a useful approximation to our updated beliefs.
\end{mdframed}

\subsection{Bayesian Linear Model with Zero-Mean Prior}\label{sec:bayesian-zero-mean}
In the linear model, given an input data matrix $\bX\in\real^{N\times D}$ and an  observation vector $\by\in\real^N$, e consider the system $\by=\bX\bbeta$, where $\bbeta\in\real^D$ is a vector of weights. When $N>D$, this system is overdetermined---that is, there are more equations than unknowns---and typically has no exact solution.
Let the column space of $\bX$ be defined as $\cspace(\bX)=\{\bX\bgamma: \forall \bgamma \in\real^D \}$. 
The absence of a solution to $\by=\bX\bbeta$ implies that  $\by\notin \cspace(\bX)$.  
To address this, we instead seek the weight vector $\bbeta$ that minimizes the mean squared error (MSE) between $\by$ and $\bX\bbeta$.

Rather than minimizing the MSE directly, we adopt a probabilistic perspective by introducing a Gaussian noise vector  $\bepsilon\in\real^N$:
$$
\rvy = \bX\bbeta + \bepsilon,
$$
where $\bepsilon \sim \normal(\bzero, \sigma^2 \boldsymbol{I})$ and $\sigma^2$ is a \textbf{fixed} variance (Here, $\normal(\ba, \bB)$ denotes a multivariate Gaussian distribution~\footnote{We defer the formal definition of the multivariate Gaussian distribution to Chapter~\ref{chapter:conjugate_models_bmf}  when we discuss regular conjugate models; see Definition~\ref{definition:multivariate_gaussian}.} with mean $\ba$ and covariance $\bB$.)
A detailed treatment of this model can be found in \citet{rasmussen2003gaussian, hoff2009first, lu2021rigorous}).

This additive Gaussian noise assumption leads naturally to a likelihood function. Given the observed inputs $\mathcalX (\bx_{1:N})= \{\bx_1, \bx_2, \ldots, \bx_N\}$ (with corresponding design matrix $\bX\in\real^{N\times D}$), the likelihood of the data is:
\begin{equation}\label{equation:linear_gaussian_addi_likelihood}
\mathrm{likelihood} = \by \mid  \bX, \bbeta, \sigma^2 \sim \normal(\bX\bbeta, \sigma^2\bI).
\end{equation}
We now place a multivariate Gaussian prior on the weight vector  $\bbeta$:
\begin{equation}
\mathrm{prior} = 	\bbeta \sim \normal(\bzero, \bSigma_0). \nonumber
\end{equation}
Applying Bayes' theorem, ``$\mathrm{posterior} \propto \mathrm{likelihood} \times \mathrm{prior} $," we obtain the posterior distribution:
\begin{equation}
\begin{aligned}
&\gap \mathrm{posterior}
= p(\bbeta\mid \by,\bX, \sigma^2) 
\propto p(\by\mid \bX, \bbeta, \sigma^2) \cdot  p(\bbeta \mid  \bSigma_0) \\
&=  \frac{1}{(2\pi \sigma^2)^{D/2}} \exp\left\{-\frac{1}{2\sigma^2} (\by-\bX\bbeta)^\top(\by-\bX\bbeta)\right\} 
\times \frac{1}{(2\pi)^{D/2}\abs{\bSigma_0}^{1/2}}\exp\left(-\frac{1}{2} \bbeta^\top\bSigma_0^{-1}\bbeta\right) \\
&\propto \exp\left\{-\frac{1}{2} (\bbeta - \bbeta_1)^\top \bSigma_1^{-1} (\bbeta - \bbeta_1)\right\}
\propto \normal(\bbeta_1, \bSigma_1)
, \nonumber
\end{aligned}
\end{equation}
where the posterior mean $\bbeta_1$ and covariance $\bSigma_1$ are given by
$$
\bbeta_1 \triangleq \left(\frac{1}{\sigma^2}\bX^\top\bX + \bSigma_0^{-1}\right)^{-1}
\left(\frac{1}{\sigma^2}\bX^\top\by\right),
\qquad 
\bSigma_1 \triangleq \left(\frac{1}{\sigma^2} \bX^\top\bX + \bSigma_0^{-1}\right)^{-1}.
$$
Thus, the posterior distribution is also Gaussian---the same family as the prior---making the Gaussian prior conjugate for this likelihood.

\paragrapharrow{A word on the notation.} 
We use $\{\bbeta_1,\bSigma_1\}$ to denote the posterior mean and covariance under the zero-mean Gaussian prior. For clarity, we will use $\{\bbeta_2,\bSigma_2\}$ and $\{\bbeta_3,\bSigma_3\}$ to represent the corresponding quantities in the semi-conjugate and fully conjugate prior settings, respectively (see subsequent sections).

\index{Ordinary least squares}
\paragrapharrow{Connection to ordinary least squares (OLS).}
The Bayesian linear model does not require $\bX$ to have full column rank. 
However, if $\bX$ is full rank (so that $\bX^\top\bX$ is invertible when $N>D$), 
then in the limit as the prior becomes non-informative---i.e., $\bSigma_0^{-1} \rightarrow \bzero$---the posterior mean converges to: $\bbeta_1 \rightarrow \widehatbbeta = (\bX^\top\bX)^{-1}\bX^\top\by$.
In this case, the maximum a posteriori (MAP) estimate coincides with the ordinary least squares (OLS) estimator. Moreover, the posterior distribution becomes:
$$
\bbeta\mid \by,\bX, \sigma^2 \sim \normal (\widehatbbeta, \sigma^2(\bX^\top \bX)^{-1}),
$$
which mirrors the sampling distribution of the OLS estimator under Gaussian errors:
$\widehatbbeta_{\text{OLS}} \sim \normal(\widehatbbeta, \sigma^2(\bX^\top \bX)^{-1})$  (see \citet{lu2021rigorous}).

\index{Ridge regression}
\index{Cross-validation}
\paragrapharrow{Ridge regression.}
In standard least squares, we approximate $\by$ by  $\bX\bbeta$. 
However, two practical issues may arise: (i) overfitting, especially when $D$ is large relative to $N$, and 
(ii)  singularity of $\bX^\top\bX$ when $\bX$ lacks full rank.
Ridge regression addresses these by penalizing large coefficients. Instead of minimizing $\norm{\by-\bX\bbeta}^2$, we minimize: $\norm{\by-\bX\bbeta}^2+\lambda\norm{\bbeta}^2$, where $\lambda>0$ is a regularization hyper-parameter, often selected via cross-validation (CV). The solution is:
\begin{equation*}
\widehatbbeta_{ridge} = \left(\bX^\top\bX + \lambda \bI\right)^{-1} \bX^\top\by.
\end{equation*}
Importantly, the matrix $(\bX^\top\bX + \lambda \bI)$ is always invertible (even when $\bX$ is rank-deficient), ensuring a unique solution. Further discussion of ridge regression is left to the reader.

\paragrapharrow{Connection to ridge regression.}
Observe that if we set the prior covariance to $\bSigma_0 = \bI$, then: $\bbeta_1 = \left(\bX^\top\bX + \sigma^2 \bI\right)^{-1} \bX^\top\by$ and $\bSigma_1 = \left(\frac{1}{\sigma^2}\bX^\top\bX+ \bI\right)^{-1}$. Since the posterior is $\bbeta\mid \by,\bX, \sigma^2 \sim \normal(\bbeta_1, \bSigma_1)$, the MAP estimate of $\bbeta$ becomes $\bbeta = \bbeta_1 =  (\bX^\top\bX + \sigma^2 \bI)^{-1} \bX^\top\by$, which matches the ridge regression estimator with $\sigma^2 = \lambda$. 
Thus, ridge regression corresponds exactly to the MAP estimate of a Bayesian linear model with a zero-mean isotropic Gaussian prior on $\bbeta$.
This provides a compelling Bayesian interpretation of ridge regression: it finds the mode of the posterior distribution, balancing data fit against prior belief in small parameter values---a natural form of regularization.


\index{Conjugate prior}
\index{Gamma distribution}
\subsection{Bayesian Linear Model with Semi-Conjugate Prior}\label{sec:semiconjugate}
We will use the \textit{Gamma distribution} as the prior for the precision (i.e., inverse variance) parameter of a Gaussian likelihood. A formal definition of the Gamma distribution is provided in Chapter~\ref{chapter:conjugate_models_bmf} (Definition~\ref{definition:gamma-distribution}) when we discuss conjugate models.
Regarding the motivation for choosing a Gamma prior on precision, we quote \citet{kruschke2014doing}:
\begin{mdframed}[hidealllines=\mdframehidelineNote,backgroundcolor=\mdframecolor]
Because of its role in conjugate priors for Gaussian likelihood functions, the Gamma distribution is routinely used as a prior for precision (i.e., inverse variance). But there is no logical necessity to do so, and modern MCMC methods permit more flexible specification of priors. Indeed, because precision is less intuitive than standard deviation, it can be more useful to give standard deviation a uniform prior that spans a wide range.
\end{mdframed}

In the same setting as Section~\ref{sec:bayesian-zero-mean}, we now consider the case where the variance $\sigma^2$ of the Gaussian likelihood is not fixed. The likelihood remains:
\begin{equation}
\mathrm{likelihood} = \by \mid  \bX, \bbeta, \sigma^2 \sim \normal(\bX\bbeta, \sigma^2\bI). \nonumber
\end{equation}
We specify a \textbf{non-zero-mean} Gaussian prior on the weight vector $\bbeta$,
and a Gamma prior on the precision $\gamma=1/\sigma^2$:
\begin{equation}
\begin{aligned}
{\color{mylightbluetext}\mathrm{prior:\,}} &\bbeta \sim \normal({\color{mylightbluetext}\bbeta_0}, \bSigma_0) \\
&{\color{mylightbluetext}\gamma = 1/\sigma^2 \sim \gammadist(a_0, b_0)}, \nonumber
\end{aligned}
\end{equation}
where the blue text highlights differences from earlier sections. 
Here, the Gamma density is defined as $\gammadist(a, b)=\frac{b^a}{\Gamma(a)} x^{a-1}\exp(-bx)$, $x>0$, and the \textit{Gamma function} is $\Gamma(a)=\int_{0}^{\infty} a^{t-1} \exp(-a)dt$.

\paragrapharrow{Step 1: Conditional posterior of $\bbeta$ given $\sigma^2$.}
Conditioning on $\sigma^2$ (or equivalently on $\gamma$), by  Bayes' theorem ``$\mathrm{posterior} \propto \mathrm{likelihood} \times \mathrm{prior} $," we get the conditional posterior density of $\bbeta$:
\begin{equation}
\begin{aligned}
\mathrm{posterior}&= p(\bbeta\mid \by,\bX, \sigma^2) 
\propto p(\by\mid \bX, \bbeta, \sigma^2) \cdot p(\bbeta \mid  \bbeta_0, \bSigma_0) \\
&=  \frac{1}{(2\pi \sigma^2)^{D/2}} \exp\left\{-\frac{1}{2\sigma^2} (\by-\bX\bbeta)^\top(\by-\bX\bbeta)\right\} \\
&\gap \times \frac{1}{(2\pi)^{D/2}\abs{\bSigma_0}^{1/2}}\exp\left\{-\frac{1}{2} (\bbeta-\bbeta_0)^\top\bSigma_0^{-1}(\bbeta-\bbeta_0)\right\} \\
&\propto \exp\left\{-\frac{1}{2} (\bbeta - \bbeta_2)^\top \bSigma_2^{-1} (\bbeta - \bbeta_2)\right\}
\propto \normal(\bbeta_2, \bSigma_2)
, \nonumber
\end{aligned}
\end{equation}
with posterior parameters
$$
\begin{aligned}
\bSigma_2 &\triangleq \left(\frac{1}{\sigma^2} \bX^\top\bX + \bSigma_0^{-1}\right)^{-1},\\
\bbeta_2 &\triangleq \bSigma_2 \left(\bSigma_0^{-1}\bbeta_0+\frac{1}{\sigma^2}\bX^\top\by\right)
= \left(\frac{1}{\sigma^2}\bX^\top\bX + \bSigma_0^{-1}\right)^{-1}
\left(\textcolor{mylightbluetext}{\bSigma_0^{-1}\bbeta_0}+\frac{1}{\sigma^2}\bX^\top\by\right).
\end{aligned}
$$
Thus, the conditional posterior is Gaussian:
\begin{equation}
\mathrm{posterior} = \bbeta\mid \by,\bX, \sigma^2 \sim \normal(\bbeta_2, \bSigma_2). \nonumber
\end{equation}

\paragrapharrow{Connection to the zero-mean prior model.}
The relationship between the zero-mean and non-zero-mean (semi-conjugate) prior models is as follows:
\begin{enumerate}
\item The posterior mean $\bbeta_1$ from Section~\ref{sec:bayesian-zero-mean} is a special case of $\bbeta_2$ when $\bbeta_0=\bzero$. 
\item If  $ \bX $  has full column rank and the prior becomes non-informative ($\bSigma_0^{-1} \rightarrow \bzero$), then $\bbeta_2 \rightarrow \widehatbbeta = (\bX^\top\bX)^{-1}\bX\by$, recovering the OLS estimate.
\item  As $\sigma^2 \rightarrow \infty$ (i.e., the data become uninformative), $\bbeta_2$ is approximately approaching  $\bbeta_0$, the prior expectation of parameter. 
In contrast, under a zero-mean prior,  $ \bbeta_1 \to \bzero $  in this limit.
\item  \textbf{Weighted average interpretation}: We can rewrite $\bbeta_2$ as 
\begin{equation}
\begin{aligned}
\bbeta_2 &=  \left(\frac{1}{\sigma^2}\bX^\top\bX + \bSigma_0^{-1}\right)^{-1}
\left(\bSigma_0^{-1}\bbeta_0+\frac{1}{\sigma^2}\bX^\top\by\right) \\
&= \left(\frac{1}{\sigma^2}\bX^\top\bX + \bSigma_0^{-1}\right)^{-1} \bSigma_0^{-1}\bbeta_0 + \left(\frac{1}{\sigma^2}\bX^\top\bX + \bSigma_0^{-1}\right)^{-1} \frac{\bX^\top\bX}{\sigma^2} (\bX^\top\bX)^{-1}\bX^\top\by \\
&=(\bI-\bA)\bbeta_0 + \bA \widehatbbeta, \nonumber
\end{aligned}
\end{equation}
where $\widehatbbeta=(\bX^\top\bX)^{-1}\bX^\top\by$ is the OLS estimate of $\bbeta$, and $\bA=(\frac{1}{\sigma^2}\bX^\top\bX + \bSigma_0^{-1})^{-1} \frac{\bX^\top\bX}{\sigma^2}$. 
We observe that the posterior mean of $\bbeta$ is a weighted average of the prior mean and the OLS estimate of $\bbeta$. Thus, if we set the prior parameter $\bbeta_0 = \widehatbbeta$, the posterior mean of $\bbeta$ becomes precisely $\widehatbbeta$.
\end{enumerate}

\paragrapharrow{Step 2: Conditional posterior of $\gamma=1/\sigma^2$ given $\bbeta$.}
Now conditioning on $\bbeta$, the conditional posterior of the precision $\gamma$ is:
\begin{equation}
\begin{aligned}
\mathrm{posterior}&= p(\gamma=\frac{1}{\sigma^2}\mid \by,\bX, \bbeta) 
\propto p(\by\mid \bX, \bbeta, \gamma) \cdot p(\gamma \mid  a_0, b_0) \\
&=  \frac{\gamma^{N/2}}{(2\pi )^{N/2}} 
\exp\left\{-\frac{\gamma}{2} (\by-\bX\bbeta)^\top(\by-\bX\bbeta)\right\} 
\times \frac{{b_0}^{a_0}}{\Gamma(a_0)} \gamma^{a_0-1} \exp(-b_0 \gamma) \\
&\propto \gamma^{(a_0+\frac{N}{2}-1)} \exp\left\{-\gamma\left[b_0+\frac{1}{2}(\by-\bX\bbeta)^\top(\by-\bX\bbeta)\right]\right\}. \nonumber
\end{aligned}
\end{equation}
This is the kernel of a Gamma distribution, so:
\begin{equation}
\mathrm{posterior\,\, of\,\,} \gamma \mathrm{\,\,given\,\,} \bbeta  = \gamma\mid \by,\bX, \bbeta \sim \gammadist\left(a_0+\frac{N}{2}, \left[b_0+\frac{1}{2}(\by-\bX\bbeta)^\top(\by-\bX\bbeta)\right]\right). \nonumber
\end{equation}

\paragrapharrow{Prior information on the noise/precision.}
The Gamma prior admits an intuitive interpretation:
\begin{enumerate}
\item  We notice that the prior mean and posterior mean of $\gamma$ are $\Exp[\gamma]=\frac{a_0}{b_0}$ and $\Exp[\gamma \mid \bbeta]={(a_0 + \frac{N}{2})}/{\big(b_0 +\frac{1}{2}(\by-\bX\bbeta)^\top(\by-\bX\bbeta)\big)}$, respectively. So the latent meaning of $2 a_0$ is the prior sample size associated with the noise variance $\sigma^2 = \frac{1}{\gamma}$. 

\item  As we assume $\rvy=\bX\bbeta +\bepsilon$, where $\bepsilon \sim \normal(\bzero, \sigma^2\bI)$, then ${(\rvy-\bX\bbeta)^\top(\rvy-\bX\bbeta)}/{\sigma^2} \sim \chi^2(N)$ and $\Exp\left[\frac{1}{2}(\by-\bX\bbeta)^\top(\by-\bX\bbeta)\right] = \frac{N}{2}\sigma^2$ \footnote{$\chi^2(N)$ is a Chi-squared distribution with $N$ degrees of freedom. See Definition~\ref{definition:chisquare_distribution}.}. 
So the latent meaning of ${b_0}/{a_0}$ is the prior guess for the noise variance $\sigma^2$.

\item  ome textbooks parameterize the prior as   $\gamma \sim  \gammadist(n_0/2, n_0\sigma_0^2/2)$ to make this explicit (in which case, $n_0$ is the prior sample size, and $\sigma_0^2$ is the prior variance). 
While this makes the interpretation explicit, the form may seem arbitrary without context.
\end{enumerate}

\index{Gibbs sampler}
\paragrapharrow{Gibbs sampler.}
Using the Gibbs sampling framework introduced in Section~\ref{section:gibbs-sampler}, we can construct a Gibbs sampler for this semi-conjugate Bayesian linear model as follows:

0. Set initial values to $\bbeta$ and $\gamma = \frac{1}{\sigma^2}$;

1. update $\bbeta$: $\mathrm{posterior} = \bbeta\mid \by,\bX, \gamma \sim \normal(\bbeta_2, \bSigma_2)$;

2. update $\gamma$: $\mathrm{posterior}  = \gamma\mid \by,\bX, \bbeta \sim \gammadist\left(a_0+\frac{N}{2}, [b_0+\frac{1}{2}(\by-\bX\bbeta)^\top(\by-\bX\bbeta)]\right)$.

\index{Conjugate prior}
\index{Normal-inverse-Gamma distribution}
\index{Bayes' theorem}
\subsection{Bayesian Linear Model with Full Conjugate Prior}\label{section:blm-fullconjugate}
Introducing a Gamma prior on the precision (i.e., inverse variance) $\gamma=1/\sigma^2$ is mathematically equivalent to placing an inverse-Gamma prior on the variance $\sigma^2$.
\footnote{We defer the formal definition of the inverse-Gamma distribution to Definition~\ref{definition:inverse_gamma_distribution} in the section on regular conjugate models.}
This setting mirrors the semi-conjugate prior model described in Section~\ref{sec:semiconjugate}, with the same likelihood:
\begin{equation}
\mathrm{likelihood} = \by \mid  \bX, \bbeta, \sigma^2 \sim \normal(\bX\bbeta, \sigma^2\bI). \nonumber
\end{equation}
However, we now specify a joint prior over ($\bbeta,\sigma^2$) as follows:
\begin{equation}
\begin{aligned}
{\color{mylightbluetext}\mathrm{prior:\,}} &\bbeta\mid \sigma^2 \sim \normal(\bbeta_0, {\color{mylightbluetext}\sigma^2} \bSigma_0) \\
&{\color{mylightbluetext}\sigma^2 \sim \inversegammadist(a_0, b_0)}, \nonumber
\end{aligned}
\end{equation}
where the blue text distinguishes this formulation from earlier ones.
This joint prior is known as the \textit{normal-inverse-Gamma (NIG)} distribution and can be written compactly as:
\begin{equation}
\begin{aligned}
\mathrm{prior:\,} &\bbeta,\sigma^2 \sim \nig(\bbeta_0, \bSigma_0, a_0, b_0) =  \normal(\bbeta_0, \sigma^2 \bSigma_0)\cdot \inversegammadist(a_0, b_0) . \nonumber
\end{aligned}
\end{equation}
Applying Bayes' theorem, ``$\mathrm{posterior} \propto \mathrm{likelihood} \times \mathrm{prior} $," we obtain the joint posterior:
\begin{equation}
\begin{aligned}
\mathrm{posterior}&= p(\bbeta,\sigma^2\mid \by,\bX)
\propto p(\by\mid \bX, \bbeta, \sigma^2)\cdot  p(\bbeta, \sigma^2 \mid  \bbeta_0, \bSigma_0, a_0, b_0) \\
&=  \frac{1}{(2\pi \sigma^2)^{N/2}} \exp\left\{-\frac{1}{2\sigma^2} (\by-\bX\bbeta)^\top(\by-\bX\bbeta)\right\} \\
&\, \times \frac{1}{(2\pi \sigma^2)^{D/2} \sqrt{\abs{\bSigma_0}}} \exp\left\{\frac{-1}{2\sigma^2} (\bbeta - \bbeta_0)^\top\bSigma_0^{-1} (\bbeta - \bbeta_0)\right\}
\cdot  \frac{{b_0}^{a_0}}{\Gamma(a_0)} \frac{1}{(\sigma^2)^{a_0+1}} \exp(\frac{-b_0}{\sigma^2}) \\
&\propto \frac{1}{(2\pi \sigma^2)^{D/2} }  \exp\left\{  \frac{1}{2\sigma^2} (\bbeta -\bbeta_3)^\top\bSigma_3^{-1}(\bbeta -\bbeta_3) \right\} \\
& \,\times \frac{1}{(\sigma^2)^{a_0 +\frac{N}{2}+1}} \exp\left\{-\frac{1}{\sigma^2} \left[b_0+\frac{1}{2} (\by^\top\by +\bbeta_0^\top\bSigma_0^{-1}\bbeta_0 -\bbeta_3^\top\bSigma_3^{-1}\bbeta_3) \right]\right\}, \nonumber
\end{aligned}
\end{equation}
where the parameters are 
$$
\begin{aligned}
\bSigma_3 &= \left( \bX^\top\bX + \bSigma_0^{-1}\right)^{-1}, \\
\bbeta_3 &= \bSigma_3\left(\bX^\top\by + \bSigma_0^{-1}\bbeta_0\right) = 
\left( \bX^\top\bX + \bSigma_0^{-1}\right)^{-1}(\bSigma_0^{-1}\bbeta_0 + \bX^\top\by).
\end{aligned}
$$ 
Let $a_N \triangleq a_0 +\frac{N}{2}$ and $b_N\triangleq b_0+\frac{1}{2} (\by^\top\by +\bbeta_0^\top\bSigma_0^{-1}\bbeta_0 -\bbeta_3^\top\bSigma_3^{-1}\bbeta_3) $. The posterior thus follows a NIG distribution with updated parameters:
\begin{equation}
\begin{aligned}
\mathrm{posterior}&=  \bbeta, \sigma^2 \mid  \by, \bX \sim \nig(\bbeta_3, \bSigma_3, a_N, b_N). \nonumber
\end{aligned}
\end{equation}

\paragrapharrow{Connection to zero-mean prior and semi-conjugate prior models.}
The full conjugate (NIG) model generalizes both the zero-mean and semi-conjugate Bayesian linear models:
\begin{enumerate}
\item If  $\bX$ has full column rank and the prior becomes non-informative ($\bSigma_0^{-1} \rightarrow \bzero$), then $\bbeta_3 \rightarrow \widehatbbeta = (\bX^\top\bX)^{-1}\bX\by$, recovering the OLS estimate.

\item  As $b_0 \rightarrow \infty$ (i.e., the prior strongly believes the noise variance is large), then $\sigma^2 \rightarrow \infty$ and $\bbeta_3$ is approximately approaching $\bbeta_0$, the prior expectation of parameter. 
This parallels the behavior in the semi-conjugate model (Section~\ref{sec:semiconjugate}), where fixing  $ \sigma^2 \to \infty $  also causes  $ \bbeta_2 \to \bbeta_0 $ .

\item\textbf{Weighted average interpretation}: We can rewrite  $\bbeta_3$ as
\begin{equation}
\begin{aligned}
\bbeta_3 &=  \left(\bX^\top\bX + \bSigma_0^{-1}\right)^{-1}(\bSigma_0^{-1}\bbeta_0+\bX^\top\by) \\
&= \left(\bX^\top\bX + \bSigma_0^{-1}\right)^{-1} \bSigma_0^{-1}\bbeta_0 + 
\left(\bX^\top\bX + \bSigma_0^{-1}\right)^{-1} (\bX^\top\bX) (\bX^\top\bX)^{-1}\bX^\top\by \\
&=(\bI-\bC)\bbeta_0 + \bC \widehatbbeta, \nonumber
\end{aligned}
\end{equation}
where $\widehatbbeta=(\bX^\top\bX)^{-1}\bX^\top\by$ is the OLS estimate of $\bbeta$, and $\bC=(\bX^\top\bX + \bSigma_0^{-1})^{-1} (\bX^\top\bX)$. We observe that the posterior mean of $\bbeta$ is a weighted average of the prior mean and the OLS estimate of $\bbeta$. Thus, if we set $\bbeta_0 = \widehatbbeta$, the posterior mean of $\bbeta$ precisely equals  $\widehatbbeta$.

\item  From the relationship of $a_N = a_0 +\frac{N}{2}$, we can interpret  $ 2a_0 $  as the prior sample size associated with the noise variance  $ \sigma^2 $ .

\item  The posterior precision matrix satisfies   $\bSigma_3^{-1} = \bX^\top\bX + \bSigma_0^{-1}$, 
showing that \textit{posterior precision = data precision ($\bX^\top\bX$) + prior precision}---a fundamental property of Gaussian conjugate updating.

\end{enumerate}
\index{Weighted average}

\section{Variational Bayesian Inference}\label{section:vbi_root}

Monte Carlo-based approximate inference algorithms evaluate the quality of an approximation by using sampling to represent the target distribution, with the goal of making the approximation as close as possible to the true posterior. These methods draw a large number of samples from a proposal distribution and use them to estimate key properties of the target distribution---such as its mean, variance, or higher-order moments. Consequently, the accuracy of the approximation depends directly on how well the sample set represents the target and on the efficiency of the sampling procedure.
In contrast, many other approximate inference techniques pursue closeness to the target distribution by minimizing a specific divergence measure that quantifies the discrepancy between the approximation and the true posterior. A prominent example of this approach is \textit{variational inference (VI)}, also referred to as \textit{variational Bayesian inference}.
In VI, the chosen measure is the \textit{Kullback--Leibler (KL) divergence}, which quantifies how one probability distribution diverges from another reference distribution. By casting inference as an optimization problem---specifically, minimizing the KL divergence---variational inference identifies the parameters of a simpler, tractable distribution (called the \textit{variational distribution}) that best approximates the intractable posterior. This formulation enables a principled trade-off between approximation accuracy and computational cost, often yielding efficient and scalable solutions for approximate inference in complex probabilistic models \citep{jordan1999introduction, wainwright2008graphical}.

\subsection{Motivating Model: Latent Variables and Alternating Methods}\label{section:lvm}
\index{Variational inference}
\index{Hidden variables}
\index{Latent variables}
\index{Global latent variables}

Consider a statistical model that jointly generates two random vectors, $\rvx\in \sX$ and $\rvz\in \sZ$, instead of one, according to a distribution from a parametric family $\mathcalF = \{f_{\btheta} = f(\cdot , \cdot \mid \btheta) :\btheta\in\bTheta\} $; that is, $(\rvx, \rvz) \sim f_{\btheta^*}$ for some true parameter $\btheta^*\in \bTheta$. 
However, we only observe realizations of  $\rvx$; the components of  $\rvz$ remain unobserved.
More precisely, although the (unknown) parameter $\btheta^*$ generates  $N$ pairs i.i.d. pairs  $(\bx_1, \bz_1), (\bx_2, \bz_2), \ldots , (\bx_N, \bz_N)$, we only have access to the observed data $\mathcalX=\{\bx_1, \bx_2, \ldots, \bx_N\}$. 
The unobserved  $\rvz$ components are therefore referred to as  \textit{latent variables} or \textit{hidden variables}. 
Examples include the component indicators in Gaussian or Bernoulli mixture models (see Problems~\ref{prob:mix_of_gauss}--\ref{prob:mix_of_bern})~\footnote{While more complex generative structures are possible, we restrict our attention here to the standard setting in which each observed data point is associated with its own latent variable. A more general framework will be introduced in later sections.}.
Latent variables are integral to the model but are not part of the observed dataset. Despite being unmeasurable, they play a crucial role in explaining the underlying structure and variability of the observed data.

In this context, the model parameter $\btheta$ is often regarded as a \textit{global latent variable}, as it governs the entire data-generating process and influences all observations and their associated  latent variables. 
In contrast, the individual hidden variables  $\{\bz_n\}_{n=1}^N$ are called \textit{local latent variables} (or simply \textit{latent variables}), since each $\bz_n$ is tied exclusively to its corresponding observation 
$\bx_n$ and accounts for the variation specific to that data point (see the graphical model in Figure~\ref{fig:lvm} for an illustration of this hierarchical relationship).

Latent variables commonly arise in real-world applications. For instance, in clustering, they may represent the unknown cluster assignments of data points \citep{beal2003variational, jain2017non, lu2021survey}. 
In topic modeling, they can indicate the latent topics underlying a collection of documents, and in image analysis, they might encode high-level features or structures not evident in raw pixel values \citep{blei2003latent}. 
\textit{Latent variable models (LVMs)} thus enable us to uncover such hidden patterns and draw more informed conclusions about the data-generating mechanism.

The most direct approach to parameter estimation in such models is to compute the \textit{maximum likelihood (ML) estimator} of the parameters $\btheta$. This is done by maximizing the marginal likelihood, which can be expressed as:
\begin{equation}\label{equation:lvm_maxmarg}
\btheta^* = \mathop{\argmax}_{\btheta\in\bTheta}\, 
\left\{\ell(\btheta; \bx_1, \bx_2, \ldots, \bx_N) 
\triangleq \prod_{n=1}^N \int p(\bz_n, \bx_n\mid \btheta) \,d\bz_n
\right\},
\end{equation}
where we assume the support set of the random vector $\rvz$ is continuous. When it is discrete, the integral is replaced by a sum: $\ell(\btheta; \bx_1, \bx_2, \ldots, \bx_N)=\prod_{n}\sum_{\bz_n\in\sZ} p(\bz_n, \bx_n\mid \btheta) \, d\bz_n$. 
However, in most practical cases, maximizing the marginal likelihood is computationally intractable; 
for example, when the support $\sZ$ is discrete, the summation over $\sZ$ involves $\abs{\sZ}^N$ terms~\footnote{$\abs{\sZ}$ denotes  the cardinality of the set $\sZ$.}, making both explicit evaluation and optimization of the likelihood function infeasible.
\footnote{This intractability is a hallmark of many latent variable models and motivates the use of approximate inference techniques---such as Monte Carlo sampling, variational inference, or the EM algorithm---which avoid direct computation of the marginal likelihood by approximating either the posterior distribution or the likelihood itself. These methods yield tractable solutions even in high-dimensional or complex latent spaces.
}

Faced with this intractability, a common strategy is to adopt an \textit{alternating maximization} approach. This iterative method alternates between two steps: (1) inferring the latent variables given the current parameter estimate, and (2) updating the model parameters based on the inferred latent variables.
More concretely, suppose the true parameter $\btheta^*$ were known. 
We could then estimate each latent variable via \textit{maximum a posteriori (MAP)} assignment:
$$
\text{(AM-Step 1): }\quad 
\widehatbz_n = \mathop{\argmax}_{\bz\in\sZ} p(\bz \mid \bx_n, \btheta^*)\gap \forall n\in\{1,2,\ldots,N\}.
$$
Conversely, if the latent variables $\{\bz_n\}$ were known, the ML estimate of $\btheta$ would be:
$$
\text{(AM-Step 2): }
\qquad 
\widehat{\btheta}_{\text{MLE}} =\mathop{\argmax}_{\btheta\in\bTheta} \, \ell(\btheta; \{(\bx_n,\bz_n)\}_{n=1}^N ).
$$
This step refines the parameter estimates based on the full set of data, including the inferred latent variables.
This alternating procedure yields a practical algorithm for fitting latent variable models; see Algorithm~\ref{alg:am_lvm}.

\paragrapharrow{Limitation and EM algorithm.}
Although intuitive, this alternating maximization scheme has notable limitations---especially when the latent space $\sZ$ is large or structured. 
At each iteration $t$, the algorithm makes a ``hard assignment," assigning the data point $\bx_n$ to just one specific value of the latent variable $\widehatbz_n^\toptzero \in \sZ$ (we use  subscript $(t)$ to denote the iteration count). 
This ignores other plausible latent configurations $\bz^\prime$  for which the posterior probability $p(\bz^\prime \mid  \bx_n, \btheta^\toptzero)$ may still be substantial, albeit lower than the maximum $ p(\bz_n^\toptzero \mid \bx_n, \btheta^\toptzero)$. 
As a result, valuable uncertainty information is discarded, potentially leading to suboptimal or unstable estimates.

The \textit{expectation-maximization (EM)} algorithm is designed to address this issue \citep{baum1970maximization, dempster1977maximum}.
Rather than committing to hard assignments, EM incorporates uncertainty by computing the expected \textit{complete-data log-likelihood} under the current posterior distribution over the latent variables. This ``soft assignment"  weights all possible latent values by their posterior probabilities, allowing the algorithm to account for the full range of plausible explanations for each observation. Consequently, EM typically yields more robust and accurate parameter estimates than simple alternating maximization.

\index{Alternating maximization}
\index{Latent variable models}
\index{EM algorithm}
\begin{algorithm}[h!] 
\caption{Alternating Maximization (AM) for Latent Variable Models}
\label{alg:am_lvm}
\begin{algorithmic}[1] 
\Require Observed data points $\{\bx_1, \bx_2, \ldots, \bx_N\}$;
\State \textbf{initialize:} $\btheta^\topone$; 
\State Choose the maximal number of iterations $C$;
\State $t=0$; \Comment{Count for the number of iterations}
\While{$t<C$} 
\State $t=t+1$;
\For{$n=1,2,\ldots, N$}
\State Step-1: $\widehatbz_n^\toptzero = \mathop{\argmax}_{\bz\in\sZ} p(\bz \mid \bx_n, \btheta^\toptzero)$;
\Comment{(AM$_1$)}
\EndFor
\State Step-2: $\btheta^\toptone =\mathop{\argmax}_{\btheta\in\bTheta}\, \ell(\btheta; \{(\bx_n,\widehatbz_n^\toptzero)\}_{n=1}^N )$;
\Comment{(AM$_2$)}
\EndWhile
\State Output $\btheta^\toptzero$;
\end{algorithmic} 
\end{algorithm}

\index{Evidence lower-bound}
\index{Variational free-energy}
\index{ELBO}
\index{VFE}
\index{KL divergence}
\index{EM algorithm}
\index{Model evidence}
\subsection{ELBO and VFE.}\label{section:elbo_cfe}

Similar to alternating maximization for latent variable models, the EM algorithm also seeks to increase the marginal likelihood over successive iterations. 
However, directly maximizing the marginal likelihood (also known as \textit{model evidence}; see Equation~\eqref{equation:lvm_maxmarg}) is often intractable or computationally prohibitive. To circumvent this, the EM algorithm introduces a \textbf{proxy function} that serves as a tractable lower bound on the marginal likelihood. 
The core idea of EM is to construct such a proxy---known as the \textit{Q-function}---which lower-bounds the log-marginal likelihood and transforms the original parameter estimation problem into a bivariate optimization over both the model parameters $\btheta$ and an auxiliary distribution over the latent variables.

To develop the EM algorithm, we assume---without loss of generality---that the latent variables have continuous support. 
Instead of maximizing the product of the likelihood  functions directly, which can be cumbersome, we take the logarithm and convert the product into a sum.
The logarithm of the {marginal likelihood} function $\ell(\btheta)\equiv \ell(\btheta; \{\bx_n\})$ in \eqref{equation:lvm_maxmarg} is given by:
\begin{equation}\label{equation:elbo_ineq}
\begin{aligned}
\mathcalL\left(\btheta\right) 
&\triangleq\ln \ell(\btheta)
=\sum_{n=1}^N \ln \int p(\bz_n, \bx_n\mid \btheta) \,d\bz_n
=\sum_{n=1}^N \ln \int q_{\bz_n}(\bz_n) \frac{p(\bz_n, \bx_n\mid \btheta)}{ q_{\bz_n}(\bz_n)} \,d\bz_n\\
&\geq \sum_{n=1}^N  \int q_{\bz_n}(\bz_n) \ln\frac{p(\bz_n, \bx_n\mid \btheta)}{ q_{\bz_n}(\bz_n)} \,d\bz_n\\
&= \underbrace{\sum_{n=1}^N  \int q_{\bz_n}(\bz_n) \ln p(\bz_n, \bx_n\mid \btheta) \,d\bz_n}_{\text{explain data/expected energy}} \underbrace{- \sum_{n=1}^N  \int q_{\bz_n}(\bz_n) \ln q_{\bz_n}(\bz_n) \,d\bz_n}_{\mathrm{entropy}}\\
&\triangleq 
\mathcalF\left(q_{\bz_1}(\bz_1), q_{\bz_2}(\bz_2), \ldots, q_{\bz_N}(\bz_N), \btheta\right)
\equiv \mathcalF\left(\{q_{\bz_n}(\bz_n)\}_{n=1}^N, \btheta\right),
\end{aligned}
\end{equation}
where the inequality follows from the Jensen's inequality, leveraging the concavity of the logarithm.
This quantity $\mathcalF\left(\{q_{\bz_n}(\bz_n)\}_{n=1}^N, \btheta\right)$ is known as the \textit{evidence lower-bound (ELBO, a.k.a., the marginal log-likelihood lower-bound, or the variational lower-bound)}.
The term $q_{\bz_n}(\bz_n)$ is an \textbf{arbitrary} probability distribution for  latent variables, which is called the \textit{variational distribution}.
As a hindsight (in the context of constrained EM optimization or variational inference), we want to maximize the ELBO while holding $\btheta$ fixed. 
The first term in the ELBO  rewards variational distributions $q_{\bz_n}(\bz_n)$ to assign high probability to configurations of the latent variables that well-explain the observations (i.e., \textit{expected energy});  the second term---the \textit{entropy}---encourages uncertainty by spreading probability mass across many configurations, preventing overconfident approximations.

\begin{SCfigure}
\centering
\includegraphics[width=0.6\textwidth]{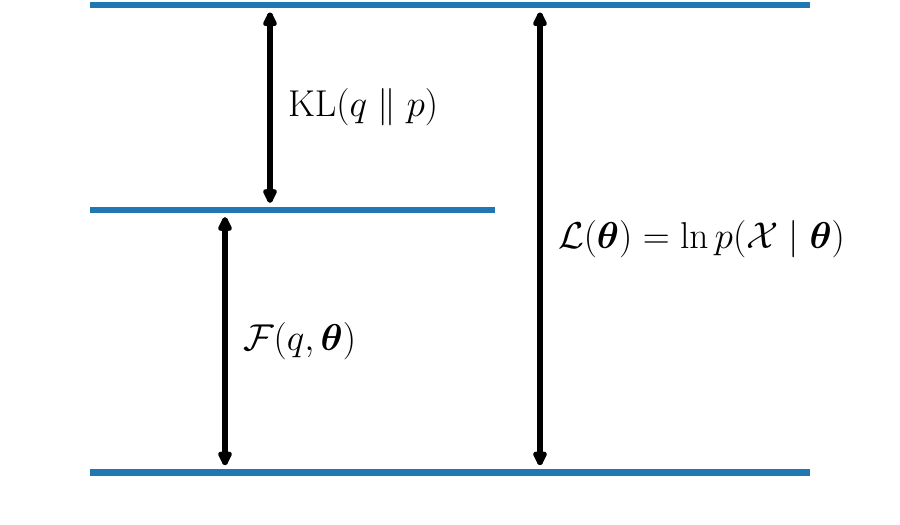}
\caption{Illustration of the ELBO decomposition
given by \eqref{equation:elbo_vfe_neg} or \eqref{equation:elbo_vfe_negv2}, which holds for
any choice of distribution $q(\bz)$.
Because  $\KL[q\parallel p] \geq 0$,
the quantity $\mathcalF(q, \btheta)$ is a lower bound on the log-marginal likelihood
function $\mathcalL(\btheta)=\ln p(\mathcalX\mid \btheta)$.}
\label{fig:ELBO_decom}
\end{SCfigure}

\paragrapharrow{ELBO decomposition.}
The ELBO admits an alternative interpretation through its relationship with the Kullback--Leibler  divergence. Specifically,
\begin{subequations}
\begin{equation}\label{equation:elbo_vfe_neg}
\begin{aligned}
\mathcalF\big(\{q_{\bz_n}&(\bz_n)\}_{n=1}^N, \btheta\big) 
= \sum_{n=1}^N  \int q_{\bz_n}(\bz_n) \ln\frac{p(\bz_n, \bx_n\mid \btheta)}{ q_{\bz_n}(\bz_n)} \,d\bz_n\\
&=\sum_{n=1}^N  \int q_{\bz_n}(\bz_n) \ln p(\bx_n \mid \btheta) \,d\bz_n  + \sum_{n=1}^N\int q_{\bz_n}(\bz_n) \ln\frac{p(\bz_n \mid \bx_n, \btheta)}{ q_{\bz_n}(\bz_n)} \,d\bz_n \\
&=\sum_{n=1}^N\ln p(\bx_n \mid \btheta) - \sum_{n=1}^N\KL\left[ q_{\bz_n}(\bz_n) \parallel p(\bz_n \mid \bx_n, \btheta)\right]\\
&\triangleq
-\mathcalL_{\text{VFE}}\left(\{q_{\bz_n}(\bz_n)\}_{n=1}^N, \btheta\right),
\end{aligned}
\end{equation}
where $\KL[P \parallel Q] \triangleq \int P(x) \ln \left( \frac{P(x)}{Q(x)} \right) dx \geq 0$ denotes the \textit{Kullback--Leibler (KL) divergence} between $P$ and $Q$ and the equality is obtained only when $P=Q$, and $\mathcalL_{\text{VFE}}\big(\left\{q_{\bz_n}(\bz_n)\right\}_{n=1}^N, \btheta\big)$ is the \textit{variational free-energy (VFE)}, a terminology from statistical  physics \citep{neal1998view}. The VFE is  the negative entropy (see Problem~\ref{problem:entropy_mgau}) of $q_{\bz}(\bz)$ minus the expected energy under $q_{\bz}(\bz)$:
$$
\mathcalL_{\text{VFE}}\left(\left\{q_{\bz_n}(\bz_n)\right\}_{n=1}^N, \btheta\right) = 
\underbrace{\sum_{n=1}^N  \int q_{\bz_n}(\bz_n)\ln{ q_{\bz_n}(\bz_n)} \,d\bz_n}_{\text{negative entropy}} 
- \underbrace{\sum_{n=1}^N  \int q_{\bz_n}(\bz_n) \ln{p(\bz_n, \bx_n\mid \btheta)}\,d\bz_n}_{\text{expected energy}}.
$$
The above derivation also shows that maximizing the ELBO is equivalent to minimizing the KL divergence between the variational distribution $q_{\bz_n}(\bz_n)$ and the true  hidden variable posterior $p(\bz_n \mid \bx_n, \btheta)$~\footnote{In deep learning community, the quantity $p(\bz_n \mid \bx_n, \btheta)$ is always known as the \textit{recognition network} or \textit{inference network} that can be characterized by neural networks \citep{hinton1995wake}. So we find the quantity $q_{\bz_n}(\bz_n)$ that acts like the recognition network.}.
On the other hand, Equation~\eqref{equation:elbo_vfe_neg} also shows that the logarithm of the marginal likelihood $\mathcalL(\btheta)$ is also the sum of the ELBO and the KL divergence $\sum_{n=1}^N\KL\left[ q_{\bz_n}(\bz_n) \parallel p(\bz_n \mid \bx_n, \btheta)\right]$:
\begin{equation}\label{equation:elbo_vfe_negv2}
\ln p(\mathcalX \mid \btheta) 
= 
\mathcalF\big(\{q_{\bz_n}(\bz_n)\}_{n=1}^N, \btheta\big) 
+\sum_{n=1}^N\KL\left[ q_{\bz_n}(\bz_n) \parallel p(\bz_n \mid \bx_n, \btheta)\right],
\end{equation}
where $\mathcalX=\{\bx_1,\bx_2,\ldots,\bx_N\}$ denotes the collection of data samples.
\end{subequations}
This, again, confirms that the ELBO is a lower bound of $\mathcalL(\btheta)$ since the KL divergence is nonnegative. And the KL distance measures the tightness of the lower bound, i.e., the gap between the ELBO and the log-marginal likelihood. 
See Figure~\ref{fig:ELBO_decom} for an illustration.

\paragrapharrow{KL divergence behavior.} 
The direction of the KL divergence has important implications for the nature of the approximation:
\begin{itemize}
\item Minimizing  $\KL[q\parallel p]$ (the \textit{reverse} or \textit{exclusive KL}) leads to \textit{mode-seeking} or \textit{zero-forcing} behavior: $q$ is forced to zero wherever $p$ is zero, often concentrating its mass around a single mode of $p$ or one of the modes of $p$.
\item In contrast, minimizing $\KL[p\parallel q]$ (the \textit{forward} or \textit{inclusive KL}) results in \textit{mass-covering} or \textit{mean-seeking} behavior: $q$ must assign non-negligible probability wherever $p$ has support, yielding broader, more inclusive approximations (see Problem~\ref{problem:forward_rever_KL}).
\end{itemize}

\index{EM algorithm}
\subsection{EM for Unconstrained Optimization}\label{section:em_uncons}
We now introduce the EM algorithm for latent variable models.
Although the general formulation is due to \citet{dempster1977maximum}, the core idea had already appeared in earlier works for specific problems \citep{baum1970maximization}.
Like alternating maximization, the EM algorithm alternates between two steps: 
the \textit{E-step}, where it computes the posterior distribution of the latent variables given the current parameter estimate, and the \textit{M-step}, where it updates the model parameters (or global latent variables) by maximizing a surrogate objective derived from the E-step. 
Using the ELBO derived in \eqref{equation:elbo_ineq}, and denoting the estimates $q_{\bz}(\bz)$ and $\btheta$  at iteration $t$  as $q^\toptzero_{\bz}(\bz)$ and $\btheta^\toptzero$, respectively, the updates of E/M steps are 
\begin{tcolorbox}[colback=white,colframe=black]
\begin{minipage}{1\textwidth}
$$
\begin{aligned}
&\textbf{E-Step: } \gap q^\toptone_{\bz_n}(\bz_n) &\leftarrow& \mathop{\argmax}_{q_{\bz_n}} \mathcalF\left(\left\{q_{\bz_n}(\bz_n)\right\}_{n=1}^N, \btheta^{\textcolor{mylightbluetext}{(t)}}\right),\gap \forall n\in\{1,2,\ldots,N\};\\
&\textbf{M-Step: } \gap \btheta^\toptone &\leftarrow& \mathop{\argmax}_{\btheta\in\bTheta} \mathcalF\left(\left\{q^{\textcolor{mylightbluetext}{(t+1)}}_{\bz_n}(\bz_n)\right\}_{n=1}^N, \btheta\right).
\end{aligned}
$$
\end{minipage}
\end{tcolorbox}
\noindent
Thus, the EM algorithm can be interpreted as follows: the E-step approximates the conditional distributions of the local latent variables  $\{\bz_n\}$,while the M-step updates the global latent variable $\btheta$.
These two steps are repeated until convergence of the sequence $\{\btheta^\toptzero\}$.
Under mild regularity conditions, convergence to a local maximum of the marginal likelihood is guaranteed \citep{gupta2011theory, jain2017non}. 
Because the log-likelihood may have multiple local maxima, it is common practice to run EM multiple times with different initializations  $\btheta^\topone$.
The final estimate is then chosen as the solution yielding the highest likelihood among all runs.

\index{Lagrangian function}
\index{Functional derivatives}
\index{Variational derivatives}
\paragrapharrow{E-Step.} 
To derive the E-step, we maximize the ELBO $\mathcalF\big(\left\{q_{\bz_n}(\bz_n)\right\}_{n=1}^N, \btheta^\toptzero\big)$ with respect to each $q_{\bz_n}(\bz_n)$ subject to the normalization constraint $\int q_{\bz_n}(\bz_n) \,d\bz_n =1$ for all $n\in\{1,2,\ldots,N\}$.
Introducing Lagrange multipliers $\{\gamma_n\}_{n=1}^N$,  the associated Lagrangian is (see, for example, \citet{boyd2004convex}):
$$
L\left(q_{\bz_n}(\bz_n), \bgamma\right)=
\mathcalF\left(\left\{q_{\bz_n}(\bz_n)\right\}_{n=1}^N, \btheta^\toptzero\right) + \sum_{n} \gamma_n \left(\int q_{\bz_n}(\bz_n) \,d\bz_n -1\right).
$$
Taking the functional (or variational) derivative of $L$ with respect to  $q_{\bz_n}(\bz_n)$ and setting it to zero yields:
\begin{equation}\label{equation:qzi_fllai}
\int \ln \frac{p(\bz_n, \bx_n \mid \btheta^\toptzero)}{q_{\bz_n}(\bz_n)} \,d\bz_n - 1+\gamma_n = 0
\implies 
q^\toptone_{\bz_n}(\bz_n) = \exp(\gamma_n-1) p(\bz_n, \bx_n\mid \btheta^\toptzero),\, \forall n.
\end{equation}
Plugging in the expressions into the Lagrangian function $L\left(q_{\bz_n}(\bz_n), \blambda\right)$ and taking maximum obtains
$$
\int \exp(\gamma_n-1)p(\bz_n, \bx_n\mid \btheta^\toptzero) \,d\bz_n=0 
\quad\implies \quad
\gamma_n = 1-\ln \int p(\bz_n, \bx_n\mid \btheta^\toptzero) \,d\bz_n.
$$
Substituting $\gamma_n$ into \eqref{equation:qzi_fllai}, we have 
$$
q^\toptone_{\bz_n}(\bz_n) = p(\bz_n \mid \bx_n, \btheta^\toptzero), \, \forall n.
$$
Hence, the optimal variational distribution in the E-step is precisely the true posterior of the latent variables given the current parameter estimate $\btheta=\btheta^\toptzero$.
Substituting this posterior  $q^\toptone_{\bz_n}(\bz_n)$ back into the ELBO $\mathcalF\big(\left\{q_{\bz_n}(\bz_n)\right\}_{n=1}^N, \btheta^\toptzero\big)$ shows that the ELBO is a \textbf{tight} bound of the logarithm of the marginal likelihood function $\mathcalL(\btheta) $, i.e., the bound becomes an equality (see Problem~\ref{problem:elbo_equa_em}): 
$$
\mathcalF\left(\left\{q_{\bz_n}(\bz_n)\right\}_{n=1}^N, \btheta^\toptzero\right)=
\mathcalL(\btheta^\toptzero).
$$
That is, the KL divergence term in the ELBO decomposition vanishes in Figure~\ref{fig:ELBO_decom} or \eqref{equation:elbo_vfe_negv2}. 

\paragrapharrow{M-Step.}
With the posterior $q^\toptone_{\bz_n}(\bz_n)= p(\bz_n \mid \bx_n, \btheta^\toptzero)$ for all $n\in\{1,2,\ldots,N\}$ fixed from the E-step, we turn to consider the global latent variable $\btheta$. Thus, the E-step in EM algorithms is considered to obtain an \textit{Q-function} in the literature \citep{gupta2011theory, jain2017non}:
$$
\begin{aligned}
&\gap \mathcalF\left(\left\{q_{\bz_n}(\bz_n) = \textcolor{mylightbluetext}{p(\bz_n \mid \bx_n, \btheta^\toptzero)}\right\}_{n=1}^N, \textcolor{mylightbluetext}{\btheta}\right)
=\sum_{n}  \int p(\bz_n \mid \bx_n, \btheta^\toptzero) \ln\frac{p(\bz_n, \bx_n\mid \btheta)}{ p(\bz_n \mid \bx_n, \btheta^\toptzero)} \,d\bz_n\\
&=\underbrace{\sum_{n}  \int p(\bz_n \mid \bx_n, \btheta^\toptzero) \ln p(\bz_n, \bx_n\mid \btheta) \,d\bz_n}_{\triangleq \text{$Q(\btheta \mid \btheta^\toptzero)$}} 
-\underbrace{\sum_{n}  \int p(\bz_n \mid \bx_n, \btheta^\toptzero)\ln p(\bz_n \mid \bx_n, \btheta^\toptzero) \,d\bz_n}_{\text{entropy of $p(\bz_n \mid \bx_n, \btheta^\toptzero$})}.
\end{aligned}
$$
Since the denominator term (the entropy term) in the $\mathcalF\Big(\left\{q_{\bz_n}(\bz_n) = \textcolor{black}{p(\bz_n \mid \bx_n, \btheta^\toptzero)}\right\}_{n=1}^N, \textcolor{black}{\btheta}\Big)$ does not depend on $\btheta$, the M-step reduces to maximizing   the Q-function:
$$
\begin{aligned}
\textbf{M-Step: } \gap \btheta^\toptone &\leftarrow \mathop{\argmax}_{\btheta\in\bTheta} 
\sum_{n}  \int p(\bz_n \mid \bx_n, \btheta^\toptzero) \ln p(\bz_n, \bx_n\mid \btheta) \,d\bz_n\\
&=\mathop{\argmax}_{\btheta\in\bTheta} Q(\btheta \mid \btheta^\toptzero).
\end{aligned}
$$
The complete EM procedure is summarized in Algorithm~\ref{alg:em_alg}:
\begin{enumerate}
\item \textit{E-step.} Compute the posterior distributions of latent variables  $  \{\bz_n\}  $  given observations  $ \{\bx_n\} $  and current parameter estimates.
\item \textit{M-step.} Update $\btheta^\toptone$ by maximizing the expected complete-data log-likelihood or the Q-function $Q(\btheta \mid \btheta^\toptzero)$.
\end{enumerate}
Because the E-step often admits a closed-form solution (as shown above), the EM algorithm effectively reduces to iteratively constructing and optimizing the Q-function. 
The Q-function can be further compactly written as 
$$
Q(\btheta \mid \btheta^\toptzero) = \sum_{n=1}^{N} \Exp_{\underbrace{\bz\sim p(\cdot \mid \bx_n, \btheta^\toptzero)}_{\text{conditional prob.}}}  \underbrace{\big[\ln p(\bz_n, \bx_n \mid \btheta)\big]}_{\text{log-joint prob.}},
$$
such that the expectation is taken over the conditional probability of hidden variables (from last iteration) for the log-joint probability for the observed and hidden variables.

\begin{algorithm}[h] 
\caption{Expectation-Maximization (EM) Algorithm}
\label{alg:em_alg}
\begin{algorithmic}[1] 
\Require Observed data points $\{\bx_1, \bx_2, \ldots, \bx_N\}$;
\State \textbf{initialize:} $\btheta^\topone$; 
\State Choose the maximal number of iterations $C$;
\State $t=0$; \Comment{Count for the number of iterations}
\While{$t<C$} 
\State $t=t+1$;
\State E-step: $q^\toptone_{\bz_n}(\bz_n) = p(\bz_n \mid \bx_n, \btheta^\toptzero), \, \forall n\in\{1,2,\ldots,N\}$;
\State M-step: $\btheta^\toptone \leftarrow \mathop{\argmax}_{\btheta\in\bTheta} Q(\btheta \mid \btheta^\toptzero)$;
\EndWhile
\State Output $\btheta^\toptzero$;
\end{algorithmic} 
\end{algorithm}

The Q-function possesses all the desirable properties of an effective proxy objective. 
Specifically, any parameter update that increases $Q(\btheta \mid \btheta^\toptzero)$  is guaranteed to increase the log-marginal likelihood $\mathcalL\left(\btheta\right)$. 
Moreover, for many important models---such as \textit{Gaussian mixture models} and \textit{mixed regression} (see Problems~\ref{prob:mix_of_gauss}--\ref{prob:mix_of_bern})---the Q-function can be both constructed and optimized efficiently in closed form \citep{jain2017non}.

\paragrapharrow{Soft assignment.}
On the other hand, for each point $\bx_n$, we can define a \textit{point-wise Q-function}:
$$
\begin{aligned}
Q(\btheta \mid \btheta^\toptzero) 
&=
\sum_{n}  \underbrace{\int p(\bz_n \mid \bx_n, \btheta^\toptzero) \ln p(\bz_n, \bx_n\mid \btheta) \,d\bz_n}_{\text{$\bQ_{\bx_n}(\btheta \mid \btheta^\toptzero)$}}
&=\sum_{n} \bQ_{\bx_n}(\btheta \mid \btheta^\toptzero),
\end{aligned}
$$
where, letting $\bw_{\bz_n}\triangleq p(\bz_n \mid \bx_n, \btheta^\toptzero)$, we have
\begin{equation}
\bQ_{\bx_n}(\btheta \mid \btheta^\toptzero) 
=
\int \bw_{\bz_n} \ln p(\bz_n, \bx_n\mid \btheta) \,d\bz_n.
\end{equation}
Therefore, upon a moment of reflexion, the key advantage of the Q-function becomes clear: 
rather than assigning each observation $\bx_n$ to a single latent configuration $\widehatbz_n^\toptzero = \mathop{\argmax}_{\bz\in\sZ} p(\bz \mid \bx_n, \btheta^\toptzero)$ (a hard assignment), 
EM uses the full posterior distribution $p(\bz_n \mid \bx_n, \btheta^\toptzero)$ to weight all possible latent values (a soft assignment).
In contrast, the alternating maximization algorithm (Algorithm~\ref{alg:am_lvm}) performs a hard assignment in Step 1: it selects the most likely latent value $\widehatbz_n$ under the current posterior  $p(\bz \mid \bx_n, \btheta^\toptzero)$ and uses only this value to update $\btheta$ in the next iteration. 
The EM algorithm, by contrast, leverages the entire posterior distribution, leading to more robust and stable updates---especially when the posterior is multimodal or uncertain.

\index{Stochastic EM}
\paragrapharrow{Stochastic EM algorithm.} 
In large-scale settings, computing the full Q-function over all $N$ data points at each iteration can be computationally prohibitive. Instead, one can apply stochastic EM, which approximates the M-step using a single randomly sampled data point (or a mini-batch). At iteration $t$, after sampling an index $n$, we perform the update:
$$
\btheta^\toptone \leftarrow \btheta^\toptzero + \eta_t \nabla Q_{\bx_n}(\btheta^\toptzero \mid \btheta^\toptzero),
$$
where $\eta_t>0$ is a step size, and $\nabla Q_{\bx_n}(\btheta^\toptzero \mid \btheta^\toptzero)$ is an \textit{ascent direction}: taking a positive step along this gradient, on average, increases the expected log-likelihood \citep{lu2022gradient}.

\index{MAP EM}
\paragrapharrow{Maximum a posteriori (MAP EM).} More generally, we can incorporate prior knowledge by placing a prior distribution $p(\btheta)$ on the parameters and maximizing the \textit{log-posterior} instead of the log-marginal likelihood:
\begin{subequations}\label{equation:mapem_all}
\begin{equation}\label{equation:mapem1}
\btheta_{\mathrm{MAP}}^* =
\mathop{\argmax}_{\btheta\in\bTheta}
\sum_{n}\ln \int p(\bz_n, \bx_n\mid \btheta) \,d\bz_n +\ln p(\btheta).
\end{equation}
The E-step remains unchanged, since the prior does not depend on the latent variables. The M-step, however, now includes the log-prior:
\begin{equation}\label{equation:mapem2}
\textbf{MAP M-Step:}\gap  \mathop{\argmax}_{\btheta\in\bTheta} Q(\btheta \mid \btheta^\toptzero)+\ln p(\btheta).
\end{equation}
\end{subequations}

We now illustrate the EM algorithm with a classic and widely used example: the Gaussian mixture model (GMM).\footnote{See, for example, \citet{jain2017non, lu2021survey} for further details.} 
GMMs are commonly applied in clustering, density estimation, and topic modeling (where they help uncover latent topics in document collections).
\begin{example}[Gaussian Mixture Model (GMM)\index{Gaussian mixture model}\index{Mixture of Gaussians}]\label{example:gmm_twoclus}
Consider a Gaussian mixture model (mixture of Gaussians) with two modes or clusters. 
The model parameters are $\btheta=\{\pi_k, \bmu_k, \bSigma_k\}_{k\in\{0,1\}}$, and the joint density is
$$
f_{\btheta} \left( \cdot,\cdot \mid \btheta \right) = \pi_0\cdot \normal_0 + \pi_1 \cdot \normal_1,
$$
where $\normal_k=\normal(\cdot \mid \bmu_k, \bSigma_k), k=\{0,1\}$ denotes a multivariate Gaussian density (Section~\ref{sec:multi_gaussian_conjugate_prior}), and $\pi_0+\pi_1=1$ are the mixture coefficients.
For simplicity, assume equal mixing proportions ($\pi_0=\pi_1=1/2$)  and shared isotropic covariance matricex ($\bSigma_0=\bSigma_1=\bI$). 
We then aim to estimate only the means: $\btheta=\{\bmu_0, \bmu_1\}$. 
Each observed data point  $\bx_n$ ($n=1,2,\ldots,N$) is associated with a discrete latent variable $z_n\in\{0,1\}$ indicating its component membership: $\normal_0$ or $\normal_1$. 
The joint and posterior distributions are:
$$
\begin{aligned}
p(\bx_n, z_n \mid \btheta) &= \normal_{z_n}(\bx_n \mid \bmu_{z_n}, \bSigma_{z_n}),\\
p(z_n \mid \bx_n, \btheta) &= \frac{p(\bx_n, z_n \mid \btheta)}{p(\bx_n \mid \btheta)}
=
\frac{p(\bx_n, z_n \mid \btheta)}{p(\bx_n,z_n \mid \btheta)+p(\bx_n,1-z_n \mid \btheta)}, \gap n\in\{1,2,\ldots,N\}.
\end{aligned}
$$
The Q-function is then constructed as:
$$
\begin{aligned}
Q(\btheta \mid \btheta^\toptzero) 
&=\sum_{n} \bQ_{\bx_n}(\btheta \mid \btheta^\toptzero)
=
\sum_{n}  \underbrace{\int p(z_n \mid \bx_n, \btheta^\toptzero) \ln p(z_n, \bx_n\mid \btheta) dz_n}_{\text{$\bQ_{\bx_n}(\btheta \mid \btheta^\toptzero)$}}\\
&=\sum_{n}  \left\{p(0 \mid \bx_n, \btheta^\toptzero) \ln p(0, \bx_n\mid \btheta) + p(1 \mid \bx_n, \btheta^\toptzero) \ln p(1, \bx_n\mid \btheta)\right\}
\end{aligned}
$$
This Q-function can be maximized analytically in the M-step, yielding closed-form updates for $\bmu_0$ and $\bmu_1$.
The resulting EM algorithm alternates between computing soft cluster responsibilities (E-step) and updating cluster means as weighted averages of the data (M-step).
Generalizations of this example---including mixtures with arbitrary numbers of components, full covariance matrices, or other distributions---are discussed in Problems~\ref{prob:mix_of_gauss}--\ref{prob:mix_of_bern}.
\end{example}

\index{Variational EM algorithm}
\subsection{EM for Constrained Optimization}\label{section:em_constrained}
In the standard (unconstrained) EM algorithm, we obtain the exact posterior in the E-step: $q^\toptone_{\bz_n}(\bz_n) = p(\bz_n \mid \bx_n, \btheta^\toptzero), \, \forall n$.
However, the real problem can be far frustrating such that the data are explained by multiple  interacting hidden variables and the hidden variable posterior $p(\bz_n \mid \bx_n, \btheta^\toptzero),  \, \forall n$ is intractable \citep{williams1991mean, ghahramani1995factorial, beal2003variational, turnertwo}. 
To address this, the variational Bayesian approach restricts the posterior to a tractable family of distributions. This leads to what is known as the \textit{constrained EM algorithm} or \textit{variational EM algorithm}.

Suppose we constrain the approximate posterior to a parametric family $\mathcalQ=\{q_{\blambda} = q_{\bz}(\bz \mid \blambda)\mid  \blambda\in{\Lambda}\}$, where $\blambda$ is called the \textit{variational parameter}. In practice, the variational distribution for each latent variable $\bz_n$ may depend on its corresponding observation  $\bx_n$, so we write $q_{\bz_n}(\bz_n \mid \blambda_n)$ with $\blambda_n=\{\bx_n, \widebarblambda_n\}$, where  only $\widebarblambda_n$ is updated during inference.
At iteration $t$, the E-step becomes (using Equation~\eqref{equation:elbo_vfe_neg}, which shows that maximizing the ELBO is equivalent to minimizing the KL divergence):
$$
\begin{aligned}
\textbf{E-Step: } \gap q^\toptone_{\bz_n}(\bz_n \mid \blambda_n) 
\leftarrow 
&\mathop{\argmax}_{\blambda \in\Lambda} \mathcalF\left(\left\{q_{\bz_n}(\bz_n \mid \blambda_n)\right\}_{n=1}^N, \btheta^\toptzero\right),\gap \forall n\in\{1,2,\ldots,N\}\\
=&\mathop{\argmin }_{\blambda \in\Lambda} \sum_{n}\KL\left[ q_{\bz_n}(\bz_n \mid \blambda_n) \parallel p(\bz_n \mid \bx_n, \btheta^\toptzero)\right].
\end{aligned}
$$
In other words, we seek the  \textit{variational posterior} $q^\toptone_{\bz_n}(\cdot \mid \blambda_n) $ that is closest---in KL divergence---to the true (but intractable) posterior $p(\bz_n \mid \bx_n, \btheta^\toptzero)$.
Because the variational family $\mathcalQ$ may not contain the true posterior, the ELBO is generally no longer tight: the gap between the ELBO and the log marginal likelihood remains positive.

The M-step proceeds as usual, but now uses the variational posterior instead of the exact one:
$$
\begin{aligned}
\textbf{M-Step: } \gap \btheta^\toptone &\leftarrow \mathop{\argmax}_{\btheta\in\bTheta} 
\sum_{n}  \int q^\toptone_{\bz_n}(\bz_n \mid \blambda_n)  \ln p(\bz_n, \bx_n\mid \btheta) \,d\bz_n
\end{aligned}
$$

One might wonder how we can minimize the KL divergence when the true posterior $p(\bz_n \mid \bx_n, \btheta^\toptzero)$  is intractable. The key insight is that if the variational family $q_{\bz}(\bz \mid \blambda)$ has tractable moments---for instance, if it is Gaussian, fully characterized by its first and second moments---then the optimization can be carried out using only expectations under the true posterior
(e.g., $\Exp[\bz_n \mid \bx_n, \btheta^\toptzero]$ and $\Exp[\bz_n\bz_n^\top \mid \bx_n, \btheta^\toptzero]$). 
We will explore this in more detail shortly.

\index{Mean-field approximation}
\index{Lagrangian function}
\index{Functional derivatives}
\index{Variational derivatives}
\subsubsection*{Mean-Field Approximation of Hidden Variables}
Let $\mathcalX=\mathcalX(\bx_{1:N})=\{\bx_1,\bx_2,\ldots,\bx_N\}$ denote the observed data, with $\bx_n\in\real^D$, and let $\mathcalZ=\mathcalZ(\bz_{1:N})=\{\bz_1,\bz_2,\ldots,\bz_N\}$ be the corresponding latent variables, where $\bz_n\in\real^Q$.
The \textit{mean-field approximation (MFA)} assumes full factorization of the variational posterior across all latent dimensions:
$$
q_{\bz_n}(\bz_n) =\prod_{q=1}^{Q} q_{z_{nq}}(z_{nq}), 
\quad \,\forall n.
$$
Under this assumption, the ELBO becomes:
$$
\begin{aligned}
\mathcalF(\{q_{\bz_n}(\bz_n)\}, \btheta) 
&= \sum_{n=1}^N  \int q_{\bz_n}(\bz_n) \ln\frac{p(\bz_n, \bx_n| \btheta)}{ q_{\bz_n}(\bz_n)} \,d\bz_n
= \sum_{n=1}^N  \int \prod_{q=1}^{Q} q_{z_{nq}}(z_{nq}) \ln\frac{p(\bz_n, \bx_n| \btheta)}{ \prod_{q=1}^{Q} q_{z_{nq}}(z_{nq})} \,d\bz_n\\
&= \sum_{n=1}^N  \int \prod_{q=1}^{Q} q_{z_{nq}}(z_{nq}) \ln p(\bz_n, \bx_n\mid \btheta) -  \sum_{q=1}^{Q} q_{z_{nq}}(z_{nq})\ln q_{z_{nq}}(z_{nq})  \,d\bz_n.
\end{aligned}
$$

To derive the E-step under the mean-field assumption at iteration $t$, we maximize the ELBO $\mathcalF\left(\left\{q_{\bz_n}(\bz_n)\right\}_{n=1}^N, \btheta^\toptzero\right)$with respect to each factor  $q_{z_{nq}}(z_{nq})$, subject to the normalization constraints: $\int q_{z_{nq}}(z_{nq}) \,dz_{nq} =1$ for all $n\in\{1,2,\ldots,N\}, q\in\{1,2,\ldots, Q\}$.
Introducing Lagrange multipliers $\{\gamma_{nq}\}$, the Lagrangian is:
$$
L\left(q_{z_{nq}}(z_{nq}), \bgamma\right)=
\mathcalF\left(\left\{q_{\bz_n}(\bz_n)\right\}_{n=1}^N, \btheta^\toptzero\right) + \sum_{n,q} \gamma_{nq} \left(\int q_{z_{nq}}(z_{nq}) dz_{nq} -1\right).
$$
Setting the functional  derivative  of the Lagrangian function with respect to $q_{z_{nq}}(z_{nq})$ to zero yields
\begin{equation}\label{equation:em_uncon_mf}
\begin{aligned}
&\ln q^\toptone_{z_{nq}}(z_{nq}) = \int \left[ \prod_{k\neq q}^{Q}q_{z_{nk}}(z_{nk}) \ln p(\bz_n, \bx_n \mid\btheta^\toptzero)\right]d\bz_{n/q} +\gamma_{nq}-1\\
&\implies q^\toptone_{z_{nq}}(v) = \frac{1}{\mathcalC_{nq}}\exp\left\{\int \left[ \prod_{k\neq q}^{Q}q_{z_{nk}}(z_{nk}) \ln p(\bz_n, \bx_n \mid\btheta^\toptzero)\right]d\bz_{n/q}\right\},
\end{aligned}
\end{equation}
where $\mathcalC_{nq}$ is the normalization constant, $d{\bz_{n/q}}$ denotes the element of integration for all elements in $\bz_n$ except $z_{nq}$, 
$\prod_{k\neq q}^{Q}$ denotes the product of all elements except the $q$-th item.
More compactly, this update can be written as:
\begin{equation}\label{equation:em_uncon_mf_comp}
q^\toptone_{z_{nq}}(z_{nq}) 
\leftarrow
\frac{1}{\mathcalC_{nq}}\exp\left\{\Exp_{q(-z_{nq})} \left[ \ln p(\bz_n, \bx_n \mid\btheta^\toptzero)\right]\right\},
\end{equation}
where the expectation can be taken over all $n^\prime \in\{1,2,\ldots,N\}, q^\prime\in\{1,2,\ldots, Q\}$ except $\{n^\prime =n, q^\prime =q\}$.
Thus, the ELBO is maximized iteratively: at each step, we update one factor $q_{z_{nq}}(z_{nq})$ while holding all others fixed (i.e., keeping the remaining factors $q_{z_{n^\prime q^\prime}}(z_{n^\prime q^\prime})$ constant, where $n^\prime \in\{1,2,\ldots,N\}$, $q^\prime\in\{1,2,\ldots, Q\}$, and $\{n^\prime \neq n, q^\prime \neq q\}$). This coordinate ascent procedure is repeated until convergence.

We will illustrate this mean-field approach with a concrete example in the next section; see Example~\ref{example:bayes_lr_mfa}.
An extension of mean-field variational inference is \textit{structured variational inference} or \textit{structured VI} \citep{saul1996mean}, which allows dependencies among subsets of latent variables. While this yields a more accurate posterior approximation, it often complicates the optimization problem, making it harder to solve analytically or computationally.

\index{ELBO}
\index{Variational Bayesian inference}
\index{Variational  inference}
\subsection{Variational Bayesian Inference}\label{section:vb_inference}
In the EM algorithm, we derive the ELBO with respect to the hidden variables only. 
In contrast, \textit{variational Bayesian (VB) inference} or simply \textit{variational inference (VI)} extends this idea by introducing a variational distribution over both the hidden variables and the model parameters. (Previously, in non-Bayesian settings, we treated $\btheta$ as deterministic and did not assign it a distribution.)

In fact, the \textbf{constrained EM algorithm} (Section~\ref{section:em_constrained}) can be viewed as a special case of variational inference that focuses exclusively on approximating the posterior over the hidden variables, while keeping the parameters fixed. 
Similarly, the  \textbf{MAP EM algorithm}  (see \eqref{equation:mapem_all}) is another special case of variational inference that incorporates the prior distribution $p(\btheta)$ directly into the EM framework, under the assumption that the posterior over the hidden variables, $p(\bz_n\mid \bx_n, \btheta^\toptzero)$, is tractable.
These perspectives clarify  the relationship between EM and VI: EM is a simplified form of VI where the variational family is restricted to delta functions over $\btheta$, effectively excluding uncertainty in the parameters.

As before, let the observed data be $\mathcalX=\{\bx_1, \bx_2, \ldots, \bx_N\}$ and the hidden variables be $\mathcalZ=\{\bz_1, \bz_2, \ldots, \bz_N\}$.
We assume that the joint distribution of the random vectors  $\rvx\in \sX$ and $\rvz\in \sZ$ belongs to a parametric family $\mathcalF = \{f_{\btheta} = f(\cdot , \cdot \mid \btheta) :\btheta\in\bTheta\} $, meaning $(\rvx, \rvz) \sim f_{\btheta^*}$ for some true parameter $\btheta^*\in \bTheta$.

Under the Bayesian framework, however, the parameters $\btheta$ themselves are treated as random variables. Before observing the data, we encode our beliefs about $\btheta$ through a prior distribution 
$p(\btheta\mid \balpha)$, where $\balpha$ denotes hyper-parameters. The full joint distribution then becomes:
$$
p(\btheta, \mathcalX, \mathcalZ \mid \balpha) = p(\btheta \mid \balpha) p(\mathcalX, \mathcalZ \mid \btheta, \balpha),
$$
where, under the i.i.d. assumption, $p(\mathcalX, \mathcalZ \mid \btheta, \balpha)=\prod_{i}^{N} p(\bx_n, \bz_n \mid \btheta, \balpha)$ (See the graphical model in Figure~\ref{fig:lvm_hyper}.)

The marginal likelihood (or model evidence) $p(\mathcalX \mid \balpha)$ can then be lower-bounded using a variational distribution $q(\mathcalZ,\btheta)$ via the ELBO:
$$
\begin{aligned}
\ln p(\mathcalX \mid \balpha) 
&= \ln \int p(\mathcalZ, \mathcalX, \btheta \mid \balpha) \,d\btheta \,d\mathcalZ 
= \ln \int q(\mathcalZ, \btheta) \frac{p(\mathcalZ, \mathcalX, \btheta \mid \balpha)}{q(\mathcalZ, \btheta)} \,d\btheta \,d\mathcalZ \\
&\geq  \int q(\mathcalZ, \btheta) \ln \frac{p(\mathcalZ, \mathcalX, \btheta \mid \balpha)}{q(\mathcalZ, \btheta)} \,d\btheta \,d\mathcalZ,
\end{aligned}
$$
where the inequality follows again from Jensen's inequality.
Following the E-step of the EM algorithm, when fixing $\balpha=\balpha^\toptzero$ for the $t$-th iteration, it leads to $q^\toptone(\mathcalZ, \btheta)=p(\mathcalZ, \mathcalX, \btheta \mid \balpha^\toptzero)$; this in turn, when substituted into the ELBO equation, attains the equality for the lower bound, i.e., a tight lower-bound.
However, the update does not simplify the problem since $p(\mathcalZ, \mathcalX, \btheta \mid \balpha^\toptzero)$ can be intractable due to the normalizing constant.
Therefore, we further assume the hidden variable and the model parameter factorize: \colorbox{\mdframecolor}{$q(\mathcalZ, \btheta) = q_{\bz}(\mathcalZ)q_{\btheta}(\btheta)$}. 
Thus, the \textit{Bayesian ELBO inequality} becomes:
$$
\begin{aligned}
\ln p(\mathcalX \mid \balpha) 
&\geq  \int q(\mathcalZ, \btheta) \ln \frac{p(\mathcalZ, \mathcalX, \btheta \mid \balpha)}{q(\mathcalZ, \btheta)} \,d\btheta \,d\mathcalZ
=\int q_{\bz}(\mathcalZ)q_{\btheta}(\btheta) \ln \frac{p(\mathcalZ, \mathcalX, \btheta \mid \balpha)}{q_{\bz}(\mathcalZ)q_{\btheta}(\btheta)} \,d\btheta \,d\mathcalZ\\
&= \int q_{\btheta}(\btheta) 
\left[\int q_{\bz}(\mathcalZ)  
\ln \frac{p(\mathcalX, \mathcalZ \mid \btheta, \balpha)}{q_{\bz}(\mathcalZ)}   d\mathcalZ
+\ln \frac{p(\btheta \mid \balpha) }{q_{\btheta}(\btheta)}   
\right]
d\btheta \\
&\triangleq \mathcalF_{\balpha} \left(q_{\bz}(\mathcalZ),q_{\btheta}(\btheta) \right)
\equiv \mathcalF_{\balpha} \left(\left\{q_{\bz_n}(\bz_n)\right\}_{n=1}^N,q_{\btheta}(\btheta) \right),
\end{aligned}
$$
where the last equality follows from the fact that the data are drawn i.i.d. from $f_{\btheta^*}$ for some $\btheta^*\in\bTheta$. 
Additionally, we assume the hidden variables  $q_{\bz}(\mathcalZ)$ factorize across data points:
\colorbox{\mdframecolor}{$q_{\bz}(\mathcalZ) = \prod_{n=1}^{N}q_{\bz_n}(\bz_n)$}. 
Under these assumptions, the Bayesian ELBO under $\balpha$ becomes:
$$
\mathcalF_{\balpha} \left(\left\{q_{\bz_n}(\bz_n)\right\}_{n=1}^N,q_{\btheta}(\btheta) \right)
= \int q_{\btheta}(\btheta) 
\bigg[
\underbrace{\sum_{n}\int q_{\bz_n}(\bz_n)   \ln \frac{p(\bx_n, \bz_n \mid \btheta, \textcolor{mylightbluetext}{\balpha})}{q_{\bz_n}(\bz_n)}   d\bz_n}_{\mathcalF\left(\left\{q_{\bz_n}(\bz_n)\right\}_{n=1}^N,q_{\btheta}(\btheta) \right)}
+\ln \frac{p(\btheta \mid \balpha) }{q_{\btheta}(\btheta)}   
\bigg]
\,d\btheta ,
$$
where we notice that the Bayesian ELBO $\mathcalF_{\balpha}(\cdot)$ (i.e., ELBO for variational Bayesian inference) generalizes the standard ELBO $\mathcalF(\cdot)$ used in latent variable models:
it reduces to $\mathcalF$ when $q_{\btheta}(\btheta)$  is a point mass (i.e., non-Bayesian inference), and it explicitly depends on the hyper-parameters $\balpha$ (see Section~\ref{section:elbo_cfe}).

\paragrapharrow{Teriminology.}  We refer to the ELBO mentioned above as the \textbf{general factorized framework}. In the following sections, we will explore alternative variational frameworks tailored to different modeling assumptions.

Analogous to the EM algorithm, we can derive coordinate ascent updates that iteratively maximize $\mathcalF_{\balpha}$ with respect to each factor in the factorized model approximation:
\begin{tcolorbox}[colback=white,colframe=black]
\begin{minipage}{1\textwidth}
$$
\begin{aligned}
&\textbf{VBE-Step: } \,\,\, q^\toptone_{\bz_n}(\bz_n) &\leftarrow& \mathop{\argmax}_{q_{\bz_n}} \mathcalF\left(\left\{q_{\bz_n}(\bz_n)\right\}_{n=1}^N, q^{\textcolor{mylightbluetext}{(t)}}_{\btheta}(\btheta)\right),\gap \forall n\in\{1,2,\ldots,N\};\\
&\textbf{VBM-Step: } \,\,\, q^\toptone_{\btheta}(\btheta) &\leftarrow& \mathop{\argmax}_{q_{\btheta}} \mathcalF\left(\left\{q^{\textcolor{mylightbluetext}{(t+1)}}_{\bz_n}(\bz_n)\right\}_{n=1}^N, q_{\btheta}(\btheta)\right).
\end{aligned}
$$
\end{minipage}
\end{tcolorbox}

\paragrapharrow{VBE-Step.}

\index{Functional derivatives}
\index{Lagrangian function}
To derive the variational Bayesian E-step (VBE-step), we maximize the Bayesian ELBO $\mathcalF\left(\left\{q_{\bz_n}(\bz_n)\right\}_{n=1}^N, \btheta^\toptzero\right)$ with respect to  $q_{\bz_n}(\bz_n)$, subject to the normalization constraints: $\int q_{\bz_n}(\bz_n) d\bz_n =1$ for all $n\in\{1,2,\ldots,N\}$.
The associated Lagrangian function  takes the form:
$$
L\left(q_{\bz_n}(\bz_n), \bgamma\right)=
\mathcalF\left(\left\{q_{\bz_n}(\bz_n)\right\}_{n=1}^N, q^\toptzero_{\btheta}(\btheta)\right) + \sum_{n} \gamma_n \left(\int q_{\bz_n}(\bz_n) d\bz_n -1\right).
$$
Taking the functional derivatives (a.k.a., variatioanl derivatives) of $L\left(q_{\bz_n}(\bz_n), \blambda\right)$ with respect to $q_{\bz_n}(\bz_n)$ and setting it to zero yields:
$$
\begin{aligned}
&\frac{\partial L\left(q_{\bz_n}(\bz_n), \blambda\right)}{\partial q_{\bz_n}(\bz_n)}
=
\int q_{\btheta}(\btheta) 
\left[
\ln p(\bx_n, \bz_n \mid \btheta, \balpha) - \ln q_{\bz_n}(\bz_n) +\gamma_n-1
\right]
\,d\btheta =0\\
&\implies
\ln q^\toptone_{\bz_n}(\bz_n) 
=
\int q^\toptzero_{\btheta}(\btheta)  \left[\ln p(\bx_n, \bz_n \mid \btheta, \balpha) \right] \,d\btheta
-
\underbrace{\int q^\toptzero_{\btheta}(\btheta) (1-\gamma_n)\,d\btheta}_{\ln \mathcalC^\toptone_{\bz_n}},
\end{aligned}
$$
where $\mathcalC^\toptone_{\bz_n}$ is a normalization constant with respect to $q^\toptone_{\bz_n}(\bz_n)$.
Therefore, the VBE-step is obtained by 
$$
\textbf{VBE-Step: } \gap q^\toptone_{\bz_n}(\bz_n) \leftarrow \frac{1}{\mathcalC^\toptone_{\bz_n}}\exp\left\{\int q^\toptzero_{\btheta}(\btheta)  \left[\ln p(\bx_n, \bz_n \mid \btheta, \balpha) \right] \,d\btheta\right\}, \gap \forall n.\\
$$

\paragrapharrow{VBM-Step.}
Similarly, for the VBM-step, we maximize $\mathcalF_{\balpha} \left(q^\toptone_{\bz}(\mathcalZ),q_{\btheta}(\btheta) \right)$ with respect to $q_{\btheta}(\btheta)$. 
The Lagrangian is:
$$
L\left(q_{\btheta}(\btheta), \gamma_{\btheta}\right)=
\mathcalF_{\balpha} \left(q^\toptone_{\bz}(\mathcalZ),q_{\btheta}(\btheta) \right) + \gamma_{\btheta} \left(\int q_{\btheta}(\btheta) \,d\btheta -1\right).
$$
Setting the functional derivative to zero gives:
$$
\begin{aligned}
&\frac{\partial L\left(q_{\btheta}(\btheta), \gamma_{\btheta}\right)}{\partial q_{\btheta}(\btheta)}
=\int q_{\bz}(\mathcalZ)  
\ln \frac{p(\mathcalX, \mathcalZ \mid \btheta, \balpha)}{q_{\bz}(\mathcalZ)}   d\mathcalZ
+
\ln p(\btheta \mid \balpha) - \ln q_{\btheta}(\btheta) -1+\gamma_{\btheta}=0\\
&\implies
\ln q_{\btheta}(\btheta) = \int q_{\bz}(\mathcalZ)  
\ln p(\mathcalX, \mathcalZ | \btheta, \balpha)   d\mathcalZ
+
\ln p(\btheta | \balpha)  -
\underbrace{(1-\gamma_{\btheta}) + \int q_{\bz}(\mathcalZ) \ln q_{\bz}(\mathcalZ) d\mathcalZ}_{\ln\mathcalC^\toptone_{\btheta}},
\end{aligned}
$$
where $\mathcalC^\toptone_{\btheta}$ is a normalization constant with respect to $q^\toptone_{\btheta}(\btheta) $.
Therefore, the VBM-step is 
$$
\textbf{VBM-Step: } 
\gap q^\toptone_{\btheta}(\btheta) \leftarrow 
\frac{1}{\mathcalC^\toptone_{\btheta}}
p(\btheta \mid \balpha) \exp\left\{ \int  q_{\bz}(\mathcalZ) \ln p(\mathcalX, \mathcalZ \mid \btheta, \balpha)  d\mathcalZ\right\}.
$$

\index{EM algorithm}
\index{Variational Bayesian inference}
\subsubsection*{Alternative Formulations}
We observe that when the model parameters $\btheta$ are treated as a point estimate---rather than as random variables drawn from a variational distribution  $q_{\btheta}(\btheta)$---the problem reduces to the standard latent variable model discussed in the context of the EM algorithm (Section~\ref{section:lvm}). In this case, the ELBO simplifies to the form given in Equation~\eqref{equation:elbo_ineq}.
Thus, the  EM algorithm can be viewed as a special case of variational inference, where uncertainty is modeled only over the latent variables, not the parameters.
More broadly, there exist several alternative formulations of the variational Bayesian inference problem. These offer different perspectives and can lead to more flexible or computationally efficient approaches, depending on the structure of the model and the nature of the inference task.

\index{Constrained VI}
\paragrapharrow{Constrained framework.}
In Section~\ref{section:em_constrained}, we introduced a constrained version of the EM algorithm, where each latent variable is approximated by a parametric distribution: $\bz_n \sim q(\bz_n \mid \blambda_n)$ for $ n\in\{1,2,\ldots,N\}$.
Similarly, in the full variational Bayesian setting, we can impose the same constraint on the latent variables while still maintaining a distribution $q_{\btheta}(\btheta)$ over the model parameters. The corresponding graphical model is shown in Figure~\ref{fig:lvm_VI}.
Under this setup, optimization of the ELBO takes the following form:
\begin{remark}[Constrained VI]\label{remark:constrained_VBi}
\begin{equation}\label{equation:constrained_VBi}
\begin{aligned}
&\gap \mathop{\argmax}_{q_{\btheta}, \blambda_n}\int \prod_{n=1}^{N}q_{\bz_n}(\bz_n\mid \blambda_n)q_{\btheta}(\btheta) \ln \frac{p(\mathcalZ, \mathcalX, \btheta \mid \balpha)}{\prod_{n=1}^{N}q_{\bz_n}(\bz_n\mid \blambda_n)q_{\btheta}(\btheta)} \,d\btheta \,d\mathcalZ\\
&\equiv\mathop{\argmax}_{q_{\btheta}, \blambda_n} \int \prod_{n=1}^{N}q_{\bz_n}(\bz_n\mid \blambda_n)q_{\btheta}(\btheta) \ln p(\mathcalZ, \mathcalX, \btheta \mid \balpha) \,d\btheta \,d\mathcalZ\\
&\gap \gap \gap\,\,-
\sum_{n=1}^{N}\int q_{\bz_n}(\bz_n\mid \blambda_n) \ln q_{\bz_n}(\bz_n\mid \blambda_n)\,d\bz_n
-
\int q_{\btheta}(\btheta) \ln q_{\btheta}(\btheta)\,d\btheta.
\end{aligned}
\end{equation}
\end{remark}

\paragrapharrow{No hidden variables, max ELBO as min of KL divergence.}
When no local latent variables are present (i.e., $\mathcalZ$ is empty), the ELBO simplifies significantly. In this case, variational inference reduces to approximating the posterior over the global parameters $\btheta$ alone:
\begin{remark}[Variational Inference Without Latent Variables]
\begin{equation}\label{equation:elbo_nohid_kl}
\begin{aligned}
\int q_{\btheta}(\btheta)  \ln \frac{p( \mathcalX, \btheta \mid \balpha)}{q_{\btheta}(\btheta) } \,d\btheta 
&=
\int q_{\btheta}(\btheta)  \ln \frac{p( \mathcalX, \btheta \mid \balpha) p(\mathcalX \mid \balpha) }{q_{\btheta}(\btheta) p(\mathcalX \mid \balpha)} \,d\btheta \\
&=\int q_{\btheta}(\btheta)  \ln \frac{p( \btheta \mid \mathcalX, \balpha)  }{q_{\btheta}(\btheta) }  \,d\btheta +\ln p(\mathcalX \mid \balpha)\\
&=-\KL\left[q_{\btheta}(\btheta) \parallel  p( \btheta \mid \mathcalX, \balpha) \right] +\ln p(\mathcalX \mid \balpha).
\end{aligned}
\end{equation}
\end{remark}
This again confirms the fundamental principle stated at the beginning of this section: variational inference seeks to minimize the KL divergence between the variational distribution and the true posterior, thereby maximizing a lower bound on the marginal likelihood.

\paragrapharrow{A unified framework.}
More generally, we can treat all unknown quantities---both global parameters and latent variables---as a single joint variable $\bomega=\{\btheta, \mathcalZ\}$. 
This leads to a unified variational inference framework that accommodates models where the number of latent variables does not correspond one-to-one with the number of observations.
For example, in topic modeling \citep{blei2003latent, blei2012probabilistic}, each document is modeled as a mixture of topics, and each topic is represented by a distribution over words in a fixed vocabulary. Here:
\begin{itemize}
\item The observed data are individual words in documents.
\item The latent variables include both the per-document topic proportions and the per-word topic assignments.
\item The total number of words (observations) and the number of topics (global latent components) are generally not equal, and the structure is hierarchical.
\end{itemize}
The goal of the model is to infer the topic structure and the distribution of topics within each document, even though the number of words and the number of topics do not align one-to-one.
Despite this complexity, the ELBO retains the same canonical form as in the parameter-only case \eqref{equation:elbo_nohid_kl}:
\begin{remark}[Unified VI]\label{remark:unified_vi}
\begin{equation}\label{equation:elbo_unified}
\begin{aligned}
\int q_{\bomega}(\bomega)  \ln \frac{p( \mathcalX, \bomega \mid \balpha)}{q_{\bomega}(\bomega) } d\bomega 
&=-\KL\left[q_{\bomega}(\bomega) \parallel  p( \bomega \mid \mathcalX, \balpha) \right] +\ln p(\mathcalX \mid \balpha).
\end{aligned}
\end{equation}
\end{remark}
This unified view emphasizes that all forms of variational inference---whether applied to parameters, local latents, or both---are instances of the same underlying optimization problem: minimizing a KL divergence to approximate a posterior distribution.

\subsubsection*{Mean-Field Approximation of Hidden Variables and Model Parameters}
To be more specific, let $\mathcalX=\{\bx_1,\bx_2,\ldots,\bx_N\}$ denote the observed data, where each $\bx_n\in\real^D$, and let $\mathcalZ=\{\bz_1,\bz_2,\ldots,\bz_N\}$ be the corresponding latent variables, with $\bz_n\in\real^Q$. The model parameters (global latent variables) are denoted by
$\btheta\in\real^P$.
Under the mean-field approximation, we assume full factorization across both the model parameters and the latent variables:
$$
q_{\btheta}(\btheta) =\prod_{p=1}^{P} q_{\theta_p}(\theta_p) 
\qquad \text{and}\qquad 
q_{\bz_n}(\bz_n) =\prod_{q=1}^{Q} q_{z_{nq}}(z_{nq}).
$$
With this factorization, the ELBO becomes:
$$
\begin{aligned}
\mathcalF_{\balpha} \left(\left\{q_{\bz_n}(\bz_n)\right\}_{n=1}^N,q_{\btheta}(\btheta) \right)
&= \int q_{\btheta}(\btheta) 
\left[
\sum_{n=1}^N\int q_{\bz_n}(\bz_n)   \ln \frac{p(\bx_n, \bz_n \mid \btheta, \balpha)}{q_{\bz_n}(\bz_n)}   \,d\bz_n
+\ln \frac{p(\btheta \mid \balpha) }{q_{\btheta}(\btheta)}   
\right]
\,d\btheta\\
&= \int \prod_{p=1}^{P} q_{\theta_p}(\theta_p) \left[\mathcalL_{\bz}(\mathcalZ,\mathcalX,\btheta,\balpha)  +\mathcalL_{\btheta}(\btheta,\balpha)  \right]  \,d\btheta,
\end{aligned}
$$
where 
$$
\begin{aligned}
\mathcalL_{\bz}(\mathcalZ,\mathcalX,\btheta,\balpha)
&\triangleq\sum_{n=1}^N  \int \prod_{q=1}^{Q} q_{z_{nq}}(z_{nq}) \ln p(\bz_n, \bx_n\mid \btheta, \balpha) -  \sum_{q=1}^{Q} q_{z_{nq}}(z_{nq})\ln q_{z_{nq}}(z_{nq})  \,d\bz_n;\\
\mathcalL_{\btheta}(\btheta,\balpha)&\triangleq \ln p(\btheta \mid \balpha) - \sum_{p=1}^{P}\ln q_{\theta_p}(\theta_p).
\end{aligned}
$$
\index{Functional derivatives}
\index{Lagrangian function}
\paragrapharrow{Update for local latent variables.}
The update for the latent variables under the mean-field approximation closely resembles the E-step in standard EM. 
For completeness and clarity, we rederive it here; key differences from the non-Bayesian case are highlighted in \textcolor{mylightbluetext}{blue}.
To obtain the E-step under the mean-field approximation, we maximize the ELBO at $t$-th iteration $\mathcalF_{\balpha}\left(\left\{q_{\bz_n}(\bz_n)\right\}_{n=1}^N, q_{\btheta^\toptzero}(\btheta^\toptzero)\right)$ with respect to each factor $q_{z_{nq}}(z_{nq})$, subject to the normalization constraints: $\int q_{z_{nq}}(z_{nq}) \,dz_{nq} =1$ for all $n\in\{1,2,\ldots,N\}, q\in\{1,2,\ldots, Q\}$.
The associated Lagrangian is:
$$
L\left(q_{z_{nq}}(z_{nq}), \bgamma\right)=
\mathcalF_{\balpha}\left(\left\{q_{\bz_n}(\bz_n)\right\}_{n=1}^N, q_{\btheta^\toptzero}(\btheta^\toptzero\right) + \sum_{n,q} \gamma_{nq} \left(\int q_{z_{nq}}(z_{nq}) dz_{nq} -1\right).
$$
Setting the functional  derivative  of the Lagrangian function with respect to $q_{z_{nq}}(z_{nq})$ to zero yields
\begin{equation}\label{equation:vbem_uncon_mf}
\begin{aligned}
&\ln q^\toptone_{z_{nq}}(z_{nq}) = \int_{\textcolor{mylightbluetext}{\btheta}} \int_{\bz} \left[ \prod_{k\neq q}^{Q}q_{z_{nk}}(z_{nk}) \ln p(\bz_n, \bx_n \mid\btheta^\toptzero, \textcolor{mylightbluetext}{\balpha})\right]d\bz_{n/q}\textcolor{mylightbluetext}{\,d\btheta} +\gamma_{nq}-1\\
&\implies q^\toptone_{z_{nq}}(z_{nq}) = \frac{1}{\mathcalC_{nq}}\exp\left\{\int_{\btheta} \int_{\bz} \left[ \prod_{k\neq q}^{Q}q_{z_{nk}}(z_{nk}) \ln p(\bz_n, \bx_n \mid\btheta^\toptzero, \textcolor{black}{\balpha})\right]d\bz_{n/q}\,d\btheta\right\},
\end{aligned}
\end{equation}
where $\mathcalC_{nq}$ is a normalization constant, $d{\bz_{n/q}}$ denotes the element of integration for all elements in $\bz_n$ except $z_{nq}$, 
$\prod_{k\neq q}^{Q}$ denotes the product of all elements except the $q$-th item.
More compactly, this can be written as:
\begin{equation}\label{equation:vbem_uncon_mf_comp}
q^\toptone_{z_{nq}}(z_{nq}) \leftarrow \frac{1}{\mathcalC_{nq}}\exp\left\{\textcolor{mylightbluetext}{\Exp_{q_{\btheta}}}
\left[
\Exp_{q(-z_{nq})} \left[ \ln p(\bz_n, \bx_n \mid\btheta^\toptzero, \textcolor{mylightbluetext}{\balpha})\right] 
\right]\right\},
\end{equation}
where the inner expectation is taken over all latent dimensions $n^\prime \in\{1,2,\ldots,N\}, q^\prime\in\{1,2,\ldots, Q\}$ except $\{n^\prime =n, q^\prime =q\}$.

\paragrapharrow{Update for model parameters (global latent variables).}
Similarly, we maximize the ELBO at $t$-th iteration $\mathcalF_{\balpha}\left(\left\{q^\toptone_{\bz_n}(\bz_n)\right\}_{n=1}^N, q_{\btheta}(\btheta)\right)$ with respect to $q_{\btheta}(\btheta)$, subject to the normalization constraints: $\int q_{\theta_p}(\theta_p) d\theta_p =1$ for all $p\in\{1,2,\ldots,P\}$. The associated Lagrangian function is
$$
L(q_{\theta_p} (\theta_p), \bnu)
=
\mathcalF_{\balpha}\left(\left\{q^\toptone_{\bz_n}(\bz_n)\right\}_{n=1}^N, q_{\btheta}(\btheta)\right) + 
\sum_{p} \nu_p \left(\int q_{\theta_p}(\theta_p) d\theta_p -1\right).
$$
Setting the functional  derivative  of the Lagrangian function with respect to $q_{\theta_p}(\theta_p)$ to zero yields
$$
\begin{aligned}
&\ln q^\toptone_{\theta_p}(\theta_p) = \int \prod_{k\neq p}^{P} q_{\theta_k}(\theta_k)
\left[   
\mathcalL_{\bz}(\mathcalZ^\toptone, \mathcalX, \btheta^\toptzero, \balpha) 
+
\mathcalL_{\btheta} (\btheta, \balpha) 
\right]
\,d\btheta_{-p}
-1+\nu_p\\
&\implies  q^\toptone_{\theta_p}(\theta_p) 
=
\frac{1}{\mathcalC_p} 
\exp
\left\{\int \prod_{k\neq p}^{P} q_{\theta_k}(\theta_k)
\left[   
\mathcalL_{\bz}(\mathcalZ^\toptone, \mathcalX, \btheta^\toptzero, \balpha) +\ln p(\btheta \mid \balpha)
\right]
\,d\btheta_{-p} 
\right\},
\end{aligned}
$$
where $d{\btheta_{-p}}$ denotes the element of integration for elements in $\btheta$ except $\theta_p$, $\mathcalC_p$ is a normalization constant.
Notably, when there are \textbf{no hidden variables} (\textbf{or the unified framework in \eqref{equation:elbo_unified}}),  the update for the model parameters reduces to 
\begin{equation}\label{equation:mf_nohidden_theta}
q^\toptone_{\theta_p}(\theta_p) 
\leftarrow
\frac{1}{\mathcalC_p} 
\exp
\left\{\int \prod_{k\neq p}^{P} q_{\theta_k}(\theta_k)
\ln p(\mathcalX,  \btheta\mid \balpha) 
\,d\btheta_{-p} 
\right\}.
\end{equation}

\begin{figure}[h!]
\centering
\includegraphics[width=0.55\textwidth]{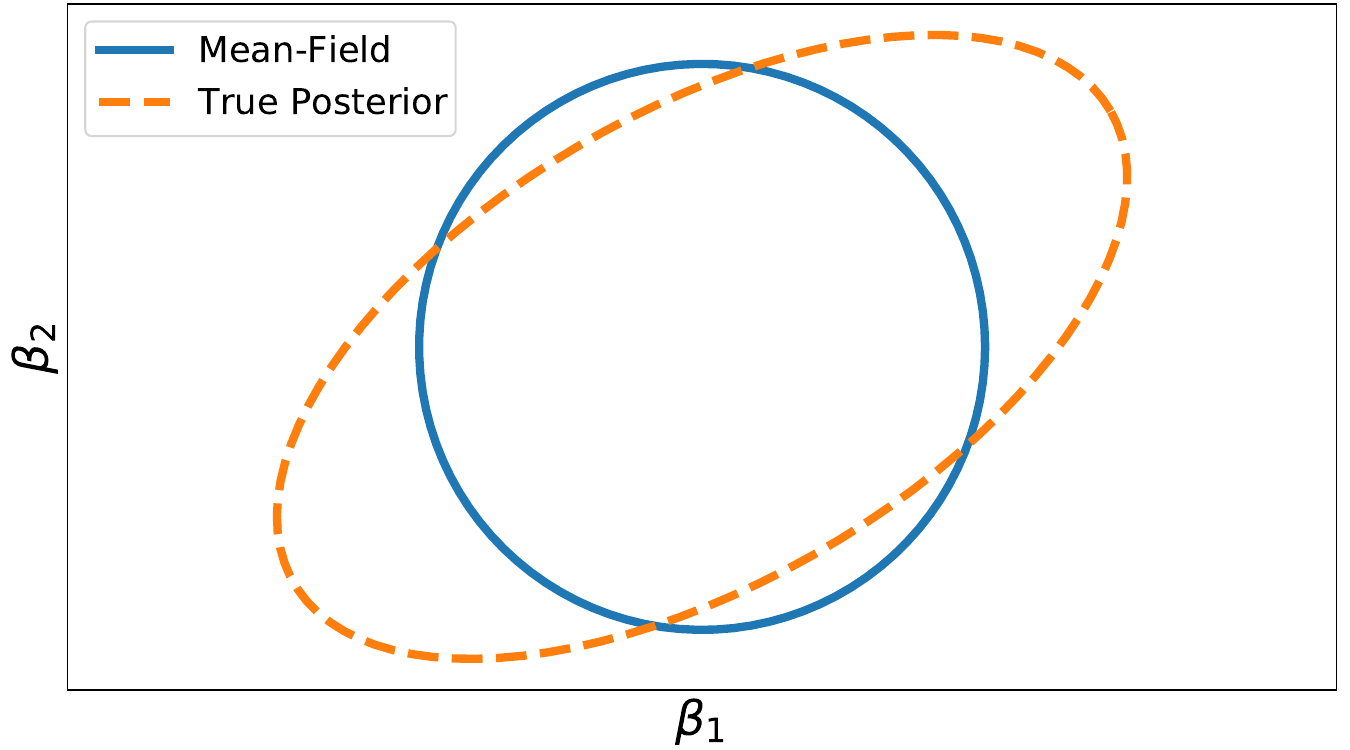}
\caption{Variational distribution (mean-field) versus true posterior in a two-dimensional Bayesian linear model. The mean-field approximation factorizes the posterior, leading to an underestimate of uncertainty (see Example~\ref{example:bayes_lr_mfa})}
\label{fig:blm_mean_field}
\end{figure}

\begin{example}[Bayesian Linear Regression]\label{example:bayes_lr_mfa}
Consider the Bayesian linear regression model from Section~\ref{sec:semiconjugate} and the notations therein; while we fix $\sigma^2$ constant.
The prior and likelihood are $\bbeta \sim \normal(\bbeta_0, \bSigma_0)$ and $\by\sim \normal(\bX\bbeta, \sigma^2\bI)$, respectively.
Therefore, the model assumes no (local) latent variables and only model parameters $\bbeta$.
The posterior distribution of $\bbeta$ is $p(\bbeta \mid \by, \bX,\sigma^2)=\normal(\bbeta_2,\bSigma_2)$, where 
$
\bSigma_2 =\left(\frac{1}{\sigma^2} \bX^\top\bX + \bSigma_0^{-1}\right)^{-1}
$ 
and 
$
\bbeta_2=\bSigma_2 \left(\bSigma_0^{-1}\bbeta_0+\frac{1}{\sigma^2}\bX^\top\by\right)
$.
Suppose model the parameters admit the following partition, within two-dimensional space:
$$
\bbeta =
\begin{bmatrix}
\beta_1\\
\beta_2
\end{bmatrix}
,\gap 
\bbeta_2 =
\begin{bmatrix}
b_1\\
b_2
\end{bmatrix}
, 
\gap 
\bSigma_2^{-1} = 
\begin{bmatrix}
\sigma_{11} & \sigma_{12}\\
\sigma_{21} & \sigma_{22}
\end{bmatrix},\,\, \sigma_{12}=\sigma_{21}.
$$
Then, under the mean-field approximation \eqref{equation:mf_nohidden_theta}, we have 
$$
\begin{aligned}
\ln q_{\beta_1}(\beta_1)
&=
\int q_{\beta_2}(\beta_2) 
\ln \underbrace{p(\by, \bX, \bbeta \mid \bbeta_0, \bSigma_0, \sigma^2)}_{\propto p(\bbeta \mid \by, \bX, \sigma^2)}
d\beta_2 +\mathcalC_1\\
&=\int q_{\beta_2}(\beta_2)   \left\{-\frac{1}{2} (\bbeta-\bbeta_2)^\top \bSigma_2^{-1}(\bbeta-\bbeta_2)\right\}  d\beta_2 +\mathcalC_2\\
&=-\frac{1}{2}\beta_1^2\sigma_{11} +\beta_1 b_1 \sigma_{11} - \beta_1\sigma_{12} (\Exp_{q_{\beta_2}(\beta_2)}[\beta_2] -b_2) +\mathcalC_3
\equiv \ln \normal(\beta_1 \mid \widehat{\mu}, \widehat{\sigma}^2),
\end{aligned}
$$ 
where $\widehat{\mu} = b_1-\frac{\sigma_{12}}{\sigma_{11}} (\Exp_{q_{\beta_2}(\beta_2)}[\beta_2] -b_2) $ and $\widehat{\sigma}^2 = \sigma_{11}^{-1}$. The variational distribution of $\beta_2$ can be computed similarly due to symmetry.
We note that the variational distribution follows a Gaussian distribution even though we have not assumed anything about the variational distribution.
However, this is not always the case; and we will provide an example that needs to assume the form of the variational distribution a priori in Section~\ref{section:mcvi}.
A caricature of variational distribution is given in Figure~\ref{fig:blm_mean_field}. 
The variational parameter $\widehat{\sigma}^2 = \sigma_{11}^{-1}$ indicates that the variance of $q_{\beta_1}(\beta_1)$ is equal to the variance of the conditional distribution $p(\beta_1\mid \beta_2, \mathcalX)$, which is smaller than the variance of the marginal distribution $p(\beta_1\mid \mathcalX)$, and therefore, mean-field VI underestimates the posterior uncertainty in this case.
\end{example}

\index{Stochastic MC gradient}
\index{Monte Carlo variational inference}
\index{Black-box variational inference}
\index{Stochastic optimization}
\subsection{Monte Carlo/Stochastic Variational Inference}\label{section:mcvi}
We should emphasize that the book primarily focuses on MCMC methods, particularly  Gibbs sampling, for Bayesian matrix factorizations. Variational inference  alternatives are discussed only briefly.
Nevertheless, we provide a concise overview of an important topic within VI: \textit{black-box VI (BBVI), also known as Monte Carlo VI (MCVI)} \citep{wingate2013automated, ranganath2014black, li2018approximate}. 
BBVI is of particular interest to practitioners because it generalizes the VI framework derived earlier to more complex models (e.g., non-conjugate or non-exponential-family models), enables inference when analytical updates are intractable, and relies on Monte Carlo estimates of gradients combined with stochastic optimization (see, e.g., \citet{lu2022gradient}).

For notational compactness, we adopt the unified framework introduced in \eqref{equation:elbo_unified}, which treats all unknowns---both model parameters and latent variables---as a single joint variable $\bomega=\{\btheta,\mathcalZ\}$.
Assume the variational distribution $q(\bomega\mid \blambda)$ is parameterized by $\blambda$. 
The ELBO then takes the form:
\begin{equation}\label{equation:mcvi_elbo}
\mathcalF(\blambda)=\int q(\bomega\mid \blambda)  \ln \frac{p( \mathcalX, \bomega \mid \balpha)}{q(\bomega\mid \blambda) } d\bomega 
=\Exp_{q(\bomega\mid \blambda)} \left[\ln p( \mathcalX, \bomega \mid \balpha) - \ln q(\bomega\mid \blambda)\right].
\end{equation}
\begin{theoremHigh}[Gradient of ELBO]\label{theorem:grad_elbo}
Consider the ELBO $\mathcalF(\blambda)$ defined in \eqref{equation:mcvi_elbo}. 
Its gradient with respect to $\blambda$ is 
\begin{equation}\label{equation:grad_elbo_theo1}
\nabla_{\blambda}\mathcalF(\blambda)
=
\Exp_{q(\bomega\mid \blambda)}
\big[
\{\nabla_{\blambda} \ln q(\bomega\mid \blambda) \}
\left\{
\ln p(\mathcalX, \bomega \mid \balpha) - \ln q(\bomega\mid \blambda)
\right\} 
\big].
\end{equation}
The quantity  $\nabla_{\blambda} \ln q(\bomega\mid \blambda)$ is known in statistics as the \textit{score function}.~\footnote{The score function measures the sensitivity of the log-likelihood to changes in the parameters. It indicates how much the log-likelihood would change if the parameters were slightly perturbed. This sensitivity is crucial for understanding how well the current parameter values fit the observed data.}
\end{theoremHigh}
\index{Score function}
\index{Sensitivity}
\index{Expected score function}
\begin{proof}[of Theorem~\ref{theorem:grad_elbo}]
To see this, using product rule and exchange derivatives with integrals via the dominated convergence theorem~\footnote{The score function exists. The score function and likelihoods are bounded.}, we have 
$$
\begin{aligned}
&\nabla_{\blambda}\mathcalF(\blambda)
=
\nabla_{\blambda} \int 
q(\bomega\mid \blambda)
\left[\ln p( \mathcalX, \bomega \mid \balpha) - \ln q(\bomega\mid \blambda)\right]
d\bomega\\
&= 
\int 
\nabla_{\blambda}q(\bomega\mid \blambda)
\left[\ln p( \mathcalX, \bomega \mid \balpha) - \ln q(\bomega\mid \blambda)\right]
d\bomega
-
\int 
q(\bomega\mid \blambda)
\left(\nabla_{\blambda} \ln q(\bomega\mid \blambda)\right)
d\bomega.
\end{aligned}
$$
Since $\nabla_{\blambda}q(\bomega \mid \blambda) = [\nabla_{\blambda}\ln q(\bomega \mid \blambda)]q(\bomega \mid \blambda)$, the first term becomes
$$
\begin{aligned}
&\gap\int
\nabla_{\blambda}[\ln q(\bomega \mid \blambda)]q(\bomega \mid \blambda)
\left[\ln p( \mathcalX, \bomega \mid \balpha) - \ln q(\bomega\mid \blambda)\right]
d\bomega\\
&=\Exp_{q(\bomega \mid \blambda)} \big[ \{\nabla_{\blambda}\ln q(\bomega \mid \blambda)\}
\left\{\ln p( \mathcalX, \bomega \mid \balpha) - \ln q(\bomega\mid \blambda)\right\} \big].
\end{aligned}
$$
The second term, the expected value of the score function, simplifies to  
$$
\begin{aligned}
\int 
q(\bomega\mid \blambda)
\left(\nabla_{\blambda} \ln q(\bomega\mid \blambda)\right)
d\bomega
&=
\Exp_{q(\bomega\mid \blambda)}
\left[ \frac{ \nabla_{\blambda}q(\bomega\mid \blambda)}{q(\bomega\mid \blambda)}\right]=\nabla_{\blambda} \int q(\bomega\mid \blambda)d\bomega = \bzero.
\end{aligned}
$$
This concludes the result.
\end{proof}

Once the gradient $\nabla_{\blambda}\mathcalF(\blambda)$ is available, we can update the variational parameter $\blambda$ via \textit{gradient ascent}:
$$
\blambda^\toptone \leftarrow  \blambda^\toptzero + \eta_t \nabla_{\blambda}\mathcalF(\blambda^\toptzero),
$$
where $\eta_t>0$ is the step size, and $\nabla_{\blambda}\mathcalF(\blambda)$ is a descent direction (again, taking a positive step in this direction leads to an increase in the function value).
However, the expectation in \eqref{equation:grad_elbo_theo1} over $\bomega\sim q(\bomega\mid \blambda)$ is often \textbf{intractable}. 
In such cases, we approximate the gradient using Monte Carlo samples (\textit{Monte Carlo gradients} or \textit{MC gradients}); hence the name Monte Carlo VI.
Specifically, drawing $S$ samples $\bomega_s \sim q(\bomega\mid \blambda)$, we compute the unbiased estimator:
\begin{equation}\label{equation:mc_grad_eq_1}
\nabla_{\blambda}\mathcalF(\blambda) \approx 
\frac{1}{S}
\sum_{s=1}^{S}\{\nabla_{\blambda} \ln q(\bomega_s\mid \blambda) \}
\left\{
\ln p(\mathcalX, \bomega_s \mid \balpha) - \ln q(\bomega_s\mid \blambda)
\right\},
\,\, \bomega_s \sim q(\bomega\mid \blambda).
\end{equation}
When the data set is large, we can further improve scalability by using stochastic mini-batch sampling.
Let $\sK\subset\{1,2,\ldots,N\}$ be a random subset of indices with $\abs{\sK}$ elements. Then the \textit{stochastic Monte Carlo gradient} is:
$$
\nabla_{\blambda}\mathcalF(\blambda) \approx 
\frac{1}{S}
\sum_{s=1}^{S}\{\nabla_{\blambda} \ln q(\bomega_s\mid \blambda) \}
\left\{
\frac{N}{\abs{\sK}}\sum_{k\in\sK}\ln p(\bx_k, \bomega_s \mid \balpha) - \ln q(\bomega_s\mid \blambda)
\right\},
\,\, \bomega_s \sim q(\bomega| \blambda),
$$
where $N$ is the total number of observed points. 
This estimator remains unbiased under random sampling of $\sK$.
The general MCVI algorithm is summarized in Algorithm~\ref{alg:mcvi}.

\begin{algorithm}[h] 
\caption{Monte Carlo Variational Inference (MCVI)}
\label{alg:mcvi}
\begin{algorithmic}[1] 
\Require Observed data points $\mathcalX=\{\bx_1, \bx_2, \ldots, \bx_N\}$;
\State \textbf{initialize:} $\blambda^\topone$; 
\State Choose the maximal number of iterations $C$;
\State $t=0$; \Comment{Count for the number of iterations}
\While{$t<C$} 
\State $t=t+1$;
\State $\blambda^\toptone \leftarrow  \blambda^\toptzero + \eta_t \nabla_{\blambda}\mathcalF(\blambda^\toptzero)$;
\EndWhile
\State Output $\blambda^\toptzero$;
\end{algorithmic} 
\end{algorithm}

\index{Variance reduction}
\index{Reparameterization trick}
\subsubsection*{Reparameterization Trick}
MCVI can be made more efficient using the \textit{reparameterization trick} \citep{williams1992simple, kingma2013auto, ho2020denoising}. 
Suppose a random variable $\btheta$ is distributed as $\btheta\sim p(\btheta)$, and that we can express 
$\btheta$ as a deterministic transformation of an auxiliary random variable: $\bepsilon \sim p(\bepsilon)$ and $\btheta=f(\bepsilon)$.
Then, for any function $g(\btheta)$, the expectation under $p(\btheta)$   can be rewritten as an expectation over $\bepsilon$:
\begin{equation}
\Exp_{p(\btheta)}\left[ g(\btheta) \right]
=
\Exp_{p(\bepsilon)} \big[ g(f(\bepsilon)) \big].
\end{equation}
This reparameterization is particularly useful when direct sampling from $p(\btheta)$ is difficult or when the Monte Carlo  gradient estimator based on the score function has high variance, which can significantly slow convergence \citep{kingma2013auto}.

In the context of variational inference, recall from Theorem~\ref{theorem:grad_elbo} that the gradient of the ELBO $\mathcalF(\blambda)$ involves an expectation over the variational distribution $q(\bomega \mid \blambda)$, where the parameter $\blambda$ appears inside the distribution $q(\bomega \mid \blambda)$.
This dependence leads to high-variance gradient estimates when using the score-function method.
The reparameterization trick addresses this issue by decoupling the randomness from the parameters. Specifically, we assume that samples from the variational distribution can be generated via a differentiable transformation:
$$
\bomega\sim q(\bomega \mid \blambda)\gap\rightarrow\gap  \bepsilon \sim p(\bepsilon), \bomega=f_{\blambda}(\bepsilon),
$$
where both the variational distribution $q(\bomega \mid \blambda)$ and the variable transformation function $ \bomega=f_{\blambda}(\bepsilon  )$ depend on the same parameter $\blambda$, but the base distribution 
$p(\bepsilon)$ is fixed and independent of $\blambda$.
Under this reparameterization, the gradient of the ELBO takes the following form.
\begin{theoremHigh}[Gradient of ELBO under Reparameterization]\label{theorem:grad_elbo_repara}
Consider the ELBO $\mathcalF(\blambda)$ defined in \eqref{equation:mcvi_elbo} and the reparameterization with $\bepsilon \sim p(\bepsilon), \bomega=f_{\blambda}(\bepsilon)$. 
The  gradient of $\mathcalF(\blambda)$  with respect to $\blambda$ is 
\begin{equation}
\nabla_{\blambda}\mathcalF(\blambda)
=
\Exp_{p(\bepsilon)} \left[ \nabla_{f}\ln p( \mathcalX, f_{\blambda}(\bepsilon) \mid \balpha) \nabla_{\blambda}f_{\blambda}(\bepsilon)\right] 
-
\Exp_{p(\bepsilon)}[\nabla_{f} \ln q(f_{\blambda}(\bepsilon)\mid \blambda)  \nabla_{\blambda}f_{\blambda}(\bepsilon) ].
\end{equation}
\end{theoremHigh}
\begin{proof}[of Theorem~\ref{theorem:grad_elbo_repara}]
The gradient of $\mathcalF(\blambda)$ w.r.t. $\blambda$ becomes, using the chain rule and the dominated convergence theorem,
$$
\begin{aligned}
\nabla_{\blambda}\mathcalF(\blambda)
&=\nabla_{\blambda}\Exp_{q(\bomega\mid \blambda)} \left[\ln p( \mathcalX, \bomega \mid \balpha) - \ln q(\bomega\mid \blambda)\right]\\
&=\nabla_{\blambda}\Exp_{p(\bepsilon)} \left[\ln p( \mathcalX, f_{\blambda}(\bepsilon) \mid \balpha)]  - \nabla_{\blambda}\Exp_{p(\bepsilon)}[\ln q(f_{\blambda}(\bepsilon)\mid \blambda)\right]\\
&=\Exp_{p(\bepsilon)} \left[ \nabla_{f}\ln p( \mathcalX, f_{\blambda}(\bepsilon) \mid \balpha)\cdot \nabla_{\blambda}f_{\blambda}(\bepsilon)\right] 
- \nabla_{\blambda}\Exp_{p(\bepsilon)}\left[\ln q(f_{\blambda}(\bepsilon)\mid \blambda)\right].
\end{aligned}
$$
The second term is
$$
\begin{aligned}
\nabla_{\blambda}\Exp_{p(\bepsilon)}[\ln q(f_{\blambda}(\bepsilon)\mid \blambda)] 
&=\Exp_{p(\bepsilon)}[\nabla_{\blambda} \ln q(\bomega\mid \blambda)|_{\bomega=f_{\blambda}(\bepsilon)}] 
+
\Exp_{p(\bepsilon)}[\nabla_{f} \ln q(f_{\blambda}(\bepsilon)\mid \blambda)  \nabla_{\blambda}f_{\blambda}(\bepsilon) ]\\
&=\Exp_{q(\bomega\mid \blambda)}[\nabla_{\blambda} \ln q(\bomega\mid \blambda)] 
+
\Exp_{p(\bepsilon)}[\nabla_{f} \ln q(f_{\blambda}(\bepsilon)\mid \blambda)  \nabla_{\blambda}f_{\blambda}(\bepsilon) ],
\end{aligned}
$$
where the first term is $\bzero$ from the proof of Theorem~\ref{theorem:grad_elbo}.
This concludes the result.
\end{proof}

Therefore, the MC gradient of the ELBO under the reparameterization trick becomes
$$
\begin{aligned}
\nabla_{\blambda}\mathcalF(\blambda)
\approx
\frac{1}{S}\sum_{s=1}^{S}
\big\{
\nabla_{f}\ln p( \mathcalX, f_{\blambda}(\bepsilon_s) | \balpha) \nabla_{\blambda}f_{\blambda}(\bepsilon_s)
-
\nabla_{f} \ln q(f_{\blambda}(\bepsilon_s)| \blambda) & \nabla_{\blambda}f_{\blambda}(\bepsilon_s) 
\big\},
\,\bepsilon\sim p(\bepsilon).
\end{aligned}
$$
A stochastic mini-batch version can be derived analogously by subsampling the data; we omit the details here for brevity.

\begin{example}[Reparameterization by Gaussian]
The most common and illustrative application of the reparameterization trick uses a Gaussian variational distribution.
We encourage the reader to briefly consult Section~\ref{sec:multi_gaussian_conjugate_prior} for background on multivariate Gaussian properties before continuing.
To be more specific, we assume $q(\bomega\mid \blambda)=\normal(\bomega \mid \bmu, \bSigma)$, where $\bSigma$ admits the Cholesky decomposition $\bSigma=\bL\bL^\top$~\footnote{See, for example, Problem~\ref{problem:cholesky} for more details.} and $\blambda=\{\bmu,\bL\}$.
Then $\bomega$ can be equivalently sampled from $f_{\blambda}(\bepsilon)=\bmu+\bL\bepsilon, \bepsilon\sim\normal(\bzero, \bI)$.
Therefore, the MC gradient w.r.t. $\bmu$ and $\bL$ are
$$
\begin{aligned}
\nabla_{\bmu}&= \frac{1}{S}\sum_{s=1}^{S} \left\{\nabla_{f}\ln p(\mathcalX, f)|_{f=\bmu+\bL\bepsilon_s} - \bSigma^{-1}(-\bL\bepsilon_s)\right\},\\
\nabla_{\bL}&= \frac{1}{S}\sum_{s=1}^{S} \left\{\nabla_{f}\ln p(\mathcalX, f)|_{f=\bmu+\bL\bepsilon_s} \cdot \bepsilon_s - \bSigma^{-1}(-\bL\bepsilon_s)\bepsilon_s\right\},\gap \bepsilon_s\sim\normal(\bzero,\bI),
\end{aligned}
$$
where we use the fact that
$\ln q(\bomega \mid \blambda)=-\frac{1}{2}\ln\abs{\bSigma}-\frac{1}{2}(\bomega-\bmu)^\top\bSigma^{-1}(\bomega-\bmu)+\mathcalC$.
\end{example}

\index{Variance reduction}
\index{Rao--Blackwell theorem}
\index{Monte Carlo variational inference}
\subsubsection*{Variance Reduction for MCVI Gradients with Rao--Blackwellization}
Although the stochastic MC gradient of $\mathcalF(\blambda)$ is an unbiased estimator of $\nabla_{\blambda}\mathcalF(\blambda)$,  its variance can be prohibitively large.
High variance may cause gradient updates to fluctuate wildly---potentially moving in unhelpful directions---and significantly slow convergence. Therefore, variance reduction is essential in stochastic gradient ascent/descent to ensure stable and efficient optimization.

The classical  \textit{Rao--Blackwell theorem} states that \citep{rao1973linear, blackwell1947conditional} provides a principled way to reduce variance. 
It states that if $\widehatbtheta$ is an unbiased estimator of a parameter $\btheta$, and $T$ is a sufficient statistic for $\btheta$, then the conditional expectation $\widehatbtheta^*=\Exp[\widehatbtheta\mid T]$ is also an unbiased estimator of $\btheta$, and satisfies
$$
\Var[\widehatbtheta^*] \leq \Var[\widehatbtheta],
$$
where the equality is attained if and only $\Pr[\widehatbtheta^*=\widehatbtheta]=1$.
In the context of variational inference, suppose the latent variables decompose as $\bomega=\{\bomega_1, \bomega_2\}$, and we wish to estimate the expectation of a function $g(\bomega_1, \bomega_2)$ under the variational distribution $q(\bomega_1, \bomega_2)$: $\Exp_{q(\bomega_1, \bomega_2)}[g(\bomega_1, \bomega_2)]$.
Define $g_2(\bomega_2) \triangleq \Exp_{q(\bomega_1\mid\bomega_2)}[g(\bomega_1, \bomega_2)]$. Then the variance can be decomposed, similar to the bias-variance decomposition~\footnote{see, for example, \citet{lu2021rigorous}.}, by 
\begin{equation}\label{equation:rao-black_deri}
\begin{aligned}
&\gap \Var_{q(\bomega_1,\bomega_2)}\left[g(\bomega_1,\bomega_2)\right]
=\Exp_{q(\bomega_1,\bomega_2)}\left[\left( g(\bomega_1,\bomega_2) - \Exp_{q(\bomega_1, \bomega_2)}[g(\bomega_1, \bomega_2)]\right)^2 \right]\\
&=\Exp_{q(\bomega_2)}\Exp_{q(\bomega_1\mid \bomega_2)}\left[\left( g(\bomega_1,\bomega_2) -g_2(\bomega_2)+g_2(\bomega_2) - \Exp_{q(\bomega_1, \bomega_2)}[g(\bomega_1, \bomega_2)]\right)^2 \right]\\
&=\Var_{q(\bomega_2)} [g_2(\bomega_2)] + \Exp_{q(\bomega_1,\bomega_2)}\left[ \left\{g_2(\bomega_2) - \Exp_{q(\bomega_1, \bomega_2)}[g(\bomega_1, \bomega_2)]\right\}^2\right]\\
&\geq \Var_{q(\bomega_2)} [g_2(\bomega_2)],
\end{aligned}
\end{equation}
where we use the fact that $\Exp_{q(\bomega_2)}[g_2(\bomega_2)] = \Exp_{q(\bomega_1, \bomega_2)}[g(\bomega_1, \bomega_2)]$.
This inequality shows that replacing $g(\bomega_1, \bomega_2)$ with its conditional expectation 
$g_2(\bomega_2)$ yields a lower-variance estimator---this is the essence of Rao--Blackwellization.

Returning to the  ELBO $\mathcalF(\blambda)$ in \eqref{equation:mcvi_elbo}, suppose we partition the latent variables as  $\bomega=\{\bomega_1, \bomega_2\}$ and assume a factorized variational distribution: \colorbox{\mdframecolor}{$q(\bomega\mid \blambda) = q(\bomega_1\mid\blambda_1)q(\bomega_2\mid\blambda_2)$}. Under this structure, we can apply Rao--Blackwellization to obtain a lower-variance gradient estimator, as stated below.

\index{Rao--Blackwellization}
\begin{theoremHigh}[Gradient of ELBO with Rao--Blackwellization]\label{theorem:grad_elbo_rao}
Consider the ELBO $\mathcalF(\blambda)$ defined in \eqref{equation:mcvi_elbo} with factorized variational distribution $q(\bomega\mid \blambda) = q(\bomega_1\mid\blambda_1)q(\bomega_2\mid\blambda_2)$ and $\bomega=\{\bomega_1,\bomega_2\}$. 
Then the gradient  with respect to $\blambda_1$ is 
\begin{equation}
\nabla_{\blambda_1}\mathcalF(\blambda)
=
\Exp_{q(\bomega_1\mid\blambda_1)}
\bigg[
\left\{\nabla_{\blambda_1}\ln q(\bomega_1 \mid \blambda_1)\right\}
\left\{
\Exp_{q(\bomega_2\mid\blambda_2)}\big[ \ln p( \mathcalX, \bomega \mid \balpha)\big]
- \ln q(\bomega_1\mid\blambda_1)\right\}
\bigg].
\end{equation}
By symmetry, an analogous expression holds for the gradient with respect to $\blambda_2$.
\end{theoremHigh}
\begin{proof}[of Theorem~\ref{theorem:grad_elbo_rao}]
The gradient w.r.t. $\blambda_1$ is 
$$
\begin{aligned}
&\nabla_{\blambda_1}\mathcalF(\blambda)
=
\nabla_{\blambda_1} \int 
q(\bomega_1\mid\blambda_1)q(\bomega_2\mid\blambda_2)
\left[\ln p( \mathcalX, \bomega \mid \balpha) - \ln q(\bomega_1\mid\blambda_1)-\ln q(\bomega_2\mid\blambda_2)\right]
d\bomega\\
&= \Exp_{q(\bomega_2\mid\blambda_2)}\left[ 
\int 
\nabla_{\blambda_1}q(\bomega_1\mid \blambda_1)
\left[\ln p( \mathcalX, \bomega \mid \balpha) - \ln q(\bomega_1\mid\blambda_1)-\ln q(\bomega_2\mid\blambda_2)\right]
d\bomega_1\right]\\
&\gap\gap \gap\gap -
\Exp_{q(\bomega_2\mid\blambda_2)}\left[ 
\int 
q(\bomega_1\mid \blambda_1)
\left(\nabla_{\blambda_1} \ln q(\bomega_1\mid \blambda_1)\right)
d\bomega_1\right].
\end{aligned}
$$
The second term is $\bzero$ (same as the one in the proof of Theorem~\ref{theorem:grad_elbo}, which is the expectation of a score function).
Since $\nabla_{\blambda_1}q(\bomega_1 \mid \blambda_1) = [\nabla_{\blambda_1}\ln q(\bomega_1 \mid \blambda_1)]q(\bomega_1 \mid \blambda_1)$, the first term becomes
$$
\begin{aligned}
&\gap\Exp_{q(\bomega_2\mid\blambda_2)}\left[ 
\int 
\nabla_{\blambda_1}q(\bomega_1\mid \blambda_1)
\left[\ln p( \mathcalX, \bomega \mid \balpha) - \ln q(\bomega_1\mid\blambda_1)-\ln q(\bomega_2\mid\blambda_2)\right]
d\bomega_1\right]\\
&=\Exp_{q(\bomega_2\mid\blambda_2)}\bigg[ 
\Exp_{q(\bomega_1\mid\blambda_1)}
\big[
\left\{\nabla_{\blambda_1}\ln q(\bomega_1 \mid \blambda_1)\right\}
\left\{\ln p( \mathcalX, \bomega \mid \balpha) - \ln q(\bomega_1\mid\blambda_1)-\ln q(\bomega_2\mid\blambda_2)\right\}
\big]
\bigg]\\
&=
\Exp_{q(\bomega_1\mid\blambda_1)}
\bigg[
\left\{\nabla_{\blambda_1}\ln q(\bomega_1 \mid \blambda_1)\right\}
\left\{\Exp_{q(\bomega_2\mid\blambda_2)}\big[ \ln p( \mathcalX, \bomega \mid \balpha)\big] - \ln q(\bomega_1\mid\blambda_1)\right\}
\bigg],
\end{aligned}
$$
where the last equality comes from that the second term above is $\bzero$.
This concludes the result.
\end{proof}

Therefore, if we can compute $\Exp_{q(\bomega_2\mid\blambda_2)}\big[ \ln p( \mathcalX, \bomega \mid \balpha)\big]$ analytically, the MC gradient w.r.t. $\blambda_1$ is 
$$
\begin{aligned}
\nabla_{\blambda_1}\mathcalF(\blambda)
\approx 
\frac{1}{S}
\sum_{s=1}^{S}
\bigg[
\left\{\nabla_{\blambda_1}\ln q(\bomega_1^s \mid \blambda_1)\right\}
\left\{
\Exp_{q(\bomega_2\mid\blambda_2)}\big[ \ln p( \mathcalX, \bomega_1^s, \bomega_2 \mid \balpha)\big]
- \ln q(\bomega_1^s\mid\blambda_1)\right\}
\bigg],
\end{aligned}
$$
where $\bomega_1^s\sim q(\bomega_1 \mid \blambda_1)$ such that the MC gradient $\nabla_{\blambda_1}\mathcalF(\blambda)$ has a smaller variance than the MC gradient $\nabla_{\blambda}\mathcalF(\blambda)$ in \eqref{equation:mc_grad_eq_1}.

\index{Control variate}
\subsubsection*{Variance Reduction for MCVI Gradients with Control Variate}
When deriving the variance reduction in Rao--Blackwellization in Equation~\eqref{equation:rao-black_deri}, we use a variate $g_2(\bomega_2)$ to decompose the total variance and thereby reduce it.
The \textit{control variate} method generalizes this idea by introducing an arbitrary auxiliary function $h(\bomega)$ with known (or easily computable) expectation under the variational distribution.
Once again, we consider the parameter $\bomega$ with its variational distribution $q(\bomega)$, and we would like to estimate the expectation for any function $g(\bomega)$: $\Exp_{q(\bomega)}[g(\bomega)]$.
We have 
\begin{equation}\label{equation:control_vari_1}
\Exp_{q(\bomega)}[g(\bomega)]
=
\Exp_{q(\bomega)}\big[\underbrace{g(\bomega) - h(\bomega) +  \Exp_{q(\bomega)}[h(\bomega)]}_{\triangleq G(\bomega)}  \big],
\end{equation}
i.e., $G(\bomega)$ is an unbiased estimator of $g(\bomega)$.
Suppose $\Var_{q(\bomega)}[h(\bomega)]<\infty$, then the variance of  $G(\bomega)$ is
$$
\Var_{q(\bomega)}\left[ G(\bomega)\right]
=
\Var_{q(\bomega)}\left[ g(\bomega)\right] + \Var_{q(\bomega)}\left[ h(\bomega)\right] - 2\Cov\left[g(\bomega),h(\bomega)\right].
$$
Therefore,  a careful choice of $h(\bomega)$ such that $\Var_{q(\bomega)}\left[ h(\bomega)\right] - 2\Cov\left[g(\bomega),h(\bomega)\right]<0$ will reduce the variance of $\Var_{q(\bomega)}\left[ G(\bomega)\right]<\Var_{q(\bomega)}\left[ g(\bomega)\right]$.

\paragrapharrow{Scaling control variate.}
To gain finer control over the magnitude of variance reduction, we introduce a scaling parameter $\mu\in\real$ for \eqref{equation:control_vari_1} and  define:
$$
\begin{aligned}
\Exp_{q(\bomega)}[g(\bomega)]
&=
\Exp_{q(\bomega)}\big[\underbrace{g(\bomega) - \mu\left(h(\bomega) -  \Exp_{q(\bomega)}[h(\bomega)]\right)  }_{\triangleq G(\bomega)}  \big];\\
\Var_{q(\bomega)}\left[ G(\bomega)\right]
&=
\Var_{q(\bomega)}\left[ g(\bomega)\right] + \mu^2\Var_{q(\bomega)}\left[ h(\bomega)\right] - 2\mu\Cov\left[g(\bomega),h(\bomega)\right].
\end{aligned}
$$
The optimal value for $\mu$ is simply $\mu^*=\Cov\left[g(\bomega),h(\bomega)\right]/ \Var_{q(\bomega)}\left[ h(\bomega)\right]$ in this sense by optimizing the quadratic function of $\mu$.

\index{Score function}
\paragrapharrow{Control variate using the  score function.}
In the context of VI, we again consider the ELBO $\mathcalF(\blambda)$ in \eqref{equation:mcvi_elbo} and notations therein. Since the expectation of the score function $\nabla_{\blambda} \ln q(\bomega\mid \blambda)$ under $q(\bomega\mid \blambda)$ is zero: $\Exp_{q(\bomega\mid \blambda)}[\nabla_{\blambda} \ln q(\bomega\mid \blambda)]=\bzero$ (see the proof of  Theorem~\ref{theorem:grad_elbo}).
Element-wise, we consider the $i$-th component $\lambda_i$ of the variational parameter $\blambda$ (see Equation~\eqref{equation:grad_elbo_theo1}), and define the following functions
$$
\begin{aligned}
g_i(\bomega) &\triangleq  \{\nabla_{\lambda_i} \ln q(\bomega\mid \lambda_i) \}
\left\{
\ln p(\mathcalX, \bomega \mid \balpha) - \ln q(\bomega\mid \lambda_i)
\right\};\\
h_i(\bomega) &\triangleq \nabla_{\lambda_i} \ln q(\bomega\mid \lambda_i).
\end{aligned}
$$
Therefore, the optical scale   $\mu_i^*$ for the $i$-th component is $\mu_i^*=\Cov\left[g_i(\bomega),h_i(\bomega)\right]/ \Var_{q_i(\bomega)}\left[ h_i(\bomega)\right]$.
sing this, the Monte Carlo gradient estimator with control variates becomes:
$$
\nabla_{\lambda_i}\mathcalF(\blambda)
\approx 
\frac{1}{S}\sum_{s=1}^{S}
\{\nabla_{\lambda_i} \ln q(\bomega_s\mid \lambda_i) \}
\left\{
\ln p(\mathcalX, \bomega_s \mid \balpha) - \ln q(\bomega_s\mid \lambda_i)
-\mu_i^* 
\right\}, 
\bomega_s\sim q(\bomega \mid \blambda).
$$
This estimator typically exhibits significantly lower variance than the naive score-function estimator, especially when $q_i(\bomega)$ and $h_i(\bomega)$ are strongly correlated. The resulting gradient is then used in standard stochastic gradient ascent updates.

\subsubsection*{Estimating the Marginal Likelihood}
After optimizing the variational distribution $q(\bomega \mid \blambda)$ as discussed, the log-marginal likelihood under the unified framework (Remark~\ref{remark:unified_vi}), estimated by the MC sample, can be approximated as:
$$
\begin{aligned}
\ln p(\mathcalX\mid \balpha) 
&= \ln \int q(\bomega \mid \blambda) \frac{p(\mathcalX, \bomega \mid \balpha)}{q(\bomega \mid \blambda)} d\bomega\\
&\approx \ln \left(\frac{1}{S}\sum_{s=1}^{S}\frac{p(\mathcalX, \bomega_s \mid \balpha)}{q(\bomega_s \mid \blambda)}\right), \gap \bomega_s\sim q(\bomega \mid \blambda).
\end{aligned}
$$
This is known as the Monte Carlo log-marginal likelihood estimator (or the log of the importance sampling estimator). 
As $S\rightarrow \infty$, the estimator converges almost surely to the true log-marginal likelihood by the law of large numbers.

\begin{figure}[h!]
\centering  
\vspace{-0.35cm} 
\subfigtopskip=2pt 
\subfigbottomskip=2pt 
\subfigcapskip=-5pt 
\subfigure[(Standard) VI.]{\label{fig:lvm_VI}
	\includegraphics[width=0.431\linewidth]{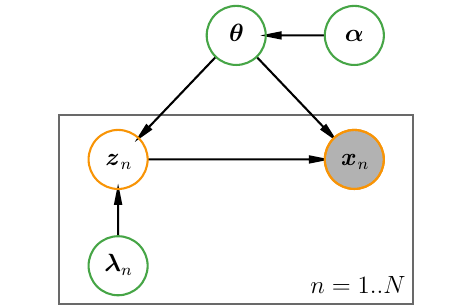}}
\subfigure[Amortized VI.]{\label{fig:lvm_VI_amortize}
	\includegraphics[width=0.431\linewidth]{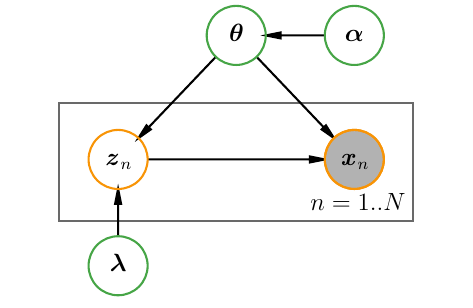}}
\caption{Graphical model representation of latent variable models under variational inference. Green circles denote prior variables, orange circles represent observed and latent variables, and plates  indicate replicated structures. }
\label{fig:lvm_and_hyp_vb}
\end{figure}

\index{Amortized variational inference}
\index{Scalability}
\subsection{Amortized Variational Inference}\label{section:amorized_VI}
We considered the constrained framework for  variation Bayesian inference  in Remark~\ref{remark:constrained_VBi}, with its graphical representation shown in Figure~\ref{fig:lvm_VI}.
The corresponding optimization problem over the ELBO is 
$$
\begin{aligned}
&\gap \mathop{\argmax}_{q_{\btheta}, \blambda_n}\mathcalF_{\balpha} \left(\left\{q_{\bz_n}(\bz_n)\right\}_{n=1}^N,q_{\btheta}(\btheta) \right)\\
&=\mathop{\argmax}_{q_{\btheta}, \blambda_n} \int \prod_{n=1}^{N}q_{\bz_n}(\bz_n\mid \blambda_n)q_{\btheta}(\btheta) \ln p(\mathcalZ, \mathcalX, \btheta \mid \balpha) \,d\btheta \,d\mathcalZ\\
&\gap \gap \gap\,\,-
\sum_{n=1}^{N}\int q_{\bz_n}(\bz_n\mid \blambda_n) \ln q_{\bz_n}(\bz_n\mid \blambda_n)\,d\bz_n
-
\int q_{\btheta}(\btheta) \ln q_{\btheta}(\btheta)\,d\btheta.
\end{aligned}
$$
A key limitation of this formulation becomes apparent when the number of observations $N$ is large: the model maintains a separate set of local variational parameters $\blambda_n$
for each data point $\bz_n$, leading to significant computational inefficiency. Moreover, there is no mechanism to share or reuse information across data points---each inference is performed in isolation, without leveraging knowledge from previously processed examples.

This issue is addressed through an ``\textit{amortization}" strategy. 
Instead of optimizing local parameters independently for every data point, \textit{amortized variational inference (amortized VI)} shares computation across the dataset by learning a global mapping from observed data to variational posteriors \citep{kingma2013auto, rezende2014stochastic, gershman2014amortized, rezende2015variational, ganguly2023amortized}.

Specifically, amortized VI introduces a stochastic encoder---often implemented as a neural network \citep{lecun2015deep, goodfellow2016deep}---that takes an observation $\bx_n$ as input and outputs the parameters of the approximate posterior $q_{\blambda}(\bz_n\mid\bx_n)$. These encoder parameters 
$\blambda$ are shared across all data points and learned jointly with any global parameters (e.g., $\btheta$; see Figure~\ref{fig:lvm_VI_amortize}).
A motivating example, known as \textit{variational autoencoder (VAE)}, will be discussed in Section~\ref{section:vae_pca}.

As a result, rather than maintaining individual parameters $\blambda_n$ for each latent variable $\bz_n$, the model predicts the variational posterior for new data using the learned encoder---eliminating the need to re-optimize from scratch for each example. This dramatically improves computational efficiency and enables generalization to unseen data by reusing patterns extracted during training.
Consequently, amortized VI avoids the per-datapoint optimization loop and allows the use of efficient stochastic gradient-based methods for end-to-end learning. The graphical model for this amortized constrained framework is illustrated in Figure~\ref{fig:lvm_VI_amortize}.
Of course, amortization also has drawbacks. Since the true global structure of the posterior is typically unknown, a poorly designed encoder (e.g., one with insufficient capacity or inappropriate architecture) can lead to amortization bias---a mismatch between the true posterior and the class of distributions representable by the encoder.
One might hope to mitigate this by using highly flexible models such as deep neural networks to parameterize  $q(\bz_n\mid \bx_n)$. However, computing the Monte Carlo approximation of the ELBO requires evaluating 
$\ln q(\bz_n\mid \bx_n)$, which must remain tractable. This imposes a practical constraint: the variational family must be both expressive and log-density computable, limiting the choice of architectures or distributional forms.

\begin{problemset}
\item \label{problem:proptionto} Suppose that both $p(x)$ and $q(x)$ are probability density functions satisfying 
$$
p(x) \propto q(x).
$$ 
Show that $p(x) = q(x)$ for all $x$.

\item When discussing Bayesian matrix decomposition methods, each model will be represented by a graphical model---for example, the GGG model in Figure~\ref{fig:bmf_ggg}. Referring to Section~\ref{section:markov-blanket}, draw the graphical representations for Bayesian linear models with (i) a zero-mean prior, (ii) a semi-conjugate prior, and (iii) a fully conjugate prior.
For each model, discuss the Markov blanket of every node.

\item \label{problem:elbo_equa_em} 
In the EM algorithm (Section~\ref{section:em_uncons}), during iteration $t$, suppose we set $q^\toptone_{\bz_n}(\bz_n) = p(\bz_n \mid \bx_n, \btheta^\toptzero)$ for all $n=1,2,\ldots,N$ in the E-step. Show that the evidence lower-bound (ELBO) becomes tight, i.e., $\mathcalF(\{q_{\bz_n}(\bz_n)\}_{n=1}^N, \btheta^\toptzero)=
\mathcalL(\btheta^\toptzero) $.

\item \label{problem:forward_rever_KL} \textbf{Exclusive/inclusive KL.}
When fitting a parametric distribution $q_{\blambda}$ to a target distribution $p$ by minimizing the KL divergence $\KL[q_{\blambda} \parallel p]$ with respect to $\blambda$, this is known as \textit{reverse/exclusive KL} minimization. 
Show that this approach exhibits \textit{mode-seeking} (or \textit{zero-forcing}) behavior: the minimization forces $q_{\blambda}(x)=0$ wherever  $p(x)=0$, often causing $q_{\blambda}$ to concentrate around a single mode of $p$. 
Conversely, when  fitting $q_{\blambda} $ to $p$ by minimizing  $\KL[p\parallel q_{\blambda} ]$ with respect to $\blambda$, this is called \textit{forward/inclusive KL}. 
Show that this leads to \textit{mass-covering} (or \textit{mean-seeking}) behavior: $q_{\blambda}$ must assign non-negligible probability mass wherever $p$ does.

\item \label{problem:kl_vae}
\textbf{KL of Gaussians.}  Let  $q(\bx)=\normal( \bmu, \diag(\bsigma^2))$ and $p(\bx)=\normal( \bzero_D, \bI_D)$, where $\bmu, \bsigma^2\in\real^D$ (see Definition~\ref{sec:multi_gaussian_conjugate_prior}). Show that $\KL[q\parallel p]=\frac{1}{2}\sum_{n=1}^{D}(\mu_n^2+\sigma^2_n-\ln \sigma^2_n-1)$. 
This expression is commonly used as the KL regularization term in variational autoencoders (VAEs); see Equation~\eqref{equation:vae_elbo_vae}.

\item \textbf{KL of Gaussians.} Let  $p(x)=\normal(\mu_1, \sigma_1^2)$ and $q(x)=\normal(\mu_2,\sigma_2^2)$. Show that 
$$
\KL[p \parallel q ] = \ln\frac{\sigma_2}{\sigma_1} +\frac{\sigma_1^2+(\mu_1-\mu_2)^2}{2\sigma_2^2}-\frac{1}{2}.
$$ 
Now consider the multivariate case: let $\normal_1(\bx)=\normal(\bmu_1, \bSigma_1)$ and $\normal_2(\bx)=\normal(\bmu_2,\bSigma_2)$ (Definition~\ref{sec:multi_gaussian_conjugate_prior}). Show that 
$$
\KL[\normal_1\parallel \normal_2] = 
\frac{1}{2}\ln\abs{\bSigma_2\bSigma_1^{-1}}+\frac{1}{2}\trace\bSigma_2^{-1}
\big( (\bmu_1-\bmu_2)(\bmu_1-\bmu_2)^\top +\bSigma_1-\bSigma_2 \big).
$$
More generally, for an arbitrary distribution $p(\bx)$ and a multivariate Gaussian $\normal(\bx)=\normal(\bmu, \bSigma)$ with $\bx\in\real^D$, show that 
$$
\KL[p\parallel \normal ]=
\int \frac{1}{2}(\bx-\bmu)^\top\bSigma^{-1}(\bx-\bmu) d\bx 
+
\frac{1}{2}\ln \abs{\bSigma}+\frac{D}{2}\ln2\pi + \int p(\bx)\ln p(\bx) d\bx.
$$

\item \textbf{KL of Bernoullis.} 
Given two Bernoulli distributions $p(x)=\bernoulli(x\mid p)$ and $q(x)=\bernoulli(x\mid q)$ (Equation~\eqref{equation:bernoulli_distribution}), show that 
$$
\KL[p \parallel q] = p\ln\frac{p}{q} + (1-p) \ln \frac{1-p}{1-q}.
$$

\item \label{problem:entropy_mgau} \textbf{Entropy of Gaussians.} 
As introduced in Section~\ref{section:elbo_cfe}, \textit{entropy} is closely related to KL divergence. The entropy of a distribution $p(\bx)$ is defined as
\begin{equation}
\entropy[p(\bx)] = -\int p(\bx)\ln p(\bx) d\bx.
\end{equation}
For a multivariate Gaussian random vector $\rvx\sim \normal(\bmu,\bSigma)$  with $\rvx\in\real^D$ (Definition~\ref{sec:multi_gaussian_conjugate_prior}), show that 
$$
\entropy[\normal(\bmu,\bSigma)] = \frac{1}{2}\ln \abs{\bSigma}+ \frac{D}{2} \ln(2\pi e).
$$

\item \label{prob:mix_of_gauss} \textbf{Mixture of Gaussians.}
Consider a dataset $\mathcalX=\{\bx_1,\bx_2, \ldots,\bx_N\}\in\real^D$ generated from a Gaussian mixture model with $K$ components:
\begin{equation}
p(\bX\mid \bpi, \{\bmu_k\}, \{\bSigma_k\}) = 
\prod_{n=1}^{N}\sum_{k=1}^{K} \pi_k \normal(\bx_n \mid \bmu_k, \bSigma_k),
\end{equation}
where $\bX\in\real^{N\times D}$ contains the data points as rows, $\bmu_k\in\real^D$,  $\bSigma\in\real^{D\times D}$, and $\bpi=[\pi_1,\pi_2,\ldots,\pi_K]^\top$ are the \textit{mixing coefficients} satisfying $0\leq \pi_k\leq 1$  and $\sum_{k=1}^{K}\pi_k=1$.
Introduce a $K$-dimensional binary latent variable $\{\bz_n\}\in\{0,1\}^K$ for each sample $n=1,2,\ldots,N$, 
such that $z_{nk}\in\{0,1\}$ and $\sum_{k=1}^{K}z_{nk}=1$.
Using the EM algorithm (Algorithm~\ref{alg:em_alg}), derive the update equations for the parameters  $\bpi, \{\bmu_k\}, \{\bSigma_k\}$.
Show that the E-step computes the posterior responsibilities as
$$
\zeta_{nk} \leftarrow
\frac{\pi_k^{\text{old}}\normal({\bx_n}\mid \bmu_k^{\text{old}}, \bSigma_k^{\text{old}})}{\sum_{\ell=1}^{K}\pi_\ell^{\text{old}}\normal({\bx_n}\mid \bmu_\ell^{\text{old}}, \bSigma_\ell^{\text{old}})}, 
\quad \forall\, n=1,2,\ldots,N, k=1,2,\ldots,K,
$$
so that  $p(z_{nk}=1\mid \bx_n) = \zeta_{nk}$.
Then show that the M-step updates the parameters as follows:
\begin{align*}
N_k &\leftarrow \sum_{n=1}^{N}\zeta_{nk}; 
&&\bmu_k^{\text{new}} \leftarrow \frac{1}{N_k}\sum_{n=1}^{N} \zeta_{nk} \bx_n;\\
\bSigma_k^{\text{new}} &\leftarrow \frac{1}{N_k}\sum_{n=1}^{N} \zeta_{nk} (\bx_n-\bmu_k)(\bx_n-\bmu_k)^\top;
&&\pi_k^{\text{new}} \leftarrow \frac{N_k}{N}, 
\qquad \forall\,n, k.
\end{align*}
\textit{Hint: See Example~\ref{example:gmm_twoclus} for a special case with two cluster.}

\item \label{prob:mix_of_bern} \textbf{Mixture of Bernoullis.}
Following the same setup as in Problem~\ref{prob:mix_of_gauss}, but now assume each observation $\bx_n\in\{0,1\}^{D}$ is binary. Consider the Bernoulli mixture model:
\begin{equation}
p(\bX\mid \bpi, \{\btheta_k\}) = 
\prod_{n=1}^{N}\sum_{k=1}^{K} \pi_k \bernoulli(\bx_n \mid \btheta_k)
=\prod_{n=1}^{N}\sum_{k=1}^{K} \pi_k \prod_{d=1}^D\left(\theta_{kd}^{x_{nd}}  (1-\theta_{kd})^{(1-x_{nd})} \right),
\end{equation}
where $\btheta_k\in\real^D$ and each component satisfies $0\leq \theta_{kd}\leq 1$ for $d=1,2,\ldots,D$.
Again, introduce a $K$-dimensional binary latent variable $\{\bz_n\}\in\{0,1\}^K$ for each sample $n=1,2,\ldots,N$, 
with $z_{nk}\in\{0,1\}$ and $\sum_{k=1}^{K}z_{nk}=1$.
Show that the E-step computes
$$
\zeta_{nk} \leftarrow
\frac{\pi_k^{\text{old}} \bernoulli({\bx_n}\mid \btheta_k^{\text{old}})}{\sum_{\ell=1}^{K}\pi_\ell^{\text{old}} \bernoulli({\bx_n}\mid \btheta_\ell^{\text{old}})}, 
\quad \forall\, n=1,2,\ldots,N, k=1,2,\ldots,K,
$$
so that  $p(z_{nk}=1\mid \bx_n) = \zeta_{nk}$.
Then show that the M-step updates the parameters as
\begin{align*}
N_k &\leftarrow \sum_{n=1}^{N}\zeta_{nk}; 
\qquad \btheta_k^{\text{new}} \leftarrow \frac{1}{N_k}\sum_{n=1}^{N} \zeta_{nk} \bx_n;
\qquad \pi_k^{\text{new}} \leftarrow \frac{N_k}{N}, 
\quad \forall\,n, k.
\end{align*}
\end{problemset}

\chapter{Regular Probability Models and Conjugacy}\label{chapter:conjugate_models_bmf}
\begingroup
\hypersetup{
	linkcolor=structurecolor,
	linktoc=page,  
}
\minitoc \newpage
\endgroup

\lettrine{\color{caligraphcolor}I}
In this chapter, we explore several specific examples of probability distributions and their conjugate properties. These distributions are not only valuable in their own right but also serve as foundational building blocks for constructing more complex models. They will be used extensively throughout the book.
One key purpose of the distributions discussed here is to model the probability distribution $p(\bx)$ of a random variable $\rvx$, given a finite set of observations $\bx_1, \bx_2, \ldots, \bx_N$. 
This task is known as \textit{density estimation}. It's important to note that density estimation is inherently ill-posed, since there are infinitely many possible probability distributions that could fit the observed data. In fact, any distribution $p(\bx)$ that assigns nonzero probability to each of the data points $\bx_1, \bx_2, \ldots, \bx_N$ is a valid candidate. The challenge of selecting an appropriate distribution is closely related to the problem of model selection, and it remains a central concern in machine learning,  statistics, or Bayesian learning.

\section{Conjugate Priors}
In Section~\ref{sec:beta-bernoulli}, we briefly discussed  conjugate priors. 
Conjugate priors play a crucial role in Bayesian statistics, providing a convenient mathematical property that simplifies the computation of posterior distributions.
We now present the formal definition as follows.

\begin{definition}[Conjugate Prior]\index{Conjugate prior}
Given a family of likelihood functions $\{p(\mathcalX \mid \btheta): \btheta \in \bTheta\}$, 
a family of prior distributions $p_\omega (\btheta)$, indexed by hyper-parameters  $\bomega \in \bOmega$, is called a \textit{conjugate prior family} if for any $\bomega$ and any observed data $\mathcalX$, the resulting posterior distribution belongs to the same family---that is, it equals $p_{\bomega^\prime} (\btheta \mid \mathcalX)$ for some updated hyper-parameters $\bomega^\prime \in \bOmega$. 
\end{definition}

In Bayesian inference, a prior distribution represents our beliefs about the parameters of a statistical model before observing any data. 
A prior is said to be conjugate to a likelihood function if their combination results in a posterior distribution that belongs to the same family of distributions as the prior. This property facilitates analytical solutions, making the computation of the posterior more tractable.
As noted previously, a simple illustration is provided by the Beta-Bernoulli model.
\begin{example}[Beta-Bernoulli]
Suppose $\mathcalX=\{x_1, x_2, ..., x_N\}$ are i.i.d. samples from a \textit{Bernoulli distribution} with parameter $\theta$, i.e., $\bernoulli(x\mid \theta)=\theta^x {(1-\theta)}^{1-x}$, where $x\in\{0,1\}$.
$\betadist(\theta \mid a,b)$ distribution~\footnote{see Definition~\ref{definition:dirichlet_dist}, which is a special case of the Dirichlet distribution.}, with $a, b >0$, is conjugate to $\bernoulli(x\mid\theta)$, since the posterior density is $p(\theta \mid \mathcalX) = \betadist(\theta \mid a+\sum_n x_n, b+N-\sum_n x_n)$. 
\end{example}

Conjugate priors enable computationally efficient Bayesian reasoning and offer an intuitive interpretation: they can be viewed as encoding real or hypothetical prior observations. Generally, models with conjugate priors are popular for two main reasons \citep{robert2007bayesian, bernardo2009bayesian, hoff2009first, gelman2013bayesian}:
\begin{itemize}
\item They often yield closed-form expressions for the posterior distribution.

\item They are easy to interpret, as the effect of observed data on the posterior is reflected directly in updated parameter values.
\end{itemize}

\index{Conditional probability}
\index{Marginal probability distribution}
\section{Random Variables and Distribution}
A random variable is a mathematical construct that models uncertain outcomes by taking on different values according to some underlying probability mechanism.
We denote a random variable using a lowercase letter in \textit{normal fonts} (e.g., $\rx$), while its possible realizations are denoted with lowercase letters in \textit{italic fonts} (e.g., $x$).  
For instance, $y_1$ and $y_2$ are  possible values of the random variable $\ry$. 
For vector-valued variables, we write the random variable as $\rvy$ (in normal fonts) and a specific realization as $\by$ (in italic fonts).  
Similarly, for matrix-valued variables, we may use $\rmY$  (in normal fonts) for the random matrix and $\bY$ (in italic fonts) for a realization.

Random variables can be either discrete or continuous. A \textit{discrete} random variable takes values in a finite or countably infinite set; these values need not be numeric (e.g., categories like ``red," ``green," ``blue"). A \textit{continuous} random variable, by contrast, takes values in an uncountable set, typically a subset of the real numbers.

Crucially, a random variable alone does not specify probabilities---it must be paired with a probability distribution that assigns likelihoods to its possible states.
A probability distribution quantifies how probability mass or density is assigned across the possible values of one or more random variables. The form of this description depends on whether the variables are discrete or continuous.

For discrete random variables, we use a \textit{probability mass function (p.m.f., PMF)}. Probability mass functions are denoted by a capital $\prob$.
The PMF maps each possible state of a random variable to its probability.
For instance, the notation $\prob(\ry = y)$ denotes the probability that the random variable $\ry$  takes the value $y$, where probabilities lie in [0,1], with 1 indicating certainty and 0 impossibility.
Alternatively, we may first define a random variable and then specify its distribution using the ``is distributed as" notation: $\ry \sim \prob(y)$.

PMFs can also describe multiple variables jointly, forming a \textit{joint probability distribution}. 
For example, $\prob(\rx = x, \ry = y)$ gives the probability that $x = x$ and $y = y$ occur  simultaneously. This is often abbreviated as $\prob(x, y)$. 
If the distribution depends on known parameters $\balpha$, we write $\prob(x, y \mid \balpha)$ for brevity.

For continuous random variables, we use a \textit{probability density function (p.d.f., PDF)} instead of  a probability mass function, typically denoted by $p$ or $f$ (lowercase). 
A function $p(y)$ qualifies as a PDF if it satisfies the following:
\begin{itemize}
\item Its domain includes all possible values of the random variable $\ry$;
\item We do not require $p(y) \leq  1$ as that in the PMF. However, it must satisfies that $\forall y\in\ry, p(y)\geq 0$.
\item It integrates to one: $\int p(y)\,dy = 1$.
\end{itemize}
Importantly, unlike a PMF, a PDF does not give probabilities directly. Instead, the probability that $\ry$ falls within an infinitesimal interval of width $\delta y$ around $y$ is approximately $p(y)\delta y$.
Moreover, if the probability density function depends on some known parameters $\balpha$, it is commonly written as $p(x\mid \balpha)$, $f(x; \balpha)$,  $f_{\rx}(x; \balpha)$, or $f_{\rx}(x)$  for brevity.

In many applications, we are interested in the probability of one event given that another has occurred. This is known as \textit{conditional probability}. 
The conditional probability that $\rx = x$ given $\ry = y$ is denoted by  $\prob(\rx = x \mid  \ry = y)$
and is computed as
$$
\prob(\rx = x \mid  \ry = y) = \frac{\prob(\rx = x, \ry = y)}{\prob(\ry = y)}.
$$
This relationship underpins Bayes' theorem (see Equation~\eqref{equation:posterior_abstract_for_mcmc}).

Conversely, when the joint distribution over a set of variables is known, we may wish to find the distribution over a subset of those variables. This is called the \textit{marginal probability distribution}.
For discrete variables $\rx$ and $\ry$, the marginal distribution of $\rx$  is obtained by summing over all possible values of $\ry$:
$$
\prob(\rx=x) = \sum_{y} \prob(\rx=x, \ry=y).
$$
For continuous variables, summation is replaced by integration:
$$
p(x) = \int p(x, y)\,dy.
$$

\section{Regular Univariate Models and Conjugacy}


\begin{table}[H]
\centering
\setlength{\tabcolsep}{5.7pt}
\begin{tabular}{l|l|l}
\hline
\hline

\hyperref[definition:gaussian_distribution]{Gaussian, p.~\pageref{definition:gaussian_distribution}} & 
\hyperref[definition:gamma-distribution]{Gamma, p.~\pageref{definition:gamma-distribution}} & 
\hyperref[equation:student_t_dist]{Student's $t$, p.~\pageref{equation:student_t_dist}}
 \\ \hline
\hyperref[{definition:inverse_gamma_distribution}]{Inverse-Gamma, p.~\pageref{definition:inverse_gamma_distribution}} & \hyperref[definition:truncated_normal]{Truncated-Normal, p.~\pageref{definition:truncated_normal}}& 
\hyperref[definition:inverse_gaussian_distribution]{Inverse-Gaussian, p.~\pageref{definition:inverse_gaussian_distribution}} 
\\ \hline
\hyperref[definition:chisquare_distribution]{Chi-Square, p.~\pageref{definition:chisquare_distribution}}  &   \hyperref[definition:normal_inverse_gamma]{Normal-Inv-Gamma, p.~\pageref{definition:normal_inverse_gamma}}& 
\hyperref[definition:inverse-chi-square]{Inverse-Chi-Squared, p.~\pageref{definition:inverse-chi-square}}
 \\ \hline
\hyperref[definition:normal_inverse_chi_square]{Nor-Inv-Chi-Squared, p.~\pageref{definition:normal_inverse_chi_square}}   & \hyperref[definition:general_truncated_normal]{General-Truncated-Nor, p.~\pageref{definition:general_truncated_normal}} & \hyperref[definition:half_normal]{Half-Normal, p.~\pageref{definition:half_normal}}    \\  \hline
\hyperref[definition:laplace_distribution]{Laplace, p.~\pageref{definition:laplace_distribution}}
& \hyperref[definition:skew_laplace_distribution]{Skew-Laplace, p.~\pageref{definition:skew_laplace_distribution}}  & \hyperref[definition:reftified_normal_distribution]{Rectified-Normal, p.~\pageref{definition:reftified_normal_distribution}} \\ \hline
\hyperref[definition:multinomial_dist]{(Bi) Multinomial, p.~\pageref{definition:multinomial_dist}}
&
\hyperref[definition:dirichlet_dist]{Dirichlet, p.~\pageref{definition:dirichlet_dist}}
&
\hyperref[definition:poisson_distribution]{Poisson, p.~\pageref{definition:poisson_distribution}} \\ \hline
\hyperref[definition:exponential_distribution]{Exponential, p.~\pageref{definition:exponential_distribution}} 
& 
\hyperref[definition:multivariate_gaussian]{Multi Gaussian, p.~\pageref{definition:multivariate_gaussian}}
& 
\hyperref[definition:multivariate-stu-t]{Multi Student's $t$, p.~\pageref{definition:multivariate-stu-t}}
\\ \hline
\hyperref[definition:wishart_dist]{Wishart, p.~\pageref{definition:wishart_dist}}
&
\hyperref[definition:multi_inverse_wishart]{Inverse-Wishart, p.~\pageref{definition:multi_inverse_wishart}}
&
\hyperref[definition:normal_inverse_wishart]{Nor-Inv-Wishart, p.~\pageref{definition:normal_inverse_wishart}} 
\\ \hline 
\hyperref[sec:beta-bernoulli]{Beta, p.~\pageref{sec:beta-bernoulli}, p.~\pageref{definition:dirichlet_dist}} & 
\hyperref[sec:beta-bernoulli]{Bernoulli, p.~\pageref{sec:beta-bernoulli}, p.~\pageref{remark:binomial_dist}} & 
\hyperref[remark:binomial_dist]{Multinoulli, p.~\pageref{remark:binomial_dist}} 
\\ \hline
\hline
\end{tabular}
\caption{Links and page references for common distributions.}
\label{table:common_distributions}
\end{table}

In most of our Bayesian matrix decomposition models, we formulate the likelihoods using univariate distributions. However, in the Gaussian case, we also make use of multivariate distributions.
In this section, we provide rigorous definitions of common univariate probability distributions and their conjugate priors.
Table~\ref{table:common_distributions} summarizes the topics covered here.
Additionally, with special attention to its unique properties, we discuss the multivariate Gaussian distribution and its conjugacy in the following section.

\begin{figure}[h]
\centering  
\vspace{-0.35cm} 
\subfigtopskip=2pt 
\subfigbottomskip=2pt 
\subfigcapskip=-5pt 
\subfigure[Gaussian PDFs.]{\label{fig:dists_gaussian}
\includegraphics[width=0.481\linewidth]{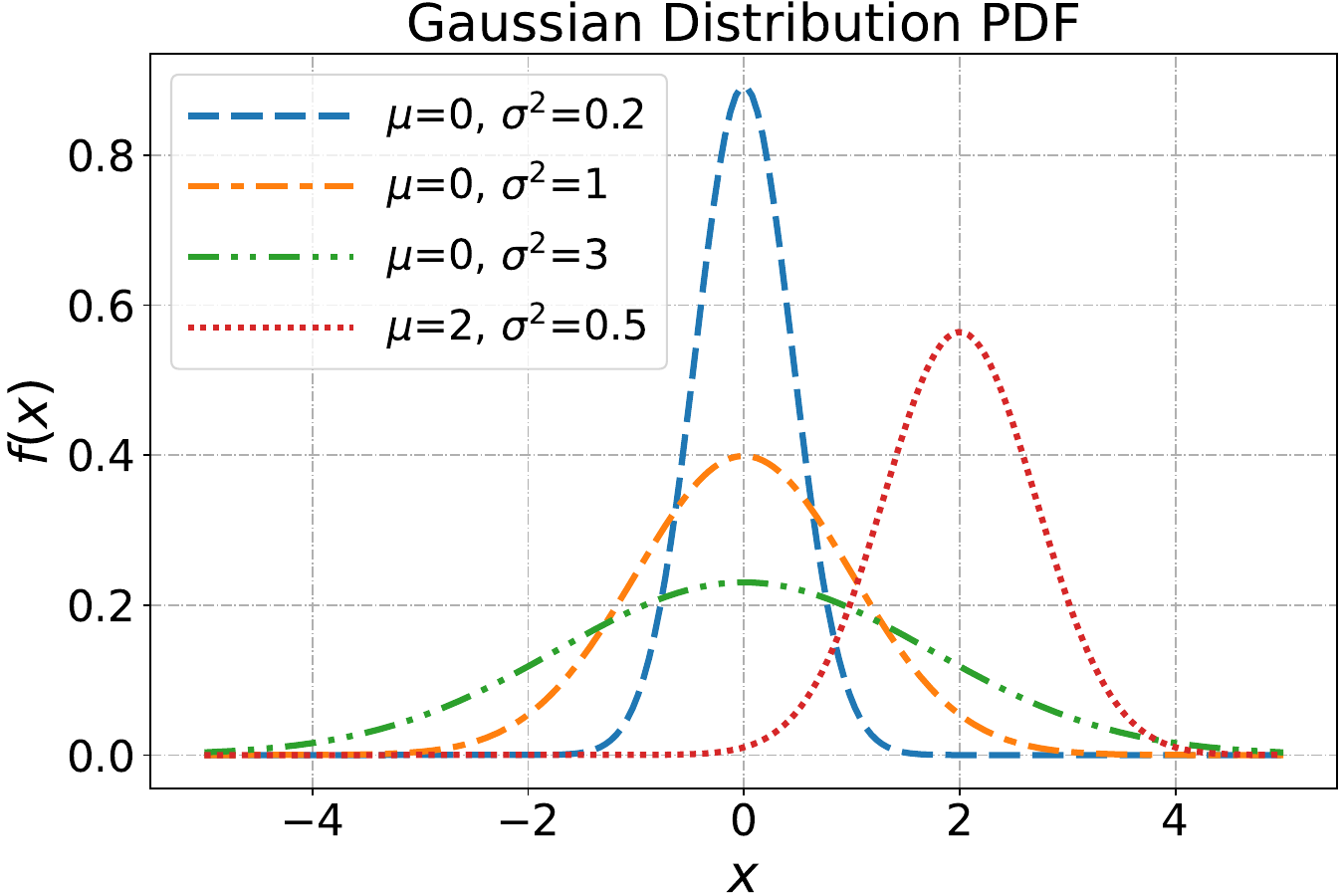}}
\subfigure[Gaussian CDFs.]{\label{fig:dists_gaussian_cdf}
\includegraphics[width=0.481\linewidth]{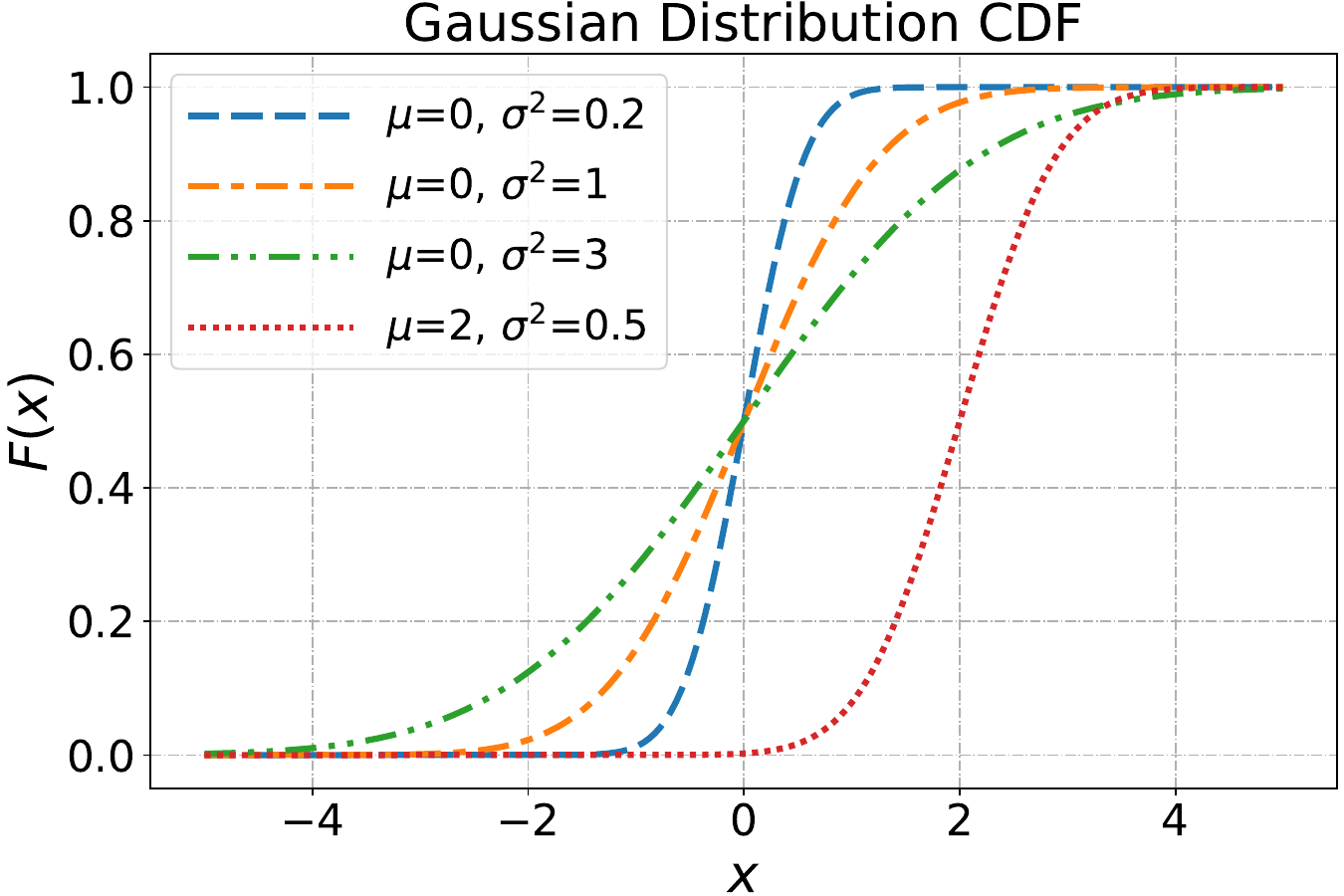}}
\caption{Gaussian probability density functions and cumulative distribution functions for different values of the mean and variance parameters $\mu$ and $\sigma^2$.}
\label{fig:dist_gaussians_all}
\end{figure}

\index{Gaussian distribution}
\begin{definition}[Gaussian or Normal Distribution]\label{definition:gaussian_distribution}
A random variable $\rx$ is said to follow a \textit{Gaussian distribution} (or a \textit{normal distribution}) with mean $\mu$ and variance  $\sigma^2>0$, denoted  $\rx \sim \normal(\mu,\sigma^2)$ \footnote{Note if two random variables $\ra$ and $\rb$ have the same distribution, then we write $\ra \sim \rb$.}, if its probability density function is
$$
f(x; \mu,\sigma^2)=\frac{1}{\sqrt{2\pi\sigma^2}} \exp \left\{-\frac{1}{2\sigma^2 }(x-\mu)^2 \right\}
=\sqrt{\frac{\tau}{2\pi}}\exp \left\{ -\frac{\tau}{2}(x-\mu)^2 \right\},
$$
where $\tau=1/\sigma^2$ is called the \textit{precision}.
The mean and variance of $\rx \sim \normal( \mu,\sigma^2)$ are given by 
$$
\Exp[\rx] = \mu, \qquad \Var[\rx] =\sigma^2=\tau^{-1}.
$$
The cumulative distribution function (c.d.f., CDF) of Gaussian is 
$$
F(x; \mu, \sigma^2) = \prob(\rx<x) = \frac{1}{\sqrt{2\pi\sigma^2}} \int_{-\infty}^x \exp\left\{-\frac{1}{2\sigma^2 }(z-\mu)^2 \right\}dz.
$$
We denote the CDF of the standard normal distribution $\normal(0,1)$ by $\Phi(y) = \int_{-\infty}^{y} \normal(u\mid 0,1)\,du= \frac{1}{\sqrt{2\pi}} \int_{-\infty}^{y} \exp(-\frac{u^2}{2}) \,du $.
Figure~\ref{fig:dist_gaussians_all} illustrates how the shape of the Gaussian distribution changes with different values of $\mu$ and $\sigma^2$.
\end{definition}

Suppose $\mathcalX=\{x_1, x_2, ..., x_N\}$ are drawn i.i.d. from a Gaussian distribution of $\normal(x\mid \mu, \sigma^2)$.
For conjugate Bayesian analysis, we may rewrite the joint likelihood as:
\begin{equation}\label{equation:uni_gaussian_likelihood}
\begin{aligned}
p(\mathcal{X} \mid \smu, \ssigma^2) &= \prod^N_{n=1} \mathcal{N} (x_n\mid\smu, \ssigma^2) \\
&= (2\pi)^{-N/2}  (\ssigma^2)^{-N/2} \exp\left\{-\frac{1}{2 \ssigma^2}  \left[ N(\widebarx - \smu)^2 + N \sum_{n=1}^N(x_n - \widebarx)^2   \right] \right\}  \\
&= (2\pi)^{-N/2}  (\ssigma^2)^{-N/2} \exp\left\{-\frac{1}{2 \ssigma^2}  \left[  N(\widebarx - \smu)^2 +  N S_{\widebarx} \right] \right\},
\end{aligned}
\end{equation}
where $S_{\widebarx}=\sum_{n=1}^N(x_n - \widebarx)^2$ and $\widebarx = (\sum_{n=1}^{N}x_n)/N$ is the sample mean.
This form is particularly useful when deriving the conditional posterior under a \textit{normal-inverse-Gamma} prior (see Equation~\eqref{equation:conjugate_nigamma_general}). 

\index{Canonical form}
With fixed  mean $\mu$ and variance $\sigma^2$ parameters, the Gaussian density can be expressed in its \textit{canonical form}:
\begin{equation}\label{equation:gaussian_form_conform}
\begin{aligned}
	p(x\mid \mu, \sigma^2) &=\normal(x \mid \mu, \sigma^2)
	\propto \exp\left\{ -\frac{1}{2\sigma^2} x^2 + \frac{\mu}{\sigma^2} x \right\},
\end{aligned}
\end{equation}
where ``$\propto$" means ``proportional to."
Thus, if a density can be written in this form, the corresponding random variable follows a Gaussian distribution $\rx\sim \normal(\mu, \sigma^2)$. 
See, for example, the posterior derivation in the Bayesian \textit{GGG} matrix decomposition model (Equation~\eqref{equation:ggg_poster_wmk1}).

While the product of two Gaussian random variables does not generally yield a standard distribution, linear combinations of Gaussian variables remain Gaussian.
\begin{remark}[Sum of Gaussians]
Let $\rx$ and $\ry$ be two Gaussian distributed variables with means $\mu_x, \mu_y$ and variance $\sigma_x^2, \sigma_y^2$, respectively.
\begin{itemize}
\item If $\rx$ and $\ry$ are uncorrelated, then:
$$
\begin{aligned}
\rx+\ry &\sim \normal(\mu_x+\mu_y, \sigma_x^2+\sigma_y^2);\\
\rx-\ry &\sim \normal(\mu_x-\mu_y, \sigma_x^2+\sigma_y^2).
\end{aligned}
$$
More generally, for independent 
$ \rx_n \sim \normal(\mu_n, \sigma_n^2),  n = 1, 2,\ldots, N, $, and constants   $ \{a_n\}  $,
$$ 
\sum_{n=1}^{N} a_n \rx_n \sim \normal\left( \sum_{n=1}^{N} a_n \mu_n, \sum_{n=1}^{N} (a_n \sigma_n)^2 \right). 
$$
\item If $\rx$ and $\ry$ have correlation $\rho $, then:
$$
\begin{aligned}
\rx+\ry &\sim \normal(\mu_x+\mu_y, \sigma_x^2+\sigma_y^2+2\rho\sigma_x\sigma_y);\\
\rx-\ry &\sim \normal(\mu_x-\mu_y, \sigma_x^2+\sigma_y^2-2\rho\sigma_x\sigma_y).
\end{aligned}
$$
\end{itemize}
Further properties of Gaussian distributions, especially in the multivariate setting, are discussed in Section~\ref{section:multi_gauss}.
\end{remark}

\index{Normal-Normal model}
\paragrapharrow{Conjugate prior for mean of a Gaussian distribution and Normal-Normal model.}
When the variance is known, the Gaussian distribution serves as a conjugate prior for the mean parameter of a Gaussian likelihood. 
Specifically, suppose $\mathcalX=\{x_1, x_2, \ldots, x_N\}$ are i.i.d. normal with mean $\theta$ and precision $\lambda$, i.e., the likelihood is $\normal(x_n \mid \theta, \lambda^{-1})$, where the variance $\sigma^2=\lambda^{-1}$ is fixed, and $\theta$ is given a $\normal(\mu_0, \lambda^{-1}_0)$ prior: $\theta \sim \normal( \mu_0, \lambda^{-1})$. Using Bayes' theorem, ``posterior $\propto$ likelihood $\times$ prior,"  the posterior density is 
$$
p(\theta \mid \mathcalX ) \propto \prod_{n=1}^{N} \normal(x_n \mid \theta, \lambda^{-1}) 
\times \normal(\theta \mid \mu_0, \lambda_0^{-1})\propto \normal(\theta \mid \widetilde{\mu}, \widetilde{\lambda}^{-1}),
$$
where, with the sample mean $\widebarx = (\sum_{n=1}^{N} x_n)/N$, 
\begin{equation}\label{equation:posterior-param-normal-normal}
\begin{aligned}
\widetilde{\mu}
&= \frac{\lambda_0 \mu_0 + \lambda \sum_{n=1}^{N} x_n}{ \lambda_0 + N\lambda}
= \frac{\lambda_0}{\lambda_0 + N\lambda} \mu_0 + \frac{N\lambda}{\lambda_0 + N\lambda} \widebarx, \\ 
\widetilde{\lambda} &= \lambda_0 + N \lambda. 
\end{aligned}
\end{equation}
Thus, the posterior mean is a weighted mean of the prior mean $\mu_0$ and the sample mean $\widebarx$; the posterior precision is the sum of the prior precision and sample precision $N\lambda$.
This setup---known as the \textit{Normal-Normal model}---demonstrates that the Gaussian distribution is self-conjugate for the mean when the variance is fixed. It has been used, for instance, to model bimodal phenomena such as human height distributions \citep{schilling2002human}.

\begin{figure}[h]
\centering  
\vspace{-0.35cm} 
\subfigtopskip=2pt 
\subfigbottomskip=2pt 
\subfigcapskip=-5pt 
\subfigure[Student's $t$ distribution by varying parameter $\nu$. When $\nu=100$, the distribution is very close to a Gaussian distribution.]{\label{fig:dists_studentt_varNU}
\includegraphics[width=0.481\linewidth]{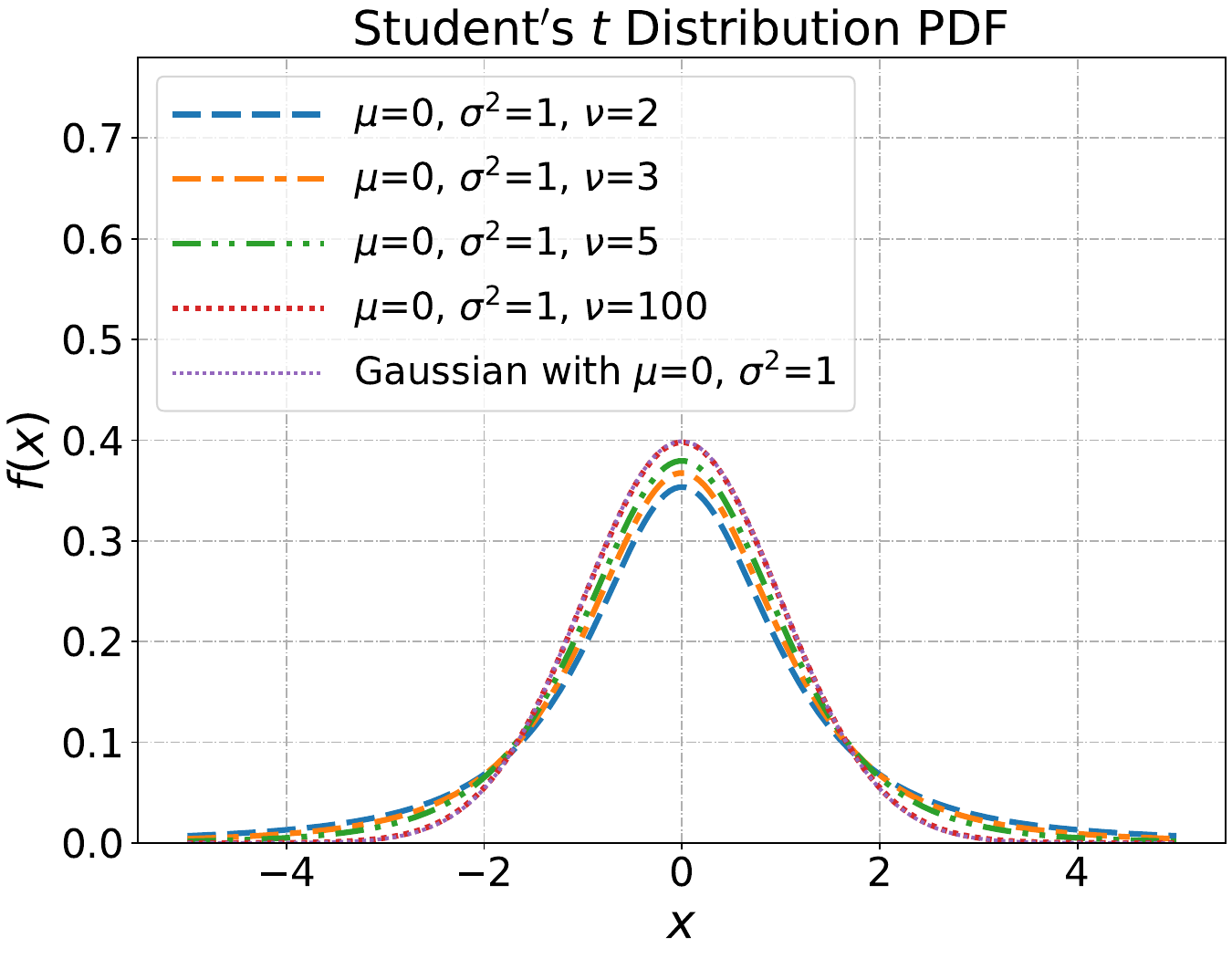}}
\subfigure[Student's $t$ distribution  by varying parameter $\sigma^2$.]{\label{fig:dists_studentt_varVar}
\includegraphics[width=0.481\linewidth]{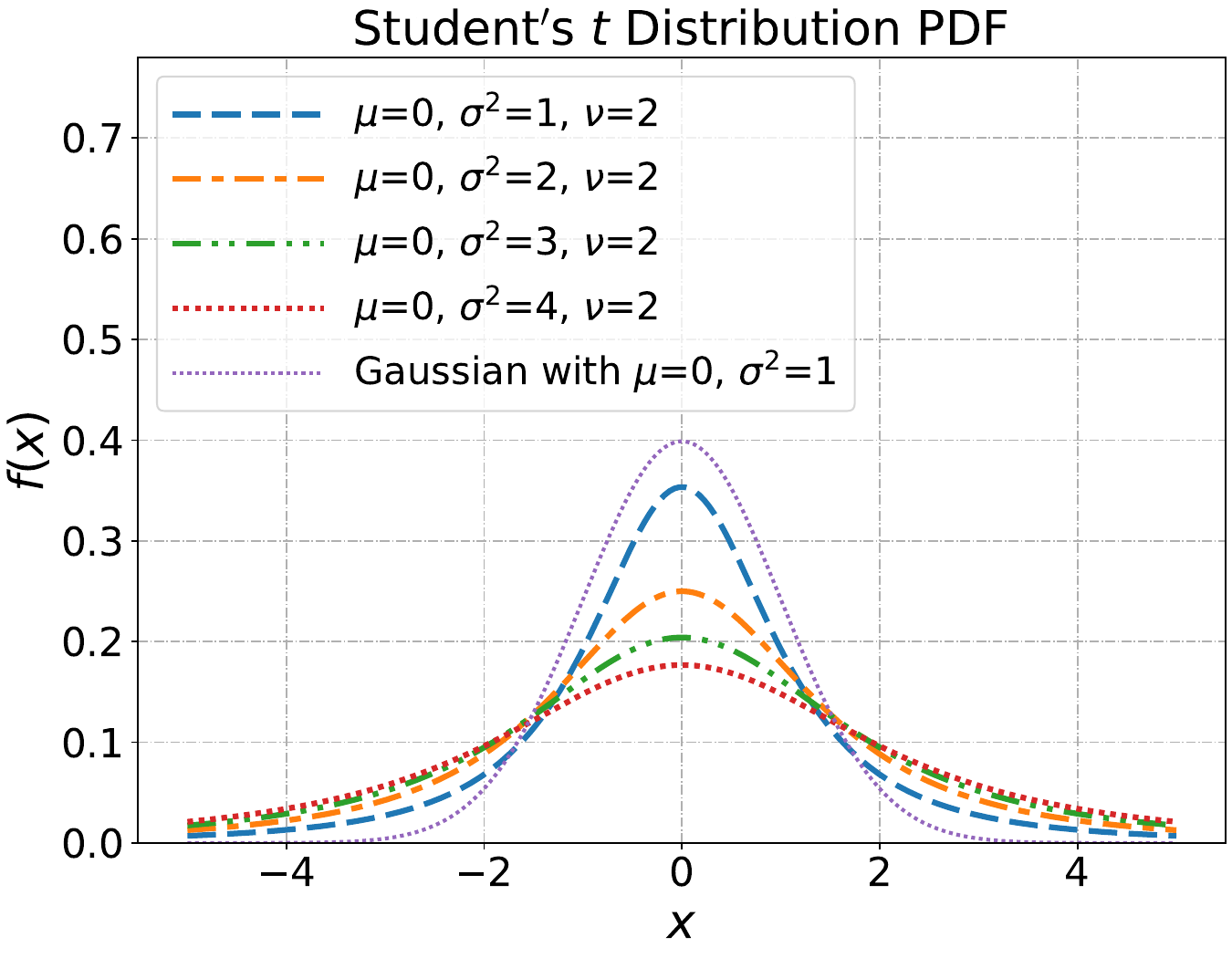}}
\caption{Student's $t$ distribution for different values of the parameters $\nu$ and $\sigma^2$.}
\label{fig:dist_students_all}
\end{figure}

\begin{definition}[Student's $t$ Distribution\index{Student's $t$ distribution}]\label{equation:student_t_dist}
A random variable $\rx$ is said to follow a \textit{Student's $t$ distribution} with parameters $\mu$, $\sigma^2>0$, and $\nu$, denoted $\rx \sim \studentt(\mu,\sigma^2, \nu)$,  if its density is
\begin{equation}\label{equation:stut_defini}
\begin{aligned}
f(x; \mu, \sigma^2, \nu)&= \frac{\Gamma(\frac{\nu+1}{2})}{\Gamma(\frac{\nu}{2})} \frac{1}{\sigma\sqrt{\nu\pi}} \times \left[ 1+ \frac{(x-\mu)^2}{\nu \sigma^2}   \right]^{-(\frac{\nu+1}{2})},
\end{aligned}
\end{equation}
where $\mu$ is the mean (location) parameter, $\sigma^2$ is called the \textit{scale parameter}, and $\nu$ is the \textit{degrees of freedom} (which controls the tail behavior).
The distribution has fatter tails than a Gaussian distribution and the degrees of freedom control the shape of the distribution. 
Smaller values of $\nu$ produce fatter tails, and as  $\nu\rightarrow \infty$, the distribution converges to $\normal(\mu, \sigma^2)$.
A notable special case  of the Student's $t$ distribution is the \textit{Cauchy distribution}, 
obtained when $\nu=1$:
\begin{equation}
\rx\sim \cauchydist(\mu, \sigma^2)
\quad\text{if}\quad \rx\sim \studentt(\mu, \sigma^2, 1).
\end{equation}
(Because of its extremely heavy tails, the Cauchy distribution has no defined mean or variance.).
For $\rx \sim \studentt(\mu,\sigma^2, \nu)$, the moments are:
$$ \Exp[\rx]=\left\{
\begin{aligned}
	&\mu , \, &\mathrm{if\,} \nu >1; \\
	&\text{undefined}, \, &\mathrm{if\,} \nu \leq 1.
\end{aligned}
\right.\qquad
\Var[\rx]=\left\{
\begin{aligned}
	&\frac{\nu}{\nu-2}\sigma^2, \, &\mathrm{if\,} \nu >2; \\
	&\infty, \, &\mathrm{if\,} 1<\nu\leq 2.
\end{aligned}
\right.
$$
Figure~\ref{fig:dist_students_all} shows how the shape of the Student's $t$ distribution varies with $\mu,\sigma^2, \nu$.
\end{definition}

\begin{remark}[Standard $t$ Distribution]
The form above is often called the \textit{non-standard $t$-distribution}. 
The \textit{standard $t$-distribution}, denoted  $\rx\sim\standardstudentt(\nu)$, corresponds to $\mu=0$ and $\sigma^2=1$:
$$
\rx\sim \standardstudentt(\nu) = \studentt(\mu=0, \sigma^2=1,\nu).
$$
\end{remark}
\begin{exercise}[Obtain Standard $t$ from Gaussians]
Suppose i.i.d. variables $\rx_n\sim \normal(\mu, \sigma^2),\forall n\in\{1,2,\ldots,N\}$, $\overline{\rx}=\frac{1}{N}\sum_{n=1}^{N}\rx_n$ (sample mean), and $\rS^2 = \frac{1}{N-1}\sum_{n=1}^{N} (\rx_n-\overline{\rx})^2$ (unbiased estimator of variance).
Show that the following $t$ variable follows the standard $\rt$ distribution with $N-1$ degrees of freedom (usually called the t-statistic):
$$
\rt=\frac{\overline{\rx}-\mu}{\rS/\sqrt{N}}\sim \standardstudentt(N-1).
$$
Note, on the other hand, the following variable is called the standardization of Gaussian variables:
$$
\rz=\frac{\overline{\rx}-\mu}{\sigma/\sqrt{N}}\sim \normal(0,1).
$$
\end{exercise}

In Figure~\ref{fig:dists_studentt_varNU}, we vary the $\nu$ parameter for the Student's $t$ distribution.
As $\nu$ decreases, the distribution becomes more spread out, leading to fatter tails compared to a Gaussian distribution. This allows for more flexibility in modeling data with greater uncertainty or \textit{outliers} since the Student's $t$ distribution has a greater probability of observing extreme values. In Bayesian modeling, the Student's $t$ distribution is often used as a prior for the mean parameter of a Gaussian likelihood, allowing for estimation of both the mean and precision of the data. This results in a \textit{Student's $t$-Normal} model.


\begin{figure}[h]
	\centering  
	\vspace{-0.35cm} 
	\subfigtopskip=2pt 
	\subfigbottomskip=2pt 
	\subfigcapskip=-5pt 
	\subfigure[Gamma distribution.]{\label{fig:dists_gamma}
		\includegraphics[width=0.481\linewidth]{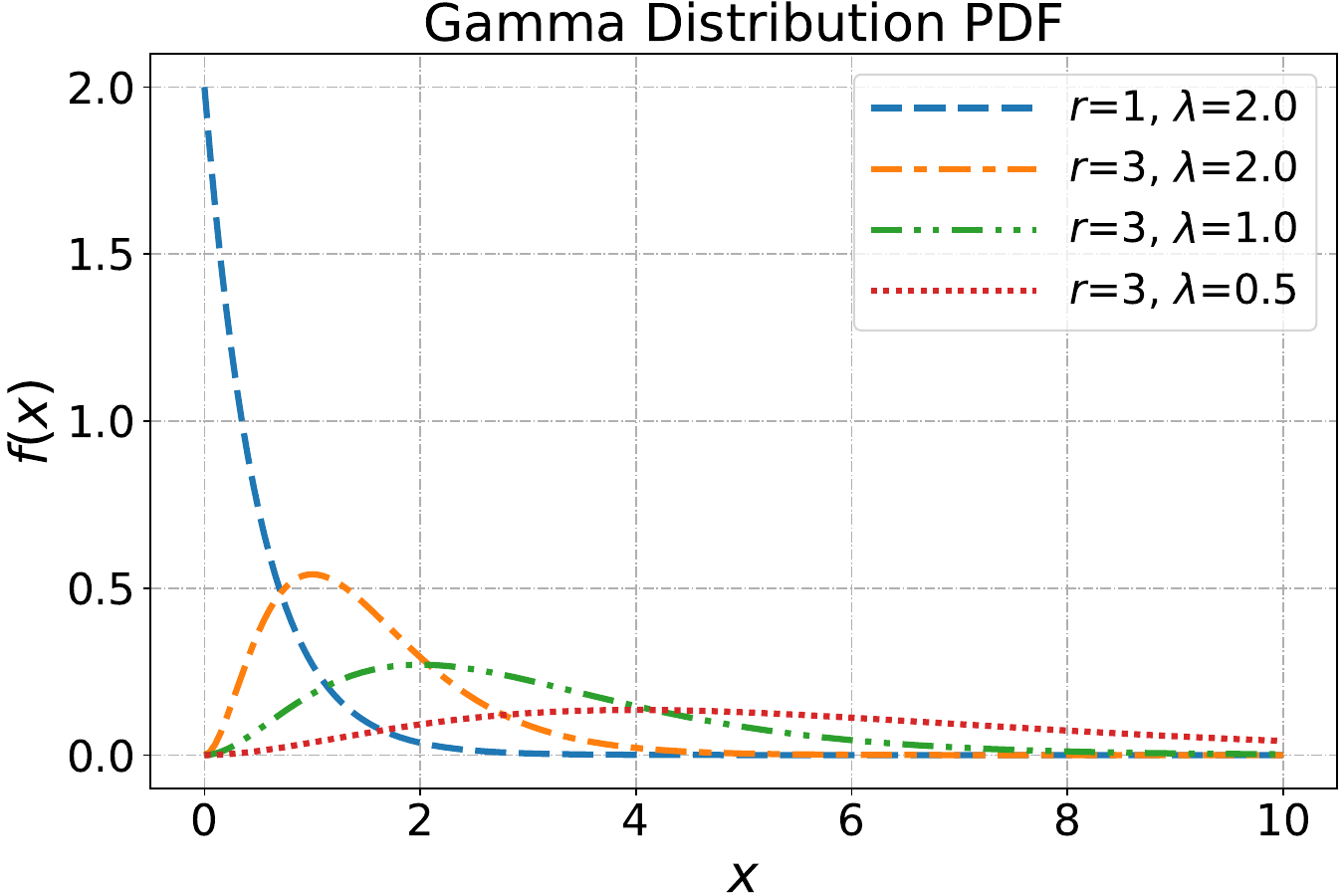}}
	\subfigure[Inverse-Gamma distribution.]{\label{fig:dists_inversegamma}
		\includegraphics[width=0.481\linewidth]{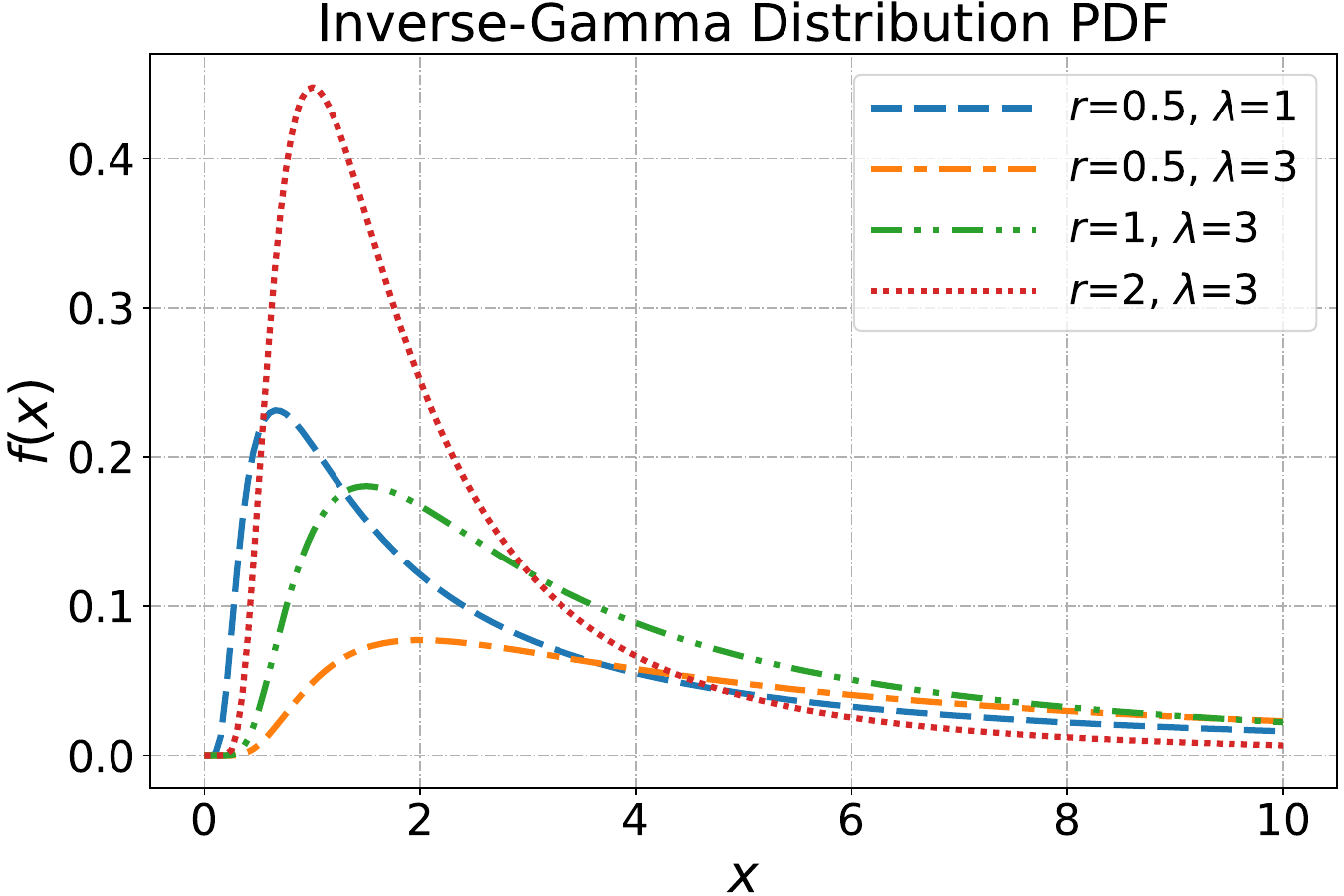}}
	\caption{Gamma and inverse-Gamma probability density functions for different values of the parameters $r$ and $\lambda$.}
	\label{fig:dists_gamma_inversegamma}
\end{figure}

\begin{definition}[Gamma Distribution\index{Gamma distribution}]\label{definition:gamma-distribution}
A random variable $\rx$ is said to follow the \textit{Gamma distribution} with shape parameter $r>0$ and rate parameter $\lambda>0$ \footnote{Note the inverse rate parameter $1/\lambda$ is called the scale parameter.
In probability theory and statistics, 
the \textbf{location} parameter shifts the entire distribution left or right, e.g., the mean parameter of a Gaussian distribution; 
the \textbf{shape} parameter compresses or stretches the entire distribution;
the \textbf{scale} parameter changes the shape of the distribution in some manner.
}, denoted $\rx \sim \gammadist(r, \lambda)$, if 
$$ f(x; r, \lambda)=\left\{
\begin{aligned}
&\frac{\lambda^r}{\Gamma(r)} x^{r-1} \exp(-\lambda x) ,& \mathrm{\,\,if\,\,} x \geq 0; \\
&0 , &\mathrm{\,\,if\,\,} x <0,
\end{aligned}
\right.
$$
where $\Gamma(x)=\int_{0}^{\infty} t^{x-1}\exp(-t)\,dt$ is the Gamma function,  and we can just take it as a function to normalize the distribution into sum to 1. In special cases when $y$ is a positive integer, $\Gamma(y) = (y-1)!$.
The mean and variance of $\rx \sim \gammadist(r, \lambda)$ are given by 
\begin{equation}
\Exp[\rx] = \frac{r}{\lambda}, \qquad \Var[\rx] = \frac{r}{\lambda^2}. \nonumber
\end{equation}
An important property of the Gamma distribution is its \textit{additivity}: let $\rx_1, \rx_2, \ldots, \rx_N$ be i.i.d. random variables drawn from $\gammadist(r_n, \lambda)$ for each $n \in \{1, 2, \ldots, N\}$. Then $\ry = \sum_{n=1}^N \rx_n$ is a random variable following from $\gammadist(\sum_{n=1}^N r_n, \lambda)$.
Figure~\ref{fig:dists_gamma} compares different parameters $r, \lambda$ for the Gamma distribution.
\end{definition}

It's crucial to note that the definition of the Gamma distribution does not restrict $r$ to be a natural number, and it allows  $r$ to take any positive number.
However, when $r$ is a positive integer, the Gamma distribution can be interpreted as a sum of $r$ exponentials of rate $\lambda$ (see Definition~\ref{definition:exponential_distribution}).
The summation property holds true more generally for Gamma variables with the same rate parameter. If $\rx_1$ and $\rx_2$ are random variables drawn from $\gammadist(r_1, \lambda)$ and $\gammadist(r_2, \lambda)$, respectively, then their sum $\rx_1+\rx_2$ is a Gamma random variable from $\gammadist(r_1+r_2, \lambda)$.

\index{Integration by parts}
In the Gamma distribution definition, we observe that the Gamma function can be defined by setting the rate parameter to 1, as follows:
$$
\Gamma(y) = \int_{0}^{\infty} x^{y-1} e^{-x} \, dx,  \qquad y\geq 0.
$$
Utilizing integration by parts $\int_{a}^{b} u(x) v^\prime(x) \, dx = u(x)v(x)|_a^b - \int_a^b u^\prime(x) v(x)\, dx$, where $u(x) =x^{y-1}$ and $v(x)=-e^{-x}$, we derive 
$$
\begin{aligned}
\Gamma(y) &= -x^{y-1}e^{-x}|_0^{\infty} - \int_{0}^{\infty} (y-1)x^{y-2}(-e^{-x})\, dx \\
&= 0+ (y-1) \int_{0}^{\infty}x^{y-2}e^{-x}\, dx = (y-1)\Gamma(y-1).
\end{aligned}
$$ 
This recurrence relation implies that for any positive integer $y$, $\Gamma(y) = (y-1)!$, as claimed.

\index{Gamma-Gamma model}
\paragrapharrow{Conjugate prior for rate  of a Gamma distribution and Gamma-Gamma model.}
The Gamma distribution is a conjugate prior of the \textit{rate} parameter of another Gamma likelihood when the shape parameter is known. 
Suppose we observe data $\mathcalX =\{x_1,x_2,\ldots,x_N\}$ drawn i.i.d. from a Gamma distribution $\rx_n\sim \gammadist(r, \lambda)$, where $r$ is fixed and $\lambda$ is unknown. 
Suppose further the rate parameter is given a Gamma prior $\lambda \sim \gammadist(a, \frac{a}{b})$. Using Bayes' theorem, the posterior density is 
$$
\begin{aligned}
p(\lambda \mid \mathcalX) &\propto 
\prod_{n=1}^{N}\gammadist(x_n \mid r, \lambda) \times 
\gammadist(\lambda \mid a, \frac{a}{b})
\propto 
\prod_{n=1}^{N} \frac{\lambda^r}{\Gamma(r)} x_n^{r-1} \exp(-\lambda x_n)
\times 
\frac{(\frac{a}{b})^a}{\Gamma(a)} 
\lambda^{a-1} \exp(-\frac{a}{b} \lambda) \\
&\propto \lambda^{Nr+a-1} \exp\left\{-\bigg(\sum_{n=1}^{N}x_n +\frac{a}{b}\bigg) \lambda\right\}
\propto \gammadist(\lambda \mid\widetildealpha, \widetildebeta),
\end{aligned}
$$
where $\widetildealpha \triangleq Nr+a$ and $\widetildebeta\triangleq\sum_{n=1}^{N}x_n +\frac{a}{b}$.
That is, the posterior density of rate $\lambda$ follows from a Gamma distribution.
We show that  the Gamma distribution is, itself, a conjugate prior for the rate parameter of a Gamma distribution when fixing the shape parameter. This is often referred to as the \textit{Gamma-Gamma model}.

\paragrapharrow{Conjugate prior for precision of a Gaussian  distribution.}

The Gamma distribution  also serves as a conjugate prior for the \textit{precision} parameter of a Gaussian distribution. To see this, suppose each entry $a_{mn}$ of matrix $\bA$ is i.i.d. normal model with mean $b_{mn}$ and precision $\tau$, i.e., the likelihood is $p(\bA \mid \bB, \tau^{-1})=\normal(\bA\mid\bB, \tau^{-1})$, the prior of $\tau$ is $p(\tau)=\gammadist(\tau \mid \alpha, \beta)$, where $\bA, \bB\in \real^{M\times N}$ are two matrices containing elements $a_{mn}$ and $b_{mn}$, respectively (the result can be applied to vector or scalar cases). Using Bayes' theorem, it can be shown that 
\begin{equation}\label{equation:gamma_conjugacy_general}
\begin{aligned}
p(\tau \mid &\bA, \bB, \alpha, \beta)
\propto \normal(\bA\mid\bB, \tau^{-1})\times \gammadist(\tau\mid \alpha, \beta)\\
&=\prod_{m,n=1}^{M,N} \normal(a_{mn}\mid b_{mn}, (\tau)^{-1}) 
\times \frac{\beta^\alpha}{\Gamma(\alpha)} \tau^{\alpha-1} \exp(-\beta \tau)\\
&\propto \tau^{\frac{MN}{2}}\exp\left\{   -\frac{\tau}{2}  \sum_{m,n=1}^{M,N}(a_{mn} - b_{mn}  )^2\right\}
\cdot \tau^{\alpha-1}\exp(-\beta \tau)\\
&=\tau^{\frac{MN}{2}+\alpha-1} \exp\left\{   -\tau  \left(\sum_{m,n=1}^{M,N}\frac{1}{2}(a_{mn} - b_{mn}  )^2 +\beta\right)\right\}
\propto \gammadist(\tau \mid \widetildealpha, \widetildebeta),\\
\end{aligned}
\end{equation}
with posterior parameters
\begin{equation}\label{equation:gamma_conju_posterior}
\widetildealpha=\frac{MN}{2}+\alpha, 
\qquad
\widetildebeta=  \sum_{m,n=1}^{M,N}\frac{1}{2}(a_{mn} - b_{mn}  )^2 +\beta.
\end{equation}
Hence, the posterior over precision $\tau$ remains Gamma-distributed.

\index{Normal-Gamma distribution}
\index{NormalGamma-Normal model}
\paragrapharrow{Joint conjugate prior for  Gaussian mean and precision.}
Going further, when the variance/precision parameter of the Gaussian distribution is not fixed with $x_1, x_2, \ldots, x_N$ drawn i.i.d. from a normal distribution with mean $\theta$ and precision $\lambda$. The \textit{normal-Gamma} distribution $\normalgamma(\alpha, \beta, \mu, c)$, with $\mu\in\real$  and $\alpha, \beta, c\in \real_+$  is a joint distribution on $(\theta, \lambda)$ by letting 
$$
\begin{aligned}
	\lambda &\sim \gammadist(\alpha, \beta); \\
	\theta \mid \lambda &\sim \normal(\mu, (c\lambda)^{-1}).
\end{aligned}
$$
That is, the joint PDF is 
$$
p(\theta, \lambda ) = \normal(\theta \mid \mu, (c\lambda)^{-1})\cdot \gammadist(\lambda \mid \alpha, \beta)
=\normalgamma(\theta, \lambda\mid \alpha, \beta, \mu, c).
$$
Given data $\mathcalX$, it turns out the posterior density is again a normal-Gamma distribution with 
$$
p(\theta, \lambda \mid\mathcalX ) \propto  
\prod_{n=1}^{N} \normal(x_n \mid \theta, \lambda^{-1}) 
\cdot 
\normalgamma(\theta, \lambda\mid \alpha, \beta, \mu, c)
\propto \normalgamma(\theta, \lambda \mid \widetildealpha,\widetildebeta, \widetilde{\mu}, \widetilde{c}),
$$
where 
$$
\begin{aligned}
\widetilde{\mu} &= \frac{c\mu +\sum_{n=1}^{N}x_n}{c+N} , \gap 
&\widetilde{c} &= c+N, \\
\widetildealpha&= \alpha+ \frac{N}{2}, \gap &\widetildebeta &= \beta + \frac{1}{2} \bigg(c\mu^2 -\widetilde{c}\widetilde{\mu}^2  + \sum_{n=1}^{N}x_n\bigg).
\end{aligned}
$$
In contrast to the Normal-Normal, this model is often referred to as the \textit{NormalGamma-Normal model}. 
The posterior mean for $\theta$ is a weighted average of the prior mean and the sample mean, 
$$
\widetilde{\mu} = \frac{c\mu +\sum_{n=1}^{N}x_n}{c+N}= \frac{c}{c+N} \mu + \frac{N}{c+N} \widebarx,
$$
where $\widebarx = \frac{1}{N}\sum_{n=1}^{N}x_n$.
From the posterior form of $\widetilde{c}$, the prior interpretation of $c$ can be described as the prior sample size for estimating the mean parameter $\theta$. The posterior shape parameter $\widetildealpha$ grows linearly with the sample size. And the posterior rate parameter $\widetildebeta$ can be written as 
$$
\widetildebeta 
= \beta + \frac{1}{2} \bigg(c\mu^2 -\widetilde{c}\widetilde{\mu}^2  + \sum_{n=1}^{N}x_n\bigg)
= \beta + \frac{1}{2}\sum_{n=1}^{N}(x_n - \widebarx)^2 + \frac{1}{2} \frac{cN}{c+N} (\widebarx-\mu)^2. 
$$
In other words, it is decomposed into the sum of a prior variation, the observed variation (sample variance), and the variation between the prior mean and sample mean:
$$
\widetildebeta=\text{(prior variation)} + \frac{1}{2}N \text{(observed variation)} + \frac{1}{2}\frac{cN}{c+N} \text{(variation between means)}.
$$

Placing  a Gamma prior over the inverse variance (i.e., precision) of a Gaussian distribution is equivalent to placing an inverse-Gamma prior on the \textit{variance}. 
We now define this distribution formally.
\begin{definition}[Inverse-Gamma Distribution\index{Inverse-Gamma distribution}]\label{definition:inverse_gamma_distribution}
A random variable $\rx$ is said to follow an \textit{inverse-Gamma distribution} with shape parameter $r>0$ and scale parameter $\lambda>0$, denoted $\rx\sim \inversegammadist(r, \lambda)$, if its density is
$$ f(x; r, \lambda)=\left\{
\begin{aligned}
	&\frac{\lambda^r}{\Gamma(r)} x^{-r-1} \exp(- \frac{\lambda}{x} ) ,& \mathrm{\,\,if\,\,} x > 0;  \\
	&0 , &\mathrm{\,\,if\,\,} x \leq 0.
\end{aligned}
\right.
$$
The mean and variance of inverse-Gamma distribution are given by 
$$ \Exp[\rx]=\left\{
\begin{aligned}
	&\frac{\lambda}{r-1}, \, &\mathrm{if\,} r\geq 1; \\
	&\infty, \, &\mathrm{if\,} 0<r<1.
\end{aligned}
\right.\qquad
\Var[\rx]=\left\{
\begin{aligned}
	&\frac{\lambda^2}{(r-1)^2(r-2)}, \, &\mathrm{if\,} r> 2; \\
	&\infty, \, &\mathrm{if\,} 0<r\leq 2.
\end{aligned}
\right.
$$
Figure~\ref{fig:dists_inversegamma} shows how the inverse-Gamma distribution behaves under different parameter settings.
\end{definition}

If $\rx$ is Gamma distributed, then $\ry=1/\rx$ is inverse-Gamma distributed.
Note that the inverse-Gamma density is not obtained by simply substituting $x=1/y$ into the Gamma density.
There is an additional factor of $y^{-2}$.~\footnote{Which follows from the \textit{Jacobian in the change-of-variables formula}. A short proof is provided here. Let $y=\frac{1}{x}$ where $y\sim \inversegammadist(r, \lambda)$ and $x\sim \gammadist(r, \lambda)$. Then, $f(y) \abs{dy} = f(x) \abs{dx}$, which results in $f(y) = f(x) \abs{\frac{dx}{dy}} = f(x)x^2 \xlongequal{ \mathrm{y}=\frac{1}{x}} \frac{\lambda^r}{\Gamma(r)} y^{-r-1} \exp(- \frac{\lambda}{y})$ for $y>0$. }  
The inverse-Gamma distribution is particularly useful as a prior for positive-valued parameters like variance. Compared to the Gamma distribution, it places more mass away from zero and has heavier tails (see Figure~\ref{fig:dists_inversegamma}), making it robust to extreme values.

\paragrapharrow{Conjugate prior for variance of a  Gaussian distribution.}
The inverse-Gamma distribution is a conjugate prior for the variance parameter of a Gaussian distribution when the mean is known. 
To see this, let the likelihood be $p(\bA \mid \bB, \sigma^2)=\normal(\bA\mid\bB, \sigma^2)$, where $\bA,\bB\in\real^{M\times N}$ are two matrices containing elements of $a_{mn}$ and $b_{mn}$, respectively (again the result can be applied to vector or scalar cases), and let the prior of $\sigma^2$ be $p(\sigma^2)=\inversegammadist(\sigma^2 \mid \alpha, \beta)$. Using Bayes' theorem, it can be shown that 
\begin{equation}\label{equation:inverse_gamma_conjugacy_general}
\begin{aligned}
p(\sigma^2 \mid &\bA, \bB, \alpha, \beta)
\propto \normal(\bA\mid \bB, \sigma^2)\times \inversegammadist(\sigma^2 \mid \alpha, \beta)\\
&=\prod_{m,n=1}^{M,N} \normal(a_{mn}\mid b_{mn}, \sigma^2) 
\times \frac{\beta^\alpha}{\Gamma(\alpha)} (\sigma^2)^{-\alpha-1} \exp(-\frac{\beta}{\sigma^2})\\
&\propto \frac{1}{\sigma^{MN}}\exp\left\{   -\frac{1}{2\sigma^2}  \sum_{m,n=1}^{M,N}(a_{mn} - b_{mn}  )^2\right\}
\cdot (\sigma^2)^{-\alpha-1}\exp(-\frac{\beta}{\sigma^2})\\
&=(\sigma^2)^{-\frac{MN}{2}-\alpha-1} \exp\left\{   -\frac{1}{\sigma^2}  \left(\sum_{m,n=1}^{M,N}\frac{1}{2}(a_{mn} - b_{mn}  )^2 +\beta\right)\right\}
\propto \inversegammadist(\sigma^2 \mid \widetildealpha, \widetildebeta),\\
\end{aligned}
\end{equation}
where the posterior parameters are
\begin{equation}\label{equation:inversegamma_conjugate_posterior}
\widetildealpha=\frac{MN}{2}+\alpha,  
\qquad
\widetildebeta=  \sum_{m,n=1}^{M,N}\frac{1}{2}(a_{mn} - b_{mn}  )^2 +\beta.
\end{equation}
That is, the posterior density of the variance $\sigma^2$ is also an inverse-Gamma distribution.
Notably, these posterior parameters match exactly those obtained when using a Gamma prior on the precision 
$\tau=1/\sigma^2$ (see Equation~\eqref{equation:gamma_conju_posterior}), confirming the duality between the two parameterizations.

\begin{figure}[htp]
\centering  
\vspace{-0.55cm} 
\subfigtopskip=2pt 
\subfigbottomskip=2pt 
\subfigcapskip=-5pt 
\subfigure[Contour plot of normal-inverse-Gamma density by varying parameter $r$ (\textcolor{brightlavender}{purple}=low, \textcolor{mydarkyellow}{yellow}=high). ]{\label{fig:dists_normalinversegamma_varyingR}
\includegraphics[width=0.955\linewidth]{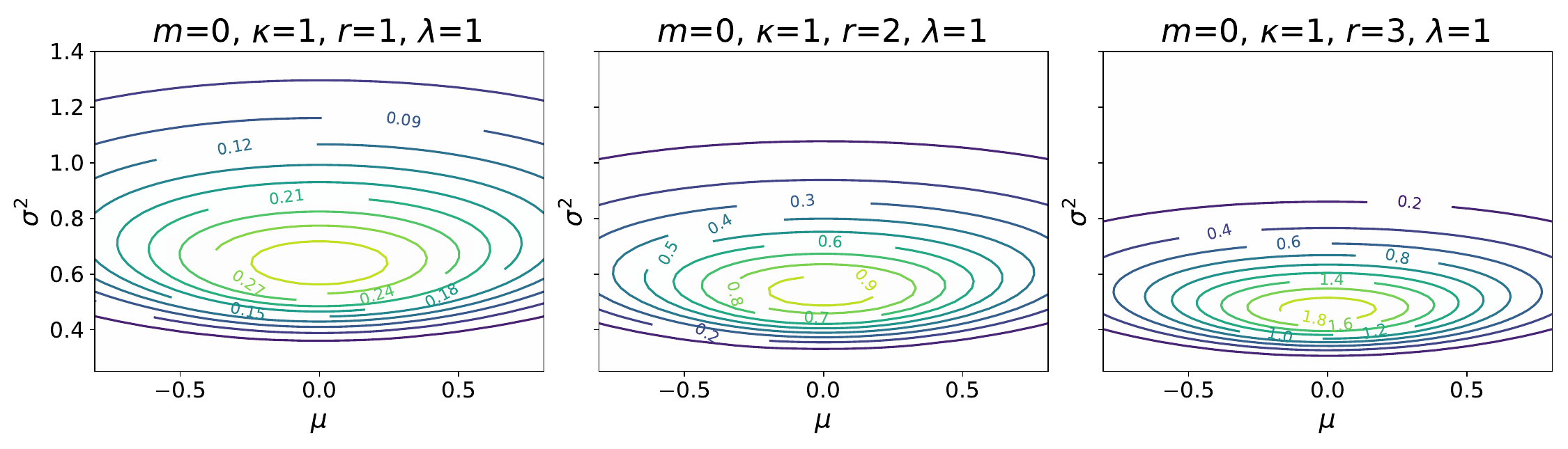}}
\subfigure[Contour plot of normal-inverse-Gamma density by varying parameter $\lambda$ (\textcolor{brightlavender}{purple}=low, \textcolor{mydarkyellow}{yellow}=high).]{\label{fig:dists_normalinversegamma_varyingLmabda}
\includegraphics[width=0.955\linewidth]{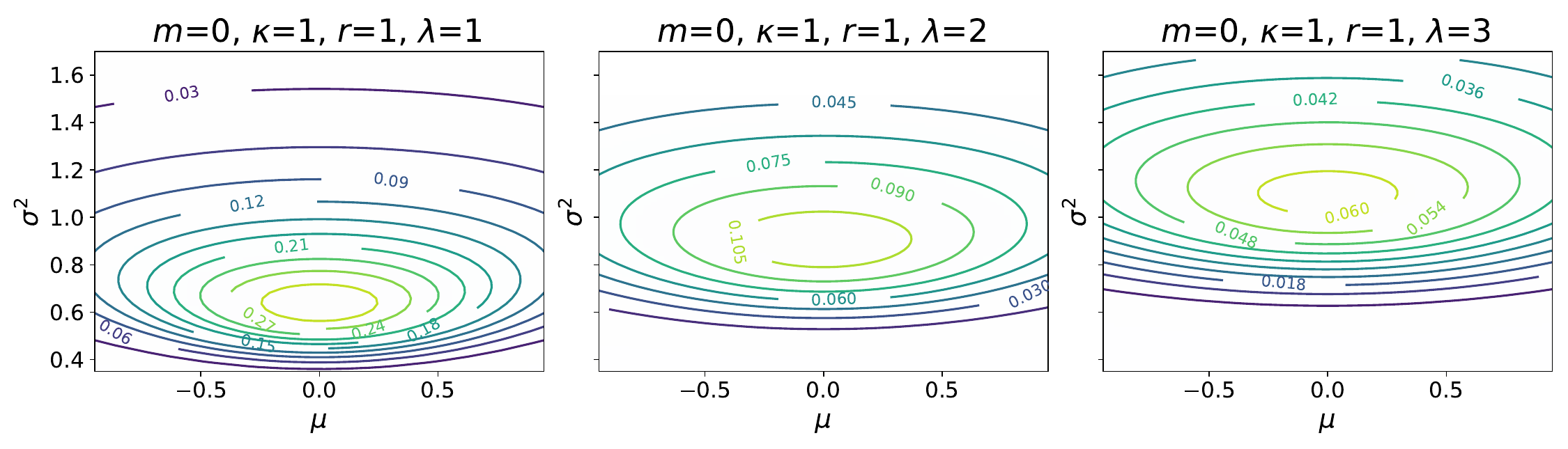}}
\subfigure[Contour plot of normal-inverse-Gamma density by varying parameter $\kappa$ (\textcolor{brightlavender}{purple}=low, \textcolor{mydarkyellow}{yellow}=high).]{\label{fig:dists_normalinversegamma_varyingKappa}
\includegraphics[width=0.955\linewidth]{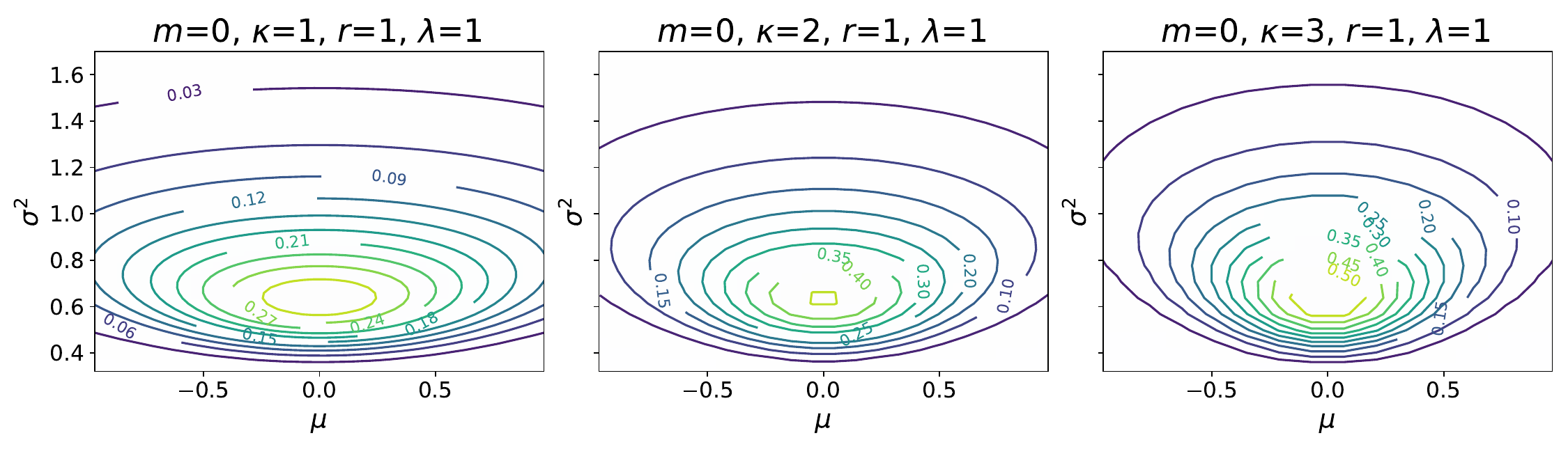}}
\subfigure[Contour plot of normal-inverse-Gamma density by varying parameter $m$ (\textcolor{brightlavender}{purple}=low, \textcolor{mydarkyellow}{yellow}=high).]{\label{fig:dists_normalinversegamma_varyingaM}
\includegraphics[width=0.955\linewidth]{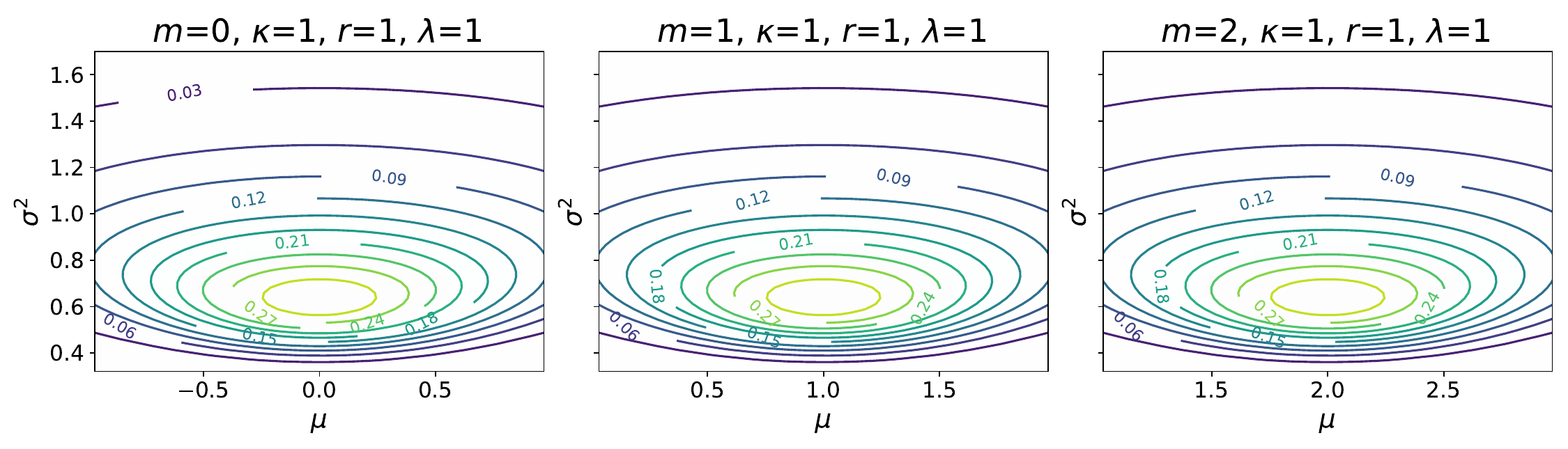}}
\caption{Normal-inverse-Gamma probability density functions by varying different parameters.}
\label{fig:dists_normalinversegamma_s}
\end{figure}

As shown earlier, the normal–Gamma distribution serves as a joint conjugate prior for the mean and precision of a Gaussian distribution. Similarly, the \textit{normal-inverse-Gamma (NIG)} distribution is a joint conjugate prior for the mean and variance of a Gaussian distribution. It is defined as follows.
\begin{definition}[Normal-Inverse-Gamma (NIG) Distribution\index{Normal-inverse-Gamma distribution}]\label{definition:normal_inverse_gamma}
The joint probability density function of the \textit{normal-inverse-Gamma distribution} is given by
\begin{equation}
\begin{aligned}
&\gap \nig (\smu, \sigma^2 \mid m, \kappa, r, \lambda) 
= \mathcal{N} (\mu\mid m, \frac{\ssigma^2}{\kappa})  \cdot \inversegammadist (\ssigma^2 \mid r,  \lambda) \\
&=\frac{1}{Z_{\nig}(\kappa, r, \lambda)}  (\ssigma^2)^{-\frac{2r +3}{2}} \exp\left\{-\frac{1}{2 \ssigma^2}\left[\kappa(m-\mu)^2 + 2\lambda \right] \right\},  \\
\end{aligned}
\label{equation:uni_gaussian_prior-nig}
\end{equation}
where $\sigma^2, r, \lambda>0$, and $Z_{\nig}(\kappa, r, \lambda)$ is a normalizing constant:
\begin{equation}\label{equation:uni_gaussian_giw_constant-nig}
Z_{\nig}(\kappa, r, \lambda) = \frac{\Gamma(r)}{\lambda^{r}} \sqrt{\frac{2\pi}{\kappa}}.
\end{equation}
Figure~\ref{fig:dists_normalinversegamma_s} displays several normal-inverse-Gamma densities under different parameter settings.
\end{definition}

\index{NIG model}
\paragrapharrow{Joint conjugate prior for the Gaussian mean and variance (NIG model).}
The normal–inverse-Gamma distribution provides an equivalent formulation to the normal–Gamma prior, but expressed in terms of variance rather than precision. This parameterization is often more convenient when working directly with variance.
Similar to the normal-Gamma prior, when the variance and mean parameters of the Gaussian distribution are not fixed with $N$ data points $\mathcalX=\{x_1, x_2, \ldots, x_N\}$ drawn i.i.d. from a normal distribution with mean $\mu$ and variance $\sigma^2$.
The normal-inverse-Gamma $\normalinversegamma(m_0, \kappa_0, r_0, \lambda_0)$ with $m_0\in\real$ and $r_0, \lambda_0, \kappa_0\in\real_+$ is a joint distribution on $\mu, \sigma^2$ by letting 
$$
\begin{aligned}
\sigma^2 &\sim \inversegammadist(r_0, \lambda_0);\\
\mu \mid \sigma^2 &\sim \normal(m_0, \frac{\sigma^2}{\kappa_0}).
\end{aligned}
$$
With this prior, $\mu$ and $\sigma^2$ decouple, and the posterior conditional densities of $\mu$ and $\sigma^2$ are Gaussian and inverse-Gamma, respectively.
The joint p.d.f of the NIG prior can be expressed as
$$
p(\mu, \sigma^2) = \normal(m_0, \frac{\sigma^2}{\kappa_0}) \cdot \inversegammadist(r_0, \lambda_0)
= \normalinversegamma(\mu, \sigma^2 \mid m_0, \kappa_0, r_0, \lambda_0).
$$

Again, by  Bayes' theorem ``$\mathrm{posterior} \propto \mathrm{likelihood} \times \mathrm{prior} $," the  posterior of the $\mu$ and $\sigma^2$ parameters under the NIG prior is
\begin{equation}\label{equation:conjugate_nigamma_general}
\begin{aligned}
&\gap p(\mu, \sigma^2\mid \mathcalX, \bbeta ) \\
&\propto \normal(\mathcalX \mid \mu, \sigma^2)
\cdot  \normalinversegamma(\mu, \sigma^2 \mid \bbeta)
\propto \prod_{n=1}^{N}\normal(x_n\mid \mu, \sigma^2)
 \cdot \normalinversegamma(\mu, \sigma^2 \mid  m_0, \kappa_0, r_0, \lambda_0)\\
&\stackrel{\star}{=} \frac{C}{(\sigma^2)^{\frac{2r_0 + 3+N}{2}}} \exp\left\{-\frac{1}{2 \ssigma^2}  \left[  N(\widebarx - \smu)^2 +  N S_{\widebarx} \right] \right\} 
 \exp\left\{ -\frac{1}{2 \ssigma^2} \left[2\lambda_0 + \kappa_0(m_0-\mu)^2\right] \right\}\\
&\propto (\sigma^2)^{-\frac{2r_N + 3}{2}}\exp\left\{ -\frac{1}{2 \ssigma^2} \left[ \lambda_N + \kappa_N(m_N-\mu)^2\right] \right\}\\
&\propto\nig(\mu, \sigma^2 \mid m_N, \kappa_{N}, r_N, \lambda_N),
\end{aligned}
\end{equation}
where $\bbeta=\{m_0, \kappa_0, r_0, \lambda_0\}$, $C=\frac{(2\pi)^{-N/2}}{Z_{\nig}(\kappa_0, r_0, \lambda_0)}$, the equality $(\star)$ follows from Equation~\eqref{equation:uni_gaussian_likelihood}, and 
$$
\begin{aligned}
m_N &= \frac{\kappa_0 m_0 + N\widebarx}{\kappa_{N}} = \frac{\kappa_0 }{\kappa_{N}}m_0 + \frac{N}{\kappa_{N}}\widebarx,\\ 
\kappa_{N}&= \kappa_{0} +N, \qquad \quad 
r_N = r_0 +\frac{N}{2},\\
\lambda_N &=\lambda_0 +\frac{1}{2}(NS_{\widebarx} + N\widebarx^2 + \kappa_{0} m_0^2 -\kappa_{N}m_N^2)\\
&= \lambda_0+\frac{1}{2}\left(NS_{\widebarx} + \frac{\kappa_0 N }{\kappa_{0}+N} (\widebarx - m_0)^2\right),
\end{aligned}
$$
where $S_{\widebarx}=\sum_{n=1}^N(x_n - \widebarx)^2$ and $\widebarx = (\sum_{n=1}^{N}x_n)/N$.
Note in the above derivation, we use the  fact about the likelihood under Gaussian in Equation~\eqref{equation:uni_gaussian_likelihood}.
The posterior mean $m_N$ is a weighted average of the prior mean $m_0$ and the sample mean $\widebarx$, with weights proportional to $\kappa_0$ and $N$, respectively.
The parameter $\kappa_0$ can be interpreted as a prior sample size for the mean.

We will discuss the posterior marginal likelihood in the \textit{normal-inverse-Chi-squared (NIX)} case.
Further discussion on the posterior marginal likelihood for the NIG prior can be found in \citet{murphy2007conjugate}. We will leave this to the readers as it is rather similar as that in the NIX prior.

\begin{figure}[h]
\centering  
\vspace{-0.35cm} 
\subfigtopskip=2pt 
\subfigbottomskip=2pt 
\subfigcapskip=-5pt 
\subfigure[Chi-squared distribution PDFs.]{\label{fig:dists_chisquared}
\includegraphics[width=0.481\linewidth]{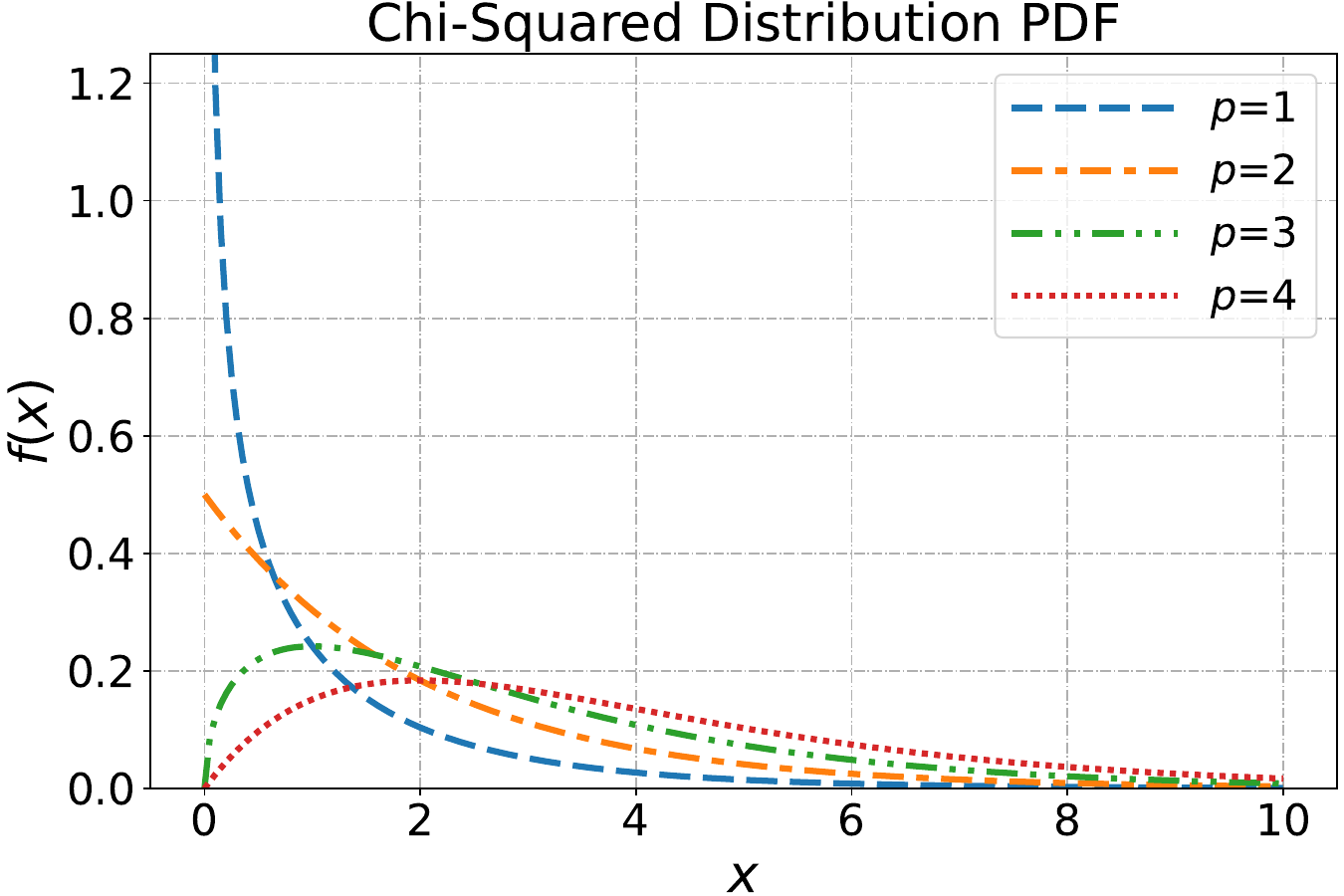}}
\subfigure[Chi-squared distribution CDFs.]{\label{fig:dists_chisquared_cdf}
\includegraphics[width=0.481\linewidth]{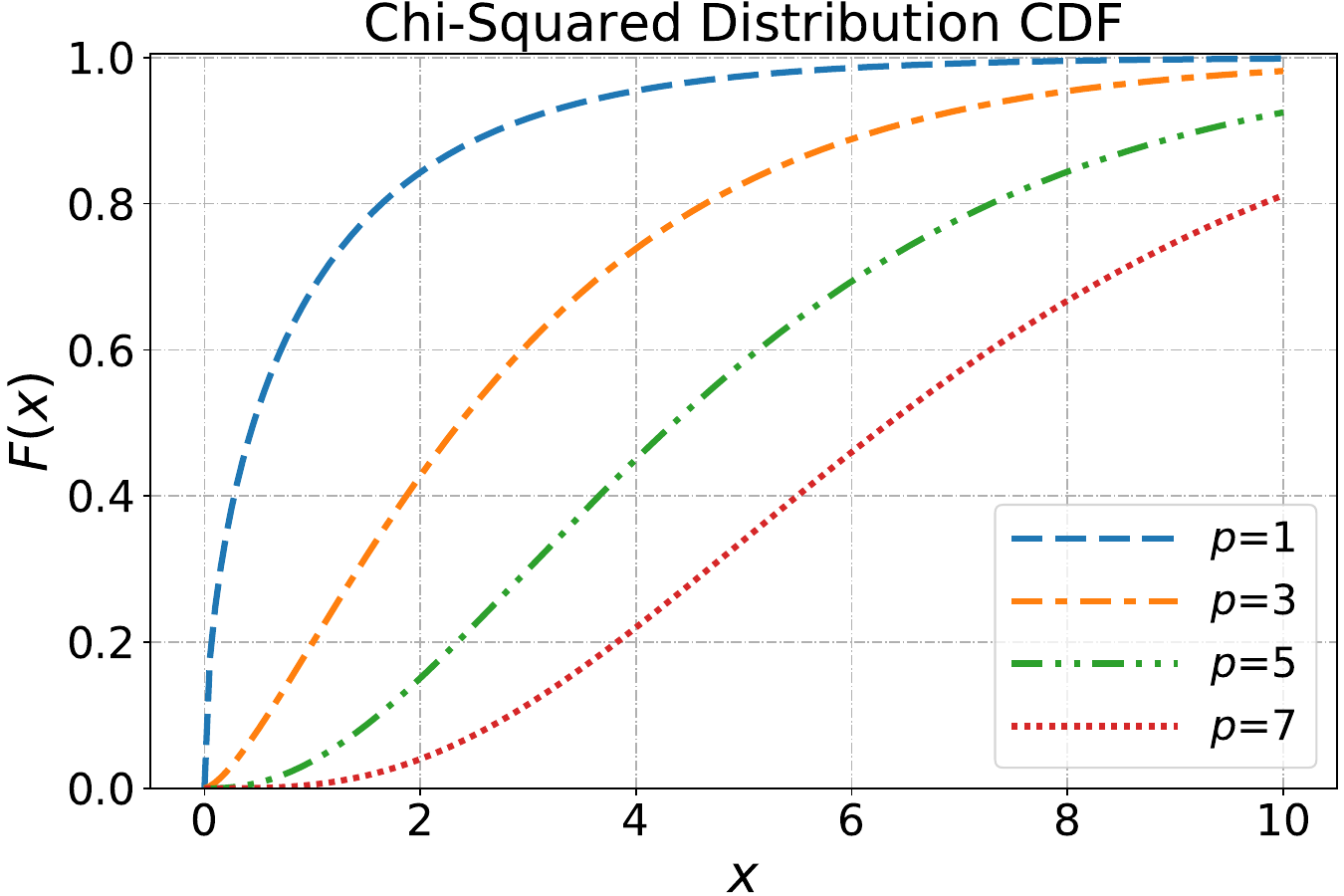}}
\caption{Chi-squared probability density functions and cumulative distribution functions for different values of  parameters.}
\label{fig:dists_chisquare_all}
\end{figure}

Another distribution that is closely related to the Gamma distribution is called the \textit{Chi-squared} distribution, which plays a central role in the distribution theory of linear models \citep{lu2021rigorous}. Its formal definition is as follows.
\begin{definition}[Chi-Squared Distribution\index{Chi-squared distribution}]\label{definition:chisquare_distribution}
Let $\rva \sim \normal(\bzero, \bI_{p})$, where $\bI_p$ is the $p\times p$ identity matrix. Then $\rx=\sum_i^p \ra_{ii}^2$ follows the \textit{Chi-squared distribution} with $p$ \textit{degrees of freedom}. We write $\rx \sim \chi^2(p)$, and we can see this is equivalent to $\rx\sim \gammadist(p/2, 1/2)$:
$$ 
f(x; p)=\left\{
\begin{aligned}
&\frac{1}{2^{p/2}\Gamma(\frac{p}{2})} x^{\frac{p}{2}-1} \exp(-\frac{x}{2}) ,& \mathrm{\,\,if\,\,} x \geq 0;  \\
&0 , &\mathrm{\,\,if\,\,} x <0.
\end{aligned}
\right.
$$
The mean and variance of $\rx\sim \chi^2(p)$ are given by 
$$
\Exp[\rx]=p, \qquad \Var[\rx]=2p. 
$$
Figure~\ref{fig:dists_chisquare_all} compares different \text{degrees of freedom} $p$ for the Chi-squared distribution.
It can be observed that the degrees of freedom parameter of the Chi-squared distribution affects the slope of its cumulative distribution function.
As the degrees of freedom increase, the distribution becomes more symmetric and shifts rightward. Conversely, smaller degrees of freedom yield highly skewed distributions concentrated near zero.
\end{definition}

\begin{exercise}[Sum of Independent Chi-Squared]
Show that $\rx=\sum_{n=1}^{N}\rx_n \sim \chisquared(p)$ given independent variables  $\rx_n\sim\chisquared(p_n)$ ($n=1,2,\ldots,N$) and $p=\sum_{n=1}^{N}p_n$.
\end{exercise}

We notice the Chi-squared distribution is a sum of i.i.d. standard normal variables. A generalization can be obtained as the sum of generally independent distributed Gaussians.
\begin{remark}[Noncentral Chi-Squared Distribution]
The sum of $p$ independently distributed $\rx_i\sim \normal(\mu_i, \sigma^2), \forall i\in\{1,2,\ldots,p\}$ is a generalization of the Chi-squared distribution, and is called the \textit{noncentral Chi-squared distribution}, denoted as $\rx\sim \nonchisquared(p)$:
$$
\rx=\frac{1}{\sigma^2}\sum_{i=1}^{p} \rx_i^2\sim \nonchisquared(p).
$$
The value $\delta = \frac{1}{\sigma^2}\sum_{i=1}^{p}\mu_i^2$ is called the \textit{noncentral parameter}.
The probability density function is 
$$
f(x;\delta, p) = 
\exp(-\delta/2) \sum_{j=0}^{\infty} \frac{(\delta/2)^j}{j!} \chisquared(x\mid 2j+k),
$$
where $\chisquared(x\mid 2j+k)$ is the PDF of a central Chi-squared distribution with degress of freedom $2j+k$.
The mean and variance of $\rx\sim \nonchisquared(p)$ are given by 
$$
\Exp[\rx]=p+\delta, \qquad \Var[\rx]=2p+4\delta. 
$$
\end{remark}

\begin{SCfigure}
\centering
\includegraphics[width=0.5\textwidth]{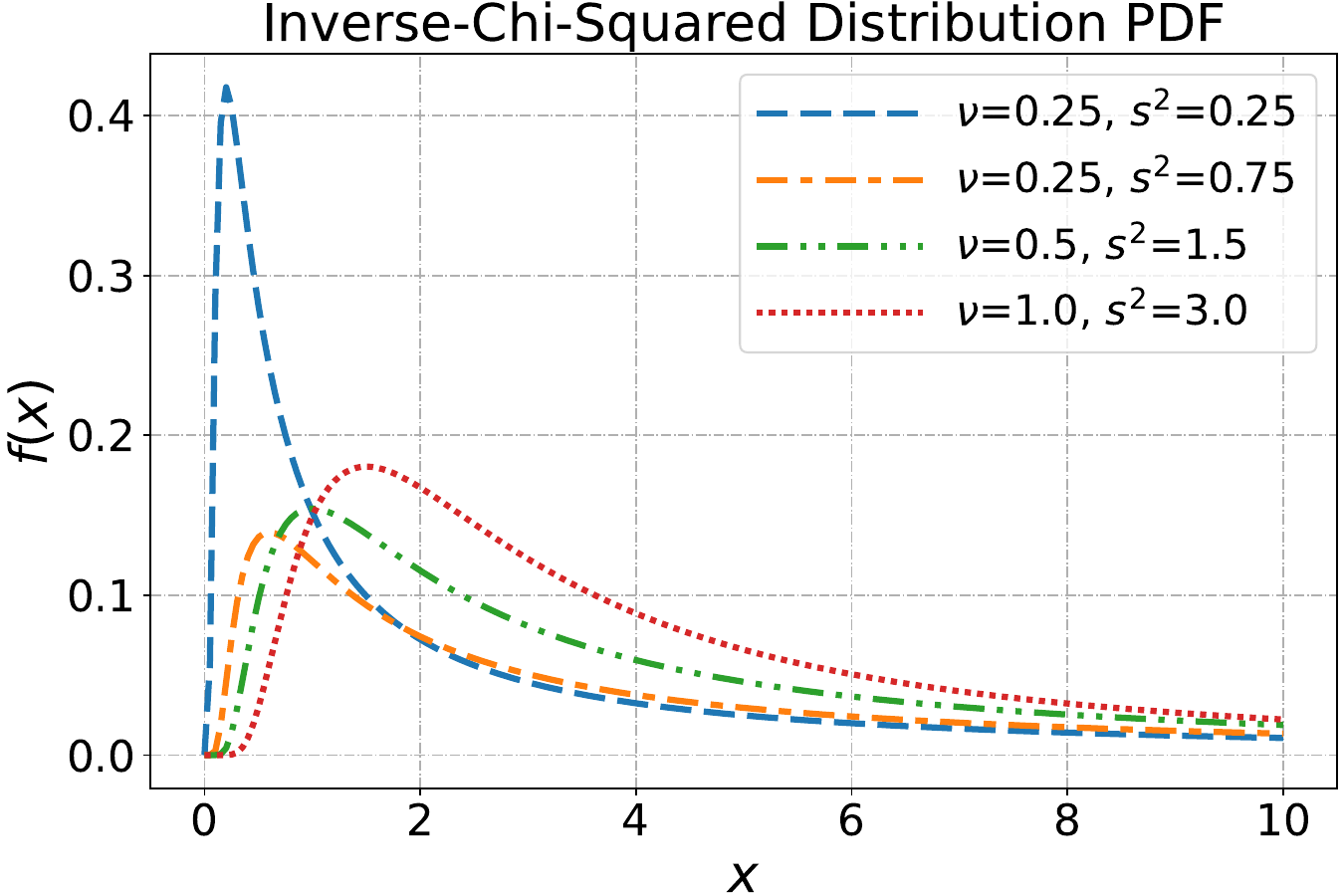}
\caption{Inverse-Chi-squared probability density functions for different values of  parameters.}
\label{fig:dist_inversechisquared}
\end{SCfigure}

An counterpart of the inverse-Gamma distribution is known as the \textit{inverse-Chi-squared} distribution. 
Following  the definition of the inverse-Gamma distribution in Definition~\ref{definition:inverse_gamma_distribution}, we provide the rigorous definition of the inverse-Chi-squared distribution as follows.
\begin{definition}[Inverse-Chi-Squared Distribution]\label{definition:inverse-chi-square}
A random variable $\rx$ is said to follow an \textit{inverse-Chi-squared distribution} with parameter $\nu>0$ and $s^2>0$, denoted $\rx\sim \inversegammadist(\frac{\nu}{2}, \frac{\nu s^2}{2})$, if
$$ f(x; \nu, s^2)=\left\{
\begin{aligned}
&\frac{{(\frac{\nu s^2}{2})}^{\frac{\nu}{2}}}{\Gamma(\frac{\nu}{2})} x^{-\frac{\nu}{2}-1} \exp(- \frac{\nu s^2}{2x} ) ,& \mathrm{\,\,if\,\,} x > 0;  \\
&0 , &\mathrm{\,\,if\,\,} x \leq 0.
\end{aligned}
\right.
$$
And it is also compactly denoted by $\rx \sim \inversechidist(\nu, s^2)$. The parameter $\nu >0$ is called the \textit{degrees of freedom}, and $s^2 > 0$ is the \textit{scale parameter}. And it is also known as the \textit{scaled} inverse-Chi-squared distribution.
The mean and variance of the inverse-Chi-squared distribution are given by 
$$ \Exp[\rx]=\left\{
\begin{aligned}
&\frac{\nu s^2}{\nu-2}, \, &\mathrm{if\,\,} \nu\geq 2; \\
&\infty, \, &\mathrm{if\,\,} 0<\nu<2.
\end{aligned}
\right.\qquad
\Var[\rx]=\left\{
\begin{aligned}
&\frac{2\nu^2 s^4}{(\nu-2)^2(\nu-4)}, \, &\mathrm{if\,\,} \nu\geq 4; \\
&\infty, \, &\mathrm{if\,\,} 0<\nu<4.
\end{aligned}
\right.
$$
To establish a connection with the inverse-Gamma distribution, we can set $S=\nu s^2$. Then the inverse-Chi-squared distribution can also be denoted by $\rx\sim \inversegammadist(\frac{\nu}{2}, \frac{S}{2})$ if $\rx \sim \inversechidist(\nu, s^2)$,  the form of which conforms to the univariate case of the inverse-Wishart distribution (Definition~\ref{definition:multi_inverse_wishart}). 
And we will observe the similarities in the posterior parameters as well.
Figure~\ref{fig:dist_inversechisquared} illustrates the impact of different parameters $\nu, s^2$ for the inverse-Chi-squared distribution.
\end{definition}

\begin{exercise}[Conjugate prior for the variance of a Gaussian distribution]
Show that the inverse-Chi-squared distribution is a conjugate prior for the Gaussian variance parameter when the mean parameter is fixed. \textit{Hint: The derivation follows the same procedure as in the inverse-Gamma case.}
\end{exercise}

\index{Joint conjugate prior}
\index{Normal-inverse-Chi-squared distribution}
As previously discussed, the normal-inverse-Gamma (NIG) distribution serves as a joint conjugate prior for the mean and variance of a Gaussian distribution. The \textit{normal-inverse-Chi-squared (NIX)} distribution defined as follows is an alternative joint conjugate prior.
\begin{definition}[Normal-Inverse-Chi-Squared (NIX) Distribution]\label{definition:normal_inverse_chi_square}
Similar to the normal-inverse-Gamma distribution, the \textit{normal-inverse-Chi-squared (NIX) distribution} is defined as (where again we set $S=\nu s^2$ as that in the inverse-Chi-square distribution to establish a connection with the normal-inverse-Gamma density)
\begin{equation}
\begin{aligned}
\nix (\smu, &\sigma^2 \mid m, \kappa, \nu, S) 
= \mathcal{N} (\mu\mid m, \frac{\ssigma^2}{\kappa})  
\cdot \inversechidist (\ssigma^2 \mid \nu,  s^2) \\
&=\frac{1}{Z_{\nix}(\kappa, \nu, s^2)} (\ssigma^2)^{-(\nu/2 + 3/2)} \exp\left\{ -\frac{1}{2 \ssigma^2} \left[\nu s^2 + \kappa(m-\mu)^2\right] \right\} \\
&\xlongequal{S = \nu s^2}
\frac{1}{Z_{\nix}(\kappa, \nu, s^2)} (\ssigma^2)^{-(\nu/2 + 3/2)} \exp\left\{ -\frac{1}{2 \ssigma^2} \left[S + \kappa(m-\mu)^2\right] \right\},
\end{aligned}
\label{equation:uni_gaussian_prior}
\end{equation}
where $\sigma^2, \nu, s^2>0$, and $Z_{\nix}(\kappa, \nu, s^2)$ is a normalizing constant:
\begin{equation}\label{equation:uni_gaussian_giw_constant}
Z_{\nix}(\kappa, \nu, s^2) = \Gamma\big(\frac{\nu}{2}\big) 
\big(\frac{2}{\nu s^2}\big)^{\nu/2} \sqrt{\frac{2\pi}{\kappa}}
=  \Gamma\big(\frac{\nu}{2}\big) \big(\frac{2}{S}\big)^{\nu/2} \sqrt{\frac{2\pi}{\kappa}}.
\end{equation}
The  normal-inverse-Chi-squared distribution can also be denoted by 
$\rx\sim \normalinversegamma(m, \kappa, \frac{\nu}{2}, \frac{S}{2})$ if 
$\rx \sim \nix(m, \kappa, \nu, s^2)$, the form of which  conforms to the univariate case of the normal-inverse-Wishart distribution (see Equation~\eqref{equation:multi_gaussian_prior}). And we will see the similarities in the posterior parameters as well.
\end{definition}

\index{NIX model}
\paragrapharrow{Joint conjugate prior for the Gaussian mean and variance (NIX model).}
Similar to the normal-inverse-Gamma prior, when the variance and mean parameters of the Gaussian distribution are not fixed with $N$ data points $\mathcalX=\{x_1, x_2, \ldots, x_N\}$ drawn i.i.d. from a normal distribution with mean $\mu$ and variance $\sigma^2$.
The normal-inverse-Chi-squared $\nix(m_0, \kappa_0, \nu_0, S_0=\nu_0\sigma_0^2)$ with $m_0\in\real$ and $\kappa_0, \mu_0, S_0\in\real_+$ is a joint distribution on $\mu, \sigma^2$ by letting 
$$
\begin{aligned}
\sigma^2 &\sim \inversechidist(\nu_0, \sigma_0^2);\\
\mu \mid \sigma^2 &\sim \normal(m_0, \frac{\sigma^2}{\kappa_0}).
\end{aligned}
$$
Again, by Bayes' theorem ``$\mathrm{posterior} \propto \mathrm{likelihood} \times \mathrm{prior} $," the conditional posterior of the $\mu$ and $\sigma^2$ parameters under the NIX prior is
\begin{equation}\label{equation:nix-posterior}
\begin{aligned}
&\gap p(\mu, \sigma^2\mid \mathcalX, \bbeta ) \\
&\propto p(\mathcalX \mid \mu, \sigma^2) p(\mu, \sigma^2 \mid \bbeta) = p(\mathcalX, \mu, \sigma^2 \mid \bbeta)\\
&= \frac{C}{(\sigma^2)^{\frac{\nu_0 + 3+N}{2}}} \exp\left\{ -\frac{1}{2 \ssigma^2}  \left[  N(\widebarx - \smu)^2 +  N S_{\widebarx} \right] \right\} 
\exp\left\{ -\frac{1}{2 \ssigma^2} \left[S_0 + \kappa_0(m_0-\mu)^2\right] \right\}\\
&= C\times  (\sigma^2)^{-\frac{\nu_N + 3}{2}}\exp\left\{ -\frac{1}{2 \ssigma^2} 
\left[ S_N + \kappa_N(m_N-\mu)^2\right] \right\}\\
&\propto\nix(\mu, \sigma^2\mid m_N, \kappa_{N}, \nu_N, \textcolor{mylightbluetext}{S_N})
= \normal (\mu\mid m_N, \frac{\ssigma^2}{\kappa_N})  \cdot \inversechidist (\ssigma^2 \mid \nu_N,  \textcolor{mylightbluetext}{\sigma^2_N}),
\end{aligned}
\end{equation}
where $\bbeta=\{m_0, \kappa_0, \nu_0, S_0=\nu_0\sigma_0^2\}$, $C=\frac{(2\pi)^{-N/2}}{Z_{\nix}(\kappa_0, \nu_0, \sigma^2_0)}$, and 

$$
\begin{aligned}
m_N &= \frac{\kappa_0 m_0 + Nv}{\kappa_{N}} = \frac{\kappa_0 }{\kappa_{N}}m_0 + \frac{N}{\kappa_{N}}\widebarx,\\ 
\kappa_{N}&= \kappa_{0} +N, \qquad\quad
\nu_N = \nu_0 +N,\\
S_N &= S_0 +NS_{\widebarx} + N\widebarx^2 + \kappa_{0} m_0^2 -\kappa_{N}m_N^2\\
&=S_0 +NS_{\widebarx} + \frac{\kappa_0 N }{\kappa_{0}+N} (\widebarx - m_0)^2,\\
\nu_N \sigma_N^2 &= S_N \leadto     \sigma_N^2 = \frac{S_N}{\nu_N} ,
\end{aligned}
$$
Thus, the posterior is again a normal-inverse-Chi-squared density. \footnote{This posterior shares the  same form as that in the multivariate case from Equation~\eqref{equation:niw_posterior_equation_1} except the $N$ in $NS_{\widebarx}$, which results from the difference between the multivariate Gaussian distribution and the univariate Gaussian distribution. Similarly, in the inverse-Chi-squared language, we can show that $\nu_N \sigma_N^2 = S_N$.}

Suppose $\nu_0\geq 2$, or $N\geq 2$ such that $\nu_N\geq 2$, the posterior expectations are given by 
$$
\Exp[\mu \mid \mathcalX, \bbeta] = m_N, \qquad
\Exp[\sigma^2 \mid \mathcalX, \bbeta] = \frac{S_N}{\nu_N-2}.
$$

\paragrapharrow{Marginal posterior of $\sigma^2$.}
Integrating out $\mu$ from the joint posterior yields:
$$
\begin{aligned}
p(\sigma^2 \mid \mathcalX, \bbeta) 
&= \int_{\mu} p(\mu, \sigma^2 \mid \mathcalX, \bbeta)\, d \mu \\
&= \int_{\mu }  \mathcal{N} (\mu\mid m_N, \frac{\ssigma^2}{\kappa_N})  
\cdot \inversechidist (\ssigma^2 \mid \nu_N,  \sigma^2_N)\,d\mu 
= \inversechidist (\ssigma^2\mid \nu_N,  \sigma^2_N),
\end{aligned}
$$
which is just an integral over a Gaussian distribution.

\paragrapharrow{Marginal posterior of $\mu$.}
Integrating out $\sigma^2$ in the posterior, we have
$$
\begin{aligned}
p(\mu \mid \mathcalX, \bbeta) 
&= \int_{\sigma^2} p(\mu, \sigma^2 \mid \mathcalX, \bbeta)\, d \sigma^2 
= \int_{\sigma^2 }  \mathcal{N} (\mu\mid m_N, \frac{\ssigma^2}{\kappa_N}) 
\cdot \inversechidist (\ssigma^2 \mid \nu_N,  \sigma^2_N)\,d\sigma^2 \\
&= \int_{\sigma^2}C(\sigma^2)^{-\frac{\nu_N + 3}{2}}
\exp\left \{-\frac{1}{2 \ssigma^2} \left[ S_N + \kappa_N(m_N-\mu)^2\right] \right\}\,d\sigma^2. 
\end{aligned}
$$
Let $\phi = \sigma^2$ and $\alpha = (\nu_N+1)/2$, $A =  S_N + \kappa_N(m_N-\mu)^2$, and $x = \frac{A}{2\phi}$, we have 
$$
\frac{d \phi}{d x} = -\frac{A}{2}x^{-2},
$$
where $A$ can be easily verified to be positive and $\phi=\sigma^2>0$.
It follows that
\begin{equation*}
\begin{aligned}
p(\mu \mid \mathcalX, \bbeta) 
&=\int_{0}^{\infty} C(\phi)^{-\alpha-1} \exp\left(-\frac{A}{2 \phi} \right)\,d\phi\\
&=\int_{\textcolor{black}{\infty}}^{\textcolor{black}{0}}  C(\frac{A}{2x})^{-\alpha-1} \exp\left( -x\right)   ( \textcolor{black}{-}\frac{A}{2}x^{-2}) \,dx \qquad &\text{(since $x=\frac{A}{2\phi}$)}\\
&=\int_{\textcolor{black}{0}}^{\textcolor{black}{\infty}}  C(\frac{A}{2x})^{-\alpha-1} \exp\left( -x\right)   ( \frac{A}{2}x^{-2}) \,dx\\
&= (\frac{A}{2})^{-\alpha} \int_{x}  Cx^{\alpha-1} \exp\left( -x\right)   \,dx \\
&= (\frac{A}{2})^{-\alpha}  (C\cdot \Gamma(1)) \int_{x}\gammadist(x\mid \alpha, 1) \,dx\qquad &\text{(see Definition~\ref{definition:gamma-distribution})}\\
&= (C\cdot \Gamma(1))\left[ \nu_N\sigma_N^2 +\kappa_{N}(m_N-\mu)^2 \right]^{-\frac{\nu_N+1}{2}}\\
&\overset{(a)}{=} (C\cdot \Gamma(1)) (\nu_N\sigma_N^2)^{-\frac{\nu_N+1}{2}} \left[ 1 +\frac{\kappa_{N}}{\nu_N\sigma_N^2}(m_N-\mu)^2 \right]^{-\frac{\nu_N+1}{2}}
\end{aligned}
\end{equation*}
We notice that $C$ is defined in Equation~\eqref{equation:nix-posterior} (in terms of $\{\kappa_N, \nu_N, \sigma^2_N\}$) with 
$$
C\overset{(b)}{=}\frac{(2\pi)^{-N/2}}{Z_{\nix}(\kappa_N, \nu_N, \sigma^2_N)} = \frac{(2\pi)^{-N/2}}{\frac{\sqrt{(2\pi)}}{\sqrt{\kappa_N}} \Gamma(\frac{\nu_N}{2}) (\frac{2}{\nu_N \sigma^2_N})^{\nu_N/2}} \propto (\nu_N \sigma^2_N)^{\nu_N/2}.
$$
Combining equalities (a) and (b) above, we obtain 
$$
p(\mu\mid \mathcalX, \bbeta)  \propto \frac{1}{\sigma_N/\sqrt{\kappa_N}} 
\left[ 1 +\frac{\kappa_{N}}{\nu_N\sigma_N^2}(\mu-m_N)^2 \right]^{-\frac{\nu_N+1}{2}} 
\propto \tau(\mu\mid m_N, \sigma_N^2/\kappa_N, \nu_N),
$$
which is a univariate Student's $t$ distribution (Definition~\ref{equation:student_t_dist}).

\paragrapharrow{Marginal likelihood of data.}
By Equation~\eqref{equation:nix-posterior}, we can obtain the marginal likelihood of data under hyper-parameters $\bbeta=( m_0, \kappa_0, \nu_0, S_0=\nu_0\sigma_0^2)$:
$$
\begin{aligned}
p(\mathcalX \mid \bbeta) 
&= \int_{\mu} 
\int_{\sigma^2}  p(\mathcalX, \mu, \sigma^2 \mid\bbeta) 
\,d\mu \,d\sigma^2\\ 
&=\frac{(2\pi)^{-N/2}}{Z_{\nix}(\kappa_0, \nu_0, \sigma^2_0)}  \int_{\mu} \int_{\sigma^2}  
(\sigma^2)^{-\frac{\nu_N + 3}{2}}\exp\left\{ -\frac{1}{2 \ssigma^2} \left[ S_N + \kappa_N(m_N-\mu)^2\right] \right\}
\,d\mu \,d\sigma^2\\
&= (2\pi)^{-N/2}\frac{Z_{\nix}(\kappa_N, \nu_N, \sigma^2_N)}{Z_{\nix}(\kappa_0, \nu_0, \sigma^2_0)} 
= (\pi)^{-N/2} \frac{\Gamma(\nu_N/2)}{\Gamma(\nu_0/2)} \sqrt{\frac{\kappa_0}{\kappa_N}} \frac{(\nu_0\sigma^2_0)^{\nu_0/2}}{(\nu_N\sigma^2_N)^{\nu_N/2}}.
\end{aligned}
$$

\paragrapharrow{Posterior predictive for new data with observations.}
Let the number of samples for data set $\{\xstar, \mathcalX\}$ be $\Nstar = N+1$, we have
\begin{equation}\label{equation:nix-posterior-new-withobser}
\begin{aligned}
&\gap p(\xstar \mid\mathcalX, \bbeta) 
= \frac{p(\xstar, \mathcalX \mid \bbeta)}{p(\mathcalX\mid \bbeta)}\\
&=\left\{(2\pi)^{-\Nstar/2}\frac{Z_{\nix}(\kappa_\Nstar, \nu_\Nstar, \sigma^2_\Nstar)}{Z_{\nix}(\kappa_0, \nu_0, \sigma^2_0)}\right\}
\bigg/\left\{(2\pi)^{-N/2}\frac{Z_{\nix}(\kappa_N, \nu_N, \sigma^2_N)}{Z_{\nix}(\kappa_0, \nu_0, \sigma^2_0)}\right\}\\
&=(2\pi)^{-1/2} \frac{Z_{\nix}(\kappa_\Nstar, \nu_\Nstar, \sigma^2_\Nstar)}{Z_{\nix}(\kappa_N, \nu_N, \sigma^2_N)}
=(\pi)^{-1/2} 
\sqrt{\frac{\kappa_N}{\kappa_{N^{\star}}}}    \frac{\Gamma(\frac{\nu_{N^{\star}}}{2})}{\Gamma(\frac{\nu_{N}}{2})}
\frac{(\nu_N \sigma_N^2)^{\frac{\nu_N}{2}}}{(\nu_{\Nstar}\sigma_{\Nstar}^2)^{\frac{\nu_{\Nstar}}{2}}}\\
&=\frac{\Gamma(\frac{\nu_{N}+1}{2})}{\Gamma(\frac{\nu_{N}}{2})}
\sqrt{\frac{\kappa_N}{(\kappa_{N}+1)}   \frac{1}{(\pi\nu_{N}\sigma_{N}^2)}}   
\left[\frac{(\nu_{\Nstar}\sigma_{\Nstar}^2)}{(\nu_N \sigma_N^2)}\right]^{-\frac{\nu_{N}+1}{2}}.
\end{aligned}
\end{equation}
We realize that 
$$
\begin{aligned}
m_N &= \frac{\kappa_{N^{\star}}m_{N^{\star}} - x^\star}{\kappa_N}=\frac{(\kappa_0 + N + 1)m_{N^{\star}} - x^\star}{\kappa_0 + N} , \\
m_{N^{\star}} &= \frac{\kappa_{N} m_N +x^{\star}}{\kappa_{N^{\star}}} 
= \frac{(\kappa_{0}+N) m_N +x^{\star}}{\kappa_{0}+N+1},\\
S_{N^{\star}} &= S_N + x^{\star} x^{\star T} - \kappa_{N^{\star}} m_{N^{\star}}^2  + \kappa_N m_N^2 \\
&= S_N + \frac{\kappa_N + 1}{\kappa_N}(m_{N^{\star}} - x^\star)^2
=S_N + \frac{\kappa_N }{\kappa_N+ 1}(m_{N} - x^\star)^2,
\end{aligned}
$$
Thus, we have 
\begin{equation}\label{equation:nix-substitute-posterior-withobser}
\begin{aligned}
\left[\frac{(\nu_{\Nstar}\sigma_{\Nstar}^2)}{(\nu_N \sigma_N^2)}\right]^{-\frac{\nu_{N}+1}{2}}&=
\left(\frac{S_\Nstar}{S_N}\right)^{-\frac{\nu_{N}+1}{2}} =1 +  \frac{\kappa_N(m_{N} - x^\star)^2 }{(\kappa_N+ 1)\nu_N\sigma_N^2}.
\end{aligned}
\end{equation}
Substituting Equation~\eqref{equation:nix-substitute-posterior-withobser} into Equation~\eqref{equation:nix-posterior-new-withobser}, it follows that 
$$
\begin{aligned}
p(\bxstar \mid\mathcalX, \bbeta) 
&=\frac{\Gamma(\frac{\nu_{N}+1}{2})}{\Gamma(\frac{\nu_{N}}{2})}
\sqrt{\frac{\kappa_N}{(\kappa_{N}+1)}   \frac{1}{(\pi\nu_{N}\sigma_{N}^2)}}   
\left[1 +  \frac{\kappa_N(m_{N} - x^\star)^2 }{(\kappa_N+ 1)\nu_N\sigma_N^2}\right]^{-\frac{\nu_{N}+1}{2}}\\
&= \tau(x^\star \mid m_N, \frac{\kappa_{N}+1}{\kappa_{N}}\sigma^2_N, \nu_N  ).
\end{aligned}
$$
That is, the predictive distribution is a Student's $t$ distribution centered at the posterior mean $m_N$, with scaled variance and $\nu_N$ degrees of freedom.

\paragrapharrow{Posterior predictive for new data without observations.}
In the absence of data ($N=0$), the prior predictive distribution is:
$$
\begin{aligned}
p(\xstar \mid \bbeta)
&=  \int_{\mu} \int_{\sigma^2} p(\xstar, \mu, \sigma^2 \mid \bbeta) d\mu d\sigma^2 
=(2\pi)^{-1/2}\frac{Z_{\nix}(\kappa_1, \nu_1, \sigma^2_1)}{Z_{\nix}(\kappa_0, \nu_0, \sigma^2_0)}\\
&=(\pi)^{-1/2} 
\sqrt{\frac{\kappa_0}{\kappa_{1}}}    \frac{\Gamma(\frac{\nu_{1}}{2})}{\Gamma(\frac{\nu_{0}}{2})}
\frac{(\nu_0 \sigma_0^2)^{\frac{\nu_0}{2}}}{(\nu_{1}\sigma_{1}^2)^{\frac{\nu_{1}}{2}}}
=\frac{\Gamma(\frac{\nu_{0}+1}{2})}{\Gamma(\frac{\nu_{0}}{2})}
\sqrt{\frac{\kappa_0}{(\kappa_{0}+1)}   \frac{1}{(\pi\nu_{0}\sigma_{0}^2)}}   
\left[\frac{(\nu_{1}\sigma_{1}^2)}{(\nu_0 \sigma_0^2)}\right]^{-\frac{\nu_{0}+1}{2}}\\
&=\tau\big(x^\star \mid m_0, \frac{\kappa_{0}+1}{\kappa_{0}}\sigma^2_0, \nu_0  \big),
\end{aligned}
$$
which is also a Student's $t$ distribution centered at the prior mean $m_0$ and $\nu_0$ degrees of freedom.

\section{Exponential and Conjugacy}
\begin{SCfigure}
\centering
\includegraphics[width=0.5\textwidth]{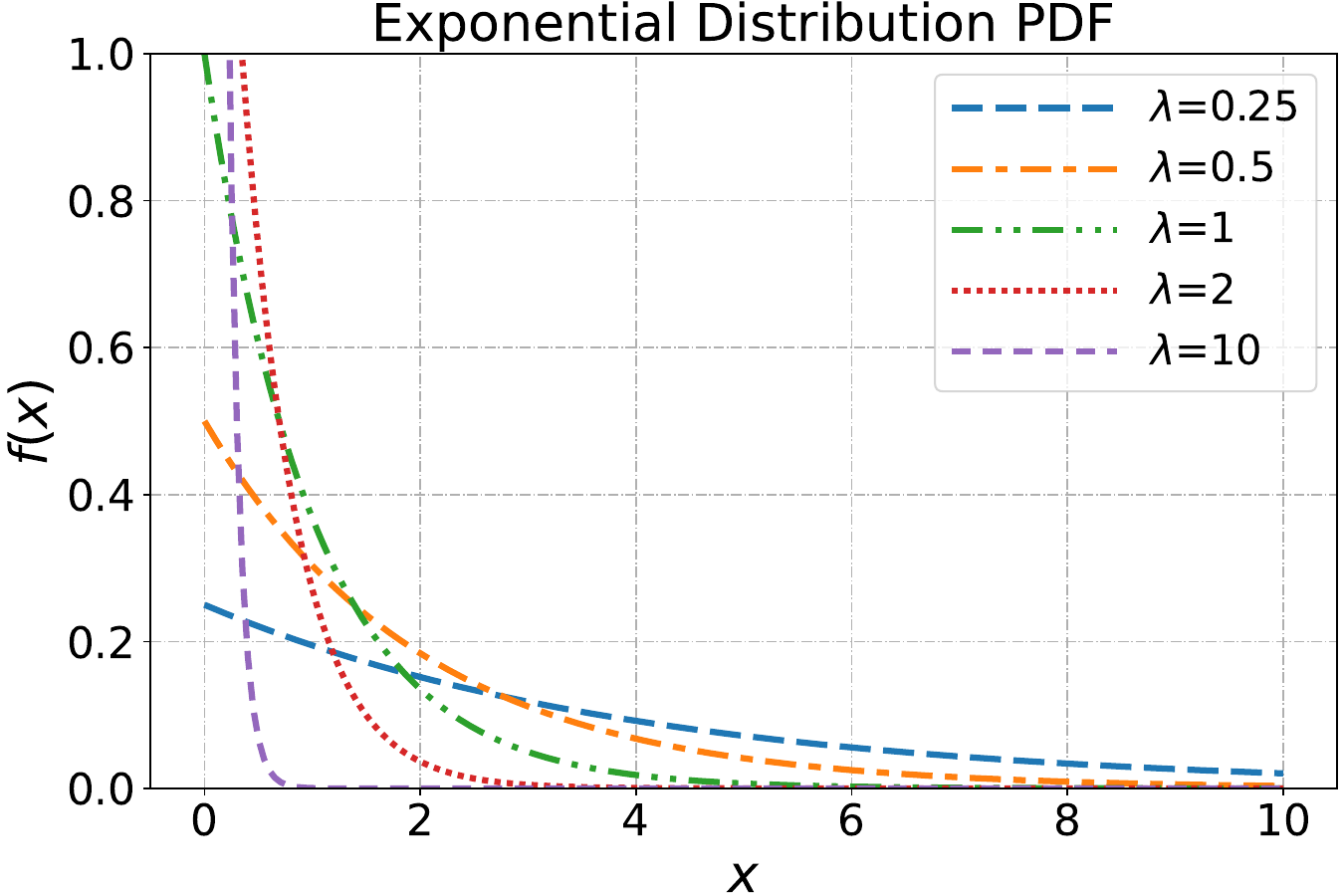}
\caption{Exponential probability density functions for different values of the rate parameter $\lambda$.}
\label{fig:dists_exponential}
\end{SCfigure}

The \textit{exponential distribution} is a continuous probability distribution commonly used to model the time until a random event occurs---such as the waiting time until the next customer arrives, the lifetime of a component, or the time between successive events in a Poisson process. It is a special case of the Gamma distribution with shape parameter equal to 1, and its support is the set of nonnegative real numbers.
\begin{definition}[Exponential Distribution\index{Exponential distribution}]\label{definition:exponential_distribution}
A random variable $\rx$ is said to follow an \textit{exponential distribution} with rate parameter $\lambda>0$, denoted $\rx \sim \exponential(\lambda)$, if its probability density function  is given by
$$ f(x; \lambda)=\left\{
\begin{aligned}
& \lambda \exp(-\lambda x)
,& \mathrm{\,\,if\,\,} x \geq 0;  \\
&0 , &\mathrm{\,\,if\,\,} x <0.
\end{aligned}
\right.
$$
This is equivalent to $\rx\sim \gammadist(1, \lambda)$, where the Gamma distribution is parameterized by shape and rate.
The mean and variance of $\rx \sim \exponential(\lambda)$ are given by 
\begin{equation}
\Exp[\rx] = \lambda^{-1}, \qquad \Var[\rx] =\lambda^{-2}. \nonumber
\end{equation}
The support of the exponential distribution is $(0,\infty)$.
Figure~\ref{fig:dists_exponential}  illustrates the PDFs of exponential distributions with different values of the rate parameter $\lambda$.
\end{definition}

Note that the mean  $\lambda^{-1}$ represents the expected waiting time until the event occurs, which justifies interpreting $\lambda$ as a rate: higher values of $\lambda$ correspond to shorter average waiting times.
An important property of the exponential distribution is \textit{memorylessness}. This means that the probability of waiting an additional amount of time $x$ does not depend on how long one has already waited. Formally:
\begin{exercise}[Memoryless of Exponential Distribution]
Let $\rx\sim \exponential(\lambda)$. Then for any $x,s\geq 0$,  show that $p(\rx\geq x + s \mid \rx \geq s) = p(\rx\geq x)$. 
That is, $\rx-s$ has an exponential distribution with parameter $\lambda$. 
This means, if $\rx$ represent the lifetime of some object under random conditions,  the memoryless property
implies that the chance that $\rx$ will ``live" longer than $x+s$ given that it has already ``lived" longer
than $s$ is the same as the chance that $\rx$ will live longer than $x$ in the first place.
\end{exercise}

\begin{exercise}[Chi-Squared and Exponential]
Show that $\rx\sim \chi^2(2)$ if and only if $\rx\sim \exponential(1/2)$.
\end{exercise}

\begin{exercise}[Transformation of Exponential]
Let  $\rx\sim\exponential(\lambda)$; show that $\ry=\frac{1}{\lambda}\rx \sim \exponential(\lambda^2)$.
Similarly,  suppose  $\rx\sim\exponential(1)$; show that $\ry=\frac{1}{\lambda}\rx \sim \exponential(\lambda)$.
\textit{Hint: use the Jacobian change of variables. Suppose $g(\cdot):\mathcalX\rightarrow \real$ is monotone and $\ry=g(\rx)$. Then $p_{\ry}(y) = \abs{\frac{\partial }{\partial y}g^{-1}(y)} p_{\rx}(g^{-1}(y)) $.}
\end{exercise}

\begin{exercise}[Sum of Independent Exponentials \citep{bibinger2013notes}]
Let $\rx_1,$$\rx_2,$ $\ldots,\rx_N$ be independent random variables with $\rx_n\sim \exponential(\lambda_n)$ and all $\lambda_n$ distinct. 
Show that 
$$
\ry=\sum_{n=1}^{N}\rx_n \sim p(y)
=
\big(\prod_{n=1}^{N}\lambda_n\big) \sum_{j=1}^{N} \frac{\exp(-\lambda_j y)}{\prod_{k\neq j}^{N}(\lambda_k-\lambda_j)}.
$$
\end{exercise}

\paragrapharrow{Conjugate prior for the exponential rate parameter.}
The Gamma distribution is a conjugate prior for the rate parameter of an exponential distribution. To see this, suppose  $\mathcalX=\{x_1, x_2, \ldots, x_N\}$ are drawn i.i.d. from an exponential distribution with rate $\lambda$, i.e., the likelihood is $\exponential(x \mid \lambda)$, and $\lambda$ is given a $\gammadist(\alpha_0, \beta_0)$ prior: $\lambda \sim \gammadist(\alpha_0, \beta_0)$. Using Bayes' theorem, the posterior is 
\begin{equation}
p(\lambda \mid \mathcalX) \propto \prod_{n=1}^{N} \exponential(x_n \mid \lambda) 
\times \gammadist(\lambda\mid \alpha_0, \beta_0)
\propto \gammadist(\theta \mid \widetildealpha, \widetildebeta),
\end{equation}
where 
\begin{equation}\label{equation:posterior-param-exponenral-gamma}
\widetildealpha= \alpha_0+N, 
\qquad 
\widetildebeta =\beta_0 + \sum_{n=1}^{N} x_n. 
\end{equation}
From this posterior form, the prior parameter $\alpha_0$ can be interpreted as the number of prior observations, and $\beta_0$ as the sum of the prior observations. 
The posterior mean of $\lambda$ is therefore
$$
\Exp[\lambda\mid\mathcalX]=
\frac{\widetildealpha}{\widetildebeta} = \frac{\alpha_0+N}{\beta_0 + \sum_{n=1}^{N} x_n},
$$
which smoothly combines prior information with observed data.

\section{Univariate Gaussian-Related Models}
We now discuss several probability distributions closely related to the univariate Gaussian distribution.

The \textit{truncated-normal (TN)} distribution is a modification of the normal distribution in which values below zero are excluded---i.e., the distribution is ``truncated" at zero. Its support is the set of nonnegative real numbers, making it suitable for applications that require nonnegativity, such as nonnegative matrix factorization.

\begin{figure}[h]
\centering  
\vspace{-0.35cm} 
\subfigtopskip=2pt 
\subfigbottomskip=2pt 
\subfigcapskip=-5pt 
\subfigure[Truncated-normal.]{\label{fig:dists_truncatednorml}
\includegraphics[width=0.481\linewidth]{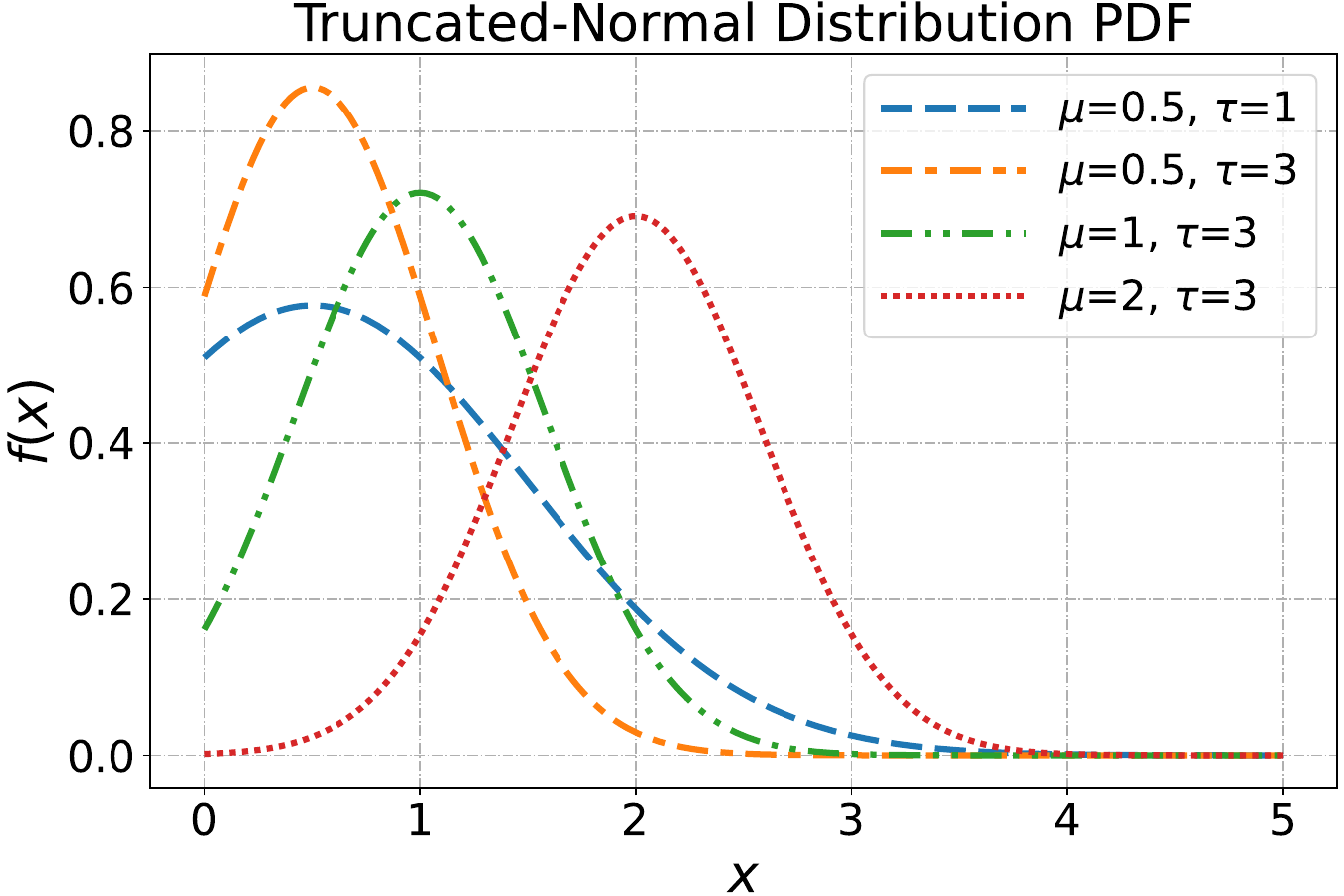}}
\subfigure[General-truncated-normal.]{\label{fig:dists_generaltruncatednorml}
\includegraphics[width=0.481\linewidth]{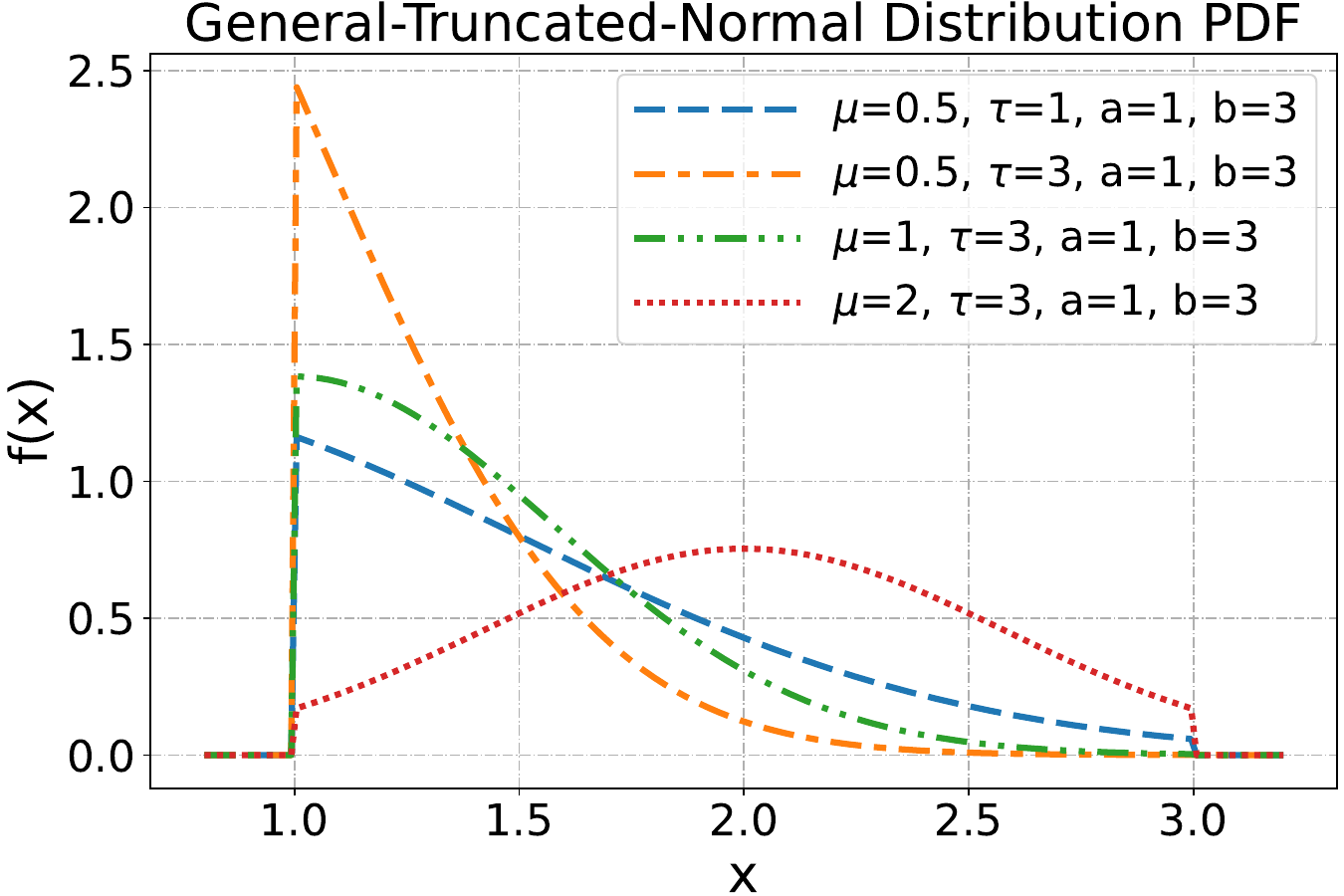}}
\caption{Truncated-normal  and general-truncated-normal  probability density functions for different values of the parameters $\mu$ and $\tau$.}
\label{fig:dists_truncatednorml_and_general}
\end{figure}

\begin{figure}[h]
\centering  
\vspace{-0.35cm} 
\subfigtopskip=2pt 
\subfigbottomskip=2pt 
\subfigcapskip=-5pt 
\subfigure[Truncated-normal.]{\label{fig:dists_truncatednorml_mean}
\includegraphics[width=0.481\linewidth]{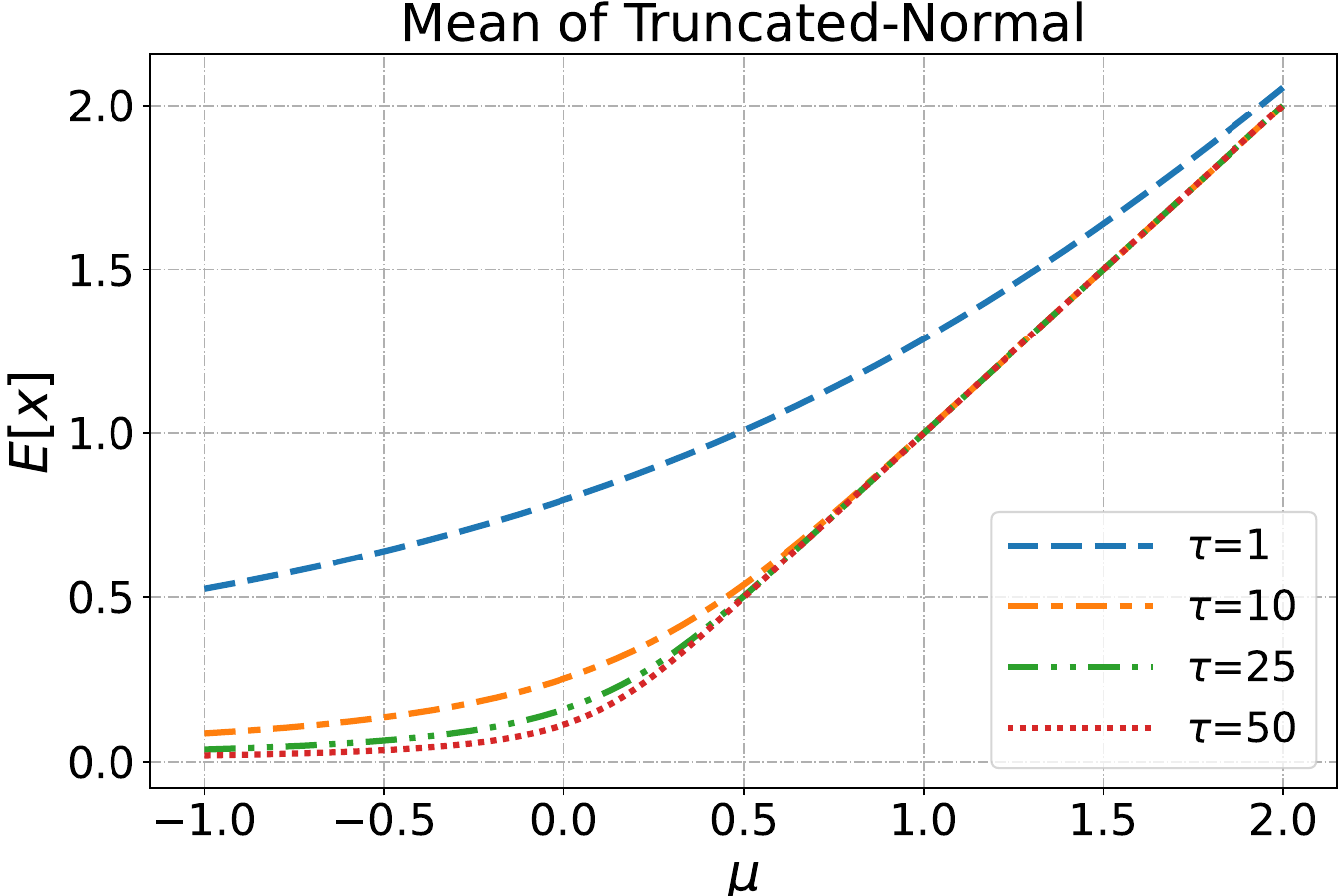}}
\subfigure[General-truncated-normal.]{\label{fig:dists_generaltruncatednorml_mean}
\includegraphics[width=0.481\linewidth]{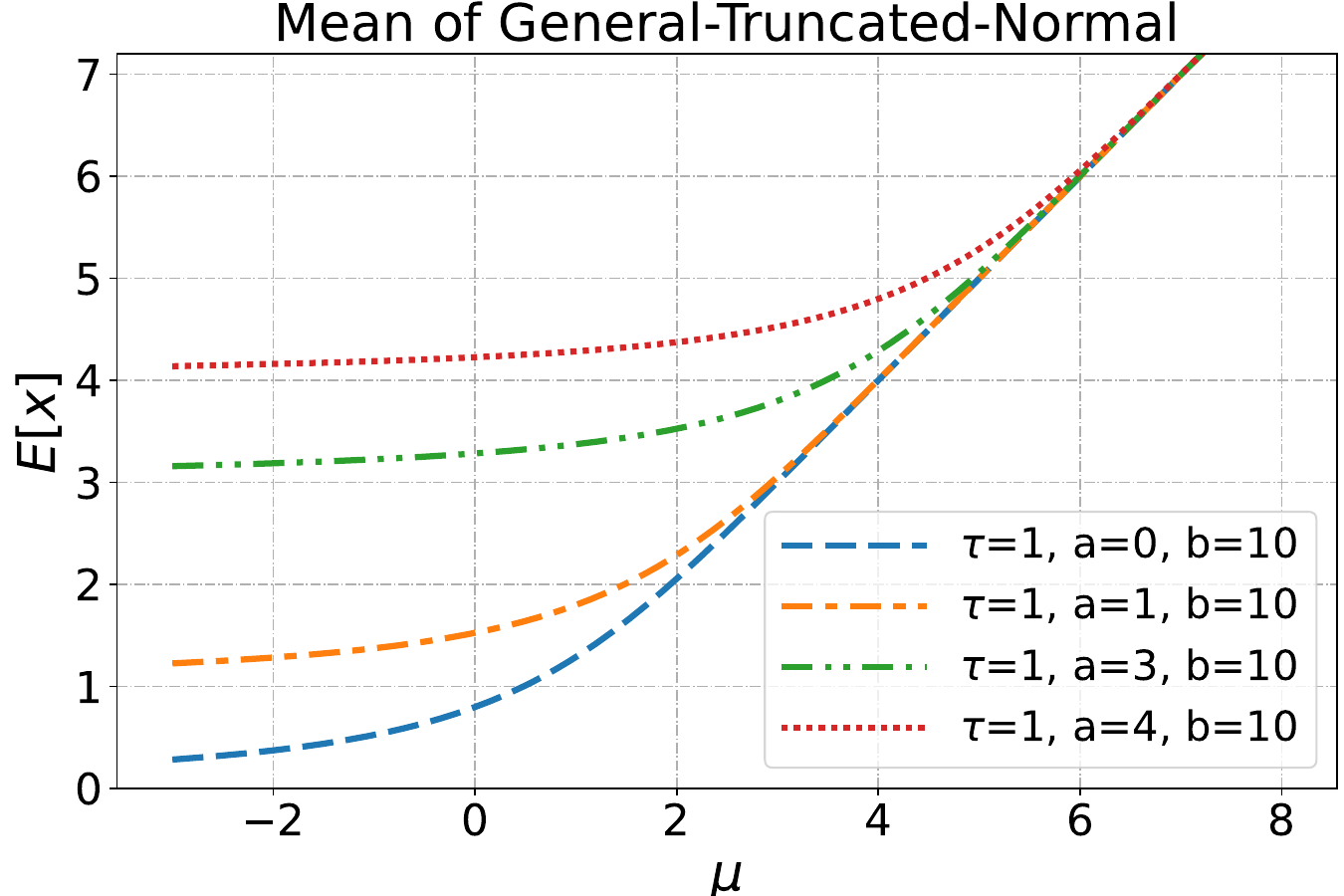}}
\caption{Mean of truncated-normal and general-truncated-normal distribution by varying $\mu, \tau, a$, and $b$ parameters.}
\label{fig:mean_truncatednorml_and_general}
\end{figure}

\begin{definition}[Truncated-Normal (TN) Distribution\index{Truncated-normal distribution}]\label{definition:truncated_normal}
A random variable $\rx$ is said to follow an \textit{truncated-normal (TN) distribution} with ``parent" mean $\mu$ and ``parent" precision $\tau>0$, denoted $\rx \sim \truncatednormal(\mu, \tau^{-1})$, if its probability density function  is given by
$$ f(x; \mu, \tau^{-1})=\left\{
\begin{aligned}
&\frac{\sqrt{\frac{\tau}{2\pi}} \exp \{-\frac{\tau}{2}(x-\mu)^2  \}  }{1-\Phi(-\mu\sqrt{\tau})}
,& \mathrm{\,\,if\,\,} x \geq 0;  \\
&0 , &\mathrm{\,\,if\,\,} x <0,
\end{aligned}
\right.
$$
where $\Phi(y) = \int_{-\infty}^{y} \normal(u\mid 0,1)\,du= \frac{1}{\sqrt{2\pi}} \int_{-\infty}^{y} \exp(-\frac{u^2}{2}) \,du $ is the cumulative distribution function of $\normal(0,1)$, the standard normal distribution. 
Generally, the cumulative density function of $\ry\sim \normal(\mu, \sigma^2)$ can be written as 
$$
F(y) = p(\ry\leq y) 
= \Phi\left(\frac{y-\mu}{\sigma}\right)
= \Phi\left((y-\mu)\cdot \sqrt{\tau}\right).~
\footnote{
Or equivalently, for general Gaussian distribution $\ry\sim\normal(\mu,\sigma^2)$, the CDF is $F(y)=\frac{1}{2}\left\{ 
1+\text{erf}\left(\frac{y-\mu}{\sigma\sqrt{2}}\right)
\right\}$, where the error function is $\text{erf}(t)=\frac{2}{\sqrt{\pi}} \int_{0}^{t}\exp(-y^2) \,dy$.}
$$
The mean and variance of $\rx \sim \truncatednormal(\mu, \tau^{-1})$ are given by 
$$
\begin{aligned}
\Exp[\rx] &= \mu - \frac{1}{\sqrt{\tau}}\cdot \frac{ - \phi(\alpha)}{1 - \Phi(\alpha)}, \\
\Var[\rx] &= \frac{1}{\tau} 
\left( 1 + 
\frac{  \alpha\phi(\alpha)}{1 - \Phi(\alpha)}+
\left(\frac{ \alpha\phi(\alpha)}{1 - \Phi(\alpha)}\right)^2
\right),
\end{aligned}
$$
where $\phi(y)=\frac{1}{\sqrt{2\pi} } \exp(-\frac{y^2}{2})$ is the PDF of the standard normal distribution, and 
$\alpha = -\mu\cdot \sqrt{\tau}$ \citep{burkardt2014truncated}.
Figure~\ref{fig:dists_truncatednorml} compares different parameters $\mu, \tau$ for the TN distribution. Figure~\ref{fig:dists_truncatednorml_mean} shows the mean value of the TN distribution by varying $\mu$ given fixed $\tau$; we can find when $\mu\rightarrow -\infty$, the mean is approaching zero.
\end{definition}

\paragrapharrow{Conjugate prior for the nonnegative mean parameter of a Gaussian.}
As previously shown, a Gaussian distribution serves as a conjugate prior for the mean of another Gaussian distribution when the variance (or precision) is known.
Similarly, the truncated-normal distribution  also acts as a conjugate prior for the \textbf{nonnegative} mean parameter of a Gaussian distribution when the variance is fixed. To see this, suppose  $\mathcalX=\{x_1, x_2, \ldots, x_N\}$ are drawn i.i.d. from a normal distribution with mean $\theta$ and precision $\tau$, i.e., the likelihood is $\normal(x \mid \theta, \tau^{-1})$ with fixing the variance $\sigma^2=\tau^{-1}$, and $\theta$ is given a $\truncatednormal(\mu_0, \tau^{-1}_0)$ prior: $\theta \sim \truncatednormal(\mu_0, \tau_0^{-1})$. Using Bayes' theorem, the posterior is 
\begin{equation}\label{equation:conjugate_truncated_nonnegative_mean}
\begin{aligned}
&\gap p(\theta \mid \mathcalX) \propto \prod_{n=1}^{N} \normal(x_n \mid \theta, \tau^{-1}) 
\times \truncatednormal(\theta \mid \mu_0, \tau_0^{-1})\\
&\propto 
\exp\left\{ -\frac{\tau_0+ N\tau}{2} \theta^2+ \bigg(\tau \sum_{n=1}^{N}x_n + \tau_0 \mu_0\bigg)\theta  \right\}\cdot u(\theta)
\propto \truncatednormal(\theta \mid \widetilde{\mu}, \widetilde{\tau}^{-1}), 
\end{aligned}
\end{equation}
where $u(y)$ is the step function with value 1 if $y\geq 0$ and value 0 if $y<0$, and
$$
\widetilde{\mu}= \frac{\tau_0 \mu_0 + \tau \sum_{n=1}^{N} x_n}{ \tau_0 + N\tau}, 
\qquad 
\widetilde{\tau} = \tau_0 + N \tau. 
$$
The posterior parameters are exactly the same as those in the Normal-Normal model (see Equation~\eqref{equation:posterior-param-normal-normal}), and the posterior ``parent" mean can also be expressed as a weighted mean of $\mu_0$ and $\widebarx$.
The only difference is the truncation at zero, which ensures the posterior remains supported on $[0,\infty)$.

Expanding beyond   the truncated-normal distribution,
the \textit{general-truncated-normal (GTN)} distribution is another variant of the normal distribution, excluding values outside a specified range. In other words, it is a normal distribution that is ``cut off" at some specified lower and/or upper bound. The range between the lower and upper bound is called the support of the distribution. 
\footnote{In some literature, the term ``truncated-normal" refers to this general form. 
Here, we distinguish the two: the truncated-normal (TN) denotes truncation at zero (i.e., $a=0,b=\infty$), while the general-truncated-normal (GTN) allows arbitrary finite or infinite bounds.}
\begin{definition}[General-Truncated-Normal (GTN) Distribution\index{General-truncated-normal distribution}]\label{definition:general_truncated_normal}
A random variable $\rx$ is said to follow a \textit{general-truncated-normal (GTN) distribution} with ``parent" mean $\mu$ and ``parent" precision $\tau>0$, denoted $\rx \sim \generaltruncatednormal(\mu, \tau^{-1}, \textcolor{mylightbluetext}{a, b})$, if its probability density function is
$$ f(x; \mu, \tau^{-1}, a, b)=\left\{
\begin{aligned}
&0,  & \mathrm{\,\,if\,\,} x < a;  \\
&\frac{\sqrt{\frac{\tau}{2\pi}} \exp \{-\frac{\tau}{2}(x-\mu)^2  \}  }{\Phi((b-\mu)\cdot \sqrt{\tau})-\Phi((a-\mu)\cdot \sqrt{\tau})}
,& \mathrm{\,\,if\,\,} a\leq x \leq b;  \\
&1 , &\mathrm{\,\,if\,\,} 0 >b,
\end{aligned}
\right.
$$
where $\Phi(\cdot)$ is the CDF of $\normal(0,1)$, the standard normal distribution.
The mean and variance of $\rx \sim \generaltruncatednormal(\mu, \tau^{-1}, a, b)$ are given by 
$$
\begin{aligned}
\Exp[\rx] &= \mu - \frac{1}{\sqrt{\tau}}\cdot \frac{\phi(\beta) - \phi(\alpha)}{\Phi(\beta) - \Phi(\alpha)}, \\
\Var[\rx] &= \frac{1}{\tau} 
\left( 1 - 
\frac{\beta \phi(\beta) - \alpha\phi(\alpha)}{\Phi(\beta) - \Phi(\alpha)}-
\left(\frac{\beta \phi(\beta) - \alpha\phi(\alpha)}{\Phi(\beta) - \Phi(\alpha)}\right)^2
  \right),
\end{aligned}
$$
where $\phi(\cdot)$ is the PDF of the standard normal distribution, and 
$$
\alpha = (a-\mu)\cdot \sqrt{\tau}, \qquad 
 \beta = (b-\mu)\cdot \sqrt{\tau}.
$$
Note that, the truncated-normal distribution is a special general-truncated-normal with $a=0$ and $b=\infty$ \citep{burkardt2014truncated}.
Figure~\ref{fig:dists_generaltruncatednorml} compares different parameters $\mu, \tau$ for the GTN distribution.
Figure~\ref{fig:dists_generaltruncatednorml_mean} shows the mean value of the GTN distribution by varying $\mu$ given fixed $\tau, a, b$; we again find when $\mu\rightarrow -\infty$, the mean is approaching zero.
\end{definition}

\paragrapharrow{Conjugate prior for the constrained mean parameter of a Gaussian.}
We have previously shown that a Gaussian distribution can serve as a conjugate prior for the mean parameter of another Gaussian distribution when the variance is fixed. 
Similarly, the general-truncated-normal distribution acts as a conjugate prior for a Gaussian mean that is \textbf{constrained} to lie within a known interval $[a,b]$. This includes the nonnegative case ($a=0,b=\infty$) as a special instance.
To see this, suppose  $\mathcalX=\{x_1, x_2, \ldots, x_N\}$ are drawn i.i.d. from a  normal distribution with mean $\theta$ and precision $\tau$, i.e., the likelihood is $\normal(x \mid \theta, \tau^{-1})$, where the variance $\sigma^2=\tau^{-1}$ is fixed, and $\theta$ is given a $\generaltruncatednormal(\mu_0, \tau^{-1}_0, a, b)$ prior: $\theta \sim \generaltruncatednormal( \mu_0, \tau_0^{-1}, a, b)$, with $a>0$ for the nonnegativity constrain. Using Bayes' theorem, the posterior is 
\begin{equation}\label{equation:conjugate_truncated_constrained_mean}
\begin{aligned}
p(\theta \mid \mathcalX) 
&\propto \prod_{n=1}^{N} \normal(x_n \mid \theta, \tau^{-1}) 
\times \generaltruncatednormal(\theta \mid \mu_0, \tau_0^{-1}, a, b)\\
&\propto \exp\left\{ -\frac{\tau_0+ N\tau}{2} \theta^2
+ \Big(\tau \sum_{n=1}^{N}x_n + \tau_0 \mu_0\Big)\theta  \right\} 
\cdot \indicator(a\leq \theta\leq b)\\
&\propto \generaltruncatednormal(\theta \mid \widetilde{\mu}, \widetilde{\tau}^{-1}, a, b), 
\end{aligned}
\end{equation}
where $\indicator(a\leq y\leq b)$ is the step function with value 1 if $ a\leq y\leq b$ and value 0 otherwise, and
$$
\widetilde{\mu}= \frac{\tau_0 \mu_0 + \tau \sum_{n=1}^{N} x_n}{ \tau_0 + N\tau}, 
\qquad 
\widetilde{\tau} = \tau_0 + N \tau. 
$$
The posterior parameters are again exactly the same as those in the Normal-Normal model (see Equation~\eqref{equation:posterior-param-normal-normal}).
The only difference is the truncation, which enforces the constraint $\theta\in[a,b]$ in both prior and posterior.

\begin{figure}[h]
\centering  
\vspace{-0.35cm} 
\subfigtopskip=2pt 
\subfigbottomskip=2pt 
\subfigcapskip=-5pt 
\subfigure[Half-normal.]{\label{fig:dists_halfnorml}
\includegraphics[width=0.481\linewidth]{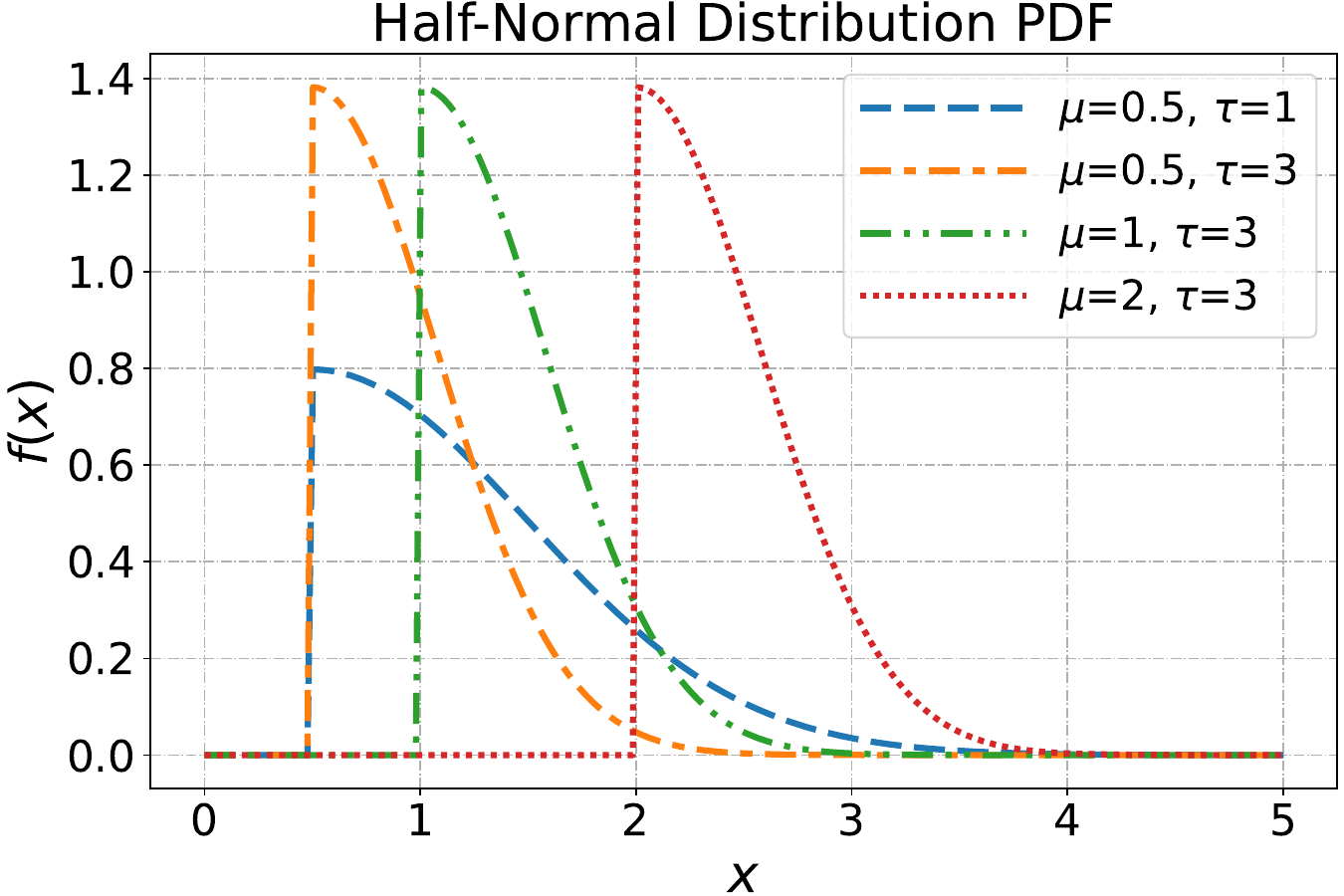}}
\subfigure[Normal, $\tau$ is the precision parameter.]{\label{fig:dists_halfnorml_gausscompare}
\includegraphics[width=0.481\linewidth]{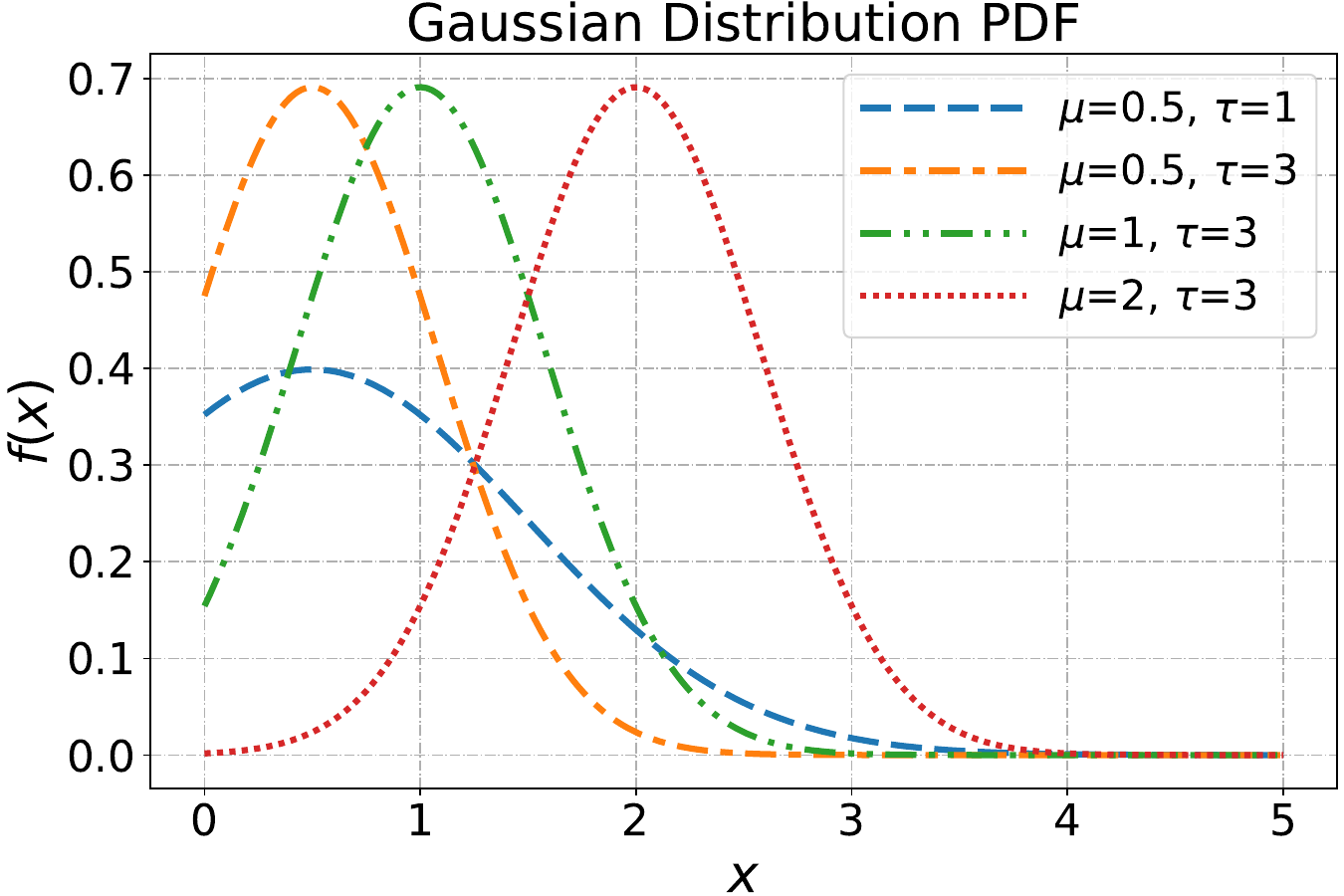}}
\caption{Half-normal and normal probability density functions for different values of the parameters $\mu$ and $\tau$.
The probability density of any value $x\geq \mu$ in the half-normal distribution is twice as that in the normal distribution with the same parameters $\mu, \tau$.
}
\label{fig:dists_half_and_general}
\end{figure}

Different to the truncated-normal distribution.
the \textit{half-normal (HN)} distribution  is another special case of the normal distribution obtained by restricting its support to values greater than or equal to the parent mean $\mu$. The resulting distribution is asymmetric, with its mode at $\mu$ and a tail extending to $+\infty$. 
It is also commonly used to model nonnegative quantities such as scales, standard deviations, or magnitudes---particularly when the underlying symmetric process is reflected or truncated at $\mu$.

\begin{definition}[Half-Normal (HN) Distribution\index{Half-normal distribution}]\label{definition:half_normal}
A random variable $\rx$ is said to follow a \textit{half-normal (HN) distribution} with ``parent" mean $\mu$ and ``parent" precision $\tau>0$, denoted $\rx \sim \halfnormal(\mu, \tau^{-1})$, if 
$$ f(x; \mu, \tau^{-1})=\left\{
\begin{aligned}
	&\sqrt{\frac{2\tau}{\pi}} \exp\left\{- \frac{\tau}{2}(x-\mu)^2 \right\}
	,& \mathrm{\,\,if\,\,} x \geq \mu;  \\
	&0 , &\mathrm{\,\,if\,\,} x <\mu.
\end{aligned}
\right.
$$
The mean and variance of $\rx \sim \halfnormal(\mu, \tau^{-1})$ are given by 
$$
\Exp[\rx] = \mu + \sqrt{\frac{2}{\pi \tau}} , \qquad
\Var[\rx] = \frac{1}{\tau}\big(1-\frac{2}{\pi} \big). 
$$
Figure~\ref{fig:dists_halfnorml} compares different parameters $\mu, \tau$ for the half-normal distribution.
\end{definition}

\index{Rectified-normal distribution}
\index{Exponentially rectified-normal distribution}

\begin{figure}[h]
\centering  
\vspace{-0.35cm} 
\subfigtopskip=2pt 
\subfigbottomskip=2pt 
\subfigcapskip=-5pt 
\subfigure[Truncated-normal. Same as Figure~\ref{fig:dists_truncatednorml}.]{\label{fig:dists_truncatednorml2}
\includegraphics[width=0.481\linewidth]{./imgs/dists_truncatednorml.pdf}}
\subfigure[Rectified-normal.]{\label{fig:dists_rectifiednorml}
\includegraphics[width=0.481\linewidth]{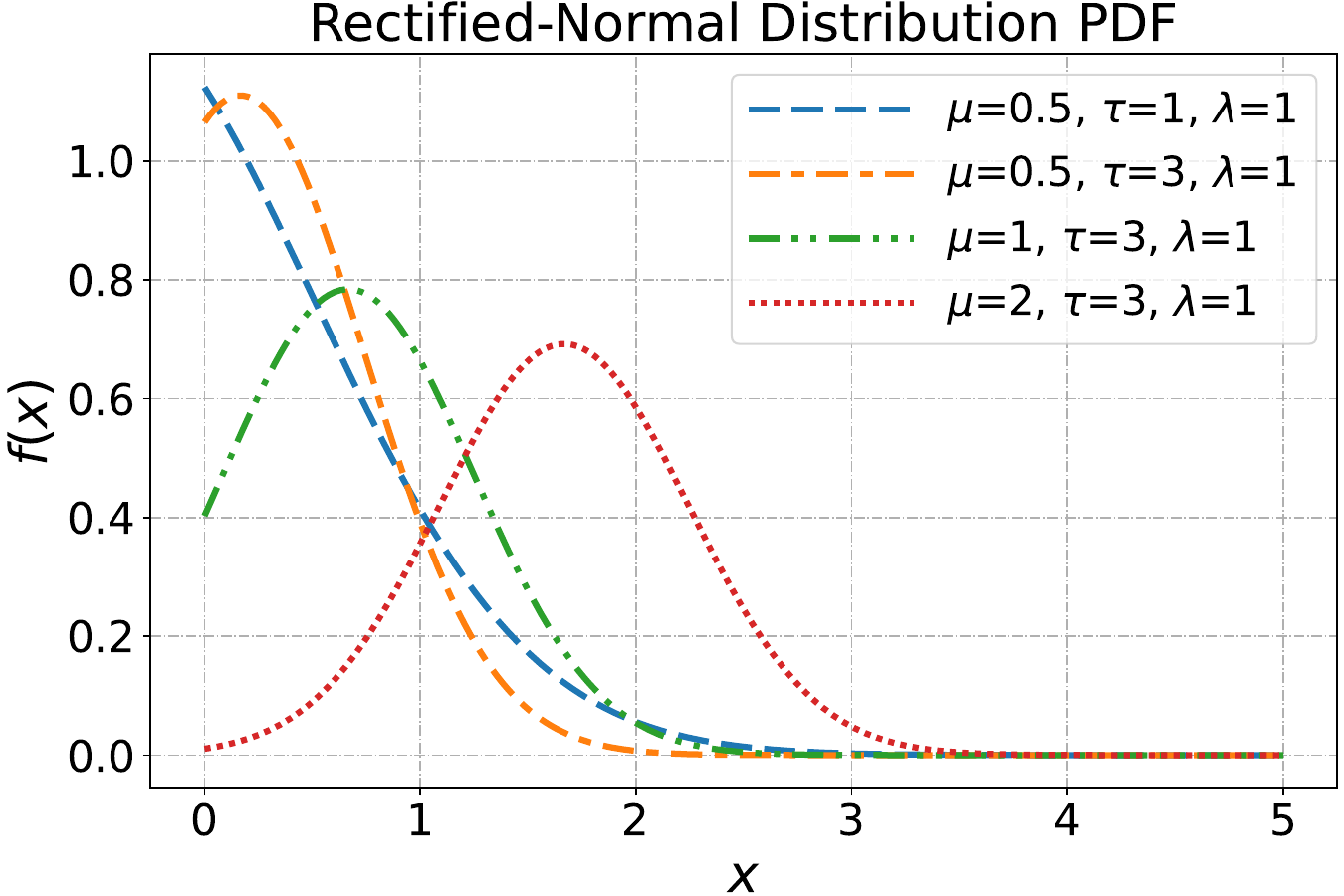}}
\caption{Truncated-normal and rectified-normal probability density functions for different values of the parameters $\mu$, $\tau$, and $\lambda$.}
\label{fig:dists_truncatednorml_and_rectified}
\end{figure}

In this text, the \textit{rectified-normal (RN)} distribution is defined as a distribution whose density is proportional to the product of a Gaussian distribution and an exponential distribution (restricted to nonnegative values). 
Thus, it is also referred to as the \textit{exponentially rectified-normal} distribution.
\begin{definition}[Rectified-Normal (RN) Distribution]\label{definition:reftified_normal_distribution}
A random variable $\rx$ is said to follow a \textit{rectified-normal (RN) distribution} (or \textit{exponentially rectified-normal} distribution) with ``parent" mean $\mu$, ``parent" precision $\tau>0$, and ``parent" rate $\lambda >0$, denoted $\rx \sim \rectifieddist(\mu, \tau^{-1}, \lambda)$, if 
$$ 
\begin{aligned}
f(x; \mu, \tau^{-1}, \lambda) 
&=\frac{1}{C}\cdot \normal(x\mid  \mu, \tau^{-1}) \cdot \exponential(x\mid \lambda)\\
&\propto \exp\left\{ -\frac{\tau}{2} \left( x - \frac{\tau\mu-\lambda}{\tau} \right)^2 \right\} \cdot u(x)
\propto \truncatednormal(x\mid\textcolor{mylightbluetext}{\frac{\tau\mu-\lambda}{\tau}}, \tau^{-1}),
\end{aligned}
$$
where $\truncatednormal(\cdot)$ is the density function of a truncated-normal distribution, and $C$ is a constant  value,
\begin{equation}\label{equation:rf_constant}
C =C^{RN}(\mu, \tau, \lambda)
= \lambda \left\{1 - \Phi\bigg(-\frac{\tau\mu-\lambda}{\sqrt{\tau}} \bigg)\right\}
\cdot \exp\left( -\mu\lambda + \frac{\lambda^2}{2\tau}\right).
\end{equation}
That is, the rectified-normal distribution is a special truncated-normal distribution with more flexibility.
The mean and variance of $\rx \sim \rectifieddist(\mu, \tau^{-1}, \lambda)$ are given by
$$
\begin{aligned}
\Exp[\rx] &= \frac{\tau\mu-\lambda}{\tau} - \frac{1}{\sqrt{\tau}}\cdot \frac{ - \phi(\alpha)}{1 - \Phi(\alpha)}, \\
\Var[\rx] &= \frac{1}{\tau} 
\left\{ 1 + 
\frac{  \alpha\phi(\alpha)}{1 - \Phi(\alpha)}+
\left(\frac{ \alpha\phi(\alpha)}{1 - \Phi(\alpha)}\right)^2
\right\},
\end{aligned}
$$
where
$
\alpha = -\frac{\tau\mu-\lambda}{\tau}\cdot \sqrt{\tau}.
$
Figure~\ref{fig:dists_rectifiednorml} compares different parameters $\mu, \tau$ for the RN distribution.
\end{definition}
The comparison between the truncated-normal and rectified-normal distributions is presented in Figure~\ref{fig:dists_truncatednorml_and_rectified}. 

\paragrapharrow{Conjugate prior for the nonnegative mean parameter of a Gaussian by RN.}
Similar to the TN distribution, the RN distribution also serves to enforce nonnegative constraint, and is conjugate to the \textbf{nonnegative} mean parameter of a Gaussian likelihood. However, due to the extra parameter $\lambda$, the RN distribution is more flexible in this sense. 
The derivation follows similarly from Equation~\eqref{equation:conjugate_truncated_nonnegative_mean}.
To see this, suppose  $\mathcalX=\{x_1, x_2, \ldots, x_N\}$ are drawn i.i.d. from a  normal distribution with mean $\theta$ and precision $\tau$, i.e., the likelihood is $\normal(x \mid \theta, \tau^{-1})$ where the variance $\sigma^2=\tau^{-1}$ is fixed, and $\theta$ is given a $\rectifieddist(\mu_0, \tau^{-1}_0, \lambda_0)$ prior: $\theta \sim \rectifieddist( \mu_0, \tau_0^{-1}, \lambda_0)$. Using Bayes' theorem, the posterior is 
\begin{equation}\label{equation:conjugate_rectified_nonnegative_mean}
\begin{aligned}
&\gap p(\theta \mid \mathcalX) \propto \prod_{n=1}^{N} \normal(x_n \mid \theta, \tau^{-1}) 
\times \rectifieddist(\theta\mid \mu_0, \tau_0^{-1}, \lambda_0) \\
&= \prod_{n=1}^{N} \normal(x_n \mid \theta, \tau^{-1})  
\times \truncatednormal(\theta\mid m_0, \tau_0^{-1})
&(m_0=\frac{\tau_0\mu_0-\lambda_0}{\tau_0})\\
&\propto \truncatednormal(\theta \mid \widetilde{\mu}, \widetilde{\tau}^{-1}), 
\end{aligned}
\end{equation}
where the updated hyper-parameters are
$$
\widetilde{\mu}= \frac{\tau_0 m_0 + \tau \sum_{n=1}^{N} x_n}{ \tau_0 + N\tau}, 
\qquad
\widetilde{\tau} = \tau_0 + N \tau. 
$$
Since the posterior is a TN distribution, it can also be interpreted as an RN distribution with appropriately updated parameters. Thus, the RN family is closed under Bayesian updating for this model---it is a valid conjugate prior for the nonnegative Gaussian mean.

\begin{SCfigure}
\centering
\includegraphics[width=0.5\textwidth]{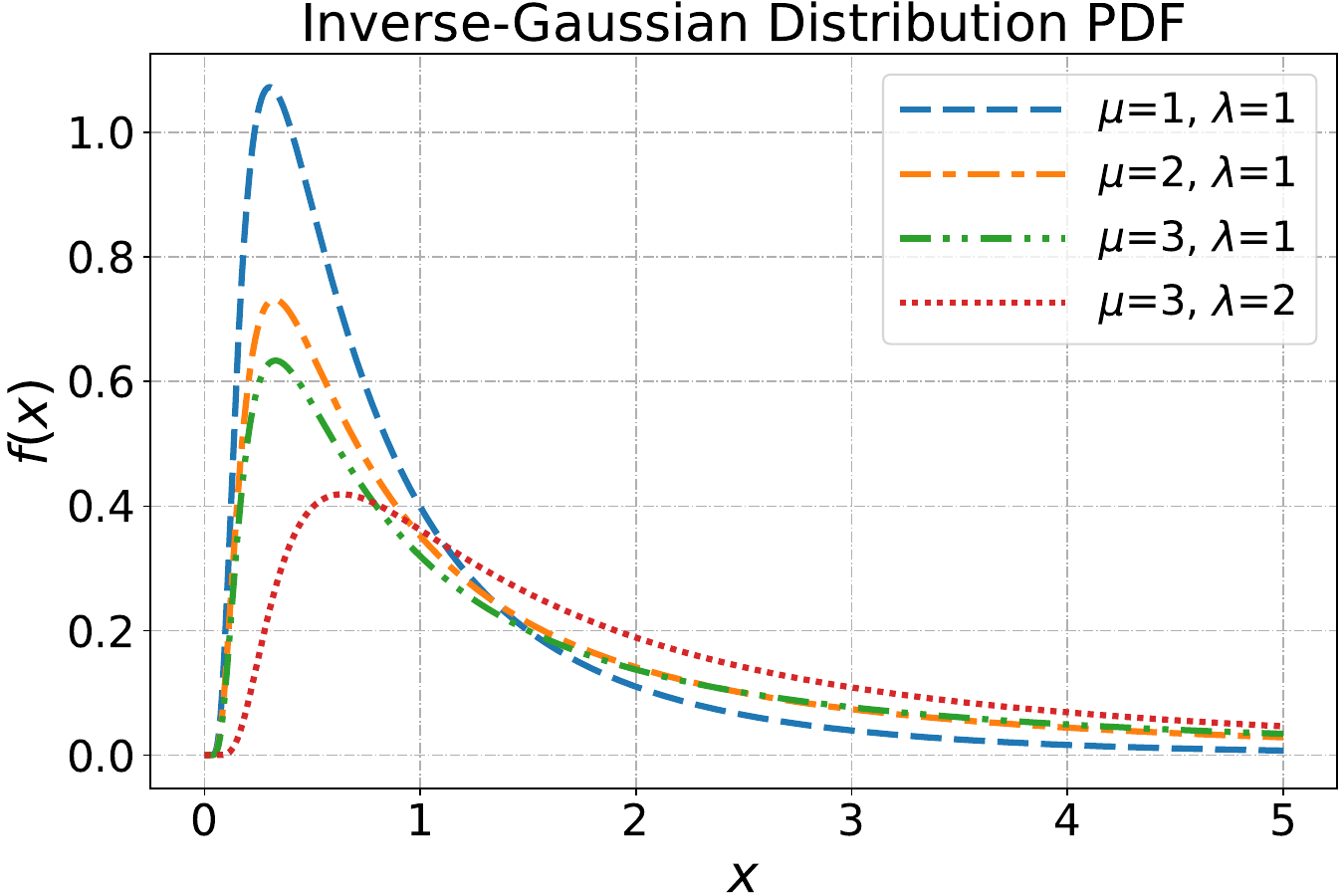}
\caption{Inverse-Gaussian probability density functions for different values of the parameters $\mu$ and $\lambda$.}
\label{fig:dists_inversegaussian}
\end{SCfigure}

The \textit{inverse-Gaussian} distribution, also known as the \textit{Wald} distribution, is a continuous probability distribution with two parameters, $\mu > 0$ and $\lambda > 0$. 
It is a versatile distribution that is used in various applications including modeling waiting times, stock prices, and lifetimes of mechanical systems.  
An important property of the distribution is that it is well-suited for modeling nonnegative, continuous, and positively skewed data with finite mean and variance.
\begin{definition}[Inverse-Gaussian Distribution\index{Inverse-Gaussian distribution}\index{Wald distribution}]\label{definition:inverse_gaussian_distribution}
A random variable $\rx$ is said to follow an \textit{inverse-Gaussian distribution} with mean parameter $\mu>0$ and shape (or precision-like) parameter $\lambda>0$, denoted $\rx \sim \inversenormaldist(\mu,\lambda)$ \footnote{Note that we use $\inversenormaldist$ to denote the inverse-Gaussian distribution and use $\inversegammadist$ to denote the inverse-Gamma distribution.}, if its probability density function  is
$$
f(x; \mu,\lambda)=\left\{
\begin{aligned}
& \sqrt{\frac{\lambda}{2\pi x^3}} \exp \left( -\frac{\lambda(x-\mu)^2}{2\mu^2 x} \right)
,& \mathrm{\,\,if\,\,} x \geq 0;  \\
&0 , &\mathrm{\,\,if\,\,} x <0.
\end{aligned}
\right.
$$
The mean and variance of $\rx \sim \inversenormaldist( \mu,\lambda)$ are given by 
\begin{equation}
\Exp[\rx] = \mu, \qquad \Var[\rx] = \frac{\mu^3}{\lambda}. \nonumber
\end{equation}
The support of an inverse-Gaussian distribution is on $(0,\infty)$.
Figure~\ref{fig:dists_inversegaussian} compares different parameters $\mu, \lambda$ for the inverse-Gaussian distribution.
\end{definition}

The \textit{Laplace} distribution, also known as  the \textit{double exponential} distribution, is named after \textit{Pierre-Simon Laplace} (1749--1827), who first derived  the distribution in 1774 \citep{kotz2001laplace, hardle2007applied}.
The Laplace distribution is useful in modeling heavy-tailed data since it has heavier tails than the normal distribution, and it is used extensively in sparse-favoring models since it expresses a high peak with heavy tails (same as the $\ell_1$-regularization term in non-probabilistic or non-Bayesian optimization methods). 
When we have a prior  belief that the parameter of interest is likely to be close to the mean  with the potential for large deviations, the Laplace distribution is then used in Bayesian modeling as a prior distribution for this context.

\begin{figure}[h]
\centering  
\vspace{-0.35cm} 
\subfigtopskip=2pt 
\subfigbottomskip=2pt 
\subfigcapskip=-5pt 
\subfigure[Laplace distribution.]{\label{fig:dists_laplace}
\includegraphics[width=0.481\linewidth]{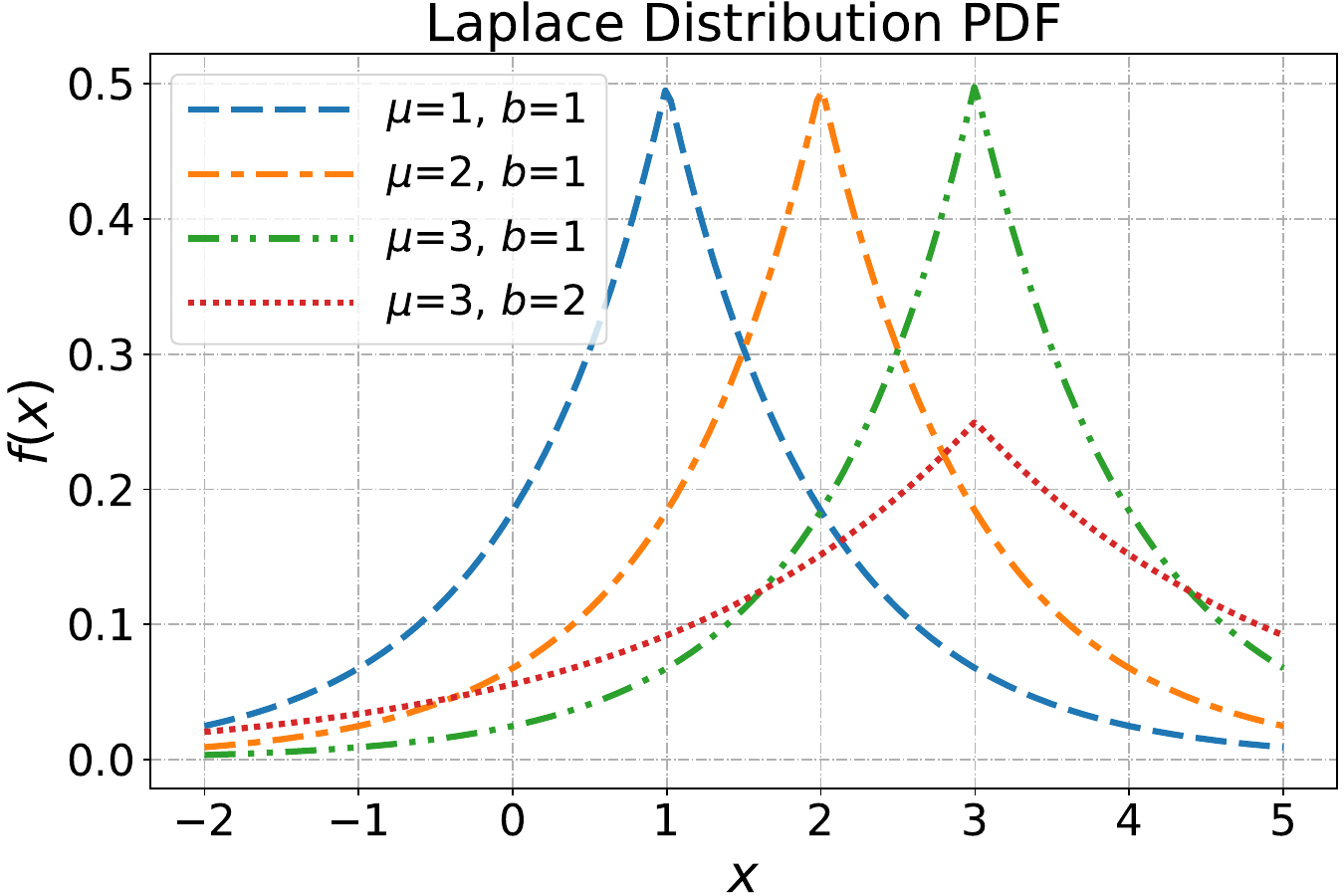}}
\subfigure[Skew-Laplace distribution.]{\label{fig:dists_skewlaplace}
\includegraphics[width=0.481\linewidth]{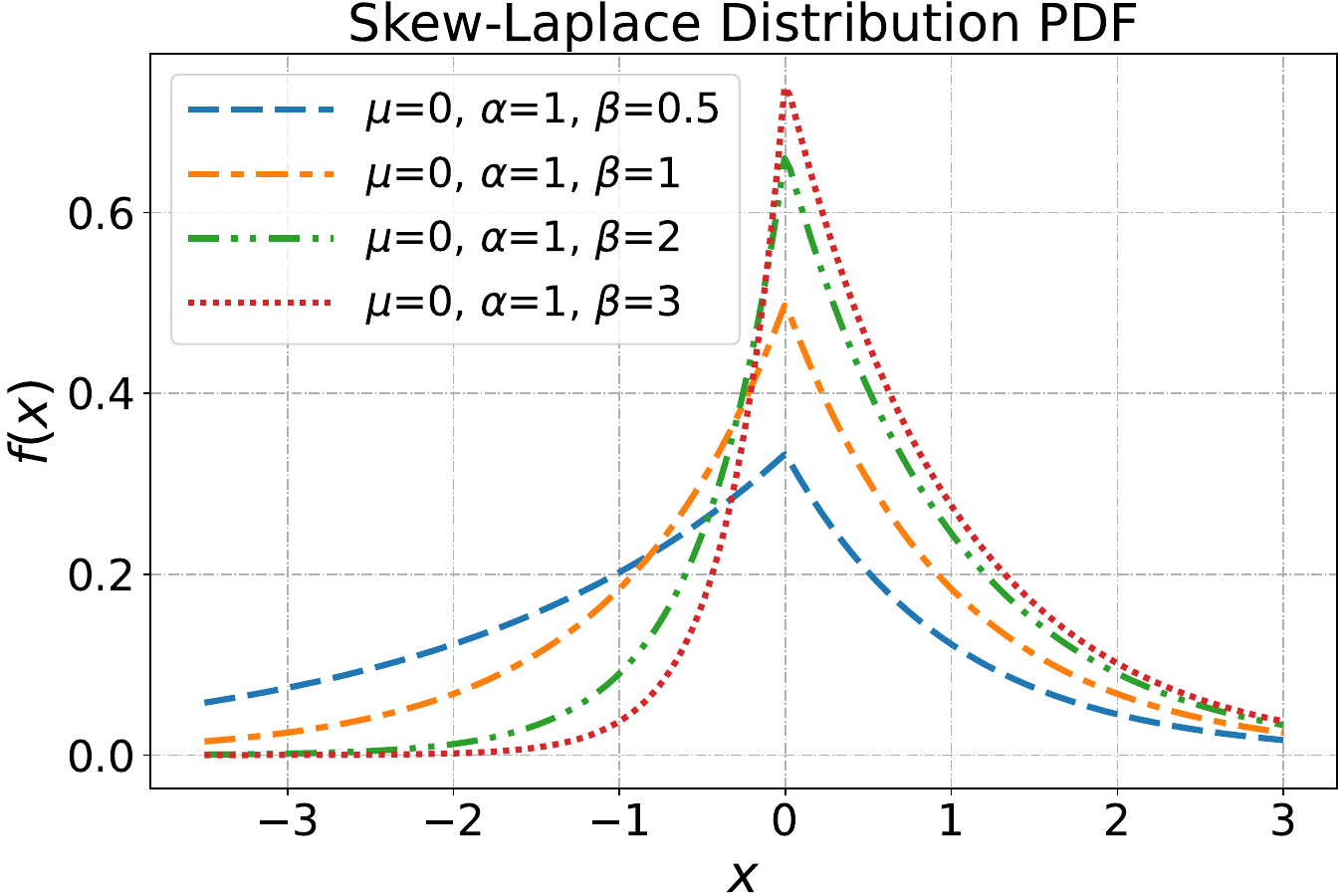}}
\caption{Laplace and skew-Laplace probability density functions for different values of the parameters.}
\label{fig:dists_lap_and_skewlap}
\end{figure}

\begin{definition}[Laplace Distribution\index{Laplace distribution}\index{Double exponential distribution}]\label{definition:laplace_distribution}
A random variable $\rx$ is said to follow a \textit{Laplace distribution} with location and scale parameters $\mu$ and $b>0$, respectively, denoted $\rx \sim \laplacedist(\mu,b)$, if its PDF is
$$
f(x; \mu,b)=\frac{1}{2b} \exp \left( -\frac{\abs{x-\mu}}{b } \right).
$$
The mean and variance of $\rx \sim \laplacedist( \mu,b)$ are given by 
\begin{equation}
\Exp[\rx] = \mu, \qquad \Var[\rx] =2b^2. \nonumber
\end{equation}
Figure~\ref{fig:dists_laplace} compares different parameters $\mu$ and $b$ for the Laplace distribution.
\end{definition}

\paragrapharrow{Laplace as a mixture of normal distributions.}
Any Laplace random variable can be thought of as the integration of a Gaussian random variable with the same mean value and a \textit{stochastic} variance that follows an exponential distribution. More formally, the Laplace distribution can be rewritten as:
\begin{equation}\label{equation:laplace_as_mixture}
\laplacedist(x\mid\mu, b) 
= \int_{0}^{\infty} \normal(x\mid \mu, \epsilon) \cdot \exponential(\epsilon \mid \frac{1}{2b^2})\,d\epsilon.
\end{equation}
To see this, we have 
$$\footnotesize
\begin{aligned}
& \int_{0}^{\infty} \normal(x\mid\mu, \epsilon) 
\cdot \exponential(\epsilon \mid \frac{1}{2b^2})\,d\epsilon 
= \int_{0}^{\infty} 
\frac{1}{\sqrt{2\pi \epsilon}} \exp\left\{-\frac{1}{2\epsilon} (x-\mu)^2\right\}
\cdot 
\frac{1}{2b^2} \exp(-\frac{1}{2b^2}\epsilon )
\,d\epsilon \\
&=
\frac{1}{2b^2} 
\int_{0}^{\infty} \frac{1}{\sqrt{2\pi \epsilon}}
\exp\left\{ - \frac{(x-\mu)^2 + \frac{\epsilon^2}{b^2}}{2\epsilon }  \right\} \,d\epsilon 
=
\frac{1}{2b^2} 
\int_{0}^{\infty} \frac{\epsilon}{\sqrt{2\pi \epsilon^3}}
\exp\left\{  
	- \frac{(x-\mu - \frac{\epsilon}{b})^2 + 2\abs{x-\mu}\frac{\epsilon}{b}   }{2\epsilon } 
\right\} \,d\epsilon \\
&\xlongequal{z\triangleq\abs{x-\mu}b}
\frac{1}{2b^2} 
\int_{0}^{\infty} \frac{\epsilon}{\sqrt{2\pi \epsilon^3}}
\exp\left\{  
- \frac{( z -\epsilon )^2  }{2\epsilon b^2 } 
\right\}  
\exp\left\{ \frac{\abs{x-\mu}}{b}\right\}
\,d\epsilon \\
&\xlongequal{\lambda\triangleq\abs{x-\mu}}
\frac{1}{\sqrt{\lambda}}
\frac{1}{2b^2} 
\exp\left\{ \frac{\abs{x-\mu}}{b}\right\}
\int_{0}^{\infty} \epsilon\frac{\sqrt{\lambda}}{\sqrt{2\pi \epsilon^3}}
\exp\left\{  
- \frac{\lambda( z -\epsilon )^2  }{2\epsilon z^2 } 
\right\}  
\,d\epsilon 
= \frac{1}{2b} \exp \left\{ -\frac{\abs{x-\mu}}{b } \right\},
\end{aligned}
$$
where the last equality follows from the mean value of an inverse-Gaussian distribution in Definition~\ref{definition:inverse_gaussian_distribution}.

\paragrapharrow{Conjugate prior for the Laplace scale parameter.}
Just as the inverse-Gamma distribution is conjugate to the variance of a Gaussian likelihood, it is also conjugate to the scale parameter $b$ of a Laplace likelihood.
Suppose we observe a matrix $\bA\in\real^{M\times N}$ 
generated from a Laplace likelihood with known location matrix $\bB\in\real^{M\times N}$
and unknown scale parameter $b>0$:
$$
p(\bA\mid \bB, b) = \prod_{m,n=1}^{M,N} \laplacedist(a_{mn}\mid b_{mn}, b).
$$
Assume an inverse-Gamma prior on $b$: $p(b)=\inversegammadist(b\mid \alpha, \beta)$. 
By Bayes' theorem, the posterior is:
\begin{equation}
\begin{aligned}
p(b \mid &\bA, \bB, \alpha, \beta)
\propto \laplacedist(\bA\mid\bB, b) \inversegammadist(b \mid \alpha, \beta)\\
&=\prod_{m,n=1}^{M,N} \laplacedist(a_{mn}\mid b_{mn}, b) 
 \frac{\beta^\alpha}{\Gamma(\alpha)} (b)^{-\alpha-1} \exp(-\frac{\beta}{b})\\
&\propto \frac{1}{b^{MN}}\exp\left\{   -\frac{1}{b}  \sum_{m,n=1}^{M,N}\abs{a_{mn} - b_{mn}}\right\}
\cdot b^{-\alpha-1}\exp(-\frac{\beta}{b})\\
&=(b)^{-{MN}-\alpha-1} \exp\left\{   -\frac{1}{b}  \left(\sum_{m,n=1}^{M,N}\frac{1}{2} \abs{a_{mn} - b_{mn} } +\beta\right)\right\}
\propto \inversegammadist(b\mid \widetildealpha, \widetildebeta),\\
\end{aligned}
\end{equation}
with updated hyper-parameters:
\begin{equation}\label{equation:inversegamma_conjugate_posterior_laplace}
\widetildealpha={MN}+\alpha,  
\qquad 
\widetildebeta=  \sum_{m,n=1}^{M,N}\frac{1}{2}\abs{a_{mn} - b_{mn} } +\beta.
\end{equation}
Thus, the inverse-Gamma distribution is indeed a conjugate prior for the scale parameter of the Laplace distribution.

The \textit{skew-Laplace} distribution (also called the \textit{asymmetric Laplace} distribution) generalizes the Laplace distribution by allowing different decay rates on either side of the location parameter, thereby introducing skewness (see Figure~\ref{fig:dists_skewlaplace}).

\begin{definition}[Skew-Laplace Distribution\index{Skew-Laplace distribution}]\label{definition:skew_laplace_distribution}
A random variable $\rx$ is said to follow a \textit{skew-Laplace} (or an \textit{asymmetric Laplace}) distribution with location and scale parameters $\mu$ and $\alpha, \beta>0$, respectively, denoted $\rx \sim \skewlaplacedist(\mu, \alpha, \beta)$, if 
$$
f(x; \mu,\alpha, \beta)=\left\{
\begin{aligned}
& \frac{\alpha\beta}{\alpha+\beta} \exp \left\{-\alpha(x-\mu)\right\}   , & \mathrm{\,\,if\,\,} x \geq \mu;  \\
& \frac{\alpha\beta}{\alpha+\beta} \exp \left\{\beta (x-\mu)\right\}   , & \mathrm{\,\,if\,\,} x< \mu.
\end{aligned}
\right.
$$
When $\alpha=\beta=\frac{1}{b}$, the skew-Laplace $\rx \sim \skewlaplacedist(\mu, \alpha, \beta)$ reduces to a Laplace density $\rx \sim \laplacedist(\mu, b)$.
The mean and variance of $\rx \sim \skewlaplacedist(\mu, \alpha, \beta)$ are given by 
\begin{equation}
	\Exp[\rx] = \mu+ \frac{\beta-\alpha}{\alpha\beta}, \qquad \Var[\rx] =\frac{\alpha^2+\beta^2}{\alpha^2\beta^2}. \nonumber
\end{equation}
Figure~\ref{fig:dists_skewlaplace} shows how the distribution changes with different  $\mu, \alpha$, and $\beta$. When $\alpha>\beta$, the right tail decays faster than the left, resulting in \textbf{left skewness} (longer left tail); conversely, $\alpha<\beta$ produces \textbf{right skewness}.
\end{definition}

\section{Multinomial Distribution and Conjugacy}
The multinomial distribution is widely used in  Bayesian mixture or  ordinal models to introduce latent categorical variables. 
It is a discrete probability distribution that describes the probabilities of observing different counts across $K$ possible outcomes in $N$ independent trials,
where each trial results in exactly one of the $K$ categories with fixed probabilities $\bpi =[\pi_1, \pi_2, \ldots, \pi_K]^\top$.
In other words, it models the distribution of counts or frequencies of events among $K$ mutually exclusive categories.

More precisely, the multinomial distribution is parameterized by:
(i) an integer $N\geq 1$ (the total number of trials), and
(ii) a probability mass function   $\bpi =[\pi_1, \pi_2, \ldots, \pi_K]^\top\in [0,1]^K$ such that $\sum_{k=1}^{K} \pi_k=1$.
It answers the following question: If we perform $N$ independent experiments, and each experiment yields outcome $k$ with probability $\pi_k$, what is the probability that outcome $k$ occurs exactly $N_k$
times (for $k=1,2,\ldots,K$), where $\sum_{k=1}^{K}N_k = N$?
Formally, we define the multinomial distribution as follows.

\begin{definition}[Multinomial Distribution\index{Multinomial distribution}]\label{definition:multinomial_dist}
A $K$-dimensional random vector $\bN=[N_1, N_2, \ldots, N_K]^\top\in \{0, 1, 2, \ldots, N\}^K$, satisfying $\sum_{k=1}^{K} N_k=N$, is said to follow a \textit{multinomial distribution} with parameters $N\in \natu$  and $\bpi =[\pi_1, \pi_2, \ldots, \pi_K]^\top\in [0,1]^K$ (with $\sum_{k=1}^{K} \pi_k=1$), denoted by $\bN \sim $ $\multinomial_K(N, \bpi)$. 
Its probability mass function  is:
\begin{equation*}
p\big(N_1, N_2, \ldots , N_K | N, \bpi = (\pi_1, \pi_2, \ldots , \pi_K)\big) = {N\choose N_1\ldots N_K} \prod^K_{k=1}\pi_k^{N_k} \cdot \indicator\left\{\sum_{k=1}^{K}N_k = N\right\},
\end{equation*}
where ${N\choose N_1\ldots N_K}=\frac{N!}{N_1! N_2!  \ldots  N_K!}$ is the \textit{multinomial coefficient}, representing the number of distinct ways to assign $N=\sum_{k=1}^{K}N_k$ trials into $K$ categories with counts $\{N_1,N_2,\ldots N_K\}$.
The mean, variance, and covariance of the multinomial distribution are
$$
\Exp[N_k] = N\pi_k, \qquad \Var[N_k] = N\pi_k(1-\pi_k), \qquad \Cov[N_k, N_m] = -N\pi_k\pi_m\, (k\neq m).
$$
When $K=2$, the multinomial distribution reduces to the \textit{binomial distribution}. 
\end{definition}

\begin{SCfigure}
\centering
\includegraphics[width=0.5\textwidth]{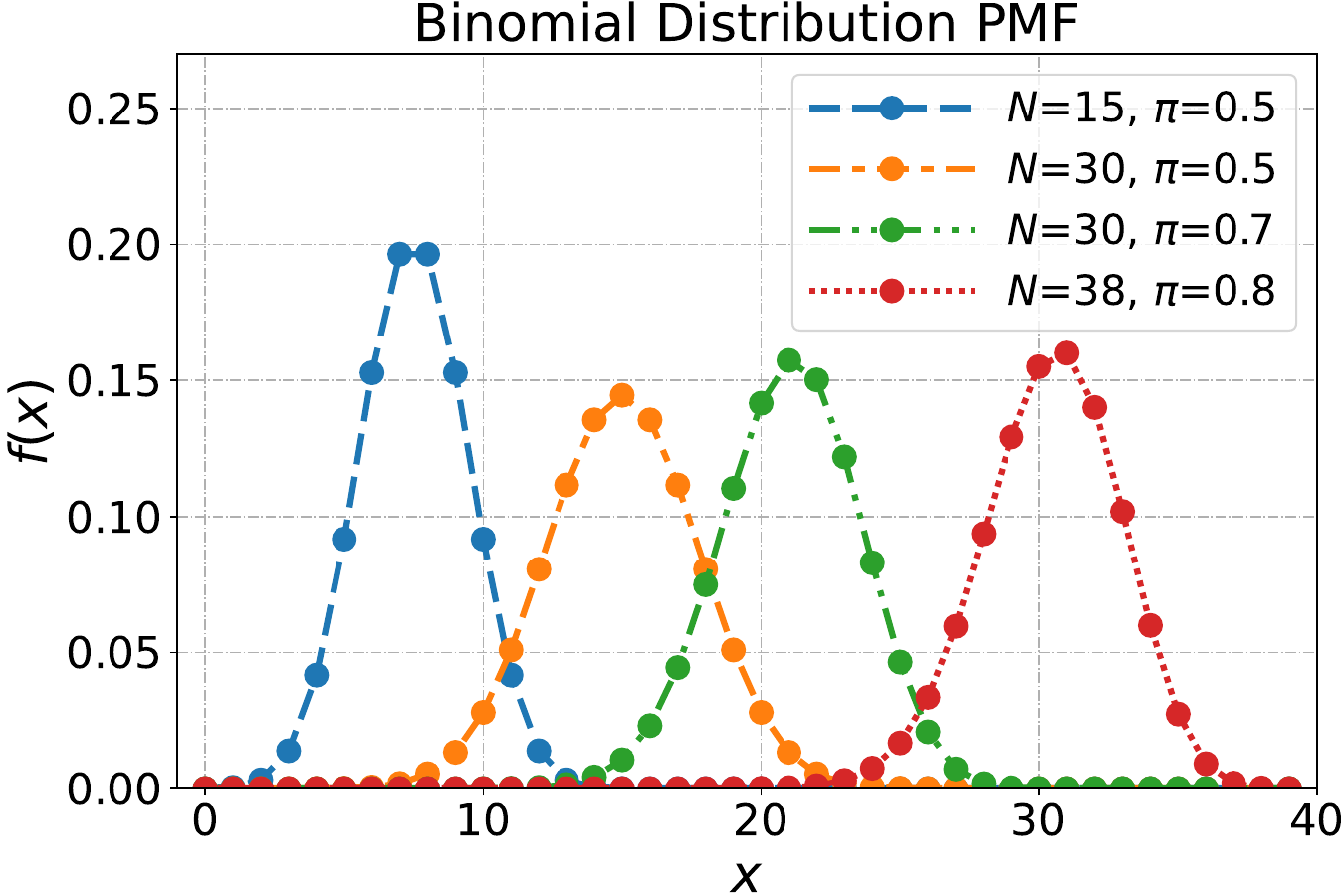}
\caption{Binomial probability mass functions for different values of the parameters $N, \pi$.}
\label{fig:dists_binom}
\end{SCfigure}

\index{Binomial distribution}
\index{One-hot encoding}
\begin{remark}[Binomial, Bernoulli, Multinoulli Distributions]\label{remark:binomial_dist}
In the multinomial distribution, when $K=2$, it is also known as a \textit{binomial} distribution. A random variable $\rx$ is said to follow the binomial distribution with parameter $\pi\in(0,1)$ and $N\in\naturalset$, denoted $\rx\sim\binomialdist(N, \pi)$, if 
$$
\prob(x\mid N, \pi) ={N\choose x} \pi^x (1-\pi)^{N-x},
$$
where ${N\choose x}=\frac{N!}{(N-x)!x!}$  is known as the \textit{binomial coefficient}, representing the number of ways to choose $x$ items from $N$.
The mean and variance of the binomial distribution are 
$$
\textbf{Binomial: }\gap \Exp[\rx] = N\pi, \qquad \Var[\rx] = N\pi(1-\pi).
$$
Figure~\ref{fig:dists_binom} compares different parameters $N$ and $\pi$ for the binomial distribution.

Apparently, the Bernoulli distribution introduced in \S~\ref{section:bayes_appetizers} is a special binomial distribution with $N=1$: 
$$
\rx\sim\bernoulli(\pi)=\binomialdist(N=1, \pi)=\pi\indicator\{\rx=1\} + (1-\pi)\indicator\{\rx=0\}, 
\,\,
\text{where }\pi\in(0,1).
$$
The mean and variance of the Bernoulli distribution are 
$$
\textbf{Bernoulli: }\gap \Exp[\rx] = \pi, \qquad \Var[\rx] = \pi(1-\pi).
$$

Suppose each element $\rx_k$ of a $K$-dimensional random vector $\rvx\in\real^K$ follows from a Bernoulli distribution $\rx_k\sim\bernoulli(\pi_k)$ with $\sum_{k=1}^K{\pi_k}=1$. Then the draws of the random  vector $\bx$ is called a \textit{one-hot encoding}, and the random vector $\rvx$ follows a \textit{Multinoulli distribution} (also known as the \textit{categorical distribution}), denoted $\rvx\sim\multinoulli(K, \bpi)$, if:
$$
p(\rvx \mid \bpi) = \prod_{k=1}^{K} \pi_k^{\indicator(\rx_i=1)}.
$$
The mean, variance, and covariance of each element in the Multinoulli distribution are 
$$
\textbf{Multinoulli: }\gap \Exp[\rx_k] = \pi_k, \gap  \Var[\rx_k] = \pi_k(1-\pi_k), 
\gap\Cov[\rx_k, \rx_m] = -\pi_k\pi_m\, (k\neq m).
$$
\end{remark}

According to the definition of the binomial distribution, it can be applied to experiments where the outcome is either ``success" or ``failure," and the experiment is repeated independently $N$ times. If $\rx\sim \binomialdist(N,\pi)$, then 
$\rx$ represents the number of successes in those $N$  trials, with  $\pi$ being the probability of success on any given trial.

\begin{exercise}[Bernoulli and Binomial]
Suppose $\rx=\sum_{n=1}^{N} \ry_n$ with $\ry_n\stackrel{i.i.d.}{\sim} \bernoullidist(\pi)$. Show that $\rx\sim \binomialdist(N, \pi)$.
\end{exercise}

\subsection{Dirichlet Distribution}\label{section:dirichlet-dist}

The Dirichlet distribution is a multivariate probability distribution defined over the probability simplex.
It takes a vector of positive real numbers as input and outputs a probability distribution over a set of values that are nonnegative and sum to 1.
In Bayesian statistics, the Dirichlet distribution is commonly used as a prior distribution---especially for discrete and categorical data---because it is the conjugate prior for the probability parameter $\bpi$ of the multinomial distribution.
\begin{definition}[Dirichlet Distribution\index{Dirichlet distribution}]\label{definition:dirichlet_dist}
A random vector $\rvx=[\rx_1, \rx_2, \ldots, \rx_K]^\top\in [0,1]^K$ is said to follow a \textit{Dirichlet distribution} with parameter $\balpha$, denoted $\rvx\sim \dirichlet( \balpha)$, if its probability density function is given by
\begin{equation} 
f(\bx ; \balpha) =  \frac{1}{D(\balpha)}  \prod_{k=1}^K x_k ^ {\alpha_k - 1},
\label{equation:dirichlet_distribution2}
\end{equation}  
such that $\sum_{k=1}^K x_k = 1$, $x_k \in$ [0, 1] and 
\begin{equation} 
D(\balpha) = \frac{\prod_{k=1}^K \Gamma(\alpha_k)}{\Gamma(\alpha_+)},
\label{equation:dirichlet_distribution3}
\end{equation}  
where $\balpha = [\alpha_1, \alpha_2, \ldots, \alpha_K]$ is a vector of positive real numbers $\alpha_k>0, \forall k$, and $\alpha_+ = \sum_{k=1}^K \alpha_k$. 
The vector  $\balpha$ is also known as the \textit{concentration parameter} of the Dirichlet distribution. Here, $\Gamma(\cdot)$ denotes the Gamma function, which is a generalization of the factorial function. 
The mean, variance, and covariance are 
$$
\Exp[\rx_k] = \frac{\alpha_k}{\alpha_+}, 
\qquad \Var[\rx_k] = \frac{\alpha_k(\alpha_+-\alpha_k)}{\alpha_+^2(\alpha_++1)}, 
\qquad \Cov[\rx_k, \rx_m]= \frac{-\alpha_k\alpha_m}{\alpha_+^2(\alpha_++1)}.
$$
When $K=2$, the Dirichlet distribution reduces to the \textit{Beta distribution}. 
The Beta distribution $\betadist(\alpha, \beta)$ is defined on $[0,1]$ with the probability density function given by 
$$
\betadist(x\mid \alpha, \beta) = \frac{\Gamma(\alpha+\beta)}{\Gamma(\alpha)\Gamma(\beta)} x^{\alpha-1}(1-x)^{\beta-1}.
$$
That is, if $\rx\sim \betadist(\alpha, \beta)$, then $\rvx=[\rx, 1-\rx] \sim \dirichlet(\balpha) $, where $\balpha=[\alpha, \beta]$.
\end{definition}

Interested readers may refer to Section~\ref{section:drive-dirichlet}  for a derivation of the Dirichlet distribution.
The sample space of the Dirichlet distribution lies on the $(K-1)$-dimensional probability simplex, denoted  $\triangle_K$, which is a subset of $\real^K$ defined as:
$$
\triangle_K \triangleq  \left\{ \bpi\mid  0 \leq \pi_k\leq 1, \,\, \sum_{k=1}^{K}\pi_k=1 \right\}.
$$
Although embedded in $\real^K$, this simplex is intrinsically $(K-1)$-dimensional because the constraint $\sum_{k=1}^{K}\pi_k=1$  reduces the degrees of freedom by one.

\begin{figure}[htp]
\centering
\subfigure[ $ \balpha=\begin{bmatrix}
10,10,10
\end{bmatrix} $, z-axis is pdf. ]{\includegraphics[width=0.4
\textwidth]{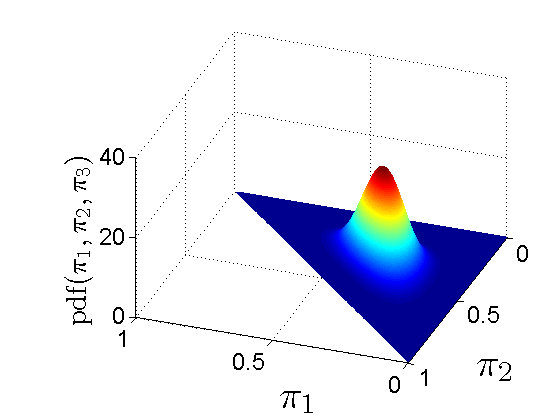} \label{fig:dirichlet_pdf}}
~
\subfigure[$\balpha=\begin{bmatrix}
10,10,10
\end{bmatrix}$, z-axis is $\pi_3$.]{\includegraphics[width=0.4
\textwidth]{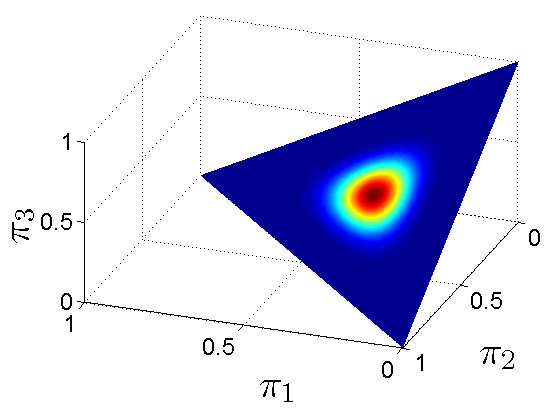} \label{fig:dirichlet_surface}}
\centering
\subfigure[ $ \balpha=\begin{bmatrix}
1,1,1
\end{bmatrix} $ ]{\includegraphics[width=0.4
\textwidth]{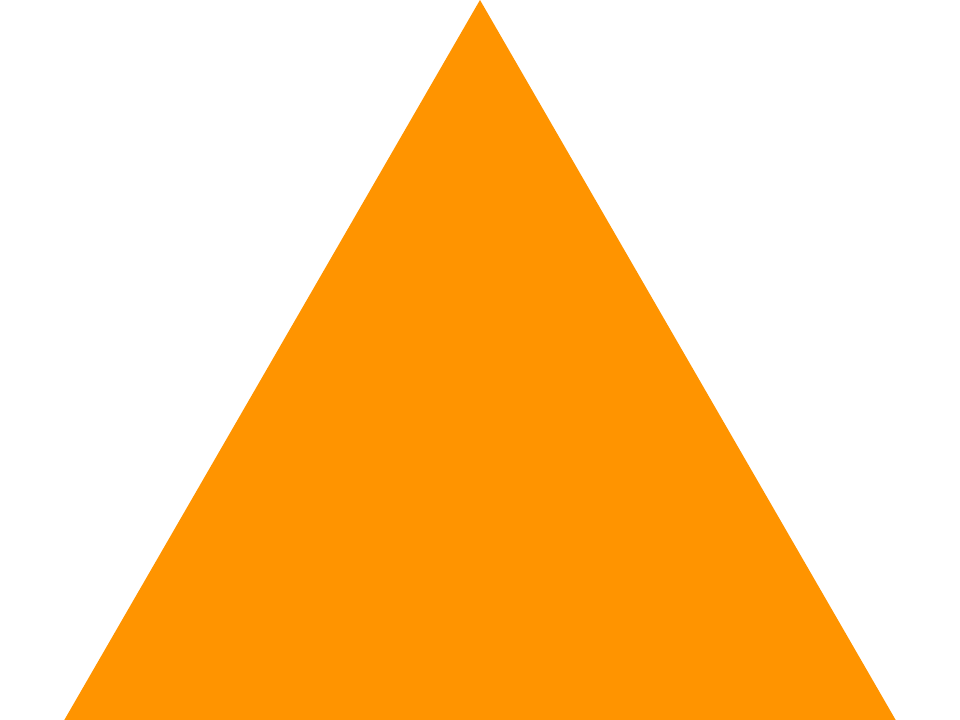} \label{fig:dirichlet_sample_111}}
~
\subfigure[$\balpha=\begin{bmatrix}
0.9,0.9,0.9
\end{bmatrix}$]{\includegraphics[width=0.4
\textwidth]{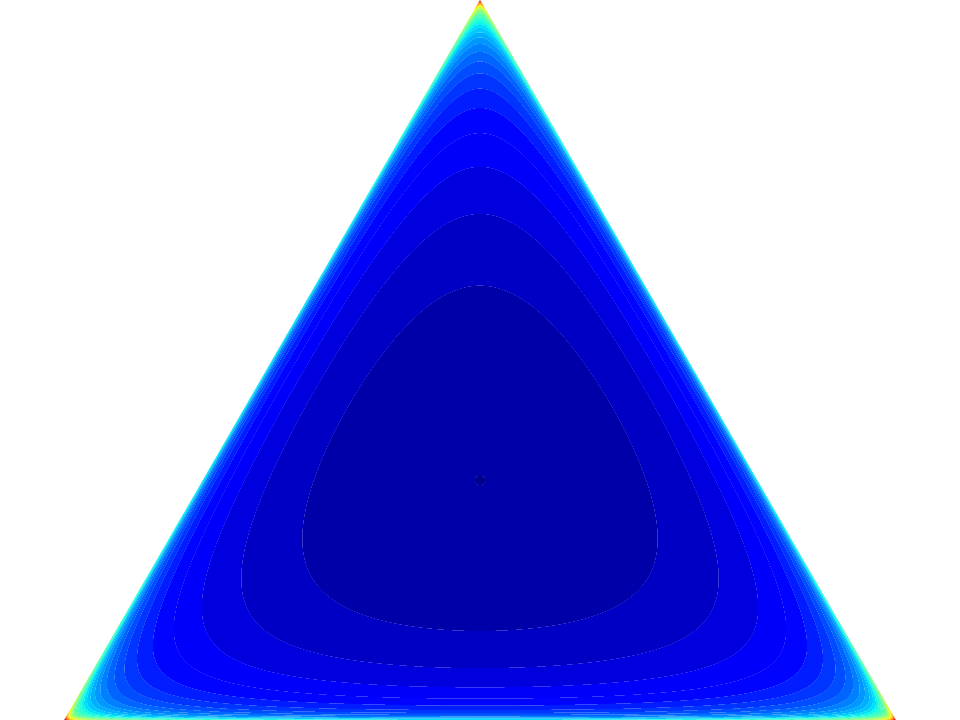} \label{fig:dirichlet_sample_090909}}
\centering
\subfigure[$\balpha=\begin{bmatrix}
10,10,10
\end{bmatrix}$]{\includegraphics[width=0.4
\textwidth]{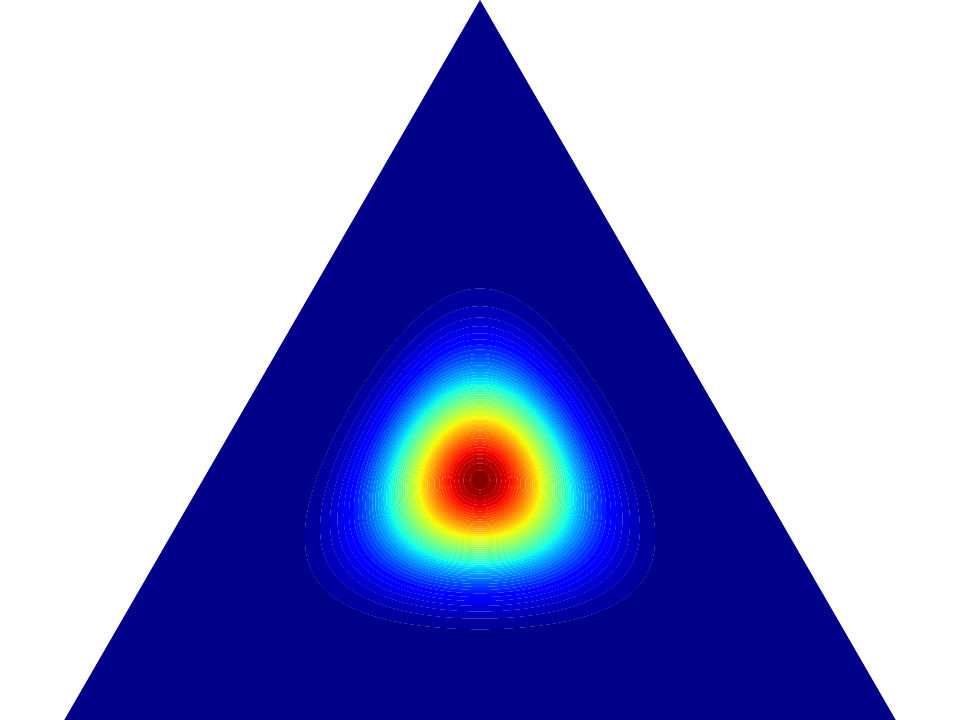} \label{fig:dirichlet_sample_101010}}
~
\subfigure[$\balpha=\begin{bmatrix}
15,5,2
\end{bmatrix}$]{\includegraphics[width=0.4
\textwidth]{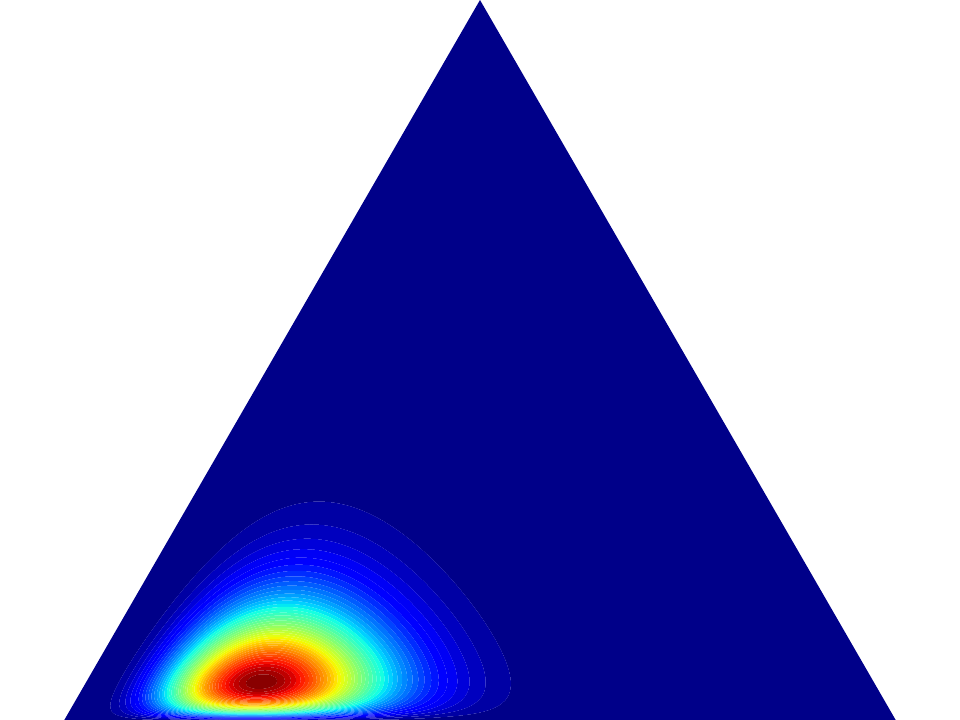} \label{fig:dirichlet_sample_1552}}
\centering
\caption{Density plots (\textcolor{mylightbluetext}{blue}=low, \textcolor{winestain}{red}=high) for the Dirichlet distribution over the probability simplex in $\real^3$ for various values of the concentration parameter $\balpha$. 
When $\balpha=[c, c, c]$, the distribution is called a \textit{symmetric Dirichlet distribution}, and the density is symmetric about the uniform probability mass function (i.e., occurs in the middle of the simplex). When $0<c<1$, there are sharp peaks of density almost at the vertices of the simplex. When $c>1$, the density becomes unimodal and concentrated in the center of the simplex. And when $c=1$, it is uniform distributed over the simplex. Finally, if $\balpha$ is not a constant vector, the density is not symmetric.}\centering
\label{fig:dirichlet_samples}
\end{figure}

\begin{figure}[h]
\center
\subfigure[ $ \balpha=\begin{bmatrix}1,1,1\end{bmatrix}$ ]{\includegraphics[width=0.4
\textwidth]{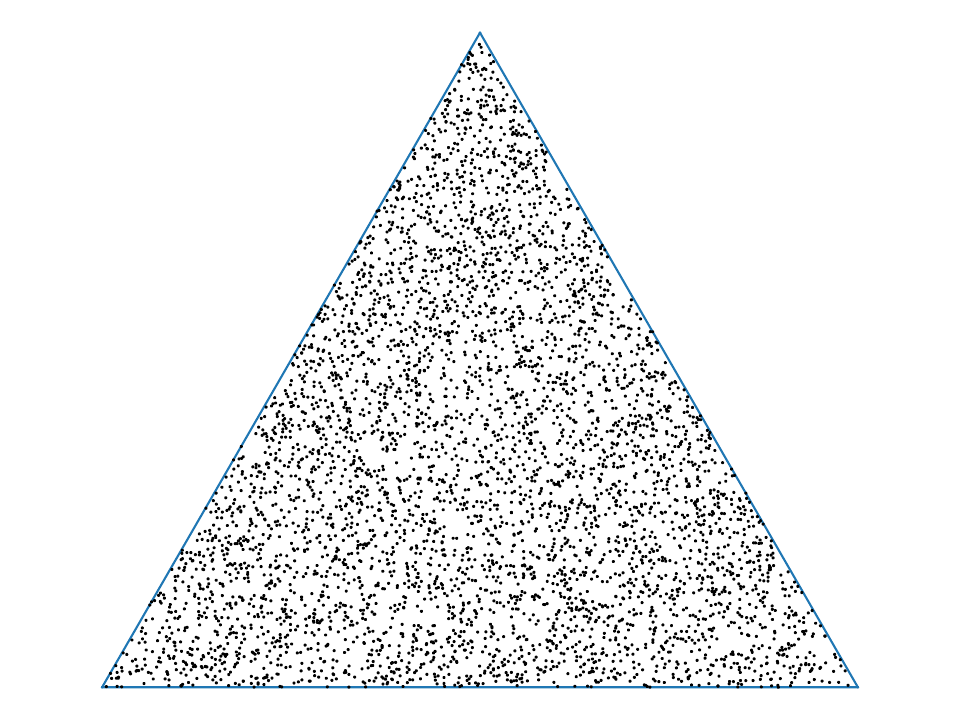} \label{fig:dirichlet_points_111}}
~
\subfigure[$\balpha=\begin{bmatrix}0.9,0.9,0.9\end{bmatrix}$]{\includegraphics[width=0.4
\textwidth]{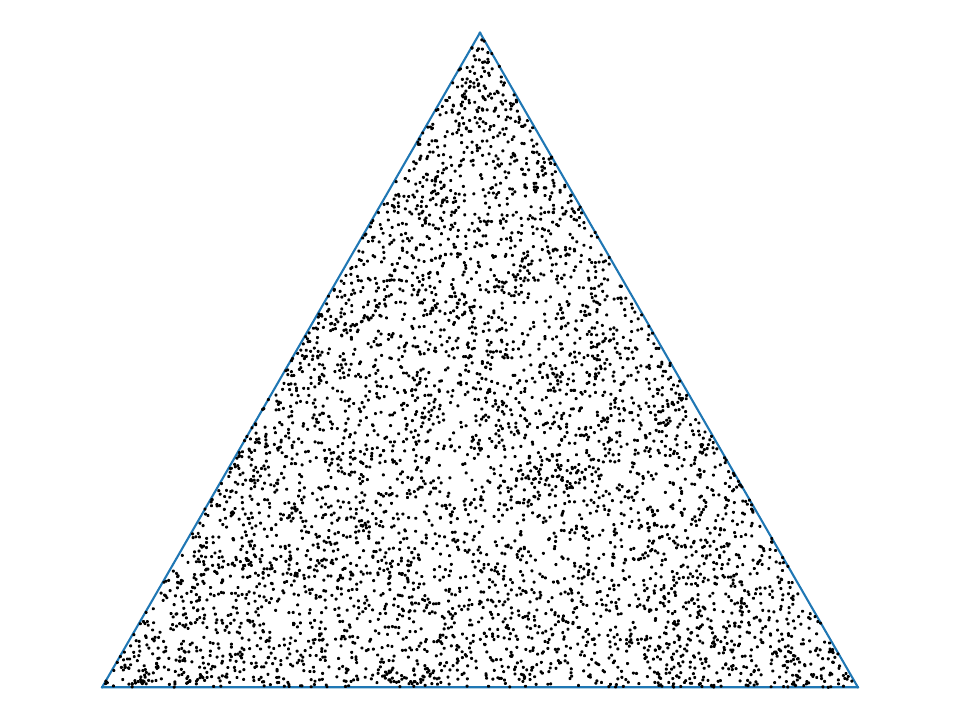} \label{fig:dirichlet_points_090909}}
\center
\subfigure[$\balpha=\begin{bmatrix}10,10,10\end{bmatrix}$]{\includegraphics[width=0.4
\textwidth]{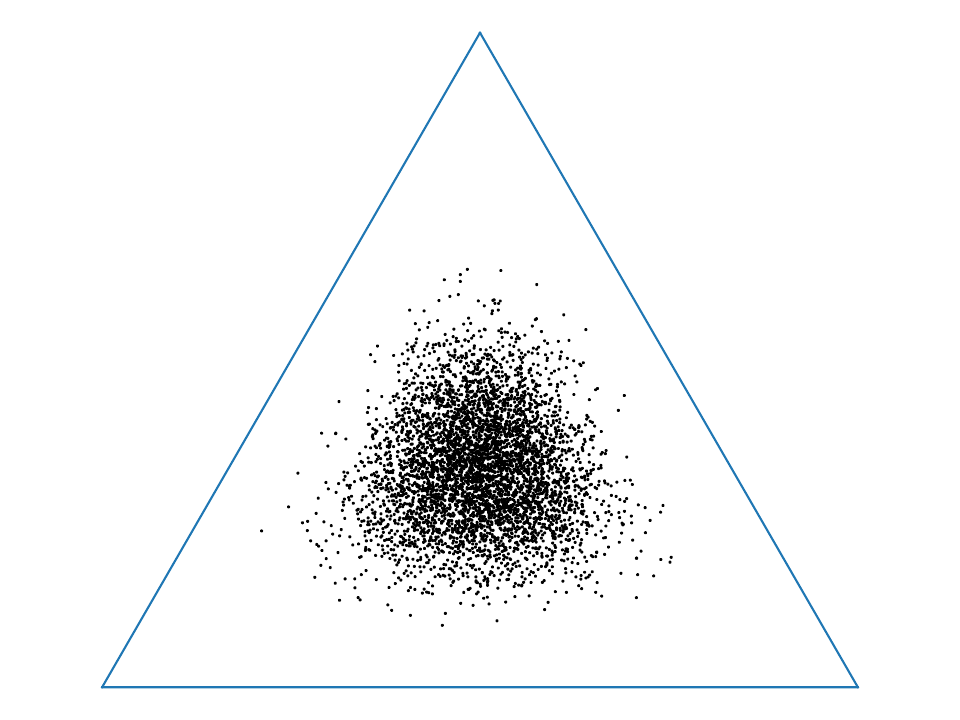} \label{fig:dirichlet_points_101010}}
~
\subfigure[$\balpha=\begin{bmatrix}15,5,2\end{bmatrix}$]{\includegraphics[width=0.4
\textwidth]{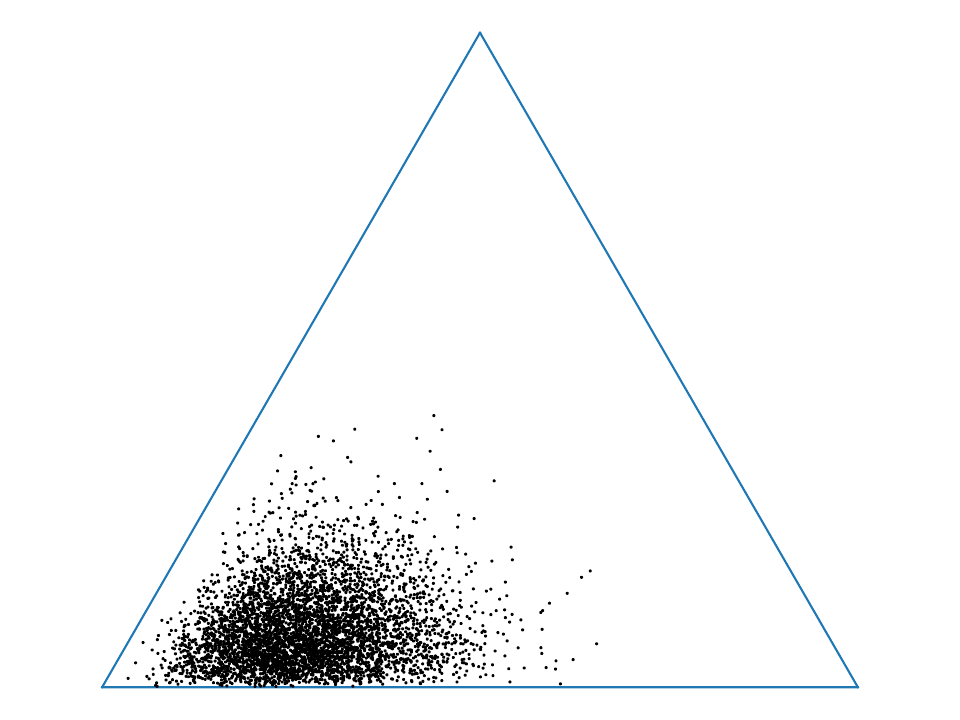} \label{fig:dirichlet_points_1552}}
\center
\caption{Draw of 5,000 points from Dirichlet distribution over the probability simplex in $\real^3$ for various values of the concentration parameter $\balpha$.}
\label{fig:dirichlet_points}
\end{figure}

Figure~\ref{fig:dirichlet_samples} displays density plots of the Dirichlet distribution over the two-dimensional simplex in $\real^3$ for several choices of $\balpha$,  
while Figure~\ref{fig:dirichlet_points} shows 5,000 random samples drawn under each setting.
Strictly speaking, the full density of a Dirichlet distribution in $\real^3$ lives in a 4-dimensional space (three coordinates plus density).
Figure~\ref{fig:dirichlet_pdf} is a projection of a surface into 3D by using the probability density as the z-axis, and Figure~\ref{fig:dirichlet_surface} uses $\pi_3$ as the z-axis instead. Figure~\ref{fig:dirichlet_sample_111} through \ref{fig:dirichlet_sample_1552} show 2D contour projections onto the simplex.

When the concentration parameter is $\balpha=[1,1,1]$, the Dirichlet distribution reduces to the uniform distribution over the simplex. This can be easily verified that $\dirichlet(\bx\mid  \balpha=[1,1,1]) =  \frac{\Gamma(3)}{(\Gamma(1))^3}= 2 $, which is a constant that does not depend on the specific value of $\bx$. When $\balpha=[c,c,c]$ with $c>1$, the density is unimodal and peaked at the center. 
This follows from the form $\dirichlet(\bx \mid \balpha=[c,c,c]) =  \frac{\Gamma(3c)}{(\Gamma(c))^3}\prod_{k=1}^3 x_k ^ {c - 1} $ such that a small value of $x_k$ will drive the density toward zero. 
On the contrary, when $\balpha=[c,c,c]$ with $c<1$, the density concentrates near the corners (vertices) of the simplex, producing sharp peaks.

Additional properties of the Dirichlet distribution are summarized in Table~\ref{table:dirichlet-property}, with proofs provided in Section~\ref{section:drive-dirichlet}. That same derivation also yields a practical method for sampling from the Dirichlet distribution: draw independent samples from Gamma distributions and normalize them.

\begin{table}[]
\begin{tabular}{l|l}
\hline
\begin{tabular}[c]{@{}l@{}}Marginal \\ Distribution\end{tabular}    & $\rx_i \sim \betadist(\alpha_i, \alpha_+-\alpha_i)$.    \\ \hline
\begin{tabular}[c]{@{}l@{}}Conditional \\ Distribution\end{tabular} & 
\begin{tabular}[c]{@{}l@{}}$\rvx_{-i} \mid \rx_i \sim (1-\rx_i)\dirichlet(\alpha_{-i})$, \\ 
where $\rvx_{-i}$ is a random vector excluding $\rx_i$. \end{tabular}                                                                                                                                                                                                                                                                                                                                                                                                                                                                      \\ \hline
\begin{tabular}[c]{@{}l@{}}Aggregation\\ Property\end{tabular}      & \begin{tabular}[c]{@{}l@{}}If $M=\rx_i+\rx_j$, then $[\rx_1, \ldots \rx_{i-1}, \rx_{i+1}, \ldots, \rx_{j-1}, \rx_{j+1}, \ldots, \rx_K, M] \sim $\\ 
\gap\gap$\dirichlet([\alpha_1, \ldots, \alpha_{i-1}, \alpha_{i+1}, \ldots, \alpha_{j-1}, \alpha_{j+1}, \ldots, \alpha_K, \alpha_i+\alpha_j])$.\\ 
In general, if $\{\sA_1, \sA_2, \ldots, \sA_r\}$ is a partition of $\{1, 2, \ldots, K\}$, then \\ 
$\left[\sum_{i\in \sA_1} \rx_i, \sum_{i\in \sA_2} \rx_i, \ldots, \sum_{i\in \sA_r} \rx_i\right] \sim$\\ \gap\gap$\dirichlet\left(\left[\sum_{i\in \sA_1} \alpha_i, \sum_{i\in \sA_2} \alpha_i, \ldots, \sum_{i\in \sA_r} \alpha_i\right]\right)$.\end{tabular} \\ \hline
\end{tabular}
\caption{Properties of the Dirichlet distribution.}
\label{table:dirichlet-property}
\end{table}

\subsection{Posterior Distribution for Multinomial Distribution}\label{section:dirichlet-dist-post}
Due to conjugacy, the Dirichlet distribution serves as a conjugate prior for the multinomial likelihood.
If $(\bN \mid \bpi) \sim $ $\multinomial_K(N, \bpi)$, and $\bpi$ is given a Dirichlet prior with $\bpi \sim$ $\dirichlet(\balpha)$,  then the posterior distribution is 
\begin{equation}
(\bpi \mid\bN) \sim \dirichlet(\balpha+\bN) = \dirichlet(\alpha_1+N_1, \ldots, \alpha_K+N_K).
\end{equation}
\begin{proof}[of conjugate prior of multinomial distribution]
Using Bayes' theorem ``$\mathrm{posterior} \propto \mathrm{likelihood} \times \mathrm{prior} $," we obtain the posterior density
$$
\begin{aligned}
\mathrm{posterior}
&= p(\bpi\mid\balpha, \bN) 
\propto \multinomial_K(\bN\mid N, \bpi) \cdot \dirichlet(\bpi\mid \balpha) \\
&= \left(\frac{N!}{N_1! N_2!  \ldots  N_K!} \prod^K_{k=1}\pi_k^{N_k}\right) \cdot  \left(\frac{1}{D(\balpha)}  \prod_{k=1}^K \pi_k ^ {\alpha_k - 1}\right)\\
&\propto   \prod_{k=1}^K \pi_k ^ {\alpha_k +N_k - 1} \propto \dirichlet(\bpi\mid \balpha+\bN).
\end{aligned}
$$
Therefore, it follows that $(\bpi \mid \bN) \sim$ $\dirichlet(\balpha+\bN)$ = $\dirichlet(\alpha_1+N_1, \ldots, \alpha_K+N_K)$.
\end{proof} 

Comparing the prior and posterior reveals that the relative magnitudes of the $\alpha_k$'s determine the prior mean of  $\bpi$, while the total $\alpha_+=\sum_k\alpha_k$ 
reflects the strength (or ``confidence") of the prior.
In fact, the Dirichlet prior $\dirichlet(\balpha)$  is mathematically equivalent to having observed $\alpha_k-1$  pseudo-counts in category $k$, for a total of $\sum_{k=1}^{K}(\alpha_k-1)$ virtual observations.

Since the Dirichlet distribution generalizes the Beta distribution to more than two categories, the Beta distribution naturally arises as the conjugate prior for the binomial likelihood \citep{hoff2009first, frigyik2010introduction}.

\section{Poisson and Multinomial}
The \textit{Poisson} distribution is a discrete probability distribution that characterizes the number of events in a fixed interval of time or space, given the average number of events in that interval.
The Poisson distribution is frequently employed for modeling  count data, such as the number of calls received by a call center in an hour or the number of emails received in a day provided that the probability of a ``success" for any given instance is ``very small." 
Other typical applications include: the number of stars in a random area of the space; the distribution of bacteria on a surface; the number of typographical errors on a typed page; the number of wrong connections to a phone number. 
These scenarios share the key feature that events are rare, independent, and occur at a roughly constant average rate.

\begin{SCfigure}
	\centering
	\includegraphics[width=0.5\textwidth]{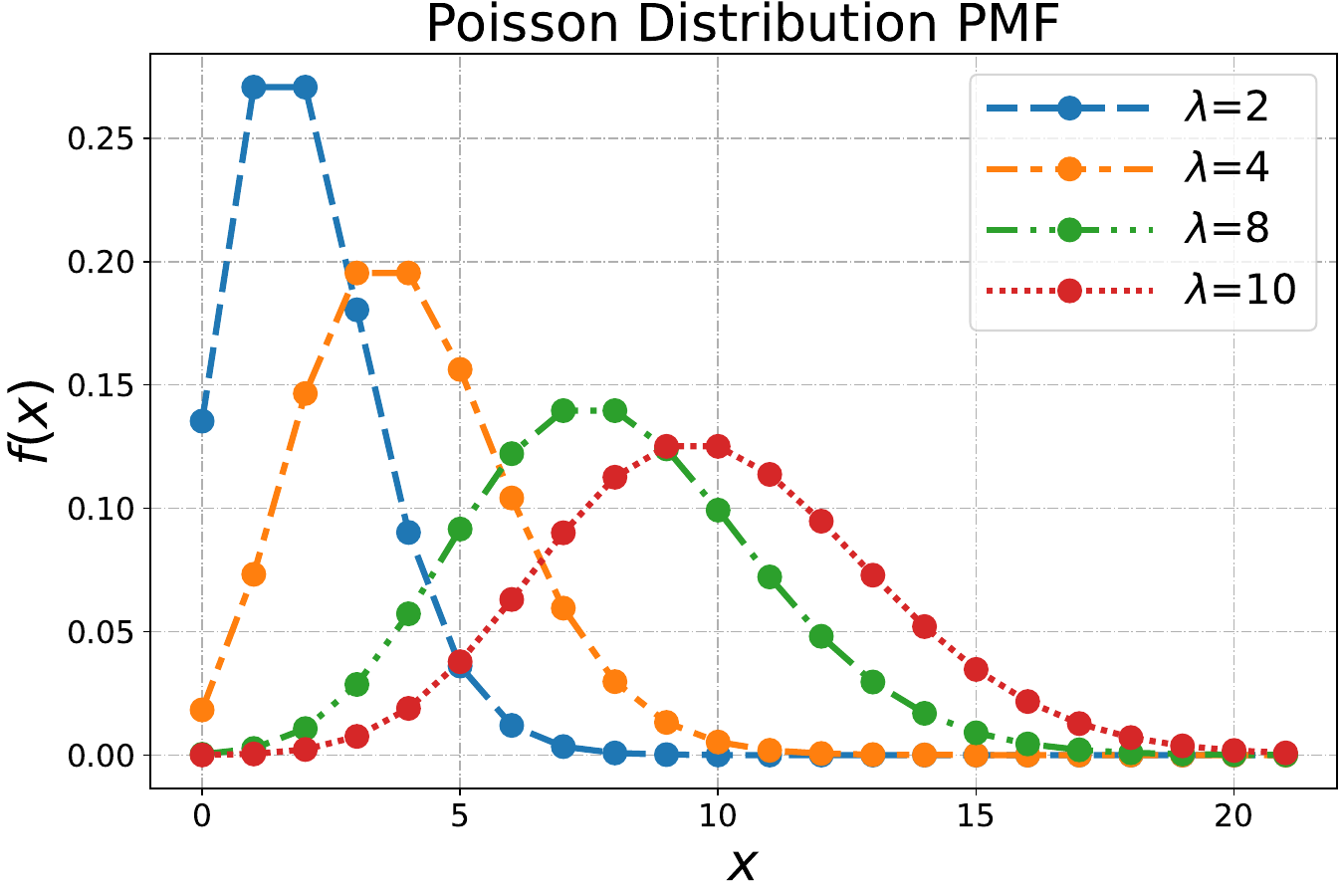}
	\caption{Poisson probability mass functions for different values of the parameter $\lambda$.}
	\label{fig:dists_poisson}
\end{SCfigure}

\begin{definition}[Poisson Distribution\index{Poisson distribution}]\label{definition:poisson_distribution}
A random variable $\rx\in\{0,1,2,3,\ldots\}$ is said to follow a \textit{Poisson distribution} with rate parameter $\lambda>0$,  denoted $\rx \sim \poissondist(\lambda)$, if its probability mass function is
$$
f(x; \lambda)= \frac{\lambda^x}{x!}  \exp(-\lambda).
$$
The mean and variance of $\rx \sim \poissondist( \lambda)$ are given by 
\begin{equation}
\Exp[\rx] = \lambda, \qquad \Var[\rx] =\lambda. \nonumber
\end{equation}
The support of the Poisson distribution is the set of nonnegative integers, $\{0,1,2,3,\ldots\} = \{0\}\cup \naturalset$.
Figure~\ref{fig:dists_poisson} illustrates the probability mass functions of the Poisson distribution for several values of $\lambda$.
\end{definition}

An important property of the Poisson distribution is that its mean equals its variance.
Intuitively, the Poisson distribution arises as the limiting case of the binomial distribution when the number of trials $N\rightarrow \infty$ and the success probability $\pi=\lambda/N \rightarrow 0$,  such that the expected number of successes $\lambda=N\pi$ remains constant. 
This limiting behavior is often referred to as the \textit{law of rare events}, which justifies using the Poisson distribution to model rare phenomena---such as radioactive decays or meteor strikes---where events are infrequent but occur over a large number of opportunities.

Another key property is that the sum of independent Poisson-distributed random variables is itself Poisson-distributed, with a rate equal to the sum of the individual rates.
\begin{theoremHigh}[Sum of Independently Distributed Poisson]\label{theorem:sum_iid_poisson}
Let $\rx_n\sim \poissondist(\lambda_n)$ for $n=1,2,\ldots,N$. Then $ \ry=\sum_{n=1}^{N} \rx_n\sim \poissondist(\sum_{n=1}^{N}\lambda_n)$.
\end{theoremHigh}
To illustrate the idea, consider two independent Poisson random variables: $\rx\sim \poissondist(\lambda_1)$ and $\ry\sim\poissondist(\lambda_2)$.
Define $\lambda=\lambda_1+\lambda_2$ and $\rz=\rx+\ry$. Then $\rz$ is a Poisson random variable with parameter $\lambda$. To see this, we have 
$$
\begin{aligned}
p(z) 
&= \prob(\rz=z) = \sum_{k=1}^{z} \prob(\rx=k) \cdot  \prob(\ry=z-k)
= \sum_{k=1}^{z} \frac{\lambda_1^k}{k!} \exp(-\lambda_1) \cdot \frac{\lambda_2^{z-k}}{(z-k)!} \exp(-\lambda_2)\\
&= \frac{\exp(-\lambda_1-\lambda_2)}{z!} \sum_{k=1}^{z} {z\choose k} \lambda_1^k\lambda_2^{z-k}
\stackrel{*}{=}\frac{\exp(-\lambda)}{z!}(\lambda_1+\lambda_2)^z = \frac{\lambda^z}{z!} \exp(-\lambda),
\end{aligned}
$$
where the equality ($*$) follows from the \textit{binomial theorem}.
By induction, this result extends to any finite number of independent Poisson variables.

Finally, the Poisson distribution has a deep connection with the multinomial distribution:
\begin{theoremHigh}[Poisson and Multinomial]\label{theorem:multinomial_poisson}
Let $\rx_k\sim \poissondist(\lambda_k)$ be independent Poisson variables for $k\in\{1,2,\ldots, K\}$. Then the conditional distribution of $\rvx=[\rx_1,\rx_2,\ldots, \rx_K]^\top$ given $\sum_{k=1}^{K}\rx_k=N$ is $\multinomial_K(N,\{ p_1, p_2, \ldots,p_K\})$ with
$$
p_k= \frac{\lambda_k}{\lambda_1+\lambda_2+\ldots+\lambda_K}, \gap \text{for all }k\in\{1,2,\ldots,K\}.
$$
\end{theoremHigh}

\section{Multivariate Gaussian Distribution and Conjugacy}\label{sec:multi_gaussian_conjugate_prior}
We have previously shown the conjugate prior for the mean parameter of a univariate Gaussian distribution when the variance (or precision) is fixed, as well as the joint conjugate prior for both the mean and variance (or precision) parameters of a univariate Gaussian distribution.
In this section, we extend this analysis to the multivariate Gaussian distribution. For further discussion, see \citet{murphy2007conjugate, murphy2012machine, teh2007exponential, kamper2013gibbs, das2014dpgmm}.

\subsection{Multivariate Gaussian Distribution}\label{section:multi_gaussian_dist}

A \textit{multivariate Gaussian distribution} (also known as a \textit{multivariate normal distribution}) is a continuous probability distribution that describes jointly normal random variables across multiple dimensions. It is fully characterized by its mean vector---whose dimension equals the number of variables---and its covariance matrix, a symmetric positive-definite square matrix of the same dimension. The covariance matrix captures the pairwise covariances between all variables, thereby encoding their linear relationships.
Widely used in machine learning, statistics, and signal processing, the multivariate Gaussian (often simply called ``Gaussian" when the context is clear) is a fundamental tool for modeling complex, high-dimensional data distributions.
We begin by providing a formal definition of the multivariate Gaussian distribution.

\begin{figure}[h]
\subfigure[Gaussian, $\bSigma =\begin{bmatrix}
1&0\\
0&1
\end{bmatrix}. $ ]{\includegraphics[width=0.31
\textwidth]{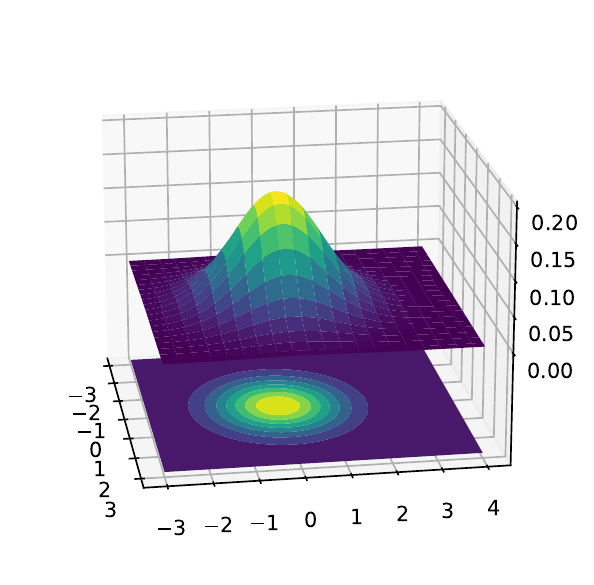} \label{fig:dists_multiGauss_sigma1}}
\subfigure[Gaussian, $\bSigma =\begin{bmatrix}
1&0\\
0&3
\end{bmatrix}.$]{\includegraphics[width=0.31
\textwidth]{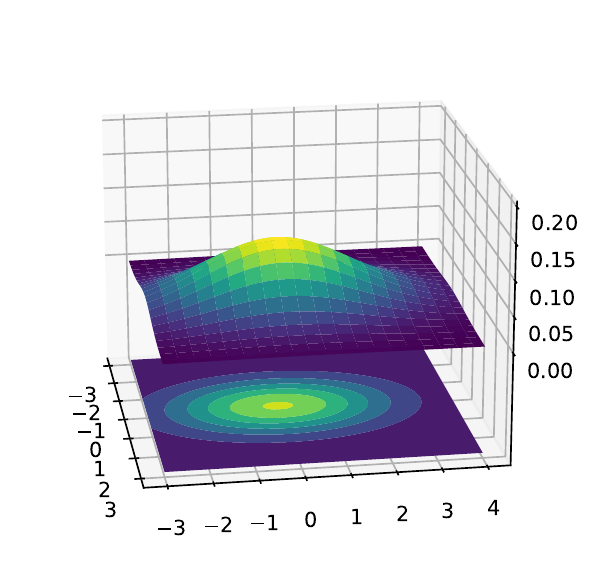} \label{fig:dists_multiGauss_sigma2}}
\subfigure[Gaussian, $\bSigma =\begin{bmatrix}
1&\textendash0.5\\
\textendash0.5&1.5
\end{bmatrix}.$]{\includegraphics[width=0.31 
\textwidth]{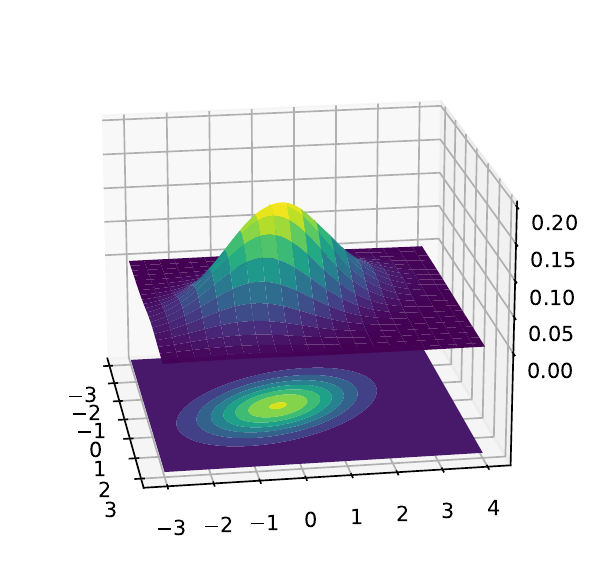} \label{fig:dists_multiGauss_sigma3}}
\subfigure[Gaussian, $\bSigma =\begin{bmatrix}
	2&0\\
	0&2
\end{bmatrix}. $ ]{\includegraphics[width=0.31
\textwidth]{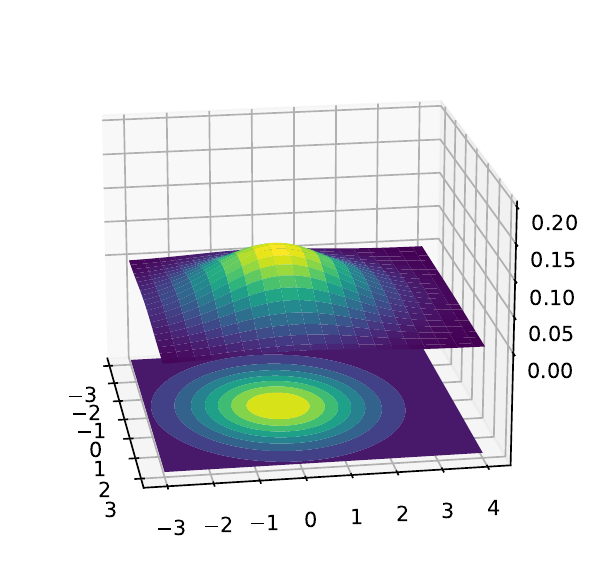} \label{fig:dists_multiGauss_sigma4}}
\subfigure[Gaussian, $\bSigma =\begin{bmatrix}
	3&0\\
	0&1
\end{bmatrix}.$]{\includegraphics[width=0.31
\textwidth]{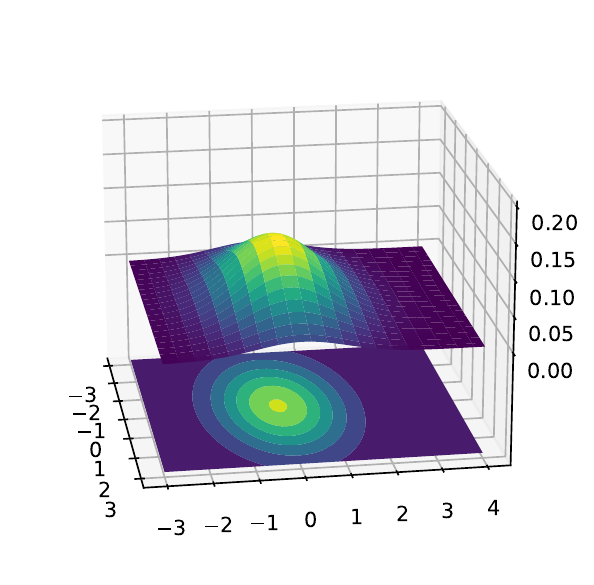} \label{fig:dists_multiGauss_sigma5}}
\subfigure[Gaussian, $\bSigma =\begin{bmatrix}
	3&\textendash0.5\\
	\textendash0.5&1.5
\end{bmatrix}.$]{\includegraphics[width=0.31
\textwidth]{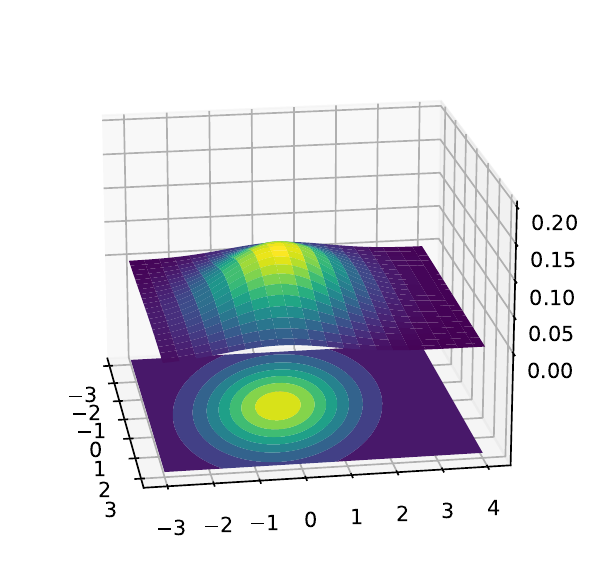} \label{fig:dists_multiGauss_sigma6}}
\centering
\caption{Density and contour plots (\textcolor{mylightbluetext}{blue}=low, \textcolor{mydarkyellow}{yellow}=high) for the multivariate Gaussian distribution over the $\mathbb{R}^2$ space for various values of the covariance/scale matrix with zero-mean vector.  Fig~\ref{fig:dists_multiGauss_sigma1} and \ref{fig:dists_multiGauss_sigma4}: A spherical covariance matrix has a circular shape; 
Fig~\ref{fig:dists_multiGauss_sigma2} and \ref{fig:dists_multiGauss_sigma5}: A diagonal covariance matrix is an \textbf{axis aligned} ellipse; 
Fig~\ref{fig:dists_multiGauss_sigma3} and \ref{fig:dists_multiGauss_sigma6}: A full covariance matrix has an  elliptical shape.}
\centering
\label{fig:multi_gaussian_density}
\end{figure}

\begin{definition}[Multivariate Gaussian Distribution\index{Multivariate Gaussian distribution}]\label{definition:multivariate_gaussian}
A random vector $\rvx \in \real^D$ is said to follow a \textit{multivariate Gaussian distribution} (multivariate normal, MVN) with parameters $\bmu\in\real^D$ and $\bSigma\in\real^{D\times D}$, denoted  $\rvx\sim \normal(\bmu, \bSigma)$, if its probability density function is given by
$$
\begin{aligned}
f(\bx; \bmu, \bSigma)&= (2\pi)^{-D/2} \abs{\bSigma}^{-1/2}\exp\left\{-\frac{1}{2}(\bx - \bmu)^\top \bSigma^{-1}(\bx - \bmu)\right\},
\end{aligned}
$$
where $\bmu \in \real^D$ is  the \textit{mean vector}, and $\bSigma\in \real^{D\times D}$ is positive definite  \textit{covariance matrix}.
The mean, mode, and covariance of the multivariate Gaussian distribution are given by 
\begin{equation*}
\begin{aligned}
\Exp [\rvx] &= \bmu, \qquad 
\mathrm{Mode}[\rvx] &= \bmu, \qquad
\Cov [\rvx] &= \bSigma. 
\end{aligned}
\end{equation*}
The covariance matrix can also be expressed as
$$
\Cov[\rvx]=\Exp[(\rvx-\bmu)(\rvx-\bmu)^\top]=\Exp[\rvx\rvx^\top]-\bmu\bmu^\top.
$$
Figure~\ref{fig:multi_gaussian_density} illustrates Gaussian density plots under different covariance structures. Additionally, a multivariate Gaussian random vector can be generated from univariate Gaussian samples; see Problem~\ref{problem:multiGauss}.
\end{definition}

Similar to the univariate case (see Equation~\eqref{equation:uni_gaussian_likelihood}), the likelihood of $N$ independent observations  $\mathcalX = \{\bx_1, \bx_2, \ldots , \bx_N \}$ drawn from a multivariate Gaussian distribution with mean $\bmu$ and covariance $\bSigma$ is given by 
\begin{equation}\label{equation:multi_gaussian_likelihood}
\begin{aligned}
&\gap p(\mathcalX \mid \bmu, \bSigma) =\prod^N_{n=1} \normal (\bxn\mid \bmu, \bSigma) \\
&\overset{(a)}{=} (2\pi)^{-ND/2} \abs{\bSigma}^{-N/2}\exp\left\{-\frac{1}{2} \sum^N_{n=1}(\bxn - \bmu)^\top \bSigma^{-1}(\bxn - \bmu)\right\} \\
&\overset{(b)}{=} (2\pi)^{-ND/2} \abs{\bSigma}^{-N/2}\exp\left\{-\frac{1}{2} \tr( \bSigma^{-1}\bS_{\bmu} )  \right\}\\
&\overset{(c)}{=} (2\pi)^{-ND/2} \abs{\bSigma}^{-N/2}\exp\left\{-\frac{N}{2}(\bmu - \widebarbx)^\top \bSigma^{-1}(\bmu - \widebarbx)\right\}  \exp\left\{-\frac{1}{2}\tr( \bSigma^{-1}\bS_{\widebarx} )\right\},
\end{aligned}
\end{equation}
where 
\begin{equation}\label{equation:mvu-sample-covariance}
\begin{aligned}
\bS_{\bmu} &\triangleq \sum^N_{n=1}(\bxn - \bmu)(\bxn - \bmu)^\top;
\quad
\bS_{\widebarx} \triangleq \sum^N_{n=1}(\bxn - \widebarbx)(\bxn - \widebarbx)^\top;
\quad
\widebarbx \triangleq\frac{1}{N}\sum^N_{n=1}\bxn.
\end{aligned}
\end{equation}
Here, $\bS_{\widebarx}$ is known as the \textit{scatter matrix} or the \textit{sum-of-squares matrix}.
The equivalence between equation (a) and equation (c)  follows from the following identity (similar reasoning applies to the equivalence between equation (a) and equation (b)):
\begin{align}
\sum^N_{n=1}(\bxn - \bmu)^\top\bSigma^{-1}(\bxn - \bmu) = \tr(\bSigma^{-1}\bS_{\widebarx}) + N \cdot (\widebarbx - \bmu)^\top\bSigma^{-1}(\widebarbx - \bmu),
\label{equation:multi_gaussian_identity}
\end{align}
where the trace of a square matrix $\bA$, denoted $\trace(\bA)$,  is the sum of its diagonal elements:  $\tr(\bA) = \sum_i a_{ii}$.
This decomposition is useful: form (b) facilitates deriving the conjugate prior for $\bSigma$ alone, while form (c) supports the joint conjugate analysis of $(\bmu, \bSigma)$, as discussed in Section~\ref{sec:niw_posterior_conjugacy}.
\begin{proof}[of Identity~\ref{equation:multi_gaussian_identity}]
A useful trick involves the cyclic invariance of the trace operator. For any vector $\bx$ and matrix $\bA$, we have
\begin{equation}
\bx^\top \bA \bx = \tr(\bx^\top \bA \bx) = \tr(\bx \bx^\top \bA) =  \tr(\bA \bx \bx^\top ),
\end{equation}
where the first equality follows from the fact that $\bx^\top \bA \bx$ is a scalar and the trace of a product is invariant under cyclical permutations of the factors \footnote{Trace is invariant under cyclical permutations: $\tr(\bA\bB\bC) = \tr(\bB\bC\bA) = \tr(\bC\bA\bB)$ whenever the products are defined.}. 

We can then rewrite $\sum^N_{n=1}(\bxn - \bmu)^\top\bSigma^{-1}(\bxn - \bmu)$ as 
\begin{equation}
\begin{aligned}
&\, \sum^N_{n=1}(\bxn - \widebarbx)^\top\bSigma^{-1}(\bxn - \widebarbx) + \sum^N_{n=1}(\widebarbx - \bmu)^\top\bSigma^{-1}(\widebarbx - \bmu) \\
&= \tr(\bSigma^{-1} \bS_{\widebarx}  ) + N \cdot (\widebarbx - \bmu)^\top\bSigma^{-1}(\widebarbx - \bmu).
\end{aligned}
\end{equation}
This concludes the proof.
\end{proof}

Although this identity \eqref{equation:multi_gaussian_identity} does not reduce computational complexity, it plays a crucial role in establishing conjugacy, as detailed in Section~\ref{sec:niw_posterior_conjugacy}.

Finally, analogous to the canonical form of the univariate Gaussian likelihood in Equation~\eqref{equation:gaussian_form_conform}, 
given fixed mean $\bmu$ and covariance $\bSigma$ parameters, the canonical form for a multivariate Gaussian distribution is:
\begin{equation}\label{equation:multi_gaussian_form_conform}
\begin{aligned}
p(\bx\mid \bmu, \bSigma) &=\normal(\bx \mid \bmu, \bSigma)
\propto \exp\left\{ -\frac{1}{2} \bx^\top\bSigma^{-1}\bx  + \bx^\top \bSigma^{-1}\bmu \right\}.
\end{aligned}
\end{equation}
Thus, if a random variable $\rvx$ has a density matching this functional form, we may conclude that $\rvx\sim \normal(\bmu, \bSigma)$. An example of this appears in the Bayesian \textit{GGGM} matrix decomposition model (see Equation~\eqref{equation:gggm_wm_post}).

\subsection{Properties of Multivariate Gaussian Distribution}\label{section:multi_gauss}
The entropy of multivariate Gaussian distributions (measured in natural units) is discussed in Problem~\ref{problem:entropy_mgau}.
Moreover, affine transformations and rotations of a multivariate Gaussian distribution remain multivariate Gaussian.
\begin{lemma}[Affine Transformation of Multivariate Gaussian Distribution]\label{lemma:affine_mult_gauss}
Let $\bA,\bB\in\real^{M\times N}$ and $\bc\in\real^M$ be fixed (non-random) matrices and vector.
Suppose  $\rvx\sim \normal(\bmu_x, \bSigma_x)$ and $\rvy\sim \normal(\bmu_y, \bSigma_y)$ are independent random vectors in $\real^N$.
Then the linear combination follows a multivariate Gaussian distribution:
$$
\rvz=\bA\rvx+\bB\rvy +\bc \sim \normal(\bA\bmu_x+\bB\bmu_y+\bc, \bA\bSigma_x\bA^\top +\bB\bSigma_y\bB^\top).
$$
Furthermore, for any fixed vector $\bd\in\real^N$, the scalar projection $\bd^\top\rvx$ is univariate Gaussian:
$$
\bd^\top\rvx \sim \normal(\bd^\top\bmu_x, \bd^\top\bSigma_x\bd).
$$
\end{lemma}
This result relies on the fact that \textit{the sum of independent Gaussian} random vectors is also Gaussian:
$$
\sum_{n=1}^N \rvx_n \sim \normal\bigg(\sum_{n=1}^N\bmu_n, \sum_{n=1}^N\bSigma_n\bigg)
\gap 
\text{if } \rvx_n\sim\normal(\bmu_n,\bSigma_n), \,\forall\, n\in\{1,2,\ldots,N\}.
$$

\begin{lemma}[Rotations on  Multivariate Gaussian Distribution]\label{lemma:rotat_multi_gauss}
Rotations preserve the form of isotropic Gaussian distributions. Specifically, if  $\rvv\sim \normal(\bzero, \sigma^2\bI)$
and  $\bQ$ is an orthogonal matrix (i.e., $\bQ\bQ^\top=\bQ^\top\bQ=\bI$), then
$$
\bQ\rvv\sim \normal(\bzero, \sigma^2\bI).
$$
More generally, let $\rvx\sim\normal(\bmu, \bSigma)$ and let $\bSigma=\bU\bLambda\bU^\top$ be the spectral decomposition of $\bSigma$ (Theorem~\ref{theorem:spectral_theorem}).
Then the rotated and centered vector follows a diagonal Gaussian distribution: 
$$
\rvy = \bU^\top(\rvx-\bmu) \sim \normal(\bzero, \bLambda).
$$
\end{lemma}

\paragrapharrow{``Standardization and decorrelation."}
The distribution $\normal(\bzero, \bI)$ is called the \textit{standard multivariate Gaussian distribution}. 
Given $\rvx\sim\normal(\bmu, \bSigma)$, then the \textit{decorrelation} (or standardized variable) of $\rvx$ follows that 
\begin{equation}\label{equation:std_mugau_recov}
\rvx\sim\normal(\bmu, \bSigma)
\quad \implies \quad
\rvz =\bSigma^{-1/2}(\rvx-\bmu) \sim \normal(\bzero,\bI).
\end{equation}
This also shows that if $\rvx\sim\normal(\bmu,\bSigma)$, then 
\begin{equation}
\rvx=\bmu+\bSigma^{1/2}\bepsilon,
\gap \text{where }\bepsilon\sim\normal(\bzero,\bI).
\end{equation}

Now suppose $\{\bx_1,\bx_2,\ldots,\bx_N\}$ are $N$  independent samples from $\normal(\bmu,\bSigma)$ and let $\widebarbx = \frac{1}{N} \sum_{n=1}^{N} \bx_n $. 
Then the sampling distribution of the sample mean satisfies
\begin{equation}\label{equation:mean_mutigau}
\sqrt{N} (\widebarbx-\bmu) \sim \normal(\bzero, \bSigma).
\end{equation}

\paragrapharrow{Quadratic of Gaussian.}
The definition of the Chi-squared distribution (Definition~\ref{definition:chisquare_distribution}) shows 
$$
\sum_{n=1}^{N} \rx_n^2\sim \chi^2(N), 
\gap \text{if }\rx_n \stackrel{i.i.d.}{\sim}  \normal(0,1).
$$
Therefore, we also have
\begin{equation}
\rvx\sim\normal(\bmu, \bSigma)
\implies 
\rvz =(\rvx-\bmu)^\top \bSigma^{-1}(\rvx-\bmu)\sim \chi^2(N), \gap \text{where }\rvx\in\real^N.
\end{equation}

More generally, consider quadratic forms $\rvx^\top\bA\rvx$ where $\bA$ is symmetric. 
The following key results hold:
\begin{theoremHigh}[Quadratic of Gaussians]
We have the following results with quadratic forms of Gaussians:
\begin{itemize}
\item 
Given $\rvx\sim\normal(\bzero,\lambda\bI)$ (of length $N$) and symmetric matrix $\bA\in\real^{N\times N}$. Then, it follows that 
$$
\frac{\rvx^\top\bA\rvx}{\lambda} \sim \chisquared(R),
$$
if and only if $\bA$ ($\bA^2=\bA$) is idempotent and has rank  $R<N$.
\item Given $\rvx\sim\normal(\bzero,\bSigma)$ (of length $N$) and symmetric matrix $\bA\in\real^{N\times N}$. Then, it follows that 
$$
\rvx^\top\bA\rvx \sim \chisquared(R),
$$
if and only if $\bA\bSigma$ is idempotent and has rank  $R<N$.
\end{itemize}
\end{theoremHigh}

\paragrapharrow{Marginal and conditional distributions.}
Let $\rvx$ and $\rvy$ be jointly Gaussian, so that
$$
\rvz=
\begin{bmatrix}
\rvx\\
\rvy 
\end{bmatrix}
\sim 
\normal\left(
\begin{bmatrix}
	\bmu_x\\
	\bmu_y 
\end{bmatrix}
,
\begin{bmatrix}
\bA & \bC\\
\bC^\top & \bB
\end{bmatrix}
\right)=
\normal\left(
\begin{bmatrix}
	\bmu_x\\
	\bmu_y 
\end{bmatrix}
,
\begin{bmatrix}
	\widetildebA & \widetildebC\\
	\widetildebC^\top & \widetildebB
\end{bmatrix}^{-1}
\right).~\footnote{
Given nonsingular $\bM$ and its inverse $\bM^{-1}$; and suppose appropriate sizes for the following partitions \citep{williams2006gaussian}:
$
\bM=
\begin{bmatrix}
\bA & \bB\\
\bC & \bD
\end{bmatrix},
\gap 
\bM^{-1}=
\begin{bmatrix}
	\widetildebA & \widetildebB\\
	\widetildebC & \widetildebD
\end{bmatrix}.
$
We have 
\begin{equation}\label{equation:mt_inv}
\begin{aligned}
&\widetildebA=\bA^{-1}+\bA^{-1}\bB\widetildebD\bC\bA^{-1}&=& (\bA-\bB\bD^{-1}\bC)^{-1} , \\
&\widetildebB=-\bA^{-1}\bB\widetildebD&=& -\widetildebA\bB\bD^{-1},\\
&\widetildebC=-\widetildebD\bC\bA^{-1}&=& -\bD^{-1}\bC\widetildebA, \\
& \widetildebD=(\bD-\bC\bA^{-1}\bB)^{-1}&=& \bD^{-1}+\bD^{-1}\bC\widetildebA\bB\bD^{-1},
\end{aligned}
\end{equation}
}
$$
where the second parametrization uses the inverse covariance (precision) matrix.
The random vectors $\rvx$ and $\rvy$ are \textit{independent} if and only if $\Cov[\rvx,\rvy]=\bC=\bzero$.
Crucially, both marginal and conditional distributions of a multivariate Gaussian are themselves Gaussian:
\begin{equation}\label{equation:marginal_multigaus}
\begin{aligned}
\rvx
\sim\normal(\bmu_x,\bA);
\qquad 
\rvx\mid \rvy=\by 
&\sim \normal(\bmu_x+\bC\bB^{-1}(\by-\bmu_y), \bA-\bC\bB^{-1}\bC^\top)\\
&=\normal(\bmu_x-\widetildebA^{-1}\widetildebC(\by-\bmu_y), \widetildebA^{-1});\\
\rvy
\sim\normal(\bmu_y,\bB);
\qquad 
\rvy\mid \rvx=\bx
&\sim \normal(\bmu_y+\bC^\top\bA^{-1}(\bx-\bmu_x), \bB-\bC^\top\bA^{-1}\bC)\\
&=\normal(\bmu_y-\widetildebB^{-1}\widetildebC^\top(\bx-\bmu_x), \widetildebB^{-1}).\\
\end{aligned}
\end{equation}
\begin{proof}[Sketch of Proof]
Suppose $\rvx^\prime =\rvx-\bC\bB^{-1}$. Then 
$$
\rvz^\prime = 
\begin{bmatrix}
\rvx^\prime \\
\rvy 
\end{bmatrix}
=
\begin{bmatrix}
\bI & -\bC\bB^{-1}\\
\bzero & \bI 
\end{bmatrix}
\rvz.
$$
Using Lemma~\ref{lemma:affine_mult_gauss}, we can show that $\rvx^\prime$ and $\rvy$ are independent.
Then, the conditional distribution of $\rvx \mid \rvy$ can be obtained by $\rvx=\rvx^\prime + \bC\bB^{-1}\rvy$  and following the distribution law.
The second result is analogous.
\end{proof}

These properties are instrumental in deriving posterior distributions in Bayesian models involving Gaussians. See the exercise below.
\begin{exercise}[Linear Gaussian Model: Affine Dependence of Gaussians\index{Linear Gaussian model}]\label{exercise:linear_gauss_model}
Suppose random vectors $\rvx\sim\normal(\bmu, \bSigma)$ and $\rvy\mid \rvx=\bx\sim\normal(\bA\bx+\bb, \bM)$.
Note $\rvy$ is not simply the affine transformation $\bA\rvx+\bb$, but it follows that $\rvy=\bA\bx+\bb+\bepsilon$ with $\bepsilon\sim \normal(\bzero, \bM)$ independent of $\rvx$.
Show that 
$$
\rvy\sim \normal(\bA\bmu+\bb, \bM+\bA\bSigma\bA^\top),
\gap 
\rvx \mid \rvy \sim \normal\big(\bL\big\{\bA^\top\bM^{-1}(\rvy-\bb)+\bSigma^{-1}\bmu  \big\}, \bL\big), 
$$
where $\bL=(\bSigma^{-1}+\bA^\top\bM^{-1}\bA)^{-1}$.
\textit{Hint: compute the cross-covariance of $\rvx$ and $\rvy$ by $\Cov[\rvx,\rvy]=\Exp[(\rvx-\bmu_x)(\rvy-\bmu_y)^\top]=\bSigma\bA^\top$ where $\bmu_x=\bmu$ and $\bmu_y=\Exp[\rvy]=\bA\bmu+\bb$, and use Woodbury matrix identity: $(\bA+\bB\bD\bC)^{-1} = \bA^{-1} - \bA^{-1} \bB(\bD^{-1} + \bC\bA^{-1}\bB)^{-1}\bC\bA^{-1}$ for conformable matrices $\bA,\bB,\bC$, and $\bD$; see, for example,  \citet{lu2021numerical}}.
\end{exercise}

\paragrapharrow{Product of Gaussians.}
The product of two Gaussian density functions is proportional to another Gaussian density (though not normalized) \citep{ahrendt2005multivariate}.
Given two Gaussians $\normal(\bmu_a, \bSigma_a)$ and $\normal(\bmu_b, \bSigma_b)$ (both of length $N$), it follows that 
\begin{equation}
\normal(\bmu_a, \bSigma_a)\cdot \normal(\bmu_b, \bSigma_b) \propto z_c\normal(\bmu_c,\bSigma_c),
\end{equation}
where 
$$
\bSigma_c=(\bSigma_a^{-1}+\bSigma_b^{-1})^{-1},
\gap
\text{and}
\gap 
\bmu_c = \bSigma_c(\bSigma_a^{-1}\bmu_a + \bSigma_b^{-1}\bmu_b).
$$
In other words, the precision (inverse covariance) of the product is the sum of the individual precisions.
The proportionality constant  $z_c$ is given by
$$
\begin{aligned}
z_c=\abs{2\pi \bSigma_a\bSigma_b\bSigma_c^{-1}}^{-\frac{1}{2}} \exp\left\{-\frac{1}{2}(\bmu_a-\bmu_b)^\top \bSigma_a^{-1}\bSigma_c\bSigma_b^{-1}(\bmu_a-\bmu_b)\right\}.
\end{aligned}
$$
 
\subsection{Multivariate Student's $t$ Distribution}
The multivariate Student's $t$-distribution is a continuous probability distribution over multiple variables that generalizes the Gaussian distribution by allowing heavier tails---i.e., it assigns higher probability to extreme values compared to a Gaussian.
The multivariate Student's $t$ distribution (often  simply called Student's $t$ distribution when the context is clear) frequently arises as the posterior predictive distribution for the parameters of a multivariate Gaussian model. We now provide a formal definition.

\begin{figure}[h]
\subfigure[Gaussian, $\bSigma =\begin{bmatrix}
1&0\\
0&1
\end{bmatrix}. $ ]{\includegraphics[width=0.31
\textwidth]{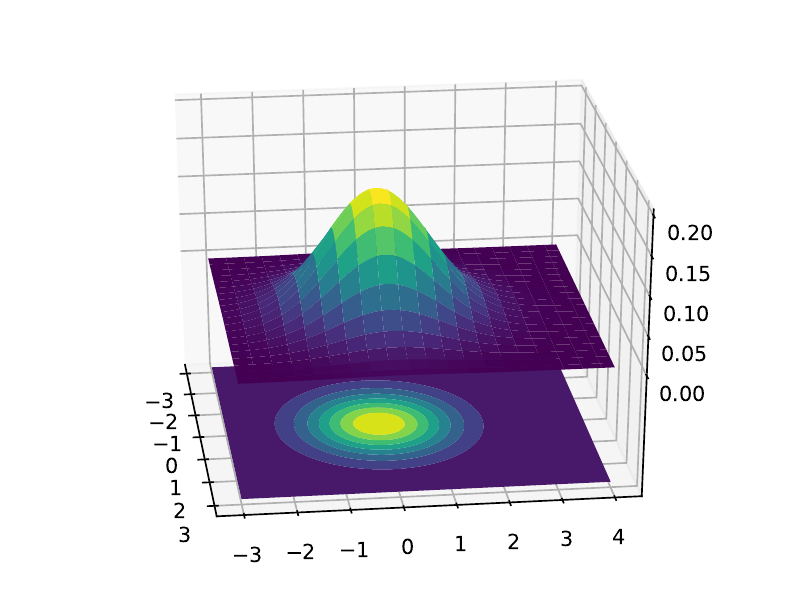} \label{fig:gauss-diagonal}}
\subfigure[Gaussian, $\bSigma =\begin{bmatrix}
1&0\\
0&3
\end{bmatrix}.$]{\includegraphics[width=0.31
\textwidth]{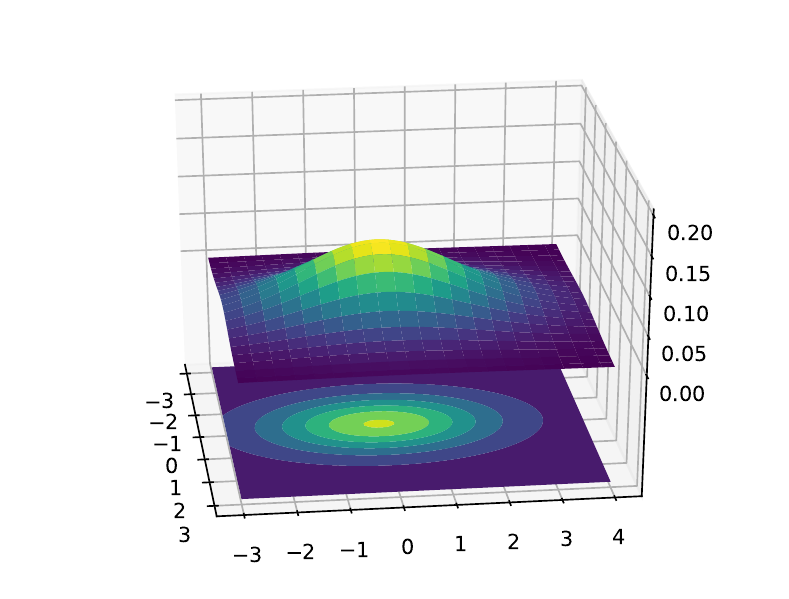} \label{fig:gauss-spherical}}
\subfigure[Gaussian, $\bSigma =\begin{bmatrix}
1&\textendash0.5\\
\textendash0.5&1.5
\end{bmatrix}.$]{\includegraphics[width=0.31 
\textwidth]{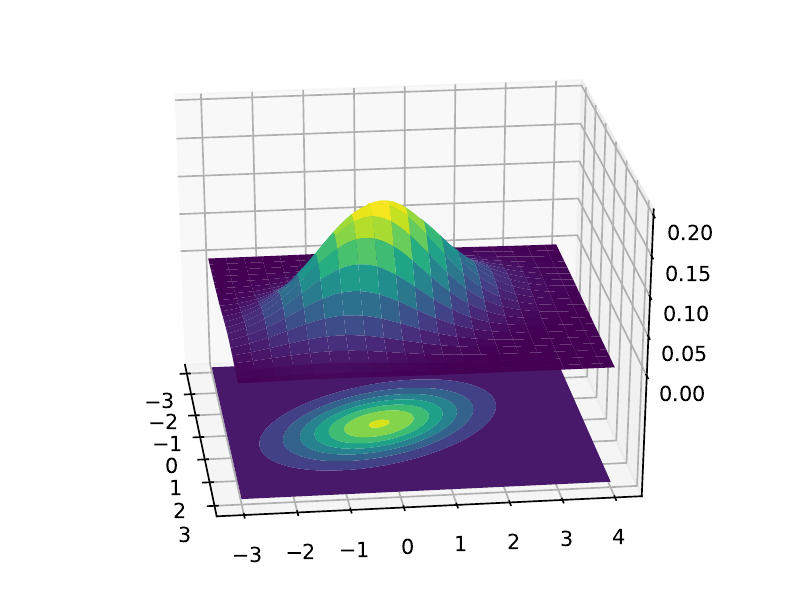} \label{fig:gauss-full}}
\subfigure[Student $t$, $\bSigma =\begin{bmatrix}
1&0\\
0&1
\end{bmatrix}, \nu=1. $ ]{\includegraphics[width=0.31
\textwidth]{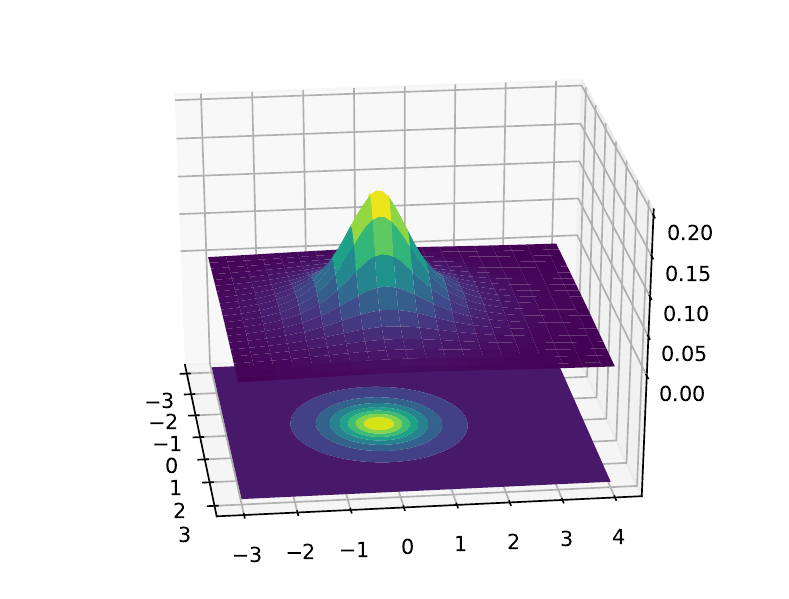} \label{fig:student-1}}
\subfigure[Student $t$, $\bSigma =\begin{bmatrix}
1&0\\
0&1
\end{bmatrix}, \nu=3. $ ]{\includegraphics[width=0.31
\textwidth]{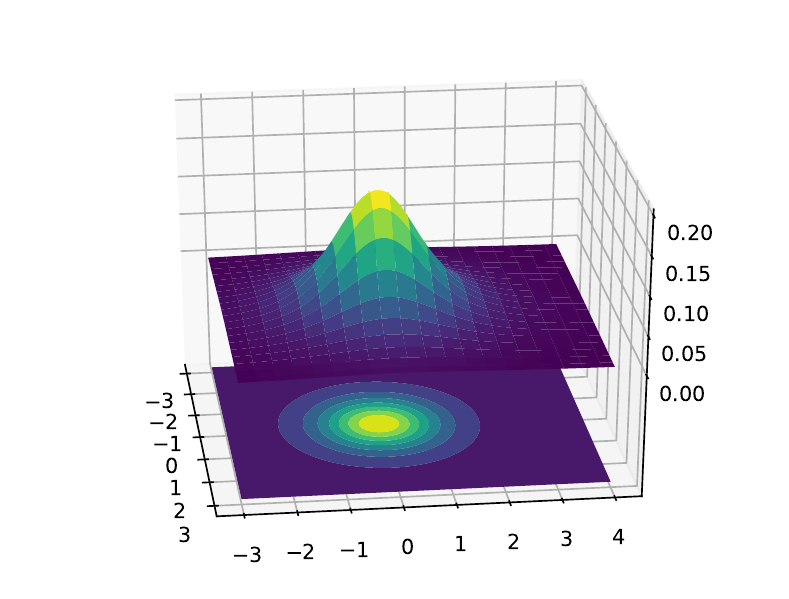} \label{fig:student-3}}
\subfigure[Stu $t$, $\bSigma =\begin{bmatrix}
1&0\\
0&1
\end{bmatrix}, \nu=200. $ ]{\includegraphics[width=0.31
\textwidth]{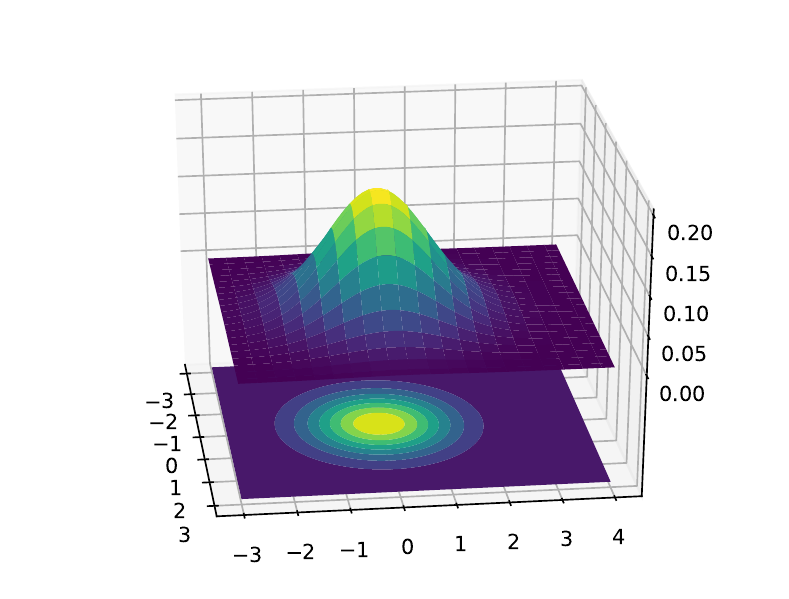} \label{fig:student200}}
\subfigure[Diff between (a) and (d)]{\includegraphics[width=0.31
\textwidth]{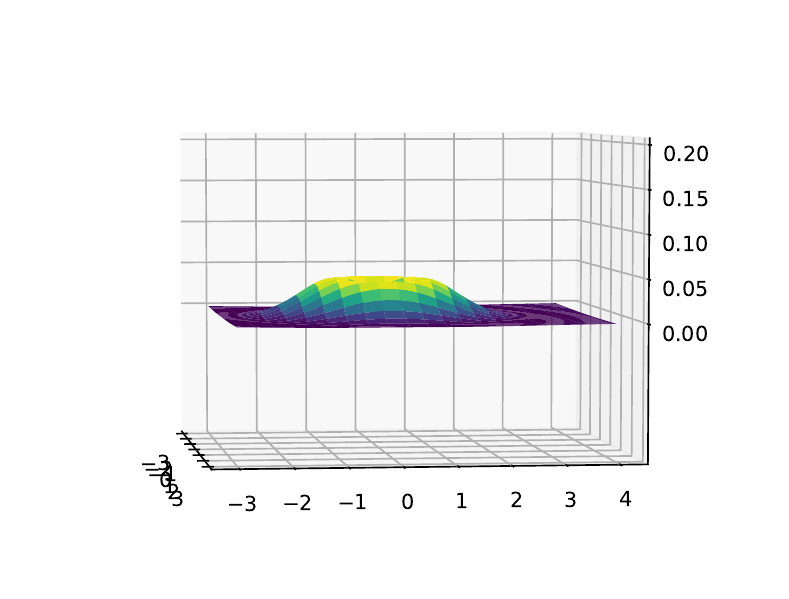} \label{fig:gauss-stu-diff1}}
\subfigure[Diff between (a) and (e)]{\includegraphics[width=0.31
\textwidth]{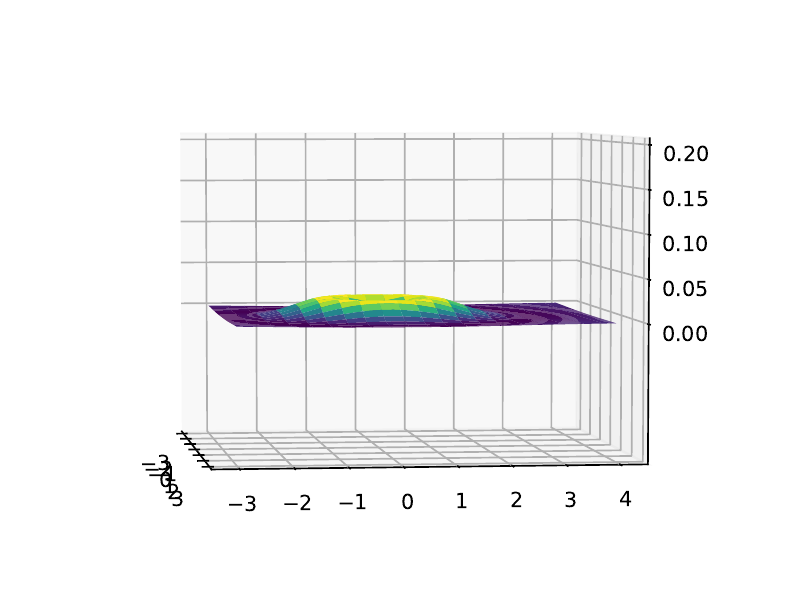} \label{fig:gauss-stu-diff3}}
\subfigure[Diff between (a) and (f)]{\includegraphics[width=0.31
\textwidth]{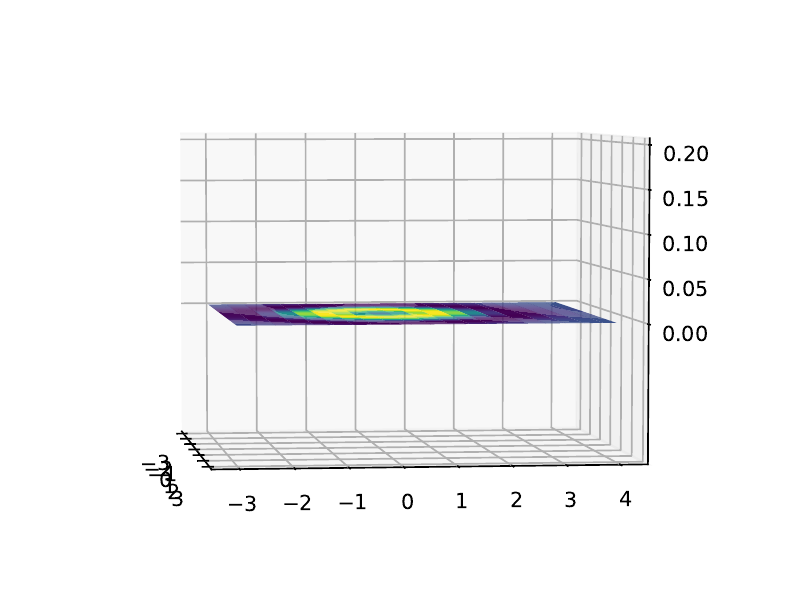} \label{fig:gauss-stu-diff200}}
\centering
\caption{Density and contour plots (\textcolor{darkblue}{blue}=low, \textcolor{mydarkyellow}{yellow}=high) for the multivariate Gaussian distribution and multivariate Student's $t$ distribution over the $\mathbb{R}^2$ space for various values of the covariance/scale matrix with zero-mean vector.  
Fig~\ref{fig:gauss-diagonal}: A spherical covariance matrix has a circular shape; 
Fig~\ref{fig:gauss-spherical}: A diagonal covariance matrix is an \textbf{axis aligned} ellipse; 
Fig~\ref{fig:gauss-full}: A full covariance matrix has a elliptical shape; \\
Fig~\ref{fig:student-1} to Fig~\ref{fig:student200} for the Student's $t$ distribution with the same scale matrix and increasing $\nu$ such that the difference between (a) and (f) in Fig~\ref{fig:gauss-stu-diff200} is approaching  zero.}\centering
\label{fig:studentt_densitys-1}
\end{figure}

\begin{definition}[Multivariate Student's $t$ Distribution\index{Multivariate Student's $t$ distribution}]\label{definition:multivariate-stu-t}
A random vector $\rvx\in\real^D$ is said to follow a \textit{multivariate Student's $t$ distribution} with location  parameters $\bmu\in\real^D$, scale matrix $\bSigma\in\real^{D\times D}$ (symmetric positive definite), and degrees of freedom  $\nu>0$, denoted $\rvx \sim \tau( \bmu, \bSigma, \nu)$, if its probability density function is
$$
\begin{aligned}
f(\bx; \bmu, \bSigma, \nu)&= \frac{\Gamma(\nu/2 + D/2)}{\Gamma(\nu/2)} \frac{\abs{\bSigma}^{-1/2}}{\nu^{D/2} \pi^{D/2}} \times \left[ 1+ \frac{1}{\nu} (\bx-\bmu)^\top \bSigma^{-1} (\bx-\bmu)  \right]^{-(\frac{\nu+D}{2})}\\
&= \frac{\Gamma(\nu/2 + D/2)}{\Gamma(\nu/2)} |\pi\bV|^{-1/2} \times \left[ 1+ \frac{1}{\nu} (\bx-\bmu)^\top \bV^{-1} (\bx-\bmu)  \right]^{-(\frac{\nu+D}{2})},
\end{aligned}
$$
where  $\bV=\nu\bSigma$. 
This distribution has heavier tails than the Gaussian. 
The smaller the value of $\nu$, the heavier the tails. 
As $\nu \rightarrow \infty$, the distribution converges towards a multivariate Gaussian.
The mean, mode, and covariance (when they exist) are given by:
\begin{equation*}
\begin{aligned}
\Exp [\bx] &= \bmu,\qquad \qquad
\mathrm{Mode}[\bx] = \bmu, \qquad\qquad
\Cov [\bx] = \frac{\nu}{\nu-2}\bSigma. 
\end{aligned}
\end{equation*}
Note that $\bSigma$ is not the covariance matrix itself (unless $\nu\rightarrow\infty$); rather, it serves as a scale parameter, which is why it is called the scale matrix.

In the univariate case ($D=1$), the density reduces to the univariate Student's $t$ distribution (see Definition~\ref{equation:student_t_dist}):
\begin{equation}\label{equation:uni-stu-nonzero}
\begin{aligned}
\tau(x\mid \mu, \sigma^2, \nu)&= \frac{\Gamma(\frac{\nu+1}{2})}{\Gamma(\frac{\nu}{2})} \frac{1}{\sigma\sqrt{\nu\pi}} \times \left[ 1+ \frac{(x-\mu)^2}{\nu \sigma^2}   \right]^{-(\frac{\nu+1}{2})}.
\end{aligned}
\end{equation}
When $D=1, \bmu=0, \bSigma=1$, then the PDF defines the \textit{univariate $t$ distribution}.
\begin{equation*}
\begin{aligned}
\tau(x\mid \nu)&= \frac{\Gamma(\frac{\nu+1}{2})}{\Gamma(\frac{\nu}{2})} \frac{1}{\sqrt{\nu\pi}} \times \left[ 1+ \frac{x^2}{\nu }   \right]^{-(\frac{\nu+1}{2})}.
\end{aligned}
\end{equation*}
\end{definition}
Figure~\ref{fig:studentt_densitys-1} compares the Gaussian and the Student's $t$ distribution for various values such that when $\nu\rightarrow \infty$, the difference between the densities is approaching zero. Given the same parameters in the densities, the Student's $t$ in general has longer ``tails" than a Gaussian, which can be seen from the comparison between Figure~\ref{fig:gauss-diagonal} and Figure~\ref{fig:student-1}. 
This heavy-tailed behavior endows the Student's $t$ distribution with an important property known as robustness: it is far less sensitive to outliers than the Gaussian distribution \citep{bishop2006pattern, murphy2012machine}.

The Student's $t$ distribution can be expressed as a \textbf{Gaussian scale mixture}
\begin{equation}\label{equation:gauss-scale-mixture}
\tau(\bx\mid \bmu, \bSigma, \nu) 
= \int_0^{\infty} \normal(\bx \mid \bmu, \bSigma/z)\cdot \gammadist\big(z\mid \frac{\nu}{2}, \frac{\nu}{2}\big)
dz.
\end{equation}
This can be thought of as an ``infinite'' mixture of Gaussians, each with a slightly different covariance matrix. In other words, a Student's $t$ distribution is obtained by adding up an
infinite number of Gaussian distributions having the same mean vector but different covariance matrices. 
From this Gaussian scale mixture view, when $\nu \rightarrow \infty$, the Gamma distribution becomes a degenerate random variable with all the nonzero mass at the point unity such that the multivariate Student's $t$ distribution converges to a multivariate Gaussian distribution.

\paragrapharrow{Affine transformations of Student's $t$.}
Similar to the multivariate Gaussian distribution, the affine transformation of a Student's $t$ also follows another Student's $t$. 
Suppose $\rvx\sim \studentt(\bmu, \bSigma, \nu)$ (of length $D$) and given a  fixed matrix $\bA\in\real^{P\times D}$ and a fixed vector $\bb\in\real^P$. Then it follows that 
\begin{equation}
\bA\rvx \sim \studentt(\bA\bmu+\bb, \bA\bSigma\bA^\top, \nu).
\end{equation}
Therefore, we can sample $\rvx\sim \studentt(\bmu, \bSigma, \nu)$ by sampling $\rvy\sim\studentt(\bzero, \bI, \nu)$ and letting $\rvx=\bmu+\bL\rvy$, where $\bSigma=\bL\bL^\top$ is the Cholesky decomposition of $\bSigma$; see Problem~\ref{problem:cholesky}.

\paragrapharrow{Marginal and conditional distributions of Student's $t$.}
Similar to the multivariate Gaussian distribution, let $\rvx$ and $\rvy$ be jointly Student's $t$ random vectors with 
$$
\rvz=
\begin{bmatrix}
\rvx\\
\rvy 
\end{bmatrix}
\sim 
\studentt\left(
\begin{bmatrix}
\bmu_x\\
\bmu_y 
\end{bmatrix}
,
\begin{bmatrix}
\bA & \bC\\
\bC^\top & \bB
\end{bmatrix} , \nu
\right)=
\studentt\left(
\begin{bmatrix}
\bmu_x\\
\bmu_y 
\end{bmatrix}
,
\begin{bmatrix}
\widetildebA & \widetildebC\\
\widetildebC^\top & \widetildebB
\end{bmatrix}^{-1}, \nu
\right),
$$
where $\rvx\in\real^{D_x}$ and $\rvy\in\real^{D_y}$.
Then every marginal distribution of a Student's $t$ distribution is itself a Student's $t$ distribution, and the conditional distribution $\rvx\mid \rvy$ also follows a Student's $t$ distribution:
\begin{equation}
\begin{aligned}
\rvx
\sim\studentt(\bmu_x,\bA, \nu),
\gap 
\rvx\mid \rvy=\by 
&\sim \studentt(\bmu_x+\bC\bB^{-1}(\by-\bmu_y), m_x(\bA-\bC\bB^{-1}\bC^\top), \nu+D_x)\\
&=\studentt(\bmu_x-\widetildebA^{-1}\widetildebC(\by-\bmu_y), m_x\widetildebA^{-1}, \nu+D_x);\\
\rvy
\sim\studentt(\bmu_y,\bB, \nu),
\gap 
\rvy\mid \rvx=\bx
&\sim \studentt(\bmu_y+\bC^\top\bA^{-1}(\bx-\bmu_x), m_y(\bB-\bC^\top\bA^{-1}\bC), \nu+D_y)\\
&=\studentt(\bmu_y-\widetildebB^{-1}\widetildebC^\top(\bx-\bmu_x), m_y\widetildebB^{-1}, \nu+D_y),\\
\end{aligned}
\end{equation}
where 
$$
\begin{aligned}
m_x &= \frac{1}{\nu+D_y}\left[\nu+ (\by-\bmu_y)^\top \bB^{-1} (\by-\bmu_y)\right];\\
m_y &= \frac{1}{\nu+D_x}\left[\nu+ (\bx-\bmu_x)^\top \bA^{-1} (\bx-\bmu_x)\right].
\end{aligned}
$$
Unlike the Gaussian case, the conditional scale matrix is scaled by a data-dependent factor ($m_x$ or $m_y$), and the degrees of freedom increase by the dimension of the conditioning variable. This reflects the \textbf{adaptive robustness} of the Student's $t$ model.

\subsection{Prior on Parameters of Multivariate Gaussian Distribution}
In Equation~\eqref{equation:inverse_gamma_conjugacy_general}, we have shown that the inverse-Gamma distribution serves as a conjugate prior for the variance parameter of a univariate Gaussian distribution. 
A natural multivariate generalization of this idea is the \textit{inverse-Wishart} distribution, which acts as a conjugate prior for the full covariance matrix of a multivariate Gaussian distribution.
Specifically, the inverse-Wishart is a probability distribution over random symmetric positive definite matrices and is commonly used to model uncertainty in covariance matrices.

Before introducing the inverse-Wishart distribution, it is helpful to recall its origin: the \textit{Wishart distribution}, a multivariate generalization of the Gamma distribution. 
As note by \citet{anderson1962introduction} in 1962,  ``\textit{The Wishart distribution ranks next to the (multivariate) normal distribution in order of importance and usefulness in multivariate statistics.}"

\begin{definition}[Wishart Distribution\index{Wishart distribution}]\label{definition:wishart_dist}
A random symmetric positive definite matrix $\bLambda\in \real^{D\times D}$ is said to follow a \textit{Wishart distribution} with scale matrix $\bM\in\real^{D\times D}$ (symmetric positive definite) and degrees of freedom $\nu\geq D$, denoted $\bLambda \sim \wishartdist(\bM, \nu)$,  if its probability density function is
$$
\begin{aligned}
&\gap f(\bLambda;\textcolor{black}{\bM}, \nu)\\
&= \abs{\bLambda}^{\textcolor{black}{\frac{\nu-D-1}{2}}} \exp\left\{-\frac{1}{2}\tr(\textcolor{black}{\bLambda} \textcolor{black}{\bM^{-1}})\right\}
\left[2^{\frac{\nu D}{2}}  \pi^{D(D-1)/4}  \textcolor{black}{\abs{\bM}^{\nu/2}  } \prod_{d=1}^D\Gamma\big(\frac{\nu+1-d}{2}\big) \right]^{-1}.
\footnote{In some texts, the density function is defined using the generalized Gamma function: $\Gamma_d(x)=\pi^{d(d-1)/4}\prod_{i=1}^{d}\Gamma(\frac{2x+1-i}{2})$, such that 
$$
f(\bLambda;\textcolor{black}{\bM}, \nu)\\
= \abs{\bLambda}^{\textcolor{black}{\frac{\nu-D-1}{2}}} \exp\left\{-\frac{1}{2}\tr(\textcolor{black}{\bLambda} \textcolor{black}{\bM^{-1}})\right\}
\left[2^{\frac{\nu D}{2}}    \textcolor{black}{\abs{\bM}^{\nu/2}  } \Gamma_D(\nu/2) \right]^{-1}.
$$
}
\end{aligned}
$$
Here, $\abs{\bLambda} = \det(\bLambda)$ denotes the determinant of matrix $\bLambda$.
The mean and element-wise variances of the Wishart distribution are given by:
$$
\begin{aligned}
\Exp [\bLambda] = \nu \bM,
\qquad \qquad
\Var[\bLambda_{ij}] = \nu (m_{ij}^2 + m_{ii}m_{jj}),
\label{equation:wishart_expectation}
\end{aligned}
$$
where $m_{ij}$ is the ($i,j$)-th entry  of $\bM$.
Moreover, as $\nu\rightarrow \infty$, the scaled matrix $\bLambda/\nu$ converges in probability to $\bM$ (using law of large numbers and the Cramer-Wold device).

When $D=1$ and $\bM=1$, the Wishart distribution reduces to the Chi-squared distribution (Definition~\ref{definition:chisquare_distribution}) such that:
$$
\wishartdist(x\mid  1, \nu) = \chisquared(x\mid \nu).
$$
\end{definition}

An intuitive interpretation of the Wishart distribution arises from sampling. 
Suppose we independently draw  vectors $\bz_1, \bz_2, \ldots, \bz_{\nu}\in\real^D$ from $\normal(\bzero, \bM)$. The sum of squares matrix of the collection of multivariate vectors is given by 
$$
\sum_{i=1}^{\nu} \bz_i\bz_i^\top = \bZ^\top\bZ,
$$
where $\bZ$ is the $\nu \times D$ matrix whose $i$-th row is $\bz_i^\top$. 
It is evident that $\bZ^\top\bZ$ is positive semidefinite (PSD). 
If $\nu >D$ and the vectors $\bz_i$ are linearly independent, then $\bZ^\top \bZ$ is positive definite (PD). 
In other words, $\bZ\bx=\bzero$ only happens when $\bx=\bzero$. We can repeat over and over again, generating matrices $\bZ_1^\top\bZ_1, \bZ_2^\top\bZ_2, \ldots, \bZ_l^\top\bZ_l$. The population distribution of these matrices follows a Wishart distribution with parameters $(\bM, \nu)$. By definition, 
$$
\begin{aligned}
\bLambda&=\bZ^\top\bZ = \sum_{i=1}^{\nu} \bz_i\bz_i^\top; \\
\Exp[\bLambda]&=\Exp[\bZ^\top\bZ] = \Exp\left[\sum_{i=1}^{\nu} \bz_i\bz_i^\top\right] = \nu \Exp[\bz_i\bz_i^\top] = \nu\bM. \\
\end{aligned}
$$

In the scalar case ($D=1$), this reduces to the well-known result: if $z$ is drawn from a zero-mean univariate normal random variable, then $z^2$ is drawn from a Gamma random variable:
$$
\text{if }   z \sim \normal(0, a) , \qquad \text{then }  z^2\sim \gammadist(a/2, 1/2).
$$

\begin{remark}[Properties of Wishart Distribution]
We list several important properties without proof:
\begin{itemize}
\item \textbf{``Decorrelation."} Suppose $\bLambda\sim\wishartdist(\bM, \nu)$ with $\bLambda\in\real^{D\times D}$. 
Then, it follows that $\bM^{-1/2}\bLambda\bM^{-1/2}\sim \wishartdist(\nu, \bI_D)$.
\item \textbf{Quadratic transformation.} Suppose $\bLambda\sim\wishartdist(\bM, \nu)$ with $\bLambda\in\real^{D\times D}$ and $\bA\in\real^{P\times D}$. 
Then, it follows that  $\bA\bLambda\bA^\top\sim \wishartdist(\bA\bM\bA^\top, \nu)$.
\item Suppose $\bLambda\sim\wishartdist(\bM, \nu)$ with $\bLambda\in\real^{D\times D}$, $\ba\in\real^D$, and $\nu>D-1$. 
Then, it follows that $\frac{\ba^\top\bM^{-1}\ba}{\ba^\top\bLambda^{-1}\ba}\sim \chisquared(\nu-D-1)$.
\item Suppose $\bLambda\sim\wishartdist(\bM, \nu)$ with $\bLambda\in\real^{D\times D}$ and $\ba\in\real^D$. 
Then, it follows that $\frac{\ba^\top\bLambda\ba}{\ba^\top\bM\ba}\sim \chisquared(\nu)$.
\item \textbf{\text{Sum of independent Wisharts}.} Given independent random matrices $\bLambda_i\sim\wishartdist(\bM, \nu_i)$ with $\nu=\sum_{i}\nu_i$. 
Then, it follows that $\sum_{i}\bLambda_i\sim \wishartdist(\bM, \nu)$.
\item \textbf{\text{Sum of independent Wisharts}.} Similarly, given independent random matrices $\bLambda\sim\wishartdist(\bM, \nu)$ and $\bLambda_1\sim\wishartdist(\bM, \nu_1)$. 
Then, it follows that $\bLambda_2=\bLambda-\bLambda_1\sim \wishartdist(\bM,\nu-\nu_1)$.
\item \textbf{``Standardization".} Suppose $\{\bx_1,\bx_2,\ldots,\bx_N\}$ are random samples of $\normal(\bmu,\bSigma)$, let $\widebarbx = \frac{1}{N} \sum_{n=1}^{N} \bx_n $ and $\bS=\frac{1}{N-1}\sum_{n=1}^{N}(\bx_n-\widebarbx)(\bx_n-\widebarbx)^\top$.
Then, it follows that $(N-1)\bS\sim \wishartdist(\bSigma,N-1)$. And it can be shown that $\widebarbx$ and $\bS$ are independent (the distribution of $\widebarbx$ is shown in \eqref{equation:mean_mutigau}). 
\end{itemize}
\end{remark}

Just as the inverse-Gamma distribution is related to the Gamma distribution---namely, if $x \sim \gammadist(r, \lambda)$, then $y={1}/{x} \sim \inversegammadist(r, \lambda)$---the \textit{inverse-Wishart (IW) distribution} is defined analogously from the Wishart distribution.

Since the inverse-Wishart is typically used as a prior for a covariance matrix, it is often useful to replace $\bM$ in the Wishart distribution with $\bS=\bM^{-1}$. This results in that
a random $D\times D$ symmetric positive definite matrix $\bSigma$ follows an inverse-Wishart $\inversewishart(\bSigma\mid \bS, \nu)$ distribution if $\bSigma^{-1}=\bLambda$ follows a Wishart $\mathrm{Wi}(\bLambda\mid \bM, \nu)$ distribution. 

\begin{definition}[Inverse-Wishart Distribution\index{Inverse-Wishart distribution}]\label{definition:multi_inverse_wishart}
A random symmetric positive definite matrix $\bSigma\in \real^{D\times D}$ is said to follow an \textit{inverse-Wishart distribution} with \textit{scale matrix} $\bS\in\real^{D\times D}$ (symmetric positive definite matrix) and \textit{degrees of freedom} $\nu> D$, denoted $\bSigma\sim \inversewishart(\bS, \nu)$, if
$$
\begin{aligned}
&\gap f(\bSigma; \textcolor{winestain}{\bS}, \nu)\\
&= \abs{\bSigma}^{\textcolor{mylightbluetext}{-\frac{\nu+D+1}{2}}} \exp\left\{-\frac{1}{2}\tr(\textcolor{mylightbluetext}{\bSigma^{-1}} \textcolor{winestain}{\bS})\right\}
\times \left[2^{\frac{\nu D}{2}}  \pi^{D(D-1)/4}  \textcolor{winestain}{\abs{\bS}^{-\nu/2}}   \prod_{d=1}^D\Gamma\big(\frac{\nu+1-d}{2}\big) \right]^{-1},
\end{aligned}
$$
where  $\abs{\bSigma} = \det(\bSigma)$ denotes the determinant.
The mean and mode of the inverse-Wishart distribution are given by 
\begin{equation}
\begin{aligned}
\Exp [\bSigma ^{-1}] = \nu \bS^{-1}=\nu \bM, \qquad
\Exp [\bSigma] = \frac{1}{\nu - D - 1} \bS, \qquad
\mathrm{Mode}[\bSigma] = \frac{1}{\nu + D + 1} \bS.
\label{equation:iw_expectation}
\end{aligned}
\end{equation}
Note that, sometimes, we replace $\bS$ by $\bM=\bS^{-1}$ such that $\Exp [\bSigma ^{-1}] = \nu \bM$, which does not involve the inverse of the matrix.

When $D=1$, the inverse-Wishart distribution reduces to the inverse-Gamma such that $\frac{\nu}{2} = r$ and $\frac{S}{2}=\lambda$ (see Definition~\ref{definition:inverse_gamma_distribution}):
$$
\mathrm{IW}(y\mid S, \nu) = \inversegammadist(y\mid r, \lambda).
$$
\end{definition}

Note that the Wishart density is not simply the inverse-Wishart density with
$\bSigma$ replaced by $\bLambda = \bSigma^{-1}$. There is an additional factor of $\abs{\bSigma}^{-(D+1)}$. See  Theorem 7.7.1 in \citet{anderson1962introduction} that the change of variables $\bLambda = \bSigma^{-1}$ introduces a Jacobian factor of $\abs{\bSigma}^{-(D+1)}$. 
Substitution of $\bSigma^{-1}$ in the definition of the Wishart distribution and multiplying by $\abs{\bSigma}^{-(D+1)}$ yield the inverse-Wishart distribution.
\footnote{Which is from the Jacobian in the change-of-variables formula. A short proof is provided here. Let $\bLambda = g(\bSigma)=\bSigma^{-1}$, where $\bSigma\sim \inversewishart(\bS, \nu)$ and $\bLambda\sim \wishartdist(\bS, \nu)$. Then, $f(\bSigma)  = f(\bLambda) \abs{J_g}$, where $J_g$ is the Jacobian matrix,  results in $f(\bSigma) = f(\bLambda) \abs{J_g} = f(\bLambda)\abs{\bSigma}^{-(D+1)} $. }

The multivariate analog of the normal-inverse-Chi-squared distribution (Definition~\ref{definition:normal_inverse_chi_square}) is the \textit{normal-inverse-Wishart (NIW) distribution} \citep{murphy2007conjugate}. 
We will see that a sample drawn from a normal-inverse-Wishart  distribution, a joint conjugate prior, provides a mean vector and a covariance matrix that can define a multivariate Gaussian distribution. Separately, we can first sample a matrix $\bSigma$ from an inverse-Wishart distribution parameterized by \{$\bS_0, \nu_0$, $\bmu$\} (this is called a \textit{semi-conjugate prior}), and then sample a mean vector from a Gaussian distribution parameterized by \{$\bmm_0, \bV_0, \bSigma$\}. \footnote{Here we use a subscript value of $0$ to indicate the parameters are used for prior density. However, in the Bayesian matrix decomposition analysis, things become more complex and the prior parameters could have other subscript values.}

\begin{definition}[Normal-Inverse-Wishart (NIW) Distribution\index{Normal-inverse-Wishart distribution}]\label{definition:normal_inverse_wishart}
Analog to the (univariate) normal-inverse-Chi-squared distribution, the multivariate counterpart, \textit{the normal-inverse-Wishart (NIW)} distribution is defined as
$$
\begin{aligned}
&\gap \niw (\bmu, \bSigma\mid \bmm, \kappa, \nu, \bS) 
= \normal(\bmu\mid \bmm, \frac{1}{\kappa}\bSigma) \cdot  \inversewishart(\bSigma\mid \bS, \nu) \\
&=\frac{1}{Z_{\niw}(D, \kappa, \nu, \bS)} \abs{\bSigma}^{\frac{-1}{2}}\exp\left\{\frac{\kappa}{2}(\bmu - \bmm)^\top\bSigma^{-1}(\bmu - \bmm)\right\} 
 \abs{\bSigma}^{-\frac{\nu+D+1}{2}} \exp\left\{-\frac{1}{2}\tr(\bSigma^{-1} \bS)\right\} \\
&= \frac{1}{Z_{\niw}(D, \kappa, \nu, \bS)} \abs{\bSigma}^{-\frac{\nu+D+2}{2}}
 \exp\left\{-\frac{\kappa}{2}(\bmu - \bmm)^\top\bSigma^{-1}(\bmu - \bmm) -\frac{1}{2}\tr(\bSigma^{-1} \bS)\right\}, 
\end{aligned}
$$
where the random vector $\bmu\in\real^D$ and the random positive definite matrix $\bSigma\in\real^{D\times D}$ are said to follow NIW, denoted $\bmu, \bSigma\sim \niw(\bmm, \kappa, \nu, \bS)$. And $Z_{\niw}(D, \kappa, \nu, \bS)$ is a normalizing constant:
\begin{equation}
Z_{\niw}(D, \kappa, \nu, \bS) =  2^{\frac{(\nu+1)D}{2}} \pi^{D(D+1)/4} \kappa^{-D/2} | \bS|^{-\nu/2}\prod_{d=1}^D\Gamma\big(\frac{\nu+1-d}{2}\big).
\label{equation:multi_gaussian_giw_constant}
\end{equation}
\end{definition}

\subsection{Posterior Distribution of $\bmu$: Separated View}\label{section:sep_mu_niw}
We now proceed to discuss the posterior distribution of a multivariate Gaussian model under NIW or inverse-Wishart priors, considering both a separated view (treating mean and covariance separately) and a unified view (using the joint NIW prior).

Consider $N$ independent observations $\mathcalX = \{\bx_1, \bx_2, \ldots , \bx_N \}$ drawn from a multivariate Gaussian distribution with unknown mean $\bmu$ and covariance $\bSigma$.
Suppose the covariance matrix $\bSigma$ is known. 
Then, from Equation~\eqref{equation:multi_gaussian_likelihood} (equality (a)), the likelihood function is:
$$
\begin{aligned}
\mathrm{\textbf{likelihood}}
&=p(\mathcalX \mid \bmu) =\normal(\mathcalX\mid \bmu, \bSigma)=\prod^N_{n=1} \normal (\bxn\mid \bmu, \bSigma)\\
&= (2\pi)^{-ND/2} \abs{\bSigma}^{-N/2}\exp\left\{-\frac{1}{2} \sum^N_{n=1}(\bxn - \bmu)^\top \bSigma^{-1}(\bxn - \bmu)\right\}\\
&\propto \exp\left\{   - \frac{1}{2}N \bmu^\top \bSigma^{-1}\bmu  + N\widebarbx^\top \bSigma^{-1} \bmu  \right\},
\end{aligned}
$$
where $\widebarbx =  (\sum_{n=1}^{N}\bx_n)/N$ is the sample mean.
The conjugate prior for the mean vector is Gaussian: $p(\bmu)= \normal(\bmu \mid \bmo, \bV_0)$. 
Its density is:
$$
\begin{aligned}
\mathrm{\textbf{prior}}
&=p(\bmu)= \normal(\bmu \mid \bmo, \bV_0)
= (2\pi)^{-D/2} \abs{\bV_0}^{-1/2}\exp\left\{-\frac{1}{2} (\bmu - \bmo)^\top \bV_0^{-1}(\bmu - \bmo)\right\} \\
&= (2\pi)^{-D/2} \abs{\bV_0}^{-1/2}\exp \left\{-\frac{1}{2}\bmu^\top\bV_0^{-1}\bmu + \bmu^\top \bV_0^{-1}\bmo - \frac{1}{2} \bmo^\top\bV_0^{-1} \bmo    \right\}\\
&\propto  \exp \left\{ -\frac{1}{2}\bmu^\top\bV_0^{-1}\bmu + \bmu^\top \bV_0^{-1}\bmo\right\}.
\end{aligned}
$$
Applying  Bayes' theorem ``$\mathrm{posterior} \propto \mathrm{likelihood} \times \mathrm{prior} $," the posterior distribution of $\bmu$ is also Gaussian:
$$
\begin{aligned}
\mathrm{\textbf{posterior}}&=p(\bmu \mid \mathcalX, \bSigma) \propto p(\mathcalX \mid \bmu, \bSigma ) \times  p(\bmu)\\
&=\exp\left( N\widebarbx^\top \bSigma^{-1} \bmu - \frac{1}{2}N \bmu^\top \bSigma^{-1}\bmu    \right) \times \exp \left( -\frac{1}{2}\bmu^\top\bV_0^{-1}\bmu + \bmu^\top \bV_0^{-1}\bmo\right)\\
&= \exp\left\{  -\frac{1}{2}\bmu^\top(\bV_0^{-1} + N\bSigma^{-1})\bmu + \bmu^\top( \bV_0^{-1}\bmo +  N\bSigma^{-1}\widebarbx )   \right\} \\
&\propto \normal(\bmu \mid \bmm_N, \bV_N),
\end{aligned}
$$
where $\bV_N^{-1} = \bV_0^{-1} + N\bSigma^{-1}$, and $\bmm_N = \bV_N ( \bV_0^{-1}\bmo +  N\bSigma^{-1}\widebarbx )$. 
Thus, the posterior precision matrix equals the sum of the prior precision $\bV_0^{-1}$ and the data precision $N\bSigma^{-1}$. 
In the limit of a \textbf{non-informative (flat) prior}---achieved by letting $\bV_0 \rightarrow \infty \bI$---the posterior simplifies to: $p(\bmu \mid \mathcalX, \bSigma) =\normal(\bmu \mid \widebarbx, \frac{1}{N}\bSigma)$.

\subsection{Posterior Distribution of $\bSigma$: Separated View}\label{section:sep_sigma_niw}
Now suppose the mean vector $\bmu$ is known. From Equation~\eqref{equation:multi_gaussian_likelihood} (equality (b)), the likelihood becomes:
\begin{equation*}
\begin{aligned}
\mathrm{\textbf{likelihood}}
= p(\mathcalX\mid \bmu, \bSigma) =\prod^N_{n=1} \normal (\bxn\mid \bmu, \bSigma) = (2\pi)^{-ND/2} \abs{\bSigma}^{-N/2}\exp\left\{-\frac{1}{2} \tr(\bSigma^{-1}\bS_{\bmu} )  \right\}.
\end{aligned}
\end{equation*}
The conjugate prior for $\bSigma$ is the inverse-Wishart distribution:
$$
\begin{aligned}
\mathrm{\textbf{prior}}=	\mathrm{IW}(\bSigma\mid \bso, \nu_0)&= \abs{\bSigma}^{-\frac{\nu_0+D+1}{2}} \exp\left\{-\frac{1}{2}\tr(\bSigma^{-1} \bso)\right\}\\
	&\gap \times \left[2^{\frac{\nu_0 D}{2}}  \pi^{D(D-1)/4}  \abs{\bso}^{-\nu_0/2}   \prod_{d=1}^D\Gamma(\frac{\nu_0+1-d}{2}) \right]^{-1}.
\end{aligned}
$$
By Bayes' theorem, the posterior is again a inverse-Wishart distribution with updated parameters:
$$
\begin{aligned}
\mathrm{\textbf{posterior}}  
&=p(\bSigma \mid \mathcalX, \bmu) \propto p(\mathcalX \mid \bmu, \bSigma ) \times  p(\bSigma)\\
&\propto \abs{\bSigma}^{-N/2}\exp\left\{-\frac{1}{2} \tr(\bSigma^{-1}\bS_{\bmu} )  \right\} \times \abs{\bSigma}^{-\frac{\nu_0+D+1}{2}} \exp\left\{-\frac{1}{2}\tr(\bSigma^{-1} \bso)\right\} \\
&= \abs{\bSigma}^{-\frac{\nu_0+N+D+1}{2}} \exp\left\{-\frac{1}{2}\tr\left(\bSigma^{-1} [\bso+\bS_{\bmu}]\right)\right\}
\propto \inversewishart(\bSigma \mid\bso+\bS_{\bmu}, \nu_0+N).
\end{aligned}
$$
The posterior degree of freedom is the sum of the prior degree of freedom $\nu_0$ and the number of observations $N$. And the posterior scale matrix is the sum of the prior scale matrix $\bso$ and the data scale matrix $\bS_{\bmu}$. 
The posterior mean of $\bSigma$ (for $\nu_0+N>D+1$) is given by
$$
\begin{aligned}
\Exp[\bSigma\mid  \mathcalX, \bmu] &= \frac{1}{\nu_0 +N- D - 1} (\bso+\bS_{\bmu}) \\
&= \frac{\nu_0 -D-1}{\nu_0 +N- D - 1} \cdot (\frac{1}{\nu_0 -D-1} \bso) +    \frac{N}{\nu_0 +N- D - 1} \cdot (\frac{1}{N}\bS_{\bmu})\\
&= \lambda \cdot \big(\frac{1}{\nu_0 -D-1} \bso \big) +   (1-\lambda) \cdot \big(\frac{1}{N}\bS_{\bmu}\big),
\end{aligned}
$$
where $\lambda=\frac{\nu_0 -D-1}{\nu_0 +N- D - 1}$, $(\frac{1}{\nu_0 -D-1} \bso)$ is the prior mean of $\bSigma$, and $(\frac{1}{N}\bS_{\bmu})$ is an unbiased estimator of the covariance. 
(As $N\rightarrow \infty$, this estimator  $(\frac{1}{N}\bS_{\bmu})$ converges to the true population covariance matrix.) 
Thus, \textbf{the posterior mean of the covariance matrix can be seen as the weighted average of the prior expectation and the unbiased estimator}. The unbiased estimator can also  be   shown to be equal to the maximum likelihood estimator (MLE) of $\bSigma$. As $N\rightarrow \infty$, it can be shown that the posterior expectation of $\bSigma$ is a consistent \footnote{An estimator $\hat{\theta}_N$ of $\theta$ constructed on the basis of a sample of size $N$ is said to be consistent if $\hat{\theta}_N \stackrel{p}{\longrightarrow} \theta$ as $N\rightarrow \infty$. See also \citet{lu2021rigorous}.} estimator of the population covariance. In particular, setting $\nu_0=D+1$ yields $\lambda=0$,  and the posterior mean reduces exactly to the MLE.

Similarly, the posterior mode of $\bSigma$ is given by
\begin{equation}\label{equation:map-covariance-multigauss}
\begin{aligned}
	\mathrm{Mode}[\bSigma] &= \frac{1}{\nu_0 +N+ D + 1} (\bso+\bS_{\bmu})\\
	&= \frac{\nu_0+D+1}{\nu_0 +N+ D + 1} (\frac{1}{\nu_0+D+1} \bso)+  \frac{N}{\nu_0 +N+ D + 1} (\frac{1}{N}  \bS_{\bmu})\\
	&=\beta (\frac{1}{\nu_0+D+1} \bso)+  (1-\beta)(\frac{1}{N}  \bS_{\bmu}),
\end{aligned}
\end{equation}
where $\beta=\frac{\nu_0+D+1}{\nu_0 +N+ D + 1}$, and $(\frac{1}{\nu_0+D+1} \bso)$ is the prior mode of $\bSigma$. \textbf{The posterior mode is a weighted average of the prior mode and the unbiased estimator}. Again, the maximum a posterior (MAP) estimator in Equation~\eqref{equation:map-covariance-multigauss} is a consistent estimator.

\subsection{Gibbs Sampling of the Mean and Covariance: Separated View}
The separated view presented here is known as a \textbf{semi-conjugate prior} on the mean and covariance of a multivariate Gaussian distribution since both conditionals, $p(\bmu\mid\mathcalX,\bSigma)$ and $p(\bSigma\mid\mathcalX,\bmu)$, are individually conjugate.
In the last two sections, we have shown 
$$
\begin{aligned}
\bmu \mid \mathcalX, \bSigma &\sim \normal(\bmm_N, \bV_N),\\
\bSigma \mid \mathcalX, \bmu &\sim \inversewishart(\bso+\bS_{\bmu}, \nu_0+N).
\end{aligned}
$$
The two full conditional distributions can be used to construct a Gibbs sampler. The Gibbs sampler generates the mean and covariance $\{\bmu^\toptone, \bSigma^\toptone\}$ for $(t+1)$-th step from $\{\bmu^\toptzero, \bSigma^\toptzero\}$ in $t$-th step  via the following two steps:
\begin{enumerate}
\item  Sample $\bmu^\toptone$ from its full conditional distribution: $\bmu^\toptone \sim \normal(\bmm_N, \bV_N)$, where $\{\bmm_N, \bV_N\}$ depend on $\bSigma^\toptzero$.

\item  Sample $\bSigma^\toptone$ from its full conditional distribution: $\bSigma^\toptone \sim \inversewishart(\bso+\bS_{\bmu}, \nu_0+N)$, where $\{\bso+\bS_{\bmu}, \nu_0+N\}$ depend on $\bmu^\toptone$.
\end{enumerate}
After sufficient burn-in and thinning iterations, the sequence $\{\bmu^\toptzero, \bSigma^\toptzero\}_{t=1}^T$ approximates draws from the true posterior distribution.

\subsection{Posterior Distribution of $\bmu$ and $\bSigma$ under NIW: Unified View}\label{sec:niw_posterior_conjugacy}

We now show that the normal-inverse-Wishart (NIW) distribution serves as a fully conjugate prior for the joint parameters---mean vector $\bmu$ and covariance matrix $\bSigma$---of a multivariate Gaussian model.

\paragrapharrow{Likelihood.}
Given $N$ independent observations $\mathcalX = \{\bx_1, \bx_2, \ldots , \bx_N \}$ drawn from the multivariate Gaussian with mean vector $\bmu$ and covariance matrix $\bSigma$, the likelihood (see equality (c) in Equation~\eqref{equation:multi_gaussian_likelihood}) can be written as:
$$
\begin{aligned}
p(\mathcalX \mid \bmu, \bSigma)=  \frac{1}{(2\pi)^{ND/2}} \abs{\bSigma}^{-N/2}\exp\left\{-\frac{N}{2}(\bmu - \widebarbx)^\top \bSigma^{-1}(\bmu - \widebarbx)-\frac{1}{2}\tr(\bSigma^{-1} \bS_{\widebarx})\right\}.
\end{aligned}
$$
\paragrapharrow{Prior.}
A naive approach might be to combine the individual conjugate priors for $\bmu$ and $\bSigma$:
$$
p(\bmu, \bSigma) = \normal(\bmu \mid \bmo, \bV_0)\cdot \inversewishart(\bSigma \mid \bso, \nu_0).
$$
However, this factorized prior is not conjugate to the full likelihood because $\bmu$ and $\bSigma$ appear together in a coupled way in the exponent.
In contrast, the NIW prior is fully conjugate. It is defined as:
\begin{equation}\label{equation:multi_gaussian_prior}
\begin{aligned}
&\gap \niw (\bmu, \bSigma\mid \bmo, \kappa_0, \nu_0, \bso) 
= \normal(\bmu\mid \bmo, \frac{1}{\kappa_0}\bSigma) 
\cdot  \inversewishart(\bSigma\mid \bso, \nu_0) \\
&= \frac{\abs{\bSigma}^{-\frac{\nu_0+D+2}{2}}}{Z_{\niw}(D, \kappa_0, \nu_0, \bso)} 
\cdot  \exp\left\{-\frac{\kappa_0}{2}(\bmu - \bmo)^\top\bSigma^{-1}(\bmu - \bmo) -\frac{1}{2}\tr(\bSigma^{-1} \bso)\right\}, 
\end{aligned}
\end{equation}
where
\begin{equation}\label{equation:multi_gaussian_giw_constant_2}
Z_{\niw}(D, \kappa_0, \nu_0, \bso) =  2^{\frac{(\nu_0+1)D}{2}} \pi^{D(D+1)/4} \kappa_0^{-D/2} | \bso|^{-\nu_0/2}\prod_{d=1}^D\Gamma\big(\frac{\nu_0+1-d}{2}\big).
\end{equation}
The specific form of the normalization term $Z_{\niw}(D, \kappa_0, \nu_0, \bso)$ will be useful to show the posterior marginal likelihood of the data in Section~\ref{section:posterior-marginal-of-data}.

\paragrapharrow{A ``prior" interpretation for the NIW prior.}
The inverse-Wishart distribution will ensure that the resulting covariance matrix is positive definite when $\nu_0 > D$. 
If we believe the true covariance is close to some target matrix $\bSigma_0$, then we can choose a large value of $\nu_0$ and set $\bso = (\nu_0 - D - 1) \bSigma_0$, making the distribution of the covariance matrix $\bSigma$ concentrated around $\bSigma_0$. On the other hand, choosing $\nu_0 = D+2$ and $\bso = \bSigma_0$ will make $\bSigma$ loosely concentrated around $\bSigma_0$.
See \citet{chipman2001practical, fraley2007bayesian, hoff2009first, murphy2012machine} for further discussion.

An intuitive interpretation of the hyper-parameters is as follows: $\bmo$ is our prior
mean for $\bmu$, $\kappa_0$ denotes how strongly we believe this prior for $\bmu$ (the larger the stronger belief in the prior mean), $\bso$ is proportional to our prior mean for $\bSigma$, and $\nu_0$ controls how strongly we believe this prior for $\bSigma$. Because the Gamma function is not defined for negative integers and zero, from Equation~\eqref{equation:multi_gaussian_giw_constant_2}, we require $\nu_0 > D - 1$ (which  can also be shown from the expectation of the covariance matrix from Equation~\eqref{equation:iw_expectation}). 
Additionally, $\bso$ must be a positive definite matrix, where an intuitive reason can be shown from Equation~\eqref{equation:iw_expectation}. See \citet{hoff2009first, murphy2012machine} for further discussion.

\paragrapharrow{Posterior.}
By  Bayes' theorem, ``$\mathrm{posterior} \propto \mathrm{likelihood} \times \mathrm{prior} $," the joint posterior of the $\bmu$ and $\bSigma$ parameters under the NIW prior is
\begin{equation}
p(\bmu, \bSigma\mid \mathcalX, \bbeta ) \propto p(\mathcalX \mid \bmu, \bSigma) p(\bmu, \bSigma \mid \bbeta) = p(\mathcalX, \bmu, \bSigma \mid \bbeta),  
\label{equation:niw_full_posterior}
\end{equation}
where $\bbeta=\{\bmo, \kappa_0, \nu_0, \bso\}$ denotes  the hyper-parameters. 
Expanding the joint density gives:
\begin{equation}
\begin{aligned}
p(\mathcalX, \bmu, \bSigma\mid \bbeta)&=p(\mathcalX \mid \bmu, \bSigma) \cdot p(\bmu, \bSigma \mid \bbeta) \\
&= C\times  \abs{\bSigma}^{- \frac{\nu_0+N+D+2}{2}} \times\\
&\gap \exp \Bigg\{ -\frac{N}{2}(\bmu - \widebarbx)^\top \bSigma^{-1} (\bmu - \widebarbx) - \frac{\kappa_0}{2} (\bmu - \bmo)^\top \bSigma^{-1}(\bmu - \bmo) \\
&\gap -\frac{1}{2} \tr(\bSigma^{-1} \bS_{\widebarx}) - \frac{1}{2} \tr(\bSigma^{-1} \bS_0)  \Bigg\},\\
\end{aligned}
\label{equation:niw_full_joint}
\end{equation}
where $C ={(2\pi)^{-ND/2}}/{Z_{\niw}(D, \kappa_0, \nu_0, \bS_0)}$ is a constant normalization term.
This expression can be rewritten as:
\begin{equation}
\begin{aligned}
p(\mathcalX, \bmu, \bSigma \mid \bbeta)
&=C\abs{\bSigma}^{- \frac{\nu_0+N+D+2}{2}} \times \\
&\gap \exp \Bigg\{-\frac{\kappa_0+N}{2} \left(\bmu - \frac{\kappa_0 \bmo+N \widebarbx}{\kappa_N} \right)^\top \bSigma^{-1} \left(\bmu - \frac{\kappa_0 \bmo+N \widebarbx}{\kappa_N} \right) \\
&\gap - \frac{1}{2} \tr \left[\bSigma^{-1} \left( \bS_0 + \bS_{\widebarx} + \frac{\kappa_0 N}{\kappa_0 + N} (\widebarbx - \bmo)(\widebarbx-\bmo)^\top \right) \right] \Bigg\}, 
\end{aligned}
\label{equation:niw_full_joint2}
\end{equation}
which is reformulated to compare with the NIW form in Equation~\eqref{equation:multi_gaussian_prior}, and we can see the reason why we rewrite the multivariate Gaussian distribution into Equation~\eqref{equation:multi_gaussian_identity} by the trace trick. It follows that the posterior is also a NIW density with updated parameters and confirms the conjugacy of the NIW prior for the multivariate Gaussian model:
\begin{equation}\label{equation:niw_posterior_equation_1}
p(\bmu, \bSigma\mid \mathcalX , \bbeta) = \niw (\bmu, \bSigma \mid \bmm_N, \kappa_N, \nu_N, \bS_N), 
\end{equation}
where 
\begin{align}
\bmm_N &= \frac{\kappa_0\bmo + N\widebarbx}{\kappa_N} = 
\frac{\kappa_0 }{\kappa_N}\bmo+\frac{N}{\kappa_N}\widebarbx , \label{equation:niw_posterior_equation_2}\\
\kappa_N  &= \kappa_0 + N,  \label{equation:niw_posterior_equation_3}\\
\nu_N    &=\nu_0 + N,  \label{equation:niw_posterior_equation_4}\\
\bS_N  &=\bS_0 + \bS_{\widebarx} + \frac{\kappa_0N}{\kappa_0 + N}(\widebarbx - \bmo)(\widebarbx - \bmo)^\top \label{equation:niw_posterior_equation_5}\\
&=\bS_0 + \sum_{n=1}^N \bx_n \bx_n^\top + \kappa_0 \bmo \bmo^\top - \kappa_N \bmm_N \bmm_N^\top . \label{equation:niw_posterior_equation_6}
\end{align}
\paragrapharrow{A ``posterior" interpretation for the NIW prior.}
An intuitive interpretation of the parameters in NIW can be obtained from the updated parameters above. 
The parameter $\nu_0$ acts like a prior sample size for the covariance; and $\nu_N =\nu_0 + N$ is the total (prior + observed) sample size. 
The posterior mean $\bmm_N$ of the model mean $\bmu$ is a weighted average of the prior mean $\bmo$ and the sample mean $\widebarbx$, with weights proportional to $\kappa_0$ and $N$.
The posterior scale matrix $\bS_N$ combines three components: the prior scale matrix $\bS_0$, the empirical covariance matrix $\bS_{\widebarx}$, and an additional term accounting for uncertainty in the mean estimate.

\subsubsection{Practical  Parameter Choice}
In practice, it is often preferable  to use a weakly informative data-dependent prior. A common choice is to set $\bS_0 = \diag(\bS_{\widebarx})/N$ and $\nu_0 =D+2$, ensuring $\Exp[\bSigma]=\bS_0$. Additionally, set  $\bmm_0 =\widebarbx$ and $\kappa_0$ to a small value, such as 0.01, where $\bS_{\widebarx}$ is the sample covariance matrix, and $\widebarbx$ is the sample mean vector as shown in Equation~\eqref{equation:mvu-sample-covariance} \citep{chipman2001practical, fraley2007bayesian, hoff2009first, murphy2012machine}. 
Alternatively, one may first standardize the data so that each feature has zero mean and unit variance. Then set: $\bS_0 = \bm{I}_D$ and $\nu_0 =D+2$ (so $\Exp[\bSigma]=\bm{I}_D$); and  set $\bmm_0 =\bm{0}$ and $\kappa_0$ to a small number, such as 0.01.

\subsubsection{Reducing Sampling Time by Maintaining Squared Sum of Customers}\label{section:reduce-sampling-sum-square}

When implementing NIW in dynamic settings---such as Gibbs sampling for Gaussian mixture models \citep{das2014dpgmm, lu2021survey} or Bayesian matrix factorization with cross-validation (see Section~\ref{section:gggw_model})---it is computationally advantageous to maintain sufficient statistics rather than recompute summaries from scratch after every update.

Note the equivalence between Equation~\eqref{equation:niw_posterior_equation_5} and Equation~\eqref{equation:niw_posterior_equation_6}. 
While the former uses the centered scatter $\bS_{\widebarx}$ and sample mean $\widebarbx$, the latter expresses $\bS_N$ in terms of raw sums:
$$
\bS_N = \bS_0 + \sum_{n=1}^N \bx_n \bx_n^\top + \kappa_0 \bmo \bmo^\top - \kappa_N \bmm_N \bmm_N^\top .
$$
This formulation is particularly useful in clustering contexts (e.g., Chinese restaurant process; see \citet{lu2021survey}). Suppose data points are dynamically added to or removed from a cluster (or “table”). If we used Equation~\eqref{equation:niw_posterior_equation_5}, we would need to recompute $\bS_{\widebarx}$ and  $\widebarbx$ over all current points each time---a costly $O(N)$ operation.
In contrast, Equation~\eqref{equation:niw_posterior_equation_6} depends only on:
the sum of outer products  $\sum_{n} \bx_n \bx_n^\top$ and the sum of vectors $\sum_n\bx_n $ (needed for $\bmm_N$).
Thus, when a point $\bx'$ is added or removed, we simply update these two aggregates by $\pm\bx'\bx'^\top $ and $\pm\bx'$, respectively---each in $O(D^2)$ and $O(D)$ time. This leads to significant computational savings, especially in iterative algorithms like Gibbs sampling.

\subsection{Posterior Marginal Likelihood of Parameters}
The marginal posterior distribution of the covariance matrix $\bSigma$ is obtained by integrating out the mean $\bmu$:
\begin{equation*}
\begin{aligned}
p(\bSigma \mid\mathcalX,\bbeta) 
&= \int_{\bmu} p(\bmu, \bSigma \mid\mathcalX,\bbeta) \,d\bmu =\inversewishart(\bSigma\mid \bS_N, \nu_N),
\end{aligned}
\end{equation*}
where $\bS_N$ and $\nu_N$ are the updated scale matrix and degrees of freedom defined in Equations~\eqref{equation:niw_posterior_equation_5}--\eqref{equation:niw_posterior_equation_4}.
Using standard results for the inverse-Wishart distribution (see Equation~\eqref{equation:iw_expectation}), its posterior mean and mode are:
$$
\begin{aligned}
\Exp[\bSigma \mid \mathcalX, \bbeta] &= \frac{\bS_N}{\nu_N - D-1} 
\qquad\text{and}\qquad
\mathrm{Mode}[\bSigma \mid \mathcalX, \bbeta]=\frac{\bS_N}{\nu_N + D+1}.
\end{aligned}
$$
provided that $\nu_N>D+1$ for the mean to exist.
Similarly, the marginal posterior of the mean vector $\bmu$ follows a multivariate Student's $t$ distribution (Definition~\ref{definition:multivariate-stu-t}).
Specifically,
\begin{equation*}
\begin{aligned}
p(\bmu \mid \mathcalX,\bbeta) 
&= \int_{\bSigma} p(\bmu, \bSigma\mid \mathcalX,\bbeta)\,d\bSigma
= \int_{\bSigma} \niw (\bmu, \bSigma \mid \bmm_N, \kappa_N, \nu_N, \bS_N) \,d\bSigma \\
&=\tau\big(\bmu \mid \bmm_N, \frac{1}{\kappa_N(\nu_N-D+1)}\bS_N, \nu_N-D+1\big),
\end{aligned}
\end{equation*}
which follows from the Gaussian scale mixture representation of the Student's $t$ distribution; see Equation~\eqref{equation:gauss-scale-mixture} and further discussion in \citet{murphy2012machine}.


\subsection{Posterior Marginal Likelihood of Data}\label{section:posterior-marginal-of-data}
The marginal likelihood (also called the evidence) of the observed data $\mathcalX$ under hyper-parameters $\bbeta=\{\bmo, \kappa_0, \nu_0, \bso\}$ is obtained by integrating out both $\bmu$ and $\bSigma$ from the full joint distribution in Equation~\eqref{equation:niw_full_joint2}:
\begin{equation}
\begin{aligned}
p(\mathcalX&\mid\bbeta) 
= \iint p(\mathcalX, \bmu, \bSigma \mid \bbeta) d\bmu d\bSigma 
= \iint \normal(\mathcalX\mid\bmu, \bSigma) \cdot \niw(\bmu, \bSigma\mid\bbeta)  d\bmu d\bSigma \\
&= \frac{(2\pi)^{-ND/2}}{Z_{\niw}(D, \kappa_0, \nu_0, \bS_0)} \int_{\bmu}\int_{\bSigma}\abs{\bSigma}^{-\frac{\nu_0+N+D+2}{2}}  \\
&\gap \times \exp \left(-\frac{\kappa_N}{2} (\bmu - \bmm_N) \bSigma^{-1} (\bmu - \bmm_N ) - \frac{1}{2} \tr(\bSigma^{-1} \bS_N )  \right)d \bmu d \bSigma \\
&\overset{(*)}{=} (2\pi)^{-\frac{ND}{2}} \frac{Z_{\niw}(D, \kappa_N, \nu_N, \bS_N)}{Z_{\niw}(D, \kappa_0, \nu_0, \bS_0)} 
= \pi^{-\frac{ND}{2}} \cdot\frac{\kappa_0^{D/2}\cdot \abs{\bS_0}^{\nu_0/2}}{\kappa_N^{D/2}\cdot \abs{\bS_N}^{\nu_N/2}} \prod_{d=1}^D \frac{\Gamma(\frac{\nu_N+1-d}{2})}{\Gamma(\frac{\nu_0+1-d}{2})}, 
\end{aligned}
\label{equation:niw_marginal_data}
\end{equation}
where the identity ($*$)  follows from the fact that the integral reduces to the normalizing constant of the NIW density given in Equation~\eqref{equation:niw_posterior_equation_1}. 

\subsection{Posterior Predictive for Data without Observations}
Suppose we wish to predict a new data point  $\bx^{\star}$ before observing any data. The prior predictive distribution is:
\begin{equation}
\begin{aligned}
p(\bx^{\star} \mid \bbeta) 
&=  \iint p(\bx^{\star}, \bmu, \bSigma \mid \bbeta) d\bmu d\bSigma 
= \iint \normal(\bx^{\star} \mid \bmu, \bSigma) \cdot \niw(\bmu, \bSigma \mid \bbeta) d\bmu d\bSigma\\
&=  \frac{\pi^{-D/2}\kappa_0^{D/2}  \abs{\bso}^{\nu_0/2}  }{(\kappa_0 + 1) ^{D/2}  \abs{\bS_1}^{\nu_1/2}}  \prod_{d=1}^D \frac{\Gamma(\frac{\nu_1+ 1-d}{2})}{\Gamma(\frac{\nu_0 + 1-d}{2})}
=  \frac{\pi^{-D/2}\kappa_0^{D/2}  \abs{\bso}^{\nu_0/2}  }{(\kappa_0 + 1) ^{D/2} \abs{\bS_1}^{\nu_1/2}}   \frac{\Gamma(\frac{\nu_0+ 2-D}{2})}{\Gamma(\frac{\nu_0 }{2})},
\end{aligned}
\label{equation:niw_prior_predictive_abstract}
\end{equation}
where $\nu_1 = \nu_0+1$, $\bS_1 = \bso+ \frac{\kappa_0 }{\kappa_0+1} (\bx^{\star}-\bmm_0)(\bx^{\star}-\bmm_0)^\top$.
Equivalently, this predictive distribution can be expressed as a multivariate Student's $t$ distribution:
\begin{equation}
p(\bx^{\star} | \bbeta) = \tau\big(\bx^{\star} \mid \bmo, \frac{\kappa_0 + 1}{\kappa_0(\nu_0 - D + 1)}\bS_0, \nu_0 - D + 1\big).
\end{equation}

\subsection{Posterior Predictive for New Data with Observations}
Now suppose we have already observed $\mathcalX$ and wish to predict a new data point $\bx^{\star}$. 
The posterior predictive distribution is:
\begin{equation}
p(\bx^{\star} \mid \mathcalX, \bbeta) = \frac{p(\bx^{\star}, \mathcalX\mid \bbeta) }{p(\mathcalX \mid \bbeta)}. 
\label{equation:niw_posterior_predictive_abstract}
\end{equation}
The denominator of Equation~\eqref{equation:niw_posterior_predictive_abstract} can be obtained directly from Equation~\eqref{equation:niw_marginal_data}. Its numerator can be obtained in a similar way from Equation~\eqref{equation:niw_marginal_data} by considering the marginal likelihood of the new set $\{\mathcalX, \bx^{\star}\}$. 
This amounts to replacing $N$ by $N^{\star}=N+1$ in Equation~\eqref{equation:niw_posterior_equation_2}, Equation~\eqref{equation:niw_posterior_equation_3}, and Equation~\eqref{equation:niw_posterior_equation_4}; and replacing $\bS_N$ by $\bS_{N^{\star}}$  in Equation~\eqref{equation:niw_posterior_equation_5}. 
Thus, the posterior predictive density becomes:
\begin{equation}
\begin{aligned}
p(\bx^{\star} \mid \mathcalX, \bbeta) &= (2\pi)^{-D/2} \frac{Z_{\niw}(D, \kappa_{N^{\star}}, \nu_{N^{\star}}, \bS_{N^{\star}})}{Z_{\niw}(D, \kappa_{N}, \nu_{N}, \bS_{N})} \\
&= \pi^{-D/2} \frac{(\kappa_{N^{\star}})^{-D/2}|\bS_{N} |^{(\nu_N)/2}}{(\kappa_{N})^{-D/2}  |\bS_{N^{\star}} |^{(\nu_{N^{\star}})/2} }  
\prod_{d=1}^D\frac{ \Gamma(\frac{\nu_{N^{\star}} + 1-d}{2})}{ \Gamma(\frac{\nu_{N} + 1-d}{2})} \\
&= \pi^{-D/2} \frac{(\kappa_{N^{\star}})^{-D/2}|\bS_{N} |^{(\nu_N)/2}}{(\kappa_{N})^{-D/2}  |\bS_{N^{\star}} |^{(\nu_{N^{\star}})/2} }  
\frac{ \Gamma(\frac{\nu_{0} + N+2-D}{2})}{ \Gamma(\frac{\nu_{0} + N}{2})} . 
\end{aligned}
\label{equation:niw_posterior_predictive_equation}
\end{equation}
As with the prior predictive, this can also be written in closed form as a multivariate Student's $t$ distribution:
\begin{equation}
p(\bx^{\star} \mid \mathcalX, \bbeta) = \tau \big(\bx^{\star}  \mid  \bmm_N, \frac{\kappa_N + 1}{\kappa_N (\nu_N - D + 1)} \bS_N, \nu_N - D + 1\big).
\end{equation}
Consequently, the predictive mean and covariance are:
\begin{equation*}
\begin{aligned}
\Exp [\bx^\star \mid\mathcalX, \bbeta] &= \bmm_N = \frac{\kappa_0 }{\kappa_0+N}\bmo+\frac{N}{\kappa_0+N}\widebarbx ,\\
\Cov [\bx^\star \mid\mathcalX, \bbeta] &= \frac{\kappa_N + 1}{\kappa_N (\nu_N - D - 1)} \bS_N
= \frac{\kappa_0+N + 1}{(\kappa_0+N) (\nu_0+N - D - 1)} \bS_N,
\end{aligned}
\end{equation*}
where the expectation is a weighted average of the prior mean and the sample mean. 
A few observations are worth noting:
\begin{itemize}
\item  As mentioned previously, $\kappa_0$ controls how strongly we believe this prior for $\bmu$. When $\kappa_0$ is large enough, $\Exp [\bx^\star \mid\mathcalX, \bbeta]$ converges to $\bmm_0$, the prior mean. 
And $\Cov [\bx^\star \mid\mathcalX, \bbeta]$ converges to ${\bS_N}/{(\nu_0+N - D - 1)} $.
\item In the meantime, when $\nu_0$ is large, the prior on $\bSigma$ becomes concentrated around $\bSigma_0=\bS_0/(\nu_0-D-1)$, and the variance behaves approximately as:
$$
\frac{\bS_N}{(\nu_0+N - D - 1)}
\rightarrow \frac{\bS_{\widebarx}}{\nu_0} + \frac{\kappa_0N}{\nu_0(\kappa_0 + N)}(\widebarbx - \bmo)(\widebarbx - \bmo)^\top ,
$$
meaning the posterior predictive uncertainty is increasingly dominated by the observed data rather than the prior hyper-parameters.
\end{itemize}


\subsection{Further Optimization via the Cholesky Decomposition}

\subsubsection*{Definition}
The \textit{Cholesky decomposition} of a symmetric positive definite matrix $\bS$ expresses it as the product of a lower triangular matrix $\bL$ and its transpose:
\begin{equation}
\bS = \bL \bL^\top, 
\end{equation}
where $\bL$ is called the \textit{Cholesky factor} of $\bS$ (see Problem~\ref{problem:cholesky}). 
An equivalent formulation uses the upper triangular matrix $\bR=\bL^\top$, yielding  $\bS = \bR^\top \bR$. 
A triangular matrix is a special type of square matrix. Specifically, a square matrix is \textit{lower triangular} if all entries above the main diagonal are zero, and \textit{upper triangular} if all entries below the main diagonal are zero.

For a matrix of dimension $D$, computing the Cholesky decomposition requires approximately $\sim \frac{1}{3}D^3$ floating-point operations (flops) \citep{lu2021numerical}.
Here, the symbol ``$\sim$" denotes asymptotic equivalence:
\begin{equation*}
	\lim_{D \to +\infty} \frac{\mathrm{number\, of\, flops}}{(1/3)D^3} = 1.
\end{equation*}

\subsubsection*{Rank-One Update}
A rank-one update of a matrix $\bS$  by a vector $\bx$ takes the form \citep{seeger2004low, lu2021numerical}:
\begin{equation*}
\bS^\prime = \bS + \bx \bx^\top. 
\end{equation*}
If the Cholesky factor $\bL$ of $\bS$ is already known, then the Cholesky factor $\bL^\prime$ of $\bS^\prime$ can be computed efficiently.  
Since $\bS^\prime$ differs from $\bS$ only by the rank-one term $\bx\bx^\top$, 
we can obtain $\bL^\prime$ from $\bL$ using a \textit{rank-one Cholesky update}. This operation costs only $\mathcalO(D^2)$ flops---significantly less than the $\mathcalO(D^3)$ required to recompute the decomposition from scratch. See \citet{lu2021numerical} for further details.

\subsubsection*{Speedup for Determinant Computation}
The determinant of a positive definite matrix $\bS$ can be efficiently computed from its Cholesky factor $\bL$:
\begin{equation*}
\abs{\bS} = \prod_{d=1}^{D} l_{dd}^2,\qquad \ln(\abs{\bS}) = 2\ln(\abs{\bL})= 2 \times \sum_{d=1}^D \ln(l_{dd}), 
\end{equation*}
where $l_{dd}$ denotes the ($d,d$)-th entry of  $\bL$; see Problem~\ref{prob:cholesky_det}. 
This computation requires only $\mathcalO(D)$ operations---once the Cholesky decomposition is available, the determinant is simply the product of the squares of the diagonal elements.

\subsubsection*{Update in NIW}
We now consider two closely related computations: the marginal likelihood of observed data in Equation~\eqref{equation:niw_marginal_data} and the posterior predictive distribution for new data in Equation~\eqref{equation:niw_posterior_predictive_abstract}. Both require efficient evaluation of determinants such as   $\abs{\bS_N}$ and $\abs{\bS_{N^{\star}}}$, where $N^*=N+1$.

For example, to compute the posterior predictive density $p(\bx^{\star} \mid \mathcalX, \bbeta)$ in Equation~\eqref{equation:niw_posterior_predictive_abstract}, we evaluate the ratio $ {p(\bx^{\star}, \mathcalX \mid \bbeta) }/{p(\mathcalX \mid \bbeta)}$, which involves the determinants $\abs{\bS_N}$ and $\abs{\bS_{N^{\star}}}$  with $N^{\star} = N+1$.
We handle these efficiently by maintaining and updating the Cholesky decompositions of  $\bS_N$ and $\bS_{N^{\star}}$.
As noted earlier, once the Cholesky factor is available, the determinant is obtained in $\mathcalO(D)$ time.
Expressing $\bS_{N^{\star}}$ in terms of $\bS_N$, we have:
\begin{align}
\bmm_N &= \frac{\kappa_{N^{\star}}\bmm_{N^{\star}} - x^\star}{\kappa_N}=\frac{(\kappa_0 + N + 1)\bmm_{N^{\star}} - x^\star}{\kappa_0 + N} , \\
\bmm_{N^{\star}} &= \frac{\kappa_{N} \bmm_N +\bx^{\star}}{\kappa_{N^{\star}}} 
= \frac{(\kappa_{0}+N) \bmm_N +\bx^{\star}}{\kappa_{0}+N+1},\\
\bS_{N^{\star}} &= \bS_N + \bx^{\star} \bx^{\star T} - \kappa_{N^{\star}} \bmm_{N^{\star}} \bmm_{N^{\star}}^\top + \kappa_N \bmm_N \bmm_N^\top \\
&= \bS_N + \frac{\kappa_0 + N + 1}{\kappa_0 + N}(\bmm_{N^{\star}} - \bx^\star)(\bmm_{N^{\star}} - \bx^\star)^\top, \label{equation:cholesky_rank_1_form}
\end{align}
where Equation~\eqref{equation:cholesky_rank_1_form} implies that Cholesky decomposition of $\bS_{N^\star}$ can be obtained from Cholesky decomposition of $\bS_N$ by a rank-one update. 
Consequently, if the Cholesky decomposition of $\bS_N$ is known, the Cholesky decomposition of  $\bS_{N^\star}$ can be updated in $\mathcalO(D^2)$ time.

\section{Deriving the Dirichlet Distribution*}\label{section:drive-dirichlet}
We derive the Dirichlet distribution and its key properties in this section.
\subsection*{Derivation}
Let $\rx_1, \rx_2, \ldots, \rx_K$ be i.i.d. random variables drawn from the Gamma distribution such that $\rx_k \sim \gammadist(\alpha_k, 1)$ for $k\in \{1, 2, \ldots, K\}$. The joint PDF  of $\rx_1, \rx_2, \ldots, \rx_K$ is given by 
$$ f_{\rx_1,\rx_2,\ldots,\rx_K}(x_1, x_2, \ldots, x_K) =\left\{
\begin{aligned}
& \prod_{k=1}^{K}\frac{1}{\Gamma(\alpha_k)} x_k^{\alpha_k-1} \exp(- x_k) ,& \text{ if all } x_k \geq 0.  \\
&0 , &\mathrm{\,\,if\,\,} \text{otherwise}.
\end{aligned}
\right.
$$
Define new random variables $\ry_1, \ry_2,\ldots,\ry_K$ and $\rz_K$ as follows:
\begin{align}
\ry_k &= \frac{\rx_k}{\sum_{k=1}^{K}\rx_k}, \, \forall\, k \in \{1, 2, \ldots, K-1\};
\qquad 
\ry_K = \frac{\rx_K}{\sum_{k=1}^{K}\rx_k} =1- \sum_{k=1}^{K-1}\ry_k;
\label{equation:dirichlet-define1}\\
\rz_K &= \sum_{k=1}^{K}\rx_k. 
\label{equation:dirichlet-define2}
\end{align}
Let  $\rvx \triangleq [\rx_1, \rx_2, \ldots, \rx_K]^\top$ denote the random vector of variables $\{\rx_k\}$, and define the transformed vector $\rvy \triangleq [\ry_1, \ry_2, \ldots, \ry_{K-1}, \rz_K]^\top$. 
Denote their realizations by $\bx \triangleq [x_1, x_2, \ldots, x_K]^\top$, $\by \triangleq [y_1, y_2, \ldots, y_{K-1}, z_K]^\top$. 
Using the method of  \textit{multidimensional transformation of variables}, the joint PDF of $\rvy$ is
$$
f_{\rvy}(\by) =  f_{\rvx}(g^{-1}(\by))
\abs{\det \left[J_{g^{-1}}(\by)  \right] },
$$
where the inverse transformation $g^{-1}$ is given by
$$
\begin{bmatrix}
x_1\\
x_2\\
\vdots\\
x_K
\end{bmatrix}
=g^{-1}(\by) = g^{-1}\left(\begin{bmatrix}
y_1\\
y_2\\
\vdots\\
z_K
\end{bmatrix}\right)=
\begin{bmatrix}
y_1\cdot z_K\\
y_2\cdot z_K\\
\vdots\\
y_K \cdot z_K
\end{bmatrix},
$$
and the Jacobian matrix is given by 
$$
\begin{aligned}
\small
J_{g^{-1}}(\by)&=
\begin{bmatrix}
\frac{\partial}{\partial y_1}g_1^{-1}(\by) & \cdots 
& \frac{\partial}{\partial z_K}g_1^{-1}(\by) \\
\vdots & \ddots 
& \vdots \\
\frac{\partial}{\partial y_1}g_K^{-1}(\by) & \cdots 
& \frac{\partial}{\partial z_K}g_K^{-1}(\by) \\
\end{bmatrix}
&=
\begin{bmatrix}
z_K & 0& \cdots & 0  &y_1 \\
0 & z_K& \cdots & 0  &y_2 \\
\vdots & \vdots& \ddots & \vdots & \vdots \\
0 & 0& \cdots & z_K  &y_{K-1} \\
-z_K & -z_k &  \cdots & -z_K  &(1-\sum_{k=1}^{K-1}y_k   ) 
\end{bmatrix}.
\end{aligned}
$$
Its determinant evaluates to $\abs{J_{g^{-1}}(\by)}=z_K^{K-1}$.
Therefore, the joint PDF of $\rvy$ becomes 
$$
f_{\rvy}(\by)  = f_{\rvx}(g^{-1}(\by)) z_K^{K-1}
=\frac{y_1^{\alpha_1-1}  y_2^{\alpha_2-1} \ldots y_{K-1}^{\alpha_{K-1}-1}  y_K^{\alpha_K-1}}{\prod_{k=1}^{K} \Gamma(\alpha_k)} \exp(-z_K) z_K^{\alpha_+-1}.
$$
where $\alpha_+=\alpha_1 + \alpha_2+\ldots+\alpha_K$.
We realize that the right-hand side of the above equation is proportional to a PDF of a Gamma distribution and 
$$
\int \exp(-z_K) z_K^{\alpha_+-1} dz_K = \Gamma(\alpha_+).
$$
Integrating out $z_K$, we obtain the marginal PDF of $(y_1, y_2, \ldots, y_{K-1})$ with $\ry_K =1- \sum_{k=1}^{K-1}\ry_k$:
$$
f(y_1, y_2, \ldots, y_{K-1}) = \frac{\Gamma(\alpha_+)}{\prod_{k=1}^{K} \Gamma(\alpha_k)} \prod_{k=1}^{K}y_k^{\alpha_k-1}.
$$
Since each  $\ry_k\in(0,1)$ for all $k\in \{1, 2, \ldots, K\}$, and $\sum_{k=1}^{K}\ry_k=1$, 
this is precisely the PDF of the Dirichlet distribution with parameter vector $\balpha=[\alpha_1, \alpha_2,\ldots,\alpha_K]^\top$ (see Definition~\ref{definition:dirichlet_dist}). 
This construction also provides a practical method for generating Dirichlet-distributed random variables: sample independent Gamma variables and normalize them.

\subsection*{Properties of the Dirichlet Distribution}
Suppose $\rvy=[\ry_1, \ry_2, \ldots, \ry_K]^\top\sim \dirichlet(\balpha)$ with $\balpha=[\alpha_1, \alpha_2, \ldots, \alpha_K]^\top$ and $\alpha_+=\sum_k \alpha_k$.
We now summarize its key properties.
\paragrapharrow{Mean of Dirichlet distribution.}
The expectation of $\ry_1$ (or any $\ry_i$) is:
$$
\begin{aligned}
\Exp[\ry_1] 
&=\int \cdots \int  
y_1 \cdot \dirichlet(\by \mid  \balpha )
\,d y_1 dy_2\cdots dy_K\\
&=\int \cdots \int  
y_1\frac{\Gamma(\alpha_+)}{\prod_{k=1}^{K} \Gamma(\alpha_k)} \prod_{k=1}^{K}y_k^{\alpha_k-1}
\,d y_1 dy_2\cdots dy_K\\
&=\frac{\Gamma(\alpha_+)}{\prod_{k=1}^{K} \Gamma(\alpha_k)} 
\int \cdots \int  
y_1^{\alpha_1+1-1} \prod_{k=2}^{K}y_k^{\alpha_k-1}
\,d y_1 dy_2\cdots dy_K\\
&=\frac{\Gamma(\alpha_+)}{\prod_{k=1}^{K} \Gamma(\alpha_k)}   
\frac{\Gamma(\alpha_1+1) \prod_{k=2}^{K} \Gamma(\alpha_k)}{\Gamma(\alpha_++1)}
=\frac{\Gamma(\alpha_+)}{\Gamma(\alpha_1)}   
\frac{\Gamma(\alpha_1+1) }{\Gamma(\alpha_++1)}
=\frac{\alpha_1}{\alpha_+},
\end{aligned}
$$
where the last equality follows from the fact that $\Gamma(x+1)=x\Gamma(x)$.
Thus, the expectation of any $\ry_i$ is $\Exp[\ry_i] =\frac{\alpha_i}{\alpha_+}$.

\paragrapharrow{Variance of Dirichlet distribution.}
The variance of $\ry_i$ is $\Var[\ry_i] = \Exp[\ry_i^2] - \Exp[\ry_i]^2$. 
Similarly, from the proof of the mean, we have
$$
\Exp[\ry_i^2] = \frac{\Gamma(\alpha_+)}{\Gamma(\alpha_++2)} \frac{\Gamma(\alpha_i+2)}{\Gamma(\alpha_i)} =\frac{(\alpha_i+1)\alpha_i}{(\alpha_++1)\alpha_+}.
$$
This implies 
$$
\Var[\ry_i] = \Exp[\ry_i^2] - \Exp[\ry_i]^2 = \frac{(\alpha_i+1)\alpha_i}{(\alpha_++1)\alpha_+} - (\frac{\alpha_i}{\alpha_+})^2 = \frac{\alpha_i (\alpha_+-\alpha_i)}{\alpha_+^2(\alpha_++1)}.
$$

\paragrapharrow{Covariance of Dirichlet distribution.}
For $i\neq j$, the covariance is  $\Cov[\ry_i \ry_j] = \Exp[\ry_i \ry_j] - \Exp[\ry_i]\Exp[\ry_j]$.
Again, similar to the proof of the mean, for $i\neq j$, we have 
$$
\Exp[\ry_i \ry_j ] = \frac{\Gamma(\alpha_+)}{\Gamma(\alpha_++2)} \frac{\Gamma(\alpha_i+1)}{\Gamma(\alpha_i)}
\frac{\Gamma(\alpha_j+1)}{\Gamma(\alpha_j)} = \frac{\alpha_i \alpha_j}{\alpha_+ (\alpha_++1)}.
$$
This implies
$$
\Cov[\ry_i \ry_j ]=\Exp[\ry_i \ry_j] - \Exp[\ry_i]\Exp[\ry_j]=\frac{\alpha_i \alpha_j}{\alpha_+ (\alpha_++1)} - \frac{\alpha_i \alpha_j}{\alpha_+^2} = \frac{-\alpha_i \alpha_j}{\alpha_+^2 (\alpha_++1) }.
$$

\paragrapharrow{Marginal distribution of $\ry_i$.}
By definitions in Equation~\eqref{equation:dirichlet-define1} and Equation~\eqref{equation:dirichlet-define2}, we have
$\rz_K - \rx_i \sim \gammadist(\alpha_+ - \alpha_i, 1)$.
This implies 
$$
\ry_i = \frac{\rx_i}{\rz_K} =\frac{\rx_i}{\rx_i +(\rz_K-\rx_i)} \sim \betadist(\alpha_i, \alpha_+-\alpha_i).
$$
which follows from the fact about the PDF of two independent Gamma random variables. \footnote{Suppose $\rx\sim \gammadist(a, \lambda)$ and $\ry\sim \gammadist(b, \lambda)$, then $\frac{\rx}{\rx+\ry} \sim \betadist(a, b)$.}

\paragrapharrow{Aggregation property.}
The Dirichlet distribution is closed under summation of components. Specifically, if $[\ry_1, \ry_2, \ldots, \ry_K]^\top\sim \dirichlet([\alpha_1, \alpha_2, \ldots, \alpha_K])$, and we define $M\triangleq \ry_i+\ry_j$, then the vector obtained by replacing $\ry_i$ and $\ry_j$ with $M$ follows a Dirichlet distribution with parameters where $\alpha_i$ and  $\alpha_j$ are replaced by $\alpha_i+\alpha_j$:
$$
\begin{aligned}
&\gap [\ry_1, \ldots \ry_{i-1}, \ry_{i+1}, \ldots, \ry_{j-1}, \ry_{j+1}, \ldots, \ry_K, M] \\
&\sim \dirichlet([\alpha_1, \ldots, \alpha_{i-1}, \alpha_{i+1}, \ldots, \alpha_{j-1}, \alpha_{j+1}, \ldots, \alpha_K, \alpha_i+\alpha_j]).
\end{aligned}
$$
This follows from the fact  that $M\sim \gammadist(\alpha_i+\alpha_j, 1)$ and the multidimensional transformation of variables as shown at the beginning of this section.

The results can be extended to a more general case. If $\{\sA_1, \sA_2, \ldots, \sA_r\}$ is a partition of $\{1, 2, \ldots, K\}$, the aggregated vector satisfies 
$$
\left[\sum_{i\in \sA_1} \ry_i, \sum_{i\in \sA_2} \ry_i, \ldots, \sum_{i\in \sA_r} \ry_i\right] \sim \dirichlet\left(\left[\sum_{i\in \sA_1} \alpha_i, \sum_{i\in \sA_2} \alpha_i, \ldots, \sum_{i\in \sA_r} \alpha_i\right]\right).
$$

\paragrapharrow{Conditional distribution.}
Consider the conditional distribution of a subset of components given the others. For example, let  $\ry_0 \triangleq \sum_{k=3}^{K}\ry_i$ and $\alpha_0\triangleq\alpha_+-\alpha_1-\alpha_2$. 
Then $[\ry_1, \ry_2, \ry_0] \sim \dirichlet([\alpha_1, \alpha_2, \alpha_0])$. Therefore, so the joint density of  ($\ry_1, \ry_2$) is
$$
f_{\ry_1,\ry_2}(y_1,y_2)=
\frac{\Gamma(\alpha_1+\alpha_2+\alpha_0)}{\Gamma(\alpha_1)\Gamma(\alpha_2)\Gamma(\alpha_0)}
y_1^{\alpha_1-1}y_2^{\alpha_2-1} (1-y_1-y_2)^{\alpha_0-1}.
$$
Similarly, the marginal of $\ry_2$ is:
$$
f_{\ry_2}(y_2) 
=  \frac{\Gamma(\alpha_1+\alpha_2+\alpha_0)}{\Gamma(\alpha_2)\Gamma(\alpha_1+\alpha_0)}
y_2^{\alpha_2-1} (1-y_2)^{\alpha_1+\alpha_0-1}= \betadist(y_1\mid  \alpha_2, \alpha_1+\alpha_0),
$$
which is a PDF of a Beta distribution. Therefore, the conditional PDF of $\ry_1\mid \ry_2 = y_2$ is given by 
$$
f_{\ry_1\mid \ry_2 = y_2}(y_1\mid y_2) = \frac{f_{\ry_1,\ry_2}(y_1,y_2)}{f_{\ry_2}(y_2)} 
=\frac{\Gamma(\alpha_1+\alpha_0)}{\Gamma(\alpha_1)\Gamma(\alpha_0)}
\left(\frac{y_1}{1-y_2}\right)^{\alpha_1-1} \left(1-\frac{y_1}{1-y_2}\right)^{\alpha_0-1} \frac{1}{1-y_2},
$$
which implies 
$$
\frac{1}{1-y_2} \ry_1\mid \ry_2=y_2 \sim \betadist(\alpha_1, \alpha_0).
$$
More generally, for any $i$, the conditional distribution of the remaining components given $\ry_i=y_i$ satisfies
$$
\rvy_{-i} \mid  \ry_i\sim (1-y_i)\dirichlet(\balpha_{-i}),
$$
where $\rvy_{-i}$ denotes the vector of all components except $\ry_i$, and $\balpha_{-i}$ is the corresponding parameter vector.

\begin{problemset}
\item \textbf{Chernoff bound for centered Gaussian.}
Let $\rx \sim \normal(0, \sigma^2)$. Show that
$$
p(\abs{\rx} \geq t) \leq 2 e^{-\frac{t^2}{2\sigma^2}}.
$$
\textit{Hint: Apply the Chernoff bound using the moment-generating function of the Gaussian distribution.}

\item \textbf{Bernoulli model.}
Given a data set of binary variables $\mathcalX=\{x_1, x_2, \ldots,x_N\}$, i.e., each $x_n\in\{0,1\}$, and consider the Bernoulli model for these variables, $\rx_n\sim \bernoulli(\theta)$.
Derive the log-likelihood of these observations $\ln p(\mathcalX\mid \theta)$. Show that the maximum likelihood estimate of $\theta$ is the sample mean:
$\theta_{\text{ML}} = (\sum_{n=1}^{N} x_n)/N$.

\item \label{problem:cholesky} \textbf{Cholesky decomposition.} 
Prove the Cholesky decomposition:
Every \textit{positive definite} (PD) matrix $\bA\in \real^{D\times D}$ can be factored as 
$$
\bA = \bR^\top\bR,
$$
where $\bR \in \real^{D\times D}$ is an upper triangular matrix \textbf{with positive diagonal elements}. This decomposition is known as the \textit{Cholesky decomposition}  of $\bA$, and $\bR$ is called the \textit{Cholesky factor} or \textit{Cholesky triangle} of $\bA$.
Specifically, the Cholesky decomposition is \textbf{unique}.
In cases where the diagonal elements of $\bR$ are not restricted to positive values, then the factorization $\bA=\bR^\top\bR$ is \textbf{not unique}.

\item \label{prob:cholesky_det} Let $\bS=\bL\bL^\top\in\real^{D\times D}$ be the Cholesky decomposition of the positive definite matrix $\bS$. Show that the determinant of $\sS$ is given by:  $\abs{\bS} = \prod_{d=1}^{D} l_{dd}^2$.

\item \textbf{Poisson and conjugacy.} Let $\rx_1, \rx_2, \ldots, \rx_N$ be i.i.d. random variables drawn from the Poisson distribution $\poissondist(\lambda)$. Suppose the prior for $\lambda$ is 
$$
\gammadist(\lambda \mid a, b) = \frac{b^a}{\Gamma(a)} \lambda^{a-1} \exp(-b \lambda) \indicator(\lambda >0).
$$
Derive the posterior distribution of $\lambda$.

\item \label{problem:multiGauss} 
Assume you can generate independent samples from the standard univariate Gaussian $\normal(0, 1)$. Show how to generate a sample from the multivariate Gaussian distribution $\normal(\bmu, \bSigma)$, where $\bSigma=\bC\bC^\top$ for some matrix $\bC$ (e.g., via Cholesky decomposition), $\bmu\in\real^N$, and $\bSigma\in\real^{N\times N}$. \textit{Hint: if $\rx_1, \rx_2, \ldots, \rx_N$ are i.i.d. from $\normal(0, 1)$, and let $\rvx=[\rx_1, \rx_2, \ldots, \rx_N]^\top$, then it follows that $\bC\rvx +\bmu\sim \normal(\bmu, \bSigma)$}.

\item \textbf{Maximum entropy.} 
The entropy of a distribution $p(\bx)$ is given in Problem~\ref{problem:entropy_mgau}.
Among all probability distributions $p(\bx)$ with a fixed mean $\bmu$ and covariance matrix $\bSigma$, show that  the Gaussian distribution $\normal(\bmu, \bSigma)$ maximizes entropy.
That is, we wish to maximize $\entropy[p(\bx)]$ over all distributions $p(\bx)$ subject to the constraints that $p(\bx)$ is normalized $\int p(\bx) \, d\bx = 1$ and that it has a specific mean and covariance, so that
\begin{align*}
\int p(\bx) \bx \, d\bx = \bmu 
\qquad\text{and}\qquad
\int p(\bx) (\bx - \bmu)(\bx - \bmu)^\top \, d\bx = \bSigma. 
\end{align*}
\textit{Hint: Use Lagrange multipliers to enforce the constraints}

\item Prove the marginal and conditional distributions of a multivariate Gaussian distribution in \eqref{equation:marginal_multigaus} rigorously.

\item Prove Lemma~\ref{lemma:affine_mult_gauss} (affine transformations of Gaussians are Gaussian)  and Lemma~\ref{lemma:rotat_multi_gauss} (rotations preserve the Gaussian form).

\item Following \eqref{equation:mt_inv}, consider a nonsingular matrix $\bA\in\real^{N\times N}$, an index set $\sI$ and its complement $\sJ=\{1,2,\ldots,N\}\setminus \sI$. Show that 
$$
\begin{aligned}
\bA^{-1}[\sI, \sI] &= \left(\bA[\sI,\sI] - \bA[\sI, \sJ]\bA[\sJ,\sJ]^{-1}\bA[\sJ,\sI] \right)^{-1};\\
\bA^{-1}[\sI,\sJ] &= \bA[\sI,\sI]^{-1} \bA[\sI,\sJ] \left(\bA[\sJ,\sI]\bA[\sI,\sI]^{-1}\bA[\sI,\sJ]-\bA[\sJ,\sJ]\right)^{-1}\\
&= \left(\bA[\sI,\sJ]\bA[\sJ,\sJ]^{-1}\bA[\sJ,\sI]-\bA[\sI,\sI]\right)^{-1} \bA[\sI,\sJ]\bA[\sJ,\sJ]^{-1},
\end{aligned}
$$
where $\bA^{-1}[\sI,\sJ]$ denotes the submatrix of $\bA^{-1}$, and $\bA[\sI,\sI]^{-1}$ denotes the inverse of $\bA[\sI,\sI]$.

\item \label{prob:part_inv} \textbf{Partitioned inverse.} 
Verify the identity \eqref{equation:mt_inv} by multiplying both sides by the matrix 
$\begin{bmatrixfoot}
\bA & \bB\\
\bC& \bD 
\end{bmatrixfoot}$.

\end{problemset}

\newpage 
\part{Non-Bayesian Matrix Decomposition}

\newpage
\chapter{Alternating Least Squares (ALS)}\label{chapter:als}
\begingroup
\hypersetup{
	linkcolor=structurecolor,
	linktoc=page,  
}
\minitoc \newpage
\endgroup

\index{Least squares}
\index{Linear models}
\index{Regression analysis}
\section{Preliminary: Least Squares Approximations}\label{section:pre_ls}
\lettrine{\color{caligraphcolor}T}
The linear model is the cornerstone of regression analysis, with the least squares approximation serving as its fundamental tool for minimizing the sum of squared errors.
This approach is a natural choice when seeking the regression function that minimizes the expected squared prediction error.
Over the past few decades, linear models have been widely applied across diverse fields, including decision-making \citep{dawes1974linear}, time series analysis \citep{christensen1991linear, lu2017machine}, quantitative finance \citep{menchero2011barra}, and numerous other disciplines such as production science, social science, and soil science \citep{fox1997applied, lane2002generalized, schaeffer2004application, mrode2014linear}.

To be more concrete, consider an overdetermined system represented by $\bb = \bA\bx $, where $\bA\in \real^{M\times N}$ is the \textit{input data matrix} (also called the \textit{predictor variables}), $\bb\in \real^M$ is the \textit{observation vector} (or \textit{target/response vector}), and the number of samples $M$ exceeds the number of features $N$. 
The vector $\bx\in\real^N$ represents the \textit{weights} (\textit{coefficients}) of the linear model.
In practice, $\bA$ typically has full column rank, as real-world predictor variables are often uncorrelated---or can be made so through preprocessing.
Moreover, a  \textit{bias term} (also known as an \textit{intercept}) is commonly introduced by augmenting $\bA$  with a column of ones. This leads to the modified system:
\begin{equation}\label{equation:ls-bias}
	\widetildebA \widetildebx = 
[\bm{1} ,\bA ] 
\begin{bmatrix}
	x_0\\
	\bx
\end{bmatrix}
 = \bb .
\end{equation}
where $x_0$ is the intercept coefficient.

However, because the system is overdetermined (i.e., there are more equations than unknowns), the equation 
$\bb=\bA\bx$ often has no exact solution---it is \textit{inconsistent}.
Let the column space of $\bA$ be denoted by $\cspace(\bA)=\{\bA\bgamma\mid  \forall \bgamma \in \real^N\}$.
When  $\bb \notin  \cspace(\bA)$, the residual error $\be=\bb-\bA\bx$ cannot be reduced to zero.
In other words, the error $\be = \bb -\bA\bx$ cannot be reduced to zero. 
In such cases, the goal becomes minimizing this error---typically measured by the mean squared error (MSE).
The resulting solution $\bx_{\text{LS}}$ that minimizes $\normtwo{\bb-\bA\bx}^2$ is known as the \textit{least squares (LS) solution} or the \textit{ordinary least squares (OLS) solution}. The least squares method is a foundational technique in the mathematical sciences, and entire textbooks are devoted to it (e.g., \citet{trefethen1997numerical, strang2019linear, strang2021every, lu2021rigorous}).

\paragrapharrow{Least squares via calculus.}
Assume that the objective function $\normtwo{\bb-\bA\bx}^2$ is differentiable and that the parameter space for $\bx$ is the entire $\real^N$ (i.e., an unconstrained optimization problem); that is, the domain of $\mathop{\min}_{\bx} \normtwo{\bb-\bA\bx}^2$ is  $\real^N$. 
Then, the least squares estimate corresponds to the point where the gradient of the objective function vanishes. This leads to the following lemma.~\footnote{Variants of the least squares problem are explored in Problems~\ref{problem:rls}--\ref{problem:twls2}.}

\index{Fermat's theorem}
\index{Normal equation}
\begin{lemma}[Least Squares via  Calculus]\label{lemma:ols}
Let  $\bA \in \real^{M\times N}$ be a  fixed data matrix with full column rank and  $M\geq N$ (i.e., its columns  are linearly independent)
\footnote{See Problems~\ref{prob:als_pseudo1}--\ref{prob:als_pseudon} for a relaxation using the pseudo-inverse.}.
For the overdetermined system $\bb = \bA\bx$, the least squares solution---obtained by setting the gradient of $\normtwo{\bb-\bA\bx}^2$ to  zero (i.e., the gradient vanishes)---is given by $\bx_{\text{LS}} = (\bA^\top\bA)^{-1}\bA^\top\bb$ 
This is known as the \textit{first-order optimality condition} for  local optima. 
The proof relies on \textit{Fermat's theorem} for multivariate functions, which itself follows from the univariate case.
\footnote{See Problem~\ref{problem:fist_opt}.}
The value, $\bx_{\text{LS}} = (\bA^\top\bA)^{-1}\bA^\top\bb$, is commonly referred to as the \textit{ordinary least squares (OLS)} estimate or simply the \textit{least squares (LS)} estimate of $\bx$.
\end{lemma}

To prove this lemma, we must confirm that $\bA^\top\bA$ is invertible. 
Under the assumption that $\bA$ has full rank and $M\geq N$, the matrix $\bA^\top\bA \in \real^{N\times N}$ is symmetric positive definite and thus invertible (see Problem~\ref{prob:rank-of-ata}).
\begin{proof}[of Lemma \ref{lemma:ols}]
From calculus, a differentiable function $f(\bx)$ attains a minimum at $\bx_{\text{LS}}$ only if $\nabla f(\bx)=\bzero$. The gradient of $\normtwo{\bb-\bA\bx}^2$ is $2\bA^\top\bA\bx -2\bA^\top\bb$. 
$\bA^\top\bA$ is invertible since we assume that $\bA$ is fixed and has full rank with $M\geq N$ (Problem~\ref{prob:rank-of-ata}). 
Consequently, the OLS solution for $\bx$ is $\bx_{\text{LS}} = (\bA^\top\bA)^{-1}\bA^\top\bb$, which completes the proof.
\end{proof}

The equation $\bA^\top\bA \bx = \bA^\top\bb$ is called the \textit{normal equation}. 
Under the assumption that $\bA$ has full rank with $M\geq N$,  $\bA^\top\bA$ is invertible, and the least squares solution is uniquely given by $\bx_{\text{LS}} = (\bA^\top\bA)^{-1}\bA^\top\bb$.

\index{Convex functions}
\begin{figure}[h!]
\centering  
\vspace{-0.35cm} 
\subfigtopskip=2pt 
\subfigbottomskip=2pt 
\subfigcapskip=-5pt 
\subfigure[A convex function.]{\label{fig:convex-1}
\includegraphics[width=0.26\linewidth]{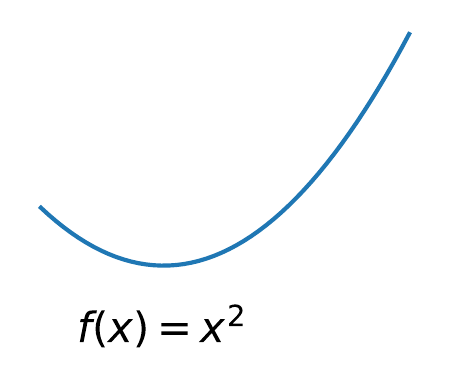}}
\subfigure[A concave function.]{\label{fig:convex-2}
\includegraphics[width=0.26\linewidth]{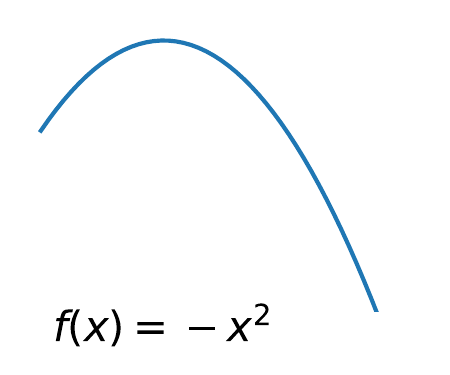}}
\subfigure[A random function.]{\label{fig:convex-3}
\includegraphics[width=0.26\linewidth]{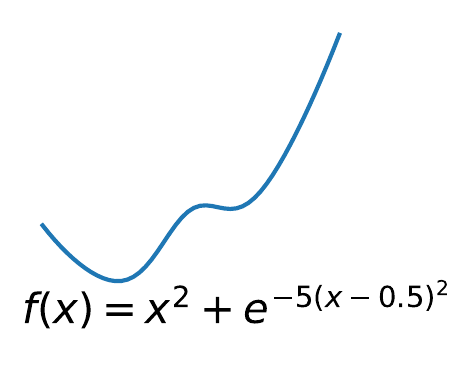}}
\caption{Three functions.}
\label{fig:convex-concave-none}
\end{figure}

However, a vanishing gradient alone does not guarantee that the solution is a global minimum---it could correspond to a maximum or a \textit{saddle point}.
\footnote{
A \textit{saddle point} is a point at which the gradient vanishes (a \textit{stationary point}), and  there exists a direction where the objective function decreases and another direction where it increases.
}
Figure~\ref{fig:convex-concave-none} illustrates this ambiguity.
What we can assert is that any local minimum must satisfy the first-order condition (zero gradient), but this condition is necessary---not sufficient---without additional assumptions.
The following remark clarifies why the OLS solution indeed minimizes the squared error.

\begin{remark}[Verification of Least Squares Solution\index{Saddle point}]
Why does a zero gradient imply minimal mean squared error?
While convexity provides a standard explanation (as we will see shortly), we can directly verify minimality.
For any $\bx \neq \bx_{\text{LS}}$, expand the squared norm: 
\begin{align*}
\normtwo{\bb - \bA\bx}^2 &= \normtwo{\bb - \bA\bx_{\text{LS}} + \bA\bx_{\text{LS}} - \bA\bx}^2
= \normtwo{\bb-\bA\bx_{\text{LS}} + \bA (\bx_{\text{LS}} - \bx)}^2 \\
&=\normtwo{\bb-\bA\bx_{\text{LS}}}^2 + \normtwo{\bA(\bx_{\text{LS}} - \bx)}^2 + 2\left(\bA(\bx_{\text{LS}} - \bx)\right)^\top(\bb-\bA\bx_{\text{LS}}) \\ 
&=\normtwo{\bb-\bA\bx_{\text{LS}}}^2 + \normtwo{\bA(\bx_{\text{LS}} - \bx)}^2 + 2(\bx_{\text{LS}} - \bx)^\top(\bA^\top\bb - \bA^\top\bA\bx_{\text{LS}}). 
\end{align*} 
The cross term vanishes due to the normal equation, and the second term is nonnegative. Hence,
$
\normtwo{\bb - \bA\bx}^2 \geq \normtwo{\bb-\bA\bx_{\text{LS}}}^2$, with equality only when $\bx=\bx_{\text{LS}}$.
Thus, the OLS estimate yields the global minimum, not a maximum or saddle point.
Indeed, this condition from the least squares estimate is also known as the \textit{sufficiency of stationarity under convexity}. When $\bx$ is defined over the entire space $\real^N$, this condition is also known as the \textit{necessity of stationarity under convexity}.
\end{remark}

One might wonder: Why does left-multiplying the system by $\bA^\top$---as in the normal equation---``magically" produce a solvable system?
Consider a simple analogy: the equation $x^2=-1$ has no real solution, but multiplying both sides by $x$ yields $x^3=-x$, which does have a real solution ($x=0$)---the value that makes $x^2$ as close as possible to 
$-1$ in the least squares sense.

\begin{example}[Altering the Solution Set by Left Multiplication]
Consider the  data matrix and target vector:
$\tiny
\bA=
\begin{bmatrix}
-3 & -4 \\
4 & 6  \\
1 & 1
\end{bmatrix}
$
and
$
\bb=
\tiny
\begin{bmatrix}
1  \\
-1   \\
0
\end{bmatrix}
.
$
It is straightforward to verify that $\bA\bx = \bb$ has no exact solution.
However, if we left-multiply both sides by
$
\bB=
\scriptsize
\begin{bmatrix}
0 & -1 & 6\\
0 & 1  & -4
\end{bmatrix},
$
then  $\bx_{\text{LS}} = [1/2, -1/2]^\top$ solves  $\bB\bA\bx= \bB\bb$. 
This illustrates how the normal equation (which uses $\bB=\bA^\top$) transforms the original inconsistent system into a consistent one in a lower-dimensional space---effectively projecting $\bb$ onto $\cspace(\bA)$ and yielding the least squares solution.
\end{example}

\paragrapharrow{Rank-deficiency.}
Our discussion so far assumes $\bA\in \real^{M\times N}$ has full column rank with $M\geq N$, ensuring  $\bA^\top\bA$ is invertible. 
However, if two or more columns of $\bA$ are perfectly correlated,  $\bA$ becomes deficient, and $\bA^\top\bA$ is singular.
In such cases, infinitely many least squares solutions may exist. A common strategy is to select the solution with the smallest Euclidean norm.
This leads naturally to the Moore--Penrose pseudo-inverse. See Problems~\ref{prob:als_pseudo1}--\ref{prob:als_pseudon} or the following paragraph for further details.

\index{Contion number}
\index{Tikhonov regularization}
\index{$\ell_2$-regularization}
\paragrapharrow{Regularizations and stability.}
Even when $\bA$ is full rank, the ordinary least squares solution can be numerically unstable if $\bA$ is \textit{nearly} singular.
Let the SVD of $\bA$ be $\bA=\bU\bSigma\bV^\top\in\real^{M\times N}$, where $\bU\in\real^{M\times M}$ and $\bV\in\real^{N\times N}$ are orthogonal, and  $\bSigma\in\real^{M\times N}$ contains the singular values $\sigma_1\geq\sigma_2\ldots\geq\sigma_N\geq 0$ on its main diagonal. 
Consequently, $\bA^\top\bA = \bV(\bSigma^\top\bSigma)\bV^\top \triangleq \bV\bS\bV^\top$, where $\bS\triangleq \bSigma^\top\bSigma  = \diag(\sigma_1^2, \sigma_2^2, \ldots,\sigma_N^2)\in\real^{N\times N}$ contains the squared singular values of $\bA$. When $\bA$ is nearly singular, $\sigma_N^2\approx 0$, making the inverse operation $(\bA^\top\bA)^{-1} = \bV\bS^{-1}\bV^\top$ numerically unstable and ill-conditioned. 
As a result, the solution $\bx_{\text{LS}} =(\bA^\top\bA)^{-1}\bA^\top\bb $ may diverge (small perturbations in $\bb$ can cause large changes in $\bx_{\text{LS}}$).
To mitigate this, we introduce \textit{$\ell_2$-regularization} (also known as \textit{Tikhonov regularization}  \citep{tikhonov1963solution}), which solves:
\begin{equation}
\bx_{\text{Tik}} = \mathop{\argmin}_{\bx} \normtwo{\bb-\bA\bx}^2 +\lambda\normtwo{\bx}^2, 
\quad \text{with }\lambda>0.
\end{equation}
The gradient is $2(\bA^\top\bA+\lambda\bI)\bx-2\bA^\top\bb$,  yielding the regularized solution:
$$
\bx_{\text{Tik}} = (\bA^\top\bA+\lambda\bI)^{-1}\bA^\top\bb.
$$
The inverse operation becomes $(\bA^\top\bA+\lambda\bI)^{-1} = \bV(\bS+\lambda\bI)^{-1}\bV^\top$, where $\widetildebS\triangleq(\bS+\lambda\bI)=\diag(\sigma_1^2+\lambda, \sigma_2^2+\lambda, \ldots,\sigma_N^2+\lambda)$. 
The solutions for OLS and Tikhonov regularized LS are given, respectively, by 
\begin{equation}\label{equation:ols_tik}
\begin{aligned}
\bx_{\text{LS}} &= (\bA^\top\bA)^{-1}\bA^\top\bb = \bV\left(\bS^{-1}\bSigma\right)\bU^\top\bb;\\
\bx_{\text{Tik}} &= (\bA^\top\bA+\lambda\bI)^{-1}\bA^\top\bb = \bV\left((\bS+\lambda\bI)^{-1}\bSigma\right)\bU^\top\bb,\\
\end{aligned}
\end{equation}
where the main diagonals of $\left(\bS^{-1}\bSigma\right)$ are $\diag(\frac{1}{\sigma_1}, \frac{1}{\sigma_2}, \ldots, \frac{1}{\sigma_N})$; and the main diagonals of $\left((\bS+\lambda\bI)^{-1}\bSigma\right)$ are $\diag(\frac{\sigma_1}{\sigma_1^2+\lambda}, \frac{\sigma_2}{\sigma_2^2+\lambda}, \ldots, \frac{\sigma_N}{\sigma_N^2+\lambda})$. The latter solution is more stable if $\lambda$ is greater than the   smallest nonzero squared singular value.
The \textit{condition number} becomes smaller if  the smallest singular value $\sigma_N$ is close to zero \citep{lu2021numerical}:
$$
\kappa(\bA^\top\bA) = \frac{\sigma_1^2}{\sigma_N^2}
\qquad \rightarrow \qquad
\kappa(\bA^\top\bA+\lambda\bI) = \frac{\lambda+\sigma_1^2}{\lambda+\sigma_N^2},
$$
which is closer to 1 when $\lambda \gg \sigma_1^2$.
Thus, Tikhonov regularization prevents divergence in nearly singular or rank-deficient settings, enhances algorithmic convergence (e.g., in alternating least squares), and resolves identifiability issues. It is now a standard tool in practice.

\index{Decomposition: ALS}
\index{Netflix}
\section{Netflix Recommender and Matrix Factorization}\label{section:als-netflix}
The rapid growth of data driven by advances in sensor technology and computing hardware has introduced new challenges in data analysis.
Large-scale datasets often contain noise and other distortions, necessitating preprocessing before deductive scientific methods can be effectively applied.
For instance, signals captured by antenna arrays are frequently corrupted by noise and other degradations.
To analyze such data meaningfully, it must be reconstructed or represented in a way that reduces inaccuracies while preserving essential structural or feasibility constraints.

Moreover, in many real-world scenarios, data collected from complex systems arises from the joint influence of multiple interrelated variables. When these variables are poorly defined or entangled, the information in the raw data becomes redundant and ambiguous.
By constructing a reduced-order model, we can approximate the original system with high fidelity.
A common strategy for denoising, model reduction, data compression, and feasibility-preserving reconstruction is to replace the original data with a lower-dimensional representation obtained via subspace approximation.
Consequently, \textit{low-rank matrix approximations (LRMA)}---or \textit{low-rank matrix decompositions}---play a pivotal role across numerous applications, including data compression, feature selection, and noise filtering.
\footnote{Strictly speaking, ``approximation" typically refers to a scenario where a matrix $\bA$ is expressed as $\bA\approx\bW\bZ$, with $\bW\bZ$ being a close but not exact estimate of $\bA$. In contrast, “decomposition” usually implies an exact factorization: $\bA=\bW\bZ$. However, in this context, we use the terms interchangeably, acknowledging that both exact and approximate factorizations may be discussed under either label.}

Low-rank matrix decomposition is a powerful technique widely used in machine learning and data mining to express a given matrix as the product of two (or more) lower-dimensional matrices.
It captures the essential structure of the data while discarding noise and redundancies.
Common approaches include singular value decomposition (SVD), principal component analysis (PCA), nonnegative matrix factorization (NMF) with multiplicative updates, and the alternating least squares (ALS) method, which we introduce in this section.

\index{Collaborative filtering}
\subsection*{Example: The Netflix Prize}
In the Netflix Prize competition \citep{bennett2007netflix}, the goal was to predict user ratings for movies based on their existing ratings for other movies (i.e., observed user-item interactions).
Let there be $M$ movies indexed by $m= 1, 2,\ldots,M$ and  $N$ users indexed by $n = 1, 2,\ldots,N$. 
(Throughout this discussion, lowercase letters such as $m,n,k$ denote running indices, while uppercase letters $M, N, K$ represent their respective upper bounds.)
Denote the rating of the $n$-th user for the $m$-th movie by $a_{mn}$, and define the rating matrix  $\bA\in\real^{M \times N}$  (also called a \textit{movie-by-user matrix} or \textit{preference matrix}) whose  $n$-th column $\ba_n \in \real^M$ contains all ratings provided by user $n$.
Crucially, most entries $\{a_{mn}\}$ are missing. 
The task is to accurately predict these unobserved ratings---i.e., to complete the matrix.

Without some underlying structure in $\bA$, there would be no meaningful relationship between observed and missing entries, rendering the completion problem ill-posed and admitting infinitely many solutions.
Thus, it is essential to impose a structural assumption. A widely adopted one is that $\bA$ is approximately low-rank: the ratings arise from a small number of latent factors (e.g., genres, user preferences), implying strong correlations among rows and columns.
This low-rank assumption makes the matrix completion problem well-posed and enables a unique, data-consistent solution. Under this model, unobserved entries are no longer independent of observed ones---they are linked through the shared low-dimensional latent space.
\footnote{This is, however, a strong assumption. Consider a rank-$R$ matrix $\bA=\sum_{r=1}^{R} \be_r\widetilde{\be}_r^\top$, where $\be_i$  and $\widetilde{\be}_r$ are standard basis vectors in $\real^M$ and $\real^N$, respectively. Such a matrix has only $R$ nonzero entries. In recommendation settings, where only a few random entries are observed, it is highly likely that some true nonzero entries remain unseen---posing significant challenges for recovery. We do not address this issue further in this book.}

It is important to note that, except in highly structured cases, any low-rank representation of $\bA$ necessarily incurs some compression error, as it is only an approximation of the original matrix.
This approach---known as \textit{collaborative filtering}---exploits recurring patterns in observed user behaviors to predict future preferences.

\subsection*{Formal Problem Statement}
Let  $\bM\in \{0,1\}^{M\times N}$ be a \textit{mask matrix}, where $m_{mn}=1$ if user $n$ has rated movie $m$, and $0$ otherwise.
The low-rank matrix completion problem can then be formulated as:
\begin{equation}
\widetildebA = \mathop{\argmin}_{\bX\in\real^{M\times N}} \sum_{m,n=1}^{M,N} (x_{mn} - a_{mn})^2\cdot m_{mn}\gap \text{s.t.} \gap \rank(\bX)\leq K.
\end{equation}
However, this problem is NP-hard  \citep{hardt2014computational}. 
An equivalent unconstrained formulation---derived via the singular value decomposition---is:
\begin{equation}
\widetildebA =\widetildebW\widetildebZ = \mathop{\argmin}_{\substack{\bW\in\real^{M\times K}\\ \bZ\in\real^{K\times N}}} \sum_{m,n=1}^{M,N} ((\bW\bZ)_{mn}- a_{mn})^2\cdot m_{mn}.
\end{equation}
This form enables practical approximation via iterative algorithms.

We now consider the general factorization problem: approximate $\bA$ as  $\bA\approx\bW\bZ$, where $\bW\in\real^{M\times K}$ and  $\bZ\in \real^{K \times  N}$. 
Typically, $K\ll\min\{M,N\}$, ensuring dimensional reduction and data compression.
The choice of $K$ is critical in practice and is usually problem-dependent (and is often selected via \textit{cross-validation (CV)}).
If we partition  $\bA=[\ba_1, \ba_2, \ldots, \ba_N]$ and $\bZ=[\bz_1, \bz_2, \ldots, \bz_N]$ into columns, then $\ba_n \approx \bW\bz_n$. 
Thus, each observed rating vector $\ba_n$ is approximated as a linear combination of the columns of $\bW$, with coefficients given by $\bz_n$. 
In this view, the columns of $\bW$ serve as a basis (or \textit{template columns}) spanning the column space of $\bA$, while $\bz_n$ encodes the coordinates (or \textit{activations}) of $\ba_n$ in that basis.

\index{Marginal convexity}
\index{Two-block coordinate descent}
\begin{algorithm}[H] 
\caption{2-Block Coordinate Descent: Framework of Most ALS and NMF Algorithms}
\label{alg:two_bcd_gen_inals}
\begin{algorithmic}[1] 
\Require A loss function for a variable with two blocks $\bX=(\bW,\bZ)$: $f(\bX)=f(\bW,\bZ)$, and data matrix $\bA$;
\Ensure Constraint on $\bW$ and $\bZ$;
\State Generate some initial matrices $\bW^{(0)}$ and $\bZ^{(0)}$;
\For{$t = 1, 2, \ldots$}
\State $\bW^\toptzero \leftarrow \text{update}\big(\bA, \bZ^{(t-1)}, \bW^{(t-1)}\big)$;
\State $\bZ^\toptzero \leftarrow \text{update}\big(\bA, \bW^\toptzero, \bZ^{(t-1)}\big)$;
\EndFor
\end{algorithmic} 
\end{algorithm}
In most cases, the factorization problem admits no closed-form solution and must be solved numerically.
A standard approach is \textit{two-block coordinate descent (2-BCD)}, as outlined in Algorithm~\ref{alg:two_bcd_gen_inals}.
To simplify the problem, let's first assume that there are no missing ratings. 
To quantify the quality of the approximation $\bA\approx\bW\bZ$, we adopt the Frobenius norm as our loss function (assume $K$ is known):
\begin{equation}\label{equation:als-per-example-loss2}
L(\bW,\bZ)\triangleq D(\bA,\bW\bZ) \triangleq \frac{1}{2}\sum_{n=1}^N \sum_{m=1}^{M} \left(a_{mn} - \bw_m^\top\bz_n\right)^2 
=\frac{1}{2} \normf{\bW\bZ-\bA}^2.~\footnote{
Note that we include a scaling factor of $\frac{1}{2}$ for easier discussion of gradients. Minimizing over $\frac{1}{2}\normf{\bW\bZ-\bA}^2$ is equivalent to minimizing over $\normf{\bW\bZ-\bA}^2$ or $\normf{\bW\bZ-\bA}$.
The choice of the Frobenius norm assumes i.i.d. Gaussian noise on the data ($\bA=\bW\bZ+\bN$, where each entry of $\bN$ follows i.i.d. Gaussian noise) and leads to a smooth optimization via least squares. When the loss is measured by the $\ell_1$ matrix norm, one obtains a robust low-rank matrix factorization; and the noise is assumed i.i.d. Laplace.}
\end{equation}
Here, $\bW=[\bw_1^\top; \bw_2^\top; \ldots; \bw_M^\top]\in \real^{M\times K}$ (rows are $\bw_m^\top$) and $\bZ=[\bz_1, \bz_2, \ldots, \bz_N] \in \real^{K\times N}$  (columns are $\bz_n$).
In \eqref{equation:als-per-example-loss2},  $L(\bW,\bZ)$ indicates it is a loss function w.r.t. $\bW$ and $\bZ$, and $D(\bA, \bW\bZ)$ implies it is a distance/divergence~\footnote{In words, the \textit{distance} $D(\bE,\bF)$ indicates $D(\bE,\bF)=D(\bF,\bE)\geq 0$ and the equality holds if and only if $\bE=\bF$; while the \textit{divergence} holds that  $D(\bE,\bF)\neq D(\bF,\bE)\geq 0$ and the equality holds if and only if $\bE=\bF$.} between $\bA$ and $\bW\bZ$ (we will use the two terms interchangeably  when necessary).

Note that $L(\bW,\bZ)=\frac{1}{2} \normf{\bW\bZ-\bA}^2$ is convex in $\bZ$ when $\bW$ is fixed, and convex in $\bW$ when $\bZ$ is fixed---a property known as \textit{marginal convexity}.
This motivates an alternating optimization scheme. 
We can first minimize the loss with respect to $\bZ$
while keeping $\bW$ fixed, and subsequently minimize it with respect to $\bW$ with $\bZ$ fixed:
$$
\left\{
\begin{aligned}
	\bZ &\leftarrow \mathop{\arg \min}_{\bZ} L(\bW,\bZ);    \qquad \text{(ALS1)} \\ 
	\bW &\leftarrow \mathop{\arg \min}_{\bW} L(\bW,\bZ). \qquad \text{(ALS2)}
\end{aligned}
\right.
$$
This is the \textit{alternating least squares (ALS)} algorithm \citep{comon2009tensor, takacs2012alternating, giampouras2018alternating}, a special case of 2-BCD.
Convergence to a local minimum is guaranteed if the loss decreases at each iteration---a property we will discuss further in the sequel.

\index{Coordinate descent algorithm}
\index{Convexity}
\index{Global minimum}
\index{ALS}
\begin{remark}[Convexity and Global Minimum]
Although it can be shown that  the loss function $L(\bW,\bZ)=\frac{1}{2}\normf{\bW\bZ-\bA}^2$ is marginally convex, it is not jointly convex in ($\bW, \bZ$). 
Therefore, the global minimum cannot generally be found.
However, the alternating scheme is guaranteed to converge to a stationary point (typically a local minimum).

More generally, let $D(\bA, \bB)$ be convex in the second argument $\bB$. Then, $D(\bA,\bW\bZ)$ is convex in $\bW$ for $\bZ$ fixed and vice versa; see Problem~\ref{prob:sep_conv}.
\end{remark}

\subsection*{Updating $\bZ$ Given $\bW$}

Fixing $\bW$, we minimize $L(\bZ|\bW) = \frac{1}{2}\normf{\bW\bZ-\bA}^2$. 
The gradient with respect to $\bZ$ is:
\begin{equation}\label{equation:givenw-update-z-allgd}
\begin{aligned}
\nabla_{\bZ} L(\bZ|\bW) &= \frac{1}{2}
\frac{\partial \,\,\trace\left((\bW\bZ-\bA)(\bW\bZ-\bA)^\top\right)}{\partial \bZ}
\stackrel{\star}{=}  \bW^\top(\bW\bZ-\bA) \in \real^{K\times N},
\end{aligned}
\end{equation}
where the first equality arises from the definition of the Frobenius norm (Definition~\ref{definition:frobernius-in-svd}) such that $\normf{\bA}=\sqrt{\trace(\bA\bA^\top)}$, and the equality ($\star$)  is a consequence of the fact that $\frac{\partial \trace(\bA\bA^\top)}{\partial \bA} = 2\bA$. 
Setting this to zero yields the  ``candidate" update:
\begin{equation}\label{equation:als-z-update}
\textbf{(``Candidate" update for $\bZ$)}: \gap {\bZ = (\bW^\top\bW)^{-1} \bW^\top \bA  \leftarrow \mathop{\arg \min}_{\bZ} L(\bZ|\bW).}
\end{equation}
To confirm this candidate is indeed a minimizer, we need to examine that the Hessian matrix is positive definite (Definition~\ref{definition:psd-pd-defini}):
$$
\nabla^2_{\bZ} L(\bZ|\bW) \succ 0.~
\footnote{In short, a twice continuously differentiable function $f$ over an open convex set $\sS$ is called \textit{convex} if and only if $\nabla^2f(\bx)\succeq  \bzero $ for any $\bx\in \sS$ (sufficient and necessary for convex); and called \textit{strictly convex} if $\nabla^2f(\bx)\succ \bzero$ for any $\bx\in \sS$ (only sufficient for strictly convex, e.g., $f(x)=x^6$ is strictly convex, but $f^{\prime\prime}(x)=30x^4$ is equal to zero at $x=0$.). 
And when the convex function $f$ is a continuously differentiable function over a convex set $\sS$, the stationary point $\nabla f(\bx^\star)=\bzero$ of $\bx^\star\in\sS$ is  a \textit{global minimizer} of $f$ over $\sS$.
In our context, when given $\bW$ and updating $\bZ$, the function is defined over the entire space $\real^{K\times N}$ \citep{lu2025practical}.
}
$$
To demonstrate this, we explicitly express the Hessian matrix as
\begin{equation}\label{equation:als-z-update_hessian}
\nabla^2_{\bZ} L(\bZ|\bW)= \widetildebW^\top\widetildebW \in \real^{KN\times KN},
~\footnote{
A block-diagonal matrix whose block matrix on the diagonal is $\bW^\top\bW$. And it can be equivalently denoted as $\nabla^2_{\bZ} L(\bZ|\bW) = \diag(\bW,\bW,\ldots,\bW)^\top\diag(\bW,\bW,\ldots,\bW)$.
}
\end{equation}
 where 
$$
\widetildebW \triangleq
\diag(\bW,\bW,\ldots,\bW)=
\begin{bmatrix}
\bW & \bzero & \ldots & \bzero\\
\bzero & \bW & \ldots & \bzero\\
\vdots & \vdots & \ddots & \vdots \\
\bzero & \bzero & \ldots & \bW
\end{bmatrix}.
$$
This Hessian is positive definite if $\bW\in \real^{M\times K}$ has full column rank $K<M$ (Problem~\ref{prob:rank-of-ata}), ensuring strict convexity and a unique global minimizer.

The key question now is: Does $\bW$ maintain full rank during iterations? 
Otherwise, we cannot claim the update of $\bZ$ in Equation~\eqref{equation:als-z-update} reduces the loss (due to convexity) so that the matrix decomposition progressively improves the approximation of the original matrix $\bA$ by $\bW\bZ$ in each iteration.
We will address the positive definiteness of the Hessian matrix shortly, relying on the following lemma.
\begin{lemma}[Rank of $\bZ$ after Updating]\label{lemma:als-update-z-rank}
	Suppose $\bA\in \real^{M\times N}$ has full rank with \textcolor{mylightbluetext}{$M\leq N$} and $\bW\in \real^{M\times K}$ has full rank with $K<M$ (i.e., $K<M\leq N$). Then the update of $\bZ=(\bW^\top\bW)^{-1} \bW^\top \bA \in \real^{K\times N}$ in Equation~\eqref{equation:als-z-update} has full rank.
\end{lemma}
\begin{proof}[of Lemma~\ref{lemma:als-update-z-rank}]
Since $\bW^\top\bW\in \real^{K\times K}$ has full rank if $\bW$ has full rank (Problem~\ref{prob:rank-of-ata}), it follows that $(\bW^\top\bW)^{-1} $ has full rank. 

Suppose $\bW^\top\bx=\bzero$. This implies that $(\bW^\top\bW)^{-1} \bW^\top\bx=\bzero$. Thus, the following two null spaces satisfy:
$
\nspace(\bW^\top) \subseteq \nspace\left((\bW^\top\bW)^{-1} \bW^\top\right).
$
Moreover, suppose $\bx$ lies in the null space of $(\bW^\top\bW)^{-1} \bW^\top$ such that $(\bW^\top\bW)^{-1} \bW^\top\bx=\bzero$. And since $(\bW^\top\bW)^{-1} $ is invertible, it implies $ \bW^\top\bx=(\bW^\top\bW)\bzero=\bzero$, leading to
$
\nspace\left((\bW^\top\bW)^{-1} \bW^\top\right)\subseteq \nspace(\bW^\top).
$
Consequently, through ``sandwiching," it follows that 
\begin{equation}\label{equation:als-z-sandiwch1}
\nspace(\bW^\top) = \nspace\left((\bW^\top\bW)^{-1} \bW^\top\right).
\end{equation}
Therefore, $(\bW^\top\bW)^{-1} \bW^\top$ has full rank $K$. Let $\bT\triangleq(\bW^\top\bW)^{-1} \bW^\top\in \real^{K\times M}$, and suppose $\bT^\top\bx=\bzero$. This implies $\bA^\top\bT^\top\bx=\bzero$, yielding 
$
\nspace(\bT^\top) \subseteq \nspace(\bA^\top\bT^\top).
$
Similarly, suppose $\bA^\top(\bT^\top\bx)=\bzero$. Since $\bA$ has full rank with the dimension of the null space being 0: $\dim\left(\nspace(\bA^\top)\right)=0$, $(\bT^\top\bx)$ must be zero. The claim follows  since $\bA$ has full rank $M$ with the row space of $\bA^\top$ being equal to the column space of $\bA$, where $\dim\left(\cspace(\bA)\right)=M$ and  $\dim\left(\nspace(\bA^\top)\right) = M-\dim\left(\cspace(\bA)\right)=0$. 
Consequently, $\bx$ is in the null space of $\bT^\top$ if $\bx$ is in the null space of $\bA^\top\bT^\top$:
$
\nspace(\bA^\top\bT^\top)\subseteq \nspace(\bT^\top).
$
By ``sandwiching" again, we obtain
\begin{equation}\label{equation:als-z-sandiwch2}
\nspace(\bT^\top) = \nspace(\bA^\top\bT^\top).
\end{equation}
Since $\bT^\top$ has full rank $K<M\leq N$, it follows that $\dim\left(\nspace(\bT^\top) \right) = \dim\left(\nspace(\bA^\top\bT^\top)\right)=0$.
Therefore,
$\bZ^\top=\bA^\top\bT^\top$ has full rank $K$.
We complete the proof.
\end{proof}

\subsection*{Updating $\bW$ Given $\bZ$}
By symmetry ($\bA=\bW\bZ$ if and only if $\bA^\top=\bZ^\top\bW^\top$),  the update for $\bW$ mirrors that for 
$\bZ$.
With $\bZ$ fixed, the gradient is:
$$
\begin{aligned}
\nabla_{\bW} L(\bW|\bZ) &= \frac{1}{2}
\frac{\partial \trace\left((\bW\bZ-\bA)(\bW\bZ-\bA)^\top\right)}{\partial \bW}
= (\bW\bZ-\bA)\bZ^\top \in \real^{M\times K}.
\end{aligned}
$$
leading to the ``candidate" update:
\begin{equation}\label{equation:als-w-update}
\textbf{(``Candidate" update for $\bW$)}:\gap{\bW^\top = (\bZ\bZ^\top)^{-1}\bZ\bA^\top  \leftarrow \mathop{\arg\min}_{\bW} L(\bW|\bZ).}
\end{equation}
Once more, we emphasize that the update is merely a ``candidate" update. 
Further validation is necessary  to ascertain the positive definiteness of the Hessian matrix.
The Hessian matrix is obtained as follows:
\begin{equation}\label{equation:als-w-update_hessian}
\begin{aligned}
\nabla_{\bW}^2 L(\bW|\bZ) =\widetildebZ\widetildebZ^\top \in \real^{KM\times KM},
\end{aligned}
\end{equation}
where $\widetildebZ\triangleq\diag(\bZ,\bZ,\ldots,\bZ)\in\real^{KM\times NM}$ is defined analogously to $\widetildebW$ in \eqref{equation:als-z-update_hessian}.
Therefore, by similar analysis, if $\bZ$ has full rank with $K<N$, the Hessian matrix is positive definite.

In Lemma~\ref{lemma:als-update-z-rank}, we proved that $\bZ$ has full rank under certain conditions, ensuring that  the Hessian matrix in Equation~\eqref{equation:als-w-update_hessian} is positive definite, and the update in Equation~\eqref{equation:als-w-update} exists.
We now prove that $\bW$ also has full rank under certain conditions, such that the Hessian in Equation~\eqref{equation:als-z-update_hessian} is positive definite, and the update in  Equation~\eqref{equation:als-z-update} exists.
\begin{lemma}[Rank of $\bW$ after Updating]\label{lemma:als-update-w-rank}
Suppose $\bA\in \real^{M\times N}$ has full rank with \textcolor{mylightbluetext}{$M\geq N$} and $\bZ\in \real^{K\times N}$ has full rank with $K<N$ (i.e., $K<N\leq M$). Then the update of $\bW^\top = (\bZ\bZ^\top)^{-1}\bZ\bA^\top$ in Equation~\eqref{equation:als-w-update} has full rank.
\end{lemma}
The proof of Lemma~\ref{lemma:als-update-w-rank} is similar to that of  Lemma~\ref{lemma:als-update-z-rank}, and we shall not repeat the details.

\paragrapharrow{Key observation.}
Combining the observations in Lemma~\ref{lemma:als-update-z-rank} and Lemma~\ref{lemma:als-update-w-rank}, 
e see that if $\bZ$ and $\bW$ are initialized with full rank, then subsequent updates preserve full rank---provided the data matrix $\bA$ satisfies compatible rank conditions.
However, note that Lemma~\ref{lemma:als-update-z-rank} requires $M\leq N$, while Lemma~\ref{lemma:als-update-w-rank} requires $M\geq N$. 
Thus, both lemmas hold simultaneously only when $M=N$---an unrealistic constraint in practice (e.g., Netflix has far more users than movies).
To overcome this limitation, we will introduce regularization in the next section, which ensures numerical stability and full-rank updates even when $M\neq N$.
(Alternatively, one may use the Moore--Penrose pseudo-inverse, as discussed in Problems~\ref{prob:als_pseudo1}--\ref{prob:als_pseudon}.)
Algorithm~\ref{alg:als} summarizes the ALS procedure.
Since the loss $\frac{1}{2}\normf{\bA-\bW\bZ}^2$ is nonincreasing and bounded below, it converges.
At convergence, the gradients vanish: $\nabla_{\bZ} L(\bZ|\bW)=\bzero$ and $\nabla_{\bW} L(\bW|\bZ)=\bzero$.

\begin{algorithm}[H] 
\caption{Alternating Least Squares (ALS)}
\label{alg:als}
\begin{algorithmic}[1] 
\Require Matrix $\bA\in \real^{M\times N}$ \textbf{with $M= N$};
\State Initialize $\bW\in \real^{M\times K}$, $\bZ\in \real^{K\times N}$ \textbf{with full rank and $K<M= N$}; 
\State Choose a stopping criterion on the approximation error $\delta$;
\State Choose the maximum number of iterations $C$;
\State $iter=0$; \Comment{Count for the number of iterations}
\While{$\normf{\bA-\bW\bZ}>\delta $ and $iter<C$} 
\State $iter=iter+1$;
\State $\bZ \leftarrow (\bW^\top\bW)^{-1} \bW^\top \bA  \leftarrow \mathop{\arg \min}_{\bZ} L(\bZ|\bW)$; 
\Comment{(ALS$_1$)}
\State $\bW^\top \leftarrow (\bZ\bZ^\top)^{-1}\bZ\bA^\top  \leftarrow \mathop{\arg\min}_{\bW} L(\bW|\bZ)$;
\Comment{(ALS$_2$)}
\EndWhile
\State Output $\bW,\bZ$;
\end{algorithmic} 
\end{algorithm}

\section{More on the Error Measure and Statistical Interpretation$^*$}\label{section:more_err_sta_als}
We briefly introduce alternative error measures for matrix factorization and approximation problems, along with their statistical interpretations.
\paragrapharrow{Frobenius norm and Gaussian noise.}
As previously noted, using the Frobenius norm as a loss function corresponds to assuming that the observed data are corrupted by i.i.d. Gaussian noise. This leads to a smooth least-squares optimization problem.
To see this, let  $\bB\triangleq\bW\bZ$ denote the low-rank approximation, and assume the noise $\epsilon$ is i.i.d. Gaussian with zero mean and variance  $\sigma^2$:
\begin{equation}\label{equation:gau_noise}
a_{mn}= b_{mn}+ \epsilon_{mn}, \gap 
\epsilon_{mn}\sim \normal(0, \sigma^2), 
\gap \forall m,n.
\end{equation}
Under this assumption, the log-likelihood of the observed matrix $\bA$ given $\bB$ and $\sigma^2$ is
$$
\ln p(\bA\mid \bB, \sigma^2) 
=
-\frac{1}{2\sigma^2}\sum_{m,n=1}^{M,N} (a_{mn}-b_{mn})^2 + C(\sigma)
=
-\frac{1}{2\sigma^2}\normf{\bA-\bW\bZ}^2 +C(\sigma),
$$
where $C(\sigma^2)$ is a constant depending only on $\sigma$.
Thus, given $\bA$, the parameters $\bW,\bZ$, and $\sigma^2$ can be estimated via maximum likelihood estimation (MLE) by maximizing the log-likelihood:
$$
\mathopmax{\bW,\bZ,\sigma^2} -\frac{1}{2\sigma^2}\normf{\bA-\bW\bZ}^2
=
\mathopmin{\bW,\bZ,\sigma^2}\frac{1}{2\sigma^2}\normf{\bA-\bW\bZ}^2.
$$
If the noise level $\sigma^2$ is constant (or treated as fixed), this reduces to the standard Frobenius-norm minimization in \eqref{equation:als-per-example-loss2}.
When the noise is not i.i.d., the resulting loss becomes a weighted $\ell_2$-norm; see Problem~\ref{prob:non_iid_gaus}.

\paragrapharrow{Matrix $\ell_1$-norm and Laplace noise.}
Similarly, if the noise follows an i.i.d. Laplace distribution (Definition~\ref{definition:laplace_distribution}) with location parameter  $\mu=0$ and scale parameter $\sigma$, the log-likelihood function becomes 
$$
\ln p(\bA\mid \bB, \sigma) = -\frac{1}{2\sigma}\sum_{m,n=1}^{M,N}\abs{a_{mn}-b_{mn}} +C(\sigma)
=
-\frac{1}{2\sigma} \normmone{\bA-\bW\bZ} +C(\sigma),
$$
where $\normmone{\cdot}$ denotes the matrix $\ell_1$-norm.
In this case, MLE is equivalent to minimizing the matrix $\ell_1$-norm if the scale is held constant.
Like Gaussian noise, Laplace noise is additive:
$$
a_{mn}= b_{mn}+ \epsilon_{mn}, \gap 
\epsilon_{mn}\sim \text{Laplace}(0, \sigma), 
\gap \forall m,n.
$$
In practice, the (matrix) $\ell_1$-norm is more robust than the Frobenius norm, making it suitable for data contaminated by sparse, large errors.

\paragrapharrow{Matrix $\ell_\infty$-norm and uniform noise.}
Alternatively, one may assume i.i.d. uniform noise:
$$
a_{mn} \sim \text{Uniform}(b_{mn} -\sigma, b_{mn}+\sigma), 
\gap \forall m,n,
$$
where $\text{Uniform}(a, b)$ denotes the uniform distribution over $[a,b]$.
This implies an additive noise model:
$$
a_{mn}= b_{mn}+ \epsilon_{mn}, \gap 
\epsilon_{mn}\sim \text{Uniform}(-\sigma, \sigma), 
\gap \forall m,n.
$$
The corresponding MLE minimizes the matrix $\ell_\infty$-norm:
$$
L(\bW,\bZ) = \normminf{\bA-\bW\bZ} = \mathopmax{m,n} \abs{a_{mn}-b_{mn}}.
$$
This formulation is useful when the goal is to control the worst-case error rather than average error.

\paragrapharrow{Poisson and nonnegative integers.}
When the entries of $\bA$ are nonnegative integers (such as word counts in text mining or event frequencies), it is natural to model them using a Poisson distribution (Definition~\ref{definition:poisson_distribution}):
$$
p(a_{mn}=x\mid b_{mn}) = \frac{b_{mn}^{x}}{x!} \exp(-b_{mn})
\,\,\text{ with }\,\,
\Exp[a_{mn}] = b_{mn}, \,\, \Var[a_{mn}] = b_{mn},
\gap \forall m,n.
$$
This indicates $a_{mn}=0$ if $b_{mn}=0$.
The MLE corresponds to minimizing the Kullback--Leibler (KL) divergence between $\bA$ and $\bW\bZ$:
\begin{equation}\label{equation:als_poi_los}
L(\bW,\bZ) = \sum_{m,n=1}^{M,N} l(a_{mn}, b_{mn}) = \sum_{m,n=1}^{M,N}\left(a_{mn} \ln(\frac{a_{mn}}{b_{mn}}) -a_{mn}+b_{mn}\right),
\end{equation}
where $l(a_{mn}, b_{mn})$ is assumed to be 0 if $b_{mn}=0$. 
The term $2l(a_{mn}, b_{mn}) =2( \ln p(a_{mn} \mid a_{mn}) - \ln p(a_{mn} \mid b_{mn}))$ is called the \textit{deviance (goodness-of-fit statistic)}, which measures the goodness of fit of the model by  comparing the log-likelihood difference between the \textit{saturated model} and the \textit{current model}.
Importantly, Poisson noise is not additive---the variance depends on the signal level $b_{mn}$.

\paragrapharrow{Multiplicative Gamma noise.}
Unlike Gaussian, Laplace, and uniform noise---which are all additive---the Gamma noise (Definition~\ref{definition:gamma-distribution}) model assumes multiplicative noise:
$$
a_{mn}= b_{mn} \cdot  \epsilon_{mn}, \gap 
\epsilon_{mn}\sim \text{Gamma}(r, \lambda)=\frac{\lambda^r}{\Gamma(r)} \epsilon_{mn}^{r-1}\exp(-\lambda \epsilon_{mn}), 
\gap \forall m,n.
$$
where $\Gamma(r)$ is the Gamma function.
When the mean of $a_{mn}$ satisfies  $\Exp[a_{mn}]=\frac{r}{\lambda}=1$, 
the MLE corresponds to minimizing the following loss (\textit{Itakura--Saito divergence, IS divergence} \citep{itakura1968analysis}):
\begin{equation}\label{equation:als_gama_los}
L(\bW,\bZ) = D(\bA,\bB) = \sum_{m,n=1}^{M,N} \frac{a_{mn}}{b_{mn}} - \ln\left(\frac{a_{mn}}{b_{mn}}\right)-1.
\end{equation}
\begin{exercise}[Scale Invariant of IS Divergence]
Show that the IS divergence is scale invariant: $D(\bA,\bB) = D(\gamma\bA, \gamma\bB)$ for any $\gamma>0$.
\end{exercise}
This scale invariance is particularly valuable in applications such as audio source separation, where low-energy frequency components can be perceptually as important as high-energy ones \citep{gillis2020nonnegative}.

\index{Regularization}
\section{Regularization and Identifiability: Extension to General Matrices}\label{section:regularization-extention-general}

\textit{Regularization} is a machine learning technique used to prevent overfitting and improve a model's generalization performance. Overfitting occurs when a model becomes overly complex and fits the training data too closely, leading to poor performance on unseen data.
To address this issue, regularization introduces a constraint or penalty term into the loss function during optimization, discouraging excessive model complexity. This creates a trade-off between fitting the training data well and maintaining a simple, generalizable model.
Common types of regularization include $\ell_1$-regularization (LASSO \citep{lu2026first}), $\ell_2$-regularization (Tikhonov or Ridge regularization; see Section~\ref{section:pre_ls}), and elastic net regularization (a combination of $\ell_1$- and $\ell_2$-penalties). 
Regularization is widely used in algorithms such as linear regression, logistic regression, and neural networks \citep{lecun2015deep, goodfellow2016deep}.

In the context of the ALS problem, we can incorporate an $\ell_2$-regularization term to minimize the following regularized loss function:
\begin{equation}\label{equation:als-regularion-full-matrix}
L(\bW,\bZ)  = \frac{1}{2}\normf{\bW\bZ-\bA}^2 +\frac{1}{2}\lambda_w \normf{\bW}^2 + \frac{1}{2}\lambda_z \normf{\bZ}^2, \qquad \lambda_w>0, \lambda_z>0,
\end{equation}
where the gradient with respect to $\bZ$ and $\bW$ are given, respectively, by 
\begin{equation}\label{equation:als-regulari-gradien}
\left\{
\begin{aligned}
\nabla_{\bZ} L(\bZ|\bW) &= \bW^\top(\bW\bZ-\bA) + \lambda_z\bZ \in \real^{K\times N};\\
\nabla_{\bW} L(\bW|\bZ)  &= (\bW\bZ-\bA)\bZ^\top + \lambda_w\bW \in \real^{M\times K}.
\end{aligned}
\right.
\end{equation}
The Hessian matrices are given, respectively, by 
$$
\left\{
\begin{aligned}
\nabla^2_{\bZ} L(\bZ|\bW) &= \widetildebW^\top\widetildebW+ \lambda_z\bI \in \real^{KN\times KN};\\
\nabla^2_{\bW} L(\bW|\bZ)  &= \widetildebZ\widetildebZ^\top + \lambda_w\bI \in \real^{KM\times KM}, \\
\end{aligned}
\right.
$$
which are positive definite due to the perturbation by the regularization:
$$
\left\{
\begin{aligned}
\bx^\top (\widetildebW^\top\widetildebW +\lambda_z\bI)\bx 
&= \underbrace{\bx^\top\widetildebW^\top\widetildebW\bx}_{\geq 0} + \lambda_z \normtwo{\bx}^2>0, \gap \text{for nonzero $\bx$};\\
\bx^\top (\widetildebZ\widetildebZ^\top +\lambda_w\bI)\bx 
&= \underbrace{\bx^\top\widetildebZ\widetildebZ^\top\bx}_{\geq 0} + \lambda_w \normtwo{\bx}^2>0,\gap \text{for nonzero $\bx$}.
\end{aligned}
\right.
$$
Regularization ensures that the Hessian matrices remain invertible, guaranteeing unique minimizers in each ALS subproblem---even when $\bW$ and $\bZ$ are rank-deficient. 
As a result, matrix factorization via ALS can be applied to any matrix, regardless of whether 
$M>N$ or $M<N$. In rare cases, one may even choose $K>\max\{M, N\}$ to obtain a high-rank approximation of $\bA$. However, in most practical settings, the goal is to find a low-rank approximation with $K<\min\{M, N\}$.
Setting the gradients to zero yields the closed-form updates:
\begin{equation}\label{equation:als-regular-final-all}
\left.
\begin{aligned}
\bZ &\leftarrow (\bW^\top\bW+ \lambda_z\bI)^{-1} \bW^\top \bA 
\qquad \text{and}\qquad 
\bW^\top \leftarrow (\bZ\bZ^\top+\lambda_w\bI)^{-1}\bZ\bA^\top .
\end{aligned}
\right.
\end{equation}
The regularization parameters $\lambda_z, \lambda_w\in \real_{++}$ control the trade-off between approximation accuracy and solution smoothness (or simplicity). Their optimal values are typically problem-dependent and are often selected via cross-validation (CV). The full procedure is summarized in Algorithm~\ref{alg:als-regularizer}.
We will also introduce the \textit{alternating direction methods of multipliers  (ADMM)} in Section~\ref{section:nmf_admm_all} for solving matrix factorization problems with $\ell_2$- or $\ell_1$-regularization. ADMM can be extended to handle additional constraints, such as nonnegativity.

While $\ell_2$ (or $\ell_1$) regularizations generalize ALS to arbitrary matrices, they are not the only options---especially in matrix completion settings where many entries of $\bA$ are missing.
For example, the \textit{nuclear norm}, defined as the sum of singular values of $\bW\bZ$, is a popular convex surrogate for rank minimization. The Soft-Impute algorithm for matrix completion guarantees exact recovery of an $N\times N$ matrix $\bA$ of rank $R$ when the number of observed entries $z$ satisfies 
$
z\geq C R N \ln N,
$
for some universal constant $C>0$ \citep{gross2011recovering, hastie2015statistical}.
Interestingly,  $\ell_2$-regularization on $\bW$ and $\bZ$ can sometimes be reformulated in terms of the nuclear norm (see Problem~\ref{problem:nuclear_equi}).

\index{Cross-validation}
\begin{algorithm}[H] 
\caption{Alternating Least Squares with Regularization}
\label{alg:als-regularizer}
\begin{algorithmic}[1] 
\Require Matrix $\bA\in \real^{M\times N}$;
\State Initialize $\bW\in \real^{M\times K}$, $\bZ\in \real^{K\times N}$ \textbf{randomly, without requiring any condition on rank or the relationship among  $M, N, K$}; 
\State Choose a stopping criterion on the approximation error $\delta$;
\State Choose regularization parameters $\lambda_w, \lambda_z$;
\State Choose the  maximal number of iterations $C$;
\State $iter=0$; \Comment{Count for the number of iterations}
\While{$\normf{\bA-\bW\bZ}>\delta $ and $iter<C$}
\State $iter=iter+1$; 
\State $\bZ \leftarrow (\bW^\top\bW+ \lambda_z\bI)^{-1} \bW^\top \bA  \leftarrow \mathop{\arg \min}_{\bZ} L(\bZ|\bW)$;
\State $\bW^\top \leftarrow (\bZ\bZ^\top+\lambda_w\bI)^{-1}\bZ\bA^\top  \leftarrow \mathop{\arg\min}_{\bW} L(\bW|\bZ)$;
\EndWhile
\State Output $\bW,\bZ$;
\end{algorithmic} 
\end{algorithm}

\index{K-means problem}
\index{Regularization}
\index{Constraint}
\index{Identifiability}
\index{Sparsity}
\index{Group sparsity}
\paragrapharrow{Regularization as constraints and identifiability.} 
Regularization terms like $\lambda_w\normf{\bW}^2$ in \eqref{equation:als-regularion-full-matrix} can be interpreted as soft constraints of the form $\normf{\bW}\leq C$, where $C$ is a constant linked via Lagrange multipliers (see, for example, \citet{boyd2004convex, lu2025practical}). 
Various constraints can be imposed on the factors $\bW$ and $\bZ$, such as: nonnegativity  (discussed in Chapter~\ref{chapter:nmf})  and  sparsity (discussed in Section~\ref{section:reg_geom_inter}).
Moreover, the factorization  $\bA=\bW\bZ$ suffers from non-identifiability; the product remains unchanged under scaling transformations: $\bW[:, k]\bZ[k,:] = (\gamma\bW[:, k])(\frac{1}{\gamma}\bZ[k,:])$ for any scalar $\gamma\neq 0$ and $k\in\{1,2,\ldots,K\}$.
Thus, while $\bW$ and $\bZ$ have $(M+N)K$ parameters, the effective degrees of freedom are only $(M+N-1)K$. Regularization helps mitigate this ambiguity by incorporating prior knowledge through penalties or constraints.
Beyond $\ell_1$ and $\ell_2$, several alternative regularizers have been proposed \citep{lee2009semi, bach2011convex, cai2010graph, iordache2012total, gillis2020nonnegative}:
\begin{itemize}
\item \textbf{Minimum-volume.} Impose the regularizer $\lambda_w \det(\bW^\top\bW)$ to encourage the columns of $\bW$ to span a small-volume simplex that tightly encloses the data points; see Section~\ref{section:gvg_model}.

\item \textbf{K-means constraint (vector quantization).}
In the context of the \textit{K-means problem}, where each column of $\bA$ is a data point, the goal is to determine a set of $K$ centroids $\bw_k, k\in\{1,2,\ldots,K\}$ (i.e., the columns of $\bW$) such that the sum of distances between each data point and its nearest centroid is minimized. 
This setup is equivalent to the low-rank matrix factorization problem where the second factor $\bZ$ must have exactly one nonzero entry per column, which is set to one, indicating the assignment of data points to their respective centroids: $\bZ\bZ^\top$ is diagonal and $\bZ\in\{0,1\}^{K\times N}$, where the diagonal values of $\bZ\bZ^\top$ indicate the number of data points associated with each cluster. The columns of $\bW$ then represent the cluster centroids \citep{zhang2017matrix, gillis2020nonnegative}.

\item \textbf{Sparsity.} 
The ``$\ell_0$-norm" (number of nonzero entries) directly measures sparsity but is non-convex and discontinuous, making optimization difficult. 
Instead, the $\ell_1$-norm is commonly used as a convex relaxation that promotes sparsity while enabling efficient optimization (e.g., via gradient-based methods) \citep{lu2025practical}.
Applying $\ell_1$-regularization to $\bZ$ induces element-wise sparsity, which is useful in applications like facial feature extraction---yielding localized, interpretable features.
The \textit{spar} operator, introduced in \citet{hoyer2004non}, is continuous and promotes sparsity based on the ratio of $\ell_1$ to $\ell_2$-norm: for $\bx\neq\bzero\in\real^N$, $spar(\bx)=\frac{\sqrt{N}-\normone{\bx}/\normtwo{\bx}}{\sqrt{N}-1}\in[0,1]$.
$spar(\bx)=0$ if and only if $\normone{\bx}=\sqrt{N}\normtwo{\bx}$ and all entries of $\bx$ are equal. 
And $spar(\bx)=1$ if and only if $\normone{\bx}=\normtwo{\bx}$, in which case $\normzero{\bx}=1$.
The higher the value, the sparser; for example $spar([1,0,0]) > spar([1,1,1])$.
In the matrix factorization context, such a sparsity constraint can be applied by ensuring $spar(\bW[:,k]) \geq r_w$ and $spar(\bZ[k,:])\geq r_z$ for all $k\in\{1,2,\ldots,K\}$, where $r_w$ and $r_z$ are constants in the interval $[0,1]$ that impose a minimal sparsity level on the columns of $\bW$ and rows of $\bZ$.

\item \textbf{Group sparsity.}
When columns of $\bZ$ naturally group into subsets $\sG$ (e.g., by semantic meaning or time segments), we may wish to enforce sparsity at the group level rather than per element.
Assuming $\cup_{\scriptsize g\in\sG}g=\{1,2,\ldots,N\}$, we penalize aggregated norms per group---e.g., using 
$\sum_{g\in\sG}\normtwo{\bZ_{:,g}}$ or $\sum_{g\in\sG}\normone{\bZ_{:,g}}$.
This encourages entire groups of features to be zero or active together, aligning with structured data assumptions.

\index{Orthogonal matrix factorization}
\item \textbf{Orthogonality.} 
To reduce redundancy among features or activations, orthogonality can be encouraged via penalties: 
$\lambda_w\normf{\bW^\top \bW - \bI_K}^2$
and 
$\lambda_z\normf{\bZ \bZ^\top - \bI_K}^2$.
These promote distinct, uncorrelated components---opposite in spirit to minimum-volume regularization.
If the hard constraint $\bZ\bZ^\top=\bI_K$ is enforced (rather than a penalty), the problem becomes \textit{orthogonal matrix factorization} (Problem~\ref{prob:ortho_mf}), where $\bZ$ has orthonormal rows. This formulation resembles a soft clustering model, where each data point is expressed as a linear combination of orthogonal basis vectors.

\item \textbf{Spatial smoothness.} When factorizing image data (vectorized and stacked as columns of $\bA$), spatial structure is lost. To preserve local coherence in the basis images (columns of $\bW$), spatial regularizers can be added.
For edge-preserving smoothness, \textit{total variation (TV)} regularization is effective:
$
\lambda_w \big(\sum_{k=1}^{K} \sum_{(i_1, i_2) \in \sS} \abs{w_{i_1,k}-w_{i_2,k}}\big),
$
where $\sS$  denotes neighboring pixel pairs. TV uses $\ell_1$ differences to maintain sharp edges. See also the denoising least squares problem in Problem~\ref{prob:denoise_rls}.

\item  \textbf{Graph regularization.} To preserve geometric relationships among data points in the latent space, we can encourage nearby points in $\bA$ to remain close in $\bZ$. Specifically, if $\normtwo{\ba_i-\ba_j}$ is small, then $\normtwo{\bz_i-\bz_j}$ should also be small.
This is achieved via the regularizer: $\lambda_z \sum_{i,j} r_{ij} \normtwo{\bz_i-\bz_j}^2$, where  $r_{ij}$ reflects similarity between points $i$ and $j$. A common choice is $r_{ij} = \exp\{-\gamma\normtwo{\ba_i-\ba_j}\}$ for $\gamma > 0$.
The matrix $\bR=[r_{ij}]\in\real^{N\times N}$ defines a weighted graph over the data, giving rise to \textit{graph-regularized matrix factorization}.
This framework also supports semi-supervised learning: if partial labels are known (e.g., certain face images belong to the same person), we can set $r_{ij} = 1$ for same-label pairs and 0 otherwise. This contrasts with purely unsupervised factorization.
\end{itemize}

\index{Missing entries}
\index{Hadamard product}
\index{Netflix recommender}
\section{Missing Entries and Rank-One Updates}\label{section:alt-columb-by-column}
As noted previously, matrix decomposition via ALS is widely used in recommender systems such as the Netflix Prize dataset, where a large fraction of entries are missing because users have not watched certain movies or have chosen not to rate them.
In this setting, the low-rank matrix factorization problem is commonly referred to as \textit{matrix completion}, which aims to recover unobserved entries from partial observations \citep{jain2017non}.
To handle missing data, we introduce a \textit{mask matrix} $\bM\in \{0,1\}^{M\times N}$, where $m_{mn}\in \{0,1\}$ indicates whether  user $n$ has rated  movie $m$ or not. 
The loss function then becomes:
$$
L(\bW,\bZ) = \frac{1}{2}\normf{\bM\hadaprod  \bA - \bM\hadaprod (\bW\bZ)}^2,
$$
where $\hadaprod$ denotes the Hadamard product (element-wise multiplication). 
This formulation concisely captures the objective: find a low-rank approximation of the rating matrix that agrees with all observed entries.
To solve this problem, we adapt the updates from Equation~\eqref{equation:als-regular-final-all} to operate column-wise (or row-wise), leading to:
\begin{equation}\label{equation:als-ori-all-wz}
	\left\{
	\begin{aligned}
		\bz_n &= (\bW^\top\bW+ \lambda_z\bI)^{-1} \bW^\top \ba_n, &\gap& \text{for $n\in \{1,2,\ldots, N\}$}  ;\\
		\bw_m &= (\bZ\bZ^\top+\lambda_w\bI)^{-1}\bZ\bb_m,  &\gap& \text{for $m\in \{1,2,\ldots, M\}$} ,
	\end{aligned}
	\right.
\end{equation}
where $\bZ=[\bz_1, \bz_2, \ldots, \bz_N]$ and $\bA=[\ba_1,\ba_2, \ldots, \ba_N]$ represent the column partitions of $\bZ$ and $\bA$, respectively. Similarly, $\bW^\top=[\bw_1, \bw_2, \ldots, \bw_M]$ and $\bA^\top=[\bb_1,\bb_2, \ldots, \bb_M]$ are the column partitions of $\bW^\top$ and $\bA^\top$, respectively. This decomposition shows that updates can be performed independently per user or per movie (i.e., \textit{rank-one update}), enabling efficient parallelization---a key advantage of ALS.

\paragrapharrow{Given $\bW$: Update user vectors.}
Let $\bo_n\in \{0,1\}^M$ indicate which movies user $n$ has rated:  $o_{nm}=1$ if user $n$ has rated movie $m$, and $o_{nm}=0$ otherwise. 
Using Matlab-style indexing, the observed ratings for user $n$ are denoted $\ba_n[\bo_n]$, 
and the corresponding rows of $\bW$ are $\bW[\bo_n, :]$.
We aim to approximate the observed entries via:  $\ba_n[\bo_n] \approx \bW[\bo_n, :]\bz_n$, 
which is a (regularized) least squares problem in $\bz_n$. 
The solution is:
\begin{equation}\label{equation:als-ori-all-wz-modif-z}
\begin{aligned}
\bz_n &= \left(\bW[\bo_n, :]^\top\bW[\bo_n, :]+ \lambda_z\bI\right)^{-1} \bW[\bo_n, :]^\top \ba_n[\bo_n], \quad \text{for $n\in \{1,2,\ldots, N\}$} .
\end{aligned}
\end{equation}
The associated loss functions are:
$$
\begin{aligned}
L(\bz_n|\bW) &=\sum_{m\in \bo_n} \left(a_{mn} - \bw_m^\top\bz_n\right)^2
\gap \text{and}\gap
L(\bZ|\bW) =\sum_{n=1}^N\ \sum_{m\in \bo_n} \left(a_{mn} - \bw_m^\top\bz_n\right)^2.
\end{aligned}
$$

\paragrapharrow{Given $\bZ$.}
Similarly, let $\bp_m \in\{0,1\}^{N}$ indicate which users have rated movie $m$:  $p_{mn}=1$ if  movie $m$ was rated by user $n$, and $p_{mn}=0$ otherwise. 
The observed ratings in the $m$-th row of $\bA$ are denoted $\bb_m[\bp_m]$, and the corresponding columns of $\bZ$ are $\bZ[:,\bp_m]$. 
We approximate: $\bb_m[\bp_m] \approx \bZ[:, \bp_m]^\top\bw_m$,  which leads to the update:
\footnote{Note that $\bZ[:, \bp_m]^\top$ is the transpose of $\bZ[:, \bp_m]$, which is equal to $\bZ^\top[\bp_m,:]$, i.e., transposing first and then selecting.}
\begin{equation}\label{equation:als-ori-all-wz-modif-w}
\begin{aligned}
\bw_m &= (\bZ[:, \bp_m]\bZ[:, \bp_m]^\top+\lambda_w\bI)^{-1}\bZ[:, \bp_m]\bb_m[\bp_m],  \quad \text{for $m\in \{1,2,\ldots, M\}$} .
\end{aligned}
\end{equation}
The corresponding loss functions are:
$$
\begin{aligned}
L(\bw_m|\bZ) &=\sum_{n\in \bp_m} \left(a_{mn} - \bw_m^\top\bz_n\right)^2 
\gap \text{and}\gap
L(\bW|\bZ) =\sum_{m=1}^M  \sum_{n\in \bp_m} \left(a_{mn} - \bw_m^\top\bz_n\right)^2 .
\end{aligned}
$$
This procedure is summarized in Algorithm~\ref{alg:als-regularizer-missing-entries}.
Other methods, such as \textit{singular value projection (SVP)}, also address matrix completion. SVP is a projected gradient descent method that iteratively applies gradient steps followed by rank truncation via singular value decomposition (SVD). However, in practice, ALS generally outperforms SVP for matrix completion tasks---particularly in large-scale recommender systems---so we focus on ALS here. For further details on SVP, see \citet{jain2017non}.

\begin{algorithm}[h] 
\caption{Alternating Least Squares with Missing Entries and Regularization}
\label{alg:als-regularizer-missing-entries}
\begin{algorithmic}[1] 
\Require Matrix $\bA\in \real^{M\times N}$;
\State Initialize $\bW\in \real^{M\times K}$, $\bZ\in \real^{K\times N}$ \textbf{randomly, without requiring any condition on rank or the relationship among  $M, N, K$}; 
\State Choose a stopping criterion on the approximation error $\delta$;
\State Choose regularization parameters $\lambda_w, \lambda_z$;
\State Compute the mask matrix $\bM$ from $\bA$;
\State Choose the maximum number of iterations $C$;
\State $iter=0$; \Comment{Count for the number of iterations}
\While{\textbf{$\normf{\bM\hadaprod  \bA- \bM\hadaprod (\bW\bZ)}^2>\delta $} and $iter<C$}
\State $iter=iter+1$; 
\For{$n=1,2,\ldots, N$}
\State $\bz_n \leftarrow \left(\bW[\bo_n, :]^\top\bW[\bo_n, :]+ \lambda_z\bI\right)^{-1} \bW[\bo_n, :]^\top \ba_n[\bo_n]$; \Comment{$n$-th column of $\bZ$}
\EndFor

\For{$m=1,2,\ldots, M$}
\State $\bw_m \leftarrow (\bZ[:, \bp_m]\bZ[:, \bp_m]^\top+\lambda_w\bI)^{-1}\bZ[:, \bp_m]\bb_m[\bp_m]$;\Comment{$m$-th column of $\bW^\top$}
\EndFor
\EndWhile
\State Output $\bW^\top=[\bw_1, \bw_2, \ldots, \bw_M],\bZ=[\bz_1, \bz_2, \ldots, \bz_N]$;
\end{algorithmic} 
\end{algorithm}

\index{Hidden features}
\index{Inner product}
\section{Vector Inner Product and Latent Representations}\label{section:als-vector-product}
The ALS algorithm seeks low-dimensional matrices $\bW\in\real^{M\times K}$ and $\bZ\in\real^{K\times N}$ such that their product approximates the observed data: $\bA\approx \bW\bZ$ in terms of the  squared loss
$\mathop{\min}_{\bW,\bZ}  \sum_{n=1}^N \sum_{m=1}^{M} \left(a_{mn} - \bw_m^\top\bz_n\right)^2$. 
Thus, each entry $a_{mn}$  is modeled as the inner product of two vectors.
Geometrically, the inner product is defined as:
$$
\bw_m^\top\bz_n = \normtwo{\bw_m} \normtwo{\bz_n} \cdot\cos (\theta),
$$
where $\theta$ is the angle between $\bw_m$ and $\bz_n$. 
For fixed vector norms, a smaller angle (i.e., greater alignment) yields a larger inner product.

In the Netflix context, ratings range from 0 to 5, with higher values indicating stronger preference. If the latent vectors $\bw_m$ (movie attributes) and $\bz_n$ (user preferences) are well-aligned, their inner product $\bw_m^\top\bz_n$ will be large---accurately predicting a high rating.
This reveals the core idea of ALS:
each movie is represented by a \textit{latent attribute vector} $\bw_m\in\real^K$, 
and each user is represented by a \textit{latent preference vector} $\bz_n\in\real^K$.
Each dimension in these vectors corresponds to a hidden feature. 
For instance:
the second component $w_{m2}$ might encode how strongly movie $m$ belongs to the ``action" genre,
While $z_{n2}$ might reflect user $n$'s affinity for action movies.
When both are large and positive, their product contributes significantly to $\bw_m^\top\bz_n$, yielding a high predicted rating.

In the factorization $\bA\approx\bW\bZ$:
the rows of $\bW$ capture hidden features of movies,
and the columns of $\bZ$ capture hidden features of users.
However, the semantic meaning of individual latent dimensions is not explicitly defined---it is learned implicitly from data. These dimensions may correspond to interpretable concepts like genre, mood, or director---but they could also represent abstract combinations with no direct real-world label. It is precisely this unobserved, inferred nature that gives rise to the terms ``latent" or ``hidden" vectors.

\index{Stochastic gradient descent}
\index{Gradient descent}
\index{Matrix inverse}
\index{LU decomposition}
\section{Gradient Descent}\label{section:als-gradie-descent}
In Algorithm~\ref{alg:als}, \ref{alg:als-regularizer}, and \ref{alg:als-regularizer-missing-entries}, the loss is minimized by solving linear systems through matrix inversion (e.g., via LU decomposition \citep{lu2021numerical}).
However, this approach becomes impractical in the era of big data. As the volume of data grows, the size of the matrices involved increases, and the computational cost of matrix inversion scales cubically with the number of samples---posing significant challenges for both memory and processing power.
This limitation has motivated the development of gradient-based optimization methods, which avoid explicit matrix inversion. Among these, \textit{gradient descent (GD)} and its variant \textit{stochastic gradient descent (SGD)}  are among the simplest, most efficient, and widely used techniques \citep{lu2022gradient}. They are particularly effective for minimizing convex loss functions. We now describe the core principles behind these methods.

Recall from Equation~\eqref{equation:als-ori-all-wz} that the column-wise update rules are derived directly from the full-matrix formulation in Equation~\eqref{equation:als-regular-final-all}, which includes regularization. To understand the connection to gradient-based methods, consider the regularized loss function:
\begin{equation}\label{als:gradient-regularization-zn}
\footnotesize
\begin{aligned}
L(\bz_n)  &=\frac{1}{2}\normf{\bW\bZ-\bA}^2 +\frac{1}{2}\lambda_w \normf{\bW}^2 + \frac{1}{2}\lambda_z \normf{\bZ}^2
= \frac{1}{2}\normtwo{\bW\bz_n-\ba_n}^2 + \frac{1}{2}\lambda_z \normtwo{\bz_n}^2 + C_{z_n},
\end{aligned}
\end{equation}
where $C_{z_n}$ is constant with respect to $\bz_n$, and $\bZ=[\bz_1, \bz_2, \ldots, \bz_N]$ and $\bA=[\ba_1,\ba_2, \ldots, \ba_N]$ represent the column partitions of $\bZ$ and $\bA$, respectively. 
Taking the gradient and  setting it to zero give the closed-form solution:
$$
\begin{aligned}
\nabla_{\bz_n} L(\bz_n) = \bW^\top\bW\bz_n - \bW^\top\ba_n + \lambda_z\bz_n
\,\,\implies \,\,
\bz_n = (\bW^\top\bW+ \lambda_z\bI)^{-1} \bW^\top \ba_n, \,\,  \forall\,n,
\end{aligned}
$$
which matches the first update rule in Equation Equation~\eqref{equation:als-ori-all-wz}.
Similarly, when minimizing with respect to a movie vector  $\bw_m$, we rewrite the loss using the transpose:
\begin{equation}\label{als:gradient-regularization-wd}
\footnotesize
\begin{aligned}
L(\bw_m )  
&
=\frac{1}{2}\normf{\bZ^\top\bW-\bA^\top}^2 +\frac{1}{2}\lambda_w \normf{\bW^\top}^2 + \frac{1}{2}\lambda_z \normf{\bZ}^2
= \frac{1}{2}\normtwo{\bZ^\top\bw_m-\bb_n}^2 + \frac{1}{2}\lambda_w \normtwo{\bw_m}^2 + C_{w_m},
\end{aligned}
\end{equation}
where $C_{w_m}$ is a constant with respect to $\bw_m$, and $\bW^\top=[\bw_1, \bw_2, \ldots, \bw_M]$ and $\bA^\top=[\bb_1,\bb_2, \ldots,$ $\bb_M]$ represent the column partitions of $\bW^\top$ and $\bA^\top$, respectively. 
Taking the gradient and  setting it to zero lead to the solution:
$$
\begin{aligned}
\nabla_{\bw_m} L(\bw_m) = \bZ\bZ^\top\bw_m - \bZ\bb_n + \lambda_w\bw_m
\,\,\implies\,\,
\bw_m = (\bZ\bZ^\top+\lambda_w\bI)^{-1}\bZ\bb_m, \,\, \forall \, m,
\end{aligned}
$$
which corresponds to the second update rule in Equation~\eqref{equation:als-ori-all-wz}:

Now, suppose we denote the iteration number by a superscript ($t=1,2,\ldots$), and aim to compute the updated variables $\{\bz^\toptone_n, \bw^\toptone_m\}$ based on the current estimates $\{\bZ^\toptzero, \bW^\toptzero\}$. 
In the exact ALS approach, we solve:
$$
\left.
\begin{aligned}
\bz^\toptone_n    &\leftarrow \mathop{\arg \min}_{\bz_n} L(\bz_n^\toptzero)
\qquad\text{and}\qquad
\bw_m^\toptone    \leftarrow \mathop{\arg\min}_{\bw_m} L(\bw_m^\toptzero).
\end{aligned}
\right.
$$
For simplicity, we focus on deriving a gradient-based update for $\bz^\toptone_n    \leftarrow \mathop{\arg \min}_{\bz_n} L(\bz_n^\toptzero)$; the derivation for $\bw_m^\toptone$ follows analogously.

\index{Linear approximation}
\index{Linear update}
\index{Greedy search}
\index{Gradient descent}
\paragrapharrow{Approximation by linear update.} 
Instead of solving the minimization exactly, we approximate the next iterate using a \textit{linear update}:
$$
\textbf{(Linear Update)}: \qquad {\bz^\toptone_n = \bz^\toptzero_n + \eta \bv},
$$
where $\eta> 0$ is a small step size and $\bv$ is a search direction to be determined.
We choose $\bv$ to minimize the loss along this direction:
$$
\bv=\mathop{\arg \min}_{\bv} L(\bz^\toptzero_n + \eta \bv) .
$$
Using a first-order Taylor expansion (Theorem~\ref{theorem:linear_approx}), we approximate:
$$
L(\bz^\toptzero_n + \eta \bv) \approx L(\bz^\toptzero_n ) + \eta \bv^\top \nabla  L(\bz^\toptzero_n ),
$$
where  $\nabla  L(\bz^\toptzero_n )$ represents the gradient of $L(\bz)$ at $\bz^\toptzero_n$. 
To ensure a meaningful direction, we constrain  $\normtwo{\bv}=1$ and solve:
$$
\bv=\mathop{\argmin}_{\normtwo{\bv}=1} L(\bz^\toptzero_n + \eta \bv) 
\approx\mathop{\argmin}_{\normtwo{\bv}=1}
 \left\{L(\bz^\toptzero_n ) + \eta \bv^\top \nabla  L(\bz^\toptzero_n )\right\}.
$$
This is known as a \textit{greedy search}. The minimum occurs when $\bv$ points in the direction opposite to the gradient:
$$
\bv = -\nabla L(\bz^\toptzero_n )\big/{\normtwobig{\nabla L(\bz^\toptzero_n )}}.
$$
Substituting back, we obtain the \textit{gradient descent (GD)} update:
$$
\bz^\toptone_n =\bz^\toptzero_n + \eta \bv = \bz^\toptzero_n - \eta {\nabla L(\bz^\toptzero_n )}\big/{\normtwobig{\nabla L(\bz^\toptzero_n )}},
$$
which is commonly referred to as   the \textit{gradient descent} (GD). 
Similarly, for the movie factors:
$$
\bw^\toptone_m =\bw^\toptzero_m + \eta \bv = \bw^\toptzero_m - \eta {\nabla L(\bw^\toptzero_m )}\big/{\normtwobig{{\nabla L(\bw^\toptzero_m )}}}.
$$
Algorithm~\ref{alg:als-regularizer-missing-stochas-gradient} presents the resulting procedure, which replaces exact ALS updates with normalized gradient steps.

It's noteworthy that the ALS without GD (Algorithm~\ref{alg:als-regularizer})  lacks explicit parameters like step size. 
This characteristic can be both advantageous and disadvantageous. 
On one hand, it absolves the  user from the time-consuming task of fine-tuning parameters, making the method more accessible and less demanding. 
On the other hand, this absence of adjustable parameters also restricts the user's control to directly influence the progression of the algorithm, leaving the convergence of ALS entirely contingent upon the inherent structure of the optimization problem at hand.

In practice, it is common to combine pure ALS iterations with gradient-based variants. The latter provide flexibility through tunable step sizes ($\eta_z, \eta_w$), enabling finer control over convergence speed and stability.

\begin{algorithm}[h] 
\caption{Alternating Least Squares with Full Entries and Gradient Descent}
\label{alg:als-regularizer-missing-stochas-gradient}
\begin{algorithmic}[1] 
\Require Matrix $\bA\in \real^{M\times N}$;
\State Initialize $\bW\in \real^{M\times K}$, $\bZ\in \real^{K\times N}$ \textbf{randomly, without requiring any condition on rank or the relationship among  $M, N, K$}; 
\State Choose a stopping criterion on the approximation error $\delta$;
\State Choose regularization parameters $\lambda_w, \lambda_z$, and step sizes $\eta_w, \eta_z$;
\State Choose the maximum number of iterations $C$;
\State $iter=0$; \Comment{Count for the number of iterations}
\While{$\normf{\bA- (\bW\bZ)}^2>\delta $ and $iter<C$}
\State $iter=iter+1$; 
\For{$n=1,2,\ldots, N$}
\State $\bz^\toptone_n \leftarrow\bz^\toptzero_n - \eta_z {\nabla L(\bz^\toptzero_n )}\big/{\normtwobig{{\nabla L(\bz^\toptzero_n )}}}$; \Comment{$n$-th column of $\bZ$}
\EndFor

\For{$m=1,2,\ldots, M$}
\State $\bw^\toptone_m  \leftarrow \bw^\toptzero_m - \eta_w {\nabla L(\bw^\toptzero_m )}\big/{\normtwobig{{\nabla L(\bw^\toptzero_m )}}}$;\Comment{$m$-th column of $\bW^\top$}
\EndFor
\EndWhile
\State Output $\bW^\top=[\bw_1, \bw_2, \ldots, \bw_M],\bZ=[\bz_1, \bz_2, \ldots, \bz_N]$;
\end{algorithmic} 
\end{algorithm}

\index{Level curves}
\index{Level surfaces}
\paragrapharrow{Geometrical interpretation of gradient descent.} 
\begin{lemma}[Direction of Gradients]\label{lemm:direction-gradients}
The gradient of a differentiable function is orthogonal to its level curves (or level surfaces in higher dimensions).
\end{lemma}
\begin{proof}[of Lemma~\ref{lemm:direction-gradients}, the informal proof]
Consider a two-dimensional level curve defined by  $f(x,y)=c$. 
Assuming sufficient smoothness, we can locally express $y$  as a function of $x$, i.e.,  $y=y(x)$,~\footnote{This is known as the implicit function  theorem, under the conditions of the nonzero partial derivative and smoothness.}, so that 
$
f(x, y(x)) = c.
$
Differentiating both sides with respect to $x$ using the chain rule gives:
$$
\frac{\partial f}{\partial x} \underbrace{\frac{dx}{dx}}_{=1} + \frac{\partial f}{\partial y} \frac{dy}{dx}=0
\gap \implies \gap 
\left\langle \frac{\partial f}{\partial x}, \frac{\partial f}{\partial y}\right\rangle
\cdot 
\left\langle \frac{dx}{dx}, \frac{dy}{dx}\right\rangle=0.
$$
Thus, the gradient is perpendicular to the tangent vector of the level curve.

In full generality, consider the level curve of a vector $\bx\in \real^N$: $f(\bx) = f(x_1, x_2, \ldots, x_N)=c$. Each variable $x_n$ can be regarded as a function of a parameter $t$ on the level curve $f(\bx)=c$: $f(x_1(t), x_2(t), \ldots, x_N(t))=c$. Differentiating the equation with respect to $t$ using the chain rule:
$$
\frac{\partial f}{\partial x_1} \frac{dx_1}{dt} + \frac{\partial f}{\partial x_2} \frac{dx_2}{dt}
+\ldots + \frac{\partial f}{\partial x_N} \frac{dx_N}{dt}
=0.
$$
Thus, the gradient is perpendicular to the tangent in the $N$-dimensional case:
$$
\left\langle \frac{\partial f}{\partial x_1}, \frac{\partial f}{\partial x_2}, \ldots, \frac{\partial f}{\partial x_N}\right\rangle
\cdot 
\left\langle \frac{dx_1}{dt}, \frac{dx_2}{dt}, \ldots \frac{dx_N}{dt}\right\rangle=0.
$$
This completes the proof.
\end{proof}

\begin{figure}[h]
\centering  
\vspace{-0.15cm}    
\subfigtopskip=2pt  
\subfigbottomskip=2pt 
\subfigcapskip=-5pt  
\subfigure[A two-dimensional convex function $L(\bz)$.]{\label{fig:alsgd1}
\includegraphics[width=0.47\linewidth]{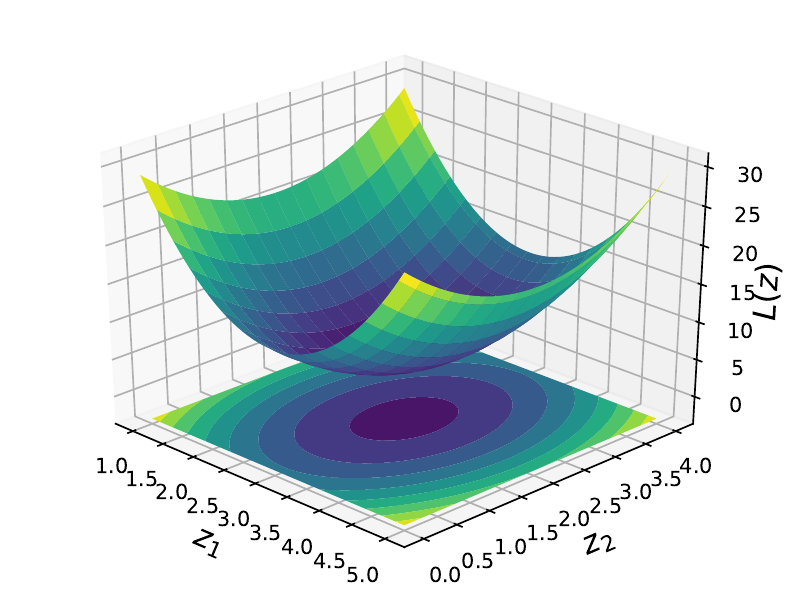}}
\subfigure[Level curve $L(\bz)=c$ and descent direction.]{\label{fig:alsgd2}
\includegraphics[width=0.44\linewidth]{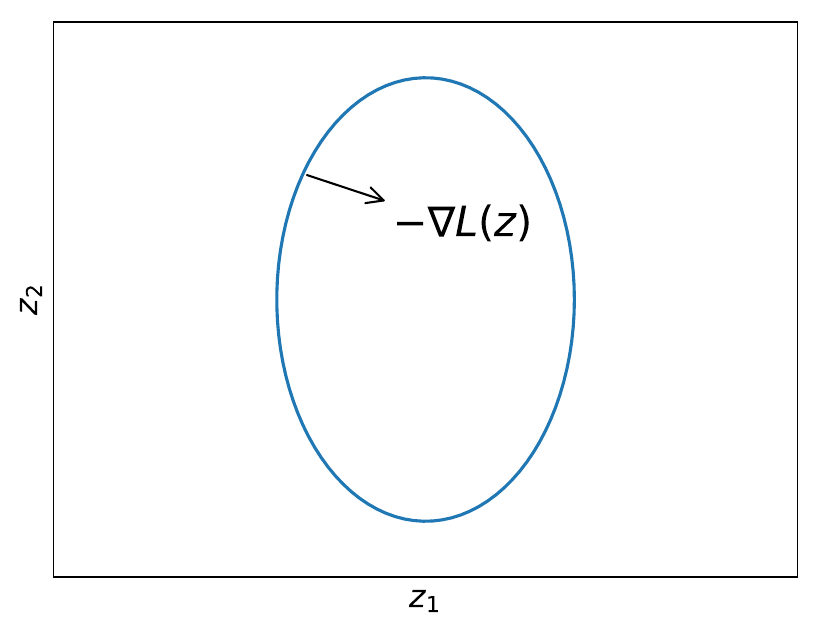}}
\caption{Figure~\ref{fig:alsgd1} shows  surface and  contour plots of a convex function (\textcolor{mydarkblue}{blue}=low, \textcolor{mydarkyellow}{yellow}=high), where the upper graph is the surface plot, and the lower one is its projection  (i.e., contour). Figure~\ref{fig:alsgd2} illustrates that the negative gradient  $-\nabla L(\bz)$ is orthogonal to the level curve and points toward decreasing values of $L(\bz)$.}
\label{fig:alsgd-geometrical}
\end{figure}

This lemma provides a geometric foundation for gradient descent. Since the gradient points in the direction of steepest ascent, moving in the opposite direction---i.e., $-\nabla L(\bz)$---ensures the fastest local decrease in the loss. Figure~\ref{fig:alsgd-geometrical} illustrates this in two dimensions: the negative gradient pushes the iterate toward lower values of the convex function $L(\bz)$.

\index{Convex function}
\index{Contour plot}
\index{Contour plot}
\index{Regularization}

\index{Geometrical interpretation}
\index{Projection gradient descent}
\index{Overfitting}
\section{Regularization: A Geometrical Interpretation}\label{section:reg_geom_inter}
\begin{figure}[h]
\centering
\includegraphics[width=0.95\textwidth]{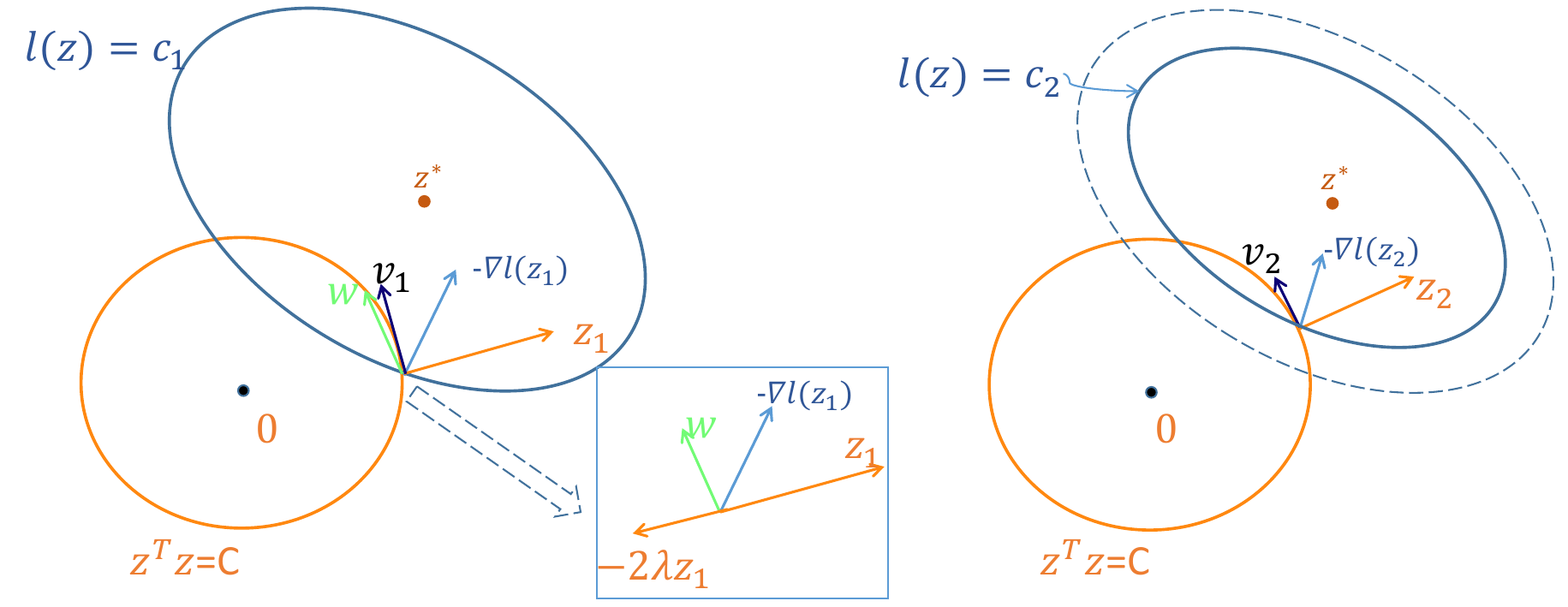}
\caption{Constrained gradient descent under the constraint $\bz^\top\bz\leq C$. 
The \textcolor{mydarkgreen}{green} vector $\bw$ represents the projection of $\bv_1$ onto the set $\bz^\top\bz\leq C$, where $\bv_1$ is the component of $-\nabla l(\bz)$ that is perpendicular to $\bz_1$. 
The right panel shows the next step after the update from the left. 
Here, $\bz^\star$ denotes unconstrained minimizer of \{$\min l(\bz)$\}.}
\label{fig:alsgd3}
\end{figure}
In Section~\ref{section:regularization-extention-general}, we discussed how regularization extends the ALS algorithm to general matrices.
Gradient descent offers a geometric interpretation of this regularization.
To avoid confusion, let: $l(\bz): \real^N\rightarrow \real$ denote the unregularized loss, and  $L(\bz) = l(\bz)+\lambda_z \normtwo{\bz}^2$ denote the regularized loss, with $\lambda_z>0$. 
When minimizing $l(\bz)$, standard gradient descent searches over the entire space $\real^N$. 
However, in machine learning, this can lead to overfitting, as the solution may fit noise rather than underlying patterns.
A common remedy is to constrain the search space---for example, by requiring $\bz^\top\bz < C$ for some constant $C>0$. This leads to the constrained optimization problem:
\begin{equation}
\mathop{\arg\min}_{\bz} \,\, l(\bz), \gap \text{s.t.,} \gap \bz^\top\bz\leq C.
\end{equation}
In unconstrained gradient descent, we update $\bz$ as: $\bz\leftarrow \bz-\eta \nabla l(\bz)$ for a small step size $\eta>0$.
Suppose the current iterate is $\bz_1$, lying at the intersection of the level curve $l(\bz)=c_1$ and the boundary $\bz^\top\bz=C$ (see the left panel of Figure~\ref{fig:alsgd3}). By Lemma~\ref{lemm:direction-gradients}, the descent direction $-\nabla l(\bz_1)$ is perpendicular to the level curve $l(\bz)=c_1$.
However, if we enforce the constraint $\bz^\top\bz\leq C$, a naive step in the direction $-\nabla l(\bz_1)$ would move the next iterate $\bz_2=\bz_1-\eta \nabla l(\bz_1)$ outside the feasible region.
To address this, we decompose the gradient into components normal and tangential to the constraint boundary:
$$
-\nabla l(\bz_1) = a\bz_1 + \bv_1,
$$ 
where $a\bz_1$ is normal (radial) to the sphere $\bz^\top\bz=C$, and $\bv_1$ is tangential (parallel) to the sphere. 
By taking only the tangential component $\bv_1$,  we stay on (or near) the constraint surface. The update becomes:
$$
\bz_2 = \text{project}(\bz_1+\eta \bv_1) = \text{project}\bigg(\bz_1 + \eta
\underbrace{(-\nabla l(\bz_1) -a\bz_1)}_{\bv_1}\bigg),~\footnote{where the operation project($\bx$) 
will project the vector $\bx$ to the closest point inside $\bz^\top\bz\leq C$. Notice here the unprojected update $\bz_2 = \bz_1+\eta \bv_1$ can still make $\bz_2$ fall outside the curve of $\bz^\top\bz\leq C$.}
$$
This method is known as \textit{projected gradient descent (PGD)}. As illustrated in Figure~\ref{fig:alsgd3} (left), the resulting update direction corresponds to a vector $\bw$ (shown in \textcolor{mydarkgreen}{green}) such that $\bz_2=\bz_1+\bw$ remains feasible.
Interestingly, this projected update is equivalent to performing standard gradient descent on the regularized loss  $L(\bz)=l(\bz)+\lambda\normtwo{\bz}^2$ for some $\lambda$.
Indeed, the gradient of $L$ is: $\nabla L(\bz)=\nabla l(\bz)+2\lambda{\bz}$, 
so the negative gradient is:
$$
\begin{aligned}
\bw=-\nabla L(\bz) &= -\nabla l(\bz) - 2\lambda \bz 
\gap \implies \gap
\bz_2 = \bz_1+ \bw =\bz_1 -  \nabla L(\bz).
\end{aligned}
$$
And in practice, a small step size $\eta$ prevents the trajectory from moving outside the curve of $\bz^\top\bz\leq C$:
$$
\bz_2  =\bz_1 -  \eta\nabla L(\bz),
$$
which aligns with the regularization term discussed in Section~\ref{section:regularization-extention-general}.

\index{Sparsity}
\index{$\ell_1$-norm}
\begin{figure}[h]
\centering
\includegraphics[width=0.95\textwidth]{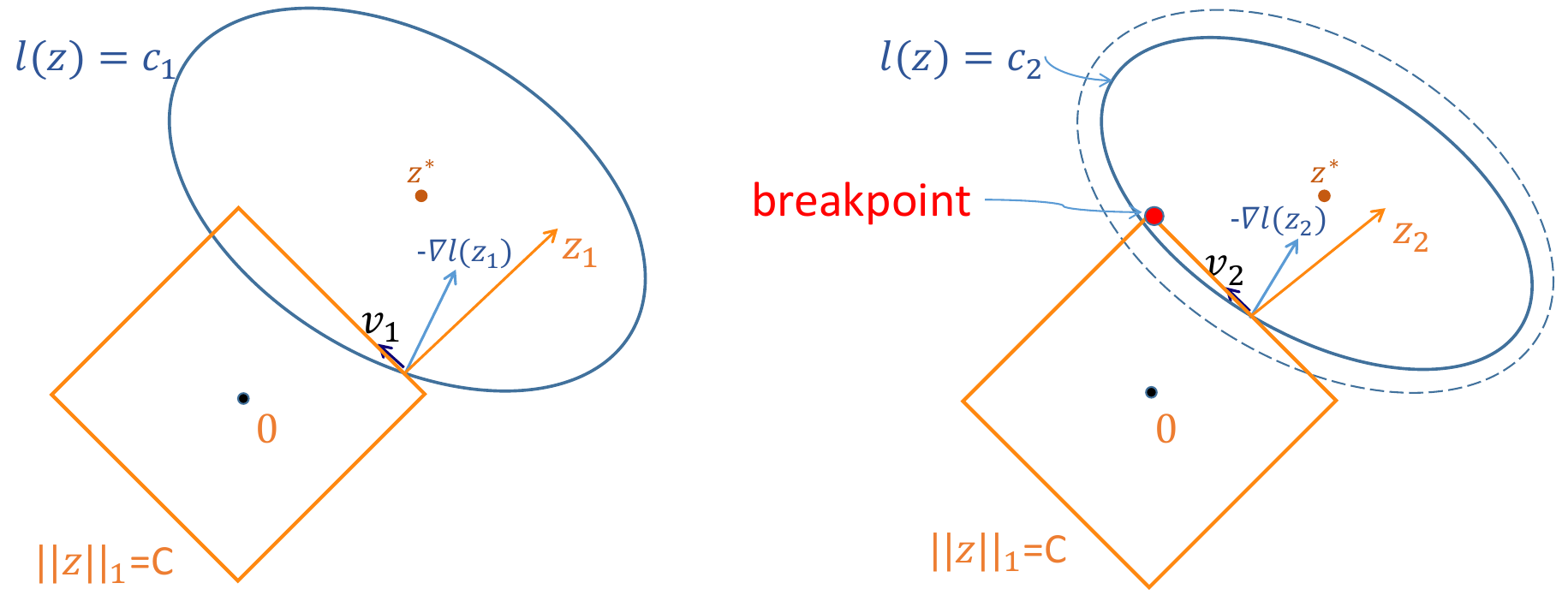}
\caption{Constrained gradient descent with $\normone{\bz}\leq C$, where the \textcolor{winestain}{red} dot denotes the breakpoint in the $\ell_1$-norm. 
The right panel shows the next step after the update from the left.
Here, $\bz^\star$ denotes the unconstrained minimizer of \{$\min l(\bz)$\}.}
\label{fig:alsgd4}
\end{figure}
\paragrapharrow{Sparsity.}
In some applications, we seek a sparse solution---i.e., a vector $\bz$ with many zero entries---that still minimizes $l(\bz)$. 
image is reconstructed using only a few active components.
To promote sparsity, we constrain the solution to the $\ell_1$-ball: $\normone{\bz} \leq C$, where $\normone{\cdot}$ is the $\ell_1$-norm of a vector or a matrix. 
As with the $\ell_2$ case, gradient descent under this constraint tends to push the iterates toward the boundary $\normone{\bz}=C$. However, unlike the smooth $\ell_2$-ball, the $\ell_1$-ball has sharp corners (breakpoint, see right panel of Figure~\ref{fig:alsgd4}). When the gradient descent trajectory hits such a corner, the solution often becomes exactly sparse---i.e., one or more components of $\bz$ become precisely zero.
In high dimensions, this effect is amplified: many coordinates are driven to zero because the geometry of the 
$\ell_1$-ball favors solutions aligned with the coordinate axes. This is why $\ell_1$-regularization  is widely used to induce sparsity in machine learning and signal processing.

\index{Stochastic gradient descent}
\index{Stochastic coordinate descent}
\section{Stochastic Gradient Descent}
The gradient descent (GD) method is a powerful optimization algorithm; however, it has notable limitations in practical settings---particularly when applied to large-scale problems.
To understand these issues, consider the mean squared error  derived from  Equation~\eqref{equation:als-per-example-loss2}:
\begin{equation}\label{equation:als-per-example-loss_mse}
\frac{1}{MN}\mathop{\min}_{\bW,\bZ}  \sum_{n=1}^N \sum_{m=1}^{M} \left(a_{mn} - \bw_m^\top\bz_n\right)^2.
\end{equation}
This objective requires computing the residual $e_{mn} = (a_{mn} - \bw_m^\top\bz_n)^2$ for each observed entry $a_{mn}$, which measures the squared difference between the predicted and actual values. 
The total sum of squared residuals is denoted by $e = \sum_{m,n=1}^{MN}e_{mn}$.
When the number of training entries $MN$ is large, evaluating the full gradient becomes computationally expensive, significantly slowing down each iteration. Moreover, gradients computed from different samples may partially cancel each other out, leading to small net parameter updates and slow convergence.
To address these challenges, researchers introduced stochastic gradient descent (SGD). Instead of computing the exact gradient over the entire dataset---which is costly---SGD approximates the gradient using a single randomly selected sample at each iteration. This estimate is then used to update the parameters in a direction that reduces the loss. Although noisy, this approximation is computationally efficient and often sufficient for convergence, especially on large datasets.

Consider the regularized  loss:
$$
L(\bW,\bZ)= \frac{1}{2} \sum_{n=1}^N \sum_{m=1}^{M} \left(a_{mn} - \bw_m^\top\bz_n\right)^2 +\frac{1}{2} \lambda_w\sum_{m=1}^{M}\normtwo{\bw_m}^2 +\frac{1}{2}\lambda_z\sum_{n=1}^{N}\normtwo{\bz_n}^2.
$$
Minimizing the overall loss $L(\bW,\bZ)$ can be achieved by iteratively reducing the per-example loss term: $l(\bw_m, \bz_n)=\frac{1}{2}\left(a_{mn} - \bw_m^\top\bz_n\right)^2+\frac{1}{2} \lambda_w\normtwo{\bw_m}^2 +\frac{1}{2}\lambda_z\normtwo{\bz_n}^2$ for all $m\in \{1,2,\ldots, M\}, n\in\{1,2,\ldots,N\}$. 
This strategy is also known as \textit{stochastic coordinate descent}, as it updates one pair 
$(\bw_m, \bz_n)$ at a time.
The gradients of $l(\cdot,\cdot)$ with respect to $\bw_m$ and $\bz_n$,  along with their closed-form solutions, are:
$$
\left\{
\begin{aligned}
\nabla_{\bz_n} l(\bz_n) &= \bw_m\bw_m^\top \bz_n +\lambda_z\bz_n  -a_{mn} \bw_m 
 &\implies&\,\, \bz_n= a_{mn}(\bw_m\bw_m^\top+\lambda_z\bI)^{-1}\bw_m;\\
\nabla_{\bw_m} l(\bw_m) &= \bz_n\bz_n^\top\bw_m +\lambda_w\bw_m - a_{mn}\bz_n &\implies&\,\, \bw_m= a_{mn}(\bz_n\bz_n^\top+\lambda_w\bI)^{-1}\bz_n.
\end{aligned}
\right.
$$
Alternatively, we can apply gradient descent using the per-example loss. Since each update is based on a single data point, this approach is referred to as \textit{stochastic gradient descent (SGD)}:
$$
\left.
\begin{aligned}
\bz_n&\leftarrow \bz_n - \eta_z \frac{\nabla_{\bz_n} l(\bz_n)}{\normtwo{\nabla_{\bz_n} l(\bz_n)}}
\qquad \text{and}\qquad 
\bw_m\leftarrow \bw_m - \eta_w \frac{\nabla_{\bw_m} l(\bw_m)}{\normtwo{\nabla_{\bw_m} l(\bw_m)}}.
\end{aligned}
\right.
$$
This SGD-based update for ALS is formalized in Algorithm~\ref{alg:als-regularizer-missing-stochas-gradient-realstoch}.
In practice, the indices $m$ and $n$ are typically chosen randomly at each step---hence the term \textit{stochastic}. If instead they are cycled through in a fixed order, the method is sometimes called \textit{incremental gradient descent}.
It is also worth noting that both GD and SGD may fail to converge if the learning rate is too large. In such cases, re-running the algorithm with a smaller step size often resolves the issue.

\begin{algorithm}[h] 
\caption{Alternating Least Squares with Full Entries and SGD}
\label{alg:als-regularizer-missing-stochas-gradient-realstoch}
\begin{algorithmic}[1] 
\Require  Matrix $\bA\in \real^{M\times N}$;
\State Initialize $\bW\in \real^{M\times K}$, $\bZ\in \real^{K\times N}$ \textbf{randomly, without requiring any condition on rank or the relationship among  $M, N, K$}; 
\State Choose a stopping criterion on the approximation error $\delta$;
\State Choose regularization parameters $\lambda_w, \lambda_z$, and step size $\eta_w, \eta_z$;
\State Choose the maximum number of iterations $C$;
\State $iter=0$; \Comment{Count for the number of iterations}
\While{$\normf{  \bA- (\bW\bZ)}^2>\delta $ and $iter<C$}
\State $iter=iter+1$; 
\For{$n=1,2,\ldots, N$}
\For{$m=1,2,\ldots, M$} \Comment{in practice, $m,n$ can be randomly produced}
\State $\bz_n\leftarrow \bz_n - \eta_z {\nabla l(\bz_n)}/{\normtwo{\nabla l(\bz_n)}}$;\Comment{$n$-th column of $\bZ$}
\State $\bw_m\leftarrow \bw_m - \eta_w {\nabla l(\bw_m)}/{\normtwo{\nabla l(\bw_m)}}$;\Comment{$m$-th column of $\bW^\top$}
\EndFor
\EndFor

\EndWhile
\State Output $\bW^\top=[\bw_1, \bw_2, \ldots, \bw_M],\bZ=[\bz_1, \bz_2, \ldots, \bz_N]$;
\end{algorithmic} 
\end{algorithm}

\section{Bias Term}

\begin{figure}[htp]
\centering
\includegraphics[width=0.95\textwidth]{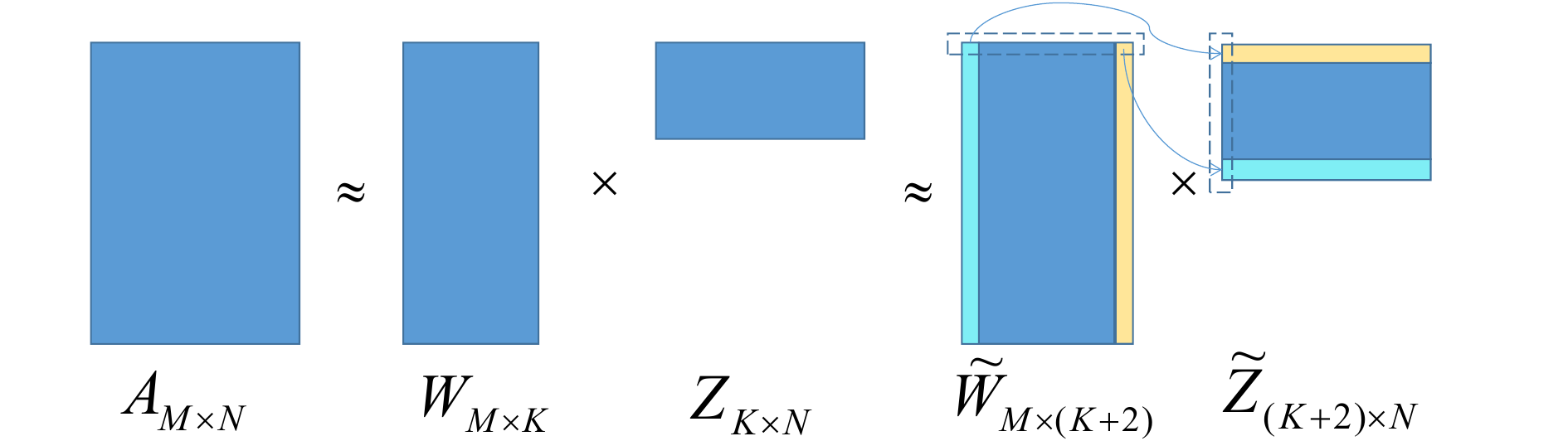}
\caption{Bias terms in alternating least squares, where the \textcolor{mydarkyellow}{yellow} entries denote ones (which are fixed), and the \textcolor{cyan}{cyan} entries denote the added features to fit the bias terms. The dotted boxes illustrate how bias terms operate in the factorization.}
\label{fig:als-bias}
\end{figure}
In ordinary least squares, a bias term (or intercept) is commonly added to improve model flexibility, as shown in Equation~\eqref{equation:ls-bias}. A similar idea applies to ALS.
Specifically, we can incorporate global, user-specific, and item-specific biases by augmenting the factor matrices:
\begin{itemize}
\item Append a fixed column of ones to the last column of $\bW$. To accommodate this, an extra row must be added to the last row of $\bZ$ to model the corresponding bias weights.
\item Similarly, prepend a fixed row of ones to the first row of $\bZ$, and add an extra column to the first column of $\bW$.
\end{itemize}
This construction is illustrated in Figure~\ref{fig:als-bias}.

Let us first consider the update for $\bz_n$. Define the augmented vector:
$$
\widetildebz_n
\triangleq
\begin{bmatrix}
	1\\
	\bz_n
\end{bmatrix}
\in\real^{K+2},
$$
which is the $n$-th column of the extended matrix $\widetildebZ$. 
Assume the extended weight matrix is partitioned as $\widetildebW=[\widebarbw_0 , \widebarbW]$, 
where $\widebarbw_0\in\real^M$ is the first column (modeling the global/user bias), and  $\widebarbW\in\real^{M\times (K+1)}$ contains the remaining latent factors (the last column is the ones vector).
Then, the loss (up to constants) becomes:
\begin{equation}
\footnotesize
\begin{aligned}
2L(\bz_n) 
&=\normf{\widetildebW\widetildebZ-\bA}^2 +\lambda_w \normf{\widetildebW}^2 + \lambda_z \normf{\widetildebZ}^2
= 
\left\Vert
\widetildebW
\begin{bmatrix}
	1 \\
	\bz_n
\end{bmatrix}-\ba_n
\right\Vert_2^2 
+ 
\underbrace{\lambda_z \normtwo{\widetildebz_n}^2}_{=\lambda_z \normtwo{\bz_n}^2+\lambda_z}
+ 
C_{z_n}\\
&= 
\left\Vert
\begin{bmatrix}
	\widebarbw_0 & \widebarbW
\end{bmatrix}
\begin{bmatrix}
	1 \\
	\bz_n
\end{bmatrix}-\ba_n
\right\Vert_2^2 
+ \lambda_z \normtwo{\bz_n}^2 + C_{z_n}
= 
\bigg\Vert
\widebarbW \bz_n - 
\underbrace{(\ba_n-\widebarbw_0)}_{\triangleq\widebarba_n}
\bigg\Vert_2^2 
+ \lambda_z \normtwo{\bz_n}^2 + C_{z_n},
\end{aligned}
\end{equation}
where  $C_{z_n}$ is a constant with respect to $\bz_n$. Let $\widebarba_n \triangleq \ba_n-\widebarbw_0$. 
The update for $\bz_n$ is just similar to the one in Equation~\eqref{als:gradient-regularization-zn}, with the gradient given by
$$
\nabla_{\bz_n} L(\bz_n) = \widebarbW^\top\widebarbW\bz_n - \widebarbW^\top\widebarba_n + \lambda_z\bz_n.
$$
Therefore, the update for $\bz_n$ is given by determining the root of the  gradient above:
$$
\textbf{(update for $\widetildebz_n$)}: \quad \bz_n = (\widebarbW^\top\widebarbW+ \lambda_z\bI)^{-1} \widebarbW^\top \widebarba_n
\gap 
\implies 
\gap 
\widetildebz_n = \begin{bmatrix}
	1\\\bz_n 
\end{bmatrix},
\,\forall n.
$$
Similarly, following the loss with respect to each row of $\bW$ in Equation~\eqref{als:gradient-regularization-wd}, let $\widetildebw_m \triangleq
\scriptsize
\begin{bmatrix}
	\bw_m \\
	1
\end{bmatrix} \in\real^{K+2}$ be the $m$-th row of $\widetildebW$ (or $m$-th column of $\widetildebW^\top$), we have 
\begin{equation}
\footnotesize
\begin{aligned}
&2L(\bw_m ) 
=\normf{\widetildebZ^\top\widetildebW^\top-\bA^\top}^2 +\lambda_w \normf{\widetildebW^\top}^2 + \lambda_z \normf{\widetildebZ}^2
= 
\normtwo{\widetildebZ^\top\widetildebw_m-\bb_m}^2 + 
\underbrace{\lambda_w \normtwo{\widetildebw_m}^2}_{=\lambda_w \normtwo{\bw_m}^2+\lambda_w}
+ 
C_{w_m} \\
&= 
\bigg\Vert
\begin{bmatrix}
\widebarbZ^\top&
\widebarbz_0
\end{bmatrix}
\scriptsize
\begin{bmatrix}
\bw_m \\
1
\end{bmatrix}
\footnotesize
-\bb_m
\bigg\Vert_2^2 
+ 
\lambda_w \normtwo{\bw_m}^2
+ 
C_{w_m}
= 
\left\Vert
\widebarbZ^\top\bw_m 
-(\bb_m-\widebarbz_0) 
\right\Vert_2^2+ 
\lambda_w \normtwo{\bw_m}^2
+ 
C_{w_m},
\end{aligned}
\end{equation}
where $\widebarbz_0$ represents the last column of $\widetildebZ^\top$, $\widebarbZ^\top$ contains the remaining $K+1$ columns of $\widetildebZ^\top$ (i.e., $\widetildebZ^\top\triangleq[\widebarbZ^\top, \widebarbz_0]$),
$C_{w_m}$ is a constant with respect to $\bw_m$.  
$\bW^\top=[\bw_1, \bw_2, \ldots, \bw_M]$ and $\bA^\top=[\bb_1,\bb_2, \ldots, \bb_M]$ are the column partitions of $\bW^\top$ and $\bA^\top$, respectively. Let $\widebarbb_m \triangleq \bb_m-\widebarbz_0$, the update for $\bw_m$ is again just similar to  the one in Equation~\eqref{als:gradient-regularization-wd}, with the gradient given by
$$
\nabla_{\bw_m} L(\bw_m) = \widebarbZ\cdot \widebarbZ^\top\bw_m - \widebarbZ\cdot \widebarbb_m + \lambda_w\bw_m.
$$
Therefore, the update for $\bw_m$ is given by the root of the  gradient above:
$$
\textbf{(update for $\widetildebw_m$)}:\quad 
\bw_m=(\widebarbZ\cdot \widebarbZ^\top+\lambda_w\bI)^{-1}\widebarbZ\cdot \widebarbb_m
\gap \implies \gap 
\widetildebw_m = \begin{bmatrix}
	\bw_m \\ 1
\end{bmatrix}, 
\forall m.
$$
These closed-form updates naturally extend ALS to include bias terms. Similar derivations can be carried out using gradient descent, and the framework readily accommodates missing entries (see Sections~\ref{section:als-gradie-descent} and \ref{section:alt-columb-by-column} for details).

\section{Convergence}\label{section:als_convergence}

We previously noted that the ALS Algorithm~\ref{alg:als} belongs to a class of optimization methods known as block coordinate descent (BCD) (see Algorithm~\ref{alg:two_bcd_gen_inals}, which is known as the alternating method in the EM algorithm context; see Algorithm~\ref{alg:am_lvm}). 
In BCD, the variables are partitioned into blocks, and the algorithm iteratively optimizes the objective function with respect to one block at a time while keeping the others fixed. This strategy is especially advantageous for large-scale problems, where jointly optimizing over all variables is computationally prohibitive or infeasible.
Although BCD does not always guarantee convergence to a global optimum---particularly for non-convex problems---it often converges to a stationary point (i.e., a point where the gradient of the objective vanishes). Moreover, when multiple blocks are independent, their updates can be performed in parallel, further improving efficiency.
In the context of ALS, the factor matrices $(\bW,\bZ)$ naturally form two blocks of variables, which is why this approach is often referred to as 2-block coordinate descent (2-BCD). The convergence properties of this method are established by the following result.

\begin{theoremHigh}[Convergence of 2-BCD Method \citep{grippo2000convergence, beck2017first, jain2017non, gillis2020nonnegative}]
Let the iterates be generated by a 2-BCD algorithm. Then every limit point of the sequence is a stationary point of the objective function, provided that:
\begin{enumerate}
\item  The objective function is continuously differentiable.
\item Each block of variables is constrained to lie in a closed convex set.
\end{enumerate}
\end{theoremHigh}
For standard ALS algorithms---and for most nonnegative matrix factorization (NMF) methods discussed in the next chapter---both conditions are satisfied. Consequently, convergence to a stationary point is guaranteed.
More generally, the convergence of multi-block BCD methods requires additional assumptions, as stated below.

\begin{theoremHigh}[Convergence of BCD Method \citep{bertsekas1997nonlinear, gillis2020nonnegative}]
Consider a BCD algorithm applied to a problem with more than two blocks. Every limit point of the generated sequence is a stationary point, provided that:
\begin{enumerate}
\item  The objective function is continuously differentiable.
\item  Each block of variables belongs to a closed convex set.
\item  For each block, the subproblem solved at every iteration has a unique minimizer.
\item  The objective function value decreases monotonically across successive iterates (i.e., after each block update).
\end{enumerate}
\end{theoremHigh}

\section{Movie Recommender}\label{section:movie_rec_als}
\begin{figure}[h]
\centering  
\vspace{-0.35cm} 
\subfigtopskip=2pt 
\subfigbottomskip=2pt
\subfigcapskip=-5pt
\subfigure[Training set.]{\label{fig:movie100k1}
\includegraphics[width=0.47\linewidth]{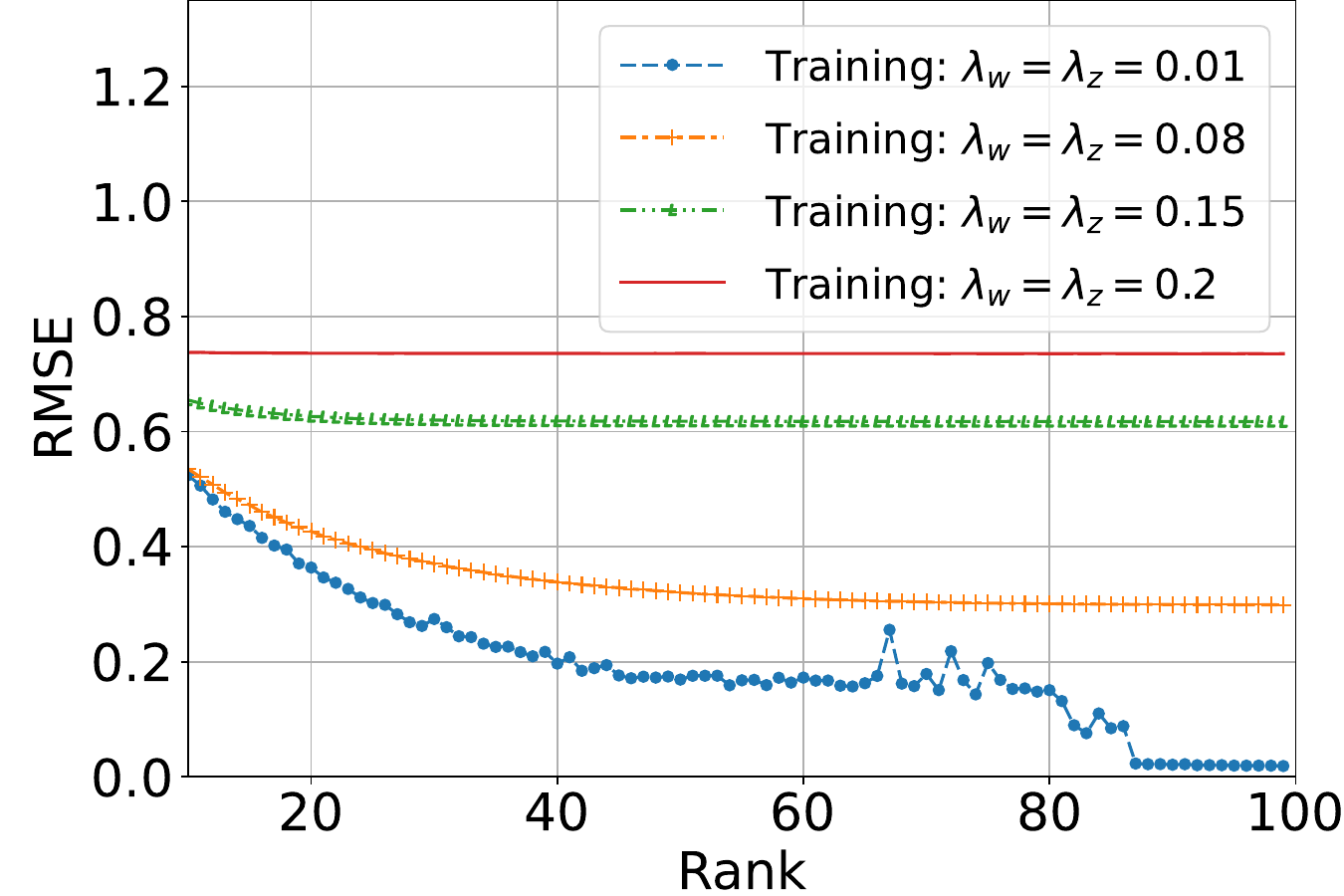}}
\quad 
\subfigure[Validation set.]{\label{fig:movie100k2}
\includegraphics[width=0.47\linewidth]{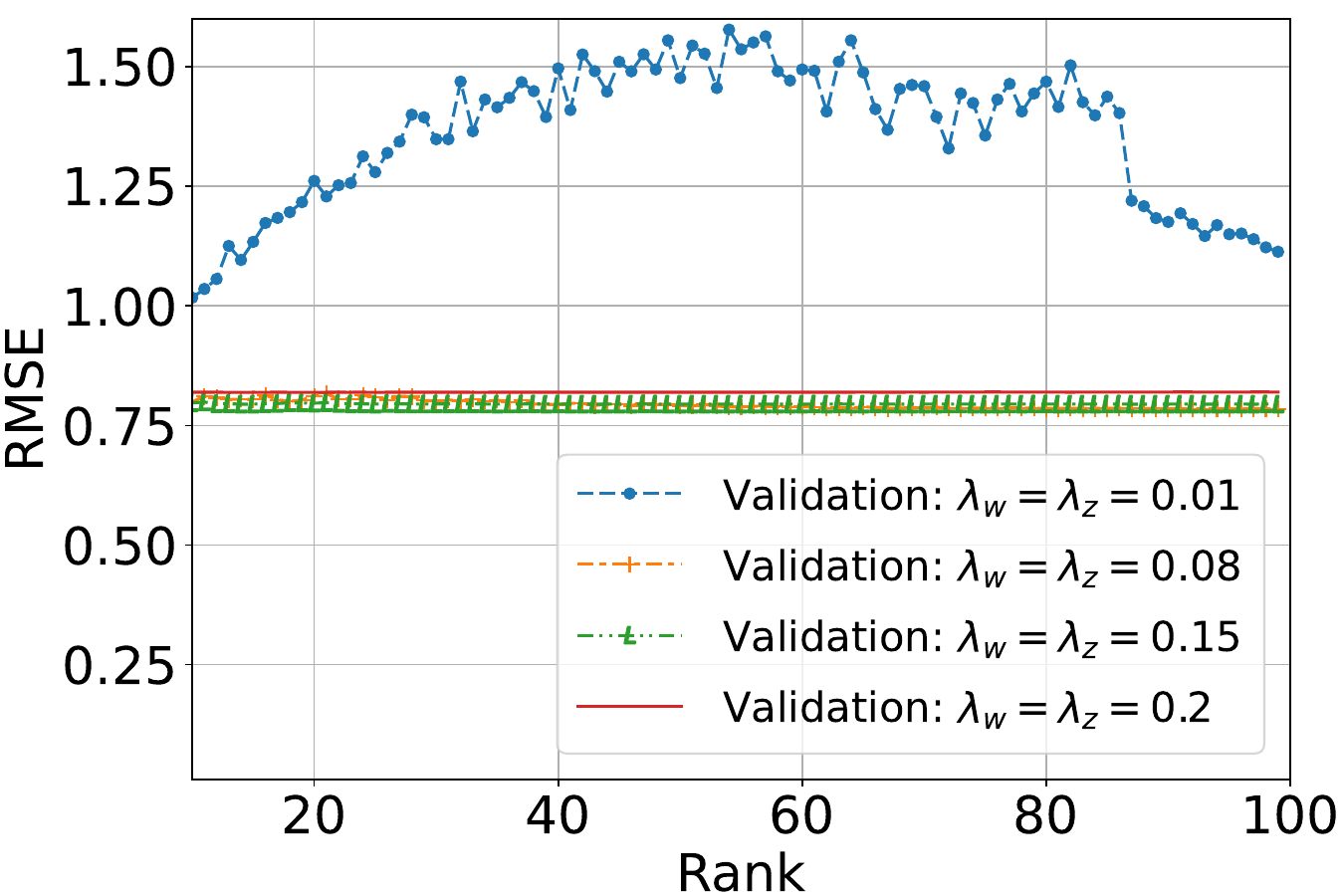}}
\caption{
Training and validation RMSE for the ``MovieLens 100K" dataset across different ranks ($K$) and regularization strengths ($\lambda$).}
\label{fig:movie100k}
\end{figure}

The ALS algorithm has been widely applied in movie recommendation systems.
To illustrate this, we use the ``MovieLens 100K" dataset from MovieLens \citep{harper2015movielens}\footnote{\url{http://grouplens.org}}---a benchmark dataset in recommender systems research due to its rich collection of user-movie ratings.
The dataset contains 100,000 ratings from 943 users on 1,682 movies, with integer ratings ranging from 1 to 5. 
The data was collected via the MovieLens website over a seven-month period, from September 19, 1997, to April 22, 1998. To ensure data quality, users with fewer than 20 ratings or incomplete demographic information were removed. While demographic attributes (age, gender, occupation, ZIP code) are available, our focus here is solely on the raw rating matrix, to evaluate how well the low-rank ALS model captures the underlying preference structure and enables accurate recommendations.
The dataset is split into a training set (95,015 ratings) and a validation set (4,985 ratings). Model performance is measured using the root mean squared error (RMSE), defined as:
$
\text{RMSE}(\bx, \widehatbx) = \sqrt{\frac{1}{N} \sum_{n=1}^{N}(x_n-\widehatx_n)^2},
$
which quantifies the average magnitude of prediction errors.
For the ALS algorithm, the lowest validation RMSE (0.806, i.e., less than 1) is achieved with rank $K=62$ and regularization parameters $\lambda_w=\lambda_z=0.15$, as shown in Figure~\ref{fig:movie100k}. Given that ratings range from 1 to 5, an RMSE below 1 indicates that the model can reliably distinguish between positive (e.g., ratings 4--5) and negative (e.g., ratings 1--2) user preferences.

\paragrapharrow{Recommender 1.}
A simple rule-based system recommends movie $m$ to user $n$ if the predicted rating $\widehata_{mn}\geq4$ and user $n$ has not yet rated movie  $m$. 
\paragrapharrow{Recommender 2.}
Alternatively, we can recommend movies similar to those the user has highly rated. Suppose user $n$ gave movie $m$ a rating of 5 ($a_{mn}=5$). Under the ALS factorization $\bA\approx\bW\bZ$, each row of $\bW$ represents the latent feature vector of a movie (see Section~\ref{section:als-vector-product}). To find recommendations, we identify movies that are most similar to movie $m$ but have not yet been rated by user $n$:
$$
\mathop{\argmax}_{\bw_i} \gap \text{similarity}(\bw_i, \bw_m), \qquad \text{for all} \gap i \notin \bo_n,
$$
where $\bw_i$ is the latent vector of movie $i$, and $\bo_n$ denotes the set of movies already rated by user $n$.
This approach relies on a vector similarity measure. The most common choice is \textit{cosine similarity}, defined as:
$$
\cos(\bx, \by) = \frac{\bx^\top\by}{\normtwo{\bx}\cdot \normtwo{\by}},
$$
Cosine similarity ranges from $-1$ (completely dissimilar) to $1$ (identical direction), and depends only on the angle between vectors---not their magnitudes---since it operates on normalized vectors.
Another widely used measure is \textit{Pearson correlation}:
$$
\text{Pearson}(\bx,\by) =\frac{\Cov(\bx,\by)}{\sigma_x \cdot \sigma_y}
= \frac{\sum_{n=1}^{N} (x_n - \bar{x} ) (y_n -\bar{y})}{ \sqrt{\sum_{n=1}^{N} (x_n-\bar{x})^2}\sqrt{ \sum_{n=1}^{N} (y_n-\bar{y})^2 }}.
$$
Like cosine similarity, Pearson correlation ranges from $-1$ to $1$, with $0$ indicating no linear relationship, $-1$ indicating perfectly dissimilarity, and 1 denoting perfectly similarity. It is commonly used to assess linear dependence in statistics and regression.

Both measures are prevalent in machine learning: Pearson correlation is often used in statistical modeling, while cosine similarity dominates in recommendation systems and information retrieval due to its robustness to magnitude differences.
In our experiments, cosine similarity yields better performance, as confirmed by \textit{precision-recall (PR)} curve analysis.

\begin{figure}[h]
\centering
\vspace{-0.35cm}
\subfigtopskip=2pt
\subfigbottomskip=2pt
\subfigcapskip=-5pt
\subfigure[Cosine Bin Plot.]{\label{fig:als-cosine}%
\includegraphics[width=0.32\linewidth]{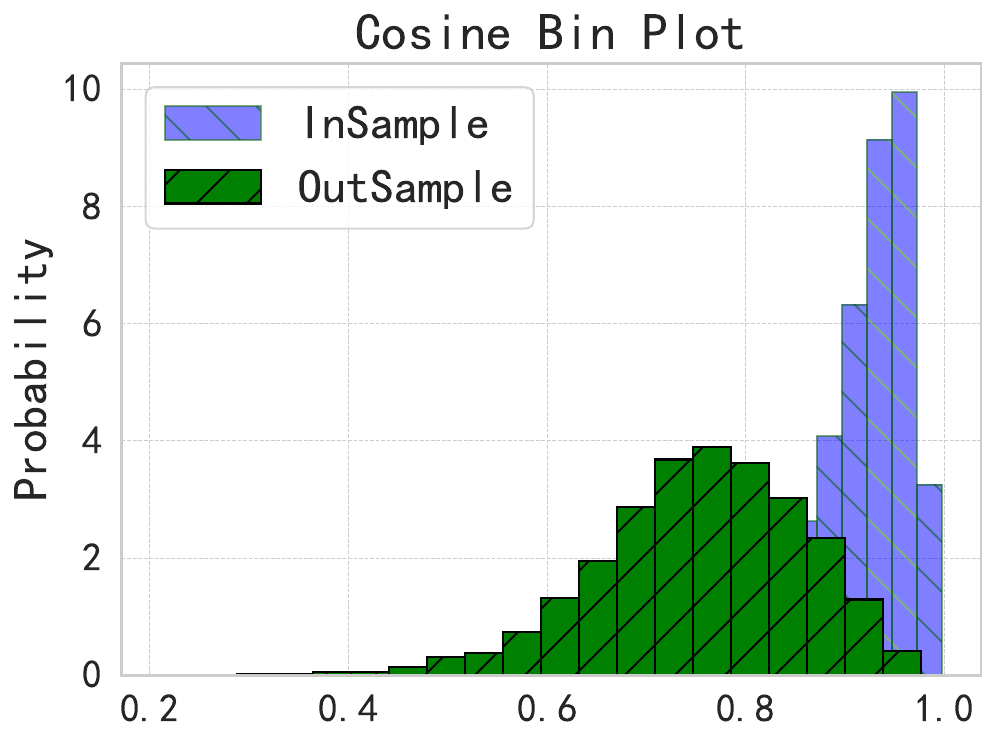}}%
\subfigure[Pearson Bin Plot.]{\label{fig:als-pearson}%
\includegraphics[width=0.32\linewidth]{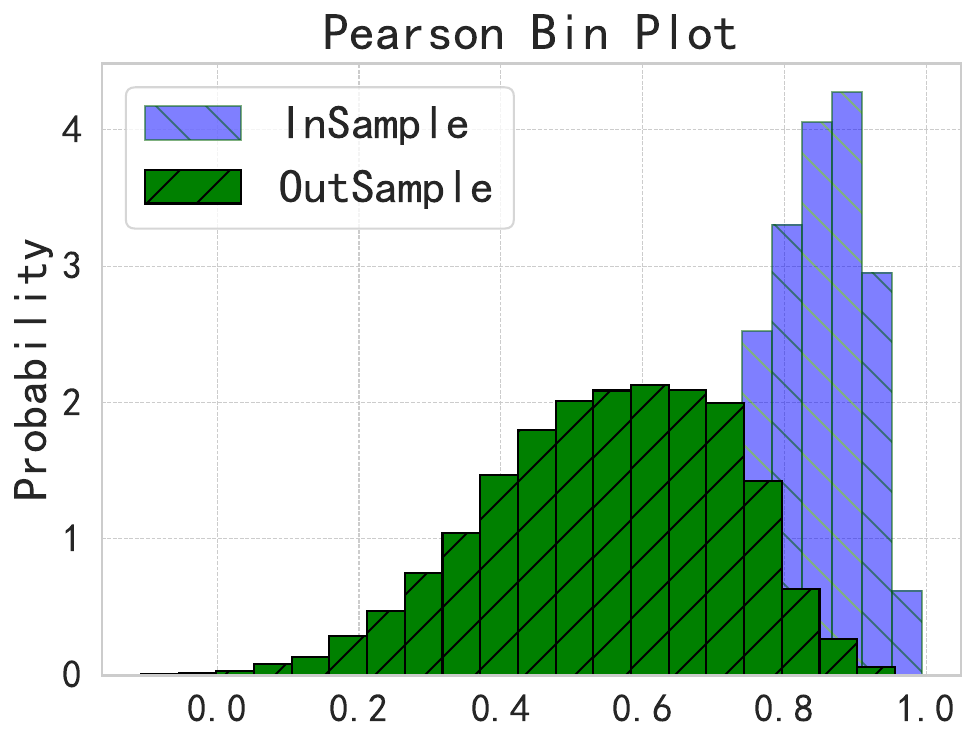}}%
\subfigure[PR Curve.]{\label{fig:als-prcurve}%
\includegraphics[width=0.35\linewidth]{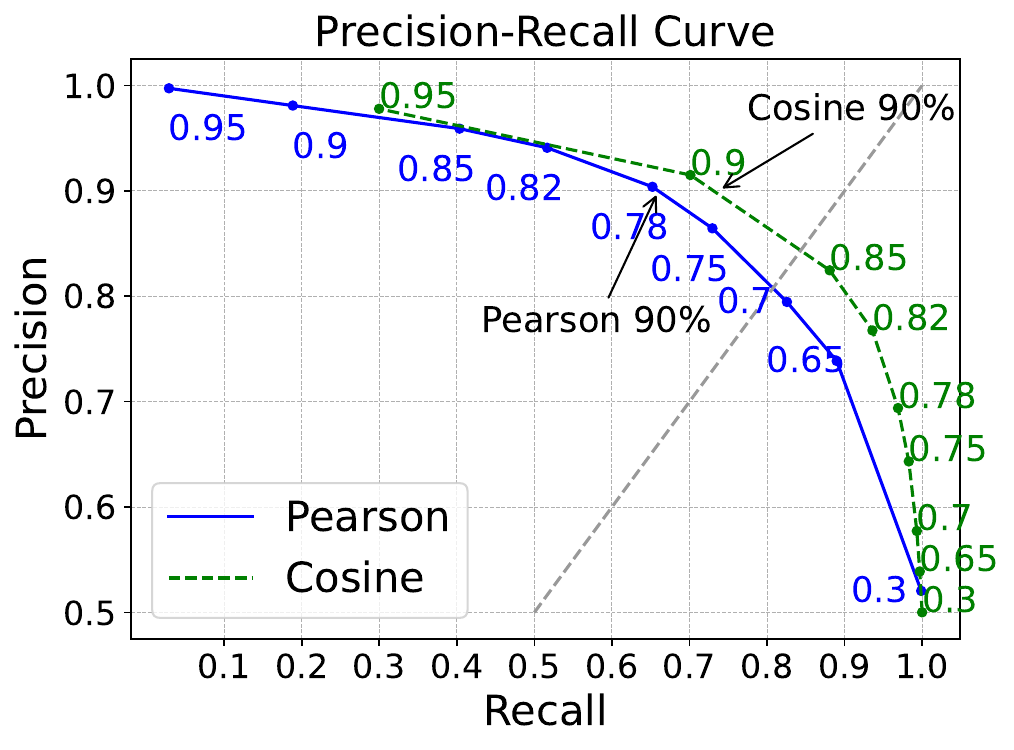}}%
\caption{Distribution of the insample and outsample similarities using cosine and Pearson similarities, and the Precision-Recall curves for both.}
\label{fig:als-prcurive-bin}
\end{figure}

Using the same MovieLens 100K setup ($K=62$, $\lambda_w=\lambda_z=0.15$), we analyze whether ALS can distinguish between movies that users rate highly versus poorly (i.e., the Recommender 2 context). 
We define
``insample" as the similarity (of movie latent vectors) between  pairs of movies both rated 5 by the same user, and
``outsample" as the similarity between a movie rated 5 and another rated 1 by the same user.
Figures~\ref{fig:als-cosine} and~\ref{fig:als-pearson} show the distributions of these similarities under cosine and Pearson measures, respectively. In both cases, the insample and outsample distributions are clearly separated---demonstrating that ALS successfully learns meaningful latent features that reflect user preferences.
Figure~\ref{fig:als-prcurve} presents the precision-recall curves. Cosine similarity achieves over 73\% recall at 90\% precision, whereas Pearson correlation reaches only about 64\% recall at the same precision level. This confirms that cosine similarity is better suited for this recommendation task.
Other similarity measures, such as \textit{negative Euclidean distance}, could also be explored. 
While Euclidean distance quantifies dissimilarity, its negative can serve as a similarity score---though it is sensitive to vector magnitudes and less common in collaborative filtering.

\paragrapharrow{Explicit vs. implicit feedback.}
The ALS method described above is designed for \textit{explicit feedback}, where user ratings carry clear semantic and hierarchical meaning (e.g., ``I like this movie"; higher value indicates more preference).
In contrast, many real-world systems rely on \textit{implicit feedback}, where preferences are inferred from user behavior---such as clicks, views, purchases, or time spent on a page. These signals are abundant but noisy, as they do not directly indicate dislike (e.g., a user may simply not see an item).
To handle implicit data, ALS can be extended in several ways:
using a dictionary-based transformation to map interactions into latent user/item representations \citep{he2017neural};
incorporating multinomial priors into variational autoencoders (VAE, see Sections~\ref{section:vae_pca} and~\ref{section:vae_multino_gen});
leveraging probabilistic models that explicitly account for uncertainty in implicit signals \citep{liang2018variational}.
These extensions enhance ALS's flexibility, enabling effective recommendations even when explicit ratings are unavailable.

\begin{problemset}

\item \label{prob:als_pseudo1} \textbf{Least squares for rank-deficiency \citep{lu2021rigorous}.} Let $\bA\in\real^{M\times N}$ and $\bb\in\real^M$. Show that the least squares problem $L(\bx)=\normtwo{\bA\bx-\bb}^2$ admits a minimizer $\bx^*\in\real^N$ if and only if there exists a vector $\by\in\real^N$ such that $\bx^*=\bA^+\bb+(\bI-\bA^+\bA)\by$, where $\bA^+$ is the \textit{pseudo-inverse} of $\bA$ \citep{lu2021numerical}. 
\begin{itemize}
\item This shows that the least squares has a \textbf{unique} minimizer of $\bx^*=\bA^+\bb$ only when $\bA^+$ is a left inverse of $\bA$ (i.e., $\bA^+\bA=\bI_N$). The solution in Lemma~\ref{lemma:ols} corresponds to this special case.
\item The minimal value of the objective is $L(\bx^*)=\bb^\top(\bI-\bA\bA^+)\bb$.
\item If $\by\neq \bzero$, then  $\normtwo{\bA^+\bb}\leq \normtwo{\bA^+\bb+(\bI-\bA^+\bA)\by}$.
\end{itemize}
\textit{Hint: Use SVD (Theorem~\ref{theorem:reduced_svd_rectangular}).}

\item \label{prob:als_pseudo2} \textbf{Least squares for rank-deficiency.} Let  $\bA\in\real^{M\times N}$ and $\bB\in\real^{M\times P}$. Show that the least squares problem $L(\bX) = \normf{\bA\bX-\bB}^2$ has a minimizer $\bX^*=\bA^+\bB\in\real^{N\times P}$. Determine all minimizers by applying the result from Problem~\ref{prob:als_pseudo1}.

\item \label{prob:als_pseudon} \textbf{Least squares for rank-deficiency.}  Let  $\bA\in\real^{M\times N}$ and $\bB\in\real^{P\times N}$. Show that the least squares problem $L(\bX) = \normf{\bX\bA-\bB}^2$ has a minimizer $\bX^*=\bB\bA^+\in\real^{p\times m}$.
	
\item Prove Lemma~\ref{lemma:als-update-w-rank}.

\item \label{prob:sep_conv} \textbf{Marginally convex.} Let $D(\bA, \bB)$ be convex in its  second argument $\bB$. Show that $D(\bA,\bW\bZ)$ is convex in $\bW$ when $\bZ$ is fixed, and convex in $\bZ$ when $\bW$ is fixed.

\item Show that any function that is jointly convex in its arguments is necessarily marginally convex in each argument.


\item \label{prob:non_iid_gaus} \textbf{Weighted $\ell_2$ loss from non i.i.d. Gaussian noise.} Suppose the Gaussian noise in \eqref{equation:gau_noise} is not i.i.d. 
Discuss the likelihood function for the problem $\bA=\bW\bZ$. 
Show that the resulting loss function takes the form of a \textit{weighted $\ell_2$-norm (or a weighted Frobenius norm)}: $L(\bW,\bZ)=\normf{\bW\hadaprod(\bA-\bW\bZ)}^2=\sum_{m,n=1}^{M,N} w_{mn}(a_{mn}-b_{mn})^2$ if $\bB\triangleq\bW\bZ=\{b_{mn}\}\in\real^{M\times N}$. 
Explain how the weight $w_{mn}$ relates to the noise variance $\sigma^2_{mn}$ at entry $(m,n)$.

\item Show that the loss function in \eqref{equation:als_gama_los} arises from the deviance defined in \eqref{equation:als_poi_los}.

\item \label{prob:ortho_mf} \textbf{Orthogonal and projective  matrix factorization.} Consider the optimization problem $\mathopmin{\bW}\normf{\bA-\bW\bZ}^2$ such that $\bZ\bZ^\top=\bI_K$, where $\bA\in\real^{M\times N}, \bW\in\real^{M\times K}, \bZ\in\real^{K\times N}$, and $K\leq \min\{M,N\}$. Show that the optimal value $\bW^*$ given $\bZ$ is $\bA\bZ^\top$.
This indicates that the matrix factorization optimization can be equivalently stated as $\mathopmin{\bZ\bZ^\top=\bI_K}\normf{\bA-\bA\bZ^\top\bZ}^2$. 
The relaxed version---dropping the orthogonality constraint---is known as \textit{projective matrix factorization} \citep{yuan2005projective, yang2010linear}:
$$
\mathopmin{\bZ}\normf{\bA-\bA\bZ^\top\bZ}^2,
$$
where each row of $\bA$ is projected onto a $K$-dimensional subspace, hence the name.
Further interpretations of orthogonal and projective factorizations are discussed in Problem~\ref{prob:ortho_nmf}.

\item \label{problem:rls} \textbf{Regularized least squares (RLS).} 
Given $\bA\in\real^{M\times N}, \bb\in\real^M, \bB\in\real^{P\times N}$, and $\lambda\in\real_{++}$, we consider the regularized least squares (RLS) problem:
$$
\mathop{\min}_{\bx\in\real^N} \normtwo{\bA\bx-\bb}^2 + \lambda\normtwo{\bB\bx}^2.
$$
Show that this problem has a unique solution if and only if $\nspace(\bA)\cap \nspace(\bB) = \{\bzero\}$.

\item \label{prob:denoise_rls} \textbf{Denoising via RLS.} Suppose we observe a noisy signal
$
\by = \bx+\be,
$
where $\bx$ is the true signal, and $\be$ is the noise vector.
We want to find an estimate $\bx$ of the observed measurement $\by$ such that $\bx \approx \by$:
$
\min \normtwo{\bx-\by}^2.
$
Apparently, the optimal solution of this optimization is given by $\bx=\by$; however, it is meaningless.
To obtain a smoother estimate, introduce a penalty on differences between consecutive entries:
$
R(\bx) = \sum_{i=1}^{n-1} (x_i - x_{i+1})^2.
$
Then, 
\begin{itemize}
\item Reformulate this as a regularized least squares problem and derive its closed-form solution.
\item Provide real-world applications. For example, in modeling the profit-and-loss trajectory of a financial asset, consecutive daily observations should vary smoothly rather than exhibit abrupt jumps.
\end{itemize}

\item \textbf{Weighted least squares (WLS).}
Building on Lemma~\ref{lemma:ols}, assume each data point $m\in\{1,2,\ldots, M\}$ (i.e., each row of $\bA$) is assigned a positive weight $w_m$. 
This means some  data points may carry greater significance than others, and we can produce approximate minimizers that reflect this.
Show that the value $\bx_{WLS} = (\bA^\top\bW^2\bA)^{-1}\bA^\top\bW^2\bb$ serves as the \textit{weighted least squares (WLS)}  estimate of $\bx$, where $\bW=\diag(w_1, w_2, \ldots, w_M)\in\real^{M\times M}$. \textit{Hint: Derive the normal equations for this weighted problem.}

\item \textbf{Positive definite  weighted least squares (PDWLS).}
Building on Lemma~\ref{lemma:ols}, we consider further  the matrix equation $\bA\bx + \be =\bb$, where $\be$ is an error vector. Define the weighted error squared sum $E_w = \be^\top \bW \be$, where the weighting matrix $\bW$ is  positive definite. 
Show that the positive definite weighted least squares solution is $\bx^* = (\bA^\top\bW\bA)^{-1}\bA^\top\bW\bb$. \textit{Hint: Compute the gradient of $E_w = (\bb-\bA\bx)^\top\bW(\bb-\bA\bx)$.}

\item \textbf{Weighted color noise least squares.}
Building on Lemma~\ref{lemma:ols}, we consider  the matrix equation $ \bA\bx + \be = \bb $, where $\be$ is an additive color noise vector satisfying the conditions $\Exp[\be] = \bzero$ and $\Exp[\be\be^\top] = \bSigma$, where $\bSigma$ is known. Use the weighting error function $E_w = \be^\top \bW \be$ as the cost function for finding the optimal estimate $\bx^*$. Show that
$\bx^* = (\bA^\top \bW \bA)^{-1} \bA^\top \bW \bb$,
where the optimal choice of the weighting matrix $\bW$ is $\bW^* = \bSigma^{-1}$.
\textit{Hint: Compute the gradient of $E_w = (\bb-\bA\bx)^\top\bW(\bb-\bA\bx)$.}

\item \label{problem:tls} \textbf{Transformed least squares (TLS).}
Building on Lemma~\ref{lemma:ols}, we consider further the restriction $\bx=\bC\bgamma+\bc$, where $\bC\in\real^{N\times K}$ is a known matrix such that $\bA\bC$ has full rank, $\bc$ is a known vector, and $\bgamma$ is an unknown vector.
Show that the value $\bx_{TLS}=\bC(\bC^\top\bA^\top\bA\bC)^{-1}(\bC^\top\bA^\top)(\bb-\bA\bc) +\bc$ serves as the \textit{transformed least squares (TLS)} estimate of $\bx$.

\item \label{problem:twls2} Derive the transformed weighted least squares estimate.

\index{First-order optimality condition}
\index{Fermat's theorem}
\item \label{problem:fist_opt} \textbf{First-order optimality condition for local optima points.} 
Consider  \textit{Fermat's theorem}: for a one-dimensional function $g(\cdot)$ defined and differentiable over an interval ($a, b$), if a point $x^*\in(a,b)$ is a local maximum or minimum, then $g^\prime(x^*)=0$. 
Prove the first-order optimality conditions for multivariate functions based on  Fermat's theorem for one-dimensional functions.
That is, let  $f: \sS\rightarrow \real$ be a function defined on a set $\sS\subseteq \real^N$. Suppose that $\bx^*\in\text{int}(\sS)$, i.e., in the interior point of the set, is a local optimum point and that all the partial derivatives (Definition~\ref{definition:partial_deri}) of $f$ exist at $\bx^*$. Then $\nabla f(\bx^*)=\bzero$, i.e., the gradient vanishes at all local optimum points. (Note that, this optimality condition is a necessary condition but not sufficient; however, there could be vanished points which are not local maximum or minimum point.)
\textit{Hint: Consider the one-dimensional function $g(t)=f(\bx^* + t\be_n)$ for $n\in\{1,2,\ldots,N\}$.}

\item \label{prob:rank-of-ata} \textbf{Rank of $\bA^\top \bA$.}
Show that the matrices $\bA^\top \bA$ and $\bA$ share the same rank.
Similarly, show that $\bA\bA^\top$ and $\bA$ share the same rank.

\item \label{problem:pos_hessian} \textbf{Global minimum point of convex functions.} Let function $f$ be a twice continuously differentiable function defined over $\real^N$. Suppose that the Hessian $\nabla^2f(\bx) \succeq 0$ for any $\bx\in\real^N$ (i.e., the Hessian is always positive semidefinite~\footnote{Instead, if we assume the Hessian is positive semidefinite at a given point, then the point is a local minimum point.}).
This property is also referred to as the \textit{convexity}.
Show that $\bx^*$ is a global minimum point of $f$ if $\nabla f(\bx^*)=\bzero$. \textit{Hint: Use the linear approximation theorem in Theorem~\ref{theorem:linear_approx}.}

\item \textbf{Two-sided matrix least squares} Let $\bA\in\real^{M\times N}$, $\bB\in\real^{M\times K}$, and $\bC\in\real^{P\times N}$. Find the $K\times P$ matrix $\bX$ such that $L(\bX)=\norm{\bA - \bB\bX\bC}_F^2$ is minimized.
\begin{itemize}
\item Derive the derivative of $L$ with respect to $\bX$ and the optimality conditions. 
\item Show that one possible solution to the optimality conditions is $\bX^*=\bB^+\bA\bC^+$, where $\bB^+$ and $\bC^+$ are the pseudo-inverses of $\bB$ and $\bC$, respectively.
\end{itemize}
Similarly, consider the optimization with $\rank(\bX)\leq p$:
$
L(\bX)=\norm{\bA - \bB\bX\bC}_F^2$, s.t. 
$\rank(\bX)\leq p$.
Show that 
\begin{itemize}
\item One possible solution to this is $\bX^*=\bB^+\bA_p\bC^+$, where $\bA_p$ a truncated SVD of $\bB\bB^+\bA\bC^+\bC$ by replacing all but the $p$ largest singular values by zero.
\item  $\bX^*$ also minimizes $\normf{\bX}$, i.e., has the smallest magnitude among all  solutions.
\item $\bX^*$ is the \textbf{unique} solution if and only if either $\rank(\bB\bB^+\bA\bC^+\bC)\leq p$ or both $\rank(\bB\bB^+\bA\bC^+\bC)\geq p$ and $\sigma_{p+1}(\bB\bB^+\bA\bC^+\bC) < \sigma_{p}(\bB\bB^+\bA\bC^+\bC)$.
\end{itemize}


\item \label{problem:nuclear_equi} \citep{rennie2005fast, mazumder2010spectral} Consider the nuclear norm (i.e., the sum of singular values of a matrix, which provides the tightest convex envelope of the rank function of a matrix) $\norm{\bA}_n$ of any matrix $\bA\in\real^{M\times N}$ with rank $R$. Show that 
$$
\norm{\bA}_n = \mathop{\min}_{\substack{\bW\in\real^{M\times R} \\ 
\bZ\in\real^{R\times N} \\
}}
\frac{1}{2} (\normf{\bW}^2 + \normf{\bZ}^2)
\gap \text{s.t.} \gap 
\bA=\bW\bZ
$$

\item Discuss the gradient descent updates corresponding to the different regularization schemes presented in Section~\ref{section:regularization-extention-general}.

\end{problemset}

\newpage
\chapter{Nonnegative Matrix Factorization (NMF)}\index{NMF}\label{chapter:nmf}
\begingroup
\hypersetup{
	linkcolor=structurecolor,
	linktoc=page,  
}
\minitoc \newpage
\endgroup

\index{Decomposition: NMF}
\index{Sparsity}
\index{Nonnegativity constraint}
\section{Nonnegative Matrix Factorization}
\lettrine{\color{caligraphcolor}I}
In the era of big data, extracting meaningful patterns and latent structures from high-dimensional datasets has become a central challenge across scientific and technological domains.
Singular value decomposition (SVD) is grounded in strong theoretical foundations and enjoys broad applicability. However, it has notable limitations---particularly when applied to nonnegative matrices~\footnote{Nonnegative matrices exhibit unique properties in linear algebra and are essential for theoretical analysis; see Problems~\ref{prob:nnga_algebra}--\ref{prob:nonn_lin_12}.}. In such cases, SVD may produce negative components, which often lack physical interpretability.

To address this issue, \textit{nonnegative matrix factorization (NMF)} has emerged as a powerful and interpretable tool for dimensionality reduction, feature extraction, and uncovering latent structure in complex data. Early work on this problem was carried out by \citet{paatero1994positive} and \citet{cohen1993nonnegative}, who referred to it as \textit{positive matrix factorization}. The method gained widespread attention following the introduction of the \textit{multiplicative update} rule by \citet{lee2001algorithms}.

Building on the alternating least squares (ALS) framework for matrix factorization, we now turn to algorithms for solving the NMF problem:
\begin{itemize}
\item Given a nonnegative matrix $\bA\in \real_+^{M\times N}$ of rank $R$, find nonnegative matrix factors $\bW\in \real_+^{M\times K}$ and $\bZ\in \real_+^{K\times N}$ such that: 
$
\bA\approx\bW\bZ.
$
\end{itemize}

As discussed in the ALS section, a core goal in linear data analysis is to represent high-dimensional data vectors as linear combinations of lower-dimensional basis vectors. These basis vectors---often called \textit{hidden vectors}, \textit{pattern vectors}, or \textit{feature vectors}---capture the essential characteristics of the data and are crucial for tasks like pattern recognition.
For such pattern vectors to be useful in modeling and interpretation, they should satisfy two key criteria:
\begin{itemize}
\item \textit{Interpretability.} Each component should correspond to a physically or physiologically meaningful quantity, enabling intuitive understanding of the underlying data.
\item \textit{Statistical fidelity.} When the data are reliable and low-noise, the pattern vectors should effectively capture the dominant modes of variation and reflect the primary distribution of information.
\end{itemize}
NMF excels at meeting these requirements across diverse applications:
\begin{itemize}
\item In document collections, each document is represented as a vector of term frequencies (often weighted, e.g., via TF-IDF). Stacking these vectors yields a nonnegative term-by-document matrix encoding the entire corpus.
\item In image collections, each image is flattened into a pixel-intensity vector with nonnegative entries. Arranging these vectors column-wise produces a nonnegative pixel-by-image matrix.
\item In gene expression analysis, measurements under different experimental conditions form a gene-by-experiment matrix, capturing how gene activity varies across conditions.
\item In recommender systems, user-item interactions (e.g., purchase counts or ratings) are stored in a large, sparse, nonnegative matrix that reflects the limited engagement of users with most items.
\end{itemize}
Unlike general linear decompositions, NMF restricts both the basis vectors (columns of $\bW$) and their combination coefficients (entries of $\bZ$) to be nonnegative. This eliminates phenomena like destructive interference, where positive and negative contributions cancel each other out. Instead, data reconstruction relies solely on additive, parts-based representations.

The nonnegativity constraint inherently promotes sparsity, allowing NMF to isolate distinct, interpretable features. This property makes it especially valuable in domains where data naturally decompose into constituent parts. 
For example, in image processing, NMF has been successfully applied to object detection, image segmentation, and facial recognition \citep{lee2001algorithms, gillis2014and, gillis2020nonnegative}, where the nonnegative components align with intuitive visual parts (e.g., eyes, noses, textures).
In topic modeling or document analysis, each column of $\bA$ represents a document. NMF yields a soft clustering where columns of $\bW$ correspond to topics, and the entries of $\bZ$ indicate the degree to which each document belongs to each topic \citep{shahnaz2006document}.
In clustering,  a nonnegative factorization $\bA\approx\bW\bZ$ can also serve as a clustering tool. Specifically, data vector $\ba_j$ is assigned to cluster $i$ if $z_{ij}$ is the largest entry in column $j$ of $\bZ$ \citep{brunet2004metagenes, gao2005improving}.
For broader context, see the survey by \citet{berry2007algorithms}. In summary, NMF's popularity stems from its ability to automatically extract sparse, nonnegative, and interpretable latent factors.

To measure the quality of the  approximation, we evaluate the loss by computing  the Frobenius norm of the difference between the original matrix and the approximation:
\begin{equation}\label{equation:frob_nmf}
L(\bW,\bZ) \triangleq D(\bA, \bW\bZ) = \frac{1}{2}\normf{\bW\bZ-\bA}^2,~\footnote{The factor ${1}/{2}$ simplifies gradient computations.}
\end{equation}
where $L(\bW,\bZ)$ indicates it is a loss function w.r.t. $\bW$ and $\bZ$, and $D(\bA, \bW\bZ)$ implies it is a distance/divergence between $\bA$ and $\bW\bZ$ (we will use the two notations interchangeably  when  necessary).
The Frobenius norm is arguably the most widely used norm for NMF because it corresponds to Gaussian additive noise, which is reasonable in many situations and allows for the design of particularly efficient algorithms; see Section~\ref{section:more_err_sta_als}.
For nonnegative data, Gaussian noise can be interpreted as a truncated version of Gaussian noise.
Later, we will generalize this framework to other loss functions based on $\beta$-divergences (Section~\ref{section:beta_div_altmu}).

When an exact factorization $\bA=\bW\bZ$ with  $\bW\in\real^{M\times R}_+$ and $\bZ\in\real_+^{R\times N}$ exist, the problem is known as \textit{exact NMF} of size $R$. However, Exact NMF is NP-hard \citep{vavasis2010complexity, gillis2020nonnegative}, so in practice we focus exclusively on approximate NMF.

\index{Overfitting}
\index{Bayesian inference}
\index{Bayesian optimization}
\index{Bayesian matrix decomposition}
In collaborative filtering, NMF trained via multiplicative updates---despite favorable convergence guarantees---can suffer from overfitting. While regularization helps mitigate this issue, out-of-sample performance often remains suboptimal. In contrast, Bayesian approaches based on generative models can effectively control overfitting in NMF; see Chapter~\ref{chapter:bnmf}.
In the following sections, we introduce several algorithms for solving NMF problems and briefly discuss their practical applications.

\begin{algorithm}[h] 
\caption{Projected Gradient Descent (PGD) Method}
\label{alg:pgd_gen}
\begin{algorithmic}[1] 
\Require A function $f(\bx)$ and a set $\sS$; 
\For{$t=1,2,\ldots$}
\State Pick a step size $\eta_t$;
\State Set $\bx^{(t+1)} \leftarrow \mathcalP_{\sS}(\bx^{(t)} - \eta_t \nabla f(\bx^{(t)}))$;
\EndFor
\State Output final  $\bx$;
\end{algorithmic} 
\end{algorithm}
\section{NMF via Alternating Projected Gradient Descent (APGD)}\label{section:nmf_apgd}
\textit{Projected gradient descent} (\textit{PGD}; Algorithm~\ref{alg:pgd_gen}) solves optimization problems of the form
$$
\mathopmin{\bx\in\sS} f(\bx),
$$
where $\sS\subseteq\real^N$ is a constraint set. The method relies on the orthogonal projection onto $\sS$, defined as
$
\mathcalP_{\sS} (\bx) \triangleq \mathop{\argmin}_{\by\in\sS} \normtwo{\by-\bx}.
$
When $\sS$ is the nonnegative orthant ($\real^N_+$), this projection simplifies to componentwise thresholding: $\mathcalP_{\sS} (\bx) = \max\{\bzero, \bx\}$.

Applying this idea to NMF yields the \textit{alternating projected gradient descent (APGD)}  approach, which updates the factors $\bW$ and $\bZ$ iteratively:
$$
\bZ\leftarrow \max\bigg\{\bzero, \mathop{\argmin}_{\bZ\in\real^{K\times N}}\normf{\bW\bZ-\bA} \bigg\}
\gap\text{and}\gap
\bW\leftarrow \max\bigg\{\bzero, \mathop{\argmin}_{\bW\in\real^{M\times K}}\normf{\bW\bZ-\bA} \bigg\}.
$$
Each subproblem is a least squares problem followed by projection onto the nonnegative orthant.
However, due to the projection step, the resulting factors may be poorly scaled. A simple remedy is to apply a closed-form scaling factor $\gamma\geq 0$ at each iteration:
$$
\gamma^* = \mathop{\argmin}_{\gamma\geq 0} \normf{\gamma\bW\bZ-\bA} 
=
\frac{\langle \bA, \bW\bZ\rangle}{\langle \bW\bZ, \bW\bZ\rangle}
=
\frac{\langle \bA\bZ^\top, \bW\rangle}{\langle \bW^\top\bW, \bZ\bZ^\top\rangle}.
$$
Although APGD is generally not recommended as a standalone solver due to slow or unstable convergence, it can be highly effective as an initialization strategy. Specifically, running a few APGD iterations before switching to a more robust NMF algorithm often yields significant improvements---especially for sparse matrices \citep{gillis2014and}.

\index{Nonnegative least squares}
\index{NNLS|see {Nonnegative least squares}}
\index{ANLS|see {Nonnegative least squares}}
\section{NMF via Alternating Nonnegative Least Squares (ANLS)}\label{section:nmf_anls}
The alternating least squares (ALS) framework hinges on solving ordinary least squares (OLS) subproblems (Lemma~\ref{lemma:ols}). For NMF, we can replace OLS with the \textit{nonnegative least squares} (NNLS) problem:
\begin{equation}
\mathopmin{\bx\geq \bzero } f(\bx) = \mathopmin{\bx\geq \bzero } \frac{1}{2}\normtwo{\bb-\bM\bx}^2
\gap
\text{with }\bM\in\real^{M\times N}, \bb\in\real^M, \bx\in\real_+^N. 
\end{equation}
The KKT conditions for this problem imply complementary slackness: $\lambda_nx_n^*=0, \forall n$, where $\lambda_n$ is the Lagrange multiplier associated with the constraint $x_n\geq 0$. 
Additionally, the stationarity condition gives $\nabla f(\bx^*) -\sum_{n}\lambda_n \be_n=\bzero$, where $\bx^*$ denotes the optimal solution of the NNLS problem \citep{lu2021numerical, lu2025practical}.
Together, the complementary slackness and the optimal condition indicate that:
$$
\nabla f(\bx^*)
=
\sum_{n:x_n^*=0} \lambda_n \be_n,
$$
which leads to the following equivalent characterization of the KKT conditions:
\begin{equation}\label{equation:kkv_nnn_raw}
(\textbf{KKT of NNLS})\gap 
	\bx^*\geq \bzero, 
	\gap
	\nabla f(\bx^*)\geq 0, 
	\gap
	\text{and}
	\gap
	x_n^* (\nabla f(\bx^*))_n=0,\, \forall n.
\end{equation}
These conditions reveal that NNLS---and by extension, NMF---naturally induces sparsity.

Suppose we are given the \textit{inactive set} $\sI\subseteq \{1,2,\ldots,N\}$, defined as 
$$
\sI = \left\{ n \mid x_{n}^{*} > 0, \,\forall n \in\{1,2,\ldots,N\} \right\}.
$$
Its complement, the \textit{active set}, contains indices where $x_n^*=0$. 
On the inactive set, the nonnegativity constraints are inactive, so the solution satisfies the unconstrained optimality condition:
$$
\begin{aligned}
[\nabla_{\bx} f(\bx)]_{\sI} = \bzero 
\gapthree\Longleftrightarrow\gapthree [\bM^\top(\bM \bx - \bb)]_{\sI} = \bzero 
\gapthree\Longleftrightarrow\gapthree \bM[:, \sI]^\top \bM[:, \sI] \bx[\sI] = \bM[:, \sI]^\top \bb.
\end{aligned}
$$
This is precisely the normal equation for the unconstrained least squares problem for $\bx[\sI]$:
$$
\min_{\bx[\sI]} \frac{1}{2}\normtwo{\bb -\bM[:, \sI] \bx[\sI]}^2.
$$
This observation underpins the \textit{active-set method}, which iteratively refines the active and inactive sets through pivoting (adding or removing variables) to ensure monotonic decrease of the objective function \citep{lawson1995solving}; see Algorithm~\ref{alg:nmf_anls}.

\paragrapharrow{Alternating nonnegative least squares (ANLS).} 
Equipped with an NNLS solver, we can adapt the ALS framework to NMF by replacing OLS with NNLS---a strategy known as \textit{alternating nonnegative least squares} (ANLS) \citep{kim2011fast}.
Given a fixed $\bW$, the NMF objective for $\bZ$ decomposes column-wise:
$$
\frac{1}{2}\normf{\bA-\bW\bZ}^2 =
\frac{1}{2}\sum_{n=1}^{N}\normtwo{\ba_n - \bW\bz_n}^2,
$$
where each subproblem $\mathopmin{\bz_n\geq \bzero}\normtwo{\ba_n - \bW\bz_n}^2$ can be updated independently by solving an NNLS problem.
By symmetry of the NMF formulation---$\bA=\bW\bZ$ if and only if $\bA^\top=\bZ^\top\bW^\top$, and  $D(\bA, \bW\bZ)=D(\bA^\top, \bZ^\top\bW^\top)$---the update for $\bW$ (given $\bZ$)  follows analogously.
It is worth noting that during early iterations, when $\bW$ and $\bZ$ provide a poor approximation of $\bA$, solving the NNLS subproblems to high accuracy is often unnecessary. A more efficient strategy is to use ANLS as a refinement step within a faster, less accurate NMF algorithm---such as APGD or multiplicative updates (MU, as discussed in later sections).

\begin{algorithm}[h] 
\caption{Nonnegative Least Squares (NNLS) via Active-Set Method}
\label{alg:nmf_anls}
\begin{algorithmic}[1] 
\Require A real-valued matrix $\bM\in\real^{M \times N}$, a real-valued vector $\bb\in\real^M$;
\State Initialize index sets $\sI = \varnothing$ and $\sJ = \{1,2, \ldots, N\}$;
\State Initialize unknown $\bx\in\real^N$ to an all-zero vector and let $\bw \leftarrow \bM^\top(\bb - \bM\bx)$;
\State Let $\bw[\sJ]$ denote the sub-vector with indices from $\sJ$;
\State Choose a stopping criterion on the approximation error $\delta$;
\State Choose the maximum number of iterations $C$;
\State $iter=0$; \Comment{Count for the number of iterations}
\While{$\sJ \neq \varnothing$ and $\max(\bw[\sJ]) > \delta$ and $iter<C$}
\State $iter=iter+1$; 
\State Let $j$ in $\sJ$ be the index of $\max(\bw[\sJ])$ in $\bw$: $j=\mathop{\argmax}_{j\in \sJ} w_j$;
\State Add $j$ to $\sI$  and remove $j$ from $\sJ$ such that $\sI\cup \sJ = \{1,2,\ldots,N\}$;
\State Let $\bM[:,\sI]$ be $\bM$ restricted to the variables/columns included in $\sI$;
\State \parbox[t]{\dimexpr\linewidth-\algorithmicindent}{Let $\bs$ be vector of same length as $\bx$;
Let $\bs[\sI]$ denote the sub-vector with indices from $\sI$, and let $\bs[\sJ]$ denote the sub-vector with indices from $\sJ$;}
\State Set $\bs[\sI] \leftarrow ((\bM[:,\sI])^\top \bM[:,\sI])^{-1} (\bM[:,\sI])^\top \bb$ and $\bs[\sJ]$ to zero;
\While{$\min(\bs[\sI]) \leq 0$}
\State Let $\alpha \leftarrow \min \frac{x_i}{x_i - s_i}$ for $i$ in $\sI$ where $s_i \leq 0$;
\State Set $\bx\leftarrow \bx + \alpha(\bs - \bx)$;
\State Move to $\sJ$ all indices $j$ in $\sI$ such that $x_j \leq 0$;
\State Set $\bs[\sI] \leftarrow ((\bM[:,\sI])^\top \bM[:,\sI])^{-1} (\bM[:,\sI])^\top \bb$;
\EndWhile
\State Set $\bs[\sJ]$ to zero;
\State Set $\bx\leftarrow \bs$;
\State Set $\bw\leftarrow \bM^\top(\bb - \bM\bx)$;
\EndWhile
\State Output $\bx$;
\end{algorithmic} 
\end{algorithm}

\index{Hierarchical ANLS}
\section{NMF via Hierarchical Alternating Nonnegative Least Squares}
Let $\ba, \bb\in\real_+^N$ be two nonnegative vectors. 
The \textit{univariate NNLS} problem is then formulated as
$$
\mathopmin{x\geq 0} \normtwo{\ba-x\bb}^2.
$$
If $\normtwo{\bb}\neq 0$, this problem admits a closed-form solution: $x=\max\big\{0, {\bb^\top\ba}/{\normtwo{\bb}^2}\big\}$.
Motivated by this simple case, consider the $k$-th row of $\bZ$ for $k\in\{1,2,\ldots,K\}$. 
In NMF, the corresponding subproblem becomes
\begin{equation}\label{equation:llipschi_hianls}
\mathopmin{\bZ[k,:]\geq \bzero} \bigg\Vert\underbrace{\big(\bA-\sum_{p\neq k}^{K} \bW[:,p]\bZ[p,:]\big)}_{\triangleq\bA_k} - \bW[:,k]\bZ[k,:]\bigg\Vert_F^2, 
\gap \forall k.~\footnote{This subproblem is convex and is $L$-Lipschitz gradient continuous ($L$-strongly smooth); see Problem~\ref{prob:llipschi_hianls}.}
\end{equation}
Equation~\eqref{equation:llipschi_hianls} reveals that the entries within a single row of $\bZ$ do not interact with one another---similarly, entries within a single column of $\bW$ are decoupled.
Consequently, the optimization over each entry in a row of $\bZ$ can be performed independently.
Defining $\bA_k\triangleq\big(\bA-\sum_{p\neq k}^{K} \bW[:,p]\bZ[p,:]\big)$, the NMF update reduces to a set of rank-one approximations of $\bA_k$, for $k\in\{1,2,\ldots,K\}$.
The optimal solution is given by
$$
\bZ^*[k,:]=\mathop{\argmin}_{\bZ[k,:]\geq \bzero} \normf{\bA_k -\bW[:,k]\bZ[k,:]}^2
=
\max\left(
\bzero, 
\frac{\bW[:,k]^\top\bA_k}{\normtwo{\bW[:,k]}^2}
\right),
\gap \forall k,
$$
where the \texttt{max} operator is applied componentwise.
This leads to the \textit{hierarchical ANLS (Hi-ANLS)} method for NMF, which iteratively solves a sequence of univariate NNLS problems.
The procedure is summarized  in Algorithm~\ref{alg:hie_anls}, where we note that $\bZ[k,:]^\top = \bZ^\top[:,k]$.
In the algorithm,  the $k$-th row of $\bZ$ and the $k$-th column of $\bW$ are updated in an interleaved fashion. As shown by \citet{gillis2012accelerated}, updating  $\bZ$ several times before updating $\bW$ can significantly improve performance, since it reuses precomputed quantities such as  $\bW^\top\bA$ and $\bW^\top\bW$.

\begin{algorithm}[h] 
\caption{NMF via Hierarchical Alternating Nonnegative Least Squares (Hi-ANLS)}
\label{alg:hie_anls}
\begin{algorithmic}[1] 
\Require Matrix $\bA\in \real_+^{M\times N}$;
\State Initialize $\bW\in \real_{++}^{M\times K}$, $\bZ\in \real_{++}^{K\times N}$ randomly with positive entries;
\State Choose a stopping criterion on the approximation error $\delta$;
\State Choose maximal number of iterations $C$;
\State $iter=0$; \Comment{Count for the number of iterations}
\While{$\normf{\bA- (\bW\bZ)}^2>\delta $ and $iter<C$}
\State $iter=iter+1$;  
\For{$k=1$ to $K$}
\State $
\bZ[k,:]\leftarrow 
\max\left(
\bzero, 
\frac{\bW[:,k]^\top\bA_k}{\normtwo{\bW[:,k]}^2}
\right)$; \Comment{$\bA_k\triangleq\big(\bA-\sum_{p\neq k}^{K} \bW[:,p]\bZ[p,:]\big)$}

\State $
\bW[:, k]\leftarrow 
\max\left(
\bzero, 
\frac{\bA_k\bZ[k,:]^\top}{\normtwo{\bZ[k,:]}^2}
\right)$;
\EndFor
\EndWhile
\State Output $\bW,\bZ$;
\end{algorithmic} 
\end{algorithm}

\section{NMF via Alternating Direction Methods of Multipliers (ADMM)}\label{section:nmf_admm_all}
We briefly introduce the \textit{alternating direction methods of multipliers (ADMM)} and then  discuss its application to matrix factorization and NMF.
\paragrapharrow{ADMM.}
ADMM solves convex optimization problems of the form
\begin{equation}\label{equation:admm_prob}
\mathopmin{\bx, \bz} f(\bx)+g(\bz) \gap \text{s.t.}\gap \bD\bx+\bE\bz=\bff.
\end{equation}
Given a penalty  parameter $\rho>0$, the \textit{augmented Lagrangian} associated with \eqref{equation:admm_prob} is 
\begin{equation}
L_\rho (\bx, \bz, \bl) = f(\bx)+g(\bz) +\langle \bl, \bD\bx+\bE\bz-\bff\rangle + \frac{\rho}{2}\normtwo{\bD\bx+\bE\bz-\bff}^2.
\end{equation}
When $\rho=0$, this reduces to the standard Lagrangian; when $\rho>0$,  it becomes a penalized version that improves numerical stability.
The classical \textit{augmented Lagrangian method} solves the problem by performing the following steps  (at the $(t+1)$-th iteration):
$$
\text{augmented Lagrangian:}
\gap 
\left\{
\begin{aligned}
(\bx^{(t+1)}, \bz^{(t+1)}) &\in \mathop{\argmin}_{\bx, \bz} L_\rho (\bx, \bz, \bl);\\
\bl^{(t+1)}&=\bl^{(t)} + \rho(\bD\bx^{(t+1)}+\bE\bz^{(t+1)} -\bff ),
\end{aligned}
\right.
$$
where the update on $\bl^{(t+1)}$ is derived from the \textit{conjugate subgradient theorem} (see, for example, \citet{bach2011convex, lu2026first}), and 
the symbol `$\in$' acknowledges that minimizers may not be unique.
A key challenge lies in the coupling between  $\bx$ and  $\bz$ through the quadratic term $\rho(\bx^\top\bD^\top\bE\bz)$.
ADMM tackles this difficulty by replacing the exact minimization of $(\bx,\bz)$ with one iteration of the alternating minimization method (see Algorithm~\ref{alg:am_lvm}).
To be more specific, for the  $(t+1)$-iteration, the solution of ADMM takes the following form:
\begin{equation}
\text{ADMM:}\gap
\left\{
\begin{aligned}
	\bx^{(t+1)}&\in\mathop{\argmin}_{\bx} \left\{ f(\bx)+ \frac{\rho}{2}\normtwo{\bD\bx+\bE\bz^{(t)} -\bff +\frac{1}{\rho}\bl^{(t)}}^2 \right\};\\
	\bz^{(t+1)}&\in\mathop{\argmin}_{\bz} \left\{ g(\bz)+ \frac{\rho}{2}\normtwo{\bD\bx^{(t+1)}+\bE\bz -\bff +\frac{1}{\rho}\bl^{(t)}}^2 \right\};\\
	\bl^{(t+1)}&=\bl^{(t)} + \rho(\bD\bx^{(t+1)}+\bE\bz^{(t+1)} -\bff ).
\end{aligned}
\right.
\end{equation}
Introducing the scaled dual variable $\widetildebl\triangleq\frac{1}{\rho}\bl$, this is equivalently expressed as (the form we adopt hereafter):
\begin{equation}\label{equation:admm_gen_up}
\text{ADMM:}\gap
\left\{
\begin{aligned}
\bx^{(t+1)}&\in\mathop{\argmin}_{\bx} \left\{ f(\bx)+ \frac{\rho}{2}\normtwo{\bD\bx+\bE\bz^{(t)} -\bff +\widetildebl^{(t)}}^2 \right\};\\
\bz^{(t+1)}&\in\mathop{\argmin}_{\bz} \left\{ g(\bz)+ \frac{\rho}{2}\normtwo{\bD\bx^{(t+1)}+\bE\bz -\bff +\widetildebl^{(t)}}^2 \right\};\\
\widetildebl^{(t+1)}&=\widetildebl^{(t)} + (\bD\bx^{(t+1)}+\bE\bz^{(t+1)} -\bff ).
\end{aligned}
\right.
\end{equation}
Thus, ADMM iteratively updates $\bx, \bz$, and the scaled dual variable $\widetildebl$.

\paragrapharrow{ADMM applied to matrix factorization.}
We return to the problem discussed in ALS (Equation~\eqref{equation:als-per-example-loss2}, i.e., matrix factorization with Frobenius norm; not necessarily a NMF problem) together with a regularization function $r(\bZ)$:
$$
\mathopmin{\bZ} \frac{1}{2}\normf{\bA-\bW\bZ}^2+r(\bZ).
$$
Introducing an auxiliary variable $\widetildebZ\in\real^{K\times N}$, we reformulate this as
\begin{equation}\label{equation:mf_admm_prob1}
\mathopmin{\bZ} \frac{1}{2}\normf{\bA-\bW\bZ}^2+r(\widetildebZ), 
\gap 
\text{s.t.}
\gap 
\bZ=\widetildebZ.
\end{equation}
Applying \eqref{equation:admm_gen_up} with either (a). \{$\bx\leftarrow \bZ$, $\bz\leftarrow \widetildebZ$, $\widetildebl\leftarrow \bL$, $\bD=-\bI$,  $\bE=\bI$\} or (b). \{$\bx\leftarrow \bZ$, $\bz\leftarrow \widetildebZ$, $\widetildebl\leftarrow \bL$, $\bD=\bI$,  $\bE=-\bI$\}, 
yields the following ADMM updates for \eqref{equation:mf_admm_prob1}:
\begin{equation}\label{equation:admm_gen_als}
\left\{
\begin{aligned}
\bZ 
&\stackrel{(a)}{\leftarrow} (\bW^\top\bW+\rho \bI)^{-1} \left[ \bW^\top\bA +\rho(\widetildebZ+\bL) \right]
&\stackrel{(b)}{\leftarrow}& (\bW^\top\bW+\rho \bI)^{-1} \left[ \bW^\top\bA +\rho(\widetildebZ-\bL) \right];\\
\widetildebZ
&\stackrel{(a)}{\leftarrow}\mathop{\argmin}_{\widetildebZ} r(\widetildebZ) + \frac{\rho}{2}\normf{-\bZ+\widetildebZ + \bL}^2
&\stackrel{(b)}{\leftarrow}&\mathop{\argmin}_{\widetildebZ} r(\widetildebZ) + \frac{\rho}{2}\normf{\bZ-\widetildebZ + \bL}^2\\
\bL&\stackrel{(a)}{\leftarrow}\bL -\bZ+\widetildebZ &\stackrel{(b)}{\leftarrow}& \bL +\bZ-\widetildebZ.
\end{aligned}
\right.
\end{equation}
In practice, the Cholesky decomposition of $(\bW^\top\bW+\rho \bI)$ can be precomputed, enabling efficient updates via forward and backward substitution \citep{lu2021numerical}.
By symmetry, the update for $\bW$ (given $\bZ$) follows analogously. We adopt formulation (a) in subsequent discussions.

%

\paragrapharrow{ADMM applied to $\ell_1$-regularization.}
We may also consider the $\ell_1$-regularization (see Section~\ref{section:regularization-extention-general}): $r(\widetildebZ)=\lambda \Vert\widetildebZ\Vert_1$. The update for each element $(k,n)$ of $\widetildebZ$ is $\widetildez_{kn}\leftarrow \max(0, 1-\frac{\lambda}{\rho} \abs{h_{kn}}^{-1}) h_{kn}$ for all $k\in\{1,2,\ldots,K\}$ and $n\in\{1,2,\ldots,N\}$, where $h_{kn} = z_{kn}-l_{kn}$ (i.e., the elements of $\bH=\bZ-\bL$). 

\paragrapharrow{ADMM applied to smoothness/denoising regularization.}
The smoothness regularization on $\bZ$ can be defined as $r(\widetildebZ)=\frac{\lambda}{2} \Vert\bT\widetildebZ^\top\Vert_F^2$, where $\bT$ is an $N\times N $ tridiagonal matrix with 2 on the main diagonal and $-1$ on the superdiagonal and subdiagonal. This regularization ensures the proximal components in each row of $\widetildebZ$ is smooth (see Problem~\ref{prob:denoise_rls}). The update  of $\widetildebZ$ becomes $\widetildebZ\leftarrow \rho\bZ(\lambda \bT^\top\bT +\rho\bI)^{-1}$ \citep{huang2016flexible}.

\paragrapharrow{ADMM applied to NMF.}
To enforce nonnegativity in NMF, we replace $r(\bZ)$ with the indicator function of the nonnegative orthant.
The update on $\widetildebZ$ becomes $\max\left(\bzero, \bZ-\bL\right)$, where the \texttt{max} operator is applied  componentwise.
However, unlike the ANLS, Hi-ANLS, or multiplicative update (MU) methods discussed later in the next section, ADMM updates for NMF are generally not guaranteed to produce a monotonically decreasing objective value.


\index{Alternating update}
\index{Kullback--Leibler divergence}
\index{Multiplicative update}
\section{NMF via Multiplicative Update (MU)}\label{section:nmf_frob_mu}
We now consider an alternative alternating update strategy for NMF. The latent factors $\bW$ and $\bZ$ are modeled as nonnegative vectors in a low-dimensional space. These hidden representations are initialized randomly and then iteratively refined using an alternating multiplicative update rule to minimize the Frobenius norm between the observed data matrix $\bA$ and its low-rank approximation $\bW\bZ$.
Following the setup in Section~\ref{section:als-netflix}, we assume a rank-$K$ factorization. Given $\bW\in \real_+^{M\times K}$, our goal is to update $\bZ\in \real_+^{K\times N}$. The gradient of the loss function 
$L(\bW,\bZ)=\frac{1}{2}\normf{\bA-\bW\bZ}^2$ with respect to $\bZ$ is (see Equation~\eqref{equation:givenw-update-z-allgd}):
$
\begin{aligned}
	\nabla_{\bZ} L(\bW, \bZ) =\bW^\top(\bW\bZ-\bA) \in \real^{K\times N}.
\end{aligned}
$
Applying standard gradient descent (as discussed in Section~\ref{section:als-gradie-descent}), a naive update for $\bZ$ would be:
$$
(\text{GD on $\bZ$})\gap \bZ \leftarrow \bZ - \eta \left(\nabla_{\bZ} L(\bW, \bZ)\right)=\bZ - \eta \nabla_{\bZ} L(\bW, \bZ),
$$
where $\eta>0$ is a small, fixed step size.

\paragrapharrow{Multiplicative update (MU).}
Instead of using a uniform step size, suppose we allow a distinct step size $\eta_{kn}>0$ for each entry $z_{kn}$ of $\bZ$. 
The update becomes:
$$
(\text{GD$^\prime$ on $\bZ$})\gap 
\begin{aligned}
	z_{kn} &\leftarrow z_{kn} - \frac{\eta_{kn}}{2} \left(\nabla_{\bZ} L(\bW, \bZ)\right)_{kn}
	=z_{kn} - \eta_{kn}(\bW^\top\bW\bZ-\bW^\top\bA)_{kn}, \,\, \forall k,n,
\end{aligned}
$$
Now choose the adaptive step size:
$$
\eta_{kn} = \frac{z_{kn}}{(\bW^\top\bW\bZ)_{kn}}.
$$
Substituting this into the update yields the \textit{multiplicative update (MU)} rule \citep{lee2001algorithms}:
\begin{equation}\label{equation:multi-update-z}
(\text{MU on $\bZ$})\gap 
\bZ \leftarrow \bZ\hadaprod \frac{[\bW^\top\bA]}{[\bW^\top\bW\bZ]}
\stackrel{*}{=}
\bZ - \frac{[\bZ]}{[\bW^\top\bW\bZ]}\hadaprod \nabla_{\bZ} L(\bW, \bZ) 
,
\end{equation}
where $\frac{[\cdot]}{[\cdot]}$ represents the componentwise division, and $\hadaprod$ is the Hadamard (elementwise) product. 
By symmetry, the corresponding update for $\bW$ is:
\begin{equation}\label{equation:multi-update-w}
(\text{MU on $\bW$})\gap
\bW \leftarrow \bW \hadaprod \frac{[\bA\bZ^\top]}{[\bW\bZ\bZ^\top]} 
\stackrel{*}{=} 
\bW - \frac{[\bW]}{[\bW\bZ\bZ^\top]}\hadaprod \nabla_{\bW} L(\bW, \bZ) .
\end{equation}
The ratios $\frac{(\bW^\top\bA)_{kn}}{(\bW^\top\bW\bZ)_{kn}}$ and $\frac{(\bA\bZ^\top)_{mk}}{(\bW\bZ\bZ^\top)_{mk}}$ for all $m,k,n$ in \eqref{equation:multi-update-z} and \eqref{equation:multi-update-w} are called \textit{multiplicative factors}.
When $\bA=\bW\bZ$, these factors equal one, and the gradients vanish---indicating a stationary point.

\paragrapharrow{MU vs. gradient descent.}
The derivation above reveals that MU is fundamentally a variant of gradient descent, differing only in how the step size is chosen. In standard gradient descent, the step size $\eta$ may be fixed or adapt globally over time, but it is shared across all entries of the variable matrix at each iteration. In contrast, MU assigns a different, entry-specific step size that scales inversely with the current magnitude of the variable and the curvature of the objective. This adaptivity often leads to faster practical convergence and automatic satisfaction of nonnegativity constraints.

\paragrapharrow{KKT conditions for NMF with Frobenius norm.}
The KKT conditions for the NMF problem (Equation~\eqref{equation:frob_nmf}) are (cf. Equation~\eqref{equation:kkv_nnn_raw}):
\begin{equation}\label{equation:nmf_fro_kkt1}
\begin{aligned}
\bZ\geq \bzero,\gap&\nabla_{\bZ} L(\bW,\bZ)&\geq& \bzero, \gap \langle \bZ, \nabla_{\bZ} L(\bW,\bZ)\rangle &=&\bzero_{K\times N}; \\
\bW\geq \bzero,\gap&\nabla_{\bW} L(\bW,\bZ)&\geq& \bzero, \gap \langle \bW, \nabla_{\bW} L(\bW,\bZ)\rangle &=&\bzero_{M\times K},
\end{aligned}
\end{equation}
where all inequalities and inner products are interpreted componentwise. Equivalently,
\begin{equation}
\begin{aligned}
\min\{\bZ, \nabla_{\bZ} L(\bW,\bZ)\} = \bzero_{K\times N}
\gap \text{and}\gap 
\min\{\bW, \nabla_{\bW} L(\bW,\bZ)\} = \bzero_{M\times K},
\end{aligned}
\end{equation}
where the \texttt{min} operator $\min\{\cdot, \cdot \}$ is applied componentwise. Any pair $(\bW,\bZ)$ satisfying these conditions is a stationary point of the NMF problem.

\paragrapharrow{Problems in MU.}
Equality ($*$) in \eqref{equation:multi-update-z} shows that MU corresponds to a rescaled gradient descent step. Moreover, observe that:
$$
\frac{[\bW^\top\bA]_{kn}}{[\bW^\top\bW\bZ]_{kn}}\geq 1 \gap\Longleftrightarrow\gap (\nabla_{\bZ} L(\bW, \bZ))_{kn}\leq 0, \gap  \forall k, n.
$$
Thus, MU implements three intuitive rules: (i) Increase $z_{kn}$ if its partial derivative is negative; (ii) Decrease it if its partial derivative is positive; (iii) Keep it unchanged if its partial derivative is zero.
However, if an entry $z_{kn}=0$, the MU update leaves it unchanged---even if the gradient is negative (i.e., decreasing the loss would require increasing $z_{kn}$). 
In such cases, the KKT conditions in \eqref{equation:nmf_fro_kkt1} are violated, since $z_{kn}=0$ but $(\nabla_{\bZ} L(\bW,\bZ)\})_{kn}<0$. Consequently, MU iterates are not guaranteed to converge to a stationary point.
Common remedies include: (i) Initializing $\bW$ and $\bZ$ with strictly positive entries and enforcing a small lower bound (e.g., $\epsilon=10^{-9}$) on all entries \citep{gillis2012accelerated};
(ii) Reinitializing any zero entry to a small positive value whenever its gradient becomes negative \citep{chi2012tensors}.

\paragrapharrow{Monotonicity of MU.}
Despite these issues, MU enjoys a crucial theoretical property: it guarantees monotonic decrease of the objective under mild conditions.
\begin{theoremHigh}[Monotonically Nonincreasing of Multiplicative Update]\label{theorem:conv_mu_fro}
The loss $L(\bW,\bZ)=\frac{1}{2}\normf{\bW\bZ-\bA}^2$ remains nonincreasing under the following multiplicative update rules:~\footnote{More general results for $\beta$-divergences are discussed in Theorem~\ref{theorem:conv_mu_beta}.}
$$
\begin{aligned}
\bZ &\leftarrow \bZ\hadaprod \frac{[\bW^\top\bA]}{[\bW^\top\bW\bZ]}
\gap\text{and}\gap 
\bW &\leftarrow \bW \hadaprod \frac{[\bA\bZ^\top]}{[\bW\bZ\bZ^\top]},
\end{aligned}
$$
where $\bA\in\real_+^{M\times N}, \bW\in\real_+^{M\times K}$, and $\bZ\in\real_+^{K\times N}$. 
The operator $\frac{[\cdot]}{[\cdot]}$ denotes componentwise division, and $\hadaprod$ is the Hadamard product.

The MU update requires that $\bZ$ and $\bW$ should be initialized with positive (nonzero) entries; otherwise, the MU will not modify any entries due to the Hadamard product.
\end{theoremHigh}

The MU approach has played a pivotal role in the development of NMF, becoming a cornerstone of the field for several reasons: (i) The update rules are extremely easy to implement; (ii) In practice, the convergence is relatively fast compared to many other methods; (iii) Nonnegativity is preserved automatically without explicit projection.
To prove Theorem~\ref{theorem:conv_mu_fro}, we use the \textit{majorization-minimization (MM)} framework, which relies on the concept of an auxiliary function.
\begin{definition}[Auxiliary Function (Majorizer)]\label{definition:aux_func}
A function $G(\bx, \widetildebx)$ is  an \textit{auxiliary function} for $F(\bx)$ (or a majorizer of $F$ at $\widetildebx$) if, for all $\bx$~\footnote{$\bx$ can be scalars, vectors, or matrices.}
$$
G(\bx, \widetildebx) \geq F(\bx)
\gap 
\text{and}
\gap 
G(\bx, \bx) = F(\bx).
$$
In other words, the auxiliary function $G(\bx, \widetildebx)$ is an upper bound of $F(\bx)$, and the bound is tight when $\widetildebx=\bx$.
\end{definition}

\begin{lemma}[Nonincreasing in Auxiliary Functions]\label{lemma:noninmuaux}
If $ G $ is an auxiliary function for $F$, then $F$ is nonincreasing under the update
\begin{equation}\label{equation:aux_update}
\bx^{(t+1)} = \mathop{\argmin}_{\bx} \, G(\bx, \bx^{(t)}).
\end{equation}
\end{lemma}
\begin{proof}[of Lemma~\ref{lemma:noninmuaux}]
By definition,  $ F(\bx^{(t+1)}) \leq G(\bx^{(t+1)}, \bx^{(t)}) \leq G(\bx^{(t)}, \bx^{(t)}) = F(\bx^{(t)})$.
\end{proof}

Note that $ F(\bx^{(t+1)}) = F(\bx^{(t)})$ only if $ \bx^{(t)}$ is a local minimum of $ G(\bx, \bx^{(t)})$ w.r.t. $\bx$. If the partial derivatives of $ F$ exist and are continuous in a small neighborhood of $ \bx^{(t)}$, this also implies that the gradient $ \nabla F(\bx^{(t)}) = \bzero$. Thus, by iterating the update in \eqref{equation:aux_update}, we obtain a sequence of estimates that converge to a local minimum $ \bx_{\min} = \argmin_{\bx} F(\bx)$ of the objective function:
\begin{equation}
 F(\bx^{(0)}) \geq 	 F(\bx^{(1)}) \geq  F(\bx^{(2)})\geq \ldots  \geq  F(\bx^{(t)}) \geq  F(\bx^{(t+1)})\geq \ldots \geq F(\bx_{\min}).
\end{equation}
Definition~\ref{definition:aux_func} finds a majorizer $G$ of $F$, and Lemma~\ref{lemma:noninmuaux} shows the minimization property in $G$, hence the algorithm is often referred to as the \textit{majorization-minimization (MM) framework}.
The update benefits when the global minimizer of $G$ admits a closed-form solution or can be computed efficiently.

Therefore, if we can find an appropriate auxiliary function $ G(\bx, \bx^{(t)})$ for both variables in $\normf{\bA-\bW\bZ}$, the update rules in  Theorem~\ref{theorem:conv_mu_fro} follow from  \eqref{equation:aux_update}.
To apply the auxiliary function to the NMF problem, we consider a column in $\bA$ or $\bZ$: $\ba\triangleq\ba_n$ and $\bz\triangleq\bz_n$  in the following lemma, where $n\in\{1,2,\ldots,N\}$.

\begin{lemma}[Auxiliary Function for NMF]\label{lemma:aux_nmf}
Let $\bW\in\real^{K\times N}, \ba\in\real^M$, and $\bz\in\real^{K}$.
Let further $\bD\in\real^{K\times K}$ be a diagonal matrix with the $(k,k)$-th entry being $d_{kk}=\frac{(\bW^\top\bW \bz)_k}{z_k}=\frac{\bw_k^\top\bW\bz}{z_k} = \frac{\sum_{j=1}^{K} (\bW^\top\bW)_{kj}z_j}{z_k}, \, \,\forall k\in\{1,2,\ldots,K\}$, where $\bw_k$ is the $k$-th column of $\bW$ and $z_k$ is the $k$-th component of $\bz$.
Then, the following function is an auxiliary function for  $F(\bz)=\frac{1}{2}\normtwo{\ba-\bW\bz}^2$:
\[ 
G(\bz, \bz^{(t)}) = F(\bz^{(t)}) + (\bz-\bz^{(t)})^\top \nabla F(\bz^{(t)}) +\frac{1}{2}(\bz-\bz^{(t)})^\top \bD (\bz-\bz^{(t)}).
\]
\end{lemma}
\begin{proof}[of Lemma~\ref{lemma:aux_nmf}]
Since the third-order partial derivatives of $F(\bz)$ vanish (see Problem~\ref{prob:third_order_nmf}), $F(\bz)$ can be factored as 
$$
F(\bz) = F(\bz^{(t)})+(\bz-\bz^{(t)})^\top \nabla F(\bz^{(t)}) +\frac{1}{2}(\bz-\bz^{(t)})^\top \bW^\top\bW (\bz-\bz^{(t)}).
$$
Apparently, $G(\bz, \bz)=F(\bz) $. To complete the proof, we need to show that $G(\bz, \bz^{(t)}) \geq F(\bz)$; that is, $\bD-\bW^\top\bW$ is positive semidefinite.
To prove this, consider the matrix $\bM\in\real^{K\times K}$ whose entries are $m_{ij}=z_i (\bD-\bW^\top\bW)_{ij}z_j$ for all $i,j\in\{1,2,\ldots,K\}$, which is a rescaling of the components of $\bD-\bW^\top\bW$. Then $\bD-\bW^\top\bW$ is positive semidefinite if and only if $\bM$ is: 
$$
\begin{aligned}
\bx^\top& \bM\bx
= \sum_{i,j=1}^{K,K} x_i m_{ij}x_j
\stackrel{*}{=}\sum_{i,j=1}^{K,K} \left\{(\bW^\top\bW)_{ij} z_i z_j  x_i^2 - (\bW^\top\bW)_{ij}z_i z_j  x_i x_j\right\}\\
&\stackrel{\dag}{=}\sum_{i,j=1}^{K,K} (\bW^\top\bW)_{ij}z_i z_j \left( \frac{1}{2}x_i^2 + \frac{1}{2}x_j^2 - x_ix_j\right)
=\sum_{i,j=1}^{K,K} (\bW^\top\bW)_{ij}z_i z_j\frac{1}{2} \left( x_i-x_j\right)^2 \geq 0,
\end{aligned}
$$
where the equality $(\dag)$ follows from the symmetry of $\bM$, and the equality ($*$) follows from the diagonality of $\bD$:
$$
\sum_{i,j=1}^{K,K} x_iz_i d_{ij}z_jx_j
=
\sum_{i=1}^{K} x_iz_i d_{ii}z_ix_i
=
\sum_{i=1}^{K} x_i^2 z_i^2 \frac{\sum_{j=1}^{K}(\bW^\top\bW)_{ij}z_j}{z_i}
=
\sum_{i,j=1}^{K,K} (\bW^\top\bW)_{ij} z_i z_j  x_i^2.
$$
This completes the proof.
\end{proof}

Theorem~\ref{theorem:conv_mu_fro} follows directly from Lemma~\ref{lemma:aux_nmf} by minimizing $G(\bz,\bz^\toptzero)$ with respect to $\bz$, which yields the MU update.
It is generally better to update $\bW$ and $\bZ$ ``simultaneously” rather than ``sequentially," i.e., updating each matrix completely before the other. In this case, after updating a row of $\bZ$, we update the corresponding column of $\bW$.  
In the implementation, it is advisable to introduce a small positive quantity, say the square root of the machine precision, to the denominators in the approximations of $\bW$ and $\bZ$
at each iteration. 
And a trivial value like $\epsilon=10^{-9}$  suffices. The full procedure is shown in Algorithm~\ref{alg:nmf-multiplicative}.
In practice, the algorithm can also be accelerated by updating $\bW$ several times before updating $\bZ$, during which process we can reuse the result of $\bA\bZ^\top$ and $\bZ\bZ^\top$, and vice versa.

\index{Machine precision}
\begin{algorithm}[h] 
\caption{NMF via Multiplicative Updates}
\label{alg:nmf-multiplicative}
\begin{algorithmic}[1] 
\Require Matrix $\bA\in \real_+^{M\times N}$;
\State Initialize $\bW\in \real_{++}^{M\times K}$, $\bZ\in \real_{++}^{K\times N}$ randomly with positive entries;
\State Choose a stopping criterion on the approximation error $\delta$;
\State Choose maximal number of iterations $C$;
\State $iter=0$; \Comment{Count for the number of iterations}
\While{$\normf{\bA- (\bW\bZ)}^2>\delta $ and $iter<C$}
\State $iter=iter+1$; 
\State $\bZ \leftarrow \bZ\hadaprod \frac{[\bW^\top\bA]}{[\bW^\top\bW\bZ]+\epsilon}$;
\State $\bW \leftarrow \bW \hadaprod \frac{[\bA\bZ^\top]}{[\bW\bZ\bZ^\top]+\epsilon}$;
\EndWhile
\State Output $\bW,\bZ$;
\end{algorithmic} 
\end{algorithm}

\index{Regularization}
\subsection*{Regularization}
As noted  in \eqref{equation:kkv_nnn_raw},  the NNLS or NMF problem implicitly imposes a \textbf{sparsity constraint}. 
However, similar to regularized ALS (Section~\ref{section:regularization-extention-general}), adding explicit regularization can improve generalization and numerical stability. Consider the regularized objective:
$$
L(\bW,\bZ)  =\frac{1}{2}\normf{\bW\bZ-\bA}^2 +\frac{1}{2}\lambda_w \normf{\bW}^2 + \frac{1}{2}\lambda_z \normf{\bZ}^2, \qquad \lambda_w>0, \lambda_z>0,
$$
where the employed matrix norm is still the Frobenius norm. The gradient with respect to $\bZ$ (given $\bW$) is the same as that in Equation~\eqref{equation:als-regulari-gradien}:
$$
\begin{aligned}
	\frac{\partial L(\bZ|\bW)}{\partial \bZ} =\bW^\top(\bW\bZ-\bA) + \textcolor{mylightbluetext}{\lambda_z\bZ}  \in \real^{K\times N}.
\end{aligned}
$$
Repeating the MU derivation with this modified gradient and the same adaptive step size 
$\eta_{kn} = \frac{z_{kn}}{(\bW^\top\bW\bZ)_{kn}}$,
we obtain the \textit{regularized MU rules}:
$$
\begin{aligned}
\bZ \leftarrow  \bZ\hadaprod \frac{[\bW^\top\bA-\lambda_z\bZ]}{[\bW^\top\bW\bZ]}
\gap\text{and}\gap 
\bW \leftarrow \bW \hadaprod \frac{[\bA\bZ^\top-\lambda_w\bW]}{[\bW\bZ\bZ^\top]}.
\end{aligned}
$$

However, the numerators may become negative, violating nonnegativity. Two common fixes are: 
\paragrapharrow{(MMU1).} Clamp the entire update:
$$
\textbf{(MMU1):}\quad  \bZ \leftarrow \left[\bZ\hadaprod \frac{[\bW^\top\bA-\lambda_z\bZ]}{[\bW^\top\bW\bZ]}\right]_+;
\gap\text{and}\gap 
\bW \leftarrow \left[\bW \hadaprod \frac{[\bA\bZ^\top-\lambda_w\bW]}{[\bW\bZ\bZ^\top]}\right]_+,
$$
where $[x]_+ = \max\{x, \epsilon\}$. The parameter $\epsilon$ is usually a very small positive number that prevents the emergence of negative update.
That is, we add a small lower bound for entries of $\bW$ and $\bZ$.
\paragrapharrow{(MMU2).} Clamp only the numerator:
$$
\begin{aligned}
\textbf{(MMU2):}\quad 
\bZ \leftarrow \bZ\hadaprod \frac{[\bW^\top\bA-\lambda_z\bZ]_+}{[\bW^\top\bW\bZ]}
\gap\text{and}\gap 
\bW \leftarrow\bW \hadaprod \frac{[\bA\bZ^\top-\lambda_w\bW]_+}{[\bW\bZ\bZ^\top]}.
\end{aligned}
$$
Both strategies ensure nonnegativity while incorporating regularization effects.

\section{NMF with Three Factors}\label{section:nmf_tree_fac}
The NMF framework can be extended to involve three nonnegative factor matrices, a formulation known as \textit{nonnegative matrix  tri-factorization (tri-NMF or NMTF)}. 
This approach approximates the data matrix as
\begin{equation}\label{equation:tri_nmv}
\bA \approx \bW\bU\bZ,
\end{equation}
where $\bW\in\real_+^{M\times K}$, $\bU\in\real_+^{K\times J}$, and $\bZ\in\real_+^{J\times N}$.
Consider a user-item interaction matrix $\bA\in\real^{M\times N}$, where each element is a binary number $\{0,1\}$---a setting commonly referred to as \textit{implicit feedback} data, in contrast to the \textit{explicit rating} data used previously. (For instance, in datasets like Netflix or MovieLens, ratings above 4 might be mapped to 1, and ratings below or equal to a threshold (e.g., 2) to 0, yielding an implicit binary dataset.)
Standard NMF decomposes $\bA$ as a sum of $K$ rank-one components: $\bA\approx\sum_{k=1}^{K} \bW[:,k]\bZ[k,:]$, 
where each term captures a latent pattern linking a subset of users (via $\bW[:,k]$) to a subset of items (via $\bZ[k,:]$).
In the context of implicit data, each rank-one matrix can be interpreted as finding a subset of users and a subset of items (e.g., movies)  that interact strongly with each other.
In contrast, tri-NMF yields a double-sum decomposition:
$$
\bA\approx\sum_{k=1}^{K}\sum_{j=1}^{J} \bW[:,k]\bU[k,j]\bZ[j,:].
$$
This formulation can be interpreted as identifying  separately $J$ subsets of movies that are watched together (the rows of $\bZ$) and $K$ subset of users that behave similarly (the columns of $\bW$); while the matrix $\bU$ tells us how these subsets interact together.
Specifically, if $u_{kj}>0$, then the $k$-th subset of users (corresponding to the positive entries of $\bW[:,k]$) watches the movies from the $j$-th subset of movies (corresponding to the positive entries of $\bZ[j,:]$).

In essence, tri-NMF simultaneously discovers clusters of similar users and clusters of similar items, and explicitly models their pairwise associations through the nonnegative interaction matrix $\bU$.
This framework also finds application in text mining: there, tri-NMF can identify groups of documents that share common words (columns of $\bW$), groups of words that co-occur across similar documents (rows of $\bZ$), and their interplay via $\bU$ \citep{brouwer2017comparative, gillis2020nonnegative}.

\section{$\beta$-Divergence, Alternative Perspectives of MU}\label{section:beta_div_altmu}

We have introduced several alternative error measures for matrix factorization problems in Section~\ref{section:more_err_sta_als}.
As noted earlier, the sum-of-squared loss---given in either \eqref{equation:als-per-example-loss2} or \eqref{equation:frob_nmf}---is convex in one factor when the other is held fixed, which facilitates a smooth optimization process.
This loss belongs to a broader family of dissimilarity measures known as \textit{$\beta$-divergences}, commonly used in NMF.
For two nonnegative scalars $x$ and $y$, the $\beta$-divergence is defined as:
\begin{equation}\label{equation:scalar_beta_div}
d_{\beta}(x, y)
=
\left\{
\begin{aligned}
&\frac{x}{y}-\ln \frac{x}{y}-1, &\text{if } \beta=0;\\
&x\ln\frac{x}{y} -x+y , &\text{if } \beta=1;\\
&\frac{1}{\beta^2-\beta}(x^\beta+(\beta-1)y^\beta - \beta xy^{\beta-1} ) , &\text{otherwise}.\\
\end{aligned}
\right.
\end{equation}
The $\beta$-divergence is continuous in $\beta$ since $\mathoplim{\beta\rightarrow 0}(x^\beta -y^\beta)/\beta=\ln(x/y)$. 
When $\beta=0, 1,$ and $2$, the $\beta$-divergences are also known as the \textit{Itakura--Saito (IS), KL, and Frobenius/Euclidean distances/divergences}, respectively.
The $\beta$-divergence between two matrices $\bB, \bC\in\real^{M\times N}$ is defined componentwise as: 
\begin{equation}\label{equation:beta_div_mat_def}
D_{\beta}(\bB,\bC)=
\sum_{n=1}^N d_{\beta}(\bb_n, \bc_n)
=
\sum_{m,n=1}^{M,N} d_{\beta}({b_{mn}, c_{mn}}).
\end{equation}

\begin{figure}[h]
\centering  
\vspace{-0.15cm}    
\subfigtopskip=2pt  
\subfigbottomskip=2pt 
\subfigcapskip=-5pt  
\subfigure[$\beta$-divergence for $d_{\beta}(1,y)$.]{\label{fig:beta_divergence_1}
\includegraphics[width=0.47\linewidth]{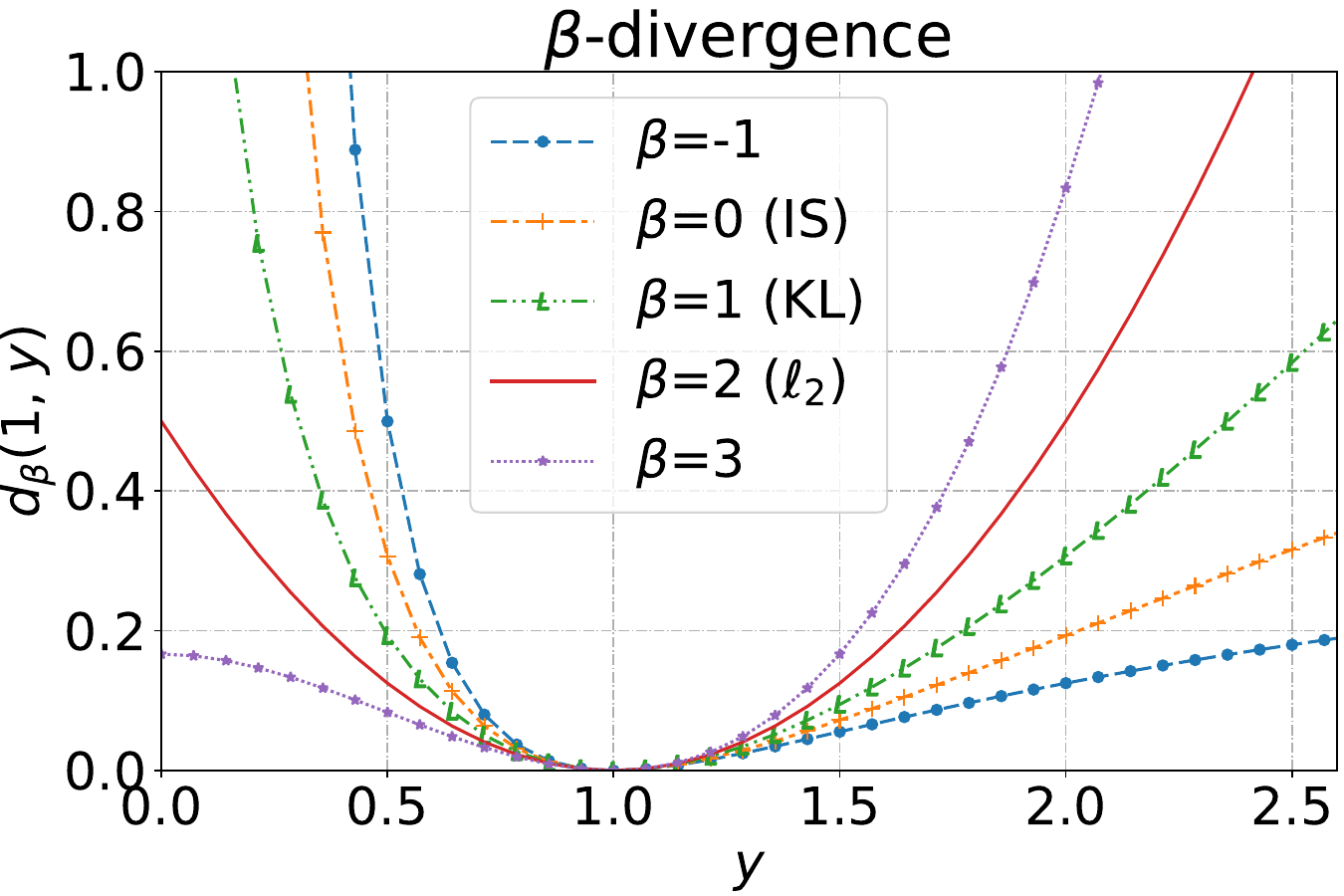}}
\subfigure[$\beta$-divergence for $d_{\beta}(2,y)$.]{\label{fig:beta_divergence_2}
\includegraphics[width=0.47\linewidth]{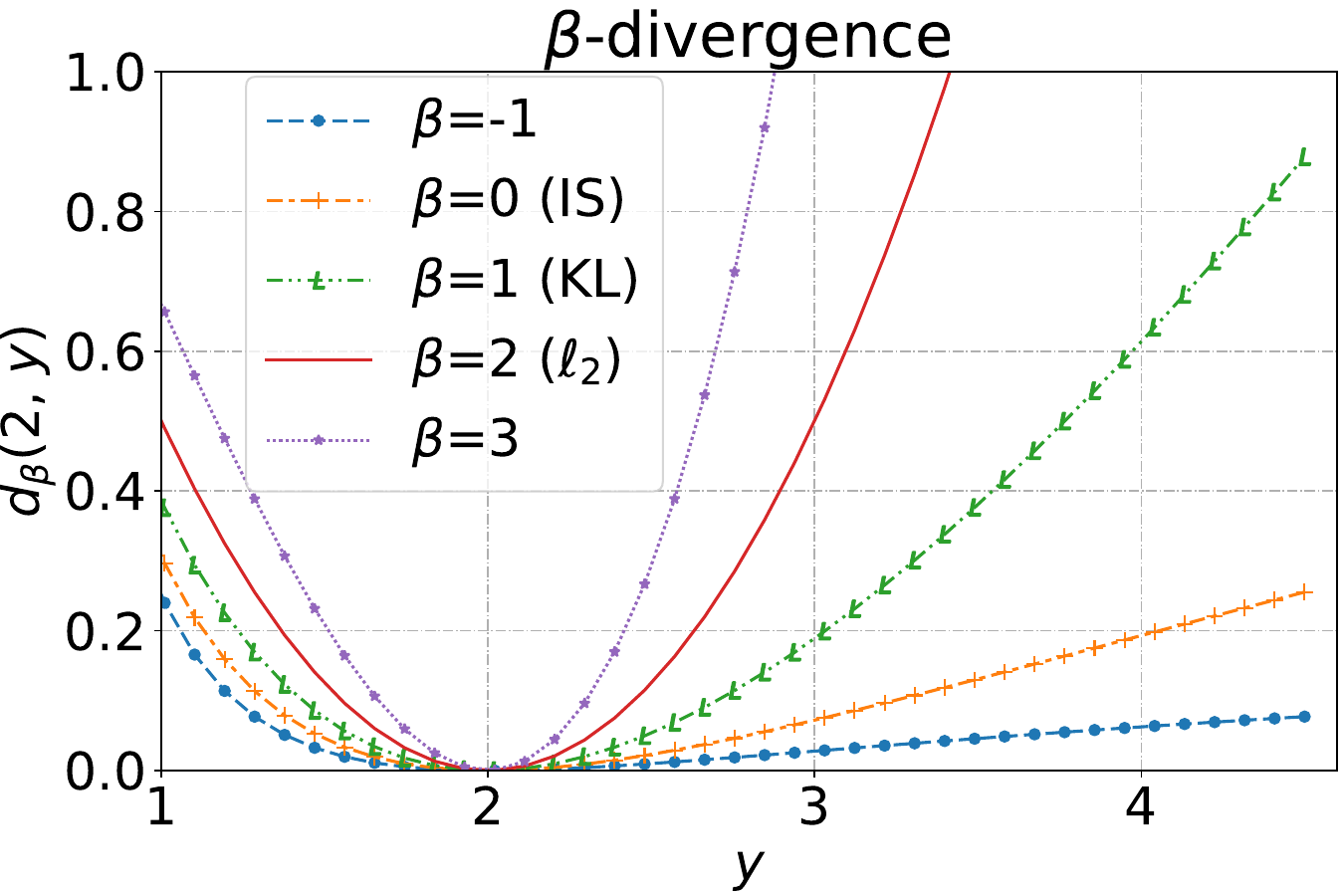}}
\caption{
The analysis of $\beta$-divergence is complex. When the first argument is fixed at 1,  smaller values are less penalized as  the $\beta$ value increases; however, when the first argument is $2$,  smaller values  are more penalized as  the $\beta$ value increases.
In both cases, larger values are more heavily penalized as the $\beta$ value increases.
}
\label{fig:beta_divergence_all}
\end{figure}

The behavior of $\beta$-divergence is nuanced. 
As illustrated in Figure~\ref{fig:beta_divergence_all}, when the first argument is fixed at 1,  smaller values are \textit{less} penalized as  the $\beta$ value increases; however, when the first argument is $2$,  smaller values  are \textit{more} penalized as  the $\beta$ value increases.
In both cases, larger values are more heavily penalized as the $\beta$ value increases. 

\paragrapharrow{Convexity of $\beta$-divergence.}
For $\beta\in[1,2]$, the scalar function $d_{\beta}(x,y)$ is convex in the second argument $y$. 
Consequently, the matrix divergence $D_{\beta}(\bA,\bW\bZ)$ is convex in $\bW$ when fixing $\bZ$, and vice versa (see Problem~\ref{prob:sep_conv}).  
This property ensures that coordinate descent algorithms are well-suited for NMF under $\beta$-divergence in this range.

\paragrapharrow{Scaling property.} Let $\gamma>0$ be a scale factor. Then:
\begin{equation}
d_{\beta}(\gamma x, \gamma y) = \gamma^\beta d_{\beta}(x,y).
\end{equation}
This indicates that the larger the $\beta$, the more sensitive the $\beta$-divergence is to large values of $x$ or $y$; on the contrary, $\beta$-divergence with small $\beta<0$ values relies more heavily on the smallest data values.
However, when $\beta=0$ (the \textit{Itakura--Saito divergence, IS divergence}, see \eqref{equation:als_gama_los}),
the divergence depends only on the ratio $x/y$, making it invariant to global scaling---a property evident from \eqref{equation:scalar_beta_div}.

\paragrapharrow{Gradient and domain considerations.} 
In NMF, the data matrix $\bA$ is nonnegative, but special care is needed when entries are zero.
Specifically, for $x=0$, the divergence $d_{\beta}(x=0, \cdot)$  is not defined for all  $\beta$:
$$
d_{\beta}(0, y)=
\left\{
\begin{aligned}
\text{not defined}, \gap &\text{if $\beta\leq 0$};\\
\frac{1}{\beta} y^\beta, \gap \gap &\text{if $\beta> 0$},
\end{aligned}
\right.
\gapthree \implies \gapthree
d'_{\beta}(0, y)=
\left\{
\begin{aligned}
\text{not defined}, \gap &\text{if $\beta\leq 0$};\\
 y^{\beta-1}, \gap \gap &\text{if $\beta> 0$},
\end{aligned}
\right.
$$
where the derivative $d'_{\beta}(0, y)$  is taken with respect to the second argument $y$.
Therefore, algorithms based on $\beta$-divergence with $\beta\leq 0$ require strictly positive input matrices 
$\bA$.
Table~\ref{table:dom_beta_1} and Table~\ref{table:dom_beta_2} summarize the domains of $d_\beta(x, \cdot)$ and $d'_{\beta}(x, \cdot)$, respectively, for different values of $\beta$ and $x$.
\noindent
\begin{table}[h]
\centering
\begin{minipage}{0.37\textwidth}
\centering
\caption{Domain of $d_\beta(x, \cdot)$.}
\label{table:dom_beta_1}
\begin{tabular}{c|c|c|c}
\hline
& $\beta \leq 0$ & $\beta \in (0,1]$ & $\beta > 1$ \\
\hline
$x = 0$ & $\varnothing$ & $\real_+$ & $\real_+$ \\
$x > 0$ & $\real_{++}$ & $\real_{++}$ & $\real_+$ \\
\hline
\end{tabular}
\end{minipage}
\quad
\begin{minipage}{0.54\textwidth}
\centering
\caption{Domain of $d'_{\beta}(x, \cdot)$.}
\label{table:dom_beta_2}
\begin{tabular}{c|c|c|c|c}
\hline
& $\beta \leq 0$ & $\beta \in (0,1)$ & $\beta \in [1,2)$ & $\beta \geq 2$ \\
\hline
$x = 0$ & $\varnothing$ & $\real_{++}$ & $\real_+$ & $\real_+$ \\
$x > 0$ & $\real_{++}$ & $\real_{++}$ & $\real_{++}$ & $\real_+$ \\
\hline
\end{tabular}
\end{minipage}
\end{table}

When the gradients exist, the partial derivatives of  $D_{\beta}(\bA,\bW\bZ)$ with respect to $\bZ$ and $\bW$ are:
\begin{equation}
\begin{aligned}
\nabla_{\bZ} D_{\beta}(\bA,\bW\bZ) &= \bW^\top \big( (\bW\bZ)^{\beta-2} \hadaprod (\bW\bZ-\bA) \big);\\
\nabla_{\bW} D_{\beta}(\bA,\bW\bZ) &=  \big( (\bW\bZ)^{\beta-2} \hadaprod (\bW\bZ-\bA) \big)\bZ^\top, 
\end{aligned}
\end{equation}
where $(\bW\bZ)^{\beta-2}$ denotes the componentwise exponentiation.
For $\beta=2$, these expressions reduce to the familiar gradient in \eqref{equation:givenw-update-z-allgd}.

\paragrapharrow{Convex-concave decomposition.}
The $\beta$-divergence can be decomposed into convex, concave, and constant components with respect to the second argument $y$:
\begin{equation}\label{equation:beta_decom_defi}
d_{\beta}(x,y) = \convd_{\beta}(x,y)+\concd_{\beta}(x,y)+\cnstd_{\beta}(x,y),
\end{equation} 
where $\convd_{\beta}(x,y)$ is convex in $y$, $\concd_{\beta}(x,y)$ is concave in $y$, and $\cnstd_{\beta}(x,y)$ is constant in $y$. 
Note that this decomposition is not unique---any affine term can be assigned to either the convex or concave part---but we follow the convention in \citet{fevotte2011algorithms}.
Table~\ref{table:beta_decom} provides the explicit forms for different ranges of $\beta$.

\begin{table}[h]
\setlength{\tabcolsep}{6.pt}    
\centering
\begin{tabular}{c|c|c|c}
\hline
& $\convd_{\beta}(x,y)/\convd'_{\beta}(x,y)$, convex & $\concd_{\beta}(x,y)/\concd'_{\beta}(x,y)$, concave & $\cnstd_{\beta}(x,y)$, constant \\ \hline\hline
$\beta < 1, \beta \neq 0$ & $-\frac{1}{\beta-1}xy^{\beta-1}/ -xy^{\beta-2}$ & $\frac{1}{\beta}y^{\beta}/y^{\beta-1}$ & $\frac{1}{\beta(\beta-1)}x^{\beta}$ \\ \hline
$\beta = 0$ & $xy^{-1}/-xy^{-2}$ & $\ln y/y^{-1}$ & $x(\ln x - 1)$ \\ \hline
$1 \leq \beta \leq 2$ & $d_{\beta}(x,y)/d'_{\beta}(x,y)$ & 0/0 & 0 \\ \hline
$\beta > 2$ & $\frac{1}{\beta}y^{\beta}/y^{\beta-1}$ & $-\frac{1}{\beta-1}xy^{\beta-1}$ & $\frac{1}{\beta(\beta-1)}x^{\beta}$ \\ \hline
\end{tabular}
\caption{Scalar convex-concave-constant decomposition  of $d_{\beta}(x,y)$ with respect to the second variable $y$, along with corresponding derivatives.}
\label{table:beta_decom}
\end{table}

\paragrapharrow{KKT conditions for NMF with $\beta$-divergence.}
The KKT conditions  and derivation in \eqref{equation:kkv_nnn_raw}) for a stationary point of the NMF problem under $\beta$-divergence are:
\begin{equation}\label{equation:nmf_beta_kkt1}
\begin{aligned}
\bZ\geq \bzero,\gap&\nabla_{\bZ} D_{\beta}(\bA, \bW\bZ) &\geq& \bzero, \gap \langle \bZ, \nabla_{\bZ} D_{\beta}(\bA, \bW\bZ)\rangle &=&\bzero_{K\times N}; \\
\bW\geq \bzero,\gap&\nabla_{\bW} D_{\beta}(\bA, \bW\bZ) &\geq& \bzero, \gap \langle \bW, \nabla_{\bW} D_{\beta}(\bA, \bW\bZ)\rangle &=&\bzero_{M\times K}.
\end{aligned}
\end{equation}
with all inequalities and inner products interpreted componentwise. 
Equivalently,
\begin{equation}
\begin{aligned}
\min\{\bZ, \nabla_{\bZ} D_{\beta}(\bA, \bW\bZ)\} = \bzero_{K\times N}
\gap \text{and}\gap 
\min\{\bW, \nabla_{\bW} D_{\beta}(\bA, \bW\bZ)\} = \bzero_{M\times K},
\end{aligned}
\end{equation}
where the \texttt{min} operator $\min\{\cdot, \cdot \}$ is applied componentwise.

\subsection{MU for $\beta$-Divergence Obtained via  Gradient Ratio Heuristic}\label{section:mu_gd_ratio}
We have shown that the MU update with Frobeius norm can be derived from rescaled gradient descent. 
For brevity, decompose the gradient  with respect to $\bZ$ as 
\begin{equation}
\nabla_{\bZ}\triangleq\nabla_{\bZ} D_{\beta}(\bA,\bW\bZ) = \nabla_{\bZ}^+ - \nabla_{\bZ}^-,
\end{equation}
where 
\begin{equation}\label{equation:mu_grati_decom}
\nabla_{\bZ}^+ = \bW^\top \big( (\bW\bZ)^{\beta-1} \big)
\gap \text{and}\gap 
\nabla_{\bZ}^- = \bW^\top \big( (\bW\bZ)^{\beta-2} \hadaprod \bA \big).
\end{equation}
When $z_{kn}>0, \forall k,n$, the KKT conditions \eqref{equation:nmf_beta_kkt1} imply that $(\nabla_{\bZ}^+)_{kn} = (\nabla_{\bZ}^-)_{kn}$. 
In standard gradient descent (i.e., $\bZ^{(t+1)}\leftarrow \bZ^{(t)} - \eta\nabla_{\bZ}$ at iteration $t$) indicates a small decrease (resp. increase) of $z_{kn}$ will lead to a  decrease of the loss function if $(\nabla_{\bZ})_{kn}>0$ (resp. $<0$).
This motivates a multiplicative update based on the componentwise ratio of the negative and positive gradient parts:
\begin{equation}\label{equation:mu_grati_z}
\bZ\leftarrow  \bZ \hadaprod \frac{[\nabla_{\bZ}^-]}{[\nabla_{\bZ}^+]},
\end{equation}
where $\frac{[\cdot]}{[\cdot]}$ denotes elementwise division.
When $\beta=2$, this recovers the Frobenius-norm MU rule in Theorem~\ref{theorem:conv_mu_fro}.
When $\beta=1$, the loss becomes the KL divergence, and the update simplifies to:
$$
\textbf{($\beta=1$)}:\gap \bZ\leftarrow \bZ \hadaprod \frac{[\bW^\top \frac{[\bA]}{[\bW\bZ]}]}{[\bW^\top \bone_{M\times N}]}.
$$
It can be shown that for $\beta\in[1,2]$, these MU updates monotonically decrease the objective $D_{\beta}(\bA, \bW\bZ)$.

\subsection{MU for $\beta$-Divergence  Obtained via Rescaled PGD}
As discussed in Section~\ref{section:nmf_apgd},
the PGD approach updates a variable by taking a gradient step and projecting back onto the feasible set. 
Consider a standard GD update on $f(\bx)$: $\bx^{(t+1)}\leftarrow \bx^{(t)}-\eta\nabla f(\bx^{(t)})$, where $\eta$ is a step size and  $-\nabla f(\bx^{(t)})$ is a \textit{descent direction} ($\bg$ is a descent direction if $\bg^\top\nabla f(\bx^{(t)})<0$.)
Consider further  a diagonal $\bD$ such that $-\eta \nabla f(\bx^{(t)}) \rightarrow -\bD \nabla f(\bx^{(t)})$ is also a descent direction (i.e., replacing the step size with a diagonal matrix)~\footnote{Indeed, $\bD$ could be any positive definite matrix, though diagonal choices preserve separability.}.
In this case, if the feasible set of $\bx$ is nonnegative, then the PGD is useful: $\bx^{(t+1)}\leftarrow \mathcalP(\bx^{(t)}-\bD\nabla f(\bx^{(t)}))$, where $\mathcalP(x)=\max\{x, 0\}$ enforces nonnegativity \citep{lu2025practical, lu2026first}.
Now decompose the gradient into its positive and negative parts: 
$$
\nabla f(\bx^{(t)}) = \nabla^+f(\bx^{(t)}) - \nabla^-f(\bx^{(t)}),
$$ 
with $\nabla^+f(\bx^{(t)})>0 $ and $\nabla^-f(\bx^{(t)})>0$ elementwise. 
Choosing the scaling matrix as $\bD=\diag\big( \frac{[\bx^{(t)}]}{[\nabla^+f(\bx^{(t)}))]}  \big)$, the rescaled PGD update becomes a MU rule:
\begin{equation}
\bx^{(t+1)}\leftarrow 
\mathcalP\bigg(\bx^{(t)}-\diag\big( \frac{[\bx^{(t)}]}{[\nabla^+f(\bx^{(t)})]}  \big)\nabla f(\bx^{(t)})\bigg)
=
\mathcalP\bigg(\bx^{(t)} \hadaprod \frac{[\nabla^-f(\bx^{(t)})]}{[\nabla^+f(\bx^{(t)})]}\bigg).
\end{equation}
Applying this to the NMF gradient decomposition in \eqref{equation:mu_grati_decom} yields exactly the MU rule in \eqref{equation:mu_grati_z}.
If we further incorporate a step size $\eta\in(0,1)$ in the rescaled PGD update, the update becomes a convex combination:
\begin{equation}
\bx^{(t+1)}
=
\mathcalP\bigg((1-\eta)\bx^{(t)} + \eta\bx^{(t)} \hadaprod \frac{[\nabla^-f(\bx^{(t)})]}{[\nabla^+f(\bx^{(t)})]}\bigg).
\end{equation}
Since $-\bD \nabla f(\bx^{(t)})$ remains a descent direction, any $\eta \in(0,1)$  ensures that the update is monotonically nonincreasing.
Moreover, because all terms are nonnegative, the projection operator $\mathcalP$ is redundant and may be omitted.

\subsection{MU for $\beta$-Divergence Obtained via MM Framework}
The $\beta$-divergence between two matrices can be defined columnwise (see Equation~\eqref{equation:beta_div_mat_def}), 
and each scalar $\beta$-divergence admits a decomposition into three parts---convex, concave, and constant---with respect to the second argument (see Equation~\eqref{equation:beta_decom_defi}).
Consequently, the NMF loss function can be decomposed as follows (note that it can also be split further componentwise):
$$
D_{\beta}(\bA, \bW\bZ)
=
\sum_{n=1}^{N}d_{\beta}(\ba_n, \bW\bz_n)
=
\sum_{n=1}^{N}\left(\convd_{\beta}(\ba_n, \bW\bz_n) + \concd_{\beta}(\ba_n, \bW\bz_n) + \cnstd_{\beta}(\ba_n, \bW\bz_n)\right).
$$
In the majorization-minimization (MM) framework, we construct an auxiliary function for each column $n$ by handling the three components separately. This approach is justified by the following lemma:
\begin{lemma}[Auxiliary Function By Parts]
Let $F(\bx)=\sum_{k=1}^K F_k(\bx)$, and  let $G_k(\bx, \widetildebx)$ be an auxiliary function for $F_k(\bx)$ at $\widetildebx$ for all $k$. 
Then, $G(\bx, \widetildebx)=\sum_{k=1}^K G_k(\bx, \widetildebx)$ is an auxiliary function for $F(\bx)$ at $\widetildebx$.
\end{lemma}
This lemma shows that constructing auxiliary functions componentwise allows us to decouple the overall optimization problem.

\paragrapharrow{Constant part.}
No auxiliary function is needed for the constant term  $\cnstd_{\beta}(\ba_n, \bW\bz_n)$  as it does not depend on $\bz_n$ and therefore has no effect on the minimization of $d_{\beta}(\ba_n, \bW\bz_n)$
with respect to $\bz_n$.

\paragrapharrow{Concave part.}
Any concave function can be upper-bounded using  linearization (the tangent plane, $f(x) +\nabla f(x)\cdot (y-x)\geq  f(y), \text{ for any $x,y\in\sS$ if $f:\sS\rightarrow \real$ is concave}$). 
Applying this to the concave component yields:
$$
\concd_{\beta}(x, y) 
\leq 
\concd_{\beta}(x, \widetildey)
+
(y-\widetildey) \concd_{\beta}'(x, \widetildey), 
$$
where $\concd_{\beta}'(x, \widetildey)$ denotes the gradient of $\concd(x, \widetildey)$ with respect to its second component $\widetildey$.
Therefore, for any $\widetildebz_n\in\real^{K}$, an auxiliary function for the concave component $\concd_{\beta}(\ba_n, \bW\bz_n)$ can be constructed by  
$$
\concG(\bz_n, \widetildebz_n) = \concd_{\beta}(\ba_n, \bW\widetildebz_n) + (\bW\bz_n-\bW\widetildebz_n) \hadaprod \concd_{\beta}' (\ba_n, \bW\widetildebz_n).
$$

\paragrapharrow{Convex part.}
For the convex component, we apply the convexity inequality {(or Jensen's inequality: $
f\big(\sum_{i=1}^{p} \lambda_ix_i\big) \leq \sum_{i=1}^{p}\lambda_i f(x_i),
$ if $\lambda_i\geq 0$ and $\sum_{i=1}^{p}\lambda_i=1$ for a convex function $f$)}.
Specifically, we construct a matrix $\bP\in\real^{M\times K}$ with entries:
\begin{equation}
p_{mk} = \frac{w_{mk}{\widetildez}_{kn}}{\sum_{j}w_{mj}{\widetildez}_{jn}}
=
\frac{w_{mk}{\widetildez}_{kn}}{\bW[m,:]\widetildebz_n}
\gap \implies\gap
\bP\geq\bzero \text{ and } \bP\bone = \bone.
\end{equation}
That is, each row of $\bP$ belongs to the unit simplex in $\real^K$.
Therefore, we have
$$
\begin{aligned}
\convd_{\beta}(a_{mn}, \bW[m,:]\bz_n ) 
&=
\convd_{\beta}\big(a_{mn}, \sum_{k=1}^{K}w_{mk}z_{kn} \big) 
=
\convd_{\beta}\big(a_{mn}, \sum_{k=1}^{K} p_{mk}\frac{w_{mk}z_{kn}}{p_{mk}} \big) \\
&\leq 
\sum_{k=1}^{K} p_{mk}\convd_{\beta}\big(a_{mn}, \frac{w_{mk}z_{kn}}{p_{mk}} \big). 
\end{aligned}
$$
Combining these constructions leads to the following auxiliary function for the full objective $D_{\beta}(\bA,\bW\bZ)$ w.r.t. $\bZ$.
\begin{theoremHigh}[Auxiliary Function for $D_{\beta}(\bA,\bW\bZ)$ w.r.t. $\bZ$]\label{theorem:aux_db_z}
Let $\widetildeba_n=\bW\widetildebz_n$ with $\widetildea_{mn}\triangleq\bW[m,:]\widetildebz_n$ for all $m, n$, where $\widetildebz_n$ is any vector in $\real^K$. Then, $G(\bZ,\widetildebZ)=\sum_{n=1}^{N} G_{n}(\bz_n, \widetildebz_n) = \sum_{n=1}^{N} \sum_{m=1}^{M} G_{mn}$ is an auxiliary function for $D_{\beta}(\bA,\bW\bZ)$ w.r.t. $\bZ$, where 
$$
\begin{aligned}
G_{mn}
&=
\cnstd_{\beta}(a_{mn}, \widetildea_{mn})
+
\concd_{\beta}(a_{mn}, \widetildea_{mn}) + \sum_{k=1}^{K}w_{mk}(z_{kn}-\widetildez_{kn}) \concd_{\beta}'(a_{mn}, \widetildea_{mn})\\
&\qquad + \sum_{k=1}^{K} \frac{w_{mk}\widetildez_{kn}}{\widetildea_{mn}}\convd_{\beta}\big(a_{mn},  \frac{\widetildea_{mn}z_{kn}}{\widetildez_{kn}} \big).
\end{aligned}
$$
\end{theoremHigh}

\begin{exercise}[Gradient and Hessian of Auxiliary Functions]\label{exercise:gra_hes_aux}
Consider the setting and notation of  Theorem~\ref{theorem:aux_db_z}. 
Let $G_n(\bz_n, \widetildebz_n)= \sum_{k=1}^{K} G_k (\bz_{kn}, \widetildebz_n) + C(\bz_n)$, where $C(\bz_n)$ is a constant w.r.t. $\bz_n$. 
That is, 
$$
G_k (\bz_{kn}, \widetildebz_n) 
= 
\sum_{m=1}^{M}w_{mk} z_{kn}\concd_{\beta}'(a_{mn}, \widetildea_{mn})
+ 
\sum_{m=1}^{M} \frac{w_{mk}\widetildez_{kn}}{\widetildea_{mn}}\convd_{\beta}\big(a_{mn},  \frac{\widetildea_{mn}z_{kn}}{\widetildez_{kn}} \big).
$$
Show that gradient of the auxiliary function is 
$$
\nabla_{z_{kn}}G_n(\bz_n, \widetildebz_n) 
=
\sum_{m=1}^{M} 
w_{mk}
\bigg(
\concd_{\beta}'(a_{mn}, \widetildea_{mn})
+
\convd_{\beta}'\big(a_{mn},  \frac{\widetildea_{mn}z_{kn}}{\widetildez_{kn}} \big)
\bigg), 
$$
and the Hessian matrix is diagonal with entries
$$
\nabla^2_{z_{kn}}G_n(\bz_n, \widetildebz_n) 
=
\sum_{m=1}^{M} 
w_{mk}
\frac{\widetildea_{mn}}{\widetildez_{kn}}
\bigg(
\convd_{\beta}{''}\big(a_{mn},  \frac{\widetildea_{mn}z_{kn}}{\widetildez_{kn}} \big)
\bigg).
$$
Note in all cases, derivatives are taken with respect to the second argument of $d_{\beta}(\cdot, \cdot)$.
\end{exercise}

Since $\convd_{\beta}(\cdot, \cdot)$ is convex in its second argument, the Hessian is positive definite. Thus, the auxiliary function is convex.
These constructions result in the following theorem by minimizing the  auxiliary function obtained in Theorem~\ref{theorem:aux_db_z}.
\begin{theoremHigh}[Nonincreasing of MU for $\beta$-Divergence \citep{fevotte2011algorithms, gillis2020nonnegative}]\label{theorem:conv_mu_beta}
Let $\bA\in\real_+^{M\times N}$, $\bW\in\real_{++}^{M\times K}$, and $\bZ\in\real_{++}^{K\times N}$. 
The  loss $D_{\beta}(\bA,\bW\bZ)$ remains nonincreasing under the following multiplicative update rules:
$$
\footnotesize
\begin{aligned}
\bZ\leftarrow \bZ\hadaprod \left(\frac{\left[\bW^\top \left\{ (\bW\bZ)^{(\beta-2)} \hadaprod \bA \right\}\right]}{[\bW^\top (\bW\bZ)^{(\beta-1)}]}\right)^{m(\beta)},
\gapthree \text{and}\gapthree
\bW\leftarrow \bW\hadaprod \left(\frac{\left[\left\{ (\bW\bZ)^{(\beta-2)} \hadaprod \bA \right\} \bZ^\top \right]}{[ (\bW\bZ)^{(\beta-1)}\bZ^\top]}\right)^{m(\beta)},
\end{aligned}
$$
where 
$$
m(\beta)\triangleq 
\left\{
\begin{aligned}
&\frac{1}{2-\beta}, \gap &\textit{if }&\beta<1; \\
&1,                  &\textit{if }& 1\leq\beta\leq 2; \\
&\frac{1}{\beta-1}, &\textit{if }& \beta>1.
\end{aligned}
\right.
$$
When $\beta=2$, this reduces to the Frobenius-norm update in Theorem~\ref{theorem:conv_mu_fro}.
For $1\leq \beta\leq 2$, the MM-derived update coincides with the gradient-ratio heuristic described in Section~\ref{section:mu_gd_ratio}.
\end{theoremHigh}
The update in Theorem~\ref{theorem:conv_mu_beta} ensures nonnegativity of the parameter updates, provided they are initialized with positive values.

\paragrapharrow{Choice of $\beta$ for NMF.}
The selection of $\beta$ depends on the application.
\citet{fevotte2009nonnegative} show that using $\beta=0$ (Itakura--Saito divergence) to decompose a piano power spectrogram accurately captures components corresponding to very low-level residual noise and hammer strikes---features that are either ignored or severely distorted when using Euclidean ($\beta=2$) or KL ($\beta=1$) divergences.
Similarly, \citet{fitzgerald2009use} demonstrate that $\beta=0.5$ is optimal for music source separation tasks.

\paragrapharrow{Convergence.}
An algorithm is said to be  \textit{convergent} if it produces a sequence of iterates $\{\bZ^{(t)}\}_{t\geq 1}$ or $\{\bW^{(t)}\}_{t\geq 1}$ that converges to a limit point $\bW^*$ or $\bZ^*$ satisfying the KKT conditions in \eqref{equation:nmf_beta_kkt1}. Monotonic nonincreasingness does not imply convergence in general, and neither is monotonicity necessary for convergence. Proving convergence of the MU methods is beyond the scope of this book; we refer the readers to \citet{gillis2020nonnegative, fevotte2011algorithms} and references therein for further details.

\section{Initialization}
Like ALS, a significant challenge in NMF is the lack of guaranteed convergence to a global minimum.
In practice, convergence can be slow, and the algorithm often settles at a suboptimal local minimum.
In the preceding discussion, we initialized $\bW$ and $\bZ$ randomly. To address these limitations, several alternative initialization strategies have been proposed to obtain better starting points, with the aim of accelerating convergence and improving solution quality \citep{boutsidis2008svd, gillis2014and}. We briefly outline these methods below:
\begin{itemize}
\item \textit{Clustering-based initialization.} Apply a clustering algorithm (e.g., $K$-means) to the columns of $\bA$. 
Set the cluster centroids of the top $K$ clusters as the initial columns of $\bW$, and initialize $\bZ$ using a scaled version of the cluster assignment matrix---i.e., $z_{kn}\neq 0$ indicates that column  $\ba_n$ belongs to cluster $k$.

\item \textit{Subset selection.} Select  $K$ representative columns of $\bA$ to form the initial $\bW$. And analogously select $K$ rows of $\bA $to initialize $\bZ$.

\item \textit{SVD-based approach.} Let the optimal rank-$K$ approximation of $\bA$ be given by its truncated SVD: $\bA=\sum_{k=1}^{K}\sigma_k\bu_k\bv_k^\top$, where each factor $\sigma_k\bu_k\bv_k^\top$ is a rank-one matrix with possible negative values in $\bu_k$ and $\bv_k$, and nonnegative $\sigma_k$. Denote $[x]_+\triangleq\max(x, 0)$, we notice 
$$
\bu_k\bv_k^\top = [\bu_k]_+[\bv_k]_+^\top+[-\bu_k]_+[-\bv_k]_+^\top-[-\bu_k]_+[\bv_k]_+^\top-[\bu_k]_+[-\bv_k]_+^\top,
$$
where the first two rank-one factors in this decomposition are nonnegative.
Then, either $[\bu_k]_+[\bv_k]_+^\top$ or $[-\bu_k]_+[-\bv_k]_+^\top$ can be selected to replace the factor $\bu_k\bv_k^\top$. \citet{boutsidis2008svd} suggest to replace each rank-one factor in $\sum_{k=1}^{K}\sigma_k\bu_k\bv_k^\top$ with  either $[\bu_k]_+[\bv_k]_+^\top$ or $[-\bu_k]_+[-\bv_k]_+^\top$, selecting the one with the larger norm and scaling it properly.
In other words, if we select $[\bu_k]_+[\bv_k]_+^\top$, then $\sigma_k\cdot [\bu_k]_+  $ can be initialized as the $k$-th column of $\bW$, and $[\bv_k]_+^\top$ can be chosen as the $k$-th row of $\bZ$.
\end{itemize}
It should be noted, however, that none of these initialization techniques are theoretically guaranteed to yield a better final solution---they are heuristic improvements aimed at practical performance. For more details, we recommend consulting the original papers cited above.

\index{Implicit hierarchy}
\section{Movie Recommender Context}
Both NMF and ALS approximate a data matrix by reconstructing its entries from a set of basis (or template) vectors. The key difference lies in the nature of these bases and how the reconstruction is performed.
In NMF, all basis vectors have nonnegative entries, and each data vector is expressed as a nonnegative linear combination of these bases, typically with relatively small coefficients along each direction.
In contrast, ALS allows basis vectors to contain both positive and negative values, and the reconstruction uses a general linear combination, where components can be large in magnitude and of either sign. This means that basis vectors can effectively be added or subtracted during reconstruction.
Consequently, depending on the application, one approach may offer a more meaningful interpretation than the other.

\paragrapharrow{Movie recommender context.}
In a movie recommendation system, the rows of $\bW$  represent latent features of movies (e.g., genres or themes), while the columns of $\bZ$ represent user preferences for those features.
For example, under NMF, a movie might be represented as: 0.5 comedy, 0.002 action, and 0.09 romantic, indicating purely additive contributions. 
In contrast, ALS might assign weights such as: 4 comedy, $-0.05$ action, and $-3$ drama, where negative values indicate that the presence of certain features reduces the relevance or rating for a user. While mathematically valid, such interpretations can be less intuitive in contexts where only positive contributions are meaningful.

\paragrapharrow{Implicit hierarchy.}
Unlike SVD, neither ALS nor NMF imposes an inherent ranking among the basis vectors. In SVD, the importance of each component is explicitly ordered by the magnitude of its corresponding singular value: $\bA=\sum_{i=1}^R\sigma_i\bu_i\bv_i^\top$.
This creates an implicit hierarchy: the first term $\sigma_1\bu_1\bv_1^\top$  captures the dominant pattern in the data, the second term refines it, and so on.
In contrast, the factors in ALS and NMF are not ordered by importance---each plays an equally weighted role unless additional constraints or post-processing are applied. Thus, SVD provides a natural notion of component significance that is absent in ALS and NMF.

\paragrapharrow{Interpretability of basis vectors.}
The basis vectors in SVD correspond to directions of maximum variance in the data and are statistically well-founded. However, due to the presence of zero, positive, and negative entries, they often lack clear semantic or visual interpretability---especially for nonnegative data such as pixel intensities in images or user ratings in recommender systems.
When reconstructing data via SVD, the combination of basis vectors involves intricate cancellations between positive and negative components, which can obscure the physical meaning of individual patterns.
Moreover, there is a fundamental tension between orthogonality and nonnegativity:
A meaningful ``pattern" in nonnegative data should itself be nonnegative.
Yet, mutually orthogonal vectors (as required in SVD) cannot all be nonnegative unless they are trivially sparse.
For instance, suppose the leading left singular vector $\bu_1$ has all nonnegative entries. Then any other left singular vector $\bu_j$ ($j\neq 1$) must satisfy $\bu_1^\top\bu_j=0$. This orthogonality condition forces $\bu_j$ to contain at least one negative entry; otherwise, the inner product would be strictly positive.
Therefore, while SVD is powerful for compression and denoising, its basis vectors are generally unsuitable as interpretable ``parts-based" representations for nonnegative data---a key advantage of NMF.

\section{Other Applications}

\begin{figure}[h]
\centering       
\subfigtopskip=2pt               
\subfigbottomskip=-2pt         
\subfigcapskip=-30pt      
\includegraphics[width=0.98\textwidth]{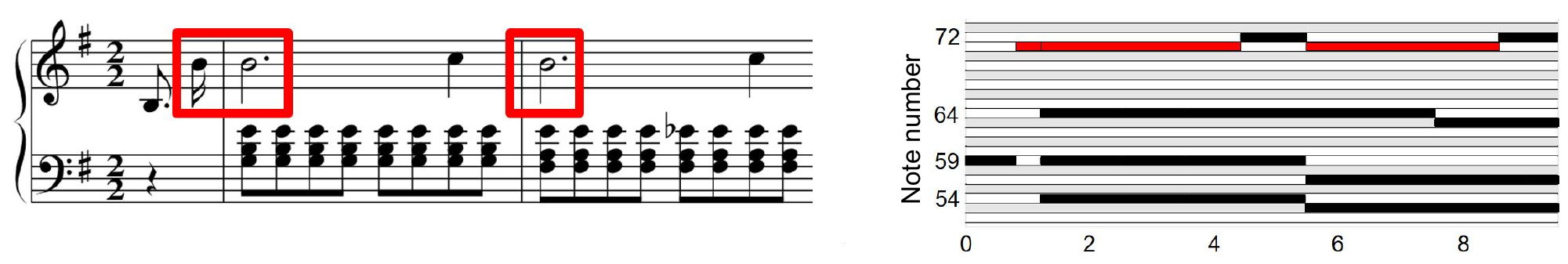}
\caption{Musical score and piano-roll representation. Figure adapted from \citet{muller2015fundamentals}.}
\label{fig:nmf_music_note}
\end{figure}

\paragrapharrow{Music spectral reconstruction.}
To illustrate the application of nonnegative matrix factorization (NMF), we demonstrate how it can decompose the spectrogram of a music recording into components that carry musical meaning \citep{muller2015fundamentals}. As an example, consider the opening measures of \textit{Frédéric Chopin's Prélude Op. 28, No. 4}. Figure~\ref{fig:nmf_music_note}  shows the musical score alongside a synchronized piano-roll visualization of an audio recording of the piece. 
For clarity, all elements related to the note with  pitch number $p=71$ are highlighted with red rectangular frames.

\begin{figure}[h]
\centering       
\vspace{-0.35cm}                 
\subfigtopskip=2pt               
\subfigbottomskip=-2pt         
\subfigcapskip=-10pt      
\includegraphics[width=0.98\textwidth]{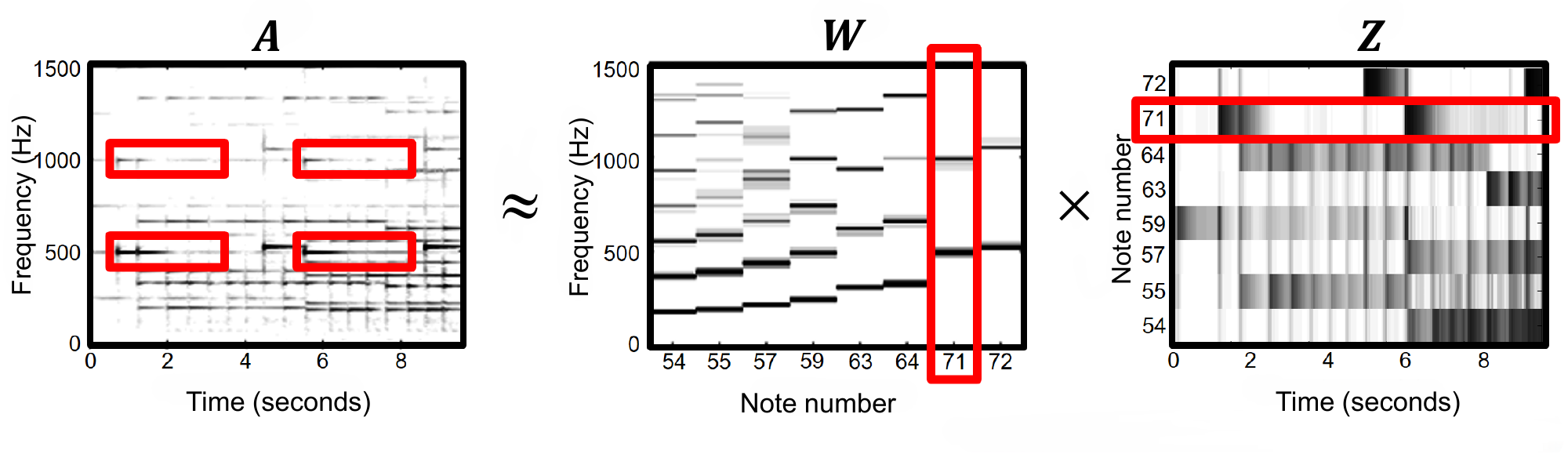}
\caption{Ideal NMF decomposition of the spectrogram. Figure adapted from \citet{muller2015fundamentals}.}
\label{fig:nmf_music_decom}
\end{figure}

The original data matrix $\bA$ is constructed from the magnitude STFT of the audio signal---a sequence of spectral vectors representing frequency content over time \citep{lopez2019nmf}.
Applying NMF factorizes $\bA$ into two nonnegative matrices $\bW$ and $\bZ$. Ideally, the columns of $\bW$ capture the spectral patterns (i.e., timbres) associated with the pitches present in the piece, while $\bZ$ encodes the temporal activation of these patterns---indicating when each note occurs in the recording. Figure~\ref{fig:nmf_music_decom} illustrates this idealized decomposition for the Chopin prelude.

In this setting, each column of $\bW$ corresponds to the spectral signature of a specific pitch, and $\bZ$ closely resembles the piano-roll representation of the musical score. This highlights two key advantages of NMF over general (unconstrained) matrix factorization:
\begin{itemize}
\item \textbf{Nonnegativity constraint.} 
NMF enforces nonnegativity on both $\bW$ and $\bZ$. This aligns naturally with physical quantities in many domains---such as energy, intensity, or amplitude---which cannot be negative. In music analysis, this ensures that the learned components correspond to meaningful, additive sound sources like individual notes or chords.

\item \textbf{Interpretability.} The factor $\bW$ represents spectral templates (timbral profiles) of musical notes, while $\bZ$ indicates their temporal activations. This yields a parts-based, interpretable decomposition. In contrast, unconstrained methods (e.g., SVD or PCA) often produce factors with mixed signs, making them difficult to interpret in terms of real-world musical events.
\end{itemize}

\begin{figure}[h]
\centering       
\vspace{-0.35cm}                 
\subfigtopskip=2pt               
\subfigbottomskip=-2pt         
\subfigcapskip=-10pt      
\includegraphics[width=0.98\textwidth]{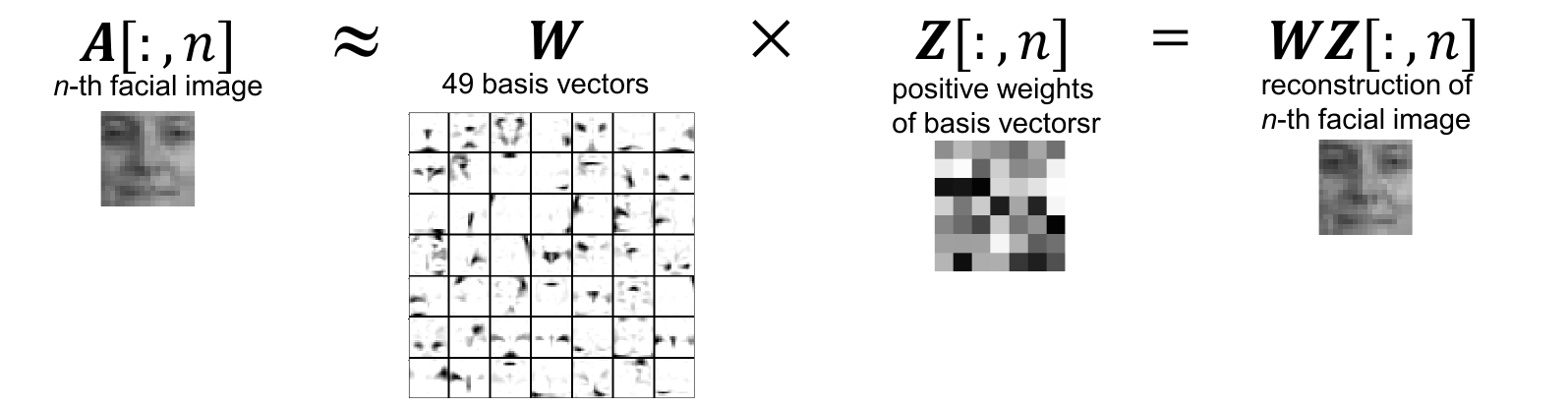}
\caption{NMF applied to the CBCL face database with $K=49$. The basis vectors in $\bW$ are reshaped into $19\times 19$ images. Localized facial features can be observed from these reshaped basis vectors, e.g., eyes, noses, nasolabial folds, and lips. Figure adapted from \citet{lee1999learning, gillis2014and}.}
\label{fig:nmf_face_decom}
\end{figure}

\paragrapharrow{Facial feature extraction and reconstruction.}
Suppose each column of the data matrix $\bA \in \real_+^{M\times N}$ represents a vectorized grayscale image of a human face, where entry $a_{mn}$ denotes the intensity of the $m$-th pixel in the $n$-th image. NMF decomposes $\bA$ into two nonnegative matrices $\bW\in\real_+^{M\times K}$ abd $\bZ\in\real_+^{K\times N}$, such that each face $\ba_n$  is approximated by a nonnegative linear combination of the columns of $\bW$.
Because $\bW$ is nonnegative, its columns can be interpreted as template images---each representing a localized facial feature (e.g., an eye, nose, or lip region). The corresponding weights in $\bZ$ combine these templates additively to reconstruct each original face. Since $K\ll N$ in typical applications, the basis images must capture recurring, sparse, and localized structures shared across the dataset. As shown in Figure~\ref{fig:nmf_face_decom} (using the CBCL face database~\footnote{{http://cbcl.mit.edu/software-datasets/FaceData2.html}}), these basis images often resemble interpretable facial parts such as eyes, noses, nasolabial folds, and lips \citep{lee1999learning, gillis2014and}. Meanwhile, each column of $\bZ$ encodes which features are present---and to what degree---in a given face image.

Moreover, when each column of $\bA$ corresponds to multiple images of the same person, NMF can be used for face recognition. Compared to methods like PCA or ALS---which produce dense, globally distributed basis vectors---NMF's sparse, parts-based representation is more robust to occlusions (e.g., sunglasses, scarves, or distortions). Even if part of a new face is occluded, the non-occluded regions (e.g., mouth or forehead) can still be accurately reconstructed using the relevant basis components \citep{jain2017non}.

\paragrapharrow{Topic recovery.} 
As introduced at the beginning of this chapter, NMF is also highly effective for topic modeling in text analysis. In this context, one constructs a term-document matrix $\bA$, where rows correspond to terms (words or phrases) and columns to documents. Each entry $a_{mn}$ reflects the weight of term $m$ in document $n$---commonly represented as binary indicators, \textit{term frequency (TF)}, or \textit{term frequency-inverse document frequency (TF-IDF)} scores \citep{shahnaz2006document}.

Under NMF, each column of $\bW$ represents a topic, defined as a nonnegative distribution over terms (i.e., a set of co-occurring words with high weights). Each column of $\bZ$ gives the topic proportions for a document---indicating how much each topic contributes to it. This formulation naturally supports soft clustering, where a document can belong to multiple topics simultaneously.

NMF is particularly well-suited for topic recovery because it respects the additive nature of textual content: documents are formed by combining topics, not by subtracting or canceling them. The resulting topics are often highly interpretable, consisting of semantically coherent word groups. However, the quality of the decomposition depends critically on:
\begin{itemize}
\item the choice of the number of topics $K$,
\item the initialization of $\bW$ and $\bZ$, and
\item appropriate preprocessing and scaling of the input matrix (e.g., normalization, TF-IDF weighting).
\end{itemize}
Careful tuning of these aspects is essential for obtaining meaningful and stable results.

\index{$L$-strongly smoothness}
\begin{problemset}
	

\item \label{prob:llipschi_hianls} \textbf{$L$-strongly smooth and PGD in Hi-ANLS problems.} A function $f(\bx): \real^N\rightarrow \real$ is said to be $L$-smooth  (i.e., its gradient is $L$-Lipschitz continuous) if it satisfies that $\normtwo{\nabla f(\bx) - \nabla f(\by)} \leq L\normtwo{\bx-\by}$ for all $\bx,\by$. Show that the subproblem~\eqref{equation:llipschi_hianls} in Hi-ANLS is $L$-strongly smooth with constant $L=\normtwo{\bW[:,k]}^2$. Therefore, the subproblem can be solved  via a \textit{projected gradient descent (PGD)} update with a step size $\eta=\frac{1}{L}$, i.e., using gradient descent update with a step size $\eta=\frac{1}{L}$ first and projecting the update onto the nonnegative orthant afterwards \citep{lu2025practical}. 


\item \label{prob:lsmooth_dslemma} \textbf{Descent lemma for $L$-strongly smooth functions.} Let $f:\sS\rightarrow (-\infty, \infty]$ be a function defined over a convex set $\sS\subseteq\real^N$ such that $\normtwo{\nabla f(\bx) - \nabla f(\by)} \leq L\normtwo{\bx-\by}$ for all $\bx$ and $\by$.
Show that 
$
f(\by)\leq f(\bx)+ \nabla f(\bx)^\top (\by-\bx) + \frac{L}{2}\normtwo{\bx-\by}^2.
$
\textit{Hint: Use fundamental theorem of calculus (Theorem~\ref{theorem:linear_approx}).}

\item \label{prob:third_order_nmf} Let $\ba\in\real^M$, $\bz\in\real^K$, and $\bW\in\real^{M\times K}$. Show that all the third-order partial derivatives of  $F(\bz)=\frac{1}{2}\normtwo{\ba-\bW\bz}^2$ vanish.

\item \textbf{MM applied to $L$-strongly smooth functions.} Let $f(\bx): \real^N\rightarrow \real$ be an $L$-smooth function such that $\normtwo{\nabla f(\bx) - \nabla f(\by)} \leq L\normtwo{\bx-\by}$ for all $\bx,\by$. Show that the function $g(\bx, \widetildebx) = f(\widetildebx)+\nabla f(\widetildebx)^\top (\bx-\widetildebx) +\frac{L}{2} \normtwo{\bx-\widetildebx}^2$ is a valid auxiliary function for $f(\bx)$.
Derive the corresponding MM update rule.

\item Derive the gradients and gradient descent updates for the tri-NMF problem in \eqref{equation:tri_nmv}.


\item \label{prob:projpro1} \textbf{Projection property-I.} Let $\sS\subset \real^N$ be a \textbf{convex set}, and let  $\by\in\real^N$ such that $\widetildeby\triangleq\mathcalP_{\sS}(\by)$. Show that for all $\bx\in\sS$, we have $\langle \bx-\widetildeby, \by-\widetildeby \rangle\leq 0$, i.e., the angle between the two vectors is greater than 90\textdegree.

\item \label{prob:projpro2} \textbf{Projection property-II.} Let $\sS\subset \real^N$ be a \textbf{convex set}, and let  $\by\in\real^N$ such that $\widetildeby\triangleq\mathcalP_{\sS}(\by)$. Show that for all $\bx\in\sS$, we have $\normtwo{\widetildeby - \bx} \leq \normtwo{\by-\bx}$ and $\normtwo{\widetildeby - \bx}^2\leq \normtwo{\by-\bx}^2 - \normtwo{\by-\widetildeby}^2$ (the latter is related to the Pythagorean theorem). \textit{Hint: Examine $\normtwo{\by-\bx}^2=\normtwo{(\widetildeby-\bx)-(\widetildeby-\by)}^2$ and Problem~\ref{prob:projpro1}.}

\item \textbf{Linear feasibility projection.} Let $\sS=\{\bx\in\real^N\mid  \bA\bx=\bb\}$, where $\bA$ has full row rank. Show that the projection satisfies $\mathcalP_{\sS}(\bx) = \bx-\bA^\top(\bA\bA^\top)^{-1}(\bA\bx-\bb)$.

\item \label{prob:ortho_nmf} \textbf{Orthogonal and projective NMF, and clustering.} Consider the same setting as the  orthogonal or projective matrix factorization in Problem~\ref{prob:ortho_mf},  and further assume that $\bA,\bW$, and $\bZ$ are nonnegative. 
Show that there is only one positive entry in each column of $\bZ$ in this case. 
How  is this related to the $K$-means problem discussed in Section~\ref{section:regularization-extention-general}?
When each column of $\bA$ represents a data point, discuss the interpretation of $z_{kn}$ (the $(k,n)$-th entry of $\bZ$) as the importance of the $k$-th cluster to the $n$-th data point in the projective NMF case; that is, each data point can belong to several clusters.

\item Show that the Poisson loss in \eqref{equation:als_poi_los} is equivalent to minimizing the $\beta$-divergence between $\bA$ and $\bW\bZ$ with $\beta=1$.

\item Show that the Gamma loss in \eqref{equation:als_gama_los} is equivalent to minimizing the $\beta$-divergence between $\bA$ and $\bW\bZ$ with $\beta=0$.


\item \label{prob:ab_diverg} \textbf{AB divergence \citep{amari2000methods}.}
The \textit{$\alpha$-$\beta$ (AB) divergence} between two positive scalars $x,y>0$ is defined as:
\begin{align*}
d_{\alpha,\beta}(x,y)&=\begin{cases}
-\frac{1}{\alpha\beta}(x^\alpha y^\beta-\frac{\alpha}{\alpha+\beta}x^{\alpha+\beta}-\frac{\beta}{\alpha+\beta}y^{\alpha+\beta}),&\alpha,\beta,\alpha+\beta\neq 0;\\
\frac{1}{\alpha^2}(x^\alpha\ln(\frac{x^\alpha}{y^\alpha})-x^\alpha+y^\alpha),&\alpha\neq 0,\beta=0;\\
\frac{1}{\alpha^2}(\ln(\frac{y^\alpha}{x^\alpha})+(\frac{y^\alpha}{x^\alpha})^{-1}-1),&\alpha=-\beta\neq 0;\\
\frac{1}{\beta^2}(y^\beta\ln(\frac{y^\beta}{x^\beta})-y^\beta+x^\beta),&\alpha=0,\beta\neq 0;\\
\frac{1}{2}(\ln (x)-\ln (y))^2,&\alpha=0,\beta=0.
\end{cases}
\end{align*}
When $\alpha+\beta=1$, this is known as the  \textit{$\alpha$-divergence}. 
Discuss under what conditions the AB divergence reduces to the $\beta$-divergence.
Furthermore, show that $d_{\alpha,\beta}(x,y)\geq 0$ for all $x,y>0$, with equality if and only if $x=y$.



\item \label{prob:nnga_algebra} \textbf{Nonnegative algebra.} 
Many useful properties follow from nonnegativity. We investigate several of them in this problem. 
For square matrices $\bA,\bB,\bC,\bD\in\real^{N\times N}$, show that 
\begin{itemize}
\item \textbf{Triangle inequality.} $\abs{\bA\bB}\leq \abs{\bA}\abs{\bB}$, where $\abs{\cdot}$ denotes the nonnegative part of the matrix.
\item \textbf{Nonexpansiveness.} $\abs{\bA^k}\leq \abs{\bA}^k$, for all $k=\{1,2,\ldots\}$.
\item \textbf{Equal norm.} $\normf{\bA}=\normf{\abs{\bA}}$.
\item If $\abs{\bB}\geq \abs{\bA}$, then $\normf{\bB}\geq \normf{\bA}$.
\item If $ \bB \geq \bA\geq \bzero$ and $ \bD \geq \bC\geq \bzero$, then $ \bB\bD \geq \bA\bC\geq \bzero$.
\item If $ \bB \geq \bA\geq \bzero$, then $\bB^k\geq \bA^k\geq \bzero$, for all $k=\{1,2,\ldots\}$.
\end{itemize}
For rectangular matrices $\bA,\bB\in\real^{M\times N}$, show that 
\begin{itemize}
\item $\abs{\bA+\bB}\leq \abs{\bA}+\abs{\bB}$.
\end{itemize}

\item$^\ast$ \label{prob:nnga_algebra2} \textbf{Eigenvalue interlacing in nonnegative matrices.}  Let   $\bB - \abs{\bA}\in\real_+^{N\times N}$. Show  that 
$$
\rho(\bA) \leq \rho(\abs{\bA}) \leq \rho(\bB),
$$
where $\rho(\bX)$ denotes the spectral radius of $\bX$ (Definition~\ref{definition:spectrum}).
\textit{Hint: Use Problem~\ref{prob:nnga_algebra}, and show that $\normf{\bA^k}\leq \normf{\abs{\bA}^k} \leq \normf{\bB^k}$.}

\item \label{prob:nnga_algebra3} Use Problem~\ref{prob:nnga_algebra2} to show that $\rho(\bB)\geq \rho(\bA)$ if $\bB\geq \bA\geq \bzero$.

\item Let $\bA\in\real_+^{N\times N}$, $\bB=\bA[1:k,1:k], \,\forall k\in\{1,2,\ldots,N\}$  (any leading principal submatrix of $\bA$), and $\bC\in\real^{k\times k}, \,\forall k\in\{1,2,\ldots,N\}$ be any principal submatrix. Show that 
\begin{itemize}
\item $\rho(\scriptsize\begin{bmatrix}
\bB & \bzero \\
\bzero & \bzero 
\end{bmatrix}
\normalsize
) 
\leq 
\rho (\bA)
$
$\implies\rho(\bB)\leq \rho(\bA)$.
\item Use the first result to prove $\rho(\bC)\leq \rho(\bA)$. \textit{Hint: Apply  permutation transformations.}
\item $\mathopmax{n=1,2,\ldots,N}a_{nn} \leq \rho(\bA)$.
\end{itemize}

\item$^\ast$  \label{prob:nonn_lin_12} Let $\bA\in\real_+^{N\times N}$. 
Show that the spectral radius satisfies the following bounds:
$$
\begin{aligned}
\text{Row sum: }\gap  \mathop{\min}_{1\leq i \leq N} \sum_{j=1}^N a_{ij} 
&\leq \rho(\bA)
\leq 
\mathop{\max}_{1\leq i \leq N} \sum_{j=1}^N a_{ij}; \\
\text{Column sum: }\gap \mathop{\min}_{1\leq j \leq N} \sum_{i=1}^N a_{ij} 
&\leq \rho(\bA)
\leq 
\mathop{\max}_{1\leq j \leq N} \sum_{i=1}^N a_{ij}. \\
\end{aligned}
$$

\end{problemset}

\part{Bayesian Matrix Decomposition}
\chapter{Principal Component Analysis (PCA)}
\begingroup
\hypersetup{
	linkcolor=structurecolor,
	linktoc=page,  
}
\minitoc \newpage
\endgroup

\lettrine{\color{caligraphcolor}P}
Principal component analysis (PCA) is one of the most widely used techniques for dimensionality reduction, data compression, and exploratory data analysis. At its core, however, PCA is fundamentally a {matrix decomposition} method. Given a data matrix  $ \bX \in \real^{N\times D} $  (with  $ N $  observations and  $ D $  features), classical PCA seeks a low-rank approximation by decomposing  $ \bX $  into the product of two lower-dimensional matrices: a score matrix  $ \bW \in \real^{N \times K} $  capturing the coordinates of the data in a reduced subspace, and a loading matrix  $ \bZ \in \real^{D \times K} $  defining the directions (principal components) of maximal variance. This yields the approximation  
$$
\bX \approx \bW\bZ^\top,
$$
which can be derived via the singular value decomposition (SVD) of the centered data matrix. 
Viewed this way, PCA is not merely a statistical tool---it is an elegant example of how structured matrix factorization can reveal the latent geometry of high-dimensional data.

While powerful, classical PCA is purely deterministic and offers no mechanism to quantify uncertainty, incorporate prior knowledge, or handle missing data in a principled way. These limitations motivate a shift from an algebraic perspective to a \textit{probabilistic} one. 
\textit{Probabilistic PCA (PPCA)}, introduced by {\citet{roweis1997algorithms,tipping1999probabilistic}}, reinterprets PCA as a latent variable model (see also the motivating example discussed in Section~\ref{section:lvm}): each observed data point is modeled as a linear transformation of a lower-dimensional latent variable, corrupted by isotropic Gaussian noise. This reformulation embeds PCA within the framework of generative models, enabling likelihood-based inference, model comparison, and seamless extension to incomplete data.

Building on PPCA, \textit{Bayesian PCA} takes the next logical step by placing prior distributions over the model parameters---typically the loading matrix (which describes the relationship between the observed variables and the latent components) and the noise variance. Through Bayesian inference, we obtain full posterior distributions rather than point estimates, naturally quantifying uncertainty in the latent structure. Moreover, hierarchical priors (e.g., \textit{automatic relevance determination} or \textit{ARD}) allow the model to infer the effective dimensionality of the latent space, effectively performing \textit{automatic complexity control}. In this chapter, we will explore this progression---from the geometric intuition of classical PCA, through the generative perspective  of probabilistic PCA, to the full inferential machinery of Bayesian PCA---highlighting how each step enriches matrix decomposition with the expressive power of probability theory.

\index{Automatic relevance determination}
\index{PCA}
\index{Principal component analysis}
\index{Unbiased estimator}
\index{Consistent estimator}
\section{Principal Component Analysis}\label{section:pca_isvd}
Principal component analysis (PCA) is frequently employed  to identify patterns in data and to uncover its underlying variance-covariance structure. In doing so, PCA serves two main purposes:
\begin{enumerate}
\item \textit{Data reduction}. The dimensionality of the data is reduced by representing it with a smaller number of \textit{principal components}.
\item \textit{Interpretation}.  PCA can reveal previously unsuspected relationships among variables or observations.
\end{enumerate}

Dimensionality reduction is also advantageous in applications that require lower-dimensional representations---such as data visualization, efficient storage, and computationally intensive tasks.
Consider a data set consisting of $N$ observations $\mathcalX=\{\bx_1,\bx_2,\ldots,\bx_N\}$, where each $\bx_n\in \real^D$ for  $n= 1,2,\ldots,N$. 
Our goal is to project this data into a lower-dimensional space of dimension $K<D$. 
We begin by defining the sample mean vector and the sample covariance matrix:
$$
\widebarbx \triangleq \frac{1}{N}\sum_{n=1}^{N}\bx_n
\qquad 
\text{and}
\qquad 
\bS \triangleq \frac{1}{N}\sum_{n=1}^{N} (\bx_n - \widebarbx)(\bx_n-\widebarbx)^\top.
$$
Here, the divisor $N$ ensures that $\bS$ is a  \textit{consistent estimator} of the true covariance matrix.~\footnote{Consistency: An estimator $\widehattheta_N $ of a parameter $\theta$, based on a sample of size $N$, is said to be consistent if $\widehattheta_N\stackrel{p}{\rightarrow} \theta$   as $N \rightarrow \infty $.}
Alternatively, one may define the covariance matrix using $N-1$ in the denominator: $\bS \triangleq \frac{1}{\textcolor{mylightbluetext}{N-1}}\sum_{n=1}^{N} (\bx_n - \widebarbx)(\bx_n-\widebarbx)^\top$, which yields an \textit{unbiased consistent estimator} of the covariance matrix \citep{lu2021rigorous}.

Each data point $\bx_n$ is then projected onto a scalar value using a vector $\bu_1$ (see discussion below) such that the projection is given by  $\bu_1^\top\bx_n$. The mean of the projected data is  $\Exp[\bu_1^\top\bx_n] = \bu_1^\top \widebarbx$, and its variance is given by
$$
\begin{aligned}
\Cov[\bu_1^\top\bx_n] &= \frac{1}{N} \sum_{n=1}^{N}( \bu_1^\top \bx_n - \bu_1^\top\widebarbx)^2=
\frac{1}{N} \sum_{n=1}^{N}\bu_1^\top  ( \bx_n -\widebarbx)( \bx_n -\widebarbx)^\top\bu_1
=\bu_1^\top\bS\bu_1.
\end{aligned}
$$

\begin{figure}[h]
\centering   
\vspace{-0.35cm}  
\subfigtopskip=2pt  
\subfigbottomskip=2pt 
\subfigcapskip=-5pt 
\subfigure[Project onto y-axis.]{\label{fig:pca_cluster1}
\includegraphics[width=0.31\linewidth]{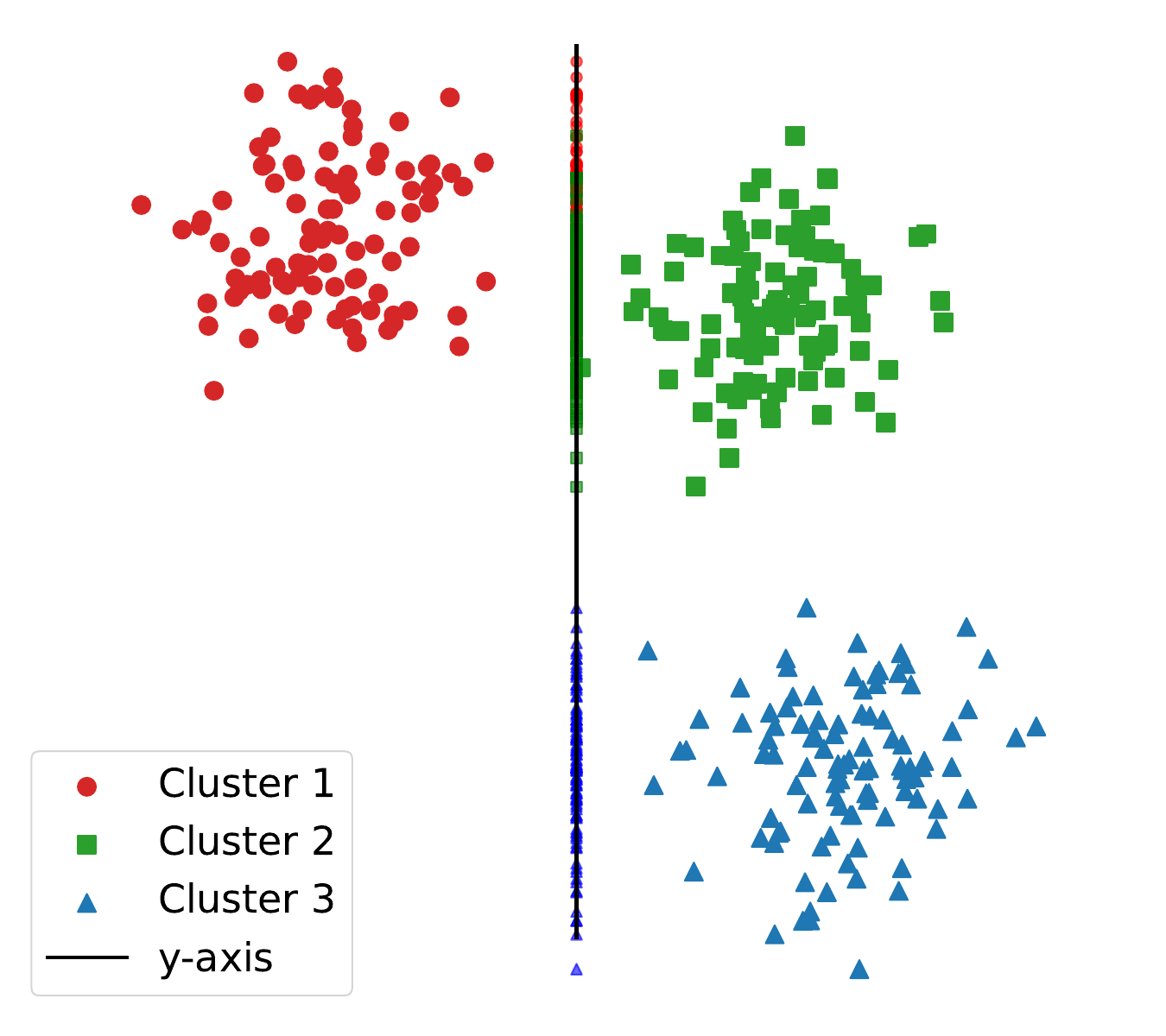}}
\subfigure[Project onto x-axis.]{\label{fig:pca_cluster2}
\includegraphics[width=0.31\linewidth]{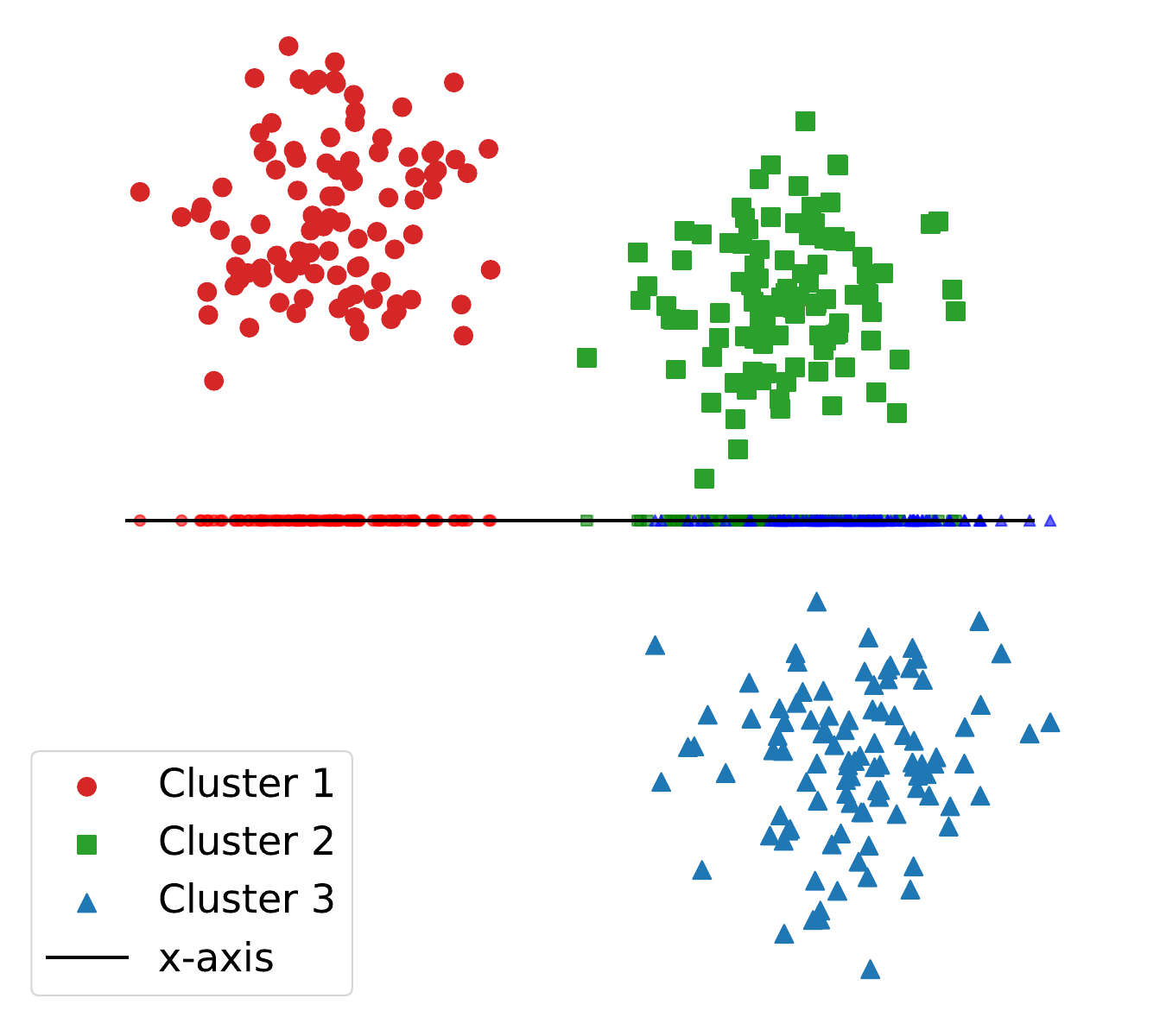}}
\subfigure[Project onto the first principal axis.]{\label{fig:pca_cluster3}
\includegraphics[width=0.31\linewidth]{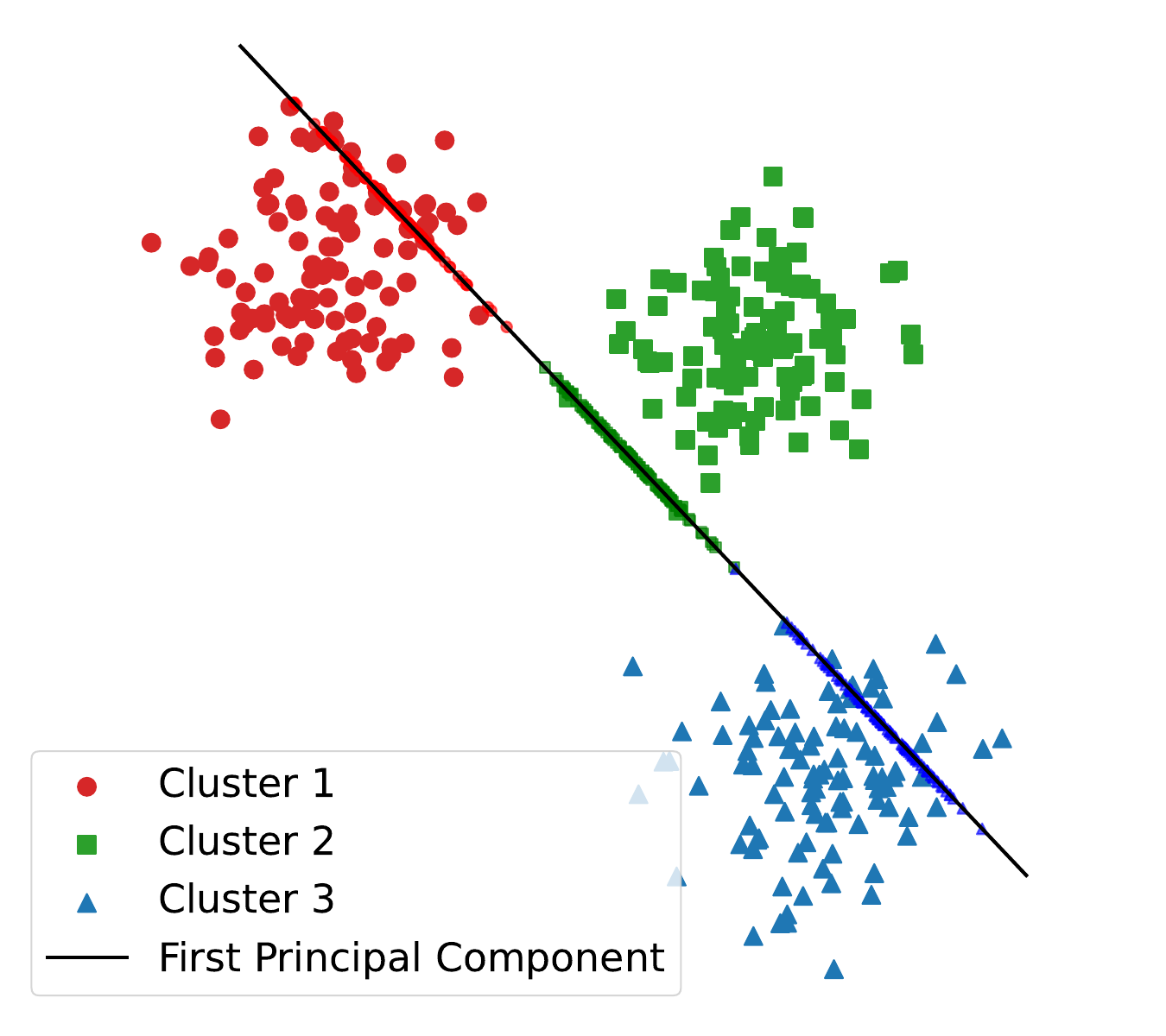}}
\caption{Dimension reduction of a two-dimensional data set that contains three clusters can lead to significant information loss when projecting onto either the x-axis or the y-axis. In contrast, projecting the data onto the first principal axis---the direction of maximal projected variance---preserves much of the cluster structure.}
\label{fig:pca_clusters}
\end{figure}

\subsection{Different Perspectives on PCA}
\paragrapharrow{Maximum-variance formulation.}
The objective of PCA is to find the direction $\bu_1$ that maximizes the projected variance $\bu_1^\top\bS\bu_1$, thereby retaining as much information as possible in the reduced representation (see visual description in Figure~\ref{fig:pca_clusters}). 
To avoid unbounded solutions, we constrain $\normtwo{\bu_1}$ to be a unit vector: $\bu_1^\top\bu_1=1$. 
Using the method of Lagrange multipliers (see, for example,  \citet{bishop2006pattern, boyd2004convex, bishop2023deep}), we maximize the following objective:
\begin{equation}\label{equation:pca_raw1}
\bu_1^\top\bS\bu_1 + \lambda_1 (1 - \bu_1^\top\bu_1).
\end{equation}
Taking the derivative with respect to $\bu_1$ and setting it to zero yields
$$
\bS\bu_1 = \lambda_1\bu_1 
\quad\implies \quad 
\bu_1^\top\bS\bu_1 = \lambda_1,
$$
which shows that $\bu_1$ is an eigenvector of $\bS$ corresponding to eigenvalue $\lambda_1$. Moreover, the projected variance equals $\lambda_1$. Thus, the direction that maximizes variance corresponds to the eigenvector associated with the \textit{largest} eigenvalue of $\bS$. This eigenvector is known as the \textit{first principal axis}.

Using the spectral decomposition (Theorem~\ref{theorem:spectral_theorem}), subsequent principal axes are obtained by selecting eigenvectors corresponding to the next largest eigenvalues, continuing until we have $K<D$ such components. This procedure achieves the desired dimensionality reduction and is referred to as the \textit{maximum-variance formulation} of PCA \citep{hotelling1933analysis, bishop2006pattern, shlens2014tutorial}.

Finally, note that the PCA framework remains valid even when  $K=D$. In this case, no dimensionality reduction occurs; instead, the data is simply rotated into a new coordinate system aligned with the principal components.

\paragrapharrow{Minimum-error formulation.}
An alternative perspective on PCA, known as the \textit{minimum-error formulation}, is discussed in \citet{pearson1901liii, bishop2006pattern}. We now review this approach.
Let  $\bU\in\real^{D\times D}$ be an orthogonal matrix whose columns $\{\bu_1, \bu_2, \ldots, \bu_D\}$ form an orthonormal basis for  $\real^D$.
Since this basis spans the entire space $\real^D$, each data point  $\bx_n$ ($n=1,2,\ldots,N$) can be expressed exactly as   a linear combination of these basis vectors:
\begin{equation}\label{equation:pca_genexp}
\bx_n = \sum_{\ell=1}^{D} \gamma_{n\ell} \bu_\ell
\equiv  \sum_{\ell=1}^{D} (\bx_n^\top \bu_\ell)  \bu_\ell,
\quad n=1,2,\ldots,N, 
\end{equation}
where the coefficients $\gamma_{n\ell} = \bx_n^\top \bu_\ell$ follow from the orthonormality of the basis.
This transformation amounts to a rotation of the coordinate system: the original coordinates   $\{x_{n1},x_{n2}, \ldots, x_{np}\}$ are replaced by new coordinates $\{\gamma_{n1},\gamma_{n2}, \ldots, \gamma_{np}\}$ in the rotated basis $\{\bu_\ell\}$.

However, our aim is not exact reconstruction but approximation using only $K<D$  dimensions---i.e., by projecting the data onto a $K$-dimensional  linear subspace. 
Without loss of generality, assume this subspace is spanned by the first $K$ basis vectors $\{\bu_1, \bu_2, \ldots,\bu_K\}$.
We then approximate each $\bx_n$ as
\begin{equation}
\widetildebx_n = \sum_{\ell=1}^{K} a_{n\ell} \bu_\ell + \sum_{\ell=K+1}^{D} b_\ell \bu_\ell 
\end{equation}
where the coefficients $\{a_{n\ell}\}$ depend on the specific data point $\bx_n$, while the offsets $\{b_\ell\}$ are shared across all  points. 
We are free to choose the basis $\{\bu_\ell\}$, the point-specific coefficients $\{a_{n\ell}\}$, and the global offsets $\{b_\ell\}$ so as to minimize the reconstruction error. 
As our error measure, we use the average squared Euclidean distance between the original and reconstructed points:
\begin{equation}
F = \frac{1}{N} \sum_{n=1}^{N} \normtwo{\bx_n - \widetildebx_n}^2. 
\end{equation}

We first minimize $F$ with respect to $\{a_{n\ell}\}$ and $\{b_\ell\}$.
Substituting the expression for $\widetildebx_n$, differentiating $F$ with respect to $a_{n\ell}$ or $b_\ell$, and using the orthonormality of $\{\bu_\ell\}$, we obtain the optimal values:
\begin{align*}
a_{n\ell} &= \bx_n^\top \bu_\ell,
\quad \ell = 1, 2,\ldots, K;\\
b_\ell &= \widebarbx^\top \bu_\ell, 
\quad \ell = K+1, \ldots, D.
\end{align*}
Substituting these back into the error expression and using the expansion in \eqref{equation:pca_genexp}, we find:
\begin{equation}
\bx_n - \widetildebx_n = \sum_{\ell=K+1}^{D} \left\{ (\bx_n - \widebarbx)^\top \bu_\ell \right\} \bu_\ell.
\end{equation}
This shows that the reconstruction error lies entirely in the subspace orthogonal to the chosen  $K$-dimensional principal subspace---i.e., it is a linear combination of $\{\bu_{K+1}, \ldots,\bu_D\}$.
This is intuitive: the best approximation within the subspace is the orthogonal projection of  $\bx_n$ onto it.

Consequently, the total error depends only on the choice of basis vectors and simplifies to:
\begin{equation}
F = \frac{1}{N} \sum_{n=1}^{N} \sum_{\ell=K+1}^{D} \left( \bx_n^\top \bu_\ell - \widebarbx^\top \bu_\ell \right)^2 = \sum_{\ell=K+1}^{D} \bu_\ell^\top \bS \bu_\ell.
\end{equation}

To minimize $F$, we must choose an orthonormal set $\{\bu_\ell\}$.
Without constraints, the trivial solution $\bu_\ell = \bzero$ would minimize $F$; hence, orthonormality is essential. 
The solution emerges naturally from the spectral decomposition of $\bS $.

To build intuition, consider the case $D = 2$ and  $K = 1$. 
We must choose a unit vector $\bu_2$ (orthogonal to the principal subspace) to minimize $F = \bu_2^\top \bS \bu_2$, subject to the normalization constraint $\bu_2^\top \bu_2 = 1$. 
Introducing a Lagrange multiplier $\lambda_2$, we minimize:
$$
F_{\lambda_2} = \bu_2^\top \bS \bu_2 + \lambda_2 \left(1 - \bu_2^\top \bu_2 \right).
$$
Setting the derivative to zero yields $\bS \bu_2 = \lambda_2 \bu_2$, so $\bu_2$ is an eigenvector of $\bS$, with eigenvalue $\lambda_2$, and  $F = \lambda_2$.
To minimize $F$, we select $\bu_2$ as the eigenvector corresponding to the {smaller} eigenvalue.
Consequently, the principal direction $\bu_1$ aligns with the eigenvector of the {larger} eigenvalue---precisely matching the maximum-variance criterion.
If the two eigenvalues are equal, all directions are equivalent, and any choice of $\bu_1$ yields the same reconstruction error.

This reasoning extends to the general case. For arbitrary $N$ and  $K<D$, the minimum of $F$ is achieved when $\{\bu_\ell\}$ are the orthonormal eigenvectors of $\bS$:
\begin{equation*}
\bS \bu_\ell = \lambda_\ell \bu_\ell , 
\quad \ell=1,2,\ldots,D,
\end{equation*}
ordered such that $\lambda_1\geq\lambda_2\geq \ldots\geq \lambda_D$.
The reconstruction error then becomes
\begin{equation*}
F = \sum_{\ell=K+1}^{D} \lambda_\ell, 
\end{equation*}
the sum of the $D-K$ smallest eigenvalues. 
Thus, to minimize reconstruction error, we retain the $K$ eigenvectors with the largest eigenvalues---exactly the same solution as in the maximum-variance formulation.

\paragrapharrow{Optimization perspectives.}
Alternatively, assume the data have been centered so that the sample mean vector $\widebarbx$ is zero (i.e., the data are mean-centered). 
If they are not already centered, we can achieve this by replacing each observation with $\widebarbx_n \leftarrow \bx_n-\widebarbx$ thereby subtracting the sample mean from every data point.
Our goal is to project the centered data points $\{\widebarbx_1, \widebarbx_2, \ldots, \widebarbx_N\}$ from  $\real^D$ into a lower-dimensional subspace $\real^K$, where $K<D$. 
Let  $\bP\in\real^{N\times K}$ be a \textit{semi-orthogonal matrix} satisfying $\bP^\top\bP=\bI_K$. 
This means the columns of $\bP$ form an orthonormal basis for a $K$-dimensional linear subspace $\mathcalV\subset \real^D$.
Then the matrix $\bH=\bP\bP^\top$ defines an \textit{orthogonal projection} (i.e., a symmetric and idempotent matrix) onto the low-dimensional subspace defined by the column space $\mathcalV$ of $\bP$ (see Problem~\ref{prob:projection-from-matrix}). 
The orthogonal projection of any centered data point $\widebarbx_n$ onto the
subspace $\mathcalV$ is 
\begin{equation}
\mathcalP_{\bP}(\widebarbx_n) =\bP\bP^\top\widebarbx_n. 
\end{equation}
PCA  seeks the projection matrix $\bP$ that   maximizes the variance of the projected data.
It can be shown that the covariance matrix of the projected data is
\begin{equation}
\frac{1}{N} \sum_{n=1}^{N}\bP\bP^\top \widebarbx_n  (\bP\bP^\top \widebarbx_n  )^\top 
=\frac{1}{N} \bP\bP^\top \bX_c^\top\bX_c \bP\bP^\top,
\end{equation}
where $\bX_c\in\real^{N\times D}$ is the centered data matrix, with each row containing one centered observation:
$$
\bX_c \triangleq 
\begin{bmatrix}
\widebarbx_1^\top\\
\widebarbx_2^\top\\
\vdots\\
\widebarbx_N^\top\\
\end{bmatrix}
\equiv 
\bX - \bone \widebarbx^\top.
$$
Here, $\bone\in\real^N$ is a column vector of ones, and $\widebarbx\in\real^D$ is the sample mean vector of the original data matrix $\bX\in\real^{N\times D}$.
Since the total variance of the projected data equals the trace of its covariance matrix, PCA can be formulated as the following optimization problem:
\begin{equation}\label{equation:pca_second}
\mathop{\max}_{\bP}\,\, \trace(\bP\bP^\top \bX_c^\top\bX_c \bP\bP^\top) \gap \text{s.t.} \gap \bP^\top\bP=\bI_K.
\end{equation}
Using the cyclic property of the trace, this objective simplifies to
$$
\trace(\bP\bP^\top \bX_c^\top\bX_c \bP\bP^\top)
=\trace(\bP^\top \bX_c^\top\bX_c \bP).
$$
It can then be shown that the optimal $\bP$ consists of the eigenvectors of $\bX_c^\top\bX_c$ corresponding to its $K$ largest eigenvalues.

\begin{SCfigure}
\centering
\includegraphics[width=0.4\textwidth]{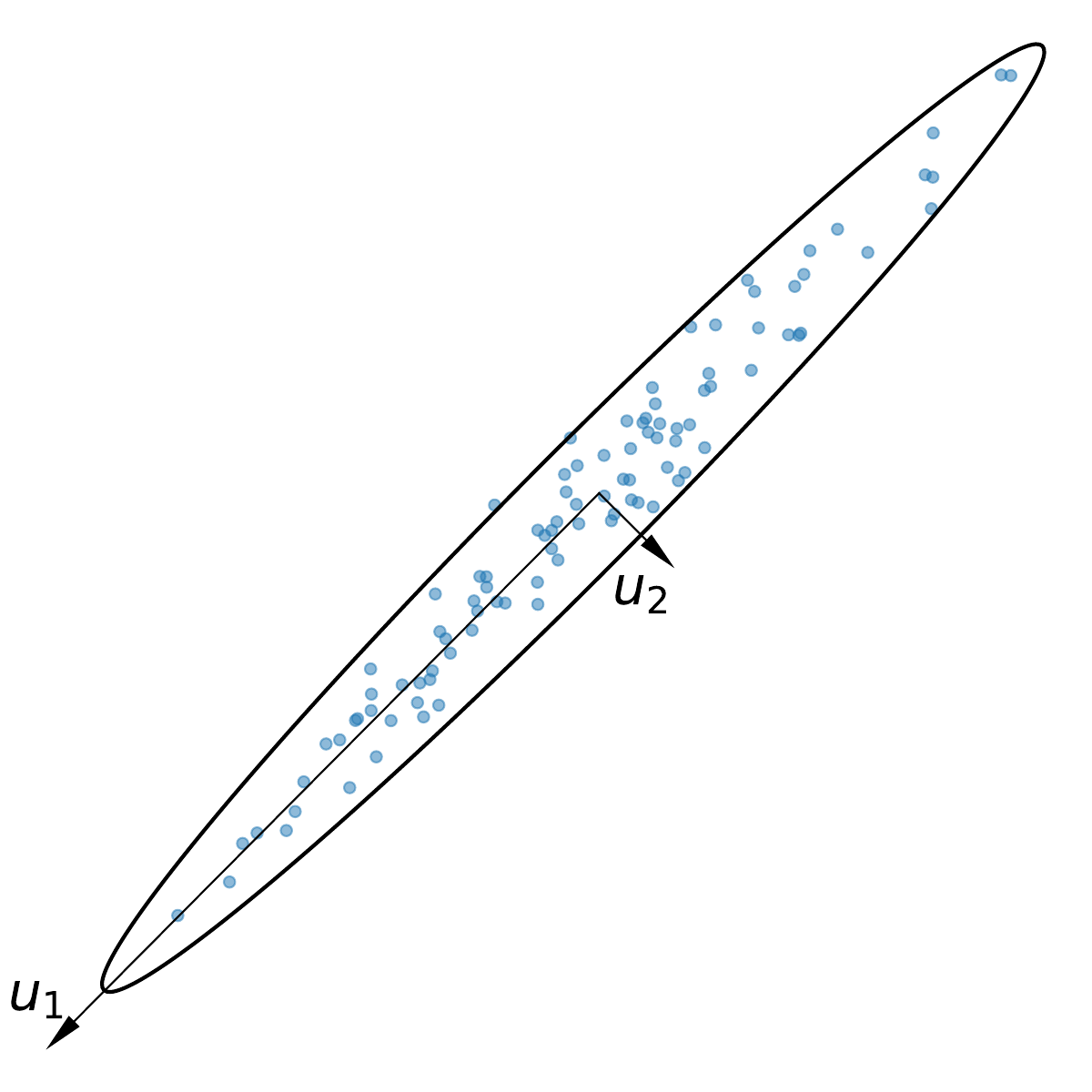}
\caption{Description of PCA in a two-dimensional case. $\bu_1$ and $\bu_2$ are the directions of corresponding eigenvectors of the covariance matrix. Therefore, $\bu_1$ encodes the first principal axis, and $\bu_2$ is the second principal axis.}
\label{fig:pca_axis_2d}
\end{SCfigure}

As noted above, the projection of $\widebarbx_n$ onto the subspace $\mathcalV$ (the column space of $\bP$) is $\bP\bP^\top \widebarbx_n$. The total squared reconstruction error---i.e., the sum of squared distances between the original points and their projections---is
$$
\sum_{n=1}^{N} \normtwo{\bP\bP^\top \widebarbx_n-\widebarbx_n}^2
=\normf{\bP\bP^\top \bX_c^\top - \bX_c^\top}^2
=-\trace(\bP^\top\bX_c^\top\bX_c\bP)+\trace(\bX_c\bX_c^\top).
$$
Thus, minimizing the reconstruction error leads to the equivalent optimization problem:
\begin{equation}\label{equation:pca_third}
\mathop{\min}_{\bP}\,\, -\trace(\bP^\top\bX_c^\top\bX_c\bP)+\trace(\bX_c\bX_c^\top)\gap \text{s.t.} \gap \bP^\top\bP=\bI_K.
\end{equation}
Since $\trace(\bX_c\bX_c^\top)$ is constant with respect to $\bP$, this minimization is equivalent to maximizing $\trace(\bP^\top\bX_c^\top\bX_c\bP)$---exactly the objective in \eqref{equation:pca_second}.
Therefore, PCA simultaneously maximizes the variance of the projected data and minimizes the reconstruction error. These two viewpoints are mathematically equivalent.
Figure~\ref{fig:pca_axis_2d} illustrates this idea in two dimensions: the first principal axis (denoted $\bu_1$) captures the direction of maximum variance, while the second principal axis (denoted $\bu_2$) is orthogonal to the first and captures the remaining variance.

\paragrapharrow{PCA via the spectral decomposition.}
Let the data matrix $\bX_c \in \real^{N\times D}$ contain the mean-centered observations as its rows. The sample covariance matrix is then given by
\begin{equation}\label{equation:pca-equ0}
\bS= \frac{1}{N}\bX_c^\top\bX_c,
\end{equation}
which is symmetric and positive semidefinite. Its spectral decomposition is
\begin{equation}\label{equation:pca-equ1}
\bS = \bU\bLambda\bU^\top,
\end{equation}
where $\bU\in\real^{D\times D}$ whose columns $\bu_1, \bu_2,\ldots,\bu_D$ are the eigenvectors of $\bS$, and $\bLambda=\diag(\lambda_1, \lambda_2,\ldots, \lambda_D)$ is a diagonal matrix of eigenvalues, ordered such that $\lambda_1 \geq \lambda_2 \geq \ldots \geq \lambda_D$. 
The eigenvectors $\{\bu_d\}$ are called the \textit{principal axes} of the data. 
They decorrelate the covariance matrix, meaning that when the data are projected onto these axes, the resulting variables are uncorrelated. 
As noted previously,  the projections of the data onto the principal axes are known as the \textit{principal components}. 
Specifically, the $k$-th principal component is the $k$-th column of the matrix $\bX_c\bU$. 
If our goal is to reduce the dimensionality from $D$ to $K<D$, we retain only the first $K$ principal components by selecting the first $K$ columns of $\bX_c\bU$:
$$
\widetildebX\triangleq
\bX_c\bU_K
=
\bX_c[\bu_1, \bu_2, \ldots, \bu_K].
$$
We make the following key observations about PCA:
\begin{itemize}
\item The matrix $\widetildebX$ is also mean-centered. 
Since $\bX_c$ has zero row-wise mean, and $\bU_K$ is a linear transformation, the reduced representation inherits this property: $\bone^\top\widetildebX = \bone^\top\bX_c\bU_K=\bzero^\top$, where $\bone\in\real^N$ is the vector of all ones.

\item The covariance matrix of $\widetildebX$ is diagonal and given by $\bLambda_K\triangleq \diag([\lambda_1, \ldots,\lambda_K])$.
Since the matrix $\widetildebX$ is mean-centered, its covariance matrix can be represented as $\widetildebX^\top\widetildebX/N$, which simplifies to:
\begin{align*}
\frac{\widetildebX^\top\widetildebX}{N} 
&= \bU_K^\top\left[\frac{1}{N}\bX_c^\top\bX_c\right]\bU_K = [\bu_1, \bu_2, \ldots, \bu_K]^\top(\bS[\bu_1,\bu_2, \ldots, \bu_K]) \\
&= [\bu_1,\bu_2, \ldots, \bu_K]^\top[\lambda_1\bu_1, \lambda_2\bu_2, \ldots, \lambda_K\bu_K] = \bLambda_K.
\end{align*}

\item  The total variance retained in the reduced representation is $\sum_{k=1}^{K} \lambda_k$.  Since the total variance in the original data is $\sum_{d=1}^{D} \lambda_d$,  the fraction of explained variance is $(\sum_{k=1}^{K} \lambda_k)/(\sum_{d=1}^{D} \lambda_d)$.
\end{itemize}
To reconstruct an approximation of the original (uncentered) data from $\widetildebX$ and $\bU_K^\top$, we must store the sample mean vector $\widebarbx$ used during centering. 
The reconstruction is then given by
\begin{equation}\label{equation:pca_proces}
\bX \approx \bX_{\text{pca}} = \underbrace{\widetildebX \bU_K^\top}_{\approx\bX_c} + \bone\widebarbx^\top.
\end{equation}
The storage overhead for $\widebarbx$ is negligible---only $D$ additional numbers---and becomes increasingly insignificant as the dataset size $N$ grows.

\index{Truncated}
\index{Truncated SVD}
\paragrapharrow{PCA via the SVD.}
Let the SVD of the centered data matrix be $\bX_c = \bP\bSigma\bQ^\top$, where $\bP\in\real^{N\times N}, \bQ\in\real^{D\times D}$ are orthogonal matrices, and $\bSigma\in\real^{N\times D}$ is a rectangular diagonal matrix with nonnegative singular values $\sigma_1\geq \sigma_2\geq \ldots\geq\sigma_R\geq 0$, on its main diagonal, where $R=\min\{N,D\}$.
Then the covariance matrix can be expressed as
\begin{equation}\label{equation:pca-equ2}
\bS= \frac{1}{N}\bX_c^\top\bX_c = \bQ \frac{\bSigma^2}{N}\bQ^\top,
\end{equation}
Comparing \eqref{equation:pca-equ2} with the spectral decomposition in \eqref{equation:pca-equ1}, we see that:
the right singular vectors $\bQ$ are precisely the eigenvectors of $\bS$ (i.e., principal axes), and the eigenvalues of $\bS$ are related to the singular values by $\lambda_d = \sigma_d^2/N$ for $d=1,2,\ldots,D$.
Thus, to reduce the dimensionality to $K$, we select the top $K$ singular values and their corresponding right singular vectors. This is equivalent to computing the truncated SVD (TSVD):
$$
\bX_K = \sum_{k=1}^{K}\sigma_k\bp_k\bq_k^\top,
$$
where $\bp_k$'s and $\bq_k$'s are the columns of $\bP$ and $\bQ$, respectively (see Problem~\ref{theorem:young-theorem_frob}).

\paragrapharrow{A computational shortcut for high-dimensional data.} 
Consider a principal axis $\bu_i$ of $\bS = \frac{1}{N}\bX_c^\top\bX_c$, we have
$$
\frac{1}{N}\bX_c^\top\bX_c \bu_i = \lambda_i \bu_i.
$$
Premultiplying both sides by $\bX_c$ yields
$$
\frac{1}{N}\bX_c\bX_c^\top (\bX_c\bu_i) = \lambda_i (\bX_c\bu_i),
$$
which shows that $\lambda_i$ is also an eigenvalue of the $N\times N$ matrix $\frac{1}{N}\bX_c\bX_c^\top \in \real^{N\times N}$, with corresponding eigenvector  $\bX_c\bu_i$. 
When the number of features greatly exceeds the number of samples ($D \gg N$), it is computationally more efficient to compute the eigenvectors of the smaller $N\times N$ matrix $\frac{1}{N}\bX_c\bX_c^\top$ rather than the $D\times D$ covariance matrix $\bS=\frac{1}{N}\bX_c^\top\bX_c$. 
This reduces the computational complexity from $\mathcalO(D^3)$ to $\mathcalO(N^3)$---a significant saving when  $D$ is very large.

Specifically, suppose $\bv_i\in\real^N$  is an eigenvector of $\frac{1}{N}\bX_c\bX_c^\top$  corresponding to a nonzero eigenvalue $\lambda_i$: 
$$
\frac{1}{N}\bX_c\bX_c^\top \bv_i = \lambda_i \bv_i.
$$
Premultiplying by $\bX_c^\top$ gives 
$$
\frac{1}{N}\bX_c^\top\bX_c (\bX_c^\top\bv_i) = \bS(\bX_c^\top\bv_i)   = \lambda_i (\bX_c^\top\bv_i).
$$
Hence,  $\bX_c^\top\bv_i$  is an eigenvector of $\frac{1}{N}\bX_c\bX_c^\top$  associated with the same eigenvalue $\lambda_i$. 
To obtain a unit-norm principal axis, we normalize: $\bu_i = \bX_c^\top\bv_i/\normtwo{\bX_c^\top\bv_i}$.
Therefore, when $D\gg N$, the principal axes can be efficiently computed via the spectral decomposition (or SVD) of the much smaller $N\times N$ matrix $\frac{1}{N}\bX_c\bX_c^\top$.

\index{Orthogonal matrix factorization}
\index{Nonnegative PCA}
\subsection{Orthogonal Matrix Factorization and Nonnegative PCA}
Consider the same centered data matrix $\bX_c$ as defined in \eqref{equation:pca-equ0}, and the following orthogonal matrix factorization problem (see Chapters~\ref{chapter:als} and \ref{chapter:nmf} for further details):
\begin{subequations}\label{equation:or_pca_all}
\begin{equation}\label{equation:or_pca1}
\mathopmin{\bW,\bZ} \normf{\bX_c-\bW\bZ}^2, 
\gap 
\text{with}\gap
\bZ\bZ^\top =\bI_K,
\end{equation} 
where $\bX_c\in\real^{N\times D}, \bW\in\real^{N\times K}$, and $\bZ\in\real^{K\times D}$ with $K\leq \min\{N, D\}$. 
For a fixed $\bZ$,  the optimal  $\bW$ is given by $\bW^*=\bX_c\bZ^\top$ (see Problem~\ref{prob:ortho_mf}; this follows from setting the gradient with respect to $\bW$ to zero).
Substituting this back into \eqref{equation:or_pca1} yields an equivalent optimization problem in terms of $\bZ$ alone:
\begin{equation}\label{equation:or_pca2}
\mathopmin{\bZ\bZ^\top =\bI_K} \normf{\bX_c-\bX_c\bZ^\top\bZ}^2
=
\mathopmin{\bZ\bZ^\top =\bI_K}\normf{\bX_c}^2 - \normf{\bX_c\bZ^\top}^2
=
\mathopmax{\bZ\bZ^\top =\bI_K}\normf{\bX_c\bZ^\top}^2.
\end{equation} 
Note that the final expression is a maximization problem with an orthogonality constraint on the rows of $\bZ$. Let $\bz_k$ denote the $k$-th row of $\bZ$. 
Then the objective function can be written as
\begin{equation}\label{equation:or_pca3}
\normf{\bX_c\bZ^\top}^2
=\sum_{k=1}^{K}\normtwo{\bX_c\bz_k}^2 
=\sum_{k=1}^{K}\bz_k^\top\bX_c^\top\bX_c\bz_k.
\end{equation}
\end{subequations}
This is a classic trace maximization problem under orthonormality constraints. By the Rayleigh--Ritz theorem (see Problems~\ref{prob:rayleigh_v2}--\ref{prob:rayleigh_v4}) the solution is obtained when the rows of $\bZ$ are the top-$K$ eigenvectors of $\bX_c^\top\bX_c$.
This result holds even for a data matrix $\bX_c$ that is not mean-centered.

When $K=1$, \eqref{equation:or_pca2} reduces to $\mathopmax{\normtwo{\bz} =1}\bz^\top(\bX_c^\top\bX_c)\bz$, which is identical to the standard PCA formulation in  \eqref{equation:pca_raw1}.
In this context,  \textit{nonnegative PCA} extends the idea by imposing a nonnegativity constraint on the loading vector:
\begin{equation}
	\mathopmax{\bz\in\real_+^D, \normtwo{\bz} =1}\bz^\top(\bX_c^\top\bX_c)\bz.
\end{equation} 
This variant is particularly useful in applications where interpretability requires nonnegative components---for example:
\begin{itemize}
\item identifying co-expressed gene sets in gene expression analysis, or
\item extracting image features that respect the nonnegative nature of pixel intensities \citep{montanari2015non}.
\end{itemize}
Note that if the original data matrix satisfies $\bX_c\in\real_+^{N\times D}$, then the unconstrained and nonnegative PCA problems may yield similar solutions---but they are not generally equivalent unless additional conditions hold.

\index{Data whitening}
\subsection{Data Whitening}

PCA is  widely used in machine learning for feature preprocessing. Beyond dimensionality reduction, PCA can also normalize the transformed features so that each has unit variance. This two-step process is known as \textit{whitening} (or \textit{sphering}).
Let $\bU_K\in\real^{D\times K}$ contain the top-$K$ eigenvectors of the sample covariance matrix $\bS=\frac{1}{N}\bX_c^\top\bX_c$. 
The first step of whitening projects the mean-centered data onto the principal subspace:
\begin{subequations}
\begin{equation}
\widetildebX = \bX_c\bU_K\in\real^{N\times K}.
\end{equation}
The second step rescales each principal component by the inverse square root of its corresponding eigenvalue.
Denote $\bLambda_K=\diag([\lambda_1,\lambda_2,\ldots,\lambda_K])$,  the whitened data matrix is:
\begin{equation}
\bY = \bX_c\bU_K\bLambda_K^{-1/2}.
\end{equation}
\end{subequations}
This transformation renders the data distribution approximately spherical: all directions in the new space have equal variance and are uncorrelated. When $K=D$ (i.e., no dimensionality reduction), the covariance of the whitened data is exactly the identity matrix:
$$
\frac{1}{N}\sum_{n=1}^{N}\by_n\by_n^\top
=\frac{1}{N}\sum_{n=1}^{N} \bLambda^{-1/2}\bU^\top\widebarbx_n\widebarbx_n^\top\bU\bLambda^{-1/2}
=\bLambda^{-1/2}\bU^\top\bS\bU\bLambda^{-1/2}
=\bI.
$$

Whitened data often leads to faster convergence in gradient-based optimization algorithms \citep{lu2025practical}. This is because large differences in feature variances create loss landscapes with highly varying curvature across dimensions, which can slow down or destabilize optimization. By equalizing the scale of all features, whitening reduces ill-conditioning and ensures that no single direction dominates the gradient updates.
Moreover, whitening prevents certain features from exerting disproportionate influence simply due to their scale---a common issue when features are measured in different units.

\begin{figure}[h]
	\centering
	\includegraphics[width=0.9\textwidth]{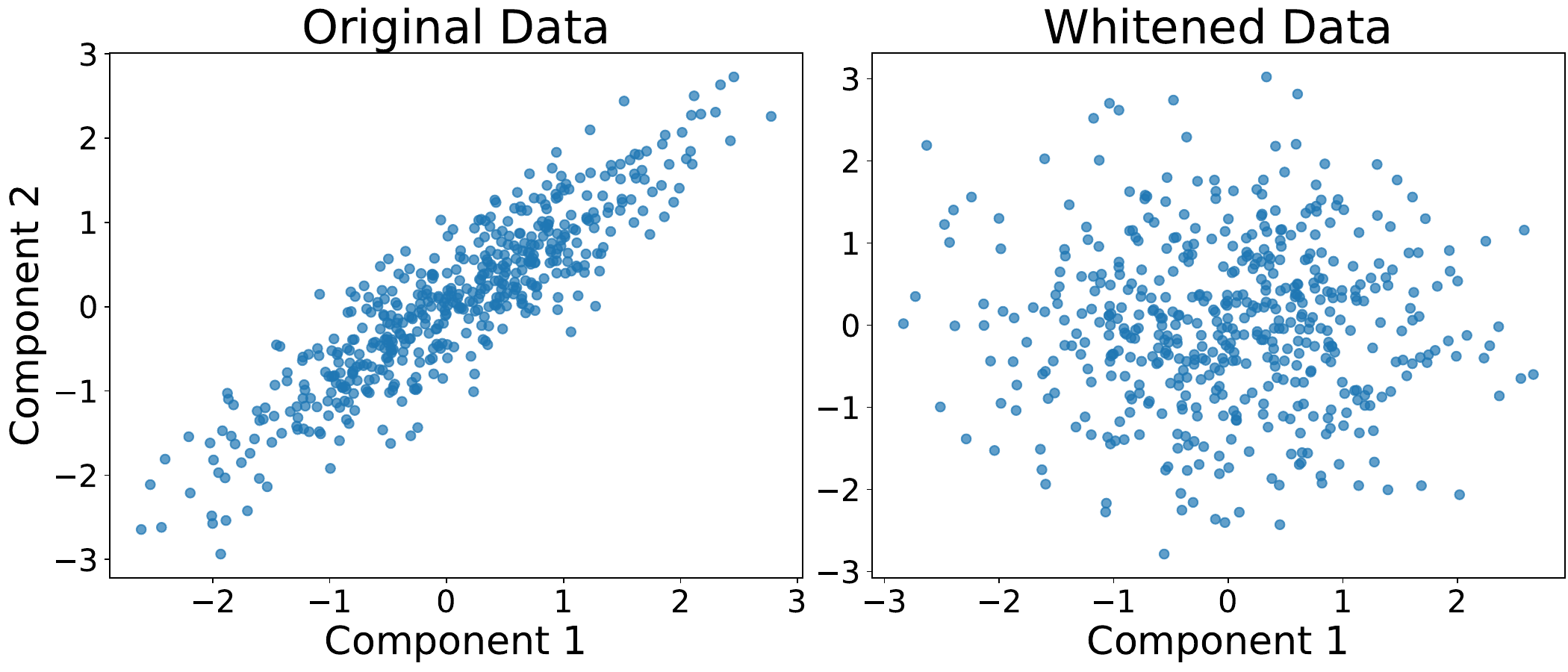} 
	\caption{An example of whitening an ellipsoidal data distribution using principal component analysis.}
	\label{fig:whitening}
\end{figure}

This preprocessing technique is especially valuable in unsupervised learning, such as anomaly or outlier detection, where there are no labels to indicate which directions in the data are important. In such settings, whitening helps ensure that distance-based methods operate on a geometrically balanced representation of the data.
An illustration of this effect is shown in Figure~\ref{fig:whitening}: an initially ellipsoidal data cloud is transformed into a spherical one through PCA-based whitening.

\section{Probabilistic and Bayesian Principal Component Analysis}\label{section:ppca_bpca_all}
In the previous section, we saw that PCA can be interpreted as a linear projection of the data onto a lower-dimensional subspace of the original $\real^D$ space.
The projected data points can be viewed as deterministic latent variables: each observation $\bx_n\in\real^D$ maps to a unique latent representation $\bz_n\in\real^K$. 
To motivate the use of probabilistic continuous latent variables, we now show that PCA can also be derived as the maximum likelihood solution of a probabilistic latent variable model (see also the motivating example for latent variable models (LVMs) in Section~\ref{section:lvm}).
This probabilistic reinterpretation of PCA is known as \textit{probabilistic PCA (PPCA)} \citep{roweis1997algorithms,tipping1999probabilistic}.
Building on this, a Bayesian treatment (BPCA) of the model parameters can also be introduced \citep{bishop1998bayesian}.
These probabilistic and Bayesian reformulations of PCA offer several advantages over standard PCA:
\begin{itemize}
\item Probabilistic/Bayesian PCA defines a constrained Gaussian distribution whose number of free parameters can be controlled while still capturing the dominant correlations in the data.
\item  An EM algorithm can be derived for PPCA that is computationally efficient---particularly when only a few leading principal components are needed---and avoids explicitly computing the full data covariance matrix.
\item {The combination of a probabilistic (or Bayesian) model with the EM algorithm provides a principled way to handle missing data.}
\item {Mixtures of probabilistic/Bayesian PCA models can be formulated and trained in a coherent, principled manner using the EM algorithm.}
\item Because PPCA or BPCA is based on a likelihood function, it enables direct comparison with other probabilistic density models. In contrast, standard PCA assigns low reconstruction error to any point near the principal subspace---even if that point lies far outside the region occupied by the training data---making it unsuitable for density modeling.
\item The models can be run generatively: once trained, they can produce synthetic samples from the learned data distribution.
\item In Bayesian PCA, the effective dimensionality $K$ of the latent subspace can be automatically inferred from the data, eliminating the need to pre-specify it.
\end{itemize}

\subsection{Probabilistic Principal Component Analysis}\label{section:ppca}
PPCA is a simple instance of the linear Gaussian framework, in which all marginal and conditional distributions are Gaussian. 
We can formulate PPCA by first introducing an explicit $K$-dimensional latent variable $\bz$, which corresponds to the principal-component subspace. We then define a Gaussian prior distribution $p(\bz)$ over this latent variable, along with a Gaussian conditional distribution $p(\bx\mid\bz)$ for the $D$-dimensional observed variable $\bx$, conditioned on $\bz$.
Specifically, the prior over $\bz$ is a zero-mean, unit-covariance Gaussian:
\begin{equation}\label{equation:probpca_pz}
	p(\bz) = \normal(\bz \mid \bzero, \bI).
\end{equation}
Similarly, the conditional distribution of the observed variable $\bx$, given $\bz$, is also Gaussian:
\begin{equation}\label{equation:probpca_pxcdz}
	p(\bx\mid\bz) = \normal(\bx \mid \bW\bz + \bmu, \sigma^2 \bI), 
\end{equation}
where the mean of  $\bx$  is a linear function of $\bz$, governed by the matrix $\bW\in \real^{D \times K}$ and the vector $\bmu\in\real^D$. 
Note that this distribution factorizes across the components of  $\bx$. 
As we will see shortly, the columns of $\bW$ span a linear subspace in the data space that corresponds to the principal subspace. The scalar parameter $\sigma^2$ controls the variance of the conditional distribution.
There is no loss of generality in assuming a zero-mean, unit-covariance Gaussian prior for $\bz$: a more general Gaussian prior would lead to an equivalent probabilistic model; see Problem~\ref{prob:ppca_gengauss}.

From a generative perspective, PPCA works as follows: to generate a sample of the observed variable $\bx$, we first draw a value for the latent variable $\bz$, and then sample $\bx$ conditioned on that latent value. Concretely, the $D$-dimensional observed variable $\bx$ is obtained via a linear transformation of the $K$-dimensional latent variable $\bz$, plus additive Gaussian noise:
\begin{equation}\label{equation:probpca_likeli_noise}
\bx = \bW\bz + \bmu + \bepsilon,
\end{equation}
where $\bz\in\real^K$ is a Gaussian latent variable, and $\bepsilon\in\real^D$ is a zero-mean Gaussian noise variable with covariance $\sigma^2 \bI$. 
Note that this formulation defines a mapping from latent space to data space---unlike the standard (non-probabilistic) view of PCA, which typically emphasizes projection from data space to a lower-dimensional subspace. The reverse mapping (from data to latent space) can be derived  using Bayes' theorem (Theorem~\ref{theorem:bayes_theo}).

\subsubsection{Reverse Mappings}
We wish to estimate the parameters $\bW$, $\bmu$, and $\sigma^2$ by maximum likelihood. 
To do so, we need the marginal distribution $p(\bx)$ of the observed variable. By the sum and product rules of probability, this is given by
\begin{equation}\label{equation:probpca_likeli_inte}
	p(\bx) = \int p(\bx\mid\bz) p(\bz) \, d\bz.
\end{equation}
Because this is a linear Gaussian model (see Exercise~\ref{exercise:linear_gauss_model}), the marginal distribution is also Gaussian:
\begin{equation}\label{equation:probpca_likeli}
p(\bx) = \normal(\bx \mid \bmu, \bM) ,
\quad \text{with }\bM \triangleq  \sigma^2 \bI + \bW\bW^\top \in\real^{D\times D}.
\end{equation}
Since $\bz$ and $\bepsilon$ are independent random variables, this result can also be derived directly by using the affine transformation of multivariate Gaussian  using Equation \eqref{equation:probpca_likeli_noise}  (see Lemma~\ref{lemma:affine_mult_gauss}):
\begin{align*}
\Exp[\bx] &= \Exp[\bW\bz + \bmu + \bepsilon] = \bmu;  \\
\Cov[\bx] &= \Exp\left[ (\bW\bz + \bepsilon)(\bW\bz + \bepsilon)^\top \right]
= \sigma^2 \bI+\bW\bW^\top.
\end{align*}

The predictive distribution $p(\bx)$ depends on the parameters $\bmu$, $\bW$, and $\sigma^2$. 
However, this parameterization contains redundancy due to rotational symmetry in the latent space.  
To see this, consider a transformed weight matrix $\widetildebW = \bW\bQ$, where $\bQ$ is an orthogonal matrix. 
Since
$\widetildebW\widetildebW^\top = \bW\bQ\bQ^\top\bW^\top = \bW\bW^\top$,
the covariance matrix $\bM$ remains unchanged. 
Thus, an entire family of matrices $\widetildebW$---differing only by rotations in latent space---yield the same predictive distribution. We will revisit the issue of parameter identifiability later.

When evaluating the predictive distribution, we require $\bM^{-1}$ (see Definition~\ref{definition:multivariate_gaussian}), which involves inverting a $D\times D$ matrix. The computational cost can be reduced using the  \textit{Woodbury identity}: $(\bX + \bY\bP^{-1}\bZ)^{-1} = \bX^{-1} - \bX^{-1}\bY(\bP + \bZ\bX^{-1}\bY)^{-1}\bZ\bX^{-1}$, which yields
\begin{equation}
\bM^{-1} = \sigma^{-2} \bI - \sigma^{-2} \bW \bN^{-1} \bW^\top 
\end{equation}
where the $K \times K$ matrix $\bN$  is defined as
\begin{equation}\label{equation:probpca_matr_wood}
\bN \triangleq   \sigma^2 \bI + \bW^\top \bW\in\real^{K\times K}. 
\end{equation}
Since $K\ll D$ in typical applications, inverting $\bN$ instead of $\bM$ reduces the computational complexity from  $\mathcalO(D^3)$ to $\mathcalO(K^3)$.

In addition to the predictive distribution  $p(\bx)$, we also need the posterior distribution $p(\bz\mid\bx)$. 
Using standard results for linear Gaussian models (Exercise~\ref{exercise:linear_gauss_model}), this posterior is Gaussian:
\begin{equation}\label{equation:probpca_postlat}
p(\bz\mid\bx) = \normal\left( \bz \mid \bN^{-1} \bW^\top (\bx - \bmu), \sigma^2 \bN^{-1} \right).
\end{equation}
Note that while the posterior mean depends on the observed data $\bx$, the posterior covariance is constant---it does not vary with $\bx$.

\begin{SCfigure}
\caption{Graphical representation of PPCA for  a data set of $N$ observations. 
Each observation $\bx_n$ is associated with a latent value $\bz_n$. 
The condition distribution of $\bx_n$ follows from \eqref{equation:probpca_pxcdz}.}
\includegraphics[width=0.5\textwidth]{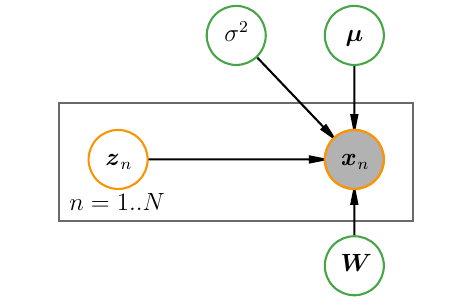}
\label{fig:PPCA_graph}
\end{SCfigure}

\subsubsection{Maximum Likelihood Estimation}
We now turn to estimating the model parameters using maximum likelihood. 
Given a data set $\mathcalX = \{\bx_1, \bx_2, \ldots,\bx_N\}$ of $N$ observed data points, the PPCA model can be represented as a directed graphical model, as shown in Figure~\ref{fig:PPCA_graph}. 
sing the marginal distribution from \eqref{equation:probpca_likeli}, the corresponding log-likelihood function is
\begin{equation}\label{equation:probpca_jointlnlike}
\begin{aligned}
\ln p(\mathcalX \mid \bmu, \bW, \sigma^2) 
&= \sum_{n=1}^{N} \ln p(\bx_n \mid \bW, \bmu, \sigma^2) \\
&= -\frac{ND}{2} \ln(2\pi) - \frac{N}{2} \ln \abs{\bM} - \frac{1}{2} \sum_{n=1}^{N} (\bx_n - \bmu)^\top \bM^{-1} (\bx_n - \bmu)\\
&= -\frac{N}{2} \left\{ D \ln(2\pi) + \ln \abs{\bM} + \trace\left( \bM^{-1} \bS_{\bmu} \right) \right\}, 
\end{aligned}
\end{equation}
where $\bS_{\bmu} \triangleq \frac{1}{N}\sum_{n=1}^{N} (\bx_n - \bmu)(\bx_n-\bmu)^\top$, and the last equality follows from using the standard form of the multivariate Gaussian log-likelihood (see \eqref{equation:multi_gaussian_likelihood}).

\begin{subequations}\label{equations:ppca_mle}
Since the log-likelihood is a quadratic function of $\bmu$, it has a unique maximum, which can be verified by examining the second derivatives.
Setting the derivative of the log-likelihood with respect to $\bmu$ to zero yields the familiar result:
\begin{equation}\label{equation:probpca_bmle}
\bmu_{\text{ML}} = \widebarbx, \quad \text{with }\widebarbx = \sum_{n=1}^{N}\bx_n/N. 
\end{equation}
Consequently, $\bS_{\bmu}$ in \eqref{equation:probpca_jointlnlike} becomes the sample covariance matrix: $\bS_{\bmu} \equiv \frac{1}{N}\sum_{n=1}^{N} (\bx_n - \widebarbx)(\bx_n-\widebarbx)^\top$.
Maximization with respect to $\bW$ and $\sigma^2$ is more involved but still admits a closed-form solution. All stationary points of the log-likelihood can be expressed as
\begin{equation}\label{equation:probpca_wmle}
\bW_{\text{ML}} = \bU_K (\bLambda_K - \sigma^2 \bI)^{1/2} \bQ, 
\end{equation}
where $\bU_K\in\real^{D\times K}$ is a matrix whose columns are any $K$ eigenvectors of the data covariance matrix $\bS$; see \eqref{equation:pca-equ1}, 
$\bLambda_K\in\real^{K\times K}$ is a  diagonal matrix containing the corresponding eigenvalues  $\{\lambda_d\}$, and $\bQ$ is an arbitrary $K \times K$ orthogonal matrix \citep{roweis1997algorithms, tipping1999probabilistic}.

Furthermore, \citet{roweis1997algorithms, tipping1999probabilistic} showed that the global maximum of the likelihood is achieved only when the selected eigenvectors correspond to the $K$ largest eigenvalues; all other stationary points are saddle points. We therefore assume the eigenvalues are ordered such that  $\lambda_1\geq \lambda_2\geq \ldots\geq\lambda_D$ in the spectral decomposition $\bS=\bU\bLambda\bU^\top$.
In this case, the columns of $\bW$ span the same principal subspace as in standard PCA. 
The corresponding maximum likelihood estimate for the noise variance is
\begin{equation}\label{equation:probpca_sigmamle}
\sigma^2_{\text{ML}} = \frac{1}{D - K} \sum_{d=K+1}^{D} \lambda_d.
\end{equation}
Thus, $\sigma^2_{\text{ML}}$ represents the average variance in the discarded dimensions.
\end{subequations}

\begin{remark}[Data Variance and Noise Level]\label{remark:pca_datavar}
It is instructive to examine the structure of the covariance matrix $\bM$ in \eqref{equation:probpca_likeli} or  \eqref{equation:probpca_jointlnlike}. Consider the variance of the predictive distribution along an arbitrary unit direction $\widetildebx$ (i.e., $\widetildebx^\top \widetildebx = 1$), given by $\widetildebx^\top \bM \widetildebx$. 
\begin{itemize}
\item If $\widetildebx$ lies orthogonal to the principal subspace (i.e., it is a linear combination of the discarded eigenvectors), then $\widetildebx^\top \bU_K = \bzero$ and hence $\widetildebx^\top \bM \widetildebx = \sigma^2$. Thus, the model predicts uniform noise variance in directions outside the principal subspace---exactly equal to the average of the discarded eigenvalues, per \eqref{equation:probpca_sigmamle}.
\item If $\widetildebx = \bu_d$, where $\bu_d$ is one of the retained eigenvectors, then $\widetildebx^\top \bM \widetildebx = (\lambda_d - \sigma^2) + \sigma^2 = \lambda_d$. 
\end{itemize}
Hence, the model exactly reproduces the data variance along the principal axes, while approximating all other directions with the single averaged noise level $\sigma^2$.
\end{remark}

\paragrapharrow{Non-identifiability.}
Because $\bQ$ is orthogonal, it acts as a rotation in the $K$-dimensional latent space. 
Substituting the solution for  $\bW_{\text{ML}}$ into the expression for $\bM$ (see \eqref{equation:probpca_likeli}) and using $\bQ\bQ^\top = \bI$, we find that  $\bM$ is independent of $\bQ$. This confirms that the predictive density remains unchanged under rotations in latent space, as noted earlier. 
This rotational freedom in latent space reflects a form of statistical \textit{non-identifiability}: there exists a continuous family of parameter settings---all related by latent-space rotations---that produce identical predictive distributions.

In the special case where $\bQ = \bI$, the columns of $\bW$ align with the principal component directions, scaled by $(\lambda_d - \sigma^2)^{1/2}$ for $d=1,2,\ldots,K$. 
This scaling has a clear interpretation: since $\bM$ arises from the convolution of two independent Gaussian sources (the unit-variance latent prior and the isotropic observation noise; see \eqref{equation:probpca_likeli_noise})---their variances add. 
Specifically, the total variance $\lambda_d$ along eigenvector $\bu_d$ decomposes into (see Remark~\ref{remark:pca_datavar}): 
\begin{itemize} 
\item a signal component $\lambda_d - \sigma^2$, contributed by the projection of the latent variable through the corresponding column of $\bW$;
\item an isotropic noise component $\sigma^2$, added uniformly in all directions.
\end{itemize}

One practical way to construct the maximum likelihood density model is to compute the spectral decomposition of the sample covariance matrix, then directly evaluate $\bW$ and $\sigma^2$ using the formulas in \eqref{equations:ppca_mle}---typically choosing $\bQ = \bI$ for simplicity. 
However, if the parameters are instead obtained via numerical optimization (e.g., using conjugate gradients \citep{nocedal1999numerical, lu2025practical} or the EM algorithm), the resulting $\bQ$ will generally be arbitrary. Consequently, the columns of $\bW$ need not be orthogonal. If an orthogonal basis is required, $\bW$ can be post-processed (e.g., via QR decomposition \citep{lu2021numerical}). Alternatively, the EM algorithm can be modified to directly yield orthogonal principal directions sorted by decreasing eigenvalue \citep{ahn2003constrained}.

Finally, consider the limiting case $K = D$, where no dimensionality reduction occurs. 
Then $\bU_K = \bU$ and $\bLambda_K = \bLambda$.  Using the orthogonality of $\bU\bU^\top = \bI$ and $\bQ\bQ^\top = \bI$, the marginal covariance becomes
\begin{equation*}
\bM =  \sigma^2 \bI = \bU (\sigma^2 \bI) \bU^\top + 
\bU (\bLambda - \sigma^2 \bI)^{1/2} \bQ \bQ^\top (\bLambda - \sigma^2 \bI)^{1/2} \bU^\top 
 = \bS.
\end{equation*}
Thus, PPCA reduces to the standard maximum likelihood estimator for a full-rank Gaussian distribution, with covariance equal to the sample covariance matrix.

\paragrapharrow{Data compression in PPCA.}
Standard PCA is typically framed as a projection from the $D$-dimensional data space onto a $K$-dimensional linear subspace. 
In contrast, PPCA is most naturally interpreted as a generative model that maps from latent space to data space via \eqref{equation:probpca_likeli_noise}. For tasks like visualization or compression, we can invert this mapping using Bayes' theorem (Theorem~\ref{theorem:bayes_theo}). Any data point $\bx$ can then be summarized by its posterior distribution over the latent variable.
From Equation~\eqref{equation:probpca_postlat}, the posterior mean is
\begin{equation*}
\Exp[\bz\mid\bx] = \bN^{-1} \bW_{\text{ML}}^\top (\bx - \widebarbx) ,
\end{equation*}
where $\bN=\sigma^2 \bI + \bW^\top \bW$ (see \eqref{equation:probpca_matr_wood}). 
Mapping this back to data space gives
\begin{equation*}
	\bW \, \Exp[\bz\mid\bx] + \bmu,
\end{equation*}
which has the same functional form as regularized linear regression---a direct consequence of the linear Gaussian structure of the model. 
Moreover,  from \eqref{equation:probpca_postlat}, the posterior covariance, 
$$
\Cov[\bz\mid\bx]  = \sigma^2 \bN^{-1},
$$ 
is constant and does not depend on $\bx$.

\paragrapharrow{Limit analysis and connection to standard PCA.}
Consider the limit $\sigma^2 \to 0$. 
In this case, the posterior mean of $\bz$ becomes
\begin{equation*}
\Exp[\bz\mid\bx] \quad\rightarrow\quad 
(\bW_{\text{ML}}^\top \bW_{\text{ML}})^{-1} \bW_{\text{ML}}^\top (\bx - \widebarbx), 
\end{equation*}
which corresponds to the orthogonal projection of the data point $\bx$ onto the latent subspace. This recovers the mapping used in standard (non-probabilistic) PCA. However, in this limit the posterior covariance vanishes, and the resulting density becomes singular (i.e., degenerate).
For any $\sigma^2 > 0$, the latent projection is shrunk toward the origin relative to the orthogonal projection---a form of regularization induced by the probabilistic model.

\paragrapharrow{Degrees of freedom.}
An important advantage of the PPCA model is that it defines a multivariate Gaussian distribution whose number of degrees of freedom---that is, the number of independent parameters---can be explicitly controlled, while still capturing the dominant correlations in the data.
Recall that a general Gaussian distribution in $D$ dimensions has $D(D+1)/2$ independent parameters in its covariance matrix, plus $D$ parameters for the mean, resulting in a total that grows quadratically with $D$. This quickly becomes impractical in high-dimensional settings.

In contrast, if we restrict the covariance to be diagonal, the number of covariance parameters drops to just $D$, yielding linear scaling with dimensionality. However, this assumption forces all variables to be independent, eliminating the ability to model any correlations.
PPCA offers an elegant compromise: it captures the $K$ most significant directions of correlation while maintaining only linear growth in the number of parameters with respect to $D$.

To see this, consider the parameter count in the PPCA model. The covariance matrix $\bM=\sigma^2\bI+\bW\bW^\top$ is determined by the $D\times K$ matrix $\bW$ and the scalar $\sigma^2$, giving a nominal total of $DK + 1$ parameters. However, this parameterization contains redundancy due to rotational invariance in the latent space: for any orthogonal $K\times K$ matrix $\bQ$, the transformation  $\bW\rightarrow \bW\bQ$ leaves $\bM$ unchanged.

The number of independent parameters in an orthogonal $K\times K$ matrix is $K(K-1)/2$. (This can be seen by noting that the first column has $K-1$ free parameters due to unit-norm constraint, the second has $K-2$ due to orthogonality and normalization, and so on.) Accounting for this redundancy, the effective number of degrees of freedom in $\bM$ is
\begin{equation}
	DK + 1 - K(K-1)/2.
\end{equation}
For fixed $K$, this expression grows linearly with $D$, making PPCA scalable to high dimensions.
Special cases illustrate the flexibility of this framework:
\begin{itemize}
\item When $K=D-1$, the model recovers the full-rank Gaussian distribution. Here, $D-1$ directions of variation are modeled explicitly via $\bW$, while the remaining direction is captured by the isotropic noise term $\sigma^2$.
\item When $K=0$, the model reduces to an isotropic Gaussian with covariance $\sigma^2\bI$, equivalent to assuming all dimensions are independent and identically distributed.
\end{itemize}

\subsubsection{EM Update}

We can now apply the EM algorithm (Algorithm~\ref{alg:em_alg}), derived by iteratively maximizing the evidence lower-bound (ELBO), to learn the parameters of the PPCA model. 
At first glance, this may seem unnecessary, since we already have a closed-form maximum likelihood solution; see \eqref{equations:ppca_mle}. 
However, in high-dimensional settings, the iterative EM approach offers practical computational advantages over explicitly forming and decomposing the sample covariance matrix. Moreover, the same EM framework extends naturally to more complex models like Bayesian PCA---which lacks a closed-form solution---and provides a principled way to handle missing data.

To derive the EM algorithm for PPCA, we follow the standard EM procedure. First, we write down the complete-data log-likelihood, then take its expectation with respect to the posterior distribution of the latent variables using the current (i.e., $t$-th iteration) parameter estimates. Maximizing this expected complete-data log-likelihood yields updated ($(t+1)$-th iteration) parameter values; see Section~\ref{section:em_uncons} for more details.
Assuming independence across data points, the complete-data log-likelihood takes the form
\begin{equation}\label{equation:ppca_complete_log_likelihood}
\ln p\left(\bX, \bZ \mid \bmu, \bW, \sigma^2\right)=\sum_{n=1}^{N}\left\{\ln p\left(\bx_n \mid \bz_n\right)+\ln p\left(\bz_n\right)\right\},
\end{equation}
where the $n$-th row of the matrix $\bZ\in\real^{N\times K}$ is given by $\bz_n^\top$, and the $n$-th row of the matrix $\bX\in\real^{N\times D}$ is given by $\bx_n^\top$.

At iteration $t$, the E-step computes the sufficient statistics of the posterior distribution $p(\bz_n\mid \bx_n)$ for each  data point $n=1,2,\ldots,N$. Specifically, using the current parameter estimates (denoted with superscript $(t)$), we evaluate:
\begin{subequations}\label{equation:probpca_estep}
\begin{align}
\widehatbz_n^\toptzero &\triangleq \Exp[\bz_n \mid \bx_n] = \bN^{-1} \bW^\toptzeroTOP (\bx_n - \bmu^\toptzero); \\
\widehatbSigma_n^\toptzero &\triangleq \Exp[\bz_n \bz_n^\top \mid \bx_n] = \sigma^2 \bN^{-1} + \Exp[\bz_n \mid \bx_n] \Exp[\bz_n \mid \bx_n]^\top,
\end{align}
\end{subequations}
where $\bN = (\sigma^2 \bI + \bW^\toptzeroTOP\bW^\toptzero)$; see \eqref{equation:probpca_postlat}. 
In the M-step, we maximize the expected complete-data log-likelihood with respect to $\bW$ and $\sigma^2$, keeping the posterior statistics fixed. (The mean $\bmu$ is updated separately.) This yields the following updates:
\begin{subequations}\label{equation:probpca_mstep}
\begin{align}
\bW^\toptone &\leftarrow 
\left[ \sum_{n=1}^N (\bx_n - \bmu^\toptzero) \widehatbz_n^\toptzeroTOP \right] 
\left[ \sum_{n=1}^N \widehatbSigma_n^\toptzero \right]^{-1}; \label{equation:probpca_emw} \\
\sigma^{2\toptone} &\leftarrow 
\sum_{n=1}^N 
\left\{\footnotesize 
\frac{\normtwo{\bx_n - \bmu^\toptzero}^2}{ND} 
- \frac{2 \widehatbz_n^\toptzeroTOP \bW^\toptoneTOP (\bx_n - \bmu^\toptzero)}{ND} 
+ \frac{\trace\big(\bW^\toptoneTOP \bW^\toptone\widehatbSigma_n^\toptzero\big)}{ND} 
\right\};
\label{equation:probpca_emsigma}
\end{align}
\end{subequations}
This derivation is a special case of Bayesian PCA (Section~\ref{section:bayespca}); the full proof is deferred to that section.
The EM algorithm for PPCA proceeds by initializing the parameters and then alternately computing the sufficient statistics of the latent space posterior distribution using \eqref{equation:probpca_estep}  in the E-step and revising the parameter values using \eqref{equation:probpca_mstep} in the M-step; see Algorithm~\ref{alg:em_probpca}.

\begin{algorithm}[h] 
\caption{Expectation-Maximization (EM) Algorithm for PPCA}
\label{alg:em_probpca}
\begin{algorithmic}[1] 
\Require Observed data points $\mathcalX=\{\bx_1, \bx_2, \ldots, \bx_N\}$;
\State \textbf{initialize:} $\bW^\topone, \sigma^{2\topone}$; 
\State Choose the maximal number of iterations $C$;
\State $t=0$; \Comment{Count for the number of iterations}
\While{$t<C$} 
\State $t=t+1$;
\State $\widehatbz_n^\toptzero \leftarrow \bN^{-1} \bW^\toptzeroTOP (\bx_n - \widebarbx)$;  
\Comment{(PPCAE$_1$)}
\State $\widehatbSigma_n^\toptzero \leftarrow  \sigma^2 \bN^{-1} + \Exp[\bz_n \mid \bx_n] \Exp[\bz_n \mid \bx_n]^\top$;  
\Comment{(PPCAE$_2$)}
\State $\bW^\toptone \leftarrow 
\left[ \sum_n (\bx_n - \widebarbx) \widehatbz_n^\toptzeroTOP \right] 
\left[ \sum_n \widehatbSigma_n^\toptzero  \right]^{-1}$; 
\Comment{(PPCAM$_1$)}
\State $\sigma^{2\toptone} \leftarrow 
\sum
\left\{\footnotesize 
\frac{\normtwo{\bx_n - \bmu^\toptzero}^2}{ND} 
- \frac{2 \widehatbz_n^\toptzeroTOP \bW^\toptoneTOP (\bx_n - \bmu^\toptzero)}{ND} 
+ \frac{\trace\big(\bW^\toptoneTOP \bW^\toptone\widehatbSigma_n^\toptzero\big)}{ND} 
\right\}$;
\Comment{(PPCAM$_2$)}
\EndWhile
\State Output $\bW^\toptzero$, $\sigma^{2\toptzero}$;
\end{algorithmic} 
\end{algorithm}

\begin{remark}[EM for Standard PCA]
An elegant feature of this framework is that it remains valid even in the limit $\sigma^2 \to 0$, which corresponds to standard (non-probabilistic) PCA \citep{roweis1997algorithms}. 
In this limit, the E-step simplifies: only the posterior mean $\Exp[\bz_n]$ is needed, and the posterior covariance vanishes. 
To highlight the simplicity of the resulting algorithm, define $\bX_c\in\real^{N \times D}$ as the centered data matrix, whose $n$-th row is $(\bx_n - \widebarbx)^\top$, and similarly define $\bG\in\real^{K \times N}$ as the matrix whose $n$-th column is  $\Exp[\bz_n]$. 
Then the E-step \eqref{equation:probpca_estep} becomes
\begin{equation*}
\bG = (\bW^\toptzeroTOP \bW^\toptzero)^{-1} \bW^\toptzeroTOP \bX_c^\top,
\end{equation*}
and the M-step update for $\bW$ is
\begin{equation*}
\bW^\toptone = \bX_c^\top  \bG^\top (\bG  \bG^\top)^{-1}.
\end{equation*}
These updates have an intuitive interpretation:
The E-step orthogonally projects each centered data point onto the current estimate of the principal subspace (i.e., the column space of $\bW^\toptzero$); see Equation~\eqref{equation:ols_tik}.
The M-step re-estimates the subspace to best reconstruct the data, assuming the projections are fixed.
\end{remark}

One key advantage of the EM algorithm for PCA is its computational efficiency in large-scale settings \citep{roweis1997algorithms}. Although standard PCA---based on spectral decomposition of the sample covariance matrix---is non-iterative, it can be expensive in high dimensions.
Specifically: (i) Computing the full covariance matrix costs $\mathcalO(ND^2)$; (ii) Its spectral decomposition costs $\mathcalO(D^3)$; (iii) Even when only the top $K$ eigenvectors are needed (using methods like Lanczos), the cost is typically $\mathcalO(KD^2)$.

In contrast, the EM algorithm never constructs the $D\times D$ covariance matrix explicitly. Its dominant operations involve sums over the data set that scale as $\mathcalO(NDK)$. When $K \ll D$, this can be substantially cheaper than $\mathcalO(ND^2)$, often outweighing the cost of iteration.
Furthermore, the EM algorithm admits an online or streaming implementation. Each data point can be processed independently in the E-step (producing a $K$-dimensional vector and a $K\times K$ matrix), and the M-step only requires accumulating running sums. Thus, data points can be read, processed, and discarded one at a time---making the approach especially suitable when both $N$ and $D$ are large.

\subsection{Bayesian Principal Component Analysis}\label{section:bayespca}
PPCA has been successfully applied to problems in data compression, density estimation, and data visualization. However, like standard PCA, the model itself provides no built-in mechanism for selecting the latent-space dimensionality $K$. When $K=D-1$, PPCA is equivalent to a full-covariance Gaussian distribution. For $K<D-1$, it corresponds to a constrained Gaussian in which the variance in the remaining $D-K$ directions is captured by a single shared parameter $\sigma^2$. Thus, choosing $K$ amounts to a model selection problem that balances complexity against fit.
If sufficient data is available, one practical approach is to use cross-validation to evaluate all candidate values of $K$. However, this quickly becomes computationally infeasible when working with mixture models of PPCA, especially if each mixture component is allowed to have its own latent dimensionality.

This issue of model complexity can be addressed naturally within a Bayesian framework \citep{bishop1998bayesian}. Building on the probabilistic formulation of PCA introduced in Section~\ref{section:ppca}, we again define the generative process as:
\begin{subequations}
\begin{align}
\rvz &\sim  \normal(\bzero, \bI);\\
\rvx\mid \bz &\sim \normal(\bW\bz + \bmu, \sigma^2 \bI), 
\end{align}
\end{subequations}

An additional advantage of this probabilistic model is its ability to handle missing data, provided the data is \textit{missing at random}---that is, the missingness mechanism does not depend on either observed or unobserved values. In such cases, the likelihood is obtained by marginalizing over the unobserved variables, and the resulting objective can be optimized using the EM algorithm.

A fully Bayesian treatment proceeds by placing a prior distribution $p(\bmu, \bW, \sigma^2)$ over the model parameters. The posterior distribution is then given by Bayes' theorem: $p(\bmu, \bW, \sigma^2\mid \mathcalX)\propto p(\mathcalX\mid \bmu, \bW, \sigma^2)p(\bmu, \bW, \sigma^2)$, where the log-likelihood $\ln p(\mathcalX\mid \bmu, \bW, \sigma^2)$ is given in \eqref{equation:probpca_jointlnlike}. 
Finally, predictions for a new data point $\bx'$ are made by marginalizing over the posterior:
\begin{equation}
p(\bx'\mid \mathcalX) = \iiint p(\bx'\mid \bmu, \bW, \sigma^2)p(\bmu, \bW, \sigma^2\mid \mathcalX)\,d\bmu\,d\bW\,d\sigma^2.
\end{equation}

\index{Automatic relevance determination}
\subsubsection{Hyperprior}
\begin{SCfigure}
\caption{Graphical representation of Bayesian PCA for  a data set of $N$ observations. 
	Each observation $\bx_n$ is associated with a latent value $\bz_n$. 
	The condition distribution of $\bx_n$ follows from \eqref{equation:probpca_pxcdz}.}
\includegraphics[width=0.5\textwidth]{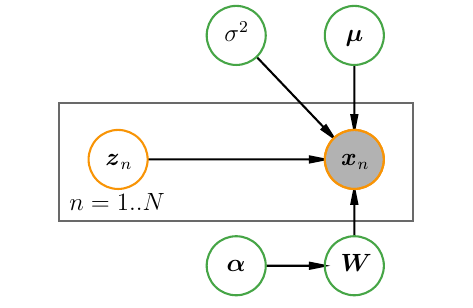}
\label{fig:BPCA_graph}
\end{SCfigure}

To implement this Bayesian framework, two key challenges must be addressed:
(1) the choice of prior distribution, and
(2) the development of a tractable inference algorithm.
For simplicity, in this book, our primary focus is on automatically controlling the effective dimensionality of the latent space---that is, determining how many principal components are truly needed---without resorting to discrete model selection. Instead, we introduce continuous hyper-parameters that adaptively prune irrelevant dimensions during inference.
This is achieved by placing a hierarchical prior over $\bW$, governed by a vector of hyper-parameters $\balpha = \{\alpha_1,\alpha_2, \ldots,\alpha_K\}$.
We set the latent dimensionality to its maximum possible value $K = D - 1$, and assign an independent Gaussian prior to each column $\bw_k$ of $\bW$:
\begin{equation}
p(\bW\mid \balpha) = \prod_{k=1}^{D-1} \left(\frac{\alpha_k}{2\pi}\right)^{D/2} \exp\left\{-\frac{1}{2}\alpha_k\normtwo{\bw_k}^2\right\},
\end{equation}
where $\bw_k$ denotes the $k$-th column of $\bW$. 
This prior is inspired by the \textit{automatic relevance determination} (ARD) framework \citep{mackay1995probable}. 
Each hyper-parameter $\alpha_k$ controls the inverse variance of $\bw_k$:
if the data provides little evidence for a particular latent direction, the posterior for $\alpha_k$ concentrates at large values, forcing $\bw_k\rightarrow \bzero$.
In effect, that dimension is ``switched off." 
The probabilistic structure of the resulting Bayesian PCA model is illustrated in Figure~\ref{fig:BPCA_graph}.

Similarly, the mode  $\bW_{\text{MAP}}$ is found by maximizing the log-posterior:
\begin{equation}
\ln p(\bW\mid \mathcalX) = L - \frac{1}{2} \sum_{k=1}^{D-1} \alpha_k \normtwo{\bw_k}^2 + \text{const.}
\end{equation}
where $L=\ln p(\mathcalX\mid \bmu, \bW, \sigma^2)$ is given by \eqref{equation:probpca_jointlnlike}. 
For the purpose of controlling latent dimensionality, we treat $\bmu$, $\sigma^2$ and $\balpha$ as fixed parameters to be estimated---not as random variables. 
This avoids the need to specify priors for them. Specifically:
(i)  $\bmu$ and $\sigma^2$   are estimated via maximum likelihood; (ii) $\balpha$ is estimated via type-II maximum likelihood, i.e., by maximizing the marginal likelihood $p(\mathcalX\mid\balpha)$, obtained by integrating out $\bW$ under a quadratic (Laplace) approximation around $\bW_{\text{MAP}}$ \citet{mackay1995probable, bishop1998bayesian}.
This leads to the following re-estimation formula for each hyper-parameter:
\begin{equation}\label{equation:bayespca_hyperalpha}
\alpha_k \leftarrow  \frac{\gamma_k}{\normtwo{\bw_k}^2},
\end{equation}
where $\gamma_k = D - \alpha_k \trace_k(\bH^{-1})$ represents the effective number of well-determined parameters in $\bw_k$, $\bH$ is the Hessian of $\ln p(\bW\mid \mathcalX)$ evaluated at  $\bW_{\text{MAP}}$, and $\trace_k(\cdot)$ denotes the trace over the block of $\bH^{-1}$ corresponding to $\bw_k$ \citep{bishop1998bayesian}.

Following \citet{bishop1998bayesian}, we adopt a further simplification: $\gamma_k\approx D$, assuming all model parameters are well-constrained by the data. This avoids the costly computation and storage of the full Hessian matrix. Under this approximation, any column $\bw_k$ that lacks sufficient support from the data will be driven to zero, causing $\alpha_k \to \infty$, thereby completely deactivating that latent dimension.
We define the \textit{effective dimensionality} $K_{\text{eff}}$ of the model as the number of columns $\bw_k$ that remain nonzero after convergence.

The effective dimensionality estimated by Bayesian PCA depends on the number $N$ of data points. As $N \to \infty$, we expect the effective dimensionality $K_{\text{eff}}$ to approach $D-1$.
In this limit, the maximum likelihood framework and the Bayesian approach yield identical results. For finite datasets, however, $K_{\text{eff}}$ may be reduced: directions in latent space that lack sufficient support from the data are automatically suppressed. The variance of the data along the remaining $D-K_{\text{eff}}$ directions is then captured by the single shared parameter $\sigma^2$.

\subsubsection{EM Update}
The MAP estimate $\bW_{\text{MAP}}$ can be computed efficiently using the EM algorithm (Algorithm~\ref{alg:em_alg}).
At iteration $t$, the E-step computes the expected sufficient statistics of the latent posterior distribution $p(\bz_n \mid \bx_n)$ for each data point  $n=1,2,\ldots,N$:
\begin{subequations}\label{equation:bayespca_estep}
\begin{align}
\widehatbz_n^\toptzero &\triangleq \Exp[\bz_n \mid \bx_n] = \bN^{-1} \bW^\toptzeroTOP (\bx_n - \bmu^\toptzero); \\
\widehatbSigma_n^\toptzero &\triangleq \Exp[\bz_n \bz_n^\top \mid \bx_n] = \sigma^2 \bN^{-1} + \Exp[\bz_n \mid \bx_n] \Exp[\bz_n \mid \bx_n]^\top,
\end{align}
\end{subequations}
where $\bN = (\sigma^2 \bI + \bW^\toptzeroTOP\bW^\toptzero)$; see \eqref{equation:probpca_postlat}. 
In the M-step, the model parameters are updated as follows:
\begin{subequations}\label{equation:bayespca_mstep}
\begin{align}
\bW^\toptone &\leftarrow 
\left[ \sum_{n=1}^N (\bx_n - \bmu^\toptzero) \widehatbz_n^\toptzeroTOP \right] 
\left[ \sum_{n=1}^N \widehatbSigma_n^\toptzero + \sigma^{2\toptzero} \diag(\balpha^\toptzero) \right]^{-1}; \label{equation:bayespca_emw} \\
\sigma^{2\toptone} &\leftarrow 
\sum_{n=1}^N 
\left\{\footnotesize 
\frac{\normtwo{\bx_n - \bmu^\toptzero}^2}{ND} 
- \frac{2 \widehatbz_n^\toptzeroTOP \bW^\toptoneTOP (\bx_n - \bmu^\toptzero)}{ND} 
+ \frac{\trace\big(\bW^\toptoneTOP \bW^\toptone\widehatbSigma_n^\toptzero\big)}{ND} 
\right\};
\label{equation:bayespca_emsigma}\\
\alpha_k^\toptone & \leftarrow 
{\gamma_k}/{\normtwobig{\bw_k^\toptone}^2}, 
\quad k=1,2,\ldots,K;
\label{equation:bayespca_emalpha}\\
\bmu^\toptone &\leftarrow \widebarbx,
\label{equation:bayespca_emmu}
\end{align}
\end{subequations}
Note that the update for $\bmu$ is identical at every iteration, i.e, the sample mean (see \eqref{equation:probpca_bmle}).
The EM algorithm for Bayesian PCA proceeds by initializing the parameters and then alternating between
the E-step (computing the posterior sufficient statistics via \eqref{equation:bayespca_estep}) and
the M-step (updating all parameters using\eqref{equation:bayespca_mstep})
until a convergence criterion is satisfied (e.g., small changes in parameter values or log-likelihood).
The updates for  $\bW$ and $\sigma^2$ are interleaved with re-estimation of the hyper-parameters $\alpha_k$ using  \eqref{equation:bayespca_hyperalpha}, where we set $\gamma_i = D$ (the data dimensionality) as a simplifying approximation. The complete procedure is summarized in Algorithm~\ref{alg:em_bayespca}.

\begin{algorithm}[h] 
\caption{Expectation-Maximization (EM) Algorithm for Bayesian PCA}
\label{alg:em_bayespca}
\begin{algorithmic}[1] 
\Require Observed data points $\mathcalX=\{\bx_1, \bx_2, \ldots, \bx_N\}$;
\State \textbf{initialize:} $\bW^\topone, \sigma^{2\topone}, \balpha^\topone$; 
\State Choose the maximal number of iterations $C$;
\State $t=0$; \Comment{Count for the number of iterations}
\While{$t<C$} 
\State $t=t+1$;
\State $\widehatbz_n^\toptzero \leftarrow \Exp[\bz_n \mid \bx_n] = \bN^{-1} \bW^\toptzeroTOP (\bx_n - \bmu^\toptzero)$;  
\Comment{(BPCAE$_1$)}
\State $\widehatbSigma_n^\toptzero \leftarrow \Exp[\bz_n \bz_n^\top \mid \bx_n] = \sigma^2 \bN^{-1} + \Exp[\bz_n \mid \bx_n] \Exp[\bz_n \mid \bx_n]^\top$;  
\Comment{(BPCAE$_2$)}
\State $\bW^\toptone \leftarrow 
\left[ \sum_n (\bx_n - \bmu) \widehatbz_n^\toptzeroTOP \right] 
\left[ \sum_n \widehatbSigma_n^\toptzero + \sigma^{2\toptzero} \diag(\balpha^\toptzero) \right]^{-1}$; \Comment{(BPCAM$_1$)}
\State $\sigma^{2\toptone} \leftarrow 
\sum
\left\{\footnotesize 
\frac{\normtwo{\bx_n - \bmu^\toptzero}^2}{ND} 
- \frac{2 \widehatbz_n^\toptzeroTOP \bW^\toptoneTOP (\bx_n - \bmu^\toptzero)}{ND} 
+ \frac{\trace\big(\bW^\toptoneTOP \bW^\toptone\widehatbSigma_n^\toptzero\big)}{ND} 
\right\}$;
\Comment{(BPCAM$_2$)}
\State $\alpha_k^\toptone  \leftarrow 
\frac{\gamma_k}{\normtwobig{\bw_k^\toptone}^2}, 
\quad k=1,2,\ldots,K$;  
\Comment{(BPCAM$_3$)}
\State $\bmu^\toptone \leftarrow \widebarbx$; \Comment{Same for each iteration, (BPCAM$_4$)}
\EndWhile
\State Output $\bW^\toptzero$, $\sigma^{2\toptzero}$, $\balpha^\toptzero$;
\end{algorithmic} 
\end{algorithm}
\begin{proof}[of EM update \eqref{equation:bayespca_mstep}]
Since the update for $\bmu$ is simply the sample mean (and thus constant across iterations), we focus on estimating  $\bW$  and  $\sigma^2$.
Because the data points are independent, we consider the joint distribution for a single observation $\{\bx, \bz\}$:
$$
p(\bx, \bz) = p(\bz) p(\bx \mid \bz) 
= \normal(\bz \mid \bzero, \bI) \cdot \normal(\bx \mid \bW \bz + \bmu, \sigma^2 \bI).
$$
In the MAP-EM framework (see Algorithm~\ref{alg:em_alg} and  Equation~\eqref{equation:mapem2}), we maximize the expected complete-data log-likelihood plus the log-prior:
\begin{align*}
Q(\btheta \mid \btheta^\toptzero) 
&= \Exp_{\bz \sim p(\bz \mid \bx, \btheta^\toptzero)} [\ln p(\bx, \bz \mid \btheta)]
+\ln p(\btheta)
=\Exp[\ln p(\bz)] + \Exp[\ln p(\bx \mid \bz)]+\ln p(\btheta)\\
&= \Exp\left[-\frac{1}{2} \bz^\top \bz\right] + \Exp\left[ -\frac{1}{2\sigma^2} \normtwo{\bx - \bW \bz - \bmu}^2 \right] - \frac{D}{2} \ln(2\pi\sigma^2) - \frac{K}{2} \ln(2\pi) +\ln p(\btheta),
\end{align*}
where we used the fact that  $  \ln p(\bz) = -\frac{1}{2} \bz^\top \bz - \frac{K}{2} \ln(2\pi)  $ 
and  $  \ln p(\bx \mid \bz) = -\frac{1}{2\sigma^2} \normtwo{\bx - \bW \bz - \bmu}^2 - \frac{D}{2} \ln(2\pi\sigma^2)  $ from \eqref{equation:probpca_pz} and \eqref{equation:probpca_pxcdz}.

\paragraph{M-step: update  $\bW$.}
We fix  $  \sigma^2=\sigma^{2\toptzero}$, $\widehatbz=\widehatbz^\toptzero$, $\widehatbSigma=\widehatbSigma^\toptzero$, $\balpha=\balpha^\toptzero$,  and optimize  $\bW$.
First, expand the quadratic term:
\begin{align*}
\Exp[\normtwo{\bx - \bW \bz - \bmu}^2]
&= \normtwo{\bx - \bmu}^2 - 2(\bx - \bmu)^\top \bW \Exp[\bz] + \Exp[\bz^\top \bW^\top \bW \bz]\\
&= \normtwo{\bx - \bmu}^2 - 2(\bx - \bmu)^\top \bW \widehatbz + \trace(\bW^\top \bW \widehatbSigma),
\end{align*}
where the last equality follows from 
$
\Exp[\bz^\top \bW^\top \bW \bz] = \trace(\bW^\top \bW \Exp[\bz \bz^\top]) = \trace(\bW^\top \bW \widehatbSigma)$.
Substituting into $  Q(\btheta \mid \btheta^\toptzero)  $:
\begin{align*}
Q(\bW \mid \bW^\toptzero) 
&= -\frac{1}{2} \sum_{n=1}^N \Exp[\bz_n^\top \bz_n] - \frac{1}{2\sigma^2} \sum_{n=1}^N \left( \normtwo{\bx_n - \bmu}^2 
- 2(\bx_n - \bmu)^\top \bW \widehatbz_n + \trace(\bW^\top \bW \widehatbSigma_n) \right)\\
&\quad  - \sum_{k=1}^{K}\frac{1}{2} \alpha_k\normtwo{\bw_k}^2 + \text{const.}
\end{align*}
Discarding terms independent of $\bW$,  we minimize the following objective:
\begin{equation}\label{equation:bayespca_emw_loss}
\sum_{n=1}^{N} \left( -2(\bx_n - \bmu)^\top \bW \widehatbz_n + \trace(\bW^\top \bW \widehatbSigma_n) \right)
+ \sigma^2 \sum_{k=1}^{K} \alpha_k\normtwo{\bw_k}^2.
\end{equation}
Taking the derivative with respect to $\bW$ and setting it to zero yields the closed-form update in \eqref{equation:bayespca_emw}.

\paragraph{M-step: update  $\sigma^2$.}
Now fix  $  \bW =\bW^\toptone$, and optimize  $\sigma^2$.
From  $  Q(\btheta\mid \btheta^\toptzero)  $:
$$
Q(\sigma^2 \mid \sigma^{2\toptzero}) 
= -\frac{1}{2\sigma^2} \sum_{n=1}^N \Exp[\normtwo{\bx_n - \bW \bz_n - \bmu}^2] - \frac{D}{2} \ln(2\pi\sigma^2)
$$
Denote $ r_n \triangleq \Exp[\normtwo{\bx_n - \bW \bz_n - \bmu}^2] $.
Then, 
$Q(\sigma^2) = -\frac{1}{2\sigma^2} \sum_{n=1}^N r_n - \frac{D}{2} \ln \sigma^2 - \text{const}$.
Taking derivative w.r.t.  $\sigma^2$ yields
$\sigma^2 = \frac{1}{ND} \sum_{n=1}^N r_n$.
This completes the derivation of the EM updates for Bayesian PCA.
\end{proof}

\subsection{Mixtures of Probabilistic and Bayesian PCA Models}
We have introduced mixture models---for example, mixtures of Gaussians or Bernoullis---in Example~\ref{example:gmm_twoclus} and Problems~\ref{prob:mix_of_gauss}--\ref{prob:mix_of_bern}.
Given a probabilistic formulation of PCA, it is straightforward to construct a mixture distribution as a linear superposition of principal component analyzers. In the case of maximum-likelihood PCA, we must choose both the number $Q$ of mixture components and the latent space dimensionality $K$ for each component \citep{tipping1999mixtures}. However, even for moderate values of $Q$ and data spaces of several dimensions, it quickly becomes computationally intractable to explore the exponentially large number of possible combinations of $K$ values across components. In this case, Bayesian PCA offers a significant advantage: it allows the effective dimensionalities of the models to be determined automatically.

Specifically, given a data set $\mathcalX=\{\bx_1, \bx_2, \ldots,\bx_n\}$ of $N$ observations, we consider a mixture of PPCA model.
For each data point  $\bx_n$:
\begin{enumerate}[label=(\roman*)]
\item Choose a component  $ q \in \{1,2,\dots,Q\} $  with probability  $ \pi_q $.
\item Draw a latent vector  $ \bz_{nq} \sim \normal(\bzero, \bI_q) $, where $\bz_{nq}\in\real^K$.
\item Generate  $ \bx_n \sim \normal\big( \bmu_q + \bW_q \bz_{nq},\; \sigma_q^2 \bI \big) $, i.e., each observed data point  $ \bx_n \in \real^D $  is generated from a lower-dimensional latent variable  $ \bz_{nq} \in \real^K $  ( $ K<D $ ) via a linear Gaussian mapping, where  $ \bW_q \in \real^{D\times K} $ denotes weight matrices,   $ \bmu_q \in \real^D $ denotes means, and $ \sigma_q^2 $ denotes isotropic noise variances.
\end{enumerate}
Using \eqref{equation:probpca_likeli}, the marginal distribution of  $\bx_n$  under component  $q$  is:
\begin{equation}
p(\bx_n \mid q) = \normal\big( \bx_n \mid \bmu_q,\; \bM_q \big), \quad \text{where } \bM_q =  \sigma_q^2 \bI+ \bW_q \bW_q^\top.
\end{equation}
In other words, we assume the data are generated from one of $Q$ such PPCA components, each with its own set of parameters. The full mixture likelihood is:
\begin{equation}
p(\mathcalX \mid \btheta) = \prod_{n=1}^N \sum_{q=1}^Q \pi_q \, \normal\big( \bx_n \mid \bmu_q, \bM_q \big),
\end{equation}
where  $ \btheta \triangleq  \{ \pi_q, \bmu_q, \bW_q, \sigma_q^2 \}_{q=1}^Q $.
This approach enables automatic clustering of the data into groups, where each group is modeled by a low-dimensional probabilistic subspace---thus unifying dimensionality reduction and clustering within a single probabilistic framework.

\paragrapharrow{EM update for mixture of PPCA.}
For such models, a mixture of maximum-likelihood PPCA components can be estimated using the EM algorithm. 
In the M-step, we apply the maximum-likelihood updates from Equation~\eqref{equations:ppca_mle}, using eigenvectors and eigenvalues derived from weighted covariance matrices, where the weights are the posterior responsibilities computed in the E-step.
The EM algorithm for a mixture of PPCA models combines two powerful ideas:
\begin{enumerate}[label=(\roman*)]
\item \textit{Probabilistic PCA (PPCA).} A latent-variable generative model in which each observed data point  $ \bx \in \real^D $  is generated from a lower-dimensional latent variable  $ \bz \in \real^K $  ($ K<D $) via a linear Gaussian mapping.
\item \textit{Mixture modeling.} Assumes the data originate from one of $Q$ such PPCA components, each with its own parameters.
\end{enumerate}

Following Problems~\ref{prob:mix_of_gauss}--\ref{prob:mix_of_bern} and given a data set of $N$ observations $\mathcalX=\{\bx_1,\bx_2,\ldots,\bx_N\}$, we introduce latent indicator variables  $ y_{nq} \in \{0,1\} $, where  $ y_{nq} = 1 $  if  $\bx_n$  belongs to component  $q$. 
The posterior probability that component $q$ generated $\bx_n$ is:
\begin{subequations}\label{equation:em_mix_ppca}
\begin{equation}
\zeta_{nq} \triangleq p(y_{nq} = 1 \mid \bx_n, \btheta^{\toptzero}) = 
\frac{ \pi_q^\toptzero \, \normal( \bx_n \mid \bmu_q^\toptzero, \bM_q^\toptzero ) }
{ \sum_{\ell=1}^Q \pi_\ell^\toptzero \, \normal( \bx_n \mid \bmu_\ell^\toptzero, \bM_\ell^\toptzero ) }.
\end{equation}
Additionally, to update $ \bW_q $  and  $ \sigma_q^2 $, we require expectations over the latent variables $ \bz_{nq} $. 
At iteration $t$, using properties of conditional Gaussians in PPCA (see \eqref{equation:probpca_estep}), we compute the sufficient statistics of the posterior distribution $p(\bz_{nq}\mid \bx_n)$ for each data point $n=1,2,\ldots,N$ and each component $q=1,2,\ldots,Q$:
\begin{align}
\widehatbz_{nq}^\toptzero &= \bN_q^{-1} \bW_q^\toptzeroTOP (\bx_n - \bmu_q^\toptzero),
\quad \text{where } \bN_q = \bW_q^\toptzeroTOP \bW_q^\toptzero + \sigma_q^{2,\toptzero} \bI_q;\\
\widehatbSigma_{nq}^\toptzero &= \sigma_q^{2,\toptzero} \bN_q^{-1} + \widehatbz_{nq}^\toptzero \widehatbz_{nq}^\toptzeroTOP.
\end{align}
These quantities are weighted by $ \zeta_{nq} $  in the M-step, yielding the following parameter updates:
\begin{align}
\pi_q^{\toptone} &\leftarrow \frac{N_q}{N};\qquad \qquad
\bmu_q^{\toptone} \leftarrow \frac{1}{N_q} \sum_{n=1}^N \zeta_{nq} \bx_n;\\
\bW_q^{\toptone} &\leftarrow \left[ \sum_{n=1}^N \zeta_{nq} (\bx_n - \bmu_q^{\toptone}) \widehatbz_{nq}^\toptzeroTOP \right] \left[ \sum_{n=1}^N \zeta_{nq} \widehatbSigma_{nq}^\toptzero \right]^{-1};\\
\footnotesize
\sigma_q^{2,\toptone} &\leftarrow  \frac{1}{N_q D}\sum_{n=1}^N \zeta_{nq} \left\{\normtwo{\bx_n - \bmu_q^{\toptone}}^2 
- A^\toptzero
+ B^\toptzero\right\},
\end{align}
where $ N_q \triangleq \sum_{n=1}^N \zeta_{nq} $  (i.e., the effective number of points assigned to component  $q$), 
$A^\toptzero\triangleq 2\big(\widehatbz_{nq}^\toptzeroTOP \bW_q^{\toptone\top} (\bx_n - \bmu_q^{\toptone})  \big) $, 
and $B^\toptzero\triangleq \trace\big( \bW_q^{\toptone\top} \bW_q^{\toptone} \widehatbSigma_{nq}^\toptzero \big)$.
\end{subequations}
The EM algorithm for a mixture of PPCA alternates between the E-step (i.e., computing component responsibilities  $ \zeta_{nq} $  and expected latent variables  $ \widehatbz_{nq}^\toptzero $,  $ \widehatbSigma_{nq}^\toptzero $) and the M-step (i.e.,  updating mixing weights, means, projection matrices  $ \bW_q $, and noise variances using weighted averages based on  $ \zeta_{nq} $); see Algorithm~\ref{alg:em_mix_ppca}.

\begin{algorithm}[h] 
\caption{Expectation-Maximization (EM) Algorithm for Mixture of PPCA}
\label{alg:em_mix_ppca}
\begin{algorithmic}[1] 
\Require Observed data points $\mathcalX=\{\bx_1, \bx_2, \ldots, \bx_N\}$;
\State \textbf{initialize:} $\{\bW_q^\topone\}, \{\bmu_q^\topone\}, \{\sigma_q^{2\topone}\}, \{y_{nq}^\topone\}, \{\pi_q^\topone\}$; 
\State Choose the maximal number of iterations $C$;
\State $t=0$; \Comment{Count for the number of iterations}
\While{$t<C$} 
\State $t=t+1$;
\State $\widehatbz_{nq}^\toptzero \leftarrow \bN_q^{-1} \bW_q^\toptzeroTOP (\bx_n - \bmu_q^\toptzero)$;  
\Comment{(MPPCAE$_1$)}
\State $\widehatbSigma_{nq}^\toptzero \leftarrow \sigma_q^{2,\toptzero} \bN_q^{-1} + \widehatbz_{nq}^\toptzero \widehatbz_{nq}^\toptzeroTOP$;  
\Comment{(MPPCAE$_2$)}
\State $\pi_q^{\toptone} \leftarrow \frac{N_q}{N}$, where $N_q = \sum_{n} \zeta_{nq}$; \Comment{(MPPCAM$_1$)}
\State $\bmu_q^{\toptone} \leftarrow \frac{1}{N_q} \sum_{n}\zeta_{nq} \bx_n$;
\Comment{(MPPCAM$_2$)}
\State $\bW_q^{\toptone} \leftarrow \left[ \sum_{n} \zeta_{nq} (\bx_n - \bmu_q^{\toptone}) \widehatbz_{nq}^\toptzeroTOP \right] \left[ \sum_{n} \zeta_{nq} \widehatbSigma_{nq}^\toptzero \right]^{-1}$;  
\Comment{(MPPCAM$_3$)}
\State $\sigma_q^{2,\toptone} \leftarrow  \frac{1}{N_q D}\sum_{n} \zeta_{nq} \left\{\normtwo{\bx_n - \bmu_q^{\toptone}}^2 
- A^\toptzero
+ B^\toptzero\right\}$; \Comment{(MPPCAM$_4$)}
\EndWhile
\State Output $\bW^\toptzero$, $\bmu^\toptzero$, $\sigma^{2\toptzero}$;
\end{algorithmic} 
\end{algorithm}

\paragrapharrow{Mixture of Bayesian PCA.} 
This mixture framework extends naturally to Bayesian PCA. In maximum-likelihood PCA, independently selecting a different latent dimensionality $K$ for each component is computationally impractical, so we typically assume all components share the same $K$.
Bayesian PCA is especially advantageous for small datasets in high-dimensional spaces, as it avoids the singularities that often plague maximum-likelihood (or standard) PCA by automatically suppressing irrelevant degrees of freedom. This benefit is particularly valuable in mixture modeling: even when the total dataset is large, the effective number of points associated with individual clusters may be small, making regularization through Bayesian inference crucial for stable estimation.

\index{Autoencoder}
\section{Autoencoder and Variational Autoencoder}
\textit{Autoencoders} and their probabilistic counterpart, the \textit{variational autoencoders (VAEs)}, are among the most influential architectures in unsupervised and self-supervised learning. They offer a deep learning–based framework for learning compressed, structured representations of data---enabling tasks ranging from denoising and anomaly detection to generative modeling. What makes these models especially compelling is that they generalize and extend ideas long familiar in classical data analysis, particularly matrix decomposition techniques such as principal component analysis (PCA) and nonnegative matrix factorization (NMF).
At its heart, an autoencoder consists of two transformations or two neural networks: an encoder that maps input data 
$\bx$ into a lower-dimensional latent code $\bz$, and a decoder that reconstructs the input from this code. When both encoder and decoder are constrained to be linear and the reconstruction loss is squared error, the optimal solution of an autoencoder with a $K$-dimensional bottleneck is mathematically equivalent to performing PCA, which is, as mentioned previously, a form of matrix decomposition that factorizes the data matrix $\bX\in\real^{N\times D}$
into orthogonal components capturing maximal variance. In this sense, the autoencoder can be viewed as a nonlinear, learnable generalization of matrix factorization, where the ``factors" are no longer restricted to linear subspaces or nonnegativity constraints, but can capture complex, hierarchical patterns through deep nonlinear mappings.
The VAE builds on this foundation by introducing a probabilistic interpretation: instead of producing a single point estimate $\bz$, the encoder outputs parameters of a distribution (typically Gaussian) over the latent space. The VAE then optimizes a lower bound on the data log-likelihood---the evidence lower bound (ELBO)---which balances reconstruction fidelity against regularization via a prior (often a standard normal distribution). This regularization encourages the latent space to be smooth and well-structured, enabling meaningful interpolation and generation. From a matrix perspective, one can view the VAE as performing a stochastic, regularized, and nonlinear decomposition of the data matrix, where uncertainty and generative capacity are explicitly modeled.
This bridge between classical linear algebra and modern deep generative modeling underscores why autoencoders and VAEs deserve careful study. They inherit the dimensionality-reduction intuition of matrix decomposition while vastly expanding its expressive power through neural networks and probabilistic inference. In doing so, they exemplify a broader theme in machine learning: the evolution of classical methods into flexible, data-driven frameworks capable of handling real-world complexity.
Thus, understanding autoencoders and VAEs not only equips us with practical tools for representation learning and generation but also deepens our appreciation of how foundational ideas---like decomposing data into interpretable parts---continue to evolve in the era of deep learning.

\begin{figure}[H]
\centering
\includegraphics[width=0.99\textwidth]{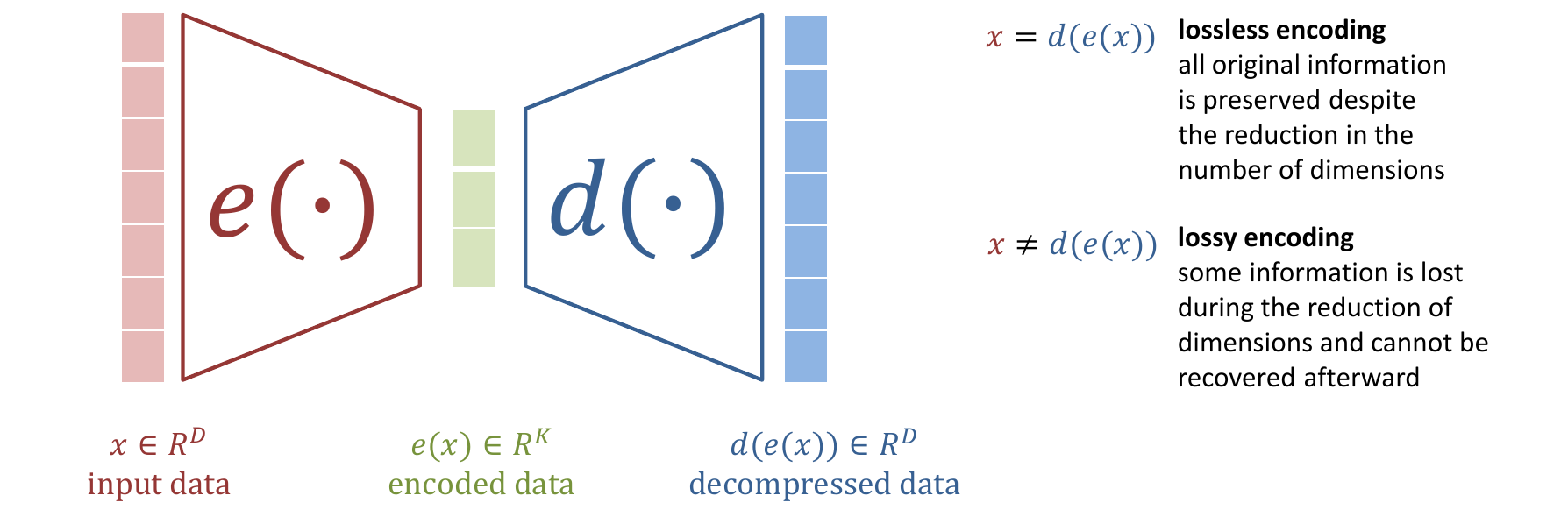}
\caption{Description of an autoencoder.}
\label{fig:autoencoder}
\end{figure}
\subsection{Autoencoder}
In machine learning, an \textit{autoencoder} performs dimensionality reduction by decreasing the number of features used to describe a data set (denoted as $\bx$ in Figure~\ref{fig:autoencoder}). 
This is achieved through an \textit{encoding process}, denoted $e(\bx)$, which either selects a subset of the original features or constructs a smaller set of new features derived from them.
The \textit{decoding process}, written as  $d(e(\bx))$, attempts to reconstruct the original input from its compressed (latent) representation. Depending on the data distribution, the dimensionality of the latent space, and the design of the encoder, this compression may be lossy---meaning some information is irreversibly lost during encoding and cannot be fully recovered during decoding.

Thus, the primary goal of an autoencoding method is to identify the optimal encoder-decoder pair from a given set of candidates. More precisely, given families of possible encoders $\gE$ and decoders $\gD$, we seek the pair that preserves as much information as possible during encoding, thereby minimizing reconstruction error during decoding. This leads to the following optimization formulation:
$$
(e, d) = \mathop{\argmin}_{(e,d)\in(\gE,\gD)} f\left(\bx, d(e(\bx))\right),
$$
where $f(\cdot, \cdot)$  is a loss function measuring the discrepancy between the original input and its reconstruction.
The families $\gE$ and $\gD$ can consist of any suitable functions---for example, multilayer perceptrons or deep neural networks \citep{lecun2015deep, goodfellow2016deep}.
However, when both the encoder and decoder are restricted to linear transformations (i.e., a \textit{linear autoencoder}), the solution aligns with PCA or SVD---provided the loss function $f$ is based on the Frobenius norm or spectral norm (see Problem~\ref{theorem:young-theorem_frob} or \citet{lu2021numerical}).\footnote{Strictly speaking, PCA (or SVD) corresponds to a special case of the linear autoencoder in which the learned basis vectors are orthonormal. In contrast, a general linear autoencoder does not require its weight vectors to be orthogonal or normalized.}

Consequently, the weight vectors defining the linear transformation in Figure~\ref{fig:autoencoder} span the principal subspace, though they need not be orthogonal or unit-length. This equivalence is unsurprising: both PCA and linear autoencoders perform linear dimensionality reduction and minimize the same sum-of-squares reconstruction error.

One might expect that the limitations of linear manifolds could be overcome by introducing nonlinear activation functions, such as those used in deep neural networks. However, even with nonlinear transformations, the minimum reconstruction error under squared loss is still achieved by projecting the data onto the principal component subspace \citep{bourlard1988auto}. Therefore, two-layer neural networks offer no advantage over PCA for linear dimensionality reduction under these conditions. In contrast, standard PCA methods based on SVD---or alternating least squares (ALS; see Chapter~\ref{chapter:als})---are guaranteed to converge to the globally optimal solution in finite time and produce an ordered set of eigenvalues with corresponding orthonormal eigenvectors.

\paragrapharrow{Solution via truncated SVD or PCA.}
Assume the data are centered so that the sample mean $\widebarbx=\bzero$. 
Let the data matrix $\bX \in \real^{N\times D}$ contain the $N$ observations as rows.
As shown in Problem~\ref{theorem:young-theorem_frob}, the \textit{truncated SVD (TSVD)}---which sets all but the top $K$ singular values to zero---provides the best rank-$K$ approximation of $\bX$ in the Frobenius norm. Denote this approximation by $\widetildebX=\bU_K\bSigma_K\bV_K^\top$, where $\bU_K\in\real^{N\times K}$, $\bV_K\in\real^{D\times K}$, and $\bSigma_K\in\real^{K\times K}$. 
From the perspective of PCA, the encoder maps each data point $\bx_n\in\real^D$ (the $n$-th row of $\bX$) to its coordinates in the principal subspace:
$$
\text{encoder: }\gap e(\bx_n) = \bV_K^\top \bx_n, \gap \bx_n\in\real^{D},\, \forall\, n\in\{1,2,\ldots, N\}.
$$
Here, the columns of $\bV_K$ are the orthonormal eigenvectors corresponding to the $K$ largest eigenvalues of the data covariance matrix. 
Since $\bV_K^\top\bV_K=\bI_K$, it follows that $\bV_K^\top\widetildebX^\top=\bSigma_K\bU_K^\top$.
The corresponding decoder reconstructs the input by projecting back into the original space:
$$
\text{decoder: }\gap d(e(\bx_n)) = \bV_K e(\bx_n), \gap  \forall\, n\in\{1,2,\ldots, N\}.
$$
Thus, the full reconstruction is $d(e({\bX}^\top))=\bV_K\bSigma_K\bU_K^\top$, which is precisely the truncated SVD of $\bX^\top$. 
This confirms that when $\gE$ and $\gD$ are linear, the autoencoder is equivalent to PCA/SVD.

\paragrapharrow{Other formulations.}
The autoencoder can also be trained by minimizing a reconstruction loss directly:
\begin{equation}
J(\btheta) = \frac{1}{2}\sum_{n=1}^{N}\normtwo{e(\bx_n; \btheta) - \bx_n}^2,
\end{equation}
where $\btheta$ denotes the model parameters (e.g., weights of a neural network).

In general, the encoder  $e(\bx_n; \btheta)$ can be implemented using various models---linear, nonlinear, or deep neural networks \citep{bishop2006pattern}. 
However, to prevent the model from simply learning the identity mapping (which would yield perfect reconstruction but no useful compression), it is essential to constrain the capacity of the latent representation. One common approach is to limit the dimensionality of the hidden (bottleneck) layer.
An alternative strategy is to encourage sparsity in the internal representation through regularization. A popular choice is the $\ell_1$ penalty, which promotes sparse activations (see Section~\ref{section:more_err_sta_als}). This leads to the regularized objective:
\begin{equation}
J(\btheta) = \frac{1}{2}\sum_{n=1}^{N}\normtwo{e(\bx_n; \btheta) - \bx_n}^2
+
\lambda\normone{\btheta},
\end{equation}
where $\lambda\in\real_+$ controls the strength of regularization.

Another effective approach is the \textit{denoising autoencoder} \citep{vincent2008extracting}, which forces the model to learn robust data representations by training it to reconstruct clean inputs from corrupted versions. Specifically, each input $\bx_n$ is artificially corrupted (e.g., by adding noise) to produce $\widetildebx_n$, which is then fed into the autoencoder. The model is trained to minimize:
\begin{equation}
J(\btheta) = \frac{1}{2}\sum_{n=1}^{N}\normtwo{e(\widetildebx_n; \btheta) - \bx_n}^2.
\end{equation}
By learning to ``undo" the corruption, the model captures meaningful structural properties of the data. For instance, in image data, it may learn that neighboring pixels are highly correlated, enabling it to correct noisy or missing pixel values.

The most common form of corruption uses additive Gaussian noise. Alternatively, when certain input dimensions are randomly masked out (set to zero), the model is referred to as a \textit{masked autoencoder}. For example, \citet{he2022masked} use deep networks to reconstruct full images from partially observed (masked) inputs.

\begin{figure}[h!]
\centering  
\vspace{-0.35cm} 
\subfigtopskip=2pt 
\subfigbottomskip=2pt 
\subfigcapskip=-5pt 
\subfigure[VAE.]{\label{fig:lvm_VAE}
\includegraphics[width=0.431\linewidth]{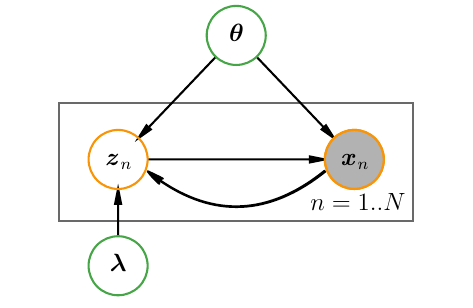}}
\subfigure[Parameter flow in VAE.]{\label{fig:lvm_VAE_flow}
\includegraphics[width=0.431\linewidth]{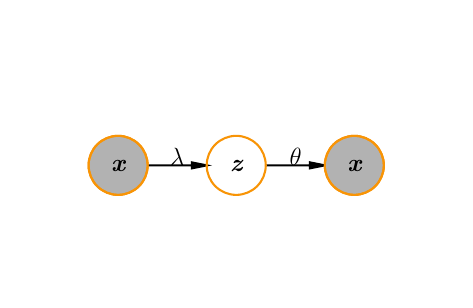}}
\caption{Graphical representation for VAE. Use the variational distribution $q_{\blambda}(\bz\mid \bx)$ to approximate the intractable posterior $p_{\btheta}(\bz\mid \bx)$.}
\label{fig:lvm_VAE_and_flow}
\end{figure}

\begin{figure}[htp]
\centering
\includegraphics[width=0.99\textwidth]{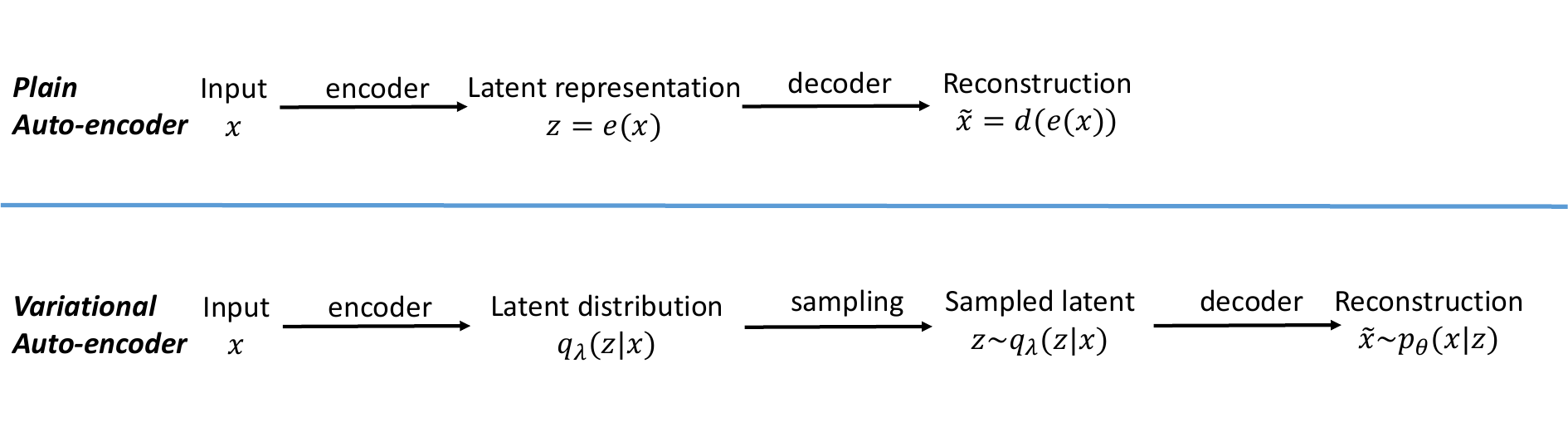}
\caption{Comparison between a standard autoencoder and a variational autoencoder.}
\label{fig:ae_vae}
\end{figure}

\index{Variational autoencoder}
\index{Reparameterization trick}
\subsection{Variational Autoencoder (VAE)}\label{section:vae_pca}
We have already seen that the likelihood function for a latent variable model (see Section~\ref{section:lvm}) is given by
\begin{equation}\label{equation:vae_likelihood}
p(\bx\mid \btheta) = \int p_{\btheta}(\bx\mid\bz) p(\bz) \, d\bz.
\end{equation}
When $p_{\btheta}(\bx\mid\bz)$ is defined by a deep neural network \citep{lecun2015deep, goodfellow2016deep},\footnote{Note that we use $p_{\btheta}(\bx\mid\bz)$ here instead of $p(\bx\mid\bz, \btheta)$ to emphasize that the parameter $\btheta$ is shared by all data points; similarly for the parameter $\blambda$ introduced later.} this integral becomes intractable, as it cannot be evaluated analytically due to the high-dimensional integration over $\bz$. 
The \textit{variational autoencoder} (VAE) \citep{kingma2013auto, rezende2014stochastic, kingma2019introduction} circumvents this issue by optimizing a tractable lower bound on the log-likelihood during training. A graphical illustration of the VAE is shown in Figure~\ref{fig:lvm_VAE}.
The VAE framework rests on three key ideas: 
(i) Use of the \textit{evidence lower-bound (ELBO)}  to approximate the intractable log-likelihood, establishing a close connection to the EM algorithm;
(ii) \textit{Amortized inference:} instead of computing a separate posterior approximation for each data point, a shared encoder network is trained to map inputs $\bx$ to approximate posterior distributions over $\bz$;
(iii) Application of the \textit{reparameterization trick} to enable efficient gradient-based optimization of the encoder parameters.

Consider a generative model where the conditional distribution $p_{\btheta}(\bx\mid\bz)$ over the observed variable $\bx\in\real^D$ is governed by a deep neural network $g(\bz,\btheta)$. For example, $g(\bz,\btheta)$ might output the mean of a Gaussian likelihood. 
Additionally, assume a standard Gaussian prior over the latent variable $\bz\in\real^K$:
\begin{equation}\label{equation:vae_latent_prior}
p(\bz) = \normal(\bz\mid\bzero,\bI).
\end{equation}
To derive the VAE objective, recall from \eqref{equation:elbo_vfe_neg}that, for any \textbf{arbitrary} probability distribution $q(\bz)$ over the latent space, the following identity holds:
\begin{equation}
\ln p(\bx\mid\btheta) = \mathcalF(\btheta) + \KL[q(\bz) \parallel p(\bz\mid\bx,\btheta)],
\label{eq:elbo_kl}
\end{equation}
where $\mathcalF$ is the evidence lower-bound (ELBO) (see also \eqref{equation:elbo_ineq}), defined as
\begin{equation}\label{equation:vae_elbo}
\mathcalF(\btheta) 
= \int q(\bz) \ln \left\{ \frac{p_{\btheta}(\bx\mid\bz) p(\bz)}{q(\bz)} \right\} d\bz,
\end{equation}
and $\KL[P \parallel Q] = \int P(x) \ln \left( \frac{P(x)}{Q(x)} \right) dx \geq 0$ denotes the {KL divergence} between $P$ and $Q$, with equality if and only if $P=Q$.
Since the KL divergence is nonnegative, it follows that
\begin{equation}\label{eq:vae_lower_bound}
	\ln p(\bx\mid\btheta) \geq \mathcalF,
\end{equation}
so $\mathcalF$ is a lower bound on the log-likelihood $\ln p(\bx\mid\btheta)$. 
Although the true log-likelihood is intractable, the ELBO can be approximated using Monte Carlo estimation---a technique known as \textit{stochastic variational inference} (see also Section~\ref{section:mcvi}). Thus, the ELBO serves as a practical surrogate for maximum likelihood training.

\paragrapharrow{Complete-data log-likelihood.}
Now consider a data set $\mathcalX = \{\bx_1, \bx_2, \ldots,\bx_N\}$, assumed to consist of independent and identically distributed samples from the model. The total log-likelihood decomposes as
\begin{equation}\label{equation:vae_log_likelihood_sum}
\ln p(\mathcalX\mid\btheta) = \sum_{n=1}^N \mathcalF_n + \sum_{n=1}^N \KL[q_{\bz_n}(\bz_n\mid \blambda_n) \parallel p(\bz_n\mid\bx_n,\btheta)],
\end{equation}
where the per-data-point ELBO is
\begin{equation}\label{equation:vae_elbo_per_data}
	\mathcalF_n = \int q_{\bz_n}(\bz_n\mid \blambda_n) \ln \left\{ \frac{p_{\btheta}(\bx_n\mid\bz_n) p(\bz_n)}{q_{\bz_n}(\bz_n\mid \blambda_n)} \right\} d\bz_n.
\end{equation}
As in models like PPCA or mixture models (see Problems~\ref{prob:mix_of_gauss}--\ref{prob:mix_of_bern}), we introduce a separate latent variable $\bz_n$ for each observation $\bx_n$.
Consequently, each $\bz_n$ has its own approximate posterior  $q_{\bz_n}(\bz_n\mid \blambda_n)$, which can---in principle---be optimized independently.

\paragrapharrow{Intractability and approximation.}
Because \eqref{equation:vae_log_likelihood_sum} holds for any choice of $q_{\bz_n}(\bz_n\mid \blambda_n)$, we can select the family of distributions that maximizes $\mathcalF_n$ (or equivalently, minimizes the KL divergence to the true posterior, $\KL[q_{\bz_n}(\bz_n\mid \blambda_n) \parallel p(\bz_n\mid\bx_n,\btheta)]$). 
In simpler models such as Gaussian mixtures or probabilistic/Bayesian PCA, the exact posterior $p(\bz_n\mid\bx_n,\btheta)$ can be computed analytically in the E-step of the EM algorithm. Setting $q_{\bz_n}$ equal to this posterior yields zero KL divergence, making the ELBO equal to the true log-likelihood (see Problems~\ref{prob:mix_of_gauss}--\ref{prob:mix_of_bern}, Algorithms~\ref{alg:em_probpca} and~\ref{alg:em_bayespca}).
By Bayes' theorem, the exact posterior is
\begin{equation}\label{equation:vae_bayes_posterior}
p(\bz_n\mid\bx_n,\btheta) = \frac{p_{\btheta}(\bx_n\mid\bz_n) p(\bz_n)}{p(\bx_n\mid\btheta)}.
\end{equation}
While the numerator is easy to evaluate---since $p(\bz_n)$ is a standard Gaussian (Equation~\eqref{equation:vae_latent_prior}) and $p_{\btheta}(\bx_n\mid\bz_n)$ is given by the decoder network---the denominator is precisely the intractable marginal likelihood $p(\bx_n\mid\btheta)$. 
Therefore, we must approximate the posterior.

In principle, one could introduce a separate set of parameters $\blambda_n$ for each $q_{\bz_n}(\cdot\mid \blambda_n)$ and optimize them individually. However, this approach would be computationally prohibitive for large datasets, and the posteriors would need to be re-estimated after every update to $\btheta$. Instead, the VAE adopts a more scalable strategy: it introduces a second neural network---the encoder---that amortizes the cost of inference by sharing parameters across all data points. This encoder maps each input $\bx_n$ to the parameters of its approximate posterior $q_{\bz_n}$, enabling efficient and joint optimization of both generative and inference models.

\subsection{VAE with Amortized Variational Inference}
In the VAE, rather than computing a separate posterior approximation $p(\bz_n\mid\bx_n,\blambda_n)$ for each data point $\bx_n$ separately, we train a single neural network---called the \textit{encoder network} or \textit{recognition model}---to approximate all of these posteriors simultaneously. 
This approach is known as \textit{amortized variational inference} (or simply \textit{amortized inference}; see Section~\ref{section:amorized_VI}). 
The encoder defines a conditional distribution $q_{\blambda}(\bz\mid\bx)$, parameterized by $\blambda$,  that maps any input $\bx$ to a distribution over the latent space; see Figure~\ref{fig:lvm_and_hyp_vb}.
The objective function---the ELBO---now depends on both the generative parameters $\btheta$ and the inference parameters $\blambda$. We maximize this bound jointly with respect to both sets of parameters using gradient-based optimization methods \citep{goodfellow2016deep, lu2022gradient}.

\paragrapharrow{Encoder and decoder phases.}
A VAE thus consists of two neural networks with independent parameters that are trained jointly: 
an \textit{encoder network} that maps an observed data vector $\bx$ to a distribution over the latent variable $\bz$, and a decoder network---the original generative model---that maps a latent vector $\bz$ back to the data space. 
This architecture resembles a standard autoencoder based on neural networks, but with a crucial difference: instead of producing a single point estimate in the latent space, the encoder outputs a probability distribution over $\bz$. As we will see, the encoder effectively learns an approximate probabilistic inverse of the decoder, consistent with Bayes' theorem; see Figure~\ref{fig:ae_vae}.

A common choice for the encoder is a Gaussian distribution with diagonal covariance, where both the mean $\bmu_{\blambda}$ and  variance $\bsigma_{\blambda}^2$ are outputs of a neural network that takes  $\bx_n$ as input:
\begin{equation}\label{equation:vae_encoder_gaussian}
\begin{aligned}
q_{\blambda}(\bz_n\mid \bx_n ) 
&= \normal(\bz_n \mid \bmu_{\blambda}(\bx_n), \diag(\bsigma_{\blambda}^2(\bx_n)))\\
&= \prod_{k=1}^K \normal\left(z_{nk} \mid 
(\bmu_\blambda(\bx_n))_k, 
(\bsigma_\blambda^2(\bx_n))_k
\right),
\quad 
n=1,2,\ldots,N,
\end{aligned}
\end{equation}
Note that the mean components can take any real value, so the corresponding output units typically use a linear activation function. In contrast, the variances must be nonnegative; therefore, the network usually outputs log-variances or applies an exponential activation (e.g., $\exp(\cdot)$) to ensure positivity.

Similarly, the decoder is often modeled as a Gaussian distribution with diagonal covariance:
$$
p_{\btheta}(\bx_n\mid \bz_n)
=\normal(\bx\mid \bmu_{\btheta}(\bz_n), \diag(\bsigma^2_{\btheta}(\bz_n))), 
\quad 
n=1,2,\ldots,N,
$$
where $\bmu_{\btheta}$ and $\bsigma_{\btheta}^2$ are given by neural network transformations of $\bz_n$ for all $n=1,2,\ldots,N$. 
In practice, the decoder variance $\bsigma^2_{\btheta}$ is often fixed (e.g., set to 1), and only the mean $\bmu_{\btheta}$ is learned.
The reconstructed output $\widetildebx_n = \bmu_{\btheta}(\bz_n)$ is then compared to the original input $\bx_n$.
For binary data (e.g., binarized MNIST), a Bernoulli likelihood is commonly used instead of a Gaussian.

\paragrapharrow{Comparison with EM algorithms.}
The ultimate goal is to estimate the generative parameters $\btheta$ via maximum likelihood::
$
\mathop{\max}_{\btheta} \sum_{n=1}^{N} \ln p_{\btheta}(\bx_n  ).
$
In standard VI, this is approached using a constrained EM algorithm (Section~\ref{section:em_constrained}), 
which iteratively maximizes the ELBO:
$$
\mathop{\max}_{\btheta,\{\blambda_n\}} 
\sum_{n=1}^{N} \Exp_{q_{\bz_n}(\bz_n\mid\blambda_n)} 
\left[ \ln \frac{p_{\btheta}(\bz_n, \bx_n  )}{q_{\bz_n}(\bz_n\mid\blambda_n)} \right].
$$ 
The VAE, however, adopts a shared, data-dependent posterior approximation: $q_{\blambda}(\bz\mid \bx )$ with variational parameter $\blambda$:
$
q_{\blambda}(\bz_n\mid \bx_n ) = \normal(\bz_n \mid \bmu_{\blambda}(\bx_n), \diag(\bsigma_{\blambda}^2(\bx_n))),
$
where $\bmu_{\blambda}$ and $\bsigma_{\blambda}^2$ are parameterized by a neural network.
This leads to an algorithm that closely resembles training a standard autoencoder---but with a stochastic encoding step that injects learnable noise (see Figure~\ref{fig:lvm_VAE_flow}), hence the name variational autoencoder.
The VAE objective is to maximize the following ELBO jointly over $\btheta$ and $\blambda$:
\begin{equation}\label{equation:vae_elbo_vae}
\begin{aligned}
&\mathop{\max}_{\btheta,\blambda} 
\left\{\sum_{n=1}^{N} \Exp_{q_{\blambda}(\bz_n\mid\bx_n)} 
\left[ \ln \frac{p_{\btheta}(\bz_n, \bx_n  )}{q_{\blambda}(\bz_n\mid\bx_n)} \right]
= \sum_{n=1}^{N}\int q_{\blambda}(\bz_n\mid\bx_n) \ln \left\{ \frac{p_{\btheta}(\bx_n\mid\bz_n) p(\bz_n)}{q_{\blambda}(\bz_n\mid\bx_n)} \right\} d\bz_n 
\right\}
\\
&=
\mathop{\max}_{\btheta,\blambda} 
\underbrace{\sum_{n=1}^{N} \Exp_{q_{\blambda}(\bz_n\mid\bx_n)} 
	\left[ \ln {p_{\btheta}(\bx_n \mid\bz_n  )} \right]}_{\text{data fit}}
\underbrace{-\sum_{n=1}^{N}\KL\left[q_{\blambda}(\bz_n\mid\bx_n) \parallel p(\bz_n  ) \right]}_{\text{regularization}},
\end{aligned}
\end{equation}
where the first term can be
interpreted as (negative) reconstruction error (the reconstruction of the observed data $\bx_n$ from the latent space), and the second KL term can be viewed as a regularization, measured by the KL divergence between the approximate posterior and the prior $p(\bz_n)$. 
By convention, the prior is chosen to be a standard multivariate Gaussian:
$$
p(\bz_n)=\normal(\bzero, \bI), 
\quad n=1,2,\ldots,N. 
$$
This choice is convenient because (i) it is easy to sample from, and (ii) the KL divergence between two Gaussians has a closed-form expression (see Problem~\ref{problem:kl_vae}). 
Moreover, it encourages the encoder to produce latent representations that are smooth, continuous, and centered around the origin with unit variance, facilitating meaningful interpolation and sampling in the latent space.
In practice, many researchers scale the KL regularization term by a hyper-parameter $\beta$ to control the trade-off between reconstruction fidelity and latent structure \citep{hoffman2016elbo}.

The ELBO is optimized using gradient-based methods, typically stochastic gradient descent (SGD) or its variants, applied to mini-batches of data \citep{goodfellow2016deep, lu2022gradient}. Although $\btheta$ and 
$\blambda$ are updated jointly, one can conceptually view the process as alternating between improving the encoder (inference) and the decoder (generative model)---analogous to the E/M-steps of the EM algorithm.

Another key distinction from classical EM is that, for a fixed $\btheta$, optimizing $\blambda$ does not generally drive the KL divergence to zero. This is because the encoder network---despite its flexibility---cannot perfectly represent the true posterior 
$p(\bz_n\mid \bx_n, \btheta)$. Several factors contribute to this residual gap:
\begin{enumerate}[label=(\roman*)]
\item The true posterior may not be Gaussian or factorized (i.e., it may have complex dependencies across latent dimensions).
\item Even deep neural networks have finite capacity and cannot represent arbitrary distributions exactly.
\item The optimization itself is approximate due to stochastic gradients, limited iterations, and local optima---limitations shared with constrained EM methods (see Section~\ref{section:em_constrained}).
\end{enumerate}
As a result, the ELBO remains a strict lower bound on the true log-likelihood, as illustrated in Figure~\ref{fig:ELBO_decom}.

\subsection{VAE with the Reparameterization Trick}
Unfortunately, even with the decomposition in \eqref{equation:vae_elbo_vae}, the ELBO remains intractable to compute exactly because the first term involves an integral over the latent variable $\bz_n$, and the integrand---due to the nonlinear decoder network---has a complex dependence on $\bz_n$. For a single data point $\bx_n$, the contribution to the ELBO can be written as (see \eqref{equation:vae_elbo_vae}):
\begin{equation}\label{equation:vae_elbo_decomposed}
\mathcalF_n(\btheta,\blambda) 
=
\Exp_{q_{\blambda}(\bz_n\mid\bx_n)} \left[ \ln {p_{\btheta}(\bx_n \mid\bz_n  )} \right] 
- \KL[q_{\blambda}(\bz_n\mid\bx_n) \parallel p(\bz_n)].
\end{equation}
The second term is a KL divergence between two Gaussian distributions and admits a closed-form expression (see Problem~\ref{problem:kl_vae}):
\begin{equation}
\KL[q_{\blambda}(\bz_n\mid\bx_n) \parallel p(\bz_n)] = \frac{1}{2} \sum_{k=1}^K \left\{  (\bmu_\blambda(\bx_n))_k + (\bsigma_\blambda^2(\bx_n))_k - \ln (\bsigma_\blambda^2(\bx_n))_k  - 1\right\}.
\end{equation}
For the first term in \eqref{equation:vae_elbo_decomposed}, a natural approach is to approximate the expectation using a Monte Carlo estimator:
\begin{equation}\label{equation:vae_monte_carlo_approx}
\int q_{\blambda}(\bz_n\mid\bx_n) \ln p_{\btheta}(\bx_n\mid\bz_n) \, d\bz_n \simeq \frac{1}{S} \sum_{s=1}^{S} \ln p_{\btheta}(\bx_n\mid\bz_n^{(s)}),
\quad \bz_n^{(s)}\sim q_{\blambda}(\bz_n\mid\bx_n).
\end{equation}
While this estimate is straightforward to differentiate with respect to $\btheta$, 
computing gradients with respect to $\blambda$ is problematic: the samples $\bz_n^{(s)}$  depend on $\blambda$ through the encoder distribution, yet once drawn, they are treated as fixed constants. Consequently, standard backpropagation cannot propagate gradients through the sampling operation.
Conceptually, fixing $\bz_n$ to sampled values blocks the flow of gradient information to the encoder parameters $\blambda$; that is, when the ELBO is estimated using fixed samples from $q_{\blambda}(\bz_n\mid \bx_n)$, the error signal cannot be backpropagated through the stochastic sampling step to update the encoder network.

This issue is resolved by the \textit{reparameterization trick} (see Section~\ref{section:mcvi}), which rewrites the sampling procedure so that randomness is decoupled from the parameters. Specifically, if $\bepsilon\sim\normal(\bzero, \bI)$, then
\begin{equation}\label{eq:reparam_trick}
\bz = \diag(\bsigma)  \bepsilon + \bmu
\end{equation}
follows a Gaussian distribution $\normal(\bmu, \diag(\bsigma^2))$  (see Lemma~\ref{lemma:affine_mult_gauss}). 
Applying this to our encoder, we replace direct sampling from $q_{\blambda}(\bz_n \mid \bx_n )$ with:
\begin{equation}\label{equation:vae_sample_reparam}
\bz_n^{(s)}\sim q_{\blambda}(\bz_n \mid \bx_n )
\quad\implies \quad
\bepsilon_n^{(s)} \sim \normal(\bzero, \bI), \bz_n^{(s)} = \bmu_{\blambda}(\bx_n)+ \bsigma_{\blambda}(\bx_n)\hadaprod \bepsilon_n^{(s)},
\end{equation}
where $\hadaprod$ denotes the \textit{Hadamard product}, and $s = 1,2,\ldots,S$ indexes the samples.
This reformulation makes the dependence on $\blambda$ explicit and differentiable, enabling gradient-based optimization via automatic differentiation.

Although the reparameterization trick applies primarily to continuous latent variables, alternative gradient estimators exist for discrete cases (e.g., the REINFORCE estimator; \citealp{williams1992simple}). However, these typically suffer from high variance. Thus, the reparameterization trick also serves as an effective variance reduction technique (see Section~\ref{section:mcvi}).

Under our modeling assumptions, the full VAE objective  (averaged over a mini-batch $\sT\subset\{1,2,\ldots,N\}$) becomes:
\begin{equation}\label{equation:vae_loss}
\mathcalF = \frac{N}{\abs{\sT}}\sum_{n\in\sT} \left\{ \frac{1}{2} \sum_{k=1}^K \left(1 + \ln \sigma_{nk}^2 - \mu_{nk}^2 - \sigma_{nk}^2 \right) + \frac{1}{S} \sum_{s=1}^S \ln p_{\btheta}(\bx_n\mid\bz_n^{(s)}) \right\}.
\end{equation}
where the latent sample $\bz_n^{(s)}$ is constructed as  $z_{nk}^{(s)} = \sigma_{nk} \bepsilon_n^{(s)} + \mu_{nk}$, in which $\mu_{nk} = (\mu_{\blambda}(\bx_n))_k$ and $\sigma_{nk} = (\sigma_{\blambda}(\bx_n))_k$.
In practice, the number of Monte Carlo samples per data point is often set to $S=1$. Although this yields a noisy estimate of the ELBO, the noise is compatible with stochastic gradient optimization and generally leads to faster and more efficient training.

VAE training proceeds as follows: for each data point in a mini-batch, (i) forward-propagate through the encoder to obtain $\bmu_{\blambda}(\bx_n)$ and $\bsigma_{\blambda}(\bx_n)$, (ii) draw $\bepsilon_n^{(s)}$ and compute $\bz_n^{(s)}$ via reparameterization, (iii) pass $\bz_n^{(s)}$ through the decoder to evaluate the reconstruction log-likelihood or ELBO \eqref{equation:vae_loss}, and (iv) compute gradients of the ELBO with respect to both $\btheta$ and  
$\blambda$ using automatic differentiation. 
The complete procedure is summarized in Algorithm~\ref{alg:amo_vae}.

\begin{algorithm}[h] 
\caption{Variational Autoencoder (VAE)}
\label{alg:amo_vae}
\begin{algorithmic}[1] 
\Require Observed data points $\mathcalX=\{\bx_1, \bx_2, \ldots, \bx_N\}$;
\State \textbf{initialize:} $\btheta^\topone,  \blambda^\topone$; 
\State \textbf{parameters:} Step size $\eta$, MC sample number $S$;
\State Choose the maximal number of iterations $C$;
\State $t=0$; \Comment{Count for the number of iterations}
\While{$t<C$} 
\State $\bepsilon^{(s)}\sim\normal(\bzero, \bI)$ for all $s=1,2,\ldots,S$;  
\Comment{(VAE$_1$)}
\State $\bz_n^{(s)} \leftarrow  \bmu_{\blambda}(\bx_n)+ \bsigma_{\blambda}(\bx_n)\hadaprod \bepsilon_n^{(s)}$  for all $n=1,2,\ldots,N$, $s=1,2,\ldots,S$;    
\Comment{(VAE$_2$)}
\State $\mathcalF \leftarrow \frac{N}{\abs{\sT}}\sum_{n\in\sT} \left\{ \frac{1}{2} \sum_{k} \left(1 + \ln \sigma_{nk}^2 - \mu_{nk}^2 - \sigma_{nk}^2 \right) + \frac{1}{S} \sum_{s} \ln p_{\btheta}(\bx_n\mid\bz_n^{(s)}) \right\}$;
\Comment{(VAE$_3$)}
\State $\btheta \leftarrow \btheta +\eta \nabla_{\btheta}\mathcalF $; \Comment{(VAE$_4$)}
\State $\blambda \leftarrow \blambda +\eta \nabla_{\blambda}\mathcalF $; \Comment{(VAE$_5$)}
\EndWhile
\State Output $\btheta$, $\blambda$;
\end{algorithmic} 
\end{algorithm}

\paragrapharrow{Evaluation and generative process.}
After training, to evaluate how well the model represents a new test point $\widetildebx$, we use the ELBO 
$\mathcalF$ as a tractable lower bound on the log-likelihood. For a tighter estimate, it is preferable to sample from the approximate posterior $q(\bz\mid\widetildebx,\blambda)$ rather than from the prior $p(\bz)$, as the former concentrates probability mass in regions relevant to $\widetildebx$.

Once the model is trained and evaluated, the encoder network is discarded and new data points are generated by sampling from the prior $p(\bz)$ and forward-propagating through the decoder network to obtain samples in the data space: $p_{\btheta}(\bx\mid \bz)$.
This contrasts with standard autoencoders, where the latent code is a deterministic function of the input. In a VAE, the encoder outputs a distribution over latent codes, and actual codes are obtained by sampling---making the model inherently probabilistic (see Figure~\ref{fig:ae_vae}).

This stochastic framework makes VAEs particularly effective for generative tasks, such as image and sequence synthesis \citep{kingma2013auto, rezende2014stochastic}. Because the latent space is regularized to be smooth and continuous, \textbf{interpolating} between two latent vectors typically yields meaningful transitions in the data space.

However, VAEs are known to sometimes produce blurrier images compared to alternatives like GANs \citep{goodfellow2020generative}, as the reconstruction objective encourages averaging over plausible outputs to minimize expected error.

\paragrapharrow{Conditional VAE.}
Several variants of the VAE exist. For image data, encoders typically use convolutional layers, while decoders employ transposed convolutions. In a \textit{conditional VAE}, both the encoder and decoder receive an additional conditioning variable $\bc$ (e.g., a class label). The prior over the latent variable can either remain the standard $p(\bz)$ or be extended to a conditional prior $p(\bz\mid\bc)$, which may be parameterized by a separate neural network. Training proceeds in the same manner as in the standard VAE.
During generation, the user can specify a particular value of $\bc$ to guide the model toward producing more relevant or targeted outputs.

\paragrapharrow{General framework.}
Instead of using the standard Gaussian assumption for $p(\bz)$, we can also analyze for general distributions. 
We still have
$$
\bz_n\sim q_{\blambda}(\bz_n \mid \bx_n )
\implies 
\bepsilon_n \sim \normal(\bzero, \bI), \bz_n = \bmu_{\blambda}(\bx_n)+ \bsigma_{\blambda}(\bx_n)\hadaprod \bepsilon_n,
$$
where $\hadaprod$ represents the {Hadamard product}. And the corresponding ELBO objective in \eqref{equation:vae_elbo} becomes
$$
\mathop{\max}_{\btheta,\blambda} 
\sum_{n=1}^{N} \Exp_{p(\bepsilon)} 
\left[ \ln \frac{p_{\btheta}\big(\left\{\bmu_{\blambda}(\bx_n)+ \bsigma_{\blambda}(\bx_n)\hadaprod \bepsilon_n\right\}, \bx_n \big)}{q_{\blambda}\big(\left\{\bmu_{\blambda}(\bx_n)+ \bsigma_{\blambda}(\bx_n)\hadaprod \bepsilon_n\right\}\mid\bx_n\big)} \right].
$$
A Monte Carlo approximation yields:
$$
\sum_{n=1}^{N} 
\frac{1}{S}
\sum_{s=1}^{S}
\left[ \ln \frac{p_{\btheta}\big(\left\{\bmu_{\blambda}(\bx_n)+ \bsigma_{\blambda}(\bx_n)\hadaprod \bepsilon_n^{(s)}\right\}, \bx_n \big)}{q_{\blambda}\big(\left\{\bmu_{\blambda}(\bx_n)+ \bsigma_{\blambda}(\bx_n)\hadaprod \bepsilon_n^{(s)}\right\}\mid\bx_n\big)} \right],
\gap
\bepsilon_n^{(s)}\sim\normal(\bzero, \bI),
$$
enabling end-to-end gradient-based learning via backpropagation.

\paragrapharrow{Other issues.}
The KL term in the ELBO \eqref{equation:vae_elbo_decomposed} regularizes the encoder to align its output with the prior $p(\bz)$, ensuring that samples from the prior yield realistic data when passed through the decoder. However, two failure modes can occur:
\begin{enumerate}[label=(\roman*)]
\item \textit{Posterior collapse.} The encoder ignores the input and outputs a distribution close to the prior, i.e., $q_{\blambda}(\bz\mid \bx)\approx p(\bz)$. This renders the latent code uninformative. Symptoms include poor reconstructions (blurry outputs) and a KL divergence near zero.
\item \textit{Poor generative quality.} Reconstructions are accurate, but samples generated from  $p(\bz)$ are unrealistic. Here, the encoder fits the data tightly, causing $q_{\blambda}(\bz\mid \bx)$ to diverge significantly from  $p(\bz)$, so prior samples fall in low-density regions of the latent space used during training.
\end{enumerate}

Both issues can be mitigated by introducing a weighting coefficient $\bbeta>0$ on the KL term:
\begin{equation}\label{equation:vae_beta_elbo_decomposed}
\mathcalF_n^{\beta}
=
\Exp_{q_{\blambda}(\bz_n\mid\bx_n)} \left[ \ln {p_{\btheta}(\bx_n \mid\bz_n  )} \right] 
- \beta\cdot \KL[q_{\blambda}(\bz_n\mid\bx_n) \parallel p(\bz_n)].
\end{equation}
This yields the \textit{$\beta$-VAE} \citep{hoffman2016elbo, higgins2017beta}. If reconstructions are poor, increase $\beta$; if generated samples are poor, decrease  $\beta$. Often, $\beta$ is scheduled to start small and increase gradually during training (\textit{KL annealing}).

\begin{problemset}


\item \textbf{Shared SVD from identical scatter matrices.}  Consider two data matrices $\bA_1$ and $\bA_2$ that have identical scatter matrices $\bA_1^\top \bA_1 = \bA_2^\top \bA_2$, but are otherwise distinct. Show that both $\bA_1$ and $\bA_2$ admit a partially shared SVD of the form: $\bA_1 = \bU_1 \bSigma \bV^\top$ and $\bA_2 = \bU_2 \bSigma \bV^\top$. Use this result to show that $\bA_2 = \bQ_{12} \bA_1$, where $\bQ_{12}$ is an orthogonal matrix.

\item \label{theorem:polar-decomposition} \textbf{Polar decomposition.}
Let $\bA\in\real^{M\times N}$. Show that  $\bA$ can be factored as 
\begin{itemize}
\item \textbf{Case $M>N$: left polar decomposition.} $\bA=\bQ_l\bS_l$, where $\bS_l^2=\bA^\top\bA$ is positive semidefinite and is \textbf{uniquely} determined. 
The factor $\bQ_l$ has orthonormal columns and is \textbf{uniquely} determined if $\rank(\bA)=N$.
\item \textbf{Case $M<N$: right polar decomposition.} $\bA=\bS_r\bQ_r$, where $\bS_r^2=\bA\bA^\top$ is positive semidefinite and is \textbf{uniquely} determined.
The factor $\bQ_r$ has orthonormal rows and is \textbf{uniquely} determined if $\rank(\bA)=M$.
\item \textbf{Case $M=N$: left/right polar decomposition.} $\bA=\bQ\bS_l=\bS_r\bQ$, where $\bS_l^2=\bA^\top\bA$ and $\bS_r^2=\bA\bA^\top$ are positive semidefinite and are \textbf{uniquely} determined. The factor $\bQ$ is orthonoal, and it is the same for both the left and right polar decompositions. $\bQ$ is \textbf{uniquely} determined if $\bA$ is nonsingular ($\rank(\bA)=N$).
\end{itemize}
\noindent
\textit{Hint: Use SVD.}

\item \label{prob:ppca_gengauss} Suppose we replace the zero-mean, unit-covariance latent distribution in Equation \eqref{equation:probpca_pz}
of the PPCA model with  a general Gaussian distribution  $\normal(\bz\mid \bmu_z, \bSigma_z)$. By redefining the model parameters appropriately, show that the resulting marginal distribution  $p(\bx)$ over the observed variables remains unchanged for any valid choice of
$\bmu_z$ and $\bSigma_z$.

\item Verify that maximizing the log-likelihood in \eqref{equation:probpca_jointlnlike} for the PPCA model with respect to the parameter $\bmu$, $\bW$, and $\sigma^2$ yields the maximum likelihood estimates given in Equations~\eqref{equation:probpca_bmle}, \eqref{equation:probpca_wmle}, and \eqref{equation:probpca_sigmamle}, respectively.

\item Verify that minimizing the objective function in  \eqref{equation:bayespca_emw_loss} during the EM update for Bayesian PCA leads to the parameter update for $\bW$ shown in \eqref{equation:bayespca_emw}.
\textit{Hint: Denote  $  \bA = \sum_{n=1}^N \widehatbSigma_n  $, and  $  \bB = \sum_{n=1}^N (\bx_n - \bmu) \widehatbz_n^\top  $, express the objective function in \eqref{equation:bayespca_emw_loss} with the trace of $\bW$, then set its derivative with respect to $\bW$ to zero.
}

\index{Truncated SVD}
\index{Eckart-Young-Mirsky theorem}
\index{Low-rank approximation}
\item \label{theorem:young-theorem_frob}\textbf{Eckart-Young-Mirsky theorem and truncated SVD (TSVD) \citep{stewart1993early, lu2021numerical}.}
Suppose we wish to approximate a rank-$R$ matrix $\bA\in \real^{M\times N}$ by a lower-rank matrix  $\bB$ of rank $K$ ($K<R$), measured in the Frobenius norm (Definition~\ref{definition:frobernius-in-svd}):
$$
\bB = \mathop{\arg\min}_{\rank(\bB)\leq K} \, \normf{\bA - \bB}.
$$
Let $\bA_K$ be the \textit{truncated SVD} (TSVD) of $\bA$ with the top $K$ terms, i.e., $\bA_K = \sum_{i=1}^{K} \sigma_i\bu_i\bv_i^\top$ from the SVD of $\bA=\sum_{i=1}^{R} \sigma_i\bu_i\bv_i^\top$ by zeroing out the $R-K$ trailing singular values of $\bA$. 
Show that $\bA_K$ is the optimal rank-$K$ approximation to $\bA$ in terms of the Frobenius norm, satisfying $\normf{\bA-\bA_K}^2 = \sum_{i\geq K+1}\sigma_i^2$.
What is the corresponding optimal approximation in the spectral norm (Definition~\ref{definition:spectral_norm})?

\item \label{prob:projection-from-matrix}\textbf{Projection matrix from a set of vectors \citep{lu2021numerical}.}
Let $\ba_1, \ba_2, \ldots, \ba_N \in\real^M$ be linearly independent vectors that span a subspace $\mathcalV$: i.e., $\cspace([\ba_1, \ba_2, \ldots, \ba_N]) = \mathcalV$, where $M\geq N$. 
Show that the orthogonal projection matrix onto $\mathcalV$ is given by
$$
\bH = \bA(\bA^\top\bA)^{-1}\bA^\top,
$$
where $\bA=[\ba_1, \ba_2, \ldots, \ba_N]\in \real^{M\times N}$.

\item \label{prob:rayleigh_v2} \textbf{Rayleigh--Ritz theorem \citep{lu2021numerical}.}
Let $\bA\in \real^{N\times N}$ be a symmetric matrix with eigenvalues $\lambda_1\leq \lambda_2\leq \ldots \leq \lambda_N$ and corresponding mutually orthonormal eigenvectors $\bq_1, \bq_2, \ldots, \bq_N$. 
For any unit vector $\bx\in \real^N$ (i.e., $\normtwo{\bx}=1$), show that the quadratic form satisfies
$$
\lambda_N \geq 
\bx^\top\bA\bx
\geq \lambda_1, 
\gap \normtwo{\bx}=1.
$$

\item \label{prob:rayleigh_v3} \textbf{Rayleigh--Ritz theorem.}
Use the Rayleigh--Ritz theorem to show that for a symmetric matrix $\bA\in\real^{N\times N}$ with eigenvalues $\lambda_1\leq \lambda_2\leq \ldots \leq \lambda_N$ and orthonormal eigenvectors $\bq_1, \bq_2, \ldots, \bq_N$, the \textit{Rayleigh quotient}
$r(\bx) = \frac{\bx^\top\bA\bx}{\bx^\top\bx}$ satisfies
$$
\begin{aligned}
&\mathop{\max}_{\bx\neq \bzero}
\frac{\bx^\top\bA\bx}{\bx^\top\bx}
=
\lambda_N
\qquad\text{and}\qquad
\mathop{\min}_{\bx\neq \bzero}
\frac{\bx^\top\bA\bx}{\bx^\top\bx}
=
\lambda_1,
\end{aligned}
$$
with the maximum value $\lambda_N$ attained at $\bx=\alpha \bq_N$ (for nonzero $\alpha$), and the minimum value $\lambda_1$ attained at $\bx=\alpha\bq_1$. Or if $\mathcalV$ is the subspace spanned by  $\{\bq_p, \bq_{p+1}, \ldots, \bq_q\}$, show that  
$$
\begin{aligned}
&\mathop{\max}_{\bx\neq \bzero, \bx\in\mathcalV}
\frac{\bx^\top\bA\bx}{\bx^\top\bx}
=
\lambda_q
\qquad\text{and}\qquad
\mathop{\min}_{\bx\neq \bzero, \bx\in\mathcalV}
\frac{\bx^\top\bA\bx}{\bx^\top\bx}
=
\lambda_p,
\end{aligned}
$$
with the maximum value $\lambda_q$ achieved at $\bx=\alpha\bq_q$, and the minimum value $\lambda_p$ achieved at $\bx=\alpha\bq_p$.
Let $\bx\triangleq\be_n$ for $n\in\{1,2,\ldots,N\}$. This implies
$$
\textbf{(Symmetric $\bA$): }\qquad \lambda_{\min}(\bA) \leq d_{\min}(\bA) \leq d_{\max}(\bA) \leq \lambda_{\max}(\bA),
$$
where $d_{\min}(\bA)$ and $d_{\max}(\bA)$ denote the smallest and largest diagonal entries of $\bA$, respectively.

\item \label{prob:rayleigh_v4} Use the Rayleigh--Ritz theorem to prove that the solution to the optimization problem in Equation \eqref{equation:or_pca_all} is given by the top-$K$ eigenvectors of $\bX_c^\top\bX_c$, where $\bX_c$ is the centered data matrix.

\item Verify the EM algorithm for the mixture of PPCA models given in \eqref{equation:em_mix_ppca}. 
Additionally, derive the EM algorithm for a mixture of Bayesian PCA models.
\textit{Consider the  EM update in \eqref{equation:bayespca_mstep}.}

\item  \label{prob:factor_ana} \textbf{Factor analysis \citep{everett2013introduction, basilevsky2009statistical}.}
Factor analysis is another linear Gaussian latent variable model closely related to PPCA. The key difference lies in the conditional distribution of the observed variable $\bx$ given the latent variable $\bz$: instead of an isotropic noise covariance, it uses a diagonal covariance matrix:
\begin{equation}
	p(\bx\mid\bz) = \normal(\bx \mid \bW\bz + \bmu, \bD),
\end{equation}
where $\bD\in\real^{D\times D}$ is a  diagonal matrix, and $\bW\in\real^{D\times K}$.
In this model, the observed covariance structure is decomposed into two parts: $\bW\bW^\top$, which captures correlations among variables (the columns of $\bW$ are called \textit{factor loadings}); and $\bD$, whose diagonal entries (called \textit{uniquenesses}) represent independent noise variances for each observed dimension.
Show that the marginal distribution of  $\bx$ is
\begin{equation}
p(\bx) = \normal(\bx \mid \bmu, \bM),
\qquad \text{with }\bM \triangleq  \bD+\bW\bW^\top .
\end{equation}
Like PPCA, this model is invariant under rotations in the latent space.
Given observations $\mathcalX=\{\bx_1, \bx_2, \ldots,\bx_N\}$, derive the EM updates for the factor analysis model.
Specifically, at iteration $t$, show that the E-step computes:
\begin{subequations}
\begin{align}
\widehatbz_n^\toptzero &= \bN^{-1} \bW^\top \bD^{-1} (\bx_n - \widebarbx);\\
\widehatbSigma_n^\toptzero &= \bN^{-1} + 	\widehatbz_n^\toptzero 	\widehatbz_n^\toptzeroTOP,
\end{align}
\end{subequations}
where now $ \bN \triangleq \bI + \bW^\top \bD^{-1} \bW$.
Note that $\bN^{-1}$ involves only $K \times K$ matrix inversions (since $K\ll D$ in typical applications), while $\bD^{-1}$ is trivial to compute because $\bD$ is diagonal.
Similarly, show that the M-step updates are:
\begin{subequations}
\begin{align}
\bW^\toptone &\leftarrow \left[ \sum_{n=1}^{N} (\bx_n - \widebarbx) 	\widehatbz_n^\toptzeroTOP \right] \left[ \sum_{n=1}^{N} \widehatbSigma_n^\toptzero \right]^{-1}; \\
\bD^\toptone &\leftarrow \diag\left\{ \bS - \bW^\toptone \frac{1}{N} \sum_{n=1}^{N} 	\widehatbz_n^\toptzero (\bx_n - \widebarbx)^\top \right\}.
\end{align}
\end{subequations}
\textit{Hint:  Leverage results from the linear Gaussian model (Exercise~\ref{exercise:linear_gauss_model}) and the EM update for Bayesian PCA in Equation \eqref{equation:bayespca_mstep}. }

\end{problemset}

\chapter{Bayesian Real Matrix Factorization}\label{chapter:bayes_rmf}
\begingroup
\hypersetup{
	linkcolor=structurecolor,
	linktoc=page,  
}
\minitoc \newpage
\endgroup

\section{Introduction}\label{section:bmf_real_intro}
\lettrine{\color{caligraphcolor}T}
The explosion of data driven by advances in sensor technology and computer hardware has introduced new challenges in data analysis.
Large datasets often contain noise and other distortions, necessitating preprocessing before deductive scientific methods can be effectively applied.
For instance, signals captured by antenna arrays are frequently contaminated by noise and other degradations.
To analyze such data effectively, it must be reconstructed or represented in a way that reduces inaccuracies while preserving essential structural properties.

Moreover, data collected from complex systems often arises from multiple interrelated variables acting simultaneously. When these variables are poorly defined or unobserved, the information in the raw data can become redundant or ambiguous.
By constructing a reduced-order model, we can approximate the behavior of the original system with high fidelity.

As previously noted, a common strategy for denoising, dimensionality reduction, and feasibility-preserving reconstruction is to replace the original data with a lower-dimensional representation obtained via subspace approximation.
Consequently, low-rank approximations---or low-rank matrix decompositions---play a pivotal role across a wide range of applications.
Low-rank matrix decomposition is a powerful technique in machine learning and data mining that expresses a given matrix as the product of two or more matrices of lower dimensionality.
This approach captures the essential structure of the data while filtering out noise and redundancies.
Well-known methods for low-rank decomposition include singular value decomposition (SVD), principal component analysis (PCA), and multiplicative-update nonnegative matrix factorization (NMF).

Bayesian low-rank decomposition extends this framework by incorporating Bayesian modeling principles.
It treats the observed data as arising from a low-rank matrix that itself is generated from prior probability distributions.
This formulation enables the integration of prior knowledge and explicit modeling of uncertainty in the factor matrices.
As a result, Bayesian methods often yield more robust, stable, and interpretable results compared to purely data-driven approaches.
Furthermore, they provide probabilistic quantification of uncertainty, offering meaningful confidence measures for the estimated factors.
These advantages help mitigate overfitting and make Bayesian low-rank decomposition particularly effective for both predictive and explanatory modeling.
Consider an observed dataset represented as an $M\times N$ real-valued matrix $\bA$, where rows correspond to observations and columns to variables of interest.
Following Chapter~\ref{chapter:als}, the \textit{real matrix factorization (RMF)} problem---a canonical bilinear decomposition---can be expressed as
$$
\bA=\bW\bZ+\bE,
$$  
where $\bA=[\ba_1, \ba_2, \ldots, \ba_N]\in \real^{M\times N}$ is approximately factorized into an $M\times K$  matrix $\bW\in \real^{M\times K}$ and a $K\times N$ matrix  $\bZ\in \real^{K\times N}$. The data set $\bA$ need not be complete; missing entries can be indicated by a binary mask matrix $\bM\in \real^{M\times N}$, where $m_{mn}=1$ denotes an observed entry and $m_{mn}=0$ a missing one.
Matrices $\bW$ and $\bZ$ represent latent explanatory factors: their product provides a predictor for the entries of $\bA$. 
When some entries of $\bA$  are missing, $\bW$ and $\bZ$ can be used to impute those values.
If either $\bW$ or $\bZ$ is known, the problem reduces to a standard regression task.

The factorization of the original data matrix $\bA$ is achieved by finding two such  real matrices:
one representing the \textit{basis} (or \textit{dictionary}) components and the other representing the \textit{activations} (or \textit{coefficients}).
Let $\bz_n$ denote the $n$-th column of $\bZ$. Then the matrix multiplication of $\bW\bZ$ can be implemented as computing each column vector $\ba_n$  of $\bA$ as a linear combination of the columns of 
$\bW$, using coefficients provided by $\bz_n$:
$$
\ba_n = \bW\bz_n.
$$

In the Netflix context, the entry $a_{mn}$---the $(m,n)$-th element of $\bA$---represents the rating given by user $n$ to movie $m$ (with higher values indicating stronger preference). In this setting, $\bw_m$ (the $m$-th row of $\bW$) can represent the latent features of movie $m$, and  $\bz_n$ (the $n$-th column of $\bZ$) encodes the preferences of user $n$ (see Section~\ref{section:als-vector-product}).

To simplify the problem, assume initially that there are no missing entries.
We seek to project each data vector $\ba_n$ into a lower-dimensional space $\bz_n \in \real^K$, with $K<M$, such that the reconstruction error, measured by the Frobenius norm, is minimized (assuming $K$ is known):
\begin{equation}\label{equation:als-per-example-loss_real}
	\mathop{\min}_{\bW,\bZ}  \sum_{n=1}^N \sum_{m=1}^{M} \left(a_{mn} - \bw_m^\top\bz_n\right)^2,
\end{equation}
where $\bW=[\bw_1^\top; \bw_2^\top; \ldots; \bw_M^\top]\in \real^{M\times K}$ and $\bZ=[\bz_1, \bz_2, \ldots, \bz_N] \in \real^{K\times N}$ contain the vectors $\bw_m$ and $\bz_n$ as \textbf{rows and columns}, respectively \footnote{Note that in some formulations, $\bZ$ is taken as an $N\times K$ matrix so that $\bA=\bW\bZ\textcolor{mylightbluetext}{^\top} +\bE$.}. 
The loss function  in Equation~\eqref{equation:als-per-example-loss_real} is known as the \textit{per-example loss}.
It can be equivalently written as
\begin{equation}\label{equation:frob_loss_brmf}
L(\bW,\bZ) = \sum_{n=1}^N \sum_{m=1}^{M} \left(a_{mn} - \bw_m^\top\bz_n\right)^2 
= \norm{\bW\bZ-\bA}^2
=\trace\left\{(\bW\bZ-\bA)^\top (\bW\bZ-\bA)\right\},
\end{equation}
where $\trace(\cdot)$ denotes the trace operator.
This matrix factorization problem resembles a standard inverse problem---except that in a typical inverse problem, one of the factors is known, allowing ordinary least squares or similar methods to recover the unknown component by minimizing residuals.
When neither $\bW$ nor $\bZ$ is known, however, the factorization becomes highly non-convex and challenging---even when the latent dimension $K$ is as small as 2 or 3.
Due to the vast number of possible solutions and the lack of an analytical method to identify them, we sample the solution space using Markov chain Monte Carlo (MCMC) procedures to characterize its properties.

We discussed the Bayesian approach in Section~\ref{section:bayes_approach}. In the context of Bayesian matrix factorization, the model leads to the following form of Bayes' rule:
\begin{equation}\label{equation:bmf_bayes}
p(\bW,\bZ \mid \bA) \propto p(\bA \mid \bW, \bZ) \times p(\bW, \bZ),
\end{equation}
where $p(\bW, \bZ)$ encodes prior beliefs about the solution independently of the data, and $p(\bA\mid \bW, \bZ)$ is the likelihood, which measures how well the model explains the observed data.

\paragrapharrow{Terminology.} In Bayesian matrix factorization modeling, three modeling choices determine the specific type of matrix decomposition:
(i) the likelihood function,
(ii) the prior distributions placed on the factor matrices $\bW$ and $\bZ$, and
(iii) whether additional hierarchical (hyper)priors are imposed on the parameters of those priors.
We name the resulting model by listing the densities in the order: likelihood---prior for $\bW$---prior for $\bZ$, optionally followed by any hyperpriors.
For example, if the likelihood is Gaussian, the prior on $\bW$ is exponential, and the prior on $\bZ$ is Gaussian, the model is called a \textit{Gaussian Exponential-Gaussian (GEG)} model.
If a Gamma hyperprior is placed on the rate parameter of the exponential prior, the model becomes a \textit{Gaussian Exponential-Gaussian Gamma (GEGA)} model (where ``A" stands for Gamma to avoid confusion with Gaussian).
Table~\ref{table:summ_real_mf} summarizes the Bayesian models for real matrix factorization presented in this chapter.

\begin{table}[tp]
\centering
\setlength{\tabcolsep}{2.7pt}
\renewcommand{\arraystretch}{1.25}
\begin{tabular}{l|llll}
\hline
Name & Likelihood &Prior $\bW$ & Prior $\bZ$ &  Hierarchical prior \\ 
\hline \hline
GGG & $\normal(a_{mn}|\bw_m^\top\bz_n, \sigma^2)$ & 
$\normal(w_{mk}|0, (\lambda_{mk}^W)^{-1})$ & $\normal(z_{kn}|0, (\lambda_{kn}^Z)^{-1})$ & \gap\gap\slash \\ \hline
GGGM & $\normal(a_{mn}|\bw_m^\top\bz_n, \sigma^2)$ & 
$\normal(\bw_m|\bzero, \lambda^{-1}\bI)$ & $\normal(\bz_{n}|\bzero, \lambda^{-1}\bI)$ & \gap\gap\slash \\ \hline
GGGA& $\normal(a_{mn}|\bw_m^\top\bz_n, \sigma^2)$ &  
$\normal(w_{mk}|0, (\lambda_k)^{-1})$ & $\normal(z_{kn}|0, (\lambda_k)^{-1})$ & $\gammadist(\lambda_k|\alpha_\lambda, \beta_\lambda)$ \\ \hline
GGGW& $\normal(a_{mn}|\bw_m^\top\bz_n, \sigma^2)$ &  
$\normal(\bw_m|\bmu_w, \bSigma_w)$& $\normal(\bz_n|\bmu_z, \bSigma_z)$ &  
\begin{tabular}{@{}c@{}}$\{\bmu_w, \bSigma_w\}$, $\{\bmu_z, \bSigma_z\}\sim$ 
	\\ $\niw(\bmm_0, \kappa_0, \nu_0, \bS_0)$\end{tabular}\\ \hline
GVG & $\normal(a_{mn}|\bw_m^\top\bz_n, \sigma^2)$ & 
\begin{tabular}{@{}c@{}}$\bW\sim$ 
	\\ $\exp\{-\gamma \bW^\top\bW\}$ \end{tabular}
& $\normal(z_{kn}|0, (\lambda_{kn}^Z)^{-1})$ &   \gap\gap\slash
\\
 \hline
\end{tabular}
\caption{Overview of Bayesian real matrix factorization models.}
\label{table:summ_real_mf}
\end{table}

\index{Decomposition: GGG}
\section{All Gaussian (GGG) Model and Markov Blanket}\label{section:markov-blanket}

\begin{figure}[h]
\centering  
\subfigtopskip=2pt 
\subfigbottomskip=2pt 
\subfigcapskip=-5pt 
\includegraphics[width=0.421\linewidth]{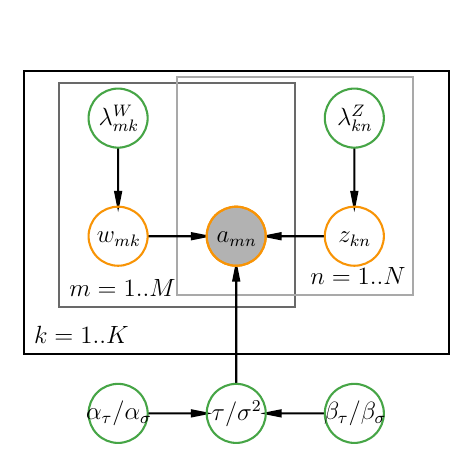}
\caption{Graphical model representation of the GGG model.  Green circles denote prior variables, orange circles represent observed and latent variables (shaded cycles denote observed variables), and plates indicate  repeated structures. The slash ``/" in the node denotes  ``or."}
\label{fig:bmf_ggg}
\end{figure}

The \textit{all-Gaussian (GGG)} model is perhaps the simplest Bayesian approach to RMF, employing Gaussian likelihood and Gaussian priors on the factor matrices \citep{salakhutdinov2008bayesian, gonen2012predicting, virtanen2011bayesian, virtanen2012bayesian}. 
\paragrapharrow{Likelihood.}
We interpret the data matrix  $\bA$ as generated by the probabilistic process depicted in Figure~\ref{fig:bmf_ggg}. 
Each observed   entry $a_{mn}$ of matrix $\bA$ is modeled using a Gaussian likelihood with   variance $\sigma^2$ and a mean determined by the latent decomposition $\bw_m^\top\bz_n$ (Equation~\eqref{equation:als-per-example-loss_real}):
\begin{equation}\label{equation:ggg_data_entry_likelihood}
	p(a_{mn} \mid \bw_m^\top\bz_n, \sigma^2) = \normal(a_{mn}\mid\bw_m^\top\bz_n, \sigma^2).
\end{equation}
Equivalently, this assumes that the residuals  $e_{mn}=a_{mn} -\bw_m^\top\bz_n$  are i.i.d. according to a zero-mean normal distribution with variance $\sigma^2$. 
This leads to the full likelihood:
\begin{equation}\label{equation:ggg_likelihood}
\begin{aligned}
p(\bA\mid \btheta) &= \prod_{m,n=1}^{M,N} \normal \left(a_{mn}\mid (\bW\bZ)_{mn}, \sigma^2 \right)\\
&= \prod_{m,n=1}^{M,N} \normal \left(a_{mn}\mid (\bW\bZ)_{mn}, \tau^{-1} \right),
\end{aligned}
\end{equation}
where $\btheta=\{\bW,\bZ,\sigma^2\}$ denotes all model parameters, $\sigma^2$ is the noise variance, and $\tau^{-1}=\sigma^2$ is the precision. 
Here,
$$
\normal(x\mid \mu,\sigma^2)=\frac{1}{(2\pi\sigma^2)^{1/2}} \exp \left\{ -\frac{1}{2\sigma^2 } (x-\mu)^2\right\}
=\sqrt{\frac{\tau}{2\pi}}\exp \left\{ -\frac{\tau}{2}(x-\mu)^2 \right\}
$$ 
is the normal density (see Definition~\ref{definition:gaussian_distribution}). 
\paragrapharrow{Prior.}
We place independent zero-mean Gaussian priors on the entries of $\bW$ and $\bZ$, with precisions $\{\lambda_{mk}^W\}$ and $\{\lambda_{kn}^Z\}$, respectively:
\begin{equation}\label{equation:ggg_prior_density_gaussian}
\begin{aligned}
\rw_{mk} &\sim \normal(w_{mk}\mid 0, (\lambda_{mk}^W)^{-1}), 
\gap 
&\rz_{kn}\sim&  \normal(z_{kn}\mid 0, (\lambda_{kn}^Z)^{-1});\\
p(\bW) &=\prod_{m,k=1}^{M,K} \normal(w_{mk}\mid 0, (\lambda_{mk}^W)^{-1}), 
\gap 
&p(\bZ) =&\prod_{k,n=1}^{K,N} \normal(z_{kn}\mid 0,( \lambda_{kn}^Z)^{-1}).
\end{aligned}
\end{equation}
For the noise variance $\sigma^2$,  we adopt an inverse-Gamma prior with shape ${\alpha_\sigma}$ and scale ${\beta_\sigma}$ (Definition~\ref{definition:inverse_gamma_distribution}),
\begin{equation}
p(\sigma^2)= \inversegammadist(\sigma^2\mid \alpha_\sigma, \beta_\sigma) = \frac{{\beta_\sigma}^{\alpha_\sigma}}{\Gamma({\alpha_\sigma})} (\sigma^2)^{-\alpha_\sigma-1} \exp\left( -\frac{{\beta_\sigma}}{\sigma^2} \right).
\end{equation}
By Bayes' rule (Equation~\eqref{equation:posterior_abstract_for_mcmc}), the posterior distribution is proportional to the product of the likelihood and the prior. This posterior can be maximized (e.g., via MAP estimation) or sampled from to obtain estimates of $\bW$ and $\bZ$, as discussed in Section~\ref{section:ppca_bpca_all}.

\index{Inverse-Gamma distribution}
\index{Markov blanket}
\paragrapharrow{Markov blanket.}
The most widely used methods for posterior inference in Bayesian models are \textit{Markov chain Monte Carlo (MCMC)} techniques, as described in Section~\ref{sec:monte_carlo_methods}. The core idea of MCMC is to construct a Markov chain over the latent variables whose stationary distribution is the target posterior \citep{andrieu2003introduction}. By simulating this chain, one eventually obtains samples from the posterior distribution.
A particularly convenient MCMC algorithm is Gibbs sampling, which iteratively samples each latent variable from its conditional posterior given all other variables and the observed data. Gibbs sampling is especially effective when these full conditional distributions are analytically tractable.

\begin{figure}[h!]
\centering
\includegraphics[width=0.35\textwidth]{imgs/markov_blanket.pdf}
\caption{The Markov blanket of a directed acyclic graphical (DAG) model. In a Bayesian network, the Markov blanket of node $A$ includes its \textbf{parents, children, and the other parents of all of its children (co-parents)}. The shaded cycle encloses all nodes in the Markov blanket of $A$. The figure is due to wikipedia page of Markov blanket.}
\label{fig:markov_blanket}
\end{figure}

To do Gibbs sampling, we need to derive the conditional posterior distributions for each parameter conditioned on all the other parameters $p(\theta_i \mid  \btheta_{-i}, \mathcalX)$, where $\mathcalX$ is again the set of  data points (here, the observed matrix $\bA$), and $\btheta_{-i}$ denotes all parameters except $\theta_i$. 
Crucially, in a graphical model, this conditional distribution depends only on the variables in the \textit{Markov blanket} of $\theta_i$.
For the GGG model in Figure~\ref{fig:bmf_ggg}---a \textit{directed acyclic graphical (DAG)} model---the Markov blanket of any node includes its \textbf{parents, children, and co-parents} \citep{jordan2004introduction}, as illustrated in Figure~\ref{fig:markov_blanket}.

\paragrapharrow{Example: Markov blanket for $w_{mk}$.} 
At first glance, the concept of a Markov blanket may seem abstract. Consider sampling the $(m,k)$-th entry 
$w_{mk}$ of $\bW$. From Figure~\ref{fig:bmf_ggg}, we identify:
(i) Parent: $\lambda_{mk}^W$; (ii)  Children: $\{a_{mn}\}_{n=1}^N$ (all observed entries in row $m$); 
(iii)  Co-parents: $\sigma^2$, the entire matrix $\bZ$, and all other entries of $\bW$ (denoted $\bW_{-mk}$).
Thus, the conditional posterior of $w_{mk}$ depends only on these variables:
$$
p(w_{mk} \mid  -)  = p(w_{mk}\mid  \bA,\bW_{-mk},\bZ,  \sigma^2,  \lambda_{mk}^W).
$$
More generally, the Markov blanket allows us to write the full conditionals for all parameters in the GGG model:
\begin{align}
p(w_{mk} \mid  -) 	
&= p(w_{mk}\mid  \bA,   \bW_{-mk}, \bZ,  \sigma^2, \lambda_{mk}^W), \label{equation:ggg_1} \\
p(z_{kn} \mid  -)   
&= p(z_{kn}\mid  \bA,  \bW, \bZ_{-kn}, \sigma^2,\lambda_{kn}^Z), \label{equation:ggg_2}\\
p(\sigma^2 \mid  -) 
&= p(\sigma^2\mid   \bA,\bW, \bZ, \alpha_\sigma, \beta_\sigma,)  \label{equation:ggg_3}.
\end{align} 
The Gibbs sampler proceeds by iteratively drawing samples from these conditional distributions:
\begin{itemize}
\item Sample each $w_{mk}$ from Equation~\eqref{equation:ggg_1} (the conditional distribution of $w_{mk}$);
\item Sample each $z_{kn}$ from Equation~\eqref{equation:ggg_2} (the conditional distribution of $z_{kn}$);
\item Sample $\sigma^2$ from Equation~\eqref{equation:ggg_3} (the conditional distribution of the noise variance).
\end{itemize}
This sequence of updates defines a Markov chain whose stationary distribution is the joint posterior 
$p(\bW,\bZ,\sigma^2\mid \bA)$. After a sufficient number of iterations (including a burn-in period), the samples provide a valid approximation to the posterior.

\paragrapharrow{Posterior.}
Given the observed matrix $\bA$, our goal is to estimate the \textit{posterior distribution} of the latent factors, $p(\bW,\bZ \mid \bA)$.
This posterior is central to matrix decomposition-based applications. For instance, in the Netflix context, we use the posterior expectations of each user's hidden preferences and each movie's latent attributes to predict which unwatched movies a user is likely to enjoy.
In this book, we employ Gibbs sampling for posterior inference because it provides highly accurate approximations of the true posterior. A key advantage of MCMC methods like Gibbs sampling is that they yield asymptotically exact results---that is, as the number of samples grows, the empirical distribution converges to the true posterior.
An alternative approach is variational Bayesian inference (see Section~\ref{section:vbi_root}), which will be shortly covered in this context.
For RMF, applying Bayes' rule together with MCMC means we must be able to sample from the following full conditional distributions (determined by the Markov blanket):
$$
\begin{aligned}
p(w_{mk} &\mid  \bA, \bW_{-mk}, \bZ, \sigma^2, \lambda_{mk}^W), \\
p(z_{kn}&\mid \bA,  \bW, \bZ_{-kn},\sigma^2,   \lambda_{kn}^Z ), \\
p(\sigma^2 &\mid  \bA,\bW, \bZ, \alpha_\sigma, \beta_\sigma), \\
\end{aligned}
$$
where $\bW_{-{mk}}$ denotes all entries  of $\bW$ except $w_{mk}$, and $\bZ_{-kn}$ denotes all entries  of $\bZ$ except $z_{kn}$. 
By Bayes' theorem, the conditional density of $w_{mk}$ depends only on its Markov blanket: parents ($\lambda_{mk}^W$), children (the observed entries $\{a_{mn}\}_{n=1}^N$, i.e., row $m$ of $\bA$), and co-parents (the noise variance $\sigma^2$, the other entries of $\bW$---$\bW_{-mk}$, and the entire matrix $\bZ$) (See Figure~\ref{fig:bmf_ggg} and Section~\ref{section:markov-blanket}.) 
The conditional posterior can be derived as follows:
\begin{equation}\label{equation:ggg_poster_wmk1}
\begin{aligned}
&\gap p(w_{mk} \mid  \bA ,   \bW_{-mk}, \bZ, \sigma^2, \lambda_{mk}^W ) 
\propto p(\bA\mid  \bW, \bZ, \sigma^2) \times p(w_{mk}\mid  \lambda_{mk}^W)\quad\\
&=\prod_{i,j=1}^{M,N} \normal \left(a_{ij}\mid  \bw_i^\top\bz_j, \sigma^2 \right)\times \normal(w_{mk}\mid 0, (\lambda_{mk}^W)^{-1}) \\
&\propto \exp\Bigg\{   -\frac{1}{2\sigma^2}  \sum_{i,j=1}^{M,N}(a_{ij} - \bw_i^\top\bz_j  )^2\Bigg\}  \times \exp\left\{ -\frac{w_{mk}^2}{2} \lambda_{mk}^W \right\} \\
\end{aligned}
\end{equation}
\begin{align*}
\qquad\quad 
&\propto \exp\Bigg\{   -\frac{1}{2\sigma^2}  \sum_{j=1}^{N}(a_{mj} - \bw_m^\top\bz_j  )^2\Bigg\}  
\times \exp\left\{ -\frac{w_{mk}^2}{2}\lambda_{mk}^W \right\} \\
& \propto \exp\Bigg\{   -\frac{1}{2\sigma^2}  \sum_{j=1}^{N}
\Bigg[ w_{mk}^2z_{kj}^2 + 2w_{mk} z_{kj}\bigg(\sum_{i\neq k}^{K}w_{mi}z_{ij} - a_{mj}\bigg)  \Bigg]
\Bigg\}  \cdot \exp\left\{ -\frac{w_{mk}^2}{2}\lambda_{mk}^W \right\}\\
&\stackrel{\star}{\propto} \exp\Bigg\{   
-\underbrace{
\left(\frac{\sum_{j=1}^{N} z_{kj}^2 }{2\sigma^2} +\frac{\lambda_{mk}^W}{2} \right) }_{\textcolor{mylightbluetext}{\triangleq 1/(2\widetilde{\sigma_{mk}^{2}})}}
w_{mk}^2 
+
w_{mk}\underbrace{\Bigg(  \frac{1}{\sigma^2} \sum_{j=1}^{N} z_{kj}\bigg( a_{mj} - \sum_{i\neq k}^{K}w_{mi}z_{ij}\bigg)  \Bigg)}_{\textcolor{mylightbluetext}{\triangleq \widetilde{\sigma_{mk}^{2}}^{-1} \widetilde{\mu_{mk}}}}
\Bigg\}  \\
&\propto   \normal(w_{mk} \mid  \widetilde{\mu_{mk}}, \widetilde{\sigma_{mk}^{2}}),
\end{align*}
where the equality ($\star$) follows from completing the square and matches the canonical form of a Gaussian density (see Equation~\eqref{equation:gaussian_form_conform}\index{Canonical form}).
Thus, the posterior variance and mean are given by
\begin{align}
\widetilde{\sigma_{mk}^{2}} 
&= 1\bigg/ \left(	\frac{1 }{\sigma^2} \sum_{j=1}^{N} z_{kj}^2 +\lambda_{mk}^W\right);
\label{equation:ggg_posterior_variance}\\
\widetilde{\mu_{mk}} &= \frac{\widetilde{\sigma_{mk}^{2}}}{\sigma^2} \cdot  \sum_{j=1}^{N} z_{kj}\bigg( a_{mj} - \sum_{i\neq k}^{K}w_{mi}z_{ij}\bigg) .
\label{equation:ggg_posterior_mean}
\end{align}
Notably, the posterior precision $1/ \widetilde{\sigma_{mk}^{2}}$ is the sum of the \textit{prior precision} $\lambda_{mk}^W$ and the \textit{data precision} $\frac{1 }{\sigma^2} \sum_{j=1}^{N} z_{kj}^2$.
In the Netflix example, this implies that the conditional distribution of a movie's latent feature vector $\bw_m$ for $m=1,2,\ldots,M$, given user features, observed ratings, and hyper-parameters, is Gaussian.

By symmetry, an analogous derivation applies to  $\{z_{kn}\}$ (for $k={1,2,\ldots, K}$ and $n={1,2,\ldots, N}$).
The resulting conditional posterior for $z_{kn}$ is also Gaussian, with updated mean and variance derived similarly; see Problem~\ref{problem:symmetry_zkn}.

The conditional density of $\sigma^2$ depends on its parents ($\alpha_\sigma$, $\beta_\sigma$), children ($\bA$), and co-parents ($\bW$, $\bZ$). 
Due to conjugacy between the Gaussian likelihood and the inverse-Gamma prior (see Equation~\eqref{equation:inverse_gamma_conjugacy_general}), the posterior is also inverse-Gamma:
\begin{equation}\label{equation:ggg_posterior_sigma2}
\begin{aligned}
&p(\sigma^2 \mid \bA, \bW,\bZ, \alpha_\sigma, \beta_\sigma ) = \inversegammadist (\sigma^2\mid  \widetilde{\alpha_{\sigma}}, \widetilde{\beta_{\sigma}}),  \\
& \widetilde{\alpha_{\sigma}} = \frac{MN}{2} +{\alpha_\sigma}, \gap 
\widetilde{\beta_{\sigma}}  =  \frac{1}{2} \sum_{m,n=1}^{M,N} (\bA-\bW\bZ)_{mn}^2 + {\beta_\sigma}.
\end{aligned}
\end{equation}

\paragrapharrow{Missing entries.} 
As noted previously, in many practical scenarios---such as the Netflix problem---some entries of $\bA$ are missing. 
Let  $\sS$ denote the set of indices $(m,n)$ corresponding to observed entries in $\bA$. Denote further $\sS_m = \{n \mid  (m,n)\in {\sS}\}$, i.e., the observed entries in the $m$-th row (with $\abs{\sS_m}\leq N$); $\sS_n = \{m \mid  (m,n)\in \sS\}$, i.e., the observed entries in the $n$-th column (with $\abs{\sS_n}\leq M$). 
When entries are missing, the posterior for $w_{mk}$ remains Gaussian, but the sums are restricted to observed data:
\begin{equation}
w_{mk}\sim \normal(w_{mk} \mid  \widetilde{\mu_{mk}}, \widetilde{\sigma_{mk}^{2}}), 
\end{equation}
where 
\begin{equation}
\begin{aligned}
\widetilde{\sigma_{mk}^{2}} &= 1\bigg/ \left(	\frac{1}{\sigma^2}\sum_{\textcolor{mylightbluetext}{j\in\sS_m}} z_{kj}^2 +\lambda_{mk}^W\right),\\
\widetilde{\mu_{mk}} &= \frac{\widetilde{\sigma_{mk}^{2}}}{\sigma^2} \cdot  \sum_{\textcolor{mylightbluetext}{j\in\sS_m}} z_{kj}\bigg( a_{mj} - \sum_{i\neq k}^{K}w_{mi}z_{ij}\bigg) .
\end{aligned}
\end{equation}
For simplicity, the following discussion assumes a fully observed data matrix $\bA$. However, extensions to handle missing entries follow the same principles and can be derived analogously.

\index{Decomposition: GGGM}
\paragrapharrow{GGG with shared prior (GGGM).}
When we set $\lambda=\lambda_{mk}^W$ for all $m\in\{1,2,\ldots, M\}$ and $k\in\{1,2,\ldots,K\}$ (and similarly $\lambda=\lambda_{kn}^Z$ for all $k,n$), the conditional posterior distributions in the Gibbs sampling algorithm can be expressed as multivariate Gaussian densities (see Definition~\ref{definition:multivariate_gaussian}). 
In this case, we place an \textit{isotropic} multivariate Gaussian prior on each row $\bw_m$ of $\bW$ and each column $\bz_n$ of $\bZ$:
\begin{equation}\label{equation:gggm_prior}
\bw_m \sim \normal(\bw_m \mid  \bzero, \lambda^{-1}\bI), \gap 
\bz_n \sim \normal(\bz_n \mid  \bzero, \lambda^{-1}\bI),
\end{equation}
meaning that each of the $M$ item factors and $N$ user factors follows a multivariate normal distribution with spherical covariance.
Let $\bW_{-m}$ denote all rows of $\bW$ except the $m$-th row. 
Again by Bayes' rule, conditional posterior for $\bw_m$ in Gibbs sampling is then:
\begin{equation}\label{equation:gggm_wm_post}
\begin{aligned}
& \gap p(\bw_m \mid  \sigma^2, \bW_{-m}, \bZ, \lambda, \bA)  
\propto p(\bA \mid  \bW, \bZ, \sigma^2) \times \normal(\bw_m \mid  \bzero ,\lambda^{-1}\bI) \\
&\propto \normal(\bA \mid  \bW\bZ, \sigma^2\bI) \times \normal(\bw_m \mid  \bzero ,\lambda^{-1}\bI) \\
&\propto \exp\Bigg\{   -\frac{1}{2\sigma^2}  \sum_{j=1}^{N}(a_{mj} - \bw_m^\top\bz_j  )^2\Bigg\}  
\times \exp\left\{ -\frac{\lambda}{2} \bw_m^\top\bw_m \right\} \\
&\stackrel{\star}{\propto} \exp\Bigg\{-\frac{1}{2} \bw_m^\top 
\underbrace{\Bigg[\lambda\bI + \frac{1}{\sigma^2}\sum_{j=1}^{N}\bz_j\bz_j^\top \Bigg]}_{\textcolor{mylightbluetext}{\triangleq \widetilde{\bSigma}^{-1}}}
\bw_m
+ \bw_m^\top 
\underbrace{\frac{1}{\sigma^2} \sum_{j=1}^{N} a_{mj}\bz_j }_{\textcolor{mylightbluetext}{\triangleq \widetilde{\bSigma}^{-1} \widetilde{\bmu} }}
\Bigg\}
\propto \normal(\bw_m \mid  \widetilde{\bmu}, \widetilde{\bSigma}),
\end{aligned}
\end{equation}
where the equality ($\star$) follows from completing the square and matches the canonical form of a multivariate Gaussian (see  Equation~\eqref{equation:multi_gaussian_form_conform}).
The posterior covariance and mean are given by:
$$
\widetilde{\bSigma}= \Bigg[\lambda\bI + \frac{1}{\sigma^2}\sum_{j=1}^{N}\bz_j\bz_j^\top \Bigg]^{-1} 
\qquad \text{and}\qquad 
\widetilde{\bmu} = \frac{1}{\sigma^2}  \widetilde{\bSigma}\cdot \sum_{j=1}^{N} a_{mj}\bz_j.
$$
In this case, the \textit{posterior precision matrix $\widetilde{\bSigma}^{-1}$ is the sum of the prior precision matrix $\lambda\bI$ and the data precision matrix $\frac{1}{\sigma^2}\sum_{j=1}^{N}\bz_j\bz_j^\top$}.
By symmetry, an analogous expression holds for each column $\bz_n$ of factor $\bZ$ (for $n\in\{1,2,\ldots, N\}$). 
See Problem~\ref{problem:symmetry_zkngggm}.

\paragrapharrow{Gibbs sampling.}
Because conjugate priors are used for both parameters and hyper-parameters in the Bayesian matrix factorization model, sampling from the full conditional distributions is straightforward.
As described in Section~\ref{section:gibbs-sampler}, we can construct a Gibbs sampler for the GGG model as outlined in Algorithm~\ref{alg:ggg_gibbs_sampler}.
Thanks to our choice of conjugate priors, all conditional distributions belong to standard families (Gaussian or inverse-Gamma), allowing direct sampling without resorting to slower methods like rejection sampling.
In practice, it is common to assume identical prior precisions across all latent dimensions, i.e.,
$\lambda=\{\lambda_{mk}^W\} = \{\lambda_{kn}^Z\}$ for all $m,k,n$. 
Default uninformative hyper-parameter settings are often:
$\alpha_\sigma=\beta_\sigma=1$, $\{\lambda_{mk}^W\} = \{\lambda_{kn}^Z\} = 0.1$.

\begin{algorithm}[htp] 
\caption{Gibbs sampler for GGG model in one iteration. 
A variance prior on $\sigma^2$ is used here; an equivalent formulation exists using precision $\tau=1/\sigma^2$. 
The algorithm is presented for clarity---not efficiency. A vectorized implementation would be significantly faster. Default hyper-parameters:  $\alpha_\sigma=\beta_\sigma=1$, $\{\lambda_{mk}^W\} = \{\lambda_{kn}^Z\} = 0.1$.} 
\label{alg:ggg_gibbs_sampler}  
\begin{algorithmic}[1] 
\Require Choose initial $\alpha_\sigma, \beta_\sigma, \lambda_{mk}^W, \lambda_{kn}^Z$;
\For{$k=1$ to $K$} 
\For{$m=1$ to $M$}
\State Sample $w_{mk}$ from $p(w_{mk}\mid  \bA ,   \bW_{-mk}, \bZ, \sigma^2, \lambda_{mk}^W )$; 
\Comment{Equation~\eqref{equation:ggg_poster_wmk1}}
\EndFor
\For{$n=1$ to $N$}
\State Sample $z_{kn}$ from $p(z_{kn}\mid \bA ,   \bW, \bZ_{-kn}, \sigma^2, \lambda_{kn}^Z)$; 
\Comment{Symmetry of Eq.~\eqref{equation:ggg_poster_wmk1}}
\EndFor
\EndFor
\State Sample $\sigma^2$ from $p(\sigma^2\mid \bA, \bW,\bZ, \alpha_\sigma, \beta_\sigma)$; 
\Comment{Equation~\eqref{equation:ggg_posterior_sigma2}}
\State Report loss in Equation~\eqref{equation:frob_loss_brmf}, stop if it converges;
\end{algorithmic} 
\end{algorithm}

\paragrapharrow{Prior by Gamma distribution.} 
Note that placing an inverse-Gamma prior on the variance $\sigma^2$ of a Gaussian density
is equivalent to assigning a Gamma prior on the precision $\tau=\sigma^{-2}$.
We model $\tau$ using a Gamma distribution with shape $\alpha_\tau>0$ and rate $\beta_\tau>0$ (Definition~\ref{definition:gamma-distribution}):
\begin{equation}\label{equation:ggg_gamma_prior}
p(\tau)\sim \gammadist (\tau\mid  \alpha_\tau, \beta_\tau) = \frac{\beta_{\tau}^{\alpha_{\tau}}}{\Gamma(\alpha_{\tau})} \tau^{\alpha_{\tau}-1} \exp({-\beta_{\tau}\cdot \tau}),
\end{equation}
Due to conjugacy (Equation~\eqref{equation:gamma_conjugacy_general}), the posterior is also Gamma:
$$
\begin{aligned}
&\gap 
p(\tau \mid  \bW,\bZ, \bA) = \gammadist (\tau; \widetilde{\alpha_\tau}, \widetilde{\beta_\tau}), \\
&\gap \widetilde{\alpha_\tau} = \frac{MN}{2} +{\alpha_\tau}, \gap 
\widetilde{\beta_\tau}  =  \frac{1}{2} \sum_{m,n=1}^{M,N} (\bA-\bW\bZ)_{mn}^2 + {\beta_\tau}.
\end{aligned}
$$
In practice, one may choose $\alpha_\tau=\alpha_\sigma$ and $\beta_\tau=\beta_\sigma$ to maintain consistency between the two parameterizations.

\index{Variational Bayesian inference}
\paragrapharrow{Variational Bayesian inference.}
Although this book primarily focuses on Gibbs sampling for Bayesian matrix factorization, it is worth noting that \textit{variational Bayesian (VB) inference} for the GGG model leads to updates that closely resemble those of the Gibbs sampler---thanks to the Gaussian likelihood and conjugate priors; see Section~\ref{section:vbi_root}.
Under the \textit{mean-field approximation}, we assume a factorized variational distribution; see Section~\ref{section:em_constrained}. For the latent variable 
$w_{mk}$, we posit a Gaussian form $q(w_{mk})=\normal(w_{mk} \mid  \widetilde{\mu_{mk}}, \widetilde{\sigma_{mk}^{2}})$. 
Following Equation~\eqref{equation:em_uncon_mf}, the optimal variational distribution satisfies:
$$
\begin{aligned}
&q_{w_{mk}}(w_{mk}) 
\propto \exp\left\{ \Exp_{q(-w_{mk})} \left[\ln p(\bA \mid \bW, \bZ) +\ln p(\bW, \cancel{\bZ}) \right] \right\}\\
&\propto \exp\Bigg\{ \Exp_{q_{\btheta(-w_{mk})}} \Bigg[\Bigg\{   -\frac{1}{2\sigma^2}  \sum_{j=1}^{N}(a_{mj} - \bw_m^\top\bz_j  )^2\Bigg\} +\ln\Bigg(\exp\left\{ -\frac{w_{mk}^2}{2} \lambda_{mk}^W \right\}\Bigg) \Bigg] \Bigg\}\\
&\propto \exp\Bigg\{ \Exp_{q_{\btheta(-w_{mk})}} \Bigg[   -\frac{1}{2\sigma^2}  \sum_{j=1}^{N}(a_{mj} - \bw_m^\top\bz_j  )^2  \Bigg] \Bigg\}\cdot \exp\left\{ -\frac{w_{mk}^2}{2} \lambda_{mk}^W \right\}\\
&\propto \exp\Bigg\{   -\frac{1}{2\sigma^2}  \sum_{j=1}^{N}
\Bigg[ w_{mk}^2z_{kj}^2 + 2w_{mk} z_{kj}\bigg(\sum_{i\neq k}^{K}w_{mi}z_{ij} - a_{mj}\bigg)  \Bigg]
\Bigg\}  \cdot \exp\left\{ -\frac{w_{mk}^2}{2}\lambda_{mk}^W \right\}\\
&\propto \exp\Bigg\{   
-\underbrace{
\left(\frac{\sum_{j=1}^{N} z_{kj}^2 }{2\sigma^2} +\frac{\lambda_{mk}^W}{2} \right) }_{\textcolor{mylightbluetext}{\triangleq 1/(2\widetilde{\sigma_{mk}^{2}})}}
w_{mk}^2 
+
w_{mk}\underbrace{\Bigg(  \frac{1}{\sigma^2} \sum_{j=1}^{N} z_{kj}\bigg( a_{mj} - \sum_{i\neq k}^{K}w_{mi}z_{ij}\bigg)  \Bigg)}_{\textcolor{mylightbluetext}{\triangleq \widetilde{\sigma_{mk}^{2}}^{-1} \widetilde{\mu_{mk}}}}
\Bigg\}  \\
&\propto   \normal(w_{mk} \mid  \widetilde{\mu_{mk}}, \widetilde{\sigma_{mk}^{2}}).
\end{aligned}
$$
These update equations are identical in form to those derived for the Gibbs sampler in Equation~\eqref{equation:ggg_poster_wmk1}. Thus, under the mean-field assumption and pointwise evaluation of expectations, variational inference yields the same parameter updates as Gibbs sampling---though the interpretations differ (deterministic optimization vs. stochastic sampling).

\begin{figure}[h]
\centering  
\vspace{-0.35cm} 
\subfigtopskip=2pt 
\subfigbottomskip=2pt 
\subfigcapskip=-5pt 
\subfigure[GGGA.]{\label{fig:bmf_ggga}
\includegraphics[width=0.421\linewidth]{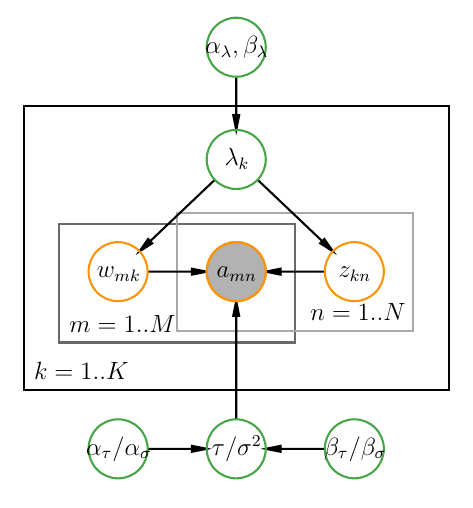}}
\subfigure[GGGW.]{\label{fig:bmf_gggw}
\includegraphics[width=0.421\linewidth]{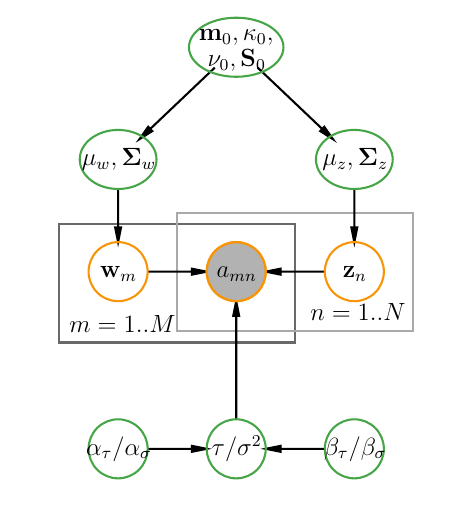}}
\caption{Graphical model representation of GGGA and GGGW models. Green circles denote prior variables, orange circles represent observed and latent variables  (shaded cycles denote observed variables), and plates represent repeated variables. 
The slash ``/" in the variable represents ``or," and the comma ``," in the variable represents ``and."}
\label{fig:bmf_ggga_gggw}
\end{figure}

\index{Decomposition: GGGA}
\section{All Gaussian Model with ARD Hierarchical Prior (GGGA)}\label{section:brmf_ggga}
The \textit{all-Gaussian   with hierarchical Gamma prior (GGGA)} model was  proposed by \citet{virtanen2011bayesian, virtanen2012bayesian} as an extension of the GGG model.
The key distinction is that GGGA places a hyperprior over the Gaussian prior on the latent factors. This enables \textit{automatic relevance determination} (ARD), a mechanism that facilitates automatic model selection by adaptively pruning irrelevant latent dimensions (see Figure~\ref{fig:bmf_ggga} and Section~\ref{section:bayespca}).

\index{Automatic relevance determination}
\index{Linear model}
\index{Linear regression}
\index{Bayesian linear regression}
\paragrapharrow{Automatic relevance determination (ARD).}
\textit{Automatic relevance determination (ARD)} is a Bayesian technique used in machine learning and statistics to automatically assess the relevance of input features or latent components  (see Section~\ref{section:bayespca}). It is commonly applied in Bayesian linear regression (see Section~\ref{section:bayespca}).

In Bayesian linear regression, the relationship between inputs and outputs is modeled probabilistically:
$$
\bb = \bA\bx +\bepsilon,
$$
where $\bA\in\real^{M\times N}$ is the input data matrix, $\bx\in\real^N$ is a vector of weights, and $\bepsilon\in\real^N$ is Gaussian noise. 
ARD extends this framework by introducing a probabilistic prior distribution over the weights. 
ARD extends this framework by assigning a separate precision parameter (i.e., inverse variance) to each component of $\bx$.

The core idea of ARD is to let the data determine feature relevance.
Irrelevant features are assigned high precision (low variance), effectively shrinking their weights toward zero.
Relevant features receive low precision (high variance), allowing them to take larger values.
This adaptive shrinkage enables the model to automatically select a subset of meaningful features or latent factors without manual intervention.

In a full Bayesian treatment, ARD introduces a hyperprior over the precision parameters. Inference then jointly estimates the weights and their precisions, accounting for uncertainty in both.

Key advantages of ARD include automatic feature/latent dimension selection, reduced overfitting, and improved model interpretability. However, careful choice of hyperpriors and hyper-parameters is essential to obtain meaningful results.

\index{Hyperprior}
\paragrapharrow{Hyperprior.}
Building on the GGG model, the GGGA model adopts an ARD-style hierarchical prior. Specifically, we place a shared Gamma hyperprior over the precision of each latent dimension $k$:
\begin{equation}\label{equation:hyper_eq_ggga}
\begin{aligned}
w_{mk} &\sim \normal( 0, (\lambda_{k})^{-1}), 
\gap 
&z_{kn}\sim&  \normal( 0, (\lambda_{k})^{-1}), \gap 
\lambda_{k} \sim \gammadist( \alpha_{\lambda}, \beta_{\lambda}),\\
\end{aligned}
\end{equation}
for all $k\in\{1,2,\ldots, K\}$.
Here, the same precision parameter $\lambda_k$ governs all entries in column $k$ of $\bW$ and the  row $k$ of $\bZ$. 
As a result, the entire latent factor $k$ is either:
(i) activated (if $\lambda_k$ is small, implying large variance), or (ii) suppressed (if $\lambda_k$  is large, forcing weights toward zero). This behavior is illustrated in Figure~\ref{fig:dists_gaussian}.

\index{Markov blanket}
\paragrapharrow{Posterior.}
For Bayesian matrix factorization, Gibbs sampling requires drawing from the full conditional distributions of all variables. Using the Markov blanket principle (Section~\ref{section:markov-blanket}), these conditionals are:
$$
\begin{aligned}
&p(\sigma^2\mid \bA, \bW, \bZ, \alpha_\sigma, \beta_\sigma ), 
&\gap& p(w_{mk}\mid \bA,   \bW_{-mk}, \bZ, \sigma^2, \blambda ), \\
&p(\lambda_k \mid \bW, \bZ, \blambda_{-k}, \sigma^2,  \alpha_\lambda, \beta_\lambda),  
&\gap& p(z_{kn}\mid \bA,  \bW, \bZ_{-kn}, \sigma^2, \blambda ), \\
\end{aligned}
$$
where $\blambda=[\lambda_1,\lambda_2,\ldots,\lambda_K]^\top\in \real_+^K$ is a vector including all $\lambda_k$ values, and $\blambda_{-k}$ denotes all elements of $\blambda$ except $\lambda_k$. 
The conditional posteriors for variables $\sigma^2$, $\{w_{mk}\}$, and $\{z_{kn}\}$ remain structurally identical to those in the GGG model, except now we replace $\lambda_{mk}^W$ and $\lambda_{kn}^Z$ for all $m\in\{1,2,\ldots, M\}$ and $n\in\{1,2,\ldots, N\}$ by shared $\lambda_k$. 
The conditional posterior density of  $\lambda_{k}$ depends on its parents ($\alpha_\lambda$, $\beta_\lambda$), children (the $k$-th column $\widetildebw_k$ of $\bW$ and the  $k$-th row $\widetildebz_k$ of $\bZ$; see definition in Equation~\eqref{equation:als-per-example-loss_real})
\footnote{Note we denote the $m$-th row of $\bW$ by $\bw_m$, and $n$-th column of $\bZ$ by $\bz_n$ previously.}, and co-parents ($\blambda_{-k}$). (See Figure~\ref{fig:bmf_ggga} and Section~\ref{section:markov-blanket}.)
The posterior is derived as follows:
\begin{equation}\label{equation:posterior-ggga_lambdak}
\begin{aligned}
&\gap p(\lambda_k\mid  \bW, \bZ, \blambda_{-k}, \sigma^2,  \alpha_\lambda, \beta_\lambda)
\propto  p(\widetildebw_k, \widetildebz_k\mid  \lambda_k) \cdot p(\lambda_k) \\
&= 
\prod_{i=1}^{M}  \normal(\widetildew_{ik}\mid 0, (\lambda_{k})^{-1}) \cdot 
\prod_{j=1}^{N} \normal(\widetildez_{kj}\mid  0,(\lambda_{k})^{-1})
\cdot \gammadist(\lambda_k \mid  \alpha_{\lambda}, \beta_\lambda)\\
&=
\prod_{i=1}^{M} \lambda_k^{1/2} \exp\left\{-\frac{\lambda_k \widetildew_{ik}^2}{2}\right\} \cdot 
\prod_{j=1}^{N} \lambda_k^{1/2} \exp\left\{-\frac{\lambda_k \widetildez_{kj}^2}{2}\right\}
\cdot \frac{\beta_{\lambda}^{\alpha_\lambda}}{\Gamma(\alpha_\lambda)} \lambda_k^{\alpha_\lambda-1} \exp(- \lambda_k \beta_\lambda)\\
&\propto\lambda_k^{\frac{M+N}{2} +\alpha_\lambda -1} 
\exp\Bigg\{-\lambda_k \cdot  \Bigg(\frac{1}{2}\sum_{i=1}^{M}\widetildew_{ik}^2 + 
\frac{1}{2}\sum_{j=1}^{N}\widetildez_{kj}^2 +\beta_\lambda\Bigg)\Bigg\}
\propto \gammadist(\lambda_k\mid  \widetilde{\alpha_\lambda}, \widetilde{\beta_\lambda}),
\end{aligned}
\end{equation}
where 
$$
\widetilde{\alpha_\lambda} = \frac{M+N}{2}+\alpha_\lambda, 
\qquad  
\widetilde{\beta_\lambda}=\frac{1}{2}\sum_{i=1}^{M}\widetildew_{ik}^2 + 
\frac{1}{2}\sum_{j=1}^{N}\widetildez_{kj}^2 +\beta_\lambda.
$$
From the properties of the Gamma distribution (Definition~\ref{definition:gamma-distribution}), the posterior mean and variance for $\lambda_k$ are:
$$
\Exp[\lambda_k] = \frac{\widetilde{\alpha_\lambda}}{\widetilde{\beta_\lambda}}, 
\qquad 
\Var[\lambda_k] = \frac{\widetilde{\alpha_\lambda}}{\widetilde{\beta_\lambda}^2}.
$$
This reveals two intuitive behaviors:
\begin{itemize}
\item \textit{Larger matrices favor stronger regularization.}
Upon observing the posterior parameters, it is evident that  with a larger  shape of the raw matrix $\bA$ (i.e., $M+N$ is larger), there is a preference for  a larger value of $\lambda_k$ (during sampling, since $\widetilde{\alpha_\lambda}$ tends to be larger). 
As indicated in Equation~\eqref{equation:hyper_eq_ggga}, this preference imposes a larger and sparser regularization over the model.
This is reasonable since a larger shape  indicates the vector product for each entry of $\bA$ through $\bW$ and $\bZ$ involves more entries to sum up.

\item \textit{Large latent values suppress shrinkage.}
If the current samples of $\widetildew_{ik}$ or $\widetildez_{kj}$ are large in magnitude, $\widetilde{\beta_\lambda}$ increases, leading to a smaller posterior mean for $\lambda_k$. This reduces shrinkage, allowing the factor to explore a wider range of values---consistent with the evidence from previous Gibbs iterations. This is reasonable in the sense that we want to explore in a larger space if the factored components have larger elements from previous Gibbs iterations.
\end{itemize}
Thus, the ARD mechanism dynamically balances model complexity and data fit through the hierarchical prior.

\paragrapharrow{Gibbs sampling.}
A Gibbs sampler for the GGGA model is summarized in Algorithm~\ref{alg:ggga_gibbs_sampler}. By default, we use uninformative hyper-parameters: $\alpha_\sigma=\beta_\sigma=1$, $\alpha_\lambda=\beta_\lambda=1$.

\begin{algorithm}[h] 
\caption{
Gibbs sampler for GGGA model in one iteration.
A variance prior on $\sigma^2$ is used; an equivalent formulation exists using precision $\tau=1/\sigma^2$. 
The algorithm prioritizes clarity over efficiency---a vectorized implementation would be faster.
By default, uninformative hyper-parameters are $\alpha_\sigma=\beta_\sigma=1$, $\alpha_\lambda=\beta_\lambda=1$.
} 
\label{alg:ggga_gibbs_sampler}  
\begin{algorithmic}[1] 
\Require Choose initial $\alpha_\sigma, \beta_\sigma, \alpha_\lambda, \beta_\lambda$;
\For{$k=1$ to $K$} 
\For{$m=1$ to $M$}
\State Sample $w_{mk}$ from $p(w_{mk}\mid \bA,    \bW_{-mk}, \bZ,\sigma^2,  \lambda_k )$; \Comment{Equation~\eqref{equation:ggg_poster_wmk1}}
\EndFor
\For{$n=1$ to $N$}
\State Sample $z_{kn}$ from $p(z_{kn}\mid \bA ,   \bW, \bZ_{-kn}, \sigma^2, \lambda_k)$; 
\Comment{Symmetry of Eq.~\eqref{equation:ggg_poster_wmk1}}
\EndFor
\State Sample $\lambda_k$ from $p(\lambda_k\mid   \bW, \bZ, \sigma^2, \alpha_\lambda, \beta_\lambda)$; 
\Comment{Equation~\eqref{equation:posterior-ggga_lambdak}}
\EndFor
\State Sample $\sigma^2$ from $p(\sigma^2\mid \bA, \bW,\bZ, \alpha_\sigma, \beta_\sigma)$;
\Comment{Equation~\eqref{equation:ggg_posterior_sigma2}}
\State Report loss in Equation~\eqref{equation:frob_loss_brmf}, stop if it converges;
\end{algorithmic} 
\end{algorithm}

\index{Decomposition: GGGW}
\section{All Gaussian Model with Wishart Hierarchical Prior (GGGW)}\label{section:gggw_model}
The hierarchical prior based on the Wishart distribution was introduced by \citet{salakhutdinov2008bayesian} to enhance flexibility and improve model calibration relative to the GGG model.
Rather than assuming independence among individual entries of the factor matrices $\bW$ and $\bZ$, the GGGW model assumes that:
each row $\bw_m$ of $\bW$ and each column $\bz_n$ of $\bZ$ follows  a multivariate Gaussian distribution  (Definition~\ref{definition:multivariate_gaussian}). 
The mean and covariance parameters of these Gaussians are themselves assigned a \textit{normal-inverse-Wishart (NIW)} hyperprior (Definition~\ref{definition:normal_inverse_wishart}).

\paragrapharrow{Prior and hyperprior.}
As in the GGG model, we assume a Gaussian likelihood for the observed data matrix $\bA$, and the variance parameter $\sigma^2$ is placed over an inverse-Gamma prior with shape $\alpha_\sigma$ and scale $\beta_\sigma$: $\inversegammadist(\sigma^2\mid \alpha_\sigma, \beta_\sigma)$.
Given the $m$-th row $\bw_m$ of $\bW$ and the $n$-th column $\bz_n$ of $\bZ$, we consider the multivariate Gaussian density and the normal-inverse-Wishart prior as follows:
\begin{align}
\bw_m&\sim \normal(\bw_m\mid  \bmu_w, \bSigma_w), \gap\gap& \bmu_w, \bSigma_w &\sim \niw(\bmu_w, \bSigma_w\mid  \bmm_0, \kappa_0, \nu_0, \bS_0);\label{equa:gggw_prior_hyper1}\\
\bz_n&\sim \normal(\bz_n\mid  \bmu_z, \bSigma_z), \gap\gap& \bmu_z, \bSigma_z &\sim \niw(\bmu_z, \bSigma_z\mid  \bmm_0, \kappa_0, \nu_0, \bS_0),\label{equa:gggw_prior_hyper2}
\end{align}
where $\niw (\bmu, \bSigma\mid  \bmo, \kappa_0, \nu_0, \bso) 
= \mathcal{N}(\bmu\mid  \bmo, \frac{1}{\kappa_0}\bSigma) \cdot  \inversewishart(\bSigma\mid  \bso, \nu_0)$ is the density of a normal-inverse-Wishart distribution, 
and $ \inversewishart(\bSigma\mid  \bso, \nu_0)$ denotes the  inverse-Wishart distribution (Definition~\ref{definition:multi_inverse_wishart}). 
Although one could alternatively place independent (semi-conjugate) priors on the mean and covariance---i.e., a normal prior on $\bmu$ and an inverse-Wishart prior on $\bSigma$---we focus here on the fully conjugate NIW formulation. Details on the semi-conjugate approach can be found in Sections~\ref{section:sep_mu_niw} and~\ref{section:sep_sigma_niw}.

By conjugacy (Section~\ref{sec:niw_posterior_conjugacy}), the posterior distribution over $\{\bmu_w, \bSigma_w\}$ remains NIW with updated hyper-parameters:
\begin{equation}\label{equation:gggw_post_niw}
\bmu_w, \bSigma_w\sim \niw(\bmu_w, \bSigma_w\mid \bmm_M, \kappa_M, \nu_M, \bS_M),
\end{equation}
where 
\begin{align}
\bmm_M &= \frac{\kappa_0\bmo + M\overline{\bw}}{\kappa_M} = 
\frac{\kappa_0 }{\kappa_M}\bmo+\frac{M}{\kappa_M}\overline{\bw} , \label{equation:niw_posterior_equation_2_gggw}\\
\kappa_M  &= \kappa_0 + M,  \label{equation:niw_posterior_equation_3_gggw}\\
\nu_M     &=\nu_0 + M,  \label{equation:niw_posterior_equation_4_gggw}\\
\bS_M     
&=\bS_0 + \bS_{\overline{w}} + \frac{\kappa_0 M}{\kappa_0 + M}(\overline{\bw} - \bmo)(\overline{\bw} - \bmo)^\top \label{equation:niw_posterior_equation_5_gggw}\\
&=\bS_0 + \sum_{m=1}^M \bw_m \bw_m^\top + \kappa_0 \bmo \bmo^\top - \kappa_M \bmm_M \bmm_M^\top,  \label{equation:niw_posterior_equation_6_gggw}\\
\overline{\bw} &= \frac{1}{M}\sum_{m=1}^{M} \bw_m, \label{equation:niw_posterior_equation_7_gggw}\\
\bS_{\overline{w}} &= \sum_{m=1}^{M} (\bw_m - \overline{\bw})(\bw_m - \overline{\bw})^\top. \label{equation:niw_posterior_equation_8_gggw}
\end{align}
An intuitive interpretation for the parameters in NIW can be obtained from the updated parameters above. 
The parameter $\nu_0$ acts as a prior pseudo-count for the covariance estimation; thus, $\nu_M =\nu_0 + M$ reflects the total (posterior) effective sample size.
The posterior mean $\bmm_M$ of the model mean $\bmu_w$ is a weighted average of the prior mean $\bmm_0$ and the empirical sample mean $\overline{\bw}$. The posterior scale matrix $\bS_M$ combines the prior scale matrix $\bS_0$, empirical covariance matrix $\bS_{\overline{w}}$, and an additional term accounting for uncertainty in the mean estimate.
By symmetry, identical updates apply to $\{\bmu_z, \bSigma_z\}$ using the $N$ columns of $\bZ$.

\paragrapharrow{Gibbs sampling.}
A Gibbs sampler for the GGGW model is outlined in  Algorithm~\ref{alg:gggw_gibbs_sampler}. By default, uninformative hyper-parameters are $\alpha_\sigma=\beta_\sigma=1$, $\bmm_0=\bzero, \kappa_0=1, \nu_0=K+1, \bS_0=\bI$.
(Note that $\nu_0=K+1$ ensures the inverse-Wishart prior is proper; see Definition~\ref{definition:multi_inverse_wishart}.)

\begin{algorithm}[h] 
\caption{
Gibbs sampler for GGGW model in one iteration.
A variance prior on $\sigma^2$ is used; an equivalent formulation exists using precision $\tau=1/\sigma^2$. 
By default, uninformative hyper-parameters are $\alpha_\sigma=\beta_\sigma=1$, $\bmm_0=\bzero, \kappa_0=1, \nu_0=K+1, \bS_0=\bI$.
} 
\label{alg:gggw_gibbs_sampler}  
\begin{algorithmic}[1] 
\Require Choose initial $\alpha_\sigma, \beta_\sigma, \bmm_0, \kappa_0, \nu_0, \bS_0$;
\For{$m=1$ to $M$}
\State Sample $\bw_m$ from $p(\bw_m\mid \bmu_w, \bSigma_w )$; \Comment{Equation~\eqref{equa:gggw_prior_hyper1}}
\EndFor
\For{$n=1$ to $N$}
\State Sample $\bz_n$ from $p(\bz_n\mid\bmu_z, \bSigma_z)$; 
\Comment{Symmetry of Eq.~\eqref{equa:gggw_prior_hyper2}}
\EndFor
\State Sample $\bmu_w, \bSigma_w$ from $p(\bmu_w, \bSigma_w\mid \bmm_M, \kappa_M, \nu_M, \bS_M )$
\Comment{Equation~\eqref{equation:gggw_post_niw}}
\State Sample $\bmu_z, \bSigma_z$ from $p(\bmu_z, \bSigma_z\mid \bmm_N, \kappa_N, \nu_N, \bS_N )$
\Comment{Symmetry of Eq.~\eqref{equation:gggw_post_niw}}
\State Sample $\sigma^2$ from $p(\sigma^2\mid \bA, \bW,\bZ, \alpha_\sigma, \beta_\sigma)$;
\Comment{Equation~\eqref{equation:ggg_posterior_sigma2}}
\State Report loss in Equation~\eqref{equation:frob_loss_brmf}, stop if it converges;
\end{algorithmic} 
\end{algorithm}

\index{Decomposition: GVG}
\section{Gaussian Likelihood with Volume and Gaussian Priors (GVG)}\label{section:gvg_model}
The Gaussian likelihood with volume and Gaussian prior (GVG) model was introduced by \citet{arngren2011unmixing}.
While the original paper applies the volume prior to unmix a set of pixels into \textit{pure spectral signatures (endmembers)} and corresponding \textit{fractional abundances} in hyperspectral image analysis, in which case the factored components are nonnegative. 
However, it can also be applied in the real-valued applications. 

\begin{figure}[h]
\centering  
\subfigtopskip=2pt 
\subfigbottomskip=6pt 
\subfigcapskip=-15pt 
\includegraphics[width=0.421\textwidth]{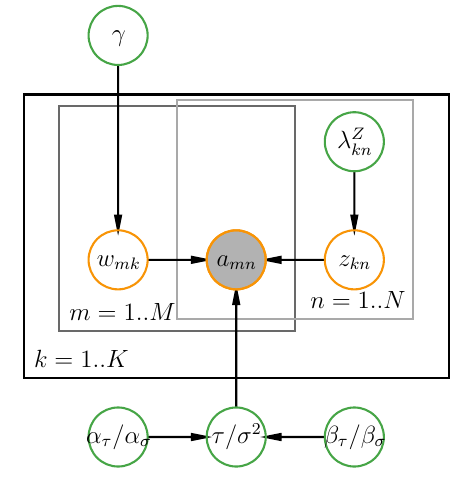}
\caption{Graphical representation of GVG model. Green circles denote prior variables, orange circles represent observed and latent variables  (shaded cycles denote observed variables), and plates indicate repeated variables.
The slash ``/" in the variable represents ``or."}
\label{fig:bmf_gvg}
\end{figure}

In the GVG model, the prior over the abundance matrix $\bZ$ remains Gaussian---identical to the one used in the GGG model. 
In contrast, instead of placing a Gaussian prior on the endmember matrix $\bW$, \citet{arngren2011unmixing} introduce a volume-promoting prior with density proportional to $\bW\propto \exp\{-\gamma \det(\bW^\top\bW)\}$, as illustrated in Figure~\ref{fig:bmf_gvg}. 
This prior encourages the columns of $\bW$ to span a large-volume simplex, which helps promote diversity among the inferred endmembers. The model includes a single hyper-parameter $\gamma>0$, which must be set manually.

\paragrapharrow{Posterior inference for $w_{mk}$.}
To derive the full conditional posterior of an individual entry  $w_{mk}$, we define the following auxiliary quantities:
let $\bw_{m,-k}\in \real^{K-1}$ denote the $m$-th row of $\bW$ excluding column $k$; 
vector $\bw_{-m, k}\in\real^{M-1}$ denote the $k$-th column of $\bW$ excluding row $m$; 
matrix $\bW_{-m,-k}\in \real^{(M-1)\times (K-1)}$ denote $\bW$ with row $m$ and column $k$ removed; 
matrix $\bW_{:,-k}\in \real^{M\times (K-1)}$ denote $\bW$ with column $k$ removed; 
scalar value $D_{-k,-k}=\det\big( \bW_{:,-k}^\top \bW_{:,-k} \big)$; 
and the \textit{matrix adjugate} of $\big( \bW_{:,-k}^\top \bW_{:,-k} \big)$ as
$\bA_{-k,-k}=\det\big( \bW_{:,-k}^\top \bW_{:,-k} \big)\big( \bW_{:,-k}^\top \bW_{:,-k} \big)^{-1}\in\real^{(K-1)\times(K-1)}$. 
Using these, the conditional posterior of $w_{mk}$ is Gaussian:
\begin{equation}\label{equation:gvg_post_w}
w_{mk} \sim \normal(w_{mk}\mid  \widetilde{\mu_{mk}}, \widetilde{\sigma_{mk}^{2}}),
\end{equation}
where 
$$
\begin{aligned}
\widetilde{\mu_{mk}}&= 
\widetilde{\sigma_{mk}^{2}} 
\Bigg\{
\gamma \bw_{m,-k}^\top \bA_{-k,-k}(\bW_{-m,-k}^\top)\bw_{-m,k} + 
\frac{1}{\sigma^2} \sum_{j=1}^{N} (a_{mj} - \sum_{i\neq k}^{K} \bw_{m,-k}^\top \bz_{j, -k})z_{kj}
\Bigg\},
\end{aligned}
$$
and
$$
\begin{aligned}
\widetilde{\sigma_{mk}^{2}}&=
1\bigg/ \Bigg(	\frac{1 }{\sigma^2} \sum_{j=1}^{N} z_{kj}^2 + 
\gamma\left(D_{-k,-k} - \bw_{m,-k}^\top \bA_{-k,-k}\bw_{m,-k}\right)
 \Bigg).
\end{aligned}
$$
This result follows from isolating $w_{mk}$ in the determinant term $\exp\{-\gamma \det(\bW^\top\bW)\}$ using the block-matrix determinant identity:
$\det(\bM) = \det(\bD)\det(\bA-\bB\bD^{-1}\bC) = \det(\bA)\det(\bD-\bC\bA^{-1}\bB)$
if the  matrix $\bM$ has the block formulation $\bM=
\begin{bmatrixfoot}
\bA & \bB\\
\bC & \bD
\end{bmatrixfoot}$.

\begin{problemset}
\item Suppose we replace the zero-mean, unit-covariance prior in  Equation~\eqref{equation:gggm_prior} of the GGGM model with a general Gaussian distribution $\normal(\bw_m\mid \bmu_m, \bSigma_m)$. 
By appropriately redefining the model parameters, does the marginal distribution $p(\bA)$ over the observed data remain unchanged for any valid choice of $\bmu_m$ and $\bSigma_m$?

\item Use the ``MovieLens 100K" dataset introduced in Section~\ref{section:movie_rec_als} to evaluate and compare the performance of the Bayesian real-valued matrix factorization methods presented in this chapter.

\item \label{problem:symmetry_zkn} Following the derivation in Equation~\eqref{equation:ggg_poster_wmk1}, derive the conditional distribution over the user feature $z_{kn}$ (for all $k\in\{1,2,\ldots, K\}$ and $n\in\{1,2,\ldots, N\}$) in  the GGG model.

\item \label{problem:symmetry_zkngggm}  Similarly, based on the derivation in Equation~\eqref{equation:gggm_wm_post}, obtain the conditional posterior for the full user feature vector $\bz_{n}$ (for all  $n\in\{1,2,\ldots, N\}$) in the GGGM model.

\item Building on the discussions in Sections~\ref{section:brmf_ggga} and~\ref{section:bayes_appetizers}, derive the automatic relevance determination (ARD) results for linear regression models.

\item We derived variational inference (VI) updates for the GGG model; see Section~\ref{section:markov-blanket}. 
Extend this derivation to obtain VI update equations for the GGGM, GGGA, GGGW, and GVG models.

\item Verify Equation~\eqref{equation:gvg_post_w} rigorously by explicitly deriving the conditional posterior of  $w_{mk}$  under the volume prior.
\end{problemset}


\chapter{Bayesian Nonnegative Matrix Factorization}\label{chapter:bnmf}
\begingroup
\hypersetup{
	linkcolor=structurecolor,
	linktoc=page,  
}
\minitoc \newpage
\endgroup

\section{Introduction}
\lettrine{\color{caligraphcolor}T}
The nonnegative matrix factorization (NMF) method is used to analyze data matrices whose entries are nonnegative---a common characteristic of datasets derived from text and images \citep{berry2007algorithms}; see Chapter~\ref{chapter:nmf}.
In cases where the entries in $\bA, \bW$, and $\bZ$ are nonnegative, NMF algorithms have frequently improved performance. 
Thus, the scope of NMF research has grown rapidly in recent years, particularly in the fields of machine learning \citep{lee1999learning, lee2000algorithms}. 

Early work on nonnegative matrix factorization was carried out in the 1990s by a Finnish research group under the name \textit{positive matrix factorization} \citep{paatero1991matrix, paatero1994positive, anttila1995source}. This body of work is seldom cited by later researchers, partly due to the misleading term ``positive matrix factorization"---despite the fact that \citet{paatero1994positive} actually developed a nonnegative matrix factorization method.
Since its popularization by \citet{lee1999learning, lee2000algorithms}, the NMF problem has attracted considerable attention, both in published and unpublished research, across diverse fields such as science, engineering, and medicine. Various authors have also proposed alternative formulations of the NMF problem \citep{schmidt2009bayesian, tan2013automatic, brouwer2017prior, lu2022flexible}.

Formally, the NMF problem can be expressed as $\bA=\bW\bZ+\bE$, where a data matrix $\bA\in \real^{M\times N}$ is approximately factorized into two nonnegative matrices: $M\times K$ nonnegative matrix $\bW\in \real_+^{M\times K}$ and a  $K\times N$ nonnegative matrix $\bZ\in \real_+^{K\times N}$. 
The data set $\bA$ need not be complete; missing entries can be indicated by a binary mask matrix $\bM\in \real^{M\times N}$, where an entry of 1 denotes an observed value and 0 denotes a missing value.
The nonnegativity constraint renders the resulting factors more interpretable and easier to inspect---especially in applications like image analysis.

To simplify the discussion, we first assume that there are no missing entries.
Handling missing data in the Bayesian NMF framework follows the same approach as in Bayesian RMF (see Section~\ref{section:markov-blanket}).
The goal of NMF is to project each data vector  $\ba_n$ into a lower-dimensional representation 
 $\bz_n \in \real^K$, where  $K<M$,
such that the reconstruction error---measured by the Frobenius norm---is minimized (assuming $K$ is known):
\begin{equation}\label{equation:als-per-example-loss_bnmf}
	\mathop{\min}_{\bW,\bZ} L(\bW,\bZ) =  \mathop{\min}_{\bW,\bZ}\sum_{n=1}^N \sum_{m=1}^{M} \left(a_{mn} - \bw_m^\top\bz_n\right)^2.
\end{equation}
Here, $\bW=[\bw_1^\top; \bw_2^\top; \ldots; \bw_M^\top]\in \real^{M\times K}$ and $\bZ=[\bz_1, \bz_2, \ldots, \bz_N] \in \real^{K\times N}$, with $\bw_m$ and $\bz_n$ representing the rows of $\bW$ and the columns of $\bZ$, respectively.

\begin{table}[h]
\centering
\small
\setlength{\tabcolsep}{1.0pt}
\renewcommand{\arraystretch}{1.25}
\begin{tabular}{l|llll}
\hline
Name & Likelihood &Prior $\bW$ & Prior $\bZ$ &  Hierarchical prior \\ 
\hline \hline
GEE & $\normal(a_{mn}|\bw_m^\top\bz_n, \sigma^2)$ & 
$\exponential(w_{mk}|\lambda_{mk}^W)$ & $\exponential(z_{kn}|\lambda_{kn}^Z)$ & \gap\gap\slash \\ \hline
GEEA & $\normal(a_{mn}|\bw_m^\top\bz_n, \sigma^2)$ & 
$\exponential(w_{mk}|\lambda_k)$ & $\exponential(z_{kn}|\lambda_k)$ & $\gammadist(\lambda_k|\alpha_\lambda, \beta_\lambda)$ \\ \hline
GTT& $\normal(a_{mn}|\bw_m^\top\bz_n, \sigma^2)$ &  
$\truncatednormal(w_{mk}|\mu_{mk}^W, \frac{1}{\tau_{mk}^W})$ & $\truncatednormal(z_{kn}|\mu_{kn}^Z, \frac{1}{\tau_{kn}^Z})$ & \gap\gap\slash \\ \hline
GTTN& $\normal(a_{mn}|\bw_m^\top\bz_n, \sigma^2)$ &  
$\truncatednormal(w_{mk}|\mu_{mk}^W, \frac{1}{\tau_{mk}^W})$ & $\truncatednormal(z_{kn}|\mu_{kn}^Z, \frac{1}{\tau_{kn}^Z})$ &
\begin{tabular}{@{}c@{}}$\{\mu_{mk}^W, \tau_{mk}^W\}$, $\{\mu_{kn}^Z, \tau_{kn}^Z\}\sim$ 
\\ $\tnsng(\mu_\mu, \tau_\mu, a, b)$\end{tabular}\\ \hline
GRR & $\normal(a_{mn}|\bw_m^\top\bz_n, \sigma^2)$ & 
$\rectifieddist(\cdot|\mu_{mk}^W, \frac{1}{\tau_{mk}^W}, \lambda_{mk}^W)$ & $\rectifieddist(\cdot|\mu_{kn}^Z, \frac{1}{\tau_{kn}^Z}, \lambda_{kn}^Z)$ & \gap\gap\slash \\ \hline
GRRN & $\normal(a_{mn}|\bw_m^\top\bz_n, \sigma^2)$ & 
$\rectifieddist(\cdot|\mu_{mk}^W, \frac{1}{\tau_{mk}^W}, \lambda_{mk}^W)$ & $\rectifieddist(\cdot|\mu_{kn}^Z, \frac{1}{\tau_{kn}^Z}, \lambda_{kn}^Z)$ & 
\begin{tabular}{@{}c@{}}$\{\mu_{mk}^W, \tau_{mk}^W, \lambda_{mk}^W\}$, \\
	$\{\mu_{kn}^Z, \tau_{kn}^Z, \lambda_{kn}^Z\}\sim$ 
\\ $\rnsng(\mu_\mu, \tau_\mu, a, b, \alpha_\lambda, \beta_\lambda)$\end{tabular}
\\\hline
GL$_1^2$& 
$\normal(a_{mn}|\bw_m^\top\bz_n, \sigma^2)$
 &  
\begin{tabular}{@{}c@{}}
	$\bW\sim\exp$  \\ 
	$\{ \frac{-\lambda^W}{2}\sum_{m}(\sum_{k}w_{mk})^2 \}$
\end{tabular} 
& \begin{tabular}{@{}c@{}}
$\bZ\sim\exp$  \\ 
$\{ \frac{-\lambda^Z}{2}\sum_{n}(\sum_{k}z_{kn})^2 \}$
\end{tabular}
 &  
\gap\gap\slash \\ \hline
GL$_2^2$& 
$\normal(a_{mn}|\bw_m^\top\bz_n, \sigma^2)$
&  
\begin{tabular}{@{}c@{}}
	$\bW\sim\exp$  \\ 
	$\{ \frac{-\lambda^W}{2}\sum_{m}(\sum_{k}w_{mk}^2) \}$
\end{tabular} 
& \begin{tabular}{@{}c@{}}
	$\bZ\sim\exp$  \\ 
	$\{ \frac{-\lambda^Z}{2}\sum_{n}(\sum_{k}z_{kn}^2) \}$
\end{tabular}
&  
\gap\gap\slash \\ \hline
GL$_\infty$& 
$\normal(a_{mn}|\bw_m^\top\bz_n, \sigma^2)$
&  
\begin{tabular}{@{}c@{}}
	$\bW\sim\exp\{ -$  \\ 
	$\lambda^W\sum_{m}(\max_k\abs{w_{mk}}) \}$
\end{tabular} 
& \begin{tabular}{@{}c@{}}
	$\bZ\sim\exp\{ -$  \\ 
$\lambda^Z\sum_{n}(\max_k\abs{z_{kn}}) \}$
\end{tabular}
&  
\gap\gap\slash \\ \hline\hline 
GEG & $\normal(a_{mn}|\bw_m^\top\bz_n, \sigma^2)$ & 
$\exponential(w_{mk}|\lambda_{mk}^W)$ & $\normal(z_{kn}|0, (\lambda_{kn}^Z)^{-1})$ & \gap\gap\slash \\ \hline
GnVG & $\normal(a_{mn}|\bw_m^\top\bz_n, \sigma^2)$ & 
\begin{tabular}{@{}c@{}}$\bW\sim$ 
	\\ $\exp\{-\gamma \bW^\top\bW\}u(\bW)$ \end{tabular}
& $\normal(z_{kn}|0, (\lambda_{kn}^Z)^{-1})$ &   \gap\gap\slash
\\ \hline
\end{tabular}
\caption{Overview of Bayesian nonnegative and semi-nonnegative matrix factorization models.}
\label{table:summ_real_nmf}
\end{table}

\paragrapharrow{Terminology.}
Following the terminology established for Bayesian RMF, Bayesian NMF models are denoted by their density functions, listed in the order: likelihood, priors, and hyperpriors (see Section~\ref{section:bmf_real_intro}).
Table~\ref{table:summ_real_nmf} summarizes the Bayesian models for nonnegative matrix factorization presented in this chapter.

\paragrapharrow{Why Bayesian NMF?}
While classical NMF provides a powerful framework for dimensionality reduction and parts-based $\bW$ and $\bZ$ as fixed parameters to be optimized. This deterministic perspective offers limited flexibility in modeling uncertainty, selecting model complexity, or incorporating prior knowledge---challenges that are especially pronounced in real-world settings with noisy, sparse, or incomplete data.
In contrast, Bayesian NMF treats $\bW$ and $\bZ$ as random variables governed by prior distributions. This probabilistic formulation offers several key advantages.

First, it enables principled handling of uncertainty in both the latent factors and predictions. Rather than returning point estimates, Bayesian inference yields full posterior distributions, which can be used to quantify confidence in reconstructions or downstream decisions---a critical feature in applications such as medical diagnostics or scientific discovery.

Second, Bayesian NMF facilitates automatic complexity control through hierarchical priors (e.g., automatic relevance determination or sparsity-inducing priors). By placing appropriate hyperpriors on model parameters, the effective rank $K$ or sparsity structure can be inferred from the data, reducing the need for ad hoc cross-validation or manual tuning.

Third, the framework naturally accommodates missing data and heterogeneous noise models. As noted earlier, missing entries can be marginalized out within the likelihood, and observation noise can be modeled more realistically (e.g., using Poisson or Bernoulli likelihoods for count or binary data), rather than relying solely on Gaussian assumptions implicit in Frobenius-norm minimization.
The Bayesian approach offers great flexibility. By choosing different prior distributions (e.g., Exponential, Truncated-Normal, Rectified-Normal), one can encode different beliefs or constraints about the structure of the factors $\bW$ and $\bZ$, such as sparsity or smoothness. This allows the model to be tailored to the specific characteristics of the data.

Finally, Bayesian NMF promotes interpretability and robustness by encoding domain knowledge through informative priors (such as smoothness, sparsity, or nonnegativity constraints) while remaining coherent under uncertainty. This makes it particularly well-suited for exploratory analysis in fields like genomics, remote sensing, and topic modeling, where understanding the structure of latent components is as important as predictive accuracy.
For these reasons, the Bayesian treatment of NMF not only generalizes the classical approach but also aligns more closely with the demands of modern data science: uncertainty-aware, adaptive, and interpretable.

\begin{figure}[h]
\centering  
\vspace{-0.35cm} 
\subfigtopskip=2pt 
\subfigbottomskip=2pt 
\subfigcapskip=-5pt 
\subfigure[GEE.]{\label{fig:bmf_gee}
\includegraphics[width=0.421\linewidth]{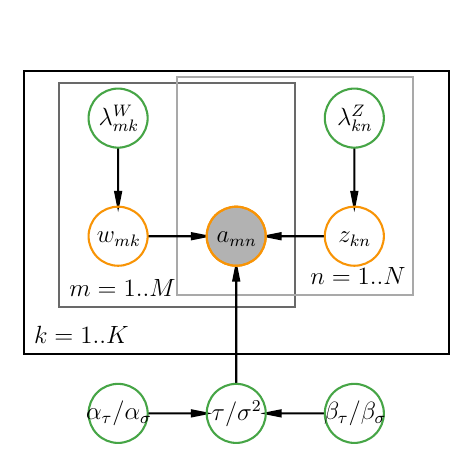}}
\subfigure[GEEA.]{\label{fig:bmf_geea}
\includegraphics[width=0.421\linewidth]{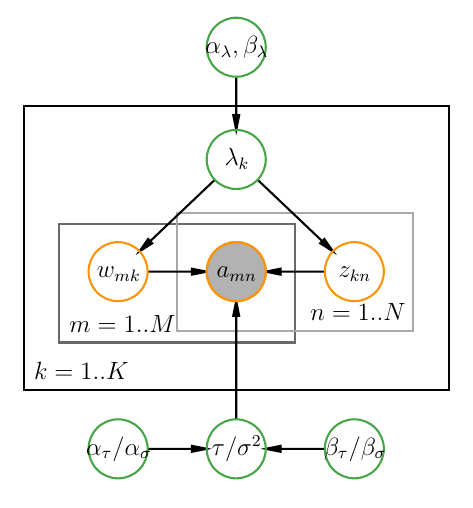}}
\caption{Graphical model representation of GEE and GEEA models. Green circles denote prior variables, orange circles represent observed and latent variables, and plates represent repeated variables. 
The slash ``/" in the variable represents ``or," and the comma ``," in the variable represents ``and."}
\label{fig:bmf_gee_geea}
\end{figure}

\index{Decomposition: GEE}
\section{Gaussian Likelihood with Exponential Priors (GEE)}\label{section:gee_model}
The \textit{Gaussian likelihood with exponential priors (GEE)} model is one of the simplest Bayesian NMF formulations, combining a Gaussian likelihood for the observed data with exponential priors on the factor matrices \citep{schmidt2009bayesian}.
\paragrapharrow{Likelihood.}
Again, we view the data $\bA$ as being produced according to the probabilistic generative process shown in Figure~\ref{fig:bmf_gee}. 
We assume the residuals, $e_{mn}$, are i.i.d. drawn from a  zero-mean Gaussian distribution with variance $\sigma^2$. 
Equivalently, each observed entry $a_{mn}$ is modeled as a Gaussian random variable with  variance $\sigma^2$ and a mean given by the latent decomposition $\bw_m^\top\bz_n$ consistent with the reconstruction error in Equation~\eqref{equation:als-per-example-loss_bnmf}. 
This leads to the following likelihood function:
\begin{equation}\label{equation:gee_likelihood}
\begin{aligned}
p(\bA\mid \btheta) 
&= \prod_{m,n=1}^{M,N} \normal \left(a_{mn}\mid (\bW\bZ)_{mn}, \sigma^2 \right)
= \prod_{m,n=1}^{M,N} \normal \left(a_{mn}\mid (\bW\bZ)_{mn}, \tau^{-1} \right),
\end{aligned}
\end{equation}
where $\btheta=\{\bW,\bZ,\sigma^2\}$ denotes all model parameters, $\sigma^2$ is the variance, $\tau^{-1}=\sigma^2$ is the precision, and  $\normal(x\mid \mu,\sigma^2)$ represents the normal density function.

\paragrapharrow{Prior.}
We treat the latent variables $\{w_{mk}\}$ (and $\{z_{kn}\}$) as random quantities and assign prior distributions to encode structural assumptions---most notably, nonnegativity.  
While other constraints are possible (e.g., semi-nonnegativity \citep{ding2008convex} or discrete support    \citep{gopalan2014bayesian, gopalan2015scalable}),
the GEE model adopts independent exponential priors (Definition~\ref{definition:exponential_distribution}) for all entries of $\bW$ and $\bZ$.
Specifically,
\begin{equation}\label{equation:gee_prior_density_exponential}
\begin{aligned}
w_{mk} &\sim \exponential(w_{mk}\mid \lambda_{mk}^W), 
\gap 
&z_{kn}\sim&  \exponential(z_{kn}\mid \lambda_{kn}^Z);\\
p(\bW) &=\prod_{m,k=1}^{M,K} \exponential(w_{mk}\mid \lambda_{mk}^W), 
\gap 
&p(\bZ) =&\prod_{k,n=1}^{K,N} \exponential(z_{kn}\mid \lambda_{kn}^Z),
\end{aligned}
\end{equation}
where $\exponential(x\mid \lambda)=\lambda\exp(-\lambda x)u(x)$ is the exponential density, and $u(x)$ is the unit step function (equal to 1 for $x\geq 0$ and 0 otherwise). This choice enforces nonnegativity by construction.
As in the GGG model, the noise variance $\sigma^2$ is assigned an inverse-Gamma prior with shape ${\alpha_\sigma}$ and scale ${\beta_\sigma}$ (Definition~\ref{definition:inverse_gamma_distribution}):
\begin{equation}\label{equation:geg_sigma_prior}
p(\sigma^2)= \inversegammadist(\sigma^2\mid \alpha_\sigma, \beta_\sigma) = \frac{{\beta_\sigma}^{\alpha_\sigma}}{\Gamma({\alpha_\sigma})} (\sigma^2)^{-\alpha_\sigma-1} \exp\left( -\frac{{\beta_\sigma}}{\sigma^2} \right).
\end{equation}
By Bayes' rule (Equation~\eqref{equation:posterior_abstract_for_mcmc}), the posterior distribution over the latent variables is proportional to the product of the likelihood and the priors. This posterior can be used for inference via optimization or sampling.

\index{Markov blanket}
\paragrapharrow{Posterior.}
In a Bayesian NMF setting, Markov Chain Monte Carlo (MCMC) methods---particularly Gibbs sampling---require sampling from the full conditional distributions of each latent variable given all others (Section~\ref{section:markov-blanket}). For the GEE model, these conditionals are:
$$
\begin{aligned}
&p(w_{mk}\mid \bA,  \bW_{-mk}, \bZ, \sigma^2,\blambda^W), \\
& p(z_{kn}\mid \bA,  \bW, \bZ_{-kn},  \sigma^2,\blambda^Z ), \\
& p(\sigma^2 \mid \bA, \bW, \bZ, \blambda^W, \blambda^Z ), \\
\end{aligned}
$$
where $\blambda^W$ is an $M\times K$ matrix containing all $\{\lambda_{mk}^W\}$ entries, $\blambda^Z$ is a $K\times N$ matrix including all $\{\lambda_{kn}^Z\}$ values, and $\bW_{-{mk}}$ denotes all elements of $\bW$ except $w_{mk}$. 
Applying Bayes' theorem, the conditional density of $w_{mk}$ depends on its parents ($\lambda_{mk}^W$), children ($a_{mn}$), and co-parents ($\tau$ or $\sigma^2$, $\bW_{-mk}, \bZ$). (See Figure~\ref{fig:bmf_gee} and Section~\ref{section:markov-blanket}.) 
It can be derived as follows:
\begin{equation}\label{equation:gee_poster_wmk1}
\begin{aligned}
&\gap p(w_{mk} \mid \bA ,   \bW_{-mk}, \bZ, \sigma^2, \lambda_{mk}^W ) \\
&\propto p(\bA\mid \bW, \bZ, \sigma^2) \times p(w_{mk}\mid \lambda_{mk}^W)
=\prod_{i,j=1}^{M,N} \normal \left(a_{ij}\mid \bw_i^\top\bz_j, \sigma^2 \right)\times \exponential(w_{mk}\mid \lambda_{mk}^W)\qquad  \\
&\propto \exp\Bigg\{   -\frac{1}{2\sigma^2}  \sum_{i,j=1}^{M,N}(a_{ij} - \bw_i^\top\bz_j  )^2\Bigg\}  \times \cancel{\lambda_{mk}^W }\exp(-\lambda_{mk}^W \cdot w_{mk})u(w_{mk})  \\
\end{aligned}
\end{equation}
\begin{align*}
&\propto \exp\Bigg\{   -\frac{1}{2\sigma^2}  \sum_{j=1}^{N}(a_{mj} - \bw_m^\top\bz_j  )^2\Bigg\}  \cdot  \exp(-\lambda_{mk}^W\cdot w_{mk})u(w_{mk})\\
&\propto \exp\Bigg\{   -\frac{1}{2\sigma^2}  \sum_{j=1}^{N}
\Bigg( w_{mk}^2z_{kj}^2 + 2w_{mk} z_{kj}\bigg(\sum_{i\neq k}^{K}w_{mi}z_{ij} - a_{mj}\bigg)  \Bigg)
\Bigg\}  \cdot \exp(-\lambda_{mk}^W\cdot w_{mk})u(w_{mk})\\
&\propto \exp\Bigg\{   
-\underbrace{\Bigg(\frac{\sum_{j=1}^{N} z_{kj}^2 }{2\sigma^2}  \Bigg) }_{\textcolor{mylightbluetext}{\triangleq 1/(2\widetilde{\sigma_{mk}^{2}})}}
w_{mk}^2 
+
w_{mk}\underbrace{\Bigg( -\lambda_{mk}^W+ \frac{1}{\sigma^2} \sum_{j=1}^{N} z_{kj}\bigg( a_{mj} - \sum_{i\neq k}^{K}w_{mi}z_{ij}\bigg)  \Bigg)}_{\textcolor{mylightbluetext}{\triangleq \widetilde{\sigma_{mk}^{2}}^{-1} \widetilde{\mu_{mk}}}}
\Bigg\}  \cdot u(w_{mk})\\
&\propto   \normal(w_{mk} \mid \widetilde{\mu_{mk}}, \widetilde{\sigma_{mk}^{2}})\cdot u(w_{mk}) 
= \truncatednormal(w_{mk} \mid \widetilde{\mu_{mk}}, \widetilde{\sigma_{mk}^{2}}),
\end{align*}
where $u(x)$ is the unit step function with value 1 if $x\geq 0$ and value 0 if $x<0$, and $\truncatednormal(x \mid \mu, \sigma^2)$ denotes  the \textit{truncated-normal (TN) density} with ``parent" mean $\mu$ and ``parent" variance $\sigma^2$ (Definition~\ref{definition:truncated_normal}).
The posterior ``parent" variance and mean are given by:
\begin{align}
\widetilde{\sigma_{mk}^{2}} &= {\sigma^2}/{(\sum_{j=1}^{N} z_{kj}^2)};
\label{equation:gee_posterior_variance}\\
\widetilde{\mu_{mk}} &= \Bigg( -\lambda_{mk}^W+ \frac{1}{\sigma^2} \sum_{j=1}^{N} z_{kj}\bigg( a_{mj} - \sum_{i\neq k}^{K}w_{mi}z_{ij}\bigg)  \Bigg)\cdot \widetilde{\sigma_{mk}^{2}}.
\label{equation:gee_posterior_mean}
\end{align}
By symmetry, an analogous expression holds for $z_{kn}$ (for all $k=1,2,\ldots,K$ and $n=1,2,\ldots,N$). 
Finally, the conditional posterior for the noise variance $\sigma^2$ is also analytically tractable due to conjugacy (see Equation~\eqref{equation:inverse_gamma_conjugacy_general})).  
It depends on its parents ($\alpha_\sigma$, $\beta_\sigma$), children ($\bA$), and co-parents ($\bW$, $\bZ$). It follows an inverse-Gamma distribution:
\begin{equation}\label{equation:gee_posterior_sigma2}
\begin{aligned}
&
p(\sigma^2 \mid \bA, \bW,\bZ, \alpha_\sigma, \beta_\sigma) = \inversegammadist (\sigma^2\mid \widetilde{\alpha_{\sigma}}, \widetilde{\beta_{\sigma}}), \qquad \\
& \widetilde{\alpha_{\sigma}} = \frac{MN}{2} +{\alpha_\sigma}, \qquad 
\widetilde{\beta_{\sigma}}  =  \frac{1}{2} \sum_{m,n=1}^{M,N} (\bA-\bW\bZ)_{mn}^2 + {\beta_\sigma}.
\end{aligned}
\end{equation}

\index{$\ell_1$-norm}
\index{$\ell_1$ constraint}
\index{Sparsity}
\index{Conjugacy}
\paragrapharrow{Interpretation of the posterior: Sparsity constraint.}
The exponential prior acts as a Bayesian analog of an $\ell_1$-norm penalty, promoting sparsity in the GEE model.
This effect arises from the negative bias term $-\lambda_{mk}^W$ in Equation~\eqref{equation:gee_posterior_mean}:
larger values of $\lambda_{mk}^W$ shift the posterior mean $\widetilde{\mu_{mk}}$ toward zero.
Consequently, samples from $\truncatednormal(w_{mk} \mid \widetilde{\mu_{mk}}, \widetilde{\sigma_{mk}^{2}})$  concentrate near zero, encouraging sparse solutions (see Figure~\ref{fig:dists_truncatednorml_mean}).

\paragrapharrow{Rectified-normal form.}
Alternatively, the posterior of $w_{mk}$ can be expressed as the product of a Gaussian and an exponential density:
$$
\small 
\begin{aligned}
&\gap p(w_{mk}\mid \bA,  \bW_{-mk}, \bZ, \sigma^2, \blambda^W )\\
&\propto \exp\Bigg\{   (\frac{-1 }{2\sigma^2}  \sum_{j=1}^{N} z_{kj}^2)  w_{mk}^2 +
w_{mk}\underbrace{\Bigg( \frac{1}{\sigma^2} \sum_{j=1}^{N} z_{kj}\bigg( a_{mj} - \sum_{i\neq k}^{K}w_{mi}z_{ij}\bigg)  \Bigg)}_{\textcolor{mylightbluetext}{\triangleq \widehat{\sigma_{mk}^{2}}^{-1} \widehat{\mu_{mk}}}}
\Bigg\}   \exp(-\lambda_{mk}^W w_{mk}) u(w_{mk})\\
&\propto \normal(w_{mk} \mid \widehat{\mu_{mk}}, \widehat{\sigma_{mk}^{2}})\cdot \exponential(w_{mk}\mid\lambda_{mk}^W) = \rectifieddist(w_{mk} \mid   \widehat{\mu_{mk}}, \widehat{\sigma_{mk}^{2}}, \lambda_{mk}^W ),
\end{aligned}
$$
where $\widehat{\sigma^2_{mk}}=\widetilde{\sigma_{mk}^{2}}= {\sigma^2}/{(\sum_{j=1}^{N} z_{kj}^2)}$ is the posterior ``parent" variance of the normal distribution with ``parent" mean $\widehat{\mu_{mk}}$, 
$$
\widehat{\mu_{mk}} =  \frac{1}{\sum_{j=1}^{N} z_{kj}^2} \cdot \sum_{j=1}^{N} z_{kj}\bigg( a_{mj} - \sum_{i\neq k}^{K}w_{mi}z_{ij}\bigg), 
$$
and $\rectifieddist(x\mid \mu, \sigma^2, \lambda)\propto \normal(x\mid \mu, \sigma^2)
\exponential(x\mid \lambda)$ denotes the \textit{rectified-normal (RN)} distribution (Definition~\ref{definition:reftified_normal_distribution}).

\begin{algorithm}[h] 
\caption{Gibbs sampler for GEE model in one iteration (prior on variance $\sigma^2$ here, similarly for the precision $\tau$). The procedure presented here may not be efficient but is explanatory. A more  efficient one can be implemented in a vectorized manner. By default,
uninformative hyper-parameters are $\alpha_\sigma=\beta_\sigma=1$, $\{\lambda_{mk}^W\}= \{\lambda_{kn}^Z\}=0.1$.} 
\label{alg:gee_gibbs_sampler}  
\begin{algorithmic}[1] 
\Require Choose initial $\alpha_\sigma, \beta_\sigma, \lambda_{mk}^W, \lambda_{kn}^Z$;
\For{$k=1$ to $K$} 
\For{$m=1$ to $M$}
\State Sample $w_{mk}$ from $p(w_{mk} \mid \bA ,   \bW_{-mk}, \bZ, \sigma^2, \lambda_{mk}^W)$; 
\Comment{Equation~\eqref{equation:gee_poster_wmk1}}
\EndFor
\For{$n=1$ to $N$}
\State Sample $z_{kn}$ from $p(z_{kn} \mid \bA ,   \bW, \bZ_{-kn},\sigma^2 \lambda_{kn}^Z )$; 
\Comment{Symmetry of Eq.~\eqref{equation:gee_poster_wmk1}}
\EndFor
\EndFor
\State Sample $\sigma^2$ from $p(\sigma^2 \mid  \bA, \bW,\bZ,\alpha_\sigma,\beta_\sigma)$; 
\Comment{Equation~\eqref{equation:gee_posterior_sigma2}}
\State Report loss in Equation~\eqref{equation:als-per-example-loss_bnmf}, stop if it converges;
\end{algorithmic} 
\end{algorithm}

\index{Inverse-Gamma distribution}
\index{Rectified-normal distribution}
\paragrapharrow{Gibbs sampling.}
Using the Gibbs sampling framework introduced in Section~\ref{section:gibbs-sampler}, Algorithm~\ref{alg:gee_gibbs_sampler} provides a straightforward (though not optimized) procedure for posterior inference in the GEE model. In practice, it is common to use a shared rate parameter across all entries, i.e.,
 $\lambda=\{\lambda_{mk}^W\} = \{\lambda_{nk}^Z\}$ for all $m,k,n$. 
By default,
uninformative hyper-parameters are $\alpha_\sigma=\beta_\sigma=1$, $\{\lambda_{mk}^W\}= \{\lambda_{kn}^Z\}=0.1$.

\index{Variational Bayesian inference}
\paragrapharrow{Variational Bayesian inference.}
Similar to the GGG model discussed in Section~\ref{section:markov-blanket}, we demonstrate that variational Bayesian inference for the GEE model aligns with the Gibbs sampling approach. This compatibility is also a result of the Gaussian likelihood on the observed elements.
Adopting the mean-field approximation with a Gaussian variational distribution, and in line with Equation~\eqref{equation:em_uncon_mf}, we can deduce the following:
$$
\begin{aligned}
&q_{w_{mk}}(w_{mk}) 
\propto \exp\left\{ \Exp_{q(-w_{mk})} \left[\ln p(\bA \mid \bW, \bZ) +\ln p(\bW, \cancel{\bZ}) \right] \right\}\\
&\propto \exp\Bigg\{ \Exp_{q_{\btheta(-w_{mk})}} \Bigg[\Bigg\{   -\frac{1}{2\sigma^2}  \sum_{j=1}^{N}(a_{mj} - \bw_m^\top\bz_j  )^2\Bigg\} +\ln(\lambda_{mk}^W\exp(-\lambda_{mk}^W \cdot w_{mk})) \Bigg] \Bigg\}u(w_{mk})\\
&\propto \exp\Bigg\{ \Exp_{q_{\btheta(-w_{mk})}} \Bigg[   -\frac{1}{2\sigma^2}  \sum_{j=1}^{N}(a_{mj} - \bw_m^\top\bz_j  )^2  \Bigg] \Bigg\}\cdot \exp(-\lambda_{mk}^W \cdot w_{mk})u(w_{mk})\\
\end{aligned}
$$
$$
\begin{aligned}
&\propto \exp\Bigg\{   -\frac{1}{2\sigma^2}  \sum_{j=1}^{N}
\Bigg( w_{mk}^2z_{kj}^2 + 2w_{mk} z_{kj}\bigg(\sum_{i\neq k}^{K}w_{mi}z_{ij} - a_{mj}\bigg)  \Bigg)
\Bigg\}  \cdot \exp(-\lambda_{mk}^W\cdot w_{mk})u(w_{mk})\\
&\propto \exp\Bigg\{   
-\underbrace{\Bigg(\frac{\sum_{j=1}^{N} z_{kj}^2 }{2\sigma^2}  \Bigg) }_{\textcolor{mylightbluetext}{\triangleq 1/(2\widetilde{\sigma_{mk}^{2}})}}
w_{mk}^2 
+
w_{mk}\underbrace{\Bigg( -\lambda_{mk}^W+ \frac{1}{\sigma^2} \sum_{j=1}^{N} z_{kj}\bigg( a_{mj} - \sum_{i\neq k}^{K}w_{mi}z_{ij}\bigg)  \Bigg)}_{\textcolor{mylightbluetext}{\triangleq \widetilde{\sigma_{mk}^{2}}^{-1} \widetilde{\mu_{mk}}}}
\Bigg\}  \cdot u(w_{mk})\\
&\propto   \normal(w_{mk} \mid \widetilde{\mu_{mk}}, \widetilde{\sigma_{mk}^{2}})\cdot u(w_{mk}) 
= \truncatednormal(w_{mk} \mid \widetilde{\mu_{mk}}, \widetilde{\sigma_{mk}^{2}}),
\end{aligned}
$$
in which the final trio of lines mirrors those found in Equation~\eqref{equation:gee_poster_wmk1}. 
Thus, both inference strategies lead to structurally identical update equations, differing primarily in how expectations are computed (exact vs. approximate).

\index{Decomposition: GEEA}
\index{Automatic relevance determination}
\index{Markov blanket}
\section{GEE Model with ARD Hierarchical Prior (GEEA)}\label{section:geea_nmf_model}
The \textit{Gaussian likelihood with exponential priors and hierarchical prior  (GEEA)} model was first introduced by \citet{tan2013automatic} as an extension of the GEE model.
The key distinction is that GEEA places a hyperprior on the rate parameters of the exponential priors. This hierarchical structure enables \textit{automatic relevance determination (ARD)}, a mechanism that facilitates automatic model selection by adaptively pruning irrelevant latent factors; see Sections~\ref{section:bayespca} and \ref{section:brmf_ggga}.

In the GEEA model, instead of assigning individual rate parameters to each entry of $\bW$ and $\bZ$, a shared rate parameter $\lambda_k$ is used for all elements in the $k$-th column of $\bW$ and the $k$-th row of $\bZ$. In other words, each latent factor $k$ is governed by a single hyper-parameter  $\lambda_k$, which controls the overall scale (or ``relevance") of that factor.

\index{Hyperprior}
\index{Gamma distribution}
\paragrapharrow{Hyperprior.}
Building on the exponential priors in Equation~\eqref{equation:gee_prior_density_exponential}, we place a Gamma hyperprior on each shared rate parameter $\lambda_k$:
$$
w_{mk}\sim \exponential(w_{mk}\mid \lambda_{k}), \gap 
z_{kn}\sim \exponential(z_{kn}\mid \lambda_{k}), \gap 
\lambda_{k} \sim \gammadist(\lambda_{k} \mid \alpha_\lambda, \beta_\lambda), 
$$
where $\lambda_k>0$ is shared across the entire $k$-th component (i.e., shared by all entries in the same column of $\bW$ and the same row of $\bZ$). 
A small value of $\lambda_k$ encourages larger values in the corresponding factor (activating it), whereas a large $\lambda_k$ shrinks the factor toward zero (effectively ``turning it off"; see Figure~\ref{fig:dists_exponential}).
The entire factor $k$ is then either activated if $\lambda_k$ has a low value or ``turned off" if $\lambda_k$ has a high value.
This mechanism allows us to specify an upper bound on the number of latent factors, $K$, without needing to predefine the exact effective rank. The graphical model for GEEA is shown in Figure~\ref{fig:bmf_geea}.

\paragrapharrow{Posterior.}
As in other Bayesian MF models, posterior inference via MCMC requires sampling from the full conditional distributions of all latent variables and hyper-parameters (Section~\ref{section:markov-blanket}). For GEEA, these conditionals are:
$$
\begin{aligned}
&p(\sigma^2 \mid \bA, \bW, \bZ, \blambda  ), &\gap& 
p(w_{mk}\mid \bA,  \bW_{-mk}, \bZ,\sigma^2,  \blambda ), \\
&p(\lambda_k \mid \bW, \bZ, \blambda_{-k},    \alpha_\lambda, \beta_\lambda),  
&\gap& p(z_{kn}\mid\bA , \bW,\bZ_{-kn},  \sigma^2, \blambda ), \\
\end{aligned}
$$
where $\blambda=[\lambda_1, \lambda_2,\ldots,\lambda_K]^\top\in \real_+^K$ is a vector including all $\lambda_k$ values, 
and $\blambda_{-k}$ denotes all components  of $\blambda$ except $\lambda_k$. The posteriors for variables $\{w_{mk}\}$ and $\{z_{kn}\}$ are identical in form to those in the GEE model (Equation~\eqref{equation:gee_poster_wmk1}), except that the individual rates $\{\lambda_{mk}^W\}$ and $\{\lambda_{kn}^Z\}$ are replaced by the shared $\lambda_k$. 
The conditional posterior for $\lambda_k$ follows again from Bayes' theorem.
It depends on its parents ($\alpha_\lambda, \beta_\lambda$), children ($k$-th column $\widetilde{\bw}_k$ of $\bW$, $k$-th row $\widetilde{\bz}_k$ of $\bZ$; note we define $\bw_m$ as the $m$-th row of $\bW$ and $\bz_n$ as the $n$-th column of $\bZ$ in Equation~\eqref{equation:als-per-example-loss_bnmf}), and co-parents (none) \footnote{See Figure~\ref{fig:bmf_geea} and Section~\ref{section:markov-blanket}.}. 
The posterior density of $\lambda_k$ is derived as follows:
\begin{equation}\label{equation:posterior-geea_lambdak}
\begin{aligned}
&\gap p(\lambda_k \mid  \bW, \bZ, \alpha_\lambda, \beta_\lambda)\\
&\propto  p(\widetilde{\bw}_k, \widetilde{\bz}_k \mid \lambda_k) \times p(\lambda_k) = 
\prod_{i=1}^{M}  \exponential(w_{ik}\mid \lambda_{k}) \cdot 
\prod_{j=1}^{N} \exponential(z_{kj}\mid \lambda_{k})
\times \gammadist(\lambda_k \mid \alpha_{\lambda}, \beta_\lambda)\\
&=
\prod_{i=1}^{M}  \lambda_k \exp(-\lambda_k w_{ik}) \cdot 
\prod_{j=1}^{N} \lambda_k \exp(-\lambda_k z_{kj})
\times \frac{\beta_{\lambda}^{\alpha_\lambda}}{\Gamma(\alpha_\lambda)} \lambda_k^{\alpha_\lambda-1} \exp(- \lambda_k \beta_\lambda)\\
&\propto\lambda_k^{M+N +\alpha_\lambda -1} \exp\left\{-\lambda_k \cdot  \left(\sum_{k=1}^{K}(w_{mk}+ z_{kn}) +\beta_\lambda\right)\right\}
\propto \gammadist(\lambda_k \mid \widetilde{\alpha_\lambda}, \widetilde{\beta_\lambda}),
\end{aligned}
\end{equation}
where the updated hyper-parameters are:
$$
\widetilde{\alpha_\lambda} = M+N+\alpha_\lambda, 
\qquad  
\widetilde{\beta_\lambda}=\sum_{k=1}^{K}(w_{mk}+ z_{kn}) +\beta_\lambda.
$$
In this formulation, the prior parameter $\alpha_\lambda$ can be interpreted as the number of pseudo-observations (prior observations), and $\beta_\lambda$ as the sum of the prior observations. Thus, weakly informative (or uninformative) priors can be set using  $\alpha_\lambda=\beta_\lambda=1$.

\paragrapharrow{Gibbs sampling.}
A Gibbs sampler for the GEEA model can be constructed as outlined in Algorithm~\ref{alg:geea_gibbs_sampler}. By default, we use uninformative hyper-parameters:$\alpha_\sigma=\beta_\sigma=1$, $\alpha_\lambda=\beta_\lambda=1$.

\begin{algorithm}[h] 
\caption{Gibbs sampler for GEEA model in one iteration (prior on variance $\sigma^2$ here, similarly for the precision $\tau$). 
The procedure presented here may not be efficient but is explanatory. A more efficient one can be implemented in a vectorized manner. 
By default,
uninformative hyper-parameters are $\alpha_\sigma=\beta_\sigma=1$, $\alpha_\lambda=\beta_\lambda=1$.} 
\label{alg:geea_gibbs_sampler}  
\begin{algorithmic}[1] 
\Require Choose initial $\alpha_\sigma, \beta_\sigma, \alpha_\lambda, \beta_\lambda$;
\For{$k=1$ to $K$} 
\For{$m=1$ to $M$}
\State Sample $w_{mk}$ from $p(w_{mk} \mid \bA,   \bW_{-mk}, \bZ, \sigma^2, \lambda_{k} )$; \Comment{Equation~\eqref{equation:gee_poster_wmk1}}
\EndFor
\For{$n=1$ to $N$}
\State Sample $z_{kn}$ from $p(z_{kn} \mid \bA,   \bW, \bZ_{-kn},\sigma^2, \lambda_{k} )$; 
\Comment{Symmetry of Eq.~\eqref{equation:gee_poster_wmk1}}
\EndFor
\State Sample $\lambda_k$ from $p(\lambda_k \mid  \bW, \bZ, \alpha_\lambda, \beta_\lambda)$; 
\Comment{Equation~\eqref{equation:posterior-geea_lambdak}}
\EndFor
\State Sample $\sigma^2$ from $p(\sigma^2 \mid  \bA, \bW,\bZ,\alpha_\sigma,\beta_\sigma)$; 
\Comment{Equation~\eqref{equation:gee_posterior_sigma2}}
\State Report loss in Equation~\eqref{equation:als-per-example-loss_bnmf}, stop if it converges;
\end{algorithmic} 
\end{algorithm}

\index{Decomposition: GTT}
\index{Truncated-normal distribution}
\section{Gaussian Likelihood with Truncated-Normal Priors (GTT)}
\begin{figure}[h]
\centering  
\vspace{-0.35cm} 
\subfigtopskip=2pt 
\subfigbottomskip=2pt 
\subfigcapskip=-5pt 
\subfigure[GTT.]{\label{fig:bmf_gtt}
\includegraphics[width=0.421\linewidth]{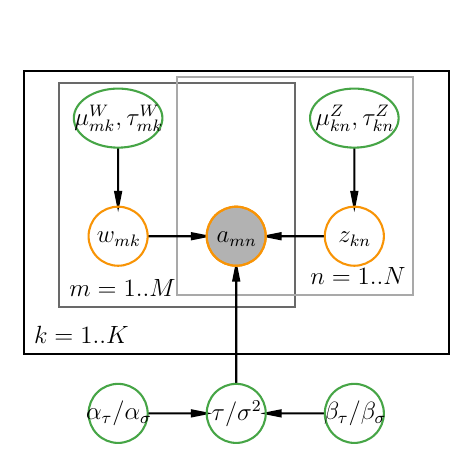}}
\subfigure[GTTN.]{\label{fig:bmf_gttn}
\includegraphics[width=0.421\linewidth]{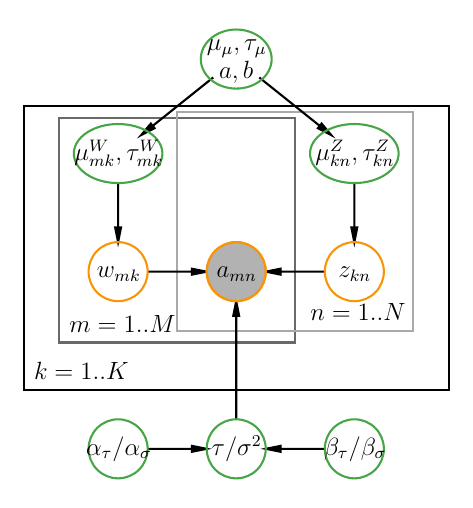}}
\caption{Graphical model representation of GTT and GTTN models. Green circles denote prior variables, orange circles represent observed and latent variables (shaded cycles denote observed variables), and plates represent repeated variables. 
The slash ``/" in the variable represents ``or," and the comma ``," in the variable represents ``and."}
\label{fig:bmf_gtt_gttn}
\end{figure}

The \textit{Gaussian likelihood with truncated-normal priors (GTT)} model was introduced in \citet{brouwer2017prior}, where truncated-normal (TN) priors  
are used over factored matrices (Figure~\ref{fig:bmf_gtt}). 
The truncated-normal distribution, a variant of the normal distribution, excludes values smaller than zero (Definition~\ref{definition:truncated_normal}), allowing it to impose nonnegativity in Bayesian models.
The likelihood function is identical to that of the GEE model (Equation~\eqref{equation:gee_likelihood}). 

\paragrapharrow{Prior.}
We assume that the entries of  $\bW$ and $\bZ$ are independently distributed according to truncated-normal distributions with means and precisions given by $\{\bmu^W ,\btau^W\}$ and  $\{\bmu^Z,\btau^Z\}$, respectively:
\begin{equation}\label{equation:gtt_prior_wmk}
	w_{mk} \sim \truncatednormal(w_{mk} \mid \mu_{mk}^{W}, (\tau_{mk}^W)^{-1} ), \gap 
	z_{kn} \sim \truncatednormal(z_{kn} \mid \mu_{kn}^Z, (\tau_{kn}^Z)^{-1} ),
\end{equation}
where $\bmu^W$ is an $M\times K$ matrix containing all $\{\mu_{mk}^W\}$ entries, $\bmu^Z$ is a $K\times N$ matrix including all $\{\mu_{kn}^Z\}$ values, $\btau^W$ is an $M\times K$ matrix containing all $\{\tau_{mk}^W\}$ entries, and $\btau^Z$ is a $K\times N$ matrix including all $\{\tau_{kn}^Z\}$ values.

\index{Markov blanket}
\paragrapharrow{Posterior.}
Again, following  Bayes' rule and MCMC, this means we need to be able to draw from distributions (by Markov blanket, Section~\ref{section:markov-blanket}):
$$
\begin{aligned}
&p(w_{mk}\mid\bA, \bW_{-mk}, \bZ, \sigma^2, \mu_{mk}^W,\tau_{mk}^W), \\
& p(z_{kn}\mid\bA, \bW, \bZ_{-kn}, \sigma^2,  \mu_{kn}^Z, \tau_{kn}^Z), \\
& p(\sigma^2 \mid \bA, \bW, \bZ, \alpha_\sigma, \beta_\sigma ), \\
\end{aligned}
$$
where $\bW_{-{mk}}$ denotes all elements of $\bW$ except $w_{mk}$, and $\bZ_{-kn}$ denotes all elements of $\bZ$ except $z_{kn}$. 
Using Bayes' theorem, 
the conditional density of $w_{mk}$ depends on its parents ($\mu_{mk}^W$, $\tau_{mk}^W$), children ($a_{mn}$), and co-parents ($\tau$ or $\sigma^2$, $\bW_{-mk}, \bZ$). 
(See Figure~\ref{fig:bmf_gtt} and Section~\ref{section:markov-blanket}.) 
And it can be obtained by (similarly to computing the conditional density of $w_{mk}$ in the GEE model, Equation~\eqref{equation:gee_poster_wmk1})
\begin{equation}\label{equation:gtt_posterior_wmk1}
\small
\begin{aligned}
&\gap p(w_{mk} \mid \sigma^2,   \bW_{-mk}, \bZ, \mu_{mk}^W, \tau_{mk}^W, \bA) 
\propto p(\bA\mid \bW, \bZ, \sigma^2) \cdot p(w_{mk}\mid  \mu_{mk}^W, (\tau_{mk}^W)^{-1})\\
&=\prod_{i,j=1}^{M,N} \normal \left(a_{ij}\mid \bw_i^\top\bz_j, \sigma^2 \right)\times
\truncatednormal(w_{mk} \mid \mu_{mk}^W, (\tau_{mk}^W)^{-1} ) \\
&\propto \exp\bigg\{ - (\frac{ \sum_{j=1}^{N} z_{kj}^2 }{2\sigma^2} +\textcolor{black}{\frac{\tau_{mk}^W}{2}})  w_{mk}^2 +
w_{mk} 
\big\{ \frac{1}{\sigma^2} \sum_{j=1}^{N} z_{kj}( a_{mj}- \sum_{i\neq k}^{K}w_{mi}z_{ij})  + \textcolor{black}{\tau_{mk}^W \mu_{mk}^W} \big\}
\bigg\}  u(w_{mk})\\
&\propto   \normal(w_{mk} \mid \widetilde{\mu_{mk}}, \widetilde{\sigma_{mk}^{2}})\cdot u(w_{mk}) 
= \truncatednormal(w_{mk} \mid \widetilde{\mu_{mk}}, \widetilde{\sigma_{mk}^{2}}),
\end{aligned}
\end{equation}
where $\widetilde{\sigma_{mk}^{2}}= {\sigma^2}/{(\sum_{j=1}^{N} z_{kj}^2 + \tau_{mk}^W\cdot \sigma^2)}$ is the posterior ``parent" variance of the normal distribution with ``parent" mean $\widetilde{\mu_{mk}}$, 
$$
\widetilde{\mu_{mk}} = \Bigg\{ \frac{1}{\sigma^2} \sum_{j=1}^{N} z_{kj}\bigg( a_{mj} - \sum_{i\neq k}^{K}w_{mi}z_{ij}  \bigg)  + \textcolor{black}{\tau_{mk}^W \mu_{mk}^W}\Bigg\}\cdot \widetilde{\sigma_{mk}^{2}}.
$$
By symmetry, an analogous expression holds for variables $\{z_{kn}\}$. 
The conditional posterior for $\sigma^2$ remains identical to that in the GEE model (Equation~\eqref{equation:gee_posterior_sigma2}).

\begin{algorithm}[h] 
\caption{Gibbs sampler for GTT model in one iteration (prior on variance $\sigma^2$ here, similarly for the precision $\tau$). 
The algorithm is explanatory; a vectorized implementation would be more efficient.
By default,
uninformative hyper-parameters are $\alpha_\sigma=\beta_\sigma=1$, $\{\mu_{mk}^W\}=\{\mu_{kn}^Z\}=0$, $\{\tau_{mk}^W\}=\{\tau_{kn}^Z\}=0.1$.} 
\label{alg:gtt_gibbs_sampler}  
\begin{algorithmic}[1] 
\Require Choose initial $\alpha_\sigma, \beta_\sigma, \mu_{mk}^W, \tau_{mk}^W,\mu_{kn}^Z, \tau_{kn}^Z$;
\For{$k=1$ to $K$} 
\For{$m=1$ to $M$}
\State Sample $w_{mk}$ from $p(w_{mk} \mid \bA, \bW_{-mk},\bZ,\sigma^2,\mu_{mk}^W,\tau_{mk}^W)$; \Comment{Equation~\eqref{equation:gtt_posterior_wmk1}}
\EndFor
\For{$n=1$ to $N$}
\State Sample $z_{kn}$ from $p(z_{kn} \mid \bA, \bW,\bZ_{-kn},\sigma^2,\mu_{kn}^Z,\tau_{kn}^Z)$; 
\Comment{Symmetry of Eq.~\eqref{equation:gtt_posterior_wmk1}}
\EndFor
\EndFor
\State Sample $\sigma^2$ from $p(\sigma^2 \mid  \bA, \bW,\bZ,\alpha_\sigma,\beta_\sigma)$;  \Comment{Equation~\eqref{equation:gee_posterior_sigma2}}
\State Report loss in Equation~\eqref{equation:als-per-example-loss_bnmf}, stop if it converges;
\end{algorithmic} 
\end{algorithm}

\paragrapharrow{Gibbs sampling.}
Algorithm~\ref{alg:gtt_gibbs_sampler} outlines a Gibbs sampler for the GTT model. In practice, it is common to use shared hyper-parameters across all entries: 
 $\mu^W=\{\mu_{mk}^W\}'s, \mu^Z=\{\mu_{nk}^Z\}'s$, $\tau^W=\{\tau_{mk}^W\}'s,  \tau^Z=\{\tau_{nk}^Z\}'s$ for all $m,k,n$. 
By default,
uninformative hyper-parameters are $\alpha_\sigma=\beta_\sigma=1$, $\{\mu_{mk}^W\}=\{\mu_{kn}^Z\}=0$, $\{\tau_{mk}^W\}=\{\tau_{kn}^Z\}=0.1$.

\index{Decomposition: GTTN}
\index{Hierarchical prior}
\index{Hyperprior}
\index{Rectified-normal distribution}
\index{TN-scaled-normal-Gamma (TNSNG) prior}
\section{GTT Model with Hierarchical Priors (GTTN)}
A hierarchical prior, known as the TN-scaled-normal-Gamma prior, was originally proposed by \citet{schmidt2009probabilistic} in the context of a rectified-normal distribution, and later adapted to the GTTN model based on the GTT model by \citet{brouwer2017prior}. The key distinction of the GTTN model is that it places a hyperprior on both parameters---the mean and precision---of the truncated-normal distribution (see Figure~\ref{fig:bmf_gttn}).

\paragrapharrow{Hyperprior.}
As shown in Equation~\eqref{equation:conjugate_truncated_nonnegative_mean}, the truncated-normal density serves as a conjugate prior for the nonnegative mean of a Gaussian likelihood, which underlies the GTT model.
If the priors over $\{w_{mk}\}$ and $\{z_{kn}\}$ were standard (untruncated) Gaussians, natural conjugate priors for their means and variances would be the normal-inverse-Gamma or normal-inverse-Chi-square distributions (Equations~\eqref{equation:conjugate_nigamma_general} and~\eqref{equation:nix-posterior}). However, these are not conjugate when the likelihood involves a truncated-normal prior.

\begin{figure}[tp]
\centering  
\vspace{-0.35cm} 
\subfigtopskip=2pt 
\subfigbottomskip=2pt 
\subfigcapskip=-5pt 
\subfigure[GTTN with same hyper-parameters. Same as Figure~\ref{fig:bmf_gttn}.]{\label{fig:bmf_gttnsss}
\includegraphics[width=0.421\linewidth]{./imgs/bmf_gttn.pdf}}
\subfigure[GTTN with different hyper-parameters.]{\label{fig:bmf_gttnddd}
\includegraphics[width=0.421\linewidth]{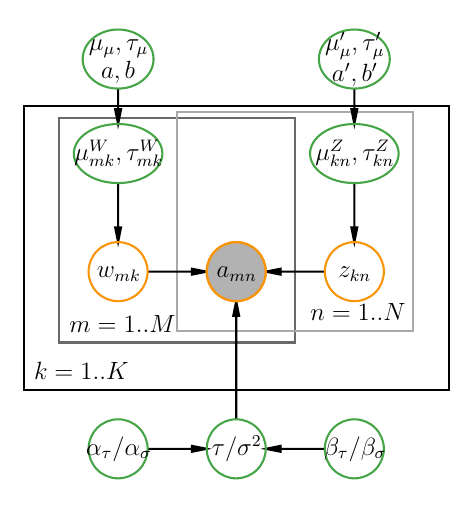}}
\caption{Graphical model representation of GTTN with same and different hyper-parameters. Green circles denote prior variables, orange circles represent observed and latent variables (shaded cycles denote observed variables), and plates represent repeated variables. The slash ``/" in the variable represents ``or," and the comma ``," in the variable represents ``and."}
\label{fig:bmf_gtt_gttn2}
\end{figure}

To address this, we adopt a tailored prior known as the \textit{TN-scaled-normal-Gamma (TNSNG)} distribution (sometimes referred to as a TN-scaled-normal-\textbf{inverse}-Gamma prior for the ``parent" mean and variance parameters)~\footnote{The original formulation in \citet{schmidt2009probabilistic} used a rectified-normal (RN) base. Here, we reinterpret it in the context of the truncated-normal density.}:
$$
\begin{aligned}
\mu_{mk}^W, \tau_{mk}^W \mid \mu_\mu, \tau_\mu, a,  b
&\sim \tnsng(\mu_{mk}^W, \tau_{mk}^W \mid \mu_\mu, \tau_\mu, a,  b)\\
&\propto  \frac{1}{\sqrt{\tau_{mk}^W}} \left(1 - \Phi\big(-\mu_{mk}^W\sqrt{\tau_{mk}^W} \big) \right) \cdot 
\normal(\mu_{mk}^W\mid \mu_\mu, (\tau_\mu)^{-1})\cdot 
\gammadist(\tau_{mk}^W \mid a, b);\\
\mu_{kn}^Z, \tau_{kn}^Z \mid \mu_\mu, \tau_\mu, a,  b
&\sim \tnsng(\mu_{kn}^Z, \tau_{kn}^Z \mid \mu_\mu, \tau_\mu, a,  b).
\end{aligned}
$$
Here, the same hyper-parameters $\{ \mu_\mu, \tau_\mu, a,  b\}$ are typically shared across all entries   $\{u_{mk}^W, \tau_{mk}^W \}$ and $\{u_{kn}^Z, \tau_{kn}^Z\}$. 
However, in certain applications---e.g., when one expects small values in $\bW$ but large values in $\bZ$---distinct hyper-parameter sets can be used for $\bW$ and $\bZ$ (see Figure~\ref{fig:bmf_gtt_gttn2} for a graphical comparison).

It is important to note that the TNSNG prior is not simply the product of independent normal and Gamma distributions. In fact, direct sampling from this joint prior is nontrivial. Nevertheless, its carefully designed form ensures that the full conditional posteriors for variables $\{\mu_{mk}^W\}$ and $\{\tau_{mk}^W\}$ remain analytically tractable---specifically, Gaussian and Gamma, respectively. This decoupling is achieved through the scaling term involving the cumulative distribution function $\Phi(\cdot)$, which compensates for the truncation.

\paragrapharrow{Posterior.}
The posterior distributions for variables $\{w_{mk}\}$, $\{z_{kn}\}$, and $\sigma^2$ are the same as those in the GTT model. The posteriors for $\{u_{mk}^W, \tau_{mk}^W\}$ can be obtained using Bayes' rule, where the conditional density of $\{\mu_{mk}^W, \tau_{mk}^W\}$ depend on their parents ($\mu_\mu, \tau_\mu, a,  b$), children ($w_{mk}$), and co-parents (none). Then it follows from the likelihood in Equation~\eqref{equation:gtt_prior_wmk} that the conditional densities of $\mu_{mk}^W$ is
\begin{equation}\label{equation:posterior_gttn_mu_tau1}
\begin{aligned}
&\gap p(\mu_{mk}^W \mid\textcolor{black}{ \tau_{mk}^W}, w_{mk}, \mu_\mu, \tau_\mu, a,  b)\\
&\propto \truncatednormal(w_{mk} \mid \mu_{mk}^W, (\tau_{mk}^W)^{-1} )
\cdot  \frac{1}{\sqrt{\tau_{mk}^W}} \left(1 - \Phi\big(-\mu_{mk}^W\sqrt{\tau_{mk}^W} \big) \right)  
\normal(\mu_{mk}^W\mid \mu_\mu, \tau_\mu) 
\gammadist(\tau_{mk}^W \mid a, b)\\
&\propto \exp\Bigg\{ -
\underbrace{\frac{\tau_{mk}^W + \tau_\mu}{2}}_{\textcolor{mylightbluetext}{\triangleq \widetilde{t}/2}}
(\mu_{mk}^W)^2 + \mu_{mk}^W 
\underbrace{(\tau_{mk}^W w_{mk} +\tau_{\mu}\mu_\mu)}_{\textcolor{mylightbluetext}{\triangleq \widetilde{m}\cdot  \widetilde{t}}} \Bigg\}
\propto \normal(\mu_{mk}^W \mid \widetilde{m}, \widetilde{t}^{-1}),
\end{aligned}
\end{equation}
where 
$
\widetilde{t}=\tau_{mk}^W + \tau_\mu, 
\widetilde{m}=(\tau_{mk}^W w_{mk} +\tau_{\mu}\mu_\mu)/\widetilde{t}.
$
And the conditional density of $\tau_{mk}^W$ is 
\begin{equation}\label{equation:posterior_gttn_mu_tau2}
\begin{aligned}
&\gap p(\tau_{mk}^W \mid\textcolor{black}{\mu_{mk}^W }, w_{mk}, \mu_\mu, \tau_\mu, a,  b)\\
&\propto \truncatednormal(w_{mk} \mid \mu_{mk}^W, (\tau_{mk}^W)^{-1} )
\cdot  \frac{1}{\sqrt{\tau_{mk}^W}} \left(1 - \Phi\big(-\mu_{mk}^W\sqrt{\tau_{mk}^W} \big) \right)  
\normal(\mu_{mk}^W\mid \mu_\mu, \tau_\mu) 
\gammadist(\tau_{mk}^W \mid a, b)\\
&\propto (\tau_{mk}^W)^{a-1} \exp\left\{-\left( b+ \frac{(w_{mk}-\mu_{mk}^W)^2}{2} \right) \tau_{mk}^W \right\}  
\propto \gammadist(\tau_{mk}^W \mid \widetilde{a}, \widetilde{b}),
\end{aligned}
\end{equation}
where $\widetilde{a} = a, \widetilde{b}= b+ {(w_{mk}-\mu_{mk}^W)^2}/{2}$.
And again due to symmetry, the expressions for $\mu_{kn}^Z$ and $\tau_{kn}^Z$ can be derived accordingly. The Gibbs sampler for the GTTN model is then formulated in Algorithm~\ref{alg:gttn_gibbs_sampler}. By default,
uninformative hyper-parameters are $\alpha_\sigma=\beta_\sigma=1$, $\mu_\mu=0, \tau_\mu=0.1$,  $a=b=1$.

\begin{algorithm}[h] 
\caption{GibbssSampler for GTTN model in one iteration (prior on variance $\sigma^2$ here, similarly for the precision $\tau$). 
The procedure presented here may not be efficient but is explanatory. A more efficient one can be implemented in a vectorized manner. 
By default,
uninformative hyper-parameters are $\alpha_\sigma=\beta_\sigma=1$, $\mu_\mu=0, \tau_\mu=0.1$,  $a=b=1$.} 
\label{alg:gttn_gibbs_sampler}  
\begin{algorithmic}[1] 
\Require Choose initial $\alpha_\sigma, \beta_\sigma, \mu_\mu, \tau_\mu, a,  b$;
\For{$k=1$ to $K$} 
\For{$m=1$ to $M$}
\State Sample $w_{mk}$ from $p(w_{mk} \mid\bA ,   \bW_{-mk}, \bZ, \sigma^2, \mu_{mk}^W, \tau_{mk}^W)$; \Comment{Equation~\eqref{equation:gtt_posterior_wmk1}}
\State Sample $\mu_{mk}^W$ from $p(\mu_{mk}^W \mid \tau_{mk}^W, w_{mk}, \mu_\mu, \tau_\mu, a,  b)$;
\Comment{Equation~\eqref{equation:posterior_gttn_mu_tau1}}
\State Sample $\tau_{mk}^W$ from $p(\tau_{mk}^W \mid\mu_{mk}^W, w_{mk} , \mu_\mu, \tau_\mu, a,  b)$;
\Comment{Equation~\eqref{equation:posterior_gttn_mu_tau2}}
\EndFor
\For{$n=1$ to $N$}
\State Sample $z_{kn}$ from $p(z_{kn} \mid \bA,   \bW, \bZ_{-kn}, \sigma^2,\mu_{kn}^Z, \tau_{kn}^Z )$; 
\Comment{Symmetry of Eq.~\eqref{equation:gtt_posterior_wmk1}}
\State Sample $\mu_{kn}^Z$ from $p(\mu_{kn}^Z \mid \tau_{kn}^Z, z_{kn}, \mu_\mu, \tau_\mu, a,  b)$;
\Comment{Symmetry of Eq.~\eqref{equation:posterior_gttn_mu_tau1}}
\State Sample $\tau_{kn}^Z$ from $p(\tau_{kn}^Z \mid\mu_{kn}^Z, z_{kn} , \mu_\mu, \tau_\mu, a,  b)$;
\Comment{Symmetry of Eq.~\eqref{equation:posterior_gttn_mu_tau2}}

\EndFor
\EndFor
\State Sample $\sigma^2$ from $p(\sigma^2 \mid  \bA, \bW,\bZ,\alpha_\sigma,\beta_\sigma)$;  \Comment{Equation~\eqref{equation:gee_posterior_sigma2}}
\State Report loss in Equation~\eqref{equation:als-per-example-loss_bnmf}, stop if it converges;
\end{algorithmic} 
\end{algorithm}

\index{Decomposition: GRR}
\index{Decomposition: GRRN}
\section{Gaussian Likelihood with RN and Hierarchical Priors (GRR, GRRN)}

\begin{figure}[h]
\centering  
\vspace{-0.35cm} 
\subfigtopskip=2pt 
\subfigbottomskip=2pt 
\subfigcapskip=-5pt 
\subfigure[GRR.]{\label{fig:bmf_grr}
\includegraphics[width=0.421\linewidth]{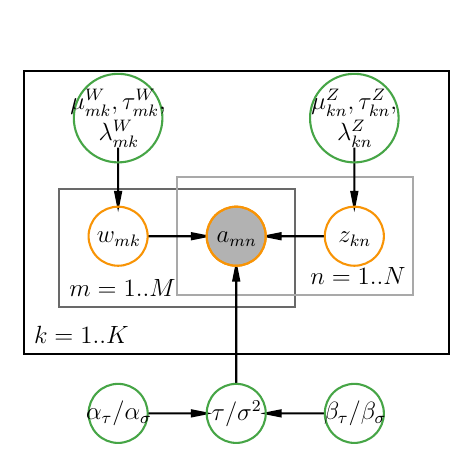}}
\subfigure[GRRN.]{\label{fig:bmf_grrn}
\includegraphics[width=0.421\linewidth]{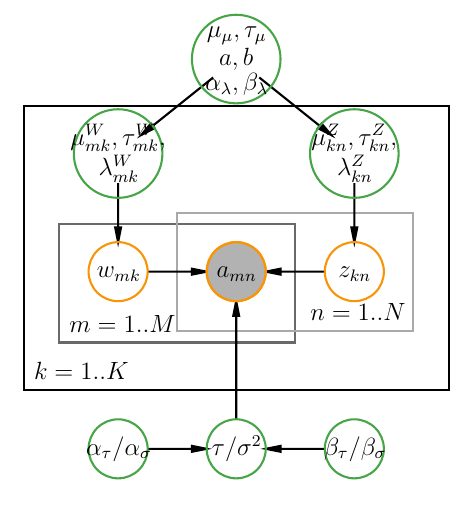}}
\caption{Graphical representation of GRR and GRRN models. Green circles denote prior variables, orange circles represent observed and latent variables (shaded cycles denote observed variables), and plates represent repeated variables. The slash ``/" in the variable represents ``or," and the comma ``," in the variable represents ``and."}
\label{fig:bmf_grr_grrn}
\end{figure}

Going further, \citet{lu2022flexible} propose the Gaussian likelihood with rectified-normal and hierarchical priors---referred to as the GRR and GRRN models---to enhance flexibility beyond the GTT and GTTN frameworks.
In this setting, we again interpret the observed data matrix $\bA$ as generated by the probabilistic process depicted in Figure~\ref{fig:bmf_grrn}. Each entry $a_{mn}$ is modeled using a Gaussian likelihood with variance $\sigma^2$ and mean given by the latent decomposition $\bw_m^\top\bz_n$ (Equation~\eqref{equation:als-per-example-loss_bnmf}). This likelihood matches that used in the GEE model (Equation~\eqref{equation:gee_likelihood}).


\index{Rectified-normal distribution}
\index{Exponentially rectified-normal distribution}
\index{Hierarchical prior}
\index{RN-scaled-normal-Gamma prior}
\paragrapharrow{Prior.}
We treat the latent variables $\{w_{mk}\}$ (and $\{z_{kn}\}$) as random quantities and assign them prior distributions to encode structural assumptions---specifically, nonnegativity in this context.
We assume that each $w_{mk}$ and $z_{kn}$ is independently drawn from a \textit{rectified-normal (RN)} prior (also known as an \textit{exponentially rectified-normal} distribution; see Definition~\ref{definition:reftified_normal_distribution}):
\begin{equation}\label{equation:rn_prior_grrn}
\begin{aligned}
p(w_{mk} \mid \cdot ) &= \rectifieddist(w_{mk} \mid \mu_{mk}^W, (\tau_{mk}^W)^{-1}, \lambda_{mk}^W);\\
p(z_{kn} \mid \cdot ) &= \rectifieddist(z_{kn} \mid \mu_{kn}^Z, (\tau_{kn}^Z)^{-1}, \lambda_{kn}^Z).
\end{aligned}
\end{equation}
This prior enforces nonnegativity on the factor matrices $\bW$ and $\bZ$ and is conjugate to the Gaussian likelihood (Equation~\eqref{equation:conjugate_rectified_nonnegative_mean}).
In principle, distinct RN priors could be used for $\bW$ and  $\bZ$---for example, to encourage sparsity in one factor but not the other. However, we do not consider such asymmetric cases here, as they lie outside the main scope of this book.
Notably, the posterior distribution for each latent variable under this model is a truncated-normal (TN), which is a special case of the rectified-normal (RN) distribution. The resulting model is called the Gaussian likelihood with rectified-normal priors (GRR). Since the RN distribution generalizes the TN, the GRR model reduces to GTT under specific choices of prior parameters. The key advantage of the RN formulation lies in its natural extension to a hierarchical model, which provides principled guidance for selecting prior hyper-parameters---as discussed next.

\paragrapharrow{Hierarchical prior.}
To increase flexibility, we place a joint hyperprior over the RN parameters $\{\mu_{mk}^W, \tau_{mk}^W, \lambda_{mk}^W\}$ in Equation~\eqref{equation:rn_prior_grrn}, namely, the \textit{RN-scaled-normal-Gamma (RNSNG)} prior,
\begin{equation}
\begin{aligned}
&\gap p(\mu_{mk}^W, \tau_{mk}^W, \lambda_{mk}^W \mid\cdot) 
= \rnsng(\mu_{mk}^W, \tau_{mk}^W, \lambda_{mk}^W\mid \mu_\mu, \tau_\mu, a, b, \alpha_\lambda, \beta_\lambda)\\
&=	C(\mu_{mk}^W, \tau_{mk}^W, \lambda_{mk}^W)
\cdot 
\normal(\mu_{mk}^W\mid \mu_\mu, (\tau_\mu)^{-1})
 \cdot \gammadist(\tau_{mk}^W \mid a, b)
\cdot \gammadist(\lambda_{mk}^W \mid \alpha_\lambda, \beta_\lambda),
\end{aligned}
\end{equation}
where $C(\mu_{mk}^W, \tau_{mk}^W, \lambda_{mk}^W)$ is a constant in terms of $\{\mu_{mk}^W, \tau_{mk}^W, \lambda_{mk}^W\}$.
This prior can decouple parameters $\mu_{mk}^W, \tau_{mk}^W$, and $\lambda_{mk}^W$, and their posterior conditional densities are Gaussian, Gamma, and Gamma respectively due to this convenient scale. 
An analogous RNSNG prior is placed over $\{\mu_{kn}^Z, \tau_{kn}^Z, \lambda_{kn}^Z \}$.

\index{Markov blanket}
\paragrapharrow{Posterior.}
Again, following  Bayes' rule and MCMC, this means we need to be able to draw from distributions (by Markov blanket, Section~\ref{section:markov-blanket}):
$$
\begin{aligned}
&p(w_{mk}\mid\bA, \bW_{-mk}, \bZ, \sigma^2, \mu_{mk}^W,\tau_{mk}^W, \lambda_{mk}^W), \\
& p(z_{kn}\mid\bA, \bW, \bZ_{-kn}, \sigma^2,  \mu_{kn}^Z, \tau_{kn}^Z, \lambda_{kn}^Z), \\
& p(\sigma^2 \mid \bA, \bW, \bZ, \alpha_\sigma, \beta_\sigma ), \\
\end{aligned}
$$
where $\bW_{-{mk}}$ denotes all elements of $\bW$ except $w_{mk}$, and $\bZ_{-kn}$ denotes all elements of $\bZ$ except $z_{kn}$. 
Using Bayes' theorem, 
the conditional density of $w_{mk}$ depends on its parents ($\mu_{mk}^W$, $\tau_{mk}^W$, $\lambda_{mk}^W$), children ($a_{mn}$), and co-parents ($\tau$ or $\sigma^2$, $\bW_{-mk}, \bZ$)~\footnote{See Figure~\ref{fig:bmf_grrn} and Section~\ref{section:markov-blanket}.}. 
The conditional density of $w_{mk}$ follows a truncated-normal density. 
And it can be obtained by (similar to computing the conditional density of $w_{mk}$ in the GEE model, Equation~\eqref{equation:gee_poster_wmk1}):
\begin{equation}\label{equation:posterior_grrn_wmk122_app}
\small
\begin{aligned}
&\gap p(w_{mk}\mid \bA, \bW_{-mk}, \bZ,\sigma^2, \mu_{mk}^W,\tau_{mk}^W, \lambda_{mk}^W) 
\propto p(\bA\mid\bW, \bZ, \sigma^2) \times p(w_{mk} \mid \mu_{mk}, \tau_{mk}, \lambda_{mk})\\
&\propto \prod_{i,j=1}^{M,N} \normal(a_{ij}\mid \bw_i^\top\bz_j, \sigma^2) 
\times \rectifieddist(\mu_{mk}, (\tau_{mk})^{-1}, \lambda_{mk})\\
&\stackrel{\star}{\propto} \prod_{i,j=1}^{M,N} \normal(a_{ij}\mid \bw_i^\top\bz_j, \sigma^2) 
\times \truncatednormal\Bigg( \underbrace{\frac{\tau_{mk}^W\mu_{mk}^W - \lambda_{mk}^W}{\tau_{mk}^W}}_{\textcolor{mylightbluetext}{\triangleq \mu^\prime}}, (\tau_{mk})^{-1} \Bigg)\\
&\propto 
\exp\Bigg\{ -\Bigg(\frac{\sum_{j=1}^{N}z_{kj}^2 }{2\sigma^2} + \frac{\tau_{mk}^W}{2}\Bigg)w_{mk}^2 + w_{mk} \Bigg(
\frac{1}{\sigma^2} \sum_{j=1}^{N}z_{kj}(a_{mj}- \sum_{i\neq k}^{K}w_{mk}z_{ij}) + \tau_{mk}^W \mu^\prime 
\Bigg)  \Bigg\} u(w_{mk})\\
& \propto \normal(w_{mk}\mid \widetilde{\mu_{mk}} , \widetilde{\sigma_{mk}^2})u(w_{mk}) 
= \truncatednormal(w_{mk}\mid \widetilde{\mu_{mk}} , \widetilde{\sigma_{mk}^2}),
\end{aligned}
\end{equation}
where the equality $(\star)$ follows from the equivalence between the RN and TN distributions (Definition~\ref{definition:reftified_normal_distribution}), $\widetilde{\sigma_{mk}^2} = {\sigma^2}/{( \sum_{j=1}^{N} z_{kj}^2 + \tau_{mk}^W \cdot \sigma^2 )}$ is the posterior ``parent" variance of the normal distribution with posterior ``parent" mean 
$$
\widetilde{\mu_{mk}} = 
\Bigg(
\frac{1}{\sigma^2} \sum_{j=1}^{N}z_{kj}(a_{mj}- \sum_{i\neq k}^{K}w_{mk}z_{ij}) + \tau_{mk}^W \mu^\prime 
\Bigg)\cdot \widetilde{\sigma_{mk}^2}.
$$
The quantity $\mu^\prime = {(\tau_{mk}^W\mu_{mk}^W - \lambda_{mk}^W)}/{\tau_{mk}^W}$ is the ``parent" mean of the truncated-normal density. 
Due to symmetry, the conditional posterior for $z_{kn}$ can be  derived similarly.


\paragrapharrow{Extra update for GRRN.}
Following   the graphical representation of the GRRN model in Figure~\ref{fig:bmf_grrn}, 
we also sample the hyper-parameters iteratively:
$$
\begin{aligned}
&p(\mu_{mk}^W \mid\textcolor{black}{ \tau_{mk}^W, \lambda_{mk}^W}, \mu_\mu, \tau_\mu, a,  b,\alpha_\lambda, \beta_\lambda, w_{mk}),\\
&p(\tau_{mk}^W \mid\textcolor{black}{ \mu_{mk}^W, \lambda_{mk}^W}, \mu_\mu, \tau_\mu, a,  b,\alpha_\lambda, \beta_\lambda, w_{mk}),\\
&p(\lambda_{mk}^W \mid\textcolor{black}{ \mu_{mk}^W,\tau_{mk}^W} , \mu_\mu, \tau_\mu, a,  b,\alpha_\lambda, \beta_\lambda, w_{mk}).
\end{aligned}
$$

The conditional density for $\mu_{mk}^W$ is a truncated-normal (a special rectified-normal):
\begin{equation}\label{equation:posterior_grrn_mu_tau1222_app}
\begin{aligned}
&\gap p(\mu_{mk}^W \mid\textcolor{black}{ \tau_{mk}^W, \lambda_{mk}^W}, \mu_\mu, \tau_\mu, a,  b,\alpha_\lambda, \beta_\lambda, w_{mk})\\
&\propto \rectifieddist(w_{mk} \mid \mu_{mk}^W, (\tau_{mk}^W)^{-1}, \lambda_{mk}^W) \cdot 
\rnsng(\mu_{mk}^W, \tau_{mk}^W, \lambda_{mk}^W\mid \mu_\mu, \tau_\mu, a, b, \alpha_\lambda, \beta_\lambda)\\
&\propto \rectifieddist(w_{mk} \mid \mu_{mk}^W, (\tau_{mk}^W)^{-1}, \lambda_{mk}^W)
\cdot  \normal(\mu_{mk}^W\mid \mu_\mu, (\tau_\mu)^{-1})
\cdot \gammadist(\tau_{mk}^W \mid a, b)
\cdot  \gammadist(\lambda_{mk}^W \mid \alpha_\lambda, \beta_\lambda)\\
&= 
\normal(w_{mk}| \mu_{mk}^W, (\tau_{mk}^W)^{-1})\cdot 
\cancel{\exponential(w_{mk}| \lambda_{mk}^W)}
\cdot  \normal(\mu_{mk}^W| \mu_\mu, (\tau_\mu)^{-1})
\cdot \cancel{\gammadist(\tau_{mk}^W | a, b)}
\cdot  \cancel{\gammadist(\lambda_{mk}^W | \alpha_\lambda, \beta_\lambda)}\\
&\propto \normal(w_{mk}\mid \mu_{mk}^W, (\tau_{mk}^W)^{-1})\normal(\mu_{mk}^W\mid \mu_\mu, (\tau_\mu)^{-1})
\propto \normal(\mu_{mk}^W \mid \widetilde{m}, \widetilde{t}^{-1}),
\end{aligned}
\end{equation}
where 
$
\widetilde{t}=\tau_{mk}^W + \tau_\mu$ and $
\widetilde{m}=(\tau_{mk}^W w_{mk} +\tau_{\mu}\mu_\mu)/\widetilde{t}
$
are the posterior precision and mean, respectively. The samples of variables $\{w_{mk}\}$ are nonnegative due to the rectification in the distribution (by exponential distribution with the density). However, this ``parent" mean parameter $\mu_{mk}^W$ is not limited to be nonnegative.

The conditional density for $\tau_{mk}^W$ is a Gamma distribution:
\begin{equation}\label{equation:posterior_grrn_tau_tau1222_app}
\begin{aligned}
&\gap p(\tau_{mk}^W \mid\textcolor{black}{ \mu_{mk}^W, \lambda_{mk}^W}, \mu_\mu, \tau_\mu, a,  b,\alpha_\lambda, \beta_\lambda, w_{mk})\\
&\propto \rectifieddist(w_{mk} \mid \mu_{mk}^W, (\tau_{mk}^W)^{-1}, \lambda_{mk}^W)
\cdot \rnsng(\mu_{mk}^W, \tau_{mk}^W, \lambda_{mk}^W\mid \mu_\mu, \tau_\mu, a, b, \alpha_\lambda, \beta_\lambda) \\
&\propto \rectifieddist(w_{mk} \mid \mu_{mk}^W, (\tau_{mk}^W)^{-1}, \lambda_{mk}^W)
\cdot  \normal(\mu_{mk}^W\mid \mu_\mu, (\tau_\mu)^{-1})
\cdot \gammadist(\tau_{mk}^W \mid a, b)
\cdot  \gammadist(\lambda_{mk}^W \mid \alpha_\lambda, \beta_\lambda)\\
&= 
\normal(w_{mk}| \mu_{mk}^W, (\tau_{mk}^W)^{-1})\cdot 
\cancel{\exponential(w_{mk}| \lambda_{mk}^W)}
\cdot  \cancel{\normal(\mu_{mk}^W| \mu_\mu, (\tau_\mu)^{-1})}
\cdot {\gammadist(\tau_{mk}^W | a, b)}
\cdot  \cancel{\gammadist(\lambda_{mk}^W | \alpha_\lambda, \beta_\lambda)}\\
&\propto \normal(w_{mk}\mid \mu_{mk}^W, (\tau_{mk}^W)^{-1})\gammadist(\tau_{mk}^W \mid a, b)\\
&\propto (\tau_{mk}^W)^{a+\frac{1}{2}-1} \exp\left\{-\left( b+ \frac{(w_{mk}-\mu_{mk}^W)^2}{2} \right) \tau_{mk}^W \right\}  
\propto \gammadist(\tau_{mk}^W \mid \widetilde{a}, \widetilde{b}),
\end{aligned}
\end{equation}
where $\widetilde{a} = a+\frac{1}{2}$ and $\widetilde{b}= b+ {(w_{mk}-\mu_{mk}^W)^2}/{2}$ are the posterior shape and rate parameters, respectively.

Furthermore, the conditional density for $\lambda_{mk}^W$ follows also a Gamma distribution:
\begin{equation}\label{equation:posterior_grrn_lambda122_app}
\begin{aligned}
&\gap p(\lambda_{mk}^W \mid\textcolor{black}{ \mu_{mk}^W,\tau_{mk}^W} , \mu_\mu, \tau_\mu, a,  b,\alpha_\lambda, \beta_\lambda, w_{mk})\\
&\propto \rectifieddist(w_{mk} \mid \mu_{mk}^W, (\tau_{mk}^W)^{-1}, \lambda_{mk}^W)\cdot  
\rnsng(\mu_{mk}^W, \tau_{mk}^W, \lambda_{mk}^W\mid \mu_\mu, \tau_\mu, a, b, \alpha_\lambda, \beta_\lambda) \\
&\propto \rectifieddist(w_{mk} \mid \mu_{mk}^W, (\tau_{mk}^W)^{-1}, \lambda_{mk}^W)
\cdot  \normal(\mu_{mk}^W\mid \mu_\mu, (\tau_\mu)^{-1})
\cdot \gammadist(\tau_{mk}^W \mid a, b)
\cdot  \gammadist(\lambda_{mk}^W \mid \alpha_\lambda, \beta_\lambda)\\
&= 
\cancel{\normal(w_{mk}| \mu_{mk}^W, (\tau_{mk}^W)^{-1})}\cdot 
{\exponential(w_{mk}| \lambda_{mk}^W)}
\cdot  \cancel{\normal(\mu_{mk}^W| \mu_\mu, (\tau_\mu)^{-1})}
\cdot \cancel{\gammadist(\tau_{mk}^W | a, b)}
\cdot  {\gammadist(\lambda_{mk}^W | \alpha_\lambda, \beta_\lambda)}\\
&\propto \exponential(w_{mk}\mid \lambda_{mk}^W)\gammadist(\lambda_{mk}^W \mid \alpha_\lambda, \beta_\lambda) \propto \gammadist(\lambda_{mk}^W \mid \widetilde{\alpha_\lambda}, \widetilde{\beta_\lambda}),
\end{aligned}
\end{equation}
where 
$
\widetilde{\alpha_\lambda}= \alpha_\lambda+1$ and $\widetilde{\beta_\lambda}= \beta_\lambda + w_{mk}.$

\paragrapharrow{Key observations.}
The {importance} of this hierarchical prior becomes evident through the interpretation of its conditional density.
Here, \textbf{the prior parameter $\alpha_\lambda$ can be interpreted as the number of prior observations, and $\beta_\lambda$ as the prior knowledge of $w_{mk}$.}
On the one hand, an uninformative choice for $\alpha_\lambda$ is $\alpha_\lambda=1$. 
On the other hand, if one prefers a sparse decomposition with a larger regularization on the model, $\beta_\lambda$ can be chosen as a small value, e.g., $\beta_\lambda=0.01$; or a large value, e.g., $\beta_\lambda=100$, can be applied since we are in the NMF context, in which case, a large value in $\bW$ will enforce the counterparts in $\bZ$ to have small values. 
While an \textit{uninformative choice} for $\beta_\lambda$ is as follows. Suppose the mean value of all entries of matrix $\bA$ is $m_0$, then $\beta_\lambda$ can be set as $\beta_\lambda=\sqrt{\frac{m_0}{K}}$, where the value $K$ is the latent dimension such that each prior  entry $a_{mn}=\bw_m^\top\bz_n$ is equal to $m_0$.  After developing this hierarchical prior, we realize its similarity with the GTTN model (first introduced in a tensor decomposition context \citep{schmidt2009probabilistic}, and further discussed in \citet{brouwer2017prior}). However, the parameters in conditional densities of the GTTN model lack interpretation and flexibility so that there are no guidelines for parameter tuning when the performance is poor. The  GRRN model, on the other hand, can work well generally when we select the uninformative prior $\beta_\lambda=\sqrt{\frac{m_0}{K}}$; moreover, one can even set $\beta_\lambda=20 \cdot \sqrt{\frac{m_0}{K}}$ or $0.1 \cdot \sqrt{\frac{m_0}{K}}$ if one prefers a larger regularization as mentioned above.

Due to symmetry, the conditional expression for $\mu_{kn}^Z$, $\tau_{kn}^Z$, and $\lambda_{kn}^Z$ can be  derived similarly; and we shall not go into the details.

\paragrapharrow{Gibbs sampling.}
The full procedure is  formulated in Algorithm~\ref{alg:grrn_gibbs_sampler}.
By default, uninformative priors are $\alpha_\sigma=\beta_\sigma=1$, $\mu_\mu =0$, $\tau_\mu=0.1, a=b=1$, $\alpha_\lambda=1, \beta_\lambda = \sqrt{\frac{m_0}{K}}$.


\begin{algorithm}[tb] 
\caption{Gibbs sampler for GRRN in one iteration (prior on variance $\sigma^2$ here, similarly for the precision $\tau$). The procedure presented here may not be efficient but is explanatory. A more efficient one can be implemented in a vectorized manner. By default, uninformative priors are $\alpha_\sigma=\beta_\sigma=1$, $\mu_\mu =0$, $\tau_\mu=0.1, a=b=1$, $\alpha_\lambda=1, \beta_\lambda = \sqrt{\frac{m_0}{K}}$. 
One can even set $\beta_\lambda=20 \cdot \sqrt{\frac{m_0}{K}}$ or $0.1 \cdot \sqrt{\frac{m_0}{K}}$ if one prefers a larger regularization.} 
\label{alg:grrn_gibbs_sampler}  
\begin{algorithmic}[1] 
\State {\bfseries Input:} Choose parameters $\alpha_\sigma, \beta_\sigma, \mu_\mu, \tau_\mu, a,  b, \alpha_\lambda, \beta_\lambda$;
\For{$k=1$ to $K$} 
\For{$m=1$ to $M$}
\State Sample $w_{mk}$ from  
$p(w_{mk} \mid \bA,   \bW_{-mk}, \bZ, \sigma^2, \mu_{mk}^W, \tau_{mk}^W,\lambda_{mk}^W )$; 
\Comment{Equation~\eqref{equation:posterior_grrn_wmk122_app}}

\State Sample $\mu_{mk}^W$ from 
$p(\mu_{mk}^W \mid \tau_{mk}^W,\lambda_{mk}^W, \mu_\mu, \tau_\mu, a,  b,\alpha_\lambda, \beta_\lambda, w_{mk})$;
\Comment{Equation~\eqref{equation:posterior_grrn_mu_tau1222_app}}

\State Sample $\tau_{mk}^W$ from 
$p(\tau_{mk}^W \mid\mu_{mk}^W ,\lambda_{mk}^W, \mu_\mu, \tau_\mu, a,  b,\alpha_\lambda, \beta_\lambda, w_{mk})$;
\Comment{Equation~\eqref{equation:posterior_grrn_tau_tau1222_app}}
\State Sample $\lambda_{mk}^W$ from 
$p(\lambda_{mk}^W \mid{ \mu_{mk}^W,\tau_{mk}^W} , \mu_\mu, \tau_\mu, a,  b,\alpha_\lambda, \beta_\lambda, w_{mk})$; 
\Comment{Equation~\eqref{equation:posterior_grrn_lambda122_app}}
\EndFor
\For{$n=1$ to $N$}
\State Sample $z_{kn}$ from 
$p(z_{kn} \mid \bA,   \bW, \bZ_{-kn},\sigma^2, \mu_{kn}^Z, \tau_{kn}^Z,\lambda_{kn}^Z )$; 
\Comment{Sytry. of  Eq.~\eqref{equation:posterior_grrn_wmk122_app}}
\State Sample $\mu_{kn}^Z$ from 
$p(\mu_{kn}^Z \mid \tau_{kn}^Z,\lambda_{kn}^Z, \mu_\mu, \tau_\mu, a,  b,\alpha_\lambda, \beta_\lambda, z_{kn})$;
\Comment{Sytry. of  Eq.~\eqref{equation:posterior_grrn_mu_tau1222_app}}
\State Sample $\tau_{kn}^Z$ from 
$p(\tau_{kn}^Z \mid\mu_{kn}^Z,\lambda_{kn}^Z , \mu_\mu, \tau_\mu, a,  b,\alpha_\lambda, \beta_\lambda, z_{kn})$;
\Comment{Sytry. of  Eq.~\eqref{equation:posterior_grrn_tau_tau1222_app}}

\State Sample $\lambda_{kn}^Z$ from 
$p(\lambda_{kn}^Z \mid{ \mu_{kn}^Z,\tau_{kn}^Z} , \mu_\mu, \tau_\mu, a,  b,\alpha_\lambda, \beta_\lambda, z_{kn})$; 
\Comment{Sytry. of  Eq.~\eqref{equation:posterior_grrn_lambda122_app}}
\EndFor
\EndFor
\State Sample $\sigma^2$ from $p(\sigma^2 \mid  \bA, \bW,\bZ,\alpha_\sigma,\beta_\sigma)$; 
\Comment{Equation~\eqref{equation:gee_posterior_sigma2}}
\State Report loss in Equation~\eqref{equation:als-per-example-loss_bnmf}, stop if it converges;
\end{algorithmic} 
\end{algorithm}

\paragrapharrow{Computational complexity.}
The adopted Gibbs sampling method for the GRRN model has a  complexity of $\mathcalO(MNK^2)$, where the most costs come from the update on the conditional density of variables $\{w_{mk}\}$ and $\{z_{kn}\}$. In the meantime, all the methods we have introduced in the above sections (GEE, GTT, GTTN) have a  complexity of $\mathcal{O}(MNK^2)$. Compared to the GTTN model, the  GRRN model only has an extra cost on the update of $\lambda_{mk}^W$, which does not constitute the bottleneck of the algorithm.

\noindent\makebox[\textwidth][c]{%
\begin{minipage}{\textwidth}
\begin{minipage}[b]{0.41\textwidth}
\centering
\begin{figure}[H]
\centering  
\subfigtopskip=2pt 
\subfigbottomskip=2pt 
\subfigcapskip=-5pt 
\subfigure{\includegraphics[width=0.451\textwidth]{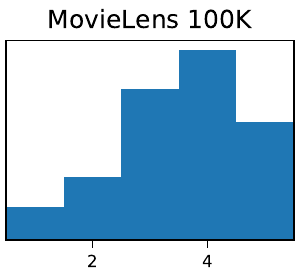} \label{fig:data_movielen100k}}
\subfigure{\includegraphics[width=0.451\textwidth]{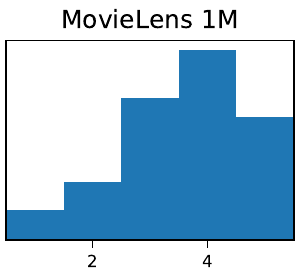} \label{fig:data_movielen1m}}
\caption{Data distribution of MovieLens 100K and MovieLens 1M datasets. 
The MovieLens 1M data set has a higher fraction of users who give a rate of 5 and a lower fraction for rates of 3.}
\label{fig:datasets_nmf}
\end{figure}
\end{minipage}
\hfill
\begin{minipage}[b]{0.57\textwidth}
\centering
\setlength{\tabcolsep}{5pt}
\renewcommand{\arraystretch}{1.34}
\begin{tabular}{llll}
\hline
Data set        & Rows & Columns & Fraction obs. \\ \hline
MovieLens 100K & 943  & 1473    & 0.072         \\ 
MovieLens 1M &6040 &3503& 0.047 \\
\hline
\end{tabular}
\captionof{table}{Data set description. 99,723 and 999,917 observed entries for MovieLens 100K and MovieLens 1M datasets, respectively (user vectors or movie vectors with less than 3 observed entries are cleaned). MovieLens 100K is relatively a small data set and the MovieLens 1M tends to be large; while both of them are sparse.}
\label{table:datadescription}
\end{minipage}
\end{minipage}
}


\subsubsection{Examples}
To demonstrate the main advantages of the introduced GRRN method, 
we conduct experiments across different analysis tasks and datasets, including MovieLens 100K and MovieLens 1M---both widely used benchmarks for movie rating prediction \citep{harper2015movielens}; see also Section~\ref{section:movie_rec_als}.

These datasets contain user ratings on a scale from 1 to 5 stars, with approximately 100,000 and 1,000,000 ratings, respectively. Our goal is to predict missing entries to enable personalized movie recommendations. To ensure data quality, we remove users or movies with fewer than three observed ratings.
A summary of the datasets is provided in Table~\ref{table:datadescription}, and their rating distributions are shown in Figure~\ref{fig:datasets_nmf}. The MovieLens 1M dataset has a higher proportion of 5-star ratings and a lower proportion of 3-star ratings compared to MovieLens 100K. While both datasets are sparse, MovieLens 100K is relatively small, whereas MovieLens 1M is significantly larger---not only in the number of users but also in the number of movies (i.e., feature dimensionality)---making it a more challenging evaluation setting.

Across all experiments, we use the same parameter initialization for fair comparison. We evaluate models in terms of convergence speed and generalization performance. In a wide range of scenarios, GRRN consistently achieves faster convergence and matches or outperforms other Bayesian NMF models in out-of-sample prediction.

\paragrapharrow{Hyper-parameters.}
We adopt the default hyper-parameter settings from \citet{brouwer2017prior}. We use $\{\lambda_{mk}^W\}=\{\lambda_{kn}^Z\}=0.1$ (GEE); $\{\mu_{mk}^Z\}=\{\mu_{kn}^Z\}=0, \{\tau_{mk}^Z\}=\{\tau_{kn}^Z\}=0.1$ (GTT); uninformative $\alpha_\sigma=\beta_\sigma=1$ (Gaussian likelihood in GEE, GTT, GTTN, GRRN); 
$\mu_\mu =0$, $\tau_\mu=0.1, a=b=1$ (hyperprior in GTTN, GRRN); $\alpha_\lambda=1, \beta_\lambda = \sqrt{\frac{m_0}{K}}$ (hyperprior in GRRN). These are very weak prior choices and the models are not sensitive to them \citep{brouwer2017prior}.
As long as the hyper-parameters are set, the observed or unobserved variables are initialized from random draws as this initialization procedure provides a better initial guess of the right patterns in the matrices.
In all experiments, we run the Gibbs sampler 500 iterations with a burn-in of 400 iterations as the convergence analysis shows the algorithm can converge in less than 200 iterations.

\begin{figure*}[h]
\centering  
\subfigtopskip=2pt 
\subfigbottomskip=2pt 
\subfigcapskip=2pt 
\subfigure[Convergence on the \textbf{MovieLens 100K} data set with increasing latent dimension $K$.]{\includegraphics[width=1\textwidth]{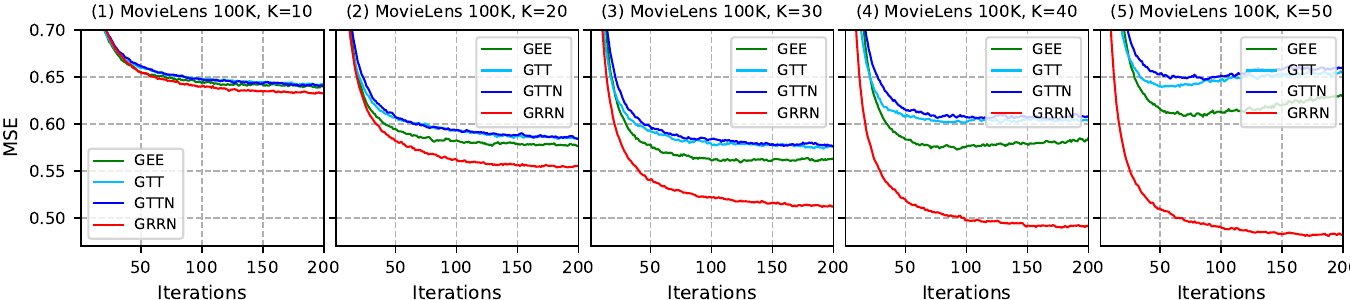} \label{fig:convergences_gdsc_20}}
\subfigure[Convergence on the \textbf{MovieLens 1M} data set with increasing latent dimension $K$.]{\includegraphics[width=1\textwidth]{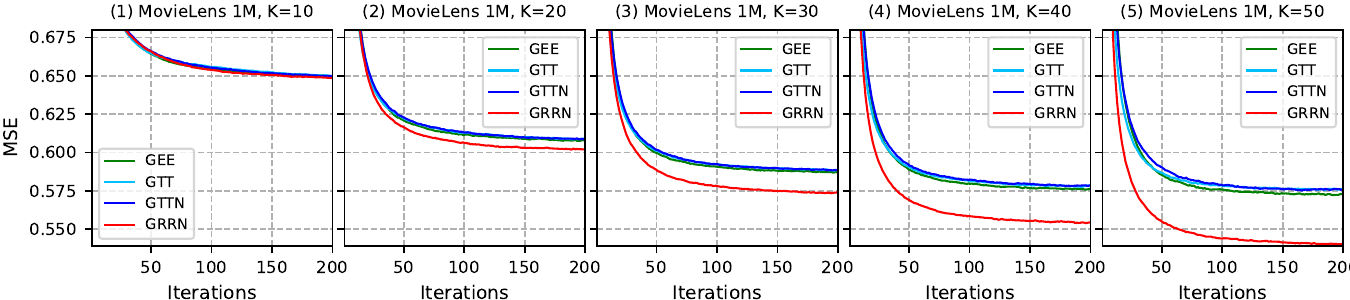} \label{fig:convergences_movielens100k_20}}
\caption{Convergence of the models on the 
MovieLens 100K (upper) and the MovieLens 1M (lower) datasets, measuring the training data fit (mean squared error). When increasing latent dimension $K$, the  GRRN continues to improve the performance; while other models start to decrease on the MovieLens 100K data set or stop increasing on the MovieLens 1M data set.}
\label{fig:convergences_gdsc_movielens100k}
\end{figure*}

\paragrapharrow{Convergence analysis.}
We first compare convergence rates on both MovieLens datasets using latent dimensions $K=10, 20, 30, 40, 50$, with performance measured by mean squared error (MSE). Figure~\ref{fig:convergences_gdsc_movielens100k} shows averaged results over ten independent runs.
On MovieLens 1M, all methods achieve lower MSE as $K$ increases, but GRRN consistently outperforms the others. The GTT and GTTN models yield similar results, as expected---their structures are closely related, with GTTN merely adding a hierarchical layer over GTT.
In contrast, on MovieLens 100K, increasing $K$ leads GEE, GTT, and GTTN to initially improve but then degrade or stagnate---indicating overfitting or optimization difficulties in smaller datasets. GRRN, however, continues to improve steadily, demonstrating superior robustness and making it a better choice for dimensionality reduction in sparse, limited-data regimes.

\paragrapharrow{Noise sensitivity.}
We further evaluate model robustness by adding Gaussian noise at varying signal-to-noise ratios: $\{0\%, 10\%,$ $20\%,$ $50\%, 100\%, 200\%, 500\%, 1000\%\}$, defined as the ratio of added noise variance to data variance. Results for MovieLens 100K with $K=50$ are shown in Figure~\ref{fig:noise_graph_movielens100k}. All models exhibit similar sensitivity to noise. Comparable behavior is observed on MovieLens 1M and across other values of $K$; thus, we omit redundant details.

\noindent\makebox[\textwidth][c]{%
\begin{minipage}{\textwidth}
\begin{minipage}[b]{0.415\textwidth}
\centering
\begin{figure}[H]
\centering  
\subfigtopskip=2pt 
\subfigbottomskip=9pt 
\subfigcapskip=-5pt 
\includegraphics[width=0.9\textwidth]{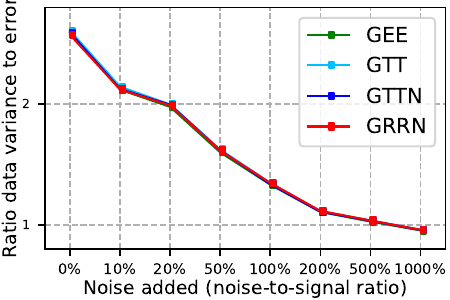}
\caption{Ratio of the variance of data to the MSE of the predictions. The higher the better. All models perform similarly. Similar results can be found on the MovieLens 1M data set and other $K$ values, and we shall not repeat the details.}
\label{fig:noise_graph_movielens100k}
\end{figure}
\end{minipage}
\hfill
\begin{minipage}[b]{0.57\textwidth}
\centering
\renewcommand{\arraystretch}{1.14}
\small
\begin{tabular}{lllll}
\hline
$K$\textbackslash{}Models & GEE & GTT   & GTTN & GRRN   \\ 
\hline
$K$=20   &  1.18   &  1.06    &   1.07  &   \textbf{ 1.02 } \\
$K$=30   &  1.43   &  1.18    &   1.20  &   \textbf{ 1.00 } \\
$K$=40   &  1.86   &  1.42    &   1.45  &   \textbf{ 0.98 } \\
$K$=50   &  2.63   &  1.84    &   1.89  &   \textbf{ 0.97 } \\
\hline
\hline
$K$=20   &  3.47   &  1.46    &   1.57  &   \textbf{ 1.10 } \\
$K$=30   &  6.86   &  2.27    &   2.52  &   \textbf{ 1.05 } \\
$K$=40   &  17056.27   &  4.07    &   4.79  &   \textbf{ 1.04 } \\
$K$=50   &  236750.39   &  2650.21    &   5452.18  &   \textbf{ 1.05 } \\
\hline
\end{tabular}
\captionof{table}{Mean squared error measure when 97\% (upper table) and 98\% (lower table) of data is unobserved for the MovieLens 100K data set. The performance of the  GRRN model exhibits only a marginal deterioration when increasing the fraction of unobserved from 97\% to 98\%.
Similar situations can be observed in the MovieLens 1M experiment.}
\label{table:movielens100k_special_sparsity_case}
\end{minipage}
\end{minipage}
}

%


\begin{figure*}[h]
\centering  
\subfigtopskip=2pt 
\subfigbottomskip=2pt 
\subfigcapskip=-2pt 
\subfigure[Predictive results on the \textbf{MovieLens 100K} data set with increasing fraction of unobserved data and increasing latent dimension $K$.]{\includegraphics[width=1\textwidth]{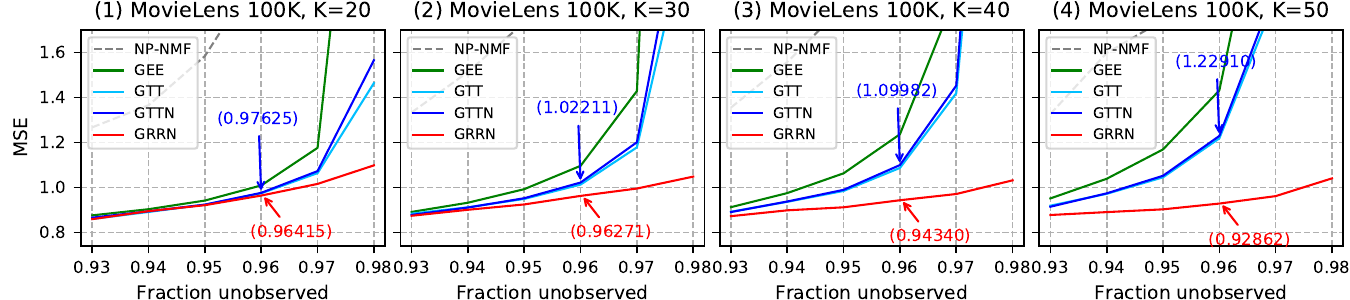} \label{fig:sparsity_movielens_100k_variousK}}
\subfigure[Predictive results on the \textbf{MovieLens 1M} data set  with increasing fraction of unobserved data and increasing latent dimension $K$.]{\includegraphics[width=1\textwidth]{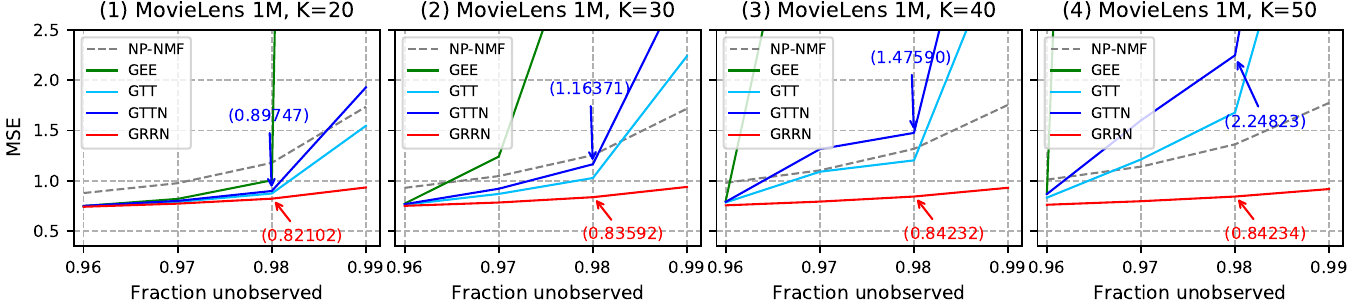} \label{fig:sparsity_movielens_1M_variousK}}
\caption{Predictive results on the MovieLens 100K (upper) and MovieLens 1M (lower) datasets, with the least fractions of unobserved data being 0.928 and 0.953, respectively (see Table~\ref{table:datadescription} for the data description).
We measure the predictive performance (mean squared error) on a held-out data set for different fractions of unobserved data. The blue and red arrows compare the MSEs of GTTN and GRRN models when the fractions of unobserved data are 0.96 and 0.98, respectively. }
\label{fig:sparsity_movielen100_1M}
\end{figure*}

\paragrapharrow{Predictive analysis.}
Motivated by GRRN's strong in-sample convergence, we assess its generalization ability under increasing data sparsity. For each sparsity level, we randomly mask a fraction of observed entries, train on the remaining data, and evaluate MSE on the held-out test set. We vary $K$ from 20 to 50 across all models.
As shown in Figure~\ref{fig:sparsity_movielen100_1M}, when sparsity is moderate (e.g., 93\% unobserved for MovieLens 100K; 96\% for MovieLens 1M) and $K=20$, all models perform similarly---GRRN offers only a slight edge. However, as sparsity increases or $K$ grows, GRRN significantly outperforms all competitors.

Table~\ref{table:movielens100k_special_sparsity_case} quantifies this advantage: when sparsity rises from 97 to 98, GRRN's MSE remains nearly constant ($\approx 1.0$), while other models suffer catastrophic degradation (e.g., GEE's MSE jumps from 2.63 to over 236,750.39 when $K=50$). This confirms that GRRN is far more robust to overfitting.
Notably, although GEE often achieves better in-sample fit (Figure~\ref{fig:convergences_gdsc_movielens100k}), this comes at the cost of poor generalization (Figure~\ref{fig:sparsity_movielen100_1M}). In contrast, GRRN excels in both in-sample and out-of-sample performance, making it a more reliable choice for real-world missing-data prediction.

Finally, we include a popular non-probabilistic NMF (NP-NMF) baseline \citep{lee2000algorithms} (see Chapter~\ref{chapter:nmf}). As shown by the grey curves in Figure~\ref{fig:sparsity_movielen100_1M}, NP-NMF overfits readily, even at low $K$ and moderate sparsity---though the effect is somewhat milder on the larger MovieLens 1M dataset. This underscores the benefit of Bayesian regularization in high-sparsity regimes.

\index{Sparsity}
\index{$\ell_1$-norm}
\index{$\ell_p$-norm}
\section{Priors as Regularization}\label{section:prior_as_reg}
Denoting the prior parameters as $\btheta=\{\bW,\bZ,\sigma^2\}$ and applying Bayes' rule, the posterior distribution is proportional to the product of the likelihood and the prior:
$$
p(\btheta \mid \bA) \propto p(\bA \mid \btheta) \cdot p(\btheta) .
$$
Taking logarithms, the log-posterior becomes
$$
\begin{aligned}
\ln p(\btheta \mid \bA) 
&= \ln p(\bA \mid \btheta) + \ln p(\btheta) + \mathcalC_1  
=\ln \prod_{m,n=1}^{M,N} \normal \left(a_{mn}\mid \bw_m^\top\bz_n, \sigma^2 \right) + \ln p(\bW,\bZ) + \mathcalC_2\\
&=-\frac{1}{2\sigma^2} \left(a_{mn} - \bw_m^\top\bz_n \right)^2  + \ln p(\bW,\bZ) + \mathcalC_3,
\end{aligned}
$$
where $\mathcalC_1,\mathcalC_2$, and $\mathcalC_3$ are constants independent of the parameters. 
The final expression consists of two key components: (1) the negative squared reconstruction error (i.e., the training loss), and (2) a regularization term derived from the prior over the factor matrices $\bW$ and $\bZ$. This prior acts as a regularizer that helps prevent overfitting and improves generalization performance.
More concretely, different choices of priors on $\bW$ correspond to different regularization penalties. In the context of NMF, common regularizers include the following:
\begin{equation}\label{equation:4norms-in-vanilla-bmf}
\begin{aligned}
\ell_1 &= \sum_{m=1}^{M} \sum_{k=1}^{K} w_{mk}, \gap\gap  
&\ell_2^{1/2} &= \sum_{m=1}^{M} \sqrt{\sum_{k=1}^{K}w_{mk}},\\
\ell_1^2 &= \sum_{m=1}^{M} \left(\sum_{k=1}^{K} w_{mk}\right)^2, \gap\gap  
&\ell_2^2 &= \sum_{m=1}^{M}{\sum_{k=1}^{K}w_{mk}^2}.
\end{aligned}
\end{equation}
We note that the $\ell_2^2$-norm 
\footnote{
Strictly speaking, this is not a norm but a squared $\ell_2$-type penalty. 
A norm should satisfy nonnegativity ($\norm{\bA}\geq 0$), positive homogeneity ($\norm{\lambda \bA}=\abs{\lambda}\cdot \norm{\bA}$), and triangle inequality ($\norm{\bA+\bB}\leq \norm{\bA}+\norm{\bB}$) for matrices $\bA,\bB$ and scalar $\lambda$; see \citet{lu2021numerical}.
} 
is equivalent to an independent Gaussian prior (GGG model); the $\ell_1$-norm is equivalent to a Laplace prior  in  real-valued decomposition and is equivalently to an exponential prior (GEE model) in nonnegative matrix factorization, which aligns with the KKT conditions derived in the NMF context (see Section~\ref{section:nmf_anls}).

In the following sections, we discuss several Bayesian NMF models derived from these different priors. The resulting differences in the conditional posterior distribution for a latent variable $\{w_{mk}\}$ are summarized in Table~\ref{table:bnmf_regularizer_posterior}. By symmetry, the conditional posteriors for the variables $\{z_{kn}\}$ take analogous forms.

\begin{table*}[t]
\setlength{\tabcolsep}{2.4pt}
\renewcommand{\arraystretch}{1.5}
\footnotesize 
\begin{tabular}{l|l|l|l}
\hline
& Conditional $w_{mk}$& $\widetilde{\mu_{mk}}$ (mean) & $\widetilde{\sigma_{mk}^{2}}$  (variance)    \\ \hline\hline
GEE         & $\truncatednormal(w_{mk} | \widetilde{\mu_{mk}}, \widetilde{\sigma_{mk}^{2}})$ & $\left( -\lambda_{mk}^W\gap \gap\gap \,\,\,\,\,+ \frac{1}{\sigma^2} \sum_{j=1}^{N} z_{kj}\big( a_{mj} - \sum_{i\neq k}^{K}w_{mi}z_{ij}\big)  \right) \widetilde{\sigma_{mk}^{2}}$                                         & $ \frac{\sigma^2}{\sum_{j=1}^{N} z_{kj}^2}$                                   
\\ \hline
GL$_1^2$    & $\truncatednormal(w_{mk} | \widetilde{\mu_{mk}}, \widetilde{\sigma_{mk}^{2}})$ & $\left( -\lambda^W\textcolor{winestain}{\sum_{j\neq k}^{K}w_{mj}}+\frac{1}{\sigma^2} \sum_{j=1}^{N} z_{kj}\big( a_{mj} - \sum_{i\neq k}^{K}w_{mi}z_{ij}\big)  \right) \widetilde{\sigma_{mk}^{2}}$ & $\frac{\sigma^2}{\sum_{j=1}^{N} z_{kj}^2 +\textcolor{winestain}{\sigma^2\lambda^W}}$ \\ 
\hline
GL$_2^2$    & $\truncatednormal(w_{mk} | \widetilde{\mu_{mk}}, \widetilde{\sigma_{mk}^{2}})$ & $\left(\gap \gap\gap\gap\gap\,\,\,\, \gap\frac{1}{\sigma^2} \sum_{j=1}^{N} z_{kj}\big( a_{mj} - \sum_{i\neq k}^{K}w_{mi}z_{ij}\big)  \right) \widetilde{\sigma_{mk}^{2}}$                                                           & $\frac{\sigma^2}{\sum_{j=1}^{N} z_{kj}^2 +\textcolor{winestain}{\sigma^2\lambda^W}}$ \\ 
\hline
GL$_\infty$ & $\truncatednormal(w_{mk} | \widetilde{\mu_{mk}}, \widetilde{\sigma_{mk}^{2}})$ & $\left(-\textcolor{winestain}{\lambda^W\cdot \indicator(w_{mk})}
\,\,\,\,\,\,\, +\frac{1}{\sigma^2} \sum_{j=1}^{N} z_{kj}\big( a_{mj} - \sum_{i\neq k}^{K}w_{mi}z_{ij}\big)  \right) \widetilde{\sigma_{mk}^{2}}$       & $\frac{\sigma^2}{\sum_{j=1}^{N} z_{kj}^2 }$                                    \\ 
\hline
GL$_{2,\infty}^2$ & $\truncatednormal(w_{mk} | \widetilde{\mu_{mk}}, \widetilde{\sigma_{mk}^{2}})$ & $\left(-\textcolor{winestain}{\lambda^W\cdot \indicator(w_{mk})}
\,\,\,\,\,\,\, +\frac{1}{\sigma^2} \sum_{j=1}^{N} z_{kj}\big( a_{mj} - \sum_{i\neq k}^{K}w_{mi}z_{ij}\big)  \right) \widetilde{\sigma_{mk}^{2}}$       & $\frac{\sigma^2}{\sum_{j=1}^{N} z_{kj}^2 +
\textcolor{winestain}{\sigma^2\lambda^W} 
}$                                    \\ 
\hline
\end{tabular}
\caption{Posterior conditional densities of $w_{mk}$'s for GEE, GL$_1^2$, GL$_2^2$, GL$_\infty$, and GL$_{2,\infty}^2$ models. 
The difference is highlighted in \textcolor{winestain}{red}.	
The conditional densities of $z_{kn}$'s are similar due to their symmetry to $w_{mk}$'s.
$\truncatednormal(x|\mu,\tau^{-1}) =\frac{\sqrt{\frac{\tau}{2\pi}} \exp\{-\frac{\tau}{2} (x-\mu)^2 \} } 
{1-\Phi(-\mu\sqrt{\tau})} u(x)$
is a truncated-normal (TN) density with zero density below $x=0$ and renormalized to integrate to one. $\mu$ and $\tau$ are known as the ``parent" mean and ``parent" precision. $\Phi(\cdot)$ is the cumulative distribution function of standard normal density $\normal(0,1)$.
}
\label{table:bnmf_regularizer_posterior}
\end{table*}

\index{Decomposition: GL$_1^2$}
\section{Gaussian $\ell_1^2$ Norm (GL$_1^2$) Model}

\begin{figure}[h]
\centering  
\subfigtopskip=2pt 
\subfigbottomskip=2pt 
\subfigcapskip=-5pt 
\includegraphics[width=0.421\linewidth]{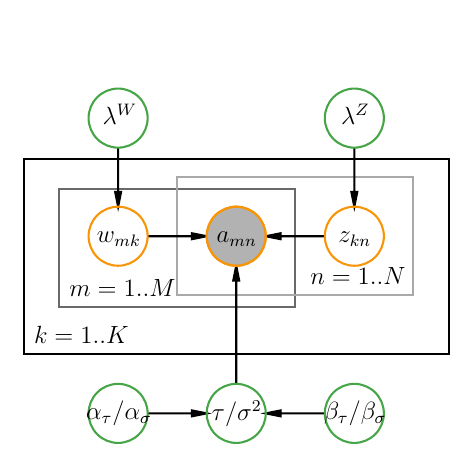}
\caption{Graphical model representation of GL$_1^2$, GL$_2^2$, GL$_\infty$, and GL$_{2,\infty}^2$ models.  
Green circles denote prior variables, orange circles represent observed and latent variables (shaded cycles denote observed variables), and plates represent repeated variables.
The slash ``/" in the variable represents ``or."}
\label{fig:bmf_gl12}
\end{figure}

The \textit{Gaussian $\ell_1^2$-norm  (GL$_1^2$) model}  was introduced by \citet{brouwer2017prior}, based on the $\ell_1^2$ regularization term defined 6in Equation~\eqref{equation:4norms-in-vanilla-bmf}, applied to both factor matrices  $\bW$ and $\bZ$.
As before, we interpret the observed data matrix $\bA$ as generated through the probabilistic graphical model depicted in Figure~\ref{fig:bmf_gl12}.
Specifically, each entry $a_{mn}$ is assumed to follow a Gaussian likelihood with variance $\sigma^2$ and mean given by the low-rank reconstruction $\bw_m^\top\bz_n$, as in Equation~\eqref{equation:als-per-example-loss_bnmf}.

\paragrapharrow{Prior.}
The  $\ell_1^2$ prior follows immediately by replacing the $\ell_1$-norm with the $\ell_1^2$-norm in the exponential prior.
We assume $\bW$ and $\bZ$ are independent, nonnegative, and proportional to an exponential function, with priors governed by hyper-parameters  $\lambda^W$ and $\lambda^Z$, respectively:
\begin{equation}\label{equation:gl12_prior_density}
\begin{aligned}
&p(\bW\mid\lambda^W) &\propto &
\left\{
\begin{aligned}
&\exp \Bigg[  -\frac{\lambda^W}{2} \sum_{m=1}^{M} \Bigg(\sum_{k=1}^{K} w_{mk}\Bigg)^2 \Bigg]
, &\gap &\text{if $w_{mk}\geq 0$ for all $m,k$ };
\\
&0, &\gap &\text{if otherwise};
\end{aligned}
\right.\\
\gap 
&p(\bZ\mid\lambda^Z) &\propto &
\left\{
\begin{aligned}
&\exp \Bigg[  -\frac{\lambda^Z}{2} \sum_{n=1}^{N} \Bigg(\sum_{k=1}^{K} z_{kn}\Bigg)^2 \Bigg]
, &\gap &\text{if $z_{kn}\geq 0$ for all $n,k$ };
\\
&0, &\gap &\text{if otherwise}.
\end{aligned}
\right.\\
\end{aligned}
\end{equation}
As in previous models, the noise variance $\sigma^2={1}/{\tau}$ is assigned  an inverse-Gamma prior  with shape ${\alpha_\sigma}$ and scale ${\beta_\sigma}$, respectively.

\index{Markov blanket}
\paragrapharrow{Posterior.}
Following Bayes' rule and the principles of MCMC, inference proceeds by sampling from the full conditional distributions of each latent variable---i.e., their Markov blankets (see Section~\ref{section:markov-blanket}). These conditionals are:
$$
\begin{aligned}
&p(w_{mk}\mid \bA, \bW_{-mk}, \bZ,\sigma^2, \lambda^W, \lambda^Z ), \\
& p(z_{kn}\mid \bA, \bW,\bZ_{-kn}, \sigma^2, \lambda^W, \lambda^Z), \\
& p(\sigma^2 \mid \bA, \bW, \bZ,\alpha_\sigma, \beta_\sigma ), \\
\end{aligned}
$$
where $\bW_{-{mk}}$ denotes all elements of $\bW$ except $w_{mk}$, and $\bZ_{-kn}$ denotes all entries  of $\bZ$ except $z_{kn}$. 
Using Bayes' theorem, 
the conditional density of $w_{mk}$ depends on its parents ($\lambda^W$), children ($a_{mn}$), and co-parents ($\tau$ or $\sigma^2$, $\bW_{-mk}, \bZ$). 
(See Figure~\ref{fig:bmf_gl12} and Section~\ref{section:markov-blanket}.)
Then, the conditional density of $w_{mk}$ can be obtained by
\begin{equation}\label{equation:gl12_poster_wmk1}
\small
\begin{aligned}
&\gap 
p(w_{mk} |  \bA,   \bW_{-mk}, \bZ, \sigma^2, \lambda^W) 
\propto p(\bA| \bW, \bZ, \sigma^2)  \cdot  p(\bW| \lambda^W)
=\prod_{i,j=1}^{M,N} \normal \left(a_{ij}| \bw_i^\top\bz_j, \sigma^2 \right)
\cdot 
p(\bW| \lambda^W)\\
&\propto \exp\Bigg\{   -\frac{1}{2\sigma^2}  \sum_{i,j=1}^{M,N}(a_{ij} - \bw_i^\top\bz_j  )^2\Bigg\}  \times 
\exp \Bigg\{ -\frac{\lambda^W}{2} \sum_{i=1}^{M} \left(\sum_{j=1}^{K} w_{ij}\right)^2 \Bigg\} \cdot u(w_{mk}) \\
&\propto \exp\Bigg\{   -\frac{1}{2\sigma^2}  \sum_{j=1}^{N}(a_{mj} - \bw_m^\top\bz_j  )^2\Bigg\}  \times   
\exp\left\{ -\frac{\lambda^W}{2} \Bigg( w_{mk} + \sum_{j\neq k }^{K}w_{mj}\Bigg)^2 \right\}\cdot u(w_{mk}) \\
&\propto \exp\Bigg\{   
-
\underbrace{
\bigg(
\frac{\sum_{j=1}^{N} z_{kj}^2 + \textcolor{winestain}{\sigma^2\lambda^W} }{2\sigma^2}  
\bigg)  
}_{\textcolor{mylightbluetext}{\triangleq 1/(2\widetilde{\sigma^2_{mk} }) }} 
w_{mk}^2 
+ w_{mk}\underbrace{
\bigg(
-\lambda^W  \textcolor{winestain}{\sum_{j\neq k}^{K}w_{mj}}+ 
 \sum_{j=1}^{N} \frac{z_{kj}}{\sigma^2}\big( a_{mj} - \sum_{i\neq k}^{K}w_{mi}z_{ij}\big) 
 \bigg)
}_{\textcolor{mylightbluetext}{\triangleq \widetilde{\sigma_{mk}^{2}}^{-1} \widetilde{\mu_{mk}}}}
\Bigg\}  u(w_{mk})\\
&\propto   \normal(w_{mk} \mid \widetilde{\mu_{mk}}, \widetilde{\sigma_{mk}^{2}})\cdot u(w_{mk}) 
= \truncatednormal(w_{mk} \mid \widetilde{\mu_{mk}}, \widetilde{\sigma_{mk}^{2}}),
\end{aligned}
\end{equation}
where $u(x)$ is the unit function equal to 1 if $x\geq 0$ and 0 otherwise;
the quantity $\widetilde{\sigma_{mk}^{2}}= {\sigma^2}/{\big(\sum_{j=1}^{N} z_{kj}^2 +\textcolor{winestain}{\sigma^2\lambda^W}\big)}$ denotes the ``parent" posterior variance of the normal distribution,  
$$
\widetilde{\mu_{mk}} = \Bigg\{ -\lambda^W
\cdot\textcolor{winestain}{\sum_{j\neq k}^{K}w_{mj}}
+
\frac{1}{\sigma^2} \sum_{j=1}^{N} z_{kj}\bigg( a_{mj} - \sum_{i\neq k}^{K}w_{mi}z_{ij}\bigg)  \Bigg\}
\cdot \widetilde{\sigma_{mk}^{2}}
$$
is the ``parent" posterior mean of the normal distribution, and $\truncatednormal(x \mid \mu, \sigma^2)$ is the {truncated-normal density} with ``parent" mean $\mu$ and ``parent" variance $\sigma^2$ Note that this posterior closely resembles that of the GEE model (Equation~\eqref{equation:gee_poster_wmk1}), with differences highlighted in \textcolor{winestain}{red}. 
A side-by-side comparison of conditional posteriors for $w_{mk}$ across models is provided in Table~\ref{table:bnmf_regularizer_posterior}.

Equivalently, the posterior of $w_{mk}$ can be described using a rectified-normal distribution (Definition~\ref{definition:reftified_normal_distribution}); we omit further details here.

By symmetry, the conditional posterior for $z_{kn}$  takes an analogous form. 
Moreover, the posterior for $\sigma^2$ in the GL$_1^2$ model is identical to that in the GEE model (Equation~\eqref{equation:gee_posterior_sigma2}).

\paragrapharrow{Connection between GEE and GL$_1^2$ models.}
Compared to GEE, the GL$_1^2$ model includes an additional term, $\sigma^2\lambda^W$, in the denominator of the posterior ``parent" variance $\widetilde{\sigma_{mk}^{2}}$. 
When all else are equal, this results in a smaller posterior variance, making the distribution more concentrated around its mean. Thus,  GL$_1^2$ imposes a stronger regularization than GEE.

Furthermore, when the hyper-parameters $\{\lambda_{mk}^W\}$ (in GEE) and $\lambda^W$ (in GL$_1^2$) are set to the same value (see Table~\ref{table:bnmf_regularizer_posterior}), the extra term $\sum_{j\neq k}^{K}w_{mj}$ in the posterior ``parent" mean $\widetilde{\mu_{mk}}$ of GL$_1^2$ plays a crucial role in controlling sparsity of factored components in the NMF context:
\begin{itemize}
\item When entries of $\bA$ are large, the sum $\sum_{j\neq k}^{K}w_{mj}$  tends to be greater than 1, which pulls 
$\widetilde{\mu_{mk}}$ toward zero or even negative values. Since the distribution is truncated at zero, this forces samples of $w_{mk}\sim \truncatednormal(w_{mk}\mid\cdot)$ to cluster near zero, promoting sparsity (see Figure~\ref{fig:dists_truncatednorml_mean}: smaller parent means yield smaller expectations under the truncated-normal density).
See also the example in Section~\ref{section:l22_linfty_models} for the experiment on the GDSC $IC_{50}$ data set.
\item Conversely, when entries of $\bA$ are small, this sum is typically less than 1, so $\lambda^W$ has little effect on $\widetilde{\mu_{mk}}$. The posterior mean remains large, leading to denser factor matrices $\bW$ and $\bZ$.
See  Section~\ref{section:l22_linfty_models} for the experiment on the Gene Body Methylation data set.
\end{itemize}
This behavior reveals a key limitation of the GL$_1^2$ model: its sensitivity to the scale of the data matrix $\bA$.
 It is neither consistent nor robust across datasets with different magnitudes (e.g., compare results on the GDSC $IC_{50}$ data in Section~\ref{section:l22_linfty_models} versus Gene Body Methylation data in the same section).
In contrast, the  GL$_2^2$ and GL$_{2,\infty}^2$ models (introduced next) exhibit consistent and robust performance across diverse data types. They also provide stronger regularization than GEE, leading to improved predictive accuracy---particularly when $\bA$ contains large values.

\begin{algorithm}[h] 
\caption{Gibbs sampler for GL$_1^2$ model in one iteration  (prior on variance $\sigma^2$ here, similarly for the precision $\tau$). The procedure presented here may not be efficient but is explanatory. A more efficient one can be implemented in a vectorized manner. By default, uninformative priors are $\alpha_\sigma=\beta_\sigma=1$, $\lambda^W =\lambda^Z=0.1$.} 
\label{alg:gl12_gibbs_sampler}  
\begin{algorithmic}[1] 
\Require Choose initial $\alpha_\sigma, \beta_\sigma, \lambda^W, \lambda^Z$;
\For{$k=1$ to $K$} 
\For{$m=1$ to $M$}
\State Sample $w_{mk}$ from $p(w_{mk} \mid \bA,   \bW_{-mk}, \bZ, \sigma^2, \lambda^W)$; 
\Comment{Equation~\eqref{equation:gl12_poster_wmk1}}
\EndFor
\For{$n=1$ to $N$}
\State Sample $z_{kn}$ from $p(z_{kn} \mid \bA,   \bW, \bZ_{-kn},\sigma^2, \lambda^Z )$; 
\Comment{Symmetry of Equation~\eqref{equation:gl12_poster_wmk1}}
\EndFor
\EndFor
\State Sample $\sigma^2$ from $p(\sigma^2 \mid  \bA, \bW,\bZ,\alpha_\sigma,\beta_\sigma)$; 
\Comment{Equation~\eqref{equation:gee_posterior_sigma2}}
\State Report loss in Equation~\eqref{equation:als-per-example-loss_bnmf}, stop if it converges;
\end{algorithmic} 
\end{algorithm}

\paragrapharrow{Gibbs sampling.}
Using the Gibbs sampling framework outlined in Section~\ref{section:gibbs-sampler}, we implement the above procedure as shown in Algorithm~\ref{alg:gl12_gibbs_sampler}. In practice, we typically use a shared hyper-parameter $\lambda=\lambda^W =\lambda^Z$. 
By default, we adopt weakly informative priors: $\alpha_\sigma=\beta_\sigma=1$, $\lambda^W =\lambda^Z=0.1$.

\begin{SCfigure}
\centering  
\vspace{-0.25cm} 
\subfigtopskip=2pt 
\subfigbottomskip=9pt 
\subfigcapskip=0pt 
\includegraphics[width=0.601\textwidth]{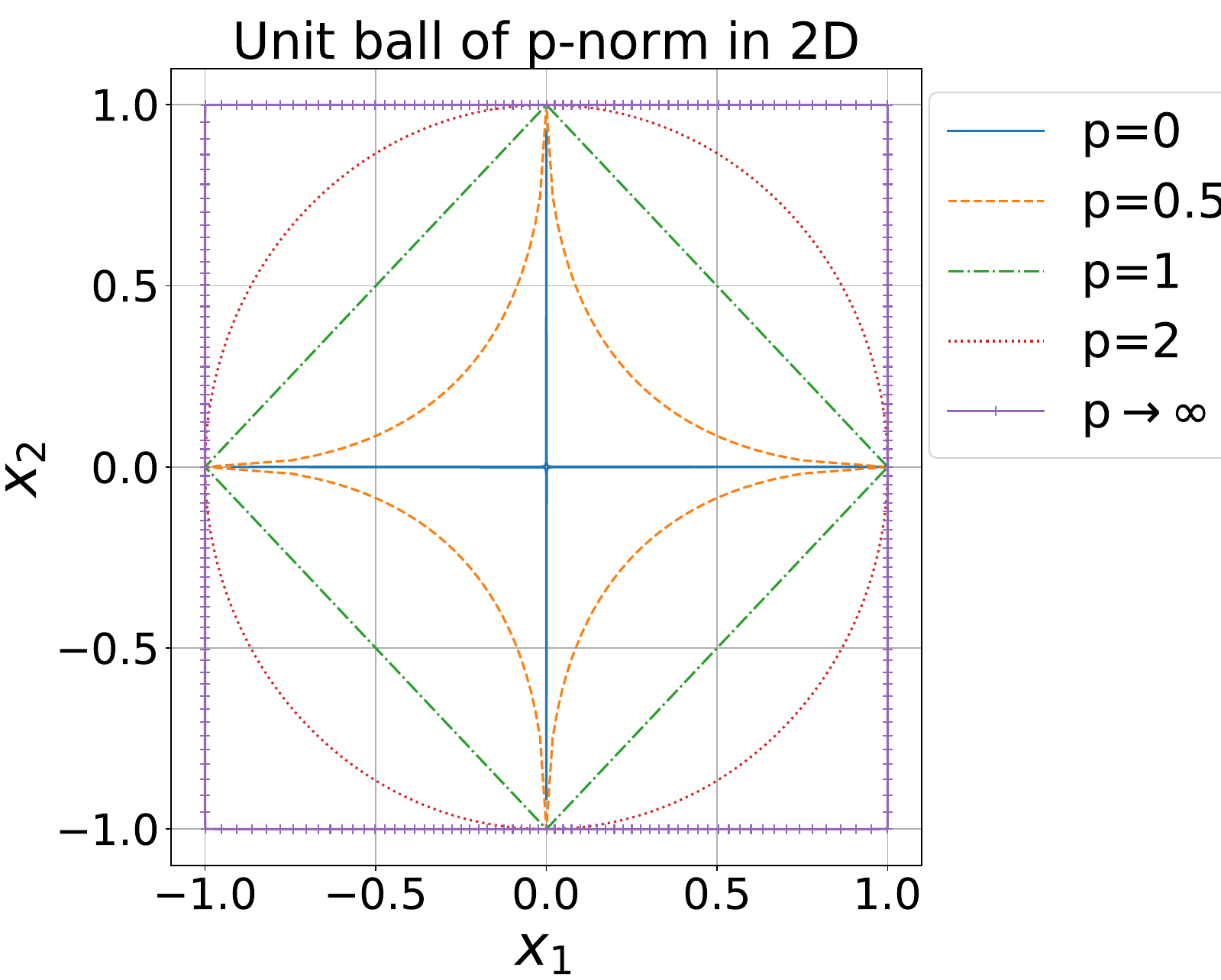}
\caption{Unit ball of $\ell_p$-norm in two-dimensional space. The $\ell_p$-norm over a vector $\bx\in\real^N$ is defined as $\ell_p(\bx) =\normp{\bx} = (\sum_{n} \abs{x_n}^p)^{1/p}$. 
For $p<1$, this does not satisfy the triangle inequality and thus is not a true norm.}
\label{fig:p-norm-2d}
\end{SCfigure}

\begin{figure}[H]
\centering  
\subfigtopskip=2pt
\subfigbottomskip=2pt 
\subfigcapskip=-5pt 
\subfigure[$p=\infty$.]{\label{fig:p-norm-3d1}
\includegraphics[width=0.18\linewidth]{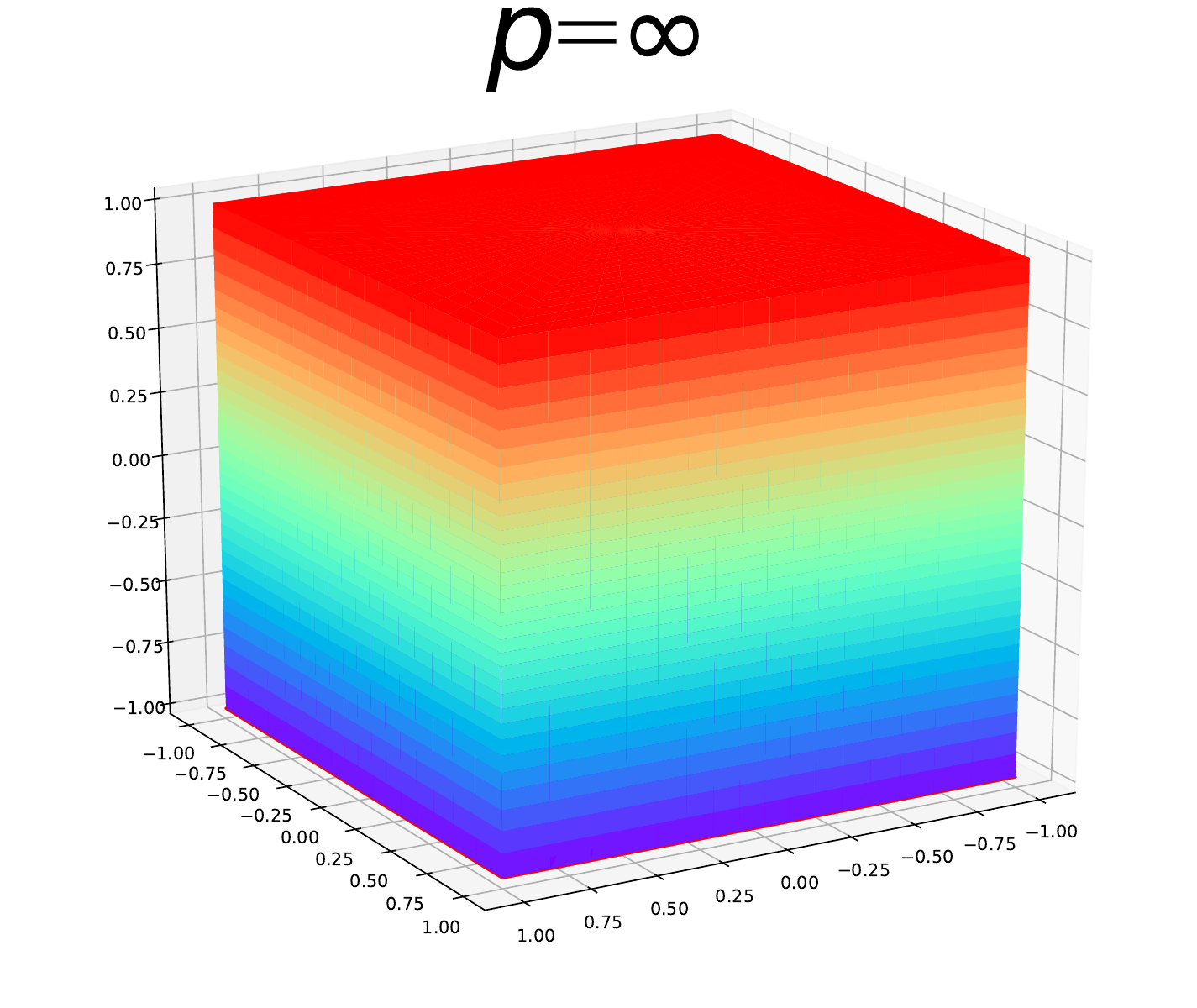}}
\subfigure[$p=2$.]{\label{fig:p-norm-3d2}
\includegraphics[width=0.18\linewidth]{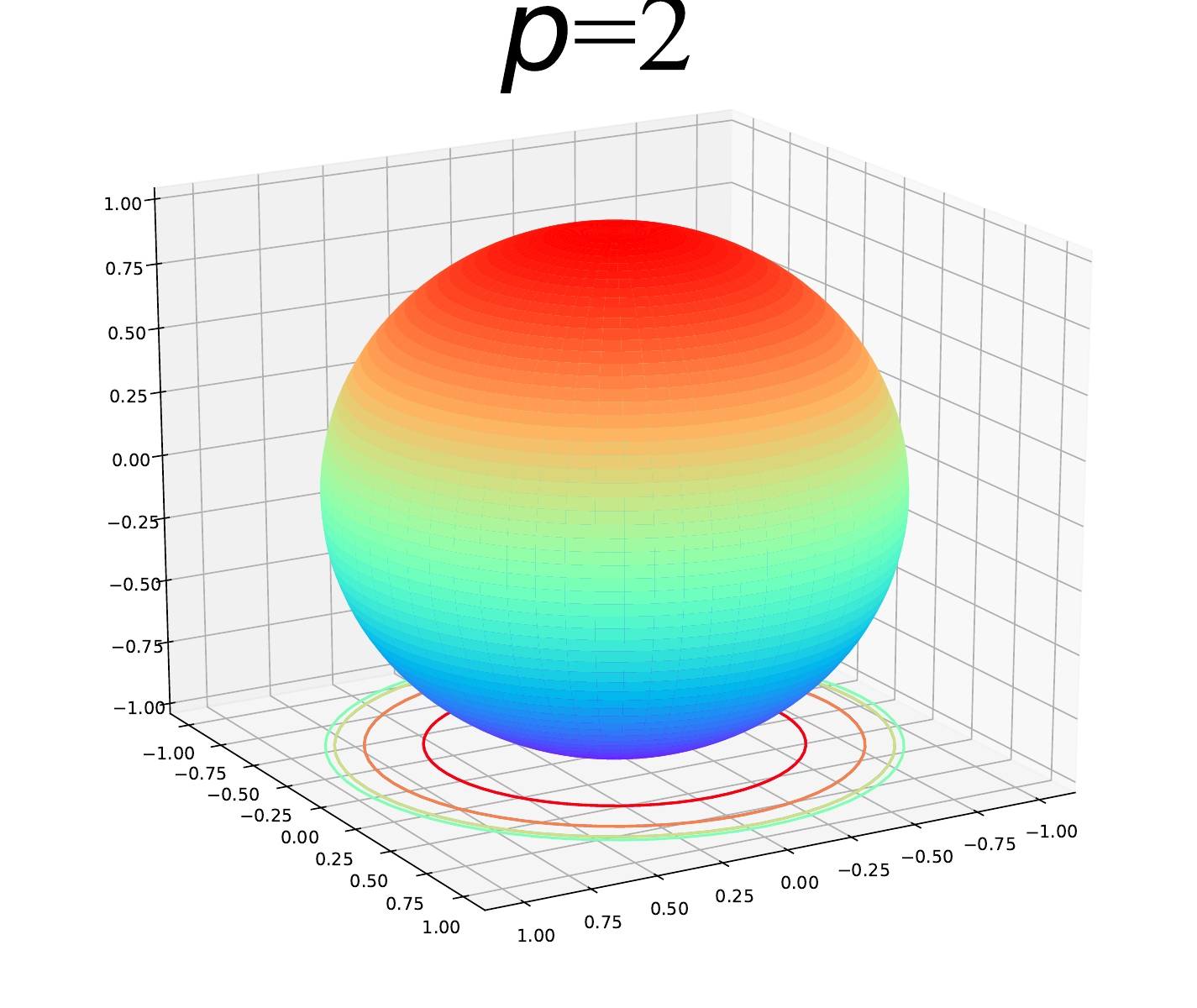}}
\subfigure[$p=1$.]{\label{fig:p-norm-3d3}
\includegraphics[width=0.18\linewidth]{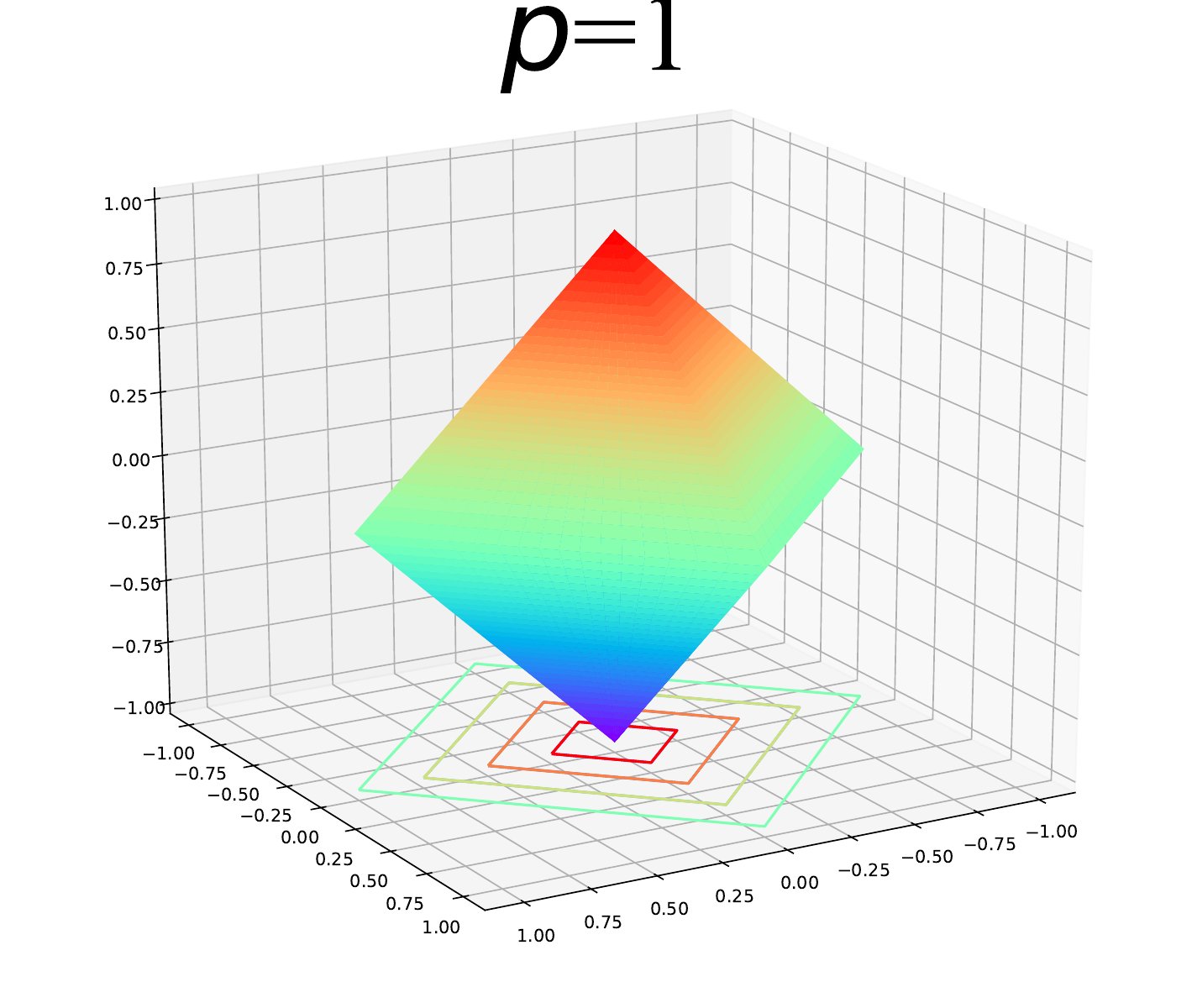}}
\subfigure[$p=0.5$.]{\label{fig:p-norm-3d4}
\includegraphics[width=0.18\linewidth]{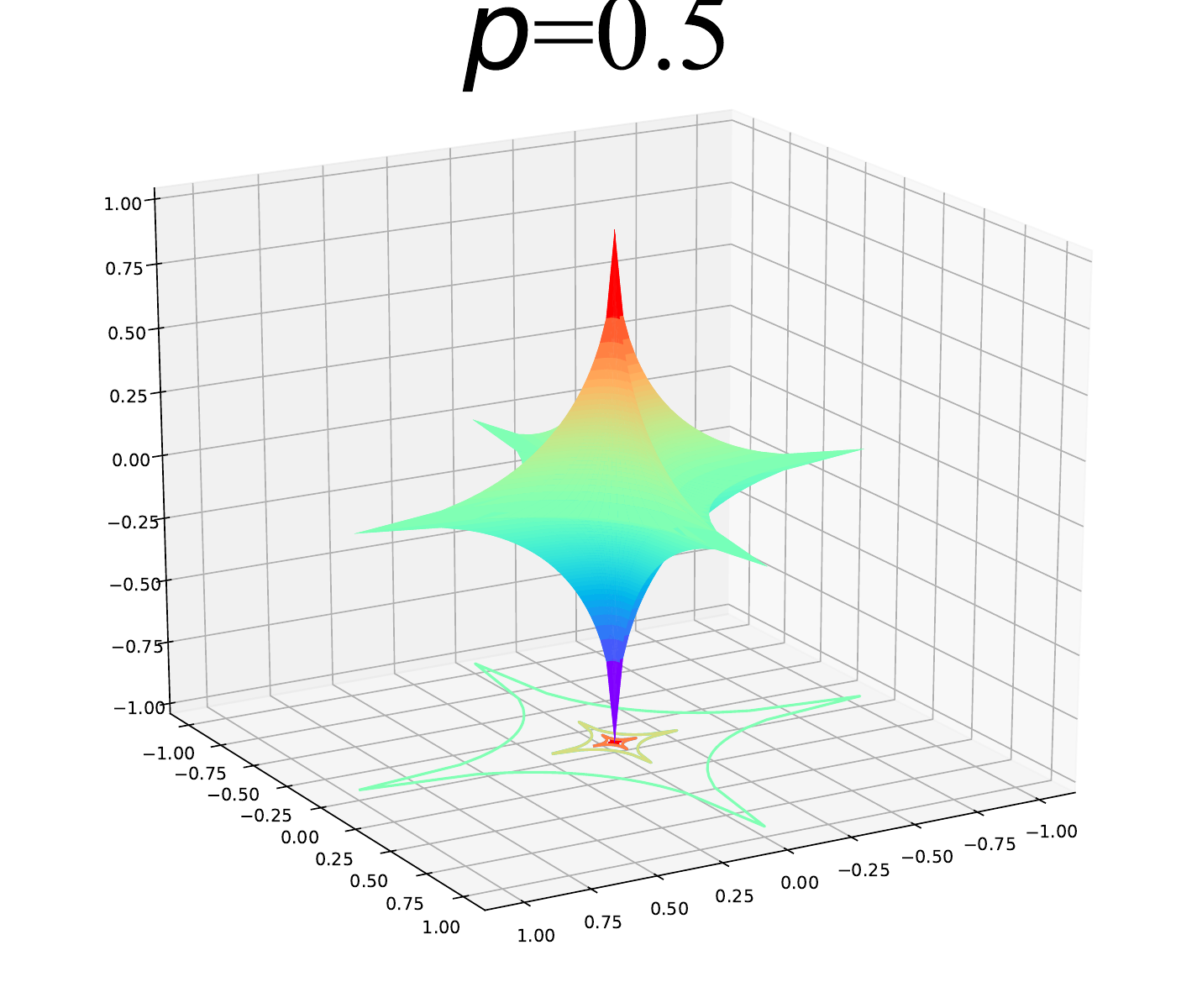}}
\subfigure[$p=0$.]{\label{fig:p-norm-3d5}
\includegraphics[width=0.18\linewidth]{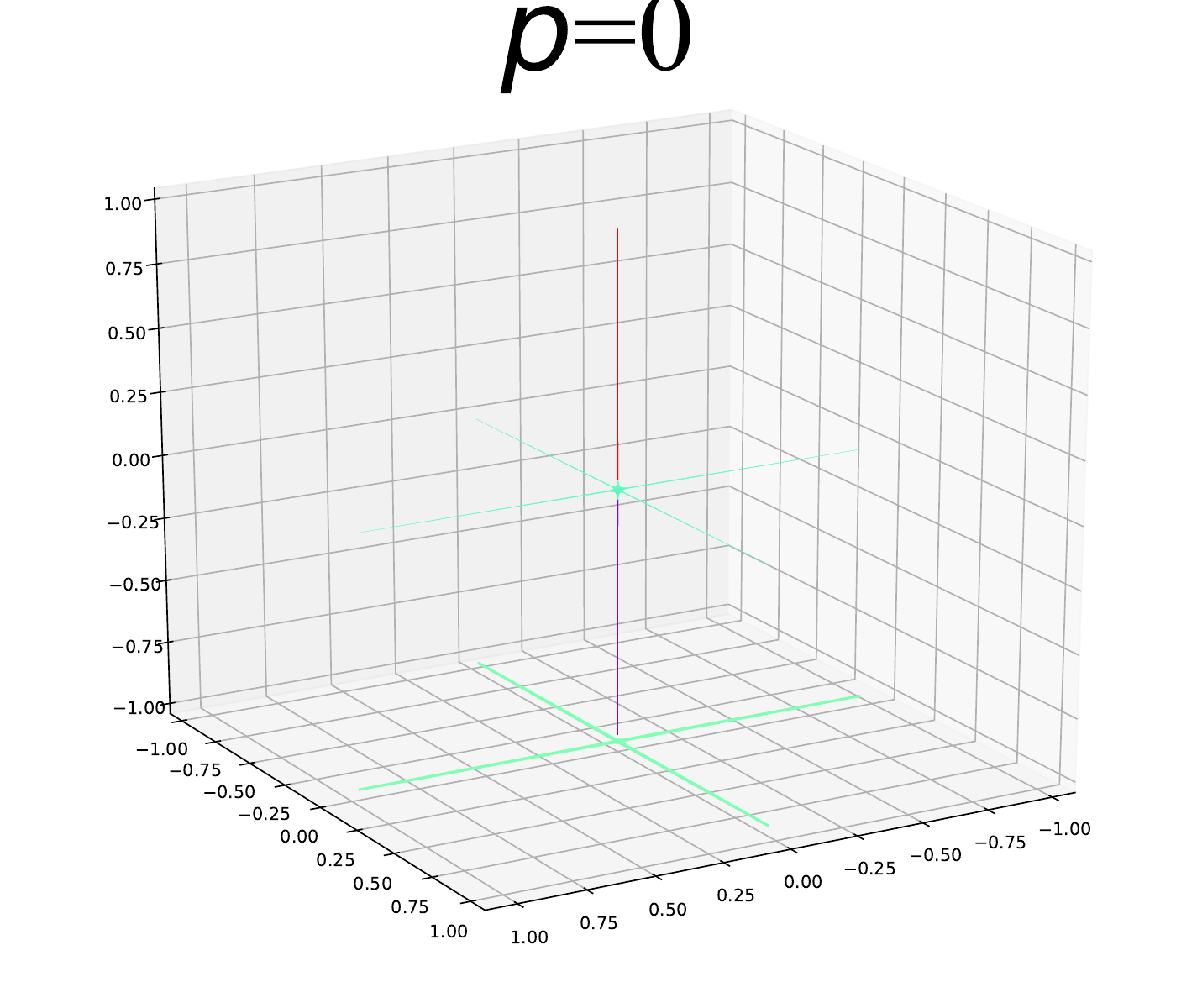}}
\caption{Unit ball of $\ell_p$-norms in three-dimensional space.}
\label{fig:p-norm-comparison-3d}
\end{figure}

\index{Decomposition: GL$_2^2$}
\index{Decomposition: GL$_\infty$}
\section{Gaussian $\ell_2^2$-Norm (GL$_2^2$) and Gaussian $\ell_\infty$-Norm (GL$_\infty$) Models}\label{section:l22_linfty_models}
Following the development of the GL$_1^2$ model, further exploration of the behavior induced by different ``norms" was conducted in \citet{lu2022robust}.
The $\ell_p$ prior builds upon the implicit regularization observed in the GL$_1^2$ model.
For any vector $\bx\in \real^N$, the $\ell_p$-norm is defined as
$$
\ell_p(\bx) = \big(\sum_{n=1}^N \abs{x_n}^p\big)^{1/p}.
$$
Figure~\ref{fig:p-norm-2d} and Figure~\ref{fig:p-norm-comparison-3d} illustrate the corresponding unit balls in two- and three-dimensional space, respectively.
Norms play a central role in machine learning. In Chapter~\ref{chapter:als}, we discussed the least squares problem, which minimizes the squared $\ell_2$ distance between an observation $\bb$ and its prediction $\bA\bx$: $\normtwo{\bA\bx-\bb}^2$. 
In contrast, minimizing the $\ell_1$-norm $\normone{\bA\bx-\bb}$ yields a more robust estimator of $\bx$  in the presence of outliers \citep{zoubir2012robust}.
For a matrix  $\bW\in\real^{M\times K}$, the $\ell_p$-norm can be extended row-wise as
\begin{equation}
\ell_p = \sum_{m=1}^{M} \Bigg(\sum_{k=1}^{K} \abs{w_{mk}}^p\Bigg)^{1/p}.
\end{equation}
In the context of NMF, where $w_{mk}\geq 0$,  the $\ell_1$-norm used in the GEE model (Equation~\eqref{equation:4norms-in-vanilla-bmf}) corresponds exactly to the $\ell_p$-norm with $p=1$. 
This norm is well known to encourage sparsity (see Section~\ref{section:gee_model}).

We now extend this Bayesian framework to models based on the $\ell_2^2$- and $\ell_\infty$-norms. As before, we assume that the data matrix $\bA$ is generated via the same probabilistic graphical model shown in Figure~\ref{fig:bmf_gl12}, with each entry $a_{mn}$ following a Gaussian likelihood
with variance $\sigma^2$ and mean given by the latent decomposition $\bw_m^\top\bz_n$ (Equation~\eqref{equation:als-per-example-loss_bnmf}).
Consequently, the posterior distribution of the noise variance $\sigma^2$, under an inverse-Gamma prior with shape $\alpha_\sigma$ and scale $\beta_\sigma$, remains identical to that in the GEE model (Equation~\eqref{equation:gee_posterior_sigma2}).

\paragrapharrow{Prior for the GL$_2^2$ model.}
Based on the squared  $\ell_2$-norm, we place independent priors on  $\bW$ and $\bZ$, , governed by hyper-parameters $\lambda^W$  and $\lambda^Z$, respectively:
\begin{equation}\label{equation:gp22_prior_density}
\begin{aligned}
&p(\bW\mid \lambda^W) &\propto &
\left\{
\begin{aligned}
&\exp \Bigg[  -\frac{\lambda^W}{2} \sum_{m=1}^{M} \Bigg(\sum_{k=1}^{K} w_{mk}^2\Bigg) \Bigg]
, &\gap &\text{if $w_{mk}\geq 0$ for all $m,k$ };
\\
&0, &\gap &\text{if otherwise};
\end{aligned}
\right.\\
\gap 
&p(\bZ\mid\lambda^Z) &\propto &
\left\{
\begin{aligned}
&\exp \Bigg[  -\frac{\lambda^Z}{2} \sum_{n=1}^{N} \Bigg(\sum_{k=1}^{K} z_{kn}^2\Bigg) \Bigg]
, &\gap &\text{if $z_{kn}\geq 0$ for all $n,k$ };
\\
&0, &\gap &\text{if otherwise}.
\end{aligned}
\right.\\
\end{aligned}
\end{equation}
\paragrapharrow{Posterior  for the GL$_2^2$ model.}
According to Bayes' rule (Equation~\eqref{equation:posterior_abstract_for_mcmc}), the posterior is proportional to the product of likelihood and prior, it can be maximized to yield an estimate of $\bW$ and $\bZ$. 
Using Bayes' theorem, 
the conditional density of $w_{mk}$ depends on its parents ($\lambda^W$), children ($a_{mn}$), and co-parents ($\tau$ or $\sigma^2$, $\bW_{-mk}, \bZ$). 
(See Figure~\ref{fig:bmf_gl12} and Section~\ref{section:markov-blanket}.) 
This yields:
\begin{equation}\label{equation:gp22_poster_wmk1}
\small
\begin{aligned}
&\gap p(w_{mk} \mid  \bA,   \bW_{-mk}, \bZ, \sigma^2,\lambda^W) \\
&\propto p(\bA\mid \bW, \bZ, \sigma^2) \times p(\bW\mid \lambda^W)
=\prod_{i,j=1}^{M,N} \normal \left(a_{ij}\mid \bw_i^\top\bz_j, \sigma^2 \right)\times 
p(\bW\mid \lambda^W)
 \cdot u(w_{mk}) \\
 &\propto \exp\Bigg\{   -\frac{1}{2\sigma^2}  \sum_{i,j=1}^{M,N}(a_{ij} - \bw_i^\top\bz_j  )^2\Bigg\}  \times 
 \exp \Bigg\{ -\frac{\lambda^W}{2} \sum_{i=1}^{M} \Bigg(\sum_{j=1}^{K} w_{ij}^2\Bigg) \Bigg\}
 \cdot u(w_{mk}) \\
\end{aligned}
\end{equation}
$$
\small
\begin{aligned}
\quad 
&\propto \exp\Bigg\{   -\frac{1}{2\sigma^2}  \sum_{j=1}^{N}(a_{mj} - \bw_m^\top\bz_j  )^2\Bigg\}  \times   
\exp\left\{ -\frac{\lambda^W}{2}w_{mk}^2 \right\} 
\cdot u(w_{mk}) \\
&\propto \exp\Bigg\{   
-
\underbrace{\Bigg(\frac{\sum_{j=1}^{N} z_{kj}^2 + \textcolor{winestain}{\sigma^2\lambda^W} }{2\sigma^2}  \Bigg)  }_{\textcolor{mylightbluetext}{\triangleq 1/(2\widetilde{\sigma^2_{mk} }) }} 
w_{mk}^2 +  w_{mk}\underbrace{\Bigg(  
	\frac{1}{\sigma^2} \sum_{j=1}^{N} z_{kj}\bigg( a_{mj} - \sum_{i\neq k}^{K}w_{mi}z_{ij}\bigg)  \Bigg)}_{\textcolor{mylightbluetext}{\triangleq \widetilde{\sigma_{mk}^{2}}^{-1} \widetilde{\mu_{mk}}}}
\Bigg\}  \cdot u(w_{mk})\\
&\propto   \normal(w_{mk} \mid \widetilde{\mu_{mk}}, \widetilde{\sigma_{mk}^{2}})\cdot u(w_{mk}) 
= \truncatednormal(w_{mk} \mid \widetilde{\mu_{mk}}, \widetilde{\sigma_{mk}^{2}}),
\end{aligned}
$$
where $\widetilde{\sigma_{mk}^{2}}= {\sigma^2}/{(\sum_{j=1}^{N} z_{kj}^2 +\textcolor{winestain}{\sigma^2\lambda^W})}$ is the posterior variance of the normal distribution,  
$$
\widetilde{\mu_{mk}} = 
\Bigg\{
\frac{1}{\sigma^2} \sum_{j=1}^{N} z_{kj}\bigg( a_{mj} - \sum_{i\neq k}^{K}w_{mi}z_{ij}\bigg)  
\Bigg\}
\cdot \widetilde{\sigma_{mk}^{2}}
$$
is the posterior mean of the normal distribution, and $\truncatednormal(x \mid \mu, \sigma^2)$ is the \textit{truncated-normal density} with ``parent" mean $\mu$ and ``parent" variance $\sigma^2$ (Definition~\ref{definition:truncated_normal}). 
Note again that this posterior closely resembles that of the GEE model (Equation~\eqref{equation:gee_poster_wmk1}), with differences highlighted in \textcolor{winestain}{red}; see also Table~\ref{table:bnmf_regularizer_posterior}.

\paragrapharrow{Connection between GEE, GL$_1^2$, and GL$_2^2$ models.}
We observe that the posterior ``parent" mean $\widetilde{\mu_{mk}}$ in the GL$_2^2$ model is larger than that in the GEE model since it does not contain the negative term $-\lambda_{mk}^W$ (see Table~\ref{table:bnmf_regularizer_posterior}). 
While the posterior ``parent" variance is smaller than that in the GEE model; therefore, the 
conditional density of GL$_2^2$ model is more clustered and it imposes a larger regularization in the sense of data/entry distribution (see Figure~\ref{fig:dists_truncatednorml_mean}, the smaller the ``parent" variance of the truncated-normal distribution, the larger the ``parent" precision, and the smaller the expectation of the truncated-normal variable). This can induce sparsity in the context of nonnegative matrix factorization.

Importantly, unlike the GL$_1^2$ model, the GL$_2^2$ model does \textbf{not} contain  the 
problematic  term $\sum_{j\neq k}^{K}w_{mj}$ appearing in the GL$_1^2$ mean, which causes the \textbf{inconsistency} across datasets with different scales of $\bA$ (as  previously discussed). 
Thus, the GL$_2^2$ model  is more robust and consistent across varying data types.

\paragrapharrow{Prior for the GL$_\infty$ model.}
As $p\rightarrow \infty$, the $\ell_p$-norm over  $\bW$  converges to 
\begin{equation}
	\ell_\infty = \sum_{m=1}^{M} \Bigg(\sum_{k=1}^{K} \abs{w_{mk}}^\infty\Bigg)^{1/\infty} 
	= \sum_{m=1}^{M} \mathop{\max}_{k} \abs{w_{mk}}.
\end{equation}
Based on the $\ell_p$-norm,
we assume $\bW$ and $\bZ$ are independently exponentially distributed with scales $\lambda^W$ and $\lambda^Z$, respectively (Definition~\ref{definition:exponential_distribution}):
\begin{equation}\label{equation:gpinfty_prior_density}
\begin{aligned}
&p(\bW\mid \lambda^W) &\propto &
\left\{
\begin{aligned}
&\exp \Bigg[  -{\lambda^W} \sum_{m=1}^{M} \mathop{\max}_{k} \abs{w_{mk}} \Bigg]
, &\gap &\text{if $w_{mk}\geq 0$ for all $m,k$ };
\\
&0, &\gap &\text{if otherwise};
\end{aligned}
\right.\\
\gap 
&p(\bZ\mid \lambda^Z) &\propto &
\left\{
\begin{aligned}
&\exp \Bigg[  -\lambda^Z \sum_{n=1}^{N} \mathop{\max}_{k} \abs{z_{kn}} \Bigg]
, &\gap &\text{if $z_{kn}\geq 0$ for all $n,k$ };
\\
&0, &\gap &\text{if otherwise}.
\end{aligned}
\right.\\
\end{aligned}
\end{equation}
Note that we omit the factor of $1/2$ in the exponent (compared to GL$_2^2$) for consistency with the resulting conditional posterior form (see Equation~\eqref{equation:gpinfty_poster_wmk1}).

\paragrapharrow{Posterior for GL$_\infty$ model.}
Applying Bayes' rule, the conditional density of $w_{mk}$ again follows a truncated normal distribution.
Let $\indicator(w_{mk})$ denote the indicator whether $w_{mk}$ is the largest one for $k=1,2,\ldots, K$ (i.e., $w_{mk}=\max_{k'}w_{mk'}$; it is the largest entry in row $m$). 
The conditional density can be obtained by
\begin{equation}\label{equation:gpinfty_poster_wmk1}
\small
\begin{aligned}
&\gap p(w_{mk} \mid \bA,   \bW_{-mk}, \bZ,\sigma^2, \lambda^W) \\
&\propto p(\bA\mid \bW, \bZ, \sigma^2) \times p(\bW\mid \lambda^W)
=\prod_{i,j=1}^{M,N} \normal \left(a_{ij}\mid \bw_i^\top\bz_j, \sigma^2 \right)\times 
p(\bW\mid \lambda^W)
\cdot u(w_{mk}) \qquad \quad \\
&\propto \exp\Bigg\{   -\frac{1}{2\sigma^2}  \sum_{i,j=1}^{M,N}(a_{ij} - \bw_i^\top\bz_j  )^2\Bigg\}  \times 
\exp \left\{ -{\lambda^W} \cdot\sum_{i=1}^{M} \mathop{\max}_{k} \abs{w_{ik}} \right\}
\cdot u(w_{mk}) \\
\end{aligned}
\end{equation}
$$
\small
\begin{aligned}
&\propto \exp\Bigg\{   -\frac{1}{2\sigma^2}  \sum_{j=1}^{N}(a_{mj} - \bw_m^\top\bz_j  )^2\Bigg\}  \times   
\exp\left\{ -{\lambda^W}\cdot w_{mk} \right\} 
\cdot u(w_{mk}) \cdot \textcolor{winestain}{\indicator(w_{mk})}\\
&\propto \exp\Bigg\{   -\frac{1}{2\sigma^2}  \sum_{j=1}^{N}
\bigg[ w_{mk}^2z_{kj}^2 + 2w_{mk} z_{kj}\bigg(\sum_{i\neq k}^{K}w_{mi}z_{ij} - a_{mj}\bigg)  \bigg]
\Bigg\} 
\cdot \exp\left\{ - w_{mk}  \textcolor{winestain}{\lambda^W\indicator( w_{mk})}\right\}  \cdot u(w_{mk})\\
&\propto \exp\Bigg\{   
-
\underbrace{\Bigg(\frac{\sum_{j=1}^{N} z_{kj}^2  }{2\sigma^2}  \Bigg)  }_{\textcolor{mylightbluetext}{ \triangleq 1/(2\widetilde{\sigma^2_{mk} }) }} 
w_{mk}^2 +
w_{mk}\underbrace{\bigg[
	-\textcolor{winestain}{\lambda^W \indicator(w_{mk})}+  
	\frac{1}{\sigma^2} \sum_{j=1}^{N} z_{kj}\bigg( a_{mj} - \sum_{i\neq k}^{K}w_{mi}z_{ij}\bigg)  \bigg]}_{\textcolor{mylightbluetext}{\triangleq \widetilde{\sigma_{mk}^{2}}^{-1} \widetilde{\mu_{mk}}}}
\Bigg\}  \cdot u(w_{mk})\\
&\propto   \normal(w_{mk} \mid \widetilde{\mu_{mk}}, \widetilde{\sigma_{mk}^{2}})\cdot u(w_{mk}) 
= \truncatednormal(w_{mk} \mid \widetilde{\mu_{mk}}, \widetilde{\sigma_{mk}^{2}}),
\end{aligned}
$$
where $u(x)$ is the unit function with value 1 if $x\geq 0$ and value 0 if $x<0$, $\widetilde{\sigma_{mk}^{2}}= {\sigma^2}/{(\sum_{j=1}^{N} z_{kj}^2 )}$ is the posterior ``parent" variance of the normal distribution,  
$$
\widetilde{\mu_{mk}} = 
\Bigg\{
-\textcolor{winestain}{\lambda^W\cdot \indicator(w_{mk})}+
\frac{1}{\sigma^2} \sum_{j=1}^{N} z_{kj}\bigg( a_{mj} - \sum_{i\neq k}^{K}w_{mi}z_{ij}\bigg)  
\Bigg\}
\cdot \widetilde{\sigma_{mk}^{2}}
$$
is the posterior ``parent" mean of the normal distribution, and $\truncatednormal(x \mid \mu, \sigma^2)$ is the \textit{truncated-normal density} with ``parent" mean $\mu$ and ``parent" variance $\sigma^2$ (Definition~\ref{definition:truncated_normal}). 

\index{Consistency}
\paragrapharrow{Connection between GEE and GL$_\infty$ models.}

The posterior ``parent" variance $\widetilde{\sigma^2_{mk}}$ in the GL$_\infty$ model matches that of GEE exactly (see Table~\ref{table:bnmf_regularizer_posterior}).
Denote $\indicator(w_{mk})$ as the indicator whether $w_{mk}$ is the largest one among $k=1,2,\ldots, K$.
Suppose further the condition $\indicator(w_{mk})$ is satisfied, parameters $\{\lambda_{mk}^W\}$ in the GEE model and $\lambda^W$ in the GL$_\infty$ model are equal, the ``parent" mean parameter $\widetilde{\mu_{mk}}$ is the same as that in the GEE model as well.
However, when $w_{mk}$ is not the maximum value among $\{w_{m1}, w_{m2}, \ldots, w_{mK}\}$, the ``parent" mean $\widetilde{\mu_{mk}}$ is larger than that in the GEE model since the GL$_\infty$ model excludes this negative term.
The GL$_\infty$ model then has the interpretation that it has a \textit{sparsity constraint} when $w_{mk}$ is the maximum value; and it has a \textit{relatively loose constraint} when $w_{mk}$ is not the maximum value. Overall, the GL$_\infty$ favors a loose regularization compared to the GEE model.

\paragrapharrow{Further extension: GL$_{2,\infty}^2$ model.}

The \textit{GL$_{2,\infty}^2$ model} takes the advantages of both GL$_2^2$ and GL$_\infty$ models.
The implicit prior of the GL$_{2,\infty}^2$ model can be obtained by 
\begin{equation}\label{equation:gp22infty_prior_density}
\begin{aligned}
&\gap p(\bW\mid \lambda^W) \propto
\exp \Bigg\{  \frac{-\lambda^W}{2} \sum_{m=1}^{M} \bigg(\sum_{k=1}^{K} w_{mk}^2+2\mathop{\max}_{k} |w_{mk}| \bigg) \Bigg\} u(\bW).
\end{aligned}
\end{equation}
The corresponding posterior parameters are summarized in Table~\ref{table:bnmf_regularizer_posterior}.

\paragrapharrow{Computational complexity and Gibbs sampler.}
All models---GEE, GL$_1^2$, GL$_2^2$, GL$_\infty$, and GL$_{2,\infty}^2$---share the same Gibbs sampling framework, with computational complexity  $\mathcalO(MNK^2)$. 
The dominant cost arises from evaluating the quadratic terms in the conditional posteriors of variables $\{w_{mk}\}$ and $\{z_{kn}\}$. The Gibbs sampler for the discussed models is formulated in Algorithm~\ref{alg:nmf_regular_12_infty_gibbs_sampler} outlines the general Gibbs sampler. 
By default, we use weakly informative priors:
$\lambda^W=\lambda^Z=0.1$ (for all regularized models; GL$_1^2$, GL$_2^2$, GL$_\infty$, GL$_{2,\infty}^2$);
$\alpha_\sigma=\beta_\sigma=1$ (for the inverse-Gamma prior in GL$_1^2$, GL$_2^2$, GL$_\infty$, GL$_{2,\infty}^2$).

\begin{algorithm}[h] 
\caption{Gibbs sampler for GL$_1^2$, GL$_2^2$, and GL$_\infty$ models (prior on  variance $\sigma^2$ here, similarly for the precision $\tau$). The procedure presented here is for explanatory purposes, and vectorization can expedite the procedure. 
By default, uninformative priors are 
$\lambda^W=\lambda^Z=0.1$ (GL$_1^2$, GL$_2^2$, GL$_\infty$, GL$_{2,\infty}^2$);
$\alpha_\sigma=\beta_\sigma=1$ (inverse-Gamma prior in GL$_1^2$, GL$_2^2$, GL$_\infty$, GL$_{2,\infty}^2$).
} 
\label{alg:nmf_regular_12_infty_gibbs_sampler}  
\begin{algorithmic}[1] 
\For{$k=1$ to $K$} 
\For{$m=1$ to $M$}
\State Sample $w_{mk}$ from $p(w_{mk} \mid \cdot )=\truncatednormal(w_{mk} \mid \widetilde{\mu_{mk}}, \widetilde{\sigma_{mk}^{2}})$ from Table~\ref{table:bnmf_regularizer_posterior}; 
\EndFor
\For{$n=1$ to $N$}
\State Sample $z_{kn}$ from $p(z_{kn} \mid\cdot )=\truncatednormal(z_{nk} \mid \widetilde{\mu_{kn}}, \widetilde{\sigma_{kn}^{2}})$; 
\Comment{symmetry of $w_{mk}$}
\EndFor
\EndFor
\State Sample $\sigma^2$ from $p(\sigma^2 \mid  \bA, \bW,\bZ,\alpha_\sigma,\beta_\sigma)$; 
\Comment{Equation~\eqref{equation:gee_posterior_sigma2}}
\State Report loss in Equation~\eqref{equation:als-per-example-loss_bnmf}, stop if it converges;
\end{algorithmic} 
\end{algorithm}

\subsubsection{Examples}
We conduct experiments across various analysis tasks to demonstrate the key advantages of the GL$_2^2$ and GL$_{2,\infty}^2$ methods. We use two datasets from bioinformatics: The first one is the Genomics of Drug Sensitivity in Cancer dataset~\footnote{\url{https://www.cancerrxgene.org/}} (GDSC $IC_{50}$) \citep{yang2012genomics}, which contains a wide range of drugs and their treatment outcomes on different cancer and tissue types (cell lines).
Following \citet{brouwer2017prior}, we preprocess the GDSC $IC_{50}$ dataset by capping high values
to 100, undoing the natural log transform, and casting them as integers. 
The second one is the Gene Body Methylation dataset \citep{koboldt2012comprehensive}, which gives the amount of methylation measured in the body region of 160 breast cancer driver genes.
We multiply the values in the Gene Body Methylation dataset by 20 and cast them as integers as well.
A summary of both datasets is provided in Table~\ref{table:datadescription_nmf_regularizer}, and their value distributions are shown in Figure~\ref{fig:datasets_nmf_regularizer}.
The GDSC $IC_{50}$ data exhibits a wide and unbalanced range, with values concentrated near 0 or capped at 100.
In contrast, the Gene Body Methylation data has a narrower and more balanced distribution.
We can see that the GDSC $IC_{50}$ is relatively a large dataset, whose matrix rank is $139$, and the Gene Body Methylation data tends to be small, possessing a matrix rank of 160.

\noindent\makebox[\textwidth][c]{%
\begin{minipage}{\textwidth}
\begin{minipage}[b]{0.41\textwidth}
\centering
\begin{figure}[H]
\centering  
\vspace{-0.35cm} 
\subfigtopskip=2pt 
\subfigbottomskip=2pt 
\subfigcapskip=-5pt 
\subfigure{\includegraphics[width=0.451\textwidth]{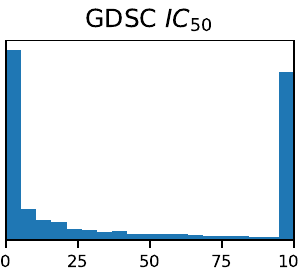} \label{fig:plot_gdsc}}
\subfigure{\includegraphics[width=0.451\textwidth]{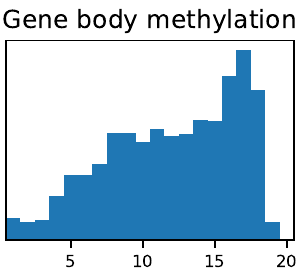} \label{fig:plot_methylation_gm}}
\caption{Data distribution of GDSC $IC_{50}$ and Gene Body Methylation datasets.}
\label{fig:datasets_nmf_regularizer}
\end{figure}
\end{minipage}
\hfill
\begin{minipage}[b]{0.575\textwidth}
\centering
\renewcommand{\arraystretch}{1.25}
\setlength{\tabcolsep}{7pt}
\small

\begin{tabular}{l|lll}
\hline
Dataset        & Rows & Columns & Fraction obs. \\ \hline
GDSC $IC_{50}$ & 707  & 139     & 0.806         \\
Gene Body Meth.  & 160  & 254     & 1.000         \\
\hline
\end{tabular}
\captionof{table}{Dataset description. Gene Body Methylation is relatively a small dataset, and the GDSC $IC_{50}$ tends to be large. The description provides the number of rows, columns, and the fraction of entries that are observed.}
\label{table:datadescription_nmf_regularizer}
\end{minipage}
\end{minipage}
}

All models use the same parameter initialization strategy. We evaluate performance in terms of convergence speed and generalization ability. Across a wide range of settings, GL$_2^2$ and GL$_{2,\infty}^2$ consistently achieve faster convergence and deliver out-of-sample performance that is as good as or better than other Bayesian NMF models with implicit regularization.

\paragrapharrow{Hyper-parameters.}
We adopt the default hyperparameter settings from \citet{brouwer2017prior}. We use $\{\lambda_{mk}^W\}=\{\lambda_{kn}^Z\}=0.1$ (GEE); 
$\lambda^W=\lambda^Z=0.1$ (GL$_1^2$, GL$_2^2$, GL$_{2,\infty}^2$);
uninformative $\alpha_\sigma=\beta_\sigma=1$ (inverse-Gamma prior in GEE, GL$_1^2$, GL$_2^2$, GL$_{2,\infty}^2$).
These represent weakly informative priors, and results are robust to small variations \citep{brouwer2017prior}.
Once hyper-parameters are fixed, all latent variables are initialized via random draws---a strategy that helps capture meaningful patterns early in inference.
In all experiments, we run the Gibbs sampler for 500 iterations, discarding the first 300 as burn-in. Convergence diagnostics confirm that the algorithm typically stabilizes within 200 iterations.

\begin{figure*}[htp]
\centering  
\vspace{-0.25cm} 
\subfigtopskip=2pt 
\subfigbottomskip=2pt 
\subfigcapskip=-1pt 
\subfigure[Convergence on the \textbf{GDSC $\boldsymbol{IC_{50}}$} dataset with increasing latent dimension $K$.]{\includegraphics[width=1\textwidth]{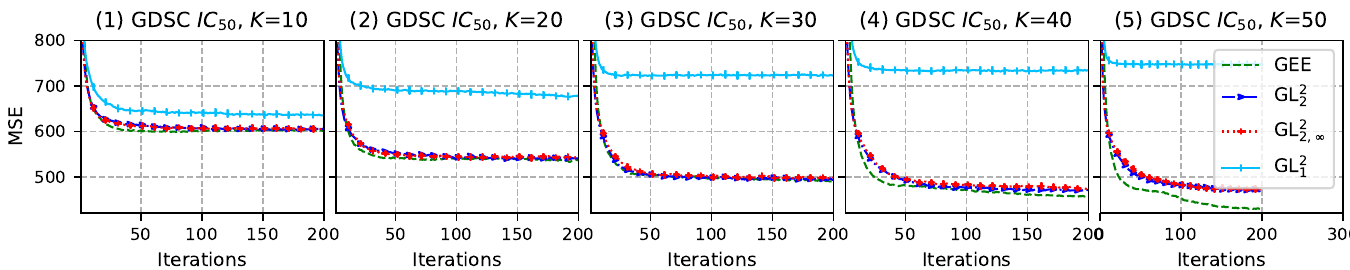} \label{fig:nmf_regularizers_convergences_gdsc}}\vspace{-0.6em}
\subfigure[Data distribution of factored component $\bW$ over the last 20 iterations for \textbf{GDSC $\boldsymbol{IC_{50}}$}. ]{\includegraphics[width=1\textwidth]{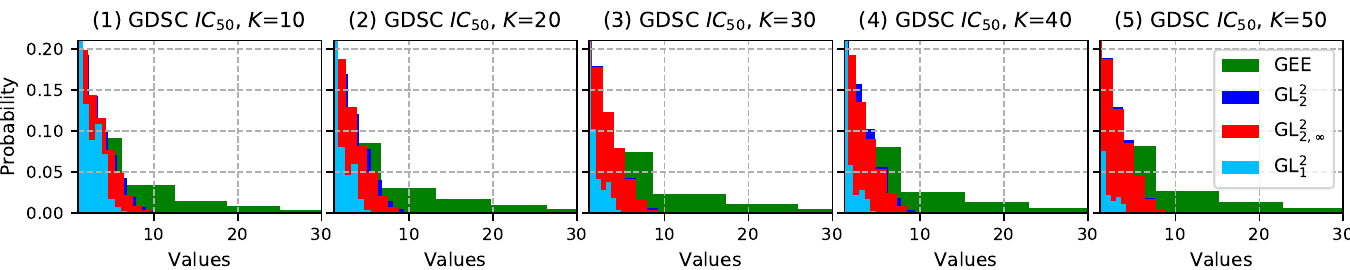} \label{fig:nmf_regularizers_distributions_gdsc}}\vspace{-0.1em}
\caption{Convergence of the models on the 
GDSC $IC_{50}$ (upper) and the distribution of factored $\bW$  (lower), measuring the training data fit (mean squared error). When we increase the latent dimension $K$, the GEE, GL$_2^2$, and GL$_{2,\infty}^2$ algorithms continue to increase the performance; while GL$_1^2$ starts to decrease. 
The results of GL$_\infty$ and GL$_{2,\infty}^2$ models are similar so we only present the results of the GL$_{2,\infty}^2$ model for brevity.}
\label{fig:convergences_gdsc_nmf_regularizer}
\end{figure*}

\begin{figure*}[h]
\centering  
\vspace{-0.25cm} 
\subfigtopskip=2pt 
\subfigbottomskip=2pt 
\subfigcapskip=2pt 
\subfigure[Convergence on the \textbf{Gene Body Methylation} dataset with increasing latent dimension $K$.]{\includegraphics[width=1\textwidth]{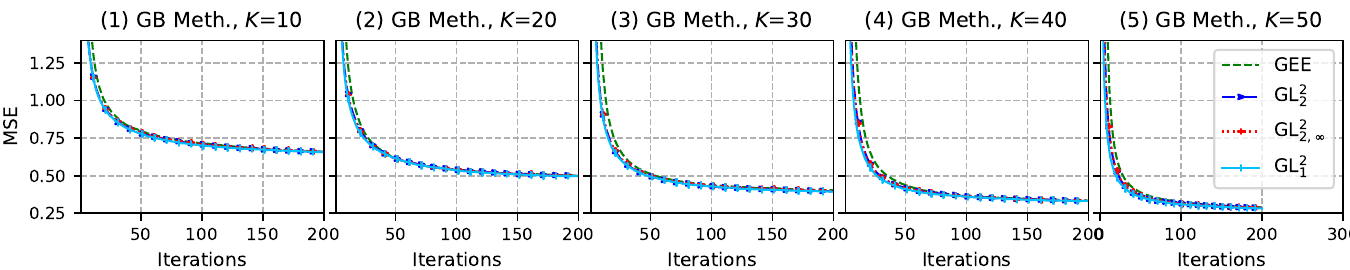} \label{fig:nmf_regularizers_convergences_gene_body}}\vspace{-0.6em}
\subfigure[Data distribution of factored component $\bW$ over the last 20 iterations for \textbf{Gene Body Methylation}.]{\includegraphics[width=1\textwidth]{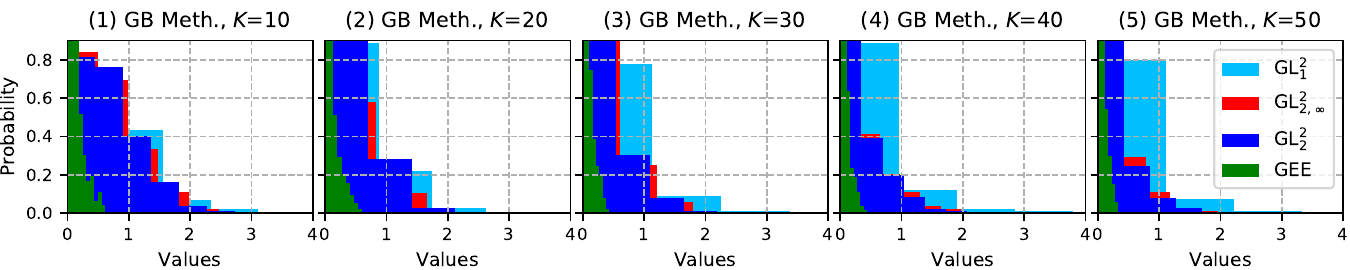} \label{fig:nmf_regularizers_distributions_gene_body}}\vspace{-0.1em}
\caption{Convergence of the models on the 
Gene Body Methylation dataset (upper) and the distribution of factored $\bW$ (lower), measuring the training data fit (mean squared error). When we increase the latent dimension $K$, all the models continue to improve the performance.}
\label{fig:convergences_ctrp_nmf_regularizer}
\end{figure*}

\paragrapharrow{Convergence analysis for GDSC $IC_{50}$ with relatively large entries.} Firstly, we compare the convergence in terms of iterations on the GDSC $IC_{50}$  and Gene Body Methylation datasets. We run each model with $K=\{10, 20, 30, 40, 50\}$, and the loss is measured by mean squared error (MSE).
Figure~\ref{fig:nmf_regularizers_convergences_gdsc}  shows the average convergence results over ten repeats, and 
Figure~\ref{fig:nmf_regularizers_distributions_gdsc} shows the distribution of entries of the factored $\bW$ for the last 20 iterations
on the GDSC $IC_{50}$ dataset. 
The result is consistent with our analysis (see the connection between different models). Since the values of the data matrix for GDSC $IC_{50}$ dataset is large, the posterior ``parent" mean $\widetilde{\mu_{mk}}$ in the GL$_1^2$ model is approaching zero or even negative; thus, it has a larger regularization than GEE model. This makes the GL$_1^2$ model converge to a worse performance. GL$_2^2$ and GL$_{2,\infty}^2$ models, on the contrary, impose a looser regularization than the GL$_1^2$ model, and the convergence performances are close to that of the GEE model.

\index{Consistency}
\paragrapharrow{Convergence analysis for Gene Body Methylation with relatively small entries.}
Figure~\ref{fig:nmf_regularizers_convergences_gene_body} further shows the average convergence results over ten repeats, and Figure~\ref{fig:nmf_regularizers_distributions_gene_body} shows the distribution of the entries of the factored $\bW$ for the last 20 iterations on the Gene Body Methylation dataset. 
The situation is \textbf{different} for the GL$_1^2$ model since the range of the entries of the Gene Body Methylation dataset is smaller than that of the GDSC $IC_{50}$ dataset (see Figure~\ref{fig:datasets_nmf_regularizer}).
This makes the $-\lambda^W\cdot\textcolor{black}{\sum_{j\neq k}^{K}w_{mj}}$ term of posterior ``parent" mean $\widetilde{\mu_{mk}}$ in the GL$_1^2$ model approach zero (see Table~\ref{table:bnmf_regularizer_posterior}), and the model then favors a looser regularization than the GEE model. 

The situation can be further presented by the distribution of the factored component $\bW$ on the GDSC $IC_{50}$ (Figure~\ref{fig:nmf_regularizers_distributions_gdsc}) and   the Gene Body Methylation (Figure~\ref{fig:nmf_regularizers_distributions_gene_body}). The GEE model has larger values of $\bW$ on the former dataset and smaller values on the latter; while GL$_1^2$ has smaller values of $\bW$ on the former dataset and larger values on the latter.
In other words, the regularization of the GEE and GL$_1^2$ is \textbf{inconsistent} on the two different data matrices.
In comparison, the  GL$_2^2$ and GL$_{2,\infty}^2$  models are \textbf{consistent} on different datasets, making them more robust algorithms to compute the nonnegative matrix factorization of the observed data.

Table~\ref{table:nmfreg_distribtuionsvalues} shows the mean values of the factored component $\bW$ over the last 20 iterations for GDSC $IC_{50}$ (upper table) and  Gene Body Methylation (lower table), where the value in the parentheses is the sparsity evaluated by taking the percentage of values smaller than 0.1.
The \textbf{inconsistency} of GEE  for different matrices can be observed (either large sparsity or small sparsity), while the results for the  GL$_2^2$ and GL$_{2,\infty}^2$ models are more \textbf{consistent}.


\begin{figure}[ht]
\centering
\noindent
\makebox[\textwidth][c]{%
\begin{minipage}{\textwidth}
\begin{minipage}[b]{0.4\textwidth}
\centering
\setlength{\tabcolsep}{3.4pt}
\renewcommand{\arraystretch}{1.35}
\footnotesize
\begin{tabular}{l|llll}
\hline
$K$ & GEE & GL$_1^2$   & GL$_2^2$ & GL$_{2,\infty}^2$  \\ \hline
10   &  8.1 (1.9)   &  1.3 (10.3)    &    2.4 (3.8) &   2.4 (4.5)  \\
20   &  8.6 (1.5)   &  0.8 (14.7)    &    2.3 (4.1) &   2.2 (4.4)  \\
30   &  8.7 (1.4)   &  0.7 (17.3)    &    2.2 (4.3) &   2.2 (4.4)  \\
40   &  8.3 (1.5)   &  0.6 (19.4)    &    2.2 (4.4) &   2.2 (4.4)  \\
50   &  8.0 (1.6)   &  0.5 (21.2)    &    2.2 (4.1) &   2.2 (4.2)  \\
\hline
\hline
10   &  0.1 (80.4)   &  0.7 (11.4)    &    0.7 (11.5) &   0.7 (12.7)  \\
20   &  0.1 (87.8)   &  0.6 (16.2)    &    0.5 (21.3) &   0.5 (21.0)  \\
30   &  0.0 (90.2)   &  0.6 (18.2)    &    0.3 (37.1) &   0.3 (36.4)  \\
40   &  0.0 (92.2)   &  0.6 (20.8)    &    0.3 (48.9) &   0.3 (49.1)  \\
50   &  0.0 (93.0)   &  0.5 (22.8)    &    0.2 (58.4) &   0.2 (58.4)  \\
\hline
\end{tabular}
\captionof{table}{Mean values of the factored component $\bW$ in the last 20 iterations, where the value in the (parentheses) is the sparsity evaluated by taking the percentage of values smaller than 0.1, for GDSC $IC_{50}$ (upper table) and  Gene Body Methylation (lower table).
The \textbf{inconsistency} of GEE and GL$_1^2$ for different matrices can be observed.
}
\label{table:nmfreg_distribtuionsvalues}
\end{minipage}
\vspace{+0.35cm}
\hfill\hfill
\begin{minipage}[b]{0.57\textwidth}
\centering
\renewcommand{\arraystretch}{1.1}
\footnotesize
\setlength{\tabcolsep}{3.8pt}
\begin{tabular}{l|l|llll}
\hline
Unobs.& $K$ & GEE & GL$_1^2$   & GL$_2^2$ & GL$_{2,\infty}^2$   \\ \hline
\parbox[t]{5.0mm}{\multirow{4}{*}{\rotatebox[origin=c]{0}{60\%}  }} 
&20   &  787.60   &  880.36    &   \textbf{ 769.24  } &   \textbf{ 768.27 } \\
&30   &  810.39   &  888.47    &   \textbf{ 774.53  } &   \textbf{ 773.27 } \\
&40   &  802.39   &  892.01    &   \textbf{ 783.26  } &   \textbf{ 784.30 } \\
&50   &  \textbf{795.72}   &  895.05    &   \textbf{ 806.14  } &   \textbf{ 807.44 } \\
\hline
\hline
\parbox[t]{5.0mm}{\multirow{4}{*}{\rotatebox[origin=c]{0}{70\%}  }} 
&20   &  841.74   &  895.77    &   \textbf{ 798.44  } &   \textbf{ 796.15 } \\
&30   &  830.45   &  902.48    &   \textbf{ 807.37  } &   \textbf{ 806.61 } \\
&40   &  \textbf{842.70}   &  907.65    &   \textbf{ 832.67  } &   \textbf{ 835.89 } \\
&50   & \textbf{846.83}   &  {1018.97}	$\uparrow$    &   \textbf{ 864.58  } &   \textbf{ 869.15 } \\
\hline
\hline
\parbox[t]{5.0mm}{\multirow{4}{*}{\rotatebox[origin=c]{0}{80\%}  }} 
&20   &  904.39   &  926.72    &   \textbf{ 842.24  } &   \textbf{ 841.84 } \\
&30   &  \textbf{887.63}   &  938.92    &   \textbf{ 879.30  } &   \textbf{ 883.57 } \\
&40   &  \textbf{942.44}   &  2634.69    &   \textbf{ 935.09  } &   \textbf{ 939.77 } \\
&50   &  \textbf{952.45}   &  {2730.30} $\uparrow$     &   \textbf{ 974.01  } &   \textbf{ 973.75 } \\
\hline
\end{tabular}
\captionof{table}{Mean squared error measure when the percentage of unobserved data is 60\% (upper table), 70\% (middle table), or 80\% (lower table) for the GDSC $IC_{50}$ dataset. The performance of the  GL$_2^2$ and GL$_{2,\infty}^2$ models is only slightly worse when we increase the fraction of unobserved from 60\% to 80\%; while the performance of GL$_1^2$ becomes extremely poor.
Similar observations occur in the Gene Body Methylation experiment. The symbol $\uparrow$ means the performance becomes extremely worse.}
\label{table:nmfregu_special_sparsity_case}
\end{minipage}
\end{minipage}
}
\end{figure}

\begin{figure*}[h]
	\centering  
	\subfigtopskip=2pt 
	\subfigbottomskip=2pt 
	\subfigcapskip=-2pt 
	\subfigure[Predictive results on the \textbf{GDSC $\boldsymbol{IC_{50}}$} dataset with increasing fraction of unobserved data and increasing latent dimension $K$.]{\includegraphics[width=1\textwidth]{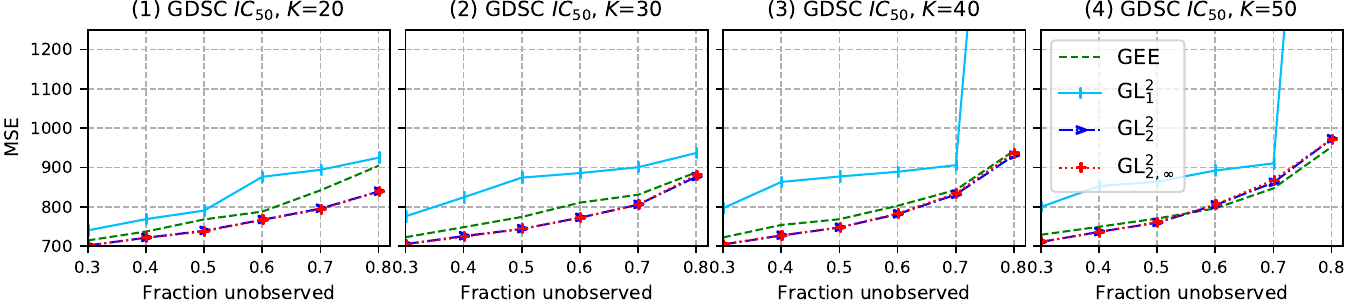} 
		\label{fig:nmf_regularizer_sparsity_movielens_gdsc}}
	\subfigure[Predictive results on \textbf{Gene Body Methylation} dataset  with increasing fraction of unobserved data and increasing latent dimension $K$.]{\includegraphics[width=1\textwidth]{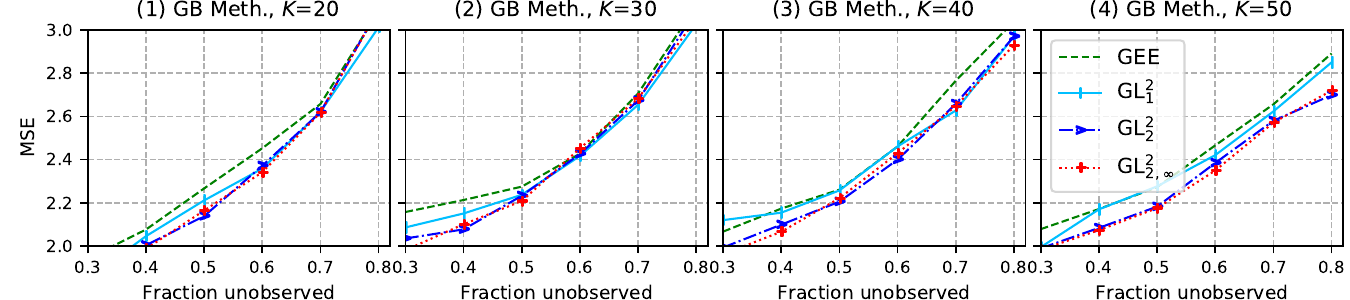} 
		\label{fig:nmf_regularizer_sparsity_movielens_genebody_meth}}
	\caption{Predictive results on the \textbf{GDSC $\boldsymbol{IC_{50}}$} (upper) and \textbf{Gene Body Methylation} (lower) datasets.
		We measure the predictive performance (mean squared error) on a held-out dataset for different fractions of unobserved data. }
	\label{fig:sparsity_nmf_regulariza}
\end{figure*}

\paragrapharrow{Predictive analysis.}

The training performances of the GEE, GL$_2^2$, and GL$_{2,\infty}^2$ models steadily improve as the model complexity grows. Inspired by this result, 
we measure the predictive performance when the sparsity of the data increases to see whether the models overfit or not. For different fractions of unobserved data, we randomly split the data based on that fraction, train the model on the observed data, and measure the performance on the held-out test data. Again, we increase the latent dimension $K$ from $K=20$ to $K=30, 40, 50$ for all models. The average MSE of ten repeats is given in Figure~\ref{fig:sparsity_nmf_regulariza}. 
We still observe the \textbf{inconsistency} issue in the GL$_1^2$ model, its predictive performance  is as good as that of the introduced GL$_2^2$ and GL$_{2,\infty}^2$ models on the Gene Body Methylation dataset; while the predictive results of the GL$_1^2$ model are extremely poor on the GDSC $IC_{50}$ dataset. 

For the GDSC $IC_{50}$ dataset, the  GL$_2^2$ and GL$_{2,\infty}^2$ models perform best when the latent dimensions are $K=20, 30, 40$; when $K=50$ and the fraction of unobserved data increases, the GEE model is slightly better. As aforementioned, the GL$_1^2$ performs the worst on this dataset; and when the fraction of unobserved data increases or $K$ increases, the predictive results of GL$_1^2$ deteriorate quickly.

For the Gene Body Methylation dataset, the predictive performance of GL$_1^2$, GL$_2^2$, and GL$_{2,\infty}^2$ models are close (GL$_1^2$ has a slightly larger error). The GEE model performs the worst on this dataset.

The comparison of the results on the two sets shows the  GL$_2^2$ and GL$_{2,\infty}^2$ models have both better in-sample and out-of-sample performance, making them a more \textbf{robust} choice in predicting missing entries.

Table~\ref{table:nmfregu_special_sparsity_case} shows the MSE predictions of different models when the fractions of unobserved data are $60\%$, $70\%$, and $80\%$, respectively. We observe that the performance of the  GL$_2^2$ and GL$_{2,\infty}^2$ models are only slightly worse when we increase the fraction of unobserved from 60\% to 80\%. This, again, indicates that the  GL$_2^2$ and GL$_{2,\infty}^2$ models are more \textbf{robust} with less overfitting. While for the GL$_1^2$ model, the performance becomes extremely poor in this scenario.

\begin{figure}[h]
\centering  
\subfigtopskip=2pt 
\subfigbottomskip=9pt 
\subfigcapskip=-5pt 
\includegraphics[width=0.421\textwidth]{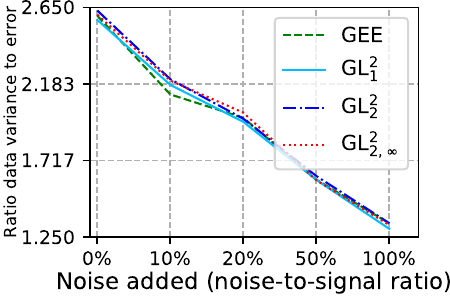}
\caption{Ratio of the variance of data to the MSE of the predictions, the higher the better.}
\label{fig:noise_graph_gdsc}
\end{figure}

\paragrapharrow{Noise sensitivity.}
Finally, we measure the noise sensitivity of different models with predictive performance when the datasets are noisy. To see this, we add different levels of Gaussian noise to the data. We add levels of 
$\{0\%, 10\%,$ $20\%,$ $50\%, 100\%\}$ 
noise-to-signal ratio noise (which is the ratio of the variance of the added Gaussian noise to the variance of the data). The results for the GDSC $IC_{50}$ with $K=10$ are shown in Figure~\ref{fig:noise_graph_gdsc}. The results are the average performance over 10 repeats. We observe that the  GL$_2^2$ and GL$_{2,\infty}^2$ models perform slightly better than other Bayesian NMF models (with implicit regularization meaning).  
The  GL$_2^2$ and GL$_{2,\infty}^2$ models perform notably better when the noise-to-signal ratio is smaller than 10\% and slightly better when the ratio is larger than 20\%.
Similar results can be found on the Gene Body Methylation dataset and other $K$ values, and we shall not repeat the details.

\index{Semi-nonnegative matrix factorization}
\index{Exponential distribution}
\section{Semi-Nonnegative Matrix Factorization}
Instead of enforcing nonnegativity constraints on both factor matrices, an alternative approach---proposed by \citet{ding2008convex} and \citet{fei2008semi}---is to apply the constraint to only one of them.
Within a Bayesian framework, this can be achieved by placing a real-valued prior on one factor matrix and a nonnegative prior on the other.
As previously discussed, standard NMF is particularly well-suited for inherently nonnegative data, such as images or text corpora. However, the key advantage of \textit{semi-nonnegative matrix factorization (semi-NMF)} is its ability to handle real-valued datasets while still preserving nonnegativity in one factor. This offers greater flexibility in capturing the underlying structure of the data without restricting all components to be nonnegative.

\index{Decomposition: GEG}
\subsection{Gaussian Likelihood with Exponential and Gaussian Priors (GEG)}
\begin{figure}[h]
\centering  
\vspace{-0.35cm} 
\subfigtopskip=2pt 
\subfigbottomskip=2pt 
\subfigcapskip=-5pt 
\subfigure[GEG.]{\label{fig:bmf_geg}
\includegraphics[width=0.421\linewidth]{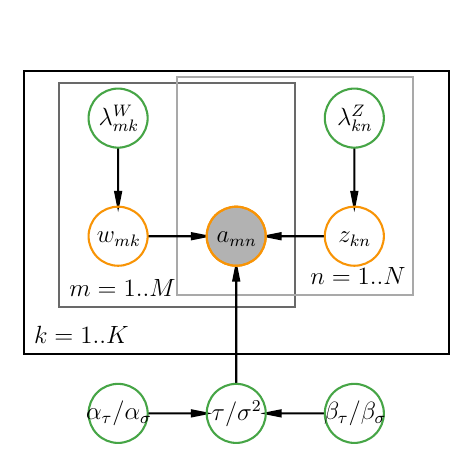}}
\subfigure[GnVG.]{\label{fig:bmf_gnvg}
\includegraphics[width=0.421\linewidth]{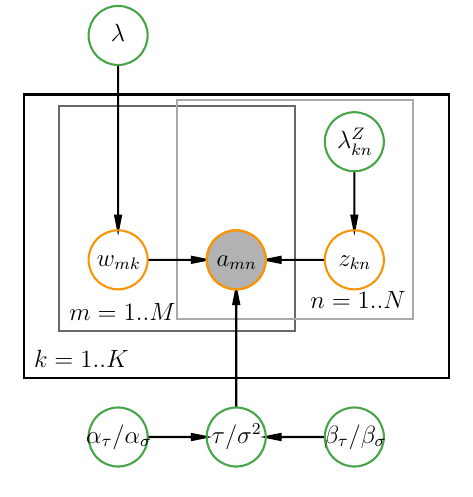}}
\caption{Graphical model representation of GEG and GnVG models. Green circles denote prior variables, orange circles represent observed and latent variables (shaded cycles denote observed variables), and plates represent repeated variables. 
The slash ``/" in the variable represents ``or."}
\label{fig:bmf_geg_gnvg}
\end{figure}

The \textit{Gaussian likelihood with exponential and Gaussian priors (GEG)} model assigns  an exponential prior to the component $\bW$ and a Gaussian prior to the component $\bZ$, following the same choices made in the GEE and GGG models, respectively (Section~\ref{section:gee_model}  and Section~\ref{section:markov-blanket}). 
The likelihood function is identical to that used in the GEE model (Equation~\eqref{equation:gee_likelihood}). A graphical representation of the GEG model is shown in  Figure~\ref{fig:bmf_geg}.

We assume that the entries of $\bW$ are independently exponentially distributed with rate parameters $\{\lambda_{mk}^W\}$, and that the entries of $\bZ$ follow independent Gaussian distributions with zero mean and precisions $\{\lambda_{kn}^Z\}$. 
Formally, 
\begin{equation}\label{equation:geg_prior_density_exponential_gaussian}
\begin{aligned}
w_{mk} &\sim \exponential(w_{mk}\mid \lambda_{mk}^W), 
\gap 
&z_{kn}\sim&  \normal(z_{kn}\mid 0, ( \lambda_{kn}^Z)^{-1});\\
p(\bW) &=\prod_{m,k=1}^{M,K} \exponential(w_{mk}\mid\lambda_{mk}^W), 
\gap 
&p(\bZ) =&\prod_{k,n=1}^{K,N} \normal(z_{kn}\mid 0, (\lambda_{kn}^Z)^{-1}),
\end{aligned}
\end{equation}
where $\exponential(x\mid  \lambda)=\lambda\exp(-\lambda x)u(x)$ is the exponential density, 
with $u(x)$ denoting  the unit step function (i.e., $u(x)=1$ if $x\geq 0$, and $0$ otherwise). 
For the noise variance $\sigma^2$, we again adopt an inverse-Gamma prior with shape parameter ${\alpha_\sigma}$ and scale parameter ${\beta_\sigma}$ (see Equation~\eqref{equation:geg_sigma_prior}).
The conditional posterior distributions for variables $\{w_{mk}\}$ and $\{z_{kn}\}$ are directly inherited from the GEE and GGG models, respectively, and thus require no rederivation here.

\index{Decomposition: GnVG}
\subsection{Gaussian Likelihood with  Volume and Gaussian Priors (GnVG)}
The volume prior introduced in the GVG model (Section~\ref{section:gvg_model}) can be adapted to enforce nonnegativity by restricting it to nonnegative values of $\bW$, as illustrated in Figure~\ref{fig:bmf_gnvg}.

As in the GGG and GEG models, we place a Gaussian prior on $\bZ$ with precisions $\{\lambda_{kn}^Z\}$:
\begin{equation}\label{equation:gnvg_prior_density_gaussian}
\begin{aligned}
z_{kn}\sim&  \normal(z_{kn}\mid 0, ( \lambda_{kn}^Z)^{-1})
\qquad\implies\qquad 
p(\bZ) =\prod_{k,n=1}^{K,N} \normal(z_{kn}\mid 0, (\lambda_{kn}^Z)^{-1}).
\end{aligned}
\end{equation}
The nonnegative volume prior over $\bW$ is defined as follows:
\begin{equation}\label{equation:gnvg_nonneg_volume_prior}
\bW\sim 
\left\{
\begin{aligned}
&\exp\{-\gamma \det(\bW^\top\bW)\} ,& \mathrm{\,\,if\,\,} w_{mk} \geq 0\text{\,\, for all }m,k;  \\
&0 , &\text{\,\,if\,\,any } w_{mk}<0.
\end{aligned}
\right.
\end{equation}
The posterior distributions for variables $\{z_{kn}\}$ are identical to that in the GEG and GGG models.
For variables $\{w_{mk}\}$, the posteriors resemble those of the GVG model, except that samples are now drawn from a truncated-normal distribution (truncated at zero) rather than an unrestricted normal distribution, to respect the nonnegativity constraint.

\index{Decomposition: NMTF}
\index{Decomposition: Tri-NMF}
\section{Nonnegative Matrix Tri-Factorization (NMTF)}\label{section:nmtf}

\begin{figure}[h]
\centering  
\vspace{-0.35cm} 
\subfigtopskip=2pt 
\subfigbottomskip=2pt 
\subfigcapskip=-5pt 
\subfigure[GEEE.]{\label{fig:bmf_geee}
\includegraphics[width=0.421\linewidth]{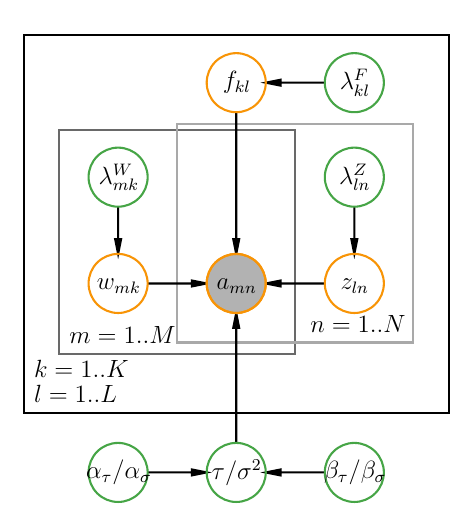}}
\subfigure[GEEEA.]{\label{fig:bmf_geeea}
\includegraphics[width=0.421\linewidth]{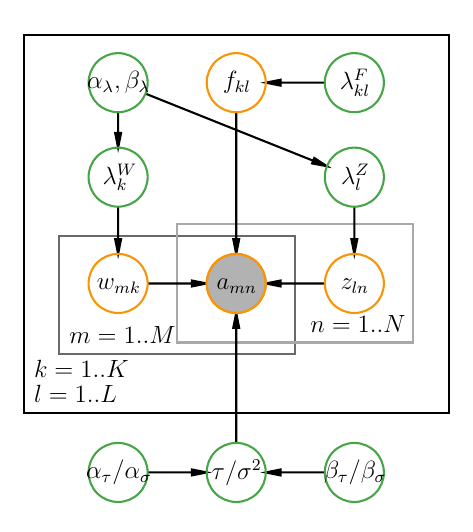}}
\caption{Graphical model representation of GEEE and GEEEA models. Green circles denote prior variables, orange circles represent observed and latent variables (shaded cycles denote observed variables), and plates represent repeated variables. The slash ``/" in the variable represents ``or," and the comma ``," in the variable represents ``and."}
\label{fig:bmf_geee_geeea}
\end{figure}

Similar to \textit{bilinear} nonnegative matrix factorization---where an observed matrix is decomposed into the product of two nonnegative factor matrices---\textit{nonnegative matrix tri-factorization (tri-NMF or NMTF)} extends this idea to three factors. Specifically, the data matrix $\bA$ is factorized as
$$
\bA = \bW\bF\bZ+\bE,
$$ 
where $\bW\in\real_+^{M\times K}, \bF\in\real_+^{K\times L}$, and $\bZ\in\real_+^{L\times N}$ are all element-wise nonnegative.
The advantages of the tri-NMF model are discussed in Section~\ref{section:nmf_tree_fac} and we shall not repeat here.

\paragrapharrow{Likelihood.}
As before, we assume the residuals $e_{mn}$ are i.i.d. zero-mean Gaussian variables with variance $\sigma^2$. This yields the following likelihood:
\begin{equation}\label{equation:geee_likelihood}
\begin{aligned}
p(\bA\mid  \btheta) 
&= \prod_{m,n=1}^{M,N} \normal \left(a_{mn}\mid  \bw_m^\top\bF\bz_n, \sigma^2 \right)
= \prod_{m,n=1}^{M,N} \normal \left(a_{mn}\mid \bw_m^\top\bF\bz_n, \tau^{-1} \right),
\end{aligned}
\end{equation}
where $\btheta=\{\bW,\bF,\bZ,\sigma^2\}$ denotes all model parameters, $\sigma^2$ is the noise variance, and $\tau^{-1}=\sigma^2$ is the precision. 
Equivalently, this corresponds to minimizing the reconstruction error measured by the Frobenius norm:
\begin{equation}\label{equation:als-per-example-loss_tnmf}
\mathop{\min}_{\bW,\bZ} L(\bW,\bZ) =  \mathop{\min}_{\bW,\bZ}\sum_{n=1}^N \sum_{m=1}^{M} \left(a_{mn} -  \bw_m^\top\bF\bz_n\right)^2.
\end{equation}

\paragrapharrow{Prior.} 
We place independent exponential priors on all entries of $\bW$, $\bF$, and $\bZ$, with scale parameters $\{\lambda_{mk}^W\}$, $\{\lambda_{kl}^F\}$, and $\{\lambda_{kn}^Z\}$, respectively (see Definition~\ref{definition:exponential_distribution}):
\begin{equation}\label{equation:geee_prior_density_exponential}
\begin{aligned}
w_{mk} &\sim \exponential(w_{mk}\mid  \lambda_{mk}^W), 
\gap 
&f_{kl}\sim&  \exponential(f_{kl}\mid  \lambda_{kl}^F),
&z_{ln}\sim&  \exponential(z_{ln}\mid  \lambda_{ln}^Z);\\
p(\bW) &=\prod_{m,k=1}^{M,K} \exponential(w_{mk}\mid  \lambda_{mk}^W), 
&p(\bF) =&\prod_{k,l=1}^{K,L} \exponential(f_{kl}\mid  \lambda_{kl}^F),
&p(\bZ) =&\prod_{l,n=1}^{L,N} \exponential(z_{ln}\mid  \lambda_{ln}^Z),
\end{aligned}
\end{equation}
where $\exponential(x\mid  \lambda)=\lambda\exp(-\lambda x)u(x)$ is the exponential density, and $u(x)$ is the unit step function.
Given this choice of priors and Gaussian likelihood, we refer to this model as the \textit{GEEE model} (Gaussian-Exponential-Exponential-Exponential). 
A graphical representation is provided in Figure~\ref{fig:bmf_geee}.

For the noise variance $\sigma^2$, we still adopt an inverse-Gamma prior with shape ${\alpha_\sigma}$ and scale ${\beta_\sigma}$ (Definition~\ref{definition:inverse_gamma_distribution}):
\begin{equation}\label{equation:geeg_sigma_prior}
p(\sigma^2)= \inversegammadist(\sigma^2\mid  \alpha_\sigma, \beta_\sigma) = \frac{{\beta_\sigma}^{\alpha_\sigma}}{\Gamma({\alpha_\sigma})} (\sigma^2)^{-\alpha_\sigma-1} \exp\left( -\frac{{\beta_\sigma}}{\sigma^2} \right).
\end{equation}
Consequently, the posterior distribution of $\sigma^2$ remains identical to that in the GEE model (Equation~\eqref{equation:gee_posterior_sigma2}).

\paragrapharrow{Posterior.}
Following Bayes' rule and using MCMC, inference proceeds by sampling from the full conditional distributions of each latent variable (via their Markov blankets; see Section~\ref{section:markov-blanket}):
$$
\begin{aligned}
&p(w_{mk}\mid \bA,  \bW_{-mk},\bF, \bZ,\sigma^2, \blambda^W,\blambda^F, \blambda^Z), \\
&p(f_{kl}\mid  \bA, \bW,\bF_{-kl}, \bZ, \sigma^2,\blambda^W,\blambda^F,  \blambda^Z), \\
& p(z_{ln}\mid \bA,  \bW, \bF,\bZ_{-ln},\sigma^2, \blambda^W,\blambda^F,  \blambda^Z), \\
& p(\sigma^2 \mid  \bA, \bW, \bZ,\alpha_\sigma, \beta_\sigma), \\
\end{aligned}
$$
where $\blambda^W$ is an $M\times K$ matrix containing all $\{\lambda_{mk}^W\}$ entries,
$\blambda^F$ is a $K\times L$ matrix of $\{\lambda_{kl}^F\}$,
 $\blambda^Z$ is an $L\times N$ matrix including all $\{\lambda_{ln}^Z\}$ values, and $\bW_{-{mk}}$ denotes all elements of $\bW$ except $w_{mk}$. 
The conditional density of $w_{mk}$ is just similar to that in the GEE model in Equation~\eqref{equation:gee_poster_wmk1}.
For simplicity, we denote the $k$-th row of $\bF$ as $\br_k$, 
and the $l$-th column of $\bF$ as $\bc_l$. The conditional density of $w_{mk}$ is the same as that in Equation~\eqref{equation:gee_poster_wmk1}, except now we replace 
$z_{kj}$ with $\br_k^\top\bz_j$ in the variance parameter of Equation~\eqref{equation:gee_posterior_variance}, and 
replace $z_{kj}$ with $\br_k^\top\bz_j$ 
and replace $z_{ij}$ with $\br_i^\top\bz_j$ in  Equation~\eqref{equation:gee_posterior_mean}.
The reason is obvious as, when considering the conditional density of $w_{mk}$, we can treat $\bF\bZ$ as a single matrix, and the problem becomes a bilinear decomposition.
By symmetry, the conditional posterior for variables $\{z_{ln}\}$ follow the same logic.

The conditional posterior for variables  $\{f_{kl}\}$, however, require explicit derivation. 
Using Bayes' theorem, 
the conditional density of $f_{kl}$ depends on its parents ($\lambda_{kl}^F$), children ($a_{mn}$), and co-parents ($\tau$ or $\sigma^2$, $\bW, \bF_{kl}, \bZ$). 
(See Figure~\ref{fig:bmf_geee} and Section~\ref{section:markov-blanket}.)
We obtain:
\begin{equation}\label{equation:geee_poster_fkl1}
\begin{aligned}
&\gap p(f_{kl}\mid  \bA , \bW,\bF_{-kl}, \bZ, \sigma^2, \cancel{\blambda^W},\cancel{\blambda^Z}, \blambda^F, \bA)=p(f_{kl} \mid \bA, \bW,   \bF_{-kl}, \bZ, \sigma^2, \lambda_{kl}^F) \\
&\propto p(\bA\mid  \bW, \bF,\bZ, \sigma^2) \times p(f_{kl}\mid  \lambda_{kl}^F)
=\prod_{i,j=1}^{M,N} \normal \left(a_{ij}\mid  \bw_i^\top\bF\bz_j, \sigma^2 \right)\times \exponential(w_{kl}\mid  \lambda_{kl}^F) \\
&\propto 
\exp\Bigg\{   -\frac{1}{2\sigma^2}  \sum_{i,j=1}^{M,N}(a_{ij} - \bw_i^\top\bF\bz_j  )^2\Bigg\}  \times \cancel{\lambda_{kl}^F }\exp(-\lambda_{kl}^F \cdot f_{kl})u(f_{kl})\\
&\propto \exp\Bigg\{   -\frac{1}{2\sigma^2}  \sum_{i,j=1}^{M,N}
\left(-2a_{ij}(\bw_i^\top\bF\bz_j) + (\bw_i^\top\bF\bz_j)^2  \right)\Bigg\}  \cdot  \exp(-\lambda_{kl}^F\cdot f_{kl})u(f_{kl}).\\
\end{aligned}
\end{equation}
To express the conditional density of $\{f_{kl}\mid \bA, \bW,\bF_{-kl},\bZ, \sigma^2,\lambda_{kl}^F \}$ in terms of $f_{kl}$, we write out $\bw_i^\top \bF\bz_j$ in the above equation as
$$
\begin{aligned}
\bw_i^\top \bF\bz_j = 
\sum_{s,t=1}^{K,L}\, w_{is} \, f_{st} z_{tj}
= f_{kl}\, (w_{ik} \, z_{lj}) +C,
\end{aligned}
$$
where 
$$
C=\sum_{(s,t)\neq (k,l)}^{K,L}\, w_{is} \, f_{st} z_{tj}
$$ 
is a constant with respect to $f_{kl}$.
Substituting this into Equation~\eqref{equation:geee_poster_fkl1} and discarding terms independent of $f_{kl}$, we find:
\begin{equation}\label{equation:geee_poster_fkl2}
\begin{aligned}
&\gap p(f_{kl} \mid \bA, \bW,   \bF_{-kl}, \bZ, \sigma^2, \lambda_{kl}^F )\\
&\propto 
\exp\Bigg\{ - 
\underbrace{\frac{\sum_{i,j=1}^{M,N} (w_{ik}z_{lj})^2}{2\sigma^2}}_{\textcolor{mylightbluetext}{\triangleq 1/(2\widetilde{\sigma_{kl}^{2}})}} f_{kl}^2 
+ f_{kl} 
\underbrace{\Bigg[ -\lambda_{kl}^F  + 
\sum_{i,j=1}^{M,N}(w_{ik}z_{lj}) \left( \frac{a_{ij}-C}{\sigma^2} \right)
\Bigg]
}
_{\textcolor{mylightbluetext}{\triangleq \widetilde{\sigma_{kl}^{2}}^{-1} \widetilde{\mu_{kl}}}}
\Bigg\}
\cdot u(f_{kl})\\
&\propto   \normal(f_{kl} \mid  \widetilde{\mu_{kl}}, \widetilde{\sigma_{kl}^{2}})\cdot u(w_{kl}) 
= \truncatednormal(f_{kl} \mid  \widetilde{\mu_{kl}}, \widetilde{\sigma_{kl}^{2}}),
\end{aligned}
\end{equation}
where $u(x)$ enforces nonnegativity, and the posterior ``parent" parameters are:
\begin{align}
\widetilde{\sigma_{kl}^{2}}&= \frac{\sigma^2}{\sum_{i,j=1}^{M,N} (w_{ik}z_{lj})^2}
\label{equation:geee_posterior_variance}\\
\widetilde{\mu_{kl}} &=
\Bigg[ -\lambda_{kl}^F  + 
\sum_{i,j=1}^{M,N}(w_{ik}z_{lj}) \left( \frac{a_{ij}-C}{\sigma^2} \right)
\Bigg]
\cdot \widetilde{\sigma_{kl}^{2}}.
\label{equation:geee_posterior_mean}
\end{align}
Once again, $\truncatednormal(x \mid  \mu, \sigma^2)$ denotes the \text{truncated-normal density} with ``parent" mean $\mu$ and ``parent" variance $\sigma^2$ (Definition~\ref{definition:truncated_normal}).

\paragrapharrow{Sparsity.} Similar to the GEE model on the factored component $w_{mk}$, the posterior parameters have a similar sparsity constraint on the component $f_{kl}$.
The sparsity comes from the negative term $-\lambda_{kl}^F$ in Equation~\eqref{equation:geee_posterior_mean}. When $\lambda_{kl}^F$ becomes larger, the posterior ``parent" mean becomes smaller, and the TN distribution will have a larger probability for smaller values (or even approaching zero) since the draws of $\truncatednormal(f_{kl} \mid \widetilde{\mu_{kl}}, \widetilde{\sigma_{kl}^{2}})$ will be around zero, thus imposing sparsity (see Figure~\ref{fig:dists_truncatednorml_mean}).


\begin{algorithm}[h] 
\caption{Gibbs sampler for GEEE Model in one iteration. 
The noise variance $\sigma^2$
is modeled with an inverse-Gamma prior (analogous formulations apply for precision $\tau$). 
While this implementation prioritizes clarity over efficiency, a vectorized version would be preferable in practice. 
Default uninformative hyper-parameters are $\alpha_\sigma=\beta_\sigma=1$, $\{\lambda_{mk}^W\}=\{\lambda_{kl}^F\}= \{\lambda_{ln}^Z\}=0.1$.} 
\label{alg:geee_gibbs_sampler}  
\begin{algorithmic}[1] 
\Require Choose initial $\alpha_\sigma, \beta_\sigma, \{\lambda_{mk}^W\}, \{\lambda_{kl}^F\}, \{\lambda_{ln}^Z\}$;
\For{$k=1$ to $K$} 
\For{$m=1$ to $M$}
\State Sample $w_{mk}$ from $p(w_{mk} \mid \bA ,   \bW_{-mk}, \bF, \bZ,\sigma^2, \lambda_{mk}^W)$; 
\Comment{Equation~\eqref{equation:gee_poster_wmk1}}
\EndFor
\For{$l=1$ to $L$}
\State Sample $f_{kl}$ from $p(f_{kl} \mid \bA ,   \bW, \bF_{-kl}, \bZ, \sigma^2,\lambda_{kl}^F)$; 
\Comment{Equation~\eqref{equation:geee_poster_fkl2}}
\EndFor
\EndFor
\For{$l=1$ to $L$}
\For{$n=1$ to $N$}
\State Sample $z_{ln}$ from $p(z_{ln} \mid \bA ,   \bW, \bF,  \bZ_{-ln},\sigma^2, \lambda_{ln}^Z)$; 
\Comment{Symmetry of Equation~\eqref{equation:gee_poster_wmk1}}
\EndFor
\EndFor
\State Sample $\sigma^2$ from $p(\sigma^2 \mid  \bA, \bW,\bZ,\alpha_\sigma,\beta_\sigma)$; 
\Comment{Equation~\eqref{equation:gee_posterior_sigma2}}
\State Report loss in Equation~\eqref{equation:als-per-example-loss_tnmf}, stop if it converges;
\end{algorithmic} 
\end{algorithm}

\paragrapharrow{Gibbs sampling.}
Once again, by this Gibbs sampling method introduced in Section~\ref{section:gibbs-sampler}, we can construct a Gibbs sampler for the GEEE model as formulated in Algorithm~\ref{alg:geee_gibbs_sampler}. And also in practice, all the parameters of the exponential distribution are set to a shared value: $\lambda=\{\lambda_{mk}^W\} =\{\lambda_{kl}^F\}= \{\lambda_{ln}^Z\}$ for all $m,k,l,n$.  
By default,
uninformative hyper-parameters are $\alpha_\sigma=\beta_\sigma=1$, $\{\lambda_{mk}^W\}=\{\lambda_{kl}^F\}= \{\lambda_{ln}^Z\}=0.1$.

\index{Automatic relevance determination}
\paragrapharrow{Automatic relevance determination.}
Similar to the GEEA model (Section~\ref{section:geea_nmf_model}), we can use ARD to share the scale parameter of exponential priors for each row of $\bW$ and each column of $\bZ$ so as to perform automatic model selection. The graphical representation is shown in Figure~\ref{fig:bmf_geeea}:
\begin{equation}\label{equation:geeea}
\begin{aligned}
\begin{aligned}
	w_{mk}&\sim \exponential(w_{mk}\mid  \lambda_{k}^W), \qquad 
	&\gap& z_{ln}\sim \exponential(z_{ln}\mid  \lambda_{l}^Z),  \\
	\lambda_{k}^W &\sim \gammadist(\lambda_{k}^W \mid  \alpha_\lambda, \beta_\lambda), \qquad
	&\gap&\lambda_{l}^Z \sim \gammadist(\lambda_{l}^Z \mid  \alpha_\lambda, \beta_\lambda).
\end{aligned}
\end{aligned}
\end{equation}
In this formulation, no parameter sharing is imposed on $\bF$. 
For brevity, we omit further details of this ARD-extended model.

\begin{problemset}
\item Following the derivation in Equation~\eqref{equation:gee_poster_wmk1}, derive the conditional distribution over the user feature $z_{kn}$, for all $k\in\{1,2,\ldots, K\}$ and $n\in\{1,2,\ldots, N\}$, 
under  the GEE model.

\item 
Similarly, following the derivation in Equation~\eqref{equation:gtt_posterior_wmk1}, obtain the conditional distribution over the user feature $z_{kn}$, for all $k\in\{1,2,\ldots, K\}$ and $n\in\{1,2,\ldots, N\}$, in the context of the GTT model.

\item We have derived the variational inference for the GEE model. 
Analogously,  derive the VI updates for the GEEA, GTT, GTTN, GRR, GRRN, GL$_1^2$, GL$_2^2$, GL$_\infty$, GEG, and GnVG models.

\item \textbf{Equivalence of matrix norms.} \label{prob:equiv_mat_norm}
Consider the norm defined in Section~\ref{section:prior_as_reg}.
Let $\norma{\cdot}$ and $\normb{\cdot}$ be two different matrix norms: $\real^{M\times N}\rightarrow \real$. Show that  there exist positive scalars $\alpha$ and $\beta$ such that for all $\bX \in \real^{M\times N}$,
$$
\alpha\norma{\bX} \leq \normb{\bX} \leq \beta\norma{\bX}.
$$
This equivalence implies that if a matrix is small in one norm, it is also small in any other norm---and vice versa.

\item \textbf{Construct  norms from  norms.}\label{prob:cons_ind}
Consider the norm defined in Section~\ref{section:prior_as_reg}.
Let $\norm{\cdot}$ be a matrix norm on $\real^{N\times N}$, and let $\bS\in\real^{N\times N}$ be nonsingular. Show that  the following function defined for $\bA\in\real^{N\times N}$ is also a matrix norm:
$$
\norm{\bA}_{\bS} = \norm{\bS\bA\bS^{-1}}.
$$

\item Following the derivation of nonnegative matrix tri-factorization in Section~\ref{section:nmtf}, formulate a real-valued matrix tri-factorization model in which the data matrix $\bA$ is decomposed as  $\bA=\bW\bF\bZ+\bE$, where $\bW\in\real^{M\times K}, \bF\in\real^{K\times L}$, and $\bZ\in\real^{L\times N}$.

\item Use the ``MovieLens 100K" dataset introduced in Section~\ref{section:movie_rec_als}, evaluate and compare the performance of Bayesian tri-NMF and non-probabilistic tri-NMF methods.

\item Derive the Gibbs sampler for the GEEEA model specified in Equation~\eqref{equation:geeea} and illustrated in Figure~\ref{fig:bmf_geeea}.

\end{problemset}

\chapter{Bayesian Poisson Matrix Factorization}
\begingroup
\hypersetup{
	linkcolor=structurecolor,
	linktoc=page,  
}
\minitoc \newpage
\endgroup

\index{Decomposition: PAA}
\section{Poisson Likelihood with Gamma Priors (PAA)}
\lettrine{\color{caligraphcolor}W}
We introduced the recommendation system problem in Section~\ref{section:movie_rec_als} using standard matrix factorization techniques.
To address this, \citet{gopalan2013scalable,gopalan2015scalable} proposed the \textit{Poisson likelihood with Gamma priors (PAA)} model in the context of recommendation systems, which is later extended to sparse  models in \citet{chang2020sparse} using Horseshoe priors. This model specifically targets challenges posed by nonnegative count data---such as those commonly found in movie recommendation datasets like the Netflix Prize challenge.
The PAA model builds upon \textit{Poisson factorization} \citep{canny2004gap,dunson2005bayesian,cemgil2009bayesian} and has been further explored in \citet{gopalan2014bayesian,hu2015zero}.

The PAA model is designed for nonnegative count data represented by a matrix $\bA\in\naturalset^{M\times N}$, where each entry captures user-item interactions.~\footnote{Items may include movies, songs, articles, or products. For concreteness, we use the term ``interaction," which can refer to any observable action such as a ``click", ``watch", ``purchase", or ``listen." 
In the Netflix context, we specifically refer to movies.}
Here, each entry $a_{mn}$ of $\bA$ (for $m=1,2,\ldots,M$ and $n=1,2,\ldots,N$) denotes the rating (or the interaction score) that user $n$ gave to item $m$, or zero if no rating was provided.
As in standard factorization approaches, we assume $\bA$ can be approximated as the product of two nonnegative matrices: $\bW\in\real_+^{M\times K}$ and $\bZ\in\real_+^{K\times N}$.

More specifically, the PAA model aims to minimize the following loss function:
\begin{equation}\label{equation:poisson_per_example}
\mathop{\min}_{\bW,\bZ} L(\bW,\bZ) =  \mathop{\min}_{\bW,\bZ}\sum_{n=1}^N \sum_{m=1}^{M} \left(a_{mn} - \bw_m^\top\bz_n\right)^2,
\end{equation}
where $\bW=[\bw_1^\top; \bw_2^\top; \ldots; \bw_M^\top]\in \real^{M\times K}$ and $\bZ=[\bz_1, \bz_2, \ldots, \bz_N] \in \real^{K\times N}$, 
with $\bw_m$ and $\bz_n$ representing the \textbf{rows} of $\bW$ and the \textbf{columns} of $\bZ$, respectively. 
Therefore, each item $m$ is represented by a vector of $K$ \textit{latent attributes} $\bw_m$ and each user $n$ by a vector of $K$ \textit{latent preferences} $\bz_n$.
In the Netflix setting, the PAA model is particularly well-suited to capture three key aspects:
(i) the \textit{heterogeneous interests of users} (some users interact with many more items than others),
(ii) the \textit{varying popularity of items} (some movies are inherently more popular), and
(iii) the realistic constraint that users have limited resources (time, attention) to consume content.
Indeed, as noted in the recommendation systems literature, an effective model should account for \textit{heterogeneity} among both users and items \citep{koren2009matrix}.

\begin{figure}[h]
\centering  
\vspace{-0.35cm} 
\subfigtopskip=2pt 
\subfigbottomskip=2pt 
\subfigcapskip=-5pt 
\subfigure[PAA.]{\label{fig:bmf_paa}
\includegraphics[width=0.421\linewidth]{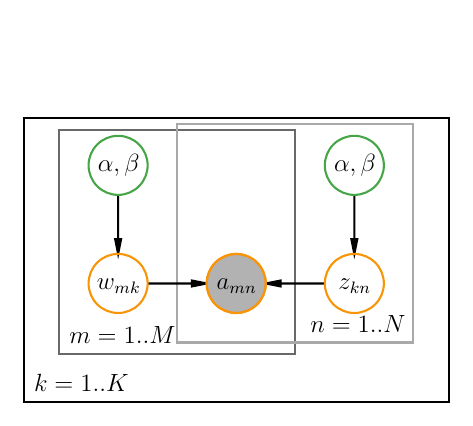}}
\subfigure[PAAA.]{\label{fig:bmf_paaa}
\includegraphics[width=0.421\linewidth]{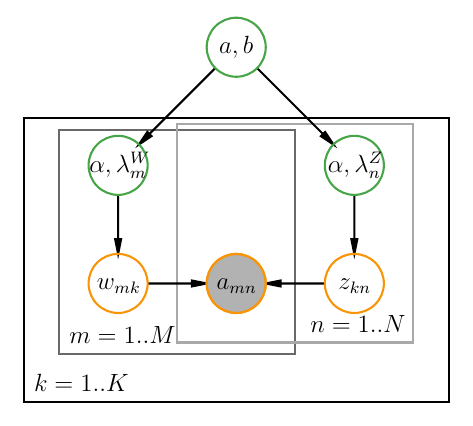}}
\caption{Graphical model representations of PAA and PAAA models. Green circles denote prior variables, orange circles represent observed and latent variables (shaded cycles denote observed variables), and plates represent repeated variables.}
\label{fig:bmf_paa_paaa}
\end{figure}

\paragrapharrow{Likelihood.}
Given a user-item interaction matrix $\bA$, 
where each user has consumed and possibly rated a subset of items, we assume each observed count $a_{mn}$ follows a Poisson distribution with mean $\bw_m^\top\bz_n$ (see Figure~\ref{fig:bmf_paa}):
$$
a_{mn}\sim \poissondist( \bw_m^\top\bz_n),
$$
where $\poissondist(\cdot)$ denotes a Poisson distribution whose parameter is the inner product of the corresponding user preference vector and item attribute vector: $\bw_m^\top\bz_n$ (see Definition~\ref{definition:poisson_distribution}).
This parallels the Gaussian likelihood used in traditional (Bayesian) matrix factorization, where the expected value of  $a_{mn}$ is also modeled as $\bw_m^\top\bz_n$.
Furthermore, suppose we decompose $a_{mn}$ into $K$ latent components:
\begin{equation}
a_{mn} = \sum_{k=1}^{K}o_{mnk}. 
\end{equation}
We then place independent Poisson priors on these components:
\begin{equation}
o_{mnk}\sim \poissondist( w_{mk}z_{kn}).
\end{equation}
By Theorem~\ref{theorem:sum_iid_poisson}, the sum of independent Poisson random variables is itself Poisson-distributed. Hence, we recover the original likelihood:
\begin{equation}\label{equation:amn_poisson}
a_{mn}=\sum_{k=1}^{K}o_{mnk}\sim  \poissondist(    \bw_m^\top\bz_n).
\end{equation}

\paragrapharrow{Prior.}
We place independent Gamma priors on the elements of $\bW$ and $\bZ$, with common shape and rate parameters $\alpha $ and $\beta$, respectively (Definition~\ref{definition:gamma-distribution}):
\begin{equation}
w_{mk}\sim\gammadist(w_{mk}\mid  \alpha, \beta), 
\qquad
z_{kn}\sim \gammadist(z_{kn} \mid  \alpha, \beta).
\end{equation}
This choice of Gamma prior encourages \textbf{sparsity} in the latent representations of users and items (see Figure~\ref{fig:dists_gamma} for illustrative examples), which better reflects real-world behavioral patterns---where users typically engage with only a small subset of available items.

\paragrapharrow{Posterior.}
Let $\bo^{mn}=[o_{mn1}, o_{mn2}, \ldots, o_{mnK}]^\top \in \real^K$,
by Theorem~\ref{theorem:multinomial_poisson}, the conditional distribution of $\bo^{mn}$, given $a_{mn} = \sum_{k=1}^{K}o_{mnk}$, is multinomial:
\begin{equation}\label{equation:paa_bomn}
\multinomial_K(\bo^{mn}\mid  a_{mn}, \bp),
\end{equation}
where $\bp=\frac{1}{\bw_m^\top\bz_n}[w_{m1}z_{1n}, w_{m2}z_{2n}, \ldots, w_{(mK)}z_{(Kn)}]^\top\in [0,1]^K$ such that $\bone^\top\bp=1$, and each element $p_i$ in $\bp$ is in the range of $[0,1]$ (see Definition~\ref{definition:multinomial_dist}).

Using Bayes' rule, the conditional posterior for  $w_{mk}$ is:
\begin{equation}\label{equation:paa_wmk}
\begin{aligned}
&\gap p(w_{mk}\mid  \bA, \bW_{-mk}, \bZ, \alpha, \beta) 
\propto 
\prod_{j=1}^{N} \poissondist(o_{mjk}\mid  w_{mk}z_{kn}) \cdot \gammadist(w_{mk}\mid \alpha, \beta)\\
&\propto 
w_{mk}^{(\sum_{j=1}^{N}o_{mjk})} \exp\Bigg\{ -w_{mk}\Bigg( \sum_{j=1}^{N}z_{kj} \Bigg) \Bigg\}
\cdot w_{mk}^{\alpha-1} \exp\left( -\beta w_{mk}\right)
\propto \gammadist(w_{mk}\mid  \widetilde{\alpha}, \widetilde{\beta}),
\end{aligned}
\end{equation}
with updated parameters
\begin{equation}\label{equation:paa_wmk_param}
\begin{aligned}
\widetilde{\alpha} = \alpha+\sum_{j=1}^{N}o_{mjk}, 
\qquad 
\widetilde{\beta} = \beta + \sum_{j=1}^{N}z_{kj} .
\end{aligned}
\end{equation}
According to the definition of the Gamma distribution (Definition~\ref{definition:gamma-distribution}), the posterior mean of $w_{mk}$ is given by 
$$
\Exp[w_{mk}\mid  \bA, \bW_{-mk}, \bZ, \alpha, \beta] =\frac{\alpha+\sum_{j=1}^{N}o_{mjk}}{\beta + \sum_{j=1}^{N}z_{kj}}.
$$
This expression is intuitive: a large total count  $\sum_{j=1}^{N}o_{mjk}$ (reflecting strong evidence for attribute $k$ in item $m$) increases the expected value of  $w_{mk}$. 
Conversely, a large sum $\sum_{j=1}^{N}z_{kj}$ (indicating high user affinity for attribute $k$) acts as a normalizing factor, tempering the estimate.
By symmetry,  analogous conditional posteriors can be derived for variables $\{z_{kn}\}$.

\begin{algorithm}[h] 
\caption{Gibbs sampler for PAA model in one iteration.  
By default,
uninformative hyper-parameters are $\alpha=\beta=1$.} 
\label{alg:paa_gibbs_sampler}  
\begin{algorithmic}[1] 
\Require Choose initial $\alpha, \beta$;
\For{$m=1$ to $M$} 
\For{$n=1$ to $N$}
\State Sample $\bo^{mn}$ from $p(\bo^{mn}\mid  a_{mn}, \bp)$; 
\Comment{Equation~\eqref{equation:paa_bomn}}
\EndFor
\EndFor
\For{$k=1$ to $K$} 
\For{$m=1$ to $M$}
\State Sample $w_{mk}$ from $p(w_{mk} \mid   \bA,\bW_{-mk}, \bZ, \alpha, \beta)$; 
\Comment{Equation~\eqref{equation:paa_wmk}}
\EndFor
\For{$n=1$ to $N$}
\State Sample $z_{kn}$ from $p(z_{kn} \mid  \bA,\bW, \bZ_{-kn}, \alpha, \beta)$; 
\Comment{Symmetry of Eq.~\eqref{equation:paa_wmk}}
\EndFor
\EndFor
\end{algorithmic} 
\end{algorithm}

\paragrapharrow{Gibbs sampling.}
Using the Gibbs sampling framework introduced in Section~\ref{section:gibbs-sampler}, we obtain the procedure outlined in Algorithm~\ref{alg:paa_gibbs_sampler}. In practice, weakly informative hyper-parameters such as $\alpha=\beta=1$ are commonly used to initialize the model.

\index{Decomposition: PAAA}
\section{PAA Model with Hierarchical Gamma Priors (PAAA)}
Extending the PAA model, \citet{gopalan2015scalable} introduced the \textit{Poisson likelihood with Gamma priors and hierarchical Gamma priors (PAAA)} model.

\paragrapharrow{Prior and diversity.}
Building on the PAA framework, the PAAA model introduces a hierarchical Gamma prior over the rate parameters of the latent factors:
\begin{equation}
\begin{aligned}
&a_{mn}\sim \poissondist(a_{mn}\mid  \bw_m^\top\bz_n), \gap\\
&w_{mk}\sim\gammadist(w_{mk}\mid  \alpha, \lambda_m^W), 
&\gap& z_{kn}\sim \gammadist(z_{kn} \mid  \alpha, \lambda_n^Z),\\
&\lambda_m^W\sim \gammadist(a, \frac{a}{b}), 
&\gap&\lambda_n^Z\sim \gammadist(a, \frac{a}{b}).
\end{aligned}
\end{equation}
This hierarchical structure enables the model to capture two key real-world phenomena:
(i) the \textit{diversity among users}---some users interact with many more items than others---and
(ii) the \textit{diversity among items}---some items are significantly more popular than others.

\paragrapharrow{Posterior.}
The conditional posteriors for $\{w_{mk}\}, \{z_{kn}\}$, and $\{\bo^{mn}\}$ are identical to those in the PAA model, except that the global rate parameter $\beta$ is replaced by the user- or item-specific rate parameters $\lambda_m^W$ or $\lambda_n^Z$ in the expression for $\widetilde{\beta}$ \eqref{equation:paa_wmk_param}.
For the hierarchical parameters, we derive the conditional posterior of  $\lambda_m^W$ using Bayes' rule:
\begin{equation}\label{equation:paaa_pos_lambda_mw}
\begin{aligned}
p(\lambda_m^W\mid&  \bW,\alpha,a,b )\propto
\prod_{k=1}^{K} \gammadist(w_{mk}\mid \alpha, \lambda_m^W) \cdot \gammadist(\lambda_m^W\mid a, \frac{a}{b})\\ 
&\propto \prod_{k=1}^{K} \frac{(\lambda_m^W)^\alpha}{\Gamma(\alpha)}
w_{mk}^{\alpha-1} \exp(-\lambda_m^W w_{mk}) 
\cdot \frac{(\frac{a}{b})^a}{\Gamma(a)} (\lambda_m^W)^{a-1}\exp(-\frac{a}{b} \lambda_m^W)\\
&\propto (\lambda_m^W)^{K\alpha+a-1} 
\exp \Bigg\{ -\lambda_m^W \Bigg( \frac{a}{b} +\sum_{k=1}^{K} w_{mk} \Bigg) \Bigg\}
\propto \gammadist(\lambda_m^W \mid  \widetilde{a}_m, \widetilde{b}_m),
\end{aligned}
\end{equation}
where the updated shape and rate parameters are
$$
\widetilde{a}_m=K\alpha+a, 
\qquad 
\widetilde{b}_m=\frac{a}{b} +\sum_{k=1}^{K} w_{mk} .
$$

\paragrapharrow{Gibbs sampling.}
A Gibbs sampler for the PAAA model is given in Algorithm~\ref{alg:paaa_gibbs_sampler}. In practice, weakly informative hyper-parameters such as  $\alpha=a=b=1$ are commonly used to initialize the model.

\begin{algorithm}[h] 
\caption{Gibbs sampler for PAAA model in one iteration. 
By default,
uninformative hyper-parameters are $\alpha=a=b=1$.} 
\label{alg:paaa_gibbs_sampler}  
\begin{algorithmic}[1] 
\Require Choose initial $\alpha,a, b$;
\For{$m=1$ to $M$} 
\For{$n=1$ to $N$}
\State Sample $\bo^{mn}$ from $p(\bo^{mn}\mid  a_{mn}, \bp)$; 
\Comment{Equation~\eqref{equation:paa_bomn}}
\EndFor
\EndFor
\For{$k=1$ to $K$} 
\For{$m=1$ to $M$}
\State Sample $w_{mk}$ from $p(w_{mk} \mid   \bA,\bW_{-mk}, \bZ, \alpha, \lambda_m^W)$; 
\Comment{Eq.~\eqref{equation:paa_wmk}, replace $\beta$ by $\lambda_m^W$}
\State Sample $\lambda_{m}^W$ from $p(\lambda_m^W\mid  \bW,\alpha,a,b )$; 
\Comment{Equation~\eqref{equation:paaa_pos_lambda_mw}}
\EndFor
\For{$n=1$ to $N$}
\State Sample $z_{kn}$ from $p(z_{kn} \mid  \bA,\bW, \bZ_{-kn}, \alpha, \lambda_n^Z)$; 
\Comment{Eq.~\eqref{equation:paa_wmk},  replace $\beta$ by $\lambda_n^Z$}
\State Sample $\lambda_{n}^Z$ from $p(\lambda_{n}^Z\mid  \bZ,\alpha,a,b )$; 
\Comment{Symmetry of Eq.~\eqref{equation:paaa_pos_lambda_mw}}
\EndFor
\EndFor
\end{algorithmic} 
\end{algorithm}

\begin{SCfigure}
\centering  
\subfigtopskip=2pt 
\subfigbottomskip=6pt 
\subfigcapskip=-15pt 
\includegraphics[width=0.621\textwidth]{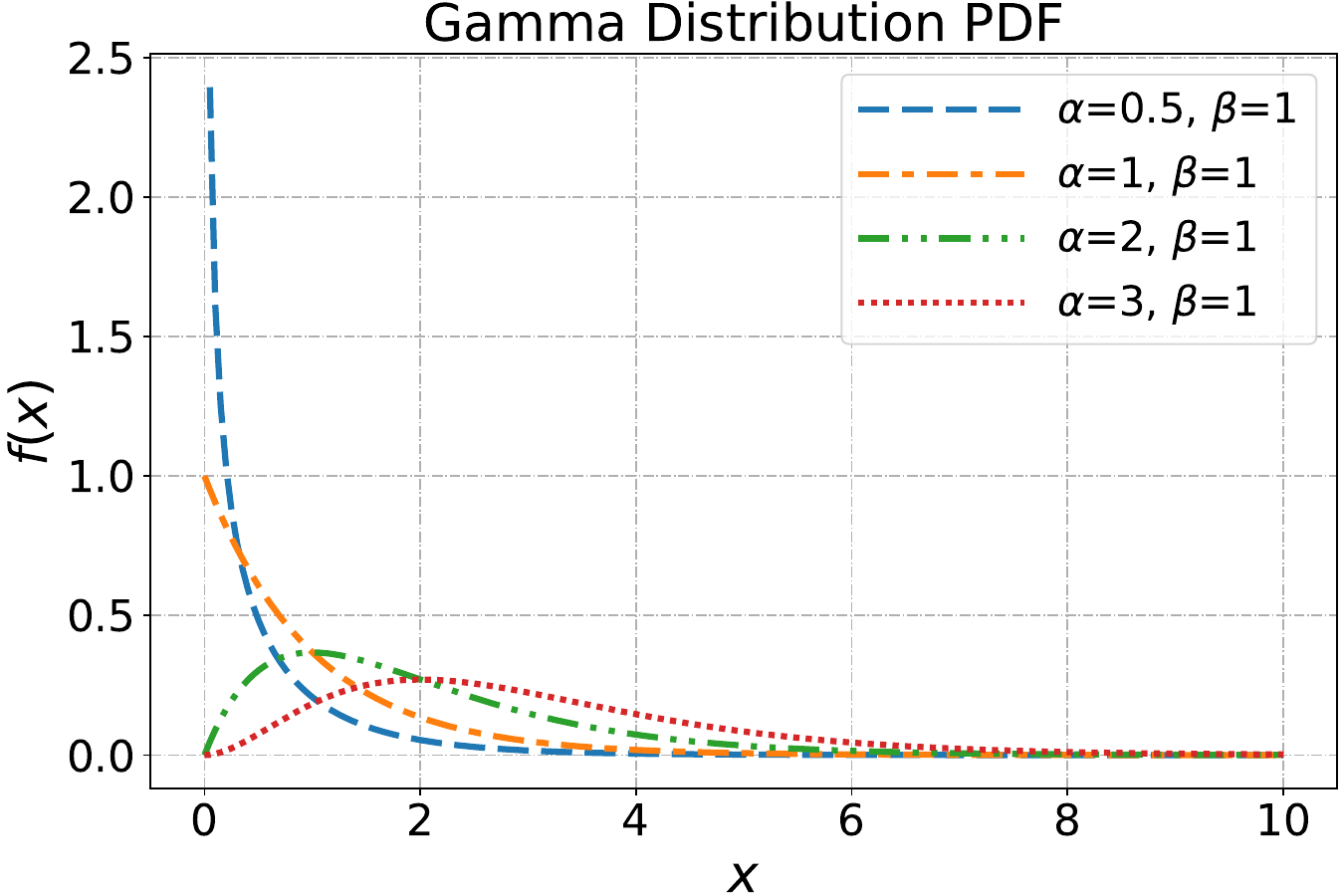}
\caption{Gamma probability density functions $\gammadist(\alpha, \beta)$ by reducing the shape parameter $\alpha$. }
\label{fig:dists_gamma_paa}
\end{SCfigure}

\section{Properties of PAA or PAAA}
Having presented the modeling details, we now highlight key statistical properties of the PAA and PAAA approaches. These features offer distinct advantages over Gaussian-likelihood matrix factorization methods---particularly in the context of implicit feedback data like that in the Netflix challenge.

\paragrapharrow{PAA or PAAA encourage sparse latent representations.}
As noted earlier, the Gamma priors placed on user preferences ($\bz_n$) and item attributes ($\bw_m$) naturally promote sparsity.
When the shape parameter $\alpha$ of the Gamma distribution is small, most latent weights are driven close to zero, with only a few taking substantial values (see Figure~\ref{fig:dists_gamma_paa}, which shows $\gammadist(\alpha, \beta=1)$ for $\alpha$=3, 2, 1, 0.5).
This yields simpler, more interpretable models where each user or item is characterized by only a handful of active latent dimensions.

\begin{figure}[h]
\centering  
\subfigtopskip=2pt 
\subfigbottomskip=2pt 
\subfigcapskip=-3pt 
\subfigure[User activity.]{\label{fig:movielen_1m_user_activity}
\includegraphics[width=0.481\linewidth]{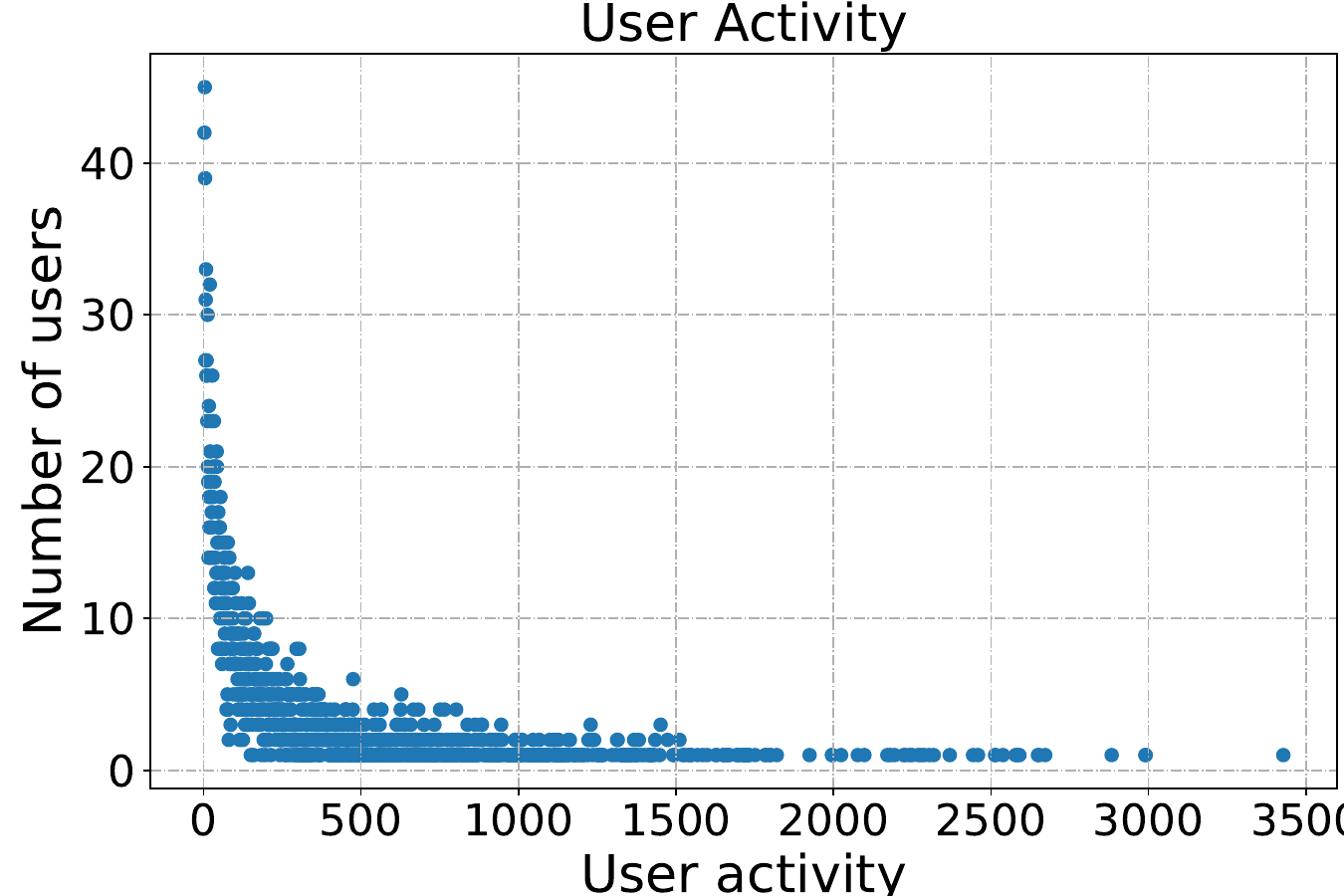}}
\subfigure[Item popularity.]{\label{fig:movielen_1m_item_activity}
\includegraphics[width=0.481\linewidth]{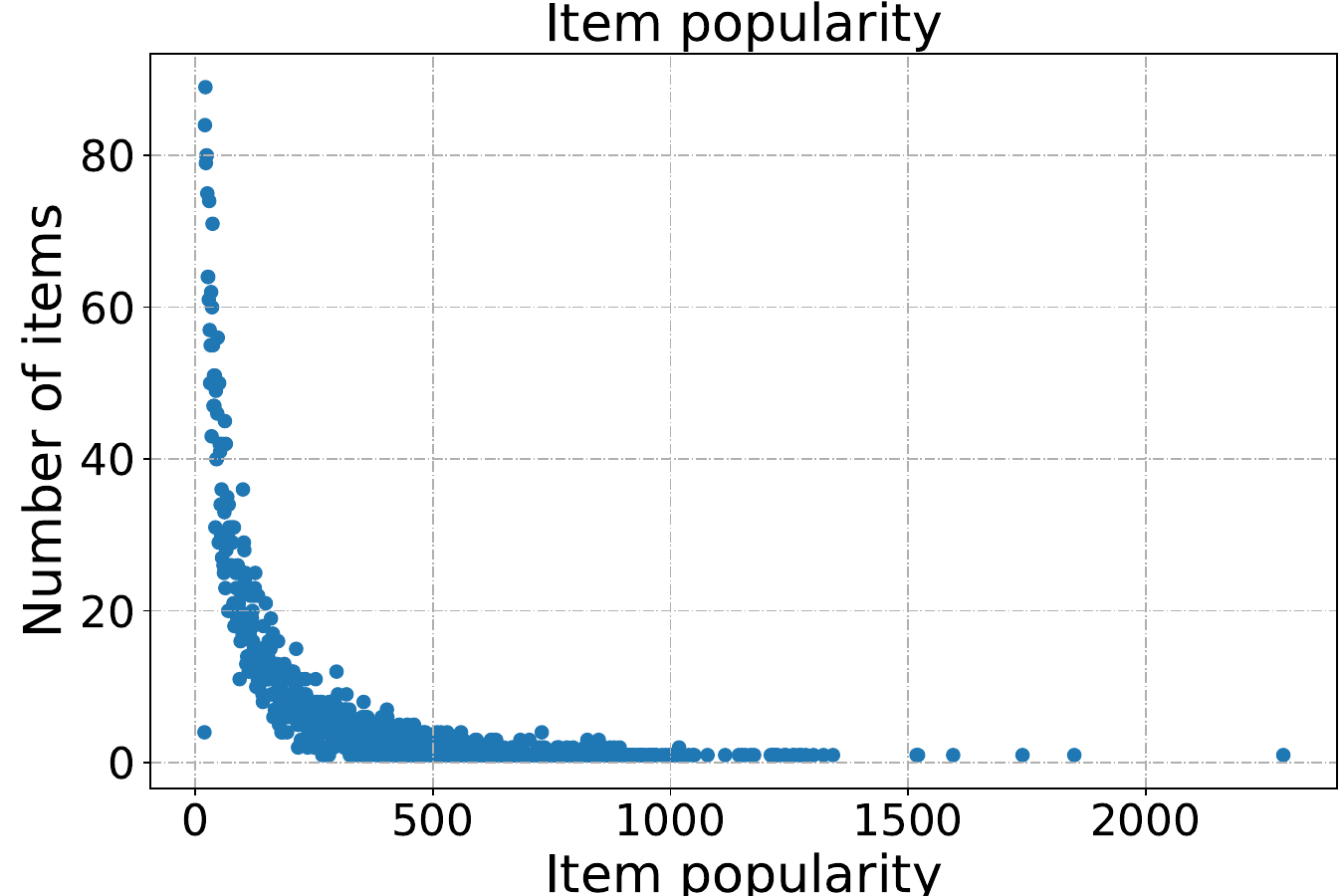}}
\caption{User activity and item popularity for the ``MovieLens 1M" data set (see data description in Table~\ref{table:datadescription}).}
\label{fig:movielen_1m_item_user}
\end{figure}

\paragrapharrow{PAA or PAAA models  long-tailed of user and item behavior.}
In implicit feedback settings---where $a_{mn}=1$ if user $n$ consumed item $m$, and $0$ otherwise~\footnote{In contrast, explicit feedback consists of explicit ratings (e.g., 1--5 stars).}---real-world interaction data typically exhibits a long-tailed distribution:
most users interact with only a few items, while a small fraction (``heavy users") interact with many; similarly, most items receive few interactions, while a few ``blockbuster" items are widely consumed.
To illustrate, consider the ``MovieLens 1M" dataset (Table~\ref{table:datadescription}), which contains ratings from 6,040 movies on 3,503 users (after pre-processing). 
Figure~\ref{fig:movielen_1m_user_activity} shows that only a small minority of users rated more than 1,500 movies, and Figure~\ref{fig:movielen_1m_item_activity} reveals that very few movies were rated by more than 500 users---clear evidence of long-tailed behavior.

The PAA and PAAA models capture this structure through a \textbf{two-stage generative process}.
From Equation~\eqref{equation:amn_poisson}, for each user   $n$, Theorems~\ref{theorem:sum_iid_poisson} and \ref{theorem:multinomial_poisson} imply:
\begin{equation}\label{equation:poisson_user_n}
\begin{aligned}
u_n = \sum_{i=1}^{M}a_{in}  &\sim \poissondist\Bigg(\sum_{i=1}^{M} \bw_i^\top\bz_n\Bigg),\\
\ba_n=[a_{1n}, a_{2n }, \ldots, a_{Mn}]^\top &\sim
\multinomial_M (u_n, \bq),
\end{aligned}
\end{equation}
where $\bq=\frac{1}{\sum_{i=1}^{M} \bw_i^\top\bz_n}[ \bw_1^\top\bz_n,  \bw_2^\top\bz_n, \ldots,  \bw_M^\top\bz_n]^\top \in [0,1]^M$ such that $\bone^\top\bq=1$. 
Thus, the PAA or PAAA  model first draws a total \textit{activity budget} $u_n$ for each user $n$, then allocates this budget across items according to $\bq$.
Learning this budget value is important for modeling the long-tail behavior of user activity.
Learning this per-user budget is crucial for modeling heterogeneous activity levels and the long tail of user behavior---something standard real-valued or nonnegative matrix factorization methods do not explicitly account for.

Similarly for each item $m$, we have
\begin{equation}\label{equation:poisson_item_m}
\begin{aligned}
v_m = \sum_{i=1}^{N}a_{mi}  &\sim \poissondist\Bigg(\sum_{i=1}^{N} \bw_m^\top\bz_i\Bigg),\\
\widehat{\ba}_m=[a_{m1}, a_{m2}, \ldots, a_{mN}]^\top &\sim
\multinomial_N (v_m, \bs),
\end{aligned}
\end{equation}
where $\bs=\frac{1}{\sum_{i=1}^{N} \bw_m^\top\bz_i}[  \bw_m^\top\bz_1,   \bw_m^\top\bz_2, \ldots,   \bw_m^\top\bz_N]^\top \in [0,1]^N$
such that $\bone^\top\bs=1$.
The PAA  or PAAA model finds the \textit{popularity} of item $m$ by $v_m$, and then learns how the popularity is distributed across users.
Again, the Poisson-Multinomial decomposition allows the model to naturally reflect the skewed, long-tailed popularity of items.

\section{Recommendation Systems}\label{section:recom_poisson}
In Section~\ref{section:movie_rec_als}, we introduced two recommendation systems based on matrix factorization.
In the following paragraphs, we briefly discuss these approaches and then propose a new recommender built upon the Bayesian matrix factorization framework.

\paragrapharrow{Recommender 1.}
A simple recommender can suggest an unconsumed movie $m$ to user $n$ by ranking items according to the posterior expected value of the Poisson rate parameter:
\begin{equation}\label{equation:recom_poisson1}
\text{score}_{mn} = \Exp[\bw_m^\top\bz_n \mid \bA].
\end{equation}
This score can be approximated by averaging $\Exp[\bw_m^\top\bz_n \mid \bA]$ over Gibbs sampling iterations after convergence.

\paragrapharrow{Recommender 2.}
After inferring the item attribute vectors $\{\bw_1, \bw_2, \ldots, \bw_M\}$, we can compute a similarity matrix between items (e.g., using Pearson correlation or cosine similarity). For each user $n$, we then recommend items that are highly similar to those they have already consumed.
A precision-recall (PR) curve can be used to select an appropriate similarity threshold for final recommendations; see Section~\ref{section:movie_rec_als}.

\paragrapharrow{Recommender 3.}
In movie recommendation, uncertainty about each unobserved entry $a_{mn}$ can be quantified by its predictive standard deviation. A practical system may choose to recommend only those items for which the prediction is highly confident.
The Bayesian framework naturally supports such uncertainty-aware recommendations. Inspired by the \textit{Sharpe ratio} from quantitative finance---which measures risk-adjusted return as the ratio of expected return to standard deviation (the higher the Sharpe ratio, the better the risk-adjusted return of the investment is considered to be)---we define an \textit{uncertainty-adjusted score}:
\begin{equation}
	\text{score}_{mn} = \frac{\Exp[\bw_m^\top\bz_n \mid \bA]}{\sqrt{\Var[\bw_m^\top\bz_n \mid \bA]}}.
\end{equation}
This score prioritizes items with high expected interaction rates relative to their prediction uncertainty. Higher values indicate more reliable recommendations.

\index{Variational autoencoder}
\index{Multinomial generation}
\section{Variational Autoencoder with  Multinomial Generation}\label{section:vae_multino_gen}
The Poisson model is closely related to the multinomial likelihood: as shown in Equation~\eqref{equation:poisson_user_n}, a user's preferences over $M$ items can be modeled via a multinomial distribution conditioned on their total activity level.
This multinomial likelihood can be directly incorporated into a variational autoencoder (VAE) framework (Section~\ref{section:vae_pca}), which provides a form of amortized inference \citep{liang2018variational}.
The graphical model is depicted in Figure~\ref{fig:lvm_VAE_and_flow}, and the generative process is defined as follows:
\begin{equation}
\begin{aligned}
&\bz_n\sim\normal(\bzero,\bI_K), 
\gap
\pi_{\btheta}(\bz_n)\propto \exp\{f_{\btheta}(\bz_n)\}, \\
\gap
&\ba_n \sim p_{\btheta}(\bz_n)= \multinomial_{M}(N_n, \pi_{\btheta}(\bz_n)),
\gap
n\in\{1,2,\ldots,N\},
\end{aligned}
\end{equation}
where  $\bz_n\in\real^{K}$ is the latent user attribute for user $n$ with $K<\min\{M,N\}$;
$\pi_{\btheta}(\bz_n)\in\real^M$ is a vector on the probability simplex, assigning higher mass to items the user is more likely to prefer; and $\multinomial_{M}(\cdot, \cdot)$ denotes  the multinomial distribution (see Definition~\ref{definition:multinomial_dist}).

The function $p_{\btheta}(\ba_n\mid\bz_n)$ is the generative model that produces the observed data based on the hidden vector $\bz_n$.
For example, let $f_{\btheta}(\bz_n)=[f_{1n},f_{2n},\ldots,f_{Mn}]\in\real^M$ where the parameter $\btheta$ can be derived  from deep neural networks, the generative function $p_{\btheta}(\cdot)$ can be taken as the Gaussian distribution
\begin{equation}
\ln p_{\btheta}(\ba_n\mid\bz_n)
= -\sum_{m=1}^{M}\frac{c_{mn}}{2}(a_{mn}-f_{mn})^2,
\end{equation}
where  $c_{mn}$ is a constant that can be used to weight the contributions of different items.
Alternatively, the generative distribution can be modeled using a logistic likelihood (under the Bernoulli likelihood):
\begin{equation}
\ln p_{\btheta}(\ba_n\mid\bz_n)
= 
\sum_{m=1}^{M} a_{mn} \ln \sigma(f_{mn})+(1-a_{mn})\ln(1-\sigma(f_{mn})),
\end{equation}
where $\sigma(x)=1/(1+\exp\{-x\})$ is the logistic sigmoid function.
In both cases, $f_{\btheta}(\bz_n)$ represents the output of the neural network, which maps the latent vector $\bz_n$ to a set of scores or probabilities for the  $M$ items. The Gaussian distribution is suitable for continuous data, while the logistic likelihood is appropriate for binary data, such as whether an item was interacted with or not (i.e., the implicit data in the recommendation context).

And $N_n=\sum_{m=1}^{M} a_{mn}$ is the total number of user interactions for user $n$, e.g., total number of clicks, watches, or purchases. 
The PAA or PAAA models first learn a budget $u_n$ for each user $n$ (Equation~\eqref{equation:poisson_user_n}) and then distribute the budget across items. However, in the VAE model, the budge $N_n$ for each user $n$ is fixed (and observed beforehand). To address this, the VAE model takes the average over multiple samples of the prediction  $\ba_n$  for each user $n$. This approach helps in handling the fixed and observed budget by incorporating the variability in the latent space.

The goal of this problem then becomes estimating the posterior distribution $p_{\btheta}(\bz_n \mid \ba_n)$.
The VAE solves this by approximating this intractable posterior distribution $p_{\btheta}(\bz_n \mid \ba_n)$ by a variational distribution $q_{\blambda}(\bz_n \mid\ba_n)=\normal(\bz_n \mid \bmu_{\blambda}(\ba_n), \diag(\bsigma^2_{\blambda}(\ba_n)))$ indexed by parameter $\blambda$.
With this reparameterization trick and stochastic gradient descent with MC gradient introduced in Section~\ref{section:vae_pca}, we can optimize the ELBO over $\btheta, \blambda$ and find the approximation $q_{\blambda}(\bz_n \mid\bz_n)$.

As of the prediction $\widehat{\ba}_n$ for each user $n$, we set the latent vector as the mean of the distribution $\bz_n=\bmu_{\blambda}(\ba_n)$ (the encoder). 
As mentioned above, we  generate the    prediction $\widehat{\ba}_n$ as a mean of a series of samples from $\widehat{\ba}_n = \frac{1}{S}\sum_{s=1}^{S}\widehat{\ba}_n^{(s)}$, where $\widehat{\ba}_n^{(s)}\sim\multinomial_{M}(N_n, \pi_{\btheta}(\bz_n))$ (the decoder). This is because the number of total interactions $N_n$ for each user $n$ is fixed a priori and measures each user's activity from the observed data.

\index{Decomposition: OGGW}
\section{Ordinal Likelihood with Gaussian  and Wishart  Priors (OGGW)}
The properties of Poisson factorization (PF) models---such as the PAA and PAAA models---show that their primary objective is to recommend items by predicting future user–item interactions.
Consequently, PF is typically applied to \textit{implicit consumer data}, where the observed data matrix $\bA\in\{0,1\}^{M\times N}$ indicates only whether a user has interacted with an item (1) or not (0).

In many real-world applications, however, the data matrix $\bA$ contains richer information.
Ordinal matrix factorization (OMF) addresses this by handling \textit{ordinal} or \textit{explicit} feedback data \citep{stevens1946theory}, where entries in $\bA$ take values from a finite, ordered set reflecting user preferences.
For example, in collaborate filtering, we seek to predict a consumer's rating of a novel item on an ordinal scale such as \textit{good} $>$ \textit{average} $>$ \textit{bad}; 
the temperature of a day  is \textit{hot} $>$ \textit{warm} $>$ \textit{cold};
a teacher always rates his/her students by giving grades on their overall performance having the ordering 
$A>B>C>D>F$
\citep{paquet2005bayesian, chu2005gaussian, gouvert2020ordinal}.
While standard real-valued or nonnegative matrix factorization methods (introduced in earlier chapters) can be adapted to such data, approaches that explicitly model the ordinal nature of the observations are generally more effective and statistically principled.

Rather than directly factorizing the observed matrix as $\bA=\bW\bZ+\bE$  (as in traditional matrix factorization), \textit{Bayesian ordinal matrix factorization} introduces an additional 
\textit{latent} (unobserved) continuous matrix $\bH=\bW\bZ+\bE\in\real^{M\times N}$. 
This hidden matrix $\bH$ serves as input to an ordinal regression model, which maps each latent value $h_{mn}$ to a probability distribution over the discrete ordinal categories, thereby generating the observed data $\bA$ (see Figure~\ref{fig:bmf_oggw}).
Because of this two-level structure---latent factors feeding into a probabilistic observation model---ordinal matrix factorization is also referred to as a hierarchical Bayesian model \citep{paquet2012hierarchical}.

\subsection{Ordinal Regression Likelihood}
\begin{figure}[htp]
\centering  
\subfigtopskip=0pt 
\subfigbottomskip=6pt 
\subfigcapskip=-15pt 
\includegraphics[width=0.421\textwidth]{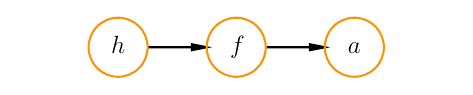}
\caption{Graphical representation of the ordinal regression model. }
\label{fig:bmf_oggw_ordireg}
\end{figure}

We now consider a data matrix $\bA\in\sA^{M\times N}$, where $\sA$ is a finite set of $A$ ordered categories. Without loss of generality, we encode these categories as consecutive integers: $\sA=\{1,2,\ldots, A\}$, preserving their inherent ordering.
The real line is partitioned into $A$ contiguous intervals using thresholds $\{b_a\}_{a=1}^{A+1}$, defined as:
$$
-\infty =b_1 <b_2 <\ldots <b_{A+1}=\infty,
$$
such that the interval $[b_a, b_{a+1})$ corresponds to the discrete category $a\in \sA$.
To model the mapping from latent variables $h$ to ordinal outcomes, we introduce an auxiliary continuous variable  $f$ (see Figure~\ref{fig:bmf_oggw_ordireg}). 
The observed category $a$ is determined by which interval $f$ falls into:
\begin{equation}
p(a\mid f)=
\left\{
\begin{aligned}
	&1, \, &\text{if } b_a\leq f < b_{a+1} \\
	&0, \, &\text{otherwise}
\end{aligned}
\right.
=u(f-b_a) - u(f-b_{a+1}),
\end{equation}
where $u(y)$ is the unit step function with a value 1 if $y\geq 0$ and  0 if $y<0$.

\begin{SCfigure}
\centering  
\subfigtopskip=2pt 
\subfigbottomskip=6pt 
\subfigcapskip=-15pt 
\includegraphics[width=0.621\textwidth]{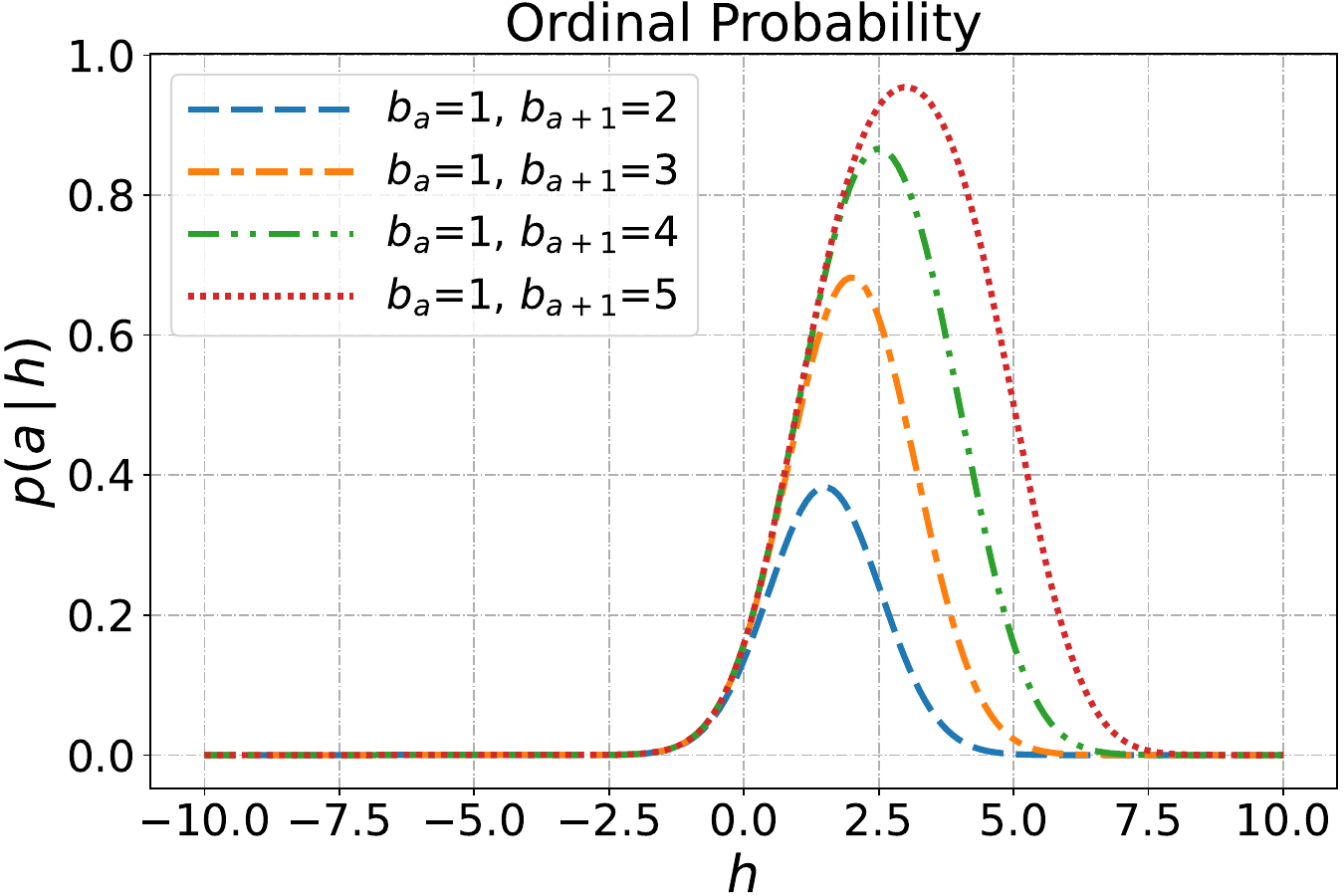}
\caption{Ordinal category probability $p(a\mid h)$ of the ordinal regression model in Equation~\eqref{equation:oggw_averaging_f}, shown as functions of the latent variable $h$. }
\label{fig:dists_oggw_ordireg}
\end{SCfigure}

Given the hidden value $h$, uncertainty about the exact location of $f$ can be modeled by a unit-variance Gaussian:
\begin{equation}
p(f\mid h) = \normal(f\mid h,1).
\end{equation}
Averaging over $f$ in $p(a, f\mid h)=p(a\mid f)p(f\mid h)$, we have  
\begin{equation}\label{equation:oggw_averaging_f}
p(a\mid h) = \int p(a, f\mid h)df = \Phi(h-b_a) - \Phi(h-b_{a+1}),
\end{equation}
where
$\Phi(y) = \int_{-\infty}^{y} \normal(u\mid 0,1)du= \frac{1}{\sqrt{2\pi}} \int_{-\infty}^{y} \exp(-\frac{u^2}{2}) du $ denotes the cumulative distribution function  of the standard normal distribution  $\normal(0,1)$.
In Equation~\eqref{equation:oggw_averaging_f}, we use the fact that 
$$
\Phi(h-b) = \int \normal(f\mid h,1)\,u(f-b) df
$$
(see \citet{albert1993bayesian}).
Figure~\ref{fig:dists_oggw_ordireg} shows the probability functions for $p(a\mid h)$ by varying $h$. 
Note that even when  $h$ lies outside the interval [$b_a, b_{a+1}$), the probability is not exactly zero due to the smoothing effect of the Gaussian noise. Moreover, when the interval width $b_{a+1}-b_a$ is small (which occurs when there are many ordinal categories), the peak probability for any single category tends to be lower---reflecting greater uncertainty in fine-grained distinctions (i.e., the probability tends to be small for falling into each interval).

\paragrapharrow{Complete likelihood.}
The ordinal regression model assigns to each entry a probability based on its latent value $h_{mn}$: maps continuous latent variables $h_{mn}$ in $\bH$ to probabilities $p(a_{mn} \mid  h_{mn})$:
\begin{equation}
p(a_{mn} \mid  h_{mn}) =
\prod_{a=1}^{A} \left[ \Phi(h_{mn} - b_a) -\Phi(h_{mn} - b_{a+1})\right]^{\indicator(a_{mn}=a)},
\end{equation}
where $\indicator(\cdot)$ is the indicator function.
Let  $\mathcalX = \{a_{mn} \mid  (m,n)\in \textit{training set}\}$ denote the set of observed entries. 
The full likelihood of the observed data under the model is then:
\begin{equation}
p(\mathcalX \mid  \bH) = \prod_{(m,n)\in\mathcalX} p(a_{mn} \mid  h_{mn}),
\end{equation}
where the product is over all the observed entries or training set entries $(m,n)$.

\subsection{Matrix Factorization Modeling on Latent Variables}
\begin{figure}[h]
\centering  
\subfigtopskip=2pt 
\subfigbottomskip=6pt 
\subfigcapskip=-15pt 
\includegraphics[width=0.421\textwidth]{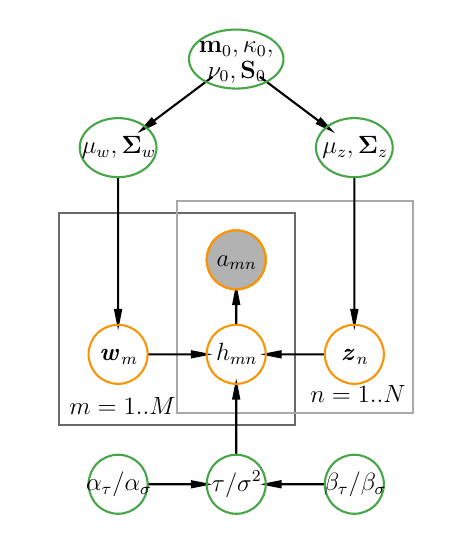}
\caption{Graphical representation of OGGW model. Green circles denote prior variables, orange circles represent observed and latent variables (shaded cycles denote observed variables), and plates represent repeated variables. The slash ``/" in the variable represents ``or," and the comma ``," in the variable represents ``and."}
\label{fig:bmf_oggw}
\end{figure}

The Bayesian treatment of the latent matrix $\bH$ mirrors that of the GGGW model (Section~\ref{section:gggw_model}), with one key distinction: here, a Gaussian likelihood is placed on the latent variables $\{h_{mn}\}$, rather than directly on the observed ordinal ratings $\{a_{mn}\}$. 
The complete graphical model---shown in Figure~\ref{fig:bmf_oggw}---is known as the \textit{ordinal likelihood with Gaussian and hierarchical normal-inverse-Wishart priors (OGGW)} model.

\paragrapharrow{Likelihood.}
We assume the residuals, $e_{mn}=h_{mn}-\bw_m^\top\bz_n$, are i.i.d. zero-mean normal with precision $\tau = {1}/{\sigma^2}$.
This yields the following likelihood:
\begin{equation}\label{equation:oggw_likelihood}
	\begin{aligned}
		p(\bH\mid  \bW,\bZ,\tau) &= \prod_{m,n=1}^{M,N} \normal \left(h_{mn}\mid  (\bW\bZ)_{mn}, \sigma^2 \right)\\
		&= \prod_{m,n=1}^{M,N} \normal \left(h_{mn}\mid  (\bW\bZ)_{mn}, \tau^{-1} \right),
	\end{aligned}
\end{equation}
where $\sigma^2$ is the noise variance and $\tau$ is the corresponding precision.

\index{Inverse-Wishart distribution}
\paragrapharrow{Prior.}
Given the $m$-th row $\bw_m$ of $\bW$ and the $n$-th column $\bz_n$ of $\bZ$, 
we place multivariate Gaussian priors with shared hyper-parameters governed by normal-inverse-Wishart prior as follows:
\begin{equation}\label{equation:oggw_prior_wm_zn}
\begin{aligned}
	\bw_m\sim \normal(\bw_m\mid  \bmu_w, \bSigma_w), &\gap\gap \bmu_w, \bSigma_w \sim \niw(\bmu_w, \bSigma_w\mid  \bmm_0, \kappa_0, \nu_0, \bS_0);\\
	\bz_n\sim \normal(\bz_n\mid  \bmu_z, \bSigma_z), &\gap\gap \bmu_z, \bSigma_z \sim \niw(\bmu_z, \bSigma_z\mid  \bmm_0, \kappa_0, \nu_0, \bS_0),\\
\end{aligned}
\end{equation}
where $\niw (\bmu, \bSigma\mid  \bmo, \kappa_0, \nu_0, \bso) 
= \mathcal{N}(\bmu\mid  \bmo, \frac{1}{\kappa_0}\bSigma) \cdot  \inversewishart(\bSigma\mid  \bso, \nu_0)$ is the density of a normal-inverse-Wishart distribution, and $ \inversewishart(\bSigma\mid  \bso, \nu_0)$ is the inverse-Wishart distribution (Equation~\eqref{equation:multi_gaussian_prior}).

Once again, the prior for the noise variance $\sigma^2$ is chosen as a conjugate inverse-Gamma density with shape ${\alpha_\sigma}$ and scale ${\beta_\sigma}$ (Definition~\ref{definition:inverse_gamma_distribution}),
$$
p(\sigma^2)= \inversegammadist(\sigma^2\mid  \alpha_\sigma, \beta_\sigma) = \frac{{\beta_\sigma}^{\alpha_\sigma}}{\Gamma({\alpha_\sigma})} (\sigma^2)^{-\alpha_\sigma-1} \exp\left( -\frac{{\beta_\sigma}}{\sigma^2} \right).
$$
Equivalently, placing an inverse-Gamma prior on the variance corresponds to placing a Gamma prior on the precision $\tau=\sigma^{-2}$. 
Thus, we may alternatively specify:
$$
p(\tau)= \gammadist (\tau\mid  \alpha_\tau, \beta_\tau) = \frac{\beta_{\tau}^{\alpha_{\tau}}}{\Gamma(\alpha_{\tau})} \tau^{\alpha_{\tau}-1} \exp({-\beta_{\tau}\cdot \tau}),
$$
with shape $\alpha_{\tau}>0$ and rate $\beta_{\tau}>0$ (Definition~\ref{definition:gamma-distribution}).

\subsection{Gibbs Sampler}
To construct a Gibbs sampler, we must derive the full conditional posterior distribution for each latent or model parameter.

\paragrapharrow{Latent variables.}
The conditional density for latent variables $h_{mn}$ is 
\begin{equation}
\begin{aligned}
p(h_{mn}\mid  a_{mn}, \bw_m, \bz_n, \tau) 
\propto p(a_{mn}\mid h_{mn})\, p(h_{mn}\mid \bw_m,\bz_n, \tau).
\end{aligned}
\end{equation}
To sample from this conditional density, we introduce back the hidden variable $f_{mn}$.
For brevity, we omit the subscript $m,n$. 
The density $f, h\mid  a, \bw,\bz,\tau$ then can be sampled from in two steps, $f\mid a, \bw, \bz, \tau$ and $h\mid f, \bw, \bz, \tau$.
The joint marginal distribution of $a, f$, and $h$, given $m=\bw^\top\bz$ and $\tau$, is 
\begin{equation}
p(a\mid f)\,p(f\mid h)\,p(h\mid m, \tau) = 
\left[ u(f-b_a) - u(f-b_{a+1})\right]\, \normal(f\mid h,1)\, \normal(h\mid m, \tau^{-1}).
\end{equation}
The conditional density of $p(f\mid a, m, \tau )$ follows from
$$
p(f\mid a, m, \tau )=\generaltruncatednormal(f\mid  m, 1+\tau^{-1}, b_a, b_{a+1}), 
$$
a general-truncated-normal density (Definition~\ref{definition:general_truncated_normal}). Therefore, the sample $h$ can be obtained by
\begin{equation}\label{equation:oggw_pos_hmn}
\begin{aligned}
p(h\mid f, m, \tau) &\propto p(f\mid h)\, p(h\mid m,\tau^{-1}) = \normal(f\mid h,1) \, \normal(h\mid m,\tau^{-1})\\
&\propto \normal\left(h \,\bigg|\, \frac{f+m\tau}{1+\tau}, (1+\tau)^{-1}\right).
\end{aligned}
\end{equation}

\paragrapharrow{Multivariate Gaussian parameters.}
Same as the GGGW model, from the discussion in Section~\ref{sec:niw_posterior_conjugacy}, the posterior density of $\{\bmu_w, \bSigma_w\}$ also follows a NIW distribution with updated parameters: 
\begin{equation}\label{equation:oggw_mu_cov}
\bmu_w, \bSigma_w\sim \niw(\bmu_w, \bSigma_w\mid \bmm_M, \kappa_M, \nu_M, \bS_M),
\end{equation}
where 
\begin{subequations}
\begin{align}
\bmm_M &= \frac{\kappa_0\bmo + M\overline{\bw}}{\kappa_M} = 
\frac{\kappa_0 }{\kappa_M}\bmo+\frac{M}{\kappa_M}\overline{\bw}  \label{equation:niw_posterior_equation_2_oggw} \qquad\qquad\quad   \\
\kappa_M  &= \kappa_0 + M  \label{equation:niw_posterior_equation_3_oggw}\\
\nu_M     &=\nu_0 + M  \label{equation:niw_posterior_equation_4_oggw}\\
\bS_M     
&=\bS_0 + \bS_{\overline{w}} + \frac{\kappa_0 M}{\kappa_0 + M}(\overline{\bw} - \bmo)(\overline{\bw} - \bmo)^\top \label{equation:niw_posterior_equation_5_oggw}\\
&=\bS_0 + \sum_{m=1}^M \bw_m \bw_m^\top + \kappa_0 \bmo \bmo^\top - \kappa_M \bmm_M \bmm_M^\top  \label{equation:niw_posterior_equation_6_oggw}\\
\overline{\bw} &= \frac{1}{M}\sum_{m=1}^{M} \bw_m. \label{equation:niw_posterior_equation_7_oggw}\\
\bS_{\overline{w}} &= \sum_{m=1}^{M} (\bw_m - \overline{\bw})(\bw_m - \overline{\bw})^\top \label{equation:niw_posterior_equation_8_oggw}
\end{align}
\end{subequations}

\paragrapharrow{Gaussian variance parameter.}

The conditional density of $\sigma^2$ depends on its parents ($\alpha_\sigma$, $\beta_\sigma$), children ($\bA$), and co-parents ($\bW$, $\bZ$). And it is an inverse-Gamma distribution (by conjugacy in Equation~\eqref{equation:inverse_gamma_conjugacy_general}) with updated parameters:
\begin{equation}\label{equation:oggw_posterior_sigma2}
\begin{aligned}
&p(\sigma^2 \mid  {\bW}, {\bZ}, \bA)=
p(\sigma^2 \mid  \bW,\bZ, \bA) = \inversegammadist (\sigma^2\mid  \widetilde{\alpha_{\sigma}}, \widetilde{\beta_{\sigma}}), \gap\gap\qquad \\
& \widetilde{\alpha_{\sigma}} = \frac{MN}{2} +{\alpha_\sigma}, 
\qquad
\widetilde{\beta_{\sigma}}  =  \frac{1}{2} \sum_{m,n=1}^{M,N} (\bA-\bW\bZ)_{mn}^2 + {\beta_\sigma}.
\end{aligned}
\end{equation}
\paragrapharrow{Gaussian precision parameter.}
Alternatively, the conditional posterior density of $\tau={1}/{\sigma^2}$ is obtained similarly (Equation~\eqref{equation:gamma_conjugacy_general}) by
\begin{equation}\label{equation:oggw_posterior_tau2}
\begin{aligned}
&p(\tau \mid  {\bW}, {\bZ},  \bA)=
p(\tau \mid  \bW,\bZ, \bA) = \gammadist (\tau\mid  \widetilde{\alpha_\tau}, \widetilde{\beta_\tau}), \gap\gap\qquad \\
&\widetilde{\alpha_\tau} = \frac{MN}{2} +{\alpha_\tau}, 
\qquad 
\widetilde{\beta_\tau}  =  \frac{1}{2} \sum_{m,n=1}^{M,N} (\bA-\bW\bZ)_{mn}^2 + {\beta_\tau}.
\end{aligned}
\end{equation}
In practice, the prior hyper-parameters are often set consistently across parameterizations, e.g. $\alpha_\tau=\alpha_\sigma$ and $\beta_\tau=\beta_\sigma$.

\paragrapharrow{Gibbs sampling.}
We can now construct a Gibbs sampler for the OGGW model, as summarized in Algorithm~\ref{alg:oggw_gibbs_sampler}. The algorithm follows the general Gibbs sampling framework introduced in Section~\ref{section:gibbs-sampler}.
In practice, the choice of hyper-parameters for the normal-inverse-Wishart prior has little impact when sufficient data are available, as the likelihood dominates weakly informative priors. A common uninformative setting is:
$\bmm_0=\bzero, \kappa_0=1, \nu_0=K+1, \bS_0=\bI$.
While the choice for $\alpha_\tau$ and $\beta_\tau$ rather depends on the datasets. A week prior choice is $\alpha_\tau=\beta_\tau=1$.

\begin{algorithm}[h] 
\caption{Gibbs sampler for OGGW model in one iteration (prior on $\tau=\frac{1}{\sigma^2}$).
By default, uninformative hyper-parameters are $\bmm_0=\bzero, \kappa_0=1, \nu_0=K+1, \bS_0=\bI$, $\alpha_\tau=\beta_\tau=1$.} 
\label{alg:oggw_gibbs_sampler}  
\begin{algorithmic}[1] 
\Require Choose initial $\alpha_\tau, \beta_\tau, \bmm_0, \kappa_0, \nu_0, \bS_0$;
\For{$m=1$ to $M$}
\State Sample $\bw_{m}$ from $p(\bw_m \mid  \bmu_m, \bSigma_m)$; 
\Comment{Equation~\eqref{equation:oggw_prior_wm_zn}}
\State Sample $h_{mn}$ from $p(h_{mn}\mid a_{mn}, \bw_m,\bz_n,\tau)$ for each $n$;
\Comment{Equation~\eqref{equation:oggw_pos_hmn}}
\EndFor
\For{$n=1$ to $N$}
\State Sample $\bz_{n}$ from $p(\bz_n \mid  \bmu_z, \bSigma_z)$; 
\Comment{Equation~\eqref{equation:oggw_prior_wm_zn}}
\State Sample $h_{mn}$ from $p(h_{mn}\mid a_{mn}, \bw_m,\bz_n,\tau)$ for each $m$;
\Comment{Equation~\eqref{equation:oggw_pos_hmn}}
\EndFor
\State Sample $\tau$ from $p(\tau \mid  \bW,\bZ, \bA)$; 
\Comment{Equation~\eqref{equation:oggw_posterior_tau2}}
\State Sample $\bmu_w, \bSigma_w$ from $p(\bmu_w, \bSigma_w\mid  \bW, M)$;
\Comment{Equation~\eqref{equation:oggw_mu_cov}}
\State Sample $\bmu_z, \bSigma_z$ from $p(\bmu_z, \bSigma_z\mid  \bZ, N)$;
\Comment{Symmetry of Eq.~\eqref{equation:oggw_mu_cov}}
\end{algorithmic} 
\end{algorithm}

\subsection{Properties of OGGW}
A key advantage of the OGGW model is that it does not merely predict a point estimate (e.g., expected rating) for missing entries in $\bA$. Instead, it provides a full predictive distribution over the discrete ordinal categories. While this richer output does not directly improve \textit{root mean squared error} (RMSE, which depends only on point predictions), it can enhance other metrics---such as \textit{mean absolute error (MAE)}---that benefit from calibrated probabilistic forecasts.

Given the hidden variables $\{h_{mn}\}$ and using the likelihood in Equation~\eqref{equation:oggw_averaging_f}, the expected value of the category value for $(m,n)$-th entry is 
$$
\begin{aligned}
\sum_{a=1}^{A} a \cdot p(a\mid h_{mn}) &=  \sum_{a=1}^{A} a \cdot \big(\Phi(h_{mn}-b_a)-\Phi(h_{mn}-b_{a+1})\big)\\
&= \sum_{a=1}^{A} \Phi(h_{mn} - b_a) - A\Phi(h_{mn}-b_{A+1}).
\end{aligned}
$$
Following the likelihood in Equation~\eqref{equation:oggw_likelihood} and integrating out $h_{mn}$, we have 
\begin{equation}\label{equation:oggw_rec_score1}
\ry_{mn} \triangleq  \sum_{a=1}^{A} a\cdot p(a \mid \bw_m,\bz_n, \tau) = \sum_{a=1}^{A} \Phi\left( \frac{\bw_m^\top\bz_n - b_a}{\sqrt{1+\tau^{-1}}}\right).
\end{equation}
Therefore, instead of using the score in Equation~\eqref{equation:recom_poisson1}, the score $\Exp[\ry_{mn} \mid \bA]$ can be obtained by averaging the values of Equation~\eqref{equation:oggw_rec_score1} during the Gibbs sampling process.

Similar to the third recommendation system introduced in Section~\ref{section:recom_poisson}, the OGGW model can also provide   uncertainty about each entry in $\bA$.
Adopting again the idea of the Sharpe ratio, we can suggest the unconsumed movie  $m$ (in the Netflix context) when $a_{mn}$ for user $n$ by the \textit{uncertainty-adjusted recommendation score}:
$$
\text{score}_{mn} = \frac{\Exp[\ry_{mn}  \mid \bA]}{\sqrt{\Var[\ry_{mn}  \mid \bA]}}.
$$

\begin{problemset}
\item Following the derivation in Equation~\eqref{equation:paa_wmk}, derive the conditional distribution over the user feature $z_{kn}$ (for all $k\in\{1,2,\ldots, K\}$ and $n\in\{1,2,\ldots, N\}$) under the PAA model.

\item Using the ``MovieLens 100K" dataset  introduced in Section~\ref{section:movie_rec_als}, evaluate and compare the performance of the PAA and PAAA models (presented in this chapter) against Bayesian real-valued or nonnegative matrix factorization methods.

\item Using the ``MovieLens 100K" dataset  introduced in Section~\ref{section:movie_rec_als}, evaluate and compare the performance of the OGGW model (introduced in this chapter) with Bayesian real-valued or nonnegative matrix factorization approaches.

\end{problemset}

\chapter{Bayesian Interpolative Decomposition}\label{chapter:bayes_id}
\begingroup
\hypersetup{
	linkcolor=structurecolor,
	linktoc=page,  
}
\minitoc \newpage
\endgroup

\section{Interpolative Decomposition (ID)}
\lettrine{\color{caligraphcolor}L}
Low-rank real-valued or nonnegative matrix factorization plays a fundamental role in modern data science.
Low-rank matrix approximation with respect to the Frobenius norm---i.e., minimizing the sum of squared differences from the target matrix---can be efficiently solved using singular value decomposition (SVD) or Bayesian real-valued/nonnegative matrix factorization methods.
However, for many applications, it is advantageous to work with a basis composed of a subset of columns directly drawn from the observed matrix itself \citep{halko2011finding, martinsson2011randomized}.
The \textit{interpolative decomposition (ID)} is one such approach that stands out by explicitly reusing actual columns from the original matrix. This property allows ID to preserve structural features such as sparsity and nonnegativity, which can significantly reduce memory requirements.

ID is widely used as a feature selection tool: it extracts the essential information from large datasets that might otherwise be too big to fit into RAM.
Moreover, it enables the removal of irrelevant components---such as noise and redundant information---through these decomposition techniques \citep{liberty2007randomized, halko2011finding, martinsson2011randomized, ari2012probabilistic, lu2022bayesian, lu2022feature}.
Identifying the indices of the spanning columns is often valuable for data interpretation and analysis. In particular, selecting a small subset of columns that captures the full informational content of the matrix can greatly simplify downstream tasks.
When the columns of the observed matrix carry specific semantic meaning---for example, representing individual transactions in a transactional dataset---the corresponding columns in the ID retain that same interpretability.

The column ID~\footnote{The term \textit{column ID} will henceforth be referred to simply as \textit{ID}, without further qualification.} factors a matrix into the product of two matrices: one consisting of selected columns from the original matrix, and the other containing an identity submatrix (possibly after column permutation), with all entries bounded in magnitude by 1.
We first state and prove the existence of the \textit{exact ID} in the following theorem, and later describe the \textit{low-rank ID} using Bayesian approaches.

\begin{theoremHigh}[Column Interpolative Decomposition]\label{theorem:interpolative-decomposition}
Any rank-$R$ matrix $\bA \in \real^{M\times N}$ can be factored as 
$$
\underset{M \times N}{\bA} = \underset{M\times R}{\bC} \gapthree  \underset{R\times N}{\bW},
$$
where $\bC\in \real^{M\times R}$ comprises $R$ linearly independent columns of $\bA$, and $\bW\in \real^{R\times N}$ is the reconstruction matrix. 
The factor $\bW$ contains an $R\times R$ identity submatrix (after a suitable column permutation), and all its entries satisfy
$$
\max \abs{w_{ij}}\leq 1, \,\, \forall \,\, i\in [1,R], j\in [1,N].
$$
The storage cost of this decomposition reduces from $MN$ floating-point numbers (for $\bA$) to $MR$ and $(N-R)R$  floats for storing $\bC$ and $\bW$, respectively, plus an additional $R$ integers to record the column indices of $\bC$ within $\bA$.
\end{theoremHigh}

\begin{figure}[htp]
\centering
\includegraphics[width=0.7\textwidth]{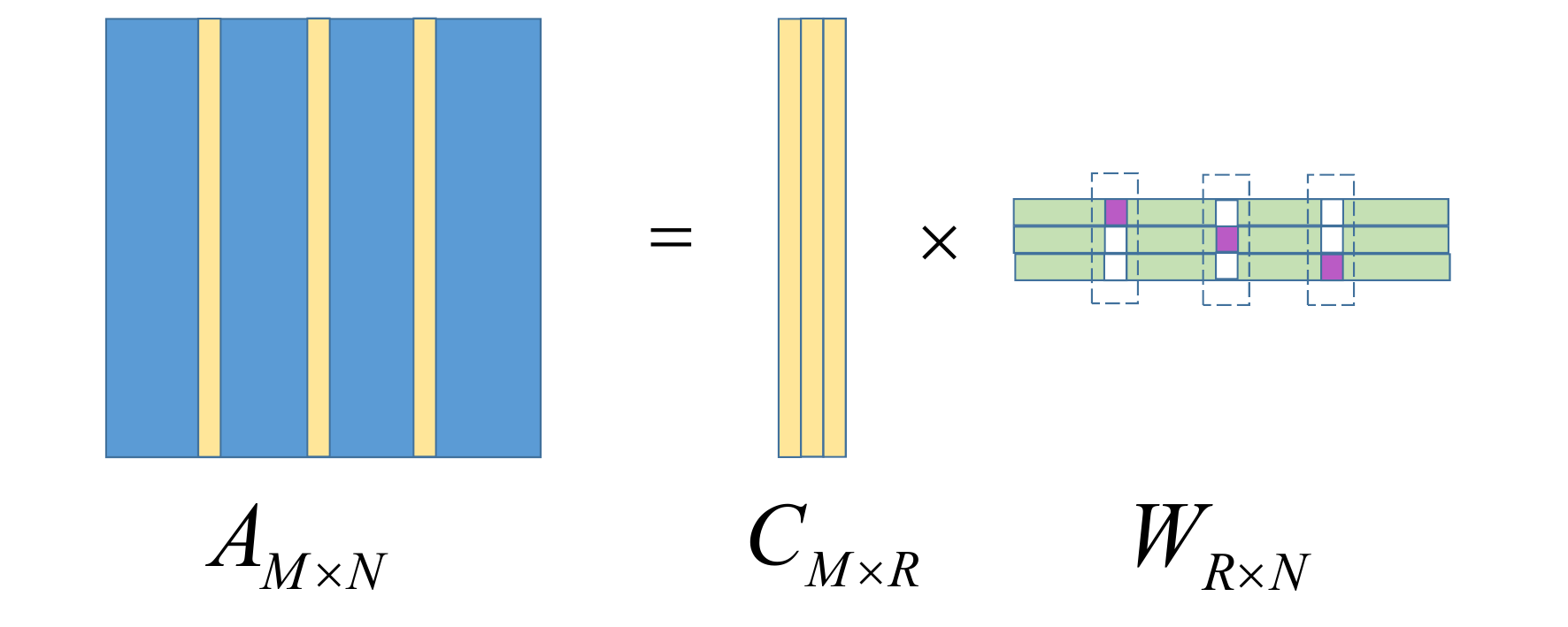}
\caption{
Illustration of the column ID of a matrix. The \textcolor{mydarkyellow}{yellow} columns denote the linearly independent (skeleton) columns of $\bA$; white entries are zero, and \textcolor{mydarkpurple}{purple} entries represent ones forming the identity submatrix in $\bW$.}
\label{fig:column-id}
\end{figure}

While we assert that all entries of $\bW$ have magnitude at most 1, some constructions guarantee only a weaker bound (e.g., $\abs{w_{ij}}\leq 2$).
Figure~\ref{fig:column-id} illustrates the column ID: the \textcolor{mydarkyellow}{yellow} columns are the selected skeleton columns of $\bA$, and the \textcolor{mydarkpurple}{purple} entries in $\bW$ form an  $R\times R$ identity submatrix.
Critically, the positions of the identity columns in $\bW$ correspond exactly to the positions of the selected skeleton columns in $\bA$.
The column ID closely resembles the \textit{CR decomposition}: both select $R$ linearly independent columns into the first factor, and both yield a second factor containing an  $R\times R$ identity submatrix \citep{strang2021every, stranglu, lu2021numerical}.
The key difference is that the CR decomposition specifically chooses the first $R$ linearly independent columns, and its second factor is derived from the \textit{reduced row echelon form (RREF)} of the matrix.
Consequently, the column ID can be applied in the same theoretical contexts as the CR decomposition---for instance, proving that the rank equals the trace for idempotent matrices \citep{lu2021numerical}, or demonstrating the fundamental result that column rank equals row rank \citep{lu2021column}.
Moreover, the column ID is a special case of \textit{rank decomposition} (see Problem~\ref{problem:rank_decom}) and, like most such decompositions, is generally not unique \citep{lu2021numerical}.

\paragrapharrow{Notations that will be extensively used in the sequel.} Following Matlab-style indexing, let $\sJ$ be an index vector of length $R$ indicating the columns of $\bA$ selected into $\bC$. 
Then we write $\bC=\bA[:,\sJ]$ (see Definition~\ref{definition:matlabnotation}).
These are the ``skeleton" columns of $\bA$. 
From the skeleton index set $\sJ$, the $R\times R$ identity submatrix in  $\bW$ is recovered as
$$
\bW[:,\sJ] = \bI_R \in \real^{R\times R}.
$$ 
Let $\sI$ denote the complementary index set of the remaining columns, so that
$$
\sJ\cap \sI=\varnothing \qquad \text{and}\qquad \sJ\cup \sI = \{1,2,\ldots, N\}.
$$
The remaining $N-R$ columns of  $\bW$  form an $R\times (N-R)$ \textit{expansion matrix}, whose entries are the expansion coefficients used to reconstruct the non-skeleton columns of $\bA$ from $\bC$:
$$
\bE = \bW[:,\sI] \in \real^{R\times (N-R)}.
$$
Finally, let  $\bP\in \real^{N\times N}$ be the column permutation matrix  (Definition~\ref{definition:permutation-matrix}) defined by $\bP=\bI_N[:,(\sJ, \sI)]$. 
Then
$$
\bA\bP = \bA[:,(\sJ, \sI)] = \left[\bC, \bA[:,\sI]\right],
$$
which implies
\begin{equation}\label{equation:interpolatibve-w-ep}
\bW\bP = \bW[:,(\sJ, \sI)] =\left[\bI_R, \bE \right] 
\qquad\implies\qquad 
\bW = \left[\bI_R, \bE \right] \bP^\top.
\end{equation}
\section{Existence of the Column Interpolative Decomposition}\label{section:proof-column-id}
\paragrapharrow{Cramer's rule.}
The proof of the existence of the column ID relies on \textit{Cramer's rule} (see Problems~\ref{prob:cramer_adj_1}--\ref{prob:cramer_adj_5}), which we briefly review below.
Cramer's rule provides an explicit formula for solving a system of linear equations with as many equations as unknowns, provided the system has a unique solution---that is, when the coefficient matrix is nonsingular.
Consider a system of $N$ linear equations in $N$ unknowns, written in matrix form as:
$$
\bG \bx = \bl,
$$
where $\bG\in \real^{N\times N}$ is nonsingular, and $\bx,\bl \in \real^N$. 
Then the system has a unique solution, with each component given by
$$
x_n = \frac{\det(\bG_n)}{\det(\bG)}, \qquad \text{for all}\gapthree n\in \{1,2,\ldots,N\},
$$
where $\bG_n$ denotes the matrix obtained by replacing the $n$-th column of $\bG$ with the vector $\bl$.
More generally, Cramer's rule applies to the matrix equation
\begin{equation}
\bG\bX = \bL,
\end{equation}
where $\bG\in \real^{N\times N}$ is nonsingular, and $\bX,\bL\in \real^{N\times M}$. 
Let $\sI=[i_1, i_2, \ldots, i_K]$ and $\sJ=[j_1,j_2,\ldots, j_K]$ be two index sets of cardinality $K<\min\{M,N\}$, with  $1\leq i_1\leq i_2\leq \ldots\leq i_K\leq N$ and $1\leq j_1\leq j_2\leq \ldots\leq j_K\leq M$. Then $\bX[\sI,\sJ]$ is the $K\times K$ submatrix of $\bX$ formed by rows $\sI$ and columns $\sJ$. 
Define $\bG_{\bL}(\sI,\sJ)$ as the $N\times N$ matrix obtained by replacing the $(i_k)$-th column of $\bG$ with the $(j_k)$-th column of $\bL$, for all $k\in \{1,2,\ldots,K\}$. Then we have 
\begin{equation}
\det(\bX[\sI,\sJ]) = \frac{\det\left(\bG_{\bL}(\sI,\sJ)\right)}{\det(\bG)}.
\end{equation}
In the special case where  $\abs{\sI}=\abs{\sJ}=1$, this reduces to
\begin{equation}\label{equation:cramer-rule-general}
x_{nm} = \frac{\det\left(\bG_{\bL}(n,m)\right)}{\det(\bG)}, 
\quad \forall \, n,m.
\end{equation}
We now use this result to prove the existence of the column ID.

\begin{proof}[of Theorem~\ref{theorem:interpolative-decomposition}]
As noted, the proof hinges on Cramer's rule. If we can express the entries of $\bW$ via \eqref{equation:cramer-rule-general} and show that the absolute value of each numerator does not exceed that of the denominator, then $\abs{w_{ij}}\leq 1$ follows immediately. However, Cramer's rule requires a square, nonsingular denominator matrix---so we must first reduce the general case to one where this holds.

\paragraph{Step 1: Column ID for full row rank matrices.}
We begin with the simpler case where $\bA$ has full row rank $R$ ($R\leq N$).
In this setting, the desired column ID takes the form $\bA=\bC\bW$, where $\bC\in\real^{R\times R}$ is a square, invertible submatrix consisting of $R$ selected columns of $\bA$.
Choose the ``skeleton" index set $\sJ$ by maximizing the absolute determinant:
\begin{equation}\label{equation:interpolative-choose-js}
\sJ = \mathop{\arg\max}_{\sJ_t} \big\{\abs{\det(\bA[:,\sJ_t])}: \text{$\sJ_t$ is a subset of $\{1,2,\ldots, N\}$ with size $R$} \big\}.
\end{equation}
Let $\sI$ denote the complementary index set, so that $\sJ\cup \sI = \{1,2,\ldots, N\}$ and  $\sJ\cup \sI=\varnothing$.
There exists a column permutation matrix $\bP$ such that
$$
\bA\bP = 
\begin{bmatrix}
\bA[:,\sJ]&\bA[:,\sI]
\end{bmatrix}.
$$
Since $\bC=\bA[:,\sJ]$ is nonsingular by construction, we can write
$$
\begin{aligned}
\bA
&=\big[\bA[:,\sJ],\, \bA[:,\sI]\big]\bP^\top
= 
\bA[:,\sJ]
\big[
\bI_R,\, \bA[:,\sJ]^{-1}\bA[:,\sI]
\big]
\bP^\top
= \bC 
\underbrace{
	\big[\bI_R , \, \bC^{-1}\bA[:,\sI]\big]
	\bP^\top}_{\bW},
\end{aligned}
$$
where the matrix $\bW$ is given by 
$
\bW=
\big[\bI_R ,\,\bC^{-1}\bA[:,\sI]\big]
\bP^\top
=
\big[\bI_R,\, \bE\big]
\bP^\top
$ 
from Equation~\eqref{equation:interpolatibve-w-ep}. To prove the claim that the magnitude of $\bW$ is no larger than 1 is equivalent to proving that entries in $\bE=\bC^{-1}\bA[:,\sI]\in \real^{R\times (N-R)}$ are no greater than 1 in absolute value.

Let  $[j_1,j_2,\ldots, j_N]$ be the permuted column indices of $[1,2,\ldots, N]$ such that 
$$
[j_1,j_2,\ldots, j_N] = [1,2,\ldots, N] \bP = [\sJ, \sI].
$$
Thus, it follows from $\bC\bE=\bA[:,\sI]$ that 
$$
\begin{aligned}
\underbrace{	[\ba_{j_1}, \ba_{j_2}, \ldots, \ba_{j_R}]}_{=\bC=\bA[:,\sJ]} \bE &=
\underbrace{[\ba_{j_{R+1}}, \ba_{j_{R+2}}, \ldots, \ba_{j_N}]}_{=\bA[:,\sI]\triangleq \bB},
\end{aligned}
$$ 
where $\ba_i$ denotes the $i$-th column of $\bA$, and we let $\bB\triangleq \bA[:,\sI]$.
Therefore, by Cramer's rule in Equation~\eqref{equation:cramer-rule-general}, we have 
\begin{equation}\label{equation:column-id-expansionmatrix}
e_{kl} = 
\frac{\det\left(\bC_{\bB}(k,l)\right)}
{\det\left(\bC\right)},
\end{equation}
where $e_{kl}$ is the entry ($k,l$) of $\bE$, and $\bC_{\bB}(k,l)$ is the $R\times R$ matrix formed by replacing the $k$-th column of $\bC$ with the $l$-th column of $\bB$. For example, 
$$
\begin{aligned}
e_{11} &= 
\frac{\det\left([\textcolor{mylightbluetext}{\ba_{j_{R+1}}}, \ba_{j_2}, \ldots, \ba_{j_R}]\right)}
{\det\left([\ba_{j_1}, \ba_{j_2}, \ldots, \ba_{j_R}]\right)},
\qquad 
&e_{12} &=
\frac{\det\left([\textcolor{mylightbluetext}{\ba_{j_{R+2}}}, \ba_{j_2},\ldots, \ba_{j_R}]\right)}
{\det\left([\ba_{j_1}, \ba_{j_2}, \ldots, \ba_{j_R}]\right)},\\
e_{21} &= 
\frac{\det\left([\ba_{j_1},\textcolor{mylightbluetext}{\ba_{j_{R+1}}}, \ldots, \ba_{j_R}]\right)}
{\det\left([\ba_{j_1}, \ba_{j_2}, \ldots, \ba_{j_R}]\right)},
\qquad 
&e_{22} &= 
\frac{\det\left([\ba_{j_1},\textcolor{mylightbluetext}{\ba_{j_{R+2}}}, \ldots, \ba_{j_R}]\right)}
{\det\left([\ba_{j_1}, \ba_{j_2}, \ldots, \ba_{j_R}]\right)}.
\end{aligned}
$$
Because $\sJ$ was chosen to maximize $\det(\bC)$ in Equation~\eqref{equation:interpolative-choose-js}, any such replacement in the numerator cannot increase the absolute determinant. Hence,
$$
\abs{e_{kl}}\leq 1, \qquad \text{for all}\gapthree k\in \{1,2,\ldots, R\}, \,\,\, l\in \{1,2,\ldots, N-R\}.
$$
	
\paragraph{Step 2: Extension to general matrices.}
Summarizing the above (and slightly abusing notation): for any matrix $\bF\in \real^{R\times N}$ of full row rank $R\leq N$, a column ID $\bF=\bC_0\bW$ exists with $\abs{w_{ij}}\leq 1$.

Now consider a general matrix $\bA\in \real^{M\times N}$ of rank $R\leq \{M,N\}$. 
It admits a \textit{rank decomposition} (see Problem~\ref{problem:rank_decom}):
$$
\underset{M\times N}{\bA} = \underset{M\times R}{\bD}\gapthree \underset{R\times N}{\bF},
$$
where $\bD$ and $\bF$ have full column rank $R$ and full row rank $R$, respectively \citep{lu2021numerical}. Apply the column ID to $\bF$: $\bF=\bC_0\bW$, where $\bC_0=\bF[:,\sJ]$ consists of $R$ linearly independent columns of $\bF$.
Then
$$
\bA[:,\sJ]=\bD\bF[:,\sJ],
$$
i.e., the columns indexed by $\sJ$ of $(\bD\bF)$ can be obtained by $\bD\bF[:,\sJ]$, which in turn are the columns of $\bA$ indexed by $\sJ$. 
Define $\bC\triangleq \bA[:,\sJ]$. 
It follows that
$$
\underbrace{\bA[:,\sJ]}_{\bC}= \underbrace{\bD\bF[:,\sJ]}_{\bD\bC_0}
$$
and 
$$
\bA = \bD\bF =\bD\bC_0\bW = \underbrace{\bD\bF[:,\sJ]}_{\bC}\bW=\bC\bW.
$$
which is the desired column ID of $\bA$. This completes the proof.
\end{proof}

The above proof suggests an intuitive algorithm for computing the optimal column ID, as shown in Algorithm~\ref{alg:column-id-intuitive}. However, any method that guarantees selection of the maximally conditioned subset of columns necessarily incurs combinatorial complexity \citep{martinsson2019randomized}. In subsequent sections, we will explore practical alternatives that yield well-conditioned (though not necessarily optimal) ID factorizations.

\begin{algorithm}[h] 
\caption{An \textcolor{mylightbluetext}{Intuitive} Method to Compute the Column ID} 
\label{alg:column-id-intuitive} 
\begin{algorithmic}[1] 
\Require 
Rank-$R$ matrix $\bA$ with size $M\times N $; 
\State Compute the rank decomposition $\underset{M\times N}{\bA} = \underset{M\times R}{\bD}\gapthree \underset{R\times N}{\bF}$, e.g., via UTV decomposition \citep{lu2021numerical};
\State Compute column ID of $\bF$: $\bF=\bF[:,\sJ]\bW = \widetildebC\bW$:
$$
\begin{aligned}
2.1. \,\,\,&\left\{
\begin{aligned}
\sJ &= \mathop{\arg\max}_{\sJ} \left\{\abs{\det(\bF[:,\sJ])}: \text{$\sJ$ is a subset of $\{1,2,\ldots, N\}$ with size $R$} \right\};&\\
\sI &= \{1,2,\ldots, N\} \setminus \sJ;&\\
\end{aligned}
\right.\\
2.2.\,\,\,&\left\{
\begin{aligned}
\widetildebC &= \bF[:,\sJ]; \\
\bM &= \bF[:,\sI];
\end{aligned}
\right.\\
2.3.\,\,\, &\bF\bP = \bF[:,(\sJ,\sI)] \text{ to obtain permutation matrix  $\bP$};\\
2.4.\,\,\, &e_{kl} = 
\frac{\det\left(\widetildebC_{\bM}(k,l)\right)}
{\det\left(\widetildebC\right)},  \qquad \text{for all}\gap k\in [1, R], l\in [1,N-R] \text{  ~(Equation~\eqref{equation:column-id-expansionmatrix})};\\
2.5.\,\,\, &\bW= [\bI_R, \bE]\bP^\top \text{ ~(Equation~\eqref{equation:interpolatibve-w-ep})}.
\end{aligned}
$$
\State 	$\bC=\bA[:,\sJ]$;
\State Output the column ID $\bA=\bC\bW$;
\end{algorithmic} 
\end{algorithm}

\begin{example}[Compute the Column ID]\label{example:column-id-a}
Consider the  matrix
$$
\bA=
\begin{bmatrix}
56 & 41 & 30\\
32 & 23 & 18\\
80 & 59 & 42
\end{bmatrix},
$$
with rank $R=2$.
The trivial process for computing the column ID of $\bA$ is shown as follows.
A rank decomposition is
$$
\bA = \bD\bF=
\begin{bmatrix}
1 & 0 \\
0 & 1 \\
2 &-1
\end{bmatrix}
\begin{bmatrix}
56 & 41 & 30 \\
32 & 23 & 18 
\end{bmatrix}.
$$
Since rank $R=2$, $\sJ$ is one of $[2,3], [1,3], [1,2]$, where the absolute determinants of $\bF[:,\sJ]$ are $48, 48, 24$, respectively. 
We may choose either $\sJ=\{2,3\}$ or $\sJ=\{1,3\}$.
We proceed by choosing $\sJ=[1,3]$:
$$
\begin{aligned}
\widetildebC &= \bF[:,\sJ]=
\begin{bmatrix}
56 & 30 \\
32 & 18 
\end{bmatrix},\qquad 
\bM = \bF[:,\sI]=\begin{bmatrix}
41 \\
23
\end{bmatrix}.
\end{aligned}
$$
And 
$$
\bF\bP = \bF[:,\{\sJ,\sI\}] = \bF[:,\{1,3,2\}]
\qquad\implies\qquad
\bP = 
\begin{bmatrix}
1 &  & \\
& &1\\
& 1 & 
\end{bmatrix}.
$$
In this example, $\bE\in \real^{2\times 1}$:
$$
\begin{aligned}
e_{11} &=
\det\left(
\begin{bmatrix}
41 & 30 \\
23 & 18
\end{bmatrix}\right)\bigg/
\det\left(
\begin{bmatrix}
56 & 30 \\
32 & 18
\end{bmatrix}\right)=1;\\
e_{21} &=
\det\left(
\begin{bmatrix}
56 & 41 \\
32 & 23
\end{bmatrix}\right)\bigg/
\det\left(
\begin{bmatrix}
56 & 30 \\
32 & 18
\end{bmatrix}\right)=-\frac{1}{2}.
\end{aligned}
$$
This makes 
$$
\bE = 
\begin{bmatrix}
1\\-\frac{1}{2}
\end{bmatrix}
\qquad\implies\qquad
\bW = [\bI_2, \bE]\bP^\top =
\begin{bmatrix}
1 & 1 & 0\\
0 & -\frac{1}{2} & 1
\end{bmatrix}.
$$
The selected skeleton columns are
$$
\bC = \bA[:,\sJ] = 
\begin{bmatrix}
56 & 30\\
32 & 18\\
80 & 42
\end{bmatrix}
\qquad\implies\qquad 
\bA=\bC\bW =
\begin{bmatrix}
56 & 30\\
32 & 18\\
80 & 42
\end{bmatrix}
\begin{bmatrix}
1 & 1 & 0\\
0 & -\frac{1}{2} & 1
\end{bmatrix},
$$
with all entries of  $\bW$ satisfying $\abs{w_{ij}}\leq 1$, as required.
\end{example}

We conclude this section by discussing the non-uniqueness of the column ID.
\begin{remark}[Non-uniqueness of the Column ID]
In the above specific Example~\ref{example:column-id-a}, we notice that both $\bF[:,\{2,3\}]$ and $\bF[:,\{1,3\}]$ achieve the maximal absolute determinant (48). 
Either choice yields a valid column ID. Moreover, once a set $\sJ$ is selected, any permutation of its indices (e.g., $\sJ=[1,3]$ vs. $[3,1]$) also produces a valid decomposition, since the identity submatrix in $\bW$ can be correspondingly permuted. These degrees of freedom---the choice among equally optimal column subsets and the ordering within a chosen subset---explain why the column ID is generally not unique.
\end{remark}

\section{Skeleton/CUR Decomposition, Row ID, and Two-Sided ID}
To delve deeper into interpolative decomposition, we first introduce a closely related factorization known as the \textit{skeleton decomposition} or \textit{CUR decomposition}.

\index{Decomposition: Skeleton}
\index{Decomposition: CUR}
\begin{theoremHigh}[Skeleton Decomposition]\label{theorem:skeleton-decomposition}
Any rank-$R$ matrix $\bA \in \real^{M \times N}$ can be decomposed as 
$$
\underset{M\times N}{\bA }= 
\underset{M\times R}{\bC} \gapthree \underset{R\times R}{\bU^{-1} }\gapthree \underset{R\times N}{\bR},
$$
where $\bC$ consists of  $R$ linearly independent columns of $\bA$, $\bR$ consists of   $R$ linearly independent rows of $\bA$, and $\bU$ is the nonsingular submatrix on the intersection of $\bC$ and $\bR$. 
Regarding storage requirements:
\begin{itemize}
\item Storing the full decomposition explicitly requires  $R(M+N)+R^2$ floating-point numbers (as opposed to $MN$ for the original matrix).
\item Alternatively, if only the positions of the selected rows and columns are recorded, one needs $MR$ floats for $\bC$, $NR$ floats for $\bR$, and $2R$ integers to store the column indices (for $\bC$) and row indices (for $\bR$ ) within $\bA$. The submatrix $\bU$ can then be reconstructed from $\bC$ and $\bR$ using these indices.
\end{itemize}
\end{theoremHigh}
\begin{proof}[of Theorem~\ref{theorem:skeleton-decomposition}]
The key ingredient is the existence of a nonsingular $R\times R$ submatrix $\bU$ within $\bA$.

\paragraph{Existence of such nonsingular matrix $\bU$.} 
Since  $\bA$ has rank $R$, it contains $R$ linearly independent columns.
Let  $\bC=[\ba_{i_1}, \ba_{i_2}, \ldots, \ba_{i_R}] \in \real^{M\times R}$ denote the matrix formed by these columns. 
Because  $\bC$ has full column rank $R$, its row space also has dimension $R$. 
Hence, there exist $R$ linearly independent rows in $\bC$. Selecting these rows yields an $R\times R$ submatrix $\bU$, which is necessarily nonsingular.

\paragraph{Main proof.}
Let $\sJ\subset \{1,2,\ldots,N\}$ and $\sS\subset \{1,2,\ldots,M\}$ be index sets of size $R$ such that $\bC=\bA[:,\sJ]$ and $\bU=\bA[\sS,\sJ]$. Since $\bU$ is invertible, any column $\ba_n$ of $\bA$ can be expressed as a linear combination of the columns of $\bC$: $\ba_n=\bC\bx_n$ for some coefficient vector $\bx_n\in\real^R$.
Now consider the restriction of $\ba_n$  to the rows indexed by $\sS$, denoted $\br_n=\bA[\sS,n]\in\real^R$. Because $\bU=\bC[\sS,:]$, we have $\br_n=\bU\bx_n$, and thus $\bx_n=\bU^{-1} \br_n$. Therefore,
$$
\ba_n = \bC \bU^{-1} \br_n, \quad \forall\,n=1,2,\ldots,N.
$$
Stacking all such columns gives
$
\bA = [\ba_1, \ba_2, \ldots, \ba_N] = \bC \bU^{-1} \bR,
$
where $\bR=\bA[\sS,:]\in\real^{R\times N}$ contains the selected rows of $\bA$. This completes the proof.

In summary: we first select $R$ linearly independent columns to form $\bC$, then identify $R$ linearly independent rows within $\bC$ to define the invertible core $\bU$, and finally use the corresponding rows of $\bA$ (i.e., $\bR$) to reconstruct the entire matrix.
\end{proof}

As a special case, if $\bA$ is square and invertible ($R=M=N$), then choosing all rows and columns yields $\bC=\bR=\bU=\bA$, and the decomposition reduces to the identity$\bA = \bA\bA^{-1}\bA$.

\begin{figure}[htp]
\centering
\includegraphics[width=0.7\textwidth]{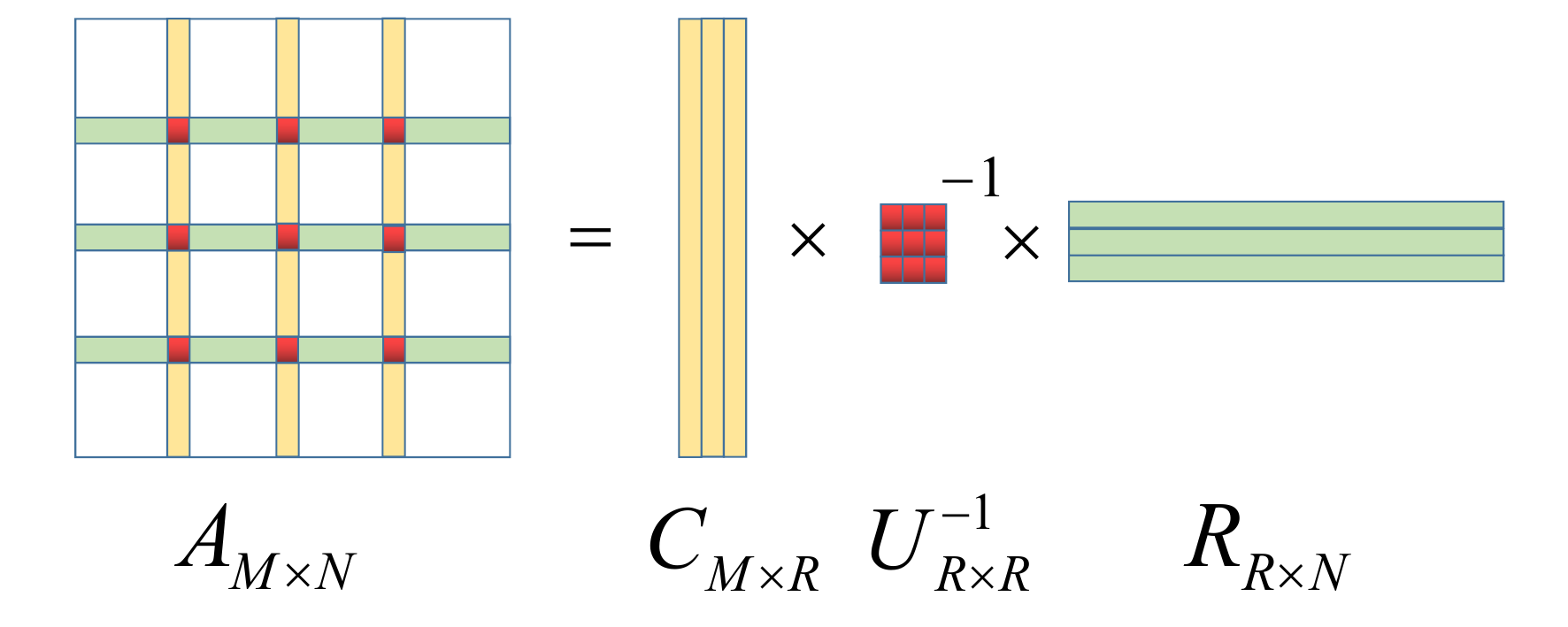}
\caption{
Illustration of the skeleton decomposition. The \textcolor{mydarkyellow}{yellow} columns represent linearly independent columns of $\bA$, and the \textcolor{mydarkgreen}{green} rows represent linearly independent rows. Their intersection forms the nonsingular core $\bU$.
}
\label{fig:skeleton}
\end{figure}

The skeleton decomposition is also commonly called the CUR decomposition, named after the three factors: 
$\bC$ (columns), $\bU$ (core), and $\bR$ (rows).
Compared to the singular value decomposition (SVD), CUR offers better interpretability: it uses actual columns and rows from the original data, whereas SVD relies on abstract singular vectors that may lack physical meaning \citep{mahoney2009cur}. Like the ID, CUR also preserves structural properties such as sparsity when the input matrix is sparse.
Moreover, similar to SVD, CUR serves as a powerful tool for data compression, feature selection, and exploratory data analysis in applications ranging from scientific computing to machine learning \citep{mahoney2009cur, an2012large, lee2008cur+}.
Figure~\ref{fig:skeleton} visualizes the decomposition: the \textcolor{mydarkyellow}{yellow} columns and \textcolor{mydarkgreen}{green} rows correspond to the selected subsets, and their overlap defines $\bU=\bA[\sS,\sJ]$, where $\sS$ and $\sJ$ are the row and column index sets, respectively.

We previously introduced the column ID. This is not an isolated concept---it belongs to a family of related decompositions.
\begin{theoremHigh}[The Whole Interpolative Decomposition]\label{theorem:interpolative-decomposition-row}
Any rank-$R$ matrix $\bA \in \real^{M \times N}$ admits the following factorizations:
$$
\begin{aligned}
\text{Column ID: }&\gap \underset{M \times N}{\bA} &=& \boxed{\underset{M\times R}{\bC}} \gap  \underset{R\times N}{\bW} ; \\
\text{Row ID: } &\gap &=&\underset{M\times R}{\bZ} \gap  \boxed{\underset{R\times N}{\bR}}; \\
\text{Two-Sided ID: } &\gap &=&\underset{M\times R}{\bZ} \gap \boxed{\underset{R\times R}{\bU}}  \gap  \underset{R\times N}{\bW}, \\
\end{aligned}
$$
where
\begin{itemize}
\item $\bC=\bA[:,\sJ]\in \real^{M\times R}$ contains  $R$ linearly independent columns of $\bA$,
and $\bW$ satisfies  $\bW[:,\sJ]=\bI_R$ (after a suitable column permutation). All entries of $\bW$ obey $\abs{w_{ij}}\leq 1$.
\item $\bR=\bA[\sS,:]\in \real^{R\times N}$ contains  $R$ linearly independent rows of $\bA$, and $\bZ$ satisfies $\bZ[\sS,:]=\bI_R$ (after a suitable row permutation). All entries of $\bZ$ obey $\abs{z_{ij}}\leq 1$.
\item $\bU=\bA[\sS,\sJ] \in \real^{R\times R}$ is the nonsingular intersection submatrix of $\bC$ and $\bR$.
\item \textbf{Skeleton decomposition:} the boxed matrices  $\bC,\bR$, and $\bU$ are identical to those in the skeleton decomposition (Theorem~\ref{theorem:skeleton-decomposition}), and indeed satisfy  $\bA=\bC\bU^{-1}\bR$.
\end{itemize}
\end{theoremHigh}
The row ID follows directly from the column ID applied to $\bA^\top$. 
If  $\bA^\top=\bC_0\bW_0$, where $\bC_0$ contains $R$ linearly independent columns of $\bA^\top$ (i.e., $R$ linearly independent rows of $\bA$), 
then transposing gives $\bA=\bW_0^\top\bC_0^\top = \bZ\bR$, with $\bR\triangleq \bC_0^\top$ and $ \bZ\triangleq\bW_0^\top$.

For the two-sided ID, recall from the skeleton decomposition that  $\bA=\bC\bU^{-1}\bR$. 
Setting $\bZ\triangleq \bC\bU^{-1}$ and noting that $\bA=\bZ\bR$, we also have from the column ID that $\bA=\bC\bW=\bZ\bU\bW$. Thus, $\bA=\bZ\bU\bW$, establishing the two-sided form.

\paragrapharrow{Storage requirements.} We summarize the memory footprint of each variant:
\begin{itemize}
\item \textit{Column ID.} It requires $MR$ and $(N-R)R$ floats to store $\bC$ and $\bW$, respectively, and $R$ integers to store the indices of the selected columns in $\bA$;
\item \textit{Row ID.} It requires $NR$ and $(M-R)R$ floats to store $\bR$ and $\bZ$, respectively, and $R$ integers to store the indices of the selected rows in $\bA$;
\item \textit{Two-Sided ID.} It requires $(M-R)R$, $(N-R)R$, and $R^2$ floats to store $\bZ,\bW$, and $\bU$, respectively. And extra $2R$ integers are required to store the indices of the selected rows and columns in $\bA$.
\end{itemize}

\paragrapharrow{Storage reduction for sparse matrices in the two-sided ID.}
Suppose we compute the column ID $\bA=\bC\bW$ with $\bC=\bA[:,\sJ]$, and further identify a set of $R$ ``spanning" rows indexed by $\sS$ such that
$$
\bA[\sS,:] = \bC[\sS,:]\bW.
$$ 
Note that $\bC[\sS,:] = \bA[\sS,\sJ]\in \real^{R\times R}$  is nonsingular (since both $\bC$ and $\bA[\sS,:]$ have rank $R$). Therefore,
$$
\bW = (\bA[\sS,\sJ])^{-1} \bA[\sS,:].
$$
This means $\bW$ need not be stored explicitly. Instead, we can store only: 
the sparse matrix $\bA[\sS,:]$ (which is cheap if $\bA$ is sparse), 
and either the inverse  $(\bA[\sS,\sJ])^{-1}$, or just the index set $\sJ$ (if the inverse can be computed on demand).
In the latter case, only $R$ integers (for $\sJ$) and the sparse row subset $\bA[\sS,:]$ are required---offering significant memory savings for large, sparse datasets.

\index{Low-Rank interpolative decomposition}
\section{Bayesian Low-Rank Interpolative Decomposition}
Instead of seeking an exact interpolative decomposition, we now consider its approximate counterpart.
The \textit{low-rank ID} problem for a given matrix $\bA$ can be formulated as
$$
\bA=\bC\bW+\bE,
$$ 
where $\bA= [\ba_1, \ba_2, \ldots, \ba_N]\in \real^{M\times N}$ is approximately factorized into an $M\times K$ matrix $\bC\in \real^{M\times K}$ containing $K$ basis columns of $\bA$ and a $K\times N$ matrix $\bW\in \real^{K\times N}$ with entries no larger than 1 in magnitude; the noise is captured by matrix $\bE\in \real^{M\times N}$. 
Here, $K<R=\rank(\bA)$, which justifies the term low-rank ID.

Several methods exist for computing low-rank ID approximations. The most widely used is the \textit{randomized ID (RID)} algorithm \citep{liberty2007randomized}. 
At a high level, the algorithm randomly samples $S > K$ columns from $\bA$, uses column-pivoted QR (CPQR) to select $K$ of those $S$ columns for basis matrix $\bC$, and then computes $\bW$ via least squares \citep{lu2021numerical}. 
Typically, the oversampling parameter is set to $S=1.2K$ to ensure that the sampled columns capture a large portion of the range (column space) of $\bA$.

However, a known drawback of randomized ID is that the resulting matrix $\bW$ may contain entries with magnitude greater than 1. While this is often tolerable in practice, it can compromise numerical stability in applications that require strict bounds on coefficient magnitudes. In fact, \citet{advani2021efficient} report cases where entries in $\bW$ exceed 167, significantly degrading stability.
In contrast, probabilistic models can naturally enforce constraints on the range of latent factors through appropriate prior distributions. Motivated by this, we focus on Bayesian ID (BID) for underlying matrices.
Bayesian ID was introduced in \citet{lu2022bayesian, lu2022comparative} and later adapted to feature selection in \citet{lu2022feature}. Training such models amounts to finding the best rank-$K$ approximation to the observed matrix $\bA$ under a specified probabilistic loss.

\begin{figure*}[h]
\centering  
\subfigtopskip=2pt 
\subfigbottomskip=9pt 
\subfigcapskip=-5pt 
\includegraphics[width=0.95\textwidth]{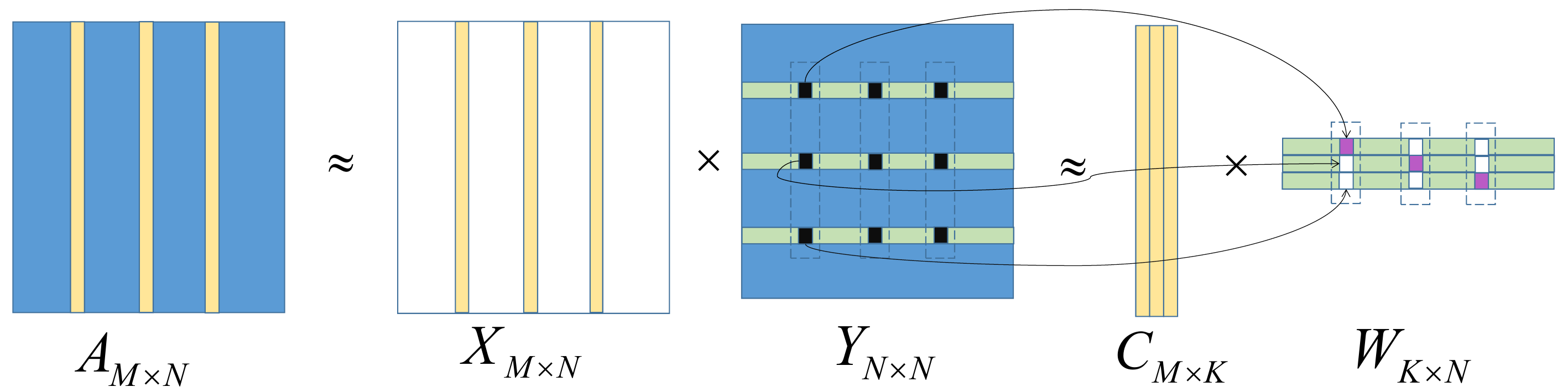}
\caption{Demonstration of the interpolative decomposition of a matrix, where the \textcolor{mydarkyellow}{yellow} vector denotes
the basis columns of matrix $\bA$, white entries denote zero, \textcolor{mydarkpurple}{purple} entries denote
one, \textcolor{mylightbluetext}{blue} and black entries denote elements that are not necessarily zero. The Bayesian ID models find the approximation $\bA\approx\bX\bY$, while the post-processing procedure calculates the approximation $\bA\approx\bC\bW$.}
\label{fig:id-column}
\end{figure*}

\paragrapharrow{Modeling the column selection process.}
Let $\br\in \{0,1\}^N$ be a \textit{state vector}  indicating the role of each column, i.e., basis column or interpolated (remaining) column: if $r_n=1$, then the $n$-th column $\ba_n$ is a basis column; if $r_n=0$, then $\ba_n$ is interpolated from the basis columns (up to noise).
Suppose further $\sJ$ is the set of the indices of the selected basis columns (with size $K$ now), $\sI$ is the set of the indices of the interpolated columns (with size $N-K$) such that 
$$
\begin{aligned}
\sJ\cap \sI = \varnothing, &\qquad \sJ \cup \sI =\{1,2,\ldots, N\}; \\
\sJ=\sJ(\br)=\{n\mid r_n=1\}_{n=1}^N, &\qquad  \sI=\sI(\br)=\{n\mid r_n=0\}_{n=1}^N.
\end{aligned}
$$
The basis matrix is then $\bC=\bA[:,\sJ]$. 
The approximation $\bA\approx \bC\bW$ can be equivalently expressed using two auxiliary matrices $\bX\in \real^{M\times N}$ and $\bY\in \real^{N\times N}$ as:
$$
\underset{M \times N}{\bA}
\approx \underset{M \times K}{\bC} \gapthree\underset{K \times N}{\bW}
=
\underset{M \times N}{\bX} \gapthree \underset{N \times N}{\bY},
$$ 
where 
$$
\begin{aligned}
\bX[:,\sJ]&=\bC\in\real^{M\times K}, \qquad &\bX[:,\sI] &= \bzero\in \real^{M\times (N-K)}; \\
\bY[\sJ,:] &= \bW\in\real^{K\times N}, \qquad &\bY[\sI, :]&=\text{random matrix }\in\real^{(N-K)\times N}.
\end{aligned}
$$ 
Crucially, the structure of $\bW$ enforces an identity submatrix corresponding to the basis columns:
\begin{equation}\label{equation:submatrix_bid_identity}
	\bI_K = \bW[:,\sJ] = \bY[\sJ,\sJ].
\end{equation}
Thus, finding a low-rank ID of  $\bA\approx\bC\bW$ is equivalent to learning $\bX$ and $\bY$  (or, implicitly, the state vector $\br$) such that $\bA\approx\bX\bY$, with the  state vector $\br$ determining which columns form $\bC$ (see Figure~\ref{fig:id-column}).

To evaluate the quality of the approximation, we minimize the reconstruction error, typically measured by the mean squared error (MSE), i.e., the squared Frobenius norm:
\begin{equation}\label{equation:idbid-per-example-loss}
	\mathop{\min}_{\bW,\bZ} \,\, \frac{1}{MN}\sum_{n=1}^N \sum_{m=1}^{M} \left(a_{mn} - \bx_m^\top\by_n\right)^2,
\end{equation}
where $\bx_m$ and $\by_n$ are the $m$-th \textbf{row} of $\bX$ and $n$-th \textbf{column}  of $\bY$, respectively.
Since $\bX$ and $\bY$ are structured via $\br$, the optimization is effectively over $\bW$ and the selection 
 $\br$.

Rather than imposing hard constraints on the entries of $\bW$  (or $\bY$), we adopt a Bayesian approach. We treat the ID as a latent factor model and place a prior distribution on the latent variables that naturally restricts their magnitude.
Specifically, we use a \textit{general-truncated-normal (GTN)} prior (see Definition~\ref{definition:general_truncated_normal}) on the entries of $\bW$. This prior ensures that sampled values lie within a bounded interval (e.g., $[-1,1]$), thereby automatically satisfying the desired magnitude constraint without explicit enforcement during optimization.
In this framework, the identity structure in Equation~\eqref{equation:submatrix_bid_identity} is preserved by fixing the corresponding entries of $\bW$ to 1 (or incorporating them as deterministic nodes in the graphical model). The remaining entries are inferred probabilistically, yielding a stable and interpretable low-rank ID.

\index{Decomposition: GBT}
\index{Decomposition: GBTN}
\index{General-truncated-normal distribution}

\subsubsection{Bayesian GBT and GBTN Models for ID}
We now introduce the Bayesian ID model termed the \textit{GBT} model.
To enhance flexibility and reduce sensitivity to hyper-parameter choices, we further propose a hierarchical extension called the \textit{GBTN} model. This variant retains simple conditional density forms while requiring only modest additional computation.
Similar to the Bayesian treatment for PCA models (Section~\ref{section:bayespca}), 
we further extend the models with automatic relevance determination (ARD). 
Therefore, the effective dimensionality $K$ of the latent subspace can be automatically inferred from the data, eliminating the need to pre-specify it.

\begin{figure}[h]
\centering  
\vspace{-0.35cm} 
\subfigtopskip=2pt 
\subfigbottomskip=6pt 
\subfigcapskip=-2pt 
\subfigure[GBT.]{\includegraphics[width=0.4\textwidth]{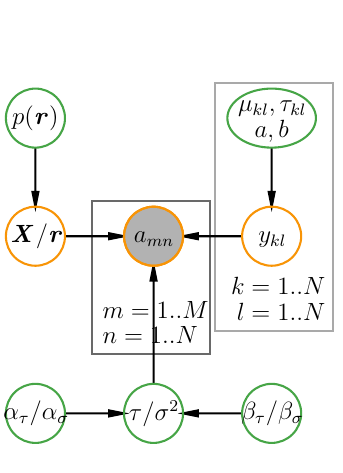} \label{fig:bmf_bid_GBT}}
\subfigure[GBTN.]{\includegraphics[width=0.4\textwidth]{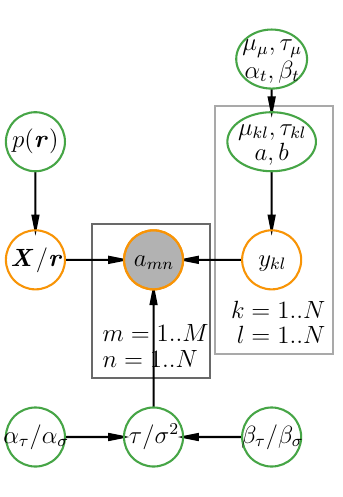} \label{fig:bmf_bid_GBTN}}
\caption{
Graphical representation of the GBT and GBTN models. Orange circles denote observed or latent variables (shaded nodes indicate observed quantities); green circles represent prior (hyper)parameters. Plates indicate replicated variables. In node labels, a slash ``/" means ``or," and a comma ``," means ``and." 
Parameters $a$ and $b$ are fixed to $a=-1$ and $b=1$ in our experiments; a weaker constraint would use $a=-2$ and $b=2$.}
\label{fig:bmf_bids}
\end{figure}

\paragrapharrow{Likelihood.}
We assume the data matrix $\bA$ is generated according to the probabilistic process depicted in Figure~\ref{fig:bmf_bids}. 
Each observed entry $a_{mn}$ of $\bA$ is modeled via a Gaussian likelihood with variance $\sigma^2$ and mean given by the low-rank reconstruction $\bx_m^\top\by_n$, consistent with the loss in \eqref{equation:idbid-per-example-loss}:
\begin{equation}\label{equation:gbt_data_entry_likelihood}
\begin{aligned}
p(a_{mn} \mid \bx_m^\top\by_n, \sigma^2) &= \normal(a_{mn}\mid \bx_m^\top\by_n, \sigma^2);\\
p(\bA\mid  \btheta) = \prod_{m,n=1}^{M,N}\normal \left(a_{mn}\mid  (\bX\bY)_{mn}, \sigma^2 \right) &
= \prod_{m,n=1}^{M,N} \normal \left(a_{mn}\mid  (\bX\bY)_{mn}, \tau^{-1} \right),
\end{aligned}
\end{equation}
where $\btheta=\{\bX,\bY,\sigma^2\}$ denotes all model parameters, $\normal(\cdot\mid \cdot)$ is the Gaussian distribution,  $\sigma^2$ is the variance, and $\tau^{-1}=\sigma^2$ is the precision.

\index{Inverse-Gamma distribution}
\paragrapharrow{Prior.}
We choose a conjugate prior over the data variance, an inverse-Gamma distribution (Definition~\ref{definition:inverse_gamma_distribution}) with shape $\alpha_\sigma$ and scale $\beta_\sigma$, 
\begin{equation}\label{equation:prior_gbt_gamma_on_variance}
p(\sigma^2 \mid  \alpha_\sigma, \beta_\sigma) = \inversegammadist(\sigma^2 \mid  \alpha_\sigma, \beta_\sigma).
\end{equation}
Equivalently, one could assign a Gamma prior $\gammadist(\tau \mid  \alpha_\tau, \beta_\tau)$ to the precision $\tau={1}/{\sigma^2}$; we omit further details here (see Equation~\eqref{equation:ggg_gamma_prior} in the GGG model).

The entries $y_{kl}$ of the latent matrix $\bY$ (with $k,l\in \{1,2,\ldots,N\}$; see Figure~\ref{fig:bmf_bids}) are treated as random variables. 
To encode our belief that their magnitudes should not exceed 1---consistent with the interpolative decomposition constraint---we assign them independent general-truncated-normal (GTN) priors (Definition~\ref{definition:general_truncated_normal}):
\begin{equation}\label{equation:rn_prior_bidd}
\begin{aligned}
p(y_{kl} \mid  \cdot ) &= \generaltruncatednormal(y_{kl} \mid  \mu_{kl}, (\tau_{kl})^{-1}, a=-1, b=1)\\
&=	
\frac{\sqrt{\frac{\tau_{kl}}{2\pi}} \exp \{-\frac{\tau_{kl}}{2}(y_{kl}-\mu_{kl})^2  \}  }
{\Phi\left((b-\mu_{kl})\cdot \sqrt{\tau_{kl}}\right)-\Phi\left((a-\mu_{kl})\cdot \sqrt{\tau_{kl}}\right)} \cdot \indicator(a\leq y_{kl}\leq b),
\end{aligned}
\end{equation}
where $\indicator(\cdot)$ is the indicator function (equal to 1 when the condition holds, and 0 otherwise).
This prior enforces the key ID constraint that no entry of $\bY$ exceeds magnitude 1. Moreover, it is conjugate to the Gaussian likelihood (see Equation~\eqref{equation:conjugate_truncated_constrained_mean}), ensuring that the posterior over $y_{kl}$ remains a GTN distribution.
In a relaxed version of the interpolative decomposition, the bound can be loosened to 2; the GTN prior accommodates this simply by setting $a=-2$ and $b=2$.

\paragrapharrow{Hierarchical prior.}
To increase model flexibility and reduce dependence on fixed hyper-parameters, we place a joint hyperprior over the GTN parameters $\{\mu_{kl}, \tau_{kl}\}$. 
Specifically, we adopt the \textit{GTN-scaled-normal-Gamma (GTNSNG)} prior:
\begin{equation}
\begin{aligned}
&\gap p(\mu_{kl}, \tau_{kl} \mid \cdot) 
=\gtnsng(\mu_{kl}, \tau_{kl}\mid  \mu_\mu, \frac{1}{\tau_\mu},\alpha_t, \beta_t)\\
&=	\left\{\Phi((b-\mu_\mu)\cdot \sqrt{\tau_\mu})-\Phi((a-\mu_\mu)\cdot \sqrt{\tau_\mu})\right\}
\cdot 
\normal(\mu_{kl}\mid  \mu_\mu, (\tau_\mu)^{-1}) \cdot \gammadist(\tau_{kl} \mid  \alpha_t, \beta_t).
\end{aligned}
\end{equation}
See Figure~\ref{fig:bmf_bid_GBTN}. This construction decouples $\{\mu_{kl}\}$ and $\{\tau_{kl}\}$, 
yielding conditionally conjugate posteriors: a normal distribution for $\{\mu_{kl}\}$ and a Gamma distribution for $\{\tau_{kl}\}$.

\paragrapharrow{Terminology.} 
Following the convention established in Section~\ref{section:bmf_real_intro} for Bayesian matrix factorization, we refer to these models as GBT and GBTN. Here, the letter ``\textit{B}" reflects the underlying Beta-Bernoulli structure inherent in the model's design.

\subsection{Gibbs Sampler}\label{section:gbt_gbtn_derivation}
We now present the derivation of the Gibbs sampler for the Bayesian ID models introduced earlier---namely, the GBT and GBTN models.

\paragrapharrow{Update of latent variables.}
The conditional posterior distribution of each latent variable $y_{kl}$ ($k,l=1,2,\ldots,N$) is a GTN distribution.
Let $\bY_{-kl}$ denote all entries of $\bY$ except $y_{kl}$. 
Based on the graphical model in Figure~\ref{fig:bmf_bids}, the full conditional posterior of $y_{kl}$ is proportional to the product of the likelihood and its GTN prior:
\begin{equation}\label{equation:posterior_gbt_ykl}
\begin{aligned}
&\gap p(y_{kl} \mid  \bA, \bX, \bY_{-kl}, \mu_{kl}, \tau_{kl}, \sigma^2) \propto  p(\bA\mid \bX,\bY, \sigma^2) \cdot p(y_{kl}\mid \mu_{kl}, \tau_{kl} )\\
&=\prod_{m,n=1}^{M,N} \normal \left(a_{mn}\mid  \bx_m^\top\by_n, \sigma^2 \right)\times
\generaltruncatednormal(y_{kl} \mid  \mu_{kl}, (\tau_{kl})^{-1},a=-1,b=1) \\
&\propto
\exp\Bigg\{  
-\frac{1}{2\sigma^2} \sum_{m,n=1}^{M,N} (a_{mn} - \bx_m^\top\by_n)^2
\Bigg\}
\exp \left\{-\frac{\tau_{kl}}{2}(y_{kl}-\mu_{kl})^2  \right\}
u(y_{kl} \mid  a,b)\\
&\propto
\exp\Bigg\{  
-\frac{1}{2\sigma^2} \sum_{m}^{M} (a_{ml} - \bx_m^\top \by_l)^2
\Bigg\}
\exp \Bigg\{-\frac{\tau_{kl}}{2}(y_{kl}-\mu_{kl})^2  \Bigg\}
u(y_{kl} \mid  a,b)\\
&{\footnotesize
\propto
\exp\Bigg\{  
-\frac{1}{2\sigma^2} \sum_{m}^{M} \bigg( x_{mk} ^2y_{kl }^2  + 2x_{mk} y_{kl } 
\big(\sum_{n\neq k}^{N}x_{mn} y_{nl}-a_{ml}\big)\bigg)
\Bigg\}
\exp \{-\frac{\tau_{kl}}{2}(y_{kl}-\mu_{kl})^2  \}
u(y_{kl} \mid  a,b)
}\\
&\propto
{\footnotesize
\exp\Bigg\{  
-y_{kl }^2
\underbrace{\Bigg(\frac{\sum_{m}^{M}  x_{mk} ^2}{2\sigma^2}+\textcolor{black}{\frac{\tau_{kl}}{2}} \Bigg)}_
{\textcolor{mylightbluetext}{\triangleq  \widetilde{\tau}/2}}
+y_{kl } 
\underbrace{\bigg(\frac{1}{\sigma^2}  \sum_{m}^{M} x_{mk}  \big(a_{ml}-\sum_{n\neq k}^{N}x_{mn}
y_{nl}\big)
+\textcolor{black}{\tau_{kl}\mu_{kl}}
\bigg)}_{\textcolor{mylightbluetext}{\triangleq\widetilde{\tau} \cdot \widetilde{\mu}}}
\Bigg\}
u(y_{kl} \mid  a,b)
}\\
&\propto \normal(y_{kl}\mid  \widetilde{\mu},( \widetilde{\tau})^{-1})u(y_{kl} \mid  a,b) 
\propto \generaltruncatednormal(y_{kl}\mid  \widetilde{\mu},( \widetilde{\tau})^{-1}, a=-1,b=1).
\end{aligned}
\end{equation}
Here, $\bx_m$ denotes the $m$-th row of $\bX$, and $\by_l$ the $l$-th column of $\bY$. 
The quantity $\widetilde{\tau} ={(\sum_{m}^{M}  x_{mk} ^2)}/{\sigma^2} +\tau_{kl}$ is the posterior ``parent" precision of the GTN density, and the posterior ``parent" mean of the GTN density is
$$
\widetilde{\mu} = \bigg(\frac{1}{\sigma^2}  \sum_{m}^{M} x_{mk}  \big(a_{ml}-\sum_{n\neq k}^{N}x_{mn}
y_{nl}\big)
+\textcolor{black}{\tau_{kl}\mu_{kl}}
\bigg) \bigg/ \widetilde{\tau}.
$$

\paragrapharrow{Update of the variance parameter.}
By conjugacy, the conditional posterior of the noise variance $\sigma^2$  is inverse-Gamma:
\begin{equation}\label{equation:posterior_gnt_sigma2}
\begin{aligned}
	&\gap p(\sigma^2 \mid  \bX, \bY, \bA)
	= \inversegammadist(\sigma^2 \mid  \widetilde{\alpha_\sigma}, \widetilde{\beta_\sigma}),
\end{aligned}
\end{equation}
with updated hyper-parameters: $\widetilde{\alpha_\sigma} = {(MN)}/{2}+\alpha_\sigma$, 
$\widetilde{\beta_\sigma}=\frac{1}{2} \sum_{m,n=1}^{M,N}(a_{mn}-\bx_m^\top\by_n)^2+\beta_\sigma$.

\paragrapharrow{Update of state vector for GBT and GBTN without ARD.}
Let  $\br\in\{0,1\}^N$ be the state vector indicating column roles: $r_n=1$ if $\ba_n$ is a basis column,   and $r_n=0$ if it is interpolated. 
Given the state vector $\br=[r_1,r_2, \ldots, r_N]^\top\in \real^N$, the relation between $\br$
and the index sets $\sJ$ is simple; 
$\sJ = \sJ(\br) = \{n\mid r_n = 1\}_{n=1}^N$ and $\sI = \sI(\br) = \{n\mid r_n = 0\}_{n=1}^N$. 

To update $\br$, we propose swapping one basis column $j\in\sJ$ with one interpolated column $i\in\sI$. Let 
$\br_{-ji}$ denote $\br$ with entries $j$ and $i$ removed. The acceptance odds for flipping $r_j=1\rightarrow 0$ and $r_i=0\rightarrow 1$ are:
\begin{equation}\label{equation:postrerior_gbt_rvector}
\begin{aligned}
j&\in \sJ; 
\qquad 
i\in \sI;\\
o_j &= 
\frac{p(r_j=0, r_i=1\mid \bA,\sigma^2, \bY, \br_{-ji})}
{p(r_j=1, r_i=0\mid \bA,\sigma^2, \bY, \br_{-ji})}\\
&=
\frac{p(r_j=0, r_i=1)}{p(r_j=1, r_i=0)}\times
\frac{p(\bA\mid \sigma^2, \bY, \br_{-ji}, r_j=0, r_i=1)}{p(\bA\mid \sigma^2, \bY, \br_{-ji}, r_j=1, r_i=0)}.
\end{aligned}
\end{equation}
Under a symmetric (uninformative) prior, we set $p(r_j=0, r_i=1)=p(r_j=1, r_i=0)$, so the ratio depends only on the likelihood. The full conditional probability becomes:
\begin{equation}\label{equation:postrerior_gbt_rvec_withoutard}
p(r_j=0, r_i=1\mid \bA,\sigma^2, \bY, \br_{-ji}) = \frac{o_j}{1+o_j}.
\end{equation}
This defines a Metropolis--Hastings step within the Gibbs sampler for updating the support of the basis matrix.

\paragrapharrow{Extra update for GBTN model.}
In the hierarchical GBTN model, the hyper-parameters $\mu_{kl}$ and $\tau_{kl}$ of the GTN prior are themselves random variables; see Figure~\ref{fig:bmf_bids}. Their conditionals are derived from the joint prior (GTNSNG) and the likelihood.
Integrating out irrelevant terms, the conditional posterior of $\mu_{kl}$ is Gaussian:
\begin{equation}\label{equation:posterior_gbt_mukl}
\begin{aligned}
&\gap p(\mu_{kl} \mid  \tau_{kl}, \mu_\mu, \tau_\mu, \alpha_t, \beta_t, y_{kl})\\
&\propto \generaltruncatednormal(y_{kl} \mid  \mu_{kl}, (\tau_{kl})^{-1}, a=-1, b=1)
 \cdot \gtnsng(\mu_{kl}, \tau_{kl}\mid  \mu_\mu, (\tau_\mu)^{-1},\alpha_t, \beta_t)\\
&\propto\generaltruncatednormal(y_{kl} \mid  \mu_{kl}, (\tau_{kl})^{-1}, a=-1, b=1)
\cdot 
\big\{\Phi((b-\mu_\mu)\cdot \sqrt{\tau_\mu})-\Phi((a-\mu_\mu)\cdot \sqrt{\tau_\mu})\big\} \\
&\gap\gap\gap\gap\gap\gap\gap\gap\gap\gap\gap\gap\gap	
\cdot{\normal(\mu_{kl}\mid  \mu_\mu, (\tau_\mu)^{-1})} \cdot 
\cancel{\gammadist(\tau_{kl} \mid  \alpha_t, \beta_t)}\\
\end{aligned}
\end{equation}
\begin{align*}
&\propto 
\sqrt{\tau_{kl}}\cdot \exp\left\{ -({\tau_{kl}}/{2}) (y_{kl}-\mu_{kl})^2\right\}
\cdot \exp\left\{ -({\tau_\mu}/{2})(\mu_\mu - \mu_{kl})^2  \right\}\\
&\propto \exp\Bigg\{  - \mu_{kl}^2 
\underbrace{{(\tau_{kl}+\tau_\mu)}/{2}}_{\textcolor{mylightbluetext}{\triangleq \widetilde{t}/2}}
+ \mu_{kl}
\underbrace{(\tau_{kl}y_{kl}+\tau_\mu\mu_\mu)}_{\textcolor{mylightbluetext}{\triangleq\widetilde{m}\cdot \widetilde{t}}}  \Bigg\}\propto 
\normal(\mu_{kl}\mid  \widetilde{m},(\,\widetilde{t}\,)^{-1}),
\qquad \qquad \quad
\end{align*}
where $\widetilde{t} = \tau_{kl}+\tau_\mu$ and $\widetilde{m} =(\tau_{kl}y_{kl}+\tau_\mu\mu_\mu)/\widetilde{t}$ are the posterior precision and mean of the normal density, respectively, and $\Phi(\cdot)$ is the cumulative distribution function  of $\normal(0,1)$. Similarly, the conditional density of $\tau_{kl}$ is,
\begin{equation}\label{equation:posterior_gbt_taukl}
\begin{aligned}
&\gap p(\tau_{kl} \mid  \mu_{kl}, \mu_\mu, \tau_\mu, \alpha_t, \beta_t, y_{kl})\\
&\propto \generaltruncatednormal(y_{kl} \mid  \mu_{kl}, (\tau_{kl})^{-1}, a=-1, b=1)
	\cdot \gtnsng(\mu_{kl}, \tau_{kl}\mid  \mu_\mu, (\tau_\mu)^{-1},\alpha_t, \beta_t)\\
&\propto\generaltruncatednormal(y_{kl} \mid  \mu_{kl}, (\tau_{kl})^{-1}, a=-1, b=1)
\cdot 
\big\{\Phi((b-\mu_\mu)\cdot \sqrt{\tau_\mu})-\Phi((a-\mu_\mu)\cdot \sqrt{\tau_\mu})\big\}\\
&\gap\gap\gap\gap\gap\gap\gap\gap\gap\gap\gap\gap\gap	
\cdot 
\cancel{\normal(\mu_{kl}\mid  \mu_\mu, (\tau_\mu)^{-1})} \cdot 
{\gammadist(\tau_{kl} \mid  \alpha_t, \beta_t)}\\
&\propto \exp\left\{  -\tau_{kl}  \frac{(y_{kl}- \mu_{kl})^2}{2}  \right\}
\tau_{kl}^{1/2} \tau_{kl}^{\alpha_t-1} \exp\left\{  -\beta_t \tau_{kl} \right\}\\
&\propto \exp\left\{   -\tau_{kl}\left[ \beta_t +  \frac{(y_{kl}- \mu_{kl})^2}{2}  \right] \right\}
\cdot \tau_{kl}^{(\alpha_t+1/2)-1}
\propto \gammadist(\tau_{kl} \mid  \widetilde{a}, \widetilde{b}),
\end{aligned}
\end{equation}
where $\widetilde{a} = \alpha_t+1/2$ and $\widetilde{b}=\beta_t +  {(y_{kl}- \mu_{kl})^2}/{2}$ are the posterior parameters of the Gamma density.

The complete Gibbs sampling procedure for both GBT and GBTN is summarized in Algorithm~\ref{alg:gbtn_gibbs_sampler}. While presented in an explanatory (element-wise) form for clarity, a vectorized implementation would significantly improve computational efficiency.

\begin{algorithm}[htb] 
\caption{Gibbs sampler for GBT and GBTN ID models. The procedure presented here may not be efficient but is explanatory. A more efficient one can be implemented in a vectorized manner. By default, uninformative priors are $a=-1, b=1,\alpha_\sigma=0.1, \beta_\sigma=1$, ($\{\mu_{kl}\}=0, \{\tau_{kl}\}=1$) for GBT, ($\mu_\mu =0$, $\tau_\mu=0.1, \alpha_t=\beta_t=1$) for GBTN.} 
\label{alg:gbtn_gibbs_sampler}  
\begin{algorithmic}[1] 
\For{$t=1$ to $T$}\Comment{$T$ iterations}
\State \algoalign{Sample state vector $\br$ from Equation~\eqref{equation:postrerior_gbt_rvec_withoutard};}
\State \algoalign{Update matrix $\bX$ by $\bA[:,\sJ]$ where index vector $\sJ$ is the index of $\br$ with value 1 and set $\bX[:,\sI]=\bzero$ where index vector $\sI$ is the index of $\br$ with value 0;}
\State Sample $\sigma^2$ from $p(\sigma^2 \mid  \bX,\bY, \bA)$ in Equation~\eqref{equation:posterior_gnt_sigma2}; 
\For{$k=1$ to $N$} 
\For{$l=1$ to $N$} 
\State Sample $y_{kl}$ from Equation~\eqref{equation:posterior_gbt_ykl};
\State (GBTN only) Sample $\mu_{kl}$ from Equation~\eqref{equation:posterior_gbt_mukl};
\State (GBTN only) Sample $\tau_{kl}$ from Equation~\eqref{equation:posterior_gbt_taukl};
\EndFor
\EndFor
\State Report loss in Equation~\eqref{equation:idbid-per-example-loss}, stop if it converges.
\EndFor
\State Report average  loss in Equation~\eqref{equation:idbid-per-example-loss} after burn-in iterations.
\end{algorithmic} 
\end{algorithm}

\index{Bayesian ID}
\subsection{Aggressive Update}
In Algorithm~\ref{alg:gbtn_gibbs_sampler}, after sampling a new state vector $\br$, we set the interpolated columns of $\bX$ to zero: $\bX[:,\sI]=\bzero$, where $\sI=\{n\mid r_n=0\}$.
However, in the next iteration, the state vector may change---specifically, an index $i\in\sI$ might switch from $r_i =0$ to $r_i=1$:
$$
r_i=0 \rightarrow r_i=1.
$$
If $\bX[:,i]$ remains zero at this point, the update of the corresponding entries in $\bY$ (via Equation~\eqref{equation:posterior_gbt_ykl}) becomes ill-defined or uninformative, since the column provides no ``gradient" signal.
To address this, we introduce an \textit{aggressive update} strategy. Instead of committing immediately to the current state vector $\br$, we maintain two candidate states:
\begin{itemize}
\item the current state $\br_1$ (with associated factor matrix $\bX_1$), and
\item a proposed state $\br_2$ (with proposal matrix $\bX_2$).
\end{itemize}
At each iteration, we sample $\br$ from $\{\br_1, \br_2\}$ according to their posterior odds (Equation~\eqref{equation:postrerior_gbt_rvec_withoutard}). Depending on the selected state, we use the corresponding $\bX$ to update $\bY$:
\begin{itemize}
\item if $\br=\br_1$, we use $\bY_1$ (computed with $\bX_1$);
\item if $\br=\br_2$, we use $\bY_2$ (computed with $\bX_2$).
\end{itemize}
This ensures that whenever a column is promoted to a basis column, its latent representation in $\bY$ is updated using a nonzero $\bX$, avoiding degenerate updates. We refer to this approach as the \textit{aggressive Gibbs sampler}.
The aggressive sampler for the GBT model is detailed in Algorithm~\ref{alg:gbtn_gibbs_sampler_aggressive}. For brevity, we omit the GBTN variant, as it follows analogously by including hyper-parameter updates.

\begin{algorithm}[h] 
\caption{\textit{Aggressive} Gibbs sampler for GBT ID model. The procedure presented here may not be efficient but is explanatory. A more efficient one can be implemented in a vectorized manner. By default, uninformative priors are $a=-1, b=1,\alpha_\sigma=0.1, \beta_\sigma=1$, ($\{\mu_{kl}\}=0, \{\tau_{kl}\}=1$) for GBT. 
} 
\label{alg:gbtn_gibbs_sampler_aggressive}  
\begin{algorithmic}[1] 
\For{$t=1$ to $T$}\Comment{$T$ iterations}
\State \algoalign{Sample state vector $\br$ from $\{\br_1, \br_2\}$ by Equation~\eqref{equation:postrerior_gbt_rvec_withoutard};}
\State Decide $\bY$: $\bY=\bY_1$ if $\br$ is $\br_1$; $\bY=\bY_2$ if $\br$ is $\br_2$;
\State Update state vector $\br_1=\br$; 
\State Sample proposal state vector $\br_2$ based on $\br$;
\State Update matrix $\bX$ by $\br=\br_1$; 
\State Update proposal $\bX_2$ by $\br_2$;
\State Sample $\sigma^2$ from $p(\sigma^2 \mid  \bX,\bY, \bA)$ in Equation~\eqref{equation:posterior_gnt_sigma2}; 
\State Sample $\bY_1 = \{y_{kl}\}$ using $\bX$;
\State Sample $\bY_2 = \{y_{kl}\}$ using $\bX_2$;
\State Report loss in Equation~\eqref{equation:idbid-per-example-loss}, stop if it converges.
\EndFor
\State Report average  loss in Equation~\eqref{equation:idbid-per-example-loss} after burn-in iterations.
\end{algorithmic} 
\end{algorithm}

\subsection{Post-Processing}\label{section:id_pprocess}
The Gibbs sampler yields an approximation $\bA\approx \bX\bY$, where $\bX\in \real^{M\times N}$ and $\bY\in \real^{N\times N}$. As described earlier, the effective low-rank factors can be extracted using the index set $\sJ=\{n\mid r_n=1\}$:
$$
\begin{aligned}
	\bC&=\bX[:,\sJ]=\bA[:,\sJ],\\
	\bW&= \bY[\sJ,:].
\end{aligned}
$$
However, in a true interpolative decomposition, the submatrix $\bY[\sJ,\sJ]=\bW[:,\sJ]$
must be the identity matrix (see Equation~\eqref{equation:submatrix_bid_identity}). The Gibbs sampling procedure does not enforce this constraint explicitly, so $\bY[\sJ,\sJ]$ may deviate from $\bI_K$.
To correct this, we perform a post-processing step: we replace $\bY[\sJ,\sJ]$ with $\bI_K$ and adjust the remaining rows of $\bW$ accordingly. This enforces the exact ID structure and typically reduces the reconstruction error further. The procedure is illustrated in Figure~\ref{fig:id-column}.

\index{Decomposition: GBT with ARD}
\index{Automatic relevance determination}
\subsection{Bayesian ID with Automatic Relevance Determination}
We further extend the Bayesian ID framework with automatic relevance determination (ARD) to eliminate the need for pre-specifying the number of basis columns $K$. Instead, the model infers $K$ automatically from the data.
Let  $\br=[r_1,r_2, \ldots, r_N]^\top\in \real^N$ be the state vector, with $\sJ = \sJ(\br) = \{n\mid r_n = 1\}_{n=1}^N$ and $\sI = \sI(\br) = \{n\mid r_n = 0\}_{n=1}^N$. 
Unlike the fixed-$K$ setting, ARD-type prior allows any index $j\in\{1,2,\ldots,N\}$ to toggle between basis and interpolated status. The full conditional for $r_j$ is:
such that 
\begin{equation}\label{equation:postrerior_gbt_rvector_ard}
\begin{aligned}
j&\in \sJ\cup \sI;\\
o_j = 
\frac{p(r_j=0\mid \bA,\sigma^2, \bY, \br_{-j})}
{p(r_j=1\mid \bA,\sigma^2, \bY, \br_{-j})}
&=
\frac{p(r_j=0)}{p(r_j=1)} \times
\frac{p(\bA\mid \sigma^2, \bY, \br_{-j}, r_j=0)}{p(\bA\mid \sigma^2, \bY, \br_{-j}, r_j=1)},
\end{aligned}
\end{equation}
where $\br_{-j}$ denotes all elements of $\br$ except the $j$-th element.
Compared to Equation~\eqref{equation:postrerior_gbt_rvector} (which swaps one basis and one interpolated column, keeping $\abs{\sJ}$ fixed), Equation~\eqref{equation:postrerior_gbt_rvector_ard} allows the size of 
$\sJ$ to vary. Under a uniform prior ($p(r_j=0)=p(r_j=1)=0.5$), the posterior probability becomes:
\begin{equation}\label{equation:postrerior_gbt_rvector222_ard}
p(r_j=0\mid \bA,\sigma^2, \bY, \br_{-j}) = \frac{o_j}{1+o_j}.
\end{equation}
The full Gibbs sampler with ARD is given in Algorithm~\ref{alg:gbtn_gibbs_sampler_withard}. A key difference is that all entries of $\br$ are updated sequentially (lines 2--4), rather than swapping pairs.
However, updating many entries of $\br$ simultaneously can cause abrupt changes in $\bX$, leading to unstable updates of $\bY$. To mitigate this, we introduce a \textit{critical update phase}: after resampling the entire state vector $\br$, we perform $\nu$ additional Gibbs sweeps over $\bY$ (and its hyper-parameters in GBTN) to allow the latent factors to adapt to the new support structure. This stabilization step is highlighted in \textcolor{mylightbluetext}{blue} in Algorithm~\ref{alg:gbtn_gibbs_sampler_withard}.

\begin{algorithm}[ht] 
\caption{Gibbs sampler for GBT and GBTN ID with \textit{ARD} models.  The procedure presented here can be inefficient but is explanatory. While a vectorized manner can be implemented to find a more efficient algorithm. By default, weak priors are $a=-1, b=1,\alpha_\sigma=0.1, \beta_\sigma=1$, ($\{\mu_{kl}\}=0, \{\tau_{kl}\}=1$) for GBT, ($\mu_\mu =0$, $\tau_\mu=0.1, \alpha_t=\beta_t=1$) for GBTN. \textcolor{mylightbluetext}{Number of critical steps: $\nu$}.} 
\label{alg:gbtn_gibbs_sampler_withard}  
\begin{algorithmic}[1] 
\For{$t=1$ to $T$} \Comment{$T$ iterations}
\For{\textcolor{mylightbluetext}{$j=1$ to $N$}}
\State \textcolor{mylightbluetext}{Sample state vector element $r_j$ from Equation~\eqref{equation:postrerior_gbt_rvector222_ard}};
\EndFor
\State \algoalign{Update matrix $\bX$ by $\bA[:,\sJ]$ where index vector $\sJ$ is the index of $\br$ with value 1 and set $\bX[:,\sI]=\bzero$ where index vector $\sI$ is the index of $\br$ with value 0;}
\State Sample $\sigma^2$ from $p(\sigma^2 \mid  \bX,\bY, \bA)$ in Equation~\eqref{equation:posterior_gnt_sigma2}; 
\For{\textcolor{mylightbluetext}{$n=1$ to $\nu$}}
\For{$k=1$ to $N$} 
\For{$l=1$ to $N$} 
\State Sample $y_{kl}$ from Equation~\eqref{equation:posterior_gbt_ykl};
\State (GBTN only) Sample $\mu_{kl}$ from Equation~\eqref{equation:posterior_gbt_mukl};
\State (GBTN only) Sample $\tau_{kl}$ from Equation~\eqref{equation:posterior_gbt_taukl};
\EndFor
\EndFor
\EndFor
\State Output loss in Equation~\eqref{equation:idbid-per-example-loss}, stop iteration if it converges;
\EndFor
\State Output averaged loss in Equation~\eqref{equation:idbid-per-example-loss} for evaluation after burn-in iterations;
\end{algorithmic} 
\end{algorithm}

\begin{table}[h]
\centering
\setlength{\tabcolsep}{4.4pt}
\begin{tabular}{lllll}
\hline
Data set        & Num. Rows & Num. Columns & Fraction observed & Matrix rank\\ \hline
CCLE $EC50$ & 502 & 48  &0.632   & 24\\
CCLE $IC50$ & 504 & 48 &0.965   & 24\\
Gene Body Methylation &160 & 254 &1.000   & 160\\
Promoter Methylation & 160  & 254    & 1.000           &160\\
\hline
\end{tabular}
\caption{Overview of the CCLE $EC50$, CCLE $IC50$, Gene Body Methylation, and Promoter Methylation data sets, giving the number of rows, columns, the fraction of entries that are observed, and the matrix rank.}
\label{table:datadescription_ard}
\end{table}
\begin{figure}[h]
\centering  
\vspace{-0.35cm} 
\subfigtopskip=2pt 
\subfigbottomskip=2pt 
\subfigcapskip=-5pt 
\subfigure{\includegraphics[width=0.231\textwidth]{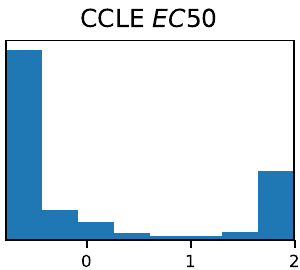} \label{fig:plot_ccle_ec}}
\subfigure{\includegraphics[width=0.231\textwidth]{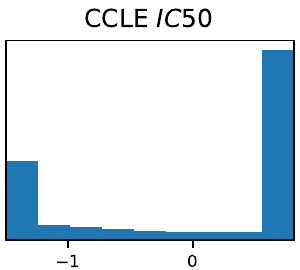} \label{fig:plot_ccle_ic}}
\subfigure{\includegraphics[width=0.231\textwidth]{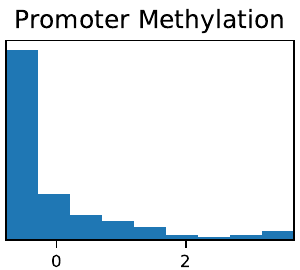} \label{fig:plot_ctrp}}
\subfigure{\includegraphics[width=0.231\textwidth]{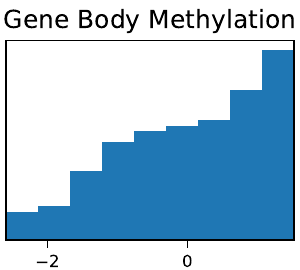} \label{fig:plot_movielens100}}
\caption{Data distribution of CCLE $EC50$, CCLE $IC50$, Gene Body Methylation, and Promoter Methylation datasets.}
\label{fig:datasets_bids_ard}
\end{figure}

\subsection{Examples for Bayesian ID}\label{section:bid_ard_experiments}
To evaluate the strategy and demonstrate the main advantages of the Bayesian ID method, 
we conduct experiments with different analysis tasks; and different datasets including 
Cancer Cell Line Encyclopedia (CCLE $EC50$ and CCLE $IC50$ datasets \citep{barretina2012cancer}), 
cancer driver genes (Gene Body Methylation \citep{koboldt2012comprehensive}), and  the promoter region (Promoter Methylation \citep{koboldt2012comprehensive}) 
from bioinformatics.
Following \citet{brouwer2017prior}, we preprocess these datasets by capping high values to 100 and undoing the natural log transform for the former three datasets.
All datasets are then standardized to have zero mean and unit variance, and missing entries are imputed with zeros.
To introduce controlled redundancy---particularly useful for evaluating column selection---we duplicate every column twice in the  CCLE $EC50$ and CCLE $IC50$ datasets. 
In contrast, the Gene Body Methylation and Promoter Methylation datasets already have more columns than their effective matrix rank, so no additional redundancy is introduced.
A summary of the four datasets is provided in Table~\ref{table:datadescription_ard}, and their empirical distributions are visualized in Figure~\ref{fig:datasets_bids_ard}.

In all experiments, we use identical parameter initialization across different tasks. Empirical results indicate that the post-processing step yields a modest performance gain, and that the GBT and GBTN models produce largely similar outcomes \citep{lu2022bayesian}. For clarity, we report only the post-processed results of the GBT model.
We compare the ARD-enhanced versions of GBT and GBTN against their vanilla (non-ARD) counterparts. Across a wide range of experiments and datasets, the ARD variants consistently achieve lower reconstruction error and match or outperform the vanilla methods in low-rank ID approximation---even when the latter are allowed to use the full matrix rank.

To quantify overall decomposition performance, we use the mean squared error (MSE), defined in Equation~\eqref{equation:idbid-per-example-loss}, which measures the discrepancy between the original and reconstructed matrices. Lower MSE indicates better performance.

\begin{figure*}[h]
\centering  
\vspace{-0.15cm} 
\subfigtopskip=2pt 
\subfigbottomskip=0pt 
\subfigcapskip=-2pt 
\subfigure[Convergence of the models on the 
CCLE $EC50$, CCLE $IC50$, Gene Body Methylation, and Promoter Methylation datasets, measured by the data fit (MSE). The algorithm almost converges in less than 50 iterations.]{\includegraphics[width=1\textwidth]{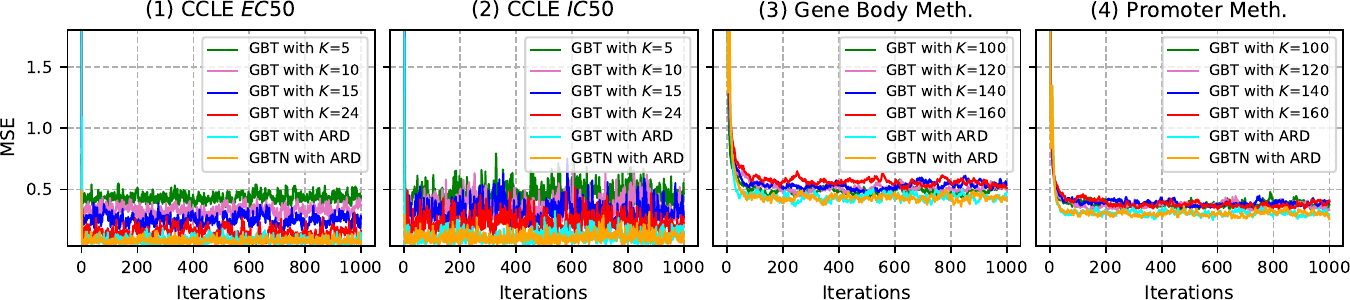} \label{fig:convergences_BIDs_ard}}
\subfigure[Averaged autocorrelation coefficients of samples of $y_{kl}$ computed using Gibbs sampling on the 
CCLE $EC50$, CCLE $IC50$, Gene Body Methylation, and Promoter Methylation datasets.]{\includegraphics[width=1\textwidth]{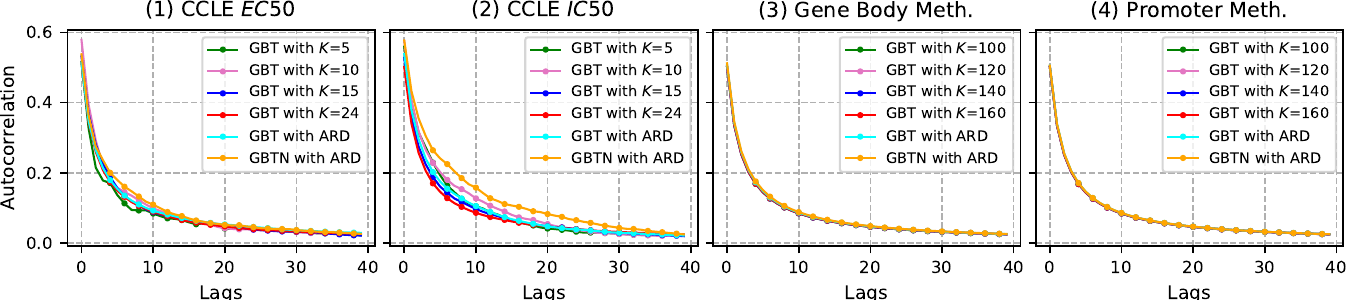} \label{fig:convergences_BIDs_autocorr_ard}}
\subfigure[Convergence of the number of selected columns on the 
CCLE $EC50$, CCLE $IC50$, Gene Body Methylation, and Promoter Methylation datasets. The algorithm almost converges in less than 100 iterations.]{\includegraphics[width=1\textwidth]{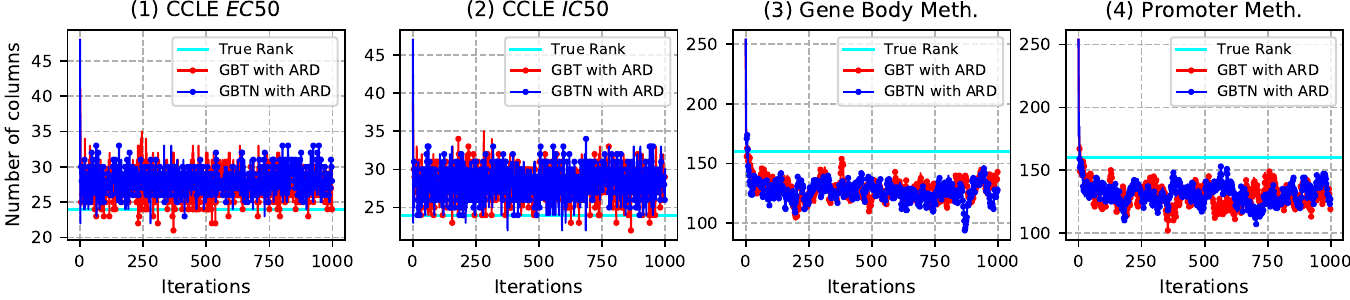} \label{fig:convergences_BIDs_numrvector_ard}}
\caption{Convergence results (upper), sampling mixing analysis (middle), and reconstructive results (lower) on the CCLE $EC50$, CCLE $IC50$, Gene Body Methylation, and Promoter Methylation datasets for various latent dimensions.
}
\label{fig:allresults_bids_ard}
\end{figure*}
\subsubsection{Hyper-parameters}
In these experiments, we use  $a=-1, b=1,\alpha_\sigma=0.1, \beta_\sigma=1$, ($\{\mu_{kl}\}=0, \{\tau_{kl}\}=1$) for GBT, ($\mu_\mu =0$, $\tau_\mu=0.1, \alpha_t=\beta_t=1$) for GBTN, and critical steps $\nu=5$ for GBT and GBTN with ARD.
These hyper-parameter choices are intentionally uninformative, and the models show little sensitivity to them. All observed and latent variables are initialized via random draws, which---given fixed hyper-parameters---provides a reasonable initial estimate of the underlying matrix structure.
In every setting, we run the Gibbs sampler for 1,000 iterations, discarding the first 100 as burn-in and applying thinning every 5 iterations. This configuration is justified by convergence diagnostics showing stabilization well before 100 iterations.

\begin{table}[]
\centering
\begin{tabular}{lllllll}
\hline
& $K_1$ & $K_2$ & $K_3$ & $K_4$ & GBT (ARD) & GBTN (ARD) \\ \hline
CCLE $EC50$    &   0.354  &   0.218  & 0.131  &  0.046  & \textbf{ 0.034 } & \textbf{ 0.031 } \\
CCLE $IC50$   &   0.301  &   0.231  & 0.161  &  0.103  & \textbf{ 0.035 } & \textbf{ 0.031 } \\
Gene Body Methylation   &   0.433  &   0.443  & 0.466  &  0.492  & \textbf{ 0.363 } & \textbf{ 0.372 } \\
Promoter Methylation   &   0.323  &   0.319  & 0.350  &  0.337  & \textbf{ 0.252 } & \textbf{ 0.263 } \\
\hline
\end{tabular}
\caption{Mean squared error (MSE) measure for varying latent dimensions $K$. 
For CCLE datasets, $K_1=5, K_2=10, K_3=15$, and $K_4=24$ (full rank); for methylation datasets, $K_1=100, K_2=120, K_3=140$, and $K_4=160$ (full rank).
ARD-based models outperform even full-rank non-ARD baselines.}
\label{table:covnergence_mse_reporte_ard}
\end{table}

\subsubsection{Convergence and Comparative Analysis}
We first examine convergence behavior across the four datasets. For the CCLE $EC50$ and CCLE $IC50$ datasets, we run GBT with $K=5,10,15,24$ (where $K=24$ equals the full matrix rank).
For the methylation datasets, we use $K=100,120,140,160$ ($K=160$ is full rank). Reconstruction error is measured by MSE.
Figure~\ref{fig:convergences_BIDs_ard} shows rapid convergence---typically within 50 iterations---across all settings. Figure~\ref{fig:convergences_BIDs_autocorr_ard} displays the averaged autocorrelation of Gibbs samples for $y_{kl}$. The autocorrelation drops below 0.1 for lags greater than 10, indicating good mixing of the Markov chain. Notably, the ARD and non-ARD variants exhibit comparable mixing properties. The sampling trajectories are smoother on the CCLE $EC50$, Gene Body Methylation, and Promoter Methylation datasets, whereas the CCLE $IC50$ dataset shows slightly noisier traces for the non-ARD GBT---suggesting that ARD may also improve numerical stability.

Comparative results (Figures~\ref{fig:convergences_BIDs_ard} and Table~\ref{table:covnergence_mse_reporte_ard}) consistently demonstrate that GBT and GBTN with ARD achieve the lowest MSE, even when compared to non-ARD models using the full matrix rank ($K=24$ for CCLE, $K=160$ for methylation data).

Figure~\ref{fig:convergences_BIDs_numrvector_ard} illustrates how the number of selected basis columns (i.e., 
$\sum_{n}r_n$) evolves during sampling. The chain stabilizes around 27 columns for the CCLE datasets and 130 columns for the methylation datasets---values close to the intrinsic ranks of the respective matrices. This confirms that the ARD mechanism automatically infers an appropriate number of basis columns, eliminating the need for manual rank selection.

\index{Decomposition: IID}
\index{Sigmoid}
\index{Intervened interpolative decomposition (IID)}
\section{Bayesian Intervened Interpolative Decomposition (IID)}\label{section:iid_main}

Expanding on the GBT model, we also introduce the  \textit{intervened interpolative decomposition (IID)} algorithm \citep{lu2022feature}.
The IID algorithm shares the same generative process as  in  Equation~\eqref{equation:gbt_data_entry_likelihood}, employing an inverse-Gamma prior over the variance parameter $\sigma^2$ (Equation~\eqref{equation:prior_gbt_gamma_on_variance})
and a GTN prior over the latent variables $\{y_{kl}\}$ (Equation~\eqref{equation:rn_prior_bidd}).
However, IID incorporates an additional assumption: some columns of the observed matrix 
$\bA$ are more important than others and should be prioritized during basis selection.

Suppose the relative importance of each column in  $\bA$ is encoded  by a \textit{raw importance vector} $\widehat{\bp}\in \real^N$,
where $\widehat{p}_n \in (-\infty, \infty)$ for all $n$ in $\{1,2,\ldots, N\}$. 
To map this into the unit interval, we apply the Sigmoid function:
$$
\bp = \text{Sigmoid}(\widehat{\bp}),
$$
where \textit{Sigmoid($\cdot$)} represents the function $f(x) = {1}/{(1+\exp\{-x\})}$ that can return a value in the range of 0 to 1.
The Sigmoid function acts as a squashing function because its domain is the set of all real numbers, and its range is (0, 1).
Then  we take the $\bp$ vector as the final \textit{importance vector} to indicate the importance of each column in   $\bA$.

Building on Equation~\eqref{equation:postrerior_gbt_rvector}, the odds ratio $o_j$ used in the Gibbs update is now modified to incorporate column importance:
\begin{equation}\label{equation:posterior_IID}
\begin{aligned}
o_j 
&=
\frac{p(r_j=0, r_i=1)}{p(r_j=1, r_i=0)} 
\times
\frac{p(\bA\mid \sigma^2, \bY, \br_{-ji}, r_j=0, r_i=1)}{p(\bA\mid \sigma^2, \bY, \br_{-ji}, r_j=1, r_i=0)}\\
&=
\textcolor{mylightbluetext}{
\frac{1-p_j }{p_j}
\frac{p_i }{1-p_i}
}
\times 
\frac{p(\bA\mid \sigma^2, \bY, \br_{-ji}, r_j=0, r_i=1)}{p(\bA\mid \sigma^2, \bY, \br_{-ji}, r_j=1, r_i=0)}.
\end{aligned}
\end{equation}
The corresponding conditional probability is then:
\begin{equation}\label{equation:postrerior_gbt_rvector222_IID}
p(r_j=0, r_i=1\mid \bA,\sigma^2, \bY, \br_{-ji}) = \frac{o_j}{1+o_j}.
\end{equation}
Because this approach intervenes in the standard Gibbs sampling procedure by biasing column selection toward higher-importance candidates, we refer to it as the \textit{intervened interpolative decomposition (IID)}.

\subsection{Quantitative Problem Statement}\label{section:iid_quantaprob}
Having introduced the IID algorithm, we now motivate its practical relevance in quantitative finance. Large hedge funds and asset managers increasingly rely on vast pools of predictive signals---often called alpha factors~\footnote{In quantitative finance, alpha factors are features designed to forecast future asset returns.}---to construct trading strategies. In industry practice, the number of such alphas can reach into the millions or even billions \citep{tulchinsky2019finding}.
Given this scale, constructing a robust meta-alpha (a composite signal that captures true trading signals) from the full alpha pool presents several challenges:
\begin{enumerate}[label=(\roman*)]
\item \textit{Trading capacity constraints.} Popular alphas may target illiquid assets. If many traders use the same signals, transaction costs and market impact can erode profitability.
\item \textit{Overfitting risk.} Using too many alphas increases model complexity, often leading to poor out-of-sample (OS) performance.
\item \textit{Multicollinearity.} Many alphas are highly correlated or functionally redundant. This can impair the performance of machine learning models (e.g., neural networks, XGBoost) that struggle with multicollinear inputs when learning meta-strategies.
\item \textit{Computational burden.} Processing billions of alphas is computationally expensive and often infeasible under real-world resource constraints.
\item \textit{Risk diversification.} To mitigate systemic risk, practitioners seek diverse subsets of alphas with low mutual correlation, enabling robust strategy testing and portfolio construction.
\end{enumerate}
For these reasons, there is a pressing need for algorithms that select a small, high-quality subset of alphas---one that avoids overfitting, scales efficiently, and delivers results in reasonable time.
A naive approach is to rank alphas by their \textit{RankIC} (defined below) and select the top performers. However, this ignores two critical issues: (i) the selected set may not be representative of the full alpha pool, and (ii) top-ranked alphas are often highly correlated, offering little diversification.

Our goal is to identify a small subset of alphas that are both high-performing and representative---meaning they can accurately reconstruct the remaining alphas with low error. Traditional ID methods (e.g., randomized ID \citep{liberty2007randomized} or Bayesian ID) can find representative columns, but they do not account for predictive performance and may select low-signal alphas.
In contrast, the IID method jointly optimizes for representativeness and desirability: it selects columns that (i) enable accurate reconstruction of the full matrix and (ii) exhibit high RankIC scores. This dual objective makes IID particularly well-suited for alpha selection in quantitative finance.

\index{Formulaic alphas}
\subsubsection{Formulaic Alphas}
Quantitative firms often design interpretable, rule-based signals known as formulaic alphas. For example, WorldQuant publicly released 101 such short-term alphas in 2016 \citep{kakushadze2016101}, and Guotai Junan Securities later published 191 alpha formulas widely adopted by practitioners \citep{guotaijunan2017}.
These alphas are derived from fundamental market data---including prices, volumes, volatility, and volume-weighted average price (VWAP)---and expressed as explicit mathematical expressions. As an illustration, a simple mean-reversion alpha might be defined as:
$$
\text{Alpha =  }- \left( \text{close(today) $-$close(5$\_$days$\_$ago ) } \right)/ 
\text{close(5$\_$days$\_$ago)}.
$$
This signal takes the opposite position of recent price movement: it goes long if the price has fallen over the past five days, and short otherwise. Intuitively, a higher alpha value suggests a greater likelihood of upward price movement in the near future.

\index{Spearman correlation}
\index{RankIC}
\subsubsection{Evaluation Metrics}
Let $r_{t}$ denote the return of a stock on day $t$, computed from closing prices $\{p_t\}_{t=1}^T$ as:
$$
r_t = \frac{p_t - p_{t-1}}{p_{t-1}}.
$$
We evaluate alpha effectiveness using the \textit{Rank information coefficient (RankIC)}:
\begin{equation}
\text{RankIC}(\ba, \br^h)=\text{Spearman}(\ba, \br^h),
\end{equation}
where $\text{Spearman}(\cdot)$ denotes the \textit{Spearman rank correlation}, $\ba$ is the vector of alpha values across stocks on a given day, and $\br^h$ is the vector of forward returns over a holding period of $h$ days (i.e., the $i$-th entry of $\br^h$ is the return realized $h$ days after day $i$). 
We must also account for the \textit{forward bias issue}: at each day $i$, the last $h$ entries in $\br^h[1:i]$ are undefined and thus set to NaN.

The RankIC measures the monotonic relationship between alpha predictions and future returns. A higher absolute RankIC indicates stronger predictive power. We use this metric directly as the importance score for each alpha, i.e., $\widehatp_n=\text{RankIC}_n$, which is then transformed via the Sigmoid function and plugged into Equation~\eqref{equation:posterior_IID} to guide column selection in IID.

\subsection{Examples for Bayesian IID}\label{section:iid_experiments}

For each stock $s\in \{1,2,\ldots, S\}$  (where $S$ is the total number of stocks), we construct an alpha matrix  $\bA_s\in \real^{N\times D}$, where $N$ is the number of alpha factors and $D$ is the number of trading days.
Each row of $\bA_s$ represents the time series of one alpha factor. Our goal is to select a subset of $M$ alphas from the full set of $N$.
We compute the RankIC between each alpha series and the forward return series with a horizon of $h=1$ day, and use this value directly as the importance score: a higher RankIC implies higher selection priority.

\begin{table}[h]
\centering
\small 
\setlength{\tabcolsep}{6pt}
\renewcommand{\arraystretch}{1.2}
\begin{tabular}{llllll}
\hline
Ticker & Type & Sector & Company &  Avg. Amount \\ \hline
SH601988   & Share & Bank               & Bank of China Limited    &  427,647,786  \\
SH601601   & Share & Public Utility    & China Pacific Insurance (Group)   & 819,382,926 \\
SH600028    & Share & Public Utility & China Petroleum \& Chemical Corporation              &  748,927,952\\
SH600016    & Share & Bank    & China Minsheng Banking Corporation &285,852,414  \\
SH601186    & Share & Public Utility    & China Railway Construction Corporation  & 594,970,588\\
SH601328    & Share & Bank & Bank of Communications Corporation &  484,445,915  \\
SH601628     & Share & Public Utility & China Life Insurance Company Limited&368,179,861  \\
SH601939    & Share   & Bank      & China Construction Bank Corporation &527,876,669  \\
\hline
SH510300    & ETF   & CSI 300         & Huatai-PineBridge CSI 300 ETF  &1,960,687,059  \\
SH510050      & ETF   & CSI 50 & ChinaAMC China CSI 50 ETF    &2,020,385,879  \\
\hline
\end{tabular}
\vspace{-0.25cm} 
\caption{Summary of the ten underlying assets in the China market. The average daily trading amount (in RMB) is computed over the test period.}
\label{table:iid_cn_data_summary}		
\end{table}

\begin{figure*}[h]
\centering  
\vspace{-0.3cm} 
\subfigtopskip=2pt 
\subfigbottomskip=9pt 
\subfigcapskip=-5pt 
\subfigure[Ten different portfolios, where we initialize each portfolio with a unitary value for clarity.]{\includegraphics[width=0.485\textwidth]{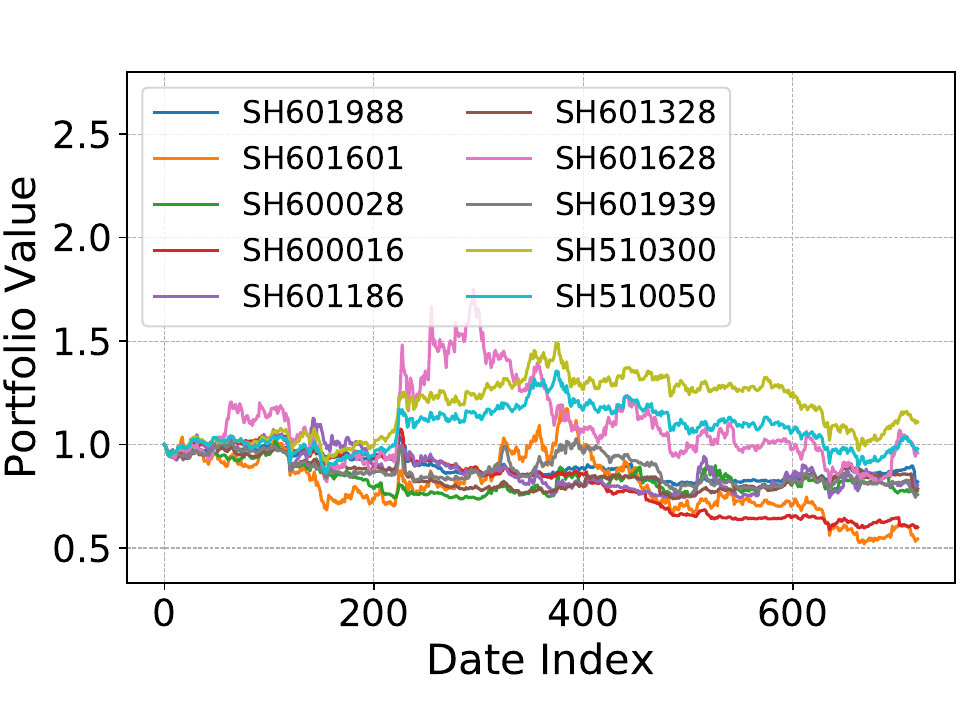} \label{fig:bid_iid_datasets_ashare}}
\subfigure[Portfolio values with the same strategy by using different alphas via comparative selection models.]{\includegraphics[width=0.485\textwidth]{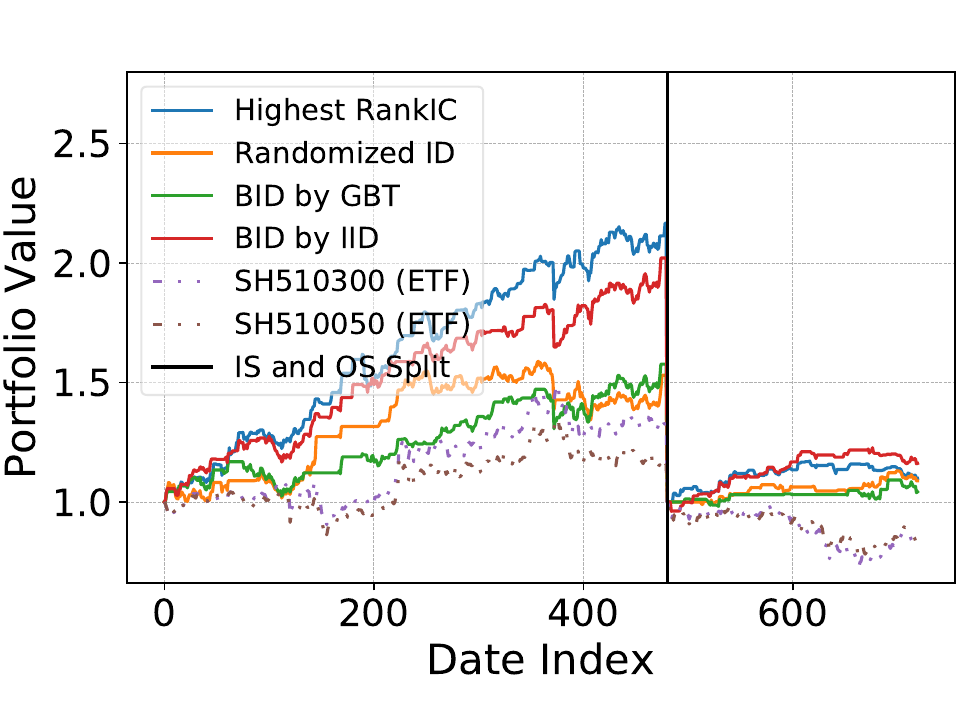} \label{fig:bid_iid_portfolio_ashare}}
\vspace{-0.7cm} 
\caption{Asset trajectories (left) and portfolio values (right), split into in-sample (IS) and out-of-sample (OS) periods. IID outperforms other methods in OS (see also Table~\ref{table:iid_selected_mean_ic}).}
\label{fig:bid_iid_portfolio_ashare_full}
\end{figure*}

\paragrapharrow{Data set.}
To evaluate the discussed algorithm and highlight the key advantages of the IID method, we conduct experiments on ten assets from the Chinese market, spanning diverse sectors including banking, public utilities, and ETFs. The data are obtained from Tushare~\footnote{\url{https://tushare.pro/}.}, and covers a three-year period, i.e., from 2018-07-18 to 2021-07-05 (720 trading days), where the data between 2018-07-18 and 2020-07-09 is considered the training set (480 calendar days); 
while data between 2020-07-10 and 2021-07-05 is taken as the test set (240 trading days).
The selected assets (summarized in Table~\ref{table:iid_cn_data_summary}) are chosen from the top 50 most liquid stocks/ETFs by average daily trading volume during the sample period, ensuring minimal trading constraints. Figure~\ref{fig:bid_iid_datasets_ashare} shows the normalized price trajectories of these assets (initialized to unit value for clarity).

We construct our alpha pool from three sources: (i) 78 alphas from the 101 formulaic alphas published by WorldQuant \citep{kakushadze2016101}; (ii) 94 alphas from the 191 formulaic alphas released by Guotai Junan Securities \citep{guotaijunan2017}; (iii) and 19 proprietary alphas.
All alphas are preprocessed to avoid extreme values. Consequently, each asset's alpha matrix $\bA_s$ has dimensions $191\times 720$.

In all experiments, we use identical parameter initialization across tasks. Empirical results show that post-processing yields a modest performance gain. For clarity, we report only the post-processed results of the GBT (without ARD) and IID models.
The IID model prioritizes columns (alphas) with high RankIC scores while maintaining low reconstruction error, consistently matching or outperforming vanilla GBT in low-rank ID approximation across datasets.
We again use MSE---defined in Equation~\eqref{equation:idbid-per-example-loss}---to evaluate decomposition quality. Lower MSE indicates better reconstruction of the original matrix.

\begin{figure*}[h]
\centering  
\vspace{-0.2cm} 
\subfigtopskip=2pt 
\subfigbottomskip=0pt 
\subfigcapskip=-2pt 
\subfigure[Convergence of the models on the 
SH510050, SH510300, SH601939, SH601628, and SH601328 datasets, as measured by MSE. 
Convergence occurs within 100 iterations.]{\includegraphics[width=1\textwidth]{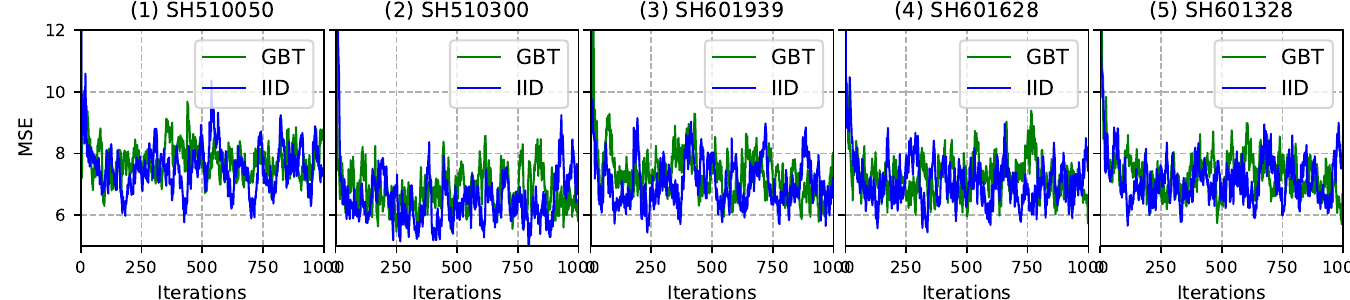} \label{fig:iid_alpha_convergence}}
\subfigure[Averaged autocorrelation coefficients of Gibbs samples for $y_{kl}$  on the 
SH510050, SH510300, SH601939, SH601628, and SH601328 datasets.]{\includegraphics[width=1\textwidth]{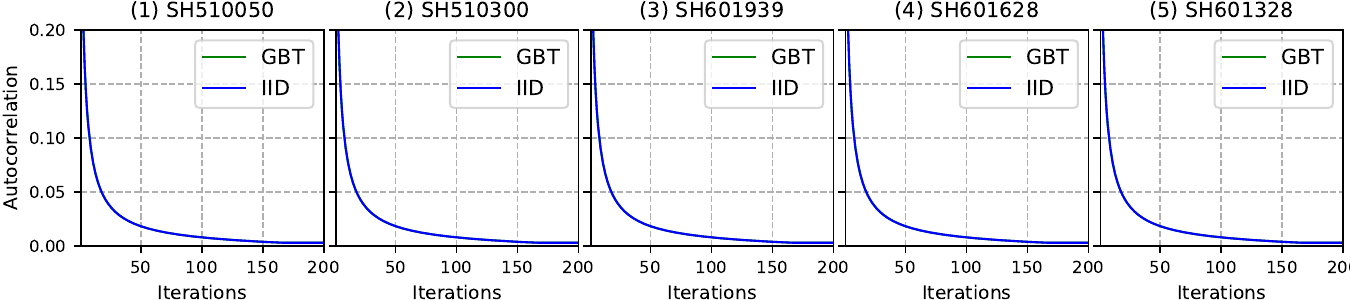} \label{fig:iid_alpha_autocorrelation}}
\vspace{-0.3cm} 
\caption{Convergence behavior (top) and sampling efficiency (bottom)  on the SH510050, SH510300, SH601939, SH601628, and SH601328 datasets for latent dimension $K=10$.
}
\label{fig:allresults_bids_ard_IID}
\end{figure*}
\paragrapharrow{Hyper-parameters.}
In these experiments, both GBT and IID use the following weakly informative priors:  $a=-1, b=1,\alpha_\sigma=0.1, \beta_\sigma=1$, ($\{\mu_{kl}\}=0, \{\tau_{kl}\}=1$).
These choices are intentionally uninformative, and the models exhibit low sensitivity to them.
All latent and observed variables are initialized via random draws, which---given fixed hyper-parameters---provides a reasonable initial estimate of the underlying structure.
We run the Gibbs sampler for 1,000 iterations, discarding the first 100 as burn-in and applying thinning every 5 iterations, consistent with convergence diagnostics showing stabilization within 100 iterations.

\begin{table}[]
\centering
\vspace{-0.35cm} 
\scriptsize
\setlength{\tabcolsep}{1.1pt}
\renewcommand{\arraystretch}{1.1}
\begin{tabular}{lllllllllll}
\hline 
& SH601988 & SH601601 & SH600028 & SH600016 & SH601186 & SH601328 & SH601628 & SH601939 & SH510300 & SH510050\\
\hline 
GBT Min. & 5.235 & 5.814 & 5.235 & \textbf{6.381}& 5.819 & 5.700 & 5.734 & 5.785 & 5.462 & 6.297 \\
IID Min. & \textbf{4.567} &\textbf{5.700}& \textbf{4.843} & 6.490 & \textbf{5.104} & \textbf{5.658} & \textbf{5.445} & \textbf{5.435} & \textbf{4.876} & \textbf{5.767} \\
GBT Mean & 6.476 & \textbf{7.367} & 6.764 & 8.053 & 7.066 & 7.250 & 7.206 & 7.242 & 6.769 & 7.776 \\
IID Mean & \textbf{6.239} & 7.449 & \textbf{6.664} & \textbf{7.831} & \textbf{6.558} & \textbf{7.081} & \textbf{7.002} & \textbf{7.031} & \textbf{6.450} & \textbf{7.492} \\
\hline
\end{tabular}
\vspace{-0.3cm} 
\caption{Minimal and mean MSE measures after burn-in across different iterations for GBT and IID models on the 10 alpha matrices from 10 assets. In all cases, $K=10$ is set as the latent dimension. 
In most cases, the results of IID converge to a smaller value than the GBT model.}
\label{table:comparis_gbt_iid_mse}
\vspace{-0.1cm} 
\end{table}

\subsubsection{Convergence and Comparative Analysis}
Due to space constraints, we present convergence results only for five assets:  SH510050, SH510300, SH601939, SH601628, and SH601328. Results for the remaining assets are qualitatively similar.

We run both GBT and IID with latent dimension $K=10$ (the full matrix rank is 191), measuring error via MSE. As shown in Figure~\ref{fig:iid_alpha_convergence}, both models converge within 100 iterations. Figure~\ref{fig:iid_alpha_autocorrelation} shows that autocorrelation drops below 0.1 for lags greater than 50, indicating good mixing of the Gibbs sampler.

Notably, despite prioritizing high-RankIC columns, the IID model does not sacrifice reconstruction accuracy---in fact, it often achieves lower MSE than vanilla GBT (see Table~\ref{table:comparis_gbt_iid_mse}).

\begin{algorithm}[!htb] 
\caption{Alpha selection for portfolio allocation. 
Select holding period $h$, number of alphas to select is $K$.
$D_{in}$ is the in-sample number of days, $D$ is the total number of days, $N$ is the total number of alphas. We then select $K$ alphas out of the $N$ alphas.
} 
\label{alg:iid_alpha_selection}  
\begin{algorithmic}[1] 
\State Split the alpha matrix for in-sample (IS) and out-of-sample (OS) periods: 
$$
\bA_{\text{in}} = \bA_s[:,0:D_{\text{in}}] \in \real^{N\times D_{\text{in}}}, \gap 
\bA_{\text{out}} = \bA_s[:,D_{\text{in}}+1:D]\in \real^{N\times (D-D_{\text{in}})};
$$
\State Apply (column) ID to $\bA_{\text{in}}^{\textcolor{mylightbluetext}{\top}}$ to select $K$ basis columns, yielding indices $\bmm$:
$$
\widehat{\bA}_{\text{in}} = \bA_s[\bmm,0:D_{\text{in}}] \in \real^{\textcolor{mylightbluetext}{K}\times D_{\text{in}}}, \gap 
\widehat{\bA}_{\text{out}} = \bA_s[\bmm,D_{\text{in}}+1:D]\in \real^{\textcolor{mylightbluetext}{K}\times (D-D_{\text{in}})};
$$
\For{$k=1$ to $K$}
\State \algoalign{Using the $k$-th IS alpha vector $\ba_k=\widehat{\bA}_{\text{in}}[k,:]\in \real^{D_{\text{in}}}$ to decide the weight $\bw_k$ and interception $b_k$ via ordinary least squares (OLS) so that the MSE between the prediction $\ba_k^\top\bw_k +b_k$ and the shifted return vector $\br^h$ is minimized, i.e., minimizing $\text{MSE}(\ba_k^\top\bw_k +b_k, \br^h)$. The weight and interception are then used in OS evaluation.}
\EndFor
\For{$d=1$ to $D-D_{\text{in}}$}
\State \algoalign{On each day in the OS period, we use the mean evaluation of each prediction from the $K$ alphas to decide to go long or not, i.e., to go long if $\sum_{k=1}^K \ba_k^\top\bw_k +b_k >0$; 
otherwise, hold cash since we restrict the analysis to  long-only portfolios.}
\State \algoalign{Though we employ a long-only portfolio, we can favor a \textit{market-neutral strategy}: we open long positions only when we anticipate that at least half of the stocks will rise on the following $h$ day, and we weight each stock equally.}

\EndFor
\end{algorithmic} 
\end{algorithm}

\begin{table}[]
\centering
\setlength{\tabcolsep}{7.3pt}
\renewcommand{\arraystretch}{1.2}
\small
\begin{tabular}{lllll}
\hline 
Methods & Highest RankIC & Randomized ID & BID with GBT & BID with IID \\
\hline
Mean RankIC & \textbf{0.1035}         & 0.0651        & 0.0553       & \textbf{0.0752}  \\
Mean Correlation & 0.2276$\downarrow$         & 0.5741$\downarrow$        & \textbf{0.1132}       & \textbf{0.1497}  \\
\hline
Sharpe Ratio (OS) & 1.0276 & 1.0544 & 0.5045 & \textbf{1.5721}  \\
Sharpe Ratio (IS) & \textbf{2.6511} & 1.3019 & 1.4965 & 2.3231  \\
\hline
Annual Return (OS) & 0.1043 & 0.0932 & 0.0484 & \textbf{0.1633}  \\
Annual Return (IS) & \textbf{0.4390} & 0.2281 & 0.2425 & 0.3805  \\
\hline
Max Drawdown (OS) & 0.0632 & \textbf{0.0373} & \textbf{0.0484} & \textbf{0.0552}  \\
Max Drawdown (IS) & \textbf{0.0892} & 0.1548 & 0.1232 & \textbf{0.0975}  \\
\hline
\end{tabular}
\vspace{-0.3cm} 
\caption{
Performance comparison across alpha selection methods. Higher RankIC and lower correlation are desirable. The symbol ``$\downarrow$" indicates poor diversification (high correlation). IID achieves the best out-of-sample risk-adjusted returns, and it balances the trade-off between the mean RankIC and the mean correlation.
In all cases, \textit{IS} means in-sample measurements, and \textit{OS} means out-of-sample measurements. 
}
\label{table:iid_selected_mean_ic}
\end{table}

\subsubsection{Quantitative Strategy}

After running GBT and IID on each asset's alpha matrix, we retain the state vector $\br$ and select the top 10 alphas with the highest average selection frequency over 1,000 Gibbs iterations (after 100 burn-in and thinning by 5).

We then implement the strategy outlined in Algorithm~\ref{alg:iid_alpha_selection} with $h=1$, $N=191$ alphas, $K=10$ alphas, $D=720$ trading days, and $D_{\text{in}}=480$ trading days. Although simple, this pipeline demonstrates how IID can be deployed in practice.

As shown in Table~\ref{table:iid_selected_mean_ic} and Figure~\ref{fig:bid_iid_portfolio_ashare}, the IID-based strategy slightly underperforms the \textit{highest-RankIC baseline} (i.e., simply selecting the highest-RankIC alphas) in-sample (in terms of \textit{Sharpe ratio}, annual return, and maximum drawdown) but  outperforms it out-of-sample---which is the ultimate objective in quantitative finance.
For comparison, we also include results from Randomized ID \citep{liberty2007randomized}, which performs worse than even vanilla BID (GBT). This highlights the value of principled selection.
Although IID does not select the absolute highest-RankIC alphas, this is by design---and beneficial---for several reasons:
\begin{enumerate}
\item \textit{Scalability.} Our alpha pool contains only 191 factors. In real-world settings with millions or billions of alphas, naive top-$K$ selection becomes unreliable, whereas IID's representativeness-aware approach scales more robustly.
\item \textit{Diversification.} The alphas selected by the highest-RankIC method exhibit high mutual correlation (0.2276 vs. 0.1497 for IID), reducing portfolio diversity and increasing vulnerability to regime shifts.
\item \textit{Model compatibility.} Simple OLS is used here, but in more complex models (e.g., neural networks, XGBoost), high multicollinearity among inputs degrades performance and interpretability. IID mitigates this by promoting diversity.
\item \textit{Risk management.} Even if top-RankIC alphas perform well historically, over-reliance on them concentrates risk. IID enables discovery of alternative, less correlated strategies, enhancing robustness.
\end{enumerate}
Thus, IID strikes a practical balance between predictive power and representational diversity, making it well-suited for real-world alpha selection.

\begin{problemset}

\item Determine the column ID for the matrix 
$$
\bA = 
\begin{bmatrix}
	1 & 3 & 2 \\
	3 & 7 & 6 \\
	4 & 5 & 8
\end{bmatrix}.
$$

\item \textbf{Magnitude matters.} Suppose you are given an $N\times N$ matrix where the absolute value of every entry is at most 1. Show that the absolute value of the determinant of this matrix is also at most $(N)^{N/2}$. Additionally, provide an example of a $2 \times 2$ matrix for which the determinant achieves this upper bound.

\item \label{prob:adjugate} \textbf{Adjugate.}
Let $\bA\in\real^{N\times N}$ be any square matrix. Then, the \textit{adjugate} of $\bA$, denoted $\adjugate(\bA)$,  is the $N\times N$ matrix whose $(i,j)$-th element is defined  by 
\begin{equation}\label{equation:adjug1}
	\adjugate(\bA)_{ij} = (-1)^{i+j} \det\left(\bA\big[\comple{\{j\}}, \comple{\{i\}}\big]\right),
\end{equation}
where $\comple{\{i\}}$ is the complementary set of $\{1,2,\ldots,N\}$: $\comple{\{i\}}=\{1,2,\ldots,N\}\setminus \{i\}$.
Show that 
\begin{equation}\label{equation:adjug2}
	\adjugate(\bA)\bA = \bA \adjugate(\bA) = \det(\bA)\bI.
\end{equation}

\index{Cramer's rule}
\item \label{prob:cramer_adj_1} \textbf{Cramer's rule.}	Consider the linear system $\bG\bx=\bl$, where $\bG\in\real^{N\times N}$, and $\bx,\bl\in\real^N$. Let $\bG_{\bl}(n)$ be the matrix formed by replacing the $n$-th column of $\bG$ with the vector $\bl$. Show that the $n$-th component  of the vector $\adjugate(\bG)\bl\in\real^N$ is given by
\begin{equation}\label{equation:cramer_adj_1}
\big(\adjugate(\bG)\bl\big)_n = \det(\bG_{\bl}(n)), \gap n\in\{1,2,\ldots,N\}.
\end{equation}
Now consider the matrix equation  $\bG\bX=\bL$, where $\bG\in\real^{N\times N}$, and $\bX,\bL\in\real^{N\times M}$. Let $\bG_{\bL}(n,m)$ be the matrix formed by replacing the $n$-th column of $\bG$ with the $m$-th column $\bl_m$ of $\bL$. Show that the $(n,m)$-th element of $\adjugate(\bG)\bL\in\real^{N\times M}$ is 
\begin{equation}\label{equation:cramer_adj_2}
\big(\adjugate(\bG)\bL\big)_{nm} = \det(\bG_{\bL}(n,m)), \gap n\in\{1,2,\ldots,N\}, m\in\{1,2,\ldots,M\}.
\end{equation}
\textit{Hint: Use the definitions of the determinant and the adjugate.}

\item \label{prob:cramer_adj_2} \textbf{Cramer's rule.} Under the same setting as Problem~\ref{prob:cramer_adj_1}, define the vector and matrix
$$
\adjugate(\bG)\bl = \big[\det(\bG_{\bl}(n))\big]_{n=1}^N\in\real^N,
\gap 
\adjugate(\bG)\bL = \big[\det(\bG_{\bl}(n,m))\big]_{n,m=1}^{N,M}\in\real^{N\times M},
$$
so that  the $n$-th element of the vector is $\det(\bG_{\bl}(n))$, and the  $(n,m)$-th element of the matrix is $\det(\bG_{\bL}(n,m))$. 
Show that
\begin{equation}\label{equation:cramer_adj_3}
\begin{aligned}
\bG \big[\det(\bG_{\bl}(n))\big]_{n=1}^N &= \bG \adjugate(\bG) \bl=\det(\bG)\bl;\\
\bG \big[\det(\bG_{\bL}(n,m))\big]_{n,m=1}^{N,M} &= \bG \adjugate(\bG) \bL=\det(\bG)\bL.
\end{aligned}
\end{equation}

\item \label{prob:cramer_adj_3} \textbf{Cramer's rule.} Assume the same setup as in Problem~\ref{prob:cramer_adj_1}, and further suppose that $\bG$ is nonsingular. Show that the $n$-th element of the solution $\bx$ is 
\begin{equation}\label{equation:cramer_adj_res1}
x_n = \frac{\det(\bG_{\bl}(n))}{\det(\bG)}, \gap \forall n\in\{1,2,\ldots,N\}.
\end{equation}
Similarly, show that the $(n,m)$-th element of the solution $\bX$ is 
\begin{equation}\label{equation:cramer_adj_res2}
x_{nm} = \frac{\det\left(\bG_{\bL}(n,m)\right)}{\det(\bG)}, \gap \forall n\in\{1,2,\ldots,N\}, m\in\{1,2,\ldots,M\}.
\end{equation}
These formulas constitute \textit{Cramer's rule}.

\item \label{prob:cramer_adj_4} \textbf{Cramer's rule: the simple way.}
Consider the same setting as Problem~\ref{prob:cramer_adj_1}, and assume further that $\bG$ is nonsingular. 
Observe that
\begin{equation}\label{equation:cramer_adj_4}
	\bG \bI_{\bl}(n) = \bG_{\bl}(n),  \gap \forall n\in\{1,2,\ldots,N\},
\end{equation}
where $\bI_{\bl}(n)$ represents the identity matrix with  the $n$-th column replaced by $\bl$. 
Taking determinants on both sides yields
\begin{equation}\label{equation:cramer_adj_5}
\det(\bG) \det(\bI_{\bl}(n)) = \det(\bG_{\bl}(n)).
\end{equation}
Show that  $\det(\bI_{\bl}(n))=x_n$, thereby recovering the result in  \eqref{equation:cramer_adj_res1}.

\item \label{prob:cramer_adj_5} \textbf{Determinant of inverses for subsets, Jacobi's equality.} 
Let  $\bG\in\real^{N\times N}$ be nonsingular, and  let $\sI,\sJ\subseteq\{1,2,\ldots,N\}$ be two index sets  with complement $\comple{\sI}$ and $\comple{\sJ}$, respectively. 
Prove that
\begin{equation}
	\det\big(\bG^{-1}[\comple{\sI},\comple{\sJ}]\big) = (-1)^{\gamma} \frac{\det(\bG[\sJ,\sI])}{\det(\bG)},
\end{equation}
where $\gamma=\sum_{i\in \sI} i +\sum_{j\in \sJ}j$ is the sum of indices.
In the special case $\sI=\sJ$, this reduces to
\begin{equation}
	\det\big(\bG^{-1}[\comple{\sI},\comple{\sI}]\big) = \frac{\det(\bG[\sI,\sI])}{\det(\bG)},
\end{equation}
which is known as \textit{Jacobi's equality}.
\textit{Hint: Examine the definitions of determinant and adjugate. Alternatively, the result can be derived via the Schur complement.}

\item \label{problem:id1} Provide a concrete example illustrating how the post-processing method described in Section~\ref{section:id_pprocess} reduces reconstruction error in a Bayesian interpolative decomposition.

\item \label{problem:rank_decom} Prove the rank decomposition:
Any rank-$R$ matrix $\bA \in \real^{M \times N}$ admits the factorization
$$
\underset{M\times N}{\bA }= \underset{M\times R}{\bD}\gapthree  \underset{R\times N}{\bF},
$$
where $\bD \in \real^{M\times R}$ has rank $R$, and $\bF \in \real^{R\times N}$ also has rank $R$, i.e., $\bD$ and $\bF$ have full rank $R$.
The storage for the decomposition is then reduced or potentially increased from $MN$ floating-point numbers to $R(M+N)$ floating-point numbers. \textit{Hint: Use elementary row and column operations or the reduced row echelon form.}

\item \label{problem:mtb_dis} \textbf{Determinantal identities via rank factorization.} 
Let  $\bA\in\real^{M\times N}$ have rank $R$, and suppose $\bA=\bD\bF$ is a  rank decomposition with $\bD\in\real^{M\times R}$ and $\bF\in\real^{R\times N}$.
Let $\sI,\sJ\subseteq \{1,2,\ldots,M\}$ and $\sK,\sL\subseteq\{1,2,\ldots,N\}$ be index sets, each of cardinality $R$.
Then, we have  $\bA[\sI,\sK]=\bD[\sI,:]\bF[:,\sK]$.
Show that 
\begin{itemize}
\item $\bA[\sI,\sK]$ is nonsingular if and only if $\rank(\bD[\sI,:])=\rank(\bF[:,\sK])=R$.
\item $\det(\bA[\sI,\sK]) \det(\bA[\sJ,\sL])=\det(\bA[\sI,\sL])\det(\bA[\sJ,\sK])$.
\end{itemize}

\end{problemset}

\newpage
\vskip 0.2in
\addcontentsline{toc}{chapter}{Bibliography}
\bibliography{bib}

\clearpage
\addcontentsline{toc}{chapter}{Index}
\printindex

\end{document}